\documentclass{amsbook}
\usepackage{amssymb}
\usepackage{mathrsfs}
\usepackage{stmaryrd} 
\usepackage{bbm}
\usepackage{oldgerm}
\usepackage[frenchb]{babel}
\usepackage[T1]{fontenc}
\usepackage[latin1]{inputenc}
\usepackage[all]{xy}
\usepackage{hyperref}
\usepackage{upref}
\usepackage{xcolor}
\usepackage{microtype}
\usepackage[Symbol]{upgreek}

\hypersetup{
    colorlinks,
    linkcolor={red!50!black},
    citecolor={blue!50!black},
    urlcolor={blue!80!black}
}

\newtheorem{teo}[subsection]{Théorème}
\newtheorem{prop}[subsection]{Proposition}
\newtheorem{cor}[subsection]{Corollaire}
\newtheorem{lem}[subsection]{Lemme}

\theoremstyle{definition}

\newtheorem{defi}[subsection]{Définition}
\newtheorem{rema}[subsection]{Remarque}
\newtheorem{remas}[subsection]{Remarques}

\numberwithin{equation}{subsection}

\newcommand{\gtimes}{\stackrel{\leftarrow}{\times}}

\mathchardef\mhyphen="2D

\newcommand{\muu}{{\mathbbm u}}
\newcommand{\mA}{{\mathbb A}}

\newcommand{\mR}{{\mathbb R}}

\newcommand{\mH}{{\mathbb H}}

\newcommand{\mQ}{{\mathbb Q}}
\newcommand{\mL}{{\mathbb L}}
\newcommand{\mN}{{\mathbb N}}
\newcommand{\mP}{{\mathbb P}}

\newcommand{\mX}{{\mathbb X}}
\newcommand{\mY}{{\mathbb Y}}
\newcommand{\mZ}{{\mathbb Z}}
\newcommand{\mF}{{\mathbb F}}
\newcommand{\mG}{{\mathbb G}}

\newcommand{\mU}{{\mathbb U}}
\newcommand{\mV}{{\mathbb V}}

\newcommand{\bA}{{\bf A}}
\newcommand{\bB}{{\bf B}}

\newcommand{\bD}{{\bf D}}

\newcommand{\Et}{{\bf \acute{E}t}}
\newcommand{\Sch}{{\bf Sch}}
\newcommand{\Ens}{{\bf Ens}}
\newcommand{\Pt}{{\bf Pt}}
\newcommand{\Mon}{{\bf Mon}}

\newcommand{\Top}{{\bf Top}}
\newcommand{\bHom}{{\bf Hom}}

\newcommand{\bRep}{{\bf Rep}}

\newcommand{\bAlg}{{\bf Alg}}
\newcommand{\bMod}{{\bf Mod}}
\newcommand{\aAlg}{{\alpha\mhyphen \bf Alg}}
\newcommand{\aMod}{{\alpha\mhyphen \bf Mod}}

\newcommand{\aI}{{\alpha\mhyphen \bf I}}
\newcommand{\DDiv}{{{\mhyphen \rm Div}}}

\newcommand{\acA}{{\alpha\mhyphen \cA}}

\newcommand{\ahcC}{{\alpha\mhyphen \hcC}}
\newcommand{\atcC}{{\alpha\mhyphen \tcC}}
\newcommand{\ahE}{{\alpha\mhyphen \hE}}
\newcommand{\atE}{{\alpha\mhyphen \tE}}
\newcommand{\eEt}{{\text -{\bf \acute{E}t}}}
\newcommand{\apf}{{\alpha{\rm pf}}}

\newcommand{\intern}{{\diamond}}

\newcommand{\et}{{\rm \acute{e}t}}
\newcommand{\fet}{{\rm f\acute{e}t}}

\newcommand{\coim}{{\rm coim}}

\newcommand{\zar}{{\rm zar}}
\newcommand{\aff}{{\rm aff}}
\newcommand{\coh}{{\rm coh}}
\newcommand{\qcoh}{{\rm qcoh}}
\newcommand{\cart}{{\rm cart}}
\newcommand{\scoh}{{\rm scoh}}

\newcommand{\longueur}{{\rm long}}

\newcommand{\atf}{{\rm atf}}

\newcommand{\proj}{{\rm proj}}

\newcommand{\rf}{{\rm f}}

\newcommand{\Spec}{{\rm Spec}}

\newcommand{\Spf}{{\rm Spf}}

\newcommand{\ob}{{\rm Ob}}

\newcommand{\fil}{{\rm fil}}

\newcommand{\tor}{{\rm tor}}

\newcommand{\coker}{{\rm coker}}

\newcommand{\Cov}{{\rm Cov}}

\newcommand{\im}{{\rm im}}

\newcommand{\cont}{{\rm cont}}

\newcommand{\gp}{{\rm gp}}
\newcommand{\id}{{\rm id}}
\newcommand{\Tr}{{\rm Tr}}
\newcommand{\rg}{{\rm rg}}

\newcommand{\Fil}{{\rm Fil}}

\newcommand{\Hom}{{\rm Hom}}

\newcommand{\End}{{\rm End}}

\newcommand{\Ext}{{\rm Ext}}
\newcommand{\Tor}{{\rm Tor}}
\newcommand{\Gal}{{\rm Gal}}

\newcommand{\AR}{{\rm AR}}

\newcommand{\rC}{{\rm C}}

\newcommand{\rE}{{\rm E}}
\newcommand{\rF}{{\rm F}}
\newcommand{\rH}{{\rm H}}
\newcommand{\rT}{{\rm T}}
\newcommand{\rK}{{\rm K}}
\newcommand{\rL}{{\rm L}}
\newcommand{\rM}{{\rm M}}

\newcommand{\rR}{{\rm R}}
\newcommand{\rS}{{\rm S}}
\newcommand{\rV}{{\rm V}}
\newcommand{\rW}{{\rm W}}

\newcommand{\ra}{{\rm a}}
\newcommand{\rc}{{\rm c}}

\newcommand{\rp}{{\rm p}}

\newcommand{\lh}{{\stackrel{\leftarrow}{h}}}

\newcommand{\lgg}{{\ttg}}

\newcommand{\oB}{{\overline{B}}}
\newcommand{\oC}{{\overline{C}}}

\newcommand{\oF}{{\overline{F}}}

\newcommand{\oK}{{\overline{K}}}

\newcommand{\oR}{{\overline{R}}}
\newcommand{\oS}{{\overline{S}}}

\newcommand{\oU}{{\overline{U}}}

\newcommand{\oW}{{\overline{W}}}
\newcommand{\oX}{{\overline{X}}}

\newcommand{\oZ}{{\overline{Z}}}
\newcommand{\oa}{{\overline{a}}}
\newcommand{\oee}{{\overline{e}}}
\newcommand{\of}{{\overline{f}}}
\newcommand{\ogg}{{\overline{g}}}

\newcommand{\ol}{{\overline{l}}}

\newcommand{\op}{{\overline{p}}}

\newcommand{\ou}{{\overline{u}}}
\newcommand{\ov}{{\overline{v}}}
\newcommand{\ow}{{\overline{w}}}
\newcommand{\ox}{{\overline{x}}}
\newcommand{\oy}{{\overline{y}}}
\newcommand{\oz}{{\overline{z}}}
\newcommand{\ocF}{{\overline{\cF}}}

\newcommand{\oGamma}{{\overline{\Gamma}}}
\newcommand{\oalpha}{{\overline{\alpha}}}

\newcommand{\oeta}{{\overline{\eta}}}
\newcommand{\ophi}{{\overline{\phi}}}
\newcommand{\opsi}{{\overline{\psi}}}

\newcommand{\oiota}{{\overline{\iota}}}

\newcommand{\ocB}{{\overline{\cB}}}

\newcommand{\ofp}{{\overline{\fp}}}
\newcommand{\otta}{{\overline{\tta}}}
\newcommand{\oupdelta}{{\overline{\updelta}}}

\newcommand{\uB}{{\underline{B}}}
\newcommand{\uE}{{\underline{E}}}
\newcommand{\uG}{{\underline{G}}}

\newcommand{\uR}{{\underline{R}}}
\newcommand{\uU}{{\underline{U}}}
\newcommand{\uX}{{\underline{X}}}
\newcommand{\uY}{{\underline{Y}}}
\newcommand{\uZ}{{\underline{Z}}}
\newcommand{\uf}{{\underline{f}}}
\newcommand{\ug}{{\underline{g}}}
\newcommand{\uh}{{\underline{h}}}

\newcommand{\upp}{{\underline{p}}}

\newcommand{\uoB}{{\underline{\oB}}}
\newcommand{\uoC}{{\underline{\oC}}}

\newcommand{\uoR}{{\underline{\oR}}}

\newcommand{\uoZ}{{\underline{\oZ}}}
\newcommand{\uhoR}{{\underline{\hoR}}}
\newcommand{\uocB}{{\underline{\ocB}}}
\newcommand{\uoX}{{\underline{\oX}}}

\newcommand{\ubeta}{{\underline{\beta}}}
\newcommand{\urho}{{\underline{\rho}}}
\newcommand{\uvarrho}{{\underline{\varrho}}}
\newcommand{\udelta}{{\underline{\delta}}}
\newcommand{\usigma}{{\underline{\sigma}}}
\newcommand{\ugamma}{{\underline{\gamma}}}
\newcommand{\upi}{{\underline{\pi}}}
\newcommand{\uPi}{{\underline{\Pi}}}
\newcommand{\uphi}{{\underline{\phi}}}
\newcommand{\utau}{{\underline{\tau}}}
\newcommand{\umu}{{\underline{\mu}}}

\newcommand{\ulambda}{{\underline{\lambda}}}
\newcommand{\uPhi}{{\underline{\Phi}}}
\newcommand{\uTheta}{{\underline{\Theta}}}
\newcommand{\uxi}{{\underline{\xi}}}
\newcommand{\unu}{{\underline{\nu}}}

\newcommand{\ujmath}{{\underline{\jmath}}}
\newcommand{\uDelta}{{\underline{\Delta}}}
\newcommand{\uGamma}{{\underline{\Gamma}}}
\newcommand{\uhbar}{{\underline{\hbar}}}
\newcommand{\uupgamma}{{\underline{\upgamma}}}
\newcommand{\ulgg}{{\underline{\lgg}}}
\newcommand{\umL}{{\underline{\mL}}}
\newcommand{\ucL}{{\underline{\cL}}}

\newcommand{\hA}{{\widehat{A}}}
\newcommand{\hB}{{\widehat{B}}}

\newcommand{\hD}{{\widehat{D}}}
\newcommand{\hE}{{\widehat{E}}}
\newcommand{\hP}{{\widehat{P}}}

\newcommand{\hcN}{{\widehat{\cN}}}

\newcommand{\hcF}{{\widehat{\cF}}}
\newcommand{\hocF}{{\widehat{\ocF}}}

\newcommand{\hRun}{{\widehat{R_1}}}

\newcommand{\huRun}{{\widehat{\uR_1}}}

\newcommand{\huRi}{{\widehat{\uR_\infty}}}

\newcommand{\hmZ}{{\widehat{\mZ}}}

\newcommand{\halpha}{\widehat{\alpha}}

\newcommand{\hsigma}{\widehat{\sigma}}

\newcommand{\cA}{{\mathscr A}}
\newcommand{\cB}{{\mathscr B}}
\newcommand{\cC}{{\mathscr C}}
\newcommand{\cD}{{\mathscr D}}
\newcommand{\cE}{{\mathscr E}}
\newcommand{\cF}{{\mathscr F}}
\newcommand{\cG}{{\mathscr G}}

\newcommand{\cK}{{\mathscr K}}
\newcommand{\cL}{{\mathscr L}}
\newcommand{\cP}{{\mathscr P}}
\newcommand{\co}{{\mathscr O}}
\newcommand{\cR}{{\mathscr R}}
\newcommand{\cS}{{\mathscr S}}
\newcommand{\cT}{{\mathscr T}}
\newcommand{\cH}{{\mathscr H}}
\newcommand{\cM}{{\mathscr M}}
\newcommand{\cN}{{\mathscr N}}
\newcommand{\cQ}{{\mathscr Q}}
\newcommand{\cU}{{\mathscr U}}
\newcommand{\cV}{{\mathscr V}}
\newcommand{\cW}{{\mathscr W}}

\newcommand{\cHom}{{\mathscr Hom}}

\newcommand{\fA}{{\mathfrak A}}
\newcommand{\fB}{{\mathfrak B}}
\newcommand{\fC}{{\mathfrak C}}
\newcommand{\fD}{{\mathfrak D}}

\newcommand{\fF}{{\mathfrak F}}
\newcommand{\fG}{{\mathfrak G}}

\newcommand{\fM}{{\mathfrak M}}
\newcommand{\fN}{{\mathfrak N}}

\newcommand{\fS}{{\mathfrak S}}

\newcommand{\fU}{{\mathfrak U}}
\newcommand{\fV}{{\mathfrak V}}
\newcommand{\fW}{{\mathfrak W}}
\newcommand{\fX}{{\mathfrak X}}

\newcommand{\fc}{{\mathfrak c}}

\newcommand{\fgg}{{\mathfrak g}}
\newcommand{\fm}{{\mathfrak m}}

\newcommand{\fp}{{\mathfrak p}}
\newcommand{\fq}{{\mathfrak q}}

\newcommand{\tta}{{\tt a}}

\newcommand{\ttg}{{\tt g}}
\newcommand{\ttt}{{\tt t}}

\newcommand{\hK}{{\widehat{K}}}
\newcommand{\hM}{{\widehat{M}}}

\newcommand{\hR}{{\widehat{R}}}

\newcommand{\hoR}{{\widehat{\oR}}}

\newcommand{\hcC}{{\widehat{\cC}}}

\newcommand{\hmX}{{\widehat{\mX}}}

\newcommand{\hotimes}{{\widehat{\otimes}}}

\newcommand{\bvg}{{\breve{g}}}
\newcommand{\bvu}{{\breve{u}}}

\newcommand{\bvR}{{\breve{R}}}

\newcommand{\bvoS}{{\breve{\oS}}}

\newcommand{\bvoX}{{\breve{\oX}}}
\newcommand{\bvocB}{{\breve{\ocB}}}
\newcommand{\bvocF}{{\breve{\ocF}}}

\newcommand{\bvtau}{{\breve{\tau}}}

\newcommand{\bvpi}{{\breve{\pi}}}

\newcommand{\bvsigma}{{\breve{\sigma}}}

\newcommand{\bvvarphi}{{\breve{\varphi}}}
\newcommand{\bvdelta}{{\breve{\delta}}}
\newcommand{\bvkappa}{{\breve{\kappa}}}

\newcommand{\bvPhi}{{\breve{\Phi}}}
\newcommand{\bvPsi}{{\breve{\Psi}}}
\newcommand{\bvpsi}{{\breve{\psi}}}
\newcommand{\bvTheta}{{\breve{\Theta}}}

\newcommand{\bvmZ}{{\breve{\mZ}}}

\newcommand{\bvogg}{{\breve{\ogg}}}
\newcommand{\bvutau}{{\breve{\utau}}}
\newcommand{\bvupgamma}{{\breve{\upgamma}}}
\newcommand{\bvlgg}{{\breve{\lgg}}}
\newcommand{\bvulgg}{{\breve{\ulgg}}}
\newcommand{\bvuoR}{{\breve{\uoR}}}

\newcommand{\coS}{{\check{\oS}}}

\newcommand{\coX}{{\check{\oX}}}

\newcommand{\tC}{{\widetilde{C}}}
\newcommand{\tD}{{\widetilde{D}}}
\newcommand{\tE}{{\widetilde{E}}}

\newcommand{\tG}{{\widetilde{G}}}
\newcommand{\tuE}{{\widetilde{\uE}}}
\newcommand{\tuG}{{\widetilde{\uG}}}

\newcommand{\tM}{{\widetilde{M}}}
\newcommand{\tN}{{\widetilde{N}}}

\newcommand{\tQ}{{\widetilde{Q}}}

\newcommand{\tS}{{\widetilde{S}}}

\newcommand{\tU}{{\widetilde{U}}}

\newcommand{\tX}{{\widetilde{X}}}
\newcommand{\tY}{{\widetilde{Y}}}

\newcommand{\ta}{{\widetilde{a}}}
\newcommand{\tb}{{\widetilde{b}}}
\newcommand{\tc}{{\widetilde{c}}}
\newcommand{\te}{{\widetilde{e}}}

\newcommand{\tx}{{\widetilde{x}}}

\newcommand{\trT}{{\widetilde{\rT}}}

\newcommand{\tOmega}{{\widetilde{\Omega}}}

\newcommand{\tbeta}{{\widetilde{\beta}}}
\newcommand{\tsigma}{{\widetilde{\sigma}}}

\newcommand{\tnu}{{\widetilde{\nu}}}

\newcommand{\tpi}{{\widetilde{\pi}}}

\newcommand{\tvarphi}{{\widetilde{\varphi}}}
\newcommand{\ttau}{{\widetilde{\tau}}}

\newcommand{\talpha}{{\widetilde{\alpha}}}

\newcommand{\tcC}{{\widetilde{\cC}}}
\newcommand{\tcA}{{\widetilde{\cA}}}

\newcommand{\tmX}{{\widetilde{\mX}}}

\newcommand{\tfD}{{\widetilde{\fD}}}
\newcommand{\tfG}{{\widetilde{\fG}}}
\newcommand{\tmZ}{{\widetilde{\mZ}}}

\begin{document}

\title{Les suites spectrales de Hodge-Tate}
\author{Ahmed Abbes et Michel Gros}
\address{A.A. Laboratoire Alexander Grothendieck, CNRS, IHES, Université Paris-Saclay,
35 route de Chartres, 91440 Bures-sur-Yvette, France}
\address{M.G. Université de Rennes, CNRS, IRMAR - UMR 6625, Campus de Beaulieu, F-35042 Rennes cedex, France}
\email{abbes@ihes.fr}
\email{michel.gros@univ-rennes1.fr}

\maketitle

\setcounter{tocdepth}{1}

\newpage{\setlength{\parindent}{0pc}}
\thispagestyle{empty}

\vspace*{13.5pc}
\begin{center}
{\Large\em À la mémoire de Jean-Marc Fontaine et Michel Raynaud}
\end{center}

\newpage\null\thispagestyle{empty}\newpage

\tableofcontents

\chapter*{Avant-Propos}

Ce livre présente deux résultats importants en théorie de Hodge $p$-adique suivant l'approche initiée par Faltings dans \cite{faltings1,faltings2},  à savoir
\begin{itemize}
\item[(i)] son principal théorème de comparaison $p$-adique;
\item[(ii)] la suite spectrale de Hodge-Tate.
\end{itemize}
Nous établissons pour chacun de ces résultats deux versions, une absolue et  une relative. 
Si les énoncés absolus peuvent raisonnablement être considérés comme bien compris, 
en particulier après leur extension aux variétés rigides par Scholze \cite{scholze1}, 
l'approche initiale de Faltings pour les variantes relatives est demeurée bien moins étudiée. 
La suite spectrale de Hodge-Tate relative est, quant à elle, nouvelle dans cette approche. 

Bien que nous suivions la même stratégie que celle utilisée par Faltings pour établir son principal théorème de comparaison 
$p$-adique \cite{faltings2}, une partie de nos démonstrations repose sur de nouveaux résultats. 
La preuve du cas absolu, que nous présentons avec tous les détails qu'elle mérite, est considérablement simplifiée par rapport 
à celle de Faltings par l'utilisation d'un résultat récent d'Achinger sur le caractère local $K(\pi,1)$ 
des schémas logarithmiques considérés \cite{achinger}. 
Faltings a formulé la version relative de son principal théorème de comparaison $p$-adique dans
\cite{faltings2} et en a très sommairement esquissé une preuve dans l'appendice. 
Certains de ses arguments doivent d'ailleurs être modifiés et la preuve donnée dans ce livre requiert bien plus de travail.
Elle est basée sur une étude fine de la structure locale de certains $\varphi$-modules presque-étales
qui est, par ailleurs, intéressante en elle-même. 

La suite spectrale de Hodge-Tate dans le cas absolu, qu'on devine en filigrane dans le travail de Faltings \cite{faltings2}, n'a été explicitement dégagée qu'ultérieurement par Scholze qui l'a généralisée aux variétés rigides \cite{scholze2}. 
Nous la déduisons du principal théorème de comparaison $p$-adique 
de Faltings et de la théorie de Kummer sur la fibre spéciale du topos de Faltings annelé. 
La suite spectrale de Hodge-Tate relative prend racine dans le topos de Faltings. Sa construction nécessite l'introduction d'une variante
relative de ce topos dont la définition et l'étude sont les principales nouveautés de notre travail.
Nous donnons également deux variantes ``locales'' qui s'en déduisent dont l'énoncé ne nécessite pas l'utilisation du topos de Faltings. 
Utilisant une approche différente, Caraiani et Scholze (\cite{cs} 2.2.4) ont construit une filtration de Hodge-Tate 
relative pour les morphismes propres et lisses entre espaces adiques. Antérieurement, Hyodo avait aussi considéré un cas particulier 
de telle filtration pour des schémas abéliens \cite{hyodo2}.  

Outre la volonté de renforcer les fondations de résultats devenus centraux en géométrie arithmétique, ce travail a été motivé par l'étude 
de la fonctorialité par image directe propre et (log-)lisse, à la Katz-Oda, de la correspondance de Simpson $p$-adique  développée dans \cite{agt}. 
Cette question est traitée dans \cite{ag2}.    
Le présent livre est, grâce aux rappels qu'il contient, essentiellement indépendant de \cite{ag2,agt}.

Au-delà de la table des matières, donnons maintenant quelques indications sur l'organisation de ce volume. 
Le chapitre \ref{Intro-intro} présente un survol des principaux résultats obtenus dans ce livre et la stratégie de leurs démonstrations. 
Par simplicité, nous avons fait le choix de nous limiter dans ce chapitre introductif aux schémas lisses, 
plutôt qu'à singularités toriques considérés dans le reste du livre.  
Nous avons aussi mis l'accent sur la situation relative qui est la plus inédite. 
Le chapitre \ref{prelim} rassemble divers résultats préliminaires de natures assez différentes utilisés ici ou là dans le texte; 
le lecteur peut ne le consulter qu'au fur et à mesure de ses besoins. 
Le chapitre \ref{toposfaltings}, après quelques préliminaires sur les produits orientés de topos, 
fournit les rappels nécessaires sur le topos de Faltings (\cite{agt} VI) et donne la construction de sa variante relative. 
Le chapitre \ref{cohtopfal} traite ensuite de la cohomologie du topos de Faltings et se conclut par la preuve du principal théorème de comparaison 
$p$-adique de Faltings. Suivant un cheminement similaire, le chapitre \ref{cohtopfalrel} aborde alors l'extension de ce type de résultats au cadre relatif. 
Les démonstrations deviennent plus techniques et un ingrédient nouveau est à mettre en \oe{}uvre 
concernant la structure de certains $\varphi$-modules presque-étales \eqref{elcpam15}. 
Enfin, le chapitre \ref{suitesspecht} exploite le formalisme des topos de Faltings et 
les théorèmes de comparaison $p$-adiques obtenus précédemment 
pour construire les suites spectrales de Hodge-Tate et établir leurs propriétés de fonctorialité et 
d'équivariance aux actions galoisiennes naturelles.

\subsection*{Remerciements} Nous tenons en premier lieu à exprimer toute notre reconnaissance à G. Faltings  
pour l'inspiration constante suscitée par ses travaux. Ce livre s'inscrit dans le prolongement direct 
de ses idées en théorie de Hodge $p$-adique. 
\`A diverses étapes de l'élaboration de ce texte, l'inestimable expertise  d'O. Gabber et de T. Tsuji  
a été  cruciale pour nous éviter de longs et inutiles détours. O. Gabber nous a de plus  aidé à corriger 
une erreur importante dans une première version de ce livre et T. Tsuji nous a fait bénéficier 
de très précieux et nombreux commentaires pour en améliorer la rédaction: 
qu'ils soient tous deux ici très chaleureusement remerciés pour leur indéfectible générosité.
Nous remercions T. He pour sa lecture attentive du manuscrit ayant conduit à de nombreuses améliorations. 
Nous remercions aussi les rapporteurs pour leurs  nombreuses suggestions et tout spécialement celui des chapitres 5 et 6 
pour sa relecture très méticuleuse et ses nombreux et utiles commentaires. 
Le premier auteur (A.A) remercie l'Université de Tokyo et l'Université Tsinghua 
pour leur hospitalité lors de plusieurs visites où des parties de ce travail ont été développées et présentées.
Il exprime sa gratitude à T.~Saito et L.~Fu pour leurs invitations.

\chapter{La suite spectrale de Hodge-Tate relative -- un survol}\label{Intro-intro}

\section{Introduction}

\subsection{}\label{Intro-intro1}
Soient $K$ un corps de valuation discrète complet de caractéristique $0$, à corps résiduel {\em algébriquement clos} de caractéristique $p>0$, 
$\co_K$ l'anneau de valuation de $K$, $\oK$ une clôture algébrique de $K$, $\co_\oK$ la clôture intégrale de $\co_K$ dans $\oK$.
On note $G_K$ le groupe de Galois de $\oK$ sur $K$,   $\co_C$ le complété $p$-adique de $\co_\oK$,  
$\fm_C$ l'idéal maximal de $\co_C$ et $C$ son corps des fractions.
On pose $S=\Spec(\co_K)$ et $\oS=\Spec(\co_\oK)$ et on note $s$ (resp.  $\eta$, resp. $\oeta$) 
le point fermé de $S$ (resp. le point générique de $S$, resp. le point générique de $\oS$). 
Pour tout entier $n\geq 0$, on pose $S_n=\Spec(\co_K/p^n\co_K)$. Pour tout $S$-schéma $X$, on pose
\begin{equation}\label{Intro-intro1a}
\oX=X\times_S\oS \ \ \ {\rm et}\ \ \  X_n=X\times_SS_n.
\end{equation}

L'énoncé suivant, appelé la {\em décomposition de Hodge-Tate}, a été conjecturé par Tate (\cite{tate} Remark page 180)
et démontré par différentes méthodes par Faltings \cite{faltings1,faltings2}, Niziol \cite{niziol1,niziol2} et Tsuji \cite{tsuji1,tsuji2}.

\begin{teo}\label{Intro-intro2}
Pour tout $\eta$-schéma propre et lisse $X$ et tout entier $n\geq 0$, 
il existe une décomposition fonctorielle $G_K$-équivariante canonique
\begin{equation}\label{Intro-intro2a}
\rH^n_\et(X_\oeta,\mQ_p)\otimes_{\mQ_p}C\stackrel{\sim}{\rightarrow}\bigoplus_{i=0}^n\rH^i(X,\Omega^{n-i}_{X/\eta})\otimes_KC(i-n).
\end{equation}
\end{teo}

La décomposition de Hodge-Tate est équivalente à l'existence d'une suite spectrale canonique, fonctorielle et 
$G_K$-équivariante, la {\em suite spectrale de Hodge-Tate},
\begin{equation}\label{Intro-intro2b}
\rE_2^{i,j}=\rH^i(X,\Omega^j_{X/\eta})\otimes_KC(-j)\Rightarrow \rH^{i+j}_\et(X_\oK,\mQ_p)\otimes_{\mQ_p}C.
\end{equation}
Les deux énoncés sont équivalents d'après un théorème de Tate (\cite{tate} theo.~2).
En effet, le groupe de cohomologie $\rH^0(G_K,C(1))$ s'annule, ce qui implique que la suite spectrale dégénère en $\rE_2$.
Le groupe de cohomologie $\rH^1(G_K,C(1))$ s'annule également, ce qui implique que la filtration aboutissement est scindée. 

Cette suite spectrale de Hodge-Tate, qu'on devine en filigrane dans le travail de Faltings \cite{faltings2}, n'a été explicitement dégagée 
qu'ultérieurement par Scholze qui l'a généralisée aux variétés rigides~\cite{scholze2}.

\subsection{}\label{Intro-intro3}
Nous donnons dans ce chapitre un survol du travail présenté dans ce livre conduisant à une  généralisation de la  suite spectrale de Hodge-Tate 
aux morphismes. Celle-ci prend racine dans le topos de Faltings. Sa construction requiert l'introduction d'une variante
relative de ce topos qui est la principale nouveauté de notre travail.
Utilisant une approche différente, Caraiani et Scholze (\cite{cs} 2.2.4) ont construit une filtration de Hodge-Tate relative pour les morphismes
propres et lisses entre espaces adiques. Antérieurement, Hyodo avait aussi considéré un cas particulier pour des schémas abéliens \cite{hyodo2}.  

Au-delà des suites spectrales de Hodge-Tate, nous donnons dans ce livre des preuves complètes des principaux théorèmes
de comparaison $p$-adiques de  Faltings. Ces derniers sont essentiels pour la construction de ces suites spectrales. Bien que la 
version absolue de ces théorèmes soit plutôt bien comprise, la version relative, seulement sommairement esquissée par Faltings
dans l'appendice de \cite{faltings2}, est restée peu étudiée.
En étendant la stratégie de Faltings, Scholze a prouvé des résultats similaires (\cite{scholze1} 1.3 and 5.12) dans le contexte des espaces adiques et des topos pro-étales.
On renvoie aux articles de He \cite{tongmu1,tongmu2} pour une analyse des liens entre les différentes approches. 

Dans un travail récent \cite{ag2} consacré à la fonctorialité de correspondance de Simpson $p$-adique \cite{agt} par image directe propre et (log-)lisse, 
nous étendons la construction de la suite spectrale de Hodge-Tate relative à des coefficients plus généraux. 

\subsection{}\label{Intro-logsch} 
Avant d'introduire la suite spectrale de Hodge-Tate relative qui requiert le formalisme du topos de Faltings, 
nous présentons tout d'abord deux versions locales qui en sont des conséquences plus faciles à énoncer. 

\vspace{2mm} 

Nous traitons dans ce livre le cas des schémas à singularités toriques
à l'aide de la géométrie logarithmique mais, par simplicité, nous considérons dans ce survol seulement le cas lisse.

\section[La suite spectrale de Hodge-Tate relative: sections globales]{La suite spectrale de Hodge-Tate relative: sections globales au-dessus d'un petit schéma affine}\label{Intro-lv}

\subsection{}\label{Intro-lv1}
Soit $X=\Spec(R)$ un $S$-schéma affine et lisse vérifiant les deux conditions suivantes: 
\begin{itemize}
\item[(i)] $X$ est {\em petit}  dans le sens de Faltings, c'est-à-dire qu'il admet un $S$-morphisme étale dans 
un $S$-tore, $X\rightarrow \mG_{m,S}^d=\Spec(\co_K[T_1^{\pm1},\dots, T_d^{\pm1}])$, pour un entier $d\geq 0$. 
Cette condition sera remplacée dans le cas logarithmique général par l'existence de cartes adéquates \eqref{cad1};
\item[(ii)] $X_s$ est non vide.
\end{itemize}

Soit $\oy$ un point géométrique de $X_\oeta$. On désigne par $X_\eta^\star$ (resp.  $X_\oeta^\ast$) la composante connexe de $X_\eta$ 
(resp. $X_\oeta$) contenant l'image de $\oy$ et par $(V_i)_{i\in I}$ le revêtement universel de $X^\ast_\oeta$ en $\oy$ (\cite{agt} VI.9.7.3). 
On pose $\Gamma=\pi_1(X^\star_\eta,\oy)$ et $\Delta=\pi_1(X^\ast_\oeta,\oy)$.
Pour tout $i\in I$, on note $X_i=\Spec(R_i)$ la normalisation de $\oX=X\times_S\oS$ dans $V_i$. 
\begin{equation}\label{Intro-lv1a}
\xymatrix{
V_i\ar[r]\ar[d]&X_i\ar[d]\\
X_\oeta\ar[r]&\oX}
\end{equation}
Les $\co_\oK$-algèbres $(R_i)_{i\in I}$ forment naturellement un système inductif. On désigne par $\oR$ sa limite inductive, 
\begin{equation}\label{Intro-lv1b}
\oR=\underset{\underset{i\in I}{\longrightarrow}}{\lim} \ R_i, 
\end{equation}
et par $\hoR$ son complété $p$-adique, qu'on munit des actions naturelles de  $\Gamma$.
La $\Gamma$-représentation $\hoR$ est un analogue de la $G_K$-représentation $\co_C$.

\begin{teo}[cf. \ref{sshtrsg6}]\label{Intro-lv2} 
Sous les hypothèses de \ref{Intro-lv1}, 
pour tout morphisme projectif et  lisse $g\colon X'\rightarrow X$, et tout entier $q\geq 0$, 
il existe $(\fil^q_r)_{0\leq r\leq q+1}$, une filtration décroissante exhaustive canonique 
de $\rH^q_\et(X'_\oy,\mZ_p)\otimes_{\mZ_p}\hoR[\frac 1 p]$ par des $\hoR[\frac 1 p]$-représentations de $\Gamma$, 
telle que $\fil^q_{q+1}$ soit nul et que pour tout entier $0\leq r\leq q$, on ait une suite exacte $\Gamma$-équivariante canonique 
\begin{equation}\label{Intro-lv2a}
0\rightarrow \fil^q_{r+1}\rightarrow \fil^q_r\rightarrow \rH^r(X',\Omega^{q-r}_{X'/X})\otimes_R\hoR[\frac 1 p](r-q)
\rightarrow 0.
\end{equation}
\end{teo}

Il revient au même de dire qu'il existe une suite spectrale canonique $\Gamma$-équivariante 
\begin{equation}\label{Intro-lv2b}
\rE_2^{i,j}=\rH^i(X',\Omega^j_{X'/X})\otimes_R\hoR[\frac 1 p](-j)
\Rightarrow \rH^{i+j}_\et(X'_\oy,\mQ_p)\otimes_{\mZ_p}\hoR.  
\end{equation}
En effet, il résulte du théorème de presque-pureté de Faltings que le groupe $\rH^0(\Gamma,\hoR[\frac 1 p](j))$ est nul pour tout $j\not=0$. 
La suite spectrale \eqref{Intro-lv2b} dégénère donc en $\rE_2$. 
Cependant, le groupe $\rH^1(\Gamma,\hoR[\frac 1 p](1))$ n'étant pas nul en général,  
la filtration aboutissement ne se scinde pas en général.

\subsection{}\label{Intro-lv20} 
Nous avons conjecturé l'existence de la suite spectrale \eqref{Intro-lv2b} dans une première version de ce travail. 
Scholze nous a immédiatement informé qu'il savait construire une telle suite spectrale 
en utilisant la filtration de Hodge-Tate relative associée à un morphisme propre et lisse entre espaces adiques 
qu'il a développée avec Caraiani (\cite{cs} 2.2.4). 
Nous n'abordons pas dans ce travail la question de la comparaison de leur filtration avec la nôtre \eqref{Intro-lv2}. 
Bhatt nous a aussi signalé qu'il a une stratégie pour déduire la suite spectrale \eqref{Intro-lv2b} du formalisme général 
de cohomologie prismatique.

He sait construire la filtration de Hodge-Tate relative \eqref{Intro-lv2} 
dans un cadre encore plus général que celui de \ref{sshtrsg6}. Il la déduit
de la variante globale de notre suite spectrale de Hodge-Tate relative \eqref{Intro-ft4} 
et de son résultat de descente cohomologique pour le topos de Faltings établi dans \cite{tongmu1}. 
Notre preuve de \ref{Intro-lv2}, similaire dans l'esprit à la sienne, a été indépendamment suggérée par le rapporteur. 
Nous démontrons pour ce faire un énoncé de descente cohomologique \ref{amtF28}, 
qui s'avère être un cas particulier de celui de He.

\section{La suite spectrale de Hodge-Tate relative: localisation}

\subsection{}\label{Intro-lv3}
Soient $X$ un $S$-schéma lisse, $\ox$ un point géométrique de $X$ au-dessus de $s$, $\uX$ le
localisé strict de $X$ en $\ox$. 
Soit $\oy\rightsquigarrow \ox$ un morphisme de spécialisation, c'est-à-dire un $X$-morphisme $u\colon \oy\rightarrow \uX$. 
Ce dernier induit un $X_\oeta$-morphisme $\oy\rightarrow \uX_\oeta$. Posons $\uGamma=\pi_1(\uX_\eta,\oy)$  
et soit $\cV_\ox$ la catégorie des $X$-schémas affine étales $\ox$-pointés.  
Pour tout objet $U$ de $\cV_\ox$, nous notons $\oR_U$ la $\co_\oK$-algèbre définie comme dans \eqref{Intro-lv1b} pour le $S$-schéma $U$ 
et le point géométrique $\oy\rightarrow U_\oeta$ induit par le morphisme canonique $\uX\rightarrow U$. 
Les $\co_\oK$-algèbres $(\oR_U)_{U\in \cV_\ox}$ forment naturellement un système inductif. Nous notons $\uoR$ sa limite inductive,
\begin{equation}\label{Intro-lv3a}
\uoR=\underset{\underset{U\in \cV_\ox}{\longrightarrow}}{\lim}\ \oR_U,
\end{equation}
et  $\uhoR$ son complété $p$-adique, que nous munissons des actions naturelles de $\uGamma$.

\begin{teo}[cf. \ref{sshtrl6}]\label{Intro-lv4}
Sous les hypothèses de \ref{Intro-lv3}, pour tout morphisme projectif et lisse $g\colon X'\rightarrow X$, posant $\uX'=X'\times_X\uX$, 
il existe une suite spectrale $\uGamma$-équivariante canonique
\begin{equation}\label{Intro-lv4a}
\rE_2^{i,j}=\rH^i(\uX',\Omega^j_{\uX'/\uX})\otimes_{\co_\uX(\uX)}\uhoR[\frac 1 p](-j)
\Rightarrow \rH^{i+j}_\et(X'_\oy,\mZ_p)\otimes_{\mZ_p}\uhoR[\frac 1 p]. 
\end{equation}
\end{teo}

Il résulte encore du théorème de presque-pureté de Faltings que le groupe $\rH^0(\uGamma,\uhoR[\frac 1 p](1))$ est nul. 
La suite spectrale \eqref{Intro-lv4a} dégénère donc en $\rE_2$. 
Cependant, le groupe $\rH^1(\uGamma,\uhoR[\frac 1 p](1))$ n'étant pas nul en général, 
la filtration aboutissement ne se scinde pas en général. 

Nous donnons deux constructions de la suite spectrale \eqref{Intro-lv4a}. L'une \eqref{sshtrl20} est une application de 
la variante globale \eqref{Intro-ft4}. L'autre \eqref{sshtrl6} est analogue à la construction globale.  

\subsection{}
Hyodo a établi dans \cite{hyodo2} le cas particulier de \ref{Intro-lv4} où $X'$ est un schéma abélien sur $X$ et $\ox$
est un point géométrique générique de la fibre spéciale de $X$. 
Il a aussi donné des exemples où la filtration aboutissement ne se scinde pas.

\subsection{}
La suite spectrale \eqref{Intro-lv4a} peut être globalisée sur un topos naturel, dont les points sont paramétrés par les  
morphismes de spécialisation $\oy\rightsquigarrow \ox$ d'un point géométrique $\oy$ de $X_\oeta$ vers un point géométrique
$\ox$ de $X$, à savoir sur le {\em topos de Faltings}. Ce dernier est au c{\oe}ur de la suite spectrale de Hodge-Tate, y compris
dans le cas absolu. Il a été largement étudié dans (\cite{agt} VI). Nous allons brièvement le passer en revue dans la section suivante.

\section{La suite spectrale de Hodge-Tate relative}\label{Intro-ft}

\subsection{}\label{Intro-ft1}
Soit $X$ un $S$-schéma lisse \eqref{Intro-logsch}. 
On note $E$ la catégorie des morphismes $V\rightarrow U$
au-dessus du morphisme canonique $X_\oeta\rightarrow X$, c'est-à-dire les diagrammes commutatifs
\begin{equation}\label{Intro-ft1a}
\xymatrix{V\ar[r]\ar[d]&U\ar[d]\\
X_\oeta\ar[r]&X}
\end{equation}
tels que $U$ soit étale au-dessus de $X$ et que le morphisme canonique $V\rightarrow U_\oeta$ soit {\em fini étale}. 
Il est utile de considérer la catégorie $E$ comme fibrée par le foncteur
\begin{equation}\label{Intro-ft1b}
\pi\colon E\rightarrow \Et_{/X}, \ \ \ (V\rightarrow U)\mapsto U,
\end{equation}
au-dessus du site étale de $X$. 
La fibre de $\pi$ au-dessus d'un objet $U$ de $\Et_{/X}$ est canoniquement équivalente à la catégorie $\Et_{\rf/U_\oeta}$ des morphismes finis 
étales au-dessus de $U_\oeta$. On la munit de la topologie étale et on note $U_{\oeta,\fet}$ le topos associé. 
Si $U_\oeta$ est connexe et si $\oy$ est un point géométrique de $U_\oeta$, alors le topos $U_{\oeta,\fet}$ est équivalent au
topos classifiant du groupe profini $\pi_1(U_\oeta,\oy)$, 
{\em i.e.}, à la catégorie des ensembles discrets munis d'une action à gauche continue de $\pi_1(U_\oeta,\oy)$.

Nous munissons $E$ de la topologie {\em co-évanescente} (\cite{agt} VI.1.10), c'est-à-dire de la topologie engendrée par les 
recouvrements  $\{(V_i\rightarrow U_i)\rightarrow (V\rightarrow U)\}_{i\in I}$ des types suivants~:
\begin{itemize}
\item[(v)] $U_i=U$ pour tout $i\in I$ et $(V_i\rightarrow V)_{i\in I}$ est un recouvrement;
\item[(c)] $(U_i\rightarrow U)_{i\in I}$ est un recouvrement et $V_i=V\times_UU_i$ pour tout $i\in I$. 
\end{itemize}
Le site ainsi obtenu est appelé {\em site de Faltings} de $X$. 
On désigne par $\tE$ et l'on appelle {\em topos de Faltings} de $X$ le topos des faisceaux d'ensembles sur $E$.  
C'est un analogue du topos co-évanescent $X_\et\gtimes_{X_\et}X_{\oeta,\et}$ (\cite{agt} VI.4).

Se donner un faisceau $F$ sur $E$ revient à se donner:
\begin{itemize}
\item[(i)] pour tout objet $U$ de $\Et_{/X}$, un faisceau $F_U$ de $U_{\oeta,\fet}$, à savoir la restriction de $F$ à la fibre 
de $\pi$ au-dessus de $U$;
\item[(ii)] pour tout morphisme $f\colon U'\rightarrow U$ de $\Et_{/X}$, un morphisme $\gamma_f\colon F_U\rightarrow f_{\oeta*}(F_{U'})$. 
\end{itemize}
Ces données sont astreintes à une condition de cocycle pour la composition des morphismes et à une  
condition de recollement pour les recouvrements de $\Et_{/X}$ (\cite{agt} VI.5.10). Un tel faisceau sera noté $\{U\mapsto F_U\}$.

Il existe trois morphismes canoniques de topos 
\begin{equation}\label{Intro-ft2a}
\xymatrix{
&{X_{\oeta,\et}}\ar[d]_-(0.4){\psi}&\\
{X_\et}&{\tE}\ar[l]_-(0.5){\sigma}\ar[r]^-(0.5)\beta&{X_{\oeta,\fet}}}
\end{equation}
tels que 
\begin{eqnarray}
\sigma^*(U)&=&(U_\oeta\rightarrow U)^a, \ \ \ \forall \ U\in \ob(\Et_{/X}),\label{Intro-intro4a}\\
\beta^*(V)&=&(V\rightarrow X)^a, \ \ \ \forall \ V\in \ob(\Et_{\rf/X_\oeta}),\label{Intro-intro4b}\\
\psi^*(V\rightarrow U)&=&V,\ \ \ \forall \ (V\rightarrow U)\in \ob(E),\label{Intro-intro4c}
\end{eqnarray}
où l'exposant $^a$ désigne le faisceau associé. 
Les morphismes $\sigma$ et $\beta$ sont les analogues de la première et de la seconde projection du topos co-évanescent 
$X_\et\gtimes_{X_\et}X_{\oeta,\et}$.  Le morphisme $\psi$ est un analogue du morphisme des cycles co-proches (\cite{agt} VI.4.13).

Chaque morphisme de spécialisation $\oy\rightsquigarrow \ox$ d'un point géométrique $\oy$ de $X_\oeta$ vers un point 
géométrique $\ox$ de $X$ détermine un point de $\tE$ noté $\rho(\oy\rightsquigarrow \ox)$ (\cite{agt} VI.10.18). 
La famille de ces points est conservative (\cite{agt} VI.10.21).

\begin{prop}[cf. \ref{acycloc2}]\label{Intro-ft6}
Pour tout faisceau abélien de torsion, localement constant et constructible $F$ de $X_{\oeta,\et}$, on a $\rR^i\psi_*(F)=0$ pour tout $i\geq 1$.
\end{prop}

Cet énoncé est une conséquence du fait que pour tout point géométrique $\ox$ de $X$ au-dessus de $s$, notant $\uX$ le localisé
strict de $X$ en $\ox$, $\uX_\oeta$ est un schéma $K(\pi,1)$ (\cite{agt} VI.9.21), {\em i.e.},  si $\oy$ est un point géométrique de $\uX_\oeta$, 
pour tout faisceau abélien de torsion, localement constant et constructible $F$ sur $\uX_\oeta$ et tout $i\geq 0$, on a un isomorphisme 
\begin{equation}
\rH^i(\uX_\oeta,F)\stackrel{\sim}{\rightarrow}\rH^i(\pi_1(\uX_\oeta,\oy),F_\oy). 
\end{equation}
Cette propriété a été prouvée par Faltings (\cite{faltings1} Lemma 2.3 page 281)  comme généralisation de résultats d'Artin (\cite{sga4} XI). 
Elle a été ensuite étendue au cas log-lisse par Achinger (\cite{achinger} 9.5).

\subsection{}\label{Intro-ft3}
Pour tout objet $(V\rightarrow U)$ de $E$, on note $\oU^V$ la clôture intégrale de $\oU$ dans $V$ et l'on pose
\begin{equation}\label{Intro-ft3a}
\ocB(V\rightarrow U)=\Gamma(\oU^V,\co_{\oU^V}).
\end{equation} 
Le préfaisceau sur $E$ ainsi défini est en fait un faisceau (\cite{agt} III.8.16). On écrira $\ocB=\{U\mapsto \ocB_U\}$ (cf. \ref{Intro-ft1}).  
Pour tout $X$-schéma étale $U$ qui est affine, la fibre du faisceau $\ocB_U$ de $U_{\oeta,\fet}$ en un point géométrique $\oy$ 
de $U_\oeta$, est la représentation $\oR_U$ de $\pi_1(U_\oeta,\oy)$ définie dans \eqref{Intro-lv1b} pour $U$.  

Pour tout morphisme de spécialisation $\oy\rightsquigarrow \ox$, on a 
\begin{equation}\label{Intro-ft3b}
\ocB_{\rho(\oy\rightsquigarrow \ox)}=\underset{\underset{U\in \cV_\ox}{\longrightarrow}}{\lim}\ \oR_U,
\end{equation}
où $\cV_\ox$ est la catégorie des $X$-schémas étales $\ox$-pointés $U$ qui sont affines.

\subsection{}\label{Intro-ft7}
Pour tout topos $T$, les systèmes projectifs d'objets de $T$ indexés par l'ensemble ordonné des entiers naturels 
$\mN$  forment un topos que nous notons $T^{\mN^\circ}$ (\cite{agt} III.7). 

Pour tout entier $n\geq 0$, on pose $\ocB_n=\ocB/p^n\ocB$. 
Afin de prendre en compte la topologie $p$-adique, on considère la $\co_C$-algèbre $\bvocB=(\ocB_n)_{n\geq 1}$
du topos $\tE^{\mN^\circ}$. 
On travaillera dans la catégorie $\bMod_{\mQ}(\bvocB)$ des $\bvocB$-modules à isogénie près (\cite{agt} III.6.1), 
qui est un analogue global de la catégorie des $\hoR[\frac 1 p]$-représentations de $\Delta$ considérée dans \ref{Intro-lv1}.

\begin{teo}[cf. \ref{sshtr5}]\label{Intro-ft4}
Soit $g\colon X'\rightarrow X$ un morphisme projectif et lisse. Notons
\begin{equation}
\xymatrix{
{X'^{\mN^\circ}_{\oeta,\et}}\ar[r]^{\bvg_\oeta}&{X^{\mN^\circ}_{\oeta,\et}}\ar[r]^\bvpsi&{\tE^{\mN^\circ}}}
\end{equation}
les morphismes induits par $g_\oeta$ et $\psi$ \eqref{Intro-ft2a}, 
et $\bvmZ_p$ la $\mZ_p$-algèbre $(\mZ/p^n\mZ)_{n\geq 1}$ de $X'^{\mN^\circ}_{\oeta,\et}$.  
Alors, il existe une suite spectrale canonique de $\bvocB_\mQ$-modules
\begin{equation}\label{Intro-ft4a}
\rE_2^{i,j}=\sigma^*(\rR^ig_*(\Omega^j_{X'/X}))\otimes_{\sigma^*(\co_X)}\bvocB_\mQ(-j)\Rightarrow \bvpsi_*(\rR^{i+j}\bvg_{\oeta*}(\bvmZ_p))\otimes_{\mZ_p}\bvocB_\mQ.
\end{equation}
\end{teo}

Ce théorème, ainsi que les énoncés \ref{Intro-lv2} et \ref{Intro-lv4}, valent en fait sous l'hypothèse plus générale que {\em $g$ soit propre}. 
En effet, l'hypothèse de projectivité sur $g$ est utilisée dans \ref{Intro-fmct5} ci-dessous
qui s'étend également aux morphismes propres (voir \ref{Intro-fmct50}). 

La suite spectrale \eqref{Intro-ft4a} est appelée la {\em suite spectrale de  Hodge-Tate relative}.
On démontre aisément qu'elle est $G_K$-équivariante pour les structures  $G_K$-équivariantes naturelles sur les topos et les objets 
intervenant dans sa formulation. On en déduit la proposition suivante. 

\begin{prop}[cf. \ref{sshtr11}]\label{Intro-ft5}
Sous les hypothèses de \ref{Intro-ft4}, la suite spectrale de Hodge-Tate relative \eqref{Intro-ft4a} dégénère en $\rE_2$.
\end{prop}

\begin{rema}
Avec une approche différente, Caraiani et Scholze ont construit une filtration de Hodge-Tate relative pour les  morphismes propres et lisses
d'espaces adiques (\cite{cs} 2.2.4). Comme pour la variante locale \eqref{Intro-lv20}, nous n'abordons pas dans ce travail la question de 
la comparaison de leur filtration avec la filtration aboutissement de notre suite spectrale de  Hodge-Tate relative. 
\end{rema}

\section{\texorpdfstring{Les principaux théorèmes de comparaison $p$-adiques de Faltings}{Les principaux théorèmes de comparaison p-adiques de Faltings}}

\subsection{}\label{Intro-fmct0}
Nous conservons dans cette section les hypothèses et notations de §~\ref{Intro-ft}.
On note $\co_{\oK^\flat}$ la limite du système projectif $(\co_\oK/p\co_\oK)_{\mN}$ 
dont les morphismes de transition sont les itérés de l'endomorphisme de Frobenius absolu de $\co_\oK/p\co_\oK$; 
\begin{equation}
\co_{\oK^\flat}= \underset{\underset{\mN}{\longleftarrow}}{\lim}\ \co_\oK/p\co_\oK.
\end{equation}
C'est un anneau de valuation non-discrète complet, parfait, de hauteur $1$ et de caractéristique $p$. 
On fixe une suite $(p_n)_{n\geq 0}$ d'éléments de $\co_{\oK}$ telle que $p_0=p$ et $p_{n+1}^p=p_n$ pour tout $n\geq 0$. 
On note $\varpi$ l'élément associé de $\co_{\oK^\flat}$ et on pose $\xi=[\varpi]-p$  
dans l'anneau $\rW(\co_{\oK^\flat})$ des vecteurs de Witt $p$-typiques de $\co_{\oK^\flat}$.
On a un isomorphisme canonique 
\begin{equation}
\co_C(1)\stackrel{\sim}{\rightarrow} p^{\frac{1}{p-1}}\xi \co_C.
\end{equation}

\begin{teo}[cf. \ref{TPCF15}, \cite{faltings2}]\label{Intro-fmct1}
Supposons que $X$ soit propre sur $S$. Soient $i,n$ deux entiers $\geq 0$, 
$F$ un faisceau localement constant constructible de $(\mZ/p^n\mZ)$-modules de $X_{\oeta,\et}$.
Alors, le noyau et le conoyau du morphisme canonique 
\begin{equation}\label{Intro-fmct1a} 
\rH^i(X_{\oeta,\et},F)\otimes_{\mZ_p}\co_C\rightarrow \rH^i(\tE,\psi_*(F)\otimes_{\mZ_p}\ocB)
\end{equation}
sont annulés par $\fm_C$. 
\end{teo}
On dira que le morphisme \eqref{Intro-fmct1a} est un {\em presque-isomorphisme}. 

C'est le {\em principal théorème de comparaison $p$-adique de  Faltings} à partir duquel se déduisent tous les théorèmes de 
comparaison entre la cohomologie étale  $p$-adique et les autres cohomologies $p$-adiques. C'est aussi le principal ingrédient 
dans la construction de la suite spectrale de Hodge-Tate absolue \eqref{Intro-intro2b}.

Nous revisitons dans ce livre la preuve de Faltings de ce résultat important en donnant plus de détails. 
Elle est basée sur la suite exacte d'Artin-Schreier pour la ``perfection'' de l'anneau $\ocB_1=\ocB/p\ocB$.
Un des ingrédients principaux est un énoncé de structure pour les $\varphi$-modules presque-étales sur $\co_{\oK^\flat}$
vérifiant certaines conditions, incluant une condition de presque-finitude au sens de Faltings.
Dans notre application à la cohomologie du topos de Faltings annelé par la ``perfection'' de $\ocB_1$,
la preuve de cette dernière condition résulte de la combinaison de trois ingrédients:
\begin{itemize}
\item[(i)] des calculs locaux de cohomologie galoisienne utilisant le théorème de presque-pureté de Faltings (\cite{faltings2}, \cite{agt} II.8.17); 
\item[(ii)] une étude fine des conditions de presque-finitude pour les faisceaux quasi-cohérents de modules sur les schémas; 
\item[(iii)] le résultat de Kiehl sur la finitude de la cohomologie pour un morphisme propre (\cite{kiehl} 2.9'a) (cf. \cite{egr1} 1.4.7).
\end{itemize}

\subsection{}\label{Intro-fmct2}
Expliquons maintenant la construction de  Faltings de la suite spectrale de Hodge-Tate absolue \eqref{Intro-intro2b}.  
Supposons que $X$ soit propre sur $S$. 
Par \ref{Intro-ft6}, pour tous $i,n\geq 0$, on a un isomorphisme canonique
\begin{equation}\label{Intro-fmct2a}
\rH^i(X_{\oeta,\et},\mZ/p^n\mZ)\stackrel{\sim}{\rightarrow}\rH^i(\tE,\psi_*(\mZ/p^n\mZ)).
\end{equation} 
Il n'est pas difficile de voir que le morphisme canonique $\mZ/p^n\mZ\rightarrow \psi_*(\mZ/p^n\mZ)$ est un isomorphisme. 
Utilisant alors le principal théorème de comparaison $p$-adique de Faltings \ref{Intro-fmct1}, on obtient un morphisme canonique  
\begin{equation}\label{Intro-fmct2b}  
\rH^i(X_{\oeta,\et},\mZ/p^n\mZ)\otimes_{\mZ_p}\co_C\rightarrow \rH^i(\tE,\ocB_n) 
\end{equation}
qui est un presque-isomorphisme. Pour calculer $\rH^i(\tE,\ocB_n)$, on utilise la suite spectrale de
Cartan-Leray pour le  morphisme $\sigma\colon \tE\rightarrow X_{\et}$ \eqref{Intro-ft2a},
\begin{equation}\label{Intro-fmct2c}
\rE_2^{i,j}=\rH^i(X_\et,\rR^j\sigma_*(\ocB_n))\Rightarrow \rH^{i+j}(\tE,\ocB_n). 
\end{equation}
On en déduit la suite spectrale de Hodge-Tate absolue \eqref{Intro-intro2b} en utilisant l'analogue global suivant du calcul de Faltings de la cohomologie galoisienne.

\begin{teo}[cf. \ref{tfkum11}]\label{Intro-fmct3} 
Il  existe un homomorphisme canonique
de $\co_{\oX_n}$-algèbres graduées de $X_{s,\et}$
\begin{equation}\label{Intro-fmct3a}
\wedge (\xi^{-1}\Omega^1_{\oX_n/\oS_n})\rightarrow \oplus_{i\geq 0}\rR^i\sigma_*(\ocB_n),
\end{equation}
où $\xi$ est l'élément de $\rW(\co_{\oK^\flat})$ défini en \ref{Intro-fmct0}, 
dont le noyau (resp. conoyau) est annulé par $p^{\frac{2d}{p-1}}\fm_\oK$ (resp. $p^{\frac{2d+1}{p-1}}\fm_\oK$), où $d=\dim(X/S)$.  
\end{teo}

Nous prouvons ce résultat en utilisant la théorie de Kummer sur la fibre spéciale du topos de Faltings annelé $(\tE,\ocB)$.

\subsection{}\label{Intro-fmct4} 
Soit $g\colon X'\rightarrow X$ un morphisme lisse \eqref{Intro-logsch}.
On associe à $X'$ des objets similaires à ceux associés à  $X$ in §~\ref{Intro-ft} et on les affecte d'un exposant $^\prime$.  
On a un diagramme commutatif 
\begin{equation}\label{Intro-fmct4a}
\xymatrix{
{X'_{\oeta,\et}}\ar[d]_{g_\oeta}\ar[r]^{\psi'}&{\tE'}\ar[d]^{\Theta}\ar[r]^{\sigma'}&{X'_\et}\ar[d]^g\\
{X_{\oeta,\et}}\ar[r]^\psi&{\tE}\ar[r]^{\sigma}&{X_\et}}
\end{equation}
dans lequel $\Theta$ est défini, pour tout objet $(V\rightarrow U)$ de $E$, par
\begin{equation}\label{Intro-fmct4b}
\Theta^*(V\rightarrow U)=(V\times_{X}X'\rightarrow U\times_XX')^a,
\end{equation}
où l'exposant $^a$ désigne la faisceau associé. 
On a un homomorphisme canonique d'anneaux 
\begin{equation}\label{Intro-fmct4c}
\ocB\rightarrow \Theta_*(\ocB').
\end{equation}

\begin{teo}[cf. \ref{TCFR18}, \cite{faltings2} §~6]\label{Intro-fmct5} 
Supposons que $g\colon X'\rightarrow X$ soit projectif. Soient $i,n$ deux entiers $\geq 0$, 
$F'$ un faisceau localement constant constructible de 
$(\mZ/p^n\mZ)$-modules de $X'_{\oeta,\et}$. Alors, le morphisme canonique 
\begin{equation}\label{Intro-fmct5a} 
\psi_*(\rR^ig_{\oeta*}(F'))\otimes_{\mZ_p}\ocB\rightarrow \rR^i\Theta_*(\psi'_*(F')\otimes_{\mZ_p}\ocB')
\end{equation}
est un presque-isomorphisme. 
\end{teo}

On observera que les faisceaux  $\rR^ig_{\oeta*}(F)$ ($i\geq 0$) sont localement constants constructibles sur $X_\oeta$ 
grâce aux théorèmes de changement de base propre et lisse. 

Faltings a formulé la {\em version relative} de son principal théorème de comparaison $p$-adique dans
\cite{faltings2} et en a très sommairement esquissé une preuve dans l'appendice. 
Certains arguments doivent être modifiés et la preuve donnée dans ce livre requiert bien plus de travail.
Elle est basée sur une étude fine de la structure locale de certains $\varphi$-modules presque-étales
qui est intéressante en elle-même (cf. \ref{elcpam15}). 

\begin{rema}[Note ajoutée le 2 décembre 2022]\label{Intro-fmct50}
Le théorème \ref{Intro-fmct5} vaut en fait sous l'hypothèse plus générale que {\em $g$ soit propre}. 
En effet, la condition de projectivité sur $g$ est utilisée pour prouver un résultat de presque-finitude pour les modules presque-cohérents \eqref{afini15}. 
Alors que dans ce texte, nous nous basons sur les résultats de finitude de \cite{sga6} (plutôt que ceux de \cite{kiehl}),  
He \cite{tongmu3} vient de reprendre la preuve de Kiehl (\cite{kiehl} 2.9'a) et  est parvenu à étendre l'énoncé \ref{afini15} aux morphismes propres.
\end{rema}

\begin{rema}
Scholze a généralisé  \ref{Intro-fmct1} au cas des variétés rigides en utilisant la même stratégie (\cite{scholze1} 1.3). 
Il a également prouvé un analogue de \ref{Intro-fmct5} dans le cadre des espaces adiques et des topos pro-étales (\cite{scholze1} 5.12).
Il le déduit du cas absolu en utilisant un théorème de changement de base dû à Huber.
Notons cependant que le résultat clé de notre preuve dans le cas relatif, à savoir l'étude fine de la structure locale de certains $\varphi$-modules presque-étales, ne semble pas résulter des arguments de Scholze. 
\end{rema}

\subsection{}\label{Intro-fmct6} 
Conservons les hypothèses de \ref{Intro-fmct5} et soit $n$ un entier $\geq 0$. 
Comme le morphisme canonique $\mZ/p^n\mZ\rightarrow \psi'_*(\mZ/p^n\mZ)$ est un isomorphisme, 
pour construire la suite spectrale de Hodge-Tate relative \eqref{Intro-ft4a}, 
nous sommes conduits par \ref{Intro-fmct5} à calculer les faisceaux de cohomologie $\rR^q\Theta_*(\ocB'_n)$ ($q\geq 0$). 
S'inspirant du cas absolu \eqref{Intro-fmct2}, 
le  problème est alors de trouver une factorisation naturelle de $\Theta$, à laquelle appliquer une suite spectrale de Cartan-Leray. 
Considérons le diagramme commutatif de morphismes de topos
\begin{equation}
\xymatrix{
\tE'\ar[rd]^{\sigma'}\ar[d]_{\tau}&\\
{\tE\times_{X_\et}X'_\et}\ar[r]^-(0.4){\pi}\ar[d]_{\Xi}&{X'_\et}\ar[d]^g\\
{\tE}\ar[r]^-(0.5){\sigma}&{X_\et}}
\end{equation}
Nous prouvons que le produit fibré de topos $\tE\times_{X_\et}X'_\et$ est en fait un {\em topos de Faltings relatif}, dont la définition est 
inspirée par celle des produits orientés de topos, au-delà du topos co-évanescent qui inspira déjà notre définition du topos de Faltings usuel.

\section{Topos de Faltings relatif}\label{Intro-rft}

\subsection{}\label{Intro-rft3} 
On conserve les hypothèses et notations de \ref{Intro-fmct4}.
On note $G$ la catégorie des morphismes $(W\rightarrow U\leftarrow V)$
au-dessus des morphismes canoniques $X'\rightarrow X\leftarrow X_\oeta$, c'est-à-dire  les diagrammes commutatifs
\begin{equation}\label{Intro-rft3a}
\xymatrix{W\ar[r]\ar[d]&U\ar[d]&V\ar[l]\ar[d]\\
X'\ar[r]&X&X_\oeta\ar[l]}
\end{equation}
tels que $W$ soit étale sur$X'$, $U$ soit étale sur$X$ et le morphisme canonique  $V\rightarrow U_\oeta$ soit {\em fini étale}. 
On la munit de la topologie engendrée par les recouvrements  
\begin{equation}
\{(W_i\rightarrow U_i\leftarrow V_i)\rightarrow (W\rightarrow U \leftarrow V)\}_{i\in I}
\end{equation}
des trois types suivants~:
\begin{itemize}
\item[(a)] $U_i=U$, $V_i=V$ pour tout  $i\in I$ et $(W_i\rightarrow W)_{i\in I}$ est un recouvrement;
\item[(b)] $W_i=W$, $U_i=U$ pour tout $i\in I$ et $(V_i\rightarrow V)_{i\in I}$ est un recouvrement;
\item[(c)]  les diagrammes
\begin{equation}\label{Intro-rft3b}
\xymatrix{
W'\ar[r]\ar@{=}[d]&U'\ar[d]\ar@{}|\Box[rd]&V'\ar[l]\ar[d]\\
W\ar[r]&U&V\ar[l]}
\end{equation}
dans lequel $U'\rightarrow U$ est un morphisme quelconque et le carré de droite est cartésien. 
\end{itemize}
Le site ainsi défini est appelé {\em site de Faltings relatif} du morphisme $g\colon X'\rightarrow X$. 
On désigne par $\tG$ et l'on appelle  {\em topos de Faltings relatif} de $g$ le topos des faisceaux d'ensembles sur $G$. 
C'est un analogue du produit orienté de topos $X'_\et\gtimes_{X_\et}X_{\oeta,\et}$ (\cite{agt} VI.3).

Il existe deux morphismes canoniques de topos
\begin{equation}\label{Intro-rft3c} 
\xymatrix{
{X'_\et}&{\tG}\ar[l]_-(0.4){\pi}\ar[r]^(0.4)\lambda&{X_{\oeta,\fet}}},
\end{equation}
définis par
\begin{eqnarray}
\pi^*(W)&=&(W\rightarrow X\leftarrow X_\oeta)^a, \ \ \ \forall \ W\in \ob(\Et_{/X'}),\\
\lambda^*(V)&=&(X'\rightarrow X\leftarrow V)^a, \ \ \ \forall \ V\in \ob(\Et_{\rf/X_\oeta}),
\end{eqnarray}
où l'exposant $^a$ désigne le faisceau associé. 
Ce sont des analogues de la première et de la seconde projection du produit orienté  $X'_\et\gtimes_{X_\et}X_{\oeta,\et}$.

Si $X'=X$, $\tG$ est canoniquement équivalent au topos de Faltings $\tE$ \eqref{Intro-ft1}. 
Ainsi, par fonctorialité du topos de Faltings relatif, on a une  factorisation naturelle de  $\Theta\colon \tE'\rightarrow \tE$ en 
\begin{equation}\label{Intro-rft3d} 
\tE'\stackrel{\tau}{\longrightarrow}\tG \stackrel{\lgg}{\longrightarrow} \tE. 
\end{equation}
Ces morphismes s'insèrent dans le diagramme commutatif suivant de morphismes de topos 
\begin{equation}
\xymatrix{
&\tE'\ar[d]^{\tau}\ar[r]^{\beta'}\ar[ld]_{\sigma'} &{X'_{\oeta,\fet}}\ar[d]^{g_\oeta}\\
X'_\et\ar[d]_{g}\ar@{}|\Box[rd]&\tG\ar[d]^-(0.5){\lgg}\ar[r]^-(0.4){\lambda}\ar[l]_-(0.4){\pi}&{X_{\oeta,\fet}}\\
X_\et&\tE\ar[ur]_{\beta}\ar[l]_{\sigma}&}
\end{equation}
On démontre que {\em le carré inférieur gauche est cartésien} \eqref{Intro-fmct6}. 

On a un morphisme canonique 
\begin{equation}
\varrho\colon X'_\et\gtimes_{X_\et}X_{\oeta,\et}\rightarrow \tG.
\end{equation}
Se donner un point de $X'_\et\gtimes_{X_\et}X_{\oeta,\et}$ revient à se donner un point géométrique $\ox'$ de $X'$, 
un point géométrique $\oy$ de $X_\oeta$ et un morphisme de spécialisation $\oy \rightsquigarrow g(\ox')$. 
On notera (abusivement) un tel point par $(\oy \rightsquigarrow \ox')$.  
On démontre que la famille des points $\varrho(\oy \rightsquigarrow \ox')$ de $\tG$ est conservative. 

\subsection{}\label{Intro-rft5} 
Soient $\ox'$ un point géométrique de $X'$, $\uX'$ le localisé strict de $X'$ en $\ox'$, $\uX$ le localisé strict de $X$ en $g(\ox')$. 
On désigne par $\tuG$ le topos de Faltings relatif du morphisme $\uX'\rightarrow \uX$ induit par $g$,  
par $\ulambda\colon \tuG\rightarrow \uX_{\oeta,\fet}$ le morphisme canonique \eqref{Intro-rft3c}
et par $\Phi\colon \tuG\rightarrow \tG$ le morphisme de fonctorialité. 
Il existe une section canonique $\theta$ de $\ulambda$,
\begin{equation}\label{Intro-rft5a}
\xymatrix{
{\uX_{\oeta,\fet}}\ar[r]^{\theta}\ar[rd]_\id&{\tuG}\ar[d]^{\ulambda}\\
&{\uX_{\oeta,\fet}}}
\end{equation}
Nous démontrons que le morphisme de changement de base induit par ce diagramme
\begin{equation}\label{Intro-rft5b}
\ulambda_*\rightarrow \theta^*
\end{equation}
est un isomorphisme. On pose
\begin{equation}\label{Intro-rft5d}
\phi_{\ox'}=\theta^*\circ \Phi^*\colon \tG\rightarrow \uX_{\oeta,\fet}.
\end{equation}

Si $\oy$ est un point géométrique de $\uX_{\oeta}$, on obtient naturellement un point $(\oy \rightsquigarrow \ox')$ de 
\mbox{$X'_\et\gtimes_{X_\et}X_{\oeta,\et}$}. Alors, pour tout faisceau $F$ de $\tG$,  on a un isomorphisme canonique et fonctoriel
\begin{equation}\label{Intro-rft5e}
F_{\varrho(\oy \rightsquigarrow \ox')}\stackrel{\sim}{\rightarrow} \phi_{\ox'}(F)_{\oy}.
\end{equation}

\begin{prop}[cf. \ref{tfr30}] \label{Intro-rft6} 
Sous les hypothèses de \ref{Intro-rft5}, pour tout faisceau abélien $F$ de $\tG$ et tout $q\geq 0$, on a un isomorphisme canonique
\begin{equation}
\rR^q\pi_*(F)_{\ox'}\stackrel{\sim}{\rightarrow}\rH^q(\uX_{\oeta,\fet},\phi_{\ox'}(F)).
\end{equation}
\end{prop}

\begin{cor}[cf. \ref{ktfr19}]
Soient $(\oy\rightsquigarrow \ox')$ un point de $X'_\et\gtimes_{X_\et}X_{\oeta,\et}$, 
$\uX'$ le localisé strict de $X'$ en $\ox'$, $\uX$ le localisé strict de $X$ en $g(\ox')$,
$\ug\colon \uX'\rightarrow \uX$ le morphisme induit par $g$, 
\begin{equation}
\varphi'_{\ox'}\colon \tE'\rightarrow \uX'_{\oeta,\fet}
\end{equation}
le morphisme canonique analogue de \eqref{Intro-rft5d}. 
Alors, pour tout groupe abélien $F$ de $\tE'$ et tout $q\geq 0$, on a un isomorphisme fonctoriel canonique 
\begin{equation}
(\rR^q\tau_*(F))_{\varrho(\oy\rightsquigarrow \ox')}\stackrel{\sim}{\rightarrow}
\rR^q\ug_{\oeta,\fet*}(\varphi'_{\ox'}(F))_{\oy}.
\end{equation}
\end{cor}

\subsection{}\label{Intro-rft7} 
Considérons l'anneau suivant de $\tG$,
\begin{equation}
\ocB^!=\tau_*(\ocB').
\end{equation}
On a des homomorphismes canoniques $\ocB\rightarrow \lgg_*(\ocB^!)$ et $\hbar'_*(\co_{\oX'})\rightarrow \pi_*(\ocB^!)$,
dans lequel $\hbar'\colon \oX'\rightarrow X'$ désigne la projection canonique. Par conséquent, on peut considérer $\lgg$ et $\pi$ comme des
morphismes de topos annelés. 

Pour tout point $(\oy \rightsquigarrow \ox')$ de  $X'_\et\gtimes_{X_\et}X_{\oeta,\et}$, nous démontrons que l'anneau
$\ocB^!_{\varrho(\oy \rightsquigarrow \ox')}$ est normal et strictement hensélien. De plus, l'homomorphisme canonique
$\co_{\oX',\ox'}\rightarrow \ocB^!_{\varrho(\oy \rightsquigarrow \ox')}$ est local et injectif. 

\subsection{}\label{Intro-rft11} 
Supposons que $X=\Spec(R)$ et $X'=\Spec(R')$ soient affines. Soient $\oy'$ un point géométrique de $X'_\oeta$, $\Delta'=\pi_1(X'_\oeta,\oy')$,
$(W_j)_{j\in J}$ le revêtement universel de $X'_\oeta$ en $\oy'$, $\oy=g_\oeta(\oy')$, $\Delta=\pi_1(X_\oeta,\oy)$,
$(V_i)_{i\in I}$ le revêtement universel de $X_\oeta$ en $\oy$.  
Pour tout $i\in I$, $(V_i\rightarrow X)$ est naturellement un objet de $E$ et pour tout $j\in J$, 
$(W_j\rightarrow X')$ est naturellement un objet de $E'$. Posons 
\begin{eqnarray}
\oR&=&\underset{\underset{i\in I}{\longrightarrow}}{\lim}\ \ocB(V_i\rightarrow X),\label{Intro-rft11a} \\
\oR'&=&\underset{\underset{j\in J}{\longrightarrow}}{\lim}\ \ocB'(W_j\rightarrow X'). \label{Intro-rft11b} 
\end{eqnarray}
On retrouve les $\co_\oK$-algèbres définis en \eqref{Intro-lv1b}. On les munit des actions naturelles de  $\Delta$ et de $\Delta'$. 
Pour tout $i\in I$, il existe un $X'$-morphisme canonique $\oy'\rightarrow X'\times_XV_i$. On désigne par 
$V'_i$ la composante irréductible de $X'\times_XV_i$ contenant $\oy'$ et par $\Pi_i$ le sous-groupe correspondant de $\Delta'$. 
Alors, $(V'_j\rightarrow X')$ est naturellement un objet de $E'$. 
Posons 
\begin{eqnarray}
\oR^!&=&\underset{\underset{i\in I}{\longrightarrow}}{\lim}\ \ocB'(V'_i\rightarrow X'),\label{Intro-rft11d} \\
\Pi&=&\bigcap_{i\in I}\Pi_i.\label{Intro-rft11c}
\end{eqnarray}
On a un homomorphisme canonique  $\oR\rightarrow \oR^! \rightarrow \oR'$.

Pour tout point géométrique $\ox'$ et tout morphisme de spécialisation $\oy\rightsquigarrow g(\ox')$, 
nous démontrons qu'il existe un isomorphisme canonique (déterminé par le choix de $\oy'$)
\begin{equation}\label{Intro-rft11e} 
\ocB^!_{\varrho(\oy\rightsquigarrow \ox')}\stackrel{\sim}{\rightarrow}
\underset{\underset{\ox'\rightarrow U'\rightarrow U}{\longrightarrow}}{\lim}\ \oR^!_{U'\rightarrow U},
\end{equation}
où la limite inductive est prise sur la catégorie des morphismes $\ox'\rightarrow U'\rightarrow U$ au-dessus de $\ox'\rightarrow X'\rightarrow X$, 
avec $U'$ affine  étale sur $X'$,  $U$ affine  étale sur $X$  et $\oR^!_{U'\rightarrow U}$ l'anneau correspondant \eqref{Intro-rft11d}.

\begin{prop}[cf. \ref{eccr51}]\label{Intro-rft12} 
Conservons les hypothèses de  \ref{Intro-rft11} et supposons de plus que $g$ s'insère dans un diagramme commutatif 
\begin{equation}
\xymatrix{
X'\ar[r]^-(0.5){\iota'}\ar[d]_g&{\mG_{m,S}^{d'}}\ar[d]^\gamma\\
X\ar[r]^-(0.5){\iota}&{\mG_{m,S}^{d}}}
\end{equation}
dans lequel les morphismes $\iota$ et $\iota'$ sont étales, $d$ et $d'$ sont des entiers $\geq 0$ et $\gamma$ 
est un homomorphisme de tores au-dessus de $S$ qui est lisse \eqref{Intro-logsch}. Soit $n$ un entier $\geq 0$. Alors,
\begin{itemize}
\item[{\rm (i)}] Il existe un homomorphisme canonique de $\oR^!$-algèbres graduées
\begin{equation}
\wedge (\Omega^1_{R'/R}\otimes_{R'}(\oR^!/p^n\oR^!)(-1))\rightarrow \oplus_{i\geq 0}\rH^i(\Pi,\oR'/p^n\oR'),
\end{equation}
qui est presque-injectif et dont le conoyau est tué par $p^{\frac{1}{p-1}}\fm_\oK$. 
\item[{\rm (ii)}] Le $\oR^!$-module $\rH^i(\Pi,\oR'/p^n\oR')$ est presque de présentation finie pour tout $i\geq 0$, et il est
presque-nul pour tout $i\geq r+1$, où $r=\dim(X'/X)$.
\end{itemize}
\end{prop}

C'est une version relative du calcul de Faltings de la cohomologie galoisienne de  $\oR$ qui repose sur son théorème de 
presque-pureté (\cite{faltings2} Theorem 4 page 192, \cite{agt} II.6.16). L'énoncé peut se globaliser de la manière suivante en 
utilisant la théorie de Kummer sur la fibre spéciale du topos annelé $(\tE',\ocB')$. 

\begin{teo}[cf. \ref{crtf4}]\label{Intro-rft8} 
Pour tout entier $n\geq 1$, il existe un  homomorphisme canonique de $\ocB^!$-algèbres graduées de $\tG$
\begin{equation}
\wedge (\pi^*(\xi^{-1}\Omega^1_{\oX'_n/\oX_n}))\rightarrow \oplus_{i\geq 0}\rR^i\tau_*(\ocB'_n),
\end{equation}
où $\pi^*$ désigne l'image inverse par le morphisme de topos annelés $\pi\colon (\tG,\ocB^!)\rightarrow (X'_\et,\hbar'_*(\co_{\oX'}))$, 
dont le noyau (resp. conoyau) est annulé par $p^{\frac{2r}{p-1}}\fm_\oK$ (resp. $p^{\frac{2r+1}{p-1}}\fm_\oK$), où $r=\dim(X'/X)$.  
\end{teo} 
 
Considérons ensuite la suite spectrale de Cartan-Leray  
\begin{equation}\label{Intro-rft13a}
\rE_2^{i,j}=\rR^i\lgg_*(\rR^j\tau_*(\ocB'_n))\Rightarrow \rR^{i+j}\Theta_*(\ocB'_n).
\end{equation} 
Tenant compte de \ref{Intro-rft8}, pour obtenir la suite spectrale de Hodge-Tate relative \eqref{Intro-ft4a}, on a besoin de démontrer un théorème de 
changement de base relativement au diagramme cartésien
\begin{equation}\label{Intro-rft13b}
\xymatrix{
\tG\ar[r]^{\pi}\ar[d]_{\lgg}&{X'_\et}\ar[d]^g\\
{\tE}\ar[r]^{\sigma}&{X_\et}}
\end{equation}

\begin{teo}[cf. \ref{ktfr25}] \label{Intro-rft14} 
Supposons que le morphisme $g$ soit propre. Alors, pour tout faisceau abélien de torsion $F$ de $X'_\et$ et tout $q\geq 0$, le morphisme 
de changement de base 
\begin{equation}
\sigma^*(\rR^qg_*(F))\rightarrow \rR^q\lgg_*(\pi^*(F))
\end{equation}
est un isomorphisme.
\end{teo}

La preuve est inspirée d'un théorème de changement de base pour les produits orientés dû à Gabber. 
Elle se ramène au théorème de changement de base propre pour le topos étale.

Nous avons en fait besoin d'une variante du théorème \ref{Intro-rft14} pour les modules quasi-cohérents. 
Pour ce faire, nous enrichissons le diagramme commutatif \eqref{Intro-rft13b} en un diagramme commutatif de morphismes de topos annelés
\begin{equation}\label{Intro-rft15b}
\xymatrix{
{(\tG,\ocB^!)}\ar[r]^-(0.5){\pi}\ar[d]_{\lgg}&{(X'_\et,\co_{X'})}\ar[d]^g\\
{(\tE,\ocB)}\ar[r]^-(0.5){\sigma}&{(X_\et,\co_X)}}
\end{equation}

\begin{prop}[cf. \ref{ktfr27}] \label{Intro-rft15} 
Pour tout entier $n\geq 0$, l'homomorphisme canonique
\begin{equation}
\ocB_n\boxtimes_{\co_X}\co_{X'}\rightarrow \ocB^!_n,
\end{equation}
où le produit tensoriel externe d'anneaux est relatif au diagramme cartésien \eqref{Intro-rft13b}, est un presque-isomorphisme.
\end{prop}

\begin{teo}[cf. \ref{ktfr31}] \label{Intro-rft16} 
Supposons le morphisme $g$ propre.
Il existe alors un entier $N\geq 0$ tel que pour tous entiers $n\geq 1$ et $q\geq 0$, 
et tout $\co_{X'_n}$-module quasi-cohérent $\cF'$ \eqref{Intro-intro1a}, 
le noyau et le conoyau du morphisme de changement de base {\rm (\cite{egr1} (1.2.3.3))} relativement au diagramme \eqref{Intro-rft15b},
\begin{equation}
\sigma^*(\rR^qg_*(\cF'))\rightarrow \rR^q\lgg_*(\pi^*(\cF')),
\end{equation}
où $\sigma^*$ et $\pi^*$ désignent les images inverses au sens des topos annelés, soient annulés par $p^N$.
\end{teo}

\begin{prop}[cf. \ref{ktfr30}] \label{Intro-rft17} 
Soient $n, q$ des entiers $\geq 0$, $\cF'$ un $\co_{X'_n}$-module cohérent qui est $X_n$-plat \eqref{Intro-intro1a}. 
Supposons que le morphisme $g$ soit propre et que pour tout entier $i\geq 0$, 
le $\co_{X_n}$-module $\rR^ig_*(\cF')$ soit localement libre (de type fini). 
Alors, le morphisme de changement de base relativement au diagramme \eqref{Intro-rft15b},
\begin{equation}
\sigma^*(\rR^q g_*(\cF'))\rightarrow \rR^q\lgg_*(\pi^*(\cF')),
\end{equation}
où $\sigma^*$ et $\pi^*$ désignent les images inverses au sens des topos annelés, est un presque-isomorphisme.
\end{prop}

\subsection{}
Soient $n,q$ des entiers $\geq 0$. 
Comme le morphisme canonique $\mZ/p^n\mZ\rightarrow \psi'_*(\mZ/p^n\mZ)$ est un isomorphisme,  
on obtient à partir de \ref{Intro-fmct5}, pour tout $q\geq 0$, un morphisme canonique
\begin{equation}
\psi_*(\rR^qg_{\oeta*}(\mZ/p^n\mZ))\otimes_{\mZ_p}\ocB\rightarrow \rR^q\Theta_*(\ocB'_n),
\end{equation}
qui est un presque-isomorphisme. 
La suite spectrale de Hodge-Tate relative \eqref{Intro-ft4a} se déduit alors de la  suite spectrale de Cartan-Leray
\eqref{Intro-rft13a} en utilisant \ref{Intro-rft8} et \ref{Intro-rft16}.\\

\chapter{Préliminaires}\label{prelim}

\section{Notations et conventions}

{\em Tous les anneaux considérés dans ce livre possèdent un élément unité~;
les homomorphismes d'anneaux sont toujours supposés transformer l'élément unité en l'élément unité.
Nous considérons surtout des anneaux commutatifs, et lorsque nous parlons d'anneau 
sans préciser, il est sous-entendu qu'il s'agit d'un anneau commutatif~; en particulier, 
il est sous-entendu, lorsque nous parlons d'un topos annelé $(X,A)$ sans préciser, que $A$ est commutatif.

On sous-entend par {\em monoïde}  un monoïde commutatif et unitaire. 
Les homomorphismes de monoïdes sont toujours supposés transformer l'élément unité en l'élément unité.}

\subsection{}\label{notconv15}
Soit $p$ un nombre premier. On munit  $\mZ_p$ de la topologie $p$-adique, ainsi que toutes les $\mZ_p$-algèbres adiques 
({\em i.e.}, les $\mZ_p$-algèbres complètes et séparées pour la topologie $p$-adique).
Soient $A$ une $\mZ_p$-algèbre adique, $i\colon A\rightarrow A[\frac 1 p]$ l'homomorphisme canonique. 
On appelle {\em topologie $p$-adique} sur $A[\frac 1 p]$ l'unique topologie 
compatible avec sa structure de groupe additif  pour laquelle les sous-groupes $i(p^nA)$,
pour $n\in \mN$, forment un système fondamental de voisinages de $0$ (\cite{tg} chap.~III §1.2, prop.~1). 
Elle fait de $A[\frac 1 p]$ un anneau topologique. 
Soient $M$ un $A[\frac 1 p]$-module de type fini, $M^\circ$ un sous-$A$-module de type fini de $M$ 
qui l'engendre sur $A[\frac 1 p]$. On appelle {\em topologie $p$-adique} sur $M$ l'unique topologie 
compatible avec sa structure de groupe additif  pour laquelle les sous-groupes $p^nM^\circ$,
pour $n\in \mN$, forment un système fondamental de voisinages de $0$. Cette topologie ne dépend pas
du choix de $M^\circ$. En effet, si $M'$ est un autre sous-$A$-module de type fini de $M$ qui l'engendre sur $A[\frac 1 p]$,
alors il existe $m\geq 0$ tel que $p^mM^\circ\subset M'$ et $p^mM'\subset M^\circ$. 
Il est clair que $M$ est un $A[\frac 1 p]$-module topologique.

\subsection{}\label{notconv16}
Soient $G$ un groupe profini, $A$ un anneau topologique muni d'une action continue de $G$ 
par des homomorphismes d'anneaux. Une {\em $A$-représentation} de $G$ est la donnée d'un $A$-module 
$M$ et d'une action $A$-semi-linéaire de $G$ sur $M$, {\em i.e.},  telle que 
pour tous $g\in G$, $a\in A$ et $m\in M$, on ait $g(am)=g(a)g(m)$.
On dit que la $A$-représentation est {\em continue} si $M$ est un $A$-module topologique 
et si l'action de $G$ sur $M$ est continue. Soient $M$, $N$ deux $A$-représentations 
(resp. deux $A$-représentations continues) de $G$. 
Un morphisme de $M$ dans $N$ est la donnée d'un morphisme $A$-linéaire et $G$-équivariant 
(resp. $A$-linéaire, continu et $G$-équivariant) de $M$ dans $N$. 
On note $\bRep_A(G)$ (resp. $\bRep_A^\cont(G)$)
la catégorie des $A$-représentations (resp. $A$-représentations continues) de $G$.
Si $M$ et $N$ sont deux $A$-représentations de $G$, 
les $A$-modules $M\otimes_AN$ et $\Hom_A(M,N)$ sont naturellement des $A$-représentations de $G$. 

\subsection{}\label{notconv1}
Soient $A$ un anneau, $p$ un nombre premier, $n$ un entier $\geq1$.  
On désigne par $\rW(A)$ (resp. $\rW_n(A)$) l'anneau des vecteurs de Witt 
(resp. vecteurs de Witt de longueur $n$) à coefficients dans $A$ relatif à $p$. 
On a un homomorphisme d'anneaux
\begin{equation}\label{notconv1a}
\Phi_n\colon 
\begin{array}[t]{clcr}
\rW_n(A)&\rightarrow& A,\\ 
(x_0,\dots,x_{n-1})&\mapsto&x_0^{p^{n-1}}+p x_1^{p^{n-2}}+\dots+p^{n-1}x_{n-1}
\end{array}
\end{equation}
appelé $n$-ième composante fantôme. 
On dispose aussi des morphismes de restriction, de décalage et de Frobenius
\begin{eqnarray}
\rR\colon \rW_{n+1}(A)&\rightarrow& \rW_n(A),\label{notconv1b}\\
\rV\colon \rW_n(A)&\rightarrow& \rW_{n+1}(A),\label{notconv1c}\\
\rF\colon \rW_{n+1}(A)&\rightarrow& \rW_n(A).\label{notconv1d}
\end{eqnarray}
Lorsque $A$ est de caractéristique $p$, $\rF$ induit un endomorphisme de $\rW_n(A)$, encore noté $\rF$.

\subsection{} \label{notconv2}
Pour tout anneau $R$ et tout monoïde $M$, on désigne par $R[M]$ la $R$-algèbre de $M$ et par 
$e\colon M\rightarrow R[M]$ l'homomorphisme canonique, où $R[M]$ est considéré comme un
monoïde multiplicatif. Pour tout $x\in M$, on notera $e^x$ au lieu de $e(x)$. 

On désigne par $\bA_M$ le schéma $\Spec(\mZ[M])$ muni de la structure logarithmique associée à 
la structure pré-logarithmique définie par $e\colon M\rightarrow \mZ[M]$ (\cite{agt} II.5.9). 
Pour tout homomorphisme de monoïdes $\vartheta\colon M\rightarrow N$, 
on note $\bA_\vartheta\colon \bA_N\rightarrow \bA_M$ le morphisme de schémas logarithmiques associé.

\subsection{}\label{notconv3}
Dans tout ce livre, on fixe un univers $\mU$ possédant un élément de cardinal infini, 
et un univers $\mV$ tel que $\mU\in \mV$. 
On appelle catégorie des $\mU$-ensembles et l'on note $\Ens$, 
la catégorie des ensembles qui se trouvent dans $\mU$. 
C'est un $\mU$-topos ponctuel (\cite{sga4} IV 2.2). 
On désigne par $\Sch$ la catégorie des schémas éléments de $\mU$. 
Sauf mention explicite du contraire, il sera sous-entendu que les anneaux et 
les schémas logarithmiques (et en particulier les schémas)
envisagés dans ce livre sont éléments de l'univers~$\mU$. 
On désigne par $\Top$ la $2$-catégorie des $\mU$-topos appartenant à $\mV$ (\cite{sga4} IV 3.3.1).

\subsection{}\label{notconv4}
Pour une catégorie $\cC$, nous notons $\ob(\cC)$ l'ensemble de ses objets,
$\cC^\circ$ la catégorie opposée, et pour $X,Y\in \ob(\cC)$, 
$\Hom_\cC(X,Y)$ (ou $\Hom(X,Y)$ lorsqu'il n'y a aucune ambiguïté) 
l'ensemble des morphismes de $X$ dans $Y$. 

Si $\cC$ et $\cC'$ sont deux catégories, nous désignons par $\Hom(\cC,\cC')$ 
l'ensemble des foncteurs de $\cC$ dans $\cC'$, et  
par $\bHom(\cC,\cC')$ la catégorie des foncteurs de $\cC$ dans $\cC'$. 

Soient $\cE$ une catégorie, $\cC$ et $\cC'$ deux catégories sur $\cE$ (\cite{sga1} VI 2). 
Nous notons $\Hom_{\cE}(\cC,\cC')$ l'ensemble des $\cE$-foncteurs de $\cC$ dans $\cC'$ 
et $\Hom_{\cart/\cE}(\cC,\cC')$ l'ensemble des foncteurs cartésiens (\cite{sga1} VI 5.2).
Nous désignons par $\bHom_{\cE}(\cC,\cC')$ la catégorie des $\cE$-foncteurs de $\cC$ dans $\cC'$ et 
par $\bHom_{\cart/\cE}(\cC,\cC')$ la sous-catégorie pleine formée des foncteurs cartésiens.

\subsection{}\label{notconv5}
Soit $\cC$ une catégorie. On désigne par $\hcC$ la catégorie des préfaisceaux 
de $\mU$-ensembles sur $\cC$, c'est-à-dire la catégorie des foncteurs
contravariants sur $\cC$ à valeurs dans $\Ens$ (\cite{sga4} I 1.2). 
Si $\cC$ est munie d'une topologie (\cite{sga4} II 1.1), on désigne par $\tcC$ le topos 
des faisceaux de $\mU$-ensembles sur $\cC$ (\cite{sga4} II 2.1). 

Pour $F$ un objet de $\hcC$, on note $\cC_{/F}$ la catégorie 
suivante (\cite{sga4} I 3.4.0). Les objets de $\cC_{/F}$
sont les couples formés d'un objet $X$ de $\cC$ 
et d'un morphisme $u$ de $X$ dans $F$. Si $(X,u)$ et $(Y,v)$ sont deux objets, 
un morphisme de $(X,u)$ vers $(Y,v)$ est un morphisme $g\colon X\rightarrow Y$
tel que $u=v\circ g$.

\subsection{}\label{notconv6}
Soient $\cC$ un $\mU$-site (\cite{sga4} II 3.0.2), $\tcC$ le topos des faisceaux de $\mU$-ensembles sur $\cC$,
$X$ un préfaisceau de $\mU$-ensembles sur $\cC$. On munit la catégorie $\cC_{/X}$ de la topologie induite par la topologie de $\cC$ au moyen 
du foncteur ``source'' $j_X\colon \cC_{/X}\rightarrow \cC$ (\cite{sga4} III 3.1).
D'après (\cite{sga4} III 5.2), $j_X$ est un foncteur continu et cocontinu. 
Il définit donc une suite de trois foncteurs adjoints~: 
\begin{equation}\label{notconv6a}
j_{X!}\colon (\cC_{/X})^\sim \rightarrow \tcC,\ \ \ \
j_X^*\colon \tcC \rightarrow (\cC_{/X})^\sim,\ \ \ \
j_{X*}\colon  (\cC_{/X})^\sim\rightarrow \tcC,
\end{equation}
dans le sens que pour deux foncteurs consécutifs de la suite, celui
de droite est adjoint à droite de l'autre. 
Le foncteur $j_{X!}$ se factorise par la catégorie $\tcC_{/X^a}$, 
où $X^a$ est le faisceau associé à $X$, 
et le foncteur induit $(\cC_{/X})^\sim \rightarrow \tcC_{/X^a}$
est une équivalence de catégories (\cite{sga4} III 5.4). 
Le couple de foncteurs $(j_X^*,j_{X*})$ définit un morphisme de topos qu'on notera aussi abusivement
$j_{X}\colon \tcC_{/X^a} \rightarrow \tcC$, dit morphisme de localisation de $\tcC$ en $X^a$ (cf. \cite{sga4} IV 5.2).
Pour tout objet $F$ de $\tcC$, on pose 
\begin{equation}\label{notconv6b}
F|X^\tta=j_X^*(F).
\end{equation}

Supposons, de plus, que les limites projectives finies soient représentables dans $\cC$ et que $X$ soit un objet de $\cC$.
Soit $e_\cC$ un objet final de $\cC$ (qui existe par hypothèse). 
On note $j_X^+\colon \cC\rightarrow \cC_{/X}$ le foncteur de changement de base par le morphisme canonique
$X\rightarrow e_\cC$. Alors, $j_X^+$ est exact à gauche et continu (\cite{sga4} III 1.6 et 3.3).
Comme $j_X^*$ prolonge $j_X^+$ (\cite{sga4} III 5.4), le morphisme de localisation
$j_X\colon \tcC_{/X^a} \rightarrow \tcC$ s'identifie canoniquement au morphisme de topos associé à $j_X^+$.

\subsection{}\label{notconv7} 
Soient $\cC$ un $\mU$-site, $\hcC$ la catégorie des préfaisceaux de $\mU$-ensembles sur $\cC$,
$\tcC$ le topos des faisceaux de $\mU$-ensembles sur $\cC$.
On désigne par $\Mon_\hcC$ (resp. $\Mon_\tcC$) la catégorie des monoïdes (commutatifs et unitaires) 
de $\hcC$ (resp. $\tcC$) (cf. \cite{sga4} I 3.2, resp. II 6.3.1). 
Les $\mU$-limites inductives et projectives dans $\Mon_\hcC$ sont représentables et se calculent terme à terme.
Les $\mU$-limites inductives et projectives dans $\Mon_\tcC$ sont représentables.  
Le foncteur d'injection canonique 
\begin{equation}
i\colon \Mon_\tcC\rightarrow \Mon_\hcC
\end{equation} 
admet un adjoint à gauche
\begin{equation}
a\colon \Mon_\hcC\rightarrow \Mon_\tcC, \ \ \ M\mapsto M^a,
\end{equation}
qui est exact à gauche (\cite{sga4} II 6.4). 
D'après (\cite{sga4} II 6.5), le foncteur ``faisceau d'ensembles sous-jacent'' $\Mon_\tcC\rightarrow \tcC$ admet un adjoint à gauche 
et commute donc aux limites projectives. De plus, pour tout préfaisceau de monoïdes $M$ sur $\cC$, 
le faisceau d'ensembles sous-jacent à $M^a$ est canoniquement isomorphe au faisceau d'ensembles associé
au préfaisceau d'ensembles sous-jacent à $M$. 

On note $1_\tcC$ le monoïde unité de $\tcC$, c'est-à-dire le faisceau de monoïdes associé au préfaisceau de monoïdes 
sur $\cC$ de valeur $1$. C'est un objet initial et final de la catégorie $\Mon_\tcC$.
On appelle {\em conoyau} d'un morphisme $u\colon N\rightarrow M$ de $\Mon_\tcC$ 
la somme amalgamée de $u$ et du morphisme canonique $N\rightarrow 1_\tcC$. 
On a un isomorphisme canonique 
\begin{equation}
\coker(u)\stackrel{\sim}{\rightarrow} a(\coker(i(u))).
\end{equation}
Si $u$ est un monomorphisme, on appelle $\coker(u)$ le {\em quotient} de $M$ par $N$ et on le note $M/N$.

\subsection{}\label{notconv8} 
Pour tout morphisme de topos $f\colon Y\rightarrow X$, les foncteurs $f^*$ et $f_*$ induisent un couple de foncteurs adjoints
$f^*\colon \Mon_X\rightarrow \Mon_Y$ et $f_*\colon \Mon_Y\rightarrow \Mon_X$. Le foncteur $f^*$ est exact à gauche
d'après \ref{notconv7}.

\subsection{}\label{notconv13}
Soit $X$ un $\mU$-topos. Les systèmes projectifs d'objets de $X$ indexés par l'ensemble ordonné 
des entiers naturels $\mN$, forment un topos que l'on note $X^{\mN^\circ}$. On renvoie à  
(\cite{agt} III.7) pour des sorites utiles sur ce type de topos. On rappelle, en particulier, qu'on a un morphisme de topos 
\begin{equation}\label{notconv13a}
\uplambda\colon X^{\mN^\circ}\rightarrow X,
\end{equation}
dont le foncteur image inverse $\uplambda^*$ associe à tout objet $F$ de $X$ le système projectif constant de valeur $F$,
et dont le foncteur image directe $\uplambda_*$ associe à tout système projectif sa limite projective (\cite{agt} III.7.4).

\subsection{}\label{notconv9}
Soit $(X,A)$ un $\mU$-topos annelé (\cite{sga4} IV 11.1.1). On note $\bMod(A)$ ou $\bMod(A,X)$ 
la catégorie des $A$-modules de $X$. Si $M$ est un $A$-module, on désigne par $\rS_A(M)$ 
(resp. $\wedge_A(M)$, resp. $\Gamma_A(M)$) l'algèbre symétrique (resp. extérieure, resp. à puissances divisées) 
de $M$ et pour tout entier $n\geq 0$, par $\rS_A^n(M)$ (resp. $\wedge_A^n(M)$,
resp. $\Gamma_A^n(M)$) sa partie homogène de degré $n$. 
Les formations de ces algèbres commutent à la localisation au-dessus d'un objet de $X$.

\begin{defi}[\cite{sga6} I 1.3.1]\label{notconv14}
Soit $(X,A)$ un topos annelé. 
On dit qu'un $A$-module $M$ de $X$ est {\em localement projectif de type fini} si les conditions 
équivalentes suivantes sont satisfaites :
\begin{itemize}
\item[{\rm (i)}] $M$ est de type fini et le foncteur $\cHom_A(M,\cdot)$ est exact~;
\item[{\rm (ii)}] $M$ est de type fini et tout épimorphisme de $A$-modules $N\rightarrow M$ admet localement une section~;
\item[{\rm (iii)}] $M$ est localement facteur direct d'un $A$-module libre de type fini. 
\end{itemize}
\end{defi}

Lorsque $X$ a suffisamment de points et que pour tout point $x$ de $X$, la fibre de $A$ en $x$
est un anneau local, les $A$-modules localement projectifs de type fini sont les $A$-modules 
localement libres de type fini (\cite{sga6} I 2.15.1).

\subsection{}\label{notconv17}
Soient $X$ un $\mU$-topos appartenant à $\mV$ \eqref{notconv3}, $G$ un groupe (du topos ponctuel). 
On note $\mG$ le groupoïde associé à $G$, {\em i.e.},
la catégorie ayant un seul objet, de classe de morphismes $G$. 
On appelle {\em action à gauche de $G$ sur $X$} la donnée d'un pseudo-foncteur normalisé  
du groupoïde $\mG$ dans la $2$-catégorie $\Top$ des $\mU$-topos appartenant à $\mV$ \eqref{notconv3},
qui fait correspondre $X$ à l'unique objet de $\mG$ (cf. \cite{sga1} VI §~8).  
Concrètement, cela revient à se donner pour tout $g\in G$, un morphisme $\gamma_g\colon X\rightarrow X$ et pour 
tout $(f,g)\in G^2$, un isomorphisme $c_{f,g}\colon \gamma_g^*\gamma_f^*\stackrel{\sim}{\rightarrow} \gamma_{fg}^*$
satisfaisant aux conditions suivantes:
\begin{itemize}
\item[(i)] notant $e$ l'élément neutre de $G$, on a $\gamma_e=\id_X$;
\item[(ii)] pour tout $f\in G$, on a $c_{f,e}=c_{e,f}=\id_{\gamma_f^*}$;
\item[(iii)] pour tout $(f,g,h)\in G^3$ et tout $F\in \ob(X)$, on a 
\begin{equation}
c_{f,gh}(F)\cdot c_{g,h}(\gamma_f^*(F))=c_{fg,h}(F)\cdot \gamma_h^*(c_{f,g}(F)).
\end{equation}
\end{itemize}
Ces conditions sont celles envisagées dans (\cite{sga1} VI 7.4). 
La condition (iii) se visualise par la commutativité du diagramme 
\begin{equation}\label{notconv17a}
\xymatrix{
{\gamma_h^*(\gamma_g^*(\gamma_f^*(F)))}\ar[rr]^-(0.5){h^*(c_{f,g}(F))}\ar[d]_{c_{g,h}(\gamma_f^*(F))}&&
{\gamma_h^*(\gamma_{fg}^*(F))}\ar[d]^{c_{fg,h}(F)}\\
{\gamma_{gh}^*(\gamma_f^*(F))}\ar[rr]^-(0.5){c_{f,gh}(F)}&&{\gamma_{fgh}^*(F)}}
\end{equation}

D'après (cf. \cite{sga1} VI §~8), se donner une action à gauche de $G$ sur $X$ revient à de se donner un foncteur fibrant clivé et normalisé
\begin{equation}\label{notconv17b}
\mX\rightarrow \mG,
\end{equation}
dont la catégorie fibre au-dessus de l'unique objet de $\mG$ est $X$. C'est un topos fibré dans le sens de (\cite{sga4} VI 7.1.1): 
pour tout $f\in G$, le foncteur image inverse par $f$ est le foncteur image inverse par le morphisme de topos $\gamma_f$.
Il est commode de noter $\gamma_f$ encore $f\colon X\rightarrow X$ ou $f_X$, ce qui n'induit aucun risque d'ambiguïté.  

\'Etant donné une action à gauche de $G$ sur $X$, on appelle {\em faisceau $G$-équivariant} (ou {\em objet $G$-équivariant}) 
de $X$ la donnée d'une section cartésienne du foncteur fibrant 
\eqref{notconv17b} (cf. \cite{agt} II.4.16). D'après (\cite{sga1} VI §12), il revient au même de se donner un objet $F$ de $X$ 
et pour tout $f\in G$, un isomorphisme 
\begin{equation}\label{notconv17c}
\tau^F_f\colon F\stackrel{\sim}{\rightarrow}f^*(F)
\end{equation}
tels que $\tau^F_\id=\id_F$ et que pour tout $(f,g)\in G^2$, on ait 
\begin{equation}\label{notconv17d}
\tau^F_{fg}=c_{g,f}\circ g^*(\tau^F_f)\circ \tau^F_{g}.
\end{equation}
On dira aussi que $F$ est muni d'une structure de faisceau $G$-équivariant
(ou abusivement qu'il est un faisceau $G$-équivariant).  
On laissera le soin au lecteur de décrire explicitement les morphismes entre faisceaux $G$-équivariants de $X$.

Un {\em groupe (resp. anneau) $G$-équivariant de $X$} est la donnée d'une section cartésienne de la catégorie fibrée 
au-dessus de $\mG$ des groupes (resp. anneaux) de $\mX$ \eqref{notconv17b}. 
Il revient au même de se donner un faisceau $G$-équivariant $F$ de $X$
qui est un groupe (resp. anneau) tel que pour tout $f\in G$, 
l'isomorphisme $\tau^F_f$ \eqref{notconv17c} soit un homomorphisme. 

Si $A$ est un anneau $G$-équivariant de $X$, un 
{\em $A$-module $G$-équivariant de $X$} est la donnée d'une section cartésienne de la catégorie fibrée 
au-dessus de $\mG$ des $A$-modules de $\mX$ \eqref{notconv17b}. 
Il revient au même de se donner un faisceau $G$-équivariant $M$ de $X$
qui soit un $A$-module tel que pour tout $f\in G$, 
l'isomorphisme $\tau^M_f$ soit $\tau_f^A$-linéaire \eqref{notconv17c}.

\subsection{}\label{notconv18}
Soient $G$ un groupe, $u\colon X\rightarrow Y$ un morphisme de $\mU$-topos munis d'actions à gauche de $G$. 
On désigne par $\mG$ le groupoïde associé à $G$ et par $\mX\rightarrow \mG$ et $\mY\rightarrow \mG$ les topos fibrés définis par 
les actions de $G$ sur $X$ et $Y$ respectivement \eqref{notconv17b}. 
On dit que $u$ est {\em $G$-équivariant} si le foncteur image inverse $u^*\colon Y\rightarrow X$ est induit par un $\mG$-foncteur cartésien $\muu^*\colon \mY\rightarrow \mX$. 
D'après (\cite{sga1} VI §12), il revient au même de demander que pour tout $f\in G$, il existe un isomorphisme 
$\varphi_f\colon f_X^*u^*\rightarrow u^*f_Y^*$ vérifiant les conditions a') et b') de {\em loc. cit.}
L'isomorphisme $\varphi_f$ rend donc commutatif le diagramme 
\begin{equation}\label{notconv18a}
\xymatrix{
X\ar[r]^u\ar[d]_{f_X}&Y\ar[d]^{f_Y}\\
X\ar[r]^u&Y}
\end{equation}
Supposons dans la suite de ce numéro que le morphisme $u$ soit $G$-équivariant. 

Pour tout faisceau $G$-équivariant $H$ de $Y$, $u^*(H)$ est canoniquement 
muni d'une structure de faisceau $G$-équivariant de $X$, à savoir l'image par $\muu^*$ de la $\mG$-section cartésienne de $\mY$ qui définit $H$. 
Concrètement, pour tout $f\in G$, $\tau_f^{u^*(H)}$ est l'isomorphisme composé de 
\begin{equation}\label{notconv18b}
u^*(H) \stackrel{\sim}{\rightarrow} u^*(f^*(H))\stackrel{\sim}{\rightarrow} f^*(u^*(H)),
\end{equation}
où la première flèche est induite par l'isomorphisme $\tau^H_f$,
et la seconde flèche est l'isomorphisme canonique \eqref{notconv18a}. 

Pour tout $f\in G$, le diagramme commutatif \eqref{notconv18a} induit un isomorphisme de changement de base 
$f_Y^*u_*\stackrel{\sim}{\rightarrow}u_*f_X^*$. On vérifie aussitôt que ces isomorphismes vérifient 
les conditions a') et b') de (\cite{sga1} VI §12) et induisent donc un $\mG$-foncteur cartésien $\muu_*\colon \mX\rightarrow \mY$
qui prolonge le foncteur $u_*\colon X\rightarrow Y$. 

Pour tout faisceau $G$-équivariant $F$ de $X$, $u_*(F)$ est canoniquement 
muni d'une structure de faisceau $G$-équivariant de $Y$, à savoir l'image par $\muu_*$ de la $\mG$-section cartésienne de $\mX$ qui définit $F$. 
Concrètement, pour tout $f\in G$, $\tau_f^{u_*(F)}$ est l'inverse de l'adjoint de l'isomorphisme composé de 
\begin{equation}\label{notconv18c}
u_*(F) \stackrel{\sim}{\rightarrow} u_*(f_*(F))\stackrel{\sim}{\rightarrow} f_*(u_*(F)),
\end{equation}
où la première flèche est induite par l'isomorphisme $F\stackrel{\sim}{\rightarrow} f_*(F)$ adjoint de $(\tau^F_f)^{-1}$,
et la seconde flèche est l'isomorphisme canonique \eqref{notconv18a}. 

Soient $F$ un faisceau abélien $G$-équivariant de $X$, $j$ un entier $\geq 0$. 
Pour tout $f\in G$, on considère le morphisme
\begin{equation}\label{notconv18d}
\rR^ju_*(F)\rightarrow f_*(\rR^ju_*(F))
\end{equation}
composé de
\begin{equation}\label{notconv18e}
\rR^ju_*(F) \stackrel{\sim}{\rightarrow} \rR^ju_*(f_*(F))\rightarrow \rR^j(uf)_*(F)\stackrel{\sim}{\rightarrow} \rR^j(fu)_*(F)
\rightarrow f_*(\rR^ju_*(F)),
\end{equation}
où la première flèche est induite par l'isomorphisme $F\stackrel{\sim}{\rightarrow} f_*(F)$ adjoint de $(\tau^F_f)^{-1}$,
la deuxième et la quatrième flèches sont les edge-homomorphismes de la suite spectrale de Cartan-Leray (\cite{sga4} V 5.4) 
et la troisième flèche est l'isomorphisme canonique \eqref{notconv18a}. 
Comme $f$ est une équivalence de topos, \eqref{notconv18d} est un isomorphisme. On note $\tau_f^{\rR^ju_*(F)}$ 
l'isomorphisme inverse de son adjoint. On vérifie aussitôt que ces morphismes munissent $\rR^ju_*(F)$ d'une structure de 
faisceau abélien $G$-équivariant sur $X$. 

Si $A$ est un anneau $G$-équivariant de $X$ et $M$ un $A$-module, $u_*(A)$ est un anneau $G$-équivariant de $Y$, 
et pour tout entier $j\geq 0$, $\rR^ju_*(M)$ est un $u_*(A)$-module $G$-équivariant de $Y$. 

\subsection{}\label{notconv19}
Soient $G$ un groupe, $X$ un $\mU$-topos muni d'une action à gauche de $G$,
$u\colon X\rightarrow \Ens$ la projection canonique dans le topos final (\cite{sga4} IV 4.3). 
Le morphisme $u$ est essentiellement déterminé par la donnée d'un objet final $T$ de $X$, 
à savoir $u^*(e)$ où $e$ est un ensemble ponctuel.
Comme $T$ est final, il est muni d'une structure, essentiellement unique, de faisceau $G$-équivariant.  
Il s'ensuit aussitôt que le morphisme $u$ est $G$-équivariant lorsque l'on munit le topos $\Ens$ de l'action triviale de $G$. 

D'après \ref{notconv18}, pour tout faisceau $G$-équivariant $F$ de $X$, $\Gamma(X,F)$ est naturellement 
un ensemble $G$-équivariant, autrement dit, un ensemble muni d'une action de $G$. 

De même, pour tout faisceau abélien $G$-équivariant $F$ de $X$ et tout entier $j\geq 0$, 
$\rH^i(X,F)$ est naturellement un groupe abélien $G$-équivariant, autrement dit, un groupe abélien muni 
d'une action $\mZ$-linéaire de $G$. 

Pour tous faisceaux $G$-équivariants $Z,F$ de $X$, l'ensemble $\Hom(Z,F)$ est naturellement muni d'une action à gauche de $G$. 
Concrètement, à tous $u\in \Hom(Z,F)$ et $f\in G$, posant $f'=f^{-1}$, 
le transformé $u^f$ de $u$ par $f$ est le morphisme composé défini comme suit:
\begin{eqnarray}
v\colon Z\stackrel{u}{\longrightarrow} F \stackrel{(\tau^F_{f'})'}{\longrightarrow} f'_*(F),\\
u^f\colon Z\stackrel{\tau_{f'}^Z}{\longrightarrow} f'^*(Z) \stackrel{v'}{\longrightarrow} F,
\end{eqnarray}
où $(\tau^F_{f'})'$ est l'adjoint de $(\tau^F_{f'})^{-1}=f'^*(\tau^F_f)$ et $v'$ est l'adjoint de $v$.

\subsection{}\label{notconv10}
Pour tout schéma $X$, on désigne par $\Et_{/X}$ le {\em site étale} de $X$, 
c'est-à-dire, la sous-catégorie pleine de $\Sch_{/X}$ \eqref{notconv3} formée des schémas étales sur $X$,
munie de la topologie étale; c'est un $\mU$-site. 
On note $X_\et$ le {\em topos étale} de $X$, c'est-à-dire le topos des faisceaux de $\mU$-ensembles sur $\Et_{/X}$.

On désigne par $\Et_{\coh/X}$ (resp.
$\Et_{\scoh/X}$) la sous-catégorie pleine de $\Et_{/X}$ formée des schémas étales 
de présentation finie sur $X$ (resp. étales, séparés et de présentation finie sur $X$), 
munie de la topologie induite par celle de $\Et_{/X}$; ce sont des sites $\mU$-petits. 
Si $X$ est quasi-séparé, le foncteur de restriction de $X_\et$ dans le topos des faisceaux de $\mU$-ensembles sur 
$\Et_{\coh/X}$ (resp. $\Et_{\scoh/X}$) est une équivalence de catégories (\cite{sga4} VII 3.1 et 3.2).

On désigne par $\Et_{\rf/X}$ la sous-catégorie pleine de $\Et_{/X}$ formée des schémas étales finis sur $X$, 
munie de la topologie induite par celle de $\Et_{/X}$; c'est un site $\mU$-petit. 
On appelle {\em topos fini étale} de $X$ et on note  $X_\fet$, le topos des faisceaux de $\mU$-ensembles sur $\Et_{\rf/X}$ 
(cf. \cite{agt} VI.9.2). L'injection canonique $\Et_{\rf/X}\rightarrow \Et_{/X}$ induit un morphisme de topos 
\begin{equation}\label{notconv10a}
\rho_X\colon X_\et\rightarrow X_\fet.
\end{equation}

\subsection{}\label{notconv12}
Soit $X$ un schéma. On désigne par $X_\zar$ le topos de Zariski de $X$ et par 
\begin{equation}\label{notconv12b}
u_X\colon X_{\et}\rightarrow X_{\zar}
\end{equation}
le morphisme canonique (\cite{sga4} VII 4.2.2). Si $F$ est un $\co_X$-module quasi-cohérent de $X_\zar$, 
on désigne par $\iota(F)$ le faisceau de $X_\et$ défini pour tout $X$-schéma étale $U$ par (\cite{sga4} VII 2 c))
\begin{equation}\label{notconv12c}
\iota(F)(U)=\Gamma(U,F\otimes_{\co_X}\co_U).
\end{equation}
Il est commode, lorsqu'il n'y aucune risque de confusion, de désigner $\iota(F)$ abusivement par $F$. 
On notera que $\iota(\co_X)$ est un anneau de $X_\et$ et que $\iota(F)$ est un $\iota(\co_X)$-module. 

Notant $\bMod^\qcoh(\co_X,X_\zar)$ la sous-catégorie pleine de $\bMod(\co_X,X_\zar)$ formée des $\co_X$-modules quasi-cohérents \eqref{notconv9},
la correspondance $F\mapsto \iota(F)$ définit un foncteur
\begin{equation}\label{notconv12a}
\iota\colon \bMod^\qcoh(\co_X,X_\zar)\rightarrow \bMod(\co_X,X_\et).
\end{equation}

Pour tout $\co_X$-module quasi-cohérent $F$ de $X_\zar$, on a un isomorphisme canonique 
\begin{equation}\label{notconv12d}
F\stackrel{\sim}{\rightarrow}u_{X*}(\iota(F)).
\end{equation}
Nous considérons donc $u_X$ comme un morphisme de topos annelés 
\begin{equation}\label{notconv12j}
u_X\colon (X_{\et},\co_X)\rightarrow (X_{\zar},\co_X). 
\end{equation}
Nous utilisons pour les modules la notation $u^{-1}_X$ pour désigner l'image
inverse au sens des faisceaux abéliens et nous réservons la notation 
$u^*_X$ pour l'image inverse au sens des modules. L'isomorphisme \eqref{notconv12d} induit par adjonction un morphisme
\begin{equation}\label{notconv12e}
u^*_X(F)\rightarrow \iota(F).
\end{equation}
Celui-ci est un isomorphisme. En effet, pour tout point géométrique $\ox$ de $X$, on a un isomorphisme canonique
\begin{equation}\label{notconv12f}
\iota(F)_\ox\stackrel{\sim}{\rightarrow}\underset{\underset{\ox\rightarrow U}{\longleftarrow}}{\lim}\ \Gamma(U,F\otimes_{\co_X}\co_U),
\end{equation}
où la limite est prise sur les voisinages de $\ox$ dans $\Et_{/X}$ (on peut évidemment se limiter aux voisinages affines au-dessus d'un ouvert affine de $X$ 
contenant l'image $x$ de $\ox$). Notant $X'$ le localisé strict de $X$ en $\ox$, on en déduit, d'après (\cite{ega4} 8.5.2(i)),
un isomorphisme canonique 
\begin{equation}\label{notconv12g}
\iota(F)_\ox\stackrel{\sim}{\rightarrow} \Gamma(X',F\otimes_{\co_X}\co_{X'}).
\end{equation}
Par ailleurs, on a un isomorphisme canonique
\begin{equation}\label{notconv12h}
u_X^*(F)_\ox\stackrel{\sim}{\rightarrow} F_x\otimes_{\co_{X,x}} \iota(\co_X)_{\ox}.
\end{equation}
La fibre du morphisme \eqref{notconv12e} en $\ox$ s'identifie au morphisme canonique
\begin{equation}\label{notconv12i}
F_x\otimes_{\co_{X,x}} \Gamma(X',\co_{X'}) \rightarrow  \Gamma(X',F\otimes_{\co_X}\co_{X'}),
\end{equation}
qui est un isomorphisme. Par suite, \eqref{notconv12e} est un isomorphisme.  On en déduit que $\iota$ est un adjoint à gauche du foncteur 
$u_{X*}$. La flèche d'adjonction $\id\rightarrow u_{X*}\circ \iota$ est un isomorphisme d'après \eqref{notconv12d}. En particulier, $\iota$ est pleinement fidèle. 

\subsection{}\label{notconv120}
Soient $f\colon Y\rightarrow X$ un morphisme de schémas, 
\begin{eqnarray}
\jmath\colon \bMod^\qcoh(\co_Y,Y_\zar)&\rightarrow &\bMod(\co_Y,Y_\et),\\
\iota\colon \bMod^\qcoh(\co_X,X_\zar)&\rightarrow& \bMod(\co_X,X_\et)
\end{eqnarray}
les foncteurs canoniques \eqref{notconv12a}. 
Pour tout $\co_X$-module quasi-cohérent $F$, il existe un morphisme canonique fonctoriel
\begin{equation}\label{notconv120a}
\iota(F)\rightarrow f_*(\jmath(f^*(F))).
\end{equation}

Le morphisme $f$ induit donc un morphisme de topos annelés que l'on note encore  
\begin{equation}\label{notconv120b}
f\colon (Y_\et,\co_Y)\rightarrow (X_\et,\co_X).
\end{equation}
Nous utilisons pour les modules la notation $f^{-1}$ pour désigner l'image
inverse au sens des faisceaux abéliens et nous réservons la notation $f^*$ pour l'image inverse au sens des modules.

Compte tenu de \eqref{notconv12g}, le morphisme \eqref{notconv120a} induit par adjonction un isomorphisme
\begin{equation}\label{notconv120c}
f^*(\iota(F))\stackrel{\sim}{\rightarrow} \jmath(f^*(F)).
\end{equation}

Supposons $f$ quasi-compact et quasi-séparé. Pour tout $\co_Y$-module quasi-cohérent $G$, 
le $\co_X$-module $f_*(G)$ est quasi-cohérent et d'après (\cite{ega1n} 9.3.2), on a un isomorphisme canonique fonctoriel 
\begin{equation}\label{notconv120d}
\iota(f_*(G))\stackrel{\sim}{\rightarrow}f_*(\jmath(G)).
\end{equation}
D'après (\cite{sga4} V 5.1 et VII 4.3), on en déduit pour tout entier $q\geq 0$, un isomorphisme canonique 
\begin{equation}\label{notconv120e}
\iota(\rR^q f_*(G))\stackrel{\sim}{\rightarrow}\rR^qf_*(\jmath(G)).
\end{equation}

\subsection{}\label{notconv11}
Soient $X$ un schéma connexe, $\ox$ un point géométrique de $X$.
On désigne par 
\begin{equation}\label{notconv11a}
\omega_\ox\colon \Et_{\rf/X}\rightarrow \Ens
\end{equation}
le foncteur fibre en $\ox$, qui à tout revêtement étale $Y$ de $X$ associe l'ensemble des points géométriques de 
$Y$ au-dessus de $\ox$, par  $\pi_1(X,\ox)$ le groupe fondamental de $X$ en $\ox$ (c'est-à-dire
le groupe des automorphismes du foncteur $\omega_\ox$) et par $\bB_{\pi_1(X,\ox)}$ 
le topos classifiant du groupe profini $\pi_1(X,\ox)$, 
c'est-à-dire la catégorie des $\mU$-ensembles discrets munis d'une action continue à gauche de $\pi_1(X,\ox)$ 
(\cite{sga4} IV 2.7). Alors $\omega_\ox$ induit un foncteur pleinement fidèle 
\begin{equation}\label{notconv11b}
\mu_\ox^+\colon \Et_{\rf/X}\rightarrow \bB_{\pi_1(X,\ox)}
\end{equation}
d'image essentielle la sous-catégorie pleine $\cC(\pi_1(X,\ox))$ de $\bB_{\pi_1(X,\ox)}$ formée des ensembles finis
(\cite{sga1} V  §~4 et §~7).
D'autre part, une famille 
$(Y_\lambda\rightarrow Y)_{\lambda\in \Lambda}$ de $\Et_{\rf/X}$ est couvrante pour la topologie étale, 
si et seulement si son image dans $\bB_{\pi_1(X,\ox)}$ est surjective, 
ou ce qui revient au même, couvrante pour la topologie canonique de $\bB_{\pi_1(X,\ox)}$. 
Par suite, la topologie étale sur $\Et_{\rf/X}$ est induite par la topologie canonique de 
$\bB_{\pi_1(X,\ox)}$ (\cite{sga4} III 3.3). 
Comme les objets de $\cC(\pi_1(X,\ox))$ forment une famille génératrice de $\bB_{\pi_1(X,\ox)}$, le foncteur 
\begin{equation}\label{notconv11e}
\mu_\ox\colon \bB_{\pi_1(X,\ox)}\rightarrow X_\fet
\end{equation}
qui à tout objet $G$ de $\bB_{\pi_1(X,\ox)}$ (vu comme faisceau représentable) associe sa restriction à $\Et_{\rf/X}$,
est une équivalence de catégories en vertu de (\cite{sga4} IV 1.2.1).

Soit $(X_i)_{i\in I}$ un système projectif sur un ensemble ordonné filtrant $I$ dans $\Et_{\rf/X}$
qui pro-représente $\omega_\ox$, normalisé par le fait que les morphismes de transition $X_i\rightarrow X_j$
$(i\geq j)$ sont des épimorphismes et que tout épimorphisme $X_i\rightarrow X'$ de $\Et_{\rf/X}$ 
est équivalent à un épimorphisme $X_i\rightarrow X_j$ $(j\leq i)$ convenable. 
Un tel pro-objet est essentiellement unique. Il est appelé {\em revêtement universel normalisé de $X$ en $\ox$} 
ou le {\em pro-objet fondamental normalisé de $\Et_{\rf/X}$ en $\ox$}. Considérons le foncteur 
\begin{equation}\label{notconv11c}
\nu_\ox\colon X_\fet\rightarrow \bB_{\pi_1(X,\ox)},
\ \ \ F\mapsto \underset{\underset{i\in I}{\longrightarrow}}{\lim}\ F(X_i).
\end{equation}
Par définition, la restriction de $\nu_\ox$ à $\Et_{\rf/X}$ est canoniquement isomorphe au foncteur $\mu_\ox^+$,
et on a un isomorphisme canonique de foncteurs
\begin{equation}\label{notconv11d}
\nu_\ox\circ \mu_\ox\stackrel{\sim}{\rightarrow} \id.
\end{equation}  
Comme $\mu_\ox$ est une équivalence de catégories, $\nu_\ox$ 
est une équivalence de catégories quasi-inverse de $\mu_\ox$. 
On l'appelle {\em le foncteur fibre} de $X_\fet$ en $\ox$.

\subsection{}\label{notconv20}
Soient $X$ un schéma connexe, $\Pi(X)$ son groupoïde fondamental. 
Les objets de $\Pi(X)$ sont les points géométriques de $X$. Pour tout point géométrique $\ox$ de $X$, on note 
$w_\ox\colon \Et_{\rf/X}\rightarrow \Ens$ le foncteur fibre correspondant \eqref{notconv11a}.  
Si $\ox$ et $\ox'$ sont deux points géométriques de $X$, 
l'ensemble $\pi_1(X,\ox,\ox')$ des morphismes  de $\ox$ vers $\ox'$ dans $\Pi(X)$ est l'ensemble des morphismes 
(ou ce qui revient au même des isomorphismes) $w_{\ox}\rightarrow w_{\ox'}$ de foncteurs fibres associés. 
Pour tout point géométrique $\ox$ de $X$, le foncteur  fibre 
$w_\ox$ est exact à gauche et transforme familles couvrantes en familles surjectives. 
Il se prolonge donc en un foncteur fibre $v_\ox\colon X_\fet\rightarrow \Ens$ (\cite{sga4} IV 6.3). 
Celui-ci se déduit du foncteur $\nu_\ox$ défini dans \eqref{notconv11c} par oubli de l'action de $\pi_1(X,\ox)$. 
D'après (\cite{agt} (VI.9.9.1)), $v_\ox$ correspond au point $\rho_X(\ox)$ de $X_\fet$, 
où $\rho_{X}\colon X_\et\rightarrow X_\fet$ est le foncteur canonique \eqref{notconv10a}. 
Interprétant $\Pi(X)$ comme la catégorie opposée à la catégorie des pro-objets fondamentaux
normalisés de $\Et_{\rf/X}$ (\cite{sga1} V 5.7), pour tous points géométriques $\ox$ et $\ox'$ de $X$, 
tout morphisme de $\pi_1(X,\ox,\ox')$ induit un morphisme 
$v_\ox\rightarrow v_{\ox'}$ des foncteurs fibres associés de $X_\fet$. On en déduit un foncteur 
\begin{equation}\label{notconv20a}
\Pi(X)\rightarrow \Pt(X_\fet)^\circ, \ \ \ \ox\mapsto v_\ox,
\end{equation}
où $\Pt(X_\fet)$ est la catégorie des points de $X_\fet$. Celui-ci est une équivalence de catégories en vertu de (\cite{sga4} IV 4.9.4 et 7.2.5). 

D'après (\cite{agt} VI.9.11), le foncteur 
\begin{equation}\label{notconv20b}
X_\fet\rightarrow \bHom(\Pi(X),\Ens), \ \ \ F\mapsto (\ox\mapsto v_\ox(F)),
\end{equation}
déduit de \eqref{notconv20a} induit une équivalence entre la catégorie $X_\fet$ et la catégorie des préfaisceaux de
$\mU$-ensembles $\varphi$ sur $\Pi(X)^\circ$ tels que pour tout point géométrique $\ox$ de $X$, 
l'action de $\pi_1(X,\ox)$ sur l'ensemble $\varphi(\ox)$ 
induite par son action canonique sur l'objet $\ox$ de $\Pi(X)$, soit continue pour la topologie discrète.

\subsection{}\label{notconv21}
Soient $X$ un schéma connexe muni d'une action à gauche d'un groupe (abstrait) $G$, $\Pi(X)$ le groupoïde fondamental de $X$. 
L'action de $G$ sur $X$ induit une action à gauche de $G$ sur la catégorie $\Pi(X)$, 
plus précisément un foncteur du groupoïde $\mG$ associé à $G$ ({\em i.e.}, la catégorie ayant un seul objet, de classe de morphismes $G$) 
dans la $2$-catégorie des catégories appartenant à $\mV$ \eqref{notconv3} qui fait correspondre $\Pi(X)$ au seul objet de $\mG$. 
Concrètement, tout $f\in G$ induit 
un foncteur 
\begin{equation}\label{notconv21a}
\gamma_f\colon \Pi(X)\rightarrow \Pi(X), \ \ \ \ox\mapsto f(\ox)
\end{equation}
satisfaisant aux conditions suivantes:
\begin{itemize}
\item[(i)] notant $e$ l'élément neutre de $G$, on a $\gamma_e=\id_{\Pi(X)}$;
\item[(ii)] pour tout $(f,f')\in G^2$, on a $\gamma_{f}\gamma_{f'}= \gamma_{ff'}$.
\end{itemize}

L'action à gauche de $G$ sur le groupoïde $\Pi(X)$ induit une action à droite de $G$ sur la catégorie  
$\bHom(\Pi(X),\Ens)$ des préfaisceaux de $\mU$-ensembles sur $\Pi(X)^\circ$. 
On appelle préfaisceau {\em $G$-équivariant} sur $\Pi(X)^\circ$, la donnée d'une section cartésienne du foncteur fibrant
\begin{equation}\label{notconv21c}
\mH\rightarrow \mG
\end{equation}
défini par l'action de $G$ sur $\bHom(\Pi(X),\Ens)$.
D'après (\cite{sga1} VI §12), il revient au même de se donner un préfaisceau $\varphi$ sur $\Pi(X)$
et pour tout $f\in G$, un isomorphisme 
\begin{equation}\label{notconv21d}
\tau^\varphi_f\colon \varphi\stackrel{\sim}{\rightarrow}\varphi\circ \gamma_f
\end{equation}
tels que $\tau^\varphi_\id=\id_\varphi$ et que pour tout $(f,f')\in G^2$, on ait 
\begin{equation}\label{notconv21e}
\tau^\varphi_{ff'}=(\tau^\varphi_{f}\circ \gamma_{f'})\circ \tau^\varphi_{f'}.
\end{equation}
On dira aussi que $\varphi$ est muni d'une structure de préfaisceau $G$-équivariant.  

L'action de $G$ sur $X$ induit une action à gauche de $G$ sur le topos $X_\fet$ \eqref{notconv17}. 
Le foncteur canonique \eqref{notconv20b}
\begin{equation}\label{notconv21f}
X_\fet\rightarrow \bHom(\Pi(X),\Ens), \ \ \ F\mapsto (\ox\mapsto v_\ox(F)),
\end{equation}
est $G$-équivariant lorsque l'on fait agir $G$ à droite sur la catégorie $X_\fet$ par image inverse;
{\em i.e.}, ce foncteur induit un $\mG$-foncteur cartésien entre les foncteurs fibrants clivés et normalisés définis par les actions de $G$ 
sur $X_\fet$ et $\bHom(\Pi(X),\Ens)$.
Il induit une équivalence entre la catégorie des faisceaux $G$-équivariants de $X_\fet$ 
et la catégorie des préfaisceaux $G$-équivariants de
$\mU$-ensembles $\varphi$ sur $\Pi(X)^\circ$ tels que pour tout point géométrique $\ox$ de $X$, 
l'action de $\pi_1(X,\ox)$ sur l'ensemble $\varphi(\ox)$ 
induite par son action canonique sur l'objet $\ox$ de $\Pi(X)$, soit continue pour la topologie discrète.

\subsection{}\label{notconv22}
Soient $K$ un corps, $\oK$ une clôture séparable de $K$, $G$ le groupe de Galois de $\oK$ sur $K$, $\eta=\Spec(K)$, $\oeta=\Spec(\oK)$, 
$X$ un $\eta$-schéma géométriquement connexe, $\oX=X\times_\eta\oeta$, $\varpi\colon \oX\rightarrow X$ la projection canonique,
$\Pi(X)$ (resp. $\Pi(\oX)$) le groupoïde fondamental de $X$ (resp. $\oX$) \eqref{notconv20}. 
Le groupe $G$ agit naturellement à gauche sur $\oX$ et $\varpi$ est $G$-équivariant lorsque l'on fait agir $G$ trivialement sur $X$.  
On en déduit une action à gauche de $G$ sur le topos $\oX_\fet$ \eqref{notconv17},  
une action à gauche de $G$ sur le groupoïde $\Pi(\oX)$ et une action à droite de $G$ sur la catégorie  
$\bHom(\Pi(\oX),\Ens)$ des préfaisceaux de $\mU$-ensembles sur $\Pi(\oX)^\circ$ \eqref{notconv21}. 

Considérons le foncteur
\begin{equation}\label{notconv22c}
\upvarpi \colon \Pi(\oX)\rightarrow \Pi(X), \ \ \ \ox\mapsto \upvarpi(\ox)=\varpi(\ox).
\end{equation}
Pour pout point géométrique $\ox$ de $\oX$, notons
$\oX_{(\ox)}$ le localisé strict de $\oX$ en $\ox$ et $X_{(\varpi(\ox))}$ le localisé strict de $X$ en $\varpi(\ox)$. 
On a alors un isomorphisme canonique 
\begin{equation}\label{notconv22d}
\bigsqcup_{f\in G}\oX_{(f(\ox))} \stackrel{\sim}{\rightarrow} X_{(\varpi(\ox))}\times_\eta\oeta.
\end{equation}
On en déduit que pour tous points géométriques $\ox$ et $\ox'$ de $\oX$, on a un isomorphisme canonique 
\begin{equation}\label{notconv22e}
\bigsqcup_{f\in G}\Hom_{\oX}(\ox',\oX_{(f(\ox))}) \stackrel{\sim}{\rightarrow}  \Hom_{X}(\varpi(\ox'),X_{(\varpi(\ox))}),
\end{equation}
compatible à la composition des flèches de spécialisation (\cite{sga4} VIII 7.4).
Compte tenu de l'équivalence \eqref{notconv20a}, le foncteur $\upvarpi$ induit un isomorphisme 
\begin{equation}\label{notconv22f}
\bigsqcup_{f\in G}\pi_1(\oX,f(\ox),\ox')\stackrel{\sim}{\rightarrow}  \pi_1(X,\varpi(\ox),\varpi(\ox')),
\end{equation}
compatible avec la composition dans les catégories $\Pi(\oX)$ et $\Pi(X)$. 
 
Tout préfaisceau d'ensembles (resp. de groupes) $\varphi$ sur $\Pi(X)^\circ$ induit par composition avec $\upvarpi$ 
un préfaisceau d'ensembles (resp. de groupes) $G$-équivariant 
$\psi=\varphi\circ \upvarpi$ sur $\Pi(\oX)^\circ$. 
Concrètement, pour tout $f\in G$ et tout point géométrique $\ox$ de $\oX$, l'isomorphisme 
\begin{equation}\label{notconv22g}
\tau_f^\psi \colon \psi(\ox)\rightarrow \psi(f(\ox))
\end{equation}
est défini par l'identité en identifiant la source et le but avec $\varphi(\upvarpi(\ox))$.

\begin{lem}\label{notconv23}
Sous les hypothèses de \ref{notconv22}, le foncteur qui à tout préfaisceau d'ensembles (resp. de groupes) $\varphi$ sur $\Pi(X)^\circ$ associe le 
préfaisceau d'ensembles (resp. de groupes) $G$-équivariant $\varphi\circ \upvarpi$ sur $\Pi(\oX)^\circ$ est une équivalence de catégories.
\end{lem}

En effet, le foncteur en question est pleinement fidèle compte tenu de \eqref{notconv22f}. 
Montrons qu'il est essentiellement surjectif. Soit $\psi$ un préfaisceau d'ensembles (resp. de groupes) $G$-équivariant sur $\Pi(\oX)^\circ$.
Pour tout point géométrique $\tx$ de $X$, il existe un ensemble (resp. groupe) $\varphi(\tx)$ unique à isomorphisme unique près, 
et pour tout point géométrique $\ox$ de $\oX$ un isomorphisme canonique 
\begin{equation}
\varphi(\varpi(\ox))\stackrel{\sim}{\rightarrow} \psi(\ox),
\end{equation}
compatible avec les isomorphismes $\tau_f^\psi \colon \psi(\ox)\rightarrow \psi(f(\ox))$ pour tout $f\in G$. 
Il résulte de \eqref{notconv22f} que la correspondence $\tx\mapsto \varphi(\tx)$ ainsi définie est un préfaisceau d'ensembles (resp. de groupes) 
sur $\Pi(X)^\circ$ muni d'un isomorphisme canonique $\varphi\circ \upvarpi\stackrel{\sim}{\rightarrow}\psi$ de préfaisceaux 
d'ensembles (resp. de groupes) $G$-équivariants sur $\Pi(\oX)^\circ$.

\section{\texorpdfstring{Schémas $K(\pi,1)$}{Schémas K(pi,1)}}

\begin{prop}\label{Kpun1}
Soient $X$ un schéma cohérent, n'ayant qu'un nombre fini de composantes connexes, 
$\mP$ un ensemble de nombres premiers, $\rho_X\colon X_\et\rightarrow X_\fet$ le morphisme canonique \eqref{notconv10a}.
Alors, les conditions suivantes sont équivalentes~:
\begin{itemize}
\item[{\rm (i)}] pour tout faisceau abélien de $\mP$-torsion $F$ de $X_\fet$ {\rm (\cite{sga4} IX 1.1)}, 
le morphisme d'adjonction $F\rightarrow \rR\rho_{X*}(\rho_X^*F)$ est un isomorphisme~;
\item[{\rm (ii)}] pour tout faisceau abélien de $\mP$-torsion, localement constant et constructible $F$ de $X_\et$ et tout entier $i\geq 1$, $\rR^i\rho_{X*}(F)=0$~;
\item[{\rm (iii)}] pour tout entier $n$ dont les diviseurs premiers appartiennent à $\mP$, 
tout $(\mZ/n\mZ)$-module localement constant et constructible $F$ de $X_\et$ et tout entier $i\geq 1$, $\rR^i\rho_{X*}(F)=0$~;
\item[{\rm (iv)}] pour tout faisceau abélien de $\mP$-torsion, localement constant et constructible $F$ de $X_\et$, tout entier $i\geq 1$
et tout $\xi\in \rH^i(X,F)$, il existe un revêtement étale surjectif $X'\rightarrow X$ tel que l'image canonique de $\xi$ dans $\rH^i(X',F)$ soit nulle~;
\item[{\rm (v)}] pour tout revêtement étale $Y\rightarrow X$, tout faisceau abélien de $\mP$-torsion, localement cons\-tant et constructible $F$ de $Y_\et$, 
tout entier $i\geq 1$ et tout $\xi\in \rH^i(Y,F)$, il existe un revêtement étale surjectif $Y'\rightarrow Y$ tel que l'image canonique de $\xi$ dans 
$\rH^i(Y',F)$ soit nulle~;
\item[{\rm (vi)}] pour tout entier $n$ dont les diviseurs premiers appartiennent à $\mP$, tout $(\mZ/n\mZ)$-module localement constant et constructible $F$ de $X_\et$,
tout entier $i\geq 1$ et tout $\xi\in \rH^i(X,F)$, il existe un revêtement étale surjectif $X'\rightarrow X$ tel que l'image canonique de $\xi$ dans $\rH^i(X',F)$ soit nulle. 
\end{itemize}
\end{prop}

En effet, le morphisme d'adjonction $\id\rightarrow \rho_{X*}\rho_X^*$ est un isomorphisme d'après (\cite{agt} VI.9.18). 
Les conditions (i) et (ii) sont donc équivalentes en vertu de (\cite{agt} VI.9.12, VI.9.14 et VI.9.20) et (\cite{sga4} VI 5.1).  
L'équivalence des conditions (ii) et (iii) est aussi une conséquence de (\cite{agt} VI.9.12 et VI.9.14) et (\cite{sga4} VI 5.1).

Pour tout faisceau abélien $G$ de $X_\et$ et tout entier $i\geq 0$, le faisceau $\rR^i\rho_{X*}(G)$ est le faisceau de $X_\fet$
associé au préfaisceau qui à tout $Y\in \ob(\Et_{\rf/X})$ associe le groupe $\rH^i(Y,G)$ (\cite{sga4} V 5.1). On a donc (iii)$\Rightarrow$(iv) et (v)$\Rightarrow$(iii). 

Montrons (iv)$\Rightarrow$(v). Supposons la condition (iv) satisfaite. Soient $f\colon Y\rightarrow X$ un revêtement étale, $F$ un faisceau abélien de $\mP$-torsion 
localement constant constructible de $Y_\et$, $i$ un entier $\geq 1$, $\xi\in \rH^i(Y,F)$. D'après (\cite{sga4} V 5.3, VIII 5.5 et XVI 2.2),
$f_*(F)$ est un faisceau abélien de $\mP$-torsion localement constant constructible de $X_\et$,  et le morphisme canonique
\begin{equation}
\rH^i(X,f_*(F))\rightarrow \rH^i(Y,F)
\end{equation}
est un isomorphisme.  Notons $\zeta$ l'image inverse de $\xi$. D'après (iv), il existe un revêtement étale surjectif $g\colon X'\rightarrow X$ tel que l'image canonique de 
$\zeta$ dans $\rH^i(X',f_*(F))$ soit nulle. Posons $Y'=Y\times_XX'$ et notons $f'\colon Y'\rightarrow X'$ la projection canonique. 
On a un diagramme commutatif
\begin{equation}
\xymatrix{
{\rH^i(X,f_*(F))}\ar[rr]^-(0.5)u\ar[d]&&{\rH^i(Y,F)}\ar[d]\\
{\rH^i(X',f_*(F))}\ar[r]^-(0.5)v&{\rH^i(X',f'_*(F|Y'))}\ar[r]^-(0.5){u'}&{\rH^i(Y',F)}}
\end{equation}
où les flèches verticales sont les morphismes canoniques, 
les morphismes $u$ et $u'$ sont induits par la suite spectrale de Cartan-Leray et le morphisme $v$ est induit par  
le morphisme de changement de base $f_*(F)|X'\rightarrow f'_*(F|Y')$.
On en déduit que l'image canonique de $\xi$ dans $\rH^i(Y',F)$ est nulle; d'où la condition (v). 

Enfin, les conditions (iv) et (vi) sont équivalentes en vertu de (\cite{sga4} VI 5.2).

\begin{defi}\label{Kpun2}
Soient $X$ un schéma cohérent, n'ayant qu'un nombre fini de composantes connexes, $\mP$ un ensemble de nombres premiers. 
On dit que $X$ est un {\em schéma $K(\pi,1)$ pour les faisceaux abéliens de $\mP$-torsion} (\cite{sga4} IX 1.1) si les conditions équivalentes
de \ref{Kpun1} sont remplies. 
Si $\mP$ est l'ensemble des nombres premiers inversibles dans $\co_X$, on dit simplement que $X$ est un {\em schéma $K(\pi,1)$}.
\end{defi}

\begin{rema}\label{Kpun3}
Soient $X$ un schéma, $F$ un faisceau abélien de torsion, localement constant et constructible de $X_\et$. 
Par descente (\cite{sga1} VIII 2.1, \cite{ega4} 2.7.1 et 17.7.3), $F$ est représentable par un schéma en groupes étale et fini $G$ au-dessus de $X$.  
Le groupe $\rH^1(X,F)$ classifie alors les $G$-fibrés principaux homogènes au-dessus de $X$ (\cite{sga4} VII 2 a)). 
Par descente, tout $G$-fibré principal homogène au-dessus de $X$ est un revêtement étale de $X$. 
Par suite, pour tout $\xi\in \rH^1(X,F)$, il existe un revêtement étale surjectif $X'\rightarrow X$ tel que l'image canonique de $\xi$ dans $\rH^1(X',F)$ soit nulle.
Notant $\rho_X\colon X_\et\rightarrow X_\fet$ le morphisme canonique \eqref{notconv10a}, on en déduit que 
$\rR^1\rho_{X*}(F)=0$ (\cite{sga4} V 5.1). Il en résulte, par la suite spectrale de Cartan-Leray, que le morphisme canonique 
\begin{equation}
\rH^1(X_\fet,\rho_{X*}(F)) \rightarrow \rH^1(X_\et,F)
\end{equation}
est un isomorphisme. 

Supposons, de plus, que $X$ soit cohérent et qu'il n'ait qu'un nombre fini de composantes connexes. Pour tout faisceau abélien de torsion $G$ de $X_\fet$, 
le morphisme d'adjonction $G\rightarrow \rho_{X*}(\rho_X^*G)$ est un isomorphisme d'après (\cite{agt} VI.9.18). On en déduit que le morphisme canonique 
\begin{equation}
\rH^1(X_\fet,G) \rightarrow \rH^1(X_\et,\rho_X^*(G))
\end{equation}
est un isomorphisme. 
\end{rema}

\begin{lem}\label{Kpun4}
Soient $X$ un schéma cohérent et étale localement connexe {\rm (\cite{agt} VI.9.7)}, $f\colon Y\rightarrow X$ un revêtement étale, 
$\mP$ un ensemble de nombres premiers. Alors, 
\begin{itemize}
\item[{\rm (i)}] Si $X$ est un schéma $K(\pi,1)$ pour les faisceaux abéliens de $\mP$-torsion,
il en est de même de $Y$.
\item[{\rm (ii)}] Supposons $f$ surjectif. Pour que $X$ soit un schéma $K(\pi,1)$ pour les faisceaux abéliens de $\mP$-torsion,
il faut et il suffit qu'il en soit de même de $Y$.
\end{itemize}
\end{lem}

Cela résulte aussitôt de \ref{Kpun1}.

\begin{prop}\label{Kpun10} 
Soient $I$ une catégorie filtrante et essentiellement $\mU$-petite {\rm (\cite{sga4} I 8.1.8)}, 
$\varphi\colon I^\circ \rightarrow \Sch$ \eqref{notconv3} un foncteur qui transforme les objets de $I$ en des schémas cohérents n'ayant qu'un nombre
fini de composantes connexes, et les morphismes de $I$ en des morphismes affines, $X$ la limite projective de $\varphi$ {\rm (\cite{ega4} 8.2.3)},
$\mP$ un ensemble de nombres premiers. Supposons que $X$ soit cohérent et qu'il n'ait qu'un nombre fini de composantes connexes. Alors,
\begin{itemize}
\item[{\rm (i)}] Si pour tout $i\in \ob(I)$, $X_i$ est un schéma $K(\pi,1)$ pour les faisceaux abéliens de $\mP$-torsion, il en est de même de $X$. 
\item[{\rm (ii)}] Supposons que $\varphi$ transforme les morphismes de $I$ en des revêtements étales surjectifs. 
Alors, pour que $X$ soit un schéma $K(\pi,1)$ pour les faisceaux abéliens de $\mP$-torsion, il faut et il suffit qu'il en soit de même de $X_i$ pour tout $i\in \ob(I)$,
et il suffit qu'il en soit de même de $X_i$ pour un seul objet $i$ de $I$. 
\end{itemize}
\end{prop}

(i) En effet, soit $F$ un faisceau abélien de $\mP$-torsion, localement constant et constructible de $X_\et$. Par descente, $F$ est représentable par un
schéma en groupes abéliens fini et étale au-dessus de $X$. 
D'après (\cite{ega4} 8.8.2, 8.10.5 et 17.7.8), il existe $\iota\in \ob(I)$ et un faisceau abélien de $\mP$-torsion localement constant constructible $F_{\iota}$
de $X_{\iota,\et}$ tels que $F$ soit l'image inverse de $F_\iota$. Quitte à remplacer $I^\circ$ par $I^\circ_{/\iota}$, 
on peut supposer que $\iota$ est un objet initial de $I$. Pour tout $i\in \ob(I)$, on note $F_i$ l'image inverse canonique de $F_\iota$ sur $X_i$. 
En vertu de (\cite{sga4} VII 5.8), pour tout entier $q\geq 0$, on a un isomorphisme canonique 
\begin{equation}\label{Kpun10a} 
\rH^q(X,F)\stackrel{\sim}{\rightarrow}\underset{\underset{i\in I}{\longrightarrow}}{\lim}\ \rH^q(X_i,F_i).
\end{equation}
La proposition s'ensuit aussitôt. 

(ii) Supposons que $X$ soit un schéma $K(\pi,1)$ pour les faisceaux abéliens de $\mP$-torsion. Soit $\iota\in \ob(I)$, $F_\iota$ un 
faisceau abélien de $\mP$-torsion, localement constant et constructible de $X_{\iota,\et}$. 
Quitte à remplacer $I^\circ$ par $I^\circ_{/\iota}$, on peut supposer 
que $\iota$ est un objet initial de $I$. Pour tout $i\in \ob(I)$, on note $F_i$ (resp. $F$) 
l'image inverse canonique de $F_\iota$ sur $X_i$ (resp. $X$). Soient $n\geq 1$, $\xi\in \rH^n(X_\iota,F_\iota)$. Par hypothèse, 
il existe un revêtement étale surjectif $X'\rightarrow X$ tel que l'image canonique de $\xi$ dans $\rH^n(X',F)$ soit nulle. 
D'après (\cite{ega4} 8.8.2, 8.10.5, 17.7.8 et 8.3.11), il existe un $i\in \ob(I)$, un revêtement étale surjectif $X'_i\rightarrow X_i$
et un $X$-isomorphisme $X'_i\times_{X_i}X \stackrel{\sim}{\rightarrow} X'$. Par un isomorphisme analogue à \eqref{Kpun10a}, on voit qu'il existe 
un morphisme $i\rightarrow j$ de $I$ tel que l'image canonique de $\xi$ dans $\rH^n(X'_i\times_{X_i}X_j,F_j)$ soit nulle. Comme le morphisme $X_j\rightarrow X_\iota$
est un revêtement étale surjectif, on en déduit que $X_\iota$ est $K(\pi,1)$ pour les faisceaux abéliens de $\mP$-torsion. 
La proposition s'ensuit compte tenu de (i) et \ref{Kpun4}(ii).

\begin{cor}\label{Kpun11}
Soient $K$ un corps, $L$ une extension algébrique de $K$, $X$ un $K$-schéma de type fini,
$\mP$ un ensemble de nombres premiers. Alors, pour que $X$ soit $K(\pi,1)$ pour les faisceaux abéliens de $\mP$-torsion,  
il faut et il suffit qu'il en soit de même de $X\otimes_KL$. 
\end{cor}

Quitte à remplacer $L$ par la clôture séparable de $K$ dans $L$ (\cite{sga4} VIII 1.1), on peut supposer $L$ séparable sur $K$. 
La proposition résulte alors de \ref{Kpun4}(ii) et \ref{Kpun10}(ii). 

\begin{prop}\label{Kpun18}
Soient $X$ un schéma, $X^\circ$ un ouvert de $X$, $\ox$ un point géométrique de $X$, $X'$ le localisé strict de $X$ en $\ox$, 
$\mP$ un ensemble de nombres premiers.  
Pour tout $X$-schéma $Y$, posons $Y^\circ=Y\times_XX^\circ$. Supposons que 
le schéma $X'^\circ$ soit cohérent et qu'il n'ait qu'un nombre fini de composantes connexes. 
Alors, les conditions suivantes sont équivalentes~:
\begin{itemize}
\item[{\rm (i)}] Le schéma $X'^\circ$ est $K(\pi,1)$ pour les faisceaux abéliens de $\mP$-torsion. 
\item[{\rm (ii)}] Pour tout $X$-schéma étale $\ox$-pointé $Y$, tout faisceau abélien de $\mP$-torsion, localement cons\-tant et constructible $\cF$ de $Y^\circ_\et$, 
tout entier $q\geq 1$ et tout $\xi\in \rH^q(Y^\circ,\cF)$, il existe un $Y$-schéma étale $\ox$-pointé $U$ et un revêtement étale surjectif $V\rightarrow U^\circ$
tels que l'image canonique de $\xi$ dans $\rH^q(V,\cF)$ soit nulle.
\end{itemize}
\end{prop}

On peut supposer $X$ affine. On désigne par $\fV_\ox$ la catégorie des $X$-schémas affines étales et $\ox$-pointés.
C'est une catégorie filtrante essentiellement $\mU$-petite (\cite{sga4} IV 6.8.2).
D'après (\cite{ega4} 8.2.3), $X'^\circ$ est la limite projective du foncteur 
\begin{equation}\label{Kpun18a}
\varphi \colon \fV_\ox\rightarrow \Sch, \ \ \  Y\mapsto Y^\circ.
\end{equation} 
Supposons d'abord la condition (ii) satisfaite. 
Soient $F$ un faisceau abélien de $\mP$-torsion localement constant constructible de $X'^\circ_\et$,
$q$ un entier $\geq 1$, $\xi\in \rH^q(X'^\circ,F)$. 
Montrons qu'il existe un revêtement étale surjectif $W\rightarrow X'^\circ$ tel que l'image canonique de $\xi$ dans $\rH^q(W,F)$ soit nulle. 
Par descente, $F$ est représentable par un schéma en groupes abéliens fini et étale au-dessus de $X'^\circ$. 
D'après (\cite{ega4} 8.8.2, 8.10.5 et 17.7.8), quitte à remplacer $X$ par un objet de $\fV_\ox$, on peut supposer qu'il existe un 
faisceau abélien de torsion localement constant constructible $\cF$ de $X^\circ_\et$ tel que $F$ soit isomorphe à l'image inverse de $\cF$. 
En vertu de (\cite{sga4} VII 5.8), on a un isomorphisme canonique 
\begin{equation}\label{Kpun18b}
\rH^q(X'^\circ,F)\stackrel{\sim}{\rightarrow}\underset{\underset{Y\in \fV^\circ_\ox}{\longrightarrow}}{\lim}\ \rH^q(Y^\circ,\cF).
\end{equation}
Il existe donc un objet $Y$ de $\fV_\ox$ tel que $\xi$ soit l'image canonique d'une classe $\zeta\in \rH^q(Y^\circ,\cF)$. 
D'après (ii), il existe un morphisme $U\rightarrow Y$ de $\fV_\ox$ et un revêtement étale surjectif $V\rightarrow U^\circ$ tels que l'image 
canonique de $\zeta$ dans $\rH^q(V,\cF)$ soit nulle. Le revêtement étale surjectif $V\times_{U^\circ}X'^\circ\rightarrow X'^\circ$ 
répond alors à la question~; d'où la condition (i).  

Inversement, supposons la condition (i) satisfaite. Soient $\cF$ un faisceau abélien de $\mP$-torsion, localement constant et constructible de $X^\circ_\et$, 
$q$ un entier $\geq 1$, $\xi\in \rH^q(X^\circ,\cF)$. D'après (i), il existe un revêtement étale surjectif $V\rightarrow X'^\circ$ tel que l'image canonique de $\xi$
dans $\rH^q(V,\cF|X'^\circ)$ soit nulle. D'après (\cite{ega4} 8.8.2, 8.10.5 et 17.7.8), il existe un objet $Y$ de $\fV_\ox$, un revêtement étale surjectif 
$W\rightarrow Y^\circ$ et un $X'^\circ$-isomorphisme $V\stackrel{\sim}{\rightarrow}W\times_{Y^\circ}X'^\circ$. 
En vertu de (\cite{sga4} VII 5.8), on a un isomorphisme canonique 
\begin{equation}\label{Kpun18c}
\rH^q(V,\cF|X'^\circ)\stackrel{\sim}{\rightarrow}\underset{\underset{U\in (\fV^\circ_\ox)_{/Y}}{\longrightarrow}}{\lim}\ \rH^q(W\times_{Y^\circ}U^\circ,\cF).
\end{equation}
Il existe donc un morphisme $U\rightarrow Y$ de $\fV^\circ_\ox$ tel que l'image canonique de $\xi$ dans $\rH^q(W\times_{Y^\circ}U^\circ,\cF)$ soit nulle~;
d'où la condition (ii).

\begin{defi}[\cite{sga4} XI 3.1, \cite{mo} 1.3.1]\label{Kpun5}
On appelle {\em courbe élémentaire} un morphisme de schémas $f\colon X\rightarrow S$ qui s'insère dans un diagramme commutatif
\begin{equation}
\xymatrix{
X\ar[r]^j\ar[rd]_f&\oX\ar[d]_-(0.4)\of&Y\ar[l]_i\ar[ld]^g\\
&S&}
\end{equation}
satisfaisant aux conditions suivantes~:
\begin{itemize}
\item[{\rm (i)}] $\of$ est une courbe relative projective et lisse, à fibres géométriquement irréductibles~; 
\item[{\rm (ii)}] $j$ est une immersion ouverte, $i$ est une immersion fermée, et $X$ est l'ouvert complémentaire de $i(Y)$ dans $\oX$;
\item[{\rm (iii)}] $g$ est un revêtement étale surjectif. 
\end{itemize}
On dit alors aussi que $X$ est une {\em courbe élémentaire au-dessus de $S$}. 
\end{defi}

\begin{defi}[\cite{sga4} XI 3.2, \cite{mo} 1.3.2]\label{Kpun6}
On appelle {\em polycourbe élémentaire} un morphisme de schémas $f\colon X\rightarrow S$ qui admet une factorisation en courbes élémentaires.
On dit alors aussi que $X$ est une {\em polycourbe élémentaire au-dessus de $S$}. 
\end{defi}

\begin{prop}[\cite{sga4} XI 3.3]\label{Kpun7} 
Soient $k$ un corps algébriquement clos, $X$ un $k$-schéma lisse, $x\in X(k)$. Alors, il existe un ouvert de Zariski de $X$ contenant $x$ qui est une 
polycourbe élémentaire au-dessus de $\Spec(k)$. 
\end{prop}

\begin{prop}\label{Kpun9} 
Soient $S$ un $\mQ$-schéma noethérien, régulier et $K(\pi,1)$, $f\colon X\rightarrow S$ une courbe élémentaire. 
Alors, $X$ est un schéma $K(\pi,1)$.
\end{prop}

On notera d'abord que $X$ est cohérent et n'a qu'un nombre fini de composantes connexes. Considérons un diagramme commutatif 
\begin{equation}\label{Kpun9a} 
\xymatrix{
X\ar[r]^j\ar[rd]_f&\oX\ar[d]_-(0.4)\of&Y\ar[l]_i\ar[ld]^g\\
&S&}
\end{equation}
vérifiant les conditions de \ref{Kpun5}. Soient $F$ un faisceau abélien de torsion, localement constant et constructible sur $X$, 
$n$ un entier $\geq 1$, $\xi\in \rH^n(X,F)$. Montrons qu'il existe un revêtement étale surjectif $X'\rightarrow X$ tel que l'image canonique
de $\xi$ dans $ \rH^n(X',F)$ soit nulle. On peut se borner au cas où $n\geq 2$ \eqref{Kpun3}. 
Il existe un revêtement étale surjectif $Z\rightarrow X$ tel que $F|Z$ soit constant. 
D'après \ref{Kpun7}, le morphisme $Z\rightarrow S$ induit par $f$ est une courbe élémentaire. On peut donc se borner au cas où $F$ est constant,
de valeur un groupe cyclique fini $\Lambda$. 

En vertu du théorème de pureté (\cite{travaux-gabber} XVI 3.1.1), 
on a un isomorphisme canonique $ i^!\Lambda_\oX\stackrel{\sim}{\rightarrow}\Lambda_Y(-1)[-2]$. Le triangle distingué de localisation induit alors 
un triangle de $\bD^+(\oX_\et,\Lambda)$
\begin{equation}\label{Kpun9b} 
\Lambda_\oX\longrightarrow \rR j_*(\Lambda_X)\longrightarrow i_*(\Lambda_Y)(-1)[-1]\stackrel{+1}{\longrightarrow}
\end{equation}
et par suite, $g$ étant fini, un triangle de $\bD^+(S_\et,\Lambda)$
\begin{equation}\label{Kpun9c} 
\rR\of_*(\Lambda_\oX)\longrightarrow \rR f_*(\Lambda_X)\longrightarrow g_*(\Lambda_Y)(-1)[-1]\stackrel{+1}{\longrightarrow}
\end{equation}
Montrons que le morphisme induit 
\begin{equation}\label{Kpun9d} 
g_*(\Lambda_Y)(-1)\rightarrow \rR^2 \of_*(\Lambda_\oX)
\end{equation}
est surjectif. En effet, la formation de ce morphisme commute à tout changement de base $S'\rightarrow S$, d'après (\cite{travaux-gabber} XVI 2.3.2). 
On peut donc se réduire au cas où $S$ est le spectre d'un corps algébriquement clos et où $Y$ est un point fermé de $\oX$. La classe de cycle définie par $Y$
engendre alors $\rH^2(\oX,\Lambda(1))$ (\cite{sga45} Cycle 2.1.5); d'où l'assertion. 

Il résulte de ce qui précède qu'on a un isomorphisme canonique $\of_*(\Lambda_\oX) \stackrel{\sim}{\rightarrow} f_*(\Lambda_X)$ et une suite exacte canonique
\begin{equation}\label{Kpun9e} 
0\rightarrow \rR^1 \of_*(\Lambda_\oX)\rightarrow \rR^1 f_*(\Lambda_X)\rightarrow g_*(\Lambda_Y)(-1) \rightarrow \rR^2\of_*(\Lambda_\oX)\rightarrow 0, 
\end{equation}
et que $\rR^q f_*(\Lambda_X)=0$ pour tout $q\geq 2$. Comme $\of$ et $g$ sont propres et lisses, pour tout $q\geq 0$, 
$\rR^q f_*(\Lambda_X)$ est localement constant constructible (\cite{sga4} XVI 2.2 et IX 2.1).

Considérons la suite spectrale de Cartan-Leray
\begin{equation}
E_2^{a,b}=\rH^a(S,\rR^bf_*(\Lambda_X))\Rightarrow \rH^{a+b}(X,\Lambda_X),
\end{equation}
et notons $(\rE^n_q)_{0\leq q\leq n}$ la filtration aboutissement sur $\rH^n(X,\Lambda_X)$ (on rappelle que $n\geq 2$).
Comme on a $E_\infty^{q,n-q}=E^n_q/E^n_{q+1}$ pour tout $0\leq q\leq n$, on en déduit une suite exacte canonique
\begin{equation}
\rH^n(S,f_*(\Lambda_X))\stackrel{u}{\rightarrow} \rH^n(X,\Lambda_X)\stackrel{v}{\rightarrow} \rH^{n-1}(S,\rR^1 f_*(\Lambda_X))\rightarrow \rH^{n+1}(S,f_*(\Lambda_X)).
\end{equation}
Celle-ci est compatible à tout changement de base $S'\rightarrow S$ dans un sens évident que nous n'explicitons pas. 
Comme le schéma $S$ est $K(\pi,1)$, il existe un revêtement étale surjectif $S'\rightarrow S$ tel que l'image canonique de $v(\xi)$ dans 
$\rH^{n-1}(S',\rR^1 f_*(\Lambda_X))$ soit nulle. On peut donc supposer qu'il existe $\zeta\in \rH^n(S,f_*(\Lambda_X))$ tel que $\xi=v(\zeta)$. 
De même, il existe un revêtement étale surjectif $S''\rightarrow S$ tel que l'image canonique de $\zeta$ dans $\rH^n(S'',f_*(\Lambda_X))$ soit nulle. 
L'image canonique de $\xi$ dans $\rH^n(X\times_SS'',\Lambda_X)$ est donc nulle~; d'où l'assertion recherchée.

\begin{cor}\label{Kpun12}
Soit $k$ un corps algébriquement clos de caractéristique $0$, $X$ un schéma lisse sur $k$, $x\in X$. 
Alors, il existe un ouvert $U$ de $X$ contenant $x$, qui est un schéma $K(\pi,1)$.
\end{cor}

En effet, quitte à remplacer $x$ par une spécialisation, on peut se réduire au cas où $x\in X(k)$. 
La proposition résulte alors de \ref{Kpun7} et \ref{Kpun9}.

\begin{lem}\label{Kpun14}
Soient $S$ un schéma noethérien régulier, $d$ un entier $\geq1$, 
$f\colon X\rightarrow \mA_S^d$ un morphisme lisse d'un schéma $X$ dans l'espace affine de dimension $d$ au-dessus de $S$. 
Notons $\mG^d_{m,S}$ l'ouvert de $\mA_S^d$ où les coordonnées ne s'annulent pas, $U=f^{-1}(\mG^d_{m,S})$ 
et $j\colon U\rightarrow X$ l'injection canonique. Soient $n$ un entier $\geq 1$ inversible dans $\co_S$,  $\Lambda=\mZ/n\mZ$,
$F$ un faisceau de $\Lambda$-modules localement constant et constructible de $X_\et$. 
Notons $\pi\colon \mA_S^d\rightarrow \mA_S^d$ le morphisme défini par l'élévation à la puissance $n$-ième des coordonnées de $\mA_S^d$,
$X'$ le changement de base de $X$ par $\pi$, $\varpi\colon X'\rightarrow X$ la projection canonique, $U'=\varpi^{-1}(U)$,
$j'\colon U'\rightarrow X'$ l'injection canonique et $F'=\varpi^*(F)$. Alors, pour tout entier $q\geq 1$, le morphisme de changement de base
\begin{equation}\label{Kpun14a}
\varpi^*(\rR^qj_*(j^*F))\rightarrow \rR^qj'_*(j'^*F')
\end{equation}
est nul. 
\end{lem}

En effet, la question étant locale pour la topologie étale de $X'$ et donc aussi pour celle de $X$, on peut se borner au cas où $F$ est constant de valeur $\Lambda$. 
Notons $(A_\alpha)_{1\leq \alpha\leq d}$ les axes de coordonnées de $\mA^d_S$ et pour tout $1\leq \alpha\leq d$, posons 
$D_\alpha=X\times_{\mA_S^d} A_\alpha$ et notons $i_\alpha\colon D_\alpha\rightarrow X$ l'injection canonique. 
En vertu de (\cite{travaux-gabber} XVI 3.1.4), on a un isomorphisme canonique
\begin{equation}
\rR^1 j_*(\Lambda_U)\stackrel{\sim}{\rightarrow}\oplus_{1\leq \alpha\leq d}\ i_{\alpha*}(\Lambda_{D_\alpha}(-1)).
\end{equation}
De plus, pour tout entier $q\geq 1$, le morphisme 
\begin{equation}
\wedge^q (\rR^1j_*(\Lambda_U))\rightarrow \rR^qj_*(\Lambda_U)
\end{equation}
défini par cup-produit est un isomorphisme. 
De même, considérons $X'$ comme un $\mA^d_S$-schéma lisse via la projection canonique 
$f'\colon X'\rightarrow \mA_S^d$ et pour tout $1\leq \alpha\leq d$, posons $D'_\alpha=X'\times_{\mA_S^d} A_\alpha$ et notons 
$i'_\alpha\colon D'_\alpha\rightarrow X'$ l'injection canonique.
On a alors des isomorphismes canoniques
\begin{eqnarray}
\rR^1 j'_*(\Lambda_{U'})&\stackrel{\sim}{\rightarrow}&\oplus_{1\leq \alpha\leq d}\ i'_{\alpha*}(\Lambda_{D'_\alpha}(-1)),\\
\wedge^q (\rR^1j'_*(\Lambda_{U'}))&\stackrel{\sim}{\rightarrow}& \rR^qj'_*(\Lambda_{U'}).
\end{eqnarray}
Pour tout $1\leq \alpha\leq d$, $D'_\alpha$ s'identifie canoniquement au sous-schéma réduit de $X'$ sous-jacent à $D_\alpha\times_XX'$. De plus, on a 
$\varpi^*(D_\alpha)=nD'_\alpha$ en tant que diviseurs de Cartier. Il résulte de (\cite{travaux-gabber} XVI 3.4.8) que le diagramme 
\begin{equation}
\xymatrix{
{\varpi^*(\rR^1j_*(\Lambda_U))}\ar[r]^-(0.5)\sim\ar[d]_\lambda&{\oplus_{1\leq \alpha\leq d}\ \varpi^*(i_{\alpha*}(\Lambda_{D_\alpha}(-1)))}\ar[d]^\gamma\\
{\rR^1j'_*(\Lambda_{U'})}\ar[r]^-(0.5)\sim&{\oplus_{1\leq \alpha\leq d}\ i'_{\alpha*}(\Lambda_{D'_\alpha}(-1))}}
\end{equation}
où $\lambda$ est le morphisme de changement de base et $\gamma$ est induit pour chaque $1\leq \alpha \leq d$ par $n$ fois le morphisme de changement de base,
est commutatif. La proposition s'ensuit puisque $\gamma$ est nul.

\begin{lem}\label{Kpun15}
Soient $f\colon Y\rightarrow X$, $g\colon Y'\rightarrow Y$ deux morphismes de schémas
tels que $f$ soit fini et que $g$ soit étale, $\ox$ un point géométrique de $X$ tels que le morphisme
\begin{equation}\label{Kpun15a}
Y'\otimes_X\kappa(\ox)\rightarrow Y\otimes_X\kappa(\ox)
\end{equation}
induit par $g$ soit surjectif. Alors, il existe un $X$-schéma étale $\ox$-pointé $X'$ et un $Y$-morphisme $Y\times_XX'\rightarrow Y'$. 
\end{lem}

Il suffit de montrer que si $X$ est le spectre d'un anneau local strictement hensélien de point fermé $\ox$, alors $g$ admet une
section (\cite{ega4} 8.8.2(i)). On a un $X$-isomorphisme $Y\simeq\amalg_{\oy\in Y\otimes_X\kappa(\ox)}Y_{(\oy)}$,
où $Y_{(\oy)}$ est le localisé strict de $Y$ en le point géométrique $\oy$. 
On peut évidemment se réduire au cas où $Y$ est connexe et non vide. Donc $Y$ est isomorphe au spectre
d'un anneau local strictement hensélien. 
D'après \eqref{Kpun15a} et (\cite{ega4} 18.5.11), $Y'$ est la somme de deux schémas, dont l'un est isomorphe à $Y$; d'où la proposition.

\begin{prop}[\cite{achinger} 8.1]\label{Kpun16}
Soient $S$ un $\mQ$-schéma noethérien régulier, $X$ un $S$-schéma lisse, $D$ un diviseur à croisements normaux sur $X$ 
relativement à $S$ {\rm (\cite{sga1} XIII 2.1)}, 
$U$ l'ouvert complémentaire de $D$ dans $X$, $\ox$ un point géométrique de $X$, 
$F$ un faisceau abélien de torsion, localement constant et constructible sur $U$,
$q$ un entier $\geq 1$, $\xi\in \rH^q(U,F)$. Alors, il existe un $X$-schéma étale $\ox$-pointé $Y$, et posant $V=U\times_XY$, un revêtement étale 
surjectif $V'\rightarrow V$ tels que l'image canonique de $\xi$ dans $\rH^q(V',F)$ soit nulle.
\end{prop}

\begin{cor}\label{Kpun19}
Soient $k$ un corps algébriquement clos de caractéristique $0$, $X$ un schéma lisse sur $k$, 
$D$ un diviseur à croisements normaux sur $X$, $X^\circ$ l'ouvert complémentaire de $D$ dans $X$, 
$\ox$ un point géométrique de $X$, $X'$ le localisé strict de $X$ en $\ox$. 
Alors, $X'\times_XX^\circ$ est un schéma $K(\pi,1)$. 
\end{cor}

Cela résulte de \ref{Kpun18} et  \ref{Kpun16} puisque $X'\times_XX^\circ$ est noethérien et intègre. 

\subsection{}\label{Kpun20}
Soient $V$ un anneau de valuation discrète, $S=\Spec(V)$,  $s$ (resp.  $\eta$, resp. $\oeta$) le point fermé 
(resp.  le point générique, resp. un point géométrique générique) de $S$. 
On suppose que le corps des fractions de $V$ est de caractéristique $0$, 
et que son corps résiduel est parfait de caractéristique $p>0$.
On munit $S$ de la structure logarithmique $\cM_S$ définie par son point fermé, 
autrement dit, $\cM_S=u_*(\co_\eta^\times)\cap \co_S$, où $u\colon \eta\rightarrow S$ est l'injection canonique.

\begin{prop}[\cite{achinger} 6.1]\label{Kpun21}
Conservons les hypothèses de \ref{Kpun20}, soient de plus $(X,\cM_X)$ un schéma logarithmique fin,
$f\colon (X,\cM_X)\rightarrow (S,\cM_S)$ un morphisme lisse tel que le schéma usuel $X_\eta$ soit lisse sur $\eta$, 
$\ox$ un point géométrique de $X$. 
Alors, il existe un voisinage étale $U$ de $\ox$ dans $X$ tel que $U_\eta$ soit un schéma $K(\pi,1)$. 
\end{prop}

\begin{cor}[\cite{achinger} 9.5]\label{Kpun22}
Conservons les hypothèses de \ref{Kpun20}, supposons de plus que $S$ soit strictement local. Soient $(X,\cM_X)$ un schéma logarithmique fin,
$f\colon (X,\cM_X)\rightarrow (S,\cM_S)$ un morphisme lisse tel que $X_\eta$ soit lisse sur $\eta$, 
$\ox$ un point géométrique de $X$ au-dessus de $s$, 
$X'$ le localisé strict de $X$ en $\ox$. Alors, $X'_\oeta$ est un schéma $K(\pi,1)$.
\end{cor}

Montrons d'abord que le schéma $X'_\oeta$ n'a qu'un nombre fini de composantes connexes. Soit $\ell$ un nombre premier différent de $p$. 
D'après (\cite{sga45} Th. finitude 3.2), le $\mF_\ell$-espace vectoriel $\rH^0(X'_\oeta,\mF_\ell)$ est de dimension finie. 
Il s'ensuit que l'ensemble $\fC$ des sous-schémas ouverts et fermés de $X'_\oeta$ est fini (\cite{sga4} VIII 6.1). 
Pour tout point $z$ de $X'_\oeta$, notons $U_z$ l'intersection de tous les sous-schémas ouverts et fermés de $X'_\oeta$ contenant $z$. 
Comme $\fC$ est fini, $U_z$ est ouvert et fermé dans $X'_\oeta$, autrement dit, c'est un objet de $\fC$. Par suite, $U_z$ est connexe,  
et est donc égal à la composante connexe de  $X'_\oeta$ contenant $z$. On en déduit que $X'_\oeta$ n'a qu'un nombre fini de composantes connexes.
La proposition résulte alors de \ref{Kpun18} et \ref{Kpun21}.

\begin{rema}\label{Kpun23}
On établira dans \ref{Kp1s6} l'énoncé le plus général de (\cite{achinger} 9.5).
\end{rema}

\subsection{}\label{Kpun30}
Soient $X,Y$ deux schémas cohérents n'ayant qu'un nombre fini de composantes connexes,
$f\colon Y\rightarrow X$ un morphisme. Le diagramme 
\begin{equation}\label{Kpun30a}
\xymatrix{
{Y_\et}\ar[d]_{f_\et}\ar[r]^{\rho_Y}&{Y_\fet}\ar[d]^{f_\fet}\\
{X_\et}\ar[r]^{\rho_X}&{X_\fet}}
\end{equation}
où $\rho_X$ et $\rho_Y$ sont les morphismes canoniques \eqref{notconv10a}
est commutatif à isomorphisme canonique près (\cite{agt} (VI.9.3.4)). 
Il induit donc un morphisme de changement de base (\cite{egr1} 1.2.2)
\begin{equation}\label{Kpun30b}
\rho_X^*f_{\fet*}\rightarrow f_{\et *}\rho_Y^*.
\end{equation}

\begin{prop}\label{Kpun31}
Conservons les hypothèses et notations de \ref{Kpun30}; supposons, de plus, que $f$ soit propre et lisse. Alors, le morphisme de changement de base 
\begin{equation}\label{Kpun31a}
\rho_X^*f_{\fet*}\rightarrow f_{\et *}\rho_Y^*
\end{equation} 
est un isomorphisme.
\end{prop}

Soit $F$ un objet de $Y_\fet$. Le morphisme composé 
\begin{equation}
f_{\fet*}(F)\rightarrow f_{\fet*}(\rho_{Y*}(\rho_Y^*(F)))\stackrel{\sim}{\rightarrow}\rho_{X*}(f_{\et*}(\rho_Y^*(F)))
\end{equation}
où la première flèche est induite par le morphisme d'adjonction $\id\rightarrow \rho_{Y*}\rho_Y^*$ et la seconde flèche est 
l'isomorphisme sous-jacent à \eqref{Kpun30a}, est l'adjoint du morphisme \eqref{Kpun31a}. C'est un isomorphisme en vertu de (\cite{agt} VI.9.18). 
Par ailleurs, le morphisme $f_{\et}\colon Y_\et\rightarrow X_\et$ est cohérent (\cite{sga4} VI 3.10).  
Il résulte alors de (\cite{agt} VI.9.20) et (\cite{sga4} VI 5.1 et XVI 2.2) qu'il existe un objet 
$G$ de $X_\fet$ et un isomorphisme 
\begin{equation}
f_{\et*}(\rho_Y^*(F))\stackrel{\sim}{\rightarrow}\rho_X^*(G).
\end{equation}
Par ailleurs, le diagramme 
\begin{equation}
\xymatrix{
{\rho_X^*(\rho_{X*}(\rho_X^*(G)))}\ar[r]^-(0.4)b&{\rho_X^*(G)}\\
{\rho_X^*(G)}\ar[u]^-(0.5)a\ar[ru]_\id&}
\end{equation}
où $a$ et $b$ sont induits par les morphismes d'adjonction est commutatif. Comme $a$ est un isomorphisme d'après (\cite{agt} VI.9.18), 
$b$ est un isomorphisme. Par suite, le morphisme d'adjonction 
\begin{equation}
\rho_X^*(\rho_{X*}(f_{\et*}(\rho_Y^*(F))))\rightarrow f_{\et*}(\rho_Y^*(F))
\end{equation}
est un isomorphisme; d'où la proposition.  

\begin{cor}\label{Kpun32}
Conservons les hypothèses et notations de \ref{Kpun30}; supposons, de plus, que $f$ soit étale et fini. Soient $F$ un objet de $Y_\fet$, 
$\ox$ un point géométrique de $X$, $\oy_1,\dots,\oy_n$ les points de $Y_\ox$, que l'on identifie à des points géométriques de $Y$. 
Alors, on a un isomorphisme canonique fonctoriel 
\begin{equation}
(f_{\fet*}(F))_{\rho_X(\ox)}\stackrel{\sim}{\rightarrow}\prod_{i=1}^nF_{\rho_Y(\oy_i)}.
\end{equation}
\end{cor}

Cela résulte de \ref{Kpun31} et (\cite{sga4} VIII 5.5).

\begin{cor}\label{Kpun33}
Conservons les hypothèses et notations de \ref{Kpun30}; supposons, de plus, que $f$ soit étale et fini.
Alors, le foncteur $f_{\fet*}$ est exact sur la catégorie des faisceaux abéliens de $Y_\fet$; 
en particulier, pour tout faisceau abélien $F$ de $Y_\fet$ et tout entier $q\geq 1$, on a $\rR^qf_{\fet*}(F)=0$.
\end{cor}

\begin{prop}[\cite{sga4} XVI 2.2 \& \cite{sga45} Th. finitude, Appendice 1.3.3]\label{Kpun34}
Soient $f\colon X\rightarrow Y$ un morphisme propre et lisse de schémas, $D$ un diviseur à croisements normaux sur $X$ relativement à $Y$ {\rm (\cite{sga1} XIII 2.1)},
$U$ l'ouvert complémentaire de $D$ dans $X$, $j\colon U\rightarrow X$ l'injection canonique, $g=f\circ j$, $\Lambda$ un anneau noethérien annulé par un entier $n$ 
inversible dans $Y$, $F$ un faisceau de $\Lambda$-modules localement constant constructible sur $U_\et$, modérément ramifié le long de $D$.  Alors, pour tout entier
$q\geq 0$, le faisceau de $\Lambda$-modules $\rR^qg_*(F)$ est localement constant constructible et sa formation commute à tout changement de base 
$Y'\rightarrow Y$ {\rm (\cite{egr1} 1.2.2)}. 
En particulier, pour tout point géométrique $\oy$ de $Y$, le morphisme de changement de base 
\begin{equation}\label{Kpun34a}
\rR^qg_*(F)_\oy\rightarrow \rH^q(U_\oy,F)
\end{equation}
est bijectif. 
\end{prop}

En effet, en vertu de (\cite{sga45} Th. finitude, Appendice 1.3.3), le complexe de $\Lambda$-modules $\rR j_*(F)$ est borné de cohomologie constructible 
et le couple $(j,F)$ est cohomologiquement propre relativement à $Y$, c'est-à-dire que la formation du complexe $\rR j_*(F)$ commute à tout changement de base $Y'\rightarrow Y$. 
Par suite, le faisceau de $\Lambda$-modules $\rR^qf_*(\rR j_*(F))$ est constructible et sa formation commute à tout changement de base $Y'\rightarrow Y$
d'après (\cite{sga4} XII 5.1 et XIV 1.1). Comme on a $\rR^qg_*(F)\stackrel{\sim}{\rightarrow} \rR^qf_*(\rR j_*(F))$, le faisceau de $\Lambda$-modules $\rR^qg_*(F)$
est constructible et sa formation commute à tout changement de base $Y'\rightarrow Y$ (\cite{egr1} 1.2.4(v)). 
Pour montrer qu'il est localement constant, il suffit grâce à (\cite{sga4} IX 2.11) de montrer que pour tout morphisme de spécialisation $\oy_1\rightarrow \oy_0$ de
points géométriques de $Y$, le morphisme de spécialisation correspondant 
\begin{equation}\label{Kpun34b}
\rR^qg_*(F)_{\oy_0}\rightarrow \rR^qg_*(F)_{\oy_1}
\end{equation}
est bijectif. Il existe un schéma strictement local intègre $Y'$ et un morphisme $Y'\rightarrow Y$ appliquant le point fermé $y'_0$ de $Y'$ en $y_0$ et 
son point générique $y'_1$ en $y_1$: il suffit de prendre d'abord le localisé strict de $Y$ en $\oy_0$, 
puis son sous-schéma fermé intègre défini par un point au-dessus de $y_1$. 
De plus, quitte à remplacer $Y'$ par sa clôture intégrale dans une clôture algébrique de son corps des fonctions $\kappa(y'_1)$, on peut supposer 
$Y'$ normal et $\kappa(y'_1)$ algébriquement clos. Comme la formation de $\rR^qg_*(F)$
commute au changement de base $Y'\rightarrow Y$, on peut donc supposer $Y$ normal et strictement local de point fermé $y_0$ et de point générique $y_1$ 
tel que le corps résiduel $\kappa(y_1)$ soit algébriquement clos. D'après (\cite{sga4} XII 5.5) appliqué à $f$ et au complexe $\rR j_*(F)$, 
le morphisme de spécialisation \eqref{Kpun34b} s'identifie au morphisme canonique 
\begin{equation}
\rH^q(U,F)\rightarrow \rH^q(U_1,i_1^*(F)),
\end{equation}
où $U_1=U\times_Yy_1$ et $i_1\colon U_1\rightarrow U$ est la projection canonique. Celui-ci est un isomorphisme d'après (\cite{sga4} XVI 2.4), d'où la proposition.

\section{Compléments sur la connexité}\label{csc}

\begin{prop}\label{csc1}
Soient $X$, $X'$ deux schémas, $f\colon X'\rightarrow X$ un morphisme plat et localement de présentation finie à fibres géométriques réduites,
$\ox'$ un point géométrique de $X'$. Alors, il existe un voisinage étale $\iota'\colon U'\rightarrow X'$ de $\ox'$,
un voisinage étale $\iota\colon U\rightarrow X$ de $f(\ox')$ et un morphisme $f'\colon U'\rightarrow U$ s'insèrant dans un diagramme commutatif 
\begin{equation}\label{csc1a}
\xymatrix{
\ox'\ar[r]\ar[rd]&U'\ar[r]^{\iota'}\ar[d]^{f'}&X'\ar[d]^f\\
&U\ar[r]^\iota&X}
\end{equation}
tel que les fibres géométriques de $f'$ soient intègres. 
\end{prop}

La question étant locale, on peut supposer $f$ quasi-compact. 
Notons $x'$ l'image de $\ox'$ dans $X'$ et posons $x=f(x')$. 
Soit $C$ une composante irréductible de la fibre $X'\otimes_X\kappa(x)$ de $X'$ au-dessus de $x$, contenant $x'$. 
Le morphisme $f$ étant lisse en le point générique de $C$ (\cite{ega4} 17.5.1 et 6.7.6), 
il est lisse sur un ouvert non vide de $X'$ qui rencontre $C$. 
Il existe donc un point fermé $z$ de $C$ tel que $f$ soit lisse en $z$ et que $\kappa(z)$ soit séparable sur $\kappa(x)$ (\cite{ega4} 17.15.10(iii)).  
En vertu de (\cite{ega4} 17.16.3(i)), il existe un morphisme étale $U\rightarrow X$
et un $X$-morphisme $\iota\colon U\rightarrow X'$ qui soit une immersion d'image contenant $z$. On peut donc considérer $z$ comme un point de $U$. 
La projection canonique $X'\times_XU\rightarrow U$ est plate de présentation finie à fibres géométriques réduites.
Elle admet une section $U\rightarrow X'\times_XU$ induite par $\iota$.
En vertu de (\cite{ega4} 15.6.5), il existe un ouvert $U'$ de $X'\times_XU$ tel que pour tout $s\in U$, 
la fibre $U'\otimes_U\kappa(s)$ de $U'$ au-dessus de $s$ soit la composante connexe de
$X'\otimes_{X}\kappa(s)$ contenant $\iota(s)$. 
En particulier, $U'\otimes_U\kappa(z)$ est la composante connexe de $X'\otimes_{X}\kappa(z)$ contenant $\iota(z)$.
Celle-ci s'envoie surjectivement sur la composante connexe de $X'\otimes_X\kappa(x)$ contenant $z$ et donc aussi $x'$. 
Par suite, l'image du morphisme canonique $U'\rightarrow X'$ contient $x'$. 
Par ailleurs, les fibres du morphisme $U'\rightarrow U$ sont géométriquement connexes d'après (\cite{ega4} 4.5.14), d'où la proposition.  

\begin{prop}\label{csc2}
Soient $X$, $X'$ deux schémas normaux et localement irréductibles {\rm (\cite{agt} III.3.1)}, 
$f\colon X'\rightarrow X$ un morphisme plat et localement de présentation finie à fibres géométriques réduites,
$\ox'$ un point géométrique de $X'$, $\uX'$ le localisé strict de $X'$ en $\ox'$, $\uX$ le localisé strict de $X$ en $f(\ox')$, 
$V\rightarrow \uX$ un morphisme plat de présentation finie tel que $V$ soit connexe. Alors, le schéma $V\times_{\uX}\uX'$ est connexe. 
\end{prop}

On désigne par $\fV$ la catégorie des $X$-schémas étales $\ox'$-pointés, 
et par $\fW$ la sous-catégorie pleine de $\fV$ formée des objets $(U,u\colon \ox'\rightarrow U)$ tels que le schéma $U$ soit affine et connexe. 
Ce sont des catégories cofiltrantes, et le foncteur d'injection canonique $\fW^\circ\rightarrow \fV^\circ$ est cofinal 
d'après (\cite{sga4} I 8.1.3(c) et \cite{agt} III.3.3). Soit $(U_0,u_0)$ un objet de $\fW$. 

Comme $X$ est normal, pour tout objet $(U,u)$ de $\fW$, 
le schéma $U$ est normal et étant connexe, il est donc intègre. Pour tout morphisme
$(U',u')\rightarrow (U,u)$ de $\fW$, comme le morphisme $U'\rightarrow U$ est étale, il est schématiquement dominant (\cite{ega4} 11.10.2).
Le schéma $\uX$ est la limite projective des schémas $U$ pour $(U,u)\in \ob(\fW_{/(U_0,u_0)})$.
Pour tout objet $(U,u)$ de $\fW_{/(U_0,u_0)}$, le morphisme $\uX\rightarrow U$ est dominant d'après (\cite{ega4} 8.3.8(i)) 
et donc schématiquement dominant (\cite{sga4} 11.10.4). 
En vertu de (\cite{ega4} 8.8.2 et 11.2.6), quitte à remplacer $(U_0,u_0)$ par un objet de $\fW_{/(U_0,u_0)}$, 
il existe un morphisme plat de présentation finie $V_0\rightarrow U_0$ et un $\uX$-isomorphisme $\tau\colon V\stackrel{\sim}{\rightarrow}V_0\times_{U_0}\uX$. 
Pour tout objet $(U,u)$ de $\fW_{/(U_0,u_0)}$, le morphisme $V\rightarrow V_0\times_{U_0}U$ induit par $\tau$
est schématiquement dominant (\cite{ega4} 11.10.5). On en déduit que $V_0\times_{U_0}U$ est connexe. 

On désigne par $\fV'$ la catégorie des $X'$-schémas étales $\ox'$-pointés, 
et par $\fW'$ la sous-catégorie pleine de $\fV'$ formée des objets $(U',u'\colon \ox'\rightarrow U')$ tels que le schéma $U'$ soit affine. 
Soit $(U',u')$ un objet de $\fV'$. En vertu de \ref{csc1}, appliqué à la projection canonique $U'\times_{X}U_0\rightarrow U_0$, et (\cite{agt} III.3.3), 
il existe un objet $(U,u)$ de $\fW_{/(U_0,u_0)}$, un objet $(U'',u'')$ de $\fV'_{/(U',u')}$
et un morphisme $f''\colon U''\rightarrow U$ qui s'insère dans un diagramme commutatif 
\begin{equation}\label{csc2a}
\xymatrix{
\ox'\ar[r]^{u}\ar[dr]_{u''}&U''\ar[r]\ar[d]^{f''}&X'\ar[d]^f\\
&U\ar[r]&X}
\end{equation}
tels que les fibres géométriques de $f''$ soient connexes. Le morphisme $f''$ étant plat de présentation finie, 
quitte à remplacer $U$ par un ouvert affine contenant $u''(\ox')$, on peut supposer $f''$ surjectif. 
Par suite, les schémas $U''$ et $V_0\times_{U_0}U''$ sont connexes d'après (\cite{ega4} 4.5.7). 

Soit $(U'_0,u'_0)$ un objet de $\fW'$ au-dessus de $(U_0,u_0)$. 
Notons $\fC$ la sous-catégorie pleine de la catégorie $\fV'_{/(U'_0,u'_0)}$ 
formée des objets $(U',u')$ tels que les schémas $U'$ et $V_0\times_{U_0}U'$ soient connexes. 
Le foncteur d'injection canonique $\fC^\circ\rightarrow \fV'^\circ_{/(U'_0,u'_0)}$ est cofinal d'après ce qui précède et (\cite{sga4} I 8.1.3(c)).
Le schéma $\uX'$ est la limite projective des schémas $U'$ pour $(U',u')\in \ob(\fW'_{/(U'_0,u'_0)})$.
Il est donc la limite projective des schémas $U'$ pour $(U',u')\in \ob(\fV'_{/(U'_0,u'_0)})$, 
et est aussi la limite projective des schémas $U'$ pour $(U',u')\in \ob(\fC)$.
Comme $X'$ est normal, pour tout objet $(U',u')$ de $\fC$, 
le schéma $U'$ est normal et étant connexe, il est donc intègre. Pour tout morphisme
$(U'',u'')\rightarrow (U',u')$ de $\fC$, comme le morphisme de schémas $U''\rightarrow U'$ est étale, il est schématiquement dominant.
Il en est donc de même du morphisme $V_0\times_{U_0}U''\rightarrow V_0\times_{U_0}U'$.
On en déduit alors que le schéma  $V_0\times_{U_0}\uX'$ est connexe d'après (\cite{ega4} 8.4.4). 
La proposition s'ensuit compte tenu de l'isomorphisme $V\times_{\uX}\uX'\stackrel{\sim}{\rightarrow}V_0\times_{U_0}\uX'$ induit par $\tau$. 

\section{Compléments de géométrie logarithmique}\label{cgl}

\begin{lem}\label{cgl1}
Soient $f\colon (X,\cM_X)\rightarrow (Y,\cM_Y)$ un morphisme lisse et saturé de schémas logarithmiques fins et saturés {\rm (\cite{agt} II.5.18)}, 
$\ox$ un point géométrique de $X$, $\oy=f(\ox)$, $\uX$ le localisé strict de $X$ en $\ox$, $\uY$ le localisé strict de $Y$ en $\oy$.
On note $X^\circ$ (resp. $Y^\rhd$)  l'ouvert maximal de $X$ (resp. $Y$) où la structure logarithmique $\cM_X$ (resp. $\cM_Y$) 
est triviale {\rm (\cite{ogus} III 1.2.8)} et on pose $\uX^\circ=\uX\times_XX^\circ$ et $\uY^\rhd=\uY\times_YY^\rhd$.  
Alors, le morphisme $\uX^\circ\rightarrow \uY^\rhd$ induit par $f$ est fidèlement plat. 
\end{lem}

En effet, $f$ est plat en vertu de (\cite{kato1} 4.5). 
Il suffit donc de montrer que le morphisme $\uX^\circ\rightarrow \uY^\rhd$ est surjectif ou, ce qui revient au même, 
que pour tout point géométrique $\oy'$ de $Y^\rhd$
se spécialisant en $\oy$, il existe un point géométrique $\ox'$ de $X^\circ$ au-dessus de $\oy'$, se spécialisant en $\ox$. 
Notons $(Y',\cM_{Y'})$ le schéma $\oy'$ muni de la structure logarithmique triviale et posons 
\begin{equation}
(X',\cM_{X'})=(X,\cM_X)\times_{(Y,\cM_Y)}(Y',\cM_{Y'}),
\end{equation} 
le produit étant indifféremment pris dans la catégorie des schémas logarithmiques ou dans celle des schémas logarithmiques saturés (\cite{agt} II.5.20),
de sorte que $X'=X\times_Y\oy'$. On voit aussitôt que $X'^\circ=X'\times_XX^\circ$ est l'ouvert maximal de $X'$ 
où la structure logarithmique $\cM_{X'}$ est triviale. 
Comme $f$ est plat, il existe un point géométrique $\ox''$ de $X$ au-dessus de $\oy'$, se spécialisant en $\ox$ (\cite{ega4} 2.3.4). 
On peut donc considérer $\ox''$ comme un point géométrique de $X'$. 
Le morphisme $(X',\cM_{X'})\rightarrow (Y',\cM_{Y'})$ étant lisse et saturé,  
le schéma logarithmique $(X',\cM_{X'})$ est saturé (\cite{tsuji4} II.2.12) et régulier (\cite{kato2} 8.2);  
cf. aussi (\cite{niziol} 2.3) et la preuve de (\cite{tsuji1} 1.5.1). 
Par suite, $X'^\circ$ est un ouvert schématiquement dense de $X'$
en vertu de (\cite{kato2} 11.6); cf. aussi (\cite{niziol} 2.6). Il existe donc un point géométrique $\ox'$ de $X'$, se spécialisant en $\ox''$. 
Considérant $\ox'$ comme point géométrique de $X$, il est clairement au-dessus de $\oy'$ et il se spécialise en $\ox$.

\begin{lem}\label{cgl2}
Soient $k$ un corps, $f\colon X\rightarrow Y$ un morphisme de $k$-schémas lisses, $D$ un diviseur à croisements normaux strict de $X$. 
On munit $X$ de la structure logarithmique $\cM_X$ définie par $D$ et $Y$ de la structure logarithmique triviale $\co^\times_Y$,
et on suppose que le morphisme de schémas logarithmiques $(X,\cM_X)\rightarrow (Y,\co^\times_Y)$, induit par $f$, est lisse. 
Alors, le morphisme de schémas $f$ est lisse et le diviseur $D$ est strictement à croisements normaux relativement à $Y$ {\rm (\cite{sga1} XIII 2.1)}.
\end{lem}

Soient $Z$ une composante irréductible de $D$, $i\colon Z\rightarrow X$ l'injection canonique,
$D'$ le diviseur à croisements normaux strict de $X$ défini par les composantes irréductibles de $D$ autres que $Z$. 
Alors, $Z$ est un $k$-schéma lisse et $i^*(D')$ est un diviseur à croisements normaux strict de $Z$. 
On désigne par $\cM'_X$ (resp. $\cM'_Z$) la structure logarithmique sur $X$ (resp. $Z$) définie par $D'$ (resp. $i^*(D')$).
Il est commode de noter le $\co_X$-module des différentielles logarithmiques de $(X,\cM_X)$ 
sur $(Y,\co_Y^\times)$ par $\Omega^1_{X/Y}(\log D)$; 
on utilisera des notations similaires en remplaçant $X$ par $Z$ ou $Y$ par $\Spec(k)$ ou $D$ par $D'$.  
Procédant par récurrence sur le nombre de composantes irréductibles de $D$,  
il suffit de montrer que les morphismes de schémas logarithmiques $(X,\cM'_X)\rightarrow (Y,\co^\times_Y)$ 
et $(Z,\cM'_Z)\rightarrow (Y,\co^\times_Y)$ induits par $f$ sont lisses. 

Considérons le diagramme commutatif de morphismes canoniques
\begin{equation}\label{cgl2a}
\xymatrix{
0\ar[r]&{f^*(\Omega^1_{Y/k})}\ar@{=}[d]\ar[r]&{\Omega^1_{X/k}(\log D')}\ar[r]\ar[d]^u&{\Omega^1_{X/Y}(\log D')}\ar[r]\ar[d]&0\\
0\ar[r]&{f^*(\Omega^1_{Y/k})}\ar[r]&{\Omega^1_{X/k}(\log D)}\ar[r]&{\Omega^1_{X/Y}(\log D)}\ar[r]&0}
\end{equation}
Comme le morphisme  $(X,\cM_X)\rightarrow (Y,\co^\times_Y)$ est lisse, la suite horizontale inférieure est exacte et localement scindée.
Il en est alors de même de la suite horizontale supérieure, ce qui implique que  le morphisme $(X,\cM'_X)\rightarrow (Y,\co^\times_Y)$ est lisse 
(\cite{ogus} IV 3.2.3). 

Considérons ensuite le diagramme commutatif de morphismes canoniques
\begin{equation}\label{cgl2b}
\xymatrix{
&&0&0&\\
0\ar[r]&{\co_Z(-Z)}\ar@{=}[d]\ar[r]&{i^*(\Omega^1_{X/Y}(\log D'))}\ar[r]\ar[u]&{\Omega^1_{Z/Y}(\log D')}\ar[r]\ar[u]&0\\
0\ar[r]&{\co_Z(-Z)}\ar[r]^-(0.5)v&{i^*(\Omega^1_{X/k}(\log D'))}\ar[r]\ar[u]&{\Omega^1_{Z/k}(\log D')}\ar[r]\ar[u]&0\\
&&{i^*(f^*(\Omega^1_{Y/k}))}\ar@{=}[r]\ar[u]^-(0.5)w&{i^*(f^*(\Omega^1_{Y/k}))}\ar[u]&\\
&&0\ar[u]&0\ar[u]&}
\end{equation}
La suite horizontale inférieure est exacte et localement scindée d'après (\cite{ogus} IV 3.2.2). Le morphisme composé 
\begin{equation}\label{cgl2c}
\xymatrix{
{\co_Z(-Z)}\ar[r]^-(0.5)v&{i^*(\Omega^1_{X/k}(\log D'))}\ar[r]^-(0.4){i^*(u)}&{i^*(\Omega^1_{X/k}(\log D))}}
\end{equation}
est nul. En effet, si $t$ est un paramètre local définissant le diviseur de Cartier $Z$ sur $X$, alors l'image de $i^*(t)\in \co_Z(-Z)$ par le composé
est $i^*(dt)=i^*(t d\log t)$ qui est nulle dans $i^*(\Omega^1_{X/k}(\log D))$.  

Par ailleurs, le morphisme $(X,\cM'_X)\rightarrow (Y,\co^\times_Y)$ étant lisse,  
la suite verticale de gauche est exacte et localement scindée d'après (\cite{ogus} IV 3.2.3). 
Le morphisme composé 
\begin{equation}\label{cgl2d}
\xymatrix{
{i^*(f^*(\Omega^1_{Y/k}))}\ar[r]^-(0.5)w&{i^*(\Omega^1_{X/k}(\log D'))}\ar[r]^-(0.4){i^*(u)}&{i^*(\Omega^1_{X/k}(\log D))}}
\end{equation}
est injectif puisque la suite horizontale inférieure de \eqref{cgl2a} est exacte et localement scindée. On a donc 
\begin{equation}
0=\co_Z(-Z)\cap i^*(f^*(\Omega^1_{Y/k}))\subset i^*(\Omega^1_{X/k}(\log D')).
\end{equation}
On en déduit que la suite horizontale supérieure et la suite verticale de droite de \eqref{cgl2b} sont exactes. Considérons un scindage local
\begin{equation}
\sigma\colon \Omega^1_{X/k}(\log D)\rightarrow f^*(\Omega^1_{Y/k})
\end{equation}
de la suite exacte inférieure de \eqref{cgl2a} et posons 
\begin{equation}
\sigma'=\sigma\circ u\colon \Omega^1_{X/k}(\log D')\rightarrow f^*(\Omega^1_{Y/k}).
\end{equation}
On a alors $i^*(\sigma')(w(x))=x$ pour toute section locale $x$ de $i^*(f^*(\Omega^1_{Y/k}))$, 
et $i^*(\sigma')(v(y))=0$ pour toute section locale $y$ de $\co_Z(-Z)$. 
Par suite, $i^*(\sigma')$ induit un scindage local de la suite verticale de droite de \eqref{cgl2b}. 
On en déduit que le morphisme $(Z,\cM'_Z)\rightarrow (Y,\co^\times_Y)$ est lisse  (\cite{ogus} IV 3.2.3).

\section{Rappel sur une construction de Fontaine-Grothendieck}\label{EIPF}

\subsection{}\label{eip1}
La construction rappelée dans ce numéro est due à Grothendieck (\cite{grot1} IV 3.3). 
On fixe un nombre premier $p$. Tous les anneaux des vecteurs de Witt considérés dans ce numéro sont relatifs à $p$ (cf. \ref{notconv1}). 
Soient $A$ une $\mZ_{(p)}$-algèbre, $n$ un entier $\geq 1$. L'homomorphisme d'anneaux 
\begin{equation}\label{eip1a}
\Phi_{n+1}\colon 
\begin{array}[t]{clcr}
\rW_{n+1}(A/p^nA)&\rightarrow& A/p^nA\\
(x_0,\dots,x_{n})&\mapsto& x_0^{p^n}+p x_1^{p^{n-1}}+\dots+p^{n}x_n
\end{array}
\end{equation}
s'annule sur $\rV^n(A/p^nA)$ et induit donc par passage au quotient un homomorphisme d'anneaux
\begin{equation}\label{eip1b}
\Phi'_{n+1}\colon 
\begin{array}[t]{clcr}
\rW_{n}(A/p^nA)&\rightarrow& A/p^nA\\
(x_0,\dots,x_{n-1})&\mapsto&x_0^{p^n}+p x_1^{p^{n-1}}+\dots+p^{n-1}x_{n-1}^p.
\end{array}
\end{equation}
Ce dernier s'annule sur 
\begin{equation}\label{eip1c}
\rW_n(pA/p^nA)=\ker(\rW_n(A/p^nA)\rightarrow \rW_n(A/pA))
\end{equation}
et se factorise à son tour en un homomorphisme d'anneaux
\begin{equation}\label{eip1d}
\uptheta_n\colon \rW_{n}(A/pA)\rightarrow A/p^nA.
\end{equation}
Il résulte aussitôt de la définition que le diagramme 
\begin{equation}\label{eip1e}
\xymatrix{
{\rW_{n+1}(A/pA)}\ar[r]^-(0.4){\uptheta_{n+1}}\ar[d]_{\rR\rF}&{A/p^{n+1}A}\ar[d]\\
{\rW_{n}(A/pA)}\ar[r]^-(0.4){\uptheta_n}&{A/p^nA}}
\end{equation}
où $\rR$ est le morphisme de restriction \eqref{notconv1b}, 
$\rF$ est le Frobenius \eqref{notconv1d} et la flèche non libellée est l'homomorphisme canonique, est commutatif. 

Pour tout homomorphisme de $\mZ_{(p)}$-algèbres commutatives $\varphi\colon A\rightarrow B$, 
le diagramme  
\begin{equation}\label{eip1f}
\xymatrix{
{\rW_n(A/pA)}\ar[r]\ar[d]_{\uptheta_n}&{\rW_n(B/pB)}\ar[d]^{\uptheta_n}\\
{A/p^nA}\ar[r]&{B/p^nB}}
\end{equation}
où les flèches horizontales sont les morphismes induits par $\varphi$, est commutatif.

\subsection{}\label{eipo3}
Soient $A$ une $\mZ_{(p)}$-algèbre, $\hA$ le séparé complété $p$-adique de $A$. 
On désigne par $A^\flat$ la limite projective du système projectif $(A/pA)_{\mN}$ 
dont les morphismes de transition sont les itérés de l'endomorphisme de Frobenius de $A/pA$;
\begin{equation}\label{eipo3a}
A^\flat=\underset{\underset{\mN}{\longleftarrow}}{\lim}A/pA.
\end{equation} 
C'est un anneau parfait de caractéristique $p$.
Pour tout entier $n\geq 1$, la projection canonique $A^\flat\rightarrow A/pA$ sur la $(n+1)$-ième composante
du système projectif $(A/pA)_{\mN}$ ({\em i.e.}, la composante d'indice $n$) induit un homomorphisme \eqref{notconv1}
\begin{equation}\label{eipo3b}
\nu_n\colon \rW(A^\flat)\rightarrow \rW_n(A/pA).
\end{equation}
Comme $\nu_n=\rF\circ\rR\circ \nu_{n+1}$, on obtient par passage à la limite projective un homomorphisme 
\begin{equation}\label{eipo3c}
\nu\colon \rW(A^\flat)\rightarrow \underset{\underset{n\geq 1}{\longleftarrow}}{\lim}\ \rW_n(A/pA),
\end{equation}
où les morphismes de transition de la limite projective sont les morphismes $\rF \rR$. On vérifie aussitôt  
qu'il est bijectif. Compte tenu de \eqref{eip1e},
les homomorphismes $\uptheta_n$ induisent par passage à la limite projective un 
homomorphisme 
\begin{equation}\label{eipo3d}
\theta\colon \rW(A^\flat)\rightarrow \hA.
\end{equation}
On retrouve l'homomorphisme défini par Fontaine (\cite{fontaine1} 2.2). 
Pour tout entier $r\geq 1$, on pose
\begin{equation}\label{eipo3e}
\cA_r(A)=\rW(A^\flat)/\ker(\theta)^r,
\end{equation}
et on note $\theta_r\colon \cA_r(A)\rightarrow \hA$ l'homomorphisme induit par $\theta$ (cf. \cite{fontaine3} 1.2.2). 
On prendra garde de ne pas confondre $\theta_r$ et $\uptheta_r$. 

Pour tout homomorphisme de $\mZ_{(p)}$-algèbres commutatives $\varphi\colon A\rightarrow B$, 
le diagramme  
\begin{equation}\label{eip3f}
\xymatrix{
{\rW(A^\flat)}\ar[r]\ar[d]_{\theta}&{\rW(B^\flat)}\ar[d]^{\theta}\\
{\hA}\ar[r]&{\hB}}
\end{equation}
où les flèches horizontales sont les morphismes induits par $\varphi$, est commutatif \eqref{eip1f}.  
La correspondance $A\mapsto \cA_r(A)$ est donc fonctorielle. 

\begin{rema}\label{eip2}
Soit $k$ un corps parfait.  
La projection canonique $\rW(k)^\flat\rightarrow k$ sur la première composante ({\em i.e.}, d'indice $0$)
est un isomorphisme. Elle induit donc un isomorphisme $\rW(\rW(k)^\flat)\stackrel{\sim}{\rightarrow}\rW(k)$,
que nous utilisons pour identifier ces deux anneaux. L'homomorphisme $\theta$ s'identifie alors à
l'endomorphisme identique de $\rW(k)$.
\end{rema}

\begin{lem}\label{eip5}
Soient $A$ une $\mZ_{(p)}$-algèbre, $\hA$ le séparé complété $p$-adique de $A$, 
$A^\flat$ l'anneau défini dans \eqref{eipo3a}, $\mA$ l'ensemble des suites $(x_n)_{\mN}$ de $\hA$ 
telles que $x_{n+1}^p= x_n$ pour tout $n\geq 0$. Alors, 
\begin{itemize}
\item[{\rm (i)}] L'application
\begin{equation}\label{eip5a}
\mA\rightarrow A^\flat, \ \ \ (x_n)_{\mN}\mapsto (\ox_n)_{\mN},
\end{equation}
où $\ox_n$ est la réduction de $x_n$ modulo $p$, est un isomorphisme de monoïdes multiplicatifs. 
\item[{\rm (ii)}] 
Pour tout $(x_0,x_1,\dots)\in \rW(A^\flat)$, on a 
\begin{equation}\label{eip5b}
\theta(x_0,x_1,\dots) = \sum_{n\geq 0}p^nx_n^{(n)},
\end{equation}
où $\theta$ est l'homomorphisme \eqref{eipo3d} et pour tout $n\geq 0$, $(x_n^{(m)})_{m\geq 0}$ est l'élément de $\mA$ associé à $x_n\in A^\flat$ \eqref{eip5a}.
\end{itemize}
\end{lem}

(i) Notant $v$ la valuation de $\mZ_{(p)}$ normalisée par $v(p)=1$, pour tous entiers $m,i\geq 1$ tels que $i\leq p^m$, on a $v(\binom{p^m}{i})=m-v(i)$
et donc $v(\binom{p^m}{i})+i\geq m$. 
Soient $(x_n)_{\mN}$ et $(y_n)_{\mN}$ deux suites d'éléments de $\hA$ qui induisent par réduction modulo $p$ la même suite 
$(\ox_n)_{\mN}=(\oy_n)_{\mN}$ de $A/pA$. Pour tous $n,m\geq 0$, on a alors
\begin{equation}\label{eip5c}
x_{n+m}^{p^m}-y_{n+m}^{p^m} \in p^m \hA.
\end{equation}
En particulier, si $(x_n)_{\mN}$ et $(y_n)_{\mN}$ sont des éléments de $\mA$, 
alors $x_n=y_n$ pour tout $n\geq 0$ puisque $\hA$ est séparé pour la topologie $p$-adique.
Par suite, l'application \eqref{eip5a} est injective. Montrons qu'elle est surjective. Soient $(z_n)_{\mN}$ un élément de $A^\flat$, 
$(y_n)_{\mN}$ une suite d'éléments de $A$ qui relève $(z_n)_{\mN}$. Appliquant \eqref{eip5c} aux suites $(y_n)_{n\geq 0}$ et $(y_{n+1}^p)_{n\geq 0}$,
on voit que pour tous $n,m\geq 0$, on a 
\begin{equation}\label{eip5d}
y_{n+m+1}^{p^{m+1}}-y_{n+m}^{p^m} \in p^m \hA.
\end{equation}
Par suite, pour tout entier $n\geq 0$,
$(y_{n+m}^{p^m})_{m\geq 0}$ converge vers un élément $x_n$ de $\hA$. La suite $(x_n)_{\mN}$ appartient clairement à $\mA$ et s'envoie 
sur $(z_n)_{\mN}$ par \eqref{eip5a}; d'où la surjectivité.  

(ii) Cela résulte aussitôt des définitions.

\begin{prop}[\cite{agt} II.9.5, \cite{tsuji1} A.1.1 et A.2.2]\label{eip4}
Soit $A$ une $\mZ_{(p)}$-algèbre commutative vérifiant les conditions suivantes:
\begin{itemize}
\item[{\rm (i)}] $A$ est $\mZ_{(p)}$-plat.
\item[{\rm (ii)}] $A$ est intégralement clos dans $A[\frac 1 p]$.
\item[{\rm (iii)}] Le Frobenius absolu  de $A/pA$ est surjectif. 
\item[{\rm (iv)}] Il existe une suite $(p_n)_{n\geq 0}$ d'éléments de $A$
tels que $p_0=p$ et $p_{n+1}^p=p_n$ pour tout $n\geq 0$.  
\end{itemize}
On désigne par $\upp$ l'élément de $A^\flat$ induit par la suite $(p_n)_{n\geq 0}$ et on pose  
\begin{equation}\label{eip4a}
\xi=[\upp]-p \in \rW(A^\flat),
\end{equation}
où $[\ ]$ est le représentant multiplicatif.
Alors, la suite
\begin{equation}\label{eip4b}
0\longrightarrow \rW(A^\flat)\stackrel{\cdot \xi}{\longrightarrow} \rW(A^\flat)\stackrel{\theta}{\longrightarrow} 
\hA \longrightarrow 0
\end{equation}
est exacte.
\end{prop}

\section{\texorpdfstring{Faisceaux de $\alpha$-modules}{Faisceaux de alpha-modules}}\label{alpha}

\subsection{}\label{alpha1}
Dans cette section, $\Lambda$ désigne un sous-groupe ordonné et dense de $\mR$, 
$\Lambda^+$ l'ensemble des éléments strictement positifs de $\Lambda$, et 
$R$ un anneau muni d'une suite d'idéaux principaux $\fm_\varepsilon$, indexée par $\Lambda^+$. 
Pour tout $\varepsilon\in \Lambda^+$, on choisit un générateur $\pi^\varepsilon\in R$ de $\fm_\varepsilon$. 
On suppose, de plus, que pour tout $\varepsilon\in \Lambda^+$, $\pi^\varepsilon$ n'est pas un diviseur de zéro dans $R$ 
et que pour tous $\varepsilon,\delta\in \Lambda^+$, 
$\pi^\varepsilon\cdot \pi^\delta=u_{\varepsilon,\delta}\cdot \pi^{\varepsilon+\delta}$, 
où $u_{\varepsilon,\delta}$ est une unité de $R$. On pose $\fm=\cup_{\varepsilon\in \Lambda^+}\fm_\varepsilon$. 
Ces hypothèses sont celles fixées dans (\cite{faltings2} page 186 et \cite{agt} V Notations). 
On notera que $\fm$ est $R$-plat et que $\fm^2=\fm$. 

\subsection{}\label{alpha2}
Soient $\cA$ une catégorie abélienne qui soit une $\mU$-catégorie \eqref{notconv3} (\cite{sga4} I 1.1), 
$\End(\id_\cA)$ l'anneau des endomorphismes du foncteur identique de $\cA$, $\varphi\colon R\rightarrow \End(\id_\cA)$ un homomorphisme.  
Pour tout objet $M$ de $\cA$ et tout $\gamma\in R$, on note $\mu_\gamma(M)$ l'endomorphisme de $M$ défini par $\varphi(\gamma)$. 
On observera que pour tout morphisme $f\colon M\rightarrow N$ de $\cA$, on a $\mu_\gamma(N)\circ f= f\circ \mu_\gamma(M)$.
En particulier, pour tous objets $M$ et $N$ de $\cA$, $\Hom(M,N)$ est naturellement muni d'une structure de $R$-module. 

Suivant (\cite{faltings2} page 187), on dit qu'un objet $M$ de $\cA$  
est {\em $\alpha$-nul} (ou {\em presque-nul}) s'il est annulé par tout élément de $\fm$, 
{\em i.e.}, si  $\mu_\gamma(M)=0$ pour tout $\gamma\in \fm$.
On vérifie aussitôt que pour toute suite exacte $0\rightarrow M'\rightarrow M\rightarrow M''\rightarrow 0$ de $\cA$, pour que 
$M$ soit $\alpha$-nul, il faut et il suffit que $M'$ et $M''$ soient $\alpha$-nuls. 
On appelle catégorie des {\em $\alpha$-objets} (ou {\em presque-objets}) de $\cA$ et l'on note 
$\acA$ le quotient de la catégorie $\cA$ par la sous-catégorie épaisse formée des objets $\alpha$-nuls (\cite{gabriel} III §~1). On note 
\begin{equation}\label{alpha2a}
\alpha\colon \cA\rightarrow \acA, \ \ \ M\mapsto \alpha(M)
\end{equation}
le foncteur canonique~; on notera aussi $M^\alpha$ au lieu de $\alpha(M)$ lorsqu'il n'y a aucun risque de confusion. 
La catégorie $\acA$ est abélienne et le foncteur $\alpha$ est exact (\cite{gabriel} III §~1 prop.~1).
Si $f\colon M\rightarrow N$ est un morphisme  de $\cA$, pour que $\alpha(f)$ soit nul (resp. un monomorphisme,
resp. un épimorphisme), il faut et il suffit que $\im(f)$ (resp. $\ker(f)$, resp. $\coker(f)$) soit $\alpha$-nul 
(\cite{gabriel} III §~1 lem.~2). On dit que $f$ est {\em $\alpha$-injectif} (resp. {\em $\alpha$-surjectif}, resp. un {\em $\alpha$-isomorphisme})  
si $\alpha(f)$ est injectif (resp. surjectif, resp. un isomorphisme), autrement dit si son noyau (resp. son conoyau, resp. son noyau et son conoyau) sont $\alpha$-nuls. 

La famille des $\alpha$-isomorphismes de $\cA$ permet un calcul de fractions bilatère (\cite{illusie1} I 1.4.2).
La catégorie $\acA$ s'identifie à la catégorie localisée de $\cA$ 
par rapport aux $\alpha$-isomorphismes et $\alpha$ est le foncteur canonique (de localisation). 
On laissera le soin au lecteur de vérifier ces propriétés, valables d'ailleurs pour tout quotient d'une catégorie 
abélienne par une sous-catégorie épaisse (\cite{gz} I 2.5(d)).

\begin{lem}\label{alpha3}
Les hypothèses étant celles de \eqref{alpha2}, soient, de plus, $f\colon M\rightarrow N$ un morphisme de $\cA$, $\gamma \in R$ tels que 
le noyau et le conoyau de $f$ soient annulés par $\gamma$. Alors, il existe un morphisme $g\colon N\rightarrow M$ de $\cA$ tel que 
$g\circ f=\mu_{\gamma^2}(M)$ et $f\circ g=\mu_{\gamma^2}(N)$.
\end{lem}

En effet, comme $\mu_\gamma(\coker(f))=0$, le morphisme $\mu_\gamma(N)\colon N\rightarrow N$ se factorise uniquement 
à travers un morphisme $\varphi_\gamma\colon N\rightarrow \im(f)$. 
Comme $\mu_\gamma(\ker(f))=0$, le morphisme $\mu_\gamma(M)\colon M\rightarrow M$ se factorise uniquement 
à travers un morphisme $\psi_\gamma\colon \im(f)\rightarrow M$. 
Il est clair que $g=\psi_\gamma\circ \varphi_\gamma\colon N\rightarrow M$ répond à la question.

\subsection{}\label{alpha33}
Soit $\cA$ une catégorie abélienne tensorielle qui soit une $\mU$-catégorie, autrement dit, $\cA$ est une catégorie abélienne
munie d'un foncteur bi-additif  $\otimes\colon \cA\times \cA\rightarrow \cA$ et d'un objet unité $A$, 
vérifiant certaines conditions (\cite{dm} 1.15), et soit $\phi\colon R\rightarrow \End(A)$ un homomorphisme. 
On a un homomorphisme canonique $\End(A)\rightarrow \End(\id_\cA)$ (cf. \cite{sr} I 1.3.3.3 et 2.3.3). 
On peut donc définir la catégorie quotient $\acA$ suivant \ref{alpha2}. Il résulte aussitôt de \ref{alpha3} que pour tout 
$\alpha$-isomorphisme $f\colon M\rightarrow M'$ de $\cA$ et tout $N\in \ob(\cA)$, $f\otimes \id_N\colon M\otimes N\rightarrow M'\otimes N$ 
est un $\alpha$-isomorphisme. Par suite, le produit tensoriel induit un foncteur 
\begin{equation}\label{alpha33a}
\acA\times \acA\rightarrow \acA, \ \ \ (M,N)\mapsto M\otimes N,
\end{equation}
qui fait de $\acA$ une catégorie abélienne tensorielle, dont $A^\alpha$ est un objet unité. 

Le foncteur $\alpha$ induit un homomorphisme $\End(A)\rightarrow \End(A^\alpha)$. Par suite, pour tous objets $M$ et $N$ de $\acA$, 
$\Hom_{\acA}(M,N)$ est canoniquement muni d'une structure de $\End(A)$-module. 
On voit aussitôt que pour tout $P\in \ob(\cC)$, l'application 
\begin{equation}\label{alpha33b}
\Hom_{\acA}(M,N)\rightarrow \Hom_{\acA}(M\otimes P,N\otimes P)
\end{equation}
définie par fonctorialité est $\End(A)$-linéaire (\cite{sr} I 2.2.6). 
 
\subsection{}\label{alpha4}
Pour toute $R$-algèbre $A$ (appartenant à $\mU$), on désigne par $\bMod(A)$ la catégorie abélienne tensorielle 
des $A$-modules qui se trouvent dans $\mU$. Prenant pour $\phi\colon R\rightarrow A=\End(A)$ l'homomorphisme structural,  
on appelle catégorie des {\em $\alpha$-$A$-modules} et l'on note $\aMod(A)$ 
le quotient de la catégorie abélienne $\bMod(A)$ par la sous-catégorie épaisse des $A$-modules $\alpha$-nuls. 
Nous utiliserons les conventions de notation de \ref{alpha2} et \ref{alpha33}. On observera en particulier que pour tous
$\alpha$-$A$-modules $M$ et $N$, $\Hom_{\aMod(A)}(M,N)$ est naturellement muni d'une structure de $A$-module.  

Pour tout homomorphisme de $R$-algèbres $A\rightarrow B$, le foncteur d'oubli  
$\bMod(B)\rightarrow \bMod(A)$ induit un foncteur exact
\begin{equation}\label{alpha4a}
\aMod(B)\rightarrow \aMod(A).
\end{equation}

\begin{lem}\label{alpha7}
Pour qu'un morphisme de $R$-modules $M\rightarrow N$ soit un $\alpha$-isomorphisme,
il faut et il suffit que le morphisme induit $\fm\otimes_RM\rightarrow \fm\otimes_RN$ soit un isomorphisme. 
\end{lem}
En effet, la condition est suffisante puisque le morphisme canonique 
$\fm\otimes_RM\rightarrow M$ est un $\alpha$-isomorphisme \eqref{alpha3},
et elle est nécessaire car $\fm$ est $R$-plat et pour tout $R$-module $\alpha$-nul $P$, $\fm\otimes_RP=0$. 

\subsection{}\label{alpha8} 
Soit $A$ une $R$-algèbre. 
D'après \ref{alpha7}, la catégorie des $\alpha$-isomorphismes de $\bMod(A)$ de but un $A$-module $M$ donné 
admet un objet initial à savoir le morphisme canonique $\fm\otimes_RM\rightarrow M$. 
Par suite, pour tous $A$-modules $M$ et $N$, on a un isomorphisme canonique
\begin{equation}\label{alpha8a}
\Hom_{\aMod(A)}(M^\alpha,N^\alpha)\stackrel{\sim}{\rightarrow}\Hom_{\bMod(A)}(\fm\otimes_RM,N).
\end{equation}
On voit aussitôt que cet isomorphisme est $A$-linéaire.

On désigne par $\sigma_*$ le foncteur 
\begin{equation}\label{alpha8b}
\sigma_*\colon \aMod(A)\rightarrow \bMod(A),\ \ \ P\mapsto \Hom_{\aMod(A)}(A^\alpha,P),
\end{equation}
et par $\sigma_!$ le foncteur
\begin{equation}\label{alpha8c}
\sigma_!\colon \aMod(A)\rightarrow \bMod(A), \ \ \ P\mapsto \fm\otimes_R\sigma_*(P).
\end{equation}
D'après \eqref{alpha8a}, pour tout $A$-module $M$, on a un isomorphisme canonique fonctoriel
\begin{equation}\label{alpha8e}
\sigma_*(M^\alpha)\stackrel{\sim}{\rightarrow}\Hom_{R}(\fm,M).
\end{equation}
On en déduit que pour tout homomorphisme de $R$-algèbres $A\rightarrow B$, les diagrammes
\begin{equation}\label{alpha8d}
\xymatrix{
{\aMod(B)}\ar[r]^{\sigma_*}\ar[d]&{\bMod(B)}\ar[d]\\
{\aMod(A)}\ar[r]^{\sigma_*}&{\bMod(A)}}\ \ \
\xymatrix{
{\aMod(B)}\ar[r]^{\sigma_!}\ar[d]&{\bMod(B)}\ar[d]\\
{\aMod(A)}\ar[r]^{\sigma_!}&{\bMod(A)}}
\end{equation}
où les flèches verticales sont les foncteurs d'oubli \eqref{alpha4a}, sont commutatifs à isomorphismes canoniques près~;
ce qui justifie l'abus d'omettre $A$ dans les notations $\sigma_*$ et $\sigma_!$.

\begin{prop}\label{alpha9}
Soit $A$ une $R$-algèbre.
\begin{itemize}
\item[{\rm (i)}] Le foncteur $\sigma_*$ \eqref{alpha8b} est un adjoint à droite du foncteur de localisation $\alpha$ \eqref{alpha2a}.
\item[{\rm (ii)}] Le morphisme d'adjonction $\alpha\circ \sigma_* \rightarrow \id$ est un isomorphisme, 
{\em i.e.}, le foncteur $\sigma_*$ est pleinement fidèle.
\item[{\rm (iii)}] Le morphisme d'adjonction $\id\rightarrow \sigma_*\circ\alpha$
induit un isomorphisme $\alpha\stackrel{\sim}{\rightarrow} \alpha\circ\sigma_*\circ\alpha$.
\item[{\rm (iv)}] Le foncteur $\sigma_!$ \eqref{alpha8c} est un adjoint à gauche du foncteur de localisation $\alpha$. 
\item[{\rm (v)}]  Le morphisme d'adjonction $\id \rightarrow \alpha\circ \sigma_!$ 
est un isomorphisme, {\em i.e.}, le foncteur $\sigma_!$ est pleinement fidèle.
\end{itemize}
\end{prop}

Soient $M$, $N$ deux $A$-modules, $P$ un $\alpha$-$A$-module.

(i) On a des isomorphismes canoniques fonctoriels
\begin{eqnarray}
\Hom_{\aMod(A)}(M^\alpha,N^\alpha)&\stackrel{\sim}{\rightarrow}&\Hom_A(\fm\otimes_RM,N)\\
&\stackrel{\sim}{\rightarrow}&\Hom_A(M,\Hom_R(\fm,N))\nonumber\\
&\stackrel{\sim}{\rightarrow}&\Hom_A(M,\sigma_*(N^\alpha)).\nonumber\label{alpha9a}
\end{eqnarray}
L'assertion s'ensuit compte tenu de (\cite{gz} I 1.2).

(ii) Le morphisme d'adjonction
$\alpha(\sigma_*(M^\alpha))\rightarrow M^\alpha$ correspond par \eqref{alpha8a} et \eqref{alpha8e} au morphisme canonique
\begin{equation}\label{alpha9b}
\fm\otimes_R\Hom_R(\fm,M)\rightarrow M,
\end{equation} 
qui est un $\alpha$-isomorphisme car les morphismes canoniques $\fm\otimes_RM\rightarrow M$ et  
$M\rightarrow \Hom_R(\fm,M)$ sont des $\alpha$-isomorphismes \eqref{alpha3}.

(iii) Le morphisme d'adjonction $M\rightarrow \sigma_*(\alpha(M))$  
s'identifie par \eqref{alpha8e} au morphisme canonique $M\rightarrow \Hom_R(\fm,M)$
qui est un $\alpha$-isomorphisme. 

(iv) D'après (ii), on a des isomorphismes canoniques fonctoriels
\begin{eqnarray}
\Hom_{\aMod(A)}(P,\alpha(M))&\stackrel{\sim}{\rightarrow}&\Hom_A(\alpha(\sigma_*(P)),\alpha(M))\\
&\stackrel{\sim}{\rightarrow}&\Hom_A(\fm\otimes_R\sigma_*(P),M)\nonumber\\
&\stackrel{\sim}{\rightarrow}&\Hom_A(\sigma_!(P),M).\nonumber\label{alpha9c}
\end{eqnarray}

(v) Le morphisme d'adjonction $M^\alpha\rightarrow \alpha(\sigma_!(M^\alpha))$ correspond 
par \eqref{alpha8a} et \eqref{alpha8e} au morphisme canonique 
\begin{equation}\label{alpha9d}
\fm\otimes_RM\rightarrow \fm\otimes_R\Hom_R(\fm,M),
\end{equation}
qui est un isomorphisme \eqref{alpha7}.

\begin{cor}\label{alpha10}
Soit $A$ une $R$-algèbre. 
\begin{itemize}
\item[{\rm (i)}] Le foncteur de localisation $\alpha$ \eqref{alpha2a} 
commute aux limites inductives (resp. projectives) représentables et
le foncteur $\sigma_*$ \eqref{alpha8b} (resp. $\sigma_!$ \eqref{alpha8c}) commute aux limites projectives (resp. inductives) 
représentables.
\item[{\rm (ii)}] Les foncteurs $\alpha$ et $\sigma_!$ sont exacts, et le foncteur $\sigma_*$ est exact à gauche. 
\item[{\rm (iii)}] Pour toute catégorie $\mU$-petite $I$ et tout foncteur $\varphi\colon I\rightarrow \aMod(A)$, les
limites projective et inductive de $\varphi$ sont représentables et les morphismes canoniques 
\begin{eqnarray}
\underset{\underset{I}{\rightarrow}}{\lim}\ \varphi \rightarrow 
\alpha(\underset{\underset{I}{\rightarrow}}{\lim}\ \sigma_*\circ \varphi)\label{alpha10a}\\
\alpha(\underset{\underset{I}{\leftarrow}}{\lim}\ \sigma_!\circ \varphi)
\rightarrow \underset{\underset{I}{\leftarrow}}{\lim}\ \varphi\label{alpha10b}
\end{eqnarray}
sont des isomorphismes. 
\end{itemize}
\end{cor}

(i) Cela résulte de \ref{alpha9}(i)-(iv) et (\cite{sga4} I 2.11). 

(ii) En effet, il résulte de (i) que $\alpha$ est exact, $\sigma_*$ est exact à gauche et $\sigma_!$ est exact à droite. 
Comme $\fm$ est $R$-plat, $\sigma_!$ est aussi exact à gauche \eqref{alpha8c}. 

(iii) Cela résulte de (i), \ref{alpha9}(ii)-(v) et du fait que les $\mU$-limites inductives et projectives sont représentables dans 
la catégorie $\bMod(A)$.  

\begin{cor}\label{alpha100}
Pour toute $R$-algèbre $A$, la catégorie abélienne $\aMod(A)$ est une $\mU$-catégorie vérifiant la propriété {\rm (AB 5)} de 
{\rm (\cite{tohoku} §~1.5)} et admettant un générateur, à savoir $A^\alpha$. 
\end{cor}

En effet, comme $\alpha$ admet un adjoint à droite, à savoir $\sigma_*$ \eqref{alpha9}, la sous-catégorie épaisse de $\bMod(A)$ formée 
des $A$-modules $\alpha$-nuls est localisante dans le sens de (\cite{gabriel} p. 372). 
Par suite, d'après ({\em loc. cit.}, III §~2 lem.~4 et §~4 prop.~9), la catégorie $\aMod(A)$ est une $\mU$-catégorie avec générateur, 
à savoir $A^\alpha$, et limites inductives exactes~; d'où la proposition compte tenu de ({\em loc. cit.}, I §~6 prop.~6).

\subsection{}\label{alpha5}
Soit $A$ une $R$-algèbre. On appelle {\em $\alpha$-$A$-algèbre} (ou {\em $A^\alpha$-algèbre})
un monoïde unitaire commutatif de $\aMod(A)$. 
On désigne par $\bAlg(A)$ la catégorie des $A$-algèbres qui se trouvent dans $\mU$ et par $\aAlg(A)$ la catégorie des 
$\alpha$-$A$-algèbres. Le foncteur de localisation $\alpha$ \eqref{alpha2a}
étant monoïdal, il induit un foncteur que l'on note encore 
\begin{equation}\label{alpha5a}
\alpha\colon \bAlg(A)\rightarrow \aAlg(A).
\end{equation}

Compte tenu de l'isomorphisme canonique $A^\alpha\stackrel{\sim}{\rightarrow}A^\alpha\otimes_{A^\alpha}A^\alpha$,
le foncteur $\sigma_*$ \eqref{alpha8b} induit un foncteur que l'on note encore
\begin{equation}\label{alpha5b}
\sigma_*\colon \aAlg(A)\rightarrow \bAlg(A), \ \ \ P\mapsto \Hom_{\aMod(A)}(A^\alpha,P). 
\end{equation} 

Pour toute $A$-algèbre $B$, l'isomorphisme $\fm\stackrel{\sim}{\rightarrow} \fm\otimes_R\fm$ induit sur
$\Hom_R(\fm,B)$ une structure canonique de $A$-algèbre. On a un isomorphisme canonique fonctoriel de $A$-algèbres
\begin{equation}\label{alpha5c}
\sigma_*(B^\alpha)\stackrel{\sim}{\rightarrow}\Hom_{R}(\fm,B).
\end{equation}

\begin{prop}\label{alpha6}
Soit $A$ une $R$-algèbre.
\begin{itemize}
\item[{\rm (i)}] Le foncteur $\sigma_*$ \eqref{alpha5b} est un adjoint à droite du foncteur de localisation $\alpha$ \eqref{alpha5a}.
\item[{\rm (ii)}] Le morphisme d'adjonction $\alpha\circ \sigma_* \rightarrow \id$ est un isomorphisme, 
{\em i.e.}, le foncteur $\sigma_*$ est pleinement fidèle.
\item[{\rm (iii)}] Le morphisme d'adjonction $\id\rightarrow \sigma_*\circ\alpha$
induit un isomorphisme $\alpha\stackrel{\sim}{\rightarrow} \alpha\circ\sigma_*\circ\alpha$.
\end{itemize}
\end{prop}

Soient $B$, $C$ deux $A$-algèbres, $u\colon B^\alpha\rightarrow C^\alpha$ un morphisme de $\aMod(A)$,
$v\colon \fm\otimes_RB\rightarrow C$ et $w\colon B\rightarrow \Hom_R(\fm,C)$
les morphismes $A$-linéaires associés \eqref{alpha8a}. On vérifie aussitôt que les conditions suivantes sont équivalentes~:
\begin{itemize}
\item[(a)] $u$ est un morphisme de $A^\alpha$-algèbres.
\item[(b)]  Le diagramme 
\begin{equation}\label{alpha6b}
\xymatrix{
{(\fm\otimes_RB)\otimes_A(\fm\otimes_RB)}\ar[r]^-(0.5){v\otimes v}\ar@{=}[d]&{C\otimes_AC}\ar[r]^-(0.5){\mu_C}&C\\
{(\fm\otimes_R\fm)\otimes_R(B\otimes_AB)}\ar[rr]^-(0.5){\mu_\fm\otimes_R\mu_B}&&{\fm\otimes_RB}\ar[u]_v}
\end{equation}
où $\mu_\fm$, $\mu_B$ et $\mu_C$ désignent les morphismes de multiplication de $\fm$, $B$ et $C$, 
respectivement, est commutatif. 
\item[(c)] $w$ est  un morphisme de $A$-algèbres.
\end{itemize}
On en déduit des isomorphismes canoniques fonctoriels
\begin{eqnarray}
\Hom_{\aAlg(A)}(B^\alpha,C^\alpha)&\stackrel{\sim}{\rightarrow}&\Hom_{\bAlg(A)}(B,\Hom_R(\fm,C))\nonumber\\
&\stackrel{\sim}{\rightarrow}&\Hom_{\bAlg(A)}(B,\sigma_*(C^\alpha)).\label{alpha6a}
\end{eqnarray}
La proposition (i) s'ensuit compte tenu de (\cite{gz} I 1.2). On notera que les morphismes d'adjonction 
$\alpha\circ \sigma_* \rightarrow \id$ et $\id\rightarrow \sigma_*\circ\alpha$ s'identifient aux morphismes d'adjonction 
pour les $A$-modules \eqref{alpha9}. Les propositions (ii) et (iii) résultent alors de \ref{alpha9}(ii)-(iii).

\subsection{}\label{alpha39}
Dans ce numéro, si $D$ est une $R$-algèbre, nous affecterons d'un indice $D$ les foncteurs 
$\alpha$ \eqref{alpha2a} et $\sigma_*$ \eqref{alpha8b} pour les $D$-modules 
ainsi que leurs variantes \eqref{alpha5a} et \eqref{alpha5b} pour les $D$-algèbres.  

Soit $A$ une $R$-algèbre. On pose $B=\alpha_R(A)$ \eqref{alpha5a}  
et on désigne par $\bMod(B)$ la catégorie des $B$-modules unitaires de $\aMod(R)$. Le foncteur $\alpha_R$ étant monoïdal, 
il induit un foncteur 
\begin{equation}\label{alpha39a}
\beta\colon \bMod(A)\rightarrow \bMod(B).
\end{equation} 
Celui-ci transforme clairement les $\alpha$-isomorphismes en des isomorphismes. Il induit donc un foncteur 
\begin{equation}\label{alpha39b}
b\colon \aMod(A)\rightarrow \bMod(B).
\end{equation} 

On pose $A'=\sigma_{R*}(B)$ \eqref{alpha5b}. 
D'après \ref{alpha6}, on a un homomorphisme canonique de $R$-algèbres $\lambda\colon A\rightarrow A'$,
qui induit un isomorphisme $\alpha_R(A)\stackrel{\sim}{\rightarrow}\alpha_R(A')$ \eqref{alpha5a}. 
Pour tous $R^\alpha$-modules $P$ et $Q$, on a un morphisme $R$-linéaire
\begin{equation}\label{alpha39c}
\Hom_{\aMod(R)}(R^\alpha,P)\otimes_R \Hom_{\aMod(R)}(R^\alpha,Q)\rightarrow \Hom_{\aMod(R)}(R^\alpha,P\otimes_{R^\alpha}Q),
\end{equation}
défini par fonctorialité et composition \eqref{alpha33b}.
Par suite, le foncteur $\sigma_{R*}$ \eqref{alpha8b} induit un foncteur 
\begin{equation}\label{alpha39d}
\tau'_*\colon \bMod(B)\rightarrow \bMod(A').
\end{equation}
Composant avec le foncteur induit par $\lambda$, on obtient un foncteur
\begin{equation}\label{alpha39e}
\tau_*\colon \bMod(B)\rightarrow \bMod(A).
\end{equation}

Soient $M$ et $N$ deux $A$-modules, $u\colon \beta(M)\rightarrow \beta(N)$ un morphisme de $\aMod(R)$,
$v\colon \fm\otimes_RM\rightarrow N$ et $w\colon M\rightarrow \Hom_R(\fm,N)$
les morphismes $R$-linéaires associés à $u$ \eqref{alpha8a}. On vérifie aussitôt que $u$ est $B$-linéaire si et seulement si $v$
(ou ce qui revient au même $w$) est $A$-linéaire.  
L'isomorphisme \eqref{alpha8a} induit donc un isomorphisme canonique 
\begin{equation}
\Hom_{\bMod(B)}(\beta(M),\beta(N))\stackrel{\sim}{\rightarrow} \Hom_{\bMod(A)}(\fm\otimes_RM,N).
\end{equation}
Calquant alors la preuve de \ref{alpha9}, on en déduit que $\tau_*$ est un adjoint à droite de $\beta$, 
que le morphisme d'adjonction $\beta\circ \tau_* \rightarrow \id$ est un isomorphisme 
et que le morphisme d'adjonction $\id\rightarrow \tau_*\circ\beta$
induit un isomorphisme $\beta\stackrel{\sim}{\rightarrow} \beta\circ\tau_*\circ\beta$. 

On pose
\begin{equation}\label{alpha39f}
t=\alpha_A\circ \tau_*\colon \bMod(B)\rightarrow \aMod(A).
\end{equation}
L'isomorphisme $\beta\circ \tau_* \stackrel{\sim}{\rightarrow} \id$ induit un isomorphisme 
\begin{equation}\label{alpha39g}
b\circ t \stackrel{\sim}{\rightarrow} \id.
\end{equation}
Par ailleurs, le morphisme d'adjonction $\id\rightarrow \tau_*\circ\beta$
induit un isomorphisme $\alpha_A\stackrel{\sim}{\rightarrow} \alpha_A\circ\tau_*\circ\beta$. 
Compte tenu de (\cite{gz} I 1.2), on en déduit un isomorphisme 
\begin{equation}\label{alpha39h}
\id\stackrel{\sim}{\rightarrow} t\circ b.
\end{equation} 
On vérifie aussitôt que les isomorphismes \eqref{alpha39g} et \eqref{alpha39h} font de $t$ un adjoint à droite de $b$. 
Par suite, $b$ et $t$ sont des équivalences de catégories. 
On a des isomorphismes canoniques
\begin{equation}
\beta\stackrel{\sim}{\rightarrow} b\circ \alpha_A \ \ \ {\rm et} \ \ \ \alpha_A\stackrel{\sim}{\rightarrow} t \circ \beta,
\end{equation}
compatibles aux isomorphismes \eqref{alpha39g} et \eqref{alpha39h}. 
On en déduit aussitôt un isomorphisme
\begin{equation}
\sigma_A \stackrel{\sim}{\rightarrow} \tau_*\circ b.
\end{equation}

\subsection{}\label{alpha13}
Soit $\mV$ un univers tel que $\mU\subset \mV$. 
On affectera d'un indice $\mU$ ou $\mV$ les catégories et foncteurs dépendants de l'univers. 
Soit $A$ une $R$-algèbre appartenant à $\mU$. 
On a un foncteur pleinement fidèle canonique
\begin{equation}\label{alpha13a}
\bMod_\mU(A)\rightarrow \bMod_\mV(A).
\end{equation}
Celui-ci induit un foncteur
\begin{equation}\label{alpha13b}
\aMod_\mU(A)\rightarrow \aMod_\mV(A)
\end{equation}
qui s'insère dans un diagramme commutatif 
\begin{equation}\label{alpha13c}
\xymatrix{
{\bMod_\mU(A)}\ar[r]\ar[d]_{\alpha_\mU}&{\bMod_\mV(A)}\ar[d]^{\alpha_\mV}\\
{\aMod_\mU(A)}\ar[r]&{\aMod_\mV(A)}}
\end{equation}
Le foncteur \eqref{alpha13b} est pleinement fidèle. En effet, comme le foncteur \eqref{alpha13a} est pleinement fidèle,
il suffit de montrer que pour tout objet $M$ de $\bMod_\mU(A)$, notant $\aI_\mU(M)$ (resp. $\aI_\mV(M)$)
la catégorie des $\alpha$-isomorphismes de $\bMod_\mU(A)$ (resp. $\bMod_\mV(A)$) de but $M$, 
le foncteur d'inclusion $\phi_M\colon \aI^\circ_\mU(M)\rightarrow \aI^\circ_\mV(M)$ est cofinal. 
Soit $f\colon N\rightarrow M$ un objet de $\aI_\mV(M)$.  D'après \ref{alpha3},
pour tout $\varepsilon \in \Lambda^+$, il existe un morphisme $A$-linéaire $g_\varepsilon\colon M\rightarrow N$ 
tel que $f\circ g_\varepsilon=\pi^\varepsilon \id_M$ et $g_\varepsilon\circ f=\pi^\varepsilon \id_N$. 
Comme l'ensemble $\Lambda^+$ est $\mU$-petit, 
$\cup_{\varepsilon\in \Lambda^+}\im(g_\varepsilon)$
est représentable par un sous-objet $N'$ de $N$ appartenant à $\bMod_\mU(A)$. 
Comme l'injection canonique $g\colon N'\rightarrow N$ est clairement un $\alpha$-isomorphisme, 
$f\circ g\colon N'\rightarrow M$ est un objet de $\aI_\mU(M)$, d'où l'assertion recherchée en vertu de (\cite{sga4} I 8.1.3(c)).

\subsection{}\label{alpha11}
Dans la suite de cette section, $\cC$ désigne une $\mU$-catégorie \eqref{notconv3} (\cite{sga4} I 1.1) et
$\hcC$ la catégorie des préfaisceaux de $\mU$-ensembles sur $\cC$ (\cite{sga4} I 1.2).
On note $R_\hcC$ le préfaisceau d'anneaux constant sur $\cC$ de valeur $R$.
On désigne par $\bMod(R_\hcC)$ la catégorie abélienne tensorielle des  $\mU$-préfaisceaux de $R$-modules sur $\cC$,
que l'on identifie à la catégorie des $(R_\hcC)$-modules de $\hcC$ (\cite{sga4} I 3.2). 
Prenant pour $\phi\colon R\rightarrow \End(R_\hcC)=R$ l'homomorphisme identique, 
on appelle catégorie des {\em $\alpha$-$(R_\hcC)$-modules} et l'on note $\aMod(R_\hcC)$ 
le quotient de la catégorie abélienne $\bMod(R_\hcC)$ par la sous-catégorie épaisse des $(R_\hcC)$-modules 
$\alpha$-nuls. Nous utiliserons les conventions de notation de \ref{alpha2} et \ref{alpha33}.

\begin{defi}\label{alpha14}
On appelle catégorie des {\em préfaisceaux de $\alpha$-$R$-modules} sur $\cC$ et l'on note $\ahcC$ 
la catégorie des préfaisceaux sur $\cC$ à valeurs dans la catégorie $\aMod(R)$, {\em i.e.}, 
la catégorie des foncteurs de $\cC^\circ$ à valeurs dans $\aMod(R)$ (\cite{sga4} II 6.0).
\end{defi}

On note $R^\alpha_\hcC$ le préfaisceau constant sur $\cC$ de valeur $R^\alpha$. 
Pour deux préfaisceaux de $\alpha$-$R$-modules $M$ et $N$, 
on définit le produit tensoriel $M\otimes_{R^\alpha_\hcC}N$  par la formule
suivante~: pour tout $X\in \ob(\cC)$, 
\begin{equation}\label{alpha14a}
(M\otimes_{R^\alpha_\hcC}N)(X)=M(X)\otimes_{R^\alpha}N(X), 
\end{equation}
où le terme de droite désigne le produit tensoriel dans la catégorie $\aMod(R)$ \eqref{alpha33a}. 
Munie de ce produit, $\ahcC$ est une catégorie abélienne tensorielle, ayant $R^\alpha_{\hcC}$ pour objet unité.

\begin{prop}\label{alpha15}
Les $\mU$-limites inductives et projectives dans $\ahcC$ sont représentables. Pour tout $X\in\ob(\cC)$, le foncteur 
\begin{equation}\label{alpha15a}
\ahcC\rightarrow \aMod(R),\ \ \ M\mapsto M(X)
\end{equation}
commute aux limites inductives et projectives.
\end{prop}
Cela résulte immédiatement de \ref{alpha10}. 

\begin{cor}\label{alpha150}
La catégorie $\ahcC$ est une catégorie abélienne vérifiant l'axiome {\rm (AB 5)} de {\rm (\cite{tohoku} §~1.5)}. 
\end{cor}

C'est une conséquence de \ref{alpha100} et \ref{alpha15}.

\subsection{}\label{alpha16}
Le foncteur $\alpha\colon \bMod(R)\rightarrow \aMod(R)$ définit un foncteur exact et monoïdal 
\begin{equation}\label{alpha16a}
\halpha\colon \bMod(R_\hcC)\rightarrow \ahcC. 
\end{equation}
Celui-ci transforme les modules $\alpha$-nuls en le préfaisceau de $\alpha$-$R$-modules nul. 
Il induit donc un foncteur exact et monoïdal
\begin{equation}\label{alpha16b}
u\colon \aMod(R_\hcC)\rightarrow \ahcC. 
\end{equation}

Les foncteurs $\sigma_*$ \eqref{alpha8b} et $\sigma_!$ \eqref{alpha8c} (pour les $R$-modules) induisent des foncteurs 
que l'on note respectivement
\begin{eqnarray}
\hsigma_*\colon \ahcC&\rightarrow& \bMod(R_\hcC),\label{alpha16c}\\
\hsigma_!\colon \ahcC&\rightarrow& \bMod(R_\hcC).\label{alpha16d}
\end{eqnarray}
D'après \ref{alpha9}, $\hsigma_*$ (resp. $\hsigma_!$) est un adjoint à droite (resp. à gauche) de $\halpha$; 
les morphismes d'adjonction $\halpha\circ \hsigma_* \rightarrow \id$ et $\id\rightarrow \halpha \circ \hsigma_!$ 
sont des isomorphismes~; et le morphisme d'adjonction $\id\rightarrow \hsigma_*\circ\halpha$
induit un isomorphisme $\halpha\stackrel{\sim}{\rightarrow} \halpha\circ\hsigma_*\circ\halpha$. 

On note $v$ le foncteur composé
\begin{equation}\label{alpha16e}
v=\alpha\circ \hsigma_*\colon  \ahcC\rightarrow \aMod(R_\hcC),
\end{equation}
où $\alpha$ désigne le foncteur de localisation pour les $R$-modules de $\hcC$. 
L'isomorphisme $\halpha\circ \hsigma_*\stackrel{\sim}{\rightarrow} \id$ induit un isomorphisme 
\begin{equation}\label{alpha16f}
u\circ v\stackrel{\sim}{\rightarrow} \id.
\end{equation}
D'après \ref{alpha9}(iii), le morphisme d'adjonction $\id\rightarrow \hsigma_*\circ \halpha$ induit un isomorphisme 
$\alpha\stackrel{\sim}{\rightarrow} v\circ u\circ \alpha$. Compte tenu de (\cite{gz} I 1.2), on en déduit un isomorphisme 
\begin{equation}\label{alpha16g}
\id\stackrel{\sim}{\rightarrow} v\circ u.
\end{equation} 
On vérifie aussitôt que les isomorphismes \eqref{alpha16f} et \eqref{alpha16g} font de $v$ un adjoint à droite de $u$. 
Par suite, $u$ et $v$ sont des équivalences de catégories abéliennes tensorielles, quasi-inverses l'une de l'autre. 
On a des isomorphismes canoniques
\begin{equation}\label{alpha16h}
\halpha\stackrel{\sim}{\rightarrow} u\circ \alpha \ \ \ {\rm et} \ \ \ \alpha\stackrel{\sim}{\rightarrow} v\circ \halpha,
\end{equation}
compatibles aux isomorphismes \eqref{alpha16f} et \eqref{alpha16g}. 
On en déduit aussitôt que le foncteur
\begin{equation}\label{alpha16i}
\hsigma_*\circ u \colon \aMod(R_\hcC)\rightarrow \bMod(R_\hcC)
\end{equation}
est un adjoint à droite du foncteur de localisation $\alpha$, et que le foncteur
\begin{equation}\label{alpha16j}
\hsigma_!\circ u \colon \aMod(R_\hcC)\rightarrow \bMod(R_\hcC)
\end{equation}
est un adjoint à gauche de $\alpha$. 

\subsection{}\label{alpha36}
La donnée d'un monoïde commutatif unitaire de $\bMod(R_\hcC)$ (resp. $\ahcC$) est équivalente à la donnée 
d'un préfaisceau sur $\cC$ à valeurs dans la catégorie $\bAlg(R)$ (resp. $\aAlg(R)$) \eqref{alpha5}. 
On note $\bAlg(R_\hcC)$ (resp. $\bAlg(\ahcC)$) la catégorie des monoïdes commutatifs unitaires de $\bMod(R_\hcC)$ 
(resp. $\ahcC$). Le foncteur $\halpha$ \eqref{alpha16a} étant monoïdal, il induit un foncteur que l'on note encore
\begin{equation}\label{alpha36a}
\halpha\colon \bAlg(R_\hcC)\rightarrow \bAlg(\ahcC).
\end{equation}
Par ailleurs, compte tenu de \ref{alpha5}, le foncteur $\hsigma_*$ \eqref{alpha16c} induit un foncteur que l'on note encore
\begin{equation}\label{alpha36b}
\hsigma_*\colon \bAlg(\ahcC)\rightarrow \bAlg(R_\hcC).
\end{equation}
D'après \ref{alpha6}, $\hsigma_*$ est un adjoint à droite de $\halpha$, 
le morphisme d'adjonction $\halpha\circ \hsigma_* \rightarrow \id$  
est un isomorphisme et le morphisme d'adjonction $\id\rightarrow \hsigma_*\circ\halpha$
induit un isomorphisme $\halpha\stackrel{\sim}{\rightarrow} \halpha\circ\hsigma_*\circ\halpha$.

\subsection{}\label{alpha17}
Lorsque la catégorie $\cC$ est $\mU$-petite, $\hcC$ est une $\mU$-catégorie.  
Mais lorsque $\cC$ est une $\mU$-catégorie, $\hcC_\mU$ n'est pas en général une $\mU$-catégorie (\cite{sga4} I 1.2). 
Soit $\mV$ un univers tel que $\cC\in \mV$ et $\mU\subset \mV$. 
On affectera d'un indice $\mU$ ou $\mV$ les catégories et foncteurs dépendants de l'univers. 
Le foncteur pleinement fidèle canonique $\bMod_\mU(R)\rightarrow \bMod_\mV(R)$ induit un foncteur pleinement fidèle 
\begin{equation}\label{alpha17a}
\bMod(R_{\hcC_\mU})\rightarrow \bMod(R_{\hcC_\mV}).
\end{equation}
De même, le foncteur pleinement fidèle canonique $\aMod_\mU(R)\rightarrow \aMod_\mV(R)$ \eqref{alpha13b} induit un foncteur pleinement fidèle 
\begin{equation}\label{alpha17b}
\ahcC_\mU\rightarrow \ahcC_\mV.
\end{equation}

Tout $(R_{\hcC_\mU})$-module $M$ définit un foncteur
\begin{equation}\label{alpha17c}
M\colon (\hcC_\mU)^\circ\rightarrow \bMod_\mV(R), \ \ \ X\mapsto M(X)=\Hom_{\hcC_\mU}(X,M).
\end{equation}
Pour tout objet $P$ de $\ahcC_\mU$, on considère le foncteur 
\begin{equation}\label{alpha17d}
\hP\colon (\hcC_\mU)^\circ\rightarrow \aMod_\mV(R), \ \ \ X\mapsto \alpha(\hsigma_*(P)(X)),
\end{equation}
où $\hsigma_*$ est le foncteur \eqref{alpha16c}. On obtient ainsi un foncteur
\begin{equation}\label{alpha17e}
\ahcC\rightarrow \bHom((\hcC_\mU)^\circ,\aMod_\mV(R)), \ \ \ P\mapsto \hP.
\end{equation}
D'après \ref{alpha9}(ii), on a un isomorphisme canonique fonctoriel $\hP|\cC^\circ\stackrel{\sim}{\rightarrow} P$. 
Il est alors commode et sans risque d'ambiguïté de noter $\hP$ encore $P$.

D'après (\cite{sga4} I 3.5), pour tout $(R_{\hcC_\mU})$-module $M$ et tout $X\in \ob(\hcC_\mU)$, on a un isomorphisme canonique et fonctoriel en $M$,
\begin{equation}\label{alpha17f}
M(X)\stackrel{\sim}{\rightarrow}\underset{\underset{(Y,u)\in (\cC_{/X})^\circ}{\longleftarrow}}{\lim}\ M(Y).
\end{equation}
Par suite, en vertu de \ref{alpha10}(i), pour tout objet $P$ de $\ahcC_\mU$ et tout $X\in \ob(\hcC_\mU)$, on a un isomorphisme canonique fonctoriel en $P$,
\begin{equation}\label{alpha17g}
P(X)\stackrel{\sim}{\rightarrow}\underset{\underset{(Y,u)\in (\cC_{/X})^\circ}{\longleftarrow}}{\lim}\ P(Y).
\end{equation}

\subsection{}\label{alpha20}
Dans la suite de cette section, on se donne une topologie sur $\cC$ qui en fait un $\mU$-site, {\em i.e.}, telle que $\cC$ admette une
$\mU$-petite famille topologiquement génératrice (\cite{sga4} II 3.0.2). 
On désigne par $\tcC$ le topos des faisceaux de $\mU$-ensembles sur $\cC$.
On note $R_\tcC$ le faisceau d'anneaux constant de valeur $R$ sur $\cC$ et $\fm_\tcC$ le faisceau d'idéaux constant de valeur 
$\fm$ sur $\cC$,  {\em i.e.}, les faisceaux associés au préfaisceaux constants sur $\cC$ de valeurs $R$ et $\fm$, respectivement 
(\cite{sga4} II 6.4). 

Soit $A$ une $(R_\tcC)$-algèbre de $\tcC$, {\em i.e.}, un $\mU$-faisceau de $R$-algèbres sur $\cC$ (\cite{sga4} II 6.3.1). 
On désigne par $\bMod(A)$ la catégorie abélienne tensorielle des $A$-modules de $\tcC$. 
Prenant pour $\phi\colon R\rightarrow \Gamma(\tcC,A)=\End(A)$ l'homomorphisme canonique, 
on appelle catégorie des {\em $\alpha$-$A$-modules} et l'on note $\aMod(A)$ 
le quotient de la catégorie abélienne $\bMod(A)$ par la sous-catégorie épaisse des $A$-modules 
$\alpha$-nuls. Nous utiliserons les conventions de notation de \ref{alpha2} et \ref{alpha33}. 
On observera en particulier que pour tous $A$-modules $M$ et $N$, $\Hom_{\aMod(A)}(M,N)$ est 
naturellement muni d'une structure de $\Gamma(\tcC,A)$-module.

\begin{lem}\label{alpha21}
Soit $A$ une $(R_\tcC)$-algèbre de $\tcC$. 
\begin{itemize}
\item[{\rm (i)}] Pour tout $A$-module $M$, les conditions suivantes sont équivalentes~:
\begin{itemize}
\item[{\rm (a)}] $M$ est $\alpha$-nul en tant qu'objet de $\bMod(A)$;
\item[{\rm (b)}] $M$ est $\alpha$-nul en tant qu'objet de $\bMod(R_\tcC)$;
\item[{\rm (c)}] $M$ est $\alpha$-nul en tant qu'objet de $\bMod(R_\hcC)$; 
\item[{\rm (d)}] pour tout $U\in \ob(\cC)$, $\fm M(U)=0$; 
\item[{\rm (e)}] $\fm_\tcC M=0$.
\end{itemize}
\item[{\rm (ii)}] Pour tout $(R_\hcC)$-module $\alpha$-nul $P$ de $\hcC$, le $(R_\tcC)$-module associé $P^\tta$ 
de $\tcC$ est $\alpha$-nul. 
\item[{\rm (iii)}]  Pour tout morphisme de $A$-modules $f\colon M\rightarrow N$, les conditions suivantes sont équivalentes~:
\begin{itemize}
\item[$(1)$] $f$ est un $\alpha$-isomorphisme en tant que morphisme de $\bMod(A)$~; 
\item[$(2)$] $f$ est un $\alpha$-isomorphisme en tant que morphisme de $\bMod(R_\tcC)$~; 
\item[$(3)$] $f$ est un $\alpha$-isomorphisme en tant que morphisme de $\bMod(R_\hcC)$~;  
\item[$(4)$] le morphisme $\fm_\tcC\otimes_{R_\tcC}M\rightarrow \fm_\tcC\otimes_{R_\tcC}N$ induit par $f$ est un isomorphisme.
\end{itemize} 
\item[{\rm (iv)}] Pour tous $A$-modules $M$ et $N$, on a un isomorphisme $\Gamma(\tcC,A)$-linéaire canonique
\begin{equation}\label{alpha21a}
\Hom_{\aMod(A)}(M^\alpha,N^\alpha)\stackrel{\sim}{\rightarrow}\Hom_{\bMod(A)}(\fm_\tcC\otimes_{R_\tcC}M,N).
\end{equation}
\end{itemize}
\end{lem}

(i) En effet, les conditions (a), (b), (c) et (d) sont clairement équivalentes. Comme $\fm_\tcC M$ est le faisceau associé au préfaisceau 
$U\mapsto \fm M(U)$, (d) implique (e). Il est clair que (e) implique (a). 

(ii) Cela résulte aussitôt de (i)

(iii) Les conditions (1) et (2) sont équivalentes en vertu de (i). 
L'implication $(2)\Rightarrow(3)$ résulte de \ref{alpha3} et l'implication $(3)\Rightarrow(2)$ est une conséquence de (ii).
L'implication $(4)\Rightarrow(2)$ est une conséquence du fait que le morphisme canonique $\fm_\tcC\otimes_{R_\tcC} M\rightarrow M$
est un $\alpha$-isomorphisme \eqref{alpha3}. Par ailleurs, le $(R_\tcC)$-module $\fm_\tcC$ est plat (\cite{sga4} V 1.7.1)
et pour tout $(R_\tcC)$-module $\alpha$-nul $F$, on a $\fm_\tcC\otimes_{R_\tcC}F=0$; d'où l'implication $(2)\Rightarrow(4)$.

(iv) D'après (iii), le morphisme canonique $\fm_\tcC\otimes_{R_\tcC}M\rightarrow M$ est un objet initial de
la catégorie des $\alpha$-isomorphismes de $\bMod(A)$ de but $M$; d'où  l'isomorphisme \eqref{alpha21a}.
On vérifie aussitôt qu'il est $\Gamma(\tcC,A)$-linéaire.

\begin{defi}\label{alpha18}
On dit qu'un préfaisceau de $\alpha$-$R$-modules $F$ sur $\cC$ est {\em séparé} 
(resp. est un {\em faisceau}) 
si  pour tout objet $X$ de $\cC$ et tout crible couvrant $\cR$ de $X$, le morphisme canonique 
\begin{equation}\label{alpha18a}
F(X)\rightarrow F(\cR)
\end{equation}
est un monomorphisme (resp. isomorphisme) (cf. \ref{alpha17}). 
On note $\atcC$ la sous-catégorie pleine de $\ahcC$ formée des faisceaux de $\alpha$-$R$-modules.
\end{defi}

D'après \eqref{alpha17g}, pour qu'un préfaisceau de $\alpha$-$R$-modules 
$F$ sur $\cC$ soit séparé (resp. un faisceau), 
il faut et il suffit que pour tout objet $X$ de $\cC$ et tout crible couvrant $\cR$ de $X$, le morphisme canonique
\begin{equation}\label{alpha18b}
F(X)\rightarrow \underset{\underset{(Y,u)\in (\cC_{/\cR})^\circ}{\longleftarrow}}{\lim}\ F(Y)
\end{equation}
soit un monomorphisme (resp. isomorphisme). 
Les faisceaux de $\alpha$-$R$-modules sont donc les faisceaux sur 
$\cC$ à valeurs dans la catégorie $\aMod(R)$ dans le sens de (\cite{sga4} II 6.1).

\begin{prop}\label{alpha182}
Soient $F$ un préfaisceau de $\alpha$-$R$-modules sur $\cC$,  $X\in \ob(\cC)$, $\cR$ un crible de $X$. 
Pour que le morphisme canonique $F(X)\rightarrow F(\cR)$
soit un monomorphisme (resp. isomorphisme), il faut et il suffit qu'il en soit de même 
du morphisme canonique \eqref{alpha16c}
\begin{equation}\label{alpha182a}
\hsigma_*(F)(X)\rightarrow \hsigma_*(F)(\cR).
\end{equation}
\end{prop} 

En effet, d'après \ref{alpha10}(i) et \eqref{alpha17g}, on a un isomorphisme canonique 
\begin{equation}\label{alpha182b}
\sigma_*(F(\cR))\stackrel{\sim}{\rightarrow} \hsigma_*(F)(\cR).
\end{equation}
L'image du morphisme $u\colon F(X)\rightarrow F(\cR)$ par le foncteur $\sigma_*$ s'identifie alors au morphisme \eqref{alpha182a}.
Si $u$ est un monomorphisme (resp. isomorphisme), il en est de même de $\sigma_*(u)$ en vertu de \ref{alpha10}(i).
Inversement, si $\sigma_*(u)$ est un monomorphisme (resp. isomorphisme), il en est de même de $u$ 
d'après \ref{alpha9}(ii) et \ref{alpha10}(i).

\begin{cor}\label{alpha183}
Pour qu'un préfaisceau de $\alpha$-$R$-modules $F$ sur $\cC$ soit séparé (resp. un faisceau), il faut et il suffit 
que le préfaisceau de $R$-modules $\hsigma_*(F)$ sur $\cC$ \eqref{alpha16c} soit séparé (resp. un faisceau).
\end{cor}

\begin{prop}\label{alpha181}
Pour tout $X\in \ob(\cC)$, soit $\cK(X)$ un ensemble de cribles de $X$. Supposons que les $\cK(X)$ soient stables par
changement de base et qu'ils engendrent la topologie de $\cC$. Alors, pour qu'un préfaisceau de $\alpha$-$R$-modules
$F$ sur $\cC$ soit séparé (resp. un faisceau), il faut et il suffit que pour tout $X\in \ob(\cC)$ et tout $\cR \in \cK(X)$, 
le morphisme canonique
\begin{equation}\label{alpha181a}
F(X)\rightarrow F(\cR)
\end{equation}
soit un monomorphisme (resp. isomorphisme). 
\end{prop} 
Cela résulte \ref{alpha182} et (\cite{sga4} II 2.3).

\begin{cor}\label{alpha19}
Si la topologie de $\cC$ est définie par une prétopologie, pour qu'un 
préfaisceau de $\alpha$-$R$-modules $F$ sur $\cC$ soit un faisceau, il faut et il suffit que pour tout objet $X$ de $\cC$
et tout recouvrement $(X_i\rightarrow X)_{i\in I}$, la suite de $\alpha$-$R$-modules
\begin{equation}\label{alpha19a}
0\rightarrow F(X)\rightarrow \prod_{i\in I}F(X_i)\rightarrow \prod_{(i,j)\in I^2}F(X_i\times_XX_j),
\end{equation}
où la dernière flèche est la différence des morphismes induits par les projections de $X_i\times_XX_j$ sur les deux facteurs, 
soit exacte. 
\end{cor}

Cela résulte de \ref{alpha181}, \eqref{alpha17g} et (\cite{sga4} I 2.12).

\subsection{}\label{alpha24}
Rappelons la définition du foncteur $L$ sur la catégorie $\bMod(R_\hcC)$ (\cite{sga4} II 3.0.5). 
Soient $\mV$ un univers tel que $\cC\in \mV$ et $\mU\subset \mV$, $G$ une $\mU$-petite famille topologiquement génératrice de $\cC$ 
(\cite{sga4} II 3.0.1). Pour tout objet $X$ de $\cC$, on désigne par $J(X)$ l'ensemble des cribles couvrants de $X$ et 
par $J_G(X)$ l'ensemble des cribles couvrants de $X$ engendrés par une famille $(X_i\rightarrow X)_{i\in I}$
telle que $X_i\in G$ pour tout $i\in I$. L'ensemble $J(X)$ est $\mV$-petit, et ordonné par l'inclusion, il est cofiltrant. 
Pour tout $\mU$-préfaisceau de $R$-modules $F$ sur $\cC$, 
\begin{equation}\label{alpha24a}
\underset{\underset{\cR\in J(X)^\circ}{\longrightarrow}}{\lim}\ F(\cR)
\end{equation}
est représentable par un $R$-module de $\mV$ \eqref{alpha17}. D'après (\cite{sga4} II 3.0.4), 
pour tout $\cR\in J_G(X)$, $F(\cR)$ est $\mU$-petit, et comme $J_G(X)$ est un ensemble $\mU$-petit cofinal dans $J(X)$ ({\em loc. cit.}), 
il résulte de (\cite{sga4} I 2.3.3) que la limite inductive \eqref{alpha24a} est représentable par un $R$-module $\mU$-petit. 
Choisissons pour tout $F$ et tout $X$ un $R$-module appartenant à $\mU$ qui représente cette limite inductive et posons 
\begin{equation}\label{alpha24b}
LF(X)=\underset{\underset{\cR\in J(X)^\circ}{\longrightarrow}}{\lim}\ F(\cR).
\end{equation}
Pour tout morphisme $f\colon Y\rightarrow X$ de $\cC$, le foncteur de changement de base $f^*\colon J(X)\rightarrow J(Y)$ définit un morphisme
$LF(f)\colon LF(X)\rightarrow LF(Y)$ faisant de $X\mapsto LF(X)$ un $\mU$-préfaisceau de $R$-modules sur $\cC$. 
Pour tout $X\in \ob(\cC)$, le morphisme identique de $X$ étant un objet de $J(X)$, 
on a une application canonique $\ell(F)(X)\colon F(X)\rightarrow LF(X)$. 
On définit ainsi un morphisme $\ell(F)\colon F\rightarrow LF$ de $\bMod(R_\hcC)$.
La correspondance $F\mapsto LF$ est clairement fonctorielle en $F$ et les $\ell(F)$ définissent un morphisme de foncteurs 
$\ell\colon \id\rightarrow L$.

Avec les notations de \ref{alpha16}, on désigne par $\cL$ le foncteur composé 
\begin{equation}\label{alpha24c}
\cL=\halpha\circ L\circ \hsigma_*\colon \ahcC\rightarrow \ahcC.
\end{equation}
Le morphisme $\ell$ et l'isomorphisme $\halpha\circ \hsigma_*\stackrel{\sim}{\rightarrow} \id$ induisent un morphisme de foncteurs 
$\lambda\colon \id\rightarrow \cL$. 
D'après \ref{alpha10}(i) et compte tenu de la définition \eqref{alpha17d}, 
pour tout objet $P$ de $\ahcC$ et tout $X\in \ob(\cC)$, on a un isomorphisme canonique fonctoriel 
\begin{equation}\label{alpha24d}
\cL P(X)\stackrel{\sim}{\rightarrow}\underset{\underset{\cR\in J(X)^\circ}{\longrightarrow}}{\lim}\ P(\cR).
\end{equation}
Le composé de cet isomorphisme et du morphisme $\lambda(P)(X)\colon P(X)\rightarrow \cL P(X)$ n'est autre que 
le morphisme induit par l'objet $\id_X$ de $J(X)$.

\begin{prop}\label{alpha26}
{\rm (i)}\  Le foncteur $\cL$ est exact à gauche. 

{\rm (ii)}\ Le foncteur $L$ transforme les $\alpha$-isomorphismes en des $\alpha$-isomorphismes. 

{\rm (iii)}\ Le diagramme 
\begin{equation}\label{alpha26a}
\xymatrix{
{\bMod(R_\hcC)}\ar[r]^L\ar[d]_\halpha&{\bMod(R_\hcC)}\ar[d]^\halpha\\
{\ahcC}\ar[r]^\cL&{\ahcC}}
\end{equation}
est commutatif à un isomorphisme canonique près
\begin{equation}\label{alpha26b}
\halpha\circ L\stackrel{\sim}{\rightarrow} \cL\circ \halpha,
\end{equation} 
induit par le morphisme d'adjonction $\id\rightarrow \hsigma_*\circ \halpha$.

{\rm (iv)}\ Pour tout préfaisceau de $\alpha$-$R$-modules $P$, le préfaisceau $\cL P$ est séparé. 

{\rm (v)}\ Pour qu'un préfaisceau de $\alpha$-$R$-modules $P$ soit séparé, 
il faut et il suffit que $\lambda(P)\colon P\rightarrow \cL P$ soit un monomorphisme.
Le préfaisceau $\cL P$ est alors un faisceau. 

{\rm (vi)}\ Pour qu'un préfaisceau de $\alpha$-$R$-modules $P$ soit un faisceau, 
il faut et il suffit que $\lambda(P)\colon P\rightarrow \cL P$ soit un isomorphisme. 
\end{prop}

(i) En effet, les foncteurs $\halpha$, $L$ et $\hsigma_*$ sont exacts à gauche en vertu de \ref{alpha10}, \ref{alpha15} et (\cite{sga4} II 3.2(i)). 

(ii) En vertu de \ref{alpha10}(i), pour tout $(R_\hcC)$-module $F$ et tout $X\in \ob(\cC)$, on a un isomorphisme canonique
\begin{equation}\label{alpha26c}
\halpha(L F)(X)\stackrel{\sim}{\rightarrow}\underset{\underset{\cR\in J(X)^\circ}{\longrightarrow}}{\lim}\ \alpha(F(\cR)).
\end{equation}
Compte tenu de \eqref{alpha17f}, pour tout crible $\cR$ de $X$, on a un isomorphisme canonique 
\begin{equation}\label{alpha26d}
\alpha(F(\cR))\stackrel{\sim}{\rightarrow}\underset{\underset{(Y,u)\in (\cC_{/\cR})^\circ}{\longleftarrow}}{\lim}\ \alpha(F(Y)).
\end{equation}
La proposition s'ensuit.

(iii) Cela résulte de (ii) et \ref{alpha9}(iii). 

(iv) D'après \ref{alpha10}(i), $\halpha$ transforme les préfaisceaux séparés de $R$-modules en des 
préfaisceaux séparés de $\alpha$-$R$-modules. 
La proposition résulte donc de (\cite{sga4} II 3.2(i)).

(v) Si $P$ est un préfaisceau séparé de $\alpha$-$R$-modules, 
$\lambda(P)$ est un monomorphisme en vertu de \eqref{alpha24d}
car une limite inductive filtrante de monomorphismes est un monomorphisme \eqref{alpha150}. 
Inversement, si $\lambda(P)$ est un monomorphisme, $P$ est un sous-préfaisceau 
d'un préfaisceau séparé de $\alpha$-$R$-modules~; il est donc séparé.  
Dans ce cas, $\hsigma_*(P)$ est un préfaisceau séparé de $R$-modules d'après \ref{alpha10}(i); donc
 $L\circ \hsigma_*(P)$ est un faisceau de $R$-modules en vertu de (\cite{sga4} II 3.2(iii)) et par suite 
$\cL P$ est un faisceau de $\alpha$-$R$-modules d'après \ref{alpha10}(i).

(vi) En effet, la condition est nécessaire compte tenu de \eqref{alpha24d} et elle est suffisante en vertu de (v).

\begin{prop}\label{alpha25}
Le foncteur d'inclusion $\iota\colon \atcC\rightarrow \ahcC$ admet un adjoint à gauche
\begin{equation}\label{alpha25a}
\otta\colon \ahcC\rightarrow \atcC
\end{equation}
tel que le composé $\iota\circ \otta$ soit canoniquement isomorphe au foncteur $\cL\circ \cL$ \eqref{alpha24c}.
Pour tout préfaisceau de $\alpha$-$R$-modules $P$, le morphisme d'adjonction $P\rightarrow \iota\circ \otta (P)$ se déduit par l'isomorphisme 
précédent du morphisme $\lambda(\cL P)\circ \lambda(P)\colon P\rightarrow \cL\circ\cL(P)$. 
\end{prop}

En effet, d'après \ref{alpha26}, il existe un foncteur $\otta\colon \ahcC\rightarrow \atcC$ tel que $\iota\circ \otta=\cL\circ \cL$. 
On a un morphisme de foncteurs 
\begin{equation}\label{alpha25b}
\id\rightarrow \iota\circ \otta,
\end{equation}
défini pour tout préfaisceau de $\alpha$-$R$-modules $P$, par le morphisme 
$\lambda(\cL P)\circ \lambda(P)\colon P\rightarrow \cL\circ\cL(P)$. 
Si $P$ est un faisceau de $\alpha$-$R$-modules, 
$\lambda(\cL P)\circ \lambda(P)$ est un isomorphisme. On en déduit un isomorphisme 
$\iota\circ \otta\circ \iota\stackrel{\sim}{\rightarrow} \iota$ et par suite un isomorphisme  
\begin{equation}\label{alpha25c}
\otta\circ \iota\stackrel{\sim}{\rightarrow} \id.
\end{equation}
On vérifie aussitôt que les morphismes \eqref{alpha25b} et \eqref{alpha25c} font de $\otta$ un adjoint à gauche de $\iota$; d'où la proposition.

\begin{defi}\label{alpha29}
Pour tout préfaisceau de $\alpha$-$R$-modules $F$, 
on appelle $\otta(F)$ le {\em faisceau associé} à $F$ \eqref{alpha25a}. 
\end{defi}

\begin{prop}\label{alpha28}
{\rm (i)} Le foncteur $\otta$ \eqref{alpha25a} commute aux limites inductives et est exact.

{\rm (ii)} Les $\mU$-limites inductives dans $\atcC$ sont représentables. Pour toute catégorie $\mU$-petite $I$ et 
tout foncteur $\varphi\colon I\rightarrow \atcC$, le morphisme canonique
\begin{equation}\label{alpha28a}
\underset{\underset{I}{\rightarrow}}{\lim}\ \varphi \rightarrow 
\otta(\underset{\underset{I}{\rightarrow}}{\lim}\ \iota\circ \varphi)
\end{equation}
est un isomorphisme.

{\rm (iii)} Les $\mU$-limites projectives dans $\atcC$ sont représentables. Pour tout objet $X$ de $\cC$, le foncteur 
$F\mapsto F(X)$ commute aux limites projectives~; {\em i.e.}, le foncteur d'inclusion $\iota\colon \atcC\rightarrow \ahcC$ commute aux limites projectives. 
\end{prop}

En effet, $\otta$ commute aux limites inductives et $\iota$ commute aux limites projectives d'après \ref{alpha25}.
Soit $\mV$ un univers tel que $\cC\in \mV$ et $\mU\subset \mV$.
Comme les limites projectives commutent aux limites projectives, il résulte de \ref{alpha15} et \eqref{alpha17g}
que pour tout $X\in \ob(\cC)$ et tout crible $\cR$ de $X$, le foncteur \eqref{alpha17e}
\begin{equation}
\ahcC_\mU\rightarrow \aMod_\mV(R), \ \ \ F\mapsto F(\cR)
\end{equation} 
commute aux limites projectives.
Les propositions (ii) et (iii) s'ensuivent compte tenu de \ref{alpha15}.  
Comme $\iota\circ \otta\simeq \cL\circ \cL$ d'après \ref{alpha25} 
et que $\cL$ est exact à gauche en vertu de \ref{alpha26}(i), $\otta$ est exact à gauche et donc exact~; d'où la proposition (i).

\begin{prop}\label{alpha30}
La catégorie $\atcC$ est une catégorie abélienne vérifiant l'axiome {\rm (AB 5)} de {\rm (\cite{tohoku} §~1.5)}. 
\end{prop}

Il résulte de \ref{alpha150} et \ref{alpha28} que $\atcC$ est une catégorie additive où les noyaux et les conoyaux sont représentables. 
Plus précisément, soit $u\colon F\rightarrow G$ un morphisme de $\atcC$. D'après \ref{alpha28}, on a des isomorphismes canoniques
\begin{eqnarray}
\iota(\ker(u))&\stackrel{\sim}{\rightarrow}&\ker(\iota(u)),\label{alpha30a}\\
\coker(u)&\stackrel{\sim}{\rightarrow}&\otta(\coker(\iota(u))).\label{alpha30b}
\end{eqnarray}
On en déduit un isomorphisme canonique
\begin{equation}\label{alpha30c}
\coim(u)\stackrel{\sim}{\rightarrow}\otta(\coim(\iota(u))).
\end{equation}
Montrons que le morphisme canonique $\coim(u)\rightarrow \im(u)$ est un isomorphisme. 
D'après \ref{alpha28}(iii), il suffit de montrer que la suite 
\begin{equation}
0\rightarrow \iota(\coim(u))\rightarrow \iota(G)\rightarrow \iota(\coker(u))
\end{equation}
est exacte. Compte tenu de \eqref{alpha30b} et \eqref{alpha30c} et comme le morphisme d'adjonction $\id\rightarrow \otta\circ \iota$ est un isomorphisme, cette suite est isomorphe à l'image par le foncteur $\iota\circ \otta$ de la suite 
\begin{equation}
0\rightarrow \coim(\iota(u))\rightarrow  \iota(G)\rightarrow \coker(\iota(u)).
\end{equation}
Or, cette suite est exacte et le foncteur $\iota\circ \otta$ est exact à gauche d'après \ref{alpha28}. 
Par suite, $\atcC$ est une catégorie abélienne. Comme $\ahcC$ est une catégorie abélienne vérifiant (AB~5) \eqref{alpha150},
il en est de même de $\atcC$ en vertu de \ref{alpha28}.

\subsection{}\label{alpha22}
On identifie la catégorie $\bMod(R_\tcC)$ à la catégorie des $\mU$-faisceaux de $R$-modules sur $\cC$ 
(\cite{sga4} II 6.3.1). 
D'après \ref{alpha10}(i), le foncteur $\halpha$ \eqref{alpha16a} transforme les faisceaux de $R$-modules en des 
faisceaux de $\alpha$-$R$-modules.
Il définit donc un foncteur 
\begin{equation}\label{alpha22a}
\talpha\colon \bMod(R_\tcC)\rightarrow \atcC,
\end{equation} 
qui s'insère dans un diagramme commutatif 
\begin{equation}\label{alpha22g}
\xymatrix{
{\bMod(R_\tcC)}\ar[r]^-(0.5)\talpha\ar[d]_i&{\atcC}\ar[d]^\iota\\
{\bMod(R_\hcC)}\ar[r]^-(0.5)\halpha&{\ahcC}}
\end{equation}
où $i$ et $\iota$ sont les foncteurs d'injection canoniques. 
D'après \ref{alpha21}(iii), $\talpha$ transforme les $\alpha$-isomor\-phismes en des isomorphismes. Il induit donc un foncteur
\begin{equation}\label{alpha22b}
\mu\colon \aMod(R_\tcC)\rightarrow \atcC.
\end{equation} 
D'autre part, le foncteur $\hsigma_*$ \eqref{alpha16c} transforme les faisceaux de 
$\alpha$-$R$-modules en des faisceaux de $R$-modules d'après \ref{alpha182}. Il induit donc un foncteur 
\begin{equation}\label{alpha22c}
\tsigma_*\colon \atcC\rightarrow \bMod(R_\tcC).
\end{equation}
On note $\nu$ le foncteur composé
\begin{equation}\label{alpha22d}
\nu=\alpha\circ \tsigma_*\colon \atcC\rightarrow \aMod(R_\tcC),
\end{equation} 
où $\alpha$ est le foncteur de localisation pour les $(R_\tcC)$-modules \eqref{alpha20}.
D'après \ref{alpha16}, le foncteur $\tsigma_*$ est un adjoint à droite de $\talpha$; 
le morphisme d'adjonction $\talpha\circ \tsigma_* \rightarrow \id$ 
est un isomorphisme~; et le morphisme d'adjonction $\id\rightarrow \tsigma_*\circ\talpha$
induit un isomorphisme $\talpha\stackrel{\sim}{\rightarrow} \talpha\circ\tsigma_*\circ\talpha$. 

L'isomorphisme $\talpha\circ \tsigma_*\stackrel{\sim}{\rightarrow} \id$ induit un isomorphisme 
\begin{equation}\label{alpha22e}
\mu\circ \nu\stackrel{\sim}{\rightarrow} \id.
\end{equation}
D'après \ref{alpha21}(iii), le morphisme d'adjonction $\id\rightarrow \tsigma_*\circ \talpha$ induit un isomorphisme 
$\alpha\stackrel{\sim}{\rightarrow} \nu\circ \mu\circ \alpha$. On en déduit un isomorphisme 
\begin{equation}\label{alpha22f}
\id\stackrel{\sim}{\rightarrow} \nu\circ \mu.
\end{equation} 
On vérifie aussitôt que les isomorphismes \eqref{alpha22e} et \eqref{alpha22f} font de $\nu$ un adjoint à droite de $\mu$. 
Par suite, $\mu$ et $\nu$ sont des équivalences de catégories, quasi-inverses l'une de l'autre. 
Ce sont donc des foncteurs exacts (\cite{gabriel} I §~1 prop.~13). Il s'ensuit que le foncteur $\talpha$ est aussi exact.

\begin{prop}\label{alpha27}
Le diagramme 
\begin{equation}\label{alpha27a}
\xymatrix{
{\bMod(R_\hcC)}\ar[r]^-(0.5){\tta}\ar[d]_{\halpha}&{\bMod(R_\tcC)}\ar[d]^{\talpha}\\
{\ahcC}\ar[r]^-(0.5)\otta&{\atcC}}
\end{equation}
où $\otta$ est le foncteur \eqref{alpha25a} et $\tta$ est le foncteur ``faisceau de $R$-modules associé'' {\rm (\cite{sga4} II 6.4)},
est commutatif à isomorphisme canonique près. 
\end{prop}

En effet, notant $i\colon \bMod(R_\tcC)\rightarrow \bMod(R_\hcC)$ le foncteur d'injection canonique, on a un isomorphisme canonique
$i\circ \tta\stackrel{\sim}{\rightarrow}L\circ L$ (\cite{sga4} II 3.4 et 6.4). 
Par suite, en vertu de \ref{alpha26}(iii), on a un isomorphisme canonique
\begin{equation}
\halpha \circ i\circ \tta\stackrel{\sim}{\rightarrow}\iota\circ \otta\circ \halpha.
\end{equation}
La proposition s'ensuit puisqu'on a $\halpha\circ i=\iota\circ \talpha$ \eqref{alpha22a}.

\begin{prop}\label{alpha32}
La catégorie abélienne $\atcC$ admet une famille de générateurs indexée par un ensemble appartenant à $\mU$. 
\end{prop}

Cela résulte de \ref{alpha22}, (\cite{sga4} II 6.7) et (\cite{gabriel} III §~2 lem.~4). 

\subsection{}\label{alpha34}
Pour tout $R^\alpha$-module $P$, on appelle {\em faisceau de $\alpha$-$R$-modules constant de valeur $P$ sur $\cC$} 
et l'on note $P_\tcC$ le faisceau associé au 
préfaisceau de $\alpha$-$R$-modules constant de valeur $P$ sur $\cC$ \eqref{alpha14}. D'après \ref{alpha27}, 
pour tout $R$-module $M$,
notant $M_\tcC$ le faisceau constant de valeur $M$ sur $\cC$, on a un isomorphisme canonique 
\begin{equation}\label{alpha34a}
M^\alpha_{\tcC}\stackrel{\sim}{\rightarrow}\talpha(M_\tcC).
\end{equation}

On munit $\atcC$ de la structure de catégorie abélienne tensorielle déduite de celle de $\aMod(R_\tcC)$ \eqref{alpha33a} 
via l'équivalence de catégories $\mu$ \eqref{alpha22b}; on note $\otimes_{R^\alpha_\tcC}$ le produit tensoriel dans $\atcC$. 
Le foncteur $\talpha$ est donc monoïdal~: pour tous $R_\tcC$-modules $M$ et $N$, on a un isomorphisme canonique 
\begin{equation}\label{alpha34b}
\talpha(M\otimes_{R_\tcC}N)\stackrel{\sim}{\rightarrow}M^\alpha\otimes_{R^\alpha_\tcC}N^\alpha.
\end{equation}
D'après \eqref{alpha34a}, $R^\alpha_\tcC$ est un objet unité de $\atcC$.

\begin{prop}\label{alpha35}
Pour tous préfaisceaux de $\alpha$-$R$-modules $M$ et $N$ sur $\cC$, on a un isomorphisme canonique fonctoriel 
\begin{equation}\label{alpha35a}
\otta(M) \otimes_{R^\alpha_\tcC}\otta(N)\stackrel{\sim}{\rightarrow}\otta(M\otimes_{R^\alpha_\hcC}N),
\end{equation}
où $M\otimes_{R^\alpha_\hcC}N$ désigne le produit tensoriel dans $\ahcC$ \eqref{alpha14a}. 
\end{prop}

Cela résulte aussitôt de \ref{alpha27} et (\cite{sga4} IV 12.10) puisque le foncteur $u$ \eqref{alpha16b} est une équivalence 
de catégories tensorielles.

\subsection{}\label{alpha37}
Considérons les foncteurs adjoints \eqref{alpha25}
\begin{equation}\label{alpha37c}
\xymatrix{{\atcC}\ar@<1ex>[r]^\iota&{\ahcC}\ar@<1ex>[l]^\otta}
\end{equation}
D'après \ref{alpha35}, pour tout objet $A$ de $\atcC$, on a un isomorphisme canonique
\begin{equation}\label{alpha37d}
\otta(\iota(A)\otimes_{R^\alpha_\hcC}\iota(A))\stackrel{\sim}{\rightarrow} A\otimes_{R^\alpha_\tcC}A.
\end{equation}
On en déduit par adjonction que la donnée d'une structure de monoïde commutatif unitaire de $\atcC$ sur $A$ 
est équivalente à la donnée sur $\iota(A)$ d'une structure de monoïde commutatif unitaire de $\ahcC$.
Par suite, la donnée d'un monoïde commutatif unitaire de $\atcC$ est équivalente à la donnée d'un 
préfaisceau sur $\cC$ à valeur dans $\aAlg(R)$ \eqref{alpha5} dont le préfaisceau de $\alpha$-$R$-modules  
sous-jacent est un faisceau \eqref{alpha36}. 
On note $\bAlg(\atcC)$ la catégorie des monoïdes commutatifs unitaires de $\atcC$. 
  
De même, la donnée d'un monoïde commutatif unitaire de $\bMod(R_\tcC)$ est équivalente à la donnée 
d'une $(R_\tcC)$-algèbre de $\tcC$ (\cite{sga4} II 6.3.1 et IV 12.10). 
On note $\bAlg(R_\tcC)$ la catégorie des monoïdes commutatifs unitaires de $\bMod(R_\tcC)$.

Le foncteur $\talpha$ \eqref{alpha22a} étant monoïdal, il induit un foncteur que l'on note encore
\begin{equation}\label{alpha37e}
\talpha\colon \bAlg(R_\tcC)\rightarrow \bAlg(\atcC).
\end{equation}
Compte tenu de \ref{alpha5}, le foncteur $\tsigma_*$ \eqref{alpha22c} induit un foncteur que l'on note encore
\begin{equation}\label{alpha37f}
\tsigma_*\colon \bAlg(\atcC)\rightarrow \bAlg(R_\tcC).
\end{equation}
D'après \ref{alpha6}, $\tsigma_*$ est un adjoint à droite de $\talpha$, 
le morphisme d'adjonction $\talpha\circ \tsigma_* \rightarrow \id$  
est un isomorphisme et le morphisme d'adjonction $\id\rightarrow \tsigma_*\circ\talpha$
induit un isomorphisme $\talpha\stackrel{\sim}{\rightarrow} \talpha\circ\tsigma_*\circ\talpha$. 

\subsection{}\label{alpha38}
Dans ce numéro, si $D$ est une $(R_\tcC)$-algèbre de $\tcC$ (resp. une $R$-algèbre), nous affecterons d'un indice $D$ le foncteur 
de localisation $\alpha$ \eqref{alpha2a} pour les $D$-modules.

Soit $A$ une $(R_\tcC)$-algèbre de $\tcC$. On pose $B=\talpha(A)$ \eqref{alpha37e}
et on désigne par $\bMod(B)$ la catégorie des $B$-modules unitaires de $\atcC$. 
Il résulte aussitôt de \ref{alpha35} 
que la donnée d'une structure de $B$-module unitaire sur un faisceau de $\alpha$-$R$-modules 
$M$ est équivalente à la donnée
pour tout $X\in \ob(\cC)$ d'une structure de $\alpha_R(A(X))$-module unitaire sur $M(X)$ dans le sens de \ref{alpha39} telle que 
pour tout morphisme $X\rightarrow Y$ de $\cC$, le morphisme $M(Y)\rightarrow M(X)$ 
soit linéaire relativement au morphisme de $R^\alpha$-algèbres $\alpha_R(A(Y))\rightarrow \alpha_R(A(X))$. 

Le foncteur $\talpha$ \eqref{alpha22a} étant monoïdal, il définit un foncteur 
\begin{equation}\label{alpha38a}
\tbeta\colon \bMod(A)\rightarrow \bMod(B).
\end{equation}
Celui-ci transforme les $\alpha$-isomorphismes en des isomorphismes d'après \ref{alpha21}(iii). Il induit donc un foncteur \eqref{alpha20}
\begin{equation}\label{alpha38b}
b\colon \aMod(A)\rightarrow \bMod(B).
\end{equation}
On pose $A'=\tsigma_*(B)$ \eqref{alpha37f}. D'après \ref{alpha37}, 
on a un homomorphisme canonique de $(R_\tcC)$-algèbres $\lambda\colon 
A\rightarrow A'$, qui induit un isomorphisme $\talpha(A)\stackrel{\sim}{\rightarrow}\talpha(A')$. 
Compte tenu de \ref{alpha39}, le foncteur $\tsigma_*$ \eqref{alpha22c} induit un foncteur 
\begin{equation}\label{alpha38c}
\ttau'_*\colon \bMod(B)\rightarrow \bMod(A').
\end{equation}
Composant avec le foncteur induit par $\lambda$, on obtient un foncteur 
\begin{equation}\label{alpha38d}
\ttau_*\colon \bMod(B)\rightarrow \bMod(A).
\end{equation}
On pose 
\begin{equation}\label{alpha38e}
t=\alpha_A\circ \ttau_*\colon \bMod(B)\rightarrow \aMod(A).
\end{equation}

D'après \ref{alpha39} et \ref{alpha16}, le foncteur $\ttau_*$ est un adjoint à droite de $\tbeta$, 
le morphisme d'adjonction $\tbeta\circ \ttau_* \rightarrow \id$ 
est un isomorphisme et le morphisme d'adjonction $\id\rightarrow \ttau_*\circ\tbeta$
induit un isomorphisme $\tbeta\stackrel{\sim}{\rightarrow} \tbeta\circ\ttau_*\circ\tbeta$. 

L'isomorphisme $\tbeta\circ \ttau_*\stackrel{\sim}{\rightarrow} \id$ induit un isomorphisme 
\begin{equation}\label{alpha38f}
b\circ t\stackrel{\sim}{\rightarrow} \id.
\end{equation}
D'après \ref{alpha21}(iii), le morphisme d'adjonction $\id\rightarrow \ttau_*\circ \tbeta$ induit un isomorphisme 
$\alpha_A\stackrel{\sim}{\rightarrow} t\circ b\circ \alpha_A$. On en déduit un isomorphisme 
\begin{equation}\label{alpha38g}
\id\stackrel{\sim}{\rightarrow} t\circ b.
\end{equation} 
On vérifie aussitôt que les isomorphismes \eqref{alpha38f} et \eqref{alpha38g} font de $t$ un adjoint à droite de $b$. 
Par suite, $b$ et $t$ sont des équivalences de catégories, quasi-inverses l'une de l'autre.

\section{\texorpdfstring{Conditions de $\alpha$-finitude}{Conditions de alpha-finitude}}\label{finita}

\subsection{}\label{finita1}
Les hypothèses et notations de \ref{alpha1} sont en vigueur dans cette section. 
On se donne, de plus, un $\mU$-topos $\cT$ \eqref{notconv3} annelé par une $R$-algèbre $A$ (\cite{sga4} IV 11.1.1). 
Nous considérons toujours $\cT$ comme muni de sa topologie canonique (\cite{sga4} II 2.5), qui en fait un $\mU$-site.
Pour tout objet $X$ de $\cT$, le topos $\cT_{/X}$ sera annelé par l'anneau $A|X$ \eqref{notconv6}. 

\subsection{}\label{finita3}
Soit $\gamma\in R$. On dit qu'un morphisme de $A$-modules est un {\em $\gamma$-isomorphisme} 
si son noyau et son conoyau sont annulés par $\gamma$. 

On dit qu'un complexe de $A$-modules à différentielle de degré $1$, le degré étant écrit en exposant,
est {\em $\gamma$-acyclique} si ses groupes  de cohomologie sont annulés par $\gamma$, 
et qu'il est {\em $\alpha$-acyclique} s'il est $\delta$-acyclique pour tout $\delta\in \fm$.

On appelle {\em suite de $A$-modules} un complexe de longueur finie de $A$-modules, 
à différentielle de degré $1$, le degré étant écrit en exposant.
On dit qu'une suite de $A$-modules est {\em $\gamma$-exacte} (resp. {\em $\alpha$-exacte}) 
si elle est $\gamma$-acyclique (resp. $\alpha$-acyclique) en tant que complexe de $A$-modules.

\begin{defi}\label{finita9}
Soient $F$ un $A$-module, $n$ un entier $\geq 0$, $\gamma\in R$.
\begin{itemize}
\item[(i)] On dit que $F$ est {\em de $n$-présentation $\gamma$-finie} si 
la sous-catégorie pleine de $\cT$ formée des objets $X$ tels qu'il existe une suite $\gamma$-exacte \eqref{finita3}
\begin{equation}
E^{-n}\rightarrow E^{-n+1}\rightarrow \dots\rightarrow E^0\rightarrow F|X\rightarrow 0,
\end{equation} 
avec $E^i$ un $(A|X)$-module libre de type fini pour tout $-n\leq i\leq 0$, est un raffinement de l'objet final de $\cT$. 
On dit que $F$ est de {\em type $\gamma$-fini} (resp. {\em présentation $\gamma$-finie}) 
s'il est de $0$-présentation $\gamma$-finie (resp. $1$-présentation $\gamma$-finie).

\item[(ii)] On dit que $F$ est de {\em $n$-présentation $\alpha$-finie} s'il est de $n$-présentation 
$\gamma$-finie pour tout $\gamma\in \fm$. 
On dit que $F$ est de {\em type $\alpha$-fini} (resp. {\em présentation $\alpha$-finie}) 
s'il est de $0$-présentation $\alpha$-finie (resp. $1$-présentation $\alpha$-finie).

\item[(iii)] On dit que $F$ est {\em $\alpha$-cohérent} s'il est de type $\alpha$-fini et si 
pour tout objet $X$ de $\cT$ et tout $(A|X)$-morphisme $u\colon E\rightarrow F|X$, où $E$ est un $(A|X)$-module libre de type fini, 
$\ker(u)$ est un $(A|X)$-module de type $\alpha$-fini. 

\item[(iv)] On dit que la $R$-algèbre $A$ est {\em $\alpha$-cohérente} si le $A$-module sous-jacent à $A$ est $\alpha$-cohérent. 
\end{itemize}
\end{defi}

Si $\gamma$ est une unité de $R$, la notion de $n$-présentation $\gamma$-finie correspond 
à la notion standard de $n$-présentation finie introduite dans (\cite{sga6} I 2.8). 
La notion de $\alpha$-cohérence est modelée sur la notion standard de cohérence (\cite{sga6} I 3.1). 

\begin{rema}\label{finita2}
(i)\ Supposons que $\cT$ soit le topos ponctuel, annelé par un anneau standard $A$. 
Pour qu'un $A$-module $F$ soit de type $\alpha$-fini (resp. présentation $\alpha$-finie), il faut et il suffit que pour tout $\gamma\in \fm$, 
il existe un $\gamma$-isomorphisme $f\colon G\rightarrow F$ avec $G$ un $A$-module de type fini (resp. de présentation finie) 
\eqref{finita3}.

(ii)\ Nous donnerons dans \ref{afini5} un cas intéressant pour lequel on dispose d'une caractérisation des modules 
de type $\alpha$-fini (resp. présentation $\alpha$-finie) similaire à celle pour le topos ponctuel (i). 
\end{rema}

\begin{lem}\label{finita15}
Soient $\gamma,\gamma'\in R$, $E\rightarrow F$ un $\gamma$-isomorphisme de $A$-modules, $n$ un entier $\geq 0$. 
\begin{itemize}
\item[{\rm (i)}] Si $E$ est de $n$-présentation $\gamma'$-finie, $F$ est de $n$-présentation $\gamma\gamma'$-finie.
\item[{\rm (ii)}] Si $F$ est de $n$-présentation $\gamma'$-finie, $E$ est de $n$-présentation $\gamma^2\gamma'$-finie.
\end{itemize}
\end{lem}

L'assertion (i) est immédiate et l'assertion (ii) résulte aussitôt de \ref{alpha3}.

\begin{lem}\label{finita16}
Soient $f\colon E\rightarrow F$ un morphisme de $A$-modules qui est un $\alpha$-isomorphisme, $n$ un entier $\geq 0$. 
Pour que $E$ soit de $n$-présentation $\alpha$-finie (resp. $\alpha$-cohérent), il faut et il suffit qu'il en soit de même de $F$.
\end{lem}

L'assertion relative à la $n$-présentation $\alpha$-finie résulte aussitôt de \ref{finita15}. Supposons que $F$ soit $\alpha$-cohérent et montrons 
que $E$ est $\alpha$-cohérent. On sait déjà que $E$ est de type $\alpha$-fini. Soient $X$ un objet de $\cT$, $E'$ un $(A|X)$-module libre de type fini,
$u\colon E'\rightarrow E|X$ un $(A|X)$-morphisme. Le morphisme canonique $\ker(u)\rightarrow \ker(f\circ u)$ étant un $\alpha$-isomorphisme, on en déduit 
que $\ker(u)$ est un $(A|X)$-module de type $\alpha$-fini. Par suite, $E$ est $\alpha$-cohérent. 
Inversement, supposons que $E$ soit $\alpha$-cohérent et montrons que $F$ est $\alpha$-cohérent. On sait que $F$ est de type $\alpha$-fini. 
Soient $X$ un objet de $\cT$, $F'$ un $(A|X)$-module libre de type fini, $v\colon F'\rightarrow F|X$ un $(A|X)$-morphisme. Pour tout $\gamma\in \fm$,
il existe un $\gamma^2$-isomorphisme $g\colon F\rightarrow E$ \eqref{alpha3}. Le morphisme canonique
$\ker(v)\rightarrow \ker(g\circ v)$ est donc un $\gamma^2$-isomorphisme. Comme $\ker(g\circ v)$ est de type $\alpha$-fini, 
on déduit que $\ker(v)$ est de type $\gamma^5$-fini d'après \ref{finita15}(ii)
et est donc de type $\alpha$-fini. Par suite, $F$ est $\alpha$-cohérent.

\begin{lem}\label{finita14}
Soient $f\colon (\cT',A')\rightarrow (\cT,A)$ un morphisme de topos annelés, $\gamma\in R$.
\begin{itemize}
\item[{\em (i)}] Si $u\colon E\rightarrow F$ est un $\gamma$-isomorphisme de $A$-modules, alors $f^*(u)\colon f^*(E)\rightarrow f^*(F)$ est un 
$\gamma^2$-isomorphisme. 
\item[{\rm (ii)}] Si $E\rightarrow F\rightarrow G\rightarrow 0$ est une suite $\gamma$-exacte de $A$-modules, 
son image inverse $f^*(E)\rightarrow f^*(F)\rightarrow f^*(G)\rightarrow 0$ est $\gamma^2$-exacte.
\item[{\rm (iii)}] Si $F$ est un $A$-module de type $\gamma$-fini 
(resp. de présentation $\gamma$-finie, resp. de type $\alpha$-fini, resp. de présentation $\alpha$-finie), son image inverse $f^*(F)$ est de type $\gamma$-fini 
(resp. de présentation $\gamma^2$-finie, resp. de type $\alpha$-fini, resp. de présentation $\alpha$-finie). 
\end{itemize}
\end{lem}

(i) Cela résulte aussitôt de \ref{alpha3}.

(ii) Notant $G'$ le conoyau du morphisme $E\rightarrow F$, il revient au même de dire que la suite 
$E\rightarrow F\rightarrow G\rightarrow 0$ est $\gamma$-exacte ou que le morphisme induit $G'\rightarrow G$ 
est un $\gamma$-isomorphisme. La proposition résulte alors de (i). 

(iii) Cela résulte aussitôt de (ii) et des définitions.

\begin{lem}\label{finita10}
Soient $X$ un objet de $\cT$, $F$ un $A$-module. 
Si $F$ est $\alpha$-cohérent, il en est de même du $(A|X)$-module $F|X$. 
\end{lem}

Cela résulte aussitôt des définitions \ref{finita9}. 

\begin{lem}\label{finita11}
Soient $(X_i)_{i\in I}$ un raffinement de l'objet final de $\cT$, $F$ un $A$-module, $n$ un entier $\geq 0$, $\gamma\in R$. 
Si pour tout $i\in I$, le $(A|X_i)$-module $F|X_i$ est de $n$-présentation $\gamma$-finie 
(resp. de $n$-présentation $\alpha$-finie, resp. $\alpha$-cohérent), il en est de même de $F$.
\end{lem}

Cela résulte aussitôt des définitions \ref{finita9}.

\begin{lem}\label{finita22}
Soient $x$ un point de $\cT$ {\rm (\cite{sga4} IV 6.1)}, $F$ un $A$-module. Si $F$ est de type $\alpha$-fini 
(resp. de présentation $\alpha$-finie, resp. $\alpha$-cohérent), il en est de même du $A_x$-module $F_x$. 
\end{lem}

Les deux premières assertions sont immédiates. Supposons que le $A$-module $F$ soit $\alpha$-cohérent. 
On sait alors que le $A_x$-module $F_x$ est de type $\alpha$-fini. Soient $s_1,\dots,s_n\in F_x$ $(n\geq 1)$. Notons
$u\colon A_x^n\rightarrow F_x$ le morphisme défini par les $s_i$. 
D'après (\cite{sga4} IV 6.8), il existe un voisinage $X$ de $x$ dans $\cT$ tel que $s_1,\dots,s_n$ soient les images canoniques de 
sections $\sigma_1,\dots,\sigma_n\in \Gamma(X,F)$. Notant $v\colon (A|X)^n\rightarrow F|X$ le morphisme défini par les $\sigma_i$, 
on a alors $u=v_x$. Comme le $(A|X$)-module $\ker(v)$ est de type $\alpha$-fini, le $A_x$-module $\ker(u)$ est de type $\alpha$-fini. 
Par suite, le $A_x$-module $F_x$ est $\alpha$-cohérent. 

\begin{lem}\label{finita21}
Soit $F$ un $A$-module tel que pour tout $\gamma\in \fm$, il existe un raffinement $(X_i)_{i\in I}$ de l'objet final de $\cT$ 
et pour tout $i\in I$, un $(A|X_i)$-module de type $\alpha$-fini 
(resp. de présentation $\alpha$-finie, resp. $\alpha$-cohérent) $G_i$ et un $\gamma$-isomorphisme $G_i\rightarrow F|X_i$. 
Alors, $F$ est de type $\alpha$-fini (resp. de présentation $\alpha$-finie, resp. $\alpha$-cohérent).
\end{lem}

En effet, les assertions sont immédiates pour les propriétés d'être de type $\alpha$-fini ou de présentation $\alpha$-finie. 
Considérons la propriété d'être $\alpha$-cohérent. 
Il est clair que $F$ est de type $\alpha$-fini. Soient $X$ un objet de $\cT$, $E$ un $(A|X)$-module libre de type fini,
$u\colon E\rightarrow F|X$ un $(A|X)$-morphisme. Montrons que $\ker(u)$ est un $(A|X)$-module de type $\alpha$-fini. 
Soit $\gamma\in \fm$. Quitte à remplacer $X$ par un raffinement \eqref{finita11}, il existe 
un $(A|X)$-module $\alpha$-cohérent $G$, un $\gamma$-isomorphisme $f\colon G\rightarrow F|X$
et un $(A|X)$-morphisme $v\colon E\rightarrow G$ tel que $f\circ v=\gamma u$. 
On a donc $\ker(u)\subset \ker(\gamma u)\supset \ker(v)$, et les quotients $\ker(\gamma u)/ \ker(u)$ et $\ker(\gamma u)/ \ker(v)$ sont annulés par $\gamma$.
Il existe donc un $\gamma^3$-isomorphisme $\ker(v)\rightarrow \ker(u)$ \eqref{alpha3}. 
Comme $G$ est $\alpha$-cohérent, $\ker(v)$ est de type $\alpha$-fini. 
On en déduit que $\ker(u)$ est de type $\alpha$-fini et par suite que $F$ est $\alpha$-cohérent.

\begin{prop}\label{finita8}
Soient $F$ un $A$-module, $n$ un entier $\geq 0$, $\gamma\in R$.
Supposons que $\cT$ soit équivalent au topos des faisceaux de $\mU$-ensembles sur un $\mU$-site $\cC$. 
Pour tout objet $X$ de $\cC$, on note $X^\tta$ le faisceau associé à $X$. Alors,
\begin{itemize} 
\item[{\rm (i)}] Pour que $F$ soit de $n$-présentation $\gamma$-finie, il faut et il suffit que pour tout $U\in \ob(\cC)$, 
la sous-catégorie pleine de $\cC_{/U}$ formée des objets $X\rightarrow U$ tels qu'il existe une suite $\gamma$-exacte 
\begin{equation}
E^{-n}\rightarrow E^{-n+1}\rightarrow \dots\rightarrow E^0\rightarrow F|X^\tta\rightarrow 0,
\end{equation} 
avec $E^i$ un $(A|X^\tta)$-module libre de type fini pour tout $-n\leq i\leq 0$, soit un raffinement de $U$. 
\item[{\rm (ii)}] Pour que $F$ soit $\alpha$-cohérent, il faut et il suffit que $F$ soit de type $\alpha$-fini et que pour tout $X\in \ob(\cC)$
et tout morphisme $u\colon E\rightarrow F|X^\tta$, où $E$ est un $(A|X^\tta)$-module libre de type fini, 
$\ker(u)$ soit un $(A|X^\tta)$-module de type $\alpha$-fini. 
\end{itemize}
\end{prop}

Cela résulte aussitôt de \ref{finita11} et (\cite{sga4} II 4.4 et 4.10).

\subsection{}\label{finita4}
Soient $\gamma\in R$, 
\begin{equation}\label{finita4a}
\xymatrix{
&E\ar[r]^u\ar[d]_e&F\ar[r]^v\ar[d]_f&G\ar[d]_g\ar[r]&0\\
0\ar[r]&E'\ar[r]^{u'}&F'\ar[r]^{v'}&G'&}
\end{equation}
un diagramme commutatif de $A$-modules tel que les lignes soient $\gamma$-exactes. 
Ce dernier induit des morphismes $E'\rightarrow \ker(v')$ et $\coker(u)\rightarrow G$ et un diagramme commutatif 
\begin{equation}\label{finita4c}
\xymatrix{
&E\ar[r]^u\ar[d]_{e'}&F\ar[r]\ar[d]_f&{\coker(u)}\ar[d]_{g'}\ar[r]&0\\
0\ar[r]&{\ker(v')}\ar[r]&F'\ar[r]^{v'}&G'&}
\end{equation}
On en déduit un diagramme commutatif 
\begin{equation}\label{finita4d}
\xymatrix{
{\ker(e)}\ar[rd]^{u_1}\ar@{^(->}[d]_a&&&{\coker(e)}\ar[d]_c\ar[rd]^{u'_1}&&\\
{\ker (e')}\ar[r]&{\ker(f)}\ar[r]\ar[rd]_{v_1}&{\ker(g')}\ar[d]^b\ar[r]&{\coker(e')}\ar[r]&{\coker(f)}\ar[rd]_{v'_1}\ar[r]&{\coker(g')}\ar@{->>}[d]^-(0.5)d\\
&&{\ker(g)}&&&{\coker(g)}}
\end{equation}
où la ligne centrale est exacte, $a$ est injectif et $d$ est surjectif. Il résulte aussitôt des hypothèses que $a$, $b$, 
$c$ et $d$ sont des $\gamma$-isomorphismes. En particulier, les suites $(u_1,v_1)$ et $(u'_1,v'_1)$ sont $\gamma^2$-exactes.

\begin{prop}\label{finita12}
Soient $\gamma,\gamma',\gamma''$ et $\lambda$ des éléments de $R$, 
$0\rightarrow F'\rightarrow F \rightarrow F''\rightarrow 0$ une suite $\lambda$-exacte de $A$-modules. 
\begin{itemize}
\item[{\rm (i)}] Si $F$ est de type $\gamma$-fini, $F''$ est de type $(\lambda\gamma)$-fini. 
\item[{\rm (ii)}] Si $F'$ est de type $\gamma'$-fini (resp. présentation $\gamma'$-finie) et
$F''$ est de type $\gamma''$-fini (resp. présentation $\gamma''$-finie), 
$F$ est de type $(\lambda^2\gamma'\gamma'')$-fini (resp. présentation $(\lambda^5\gamma'^3\gamma'')$-finie). 
\item[{\rm (iii)}] Si $F$ est de type $\gamma$-fini et $F''$ est de présentation $\gamma''$-finie, 
$F'$ est de type $(\lambda^8\gamma\gamma''^4)$-fini. 
\item[{\rm (iv)}] Si $F$ est de présentation $\gamma$-finie et $F'$ est de type $\gamma'$-fini, 
$F''$ est de présentation $(\lambda\gamma^3\gamma')$-finie.
\end{itemize}
\end{prop} 

(i) C'est immédiat. 

(ii) La question étant locale, on peut supposer qu'il existe des morphismes $A$-linéaires    
$f'\colon E'\rightarrow F'$ et $f''\colon E''\rightarrow F''$ avec $E'$ et $E''$ des $A$-modules libres de type fini
dont les conoyaux sont annulés par $\gamma'$ et $\gamma''$ respectivement. Notons $u'\colon E'\rightarrow F$ le morphisme
induit par $f'$. D'autre part, il existe un morphisme $u''\colon E''\rightarrow F$ qui relève $\lambda f''$. On voit aussitôt que 
le conoyau du morphisme $u'+u''\colon E'\oplus E''\rightarrow F$ est annulé par $\lambda^2\gamma'\gamma''$.
Par suite, $F$ est de type $(\lambda^2\gamma'\gamma'')$-fini. D'après \ref{finita4}, on a une suite $(\lambda^2\gamma')$-exacte
\begin{equation}\label{finita12a}
0\rightarrow \ker(f')\stackrel{a}{\rightarrow} \ker(u'+u'')\stackrel{b}{\rightarrow} \ker(\lambda f'')\rightarrow 0.
\end{equation}
En effet, la suite privée du zéro de droite est $\lambda^2$-exacte d'après \ref{finita4}; on vérifie aussitôt que le conoyau de $b$ 
est annulé par $\lambda^2\gamma'$. 
Supposons $F'$ de présentation $\gamma'$-finie et $F''$ de présentation $\gamma''$-finie. Montrons que $F$ est de présentation 
$(\lambda^4\gamma'^3\gamma'')$-finie. La question étant locale, on peut supposer de plus  
$\ker(f')$ de type $\gamma'$-fini et  $\ker(f'')$ de type $\gamma''$-fini. Par suite, $\ker(\lambda f'')$ est de type $(\lambda\gamma'')$-fini.
Il résulte alors de \eqref{finita12a} et de ce qui précède que $\ker(u'+u'')$ est de type 
$(\lambda^5\gamma'^3\gamma'')$-fini~; d'où l'assertion. 

(iii) La question étant locale, on peut se borner au cas où il existe une suite $\gamma''$-exacte 
\begin{equation}
G''\stackrel{g}{\longrightarrow} E''\stackrel{f}{\longrightarrow} F''\longrightarrow 0
\end{equation} 
avec $E''$ et $G''$ des $A$-modules libres de type fini. Il existe alors un diagramme commutatif 
\begin{equation}
\xymatrix{
&G''\ar[r]^{\lambda g}\ar[d]_v&E''\ar[r]^{\lambda f}\ar[d]_{u}&F''\ar@{=}[d]\ar[r]&0\\
0\ar[r]&F'\ar[r]&F\ar[r]&F''\ar[r]&0}
\end{equation}
Les lignes étant $(\lambda^2\gamma'')$-exactes,
le morphisme canonique $\coker(v)\rightarrow \coker(u)$ est un $(\lambda^4\gamma''^2)$-isomorphisme en vertu de \ref{finita4}. 
D'après (i), $\coker(u)$ est de type $\gamma$-fini. Donc $\coker(v)$ est  de type $(\lambda^8\gamma''^4\gamma)$-fini
en vertu de (ii) ou de \ref{alpha3}. Compte tenu de la suite exacte $0\rightarrow \im(v)\rightarrow F'\rightarrow \coker(v)\rightarrow 0$,
on en déduit  par (ii) que $F'$ est de type $(\lambda^8\gamma\gamma''^4)$-fini.

(iv) La question étant locale, on peut se borner au cas où il existe un $A$-module libre de type fini $E$ 
et un morphisme $A$-linéaire $f\colon E\rightarrow F$ dont le conoyau est annulé par $\gamma$ et dont le noyau est de type
$\gamma$-fini. Notons $u\colon F\rightarrow F''$ le morphisme donné et posons $f''=u\circ f\colon E\rightarrow F''$. 
On a alors une suite exacte 
\begin{equation}
0\rightarrow \ker(f)\rightarrow \ker(f'')\rightarrow \ker(u)\rightarrow \coker(f).
\end{equation}
Le conoyau du morphisme canonique $F'\rightarrow \ker(u)$ est annulé par $\lambda$. Par suite, $\ker(u)$ est de type 
$(\lambda\gamma')$-fini. On en déduit que $\ker(f'')$ est de type $(\lambda\gamma^3\gamma')$-fini en vertu de (ii). 
Comme $\coker(f'')$ est annulé par $\lambda\gamma$, $F''$ est de présentation $(\lambda\gamma^3\gamma')$-finie.

\begin{cor}\label{finita13}
Soit $0\rightarrow F'\rightarrow F \rightarrow F''\rightarrow 0$ une suite $\alpha$-exacte de $A$-modules. Alors, 
\begin{itemize}
\item[{\rm (i)}] Si $F$ est de type $\alpha$-fini, il en est de même de $F''$. 
\item[{\rm (ii)}] Si $F'$ et $F''$ sont de type $\alpha$-fini (resp. présentation $\alpha$-finie), il en est de même de $F$. 
\item[{\rm (iii)}] Si $F$ est de type $\alpha$-fini et $F''$ est de présentation $\alpha$-finie, $F'$ est de type $\alpha$-fini. 
\item[{\rm (iv)}] Si $F$ est de présentation $\alpha$-finie et $F'$ est de type $\alpha$-fini, $F''$ est de présentation $\alpha$-finie.
\end{itemize}
\end{cor} 

\begin{prop}\label{finita17}
Soit $0\rightarrow F'\stackrel{u}{\rightarrow} F \stackrel{v}{\rightarrow} F''\rightarrow 0$ une suite $\alpha$-exacte de $A$-modules. 
Si deux $A$-modules sont $\alpha$-cohérents, il en est de même du troisième.
\end{prop}

Compte tenu de \ref{finita16}, quitte à remplacer $F'$ par $\ker(v)$ et $F''$ par $\im(v)$, on peut supposer la suite $0\rightarrow F'\rightarrow F \rightarrow F''\rightarrow 0$ exacte. 

Supposons d'abord que $F$ et $F''$ soient $\alpha$-cohérents. D'après \ref{finita13}(iii), $F'$ est de type $\alpha$-fini. Il est donc $\alpha$-cohérent en tant que sous-$A$-module de $F$. 

Supposons ensuite que $F$ et $F'$ soient $\alpha$-cohérents. D'après \ref{finita13}(i), $F''$ est de type $\alpha$-fini. 
Soient $X$ un objet de $\cT$, $f\colon (A|X)^n\rightarrow F''|X$ un morphisme $(A|X)$-linéaire. Montrons que pour tout $\gamma\in \fm$, $\ker(f)$ est de type $\gamma^2$-fini.
La question étant locale, on peut supposer qu'il existe un morphisme $(A|X)$-linéaire $f'\colon (A|X)^n\rightarrow F|X$ tel que $f=v\circ f'$.
Comme $F'$ est de type $\alpha$-fini, on peut supposer qu'il existe un morphisme $(A|X)$-linéaire et $\gamma$-surjectif $h\colon (A|X)^m\rightarrow F'|X$. 
On en déduit un morphisme $g\colon (A|X)^{m+n}\rightarrow F|X$ qui s'insère dans un diagramme commutatif 
\begin{equation}
\xymatrix{
{F'|X}\ar[r]^-(0.5)u&{F|X}\ar[r]^-(0.5)v&{F''|X}\\
(A|X)^m\ar[r]^-(0.5)i\ar[u]^h&(A|X)^{m+n}\ar[r]^-(0.5)\pi\ar[u]^g&(A|X)^n\ar[u]_f}
\end{equation}
où, identifiant $(A|X)^{m+n}$ à $(A|X)^m\oplus (A|X)^n$, $i(x)=x+0$ et $\pi(x+x')=x'$. On en déduit une suite exacte
\begin{equation}
\ker(g)\rightarrow \ker(f)\rightarrow \coker(h).
\end{equation} 
Par hypothèse, $\ker(g)$ est de type $\alpha$-fini et $\coker(h)$ est $\gamma$-nul; donc $\ker(f)$ est de type $\gamma^2$-fini d'après \ref{finita12}(i), ce qui prouve que 
$F''$ est $\alpha$-cohérent. 

Supposons enfin que $F'$ et $F''$ soient $\alpha$-cohérents. 
D'après \ref{finita13}(ii), $F$ est de type $\alpha$-fini. Soient $X$ un objet de $\cT$, $f\colon (A|X)^n\rightarrow F|X$ un morphisme $(A|X)$-linéaire. 
Montrons que pour tout $\gamma\in \fm$, $\ker(f)$ est de type $\gamma^2$-fini. Comme $F''$ est $\alpha$-cohérent, $\ker(v\circ f)$ est de type $\alpha$-fini. 
La question étant locale, on peut supposer qu'il existe un morphisme $(A|X)$-linéaire et $\gamma$-surjectif $w\colon (A|X)^m\rightarrow \ker(v\circ f)$. 
Il existe un morphisme $(A|X)$-linéaire $f'\colon (A|X)^m\rightarrow F'|X$ qui s'insère dans un diagramme commutatif 
\begin{equation}
\xymatrix{
{F'|X}\ar[r]^u&{F|X}\ar[r]^v&{F''|X}\\
{(A|X)^m}\ar[r]^{w'}\ar[u]^{f'}&{(A|X)^n}\ar[u]_f&}
\end{equation}
où $w'$ est induit par $w$. On en déduit un morphisme $(A|X)$-linéaire et $\gamma$-surjectif  $\ker(f')\rightarrow \ker(f)$. 
Comme $F'$ est cohérent, $\ker(f')$ est de type $\alpha$-fini~; donc $\ker(f)$ est de type $\gamma^2$-fini d'après \ref{finita12}(i), ce qui prouve que 
$F$ est $\alpha$-cohérent.

\begin{cor}\label{finita18}
Soient $F$, $G$ deux $A$-modules $\alpha$-cohérents, $u\colon F\rightarrow G$ un morphisme $A$-linéaire.
Alors, $\ker(u)$, $\im(u)$ et $\coker(u)$ sont $\alpha$-cohérents.
\end{cor} 

En effet, $\im(u)$ est clairement de type $\alpha$-fini~; étant un sous-module d'un $A$-module $\alpha$-cohérent, il est $\alpha$-cohérent.  
La proposition résulte alors de \ref{finita17} appliqué aux suites exactes
$0\rightarrow \ker(u)\rightarrow F\rightarrow \im(u)\rightarrow 0$ et $0\rightarrow \im(u)\rightarrow G\rightarrow \coker(u)\rightarrow 0$. 

\begin{cor}\label{finita180}
Soient $F$, $G$ deux $A$-modules $\alpha$-cohérents. Alors, $F\otimes_AG$ et $\cHom_A(F,G)$ sont $\alpha$-cohérents.
\end{cor} 

Il suffit de montrer que pour tout $\gamma\in \fm$, il existe un raffinement $(X_i)_{i\in I}$ de l'objet final de $\cT$ et pour tout $i\in I$,
deux $(A|X_i)$-modules $\alpha$-cohérents $M_i$ et $N_i$ de $\cT_{/X_i}$ et deux $\gamma^2$-isomorphismes 
$M_i\rightarrow (F\otimes_AG)|X_i$ et $N_i\rightarrow \cHom_A(F,G)|X_i$ \eqref{finita21}. 
On peut clairement se borner au cas où il existe un $A$-module de présentation finie $F'$ et un $\gamma$-isomorphisme $F'\rightarrow F$. 
Par hypothèse, il existe deux entiers $m,n\geq 1$ et une suite exacte $A^m\rightarrow A^n\rightarrow F'\rightarrow 0$. On en déduit deux suites exactes 
\begin{eqnarray}
G^m\rightarrow G^n\rightarrow F'\otimes_AG\rightarrow 0,\\
0\rightarrow \cHom_A(F',G)\rightarrow \cHom_A(A^n,G)\rightarrow \cHom_A(A^m,G).
\end{eqnarray}
Il résulte alors de \ref{finita18} que $F'\otimes_AG$ et $\cHom_A(F',G)$ sont $\alpha$-cohérents, d'où la proposition recherchée \eqref{alpha3}.

\begin{prop}\label{finita19}
Pour que la $R$-algèbre $A$ soit $\alpha$-cohérente \eqref{finita9}, il faut et il suffit que tout $A$-module de présentation $\alpha$-finie soit $\alpha$-cohérent.
\end{prop}

Il est clair que la condition est suffisante. Montrons qu'elle est nécessaire. Supposons que la $R$-algèbre $A$ soit $\alpha$-cohérente
et montrons que tout $A$-module de présentation $\alpha$-finie $M$ est $\alpha$-cohérent. 
Soient $E$ un $A$-module libre de type fini, $u\colon E\rightarrow M$ un morphisme $A$-linéaire, $\gamma\in \fm$. 
Il suffit de montrer que $\ker(u)$ est de type $\gamma^3$-fini. La question étant locale, on peut supposer qu'il existe 
deux $A$-modules libres de type fini $E_0$ et $E_1$, et une suite $\gamma$-exacte
\begin{equation}\label{finita19a}
E_1\stackrel{v}{\rightarrow} E_0\rightarrow M\rightarrow 0.
\end{equation}
Notons $M'$ le conoyau de $v$, et $w\colon M'\rightarrow M$ le morphisme induit par \eqref{finita19a}. Il existe un morphisme $w'\colon M\rightarrow M'$
tel que $w'\circ w=\gamma^2\id _{M'}$ et $w\circ w'=\gamma^2\id_M$ \eqref{alpha3}. 
On a $\ker(u)\subset \ker(w'\circ u)\subset \ker(\gamma^2 u)$. Par suite, la multiplication par $\gamma^2$ sur $E$
induit un morphisme $\ker(w'\circ u)\rightarrow \ker(u)$ dont le conoyau est annulé par $\gamma^2$. 
Comme $M'$ est $\alpha$-cohérent en vertu de \ref{finita18}, $\ker(w'\circ u)$  est de type $\alpha$-fini. 
Par suite, $\ker(u)$ est de type $\gamma^3$-fini \eqref{finita12}, d'où l'assertion recherchée.

\begin{prop}\label{finita20}
Supposons que $A$ soit une $R$-algèbre $\alpha$-cohérente et soit $J$ un idéal $\alpha$-cohérent de $A$. 
Pour qu'un $(A/J)$-module $F$ soit $\alpha$-cohérent, il faut et il suffit qu'en tant que $A$-module, $F$ soit $\alpha$-cohérent.
En particulier, $A/J$ est une $R$-algèbre $\alpha$-cohérente. 
\end{prop}

On notera d'abord que $A/J$ est un $A$-module $\alpha$-cohérent \eqref{finita18}. Si $F$ est un $(A/J)$-module $\alpha$-cohérent, 
il existe un raffinement $(X_i)_{i\in I}$ de l'objet final de $\cT$ tel que pour tout $i\in I$, $F|X_i$ soit le conoyau d'un morphisme 
$(A/J)^{n_i}|X_i\rightarrow (A/J)^{m_i}|X_i$. Par suite, $F$ est un $A$-module $\alpha$-cohérent en vertu \ref{finita18}. 

Inversement, supposons qu'en tant que $A$-module, $F$ soit $\alpha$-cohérent. \'Etant un $A$-module de type $\alpha$-fini, 
$F$ est un $(A/J)$-module de type $\alpha$-fini. Soient $X$ un objet de $\cT$, $f\colon (A/J)^n|X\rightarrow F|X$ un morphisme $(A/J)$-linéaire.
Le morphisme canonique $u\colon A^n|X\rightarrow(A/J)^n|X$ induit un morphisme surjectif $\ker(f\circ u)\rightarrow \ker(f)$. Comme $\ker(f\circ u)$
est un $(A|X)$-module de type $\alpha$-fini, $\ker(f)$ est un $(A/J)|X$ de type $\alpha$-fini.

\section{\texorpdfstring{Modules $\alpha$-cohérents sur un schéma}{Modules alpha-cohérents sur un schéma}}\label{afini}

\subsection{}\label{afini1}
Les hypothèses et notations de \ref{alpha1} sont en vigueur dans cette section. 
Pour tout $R$-schéma $X$, on désigne par $\bMod(\co_X)$ la catégorie des $\co_X$-modules de $X_\zar$ \eqref{notconv12}
et  par $\aMod(\co_X)$ la catégorie des $\alpha$-$\co_X$-modules, c'est-à-dire 
le quotient de $\bMod(\co_X)$ par la sous-catégorie pleine formée des $\co_X$-modules $\alpha$-nuls (cf. \ref{alpha20}).
Lorsque nous parlons de $\co_X$-module sans préciser, il est sous-entendu qu'il s'agit d'un $\co_X$-module de $X_\zar$.

\begin{prop}\label{afini2}
Soient $X$ un $R$-schéma, $F$ un $\co_X$-module. Les conditions suivantes sont équivalentes~:
\begin{itemize}
\item[{\rm (i)}] Il existe un $\co_X$-module quasi-cohérent $G$ et un isomorphisme de $\alpha$-$\co_X$-modules
$\alpha(F)\stackrel{\sim}{\rightarrow} \alpha(G)$. 
\item[{\rm (ii)}] Pour tout $x\in X$, il existe un voisinage ouvert $U$ de $x$ dans $X$, un $\co_U$-module 
quasi-cohérent $G$ et un isomorphisme de $\alpha$-$\co_U$-modules
$\alpha(F|U)\stackrel{\sim}{\rightarrow} \alpha(G)$. 
\item[{\rm (iii)}] Le $\co_X$-module $\fm\otimes_RF$ est quasi-cohérent. 
\end{itemize}
\end{prop}

On a clairement (i)$\Rightarrow$(ii). Montrons (ii)$\Rightarrow$(iii). 
La question étant locale, on peut supposer la condition (i) remplie. D'après \ref{alpha21}(iv),
il existe donc un $\co_X$-module quasi-cohérent $G$ et un $\alpha$-isomorphisme $u\colon\fm\otimes_RF\rightarrow G$. 
Comme $\fm\otimes_R\fm\simeq \fm$, le morphisme $\fm\otimes_RF\rightarrow \fm\otimes_RG$ induit par $u$
est un isomorphisme en vertu de \ref{alpha21}(iii); d'où la condition (iii). 
Enfin, l'implication (iii)$\Rightarrow$(i) est immédiate puisque le morphisme canonique $\fm\otimes_RF\rightarrow F$
est un $\alpha$-isomorphisme. 

\begin{defi}\label{afini3}
Soient $X$ un $R$-schéma, $F$ un $\co_X$-module. 
On dit que $F$ est {\em $\alpha$-quasi-cohérent} s'il remplit les conditions de \ref{afini2}.
\end{defi}

\begin{lem}\label{afini4}
Soient $X$ un $R$-schéma cohérent ({\em i.e.}, quasi-compact et quasi-séparé), $\gamma\in R$, 
$F$ un $\co_X$-module quasi-cohérent de type $\gamma$-fini \eqref{finita9}.
Alors, il existe un $\co_X$-module quasi-cohérent de présentation finie $G$ et un morphisme $\co_X$-linéaire
$f\colon G\rightarrow F$ de conoyau annulé par $\gamma$.
\end{lem}

En effet, d'après (\cite{ega1n} 6.9.12), $F$ est isomorphe à la limite inductive d'un système inductif filtrant 
de $\co_X$-modules de présentation finie $(F_i)_{i\in I}$. 
Pour tout $i\in I$, notons $f_i\colon F_i\rightarrow F$ le morphisme canonique. 
Montrons qu'il existe $i\in I$ tel que le conoyau de $f_i$ soit annulé par $\gamma$. 
Comme $X$ est quasi-compact et que la catégorie $I$ est filtrante, on peut supposer qu'il existe
un morphisme $\co_X$-linéaire $u\colon E\rightarrow F$ avec $E$ un $\co_X$-module libre de type fini, 
dont le conoyau est annulé par $\gamma$. Il existe $i\in I$ et un morphisme $\co_X$-linéaire $u_i\colon E\rightarrow F_i$
tel que $u=f_i\circ u_i$. Par suite, le conoyau de $f_i$ est annulé par $\gamma$; d'où la proposition.

\begin{prop}\label{afini5}
Soient $X$ un $R$-schéma cohérent ({\em i.e.}, quasi-compact et quasi-séparé), $F$ un $\co_X$-module $\alpha$-quasi-cohérent.
Pour que $F$ soit de type $\alpha$-fini (resp. présentation $\alpha$-finie) \eqref{finita9}, 
il faut et il suffit que pour tout $\gamma\in \fm$, 
il existe un $\gamma$-isomorphisme $f\colon G\rightarrow F$ avec $G$ un $\co_X$-module quasi-cohérent de type fini 
(resp. de présentation finie).
\end{prop}

Considérons d'abord l'assertion non-respée.  
La condition est clairement suffisante. Inversement, supposons $F$ de type $\alpha$-fini et 
montrons que la condition est satisfaite. Soit $\gamma\in \fm$. 
On peut se borner au cas où $F$ est quasi-cohérent \eqref{afini2}. D'après \ref{afini4}, 
il existe un $\co_X$-module quasi-cohérent de présentation finie $G'$ et un morphisme $\co_X$-linéaire
$f'\colon G'\rightarrow F$ de conoyau annulé par $\gamma$. 
Le morphisme $f\colon \im(f')\rightarrow F$ induit par $f'$ satisfait alors à la condition pour $\gamma$.  

Considérons ensuite l'assertion respée. La condition est clairement suffisante. 
Inversement, supposons $F$ de présentation $\alpha$-finie et montrons que la condition est satisfaite. Soit $\gamma\in \fm$. 
On peut se borner au cas où $F$ est quasi-cohérent \eqref{afini2}. D'après \ref{afini4}, 
il existe un $\co_X$-module quasi-cohérent de présentation finie $G'$ et un morphisme $\co_X$-linéaire
$f'\colon G'\rightarrow F$ de conoyau annulé par $\gamma$. En vertu de \ref{finita12}(iii), 
$\ker(f')$ est quasi-cohérent de type $\gamma^{12}$-fini. D'après \ref{afini4}, il existe 
un $\co_X$-module quasi-cohérent de présentation finie $G''$ et un morphisme $f''\colon G''\rightarrow G'$ tels que 
$f'\circ f''=0$ et que $\ker(f')/\im(f'')$ soit annulé par $\gamma^{12}$. Le morphisme $f\colon G'/\im(f'')\rightarrow F$ induit par $f'$
répond alors à la condition pour $\gamma^{12}$.

\begin{cor}\label{afini102}
Soit $X$ un $R$-schéma affine, $\cA$ une $\co_X$-algèbre quasi-cohérente, $\cF$ un $\cA$-module quasi-cohérent (comme $\cA$-module ou comme 
$\co_X$-module). Pour que le $\cA$-module $\cF$ soit de type $\alpha$-fini (resp. de présentation $\alpha$-finie), 
il faut et il suffit que le $\cA(X)$-module $\cF(X)$ soit de type $\alpha$-fini (resp. de présentation $\alpha$-finie).
\end{cor}

En effet, posons $Y=\Spec(\cA(X))$ et notons $\pi\colon Y\rightarrow X$ le morphisme canonique. Le foncteur $\pi_*$ induit 
une équivalence entre la catégorie des $\co_Y$-modules quasi-cohérents de $Y_\zar$ et celle des $\cA$-modules quasi-cohérents 
de $X_\zar$ (\cite{ega1n} 9.2.1).
Notons $F$ le $\co_Y$-module quasi-cohérent de $Y_\zar$ tel que $\pi_*(F)=\cF$. Nous allons montrer que les trois conditions suivantes sont équivalentes~:
\begin{itemize}
\item[(i)] le $\cA$-module $\cF$ est de type $\alpha$-fini (resp. de présentation $\alpha$-finie, resp. $\alpha$-cohérent);
\item[(ii)] le $\co_Y$-module $F$ est de type $\alpha$-fini (resp. de présentation $\alpha$-finie, resp. $\alpha$-cohérent);
\item[(iii)] le $\cA(X)$-module $\cF(X)$ est de type $\alpha$-fini (resp. de présentation $\alpha$-finie, resp. $\alpha$-cohérent).
\end{itemize}

L'implication (i)$\Rightarrow$(ii) résulte immédiatement du fait que pour tout ouvert $U$ de $X$,  notant $\pi\colon Y_U\rightarrow U$ le changement de base de 
$\pi$, le foncteur $\pi_{U*}$ induit une équivalence entre la catégorie des $\co_{Y_U}$-modules quasi-cohérents et celle des $\cA|U$-modules quasi-cohérents. 
L'implication (ii)$\Rightarrow$(iii) est une conséquence de \ref{afini5}.
L'implication (iii)$\Rightarrow$(i) est immédiate.

\begin{cor}\label{afini103}
Soit $X$ un $R$-schéma affine, $\cA$ une $\co_X$-algèbre quasi-cohérente, $\cF$ un $\cA$-module quasi-cohérent (comme $\cA$-module ou comme 
$\co_X$-module). Pour que le $\cA$-module $\cF$ soit $\alpha$-cohérent, 
il faut et il suffit que le $\cA(X)$-module $\cF(X)$ soit $\alpha$-cohérent.
\end{cor}

Supposons d'abord que le $\cA(X)$-module $\cF(X)$ 
soit $\alpha$-cohérent. On vérifie aussitôt que pour tout $f\in \co_X(X)$, le $\cA(X)_f$-module $\cF(X)_f$ est $\alpha$-cohérent. On en déduit 
que le $\cA$-module $\cF$ est $\alpha$-cohérent. Inversement, si le $\cA$-module $\cF$ est $\alpha$-cohérent, le $\cA(X)$-module $\cF(X)$ 
est $\alpha$-cohérent d'après \ref{afini102}. 

\begin{cor}\label{afini104}
Soit $X$ un $R$-schéma affine d'anneau $A$. 
Pour qu'un $\co_X$-module $\alpha$-quasi-cohérent $F$ soit de type $\alpha$-fini (resp. de présentation $\alpha$-finie, resp. $\alpha$-cohérent), 
il faut et il suffit que le $A$-module $F(X)$ soit  de type $\alpha$-fini (resp. de présentation $\alpha$-finie, resp. $\alpha$-cohérent).
\end{cor}

Cela résulte de  \ref{finita16}, \ref{afini102} et \ref{afini103} en considérant $\fm\otimes_RF$.

\begin{prop}\label{afini11}
Soient $f\colon X'\rightarrow X$ un $R$-morphisme fidèlement plat et quasi-compact
de $R$-schémas, $F$ un $\co_X$-module $\alpha$-quasi-cohérent.
Pour que $F$ soit de type $\alpha$-fini (resp. de présentation $\alpha$-finie), il faut et il suffit
que son image inverse sur $X'$ le soit. 
\end{prop}

En effet, la condition est nécessaire en vertu de \ref{finita14}. Montrons qu'elle est suffisante. 
On peut supposer $X$ affine. Remplaçant $X'$ par une somme d'ouverts affines qui le recouvrent, 
on est ramené au cas où $X'$ est également affine. L'assertion résulte alors de \ref{afini104} et (\cite{agt} V.8.1 et V.8.5).

\begin{cor}\label{afini12}
Soient $f\colon X'\rightarrow X$ un $R$-morphisme fidèlement plat et quasi-compact
de $R$-schémas, $F$ un $\co_X$-module $\alpha$-quasi-cohérent.
Supposons que les $R$-algèbres $\co_X$ et $\co_{X'}$ soient $\alpha$-cohérentes \eqref{finita9}. 
Pour que $F$ soit $\alpha$-cohérent, il faut et il suffit que son image inverse sur $X'$ le soit.
\end{cor}

Cela résulte de \ref{finita19} et \ref{afini11}.

\begin{prop}\label{afini6}
Soient $X$ un $R$-schéma cohérent, 
$F$ un $\co_X$-module $\alpha$-quasi-cohérent. Supposons que l'anneau $\co_X$ soit cohérent.
Alors, les conditions suivantes sont équivalentes~:
\begin{itemize}
\item[{\rm (i)}] $F$ est $\alpha$-cohérent \eqref{finita9};
\item[{\rm (ii)}] $F$ est de présentation $\alpha$-finie~;
\item[{\rm (iii)}] pour tout $\gamma\in \fm$, il existe un $\gamma$-isomorphisme 
$G\rightarrow F$ avec $G$ un $\co_X$-module cohérent.
\end{itemize}
\end{prop}

On a clairement (i)$\Rightarrow$(ii). L'implication (ii)$\Rightarrow$(iii) résulte de \ref{afini5}. Montrons (iii)$\Rightarrow$(i). 
Supposons la condition (iii) satisfaite. D'après \ref{afini5}, $F$ est de type $\alpha$-fini. Il suffit donc de montrer que pour tout
morphisme $\co_X$-linéaire $f\colon E\rightarrow F$ où $E$ est un $\co_X$-module libre de type fini, $\ker(f)$ est de type $\alpha$-fini. 
Soient $\gamma\in \fm$, $g\colon G\rightarrow F$ un $\gamma$-isomorphisme tel que $G$ soit un $\co_X$-module cohérent.
Il existe un morphisme $\co_X$-linéaire $f'\colon E\rightarrow G$ tel que $g\circ f'=\gamma f$. 
Le conoyau de l'injection canonique $\ker(f')\rightarrow \ker(\gamma f)$ est annulé par $\gamma$. 
Par ailleurs, la multiplication par $\gamma$ dans $E$ induit un morphisme $\ker(\gamma f)\rightarrow \ker(f)$
dont le conoyau est clairement annulé par $\gamma$. On en déduit un morphisme $\co_X$-linéaire 
$\ker(f')\rightarrow \ker(f)$ dont le conoyau est annulé par $\gamma^2$. Comme $G$ est cohérent, $\ker(f')$ est de type fini. 
Par suite, $\ker(f)$ est de type $\gamma^2$-fini et donc de type $\alpha$-fini.

\begin{prop}\label{afini13}
Soient $Y$ un $R$-schéma affine, $X$ un $Y$-schéma projectif, $F$ un $\co_X$-module $\alpha$-quasi-cohérent. Supposons que l'anneau $\co_X$ soit $\alpha$-cohérent.
Pour que $F$ soit $\alpha$-cohérent, il faut et il suffit que pour tout $\gamma\in \fm$, il existe 
un morphisme de complexes de $\co_X$-modules $G^\bullet \rightarrow F[0]$ dont les noyau et conoyau sont $\gamma$-acycliques, 
tel que $G^i$ soit localement libre de type fini pour tout $i\in \mZ$ et soit nul pour $i>0$. 
\end{prop}

Montrons d'abord que la condition est suffisante. Soient $\gamma\in \fm$, $u\colon G^\bullet \rightarrow F[0]$ 
un morphisme de complexes de $\co_X$-modules dont les noyau et conoyau sont $\gamma$-acycliques et 
tel que $G^i$ soit localement libre de type fini pour tout $i\in \mZ$ et soit nul pour $i>0$. 
Les noyau et conoyau du morphisme
$\cH^0(G^\bullet)\rightarrow F$ induit par $u$ sont annulés par $\gamma^2$. Comme $\cH^0(G^\bullet)$ est $\alpha$-cohérent en vertu de \ref{finita11} et \ref{finita18}, 
on en déduit que $F$ est de présentation $\gamma^3$-finie \eqref{finita9}. Par suite, $F$ est de présentation $\alpha$-finie et est donc $\alpha$-cohérent d'après \ref{finita19}.
Montrons ensuite que la condition est nécessaire. Supposons que $F$ soit $\alpha$-cohérent. 
Soient $\cL$ un $\co_X$-module inversible très ample pour $Y$, $\gamma\in \fm$. 
En vertu de \ref{afini5}, il existe un morphisme $u^0\colon G^0 \rightarrow F$ dont le conoyau est annulé par $\gamma$,
avec $G^0$ un $\co_X$-module somme directe d'un nombre fini de copies de $\cL^n$ pour un entier $n<0$. 
Comme $\ker(u^0)$ est $\alpha$-cohérent, procédant par récurrence, on construit un morphisme de complexes de $\co_X$-modules $u\colon G^\bullet \rightarrow F[0]$ 
dont les noyau et conoyau sont $\gamma$-acycliques tel que pour tout $i>0$, $G^i=0$ et pour tout $i\leq 0$,
$G^i$ soit somme directe d'un nombre fini de copies de $\cL^{n_i}$ pour un entier $n_i<0$.

\begin{defi}[\cite{egr1} 1.4.1]\label{afini7}
On dit qu'un anneau $A$ est {\em universellement cohérent} 
si l'anneau de polynômes $A[t_1,\dots,t_n]$ est cohérent pour tout entier $n\geq 0$. 
\end{defi}

\begin{prop}\label{afini8}
Soient $Y$ un $R$-schéma, $f\colon X\rightarrow Y$ un morphisme propre de présentation finie,
$F$ un $\co_X$-module $\alpha$-quasi-cohérent et $\alpha$-cohérent.
Supposons que $Y$ admette un recouvrement par des ouverts affines $(U_i)_{i\in I}$
tels que $\Gamma(U_i,\co_Y)$ soit universellement cohérent pour tout $i\in I$. 
Alors, pour tout entier $q\geq 0$, $\rR^qf_*(F)$ est un $\co_Y$-module $\alpha$-quasi-cohérent et $\alpha$-cohérent. 
\end{prop}

On notera d'abord que les anneaux $\co_X$ de $X_\zar$ et $\co_Y$ de $Y_\zar$ sont cohérents d'après (\cite{egr1} 1.4.2 et 1.4.3).
Par ailleurs, la question étant locale sur $Y$, on peut le supposer affine, de sorte que $X$ est cohérent. 
Le morphisme canonique $\rR^qf_*(\fm\otimes_RF)\rightarrow \rR^qf_*(F)$ est un $\alpha$-isomorphisme \eqref{alpha3}. 
Comme $\fm\otimes_RF$ est quasi-cohérent, $\rR^qf_*(\fm\otimes_RF)$ 
est quasi-cohérent et donc $\rR^qf_*(F)$ est $\alpha$-quasi-cohérent.
Par ailleurs, pour tout $\gamma\in \fm$, il existe un $\gamma$-isomorphisme 
$u\colon G\rightarrow F$ avec $G$ un $\co_X$-module cohérent d'après \ref{afini6}. 
Il existe donc un morphisme $\co_X$-linéaire $v\colon F\rightarrow G$ tel que $u\circ v=\gamma^2\cdot \id_F$ et 
$v\circ u=\gamma^2\cdot \id_G$ \eqref{alpha3}. En vertu  de (\cite{egr1} 1.4.8), le $\co_Y$-module
$\rR^qf_*(G)$ est cohérent. Par suite, $\rR^qf_*(F)$ est $\alpha$-cohérent \eqref{afini6}.

\begin{cor}\label{afini9}
Soient $Y$ un $R$-schéma affine d'anneau $A$, $f\colon X\rightarrow Y$ un morphisme propre de présentation finie,
$F$ un $\co_X$-module $\alpha$-quasi-cohérent et $\alpha$-cohérent.
Supposons que $Y$ admette un recouvrement par des ouverts affines $(U_i)_{i\in I}$
tels que $\Gamma(U_i,\co_Y)$ soit universellement cohérent pour tout $i\in I$. 
Alors, pour tout entier $q\geq 0$, le $A$-module $\rH^q(X,F)$ est $\alpha$-cohérent. 
\end{cor}

On notera d'abord que l'anneau $\co_X$ de $X_\zar$ est cohérent (\cite{egr1} 1.4.2 et 1.4.3).
Pour tout $\gamma\in \fm$, il existe un $\gamma$-isomorphisme $u\colon G\rightarrow F$ avec $G$ un $\co_X$-module cohérent d'après \ref{afini6}. 
Il existe donc un morphisme $\co_X$-linéaire $v\colon F\rightarrow G$ tel que $u\circ v=\gamma^2\cdot \id_F$ et 
$v\circ u=\gamma^2\cdot \id_G$ \eqref{alpha3}. En vertu  de (\cite{egr1} 1.4.8 et 1.4.3), le $A$-module
$\rH^q(X,G)$ est cohérent. Par suite, $\rH^q(X,F)$ est $\alpha$-cohérent \eqref{afini6}. 

\begin{rema}\label{afini99}
On notera que (\cite{egr1} 1.4.8), utilisé dans les preuves de \ref{afini8} et \ref{afini9}, 
est un corollaire d'un théorème de Kiehl (\cite{kiehl} 2.9'a) (cf. \cite{egr1} 1.4.7)
\end{rema}

\begin{defi}\label{afini14}
On dit qu'une $R$-algèbre $A$ est {\em universellement $\alpha$-cohérente} 
si la $R$-algèbre de polynômes $A[t_1,\dots,t_n]$ est $\alpha$-cohérente \eqref{finita9} pour tout entier $n\geq 0$. 
\end{defi}

\begin{prop}\label{afini15}
Soient $Y$ un $R$-schéma, $f\colon X\rightarrow Y$ un morphisme projectif de présentation finie,
$F$ un $\co_X$-module $\alpha$-quasi-cohérent et $\alpha$-cohérent.
Supposons que $f$ soit pseudo-cohérent {\rm (\cite{sga6} III 1.2)} et que $Y$ admette un recouvrement par des ouverts affines $(U_i)_{i\in I}$
tels que $\Gamma(U_i,\co_Y)$ soit universellement $\alpha$-cohérent pour tout $i\in I$. 
Alors, pour tout entier $q\geq 0$, $\rR^qf_*(F)$ est un $\co_Y$-module $\alpha$-quasi-cohérent et $\alpha$-cohérent. 
\end{prop}

Le problème étant local sur $Y$, on peut se borner au cas où $Y=\Spec(A)$ est affine d'anneau $A$ universellement $\alpha$-cohérent. 
Il résulte de \ref{afini103} que $\co_Y$ est un anneau $\alpha$-cohérent. 
Montrons ensuite que l'anneau $\co_X$ est $\alpha$-cohérent. La question étant locale, on peut supposer que $X=\Spec(B)$ est affine. 
Il existe alors un $A$-homomorphisme surjectif $A'=A[t_1,\dots,t_n]\rightarrow B$ de noyau $I$ un idéal de type fini. 
Comme $A'$ est un anneau $\alpha$-cohérent par hypothèse, $I$ est un idéal $\alpha$-cohérent.
Par suite, $B$ est un anneau $\alpha$-cohérent d'après \ref{finita20} et il en est alors de même de $\co_X$ \eqref{afini103}.

Si $C$ est un complexe borné $\gamma$-acyclique de $\co_X$-modules de longueur $\ell$, 
$\rR^qf_*(C)$ est annulé par $\gamma^{\ell+1}$. Ceci se démontre par récurrence sur la longueur de $C$ en utilisant la filtration canonique. 
Par ailleurs, d'après (\cite{ega3} 1.4.12), il existe un entier $r\geq 1$ tel que pour tout entier $q\geq r$ et tout $\co_X$-module quasi-cohérent $M$, 
on ait $\rR^qf_*(M)=0$.  

Pour démontrer la proposition, on peut se réduire au cas où $F$ est un $\co_X$-module quasi-cohérent \eqref{alpha3}. 
En vertu de \ref{afini13}, pour tout $\gamma\in \fm$, il existe 
un morphisme de complexes de $\co_X$-modules $G^\bullet \rightarrow F[0]$ dont les noyau et conoyau sont $\gamma$-acycliques 
et tel que $G^i$ soit localement libre de type fini pour tout $i\in \mZ$ et soit nul pour $i>0$. Notant $H$ le noyau de la différentielle 
$G^{-r}\rightarrow G^{-r+1}$, on a une suite $\gamma$-exacte de complexes bornés de $\co_X$-modules 
\begin{equation}
0\rightarrow H[r]\rightarrow \tau_{\geq -r}(G^\bullet)\rightarrow F[0]\rightarrow 0,
\end{equation}
où $\tau_{\geq -r}(G^\bullet)$ est le $(-r)$-ième cran la filtration bête de $G^\bullet$. 
On en déduit, pour tout entier $q\geq 0$, un morphisme 
\begin{equation}
\rR^qf_*(\tau_{\geq -r}(G^\bullet))\rightarrow \rR^qf_*(F)
\end{equation}
dont les noyau et conoyau sont annulés par $\gamma^{2(r+1)}$.

Le complexe de $\co_X$-modules $\tau_{\geq -r}(G^\bullet)$ 
étant pseudo-cohérent (\cite{sga6} I 2.3; cf. aussi \cite{egr1} 1.3.8), il est pseudo-cohérent relativement à $f$ en vertu de (\cite{sga6} III 1.12). 
Comme il est borné, le complexe de $\co_Y$-modules $\rR f_*(\tau_{\geq -r}(G^\bullet))$ est pseudo-cohérent en vertu de (\cite{sga6} III 2.2). 
Il s'ensuit que le $\co_Y$-module $\rR^qf_*(\tau_{\geq -r}(G^\bullet))$ est $\alpha$-cohérent pour tout entier $q\geq 0$ \eqref{finita18}. 
On en déduit que $\rR^qf_*(F)$ est de présentation $\gamma^{2r+3}$-finie. Par suite, $\rR^qf_*(F)$ est de présentation $\alpha$-finie 
et donc $\alpha$-cohérent \eqref{finita19}.

\begin{cor}\label{afini16}
Soient $Y$ un $R$-schéma affine d'anneau $A$, $f\colon X\rightarrow Y$ un morphisme projectif de présentation finie,
$F$ un $\co_X$-module $\alpha$-quasi-cohérent et $\alpha$-cohérent.
Supposons que $f$ soit pseudo-cohérent {\rm (\cite{sga6} III 1.2)} et que $Y$ admette un recouvrement par des ouverts affines $(U_i)_{i\in I}$
tels que $\Gamma(U_i,\co_Y)$ soit universellement $\alpha$-cohérent pour tout $i\in I$. 
Alors, pour tout entier $q\geq 0$, le $A$-module $\rH^q(X,F)$ est $\alpha$-cohérent. 
\end{cor}

En effet, considérons la suite spectrale de Cartan-Leray 
\begin{equation}
\rE_2^{i,j}=\rH^i(Y,\rR^jf_*(F))\Rightarrow \rH^{i+j}(X,F).
\end{equation}
Il résulte de \ref{afini15} et \ref{alpha3} que $\rE_2^{i,j}$ est $\alpha$-nul pour tout $i\geq 1$.  
Par suite, pour tout $q\geq 0$, le morphisme canonique
\begin{equation}
\rH^q(X,F)\rightarrow \rH^0(Y,\rR^qf_*F)
\end{equation}
est un $\alpha$-isomorphisme. La proposition résulte alors de \ref{afini15}, \ref{afini104} et \ref{finita16}.

\begin{lem}\label{afini17}
Soient $X$ un $R$-schéma affine, $\ox$ un point géométrique de $X$ de support $x$, $A$ une $\co_X$-algèbre quasi-cohérente,
$F$ un $A$-module quasi-cohérent (comme $A$-module ou comme $\co_X$-module). 
Notons encore $A$ la $\co_X$-algèbre de $X_\et$ associée à $A$ par le foncteur \eqref{notconv12a} et
$F$ le $A$-module de $X_\et$ associé à $F$. Supposons que la $R$-algèbre $A(X)$ soit universellement $\alpha$-cohérente
et que le $A$-module $F$ soit $\alpha$-cohérent. Alors,
\begin{itemize}
\item[{\rm (i)}] Les $R$-algèbres $A_x$ et $A_\ox$ sont universellement $\alpha$-cohérentes. 
\item[{\rm (ii)}] Le $A_\ox$-module $F_\ox$ est $\alpha$-cohérent.  
\end{itemize}
\end{lem}

(i) En effet, la $R$-algèbre $A$ de $X_\zar$ est $\alpha$-cohérente en vertu de \ref{afini103}.
La $R$-algèbre $A_x$ est donc universellement $\alpha$-cohérente d'après \ref{finita22}. 
On a $A_\ox\simeq A_x\otimes_{\co_{X,x}}\co_{X,\ox}$ \eqref{notconv12g}. 
Considérons $\co_{X,\ox}$ comme une limite inductive filtrante de $\co_{X,x}$-algèbres étales $(B_i)_{i\in I}$. 
Soient $n$ un entier $\geq 1$, $u\colon  A_{\ox}^n\rightarrow A_\ox$ une forme $A_\ox$-linéaire. 
Il existe $i\in \ob(I)$ et un morphisme $B_i$-linéaire $u_i\colon  (A_x\otimes_{\co_{X,x}}B_i)^n\rightarrow A_x\otimes_{\co_{X,x}}B_i$ tels que 
$u=u_i\otimes_{B_i}\co_{X,\ox}$.  La $R$-algèbre $A_x\otimes_{\co_{X,x}}B_i$ étant $\alpha$-cohérente, le $B_i$-module $\ker(u_i)$ est de type $\alpha$-fini. 
Comme pour tout morphisme $i\rightarrow j$ de $I$, $B_j$ est $B_i$-plat, le $B_j$-module $\ker(u_i\otimes_{B_i}B_j)$ est de type $\alpha$-fini.  
Par suite, le $A_{X,\ox}$-module $\ker(u)$ est de type $\alpha$-fini. On en déduit que la $R$-algèbre $A_{X,\ox}$ est $\alpha$-cohérente. 

Pour tout entier $g\geq 1$, on note $A[t_1,\dots,t_g]$ la $A$-algèbre de polynômes en $g$ variables de $X_\zar$.
Remplaçant $A$ par $A[t_1,\dots,t_g]$, on déduit de ce qui précède que  
la $R$-algèbre $A_{X,\ox}$ est universellement $\alpha$-cohérente.

(ii) En effet, on a $F_\ox\simeq F_x\otimes_{\co_{X,x}}\co_{X,\ox}$ \eqref{notconv12g}. 
Le $A_x$-module $F_x$ est $\alpha$-cohérent d'après \ref{finita22}. La proposition résulte alors de (ii), \ref{finita14} et \ref{finita19}

\section{\texorpdfstring{Rappel sur les algèbres $\alpha$-étales}{Algèbres alpha-étales}}\label{aet}

\subsection{}\label{aet1}
Les hypothèses et notations de \ref{alpha1} sont en vigueur dans cette section. Nous renvoyons à (\cite{agt} V) pour 
les définitions et les principales propriétés des notions utilisées en {\em $\alpha$-algèbre} 
(ou presque-algèbre dans {\em loc. cit.})~: 
modules {\em $\alpha$-projectifs} ({\em loc. cit.} V.3.2),  {\em $\alpha$-projectifs de type $\alpha$-fini et de rang fini} ({\em loc. cit.} V.5.6 \& V.5.14), 
{\em $\alpha$-plats} ({\em loc. cit.} V.6.1) et {\em $\alpha$-fidèlement plats} ({\em loc. cit.} V.6.4).

\begin{defi}[\cite{gr} 3.1.1; \cite{agt} V.7.1]\label{aet2}
Soit $A\rightarrow B$ un homomorphisme de $R$-algèbres. Considérons les 3 conditions suivantes:
\begin{itemize}
\item[(i)] $B$ est $\alpha$-plat en tant que $A$-module;
\item[(i')] $B$ est $\alpha$-projectif de type $\alpha$-fini et de rang fini en tant que $A$-module;
\item[(ii)] $B$ est $\alpha$-projectif en tant que $B\otimes_AB$-module.
\end{itemize}  
On dit que la $A$-algèbre $B$ est {\em $\alpha$-étale} (ou que $B$ est {\em $\alpha$-étale sur} $A$)
si elle vérifie les conditions (i) et (ii), et qu'elle est {\em $\alpha$-finie-étale} (ou que $B$ est {\em $\alpha$-finie-étale sur} $A$) 
si elle vérifie les conditions (i') et (ii).
\end{defi}

Nous renvoyons à (\cite{agt} V.7) pour les principales propriétés des algèbres {\em $\alpha$-finies-étales}, 
appelées {\em revêtements presque-étales} dans {\em loc. cit.}
Nous n'utiliserons que très peu la notion d'algèbre $\alpha$-étale (cf. \cite{gr} 3.1). 
On notera qu'une algèbre $\alpha$-étale et uniformément presque finie projective au sens de (\cite{gr} 3.1.1) est $\alpha$-finie-étale.

\begin{prop}[\cite{agt} V.7.11]\label{aet3}
Soit $A$ une $R$-algèbre telle que pour tout (ou ce qui revient au même pour un) $\alpha\in \Lambda^+$, $\pi^\alpha$ ne soit un diviseur de zéro dans $A$,  
et que $\Spec(A)$ soit normal et localement irréductible {\rm (\cite{agt} III.3.1)}. Soit $ B$ une $A$-algèbre $\alpha$-finie-étale. 
Notons $A_\pi$ (resp. $B_\pi$) l'anneau des fractions 
de $A$ (resp. $B$) dont les dénominateurs sont de la forme $\pi^\alpha$ pour $\alpha\in \Lambda^+$, et $B'$ la clôture 
intégrale de $A$ dans $B_\pi$. Alors, l'homomorphisme canonique $B\rightarrow B_\pi$ se factorise en
\begin{equation}
B\rightarrow \Hom_R(\fm,B')\rightarrow \Hom_R(\fm,B_\pi) \stackrel{\sim}{\rightarrow} B_\pi,
\end{equation}
où la première flèche est un $\alpha$-isomorphisme et les autres flèches sont les homomorphismes canoniques  \eqref{alpha5c}. En particulier, l'homomorphisme $A\rightarrow B'$ est $\alpha$-fini-étale.  
\end{prop} 

En effet, on peut se borner au cas où $A$ est intègre (\cite{agt} V.7.4(2) et V.7.8), auquel cas l'assertion est démontrée dans
(\cite{agt} V.7.11).

\begin{defi}\label{aet4}
Soient $A$ une $R$-algèbre, $G$ un groupe fini, $B$ une $A[G]$-algèbre.
On dit que $B$ est un {\em $\alpha$-$G$-torseur sur $A$} s'il existe un homomorphisme $\alpha$-fidèlement
plat $A\rightarrow A'$ et un homomorphisme de $A'[G]$-algèbres 
\begin{equation}
B\otimes_AA'\rightarrow \prod_GA',
\end{equation}
où $G$ agit sur le terme de droite par multiplication à droite sur lui même, 
qui est un $\alpha$-isomorphisme.
\end{defi}

On renvoie à (\cite{agt} V.12) pour les principales propriétés de cette notion. 

\begin{prop}[\cite{agt} V.12.5] \label{aet5}
Soient $A$ une $R$-algèbre, $G$ un groupe fini, $B$ un $\alpha$-$G$-torseur sur $A$, $M$ un $B$-module
muni d'une action semi-linéaire de $G$. Alors, le morphisme canonique 
$B\otimes_AM^G\rightarrow M$ est un $\alpha$-isomorphisme.
\end{prop}

\begin{cor}[\cite{agt} V.12.6] \label{aet6}
Soient $A$ une $R$-algèbre, $G$ un groupe fini, $B$ un $\alpha$-$G$-torseur sur $A$. 
Alors, l'homomorphisme naturel $A\rightarrow B^G$ est un $\alpha$-isomorphisme.
\end{cor}

\begin{prop}[\cite{agt} V.12.8]\label{aet7}
Soient $A$ une $R$-algèbre, $G$ un groupe fini, $B$ un $\alpha$-$G$-torseur sur $A$, $M$
un $B$-module muni d'une action semi-linéaire de $G$, $\Tr_G$ le morphisme $A$-linéaire
de $M$ défini par $\sum_{\sigma\in G}\sigma$. Alors, pour tout $q\geq 1$,
$\rH^q(G,M)$ et $M^G/\Tr_G(M)$ sont $\alpha$-nuls. 
\end{prop}

\begin{prop}[\cite{agt} V.12.9]\label{aet8}
Soit $A$ une $R$-algèbre telle que pour tout
(ou ce qui revient au même pour un) $\alpha\in \Lambda^+$, $\pi^\alpha$ ne soit un diviseur de zéro dans $A$,  
et que $\Spec(A)$ soit normal et localement irréductible {\rm (\cite{agt} III.3.1)}. 
Notons $A_\pi$  l'anneau des fractions de $A$ dont les dénominateurs sont de la forme $\pi^\alpha$ pour $\alpha\in \Lambda^+$. Soient $U$ un ouvert de $\Spec(A_\pi)$, $V\rightarrow U$ un morphisme étale fini et galoisien de groupe 
$G$, $Y$ la fermeture intégrale de $\Spec(A)$ dans $V$, $B=\Gamma(Y,\co_Y)$. 
Si $B$ est une $A$-algèbre $\alpha$-finie-étale, alors $B$ est un $\alpha$-$G$-torseur au-dessus de $A$. 
\end{prop}

On peut se borner au cas où $U$ est irréductible (\cite{agt} V.7.4(3)). 
Il suffit alors de calquer la preuve de (\cite{agt} V.12.9) en remplaçant (\cite{agt} V.7.11) par \ref{aet3}.

\begin{lem}\label{aet10}
Soient  $(B_i)_{i\in I}$ un système inductif filtrant de $R$-algèbres $\alpha$-cohérentes \eqref{finita9} 
(resp. universellement $\alpha$-cohérentes \eqref{afini14}), $B_\infty$ sa limite inductive. 
Supposons que pour tout morphisme $i\rightarrow j$ de $I$, $B_j$ est un $B_i$-module $\alpha$-plat. 
Alors, la $R$-algèbre $B_\infty$ est $\alpha$-cohérente (resp. universellement $\alpha$-cohérente).  
\end{lem} 

Considérons un entier $n\geq 1$ et une forme $B_\infty$-linéaire $u_\infty \colon B_\infty^n\rightarrow B_\infty$.
Il existe $i\in \ob(I)$ et une forme $B_i$-linéaire $u_i\colon B_i^n\rightarrow B_i$ tels que $u_\infty=u_i\otimes_{B_i}B_\infty$. 
Le $B_i$-module $\ker(u_i)$ étant de type $\alpha$-fini, pour tout $\gamma \in \fm$, 
il existe une suite $\gamma$-exacte de $B_i$-modules 
\begin{equation}\label{aet10a}
B_i^m\rightarrow B_i^n\stackrel{u_i}{\rightarrow}B_i.
\end{equation}
Pour tout morphisme $i\rightarrow j$ de $I$, la suite 
\begin{equation}\label{aet10b}
B_j^m\longrightarrow B_j^n\stackrel{u_j}{\longrightarrow} B_j
\end{equation}
déduite de \eqref{aet10b}  par extension des scalaires de $B_i$ à $B_j$,
est $\gamma^2$-exacte. Il en est donc de même de la suite
\begin{equation}\label{aet10c}
B_\infty^m\longrightarrow B_\infty^n\stackrel{u_\infty}{\longrightarrow} B_\infty
\end{equation}
déduite de \eqref{aet10b}  par extension des scalaires de $B_i$ à $B_\infty$. 
Par suite, la $R$-algèbre $B_\infty$ est $\alpha$-cohérente.

Supposons les $R$-algèbres $B_i$ universellement $\alpha$-cohérentes. Pour entier $g\geq 0$, 
la $R$-algèbre $B_\infty[t_1,\dots,t_g]$ est $\alpha$-cohérente d'après ce qui précède et (\cite{agt} V.6.2). 
La $R$-algèbre $B_\infty$ est donc universellement $\alpha$-cohérente. 

\begin{lem}\label{aet100}
Soient  $(B_i)_{i\in I}$ un système inductif filtrant d'anneaux cohérents 
(resp. universellement cohérents \eqref{afini7}), $B_\infty$ sa limite inductive. 
Supposons que pour tout morphisme $i\rightarrow j$ de $I$, $B_j$ soit un $B_i$-module plat. 
Alors, l'anneau $B_\infty$ est cohérent (resp. universellement cohérent).  
\end{lem} 

La preuve est similaire à celle de \ref{aet10}, mais plus facile. 

\begin{lem}\label{aet9}
Soient $B$ une $R$-algèbre $\alpha$-cohérente \eqref{finita9} (resp. universellement $\alpha$-cohérente \eqref{afini14}),  
$(B_i)_{i\in I}$ un système inductif filtrant de $B$-algèbres $\alpha$-finies-étales,
$B_\infty$ sa limite inductive. Alors, la $R$-algèbre $B_\infty$ et chacune des  $R$-algèbres $B_i$ ($i\in I$) 
sont $\alpha$-cohérentes (resp. universellement $\alpha$-cohérentes).  
\end{lem} 

En effet, pour tout morphisme $i\rightarrow j$ de $I$, $B_j$ est une $B_i$-algèbre $\alpha$-finie-étale (\cite{agt} V.7.8) et en particulier,
un $B_i$-module $\alpha$-plat (\cite{agt} V.6.3). 
Pour tout $i\in \ob(I)$, $B_i$ est un $B$-module de presentation $\alpha$-finie (\cite{agt} V.8.8) et donc $\alpha$-cohérent \eqref{finita19}.  
Considérons un entier $n\geq 1$ et une forme $B_i$-linéaire $u \colon B_i^n\rightarrow B_i$. 
Le $B$-module $\ker(u)$ est $\alpha$-cohérent d'après \ref{finita18}
et est a fortiori de type $\alpha$-fini sur $B_i$. Par suite, la $R$-algèbre $B_i$ est $\alpha$-cohérente.

Si la $R$-algèbre $B$ est universellement $\alpha$-cohérente, pour tout $i\in \ob(I)$ et tout entier $g\geq 0$, 
la $R$-algèbre $B_i[t_1,\dots,t_g]$ est $\alpha$-cohérente d'après ce qui précède et (\cite{agt} V.7.4). 
Par suite, la $R$-algèbre $B_i$ est universellement $\alpha$-cohérente.

La proposition résulte donc de \ref{aet10}.

\begin{lem}\label{aet11}
Soient $A$ une $R$-algèbre, $(A_i)_{i\in I}$ un système inductif filtrant de $A$-algèbres $\alpha$-fidèlement plates, $A_\infty$ sa limite inductive.
Alors, la $A$-algèbre $A_\infty$ est $\alpha$-fidèlement plate.
\end{lem}

En effet, d'après (\cite{gr} 3.1.2(vi)), pour tout $i\in I$, l'homomorphisme canonique $A\rightarrow A_i$ est $\alpha$-injectif et le $A$-module
$A_i/A$ est $\alpha$-plat. Par suite, le morphisme canonique $\fm\otimes_RA\rightarrow \fm\otimes_RA_i$ est injectif.  
On en déduit aussitôt que le morphisme canonique $\fm\otimes_RA\rightarrow \fm\otimes_RA_\infty$ est injectif,
et par suite que l'homomorphisme canonique $A\rightarrow A_\infty$ est $\alpha$-injectif.  
Par ailleurs, il résulte aussitôt que le $A$-module $A_\infty/A$, qui est la limite inductive du système inductif $(A_i/A)_{i\in I}$, est $\alpha$-plat. 
La $A$-algèbre $A_\infty$ est donc $\alpha$-fidèlement plate en vertu de  (\cite{gr} 3.1.2(vi)).

\section{\texorpdfstring{Modules de type $\alpha$-fini sur un anneau de valuation non discrète de hauteur $1$}
{Modules de type alpha-fini sur un anneau de valuation non discrète de hauteur 1}}\label{mptf}

\subsection{}\label{mptf1}
Dans cette section, $R$ désigne un anneau muni d'une valuation non discrète de hauteur $1$, 
$v\colon R\rightarrow \mR\cup \{\infty\}$. On note $\fm$ l'idéal maximal de $R$ et on pose $\Lambda=v(R-\{0\})$ et $\Lambda^+=v(\fm-\{0\})$. 
Pour tout $\varepsilon\in \Lambda^+$, on note $\fm_\varepsilon$ l'idéal de $R$ formé des éléments $x\in R$ tels que $v(x)\geq \varepsilon$.
On observera que les hypothèses de \ref{alpha1} sont satisfaites. Il est donc loisible de considérer dans ce contexte
les notions introduites dans les sections \ref{alpha}, \ref{finita} et \ref{afini}. Nous considérons toujours $R$ 
comme muni de la topologie définie par l'idéal $\fm_\varepsilon$ pour un élément quelconque $\varepsilon\in \Lambda^+$; 
c'est un anneau prévaluatif séparé (\cite{egr1} 1.9.8). 

\begin{prop}[\cite{egr1} 1.12.15]\label{mptf2}
Si $R$ est complet et séparé, toute $R$-algèbre topologiquement de présentation finie est un anneau 
universellement cohérent \eqref{afini7}. 
\end{prop}

\subsection{}\label{mptf3}
Suivant (\cite{scholze1} 2.2), pour tous $R$-modules $M$ et $N$ et tout $\gamma \in \fm$, 
on dit que $M$ et $N$ sont {\em $\gamma$-équivalents} et on note $M\approx_\gamma N$,
s'il existe deux $R$-morphismes $f\colon M\rightarrow N$ et $g\colon N\rightarrow M$ 
tels que $f\circ g=\gamma\cdot \id_N$ et $g\circ f=\gamma\cdot \id_M$. 
On dit que $M$ et $N$ sont {\em $\alpha$-équivalents} et on note $M\approx N$
si $M$ et $N$ sont $\gamma$-équivalents pour tout $\gamma\in \fm$. 
Les relations $\approx_\gamma$ et $\approx$ sont clairement symétriques. 
Pour tous $R$-modules $M$, $N$ et $P$ et tous $\gamma,\delta \in \fm$, si $M\approx_\gamma N$
et $N\approx_\delta P$, alors $M\approx_{\gamma\delta} P$. En particulier, la relation $\approx$ est transitive. 

\subsection{}\label{mptf4}
Pour tout $\varepsilon\in \mR_{\geq 0}$, on désigne par $I_\varepsilon$ l'idéal de $R$ formé des éléments 
$x\in R$ tels que $v(x)>\varepsilon$. Alors $I_\varepsilon\approx R$. 
En effet, pour tout $\gamma \in \fm$, il existe $\gamma_1,\gamma_2\in R$ tels que $\gamma\gamma_1=\gamma_2$
et $v(\gamma_1)\leq \varepsilon <v(\gamma_2)$. Les morphismes 
$R\rightarrow I_\varepsilon$ et $I_\varepsilon\rightarrow R$ définis par les multiplications par $\gamma_2$ et 
$\gamma_1^{-1}$, respectivement, montrent que $I_\varepsilon\approx_\gamma R$. 
Par ailleurs, comme $\fm\cdot I_\varepsilon=I_\varepsilon$,
si $\varepsilon\not\in \Lambda$, $R$ et $I_\varepsilon$ ne sont pas $\alpha$-isomorphes en vertu de \ref{alpha7} et \eqref{alpha8a}. 

\begin{remas}\label{mptf5}  
(i)\ Deux $R$-modules $\alpha$-isomorphes sont $\alpha$-équivalents \eqref{alpha3}. 
La réciproque n'est pas vraie \eqref{mptf4}. 

(ii)\ Pour qu'un $R$-module $M$ soit {\em de type $\alpha$-fini} (resp. {\em de présentation $\alpha$-finie}) 
il faut et il suffit que pour tout $\gamma\in \fm$, il existe un $R$-module de type fini (resp. de présentation finie) 
$N$ tel que $M\approx_\gamma N$, d'après \ref{finita2}(i) et \ref{alpha3}.
\end{remas}

\begin{prop}[\cite{scholze1} 2.6]\label{mptf62}
Tout $R$-module de type $\alpha$-fini est de présentation $\alpha$-finie. 
\end{prop}

\begin{cor}\label{mptf63}
Tout sous-$R$-module d'un $R$-module de type $\alpha$-fini est de type $\alpha$-fini. 
\end{cor}

Cela résulte de \ref{mptf62} et \ref{finita13}.

\begin{teo}[\cite{scholze1} 2.5]\label{mptf7}
Soit $M$ un $R$-module de type $\alpha$-fini. Alors, il existe une unique suite décroissante de nombres réels positifs 
$(\varepsilon_i)_{i\geq 1}$, tendant vers $0$, et un unique entier $n\geq 0$ tels que 
\begin{equation}\label{mptf7a}
M\approx R^n\oplus (\oplus_{i\geq 1}R/I_{\varepsilon_i}). 
\end{equation}
\end{teo}

\subsection{}\label{mptf8}
On désigne par $\fS$ le $\mR$-espace vectoriel des suites de nombres réels 
$(\varepsilon_i)_{i\geq 1}$, tendant vers $0$, et par $\fS_{\geq}$ le sous-ensemble formé des suites décroissantes.
Pour tout $\varepsilon=(\varepsilon_i)_{i\geq 1}\in \fS$, on pose 
\begin{equation}\label{mptf8a}
\|\varepsilon\|=\sup_{i\geq 1}(|\varepsilon_i|).
\end{equation}
On munit $\fS_{\geq}$ de la relation d'ordre définie pour toutes suites
$\varepsilon=(\varepsilon_i)_{i\geq 1}$ et $\varepsilon'=(\varepsilon'_i)_{i\geq 1}$ de $\fS_\geq$, par $\varepsilon\geq \varepsilon'$
si pour tout $i\geq 1$, $\varepsilon_1+\dots+\varepsilon_i\geq \varepsilon'_1+\dots+\varepsilon'_i$.

\subsection{}\label{mptf9}
En vertu de \ref{mptf7}, on peut associer à tout $R$-module de torsion de type $\alpha$-fini $M$ 
une unique suite $\varepsilon_{M}=(\varepsilon_{M,i})_{i\geq 1}$ de $\fS_{\geq}$ telle que 
\begin{equation}\label{mptf9a}
M\approx \oplus_{i\geq 1}R/I_{\varepsilon_{M,i}}. 
\end{equation}
On appelle {\em longueur} de $M$ et l'on note $\lambda(M)$ l'élément 
\begin{equation}\label{mptf9b}
\lambda(M)=\sum_{i\geq 1}\varepsilon_{M,i}\in \mR_{\geq 0}\cup\{\infty\}.
\end{equation}
D'après (\cite{scholze1} 2.11), la suite $\varepsilon_M$ et la longueur $\lambda(M)$ vérifient les propriétés suivantes~:

\begin{itemize}
\item[(i)] Pour tout $R$-module de torsion et de présentation finie $M$, il existe des éléments  
$\gamma_i\in R$ $(1\leq i\leq n)$ tels que $v(\gamma_1)\geq v(\gamma_2)\geq \dots\geq v(\gamma_n)$ et un $R$-isomorphisme
\begin{equation}
M\stackrel{\sim}{\rightarrow} R/\gamma_1R\oplus R/\gamma_2R\oplus \dots \oplus R/\gamma_nR.
\end{equation}
On a alors $\varepsilon_{M,i}=v(\gamma_i)$ pour tout $1\leq i\leq n$ et $\varepsilon_{M,i}=0$ pour tout $i>n$.
\item[(ii)] Pour que deux $R$-modules de torsion de type $\alpha$-fini $M$ et $N$ 
soient $\gamma$-équivalents (pour $\gamma\in \fm$), il faut et il suffit que 
$\|\varepsilon_M-\varepsilon_N\|\leq v(\gamma)$.
En particulier, pour qu'un $R$-module de torsion de type $\alpha$-fini $M$ soit $\alpha$-nul, il faut et il suffit que 
la suite $\varepsilon_M$ soit nulle. 
\item[(iii)] Pour tous $R$-modules de torsion de type $\alpha$-fini $M$ et $M'$, si $M'$ est un sous-quotient de $M$
alors $\varepsilon_{M',i}\leq \varepsilon_{M,i}$ pour tout $i\geq 1$. 
\item[(iv)] Pour toute suite exacte $0\rightarrow M'\rightarrow M\rightarrow M''\rightarrow 0$ de $R$-modules de torsion 
de type $\alpha$-fini, on a $\varepsilon_M\leq \varepsilon_{M'}+\varepsilon_{M''}$ et $\lambda(M)= \lambda(M')+\lambda(M'')$. 
De plus, si $\varepsilon_M=\varepsilon_{M'}$ (resp. $\varepsilon_M=\varepsilon_{M''}$), alors $M''$ (resp. $M'$) est $\alpha$-nul.
\end{itemize}

\begin{lem}\label{mptf12}
Soient $M$ un $R$-module, $\varpi$ un élément non nul de $\fm$.  On pose $\hR=\underset{\longleftarrow}{\lim} \ R/\varpi^nR$ 
et $\hM=\underset{\longleftarrow}{\lim} \ M/\varpi^nM$, et on suppose que le $R$-module $M/\varpi M$ est de type $\alpha$-fini. Alors, 
le $\hR$-module $\hM$ est de type $\alpha$-fini.
\end{lem}

On notera d'abord que $\hR$ et $\hM$ sont complets et séparés pour les topologies $\varpi$-adiques et pour tout entier $n\geq1$, que le morphisme canonique 
$h_n\colon \hM\rightarrow M/\varpi^n M$ est surjectif (\cite{ac} III §~2.11 prop.~14) et que le morphisme canonique $\hR/\varpi^n \hR\rightarrow R/\varpi^n R$
est un isomorphisme (\cite{ac} III §~2.11 cor.~1 de prop.~14). Soient $d$ un entier $\geq 1$, $\varphi\colon \hR^d\rightarrow \hM$
un morphisme $\hR$-linéaire, $\varphi_1\colon (R/\varpi R)^d\rightarrow M/\varpi M$ le morphisme induit par $\varphi$. 
Supposons que $\coker(\varphi_1)$ soit annulé par un élément $\gamma\in \fm$ tel que $v(\gamma)<v(\varpi)$. Alors, la suite 
\begin{equation}\label{mptf12a}
\hR^d\stackrel{\varphi}{\rightarrow} \hM \stackrel{u}{\rightarrow}\coker(\varphi_1)\rightarrow 0,
\end{equation}
où $u$ est le morphisme canonique, est exacte.
En effet, $u$ est surjectif puisque $h_1$ est surjectif. Soit $x\in \hM$ tel que $u(x)=0$.
On montre par récurrence que pour tout entier $n\geq 1$, qu'il existe $y_n\in \hR^d$ et $x_n\in \hM$ tel que $x=\varphi(y_n)+\gamma^{1-n}\varpi^n x_n$
et $y_{n+1}-y_n\in \gamma^{-n}\varpi^n \hR^d$. On en déduit qu'il existe $y\in \hR^d$ tel que $x=\varphi(y)$; d'où l'assertion recherchée.
Par suite, $\coker(\varphi)$ est annulé par $\gamma$, ce qui implique la proposition.

\subsection{}\label{mptf10}
On suppose dans la suite de cette section que $R$ est de caractéristique $p$, et on note $\varphi$ son endomorphisme de Frobenius. 
Soit $\varpi$ un élément non nul de $\fm$. Pour tout entier $n\geq 1$, on pose $R_n=R/\varpi^n R$, $\bvR=(R_n)_{n\geq 1}$
et $\bvR^{(p)}=(R_{pn})_{n\geq 1}$. On considère $\bvR$ et $\bvR^{(p)}$ comme des anneaux du topos $\Ens^{\mN^\circ}$ (\cite{agt} III.7.7). 
L'endomorphisme $\varphi$ induit un homomorphisme d'anneaux que l'on note
\begin{equation}\label{mptf10a}
\bvvarphi\colon \bvR\rightarrow \bvR^{(p)}.
\end{equation}
Pour tout $\bvR$-module $M=(M_n)_{n\geq 1}$, on pose $M_0=0$, $M[-1]=(M_{n-1})_{n\geq 1}$ qui est naturellement un $\bvR$-module, 
et $M^{(p)}=(M_{pn})_{n\geq 1}$ qui est naturellement un $\bvR^{(p)}$-module. Les correspondances $M\mapsto M[-1]$ et $M\mapsto M^{(p)}$ sont clairement
fonctorielles, et on a un morphisme $\bvR$-linéaire canonique fonctoriel $M\rightarrow M[-1]$.

\begin{lem}[\cite{faltings2} page 224, \cite{scholze1} 2.12]\label{mptf11}
Conservons les hypothèses de \ref{mptf10}, soient de plus $M=(M_n)_{n\geq 1}$ un $\bvR$-module, 
$u\colon M[-1]\rightarrow M$ un morphisme $\bvR$-linéaire, 
\begin{equation}\label{mptf11a}
\phi\colon M\otimes_{\bvR}\bvR^{(p)}\stackrel{\sim}{\rightarrow} M^{(p)}
\end{equation}
un isomorphisme $\bvR^{(p)}$-linéaire. On pose $\hR=\underset{\longleftarrow}{\lim} \ R_n$ et $\hM=\underset{\longleftarrow}{\lim} \ M_n$. On note encore $\varphi$ 
l'endomorphisme de Frobenius de $\hR$ et on désigne par
\begin{equation}\label{mptf11b}
\uphi\colon \hM\rightarrow\hM
\end{equation}
le morphisme $\hR$-semi-linéaire induit par $\phi$. On suppose que 
\begin{itemize}
\item[{\rm (a)}] le $R$-module $M_1$ est de type $\alpha$-fini~;
\item[{\rm (b)}] le composé $M[-1]\stackrel{u}{\rightarrow}M\rightarrow M[-1]$, où la seconde flèche est le morphisme canonique, 
est induit par la multiplication par $\varpi$;
\item[{\rm (c)}] pour tout entier $n\geq 0$, la suite de $R$-modules
\begin{equation}\label{mptf11c}
\xymatrix{M_1\ar[r]\ar[r]^-(0.5){u^{\circ n}}&{M_{n+1}}\ar[r]&{M_n}}
\end{equation}
est exacte au centre.  
\end{itemize}
Alors,
\begin{itemize}
\item[{\rm (i)}] Pour tout entier $n\geq 1$, le morphisme canonique $M_{n+1}\rightarrow M_n$ est $\alpha$-surjectif
et le morphisme canonique $\hM/\varpi^n \hM \rightarrow M_n$ est un $\alpha$-isomorphisme. 
\item[{\rm (ii)}] Le $\hR$-module $\hM$ est de type $\alpha$-fini. 
\item[{\rm (iii)}] Si le corps des fractions de $R$ est algébriquement clos, il existe un entier $d\geq 0$ et un morphisme $\hR$-linéaire 
\begin{equation}\label{mptf11d}
\fm\otimes_R\hM\rightarrow \hR^d,
\end{equation}
compatible aux morphismes $\uphi$ et $\varphi$, et qui est un $\alpha$-isomorphisme.
\end{itemize}
\end{lem}

En effet, il résulte de (c), \ref{finita13} et \ref{mptf9}(iii)-(iv) que pour tout entier $n\geq 1$, le $R$-module $M_n$ est de type $\alpha$-fini  et on a 
$\varepsilon_{M_{n+1}}\leq \varepsilon_{M_n}+\varepsilon_{M_1}$. Par suite, $\varepsilon_{M_n}\leq n \varepsilon_{M_1}$. 
Par ailleurs, on a $\varepsilon_{M_{pn}}=p\varepsilon_{M_n}$ d'après \eqref{mptf11a}. On en déduit que $\varepsilon_{M_n}= n \varepsilon_{M_1}$.
Notons $M^*_{n+1}$ et $M'_n$ les images des première et seconde flèches de \eqref{mptf11c}, respectivement, de sorte qu'on a une suite exacte
$0\rightarrow M^*_{n+1}\rightarrow M_{n+1}\rightarrow M'_n\rightarrow 0$. D'après \ref{mptf9}, on a 
\begin{equation}\label{mptf11e}
(n+1) \varepsilon_{M_1}= \varepsilon_{M_{n+1}}\leq \varepsilon_{M^*_{n+1}}+\varepsilon_{M'_n}\leq \varepsilon_{M_1}+\varepsilon_{M_n}=(n+1) \varepsilon_{M_1}.
\end{equation}
Par suite, on a $\varepsilon_{M^*_{n+1}}=\varepsilon_{M_1}$ et $\varepsilon_{M'_n}=\varepsilon_{M_n}$. Il résulte alors de \ref{mptf9}(iv) que 
les morphismes canoniques $M_1\rightarrow M^*_{n+1}$ et $M'_n\rightarrow M_n$ sont des $\alpha$-isomorphismes, en particulier, la suite 
\begin{equation}\label{mptf11f}
0\longrightarrow M_1\stackrel{u^{\circ n}}{\longrightarrow} M_{n+1}\longrightarrow M_n\longrightarrow 0
\end{equation}
est $\alpha$-exacte. Pour tout entier $m\geq 1$, le diagramme suivant 
\begin{equation}\label{mptf11g}
\xymatrix{
M_1\ar[r]^-(0.5){u^{\circ m}}\ar@{=}[d]&M_{m+1}\ar[r]\ar[d]^-(0.5){u^{\circ n}}&{M_m}\ar[d]^-(0.5){u^{\circ n}}\\
M_1\ar[r]^-(0.5){u^{\circ (m+n)}}&M_{m+n+1}\ar[r]\ar[d]&{M_{m+n}}\ar[d]\\
&M_n\ar@{=}[r]&M_n}
\end{equation}
permet de déduire (par récurrence sur $m$) que la suite 
\begin{equation}\label{mptf11h}
0\longrightarrow M_m\stackrel{u^{\circ n}}{\longrightarrow} M_{n+m}\longrightarrow M_n\longrightarrow 0
\end{equation}
est $\alpha$-exacte. Le système projectif $(\alpha(M_n))_{n\geq 1}$ de $\aMod(R)$ vérifie donc la condition de Mittag-Leffler. Compte tenu de \ref{alpha10}
et (\cite{ega3} 0.13.2.2), on en déduit que la suite de $\hR$-modules
\begin{equation}\label{mptf11i}
0\longrightarrow \hM\stackrel{\varpi^n}{\longrightarrow} \hM\longrightarrow M_n\longrightarrow 0
\end{equation}
est $\alpha$-exacte~; d'où la proposition (i). Le morphisme
\begin{equation}
M^\rc=\underset{\longleftarrow}{\lim} \ \hM/\varpi^n\hM \rightarrow \hM,
\end{equation}
déduit par passage à la limite projective, est donc un $\alpha$-isomorphisme d'après \ref{alpha10}.
Il résulte de (i) et \ref{mptf12} que le $\hR$-module $M^\rc$ est de type $\alpha$-fini~; d'où la proposition (ii).

Supposons enfin le corps des fractions $K$ de $R$ algébriquement clos. Il résulte du lemme de Krasner que le corps des fractions $\hK$ de $\hR$ 
est aussi algébriquement clos. Par suite, $\varphi\colon \hR\rightarrow \hR$ est un isomorphisme~; 
il en est alors de même de $\bvvarphi\colon \bvR\rightarrow \bvR^{(p)}$ \eqref{mptf10a} et donc de $\uphi\colon \hM\rightarrow \hM$ \eqref{mptf11b}. 
On note encore $\varphi$ l'endomorphisme de Frobenius de $\hK$. 
D'après (ii) et (\cite{katz} 4.1.1), il existe un entier $d\geq 0$ et un isomorphisme $\hK$-linéaire
\begin{equation}
f\colon \hM\otimes_\hR\hK\stackrel{\sim}{\rightarrow}\hK^d
\end{equation}
compatible à $\uphi$ et $\varphi$. Notons $N$ l'image canonique de $\hM$ dans $\hM\otimes_\hR\hK$.
Il résulte de la suite $\alpha$-exacte \eqref{mptf11i} que le morphisme canonique $\hM\rightarrow N$ est un $\alpha$-isomorphisme.
Par ailleurs, on a $(\uphi\otimes \varphi)(N)=N$. Comme $N$ est de type $\alpha$-fini d'après (ii), 
il existe $\gamma\in \fm$ tel que l'on ait $\gamma \hR^d\subset f(N)\subset \gamma^{-1} \hR^d$. Appliquant les puissances négatives de $\varphi$ à ces inclusions, 
on en déduit que $f(\fm\otimes_R N)\subset \hR^d$ et que le morphisme 
\begin{equation}
\fm\otimes_R N\rightarrow \hR^d
\end{equation}
induit par $f$ est un $\alpha$-isomorphisme~; d'où la proposition (iii).

\section{Modules dérivé-complets}\label{mdc}

Nous rappelons dans cette section la notion de modules {\em dérivé-complets} dits aussi {\em faiblement complets} 
et quelques-unes de leurs propriétés, suivant (\cite{jannsen} §4) et (\cite{sp} \href{https://stacks.math.columbia.edu/tag/091N}{Tag 091N}). Nous nous limitons au cadre nécessaire aux besoins
de ce travail. 

\subsection{}\label{mdc1}
Dans cette section, $R$ désigne un anneau et $\pi$ un élément de $R$. Pour tout $R$-module $M$, on note $(M,\pi)$ le système projectif 
indexé par l'ensemble $\mZ_{\geq 0}$, dont toutes les composantes sont égales à $M$ et dont les morphismes de transition sont 
les multiplications par les puissances de $\pi$. On désigne par $\pi\DDiv(M)$ le sous-$R$-module $\pi$-divisible maximal de $M$.  
On note $\hR$ (resp. $\hM$) le séparé complété de $R$ (resp. $M$) pour la topologie $\pi$-adique,
\[
\hR=\underset{\underset{n\geq 0}{\longleftarrow}}{\lim}\ R/\pi^n R \ \ \ {\rm et}\ \ \
\hM=\underset{\underset{n\geq 0}{\longleftarrow}}{\lim}\ M/\pi^n M.
\]
Pour tout entier $n\geq 0$, on désigne par $M[\pi^n]$ le noyau de la multiplication par $\pi^n$ dans $M$. 
On note $\rT_\pi(M)$ le $\pi$-module de Tate de $M$, c'est-à-dire le $\hR$-module limite projective du système projectif 
$(M[\pi^n])_{n\geq 0}$ où les morphismes de transition sont induits par la multiplication par $\pi$ dans $M$;
\begin{equation}\label{mdc1a}
\rT_\pi(M)=\underset{\underset{n\geq 0}{\longleftarrow}}{\lim}\ M[\pi^n].
\end{equation}

La suite exacte canonique $0\rightarrow \pi^n M\rightarrow M\rightarrow M/\pi^nM\rightarrow 0$ induit par passage à la limite projective  
la suite exacte
\begin{equation}\label{mdc1b}
0\rightarrow \bigcap_{n\geq 0}\pi^nM\rightarrow M\rightarrow \hM\rightarrow 
\rR^1\underset{\underset{n\geq 0}{\longleftarrow}}{\lim}\ \pi^nM \rightarrow 0.
\end{equation}

Par ailleurs, pour tout entier $n\geq 0$, on a un diagramme commutatif 
\[
\xymatrix{
0\ar[r]&M[\pi^{n+1}]\ar[r]\ar[d]_{\pi}&M\ar[r]^{\pi^{n+1}}\ar[d]^{\pi}&{\pi^{n+1}M}\ar[r]\ar@{^(->}[d]&0\\
0\ar[r]&M[\pi^{n}]\ar[r]&M\ar[r]^{\pi^{n}}&{\pi^{n}M}\ar[r]&0}
\]
On en déduit par passage à la limite projective la suite exacte longue 
\begin{eqnarray}\label{mdc1c}
\lefteqn
{0\rightarrow \rT_\pi(M)\rightarrow \underset{\longleftarrow}{\lim}\ (M,\pi)\rightarrow \bigcap_{n\geq 0}\pi^nM\rightarrow}\\
&&
\rR^1\underset{\underset{n\geq 0}{\longleftarrow}}{\lim}\ M[\pi^n]\rightarrow \rR^1\underset{\longleftarrow}{\lim}\ (M,\pi) \rightarrow 
\rR^1\underset{\underset{n\geq 0}{\longleftarrow}}{\lim}\ \pi^nM\rightarrow 0.\nonumber
\end{eqnarray}
On vérifie aussitôt que l'image du morphisme $\underset{\longleftarrow}{\lim}\ (M,\pi)\rightarrow \bigcap_{n\geq 0}\pi^nM$
est le sous-$R$-module $\pi$-divisible maximal $\pi\DDiv(M)$ de $M$.

\begin{lem}[\cite{jannsen} 4.3]\label{mdc2} 
Pour tout $R$-module $M$, on a $\pi\DDiv(\rR^1\underset{\underset{n\geq 0}{\longleftarrow}}{\lim}\ M[\pi^n])=0$. 
\end{lem}

En effet, on a une suite exacte canonique 
\begin{equation}\label{mdc2a}
0\rightarrow \underset{\underset{n\geq 0}{\longleftarrow}}{\lim}\ M[\pi^n] \rightarrow \prod_{n\geq 0}M[\pi^n]\stackrel{h}{\rightarrow} \prod_{n\geq 0}M[\pi^n]
\rightarrow \rR^1\underset{\underset{n\geq 0}{\longleftarrow}}{\lim}\ M[\pi^n]\rightarrow 0,
\end{equation}
où $h=\id-(d_n)$ et $d_n\colon M[\pi^{n+1}]\rightarrow M[\pi^n]$ est le morphisme induit par la multiplication par $\pi$ (\cite{jannsen} (1.4)). 
Soient $(a^{(m)})_{m\geq 0}$ et $(u^{(m)})_{m\geq 0}$ deux suites d'éléments de $\prod_{n\geq 0}M[\pi^n]$ telles que 
\begin{equation}
a^{(m)}- \pi a^{(m+1)}=h(u^{(m)}).
\end{equation}
Pour tout $m\geq 0$, posons $a^{(m)}=(a^{(m)}_n)_{n\geq 0}$ et $u^{(m)}=(u^{(m)}_n)_{n\geq 0}$. Pour tous $m,n\geq 0$, on a 
\begin{equation}
a_n^{(m)}-\pi a_n^{(m+1)}=u_n^{(m)}-d_n(u_{n+1}^{(m)}).
\end{equation}
Par suite, pour tout $n\geq 0$, on a
\begin{equation}
a_n^{(0)}=\sum_{m\geq 0} \pi^mu_n^{(m)}-d_n(\sum_{m\geq 0} \pi^mu_{n+1}^{(m)}),
\end{equation}
les sommes étant en fait finies car $\pi^n u_n^{(m)}=0$, et en particulier, on a $a^{(0)}\in \im(h)$. 

\begin{lem}[\cite{jannsen} 4.4] \label{mdc3} 
Pour tout $R$-module $M$, les propriétés suivantes sont équivalentes:
\begin{itemize}
\item[{\rm (i)}] $M$ est complet et séparé pour la topologie $\pi$-adique;
\item[{\rm (ii)}] $M$ est isomorphe à la limite projective d'un système projectif $(M_n)_{n\geq 0}$ de $R$-modules 
annulés par des puissances de $\pi$;
\item[{\rm (iii)}] tous les $R$-modules qui apparaissent dans la suite \eqref{mdc1c} sont nuls.
\end{itemize}
\end{lem}

On a  clairement (i)$\Rightarrow$(ii). Montrons (ii)$\Rightarrow$(i). Supposons (ii) satisfait. 
Pour tout entier $n\geq 0$, notons  $e(n)$ le plus petit entier $e\geq 0$ tel que $\pi^e M_n=0$ et posons 
\begin{equation}
m(n)=\max(m(n-1),e(n)),
\end{equation} 
où $m(-1)=0$. On a $\bigcap_{n\geq 0}\pi^nM=0$ car $M\subset \prod_{n\geq 0} M_n$. Soit $(a^{(m)})_{m\geq 0}$ une suite d'éléments de $M$ 
et pour tout $m\geq 0$, posons $a^{(m)}=(a_n^{(m)})\in \prod_{n\geq 0} M_n$. Supposons que pour tout $m\geq 0$, on ait   
\begin{equation}
a^{(m+1)}-a^{(m)}\in \pi^m M.
\end{equation}
Par suite, pour tous $m,n\geq 0$, on a $a_n^{(m+1)}-a_n^{(m)}\in \pi^m M_n$. Il s'ensuit que pour tout $n\geq 0$ et tout $m\geq m(n)$, on a 
$a_n^{(m)}=a_n^{(m(n))}$. Posons $a=(a_n^{(m(n))})\in \prod_{n\geq 0} M_n$. Pour tout $n\geq 0$, notant $d_n\colon M_{n+1}\rightarrow M_n$
le morphisme canonique, on a 
\begin{equation}
d_n(a_{n+1}^{(m(n+1))})=a_{n}^{(m(n+1))}=a_n^{(m(n))},
\end{equation}
de sorte que $a\in M$. Pour tout $m,n\geq 0$, on a $a_n^{(m(n))}-a_n^{(m)}\in \pi^mM_n$. 
Par suite, $a$ est la limite de la suite $(a^{(m)})_{m\geq 0}$
dans $M$ pour la topologie $\pi$-adique. 

Il est clair que (iii)$\Rightarrow$(i) \eqref{mdc1b}. Montrons (i)$\Rightarrow$(iii). Supposons (i) satisfait. 
On a alors $\bigcap_{n\geq 0}\pi^nM=0$ et  $\rR^1\underset{\underset{n\geq 0}{\longleftarrow}}{\lim}\ \pi^nM=0$ \eqref{mdc1b}.
On en déduit un isomorphisme 
\begin{equation}
\rR^1\underset{\underset{n\geq 0}{\longleftarrow}}{\lim}\ M[\pi^n]\stackrel{\sim}{\rightarrow} \rR^1\underset{\longleftarrow}{\lim}\ (M,\pi).
\end{equation}
Le but de cet isomorphisme étant $\pi$-divisible, on en déduit qu'il est nul en vertu de \ref{mdc2}. Par ailleurs, comme 
$\bigcap_{n\geq 0}\pi^nM=0$, on a $\underset{\longleftarrow}{\lim}\ (M,\pi)=0$, d'où (iii).

\begin{defi}[\cite{jannsen} 4.6]\label{mdc4}
On dit qu'un $R$-module $M$ est {\em dérivé-$\pi$-complet} (ou {\em faiblement $\pi$-complet}) 
si les $R$-modules $\underset{\longleftarrow}{\lim}\ (M,\pi)$ et $\rR^1\underset{\longleftarrow}{\lim}\ (M,\pi)$ sont nuls.
\end{defi}

\begin{prop}[\cite{jannsen} 4.8]\label{mdc5}
\
\begin{itemize}
\item[{\rm (i)}] Si deux modules d'une suite exacte de $R$-modules $0\rightarrow M\rightarrow N\rightarrow P\rightarrow 0$ sont dérivé-$\pi$-complets,
il en est de même du troisième.
\item[{\rm (ii)}] Les $R$-modules dérivé-$\pi$-complets forment une sous-catégorie abélienne pleine de la catégorie des $R$-modules.
\item[{\rm (iii)}]  Pour tout système projectif de $R$-modules $(M_n)_{n\geq 0}$ tel que $M_n$ soit annulé par une puissance de $\pi$, 
le $R$-module $\rR^1\underset{\underset{n\geq 0}{\longleftarrow}}{\lim}\ M_n$ est dérivé-$\pi$-complet.
\end{itemize}
\end{prop}

(i) Cela résulte aussitôt de la définition en considérant la suite exacte de systèmes projectifs $0\rightarrow (M,\pi) 
\rightarrow (N,\pi)\rightarrow (P,\pi)\rightarrow 0$.

(ii) Pour tout morphisme de $R$-modules dérivé-$\pi$-complets $f\colon M\rightarrow N$, les suites exactes 
\begin{eqnarray*}
0\rightarrow \ker(f)\rightarrow M\rightarrow \im(f)\rightarrow 0,\\
0\rightarrow \im(f)\rightarrow N\rightarrow \coker(f)\rightarrow 0,
\end{eqnarray*} 
montrent que $\underset{\longleftarrow}{\lim}\ (\im(f),\pi)$ et $\rR^1\underset{\longleftarrow}{\lim}\ (\coker(f),\pi)$ sont nuls. La proposition résulte alors de (i). 

(iii) Cela résulte de (ii) compte tenu de la suite exacte canonique (\cite{jannsen} (1.4))
\begin{equation}
0\rightarrow \underset{\underset{n\geq 0}{\longleftarrow}}{\lim}\ M_n \rightarrow \prod_{n\geq 0}M_n\rightarrow \prod_{n\geq 0}M_n
\rightarrow \rR^1\underset{\underset{n\geq 0}{\longleftarrow}}{\lim}\ M_n\rightarrow 0,
\end{equation}
et du fait que les modules $\underset{\underset{n\geq 0}{\longleftarrow}}{\lim}\ M_n$ et $\prod_{n\geq 0}M_n$ sont dérivé-$\pi$-complets en vertu de
\ref{mdc3}.

\begin{cor}\label{mdc6}
Soient $X$ un $\mU$-topos, $X^{\mN^\circ}$ le topos des systèmes 
projectifs d'objets de $X$ indexés par $\mN$ {\rm (\cite{agt} III.7.7)},
$A=(A_n)$ une $R$-algèbre de $X^{\mN^\circ}$, $M=(M_n)$ un $A$-module de $X^{\mN^\circ}$. 
Supposons que pour tout $n\in \mN$, $M_n$ soit annulé par une puissance de $\pi$. 
Alors, le $R$-module $\rH^q(X^{\mN^\circ},M)$ est dérivé-$\pi$-complet pour tout $q\geq 0$.
\end{cor}

D'après (\cite{agt} III.7.11), on a une suite exacte canonique et fonctorielle
\begin{equation}
0\rightarrow \rR^1\underset{\underset{n\in \mN^\circ}{\longleftarrow}}{\lim}\ \rH^{q-1}(X,M_n)\rightarrow 
\rH^q(X^{\mN^\circ},M)\rightarrow \underset{\underset{n\in \mN^\circ}{\longleftarrow}}{\lim}\ \rH^q(X,M_n)\rightarrow 0,
\end{equation}
où l'on a posé $\rH^{-1}(X,M_n)=0$ pour tout $n\in \mN$. La proposition résulte alors de \ref{mdc3} et \ref{mdc5}.

\begin{lem}\label{mdc7}
Tout $R$-module dérivé-$\pi$-complet $M$ tel que $M=\pi M$, est nul.
\end{lem}

En effet, la projection canonique $\underset{\longleftarrow}{\lim}\ (M,\pi)\rightarrow M$ sur le premier module du système projectif $(M,\pi)$ est surjective.

\begin{lem}\label{mdc8}
Supposons que $R$ soit un anneau de valuation non-discrète de hauteur $1$
et que $\pi$ soit un élément de l'idéal maximal $\fm$ de $R$. 
Soit $M$ un $R$-module dérivé-$\pi$-complet tel que $M/\pi M$ soit $\alpha$-nul. Alors, $M$ est $\alpha$-nul.
\end{lem}

Notons $v\colon R\rightarrow \mR\cup \{\infty\}$ la valuation de $R$. Soient $x\in M$, $\gamma\in \fm-\{0\}$, 
$(\gamma_n)_{n\geq 1}$ une suite d'éléments de $\fm$ telle que $\sum_{n\geq 1}v(\gamma_n)\leq v(\gamma)$. Posons $\gamma_0=1$. 
On peut alors construire une suite $(x_n)_{n\geq 0}$ d'éléments de $M$ telle que $x_0=x$ et que pour tout $n\geq 0$, on ait 
$\pi x_{n+1}=\gamma_{n+1}x_n$. Posant $y_n= \frac{\gamma}{\gamma_0\dots \gamma_n}x_n$,
on a $\pi y_{n+1}=y_n$. Par suite, $\gamma x$ appartient à l'image de la projection canonique 
$u_0\colon \underset{\longleftarrow}{\lim}\ (M,\pi)\rightarrow M$ sur le premier module du système projectif $(M,\pi)$. 
Il s'ensuit que $u_0$ est $\alpha$-surjectif et donc $M$ est $\alpha$-nul.

\chapter{Topos de Faltings}\label{toposfaltings}

\section{Produits orientés de topos}

La notion de produit orienté de topos est due à Deligne. Elle a été étudiée par Gabber, Illusie et Orgogozo (\cite{travaux-gabber} XI, \cite{orgogozo}). 
On renvoie à (\cite{agt} VI) pour la définition et les principales propriétés, rappelées et complétées dans cette section. 

\subsection{}\label{topfl1}
Soient $X$, $Y$, $S$ trois $\mU$-sites \eqref{notconv3}
dans lesquels les limites projectives finies sont représentables, 
\begin{equation}
f^+\colon S\rightarrow X,\ \ \  g^+\colon S\rightarrow Y
\end{equation}
deux foncteurs continus et exacts à gauche. On désigne par 
$\tX$, $\tY$ et $\tS$ les topos des faisceaux de $\mU$-ensembles sur $X$, $Y$ et $S$, respectivement, par 
\begin{equation}
f\colon \tX\rightarrow \tS \ \ \ {\rm et}\ \ \ g\colon \tY\rightarrow \tS
\end{equation}
les morphismes de topos définis par $f^+$ et $g^+$, respectivement (\cite{sga4} IV 4.9.2), et par   
$e_X$, $e_Y$ et $e_S$ des objets finaux de $X$, $Y$ et $S$, respectivement, qui
existent par hypothèse. 

On désigne par $C$ la catégorie des triplets 
\begin{equation}
(W, U\rightarrow f^+(W), V\rightarrow g^+(W)),
\end{equation}
où $W$ est un objet de $S$, 
$U\rightarrow f^+(W)$ est un morphisme de $X$ et $V\rightarrow g^+(W)$ est un morphisme 
de $Y$; un tel objet sera noté $(U\rightarrow W\leftarrow V)$. 
Soient $(U\rightarrow W\leftarrow V)$ et $(U'\rightarrow W'\leftarrow V')$
deux objets de $C$. Un morphisme de $(U'\rightarrow W'\leftarrow V')$ dans $(U\rightarrow W\leftarrow V)$
est la donnée de trois morphismes $U\rightarrow U'$, $V\rightarrow V'$ et $W\rightarrow W'$ de
$X$, $Y$ et $S$, respectivement, tels que les diagrammes 
\begin{equation}
\xymatrix{
U'\ar[r]\ar[d]&{f^+(W')}\ar[d]\\
{U}\ar[r]&{f^+(W)}}\ \ \ \ \ \ \ \ 
\xymatrix{
V'\ar[d]\ar[r]&{g^+(W')}\ar[d]\\
V\ar[r]&{g^+(W)}}
\end{equation}
soient commutatifs. 

Il résulte aussitôt de la définition et du fait que les foncteurs $f^+$ 
et $g^+$ sont exacts à gauche que les limites projectives finies dans $C$ sont représentables.

On munit $C$ de la topologie engendrée par les recouvrements 
\begin{equation}
\{(U_i\rightarrow W_i\leftarrow V_i)\rightarrow (U\rightarrow W\leftarrow V)\}_{i\in I}
\end{equation}
des trois types suivants~:
\begin{itemize}
\item[(a)] $V_i=V$, $W_i=W$ pour tout $i\in I$, et $(U_i\rightarrow U)_{i\in I}$ est une famille couvrante.
\item[(b)] $U_i=U$, $W_i=W$ pour tout $i\in I$, et $(V_i\rightarrow V)_{i\in I}$ est une famille couvrante.
\item[(c)] $I=\{'\}$, $U'=U$ et le morphisme $V'\rightarrow V\times_{g^+(W)}g^+(W')$ est un isomorphisme
(il n'y a aucune condition sur le morphisme $W'\rightarrow W$). 
\end{itemize}
On notera que chacune de ces familles est stable par changement de base. 
On désigne par $\tC$ le topos des faisceaux de $\mU$-ensembles sur $C$. 

\begin{lem}[\cite{agt} VI.3.2] \label{topfl2}
Sous les hypothèses de \ref{topfl1}, pour qu'un  
préfaisceau $F$ sur $C$ soit un faisceau, il faut et il suffit que les conditions suivantes soient remplies~:
\begin{itemize}
\item[{\rm (i)}] Pour tout famille couvrante $(Z_i\rightarrow Z)_{i\in I}$ de $C$ du type {\rm (a)} ou {\rm (b)}, 
la suite 
\begin{equation}
F(Z)\rightarrow \prod_{i\in I}F(Z_i)\rightrightarrows \prod_{(i,j)\in I\times J}F(Z_i\times_ZZ_j)
\end{equation}
est exacte. 
\item[{\rm (ii)}] Pour tout recouvrement $(U\rightarrow W'\leftarrow V')\rightarrow (U\rightarrow W\leftarrow V)$ 
du type {\rm (c)}, l'application 
\begin{equation}\label{topfl2a}
F(U\rightarrow W\leftarrow V) \rightarrow F(U\rightarrow W'\leftarrow V')
\end{equation}
est bijective.
\end{itemize}
\end{lem}

\subsection{}\label{topfl4}
Conservons les hypothèses de \ref{topfl1}. Les foncteurs 
\begin{eqnarray}
\rp_1^+\colon X&\rightarrow& C,\ \ \ U\mapsto (U\rightarrow e_Z\leftarrow e_Y),\label{topfl4a}\\
\rp_2^+\colon Y&\rightarrow& C,\ \ \ V\mapsto (e_X\rightarrow e_Z\leftarrow V),\label{topfl4b}
\end{eqnarray}
sont exacts à gauche et continus. Ils définissent donc deux morphismes de topos (\cite{sga4} IV 4.9.2)
\begin{eqnarray}
\rp_1\colon \tC&\rightarrow& \tX,\label{topfl4c}\\
\rp_2\colon \tC&\rightarrow& \tY.\label{topfl4d}
\end{eqnarray}
Par ailleurs, on a un $2$-morphisme 
\begin{equation}\label{topfl4e}
\tau\colon g \rp_2\rightarrow f \rp_1,
\end{equation}
donné par le morphisme de foncteurs $(g \rp_{2})_*\rightarrow (f \rp_{1})_*$ suivant~: 
pour tout faisceau $F$ sur $C$ et tout $W\in \ob(S)$, 
\begin{equation}\label{topfl4f}
g_*(\rp_{2*}(F))(W)\rightarrow f_*(\rp_{1*}(F))(W)
\end{equation}
est l'application composée 
\[
F(e_X\rightarrow e_Z\leftarrow g^+(W))\rightarrow F(f^+(W)\rightarrow W\leftarrow g^+(W))
\rightarrow F(f^+(W)\rightarrow e_Z\leftarrow e_Y)
\]
où la première flèche est le morphisme canonique et la seconde flèche est l'inverse de l'isomorphisme \eqref{topfl2a}.

D'après (\cite{agt} VI.3.7), le quadruplet $(\tC,\rp_1,\rp_2,\tau)$ est universel dans le sens suivant:
pour tout $\mU$-topos $T$ muni de deux morphismes de topos $a\colon T\rightarrow \tX$ et $b\colon T\rightarrow \tY$ 
et d'un $2$-morphisme $t\colon gb\rightarrow fa$, il existe un triplet 
\begin{equation}
(h\colon T\rightarrow \tC, \alpha\colon \rp_1h\stackrel{\sim}{\rightarrow}a,  
\beta\colon \rp_2h\stackrel{\sim}{\rightarrow}b),
\end{equation}
unique à isomorphisme unique près, formé d'un morphisme de topos $h$ et de deux isomorphismes 
de morphismes de topos $\alpha$ et $\beta$, tel que le diagramme 
\begin{equation}\label{topfl6a}
\xymatrix{
{g\rp_2h}\ar[r]^{\tau*h}\ar[d]_{g*\beta}&{f\rp_1h}\ar[d]^{f*\alpha}\\
{gb}\ar[r]^t&{fa}}
\end{equation}
soit commutatif. Par suite, le topos $\tC$ ne dépend que du couple de morphismes de topos
$(f,g)$, à équivalence canonique près (cf. \cite{agt} VI.3.9). On l'appelle {\em produit orienté}  
de $\tX$ et $\tY$ au-dessus de $\tS$ et on le note $\tX\gtimes_\tS\tY$.

\subsection{}\label{topfl5}
\'Etant donnés deux morphismes de $\mU$-topos $f\colon \tX\rightarrow \tS$ et $g\colon \tY\rightarrow \tS$, 
on considère $\tX$, $\tY$ et $\tS$ comme des sites munis de leurs topologies canoniques (\cite{sga4} II 2.5 et VI 1.2). 
On définit alors le produit orienté $\tX\gtimes_\tS\tY$ comme le topos des faisceaux de $\mU$-ensembles sur le site défini dans \ref{topfl1}
associé au couple de foncteurs $(f^*\colon \tS\rightarrow \tX, g^*\colon \tS\rightarrow \tY)$, lesquels vérifient clairement les 
hypothèses requises (\cite{sga4} II 4.4 et III 1.6). 

\subsection{}\label{topfl10}
Considérons un diagramme de morphismes de topos 
\begin{equation}\label{topfl10a}
\xymatrix{
{\tX'}\ar[d]_u\ar[r]^-(0.5){f'}&{\tS'}\ar[d]^w&{\tY'}\ar[d]^v\ar[l]_-(0.5){g'}\\
{\tX}\ar[r]^-(0.5){f}&{\tS}&{\tY}\ar[l]_-(0.5){g}}
\end{equation}
et deux $2$-morphismes 
\begin{equation}\label{topfl10b}
a\colon wf'\rightarrow fu\ \ \ {\rm et}\ \ \
b\colon gv\rightarrow wg'.
\end{equation}
Notons $\rp_1\colon \tX\gtimes_{\tS}\tY\rightarrow \tX$,  $\rp_2\colon \tX\gtimes_\tS\tY\rightarrow \tY$,
$\rp'_1\colon \tX'\gtimes_{\tS'}\tY'\rightarrow \tX'$ et $\rp'_2\colon \tX'\gtimes_{\tS'}\tY'\rightarrow \tY'$ les projections
canoniques \eqref{topfl4}. 
Les morphismes de topos $u\rp'_1$ et $v\rp'_2$ et le $2$-morphisme composé
\begin{equation}\label{topfl10c}
t\colon \xymatrix{
{gv\rp'_2}\ar[r]^-(0.4){b*\rp'_2}&{wg'\rp'_2}\ar[r]^-(0.5){w*\tau'}&{wf'\rp'_1}\ar[r]^-(0.4){a*\rp'_1}&{fu\rp'_1}},
\end{equation}
définissent un morphisme de topos 
\begin{equation}\label{topfl10d}
h\colon \tX'\gtimes_{\tS'}\tY'\rightarrow \tX\gtimes_{\tS}\tY,
\end{equation}
que l'on notera aussi $u\gtimes_wv$, et des $2$-isomorphismes $\alpha\colon \rp_1 h\stackrel{\sim}{\rightarrow} u\rp'_1$ et 
$\beta\colon \rp_2 h\stackrel{\sim}{\rightarrow} v\rp'_2$ rendant commutatif les carrés du diagramme 
\begin{equation}
\xymatrix{
{\tX'}\ar[d]_u&{\tX'\times_{\tS'}\tY'}\ar[l]_-(0.5){\rp'_1}\ar[r]^-(0.5){\rp'_2}\ar[d]^{h}&{\tY'}\ar[d]^v\\
{\tX}&{\tX\times_{\tS}\tY}\ar[l]_-(0.5){\rp_1}\ar[r]^-(0.5){\rp_2}&{\tY}}
\end{equation}
De plus, le diagramme 
\begin{equation}\label{topfl10e}
\xymatrix{
{g\rp_2h}\ar[r]^{\tau*h}\ar[d]_{g*\beta}&{f\rp_1h}\ar[d]^{f*\alpha}\\
{gv\rp'_2}\ar[r]^t&{fu\rp'_1}}
\end{equation}
est commutatif. 

\subsection{}\label{topfl11}
Soient $X$, $Y$, $S$, $X'$, $Y'$, $S'$ des $\mU$-sites dans lesquels les limites projectives finies sont représentables,  
$\tX$, $\tY$, $\tS$, $\tX'$, $\tY'$, $\tS'$ les topos de faisceaux de $\mU$-ensembles associés, 
\begin{equation}\label{topfl11a}
\xymatrix{
{X}\ar[d]_{u^+}&{S}\ar[l]_-(0.5){f^+}\ar[d]^{w^+}\ar[r]^{g^+}&{Y}\ar[d]^{v^+}\\
{X'}&{S'}\ar[l]_-(0.5){f'^+}\ar[r]^-(0.5){g'^+}&Y}
\end{equation}
un diagramme de foncteurs continus et exacts à gauche,
\begin{equation}\label{topfl11b}
a^+\colon u^+\circ f^+\stackrel{\sim}{\rightarrow} f'^+\circ w^+\ \ \ {\rm et}\ \ \
b^+\colon v^+\circ g^+ \stackrel{\sim}{\rightarrow} g'^+\circ w^+
\end{equation}
deux isomorphismes de foncteurs.
On désigne par $C$ (resp. $C'$) le site  associé au couple de foncteurs $(f^+, g^+)$ (resp. $(f'^+, g'^+)$) introduit dans \ref{topfl1},
et on considère le foncteur
\begin{equation}
\Phi^+\colon C\rightarrow C', \ \ \ (U\rightarrow W\leftarrow V)\mapsto (u^+(U)\rightarrow w^+(W)\leftarrow v^+(V)),
\end{equation}
où le morphisme  $u^+(U)\rightarrow f'^+(w^+(W))$ est le composé $u^+(U)\rightarrow u^+(f^+W))\stackrel{\sim}{\rightarrow} f'^+(w^+(W))$,
la second flèche étant induite par $a^+$, 
et le morphisme $v^+(V)\rightarrow g'^+(w^+(W))$ est le composé $v^+(V)\rightarrow v^+(g^+(W))\stackrel{\sim}{\rightarrow} g'^+(w^+(W))$,
la second flèche étant induite par $b^+$. Celui-ci est clairement exact à gauche. Il transforme les familles couvrantes 
de $C$ de type (a) (resp. (b), resp. (c)) en familles couvrantes de $C'$ du même type. Pour tout faisceau $F$ sur $C'$,
$F\circ \Phi^+$ est donc un faisceau sur $C$ en vertu de \ref{topfl2}. Par suite, le foncteur $\Phi$ est continu. Il définit donc un 
morphisme de topos
\begin{equation}
\Phi\colon \tX'\gtimes_{\tS'}\tY'\rightarrow \tX\times_\tS\tY. 
\end{equation} 

Considérons le diagramme de morphismes de topos 
\begin{equation}\label{topfl11c}
\xymatrix{
{\tX'}\ar[d]_u\ar[r]^-(0.5){f'}&{\tS'}\ar[d]^w&{\tY'}\ar[d]^v\ar[l]_-(0.5){g'}\\
{\tX}\ar[r]^-(0.5){f}&{\tS}&{\tY}\ar[l]_-(0.5){g}}
\end{equation}
déduit de \eqref{topfl11a}, dont les carrés sont commutatifs aux $2$-isomorphismes 
\begin{equation}
a\colon wf'\stackrel{\sim}{\rightarrow} fu\ \ \ {\rm et}\ \ \
b\colon wg'\stackrel{\sim}{\rightarrow} gv,
\end{equation}
induits par $a^+$ et $b^+$, près. Il résulte de (\cite{agt} VI.3.6) que le morphisme $\Phi$ est celui associé 
dans \eqref{topfl10d} au diagramme \eqref{topfl11c}
et aux $2$-morphismes $a$ et $b^{-1}$ (cf. la preuve de \cite{agt} VI.3.7).

\subsection{}\label{topfl12}
Reprenons les hypothèses de \ref{topfl1}, soient, de plus, $(U\rightarrow W\leftarrow V)$ un objet de $C$,
$F$ (resp. $G$, resp. $H$) le faisceau associé à $U$ (resp. $V$, resp. $W$). 
Rappelons (\cite{sga4} IV 5.1) que la catégorie $\tX_{/F}$ est un $\mU$-topos, dit topos induit 
sur $F$ par $\tX$, et qu'on a un morphisme canonique $j_F\colon \tX_{/F}\rightarrow \tX$,
dit morphisme de localisation de $\tX$ en $F$. En fait, $\tX_{/F}$ est canoniquement équivalent au topos des faisceaux de $\mU$-ensembles 
sur la catégorie $X_{/U}$ munie de la topologie induite par celle de $X$ via le foncteur canonique $X_{/U}\rightarrow X$ (\cite{sga4} III 5.4). 
De même, on a des morphismes de localisation
$j_G\colon \tY_{/G}\rightarrow \tY$ et $j_H\colon \tS_{/H}\rightarrow \tS$. 
Le foncteur 
\begin{equation}\label{topfl12a}
f'^+\colon S_{/W}\rightarrow X_{/U},\ \ \ W'\mapsto f^+(W')\times_{f^+(W)}U
\end{equation}
est clairement exact à gauche et continu. Le morphisme de topos qui lui est associé est le composé 
\begin{equation}\label{topfl12b}
\xymatrix{
f'\colon {\tX_{/F}}\ar[r]&{\tX_{/f^*(H)}}\ar[r]^-(0.5){f_{/H}}&{\tS_{/H}}}
\end{equation}
où la première flèche est le morphisme de localisation de $\tX_{/f^*(H)}$
en l'objet $F\rightarrow f^*(H)$ déduit de $U\rightarrow f^+(W)$ 
(\cite{sga4} IV 5.5), et la seconde flèche est le morphisme induit par $f$ (\cite{sga4} IV 5.10).
On définit de même le foncteur $g'^+\colon S_{/W}\rightarrow Y_{/V}$
et le morphisme $g'\colon \tY_{/G}\rightarrow \tS_{/H}$. Les carrés du diagramme 
\begin{equation}\label{topfl12c}
\xymatrix{
{\tX_{/F}}\ar[d]_{j_F}\ar[r]^{f'}&{\tS_{/H}}\ar[d]^{j_H}&{\tY_{/G}}\ar[l]_{g'}\ar[d]^{j_G}\\
{\tX}\ar[r]^f&{\tS}&{\tY}\ar[l]_g}
\end{equation}
sont commutatifs à isomorphismes canoniques près. On désigne par $\tX_{/F}\gtimes_{\tS_{/H}}\tY_{/G}$
le produit orienté de $\tX_{/F}$ et $\tY_{/G}$ au-dessus de $\tS_{/H}$.
Le diagramme \eqref{topfl12c} induit alors un morphisme canonique \eqref{topfl10d}
\begin{equation}\label{topfl12d}
j_{F}\gtimes_{j_{H}}j_{G}\colon \tX_{/F}\gtimes_{\tS_{/H}}\tY_{/G}\rightarrow \tX\gtimes_{\tS}\tY.
\end{equation}

On note $F\gtimes_HG$
le faisceau de $\tX\gtimes_{\tS}\tY$ associé à $(U\rightarrow W\leftarrow V)$. 
D'après (\cite{agt} VI.3.6), on a un isomorphisme canonique
\begin{equation}\label{topfl12e}
F\gtimes_HG\stackrel{\sim}{\rightarrow}\rp_1^*(F)\times_{(g\rp_2)^*(H)} \rp_2^*(G),
\end{equation}
où le morphisme $\rp_1^*(F)\rightarrow (g\rp_2)^*(H)$ est le composé du morphisme $\rp_1^*(F)\rightarrow (f\rp_1)^*(H)$ 
et du morphisme $\tau\colon (f\rp_1)^*(H)\rightarrow (g\rp_2)^*(H)$ \eqref{topfl4e}. On note 
\begin{equation}\label{topfl12f}
j_{F\gtimes_HG}\colon (\tX\gtimes_\tS\tY)_{/(F\gtimes_HG)} \rightarrow \tX\gtimes_\tS\tY
\end{equation}
le morphisme de localisation de $\tX\gtimes_\tS\tY$ en $F\gtimes_HG$. 
En vertu de (\cite{agt} VI.3.15), il existe une équivalence de topos canonique  
\begin{equation}\label{topfl12g}
m\colon (\tX\gtimes_\tS\tY)_{/(F\gtimes_HG)}\stackrel{\sim}{\rightarrow} \tX_{/F}\gtimes_{\tS_{/H}}\tY_{/G}
\end{equation}
et un isomorphisme canonique 
\begin{equation}\label{topfl12h}
(j_{F}\gtimes_{j_H}j_G)\circ m\stackrel{\sim}{\rightarrow} j_{F\gtimes_HG}.
\end{equation}

\begin{prop}\label{topfl15}
Reprenons les hypothèses de \ref{topfl1}, supposons de plus que tous les objets de $X$, $Y$ et $S$ soient quasi-compacts. Alors, 
\begin{itemize}
\item[{\rm (i)}] Les topos $\tX$, $\tY$ et $\tS$ sont cohérents et les morphismes $f\colon \tX\rightarrow \tS$ 
et $g\colon \tY\rightarrow \tS$ sont cohérents.
\item[{\rm (ii)}] La topologie de $C$ est engendrée par la prétopologie définie pour chaque objet $A$ de $C$
par la donnée de l'ensemble $\Cov(A)$ des familles de morphismes $(A_i\rightarrow A)_{i\in I}$ obtenues par composition 
d'un nombre fini de familles de type (c) et de familles finies de type (a) et (b); en particulier, l'ensemble $I$ est fini.
\item[{\rm (iii)}] Tout objet de $C$ est cohérent. 
\end{itemize}
\end{prop}

(i) Cela résulte de (\cite{sga4} VI 2.4.5 et 3.3).

(ii) En effet, les familles couvrantes de type (c) et les familles couvrantes finies de type (a) et (b) de $C$ 
étant stables par changement de base, les ensembles $\Cov(A)$ pour $A\in \ob(C)$ vérifient les axiomes d'une prétopologie 
(\cite{sga4} II 1.3). Chaque élément de $\Cov(A)$ étant clairement couvrant pour la topologie définie dans \ref{topfl1},
la topologie engendrée par cette prétopologie est moins fine que la topologie définie dans \ref{topfl1}.  
Comme tout recouvrement de $C$ de type (a) (resp. (b)) est raffiné par un recouvrement fini de même type, 
la topologie engendrée par cette prétopologie est plus fine que la topologie définie dans \ref{topfl1}. 
Les deux topologies sont donc égales. 

(iii) En effet, tous les recouvrements de la prétopologie définie dans (ii) sont finis. 
Tout objet de $C$ est donc quasi-compact, et par suite cohérent puisque $C$ est stable par produits fibrés.

\begin{cor}\label{topfl16}
Soient $\tX$, $\tY$ et $\tS$ trois topos cohérents, $f\colon \tX\rightarrow \tS$ et $g\colon \tY\rightarrow \tS$ deux morphismes cohérents. 
Alors, 
\begin{itemize}
\item[{\rm (i)}] Le topos $\tX\gtimes_\tS\tY$ est cohérent~; en particulier, il a suffisamment de points.
\item[{\rm (ii)}] Les morphismes $\rp_1$ \eqref{topfl4c} et $\rp_2$ \eqref{topfl4d} sont cohérents. 
\end{itemize}
\end{cor}

En effet, le topos $\tX$ (resp. $\tY$, resp. $\tS$) admet une sous-catégorie pleine génératrice $X$ (resp. $Y$, resp. $S$)
formée d'objets cohérents, qui est stable par limites projectives finies (\cite{sga4} VI 2.4.5). 
Comme les morphismes $f$ et $g$ sont cohérents, on peut supposer que $f^*(S)\subset X$ et $g^*(S)\subset Y$ (\cite{sga4} VI 3.2). 
La proposition résulte alors de \ref{topfl15}(iii) et (\cite{sga4} VI 2.4.5, 3.3 et 9.0).

\subsection{}\label{topfl17}
Reprenons les hypothèses de \ref{topfl1}.  
Il résulte aussitôt de la propriété universelle de $\tX\gtimes_\tS\tY$ que la donnée d'un point de ce topos est équivalente à 
la donnée d'une paire de points $x\colon \Ens\rightarrow \tX$ et $y\colon \Ens\rightarrow \tY$ 
et d'un $2$-morphisme $u\colon gy\rightarrow fx$ \eqref{notconv3}. 
Un tel point sera noté abusivement $(y\rightsquigarrow x)$ ou encore $(u\colon y\rightsquigarrow x)$.  
Pour tous $F\in \ob(\tX)$ et $G\in \ob(\tY)$, on a des isomorphismes canoniques fonctoriels 
\begin{eqnarray}
(\rp_1^*F)_{(y\rightsquigarrow x)} &\stackrel{\sim}{\rightarrow}& F_{x}, \label{topfl17a}\\
(\rp_2^*G)_{(y\rightsquigarrow x)} &\stackrel{\sim}{\rightarrow}& G_{y}. \label{topfl17b}
\end{eqnarray}
D'après  (\cite{agt} VI.3.6), pour tout objet $(U\rightarrow W\leftarrow V)$ de $C$, on a un isomorphisme canonique fonctoriel
\begin{equation}\label{topfl17c}
(U\rightarrow W\leftarrow V)^a_{(y\rightsquigarrow x)} \stackrel{\sim}{\rightarrow} U^a_x\times_{W^a_{gy}}V^a_y,
\end{equation}
où l'exposant $^a$ désigne les faisceaux associés, 
l'application $V^a_y\rightarrow W^a_{gy}$ est induite par $V^a\rightarrow g^*(W^a)$, et 
l'application $U^a_{x}\rightarrow W^a_{gy}$ est composée de l'application $U^a_{x}\rightarrow W^a_{fx}$ induite  par $U^a\rightarrow f^*(W^a)$
et du morphisme de spécialisation $W^a_{fx}\rightarrow W^a_{gy}$ défini par $u$. 

La catégorie $\cV_{(y\rightsquigarrow x)}$ des voisinages de $(y\rightsquigarrow x)$ dans $C$ (\cite{sga4} IV 6.8.2) est donc 
équivalente à la catégorie des triplets $((W\rightarrow U\leftarrow V),\xi,\zeta)$ formés d'un objet $(W\rightarrow U\leftarrow V)$
de $C$ et de deux éléments $\xi\in U^a_x$ et $\zeta\in V^a_y$ ayant même image dans $W^a_{gy}$ \eqref{topfl17c}.
Elle est cofiltrante, et pour tout préfaisceau 
$P$ sur $C$, on a un isomorphisme canonique fonctoriel (\cite{sga4} IV (6.8.4))
\begin{equation}\label{topfl17d}
P^a_{(y \rightsquigarrow x)}\stackrel{\sim}{\rightarrow} \underset{\underset{((U\rightarrow W\leftarrow V),\xi,\zeta)\in 
\cV^\circ_{(y \rightsquigarrow x)}}{\longrightarrow}}\lim\ P(U\rightarrow W\leftarrow V).
\end{equation}

\subsection{}\label{topfl18}
La construction suivante est due à Gabber (\cite{travaux-gabber} XI 2.2). 
Soit $(\tS,s)$ un topos local de centre $s\colon \Ens\rightarrow \tS$ (\cite{sga4} VI 8.4.6).
Notons $\varepsilon\colon \tS\rightarrow \Ens$ la projection canonique (\cite{sga4} IV 4.3). 
L'isomorphisme canonique $\varepsilon s\stackrel{\sim}{\rightarrow}\id_{\Ens}$ 
induit un morphisme de changement de base $\varepsilon_*\rightarrow s^*$ qui est un isomorphisme puisque $(\tS,s)$ est local. 
On en déduit un isomorphisme $\varepsilon^*\varepsilon_*\stackrel{\sim}{\rightarrow} \varepsilon^*s^*$. 
Celui-ci permet d'identifier la flèche d'adjonction 
$\varepsilon^*\varepsilon_*\rightarrow \id$ à un morphisme $(s\varepsilon)^*\rightarrow \id$. On obtient ainsi un $2$-morphisme
\begin{equation}\label{topfl18a}
c_s\colon \id\rightarrow s\varepsilon.
\end{equation}

Soient $f\colon \tX\rightarrow \tS$, $g\colon \tY\rightarrow \tS$ deux morphismes de topos, $x$ un point de $\tX$ tel que $s=fx$. 
Les morphismes $\id_{\tY}$ et $x\varepsilon g$ et le $2$-morphisme 
\begin{equation}\label{topfl18b}
c_sg\colon g\rightarrow f x \varepsilon g
\end{equation}
induit par \eqref{topfl18a} et la relation $s=fx$,
définissent un morphisme $\gamma\colon \tY\rightarrow \tX\gtimes_{\tS} \tY$ qui s'insère dans un diagramme non-commutatif 
\begin{equation}\label{topfl18c}
\xymatrix{
{\Ens}\ar[d]_{x}&{\tY}\ar[d]_{\gamma}\ar[l]_{\varepsilon g}\ar[rd]^\id&\\
{\tX}\ar[rd]_{f}&{\tX\gtimes_{\tS}\tY}\ar[l]_-(0.4){\rp_1}\ar[r]^-(0.4){\rp_2}&{\tY}\ar[ld]^{g}\\
&{\tS}&}
\end{equation}
En particulier, $\gamma$ est une section de $\rp_2$ telle que $\rp_1\gamma=x\varepsilon g$. 
On dit que $\gamma$ est la {\em section canonique} définie par le point $x$.

De la relation $\rp_2\gamma=\id_\tY$ on déduit un morphisme de changement de base 
\begin{equation}\label{topfl18d}
\rp_{2*}\rightarrow \gamma^*,
\end{equation}
composé de $\rp_{2*}\rightarrow \rp_{2*}\gamma_*\gamma^*\stackrel{\sim}{\rightarrow}\gamma^*$,
où la première flèche est induite par le morphisme d'adjonction $\id\rightarrow \gamma_*\gamma^*$,
et la seconde flèche par la relation $\rp_2\gamma=\id_\tY$.

\begin{prop}[Gabber, \cite{travaux-gabber} XI 2.3]\label{topfl19}
Soient $f\colon (\tX,x)\rightarrow (\tS,s)$ un morphisme local de topos locaux (i.e., $s=fx$), 
$g\colon \tY\rightarrow \tS$ un morphisme de topos, $y$ un point de $\tY$. 
Notons $\gamma\colon Y\rightarrow \tX\gtimes_\tS\tY$ la section canonique définie par le point $x$ \eqref{topfl18}. 
Alors, pour tout objet $F$ de $\tX\gtimes_\tS\tY$, l'application 
\begin{equation}\label{topfl19a}
(\rp_{2*}F)_y\rightarrow (\gamma^*F)_y
\end{equation}
déduite du morphisme de changement de base \eqref{topfl18d} est bijective. 
\end{prop}

En effet, notant $\varepsilon \colon \tS\rightarrow \Ens$ la projection canonique,
le composé $x\varepsilon g y\colon \Ens\rightarrow \tX$ est canoniquement isomorphe à $x$. 
Le $2$-morphisme $c_s\colon \id\rightarrow s\varepsilon$ \eqref{topfl18a}
induit donc un $2$-morphisme $u\colon gy\rightarrow fx$ définissant 
un point de $\tX\gtimes_\tS\tY$ \eqref{topfl17}, qui n'est autre que $\gamma y$. 
Considérons $\tX$, $\tY$ et $\tS$ comme des sites munis de leurs topologies canoniques, 
et notons $C$ le site défini dans \ref{topfl1} relativement 
au couple de foncteurs $(f^*\colon \tS\rightarrow \tX, g^*\colon \tS\rightarrow \tY)$ (cf. \ref{topfl5}). 
On désigne par $\cV_{\gamma y}$ la catégorie des voisinages de $\gamma y$ dans $C$ \eqref{topfl17}, 
et par $\cW_y$ la catégorie des voisinages de $y$ dans $\tY$ (\cite{sga4} IV 6.8), 
c'est-à-dire la catégorie des couples $(V,\zeta)$ formés d'un objet $V$ de $\tY$ et d'un élément $\zeta\in V_y$.

Soient $e_\tS$ (resp. $e_\tX$) un objet final de $\tS$ (resp. $\tX$), $\xi_e$ l'unique élément de $x^*(e_\tX)$. 
On a le foncteur
\begin{equation}\label{topfl19b}
\iota\colon \cW_y\rightarrow \cV_{\gamma y}, \ \ \ (V,\zeta)\mapsto ((e_\tX\rightarrow e_\tS\leftarrow V),\xi_e,\zeta).
\end{equation}
Les catégories $\cV_{\gamma y}$ et $\cW_y$ sont cofiltrantes (\cite{sga4} IV 6.8.2) et le foncteur 
$\iota$ est pleinement fidèle. On vérifie aussitôt que \eqref{topfl19a} est l'application
\begin{equation}
\underset{\underset{(V,\zeta)\in \cW^\circ_y}{\longrightarrow}}{\lim}F(e_\tX\rightarrow e_\tS\leftarrow V)\longrightarrow 
\underset{\underset{((U\rightarrow W\leftarrow V),\xi,\zeta)\in \cV^\circ_{\gamma y}}{\longrightarrow}}{\lim}F(U\rightarrow W\leftarrow V)
\end{equation}
induite par $\iota^\circ$. Il suffit donc de montrer que $\iota^\circ$ est cofinal. 
Soit $((U\rightarrow W\leftarrow V),\xi,\zeta)\in \ob(\cV_{\gamma y})$.
Comme le topos $(\tX,x)$ est local, on a un isomorphisme canonique fonctoriel $U_x\stackrel{\sim}{\rightarrow} \Gamma(\tX,U)$. 
On peut donc identifier $\xi\in U_x$ à une section que l'on note encore $\xi\in U(e_\tX)$.
Notons $\kappa$ l'image de $\xi$ par l'application $U_x\rightarrow W_s$ déduite du morphisme $U\rightarrow f^*(W)$. 
Le topos $(\tS,s)$ étant local, on peut identifier $\kappa$ à une section que l'on note encore $\kappa\in W(e_\tS)$. On obtient alors
un morphisme de $C$
\begin{equation}
\xymatrix{
e_\tX\ar[r]\ar[d]_{\xi}&e_\tS\ar[d]^\kappa&V\ar[l]\ar@{=}[d]\\
U\ar[r]&W&V\ar[l]}
\end{equation}
Celui-ci induit clairement un morphisme de $\cV_{\gamma y}$ 
\begin{equation}
\iota(V,\zeta)\rightarrow ((U\rightarrow W\leftarrow V),\xi,\zeta).
\end{equation}
On en déduit que $\iota^\circ$ est cofinal en vertu de (\cite{sga4} IV 8.1.3(c)); d'où la proposition.

\begin{cor}\label{topfl20}
Sous les hypothèses de \ref{topfl19}, si, de plus, $\tY$ a assez de points, alors
le morphisme de changement de base $\rp_{2*}\rightarrow \gamma^*$ \eqref{topfl18d}
est un isomorphisme~; en particulier, le foncteur $\rp_{2*}$ est exact. 
\end{cor}

\begin{cor}\label{topfl21}
Sous les hypothèses de \ref{topfl19}, si, de plus, $\tY$ est local de centre $y$, alors 
le topos $\tX\gtimes_\tS\tY$ est local de centre $\gamma y$.
\end{cor}

En effet, pour tout faisceau $F$ de $\tX\gtimes_\tS\tY$, on a $\Gamma(\tX\gtimes_\tS\tY,F)=\Gamma(\tY,\rp_{2*}(F))=(\rp_{2*}(F))_{y}$. 
Par ailleurs, le morphisme 
\begin{equation}
\Gamma(\tX\gtimes_\tS\tY,F)=\Gamma(\tY,\rp_{2*}(F))\rightarrow \Gamma(\tY,\gamma^*(F))
\end{equation}
induit par le morphisme \eqref{topfl18d} s'identifie au morphisme défini par image inverse par $\gamma$. Le morphisme canonique 
\begin{equation}
\Gamma(\tX\gtimes_\tS\tY,F)\rightarrow F_{\gamma y}
\end{equation}
s'identifie donc à la fibre du morphisme $\rp_{2*}(F)\rightarrow \gamma^*(F)$ \eqref{topfl18d} en $y$, 
qui est un isomorphisme en vertu de \ref{topfl19}, d'où la proposition. 

\begin{cor}\label{topfl22}
Conservons les hypothèses de \ref{topfl19}, supposons, de plus, que $\tY$ ait assez de points. Alors~:
\begin{itemize}
\item[{\rm (i)}] Pour tout faisceau $F$ de $\tX\gtimes_\tS\tY$, l'application canonique 
\begin{equation}\label{topfl22a}
\Gamma(\tX\gtimes_\tS\tY,F)\rightarrow (\tY, \gamma^*F)
\end{equation}
est bijective. 
\item[{\rm (ii)}] Pour tout faisceau abélien $F$ de $\tX\gtimes_\tS\tY$ et tout entier $i\geq 0$, l'application canonique 
\begin{equation}\label{topfl22b}
\rH^i(\tX\gtimes_\tS\tY,F)\rightarrow \rH^i(\tY, \gamma^*F)
\end{equation}
est bijective.
\end{itemize}
\end{cor}

Cela résulte aussitôt de \ref{topfl19} (cf. la preuve de l'énoncé analogue \cite{agt} VI.10.28).

\subsection{}\label{rec1}
Soient $X$ et $Y$ deux $\mU$-sites dans lesquels les 
limites projectives finies sont représentables, $f^+\colon X\rightarrow Y$ un foncteur continu et exact à gauche. 
On désigne par $\tX$ et $\tY$ les topos des faisceaux de $\mU$-ensembles sur $X$ et $Y$, respectivement, 
par $f\colon \tY\rightarrow \tX$ le morphisme de topos associé à $f^+$, et par 
$e_X$ et $e_Y$ des objets finaux de $X$ et $Y$, respectivement, qui existent par hypothèse. 
Le topos $\tX\gtimes_\tX\tY$ est appelé topos {\em co-évanescent} de $f$. 
On peut le définir par un site un peu plus simple que le site $C$ associé au couple de foncteurs $(\id_X,f^+)$ 
introduit dans \ref{topfl1}, comme suit. 
 
On désigne par $D_{f^+}$ la catégorie des paires $(U, V\rightarrow f^+(U))$, où $U$ est un objet de $X$ et  
$V\rightarrow f^+(U)$ est un morphisme de $Y$; un tel objet sera noté $(V\rightarrow U)$. 
Soient $(V\rightarrow U)$, $(V'\rightarrow U')$ deux objets de $D_{f^+}$. 
Un morphisme de $(V'\rightarrow U')$ dans $(V\rightarrow U)$
est la donnée de deux morphismes $V'\rightarrow V$ de $Y$ et $U'\rightarrow U$ de $X$, tels que le diagramme
\begin{equation}
\xymatrix{
V'\ar[r]\ar[d]&{f^+(U')}\ar[d]\\
{V}\ar[r]&{f^+(U)}}
\end{equation}
soit commutatif. Il résulte aussitôt des hypothèses que les limites projectives finies dans $D_{f^+}$ sont représentables.
On munit $D_{f^+}$ de la topologie {\em co-évanescente} (\cite{agt} VI.4.1), c'est-à-dire, 
la topologie engendrée par les recouvrements 
\begin{equation}
\{(V_i\rightarrow U_i)\rightarrow (V\rightarrow U)\}_{i\in I}
\end{equation}
des deux types suivants~:
\begin{itemize}
\item[$(\alpha)$] $U_i=U$ pour tout $i\in I$, et $(V_i\rightarrow V)_{i\in I}$ est une famille couvrante.
\item[$(\beta)$] $(U_i\rightarrow U)_{i\in I}$ est une famille couvrante, 
et pour tout $i\in I$, le morphisme canonique $V_i\rightarrow V\times_{f^+(U)}f^+(U_i)$ est un isomorphisme. 
\end{itemize}
Le site ainsi défini est appelé site {\em co-évanescent} associé au foncteur $f^+$; c'est un $\mU$-site. 
On désigne par $\hD_{f^+}$ (resp. $\tD_{f^+}$) la catégorie des préfaisceaux (resp. le topos des faisceaux) de 
$\mU$-ensembles sur $D_{f^+}$. 

Considérons le site $C$ associé au couple de foncteurs $(\id_X,f^+)$ défini dans \ref{topfl1} et les foncteurs
\begin{eqnarray}
\iota^+\colon D_{f^+}&\rightarrow &C,\ \ \ (V\rightarrow U)\mapsto (U\rightarrow U\leftarrow V),\label{co-ev9a}\\
\jmath^+\colon C&\rightarrow &D_{f^+},\ \ \  (U\rightarrow W\leftarrow V)\mapsto (V\times_{f^+(W)}f^+(U)\rightarrow U).\label{co-ev9b}
\end{eqnarray}
Il est clair que $\iota^+$ est un adjoint à gauche de $\jmath^+$ et que $\iota^+$ et $\jmath^+$ sont exacts à gauche.
D'après (\cite{agt} VI.4.8), les foncteurs $\iota^+$ et $\jmath^+$ sont continus. Ils définissent donc 
des morphismes de topos 
\begin{eqnarray}
\iota\colon \tC\rightarrow \tD_{f^+},\label{co-ev100a}\\
\jmath\colon \tD_{f^+}\rightarrow \tC.\label{co-ev100b}
\end{eqnarray}
Les morphismes d'adjonction $\id \rightarrow \jmath^+\circ \iota^+$ et $\iota^+ \circ \jmath^+\rightarrow \id$ 
induisent des morphismes $\iota_*\circ \jmath_*\rightarrow \id$ et $\id\rightarrow \jmath_*\circ \iota_*$
qui font de $\iota_*$ un adjoint à droite de $\jmath_*$. 

\begin{prop}[\cite{agt} VI.4.10]\label{co-ev101}
Sous les hypothèses de \ref{rec1}, 
les morphismes d'adjonction $\iota_*\circ \jmath_*\rightarrow \id$ et $\id\rightarrow \jmath_*\circ \iota_*$ 
sont des isomorphismes. En particulier, $\iota$ \eqref{co-ev100a} et $\jmath$ \eqref{co-ev100b} 
sont des équivalences de topos quasi-inverses l'une de l'autre.
\end{prop}

\subsection{}\label{co-ev102}
Conservons les hypothèses de \ref{rec1}. Les foncteurs 
\begin{eqnarray}
\rp_1^+\colon X&\rightarrow& D_{f^+},\ \ \ U\mapsto (f^+(U)\rightarrow U),\label{rec1a}\\
\rp_2^+\colon Y&\rightarrow& D_{f^+},\ \ \ V\mapsto (V\rightarrow e_X),\label{rec1b}
\end{eqnarray}
sont exacts à gauche et continus. Ils définissent donc deux morphismes de topos 
\begin{eqnarray}
\rp_1\colon \tD_{f^+}\rightarrow \tX,\label{rec1c}\\
\rp_2\colon \tD_{f^+}\rightarrow \tY.\label{rec1d}
\end{eqnarray}
Les morphismes $\rp_1\circ \iota$  et $\rp_2\circ \iota$ s'identifient aux morphismes 
$\rp_1\colon \tC\rightarrow \tX$ et $\rp_2\colon \tC\rightarrow \tY$ définis dans \ref{topfl4}. 
Le $2$-morphisme \eqref{topfl4e}
\begin{equation}\label{rec1e}
\tau\colon f\rp_2\rightarrow \rp_1
\end{equation}
est alors défini par le morphisme de foncteurs $(f\rp_2)_*\rightarrow \rp_{1*}$ suivant: pour tout faisceau $F$ sur $D_{f^+}$ et 
tout $U\in \ob(X)$, 
\begin{equation}\label{rec1f}
f_*(\rp_{2*}(F))(U)\rightarrow \rp_{1*}(F)(U)
\end{equation}
est l'application canonique
\begin{equation}
F(f^+(U)\rightarrow e_X)\rightarrow F(f^+(U)\rightarrow U).
\end{equation}
Dans la suite de ce livre, nous identifierons les topos $\tD_{f^+}$ et $\tX\gtimes_{\tX}\tY$ au moyen des équivalences quasi-inverses l'une de l'autre
$\iota$ et $\jmath$.

Le foncteur 
\begin{equation}\label{rec1h}
\Psi^+\colon D_{f^+}\rightarrow Y, \ \ \ (V\rightarrow U)\mapsto V
\end{equation}
est continu et exact à gauche (\cite{agt} VI.4.13). Il définit donc un morphisme de topos
\begin{equation}\label{rec1i}
\Psi\colon \tY\rightarrow \tX\gtimes_{\tX}\tY
\end{equation}
tel que $\rp_1\Psi=f$, $\rp_2\Psi=\id_\tY$ et $\tau*\Psi=\id_f$. Ces relations déterminent $\Psi $
compte tenu de la propriété universelle du topos co-évanescent.
\begin{equation}\label{rec1j}
\xymatrix{
&\tY\ar[dl]_f\ar[rd]^{\id_\tY}\ar[d]_{\Psi}&\\
{\tX}\ar[dr]&{\tX\gtimes_\tX\tY}\ar[r]^-(0.4){\rp_2}\ar[l]_-(0.4){\rp_1}&{\tY}\ar[ld]^f\\
&{\tX}&}
\end{equation}
Le morphisme $\Psi$ est appelé morphisme des {\em cycles co-proches}.
De la relation $\rp_{2*}\Psi_*=\id_{\tY}$, on obtient par adjonction un morphisme 
\begin{equation}\label{rec1k}
\rp_2^*\rightarrow \Psi_*.
\end{equation}
Celui-ci est un isomorphisme d'après (\cite{agt} VI.4.14); en particulier, le foncteur $\Psi_*$ est exact.

\subsection{}\label{rec3}
Soient $X$, $Y$, $X'$ et $Y'$ des $\mU$-sites dans lesquels les limites projectives sont représentables,
\begin{equation}\label{rec3a}
\xymatrix{
X\ar[r]^{f^+}\ar[d]_{u^+}&Y\ar[d]^{v^+}\\
X'\ar[r]^{f'^+}&Y'}
\end{equation}
un diagramme de foncteurs continus et exacts à gauche, commutatif à isomorphisme près~:
\begin{equation}\label{rec3b}
a\colon v^+\circ f^+\stackrel{\sim}{\rightarrow} f'^+\circ u^+.
\end{equation}
On note $\tX$, $\tY$, $\tX'$ et $\tY'$ les topos des faisceaux de $\mU$-ensembles sur $X$, $Y$, $X'$ et $Y'$,
respectivement,  
\begin{equation}\label{rec3c}
\xymatrix{
{\tY'}\ar[r]^{f'}\ar[d]_v&{\tX'}\ar[d]^u\\
{\tY}\ar[r]^f&{\tX}}
\end{equation}
le diagramme de morphismes de topos déduit de \eqref{rec3a} et 
\begin{equation}\label{rec3d}
\alpha\colon uf'\stackrel{\sim}{\rightarrow} fv
\end{equation}
l'isomorphisme induit par $a$. On désigne par $D_{f^+}$ (resp. $D_{f'^+}$) le site co-évanescent associé à $f^+$ (resp. $f'^+$) \eqref{rec1}. 
On a un foncteur 
\begin{equation}\label{rec3e}
\Phi^+\colon D_{f^+}\rightarrow D_{f'^+}, 
\end{equation}
qui à tout objet $(V\rightarrow U)$ de $D_{f^+}$ associe l'objet $(v^+(V)\rightarrow u^+(U))$
de $D_{f'^+}$ défini par le composé
\begin{equation}\label{rec3f}
\xymatrix{
{v^+(V)}\ar[r]&{v^+(f^+(U))}\ar[r]^{a(U)}&{f'^+(u^+(U))}}.
\end{equation}
Il est clairement exact à gauche. Il est continu en vertu de (\cite{agt} VI.4.4). 
Il définit donc un morphisme de topos 
\begin{equation}\label{rec3j}
\Phi\colon \tX'\gtimes_{\tX'}\tY'\rightarrow \tX\gtimes_{\tX}\tY,
\end{equation}
qui n'est autre que le morphisme associé dans \eqref{topfl10d} au diagramme commutatif \eqref{rec3c}. 

Le diagramme 
\begin{equation}\label{rec3o}
\xymatrix{
{\tY'}\ar[r]^-(0.5){\Psi'}\ar[d]_v&{\tX'\gtimes_{\tX'}\tY'}\ar[d]^\Phi\\
{\tY}\ar[r]^-(0.5)\Psi&{\tX\gtimes_{\tX}\tY}}
\end{equation}
où $\Psi$ et $\Psi'$ sont les morphismes des cycles co-proches \eqref{rec1i}, 
est clairement commutatif à isomorphisme canonique près.

\section{Limite projective de topos co-évanescents}

\subsection{}\label{lptce1} 
On désigne par $\fM$ la catégorie des morphismes de schémas et par $\fC$ la catégorie des morphismes 
de $\fM$. Les objets de $\fC$ sont donc des diagrammes commutatifs de morphismes de schémas 
\begin{equation}\label{lptce1a} 
\xymatrix{
V\ar[r]\ar[d]&U\ar[d]\\
Y\ar[r]&X}
\end{equation} 
où l'on considère les flèches horizontales comme des objets de $\fM$
et les flèches verticales comme des morphismes de $\fM$; un tel objet sera noté $(V,U,Y,X)$.   
On désigne par $\fD$ la sous-catégorie pleine de $\fC$ formée des objets $(V,U,Y,X)$ tels que les morphismes
$U\rightarrow X$ et $V\rightarrow Y$ soient étales de présentation finie. Le ``foncteur but''
\begin{equation}\label{lptce1b} 
\fD\rightarrow \fM, \ \ \ (V,U,Y,X)\mapsto (Y\rightarrow X),
\end{equation}
fait de $\fD$ une catégorie fibrée, clivée et normalisée. Avec les notations de \ref{notconv10}, tout morphisme de schémas $f\colon Y\rightarrow X$ induit un foncteur
\begin{equation}
f^+_\coh\colon \Et_{\coh/X}\rightarrow \Et_{\coh/Y}.
\end{equation}
La catégorie fibre $\cD_f$ de \eqref{lptce1b} au-dessus de $f$ est la catégorie sous-jacente au site co-évanescent $D_{f^+_\coh}$ associé 
au foncteur $f^+_\coh$ \eqref{rec1}. Pour tout diagramme commutatif de morphismes de schémas
\begin{equation}\label{lptce1c}
\xymatrix{
Y'\ar[r]^{f'}\ar[d]_{g'}&X'\ar[d]^g\\
Y\ar[r]^f&X}
\end{equation}
le foncteur image inverse de \eqref{lptce1b} associé au morphisme $(g',g)$ de $\fM$ 
est le foncteur 
\begin{equation}\label{lptce1d}
\Phi^+\colon \fD_f\rightarrow \fD_{f'}, \ \ \ (V\rightarrow U)\mapsto (V\times_YY'\rightarrow U\times_XX').
\end{equation}
Celui-ci est continu et exact à gauche d'après \ref{rec3}. 
Munissant chaque fibre de $\fD/\fM$ de la topologie co-évanescente \eqref{rec1}, 
$\fD$ devient un $\mU$-site fibré (\cite{sga4} VI 7.2.4). On désigne par 
\begin{equation}\label{lptce1e}
\fF\rightarrow \fM
\end{equation}
le $\mU$-topos fibré associé à $\fD/\fM$ (\cite{sga4} VI 7.2.6).  
La catégorie fibre de $\fF$ au-dessus d'un objet $f\colon Y\rightarrow X$
de $\fM$ est le topos $\tfD_f$ des faisceaux de $\mU$-ensembles sur le site co-évanescent $\fD_f$, 
et le foncteur image inverse relatif au morphisme défini par 
le diagramme \eqref{lptce1c} est le foncteur $\Phi^*\colon \tfD_{f}\rightarrow \tfD_{f'}$ image inverse par le morphisme 
de topos $\Phi\colon \tfD_{f'}\rightarrow \tfD_{f}$ associé au morphisme de sites $\Phi^+$ \eqref{lptce1d}. 
On note 
\begin{equation}\label{lptce1f}
\fF^\vee\rightarrow \fM^\circ
\end{equation}
la catégorie fibrée obtenue en associant à tout objet $f\colon Y\rightarrow X$
de $\fM$ la catégorie $\fF_{f}=\tfD_{f}$, et à tout morphisme défini par 
un diagramme \eqref{lptce1c} le foncteur $\Phi_*\colon \tfD_{f'}\rightarrow \tfD_{f}$ image directe par le morphisme 
de topos $\Phi\colon \tfD_{f'}\rightarrow \tfD_{f}$. 

\subsection{}\label{lptce2}
Soient $I$ une catégorie cofiltrante essentiellement petite (\cite{sga4} I 2.7 et 8.1.8), 
\begin{equation}\label{lptce2a} 
\varphi\colon I\rightarrow \fM, \ \ \ i\mapsto (f_i\colon Y_i\rightarrow X_i)
\end{equation}
un foncteur tel que pour tout morphisme $j\rightarrow i$ de $I$, 
les morphismes $Y_j\rightarrow Y_i$ et $X_j\rightarrow X_i$ soient affines. 
On suppose qu'il existe $i_0\in \ob(I)$ tel que $X_{i_0}$ et $Y_{i_0}$ soient cohérents.
On désigne par 
\begin{eqnarray}
\fD_\varphi&\rightarrow& I\label{lptce2b}\\
\fF_\varphi&\rightarrow& I\label{lptce2c}\\
\fF^\vee_\varphi&\rightarrow& I^\circ\label{lptce2d}
\end{eqnarray}
les site, topos et catégorie fibrés déduits de $\fD$ \eqref{lptce1b}, $\fF$ \eqref{lptce1e} et $\fF^\vee$ \eqref{lptce1f}, respectivement,
par changement de base par le foncteur $\varphi$. On notera que $\fF_\varphi$ est le topos fibré associé à $\fD_\varphi$
(\cite{sga4} VI 7.2.6.8). D'après (\cite{ega4} 8.2.3), les limites projectives 
\begin{equation}\label{lptce2e}
X=\underset{\underset{i\in \ob(I)}{\longleftarrow}}{\lim}\ X_i \ \ {\rm et}\ \ 
Y=\underset{\underset{i\in \ob(I)}{\longleftarrow}}{\lim}\ Y_i
\end{equation}
sont représentables dans la catégorie des schémas. Les morphismes $(f_i)_{i\in I}$ induisent 
un morphisme $f\colon Y\rightarrow X$, qui représente la limite projective du foncteur \eqref{lptce2a}. 

Pour tout $i\in \ob(I)$, on a un diagramme commutatif canonique
\begin{equation}\label{lptce2f}
\xymatrix{
Y\ar[r]^f\ar[d]&X\ar[d]\\
Y_i\ar[r]^{f_i}&X_i}
\end{equation}
Il lui correspond un foncteur image inverse \eqref{lptce1d}
\begin{equation}\label{lptce2g}
\Phi_i^+\colon \fD_{f_i}\rightarrow \fD_f,
\end{equation}
qui est continu et exact à gauche, et par suite un morphisme de topos 
\begin{equation}\label{lptce2h}
\Phi_i\colon \tfD_{f}\rightarrow \tfD_{f_i}.
\end{equation}

On a un foncteur naturel 
\begin{equation}\label{lptce2i}
\fD_\varphi \rightarrow \fD_f,
\end{equation}
dont la restriction à la fibre au-dessus de tout $i\in \ob(I)$ est le foncteur $\Phi_i^+$ \eqref{lptce2g}. 
Ce foncteur transforme morphisme cartésien en isomorphisme. Il se factorise donc de façon unique à travers 
un foncteur (\cite{sga4} VI 6.3)
\begin{equation}\label{lptce2j}
\underset{\underset{I^\circ}{\longrightarrow}}{\lim}\ \fD_\varphi \rightarrow \fD_{f}.
\end{equation}
Le $I$-foncteur $\fD_\varphi\rightarrow \fD_f\times I$ déduit de \eqref{lptce2i} est un morphisme cartésien
de sites fibrés (\cite{sga4} VI 7.2.2). Il induit donc un morphisme cartésien de topos fibrés (\cite{sga4} VI 7.2.7)
\begin{equation}\label{lptce2k}
\tfD_f\times I\rightarrow \fF_\varphi.
\end{equation}

\begin{prop}\label{lptce3}
Le couple formé du topos $\tfD_f$ et du morphisme \eqref{lptce2k} 
est une limite projective du topos fibré $\fF_\varphi/I$ {\rm (\cite{sga4} VI 8.1.1)}.
\end{prop}

On notera d'abord que le foncteur \eqref{lptce2j} est une équivalence de catégories en vertu de 
(\cite{ega4} 8.8.2 et 17.7.8). Soient $T$ un $\mU$-topos, 
\begin{equation}
h\colon T\times I\rightarrow \fF_\varphi
\end{equation}
un morphisme cartésien de topos fibrés au-dessus de $I$. Notons $\varepsilon_I\colon \fD_\varphi\rightarrow \fF_\varphi$
le foncteur cartésien canonique (\cite{sga4} VI (7.2.6.7)), et posons 
\begin{equation}
h^+= h^*\circ \varepsilon_I\colon \fD_\varphi\rightarrow T\times I.
\end{equation}
Pour tout $i\in \ob(I)$, on désigne par
\begin{equation}
h_i^+\colon \fD_{f_i}\rightarrow T
\end{equation}
la restriction de $h^+$ aux fibres au-dessus de $i$. 
Compte tenu de l'équivalence de catégories \eqref{lptce2j} et de  (\cite{sga4} VI 6.2), il existe un et un essentiellement unique foncteur 
\begin{equation}
g^+\colon \fD_{f}\rightarrow T
\end{equation}
tel que $h^+$ soit isomorphe au composé 
\begin{equation}
\xymatrix{
{\fD_\varphi}\ar[r]&{\fD_{f}\times I}\ar[rr]^{g^+\times \id_I}&&{T\times I}},
\end{equation}
où la première flèche est le foncteur déduit de \eqref{lptce2i}. Montrons que $g^+$ est un morphisme de sites. 
Pour tout objet $(V\rightarrow U)$ de $\fD_f$, il existe $i\in \ob(I)$, un objet $(V_i\rightarrow U_i)$ de $\fD_{f_i}$ et un isomorphisme de $\fD_f$
\begin{equation}
(V\rightarrow U)\stackrel{\sim}{\rightarrow} \Phi_i^+(V_i\rightarrow U_i).
\end{equation}
Comme les foncteurs $h_i^+$ et $\Phi_i^+$ sont exacts à gauche, on en déduit que $g^+$ est exact à gauche. 
D'autre part, tout recouvrement {\em fini} de type $(\alpha)$ (resp. $(\beta)$) de $\fD_f$ \eqref{rec1}  
est l'image inverse d'un recouvrement de type $(\alpha)$ (resp. $(\beta)$) de $\fD_{f_i}$ pour un objet $i\in I$,  
en vertu de (\cite{ega4} 8.10.5(vi)). Comme les schémas $X$ et $Y$ sont cohérents, on en déduit 
que $g^+$ transforme les recouvrements de type $(\alpha)$ (resp. $(\beta)$) de $\fD_f$ en familles épimorphiques de $T$. 
Par suite, $g^+$ est continu en vertu de (\cite{agt} VI.4.4).  Il définit donc un morphisme de topos 
\begin{equation}
g\colon T\rightarrow \tfD_f
\end{equation}
tel que $h$ soit isomorphe au composé 
\begin{equation}
\xymatrix{
{T\times I}\ar[rr]^{g\times \id_I}&&{\tfD_f\times I}\ar[r]&{\fF_\varphi}},
\end{equation}
où la seconde flèche est le morphisme \eqref{lptce2k}. Un tel morphisme $g$ est essentiellement unique car la ``restriction'' 
$g^+\colon \fD_f\rightarrow T$ du foncteur $g^*$ est essentiellement unique d'après ce qui précède, d'où la proposition.

\subsection{}\label{lptce4}
Munissons $\fD_\varphi$ de la topologie totale (\cite{sga4} VI 7.4.1) 
et notons $\Top(\fD_\varphi)$ le topos des faisceaux de $\mU$-ensembles sur $\fD_\varphi$. 
D'après (\cite{sga4} VI 7.4.7), on a une équivalence canonique de catégories \eqref{notconv4}
\begin{equation}\label{lptce4a}
\Top(\fD_\varphi)\stackrel{\sim}{\rightarrow}\bHom_{I^\circ}(I^\circ, \fF^\vee_\varphi). 
\end{equation}
D'autre part, le foncteur naturel $\fD_\varphi\rightarrow \fD_f$ \eqref{lptce2i} est un morphisme de sites (\cite{sga4} VI 7.4.4)
et définit donc un morphisme de topos 
\begin{equation}\label{lptce4b}
\varpi\colon \tfD_f\rightarrow \Top(\fD_\varphi). 
\end{equation}
En vertu de \ref{lptce3} et (\cite{sga4} VI 8.2.9), il existe une équivalence de catégories $\Theta$
qui s'insère dans un diagramme commutatif 
\begin{equation}\label{lptce4c}
\xymatrix{
{\tfD_f}\ar[r]^-(0.5){\Theta}_-(0.5)\sim\ar[d]_{\varpi_*}&{\bHom_{\cart/I^\circ}(I^\circ, \fF^\vee_\varphi)}\ar@{^(->}[d]\\
{\Top(\fD_\varphi)}\ar[r]^-(0.5)\sim&{\bHom_{I^\circ}(I^\circ, \fF^\vee_\varphi)}}
\end{equation}
où la flèche horizontale inférieure est l'équivalence de catégories \eqref{lptce4a} et 
la flèche verticale de droite est l'injection canonique. 

Pour tout objet $F$ de $\Top(\fD_\varphi)$, si 
$\{i\mapsto F_i\}$ est la section correspondante de $\bHom_{I^\circ}(I^\circ, \fF^\vee_\varphi)$, on a un 
isomorphisme canonique fonctoriel (\cite{sga4} VI 8.5.2)
\begin{equation}\label{lptce4d}
\varpi^*(F)\stackrel{\sim}{\rightarrow} \underset{\underset{i\in I^\circ}{\longrightarrow}}{\lim}\ \Phi_i^*(F_i).
\end{equation}

\begin{cor}\label{lptce5}
Soient $F$ un faisceau de $\Top(\fD_\varphi)$, 
$\{i\mapsto F_i\}$ la section de $\bHom_{I^\circ}(I^\circ, \fF^\vee_\varphi)$ qui lui est 
associée par l'équivalence de catégories \eqref{lptce4a}.
Alors, on a un isomorphisme canonique fonctoriel
\begin{equation}\label{lptce5a}
\underset{\underset{i\in I^\circ}{\longrightarrow}}{\lim}\ \Gamma(\tfD_{f_i},F_i)\stackrel{\sim}{\rightarrow}
\Gamma(\tfD_{f},\underset{\underset{i\in I^\circ}{\longrightarrow}}{\lim}\ \Phi_i^*(F_i)).
\end{equation}
\end{cor}

\begin{cor}\label{lptce6}
Soit $F$ un faisceau abélien de $\Top(\fD_\varphi)$ et soit
$\{i\mapsto F_i\}$ la section de $\bHom_{I^\circ}(I^\circ, \fF^\vee_\varphi)$ qui lui est 
associée par l'équivalence de catégories \eqref{lptce4a}.
Alors, pour tout entier $q\geq 0$, on a un isomorphisme canonique fonctoriel
\begin{equation}\label{lptce6a}
\underset{\underset{i\in I^\circ}{\longrightarrow}}{\lim}\ \rH^q(\tfD_{f_i},F_i)\stackrel{\sim}{\rightarrow}
\rH^q(\tfD_{f},\underset{\underset{i\in I^\circ}{\longrightarrow}}{\lim}\ \Phi_i^*(F_i)).
\end{equation}
\end{cor}

Les corollaires \ref{lptce5} et \ref{lptce6} résultent de \ref{lptce3} et (\cite{sga4} VI 8.7.7). 
On notera que les conditions requises dans 
(\cite{sga4} VI 8.7.1 et 8.7.7) sont satisfaites en vertu de \ref{topfl15}(iii), \ref{topfl16}(i) et (\cite{sga4} VI 3.3, 5.1 et 5.2).

\section{Topos de Faltings}

On rappelle dans cette section la définition du topos de Faltings suivant l'approche développée dans (\cite{agt} VI.10).  
On établit aussi quelques propriétés supplémentaires. 

\subsection{}\label{tf0}
On désigne par $\cR$ la catégorie des revêtements étales ({\em i.e.}, la sous-catégorie pleine de la catégorie 
des morphismes de schémas, formée des revêtements étales) et par 
\begin{equation}\label{tf0a}
\cR\rightarrow \Sch
\end{equation}
le ``foncteur but'', qui fait de $\cR$ une catégorie fibrée clivée et normalisée au-dessus de $\Sch$ \eqref{notconv3} (\cite{sga1} VI):  
la catégorie fibre au-dessus de tout schéma $X$ est canoniquement équivalente à
la catégorie $\Et_{\rf/X}$, et pour tout morphisme de schémas 
$f\colon Y\rightarrow X$, le foncteur image inverse $f^+\colon \Et_{\rf/X}\rightarrow \Et_{\rf/Y}$
n'est autre que le foncteur changement de base par $f$ \eqref{notconv10}. 
Munissant chaque fibre de la topologie étale, $\cR/\Sch$ devient un $\mU$-site fibré (\cite{sga4} VI 7.2.1).

\subsection{}\label{tf1}
Dans cette section, $f\colon Y\rightarrow X$ désigne un morphisme de schémas et 
\begin{equation}\label{tf1a}
\pi\colon E\rightarrow \Et_{/X}
\end{equation}
le $\mU$-site fibré déduit du site fibré des revêtements étales $\cR/\Sch$ \eqref{tf0a} 
par changement de base par le foncteur 
\begin{equation}\label{tf1aa}
\Et_{/X}\rightarrow \Sch, \ \ \ U\mapsto U\times_XY.
\end{equation} 
On dit que $\pi$ est le {\em site fibré de Faltings} associé à $f$ (\cite{agt} VI.10.1). 
On peut décrire explicitement la catégorie $E$ de la façon suivante. Les objets de $E$ 
sont les morphismes de schémas $V\rightarrow U$ au-dessus de $f\colon Y\rightarrow X$ tels que le morphisme
$U\rightarrow X$ soit étale et que le morphisme $V\rightarrow U_Y=U\times_XY$ soit étale fini. 
Soient $(V'\rightarrow U')$, $(V\rightarrow U)$ deux objets de $E$. Un morphisme 
de $(V'\rightarrow U')$ dans $(V\rightarrow U)$ est la donnée d'un $X$-morphisme $U'\rightarrow U$ et 
d'un $Y$-morphisme $V'\rightarrow V$ tels que le diagramme
\begin{equation}\label{tf1ab}
\xymatrix{
V'\ar[r]\ar[d]&U'\ar[d]\\
V\ar[r]&U}
\end{equation}
soit commutatif. Le foncteur $\pi$ est alors défini pour tout objet $(V\rightarrow U)$ de $E$, par
\begin{equation}\label{tf1b}
\pi(V\rightarrow U)=U.
\end{equation}
La fibre de $E$ au-dessus d'un objet $U$ de $\Et_{/X}$ s'identifie canoniquement au site fini étale de $U_Y$. On note 
\begin{equation}\label{tf1c}
\iota_{U!}\colon \Et_{\rf/U_Y}\rightarrow E, \ \ \ V\mapsto (V\rightarrow U)
\end{equation}
le foncteur canonique (\cite{agt} (VI.5.1.2)). 

On munit $E$ de la topologie co-évanescente associée à $\pi$ (\cite{agt} VI.5.3), autrement dit,  
la topologie engendrée par les recouvrements $\{(V_i\rightarrow U_i)\rightarrow (V\rightarrow U)\}_{i\in I}$
des deux types suivants: 
\begin{itemize}
\item[(v)] $U_i=U$ pour tout $i\in I$, et $(V_i\rightarrow V)_{i\in I}$ est un recouvrement étale. 
\item[(c)] $(U_i\rightarrow U)_{i\in I}$ est un recouvrement étale et $V_i=U_i\times_UV$ pour tout $i\in I$. 
\end{itemize}
Le site co-évanescent $E$ ainsi défini est encore appelé {\em site de Faltings} associé à $f$;  c'est un $\mU$-site.  
On désigne par $\hE$ (resp. $\tE$) la catégorie des préfaisceaux (resp. le topos des faisceaux) 
de $\mU$-ensembles sur $E$. On appelle $\tE$ le {\em topos de Faltings} associé à $f$ (\cite{agt} VI.10.1).

On désigne par 
\begin{equation}\label{tf1d}
\fF\rightarrow  \Et_{/X}
\end{equation}
le $\mU$-topos fibré associé à $\pi$. La catégorie fibre de $\fF$ au-dessus de tout $U\in \ob(\Et_{/X})$ est 
canoniquement équivalente au topos fini étale $(U_Y)_\fet$ de $U_Y$ \eqref{notconv10}
et le foncteur image inverse 
pour tout morphisme $t\colon U'\rightarrow U$ de $\Et_{/X}$ s'identifie au foncteur 
$(t_Y)_{\fet}^*\colon (U_Y)_\fet\rightarrow (U'_Y)_\fet$ image inverse par le morphisme 
de topos $(t_Y)_\fet\colon (U'_Y)_\fet\rightarrow (U_Y)_\fet$ (\cite{agt} VI.9.3). On désigne par
\begin{equation}\label{tf1e}
\fF^\vee\rightarrow (\Et_{/X})^\circ
\end{equation}
la catégorie fibrée obtenue en associant à tout $U\in \ob(\Et_{/X})$ la catégorie $(U_Y)_\fet$, et à tout morphisme 
$t\colon U'\rightarrow U$ de $\Et_{/X}$ le foncteur 
$(t_Y)_{\fet*}\colon (U'_Y)_\fet\rightarrow (U_Y)_\fet$ 
image directe par le morphisme de topos $(t_Y)_\fet$. On désigne par
\begin{equation}\label{tf1f}
\cP^\vee\rightarrow (\Et_{/X})^\circ
\end{equation}
la catégorie fibrée obtenue en associant à tout $U\in \ob(\Et_{/X})$ la catégorie $(\Et_{\rf/U_Y})^\wedge$ 
des préfaisceaux de $\mU$-ensembles sur $\Et_{\rf/U_Y}$, et à tout morphisme $t\colon U'\rightarrow U$ de $\Et_{/X}$ 
le foncteur 
\begin{equation}\label{tf1g}
(t_Y)_{\fet*}\colon (\Et_{\rf/U'_Y})^\wedge\rightarrow (\Et_{\rf/U_Y})^\wedge
\end{equation} 
obtenu en composant avec le foncteur image inverse $t^+_Y\colon \Et_{\rf/U_Y}\rightarrow \Et_{\rf/U'_Y}$.  

On a une équivalence de catégories (\cite{agt} VI.5.2)
\begin{eqnarray}\label{tf1h}
\hE&\rightarrow& \bHom_{(\Et_{/X})^\circ}((\Et_{/X})^\circ,\cP^\vee)\\
F&\mapsto &\{U\mapsto F\circ \iota_{U!}\}.\nonumber
\end{eqnarray}
On identifiera dans la suite $F$ à la section $\{U\mapsto F\circ \iota_{U!}\}$ qui lui est associée par cette équivalence.

D'après (\cite{agt} VI.5.11), le foncteur \eqref{tf1h} induit un foncteur pleinement fidèle 
\begin{equation}\label{tf1i}
\tE\rightarrow \bHom_{(\Et_{/X})^\circ}((\Et_{/X})^\circ,\fF^\vee),
\end{equation}
d'image essentielle les sections $\{U\mapsto F_U\}$ vérifiant une condition de recollement. 

Les foncteurs 
\begin{eqnarray}
\sigma^+\colon \Et_{/X}\rightarrow E,&& U\mapsto (U_Y\rightarrow U),\label{tf1jj}\\
\alpha_{X!}\colon \Et_{\rf/Y}\rightarrow E,&& V\mapsto (V\rightarrow X).\label{tf1j}
\end{eqnarray}
sont exacts à gauche et continus (\cite{agt} VI.10.6). Ils définissent donc deux morphismes de topos 
\begin{eqnarray}
\sigma\colon \tE\rightarrow X_\et,\label{tf1kk}\\
\beta\colon \tE\rightarrow Y_\fet.\label{tf1k}
\end{eqnarray} 
Pour tout faisceau $F=\{U\mapsto F_U\}$ sur $E$, on a $\beta_*(F)=F_X$. 

Le foncteur 
\begin{equation}\label{tf1l}
\psi^+\colon E\rightarrow \Et_{/Y},\ \ \ (V\rightarrow U)\mapsto V
\end{equation}
est continu et exact à gauche (\cite{agt} VI.10.7). Il définit donc un morphisme de topos 
\begin{equation}\label{tf1m}
\psi\colon Y_\et\rightarrow \tE.
\end{equation}
Nous changeons ici de notations par rapport à {\em loc. cit.}, 
en réservant la notation $\Psi$ aux foncteurs des cycles co-proches dans le sens 
strict \eqref{rec1i}. Pour tout faisceau $F$ de $Y_\et$, on a un isomorphisme canonique de $\tE$
\begin{equation}\label{tf1n}
\psi_*(F)\stackrel{\sim}{\rightarrow} \{U\mapsto \rho_{U_Y*}(F|U_Y)\},
\end{equation}
où pour tout objet $U$ de $\Et_{/X}$, $\rho_{U_Y}\colon (U_Y)_\et\rightarrow (U_Y)_\fet$ est le morphisme canonique \eqref{notconv10a}.

\subsection{}\label{tf2}
On désigne par $\Et_{\coh/X}$ la sous-catégorie pleine de $\Et_{/X}$ formée des schémas étales 
de présentation finie sur $X$, munie de la topologie induite par celle de $\Et_{/X}$ \eqref{notconv10}, et par 
\begin{equation}
\pi_\coh\colon E_\coh\rightarrow \Et_{\coh/X}
\end{equation}
le site fibré déduit de $\pi$ \eqref{tf1a} par changement de base par le foncteur d'injection canonique 
\begin{equation}
\Et_{\coh/X}\rightarrow \Et_{/X}. 
\end{equation} 
On munit $E_\coh$ de la topologie co-évanescente définie par $\pi_\coh$ (\cite{agt} VI.5.3).
D'après (\cite{agt} VI.10.4), si $X$ est quasi-séparé, la projection canonique $E_\coh\rightarrow E$ induit par restriction une équivalence entre le topos 
$\tE$ et le topos des faisceaux de $\mU$-ensembles sur $E_\coh$. De plus, sous la même hypothèse,
la topologie co-évanescente de $E_\coh$ est induite par celle de $E$.

\subsection{}\label{tf3}
On désigne par $D_{f^+}$ le site co-évanescent du foncteur $f^+\colon \Et_{/X}\rightarrow \Et_{/Y}$ de changement de base par $f$, 
et par $X_\et\gtimes_{X_\et}Y_\et$ le topos co-évanescent du morphisme $f_\et\colon Y_\et\rightarrow X_\et$ (cf. \ref{rec1}). 
Tout objet de $E$ est naturellement un objet de $D_{f^+}$. 
On définit ainsi un foncteur pleinement fidèle et exact à gauche
\begin{equation}\label{tf3a}
\rho^+\colon E\rightarrow D_{f^+}.
\end{equation} 
Celui-ci est continu d'après (\cite{agt} VI.10.15). Il définit donc un morphisme de topos 
\begin{equation}\label{tf3b}
\rho\colon X_\et\gtimes_{X_\et}Y_\et\rightarrow \tE.
\end{equation}

Il résulte aussitôt des définitions \eqref{co-ev102} que les carrés du diagramme 
\begin{equation}\label{tf3c}
\xymatrix{
{X_\et}\ar@{=}[d]&{X_\et\gtimes_{X_\et}Y_\et}\ar[l]_-(0.5){\rp_1}\ar[d]^{\rho}\ar[r]^-(0.5){\rp_2}&
{Y_\et}\ar[d]^{\rho_Y}\\
{X_\et}&{\tE}\ar[l]_{\sigma}\ar[r]^{\beta}&{Y_\fet}}
\end{equation}
où $\rho_Y$ est le morphisme canonique \eqref{notconv10a}, et le diagramme 
\begin{equation}\label{tf3d}
\xymatrix{
{Y_\et}\ar[r]^-(0.5){\Psi}\ar[dr]_{\psi}&{X_\et\gtimes_{X_\et}Y_\et}\ar[d]^{\rho}\\
&{\tE}}
\end{equation}
sont commutatifs à isomorphismes canoniques près.

\begin{prop}\label{tf4}
Supposons $X$ et $Y$ cohérents. Alors,
\begin{itemize}
\item[{\rm (i)}] Pour tout objet $(V\rightarrow U)$ de $E_{\coh}$ \eqref{tf2}, le faisceau associé $(V\rightarrow U)^a$ est un objet cohérent de $\tE$. 
\item[{\rm (ii)}] Le topos $\tE$ est cohérent; en particulier, il a suffisamment de points. 
\item[{\rm (iii)}] Le morphisme $\rho$ \eqref{tf3b} est cohérent. 
\end{itemize}
\end{prop}

Les propositions (i) et (ii) sont mentionnées à titre de rappel (\cite{agt} VI.10.5). Montrons la proposition (iii). 
Avec les notations de \ref{notconv10}, le morphisme $f$ induit un foncteur 
\begin{equation}
f^+_\coh\colon \Et_{\coh/X}\rightarrow \Et_{\coh/Y}.
\end{equation}
On note $D_{f^+_\coh}$ le site co-évanescent associé à $f^+_\coh$ \eqref{rec1}. 
D'après \ref{rec3}, le foncteur d'injection canonique $D_{f^+_\coh}\rightarrow D_{f^+}$ est continu et exact à gauche. 
Il induit une équivalence entre les topos associés. 
Pour tout objet $(V\rightarrow U)$ de $E_\coh$, $\rho^+(V\rightarrow U)$ est un objet de $D_{f^+_\coh}$.
La proposition (iii) résulte alors de (i), \ref{topfl15}(iii) et (\cite{sga4} VI 3.2) 
appliqué à la famille topologiquement génératrice $E_\coh$ de $E$.

\begin{prop}[\cite{agt} VI.10.21]\label{tf5} Si les schémas $X$ et $Y$ sont cohérents,
lorsque $(\oy \rightsquigarrow \ox)$ décrit la famille des points de $X_\et\gtimes_{X_\et}Y_\et$ \eqref{topfl17}, 
la famille des foncteurs fibres de $\tE$ associés aux points $\rho(\oy \rightsquigarrow \ox)$ est conservative
\eqref{tf3b} {\rm (\cite{sga4} IV 6.4.0)}. 
\end{prop}

On notera que sous les mêmes hypothèses, la famille des foncteurs fibres de $X_\et\gtimes_{X_\et}Y_\et$ 
associés aux points $(\oy \rightsquigarrow \ox)$ est conservative en vertu de \ref{topfl16}(i).

\subsection{}\label{tf6}
Soient $\uX\rightarrow X$ un morphisme, $\uY=Y\times_X\uX$, $\uf\colon \uY\rightarrow \uX$ la projection canonique,
$\uX_\et\gtimes_{\uX_\et}\uY_\et$ le topos co-évanescent et $\tuE$ le topos de Faltings associés à $\uf$. On vérifie aussitôt que le diagramme
\begin{equation}\label{tf6a}
\xymatrix{
{\uX_\et\gtimes_{\uX_\et}\uY_\et}\ar[d]_{\urho}\ar[r]^{\nu}&{X_\et\gtimes_{X_\et}Y_\et}\ar[d]^\rho\\
{\tuE}\ar[r]^{\Phi}&\tE}
\end{equation}
où les morphismes $\nu$ et $\Phi$ sont définis par fonctorialité \eqref{rec3j} et (\cite{agt} (VI.10.12.5)), 
et $\rho$ et $\urho$ sont les morphismes canoniques \eqref{tf3b}, est commutatif à isomorphisme canonique près. 
On en déduit un morphisme de changement de base (cf. \cite{egr1} 1.2.2)
\begin{equation}\label{tf6b}
\fc\colon \Phi^*\rho_*\rightarrow \urho_*\nu^*.
\end{equation}
En restreignant aux faisceaux abéliens, on déduit un morphisme pour tout entier $q\geq 0$, 
\begin{equation}\label{tf6c}
\fc^q\colon \Phi^*\rR^q\rho_*\rightarrow \rR^q\urho_*\nu^*.
\end{equation}

\begin{prop}\label{tf7}
Conservons les hypothèses et notations de \ref{tf6}, supposons, de plus, que l'une des deux hypothèses suivantes soit satisfaite~:
\begin{itemize}
\item[{\rm (a)}] $\uX$ est étale au-dessus de $X$;
\item[{\rm (b)}] le morphisme $f\colon Y\rightarrow X$ est cohérent, et $\uX$ est le localisé strict de $X$ en un point géométrique $\ox$ de $X$. 
\end{itemize} 
Alors, 
\begin{itemize}
\item[{\rm (i)}] Pour tout faisceau $F$ de $X_\et\gtimes_{X_\et}Y_\et$, le morphisme de changement de base \eqref{tf6b}
\begin{equation}\label{tf7a}
\fc(F)\colon \Phi^*(\rho_*F)\rightarrow \urho_*(\nu^*F)
\end{equation}
est un isomorphisme. 
\item[{\rm (ii)}] Pour tout faisceau abélien $F$ de $X_\et\gtimes_{X_\et}Y_\et$ et tout entier $q\geq 0$, le morphisme de changement de base \eqref{tf6c}
\begin{equation}\label{tf7b}
\fc^q(F)\colon \Phi^*(\rR^q\rho_*F)\rightarrow \rR^q\urho_*(\nu^*F)
\end{equation}
est un isomorphisme. 
\end{itemize}
\end{prop}

Considérons d'abord le cas (a). Notons $(\uY\rightarrow \uX)^a$ le faisceau de $\tE$ associé à 
l'objet $(\uX\rightarrow \uY)$ de $E$. D'après (\cite{agt} VI.10.14), on a une équivalence canonique
\begin{equation}\label{tf7c}
\tuE\stackrel{\sim}{\rightarrow} \tE_{/(\uY\rightarrow \uX)^a},
\end{equation}
dont le composé avec le morphisme de localisation de $\tE$ en  $(\uY\rightarrow \uX)^a$ est égal à $\Phi$. D'après (\cite{agt} VI.4.18), 
on a une équivalence canonique
\begin{equation}\label{tf7d}
\uX_\et\gtimes_{\uX_\et}\uY_\et\stackrel{\sim}{\rightarrow} (X_\et\gtimes_{X_\et}Y_\et)_{/\rho^*((\uY\rightarrow \uX)^a)},
\end{equation}
dont le composé avec le morphisme de localisation de $X_\et\gtimes_{X_\et}Y_\et$ en  $\rho^*((\uY\rightarrow \uX)^a)$ est égal à $\nu$.
De plus, $\urho$ s'identifie à la restriction de $\rho$ au-dessus de $(\uY\rightarrow \uX)^a$ (\cite{sga4} IV (5.10.3)). La proposition (i) s'ensuit aussitôt,
et la proposition (ii) est une conséquence de (\cite{sga4} V 5.1(3)).   

Considérons ensuite le cas (b). 
Choisissons un voisinage étale affine $X'_0$ de $\ox$ dans $X$ et notons $\fV_\ox$ la catégorie des $X'_0$-schémas étales $\ox$-pointés
qui sont affines au-dessus de $X'_0$ (cf. \cite{sga4} VIII 3.9 et 4.5). Tout objet $X'$ de $\fV_\ox$ donne lieu à un diagramme à carrés 
commutatifs à isomorphismes canoniques près
\begin{equation}\label{tf7e}
\xymatrix{
{\uX_\et\gtimes_{\uX_\et}\uY_\et}\ar[d]_{\urho}\ar[r]^-(0.5){\unu_{X'}}&{X'_\et\gtimes_{X'_\et}(Y\times_XX')_\et}
\ar[d]^{\rho_{X'}}\ar[r]^-(0.5){\nu_{X'}}&
{X_\et\gtimes_{X_\et}Y_\et}\ar[d]^\rho\\
{\tuE}\ar[r]^{\uPhi_{X'}}&{\tE_{X'}}\ar[r]^{\Phi_{X'}}&\tE}
\end{equation}
où $X'_\et\gtimes_{X'_\et}(Y\times_XX')_\et$ est le topos co-évanescent et $\tE_{X'}$ est le topos de Faltings associés à la projection canonique 
$f'\colon Y\times_XX'\rightarrow X'$, $\rho_{X'}$ est le morphisme canonique \eqref{tf3b} et les flèches horizontales sont les morphismes de fonctorialité. 

Suivant \ref{lptce1}, on désigne par 
\begin{equation}\label{tf7f}
\fF\rightarrow \fV_\ox
\end{equation}
le topos fibré obtenu en associant à tout objet $X'$ de $\fV_\ox$ le topos $X'_\et\gtimes_{X'_\et}(Y\times_XX')_\et$, 
et à tout morphisme $X''\rightarrow X'$ de $\fV_\ox$ le foncteur 
\[
X'_\et\gtimes_{X'_\et}(Y\times_XX')_\et\rightarrow X''_\et\gtimes_{X''_\et}(Y\times_XX'')_\et
\] 
image inverse par le morphisme de fonctorialité. 
En vertu de \ref{lptce3}, les morphismes $\unu_{X'}$ 
identifient le topos $\uX_\et\gtimes_{\uX_\et}\uY_\et$ à la limite projective du topos fibré $\fF$.

Suivant (\cite{agt} VI.11.1), on désigne par  
\begin{equation}\label{tf7g}
\fG\rightarrow \fV_\ox
\end{equation}
le topos fibré obtenu en associant à tout objet $X'$ de $\fV_\ox$ le topos $\tE_{X'}$, et à tout morphisme $X''\rightarrow X'$ de $\fV_\ox$ le foncteur 
$\tE_{X'}\rightarrow \tE_{X''}$ image inverse par le morphisme de fonctorialité. 
En vertu de (\cite{agt} VI.11.3), les morphismes $\uPhi_{X'}$ identifient le topos $\tuE$ à la limite projective du topos fibré $\fG$.

Les morphismes $\rho_{X'}$, pour $X'\in \ob(\fV_\ox)$, définissent un morphisme de topos fibrés (\cite{sga4} VI 7.1.6)
\begin{equation}\label{tf7h}
\varrho\colon \fF\rightarrow \fG.
\end{equation}
Le morphisme $\urho$ se déduit de $\varrho$ par passage à la limite projective (\cite{sga4} VI 8.1.4). 

D'après (\cite{sga4} VI 8.5.10), pour tout objet $F$ de $X_\et\gtimes_{X_\et}Y_\et$, on a un isomorphisme canonique
\begin{equation}\label{tf7i}
\urho_*(\nu^*(F))\stackrel{\sim}{\rightarrow} \underset{\underset{X'\in \ob(\fV_\ox)}{\longrightarrow}}{\lim}\ 
\uPhi_{X'}^*(\rho_{X'*}(\nu_{X'}^*(F))).
\end{equation}
On notera que les conditions requises dans 
(\cite{sga4} VI 8.5.10) sont satisfaites en vertu de \ref{topfl15}(iii), \ref{topfl16}(i), \ref{tf4} et (\cite{sga4} VI 3.3 et 5.1).
Par ailleurs, compte tenu de (\cite{egr1} 1.2.4(i)), le morphisme $\fc(F)$ \eqref{tf7a} s'identifie à la limite inductive des morphismes 
\begin{equation}\label{tf7j}
\uPhi_{X'}^*(\Phi_{X'}^*(\rho_*(F)))\rightarrow \uPhi_{X'}^*(\rho_{X'*}(\nu_{X'}^*(F)))
\end{equation}
déduits des morphismes de changement de base relativement au carré de droite de \eqref{tf7e}. 
Ces derniers sont des isomorphismes d'après le cas (a); d'où la proposition (i) dans le cas (b).

On démontre de même la proposition (ii) en utilisant (\cite{sga4} VI 8.7.5).

\subsection{}\label{tf8}
Supposons $X$ strictement local de point fermé $x$. 
D'après (\cite{agt} VI.10.23 et VI.10.24), on a un morphisme canonique de topos 
\begin{equation}\label{tf8a}
\theta\colon Y_\fet\rightarrow \tE,
\end{equation}
qui est une section de $\beta$ \eqref{tf1k}, {\em i.e.}, on a un isomorphisme canonique 
\begin{equation}\label{tf8g}
\beta\theta\stackrel{\sim}{\rightarrow} \id_{Y_\fet}.
\end{equation} 
On obtient un morphisme de changement de base (\cite{agt} (VI.10.24.4))
\begin{equation}\label{tf8b}
\beta_*\rightarrow \theta^*.
\end{equation}

D'après \ref{topfl18}, le point $x$ définit une section canonique 
\begin{equation}
Y_\et\rightarrow X_\et\gtimes_{X_\et}Y_\et.
\end{equation}
En vertu de (\cite{agt} VI.10.25), le diagramme 
\begin{equation}\label{tf8c}
\xymatrix{
{Y_\et}\ar[r]^-(0.5){\gamma}\ar[d]_{\rho_Y}&{X_\et\gtimes_{X_\et}Y_\et}\ar[d]^\rho\\
{Y_\fet}\ar[r]^{\theta}&{\tE}}
\end{equation}
où $\rho_Y$ est le morphisme canonique \eqref{notconv10a}, est commutatif à isomorphisme canonique près. 
Il induit donc un morphisme de changement de base
\begin{equation}\label{tf8d}
\theta^*\rho_*\rightarrow \rho_{Y*}\gamma^*.
\end{equation}
En restreignant aux faisceaux abéliens, on déduit un morphisme pour tout entier $q\geq 0$ (cf. \cite{egr1} 1.2.2)
\begin{equation}\label{tf8e}
\theta^*\rR^q\rho_*\rightarrow \rR^q\rho_{Y*}\gamma^*.
\end{equation}

\begin{prop}\label{tf9}
Conservons les hypothèses et notations de \ref{tf8}. Alors, 
\begin{itemize}
\item[{\rm (i)}] Le morphisme de changement de base $\beta_*\rightarrow \theta^*$ \eqref{tf8b} est un isomorphisme~; en particulier, 
le foncteur $\beta_*$ est exact.   
\item[{\rm (ii)}] Pour tout faisceau $F$ de $X_\et\gtimes_{X_\et}Y_\et$, le morphisme de changement de base \eqref{tf8d}
\begin{equation}\label{tf9a}
\theta^*(\rho_*(F))\rightarrow \rho_{Y*}(\gamma^*(F))
\end{equation}
est un isomorphisme. 
\item[{\rm (iii)}] Pour tout faisceau abélien $F$ de $X_\et\gtimes_{X_\et}Y_\et$ et tout entier $q\geq 0$, le morphisme de changement de base \eqref{tf8e}
\begin{equation}\label{tf9b}
\theta^*(\rR^q\rho_*(F))\rightarrow \rR^q\rho_{Y*}(\gamma^*(F))
\end{equation}
est un isomorphisme. 
\end{itemize}
\end{prop}

La proposition (i) est mentionnée à titre de rappel (\cite{agt} VI.10.27). Considérons le diagramme 
\begin{equation}\label{tf9f}
\xymatrix{
{Y_\et}\ar[r]^-(0.5){\gamma}\ar[d]_{\rho_Y}&{X_\et\gtimes_{X_\et}Y_\et}\ar[d]^\rho\ar[r]^-(0.5){\rp_2}&{Y_\et}\ar[d]^{\rho_Y}\\
{Y_\fet}\ar[r]^{\theta}&{\tE}\ar[r]^-(0.5){\beta}&{Y_\fet}}
\end{equation}
dont les carrés sont commutatifs à isomorphismes canoniques près \eqref{tf3c} et \eqref{tf8c}.

Le diagramme de morphismes de foncteurs
\begin{equation}\label{tf9c}
\xymatrix{
{\theta^*\rho_*}\ar[r]^-(0.5)c&{\rho_{Y*}\gamma^*}\\
{\beta_*\rho_*}\ar[u]^a\ar[r]^-(0.5)d&{\rho_{Y*}\rp_{2*}}\ar[u]_b}
\end{equation}
où $a$ est induit par \eqref{tf8b}, $b$ est induit par \eqref{topfl18d}, $c$ est le morphisme de changement de base \eqref{tf8d} 
et $d$ est l'isomorphisme de commutativité du carré de droite de \eqref{tf9f}, est commutatif. 
En effet, d'après (\cite{egr1} 1.2.4(iii)), $c\circ a$ est le morphisme de changement de base relativement au carré commutatif 
\begin{equation}\label{tf9d}
\xymatrix{
{Y_\et}\ar[r]^-(0.5){\gamma}\ar[d]_{\rho_Y}&{X_\et\gtimes_{X_\et}Y_\et}\ar[d]^{\beta\rho}\\
{Y_\fet}\ar@{=}[r]&{Y_\fet}}
\end{equation}
La proposition (ii) résulte de \eqref{tf9c}, (i) et \ref{topfl19}. 

Restreignant aux faisceaux abéliens, pour tout entier $q\geq 0$, le diagramme de morphismes de foncteurs
\begin{equation}\label{tf9e}
\xymatrix{
{\theta^*\rR^q\rho_*}\ar[r]^-(0.5){c^q}&{\rR^q\rho_{Y*}\gamma^*}\\
{\beta_*\rR^q\rho_*}\ar[u]^{a^q}\ar[r]^-(0.5){d^q}&{\rR^q\rho_{Y*}\rp_{2*}}\ar[u]_{b^q}}
\end{equation}
où $a^q$ est induit par \eqref{tf8b}, $b^q$ est induit par \eqref{topfl18d}, $c^q$ est le morphisme de changement de base \eqref{tf8e} 
et $d^q$ est l'isomorphisme déduit de la relation $\beta\rho=\rho_Y\rp_2$ \eqref{tf9f} et de l'exactitude des foncteurs $\beta_*$ et $\rp_{2*}$ 
(cf. (i) et \ref{topfl19}), 
est commutatif. En effet, d'après (\cite{egr1} 1.2.4(v)), $c^q\circ a^q$ est le morphisme de changement de base relativement au carré commutatif \eqref{tf9d}.
La proposition (iii) résulte de \eqref{tf9e}, (i) et \ref{topfl19}.

\subsection{}\label{tf10}
Soient $\ox$ un point géométrique de $X$, $\uX$ le localisé strict de $X$ en $\ox$, 
$\uY=Y\times_X\uX$, $\uf\colon \uY\rightarrow \uX$ la projection canonique,
$\uX_\et\gtimes_{\uX_\et}\uY_\et$ le topos co-évanescent et $\tuE$ 
le topos de Faltings associés à $\uf$. Les deux carrés du diagramme
\begin{equation}\label{tf10a}
\xymatrix{
{\uY_\et}\ar[r]^-(0.5)\gamma\ar[d]_{\rho_\uY}&{\uX_\et\gtimes_{\uX_\et}\uY_\et}\ar[d]^\urho\ar[r]^\nu&
{X_\et\gtimes_{X_\et}Y_\et}\ar[d]^\rho\\
{\uY_\fet}\ar[r]^-(0.5)\theta&{\tuE}\ar[r]^\Phi&\tE}
\end{equation}
où $\gamma$ est la section canonique définie par le point géométrique $\ox$  \eqref{topfl18}, 
$\theta$ est le morphisme \eqref{tf8a},
$\Phi$ et $\nu$ sont induits par fonctorialité par le morphisme canonique $\uX\rightarrow X$, 
$\urho$ et $\rho$ sont les morphismes canoniques \eqref{tf3b}
et $\rho_\uY$ est le morphisme canonique \eqref{notconv10a}, 
sont commutatifs à isomorphismes canoniques près \eqref{tf6a} et \eqref{tf8c}. On pose
\begin{eqnarray}
\phi_\ox=\gamma^*\circ \nu^* \colon X_\et\gtimes_{X_\et}Y_\et&\rightarrow& \uY_\et,\label{tf10c}\\
\varphi_\ox=\theta^*\circ \Phi^*\colon \tE&\rightarrow& \uY_\fet.\label{tf10b}
\end{eqnarray}
On note 
\begin{equation}\label{tf10d}
\fc\colon \varphi_\ox\rho_*\rightarrow \rho_{\uY*}\phi_\ox
\end{equation}
le morphisme de changement de base relativement au rectangle extérieur du diagramme \eqref{tf10a}. 
En restreignant aux faisceaux abéliens, on déduit un morphisme pour tout entier $q\geq 0$ 
\begin{equation}\label{tf10e}
\fc^q\colon \varphi_\ox\rR^q\rho_*\rightarrow \rR^q\rho_{\uY*}\phi_\ox.
\end{equation}
On renvoie à (\cite{egr1} 1.2.2) pour la définition du morphisme de changement de base.

\begin{prop}\label{tf11}
Conservons les hypothèses et notations de \ref{tf10}; supposons de plus que le morphisme $f\colon Y\rightarrow X$ soit cohérent. Alors,
\begin{itemize}
\item[{\rm (i)}] Le diagramme 
\begin{equation}\label{tf11b}
\xymatrix{
{\varphi_\ox\rho_*\rho^*} \ar[r]^-(0.5){\fc*\rho^*}&{\rho_{\uY*}\phi_\ox\rho^*}\ar@{=}[d]\\ 
{\varphi_\ox}\ar[u]^a\ar[r]^-(0.5)b&{\rho_{\uY*}\rho_\uY^*\varphi_\ox }}
\end{equation}
où $a$ et $b$ sont induits par les morphismes d'adjonction $\id\rightarrow \rho_*\rho^*$
et $\id\rightarrow \rho_{\uY*}\rho^*_\uY$, est commutatif. 

\item[{\rm (ii)}]  Pour tout faisceau $F$ de $X_\et\gtimes_{X_\et}Y_\et$, le morphisme de changement de base \eqref{tf10d}
\begin{equation}\label{tf11a}
\fc(F)\colon \varphi_\ox(\rho_*(F))\rightarrow \rho_{\uY*}(\phi_\ox (F))
\end{equation}
est un isomorphisme.

\item[{\rm (iii)}] Pour tout faisceau abélien $F$ de $X_\et\gtimes_{X_\et}Y_\et$ et tout entier $q\geq 0$, le morphisme de changement de base \eqref{tf10e}
\begin{equation}\label{tf11c}
\fc^q\colon \varphi_\ox(\rR^q\rho_*(F))\rightarrow \rR^q\rho_{\uY*}(\phi_\ox (F))
\end{equation}
est un isomorphisme. 
\end{itemize}
\end{prop}

(i) En effet, les triangles du diagramme
\begin{equation}\label{tf11g}
\xymatrix{
{\rho_\uY^*\varphi_\ox\rho_*\rho^*}\ar@{=}[r]\ar[rd]^{\fc'}&{\phi_\ox\rho^*\rho_*\rho^*}\ar[d]^d\\
{\rho_\uY^*\varphi_\ox}\ar@{=}[r]\ar[u]^{\rho_\uY^**a}&{\phi_\ox\rho^*}}
\end{equation}
où $\fc'$ est le morphisme adjoint de $\fc*\rho^*$ et $d$ est induit par le morphisme d'adjonction $\rho^*\rho_*\rightarrow \id$,
sont commutatifs en vertu de (\cite{sga4} XVII 2.1.3). 

(ii) Notons 
\begin{eqnarray}
\fc_1\colon \Phi^*\rho_*&\rightarrow& \urho_*\nu^*,\label{tf11d}\\
\fc_2\colon \theta^*\urho_*&\rightarrow& \rho_{\uY*}\gamma^*,\label{tf11e}
\end{eqnarray}
les morphismes de changement de base relativement au carré de droite et au carré de gauche du diagramme \eqref{tf10a}, respectivement. 
D'après (\cite{egr1} 1.2.4(i)), on a 
\begin{equation}\label{tf11f}
\fc(F)= \fc_2(\nu^*(F))\circ \theta^*(\fc_1(F)).
\end{equation}
Comme $\fc_1(F)$ est un isomorphisme en vertu de \ref{tf7}(i), et que $\fc_2(\nu^*(F))$ est un isomorphisme en vertu de \ref{tf9}(ii), $\fc(F)$ est un isomorphisme. 

(iii) Il suffit de calquer la preuve de (ii) en utilisant \ref{tf7}(ii), \ref{tf9}(iii) et (\cite{egr1} 1.2.4(ii)).

\begin{cor}\label{tf12}
Supposons que les schémas $X$ et $Y$ soient cohérents, et que pour tout point géométrique $\ox$ de $X$,
notant $X_{(\ox)}$ le localisé strict de $X$ en $\ox$, le schéma $Y_{(\ox)}=Y\times_XX_{(\ox)}$ ait un nombre
fini de composantes connexes.
Alors, le morphisme d'adjonction $\id\rightarrow \rho_*\rho^*$ est un isomorphisme
\eqref{tf3b}~; en particulier, le foncteur 
\begin{equation}\label{tf12a}
\rho^*\colon \tE\rightarrow X_\et\gtimes_{X_\et}Y_\et
\end{equation}
est pleinement fidèle. 
\end{cor}

En effet, soient $\ox$ un point géométrique de $X$, $\rho_{Y_{(\ox)}}\colon Y_{(\ox),\et}\rightarrow Y_{(\ox),\fet}$ le foncteur canonique \eqref{notconv10a}, 
$\varphi_\ox\colon \tE\rightarrow Y_{(\ox),\fet}$ le foncteur défini dans \eqref{tf10b}. 
En vertu de (\cite{agt} VI.9.18), le morphisme d'adjonction $\id\rightarrow \rho_{Y_{(\ox)}*}\rho_{Y_{(\ox)}}^*$ est un isomorphisme. 
On en déduit, d'après \ref{tf11}(i)-(ii), que le morphisme 
$\varphi_\ox\rightarrow \varphi_\ox \rho_*\rho^*$ induit par le morphisme d'adjonction $\id\rightarrow \rho_*\rho^*$ est un isomorphisme. 
La proposition s'ensuit puisque la famille des foncteurs $\varphi_\ox$ lorsque $\ox$ décrit l'ensemble des points géométriques de $X$, 
est conservative d'après (\cite{agt} VI.10.32).

\begin{cor}\label{tf13}
Soient $\mP$ un ensemble de nombres premiers, $F$ un faisceau abélien de $\mP$-torsion de $\tE$, $q$ un entier $\geq 1$.  
Supposons que  les schémas $X$ et $Y$ soient cohérents, et que pour tout point géométrique $\ox$ de $X$,
notant $X_{(\ox)}$ le localisé strict de $X$ en $\ox$, le schéma $Y_{(\ox)}=Y\times_XX_{(\ox)}$ soit $K(\pi,1)$ pour les faisceaux abéliens de $\mP$-torsion  \eqref{Kpun2}.
Alors, on a 
\begin{equation}\label{tf13a}
\rR^q\rho_*(\rho^*(F))=0.
\end{equation}
\end{cor}

En effet, soient $\ox$ un point géométrique de $X$, $\rho_{Y_{(\ox)}}\colon Y_{(\ox),\et}\rightarrow Y_{(\ox),\fet}$ le foncteur canonique \eqref{notconv10a}, 
\begin{eqnarray}
\phi_\ox \colon X_\et\gtimes_{X_\et}Y_\et&\rightarrow& Y_{(\ox),\et},\label{tf13b}\\
\varphi_\ox \colon \tE&\rightarrow& Y_{(\ox),\fet},\label{tf13c}
\end{eqnarray}
les foncteurs définis dans \eqref{tf10c} et \eqref{tf10b}. En vertu de \ref{tf11}(iii), on a un isomorphisme canonique
\begin{equation}
\varphi_\ox(\rR^q\rho_*(\rho^*(F)))\stackrel{\sim}{\rightarrow} \rR^q\rho_{Y_{(\ox)}*}(\phi_\ox (\rho^*(F))).
\end{equation}
Par ailleurs, on a un isomorphisme canonique $\phi_\ox (\rho^*(F))\stackrel{\sim}{\rightarrow}\rho_{Y_{(\ox)}}^*(\varphi_\ox (F))$ \eqref{tf10a}.
On en déduit que $\varphi_\ox(\rR^q\rho_*(\rho^*(F)))=0$. La proposition s'ensuit compte tenu de (\cite{agt} VI.10.32).

\subsection{}\label{tf14}
Considérons un diagramme commutatif de morphismes de schémas
\begin{equation}\label{tf14a}
\xymatrix{
Y'\ar[r]^{g'}\ar[d]_{f'}&Y\ar[d]^f\\
X'\ar[r]^g&X}
\end{equation}
On désigne par $E'$ et $\tE'$ le site et le topos de Faltings associés à $f'$, par 
\begin{eqnarray}
\sigma'\colon \tE'\rightarrow X'_\et,\label{tf14c}\\
\beta'\colon \tE'\rightarrow Y'_\fet,\label{tf14d}\\
\psi'\colon Y'_\et\rightarrow \tE'\label{tf14e},
\end{eqnarray} 
les morphismes canoniques \eqref{tf1kk}, \eqref{tf1k} et \eqref{tf1m} relatifs à $f'$ et par  
\begin{equation}\label{tf14b}
\Phi\colon \tE'\rightarrow \tE
\end{equation}
le morphisme de topos induit par le diagramme \eqref{tf14a} (\cite{agt} VI.10.12).

\begin{lem}\label{tf24}
Conservons les hypothèses et notations de \ref{tf14}; soient, de plus, $\ox'$ un point géométrique de $X'$, $\ox=g(\ox')$, 
$\uX$ (resp. $\uX'$) le localisé strict de $X$ en $\ox$ (resp. $X'$ en $\ox'$). On pose $\uY=\uX\times_XY$ et $\uY'=\uX'\times_{X'}Y'$ et on note 
$\ug'\colon \uY'\rightarrow \uY$ le morphisme induit par $g'$ \eqref{tf14a} et
\begin{eqnarray}
\varphi_{\ox}\colon \tE\rightarrow \uY_\fet,\label{tf24a}\\
\varphi_{\ox'}\colon \tE'\rightarrow \uY'_\fet,\label{tf24b}
\end{eqnarray}
les foncteurs définis dans \eqref{tf10b}. Alors, le diagramme de foncteurs
\begin{equation}\label{tf24c}
\xymatrix{
{\tE}\ar[r]^-(0.5){\varphi_\ox}\ar[d]_{\Phi^*}&{\uY_\fet}\ar[d]^{\ug'^*}\\
{\tE'}\ar[r]^-(0.5){\varphi_{\ox'}}&{\uY'_\fet}}
\end{equation}
est commutatif à isomorphisme canonique près.
\end{lem}

Considérons le diagramme commutatif
\begin{equation}\label{tf24d}
\xymatrix{
\uY'\ar[r]^{\ug'}\ar[d]_{\uf'}&\uY\ar[d]^\uf\\
\uX'\ar[r]^\ug&\uX}
\end{equation}
où $\uf$ et $\uf'$ sont les projections canoniques et $\ug$ est induit par $g$. 
On désigne par $\tuE$ et $\tuE'$ les topos de Faltings associés aux morphismes $\uf$ et $\uf'$ respectivement, et par
\begin{equation}\label{tf24e}
\uPhi\colon \tuE'\rightarrow \tuE
\end{equation}
le morphisme de topos induit par le diagramme \eqref{tf24d}. Le diagramme de morphismes de topos 
\begin{equation}\label{tf24f}
\xymatrix{
{\tuE'}\ar[r]^{\iota'}\ar[d]_{\uPhi}&{\tE'}\ar[d]^\Phi\\
{\tuE}\ar[r]^\iota&{\tE}}
\end{equation}
où $\iota$ et $\iota'$ sont les morphismes définis par fonctorialité, est commutatif à isomorphisme canonique près (\cite{agt} VI.10.12).
On peut donc se réduire au cas où $X$ et $X'$ sont strictement locaux de points fermés $\ox$ et $\ox'$, respectivement. 

Soit $(V\rightarrow U)$ un objet de $E$ tel que le morphisme $U\rightarrow X$ soit étale, séparé et de présentation finie. On pose 
$U'=U\times_XX'$ et $V'=V\times_YY'$, de sorte que $(V'\rightarrow U')$ est un objet de $E'$. 
On désigne par $U^\rf$ la somme disjointe des localisés stricts de $U$ en les points de $U_\ox$; c'est un sous-schéma ouvert et fermé de $U$, qui est fini 
sur $X$ (\cite{ega4} 18.5.11). On voit aussitôt que $U^\rf\times_XX'$ est la somme disjointe des localisés stricts de $U'$ en les points de $U'_{\ox'}$. 
Il résulte alors des définitions (\cite{agt} VI.10.22 et VI.10.23) que le diagramme de morphisme de topos 
\begin{equation}
\xymatrix{
{Y'_\fet}\ar[r]^{\theta'}\ar[d]_{g'}&{\tE'}\ar[d]^{\Phi}\\
{Y_\fet}\ar[r]^{\theta}&{\tE}}
\end{equation}
où $\theta$ et $\theta'$ sont définis dans \eqref{tf8a}, est commutatif à isomorphisme canonique près, d'où la proposition.

\subsection{}\label{tf15}
Conservons les hypothèses et notations de \ref{tf14}; soit, de plus, $u\colon \uX\rightarrow X$ un morphisme de schémas. On note
\begin{equation}\label{tf15a}
\Sch_{/X}\rightarrow \Sch_{/\uX}, \ \ \ U\mapsto \uU
\end{equation}
le foncteur de changement de base par $u$ \eqref{notconv3}. Le diagramme \eqref{tf14a} induit donc un diagramme commutatif
\begin{equation}\label{tf15b}
\xymatrix{
\uY'\ar[r]^{\ug'}\ar[d]_{\uf'}&\uY\ar[d]^\uf\\
\uX'\ar[r]^\ug&\uX}
\end{equation}
On désigne par $\tuE$ et $\tuE'$ les topos de Faltings associés aux morphismes $\uf$ et $\uf'$ respectivement, et par
\begin{equation}\label{tf15c}
\uPhi\colon \tuE'\rightarrow \tuE
\end{equation}
le morphisme de topos induit par le diagramme \eqref{tf15b}. Le diagramme de morphismes de topos 
\begin{equation}\label{tf15d}
\xymatrix{
{\tuE'}\ar[r]^{\iota'}\ar[d]_{\uPhi}&{\tE'}\ar[d]^\Phi\\
{\tuE}\ar[r]^\iota&{\tE}}
\end{equation}
où $\iota$ et $\iota'$ sont les morphismes induits par fonctorialité par $u$, est commutatif à isomorphisme canonique près (\cite{agt} VI.10.12). 
Il induit donc un morphisme de changement de base  (\cite{egr1} 1.2.2)
\begin{equation}\label{tf15e}
\iota^*\Phi_*\rightarrow \uPhi_*\iota'^*.
\end{equation}
En restreignant aux faisceaux abéliens, on déduit un morphisme pour tout entier $q\geq 0$ 
\begin{equation}\label{tf15f}
\iota^*\rR^q\Phi_*\rightarrow \rR^q\uPhi_*\iota'^*.
\end{equation}

\begin{prop}\label{tf16}
Conservons les hypothèses et notations de \ref{tf15}, supposons, de plus, que l'une des deux hypothèses suivantes soit satisfaite~:
\begin{itemize}
\item[{\rm (a)}] $\uX$ est étale au-dessus de $X$;
\item[{\rm (b)}] les morphismes $f$, $g$ et $f'$ sont cohérents, et $\uX$ est le localisé strict de $X$ en un point géométrique $\ox$ de $X$. 
\end{itemize} 
Alors, 
\begin{itemize}
\item[{\rm (i)}] Pour tout faisceau $F$ de $\tE'$, le morphisme de changement de base \eqref{tf15e}
\begin{equation}\label{tf16a}
\iota^*(\Phi_*F)\rightarrow \uPhi_*(\iota'^*F)
\end{equation}
est un isomorphisme. 
\item[{\rm (ii)}] Pour tout faisceau abélien $F$ de $\tE'$ et tout entier $q\geq 0$, le morphisme de changement de base \eqref{tf15f}
\begin{equation}\label{tf16b}
\iota^*(\rR^q\Phi_*(F))\rightarrow \rR^q\uPhi_*(\iota'^*F)
\end{equation}
est un isomorphisme. 
\end{itemize}
\end{prop}

Considérons d'abord le cas (a). Notons $(\uY\rightarrow \uX)^a$ le faisceau de $\tE$ associé à 
l'objet $(\uY\rightarrow \uX)$ de $E$. D'après (\cite{agt} VI.10.14), on a une équivalence canonique
\begin{equation}\label{tf16c}
\tuE\stackrel{\sim}{\rightarrow} \tE_{/(\uY\rightarrow \uX)^a},
\end{equation}
dont le composé avec le morphisme de localisation de $\tE$ en  $(\uY\rightarrow \uX)^a$ est égal à $\iota$. 
De même, notant $(\uY'\rightarrow \uX')^a$ le faisceau de $\tE'$ associé à 
l'objet $(\uY'\rightarrow \uX')$ de $E'$, on a une équivalence canonique
\begin{equation}\label{tf16d}
\tuE'\stackrel{\sim}{\rightarrow} \tE'_{/(\uY'\rightarrow \uX')^a}
\end{equation}
dont le composé avec le morphisme de localisation de $\tE'$ en  $(\uY'\rightarrow \uX')^a$ est égal à $\iota'$.
De plus, on a 
\begin{equation}
\Phi^*((\uY\rightarrow \uX)^a)=(\uY'\rightarrow \uX')^a
\end{equation}
et $\uPhi$ s'identifie à la restriction de $\Phi$ au-dessus de $(\uY\rightarrow \uX)^a$ (\cite{sga4} IV (5.10.3)). 
La proposition (i) s'ensuit aussitôt, et la proposition (ii) est une conséquence de (\cite{sga4} V 5.1(3)).

Considérons ensuite le cas (b). 
Choisissons un voisinage étale affine $X_0$ de $\ox$ dans $X$ et notons $\fV_\ox$ la catégorie des $X_0$-schémas étales $\ox$-pointés
qui sont affines au-dessus de $X_0$ (cf. \cite{sga4} VIII 3.9 et 4.5). Tout objet $X_1$ de $\fV_\ox$ donne lieu à un diagramme
à carrés commutatifs à isomorphismes canoniques près
\begin{equation}\label{tf16e}
\xymatrix{
{\tuE'}\ar[r]^-(0.5){\jmath'_{X_1}}\ar[d]_{\uPhi}&{\tE'_{X_1}}\ar[r]^{\iota'_{X_1}}\ar[d]_{\Phi_{X_1}}&{\tE'}\ar[d]^\Phi\\
{\tuE}\ar[r]^-(0.5){\jmath_{X_1}}&{\tE_{X_1}}\ar[r]^{\iota_{X_1}}&{\tE}}
\end{equation}
où $\tE_{X_1}$ et $\tE'_{X_1}$ sont les topos de Faltings associés à $f\times_XX_1$ et $f'\times_XX_1$ respectivement, 
et $\Phi_{X_1}$, $\iota_{X_1}$, $\iota_{X_1}'$, $\jmath_{X_1}$ et $\jmath'_{X_1}$ sont les morphismes de fonctorialité. 

Suivant (\cite{agt} VI.11.1), on désigne par  
\begin{equation}\label{tf16f}
\fG\rightarrow \fV_\ox,
\end{equation}
le topos fibré obtenu en associant à tout objet $X_1$ de $\fV_\ox$ le topos $\tE_{X_1}$, 
et à tout morphisme $X_2\rightarrow X_1$ de $\fV_\ox$ le foncteur 
$\tE_{X_1}\rightarrow \tE_{X_2}$ image inverse par le morphisme de fonctorialité. 
En vertu de (\cite{agt} VI.11.3), les morphismes $\jmath_{X_1}$ identifient le topos $\tuE$ à la limite projective du topos fibré $\fG$.

De même, on désigne par  
\begin{equation}\label{tf16g}
\fG'\rightarrow \fV_\ox,
\end{equation}
le topos fibré obtenu en associant à tout objet $X_1$ de $\fV_\ox$ le topos $\tE'_{X_1}$, 
et à tout morphisme $X_2\rightarrow X_1$ de $\fV_\ox$ le foncteur 
$\tE'_{X_1}\rightarrow \tE'_{X_2}$ image inverse par le morphisme de fonctorialité. 
Comme $g$ est cohérent, les morphismes $\jmath'_{X_1}$ identifient le topos $\tuE'$ à la limite projective du topos fibré $\fG'$ (\cite{agt} VI.11.3).

Les morphismes $\Phi_{X_1}$, pour $X_1\in \ob(\fV_\ox)$, définissent un morphisme de topos fibrés (\cite{sga4} VI 7.1.6)
\begin{equation}\label{tf16h}
\Lambda\colon \fG'\rightarrow \fG.
\end{equation}
Le morphisme $\uPhi$ se déduit de $\Lambda$ par passage à la limite projective (\cite{sga4} VI 8.1.4). 

D'après (\cite{sga4} VI 8.5.10), pour tout objet $F$ de $\tE'$, on a un isomorphisme canonique
\begin{equation}\label{tf16i}
\uPhi_*(\iota'^*(F))\stackrel{\sim}{\rightarrow} \underset{\underset{X_1\in \ob(\fV_\ox)}{\longrightarrow}}{\lim}\ \jmath_{X_1}^*(\Phi_{X_1*}(\iota'^*_{X_1}(F))).
\end{equation}
On notera que les conditions requises dans (\cite{sga4} VI 8.5.10) sont satisfaites en vertu de \ref{tf4} et (\cite{sga4} VI 3.3 et 5.1).
Par ailleurs, compte tenu de (\cite{egr1} 1.2.4(i)), le morphisme \eqref{tf16a} s'identifie à la limite inductive des morphismes 
\begin{equation}\label{tf16j}
\jmath_{X_1}^*(\iota_{X_1}^*(\Phi_*(F)))\rightarrow \jmath_{X_1}^*(\Phi_{X_1*}(\iota'^*_{X_1}(F)))
\end{equation}
déduits des morphismes de changement de base relativement au carré de droite de \eqref{tf16e}. 
Ces derniers sont des isomorphismes d'après le cas (a); d'où la proposition (i) dans le cas (b).

On démontre de même la proposition (ii) en utilisant (\cite{sga4} VI 8.7.5).

\subsection{}\label{tf20}
Soient $X_0$ un sous-schéma fermé de $X$, $Z$ l'ouvert complémentaire de $X_0$ dans $X$, $\varpi\colon Y\times_XZ\rightarrow Z$ la projection canonique. 
Comme $Z$ est un ouvert de $X_\et$, {\em i.e.}, un sous-objet de l'objet final $X$ (\cite{sga4} IV 8.3),
$\sigma^*(Z)=(Y\times_XZ\rightarrow Z)^a$ est un ouvert de $\tE$ \eqref{tf1kk}. 
D'après (\cite{agt} VI.10.14), le topos $\tE_{/\sigma^*(Z)}$ est canoniquement équivalent au topos de Faltings 
associé au morphisme  $\varpi$. On note 
\begin{equation}\label{tf20a}
\gamma\colon \tE_{/\sigma^*(Z)}\rightarrow \tE
\end{equation}
le morphisme de localisation de $\tE$ en $\sigma^*(Z)$, que l'on identifie au morphisme de fonctorialité 
induit par le diagramme canonique (\cite{agt} VI.10.12)
\begin{equation}
\xymatrix{
Y\times_XZ\ar[r]\ar[d]&Z\ar[d]\\
Y\ar[r]&X}
\end{equation}
On a alors une suite de trois foncteurs adjoints 
\begin{equation}\label{tf20b}
\gamma_!\colon \tE_{/\sigma^*(Z)}\rightarrow \tE, \ \ \ 
\gamma^* \colon \tE\rightarrow \tE_{/\sigma^*(Z)}, \ \ \ 
\gamma_*\colon \tE_{/\sigma^*(Z)}\rightarrow \tE, 
\end{equation}
dans le sens que pour deux foncteurs consécutifs de la suite, celui de droite est adjoint à droite
de l'autre. Les foncteurs $\gamma_!$ et $\gamma_*$ sont pleinement fidèles (\cite{sga4} IV 9.2.4). 

On désigne par $\tE_0$ le sous-topos fermé de $\tE$ complémentaire de l'ouvert $\sigma^*(Z)$, 
c'est-à-dire la sous-catégorie pleine de $\tE$ formée des faisceaux $F$ tels que $\gamma^*(F)$
soit un objet final de $\tE_{/\sigma^*(Z)}$  (\cite{sga4} IV 9.3.5), et par 
\begin{equation}\label{tf20c}
\delta\colon \tE_0\rightarrow \tE
\end{equation} 
le plongement canonique, c'est-à-dire le morphisme de topos tel que  
$\delta_*\colon \tE_0\rightarrow \tE$ soit le foncteur d'injection canonique. 

On désigne par $\Pt(\tE)$, $\Pt(\tE_{/\sigma^*(Z)})$ et $\Pt(\tE_0)$ les catégories des points de $\tE$,
$\tE_{/\sigma^*(Z)}$ et $\tE_0$, respectivement, et par 
\begin{equation}\label{tf20d}
u\colon \Pt(\tE_{/\sigma^*(Z)})\rightarrow \Pt(\tE) \ \ \ {\rm et}\ \ \ v\colon \Pt(\tE_0)\rightarrow \Pt(\tE)
\end{equation}
les foncteurs induits par $\gamma$ et $\delta$, respectivement. Ces foncteurs sont pleinement fidèles,
et tout point de $\tE$ appartient à l'image essentielle de l'un ou l'autre de ces foncteurs exclusivement 
(\cite{sga4} IV 9.7.2).

\begin{prop}\label{tf21}
Conservons les hypothèses de \ref{tf20}. 
\begin{itemize}
\item[{\rm (i)}] Soit $(\oy\rightsquigarrow \ox)$ un point de $X_\et\gtimes_{X_\et}Y_\et$ \eqref{topfl17}. 
Pour que $\rho(\oy\rightsquigarrow \ox)$ \eqref{tf3b} appartienne à l'image essentielle de $u$ (resp. $v$) \eqref{tf20d}, 
il faut et il suffit que $\ox$ soit à support dans $Z$ (resp. $X_0$). 
\item[{\rm (ii)}] La famille des points de $\tE_{/\sigma^*(Z)}$ (resp. $\tE_0$) définie par la famille des points 
$\rho(\oy\rightsquigarrow \ox)$ de $\tE$ tels que $\ox$ soit à support dans $Z$ (resp. $X_0$)
est conservative. 
\end{itemize}
\end{prop}

(i) En effet, pour que $\rho(\oy\rightsquigarrow \ox)$ appartienne à l'image essentielle de $u$ (resp. $v$), 
il faut et il suffit que $(\sigma^*(Z))_{\rho(\oy\rightsquigarrow \ox)}$ soit un singleton (resp. vide). 
Par ailleurs, on a un isomorphisme
canonique (\cite{agt} (VI.10.18.1))
\begin{equation}\label{tf21a}
(\sigma^*(Z))_{\rho(\oy\rightsquigarrow \ox)}\stackrel{\sim}{\rightarrow} Z_{\ox},
\end{equation}
d'où la proposition.  

(ii) Cela résulte de (i), \ref{tf5} et (\cite{sga4} IV 9.7.3).

\begin{lem}\label{tf22}
Conservons les hypothèses de \ref{tf20}, soit, de plus, $F=\{U\mapsto F_U\}$ un faisceau de $\tE$.
Alors, les propriétés suivantes sont équivalentes~:
\begin{itemize}
\item[{\rm (i)}] $F$ est un objet de $\tE_0$. 
\item[{\rm (ii)}] Pour tout $U\in \Et_{/Z}$, $F_U$ est un objet final de $(U_Y)_\fet$, 
i.e., est représentable par $U_Y$. 
\item[{\rm (iii)}] Pour tout point $(\oy\rightsquigarrow \ox)$ de $X_\et\gtimes_{X_\et}Y_\et$ \eqref{topfl17}
tel que $\ox$ soit à support dans $Z$, la fibre $F_{\rho(\oy\rightsquigarrow \ox)}$ de $F$ en $\rho(\oy\rightsquigarrow \ox)$
est un singleton \eqref{tf3b}. 
\end{itemize}
\end{lem}

En effet, d'après (\cite{agt} VI.5.38), on a un isomorphisme canonique 
\begin{equation}\label{tf22a}
\gamma^*(F)\stackrel{\sim}{\rightarrow} \{W\mapsto F_W\}, \ \ \ (W\in \ob(\Et_{/Z})). 
\end{equation}
Comme $\{W\mapsto W_Y\}$, pour $W\in \ob(\Et_{/Z})$, est un objet final de $\tE_{/\sigma^*(Z)}$ (\cite{agt} III.8.4),
les conditions (i) et (ii) sont équivalentes. Par ailleurs, les conditions (i) et (iii) sont équivalentes en vertu de \ref{tf21}(ii).

\subsection{}\label{tf23}
Conservons les hypothèses et notations de \ref{tf20}, notons, de plus, $i\colon X_0\rightarrow X$ et $j\colon Z\rightarrow X$ les injections canoniques. 
Le topos $\tE_{/\sigma^*(Z)}$ étant canoniquement équivalent au topos de Faltings associé au morphisme $\varpi\colon Y\times_XZ\rightarrow Z$, notons 
\begin{equation}\label{tf23a}
\sigma_Z\colon \tE_{/\sigma^*(Z)}\rightarrow Z_\et
\end{equation}
le morphisme canonique \eqref{tf1kk}. D'après (\cite{agt} (VI.10.12.6)), le diagramme 
\begin{equation}\label{tf23b}
\xymatrix{
{\tE_{/\sigma^*(Z)}}\ar[r]^{\sigma_Z}\ar[d]_{\gamma}&{Z_\et}\ar[d]^{j}\\
{\tE}\ar[r]^\sigma&{X_\et}}
\end{equation}
est commutatif à isomorphisme canonique près. En vertu de (\cite{sga4} IV 9.4.3), il existe un morphisme 
\begin{equation}\label{tf23c}
\sigma_0\colon \tE_0\rightarrow X_{0,\et}
\end{equation}
unique à isomorphisme près tel que le diagramme 
\begin{equation}\label{tf23d}
\xymatrix{
{\tE_0}\ar[r]^-(0.5){\sigma_0}\ar[d]_{\delta}&{X_{0,\et}}\ar[d]^i\\
{\tE}\ar[r]^-(0.5)\sigma&{X_\et}}
\end{equation}
soit commutatif à isomorphisme près. Par définition, on a un isomorphisme canonique 
\begin{equation}\label{tf23f}
\sigma^*\circ i_*\stackrel{\sim}{\rightarrow}\delta_*\circ \sigma_0^*.
\end{equation}
Les foncteurs $i_*$ et $\delta_*$ étant exacts, 
pour tout groupe abélien $F$ de $\tE_0$ et tout entier $q\geq 0$, on a un isomorphisme canonique 
\begin{equation}\label{tf23g}
i_*(\rR^q\sigma_{0*}(F))\stackrel{\sim}{\rightarrow}\rR^q\sigma_*(\delta_*F). 
\end{equation}

\begin{lem}\label{tf25}
Conservons les hypothèses et notations de \ref{tf20}; soient, de plus, $\ox$ un point géométrique de $X$, 
$\uX$ le localisé strict de $X$ en $\ox$, $\uY=Y\times_X\uX$, 
\begin{equation}\label{tf25a}
\varphi_{\ox}\colon \tE\rightarrow \uY_\fet
\end{equation}
le foncteur défini dans \eqref{tf10b}. 
\begin{itemize}
\item[{\rm (i)}] Si $\ox$ est à support dans $Z$, alors pour tout faisceau $F$ de $\tE_0$, $\varphi_\ox(\delta_*(F))$ est un objet final de $\uY_\fet$.
\item[{\rm (ii)}] Si $\ox$ est à support dans $X_0$, alors pour tout faisceau $F$ de $\tE$, le morphisme 
\begin{equation}\label{tf25b}
\varphi_\ox(F)\rightarrow \varphi_\ox(\delta_*(\delta^*(F)))
\end{equation}
induit par le morphisme d'adjonction $F\rightarrow \delta_*(\delta^*(F))$ est un isomorphisme. 
\item[{\rm (iii)}] La famille des foncteurs $(\varphi_\ox\circ \delta_*\colon \tE_0\rightarrow \uY_\fet)$, 
lorsque $\ox$ décrit l'ensemble des points géométriques de $X_0$, est conservative {\rm (\cite{sga4} I 6.1)}.
\end{itemize}
\end{lem}

En effet, tout point géométrique $\oy$ de $\uY$ définit naturellement un point $(\oy\rightsquigarrow \ox)$ de 
$X_\et\gtimes_{X_\et}Y_\et$ \eqref{topfl17}. 
D'après (\cite{agt} VI.10.31), pour tout faisceau $F$ de $\tE$, on a un isomorphisme canonique fonctoriel
\begin{equation}\label{tf25d}
F_{\rho(\oy\rightsquigarrow \ox)}\stackrel{\sim}{\rightarrow}\varphi_\ox(F)_{\rho_{\uY}(\oy)},
\end{equation}
où la source est la fibre de $F$ en $\rho(\oy\rightsquigarrow \ox)$ \eqref{tf3b} et le but est la fibre de $\varphi_\ox(F)$ en le  point 
$\rho_{\uY}(\oy)$ de $\uY_\fet$, $\rho_\uY\colon \uY_\et\rightarrow \uY_\fet$ étant le morphisme canonique \eqref{notconv10a}.

(i) Cela résulte de \ref{tf22}, \eqref{tf25d} et (\cite{agt} VI.9.6). 

(ii) D'après \ref{tf21}(i), pour tout point géométrique $\oy$ de $\uY$,
$\rho(\oy\rightsquigarrow \ox)$ appartient à l'image essentielle du foncteur $\Pt(\tE_0)\rightarrow \Pt(\tE)$ défini par $\delta$.  
Par suite, le morphisme 
\begin{equation}\label{tf25e}
F_{\rho(\oy\rightsquigarrow \ox)}\rightarrow \delta_*(\delta^*(F))_{\rho(\oy\rightsquigarrow \ox)}
\end{equation}
induit par le morphisme d'adjonction $F\rightarrow \delta_*(\delta^*(F))$ est un isomorphisme. La proposition s'ensuit compte tenu de \eqref{tf25d} et (\cite{agt} VI.9.6). 

(iii) Cela résulte de \ref{tf21}(ii) et \eqref{tf25d}.

\begin{prop}\label{tf26}
Conservons les hypothèses et notations de \ref{tf23}. 
\begin{itemize}
\item[{\rm (i)}] Pour tout faisceau $F$ de $\tE$, le morphisme de changement de base relativement 
au diagramme commutatif de morphismes de topos \eqref{tf23d},
\begin{equation}\label{tf26a}
i^*(\sigma_*(F))\rightarrow \sigma_{0*}(\delta^*(F))
\end{equation}
est un isomorphisme. 
\item[{\rm (ii)}] Pour tout faisceau abélien $F$ de $\tE$ et tout entier $q\geq 0$, le morphisme de changement de base 
relativement au diagramme commutatif de morphismes de topos \eqref{tf23d}, 
\begin{equation}\label{tf26b}
i^*(\rR^q\sigma_*(F))\rightarrow \rR^q\sigma_{0*}(\delta^*(F))
\end{equation}
est un isomorphisme. 
\end{itemize}
\end{prop}

On renvoie à (\cite{egr1} 1.2.2) pour la définition du morphisme de changement de base. 

(i) Le morphisme \eqref{tf26a} est l'adjoint du morphisme composé 
\begin{equation}\label{tf26c}
\sigma_*(F)\rightarrow \sigma_*(\delta_*(\delta^*(F)))\stackrel{\sim}{\rightarrow} i_*(\sigma_{0*}(\delta^*(F))),
\end{equation}
où la première flèche est induite par le morphisme d'adjonction $\id\rightarrow \delta_*\delta^*$. Comme $i$ est une immersion, 
il suffit de montrer que la première flèche induit un isomorphisme sur $X_0$. Soient $\ox$ un point géométrique de $X_0$, 
$\uX$ le localisé strict de $X$ en $\ox$, $\uY=Y\times_X\uX$. On note 
\begin{equation}\label{tf26d}
\varphi_\ox\colon \tE\rightarrow \uY_\fet
\end{equation}
le foncteur défini dans \eqref{tf10b}. En vertu de (\cite{agt} VI.10.30(i)), on a un isomorphisme canonique et fonctoriel
\begin{equation}\label{tf26e}
\sigma_*(F)_\ox\stackrel{\sim}{\rightarrow}\Gamma(\uY_\fet,\varphi_\ox(F)). 
\end{equation}
Par ailleurs, d'après \ref{tf25}(ii), le morphisme d'adjonction
\begin{equation}\label{tf26g}
\varphi_\ox(F)\rightarrow \varphi_\ox(\delta_*(\delta^*(F)))
\end{equation}
est un isomorphisme. Le morphisme d'adjonction $\sigma_*(F)_\ox\rightarrow \sigma_*(\delta_*(\delta^*(F)))_\ox$ est donc un isomorphisme en $\ox$; d'où l'assertion recherchée.

(ii) Le morphisme \eqref{tf26b} est l'adjoint du morphisme composé 
\begin{equation}\label{tf26h}
\rR^q\sigma_*(F)\rightarrow \rR^q \sigma_*(\delta_*(\delta^*(F)))\stackrel{\sim}{\rightarrow} i_*(\rR^q\sigma_{0*}(\delta^*(F))),
\end{equation}
où la première flèche est induite par le morphisme d'adjonction $\id\rightarrow \delta_*\delta^*$, et la seconde flèche est l'inverse de \eqref{tf23g}. 
Il suffit alors de calquer la preuve de (i), compte tenu de (\cite{agt} VI.10.30(ii)).

\begin{prop}\label{tf27}
Soient $S$ un schéma strictement local de point fermé $s$, $g\colon X\rightarrow S$ un morphisme propre, 
$Z$ un sous-schéma ouvert de $X$ tel que $Z\cap X_s=\emptyset$. 
On désigne par $\tE_0$ le sous-topos fermé de $\tE$ complémentaire de l'ouvert $\sigma^*(Z)$ et par $\delta\colon \tE_0\rightarrow \tE$
le plongement canonique \eqref{tf20}. Alors,
\begin{itemize}
\item[{\rm (i)}] Pour tout faisceau $F$ de $\tE$, le morphisme canonique 
\begin{equation}\label{tf27a}
\Gamma(\tE,F)\rightarrow \Gamma(\tE_0,\delta^*(F))
\end{equation}
est un isomorphisme. 
\item[{\rm (ii)}] Pour tout faisceau abélien de torsion $F$ de $\tE$ et tout entier $q\geq 0$, 
le morphisme canonique 
\begin{equation}\label{tf27b}
\rH^q(\tE,F)\rightarrow \rH^q(\tE_0,\delta^*(F))
\end{equation}
est un isomorphisme. 
\end{itemize}
\end{prop}

On notera d'abord que comme le foncteur $\delta_*$ est exact, le morphisme \eqref{tf27b} n'est autre que l'image du morphisme d'adjonction 
$F\rightarrow \delta_*(\delta^*(F))$ par le foncteur cohomologique $\rH^q(\tE,-)$. 

Notons $\eta$ l'ouvert complémentaire de $s$ dans $S$, $X_\eta=X\times_S\eta$, $\tE_s$ le sous-topos fermé de $\tE$ complémentaire de l'ouvert 
$\sigma^*(X_\eta)$ et $\iota\colon \tE_s\rightarrow \tE$ le plongement canonique. Comme $Z\subset X_\eta$, il existe un morphisme 
$\lambda\colon \tE_s\rightarrow \tE_0$ tel que $\iota=\delta\circ \lambda$ (\cite{sga4} IV 9.4.2). 
Pour tout faisceau $F$ de $\tE$, le morphisme d'adjonction $F\rightarrow \delta_*(\delta^*(F))$ induit un isomorphisme
\begin{equation}
\iota^*(F)\stackrel{\sim}{\rightarrow} \iota^*(\delta_*(\delta^*(F))).
\end{equation}
Par suite, la proposition relative à l'ouvert $X_\eta$ de $X$, appliquée d'une part au faisceau $F$ et d'autre part au faisceau $\delta_*(\delta^*(F))$, 
implique la proposition recherchée pour l'ouvert $Z$ de $X$ et le faisceau $F$. 
On peut donc se borner au cas où $Z=X_\eta$.  On note $i\colon X_s\rightarrow X$ l'injection canonique et 
\begin{equation}
\sigma_s\colon \tE_s\rightarrow X_{s,\et}
\end{equation}
le morphisme induit par $\sigma$ (cf. \ref{tf23}). 

(i) Compte tenu de \eqref{tf23d}, on a un diagramme commutatif 
\begin{equation}\label{tf27c}
\xymatrix{
{\Gamma(\tE,F)}\ar[r]^-(0.5)\sim\ar[d]_u&{\Gamma(X_\et,\sigma_*(F))}\ar[d]^v\\
{\Gamma(\tE_s,\delta^*(F))}\ar[r]^-(0.5)\sim&{\Gamma(X_{s,\et},\sigma_{s*}(\delta^*(F)))}}
\end{equation}
où les flèches horizontales sont les isomorphismes canoniques, $u$ est induit par le morphisme d'adjonction $F\rightarrow \delta_*(\delta^*(F))$
et $v$ est induit par le morphisme de changement de base \eqref{tf26a}
\begin{equation}
i^*(\sigma_*(F))\rightarrow \sigma_{s*}(\delta^*(F)). 
\end{equation}
Comme ce dernier est un isomorphisme
d'après \ref{tf26}(i), $v$ est un isomorphisme en vertu du théorème de changement de base propre en cohomologie étale (\cite{sga4} XII 5.5). 
Par suite, $u$ est un isomorphisme. 

(ii) Considérons les suites spectrales de Cartan-Leray (\cite{sga4} V 5.4)
\begin{eqnarray}
\rE_2^{a,b}=\rH^a(X_\et,\rR^b \sigma_*(F))&\Rightarrow& \rH^{a+b}(\tE,F),\label{tf27d}\\
{_s\rE}_2^{a,b}=\rH^a(X_{s,\et},\rR^b \sigma_{s*}(\delta^*(F)))&\Rightarrow& \rH^{a+b}(\tE_s,\delta^*(F)).\label{tf27e}
\end{eqnarray}
Le morphisme de changement de base \eqref{tf26b} 
\begin{equation}\label{tf27f}
i^*(\rR^b \sigma_*(F))\rightarrow \rR^b\sigma_{s*}(\delta^*(F))
\end{equation} 
induit un morphisme de suites spectrales
\begin{equation}\label{tf27g}
\rE_2^{a,b}\rightarrow {_s\rE}_2^{a,b}
\end{equation}
dont le morphisme induit entre aboutissements 
\begin{equation}\label{tf27h}
\rH^{a+b}(\tE,F)\rightarrow \rH^{a+b}(\tE_0,\delta^*(F))
\end{equation}
est aussi induit par le morphisme d'adjonction $F\rightarrow \delta_*(\delta^*(F))$.
Comme \eqref{tf27f} est un isomorphisme d'après \ref{tf26}(ii), 
\eqref{tf27g} est un isomorphisme en vertu du théorème de changement de base propre en cohomologie étale (\cite{sga4} XII 5.5). 
Par suite, \eqref{tf27h} est un isomorphisme.

\subsection{}\label{tf28}
Considérons un diagramme commutatif de morphismes de schémas
\begin{equation}\label{tf28a}
\xymatrix{
Y'\ar[r]^{g'}\ar[d]_{f'}&Y\ar[d]^f\\
X'\ar[r]^g&X}
\end{equation}
Soient $X_0$ un sous-schéma fermé de $X$, $Z$ l'ouvert complémentaire de $X_0$ dans $X$, $Z'=g^{-1}(Z)$.  
On note $E'$ et $\tE'$ le site et le topos de Faltings associés à $f'$ et
\begin{equation}\label{tf28b}
\sigma'\colon \tE'\rightarrow X'_\et
\end{equation} 
le morphisme canonique \eqref{tf1kk}.  
On désigne par $\tE_0$ (resp. $\tE'_0$) le sous-topos fermé de $\tE$ (resp. $\tE'$)
complémentaire de l'ouvert $\sigma^*(Z)$ (resp. $\sigma'^*(Z')$), $\delta\colon \tE_0\rightarrow \tE$ et $\delta'\colon \tE'_0\rightarrow \tE'$
les plongements canoniques (cf. \ref{tf20}). Le foncteur 
\begin{equation}\label{tf28c}
\Theta^+\colon E\rightarrow E', \ \ \ (V\rightarrow U)\mapsto (V\times_{Y}Y'\rightarrow U\times_XX')
\end{equation}
est continu et exact à gauche (\cite{agt} VI.10.12). Il définit donc un morphisme de topos 
\begin{equation}\label{tf28d}
\Theta\colon \tE'\rightarrow \tE.
\end{equation}
Il résulte aussitôt des définitions que le diagramme
\begin{equation}\label{tf28e}
\xymatrix{
{\tE'}\ar[r]^-(0.5){\sigma'}\ar[d]_{\Theta}&{X'_\fet}\ar[d]^{g}\\
{\tE}\ar[r]^{\sigma}&{X_\fet}}
\end{equation}
est commutatif à isomorphisme canonique près.  
On en déduit un isomorphisme 
\begin{equation}\label{tf28f}
\Theta^*(\sigma^*(Z))\stackrel{\sim}{\rightarrow} \sigma'^*(Z').
\end{equation}
En vertu de (\cite{sga4} IV 9.4.3), il existe donc un morphisme de topos
\begin{equation}\label{tf28g}
\Theta_0\colon \tE'_0\rightarrow \tE_0
\end{equation}
unique à isomorphisme canonique près tel que le diagramme 
\begin{equation}\label{tf28h}
\xymatrix{
{\tE'_0}\ar[r]^{\Theta_0}\ar[d]_{\delta'}&{\tE_0}\ar[d]^{\delta}\\
{\tE'}\ar[r]^\Theta&{\tE}}
\end{equation}
soit commutatif à isomorphisme près, et même $2$-cartésien. Les foncteurs $\delta_*$ et $\delta'_*$ étant exacts, 
pour tout groupe abélien $F$ de $\tE_0$ et tout entier $q\geq 0$, on a un isomorphisme canonique 
\begin{equation}\label{tf28i}
\rR^q\Theta_*(\delta'_*F)\stackrel{\sim}{\rightarrow} \delta_*(\rR^q\Theta_{0*}(F)). 
\end{equation}

\begin{prop}\label{tf30}
Conservons les hypothèses et notations de \ref{tf28}; supposons, de plus, $g$ étale et le morphisme canonique $Y'\rightarrow Y\times_XX'$ fini étale, 
de sorte que $(Y'\rightarrow X')$ s'identifie à un objet de $E$. Alors, 
\begin{itemize}
\item[{\rm (i)}] Le morphisme $\Theta$ s'identifie au morphisme de localisation de $\tE$ en $(Y'\rightarrow X')^a$.
\item[{\rm (ii)}] Le morphisme $\Theta_0$ s'identifie au morphisme de localisation de $\tE_0$ en $\delta^*(Y'\rightarrow X')^a$. 
\item[{\rm (iii)}] Pour tout faisceau $F$ de $\tE_0$, le morphisme de changement de base relativement 
au diagramme commutatif \eqref{tf28h},
\begin{equation}\label{tf30a}
\Theta^*(\delta_*(F))\rightarrow \delta'_*(\Theta_0^*(F))
\end{equation}
est un isomorphisme. 
\end{itemize}
\end{prop}

(i) Cela résulte de (\cite{agt} VI.10.14). 

(ii) D'après (\cite{sga4} IV 5.11), le diagramme 
\begin{equation}
\xymatrix{
{(\tE_0)_{/\delta^*(Y'\rightarrow X')^a}}\ar[r]^-(0.5){\jmath}\ar[d]_{\delta^\dagger}&{\tE_0}\ar[d]^{\delta}\\
{\tE'}\ar[r]^\Theta&{\tE}}
\end{equation}
où $\jmath$ est le morphisme de localisation de $\tE_0$ en $\delta^*(Y'\rightarrow X')^a$ et $\delta^\dagger$ est le 
morphisme induit par $\delta$ (\cite{sga4} IV 5.10) est $2$-cartésien. Comme le diagramme \eqref{tf28h} est aussi $2$-cartésien, 
on peut identifier $\Theta_0$ à $\jmath$, et $\delta'$ à $\delta^\dagger$, d'où la proposition.

(iii) Cela résulte aussitôt de la preuve de (ii) et (\cite{sga4} V 5.1(3)).

\begin{cor}\label{tf31}
Soient $S$ un schéma strictement local de point fermé $s$, $X\rightarrow S$ un morphisme propre, 
$Z$ un sous-schéma ouvert de $X$ tel que $Z\cap X_s=\emptyset$, $g'\colon Y'\rightarrow Y$ un morphisme fini étale de sorte que 
$(Y'\rightarrow X)$ est un objet de $E$. 
On désigne par $\tE_0$ le sous-topos fermé de $\tE$ complémentaire de l'ouvert $\sigma^*(Z)$ et par $\delta\colon \tE_0\rightarrow \tE$
le plongement canonique \eqref{tf20}. Alors,
\begin{itemize}
\item[{\rm (i)}] Pour tout faisceau $F$ de $\tE$, le morphisme d'adjonction 
\begin{equation}\label{tf31a}
\Gamma((Y'\rightarrow X),F)\rightarrow \Gamma((Y'\rightarrow X),\delta_*(\delta^*(F)))
\end{equation}
est un isomorphisme. 
\item[{\rm (ii)}] Pour tout faisceau abélien de torsion $F$ de $\tE$ et tout entier $q\geq 0$, 
le morphisme d'adjonction
\begin{equation}\label{tf31b}
\rH^q((Y'\rightarrow X),F)\rightarrow \rH^q((Y'\rightarrow X),\delta_*(\delta^*(F)))
\end{equation}
est un isomorphisme. 
\end{itemize}
\end{cor}

Reprenons les notations de \ref{tf28} relativement au diagramme commutatif 
\begin{equation}\label{tf31c}
\xymatrix{
Y'\ar[r]^{g'}\ar[d]_{f'}&Y\ar[d]^f\\
X\ar[r]^{\id}&X}
\end{equation} 
D'après \ref{tf30}(i)-(ii), le morphisme $\Theta$ (resp. $\Theta_0$) s'identifie au morphisme de localisation de $\tE$ en $(Y'\rightarrow X)^a$
(resp. $\tE_0$ en $\delta^*((Y'\rightarrow X)^a)$). D'après (\cite{sga4} XVII 2.1.3), le morphisme de changement de base 
$u\colon \Theta^*\delta_*\rightarrow \delta'_*\Theta^*_0$ relativement à \eqref{tf28h} est l'adjoint du morphisme composé 
\begin{equation}
\delta'^*\Theta^*\delta_*\stackrel{\sim}{\rightarrow} \Theta_0^*\delta^*\delta_*\rightarrow \Theta_0^*,
\end{equation}
où la première flèche est l'isomorphisme qui fait commuter le diagramme \eqref{tf28h} et la deuxième flèche est le morphisme d'adjonction. 
Par suite, pour tout faisceau $F$ de $\tE$, le diagramme 
\begin{equation}
\xymatrix{
{\Theta^*(F)}\ar[rr]^-(0.5){w}\ar[d]_v&&{\delta'_*(\delta'^*(\Theta^*(F)))}\ar[d]^{\delta'_*(t)}\\
{\Theta^*(\delta_*(\delta^*(F)))}\ar[rr]^-(0.5){u(\delta^*(F))}&&{\delta'_*(\Theta^*_0(\delta^*(F)))}}
\end{equation}
où les morphismes $v$ et $w$ sont induits par les flèches d'adjonction, $t$ est l'isomorphisme qui fait commuter le diagramme \eqref{tf28h}, est commutatif.
En vertu de \ref{tf30}(iii), $u$ est un isomorphisme. 
D'après \ref{tf27}(i),  $\Gamma(\tE',w)$ est un isomorphisme. Il en est donc de même de  $\Gamma(\tE',v)$, d'où la proposition (i). 
La preuve de la proposition (ii) est similaire.

\begin{prop}\label{tf29}
Conservons les hypothèses et notations de \ref{tf28}; supposons, de plus, $f$ et $f'$ cohérents et $g$ propre. Alors, 
\begin{itemize}
\item[{\rm (i)}] Pour tout faisceau $F$ de $\tE'$, le morphisme de changement de base relativement 
au diagramme commutatif \eqref{tf28h},
\begin{equation}\label{tf29a}
\delta^*(\Theta_*(F))\rightarrow \Theta_{0*}(\delta'^*(F))
\end{equation}
est un isomorphisme. 
\item[{\rm (ii)}] Pour tout faisceau abélien $F$ de $\tE'$ et tout entier $q\geq 0$, le morphisme de changement de base 
relativement au diagramme commutatif \eqref{tf28h}, 
\begin{equation}\label{tf29b}
\delta^*(\rR^q\Theta_*(F))\rightarrow \rR^q\Theta_{0*}(\delta'^*(F))
\end{equation}
est un isomorphisme. 
\end{itemize}
\end{prop}

Nous démontrons seulement (ii); la preuve de (i) est similaire et plus simple. 

Soient $\ox$ un point géométrique de $X_0$, $\uX$ le localisé strict de $X$ en $\ox$, $u\colon \uX\rightarrow X$ le morphisme canonique. On note
\begin{equation}\label{tf29c}
\Sch_{/X}\rightarrow \Sch_{/\uX}, \ \ \ U\mapsto \uU
\end{equation}
le foncteur de changement de base par $u$ \eqref{notconv3}. On a un diagramme commutatif canonique
\begin{equation}\label{tf29d}
\xymatrix{
{\uY'}\ar[r]^{\ug'}\ar[d]_{\uf'}&{\uY}\ar[d]^{\uf}\\
{\uX'}\ar[r]^{\ug}&{\uX}}
\end{equation} 
On désigne par $\uE$ et $\uE'$ (resp. $\tuE$ et $\tuE'$) les sites (resp. topos) 
de Faltings associés aux morphismes $\uf$ et $\uf'$ respectivement \eqref{tf1}, par 
\begin{eqnarray}
\usigma\colon \tuE\rightarrow \uX_\et,\label{tf29e}\\
\usigma'\colon \tuE'\rightarrow \uX'_\et,\label{tf29ee}\\
\urho\colon \uX_\et\gtimes_{\uX_\et}\uY_\et\rightarrow \tuE, \label{tf29eee}
\end{eqnarray}
les morphismes canoniques \eqref{tf1kk} et \eqref{tf3b}, par $\tuE_0$ (resp. $\tuE'_0$) le sous-topos fermé de $\tuE$ (resp. $\tuE'$)
complémentaire de l'ouvert $\sigma^*(\uZ)$ (resp. $\sigma'^*(\uZ')$) et par $\udelta\colon \tuE_0\rightarrow \tE$ et $\udelta'\colon \tuE'_0\rightarrow \tuE'$
les plongements canoniques (cf. \ref{tf20}). On note
\begin{eqnarray}
\uTheta\colon \tuE'&\rightarrow& \tuE,\label{tf29f}\\
\iota\colon \tuE&\rightarrow& \tE,\label{tf29ff}\\
\iota'\colon \tuE'&\rightarrow& \tE',\label{tf29fff}
\end{eqnarray}
les morphismes de topos définis par fonctorialité (\cite{agt} VI.10.12). On a des isomorphismes canoniques
\begin{eqnarray}
\uTheta^*(\usigma^*(\uZ))&\stackrel{\sim}{\rightarrow}& \usigma'^*(\uZ'),\label{tf29g}\\
\iota^*(\sigma^*(Z))&\stackrel{\sim}{\rightarrow}& \usigma^*(\uZ),\label{tf29gg}\\
\iota'^*(\sigma'^*(Z'))&\stackrel{\sim}{\rightarrow}& \usigma'^*(\uZ').\label{tf29ggg}
\end{eqnarray}
En vertu de (\cite{sga4} IV 9.4.3), il existe donc des morphismes de topos
\begin{eqnarray}
\uTheta_0\colon \tuE'_0&\rightarrow& \tuE_0,\label{tf29h}\\
\iota_0\colon \tuE_0&\rightarrow& \tE_0,\label{tf29hh}\\
\iota'_0\colon \tuE'_0&\rightarrow& \tE_0,\label{tf29hhh}
\end{eqnarray}
uniques à isomorphismes canoniques près tels que dans le diagramme 
\begin{equation}\label{tf29i}
\xymatrix{
{\tuE'_0}\ar[rrr]^-(0.5){\udelta'}\ar[rd]^{\iota'_0}\ar[ddd]_{\uTheta_0}\ar@{}[drrr]|*+[o][F-]{2}&\ar@{}[dddl]|*+[o][F-]{1}&\ar@{}[dddr]|*+[o][F-]{3}&{\tuE'}\ar[ddd]^{\uTheta}\ar[ld]_{\iota'}\\
&{\tE'_0}\ar[r]^-(0.5){\delta'}\ar[d]_{\Theta_0}\ar@{}[dr]|*+[o][F-]{5}&{\tE'}\ar[d]^{\Theta}&\\
\ar@{}[drrr]|*+[o][F-]{4}&{\tE_0}\ar[r]^-(0.5)\delta&{\tE}&\\
{\tuE_0}\ar[ru]_{\iota_0}\ar[rrr]^-(0.5){\udelta}&&&{\tuE}\ar[lu]^{\iota}}
\end{equation}
les faces $(2)$ et $(4)$ et la face représentée par le grand carré extérieur, notée dans la suite $(6)$, 
soient commutatives à isomorphismes près, et même $2$-cartésiennes. Par ailleurs, la face $(5)$ est commutative à isomorphismes près \eqref{tf28h} 
et la face $(3)$ est commutative à isomorphisme près (\cite{agt} VI.10.12). 
On en déduit que la face (1) est commutative à isomorphisme près (\cite{sga4} IV 9.4.2).

D'après (\cite{egr1} 1.2.4(ii)), pour tout faisceau abélien $F$ de $\tE'$ et tout entier $q\geq 0$, le diagramme 
\begin{equation}\label{tf29j}
\xymatrix{
{\iota_0^*(\delta^*(\rR^q\Theta_*(F)))}\ar[r]^{w_5}\ar[d]_{w_4}&{\iota_0^*(\rR^q\Theta_{0*}(\delta'^*(F)))}\ar[r]^{w_1}&
{\rR^q\uTheta_{0*}(\iota'^*_0(\delta'^*(F)))}\ar[d]^{w_2}\\
{\udelta^*(\iota^*(\rR^q\Theta_*(F)))}\ar[r]^{w_3}&{\udelta^*(\rR^q\uTheta_*(\iota'^*(F)))}\ar[r]^{w_6}&{\rR^q\uTheta_{0*}(\delta'^*(\iota'^*(F)))}}
\end{equation}
où $w_n$, pour $n\in \{1,3,5,6\}$, est le morphisme de changement de base relativement à la face $(n)$, et $w_n$, pour $n\in \{2,4\}$, est 
l'isomorphisme qui fait commuter la face $(n)$, est commutatif. 
Le morphisme $w_3$ est un isomorphisme d'après \ref{tf16}(ii). Le morphisme $w_1$ est un isomorphisme en vertu de \ref{tf32}(iii). 
Pour que $w_5$ soit un isomorphisme, il faut et il suffit donc que $w_6$ le soit. 

Soient $\oy$ un point géométrique de $Y$, $(\oy\rightsquigarrow \ox)$ 
un point de $X_\et\gtimes_{X_\et}Y_\et$ défini par un $X$-morphisme $\tau\colon \oy\rightarrow \uX$ \eqref{topfl17}. 
D'après \ref{tf21}, il suffit de montrer que la fibre du morphisme \eqref{tf29b} 
en le point de $\tE_0$ d'image canonique le point $\rho(\oy\rightsquigarrow \ox)$ de $\tE$ \eqref{tf3b} est un isomorphisme. 
Le morphisme $\tau$ induit un $\uX$-morphisme $\oy\rightarrow \uY$. Notant encore (abusivement) $\oy$ le point géométrique de $\uY$ 
ainsi défini et $(\oy\rightsquigarrow \ox)$ le point de $\uX_\et\gtimes_{\uX_\et}\uY_\et$ défini par $\tau$, 
on a $\iota(\urho(\oy\rightsquigarrow \ox))=\rho(\oy\rightsquigarrow \ox)$ (\cite{agt} (VI.10.31.2)). Compte tenu de \ref{tf21}(i), 
considérant $\urho(\oy\rightsquigarrow \ox)$ et $\rho(\oy\rightsquigarrow \ox)$ comme des points de $\tuE_0$ et $\tE_0$ respectivement, 
on a $\iota_0(\urho(\oy\rightsquigarrow \ox))=\rho(\oy\rightsquigarrow \ox)$. En vertu de \eqref{tf29j}, 
on peut donc se réduire au cas où $X$ est strictement local de point fermé $\ox$.

Soient $\rho_Y\colon Y_\et\rightarrow Y_\fet$ le morphisme canonique \eqref{notconv10a}, 
$(V_i)_{i\in I}$ un revêtement universel de $Y$ en $\oy$ \eqref{notconv11},
$G=\{U\mapsto G_U\}$ $(U\in \ob(\Et_{/X}))$ un préfaisceau d'ensembles sur $E$ \eqref{tf1h}, 
$G^a$ le faisceau de $\tE$ associé à $G$,
$G^a_X$ le faisceau de $Y_\fet$ associé à $G_X$. 
D'après (\cite{agt} VI.10.36), $X$ étant strictement local de point fermé $\ox$, on a un isomorphisme canonique fonctoriel
\begin{equation}\label{tf29l}
(G^a)_{\rho(\oy\rightsquigarrow \ox)}\stackrel{\sim}{\rightarrow} (G^a_X)_{\rho_Y(\oy)}.
\end{equation}
On en déduit un isomorphisme canonique fonctoriel (\cite{sga4} IV (6.8.4))
\begin{equation}\label{tf29m}
(G^a)_{\rho(\oy\rightsquigarrow \ox)}\stackrel{\sim}{\rightarrow}\underset{\underset{i\in I^\circ }{\longrightarrow}}{\lim}\ G(V_i\rightarrow X).
\end{equation}

Soient $F$ un groupe abélien de $\tE$, $q$ un entier $\geq 0$. 
Le morphisme \eqref{tf29b} est l'adjoint du morphisme composé 
\begin{equation}\label{tf29p}
\rR^q\Theta_*(F)\rightarrow \rR^q\Theta_*(\delta'_*(\delta'^*(F)))\stackrel{\sim}{\rightarrow} \delta_*(\rR^q\Theta_{0*}(\delta'^*(F))),
\end{equation}
où la première flèche est induite par le morphisme d'adjonction $\id\rightarrow \delta'_*\delta'^*$ et la seconde flèche est l'isomorphisme \eqref{tf28i}. 
D'après \eqref{tf29m} et (\cite{sga4} V 5.1), on a un isomorphisme canonique fonctoriel
\begin{equation}\label{tf29q}
\rR^q\Theta_*(F)_{\rho(\oy\rightsquigarrow \ox)}\stackrel{\sim}{\rightarrow}\underset{\underset{i\in I^\circ }{\longrightarrow}}{\lim}\ 
\rH^q((V_i\times_{Y}Y'\rightarrow X'), F).
\end{equation}

En vertu de \ref{tf31}(ii), pour tout $i\in I$, le morphisme d'adjonction 
\begin{equation}
\rH^q((V_i\times_{Y}Y'\rightarrow X'), F)\rightarrow \rH^q((V_i\times_{Y}Y'\rightarrow X'), \delta'_*(\delta'^*(F)))
\end{equation}
est un isomorphisme. Il s'ensuit que la fibre de la première flèche de \eqref{tf29p} en $\rho(\oy\rightsquigarrow \ox)$ est un isomorphisme. 
Par suite, la fibre du morphisme \eqref{tf29b} 
en le point de $\tE_0$ d'image canonique le point $\rho(\oy\rightsquigarrow \ox)$ de $\tE$ est un isomorphisme.

\begin{lem}\label{tf32}
Conservons les hypothèses et notations de la preuve de \ref{tf29}.
\begin{itemize}
\item[{\rm (i)}] Pour tout faisceau $F$ de $\tE_0$, le morphisme de changement de base relativement 
à la face $(4)$ du diagramme \eqref{tf29i},
\begin{equation}\label{tf32a}
\iota^*(\delta_*(F))\rightarrow \udelta_*(\iota^*_0(F))
\end{equation}
est un isomorphisme. 
\item[{\rm (ii)}] Pour tout faisceau $F$ de $\tE'_0$, le morphisme de changement de base relativement 
à la face $(1)$ du diagramme \eqref{tf29i},
\begin{equation}\label{tf32b}
\iota_0^*(\Theta_{0*}(F))\rightarrow \uTheta_{0*}(\iota'^*_0(F))
\end{equation}
est un isomorphisme. 
\item[{\rm (iii)}] Pour tout faisceau abélien $F$ de $\tE'$ et tout entier $q\geq 0$, le morphisme de changement de base 
relativement à la face $(1)$ du diagramme  \eqref{tf29i}, 
\begin{equation}\label{tf32c}
\iota_0^*(\rR^q\Theta_{0*}(F))\rightarrow \rR^q\uTheta_{0*}(\iota'^*_0(F))
\end{equation}
est un isomorphisme. 
\end{itemize}
\end{lem}

(i) D'après (\cite{sga4} XVII 2.1.3), le morphisme \eqref{tf32a} est l'adjoint du morphisme composé
\begin{equation}\label{tf32d}
\udelta^*(\iota^*(\delta_*(F)))\stackrel{\sim}{\rightarrow} \iota_0^*(\delta^*(\delta_*(F)))\rightarrow \iota_0^*(F),
\end{equation}
où la première flèche est l'isomorphisme qui fait commuter la face $(4)$ du diagramme  \eqref{tf29i} et la seconde flèche est induite 
par le morphisme d'adjonction $\delta^*\delta_*\rightarrow \id$. Comme $\delta$ est un plongement, cette dernière flèche est un isomorphisme et il en est donc de
même du composé \eqref{tf32d}. 

Par suite, le morphisme \eqref{tf32a} est le composé
\begin{equation}\label{tf32e}
\iota^*(\delta_*(F)) \rightarrow \udelta_*(\udelta^*(\iota^*(\delta_*(F))))\stackrel{\sim}{\rightarrow} \udelta_*(\iota_0^*(F)),
\end{equation}
où la première flèche est induite par le morphisme d'adjonction $\id\rightarrow \udelta_*\udelta^*$ et la seconde flèche est l'image par
$\udelta_*$ de l'isomorphisme composé \eqref{tf32d}. D'après (\cite{sga4} IV 5.10) et compte tenu de \eqref{tf29gg}, 
le morphisme $\iota$ induit un morphisme de topos $\iota_{/\sigma^*(Z)}$
qui s'insère dans un diagramme commutatif à isomorphisme canonique près
\begin{equation}
\xymatrix{
{\tuE_{/\usigma^*(\uZ)}}\ar[rr]^{\iota_{/\sigma^*(Z)}}\ar[d]_{\ugamma}&&{\tE_{/\sigma^*(Z)}}\ar[d]^\gamma\\
{\tuE}\ar[rr]^{\iota}&&{\tE}}
\end{equation}
où les flèches verticales sont les morphismes de localisation. 
On a donc un isomorphisme canonique 
\begin{equation}
\ugamma^*(\iota^*(\delta_*(F)))\stackrel{\sim}{\rightarrow}\iota_{/\sigma^*(Z)}^*(\gamma^*(\delta_*(F))).
\end{equation}
Par suite, $\ugamma^*(\iota^*(\delta_*(F)))$ est un objet final de $\tuE_{/\usigma^*(\uZ)}$. La première flèche de \eqref{tf32e} est donc un isomorphisme;
d'où la proposition.

(ii) La preuve est similaire à celle de (iii) ci-dessous et est laissée au lecteur. 

(iii) D'après (\cite{egr1} 1.2.4(v)), le diagramme
\begin{equation}
\xymatrix{
{\iota^*(\rR^q\Theta_*(\delta'_*(F)))}\ar[r]^-(0.5){c_3}\ar@{=}[d]&{\rR^q\uTheta_*(\iota'^*(\delta'_*(F)))}\ar[r]^-(0.5){c_2}&{\rR^q\uTheta_*(\udelta'_*(\iota'^*_0(F)))}\ar@{=}[d]\\
{\iota^*(\delta_*(\rR^q\Theta_{0*}(F)))}\ar[r]^-(0.5){c_4}&{\udelta_*(\iota^*_0(\rR^q\Theta_{0*}(F)))}\ar[r]^-(0.5){c_1}&{\udelta_*(\rR^q\uTheta_{0*}(\iota'_{0*}(F)))}}
\end{equation}
où $c_n$, pour $n\in \{1,2,3,4\}$, est le morphisme de changement de base relativement à la face $(n)$, et les identifications verticales sont induites
par \eqref{tf28i} et son analogue pour $\uTheta$, est commutatif. Le morphisme $c_3$ est un isomorphisme en vertu de \ref{tf16}(ii). 
Les morphismes $c_2$ et $c_4$ sont des isomorphismes d'après (i).  Par suite, $c_1$ est un isomorphisme, d'où la proposition.

\begin{rema}
La proposition \ref{tf29} qui généralise \ref{tf27} dans un cadre relatif, ne sera pas utilisée dans la suite de ce livre.
\end{rema}

\section{Topos de Faltings relatif}\label{tfr}

\subsection{}\label{tfr1}
Dans cette section, $f\colon Y\rightarrow X$ et $g\colon Z\rightarrow X$ désignent deux morphismes de schémas. On désigne par $G$ la catégorie 
des triplets $(U, V\rightarrow U, W\rightarrow U)$ formés d'un $X$-schéma $U$ et de deux 
morphismes $V\rightarrow U$ et $W\rightarrow U$ au-dessus de $f$ et $g$ respectivement, c'est-à-dire des diagrammes commutatifs 
de morphismes de schémas 
\begin{equation}\label{tfr1a}
\xymatrix{
W\ar[r]\ar[d]&U\ar[d]&V\ar[l]\ar[d]\\
Z\ar[r]&X&Y\ar[l]}
\end{equation}
tels que les morphismes $U\rightarrow X$ et $W\rightarrow Z$ soient étales et que le morphisme $V\rightarrow U\times_XY$ soit fini étale; un tel objet sera 
noté $(W\rightarrow U\leftarrow V)$. Soient $(W\rightarrow U\leftarrow V)$ et $(W'\rightarrow U'\leftarrow V')$ deux objets de $G$. 
Un morphisme de $(W'\rightarrow U'\leftarrow V')$ dans $(W\rightarrow U\leftarrow V)$ est la donnée de trois morphismes 
$U'\rightarrow U$, $V'\rightarrow V$ et $W'\rightarrow W$ au-dessus de $X$, $Y$ et $Z$ respectivement, qui rendent le diagramme 
\begin{equation}\label{tfr1b}
\xymatrix{
W'\ar[r]\ar[d]&U'\ar[d]&V'\ar[l]\ar[d]\\
W\ar[r]&U&V\ar[l]}
\end{equation}
commutatif. 

On appelle topologie {\em co-évanescente} de $G$ la topologie engendrée par les recouvrements 
\[
\{(W_i\rightarrow U_i\leftarrow V_i)\rightarrow (W\rightarrow U\leftarrow V)\}_{i\in I}
\] 
des trois types suivants~:
\begin{itemize}
\item[(a)] $V_i=V$, $U_i=U$ pour tout $i\in I$, et $(W_i\rightarrow W)_{i\in I}$ est une famille couvrante.
\item[(b)] $W_i=W$, $U_i=U$ pour tout $i\in I$, et $(V_i\rightarrow V)_{i\in I}$ est une famille couvrante.
\item[(c)] $I=\{'\}$, $W'=W$ et le morphisme $V'\rightarrow V\times_{U}U'$ est un isomorphisme
(il n'y a aucune condition sur le morphisme $U'\rightarrow U$). 
\end{itemize}
Le site ainsi défini est appelé {\em site de Faltings relatif associé à $(f,g)$}; c'est un $\mU$-site.
On appelle {\em topos de Faltings relatif associé à $(f,g)$}, et l'on note $\tG$, 
le topos des faisceaux de $\mU$-ensembles sur $G$. 
Si $F$ est un préfaisceau de $G$, on note $F^a$ le faisceau associé. 

\subsection{}\label{tfr2}
Les limites projectives finies sont représentables dans $G$. En effet, $(Z\rightarrow X\leftarrow Y)$ est un objet final de $G$.
Par ailleurs, la limite projective d'un diagramme  
\begin{equation}\label{tfr2a}
\xymatrix{
&{(W''\rightarrow U''\leftarrow V'')}\ar[d]\\
{(W'\rightarrow U'\leftarrow V')}\ar[r]&{(W\rightarrow U\leftarrow V)}}
\end{equation} 
de $G$ est représentable par le morphisme $(W'\times_WW''\rightarrow U'\times_UU''\leftarrow V'\times_VV'')$. 
En effet, cet objet représente clairement la limite projective du diagramme \eqref{tfr2a}
dans la catégorie des morphismes de schémas au-dessus de $f$ et $g$, et c'est un objet de $G$ d'après (\cite{agt} VI.10.3). 

On notera que chacune des familles de recouvrements de type (a), (b) et (c) est stable par changement de base.

\begin{lem}\label{tfr3}
Pour qu'un préfaisceau $F$ sur $G$ soit un faisceau, il faut et il suffit que les conditions suivantes soient remplies~:
\begin{itemize}
\item[{\rm (i)}] Pour toute famille couvrante $(A_i\rightarrow A)_{i\in I}$ de $G$ de type {\rm (a)} ou {\rm (b)}, 
la suite 
\begin{equation}\label{tfr3a}
F(A)\rightarrow \prod_{i\in I}F(A_i)\rightrightarrows \prod_{(i,j)\in I\times J}F(A_i\times_AA_j)
\end{equation}
est exacte. 
\item[{\rm (ii)}] Pour tout recouvrement $(W\rightarrow U'\leftarrow V')\rightarrow (W\rightarrow U\leftarrow V)$ de $G$
de type {\rm (c)}, l'application 
\begin{equation}\label{tfr3b}
F(W\rightarrow U\leftarrow V) \rightarrow F(W\rightarrow U'\leftarrow V')
\end{equation}
est bijective.
\end{itemize}

En particulier, tout recouvrement $A'\rightarrow A$ de $G$ de type {\rm (c)} induit un isomorphisme entre les faisceaux associés 
$A'^a\stackrel{\sim}{\rightarrow} A^a$. 
\end{lem}

On notera d'abord que la catégorie $G$ est $\mV$-petite \eqref{notconv3} et qu'il suffit de montrer la proposition pour un préfaisceau $F$ 
de $\mV$-ensembles (\cite{sga4} II 2.7(2)). 
Pour tout recouvrement $(W\rightarrow U'\leftarrow V')\rightarrow (W\rightarrow U\leftarrow V)$ 
du type (c), le morphisme diagonal 
\begin{equation}\label{tfr3c}
(W\rightarrow U'\leftarrow V')\rightarrow (W\rightarrow U'\times_UU'\leftarrow V'\times_VV')
\end{equation}
est un recouvrement du type (c) qui égalise les deux projections canoniques
\begin{equation}\label{tfr3d}
(W\rightarrow U'\times_UU'\leftarrow V'\times_VV')\rightrightarrows 
(W\rightarrow U'\leftarrow V').
\end{equation}
La proposition résulte donc de (\cite{sga4} II 2.3, I 3.5 et I 2.12).

\subsection{}\label{tfr4}
Les foncteurs 
\begin{eqnarray}
\pi^+\colon \Et_{/Z}&\rightarrow& G,\ \ \ W\mapsto (W\rightarrow X\leftarrow Y),\label{tfr4a}\\
\lambda^+\colon \Et_{\rf/Y}&\rightarrow& G,\ \ \ V\mapsto (Z\rightarrow X\leftarrow V),\label{tfr4b}
\end{eqnarray}
sont exacts à gauche et continus (\cite{sga4} III 1.6). Ils définissent donc deux morphismes de topos (\cite{sga4} IV 4.9.2)
\begin{eqnarray}
\pi\colon \tG&\rightarrow& Z_\et,\label{tfr4c}\\
\lambda\colon \tG&\rightarrow& Y_\fet.\label{tfr4d}
\end{eqnarray}

\begin{rema}\label{tfr43}
Pour tout objet $(W\rightarrow U\leftarrow V)$ de $G$ tel que l'un des schémas $V$ ou $W$ soit vide, 
le faisceau associé $(W\rightarrow U\leftarrow V)^a$ est l'objet initial de $\tG$ (\cite{sga4} II 4.5). 
En effet, le schéma vide $\emptyset$ est l'objet
initial de $Z_\et$ (resp. $Y_\fet$). Par suite, $\pi^*(\emptyset)=(\emptyset\rightarrow X\leftarrow Y)^a$ et 
$\lambda^*(\emptyset)=(Z\rightarrow X\leftarrow \emptyset)^a$ sont égaux à l'objet initial de $\tG$.   
L'assertion recherchée résulte alors du fait que l'objet initial de $\tG$ est strict (\cite{sga4} II 4.5.1).
\end{rema}

\begin{lem}\label{tfr44}
Soit $(Z_i\rightarrow X_i\leftarrow Y_i)_{i\in I}$ une famille d'objets de $G$ telle que le schéma 
$Z$ soit la somme disjointe des schémas $(Z_i)_{i\in I}$, 
et que pour tout $i\in I$, le morphisme canonique $Y_i\rightarrow X_i\times_XY$ soit un isomorphisme. 
Pour tout $i\in I$, soit $(Y_{ij}\rightarrow Y_i)_{j\in J_i}$ une famille de morphismes étales finis telle que le schéma $Y_i$ 
soit la somme disjointe des schémas $(Y_{ij})_{j\in J_i}$. Alors, l'objet final $(Z\rightarrow X\leftarrow Y)^a$ de $\tG$
est la somme disjointe des faisceaux $(Z_i\rightarrow X_i\leftarrow Y_{ij})^a$ de $\tG$ associés aux objets 
$(Z_i\rightarrow X_i\leftarrow Y_{ij})$ de $G$, pour $i\in I$ et $j\in J_i$ {\rm (\cite{sga4} II 4.5)}.
\end{lem}

En effet, $(Z\rightarrow X\leftarrow Y)^a$ est la somme disjointe des $(Z_i\rightarrow X\leftarrow Y)^a_{i\in I}$ dans $\tG$
d'après \ref{tfr43} et (\cite{sga4} II 4.6 2)). 
Pour tout $i\in I$, le morphisme canonique $(Z_i\rightarrow X_i\leftarrow Y_i)^a\rightarrow (Z_i\rightarrow X\leftarrow Y)^a$
est un isomorphisme en vertu de \ref{tfr3}. 
Par ailleurs, $(Z_i\rightarrow X_i\leftarrow Y_i)^a$ est la somme disjointe des $(Z_i\rightarrow X_i\leftarrow Y_{ij})^a_{j\in J_i}$ dans $\tG$
d'après \ref{tfr43} et (\cite{sga4} II 4.6 2)), d'où la proposition. 

\begin{lem}\label{tfr32}
Soit $(A_i\rightarrow A)_{i\in I}$ un recouvrement fini de $G$ de type {\rm (b)}, et pour tout $i\in I$,
soit $(A'_i\rightarrow A_i)$ un recouvrement de $G$ de type {\rm (c)} \eqref{tfr1}. Alors, il existe un recouvrement de type {\rm (c)},
$(A'\rightarrow A)$ de $G$ et pour tout $i\in I$, un $A_i$-morphisme $A'\times_AA_i\rightarrow A'_i$; en particulier, 
le recouvrement $(A'\times_AA_i\rightarrow A)_{i\in I}$ raffine le recouvrement $(A'_i\rightarrow A)_{i\in I}$.
\end{lem}

Posons $A=(W\rightarrow U\leftarrow V)$ et pour tout $i\in I$, $A_i=(W\rightarrow U\leftarrow V_i)$
et $A'_i=(W\rightarrow U'_i\leftarrow V'_i)$, de sorte qu'on a un diagramme commutatif à carré supérieur droit cartésien
\begin{equation}
\xymatrix{
W\ar@{=}[d]\ar[r]&{U'_i}\ar[d]\ar@{}[rd]|\Box&V'_i\ar[l]\ar[d]\\
W\ar@{=}[d]\ar[r]&{U}\ar@{=}[d]&V_i\ar[l]\ar[d]\\
W\ar[r]&U&V\ar[l]}
\end{equation}
Notant $U'$ le produit des schémas $U'_i$ au-dessus de $U$, il existe alors un unique morphisme  
$W\rightarrow U'$ au-dessus de chacun des morphismes $W\rightarrow U'_i$ $(i\in I)$, 
de sorte que $A'=(W\rightarrow U'\leftarrow U'\times_UV)$ soit un objet de $G$ et que le morphisme canonique $A'\rightarrow A$
soit un recouvrement de type (c).  Pour tout $i\in I$, on a un isomorphisme canonique \eqref{tfr2}
\begin{equation}
A'\times_AA_i\stackrel{\sim}{\rightarrow} (W\rightarrow U'\leftarrow U'\times_UV_i).
\end{equation}
On en déduit un $A_i$-morphisme $A'\times_AA_i\rightarrow A'_i$, d'où la proposition.

\subsection{}\label{tfr5}
On désigne par $C$ le site défini dans \ref{topfl1} relativement aux foncteurs de changement de base par $f$ et $g$ respectivement, 
\begin{equation}
f^+\colon \Et_{/X}\rightarrow \Et_{/Y}\ \ \  {\rm et} \ \ \ g^+\colon \Et_{/X}\rightarrow \Et_{/Z}.
\end{equation}
Concrètement, la catégorie sous-jacente à $C$ est 
la catégorie des triplets $(U, V\rightarrow U, W\rightarrow U)$ formés d'un $X$-schéma $U$ et de deux 
morphismes $V\rightarrow U$ et $W\rightarrow U$ au-dessus de $f$ et $g$ respectivement 
tels que les morphismes $U\rightarrow X$, $V\rightarrow Y$ et $W\rightarrow Z$ soient étales; un tel objet sera 
noté $(W\rightarrow U\leftarrow V)$. Soient $(W\rightarrow U\leftarrow V)$ et $(W'\rightarrow U'\leftarrow V')$ deux objets de $C$. 
Un morphisme de $(W'\rightarrow U'\leftarrow V')$ dans $(W\rightarrow U\leftarrow V)$ et la donnée de trois morphismes 
$U'\rightarrow U$, $V'\rightarrow V$ et $W'\rightarrow W$ au-dessus de $X$, $Y$ et $Z$ respectivement, qui rendent le diagramme 
\begin{equation}\label{tfr5a}
\xymatrix{
W'\ar[r]\ar[d]&U'\ar[d]&V'\ar[l]\ar[d]\\
W\ar[r]&U&V\ar[l]}
\end{equation}
commutatif. On munit $C$ de la topologie engendrée par les recouvrements 
\begin{equation}
\{(W_i\rightarrow U_i\leftarrow V_i)\rightarrow (W\rightarrow U\leftarrow V)\}_{i\in I}
\end{equation}
des trois types suivants~:
\begin{itemize}
\item[(a)] $V_i=V$, $U_i=U$ pour tout $i\in I$, et $(W_i\rightarrow W)_{i\in I}$ est une famille couvrante.
\item[(b)] $W_i=W$, $U_i=U$ pour tout $i\in I$, et $(V_i\rightarrow V)_{i\in I}$ est une famille couvrante.
\item[(c)] $I=\{'\}$, $W'=W$ et le morphisme $V'\rightarrow V\times_{U}U'$ est un isomorphisme
(il n'y a aucune condition sur le morphisme $U'\rightarrow U$). 
\end{itemize}
On rappelle que le topos des faisceaux de $\mU$-ensembles sur $C$ est le produit orienté $Z_\et\gtimes_{X_\et} Y_\et$ 
des morphismes de topos $g_\et\colon Z_\et\rightarrow X_\et$ et $f_\et\colon Y_\et\rightarrow X_\et$ (cf. \ref{topfl4}).

\subsection{}\label{tfr6}
Tout objet de $G$ est naturellement un objet de $C$. 
On définit ainsi un foncteur pleinement fidèle et exact à gauche
\begin{equation}\label{tfr6a}
\rho^+\colon G\rightarrow C.
\end{equation} 
Celui-ci est continu en vertu de \ref{tfr3} et (\cite{agt} VI.3.2). Il définit donc un morphisme de topos 
\begin{equation}\label{tfr6b}
\rho\colon Z_\et\gtimes_{X_\et}Y_\et\rightarrow \tG.
\end{equation}
Il résulte aussitôt des définitions que les carrés du diagramme
\begin{equation}\label{tfr6c}
\xymatrix{
{Z_\et}\ar@{=}[d]&{Z_\et\gtimes_{X_\et}Y_\et}\ar[l]_-(0.5){\rp_1}\ar[d]^{\rho}\ar[r]^-(0.5){\rp_2}&
{Y_\et}\ar[d]^{\rho_Y}\\
{Z_\et}&{\tG}\ar[l]_{\pi}\ar[r]^{\lambda}&{Y_\fet}}
\end{equation}
où $\rp_1$ et $\rp_2$ sont les projections canoniques (\cite{agt} VI.3.4),
sont commutatifs à isomorphismes canoniques près. 

\subsection{}\label{tfr7}
Considérons un diagramme commutatif de morphismes de schémas 
\begin{equation}\label{tfr7a}
\xymatrix{
Z'\ar[r]^{g'}\ar[d]_w&X'\ar[d]^u&Y'\ar[l]_{f'}\ar[d]^v\\
Z\ar[r]^g&X&\ar[l]_fY}
\end{equation}
On désigne par $G'$ (resp. $\tG'$) le site (resp. topos) de Faltings relatif associé à $(f',g')$ \eqref{tfr1}. 
Le foncteur 
\begin{equation}\label{tfr7d}
\Phi^+\colon G\rightarrow G', \ \ \ (W\rightarrow U\leftarrow V)\mapsto (W\times_ZZ'\rightarrow U\times_XX'\leftarrow V\times_YY')
\end{equation}
est clairement exact à gauche. Celui-ci transforme les recouvrements de $G$ du type (a) (resp. (b), resp. (c)) 
en des recouvrements de $G'$ du même type. Il résulte alors de \ref{tfr3} 
que pour tout faisceau $F$ sur $G'$, $F\circ \Phi^+$ est un faisceau sur $G$. Par suite, $\Phi^+$
est continu. Il définit donc un morphisme de topos 
\begin{equation}\label{tfr7b}
\Phi\colon \tG'\rightarrow \tG.
\end{equation}
Il résulte aussitôt des définitions que les carrés du diagramme
\begin{equation}\label{tfr7c}
\xymatrix{
{Z'_\et}\ar[d]_{w_\et}&{\tG'}\ar[l]_-(0.5){\pi'}\ar[d]^{\Phi}\ar[r]^-(0.5){\lambda'}&
{Y'_\fet}\ar[d]^{v_\fet}\\
{Z_\et}&{\tG}\ar[l]_-(0.5){\pi}\ar[r]^-(0.5){\lambda}&{Y_\fet}}
\end{equation}
où $\pi$, $\pi'$ \eqref{tfr4c}, $\lambda$ et $\lambda'$ \eqref{tfr4d} sont les morphismes canoniques, 
sont commutatifs à isomorphismes canoniques près. 

Le diagramme de morphismes de topos
\begin{equation}\label{tfr7e}
\xymatrix{
{Z'_\et\gtimes_{X'_\et}Y'_\et}\ar[r]^-(0.5){\rho'}\ar[d]_{\nu}&{\tG'}\ar[d]^\Phi\\
{X_\et\gtimes_{X_\et}Y_\et}\ar[r]^-(0.5)\rho&{\tG}}
\end{equation}
où $\rho$ et $\rho'$ sont les morphismes canoniques \eqref{tfr6b} et $\nu$ est le morphisme de fonctorialité \eqref{topfl10d}, 
est commutatif à isomorphisme canonique près, d'après \ref{topfl11}.

\subsection{}\label{tfr40}
Soit $(Z'\rightarrow X'\leftarrow Y')$ un objet de $G$. On désigne par $G'$ (resp. $\tG'$) le site (resp. topos) de Faltings relatif associé au couple
de morphisme $(Y'\rightarrow X', Z'\rightarrow X')$ \eqref{tfr1}. Tout objet de $G'$ est naturellement un objet de $G$. 
On définit ainsi un foncteur 
\begin{equation}\label{tfr40a}
\Phi\colon G'\rightarrow G.
\end{equation}
On vérifie aussitôt que $\Phi$ se factorise à travers une équivalence de catégories 
\begin{equation}\label{tfr40b}
G'\stackrel{\sim}{\rightarrow} G_{/(Z'\rightarrow X'\leftarrow Y')}. 
\end{equation}

\begin{prop}\label{tfr41}
Conservons les hypothèses de \ref{tfr40}. Alors,  
\begin{itemize}
\item[{\rm (i)}] Le foncteur $\Phi$ \eqref{tfr40a} est continu et cocontinu {\rm (\cite{sga4} III 2.1)}.
\item[{\rm (ii)}] La topologie co-évanescente de $G'$ \eqref{tfr1} est induite par celle de $G$ via le foncteur $\Phi$.
\end{itemize}
\end{prop}

(i) Il est clair que $\Phi$ commute aux produits fibrés \eqref{tfr2} et qu'il 
transforme les familles couvrantes de $G'$ de types (a), (b) et (c) en des familles couvrantes de $G$ de même types \eqref{tfr1}. 
Par suite, pour tout faisceau $F$ de $\tG$, $F\circ \Phi$ est un faisceau sur $G'$ en vertu de \ref{tfr3}. 
Le foncteur $\Phi$ est donc continu.

Le foncteur $\Phi$ est adjoint à gauche du foncteur 
\begin{equation}\label{tfr41a}
\Phi^+\colon G\rightarrow G'
\end{equation}
défini pour tout objet $(W\rightarrow U\leftarrow V)$ de $G$, par 
\begin{equation}\label{tfr41b}
\Phi^+((W\rightarrow U\leftarrow V)) = (W\times_ZZ'\rightarrow U\times_XX'\leftarrow V\times_YY').
\end{equation}
Comme $\Phi^+$ est continu d'après \ref{tfr7}, $\Phi$ est cocontinu en vertu de (\cite{sga4} III 2.5). 

(ii) Comme $\Phi$ est continu, la topologie co-évanescente de $G'$
est moins fine que sa topologie $\cT$ induite par la topologie co-évanescente de $G$ via le foncteur $\Phi$ (\cite{sga4} III 3.1).
Soit $(A_i\rightarrow A)_{i\in I}$ une famille de $G'$ couvrante pour la topologie $\cT$.
La famille $(\Phi(A_i)\rightarrow \Phi(A))_{i\in I}$ est alors couvrante pour la topologie co-évanescente de $G$ d'après (\cite{sga4} III 3.3). 
Comme $\Phi$ est cocontinu, la famille $(A_i\rightarrow A)_{i\in I}$ est couvrante pour la topologie co-évanescente de $G$ (\cite{sga4} III 2.1), 
d'où la proposition.

\subsection{}\label{tfr42}
Conservons les hypothèses de \ref{tfr40}. 
En vertu de \ref{tfr41}, le foncteur $\Phi$ \eqref{tfr40a} définit une suite de trois foncteurs adjoints~:
\begin{equation}\label{tfr42a}
\Phi_!\colon \tG'\rightarrow \tG, \ \ \ \Phi^*\colon \tG\rightarrow \tG', \ \ \ \Phi_*\colon \tG'\rightarrow \tG,
\end{equation}
dans le sens que pour deux foncteurs consécutifs de la suite, celui de droite est
adjoint à droite de l'autre. D'après (\cite{sga4} III 5.4), le foncteur $\Phi_!$ se factorise à travers 
une équivalence de catégories 
\begin{equation}\label{tfr42b}
\tG'\stackrel{\sim}{\rightarrow} \tG_{/(Z'\rightarrow X'\leftarrow Y')^a}.
\end{equation}
Le couple de foncteurs $(\Phi^*,\Phi_*)$ définit le morphisme de localisation de $\tG$ en $(Z'\rightarrow X'\leftarrow Y')^a$
\begin{equation}\label{tfr42c}
\Phi\colon \tG'\rightarrow \tG.
\end{equation}
Comme $\Phi\colon G'\rightarrow G$ est un adjoint à gauche du foncteur $\Phi^+\colon G\rightarrow G'$ 
défini dans \eqref{tfr41a}, le morphisme \eqref{tfr42c}
s'identifie au morphisme défini dans \eqref{tfr7b}, en vertu de (\cite{sga4} III 2.5).

\subsection{}\label{tfr9}
Supposons $Z=X$ et $g=\id_X$ et notons $E$ (resp. $\tE$) le site (resp. topos) de Faltings associé à $f$ \eqref{tf1}. 
Considérons les foncteurs
\begin{eqnarray}
\iota^+\colon E&\rightarrow &G,\ \ \ (V\rightarrow U)\mapsto (U\rightarrow U\leftarrow V),\label{tfr9a}\\
\jmath^+\colon G&\rightarrow &E,\ \ \  (W\rightarrow U\leftarrow V)\mapsto (V\times_{U}W\rightarrow W).\label{tfr9b}
\end{eqnarray}
Il est clair que $\iota^+$ est un adjoint à gauche de $\jmath^+$, 
que le morphisme d'adjonction $\id\rightarrow \jmath^+\circ \iota^+$ 
est un isomorphisme ({\em i.e.}, $\iota^+$ est pleinement fidèle), et que $\iota^+$ et $\jmath^+$ sont exacts à gauche.

\begin{prop}\label{tfr10}
Conservons les hypothèses et notations de \ref{tfr9}. 
\begin{itemize}
\item[{\rm (i)}] Les foncteurs $\iota^+$ et $\jmath^+$ sont continus.
\item[{\rm (ii)}] La topologie de $E$ est induite par celle de $G$ au moyen du foncteur $\iota^+$.
\end{itemize}
\end{prop}

(i) Avec les notations de \ref{tfr1} et (\cite{agt} VI.5.3), le foncteur $\iota^+$ transforme les familles couvrantes
de $E$ du type (v) en familles couvrantes de $G$ du type (b), 
et les familles couvrantes de $E$ du type (c) en familles couvrantes de $G$~:
\begin{equation}\label{tfr10a}
\xymatrix{
{V_i}\ar[r]\ar[d]\ar@{}[rd]|{\Box}&{U_i}\ar[d]\\
V\ar[r]&U}
\ \ \ \mapsto \ \ \ 
\xymatrix{
{U_i}\ar[r]\ar[d]&{U_i}\ar[d]&{V_i}\ar[l]\ar[d]\ar@{}[ld]|{\Box}\\
{U_i}\ar[r]\ar[d]&U\ar[d]&V\ar[l]\ar[d]\\
U\ar[r]&U&V\ar[l]}
\end{equation}
Soient $P$ un préfaisceau sur $G$, $F=\{U\mapsto F_U\}=P\circ \iota^+$ \eqref{tf1h}. 
Pour tout $(V\rightarrow U)\in \ob(E)$, on a 
\begin{equation}\label{tfr10b}
F_U(V)=P(U\rightarrow U\leftarrow V).
\end{equation}
Par suite, si $P$ est un faisceau sur $G$, 
$F$ est un faisceau sur $E$ en vertu de \ref{tfr3}, \eqref{tfr10a} et  (\cite{agt} VI.5.10); donc $\iota^+$ est continu. 

Le foncteur $\jmath^+$ transforme les recouvrements de $G$ de type (a) (resp. (b))
en recouvrements de $E$ de type (c) (resp. (v)), et les recouvrements de $G$ de type (c)
en isomorphismes. Par suite, pour tout faisceau $F$ sur $E$, $F\circ \jmath^+$
est un faisceau sur $G$ en vertu de \ref{tfr3}; donc $\jmath^+$ est continu. 

(ii) On sait que la topologie de $E$ est induite par la topologie canonique de $\tE$ (\cite{sga4} III 3.5),
autrement dit, la topologie de $E$ est la plus fine telle que tout $F\in \ob(\tE)$ soit un faisceau.
D'après (i), on peut considérer les foncteurs 
\begin{eqnarray}
\iota_s\colon \tG&\rightarrow& \tE,\ \ \  P\mapsto P\circ \iota^+,\\
\jmath_s\colon \tE&\rightarrow& \tG,\ \ \ F\mapsto F\circ \jmath^+.
\end{eqnarray}
L'isomorphisme d'adjonction $\id\stackrel{\sim}{\rightarrow} \jmath^+\circ \iota^+$ induit un isomorphisme 
$\iota_s\circ \jmath_s\stackrel{\sim}{\rightarrow} \id$. Le foncteur $\iota_s$ est donc essentiellement surjectif. 
Par suite, la topologie de $E$ est la plus fine telle que, pour tout $P\in \ob(\tG)$, 
$\iota_s(P)$ soit un faisceau sur $E$; d'où la proposition.

\subsection{}\label{tfr11}
Conservons les hypothèses et notations de \ref{tfr9}. 
Les foncteurs $\iota^+$ et $\jmath^+$ étant exacts à gauche et continus \eqref{tfr10}, ils définissent 
des morphismes de topos (\cite{sga4} IV 4.9.2)
\begin{eqnarray}
\iota\colon \tG&\rightarrow& \tE,\label{tfr11a}\\
\jmath\colon \tE&\rightarrow& \tG.\label{tfr11b}
\end{eqnarray}
Les morphismes d'adjonction $\id \rightarrow \jmath^+\circ \iota^+$ et $\iota^+ \circ \jmath^+\rightarrow \id$ 
induisent des morphismes $\iota_*\circ \jmath_*\rightarrow \id$ et $\id\rightarrow \jmath_*\circ \iota_*$
qui font de $\iota_*$ un adjoint à droite de $\jmath_*$. 
Par ailleurs, les diagrammes de morphismes de topos
\begin{equation}\label{tfr11c}
\xymatrix{
&{\tG}\ar[ld]_{\pi}&{\tG}\ar[rd]^{\lambda}\ar[dd]_{\iota}&\\
Z_\et&&&Y_\fet\\
&\tE\ar[lu]^{\sigma}\ar[uu]_{\jmath}&\tE\ar[ru]_{\beta}&}
\end{equation}
où $\pi$ \eqref{tfr4c}, $\lambda$ \eqref{tfr4d}, $\sigma$ \eqref{tf1kk} et $\beta$ \eqref{tf1k} sont les morphismes canoniques,
sont commutatifs à isomorphismes canoniques près.

\begin{prop}\label{tfr12}
Les morphismes d'adjonction $\iota_*\circ \jmath_*\rightarrow \id$ et $\id\rightarrow \jmath_*\circ \iota_*$ 
sont des isomorphismes. En particulier, $\iota$ \eqref{tfr11a} et $\jmath$ \eqref{tfr11b} 
sont des équivalences de topos quasi-inverses l'une de l'autre.
\end{prop}

En effet, comme le morphisme d'adjonction $\id \rightarrow \jmath^+\circ \iota^+$ est un isomorphisme, 
$\iota_*\circ \jmath_*\rightarrow \id$ est un isomorphisme.  D'autre part, le morphisme 
d'adjonction $\iota^+ \circ \jmath^+\rightarrow \id$ est défini, pour tout objet $(W\rightarrow U\leftarrow V)$ de $G$, 
par le morphisme canonique
\begin{equation}\label{tfr12a}
(W\rightarrow W\leftarrow V\times_{U}W)\rightarrow (W\rightarrow U\leftarrow V),
\end{equation}
qui est un recouvrement de type (c). On en déduit d'après \ref{tfr3}(ii) que $\id\rightarrow \jmath_*\circ \iota_*$ est un isomorphisme.

\subsection{}\label{tfr13}
Considérons un diagramme commutatif de morphismes de schémas 
\begin{equation}\label{tfr13a}
\xymatrix{
Z\ar[d]_{g}&T\ar[l]_-(0.4){f'}\ar[d]^{g'}\\
X&Y\ar[l]_-(0.4)f}
\end{equation}
et notons $\tE$ (resp. $\tE'$) le topos de Faltings associé à $f$ (resp. $f'$) \eqref{tf1}. 
Considérons le diagramme commutatif de morphismes de schémas 
\begin{equation}\label{tfr13b}
\xymatrix{
Z\ar@{=}[r]\ar@{=}[d]&Z\ar[d]^g&T\ar[l]_{f'}\ar[d]^{g'}\\
Z\ar[r]^{g}\ar[d]_{g}&X\ar@{=}[d]&Y\ar[l]_{f}\ar@{=}[d]\\
X\ar@{=}[r]&X&Y\ar[l]_{f}}
\end{equation}
et identifions $\tE'$ au topos de Faltings relatif à $(f',\id_{Z})$ par l'équivalence de topos $\jmath'$ \eqref{tfr11b}
et $\tE$ au topos de Faltings relatif à $(f,\id_{X})$ par l'équivalence de topos $\iota$ \eqref{tfr11a}. 
On en déduit par fonctorialité \eqref{tfr7} des morphismes de topos  
\begin{equation}\label{tfr13c}
\xymatrix{
\tE'\ar[r]^{\tau}&\tG\ar[r]^{\lgg}&\tE}.
\end{equation}
dont le composé est le morphisme de fonctorialité induit par le diagramme \eqref{tfr13a} (cf. \cite{agt} VI.10.12). 

Il résulte de \eqref{tfr7c}, \eqref{tfr11c} et \ref{tfr12} que les triangles et carrés du diagramme de morphismes de topos 
\begin{equation}\label{tfr13d}
\xymatrix{
&\tE'\ar[d]^{\tau}\ar[r]^{\beta'}\ar[ld]_{\sigma'} & T_\fet\ar[d]^{g'_\fet}\\
Z_\et\ar[d]_{g_\et}&\tG\ar[d]_-(0.5){\lgg}\ar[r]^-(0.4){\lambda}\ar[l]_-(0.4){\pi}&Y_\fet\\
X_\et&\tE\ar[ur]_{\beta}\ar[l]_{\sigma}&}
\end{equation}
où $\pi$ \eqref{tfr4c} et $\lambda$ \eqref{tfr4d}, $\sigma$, $\sigma'$ \eqref{tf1kk}, 
$\beta$ et $\beta'$ \eqref{tf1k} sont les morphismes canoniques, 
sont commutatifs à isomorphismes canoniques près. 

\begin{prop}\label{tfr20}
Notant $\tE$ le topos de Faltings associé à $f$, le diagramme de morphismes de topos
\begin{equation}\label{tfr20a}
\xymatrix{
\tG\ar[d]_{\pi}\ar[r]^{\lgg}&\tE\ar[d]^{\sigma}\\
Z_\et\ar[r]^{g_\et}&X_\et}
\end{equation}
où $\lgg$, $\pi$ et $\sigma$ sont les morphismes canoniques \eqref{tfr13}, est commutatif à isomorphisme canonique près.
Il identifie $\tG$ au produit fibré des deux morphismes de topos 
$\sigma$ et $g_\et$ {\rm (\cite{travaux-gabber} XI §3)}. 
\end{prop}

La première assertion est immédiate et elle a déjà été vue dans \eqref{tfr13d}.
Considérons un diagramme de morphismes de topos
\begin{equation}\label{tfr20b}
\xymatrix{
\cT\ar[d]_{u}\ar[r]^{v}&\tE\ar[d]^{\sigma}\\
Z_\et\ar[r]^{g_\et}&X_\et}
\end{equation}
commutatif à isomorphismes près 
\begin{equation}\label{tfr20c}
g_\et\circ u\stackrel{\sim}{\rightarrow} \sigma\circ v.
\end{equation}
Montrons qu'il existe un unique triplet $(t,\alpha,\beta)$, où  
\begin{equation}\label{tfr20d}
t\colon \cT\rightarrow \tG
\end{equation}
est un morphisme de topos et $\alpha\colon v\stackrel{\sim}{\rightarrow} \lgg\circ t$ et 
$\beta\colon u\stackrel{\sim}{\rightarrow} \pi\circ t$ sont des isomorphismes.

Pour tout objet $(W\rightarrow U\leftarrow V)$ de $G$, le diagramme de morphismes de $G$
\begin{equation}\label{tfr20e}
\xymatrix{
{(W\rightarrow U\leftarrow V)}\ar[r]\ar[d]&{(U\times_XZ\rightarrow U\leftarrow V)}\ar[d]\\
{(W\rightarrow X\leftarrow Y)}\ar[r]&{(U\times_XZ\rightarrow X\leftarrow Y)}}
\end{equation}
est cartésien. De plus, on a 
\begin{eqnarray}
\lgg^*((V\rightarrow U)^a)&\stackrel{\sim}{\rightarrow}&(U\times_XZ\rightarrow U\leftarrow V)^a,\label{tfr20f1}\\
\pi^*(W^a)&\stackrel{\sim}{\rightarrow}&(W\rightarrow X\leftarrow Y)^a,\label{tfr20f2}\\
\pi^*(g^*(U^a))&\stackrel{\sim}{\rightarrow}&(U\times_XZ\rightarrow X\leftarrow Y)^a,\label{tfr20f3}
\end{eqnarray}
où l'exposant $^a$ désigne les faisceaux associés. Tenant compte de l'isomorphisme canonique 
$u^*\circ g^*\stackrel{\sim}{\rightarrow} v^*\circ \sigma^*$ \eqref{tfr20c}, on définit le foncteur 
\begin{equation}\label{tfr20h}
t^+\colon 
\begin{array}[t]{clcr}
G&\rightarrow &\cT\\
(W\rightarrow U\leftarrow V)&\mapsto& u^*(W^a)\times_{u^*(g^*(U^a))}v^*((V\rightarrow U)^a)
\end{array}
\end{equation}
Il est clair que si $(t,\alpha,\beta)$ existe, le composé de $t^*$ et du foncteur canonique $G\rightarrow \tG$ est canoniquement isomorphe à $t^+$. 
Ce dernier étant exact à gauche, il suffit donc de montrer qu'il est continu. 

Le foncteur $t^+$ transforme clairement les familles couvrantes de $G$ de type (a) et (b) en des familles épimorphiques de  $\cT$. 
Soit $(W\rightarrow U'\leftarrow V')\rightarrow (W\rightarrow U\leftarrow V)$ un recouvrement de $G$ de type (c). 
Le diagramme canonique de morphismes de $E$
\begin{equation}
\xymatrix{
{(V'\rightarrow U')}\ar[r]\ar[d]&{(V\rightarrow U)}\ar[d]\\
{(U'\times_XY\rightarrow U')}\ar[r]&{(U\times_XY\rightarrow U)}}
\end{equation}
est cartésien. On en déduit un diagramme cartésien de $\cT$
\begin{equation}
\xymatrix{
{v^*((V'\rightarrow U')^a)}\ar[r]\ar[d]&{v^*((V\rightarrow U)^a)}\ar[d]\\
{v^*(\sigma^*(U'^a))}\ar[r]&{v^*(\sigma^*(U^a))}}
\end{equation}
On en déduit que le morphisme 
\begin{equation}
t^+(W\rightarrow U'\leftarrow V')\rightarrow t^+(W\rightarrow U\leftarrow V)
\end{equation}
est un isomorphisme. Par suite, le foncteur $t^+$ est continu en vertu de \ref{tfr3}, d'où la proposition.

\subsection{}\label{tfr14}
Soit $\star$ l'un des deux symboles ``$\coh$'' pour cohérent, ou ``$\scoh$'' pour séparé et cohérent, introduits dans \ref{notconv10}. 
On désigne par $G_\star$ la sous-catégorie pleine de $G$ formée des objets $(W\rightarrow U\leftarrow V)$ tels que 
$U$ soit un objet de $\Et_{\star/X}$ et $W$ soit un objet de $\Et_{\star/Z}$. Les limites projectives finies  
sont représentables dans $G_\star$ \eqref{tfr2}. 
On munit $G_\star$ de la topologie engendrée par les familles couvrantes de type (c) 
et les familles couvrantes {\em finies} de type (a), (b) (cf. \ref{tfr1}) et on désigne par $\tG_\star$ le topos des faisceaux de $\mU$-ensembles sur $G_\star$. 
Calquant la preuve de \ref{tfr3}, on montre que pour qu'un préfaisceau $F$ sur $G_\star$ soit un faisceau, 
il faut et il suffit que les conditions suivantes soient remplies~:
\begin{itemize}
\item[{\rm (i)}] Pour toute famille couvrante {\em finie} $(A_i\rightarrow A)_{i\in I}$ de $G_\star$ de type {\rm (a)} ou {\rm (b)}, 
la suite 
\begin{equation}\label{tfr14a}
F(A)\rightarrow \prod_{i\in I}F(A_i)\rightrightarrows \prod_{(i,j)\in I\times J}F(A_i\times_AA_j)
\end{equation}
est exacte. 
\item[{\rm (ii)}] Pour tout recouvrement $(W\rightarrow U'\leftarrow V')\rightarrow (W\rightarrow U\leftarrow V)$ de $G_\star$
de type {\rm (c)}, l'application 
\begin{equation}\label{tfr14b}
F(W\rightarrow U\leftarrow V) \rightarrow F(W\rightarrow U'\leftarrow V')
\end{equation}
est bijective.
\end{itemize}

Le foncteur d'injection canonique $\phi\colon G_\star\rightarrow G$ est exact à gauche et continu \eqref{tfr3}.
Il définit donc un morphisme de topos (\cite{sga4} IV 4.9.2)
\begin{equation}\label{tfr14c}
\Phi\colon \tG\rightarrow \tG_\star.
\end{equation}

\begin{prop}\label{tfr15}
Les notations étant celles de \ref{tfr14}, supposons de plus les schémas $X$, $Y$ et $Z$ cohérents.  Alors, 
\begin{itemize}
\item[{\rm (i)}] La famille $G_\star$ est topologiquement génératrice dans $G$. 
\item[{\rm (ii)}] La topologie de $G_\star$ est engendrée par la prétopologie définie pour chaque objet $A$ de $G_\star$
par la donnée de l'ensemble $\Cov(A)$ des familles de morphismes $(A_i\rightarrow A)_{i\in I}$ obtenues par composition 
d'un nombre fini de familles de type (c) et de familles (finies) de type (a) et (b); en particulier, l'ensemble $I$ est fini.
\item[{\rm (iii)}] Le morphisme $\Phi$ \eqref{tfr14c} est une équivalence de topos. Nous identifierons dans la suite de cette section  
les topos $\tG$ et $\tG_\star$ au moyen du morphisme $\Phi$.   
\end{itemize}
\end{prop}

On notera d'abord que la condition que les schémas $X$, $Y$ et $Z$ sont cohérents est équivalente à la condition que le schéma $X$ est cohérent 
que les morphismes $f$ et $g$ sont cohérents. 

(i) Soit $(W\rightarrow U\leftarrow V)$ un objet de $G$. Comme le schéma $X$ est quasi-séparé, il existe une famille couvrante 
$(U_i\rightarrow U)_{i\in I}$ telle que $U_i\in \ob(\Et_{\star/X})$. Pour tout $i\in I$, posons $W_i=W\times_UU_i$ et $V_i=V\times_UV_i$, 
de sorte que $(W_i\rightarrow U_i\leftarrow V_i)$ est un objet de $G$. La famille 
\[
\{(W_i\rightarrow U_i\leftarrow V_i)\rightarrow (W\rightarrow U\leftarrow V)\}_{i\in I}
\] 
de morphismes de $G$ est couvrante, en tant que composée de familles de type (a) et (c). 
\[
\xymatrix{
W_{i}\ar[r]\ar@{=}[d]&U_i\ar[d]\ar@{}[rd]|{\Box}&V_i\ar[l]\ar[d]\\
W_{i}\ar[r]\ar[d]&U\ar@{=}[d]&V\ar[l]\ar@{=}[d]\\
W\ar[r]&U&\ar[l]V}
\]
Comme le schéma $Z$ est quasi-séparé, 
pour chaque $i\in I$, il existe une famille couvrante $(W_{i,j}\rightarrow W_i)_{j\in J_i}$ telle que $W_{i,j}\in \ob(\Et_{\star/Z})$.
La famille 
\[
\{(W_{i,j}\rightarrow U_i\leftarrow V_i)\rightarrow (W\rightarrow U\leftarrow V)\}_{i\in I,j\in J_i}
\] 
de morphismes de $G$ est couvrante, et chaque $(W_{i,j}\rightarrow U_i\leftarrow V_i)$ est un objet de $G_\star$. 

(ii) En effet, les familles couvrantes de type (c) et les familles couvrantes finies de type (a) et (b) de $G_\star$ étant stables par changement de base, 
les ensembles $\Cov(A)$ pour $A\in \ob(G_\star)$ vérifient les axiomes d'une prétopologie (\cite{sga4} II 1.3). 

(iii) Pour tout objet $A$ de $G$, on désigne par $G_{\star/A}$ la catégorie des couples $(B,u)$ où $B$ est un objet de $G_\star$
et $u\colon \phi(B)\rightarrow A$ est un morphisme de $G$. Si $(B,u)$ et $(B',u')$ sont deux objets de $G_{\star/A}$, un morphisme
de $(B,u)$ dans $(B',u')$ est un morphisme $t\colon B\rightarrow B'$ de $G_\star$ tel que $u'=\phi(t)\circ u$. 
Montrons que pour tout faisceau $F$ sur $G$ et tout $A\in \ob(G)$, le morphisme canonique 
\begin{equation}\label{tfr15a}
F(A)\rightarrow \underset{\underset{A'\in G_{\star/A}^\circ}{\longleftarrow}}{\lim} F(A')
\end{equation}
est un isomorphisme. En effet, il résulte de (i) que la suite 
\[
F(A)\rightarrow \prod_{A'\in \ob(G_{\star/A})}F(A')\rightrightarrows \prod_{(A',A'')\in \ob(G_{\star/A})^2}\ 
\prod_{B\in \ob(G_{\star/A'\times_AA''})}F(B)
\]
est exacte. D'autre part, le morphisme canonique 
\begin{eqnarray*}
\lefteqn{\underset{\underset{A'\in G^\circ_{\star/A}}{\longleftarrow}}{\lim}\ F(A')\rightarrow}\\
&& \ker\left(\prod_{A'\in \ob(G_{\star/A})}F(A')\rightrightarrows \prod_{(A',A'')\in \ob(G_{\star/A})^2}\ 
\prod_{B\in \ob(G_{\star/A'\times_AA''})}F(B)\right)
\end{eqnarray*}
est un monomorphisme, d'où l'assertion. 

L'isomorphisme fonctoriel \eqref{tfr15a} montre que le foncteur $\Phi_*\colon \tG\rightarrow \tG_\star$ est pleinement fidèle. Montrons qu'il  
est essentiellement surjectif. Soit $F'$ un $\mU$-faisceau sur $G_\star$. Pour tout objet $A$ de $G$, la catégorie $G_{\star/A}$ étant $\mU$-petite, 
on pose 
\begin{equation}\label{tfr15b}
F(A)= \underset{\underset{A'\in G_{\star/A}^\circ}{\longleftarrow}}{\lim} F'(A').
\end{equation}
On définit ainsi un préfaisceau sur $G$. On a clairement $\phi^*(F)=F\circ \phi=F'$. Il suffit donc de montrer que $F$ est un faisceau sur $G$,
ou encore que pour tout objet $A$ de $G$ et tout crible couvrant $\cR$ de $A$ dans $G$, le morphisme canonique 
\begin{equation}\label{tfr15c}
F(A)\rightarrow \underset{\underset{A'\in \cR^\circ}{\longleftarrow}}{\lim} F(A')
\end{equation}
est un isomorphisme. On peut se borner aux cribles $\cR$ engendrés par une famille couvrante de morphismes de but $A$ et de type (a), (b) ou (c) (\cite{sga4} II 2.3).

Pour tout objet $A'$ de $G_{/A}$, on désigne par $\cR_{\star/A'}$ la sous catégorie pleine de $G_{\star/A'}$ formée des couples $(B,u)$ tels que 
$\phi(B)$ soit un objet de $\cR$. Si $A'\in \ob(\cR)$, on a $\cR_{\star/A'}=G_{\star/A'}$. Il suffit alors de montrer que le morphisme canonique 
\begin{equation}\label{tfr15d}
F(A)\rightarrow \underset{\underset{A'\in \cR^\circ}{\longleftarrow}}{\lim} \ \underset{\underset{A''\in \cR_{\star/A'}^\circ}{\longleftarrow}}{\lim} F'(A'')
\end{equation}
est un isomorphisme.

Montrons que pour tout objet $A'$ de $G_{\star/A}$, 
$\cR_{\star/A'}$ est un crible couvrant de $A'$ dans $G_\star$. Comme les familles couvrantes de type (a), (b), (c) sont stables
par changement de base, on peut se borner au cas où $A=A'=(W\rightarrow U\leftarrow V)$ est un objet de $G_\star$. 
Comme $W$ est quasi-compact, l'assertion est immédiate si $\cR$ est engendré par une famille couvrante de morphismes de but $A$ et de type (a). 
Par ailleurs, $f$ étant cohérent, $U_Y$ est un schéma cohérent. Comme le morphisme $V\rightarrow U_Y$ est fini étale, $V$ est un schéma cohérent. 
L'assertion recherchée vaut donc aussi si $\cR$ est engendré par une famille couvrante de morphismes de but $A$ et de type (b). 
Supposons enfin que $\cR$ soit engendré par un recouvrement de type (c)
\[
(W\rightarrow U'\leftarrow V')\rightarrow (W\rightarrow U\leftarrow V),
\]
où $U'\rightarrow U$ est un morphisme étale quelconque et $V'=V\times_UU'$. 
Il existe une famille couvrante $(U'_i\rightarrow U')_{i\in I}$ telle que $U'_i\in \ob(\Et_{\star/X})$. Pour tout $i\in I$, posons $W'_i=W\times_{U'}U'_i$ et $V'_i=V'\times_{U'}U'_i$, de sorte que $(W'_i\rightarrow U'_i\leftarrow V'_i)$ est un objet de $G$. 
Il existe une famille couvrante $(W'_{i,j}\rightarrow W'_i)_{j\in J_i}$ telle que $W'_{i,j}\in \ob(\Et_{\star/Z})$. La famille de morphismes de $G$ 
\[
\{(W'_{i,j}\rightarrow U'_i\leftarrow V'_i)\rightarrow (W\rightarrow U\leftarrow V)\}_{i\in I,j\in J_i}
\] 
est alors couvrante en tant que composée de familles de type (a) et (c)
\[
\xymatrix{
W'_{i,j}\ar[r]\ar@{=}[d]&U'_i\ar[d]\ar@{}[rd]|{\Box}&V'_i\ar[l]\ar[d]\\
W'_{i,j}\ar[r]\ar[d]&U\ar@{=}[d]&V\ar[l]\ar@{=}[d]\\
W\ar[r]&U&\ar[l]V}
\]
Comme $W$ est quasi-compact, on peut extraire une sous famille finie $K$ de $\{(i,j) \ | \ i\in I, j\in J_i\}$ telle que $(W'_{i,j}\rightarrow W)_{(i,j)\in K}$ 
soit couvrante. Il s'ensuit que la famille de morphismes de $G_\star$ 
\[
\{(W'_{i,j}\rightarrow U'_i\leftarrow V'_i)\rightarrow (W\rightarrow U\leftarrow V)\}_{(i,j)\in K}
\] 
est couvrante en tant que composée d'une famille {\em finie} de type (a) et de familles de type (c), d'où l'assertion recherchée. 

Il résulte de ce qui précède que pour tout $A'\in \ob(G_{\star/A})$, le morphisme canonique 
\begin{equation}\label{tfr15e}
F'(A')\rightarrow \underset{\underset{A''\in \cR^\circ_{\star/A'}}{\longleftarrow}}{\lim} F'(A'')
\end{equation}
est un isomorphisme. Prenons la limite projective sur la catégorie $G_{\star/A}^\circ$, on en déduit un isomorphisme 
\begin{equation}\label{tfr15f}
F(A)\rightarrow  \underset{\underset{A'\in G^\circ_{\star/A}}{\longleftarrow}}{\lim} \ \underset{\underset{A''\in \cR^\circ_{\star/A'}}{\longleftarrow}}{\lim} F'(A'').
\end{equation}
Les morphismes canoniques  
\begin{equation}\label{tfr15g}
\underset{\underset{A'\in \cR^\circ}{\longleftarrow}}{\lim} \ \underset{\underset{A''\in \cR_{\star/A'}^\circ}{\longleftarrow}}{\lim} F'(A'')
\leftarrow \underset{\underset{A'\in \cR^\circ_{\star/A}}{\longleftarrow}}{\lim} F'(A')\rightarrow \underset{\underset{A'\in G^\circ_{\star/A}}{\longleftarrow}}{\lim} \ \underset{\underset{A''\in \cR^\circ_{\star/A'}}{\longleftarrow}}{\lim} F'(A'')
\end{equation}
étant bijectifs, on en déduit que le morphisme \eqref{tfr15d} est un isomorphisme, d'où la proposition.

\begin{cor}\label{tfr16}
Les notations étant celles de \ref{tfr14}, supposons de plus les schémas $X$, $Y$ et $Z$ cohérents.  Alors, 
\begin{itemize}
\item[{\rm (i)}]  Pour tout objet $(W\rightarrow U\leftarrow V)$ de $G_\star$, le faisceau associé 
$(W\rightarrow U\leftarrow V)^a$ de $\tG$ est cohérent. 
\item[{\rm (ii)}] Le topos $\tG$ est cohérent~; en particulier, il a suffisamment de points.
\item[{\rm (iii)}] Les morphismes $\pi$ \eqref{tfr4c} et $\lambda$ \eqref{tfr4d} sont cohérents. 
\end{itemize}
\end{cor}

Les propositions (i) et (ii) résultent de \ref{tfr15}(ii)-(iii) et (\cite{sga4} VI 2.1) et la proposition (iii) est une conséquence de (i) et (\cite{sga4} VI 3.2).

\subsection{}\label{tfr21}
D'après \ref{topfl17} et (\cite{sga4} VIII 7.9), la donnée d'un point de $Z_\et\gtimes_{X_\et}Y_\et$ 
est équivalente à la donnée d'une paire de points géométriques $\oz$ de $Z$ et $\oy$ de $Y$
et d'une flèche de spécialisation $u$ de $f(\oy)$ vers $g(\oz)$, c'est-à-dire, d'un $X$-morphisme 
$u\colon \oy\rightarrow X_{(g(\oz))}$, où $X_{(g(\oz))}$ désigne le localisé strict de $X$ en $g(\oz)$. 
Un tel point sera noté $(\oy\rightsquigarrow \oz)$ ou encore $(u\colon \oy\rightsquigarrow \oz)$. 
On désigne par $\rho(\oy\rightsquigarrow \oz)$ son image par  le morphisme $\rho\colon Z_\et\gtimes_{X_\et}Y_\et\rightarrow \tG$ \eqref{tfr6b}, 
qui est donc un point de $\tG$.  
Il résulte de \eqref{tfr6c} que pour tous $F\in \ob(Z_\et)$ et $H\in \ob(Y_\fet)$, on a des isomorphismes canoniques fonctoriels 
\begin{eqnarray}
(\pi^*F)_{\rho(\oy\rightsquigarrow \oz)} &\stackrel{\sim}{\rightarrow}& F_{\oz}, \label{tfr21a}\\
(\lambda^*H)_{\rho(\oy\rightsquigarrow \oz)} &\stackrel{\sim}{\rightarrow}& (\rho_Y^*H)_{\oy}, \label{tfr21b}
\end{eqnarray}
où $\rho_Y\colon Y_\et\rightarrow Y_\fet$ est le morphisme canonique \eqref{notconv10a}. 
Pour tout objet $(W\rightarrow U\leftarrow V)$ de $G$, on a un isomorphisme canonique fonctoriel
\begin{equation}\label{tfr21c}
(W\rightarrow U\leftarrow V)^a_{\rho(\oy\rightsquigarrow \oz)} \stackrel{\sim}{\rightarrow} W_{\oz}\times_{U_{f(\oy)}}V_{\oy},
\end{equation}
où l'application $V_{\oy}\rightarrow U_{f(\oy)}$ est induite par $V\rightarrow U$, et 
l'application $W_{\oz}\rightarrow U_{f(\oy)}$ est composée de l'application $W_{\oz}\rightarrow U_{g(\oz)}$ induite  par $W\rightarrow U$
et du morphisme de spécialisation $U_{g(\oz)}\rightarrow U_{f(\oy)}$ défini par $u$. 
En effet, cela résulte de (\cite{agt} VI.3.6) et du fait que $\rho^*$ prolonge $\rho^+$ (\cite{sga4} III 1.4).

\subsection{}\label{tfr22}
Soient $\oz$ un point géométrique de $Z$, $\oy$ un point géométrique de $Y$, 
$Z_{(\oz)}$ le localisé strict de $Z$ en $\oz$, 
$X_{(g(\oz))}$ le localisé strict de $X$ en $g(\oz)$, $u\colon \oy\rightarrow X_{(g(\oz))}$ un $X$-morphisme,
de sorte que $(\oy \rightsquigarrow \oz)$ est un point de $Z_\et\gtimes_{X_\et}Y_\et$ \eqref{tfr21}.
On désigne par $\cP_{\rho(\oy \rightsquigarrow \oz)}$ la catégorie des objets $\rho(\oy \rightsquigarrow \oz)$-pointés de $G$,
définie comme suit. Les objets de $\cP_{\rho(\oy \rightsquigarrow \oz)}$ sont les quadruplets $((W\rightarrow U\leftarrow V),\kappa,\xi,\zeta)$ 
formés d'un objet $(W\rightarrow U\leftarrow V)$ de $G$, d'un $Z$-morphisme $\kappa\colon \oz\rightarrow W$,
d'un $X$-morphisme $\xi\colon \oz\rightarrow U$ et 
d'un $Y$-morphisme $\zeta\colon \oy\rightarrow V$ tels que, notant encore $\kappa\colon Z_{(\oz)}\rightarrow W$
le $Z$-morphisme induit par $\kappa$ (\cite{sga4} VIII 7.3) et 
$\xi\colon X_{(g(\oz))}\rightarrow U$ le $X$-morphisme induit par $\xi$, les carrés du diagramme
\begin{equation}\label{tfr22a}
\xymatrix{
{Z_{(\oz)}}\ar[d]_\kappa\ar[r]^-(0.5)g&{X_{(g(\oz))}}\ar[d]^{\xi}&\oy\ar[l]_-(0.5)u\ar[d]^{\zeta}\\
W\ar[r]&U&V\ar[l]}
\end{equation}
où on a encore noté $g$ le morphisme induit par $g$, soient commutatifs. 
Soient $((W\rightarrow U\leftarrow V),\kappa,\xi,\zeta)$, $((W'\rightarrow U'\leftarrow V'),\kappa',\xi',\zeta')$  deux objets de $\cP_{\rho(\oy \rightsquigarrow \oz)}$.
Un morphisme de $((W'\rightarrow U'\leftarrow V'),\kappa',\xi',\zeta')$ dans $((W\rightarrow U\leftarrow V),\kappa,\xi,\zeta)$ est la donnée d'un 
morphisme $(c\colon W'\rightarrow W, a\colon U'\rightarrow U, b\colon V'\rightarrow V)$ de $G$ tel que $c\circ \kappa'=\kappa$, 
$a\circ \xi'=\xi$ et $b\circ \zeta'=\zeta$. 
Il résulte de \eqref{tfr21c} et (\cite{agt} VI.10.19(i)) que $\cP_{\rho(\oy \rightsquigarrow \oz)}$ est canoniquement équivalente à 
la catégorie des voisinages de 
$\rho(\oy \rightsquigarrow \oz)$ dans $G$ (\cite{sga4} IV 6.8.2). Elle est donc cofiltrante et pour tout préfaisceau 
$F$ sur $G$, on a un isomorphisme canonique fonctoriel (\cite{sga4} IV (6.8.4))
\begin{equation}\label{tfr22b}
(F^a)_{\rho(\oy \rightsquigarrow \oz)}\stackrel{\sim}{\rightarrow} \underset{\underset{((W\rightarrow U\leftarrow V), \kappa,\xi,\zeta)\in 
\cP^\circ_{\rho(\oy \rightsquigarrow \oz)}}{\longrightarrow}}\lim\ F(W\rightarrow U\leftarrow V).
\end{equation}
Si $X$, $Y$ et $Z$ sont cohérents, on peut remplacer dans la limite ci-dessus $\cP_{\rho(\oy \rightsquigarrow \oz)}$ 
par la sous-catégorie pleine $\cP^\coh_{\rho(\oy \rightsquigarrow \oz)}$ 
formée des objets $((W\rightarrow U\leftarrow V),\kappa,\xi,\zeta)$ tels que $W$ soit un objet de $\Et_{\coh/Z}$ 
et $U$ soit un objet de $\Et_{\coh/X}$, qui est aussi cofiltrante \eqref{tfr15}.

\subsection{}\label{tfr23}
Supposons le schéma $X$ strictement local, de point fermé $x$. Soit $\oz$ un point géométrique de $Z$ au-dessus de $x$. 
Le topos $X_\et$ étant local de centre $x$, d'après \ref{topfl18}, le point géométrique $\oz$ de $Z$ définit une section canonique
\begin{equation}\label{tfr23a}
\gamma\colon Y_\et\rightarrow Z_\et\gtimes_{X_\et} Y_\et.
\end{equation}
De même, le topos $X_\fet$ étant local de centre $x$, le point géométrique $\oz$ de $Z$ définit une section canonique
\begin{equation}\label{tfr23b}
\mu\colon Y_\fet\rightarrow Z_\fet\gtimes_{X_\fet} Y_\fet.
\end{equation}
Notons $\rho_X\colon X_\et\rightarrow X_\fet$, $\rho_Y\colon Y_\et\rightarrow Y_\fet$ et $\rho_Z\colon Z_\et\rightarrow Z_\fet$
les morphismes canoniques \eqref{notconv10a}. On vérifie aussitôt que le diagramme 
\begin{equation}\label{tfr23g}
\xymatrix{
{Y_\et}\ar[r]^-(0.5)\gamma\ar[d]_{\rho_Y}&{Z_\et\gtimes_{X_\et}Y_\et}\ar[d]^{\rho_Z\gtimes_{\rho_X}\rho_Y}\\
{Y_\fet}\ar[r]^-(0.5)\mu&{Z_\fet\gtimes_{X_\fet}Y_\fet}}
\end{equation}
où la flèche verticale de droite est définie par fonctorialité \eqref{topfl10d}, est commutatif à isomorphisme canonique près. 

\subsection{}\label{tfr24}
Supposons les schémas $X$ et $Z$ strictement locaux, de points fermés $x$ et $z$ respectivement, 
le morphisme $g\colon Z\rightarrow X$ local, {\em i.e.}, $g(z)=x$, et le morphisme $f\colon Y\rightarrow X$ cohérent.  
Pour tout $U\in \ob(\Et_{\scoh/X})$ \eqref{notconv10}, 
on désigne par $U^\rf$ la somme disjointe des localisés stricts de $U$ en les points de $U_x$;
c'est un sous-schéma ouvert et fermé de $U$, qui est fini sur $X$ (\cite{ega4} 18.5.11). 
De même, pour tout $W\in \ob(\Et_{\scoh/Z})$, 
on désigne par $W^\rf$ la somme disjointe des localisés stricts de $W$ en les points de $W_z$.
Pour tout $U\in \ob(\Et_{\scoh/X})$, on a un isomorphisme canonique 
\begin{equation}\label{tfr24a}
Z\times_XU^\rf\stackrel{\sim}{\rightarrow} (Z\times_XU)^\rf.
\end{equation}

La correspondance $U\mapsto U^\rf$ définit un foncteur 
\begin{equation}\label{tfr24b}
\jmath_x^+\colon \Et_{\scoh/X}\rightarrow \Et_{\rf/X}
\end{equation}
qui est clairement exact à gauche et continu. Notons  
\begin{equation}\label{tfr24c}
\jmath_x\colon X_\fet\rightarrow X_\et
\end{equation} 
morphisme de topos associé (cf. \cite{agt} VI.10.22). 
De même, la correspondance $W\mapsto W^\rf$ définit un foncteur 
\begin{equation}\label{tfr24d}
\jmath_z^+\colon \Et_{\scoh/Z}\rightarrow \Et_{\rf/Z}
\end{equation}
qui est exact à gauche et continu. Notons
\begin{equation}\label{tfr24e}
\jmath_z\colon Z_\fet\rightarrow Z_\et
\end{equation} 
le morphisme de topos associé.

On désigne par $G_\rf$ la sous-catégorie pleine de $G$ formée des objets $(W\rightarrow U\leftarrow V)$ 
tels que les morphismes $U\rightarrow X$ et $W\rightarrow Z$ soient finis étales. Il s'ensuit que le morphisme $V\rightarrow Y$ est aussi fini étale. 
On munit $G_\rf$ de la topologie engendrée par les recouvrements de type (a), (b) et (c) \eqref{tfr1}. 
Le topos des faisceaux de $\mU$-ensembles sur $G_\rf$ est le produit orienté $Z_\fet\gtimes_{X_\fet} Y_\fet$ 
des morphismes de topos $g_\fet\colon Z_\fet\rightarrow X_\fet$ et $f_\fet\colon Y_\fet\rightarrow X_\fet$ (cf. \ref{topfl4}). 

Le foncteur 
\begin{equation}\label{tfr24g}
\nu^+\colon 
\begin{array}[t]{clcr}
G_\scoh&\rightarrow& G_\rf\\
(W\rightarrow U\leftarrow V)&\mapsto&  (W^\rf\rightarrow U^\rf\leftarrow U^\rf\times_UV).
\end{array}
\end{equation}
est exact à gauche et continu. En effet, les foncteurs $\jmath_x$ et $\jmath_z$ étant continus et exacts à gauche,
le foncteur $\vartheta^+$ est exact à gauche, et il transforme les recouvrements finis de $G_\scoh$ de type (a), (b), (c)
en des recouvrements de $G_\rf$ de même type. Il est donc continu en vertu de \ref{topfl2} et \ref{tfr14}(i)-(ii). 
Compte tenu de \ref{tfr15}(iii), le foncteur $\nu^+$ définit donc un morphisme de topos 
\begin{equation}\label{tfr24h}
\nu\colon Z_\fet\gtimes_{X_\fet} Y_\fet\rightarrow \tG.
\end{equation}
Par ailleurs, les topos $X_\fet$ et $Z_\fet$ étant équivalents au topos final, la section canonique définie par le point $z$ de $Z$ \eqref{tfr23b}
\begin{equation}\label{tfr24j}
\mu\colon Y_\fet \rightarrow Z_\fet\gtimes_{X_\fet} Y_\fet
\end{equation}
est une équivalence de topos. 
On désigne le morphisme composé $\nu \mu$ par  
\begin{equation}\label{tfr24i}
\theta\colon Y_\fet\rightarrow \tG.
\end{equation}

\begin{prop}\label{tfr25}
Sous les hypothèses de \ref{tfr24}, le diagramme 
\begin{equation}\label{tfr25a}
\xymatrix{
{Y_\et}\ar[r]^-(0.5)\gamma\ar[d]_{\rho_Y}&{Z_\et\gtimes_{X_\et}Y_\et}\ar[d]^{\rho}\\
{Y_\fet}\ar[r]^-(0.5)\theta&{\tG}}
\end{equation}
où $\gamma$ est la section canonique définie par le point $z$ de $Z$ \eqref{tfr23a}, 
$\rho_Y$ est le morphisme \eqref{notconv10a} et $\rho$ est le morphisme \eqref{tfr6b}, 
est commutatif à isomorphisme canonique près. 
\end{prop}

En effet, notant $i_z\colon \Ens\rightarrow Z_\et$ le point défini par $z$ et $\varepsilon \colon X_\et\rightarrow \Ens$ la projection canonique,  
on a un diagramme non-commutatif \eqref{topfl18c}
\begin{equation}
\xymatrix{
{\Ens}\ar[d]_{i_z}&{Y_\et}\ar[d]_{\gamma}\ar[l]_{\varepsilon f_\et}\ar[rd]^\id&\\
{Z_\et}\ar[rd]_{g_\et}&{Z_\et\gtimes_{X_\et}Y_\et}\ar[l]_-(0.4){\rp_1}\ar[r]^-(0.4){\rp_2}&{Y_\et}\ar[ld]^{f_\et}\\
&{X_\et}&}
\end{equation}
Montrons d'abord que pour tout objet $(W\rightarrow U\leftarrow V)$ de $G$ que l'on considère aussi comme un objet de $C$ \eqref{tfr5}, 
le morphisme canonique 
\begin{equation}\label{tfr25b}
\gamma^*((W^\rf\rightarrow U^\rf\leftarrow V\times_UU^\rf)^a)\rightarrow \gamma^*((W\rightarrow U\leftarrow V)^a)
\end{equation}
est un isomorphisme. Compte tenu de (\cite{agt} (VI.3.7.2)), on a un isomorphisme canonique fonctoriel
\begin{equation}\label{tfr25d}
\gamma^*((W\rightarrow U\leftarrow V)^a) \stackrel{\sim}{\rightarrow} f^*_\et\varepsilon^*i_z^*(W) \times_{f^*_\et\varepsilon^*i_z^*(U\times_XZ)} V.
\end{equation}
Le morphisme canonique $i_z^*(W^\rf)\rightarrow i_z^*(W)$ étant un isomorphisme, il en est donc de même du morphisme 
\begin{equation}\label{tfr25e}
\gamma^*((W^\rf\rightarrow U\leftarrow V)^a)\rightarrow \gamma^*((W\rightarrow U\leftarrow V)^a).
\end{equation}
Comme le morphisme \eqref{tfr25b} est induit par le composé des morphismes canoniques 
\begin{equation}\label{tfr25f}
(W^\rf\rightarrow U^\rf\leftarrow V\times_UU^\rf)\rightarrow (W^\rf\rightarrow U\leftarrow V) \rightarrow (W\rightarrow U\leftarrow V),
\end{equation}
dont le premier est un recouvrement de type (c), on en déduit que c'est un isomorphisme.  

On déduit de l'isomorphisme \eqref{tfr25b} un isomorphisme canonique fonctoriel
\begin{equation}\label{tfr25c}
\gamma^*((\rho_Z\gtimes_{\rho_X}\rho_Y)^*(\nu^*((W\rightarrow U\leftarrow V)^a)))\stackrel{\sim}{\rightarrow} 
\gamma^*(\rho^*(W\rightarrow U\leftarrow V)^a).
\end{equation}
Compte tenu de \eqref{tfr23g}, on en déduit un isomorphisme canonique fonctoriel
\begin{equation}\label{tfr25g}
\rho_Y^*(\theta^*((W\rightarrow U\leftarrow V)^a))\stackrel{\sim}{\rightarrow} \gamma^*(\rho^*(W\rightarrow U\leftarrow V)^a),
\end{equation}
d'où la proposition. 

\begin{cor}\label{tfr250}
Conservons les hypothèses de \ref{tfr25} et soient $\oy$ un point géométrique de $Y$,   
$(\oy \rightsquigarrow z)$ le point de $Z_\et\gtimes_{X_\et}Y_\et$ défini par l'unique spécialisation de $f(\oy)$ dans $x$, 
{\em i.e.}, le $X$-morphisme canonique $u\colon \oy\rightarrow X$. 
Alors, les points $\rho(\oy \rightsquigarrow z)$ et $\theta(\rho_Y(\oy))$ de $\tG$
sont canoniquement isomorphes.
\end{cor}

Cela résulte de \ref{tfr25} quand on observe que les points $\gamma(\oy)$ et $(\oy \rightsquigarrow z)$ de $Z_\et\gtimes_{X_\et}Y_\et$ sont canoniquement isomorphes \eqref{tfr23a}.

\subsection{}\label{tfr26}
Conservons les hypothèses de \ref{tfr24}. On vérifie aussitôt que les carrés du diagramme 
\begin{equation}\label{tfr26a}
\xymatrix{
{Z_\fet}\ar[d]_{\jmath_z}&{Z_\fet\gtimes_{X_\fet} Y_\fet}\ar[l]_-(0.4){\rp_1}\ar[r]^-(0.4){\rp_2}\ar[d]^{\nu}&{Y_\fet}\ar[d]^{\id}\\
{Z_\et}&{\tG}\ar[r]^-(0.5){\lambda}\ar[l]_-(0.5){\pi}&{Y_\fet}}
\end{equation}
sont commutatifs à isomorphismes canoniques près. On en déduit des isomorphismes
\begin{eqnarray}
\pi \theta&\stackrel{\sim}{\rightarrow}& \jmath_z\iota_z\epsilon f_\fet,\label{tfr26b}\\
\lambda \theta&\stackrel{\sim}{\rightarrow}& \id_{Y_\fet},\label{tfr26c}
\end{eqnarray} 
où $\epsilon \colon X_\fet\rightarrow \Ens$ est la projection canonique et 
$\iota_z\colon \Ens\rightarrow Z_\fet$ est le morphisme défini par $z$ \eqref{topfl18c}. 
On obtient du second morphisme un morphisme de changement de base 
\begin{equation}\label{tfr26d}
\lambda_*\rightarrow \theta^*,
\end{equation}
composé de $\lambda_*\rightarrow \lambda_* \theta_*\theta^*\stackrel{\sim}{\rightarrow}\theta^*$,
où la première flèche est induite par le morphisme d'adjonction $\id\rightarrow \theta_*\theta^*$,
et la seconde flèche par l'isomorphisme \eqref{tfr26c}.

\begin{prop}\label{tfr27}
Supposons les schémas $X$ et $Z$ strictement locaux, de points fermés $x$ et $z$ respectivement, 
le morphisme $g\colon Z\rightarrow X$ local, {\em i.e.}, $g(z)=x$, et le morphisme $f\colon Y\rightarrow X$ cohérent,   
et notons $\theta\colon Y_\fet\rightarrow \tG$ le morphisme 
de topos défini dans \eqref{tfr24i}. Alors, le morphisme de changement de base
$\lambda_*\rightarrow \theta^*$ \eqref{tfr26d} est un isomorphisme; en particulier, le foncteur $\lambda_*$ est exact.  
\end{prop} 

Notons $\rho_Y\colon Y_\et\rightarrow Y_\fet$ le morphisme canonique \eqref{notconv10a}. 
Compte tenu de (\cite{agt} VI.9.6), il suffit de montrer que pour tout faisceau $F$ de $\tG$ et tout point géométrique $\oy$ de $Y$, l'application 
\begin{equation}\label{tfr27a}
(\lambda_*F)_{\rho_Y(\oy)}\rightarrow (\theta^*F)_{\rho_Y(\oy)}
\end{equation}
induite par le morphisme de changement de base \eqref{tfr26d} est bijective.
Il existe une unique spécialisation de $f(\oy)$ dans $x$, à savoir le $X$-morphisme canonique $u\colon \oy\rightarrow X$; 
on peut donc considérer le point $(\oy \rightsquigarrow z)$ de $Z_\et\gtimes_{X_\et}Y_\et$ \eqref{tfr21}.
D'après \ref{tfr250}, les points $\rho(\oy \rightsquigarrow z)$ et $\theta(\rho_Y(\oy))$ de $\tG$ sont canoniquement isomorphes.
On désigne par $\cC_\oy$ la catégorie des $Y$-schémas étales et finis, $\oy$-pointés, 
que l'on identifie à la catégorie des voisinages de $\rho_Y(\oy)$ dans le site $\Et_{\rf/Y}$ (\cite{sga4} IV 6.8.2), 
et par $\cP^\scoh_{\rho(\oy\rightsquigarrow z)}$ la sous-catégorie pleine de la catégorie 
$\cP_{\rho(\oy\rightsquigarrow z)}$ \eqref{tfr22} formée des objets $\rho(\oy \rightsquigarrow z)$-pointés $((W\rightarrow U\leftarrow V),\kappa,\xi,\zeta)$ 
de $G$ tels que $U$ soit un objet de $\Et_{\scoh/X}$ et $W$ soit un objet de $\Et_{\scoh/Z}$. 
Rappelons que $\kappa\colon z\rightarrow W$ est un $Z$-morphisme, $\xi\colon x\rightarrow U$ un $X$-morphisme
et $\zeta\colon \oy\rightarrow V$ un $Y$-morphisme tels que notant encore $\kappa\colon Z\rightarrow W$
la $Z$-section induite par $\kappa$ et $\xi\colon X\rightarrow U$ la $X$-section induite par $\xi$, les carrés du diagramme
\begin{equation}\label{tfr27b}
\xymatrix{
{Z}\ar[d]_\kappa\ar[r]^-(0.5)g&{X}\ar[d]^{\xi}&\oy\ar[l]_-(0.5)u\ar[d]^{\zeta}\\
W\ar[r]&U&V\ar[l]}
\end{equation}
soient commutatifs. Il résulte de  \eqref{tfr21c} et (\cite{agt} VI.10.19(i)) 
que $\cP^\scoh_{\rho(\oy \rightsquigarrow z)}$ est canoniquement équivalente 
à la catégorie des voisinages du point $\rho(\oy \rightsquigarrow z)$ dans le site $G_{\scoh}$. Elle est donc cofiltrante. 
Pour tout objet $V$ de $\Et_{\rf/Y}$, on a un isomorphisme canonique $\theta^*(Z\rightarrow X\leftarrow V)\stackrel{\sim}{\rightarrow}V$ \eqref{tfr26c}. 
Compte tenu de l'égalité $\theta(\rho_Y(\oy))=\rho(\oy \rightsquigarrow z)$, on en déduit un isomorphisme 
\begin{equation}\label{tfr27e}
(Z\rightarrow X\leftarrow V)_{\rho(\oy\rightsquigarrow z)}\stackrel{\sim}{\rightarrow}V_\oy.
\end{equation}
Celui-ci est clairement compatible avec l'isomorphisme \eqref{tfr21c}.
Notant $\kappa_0 \colon z\rightarrow Z$ et $\xi_0 \colon x\rightarrow X$ les injections canoniques, on a un foncteur pleinement fidèle
\begin{equation}\label{tfr27c}
\lambda_\oy\colon \cC_\oy\rightarrow \cP^\scoh_{\rho(\oy\rightsquigarrow z)}, \ \ \ (V,\zeta\colon \oy\rightarrow V)\mapsto 
((Z\rightarrow X\leftarrow V),\kappa_0,\xi_0,\zeta).
\end{equation}
Le morphisme d'adjonction $F\rightarrow \theta_*(\theta^*F)$ est défini pour tout objet $(W\rightarrow U\leftarrow V)$ de $G$ par le morphisme
canonique 
\begin{equation}\label{tfr27d}
F(W\rightarrow U\leftarrow V)\rightarrow \theta^*(F)(\theta^*(W\rightarrow U\leftarrow V)).
\end{equation}
Pour tout objet $(V,\zeta\colon \oy\rightarrow V)$ de $\cC_\oy$, le diagramme 
\begin{equation}\label{tfr27f}
\xymatrix{
{F(Z\rightarrow X\leftarrow V)}\ar[r]\ar[dd]&{\theta^*(F)(V)}\ar[d]\\
&{(\theta^*(F)_{\rho_Y(\oy)})(V_\oy)}\ar[d]^\zeta\\
{F_{\rho(\oy\rightsquigarrow z)}((Z\rightarrow X\leftarrow V)^a_{\rho(\oy\rightsquigarrow z)})}\ar[r]^-(0.5){(\kappa_0,\xi_0,\zeta)}&{F_{\rho(\oy\rightsquigarrow z)}}}
\end{equation}
où la flèche notée $\zeta$ (resp. $(\kappa_0,\xi_0,\zeta)$) est induite par $\zeta\in V_\oy$ et la relation $\theta(\rho_Y(\oy))=\rho(\oy \rightsquigarrow z)$
(resp. $(\kappa_0,\xi_0,\zeta)\in (Z\rightarrow X\leftarrow V)^a_{\rho(\oy\rightsquigarrow z)}$), est commutatif \eqref{tfr27e}. 
Par suite, \eqref{tfr27a} s'identifie à l'application
\begin{equation}\label{tfr27g}
\underset{\underset{(V,\zeta)\in \cC_\oy^\circ}{\longrightarrow}}{\lim}\ F(Z\rightarrow X\leftarrow V)
\rightarrow F_{\rho(\oy\rightsquigarrow z)}
\end{equation}
induite par le foncteur $\lambda^\circ_\oy$. Il suffit donc de montrer que  $\lambda_\oy^\circ$ est cofinal.

Soit $((W\rightarrow U\leftarrow V),\kappa,\xi,\zeta)$ un objet de $\cP^\scoh_{\rho(\oy\rightsquigarrow z)}$.
Notons encore $\kappa\colon Z\rightarrow W$
la $Z$-section induite par $\kappa$ et $\xi\colon X\rightarrow U$ la $X$-section induite par $\xi$, de sorte que les carrés du diagramme
\begin{equation}\label{tfr27h}
\xymatrix{
{Z}\ar[d]_\kappa\ar[r]^-(0.5)g&{X}\ar[d]^{\xi}&\oy\ar[l]_-(0.5)u\ar[d]^{\zeta}\\
W\ar[r]&U&V\ar[l]}
\end{equation}
sont commutatifs. 
On en déduit un $Y$-morphisme $\zeta'\colon \oy\rightarrow V\times_{U,\xi}X$ qui s'insère dans 
un diagramme commutatif 
\begin{equation}\label{tfr27i}
\xymatrix{
{Z}\ar[d]_\kappa\ar[r]^-(0.5)g&{X}\ar[d]^{\xi}&{V\times_{U,\xi}X}\ar[l]\ar[d]&\oy\ar[l]_-(0.5){\zeta'}\ar[ld]^-(0.5){\zeta}\\
W\ar[r]&U&V\ar[l]}
\end{equation}
Par suite, $(V\times_{U,\xi}X,\zeta')$ est un objet de $\cC_\oy$, et le diagramme \eqref{tfr27i} induit  un morphisme 
\[
\lambda_\oy(V\times_{U,\xi}X,\zeta')=((Z\rightarrow X \leftarrow V\times_{U,\xi}X),\kappa_0,\xi_0,\zeta')\rightarrow 
((W\rightarrow U\leftarrow V),\kappa,\xi,\zeta)
\] 
de $\cP^\scoh_{\rho(\oy\rightsquigarrow z)}$. On en déduit que  $\lambda_\oy^\circ$ est cofinal d'après (\cite{sga4} I 8.1.3(c));
d'où la proposition.

\begin{cor}\label{tfr280} 
Les hypothèses étant celles de \ref{tfr27}, supposons de plus que le schéma $Y$ soit strictement local de point fermé $y$
et considérons le point $(y \rightsquigarrow z)$ de $Z_\et\gtimes_{X_\et}Y_\et$ défini par l'unique spécialisation de $f(y)$ dans $x$. 
Alors, le topos $\tG$ est local de centre $\rho(y \rightsquigarrow z)$.
\end{cor}

En effet, pour tout faisceau $F$ de $\tG$, on a $\Gamma(\tG,F)=\Gamma(Y_\fet,\lambda_*(F))=(\lambda_*(F))_{\rho_Y(y)}$, où 
$\rho_Y\colon Y_\et\rightarrow Y_\fet$ est le morphisme canonique \eqref{notconv10a}. Par ailleurs, le morphisme 
\begin{equation}
\Gamma(\tG,F)=\Gamma(Y_\fet,\lambda_*(F))\rightarrow \Gamma(Y_\fet,\theta^*(F))
\end{equation}
induit par le morphisme \eqref{tfr26d} s'identifie au morphisme défini par image inverse par $\theta$. 
Compte tenu de \ref{tfr250}, le morphisme canonique 
\begin{equation}
\Gamma(\tG,F)\rightarrow F_{\rho(y \rightsquigarrow z)}
\end{equation}
s'identifie donc à la fibre du morphisme $\lambda_*(F)\rightarrow \theta^*(F)$ \eqref{tfr26d} en $\rho_Y(y)$, 
qui est un isomorphisme en vertu de \ref{tfr27}, d'où la proposition. 

\begin{cor}\label{tfr28}
Supposons les schémas $X$ et $Z$ strictement locaux, de points fermés $x$ et $z$ respectivement, 
le morphisme $g\colon Z\rightarrow X$ local, {\em i.e.}, $g(z)=x$, et le morphisme $f\colon Y\rightarrow X$ cohérent,   
et notons $\theta\colon Y_\fet\rightarrow \tG$ le morphisme 
de topos défini dans \eqref{tfr24i}.  Alors~:
\begin{itemize}
\item[{\rm (i)}] Pour tout faisceau $F$ de $\tG$, l'application canonique 
\begin{equation}\label{tfr28a}
\Gamma(\tG,F)\rightarrow \Gamma(Y_\fet, \theta^*F)
\end{equation}
est bijective. 
\item[{\rm (ii)}] Pour tout faisceau abélien $F$ de $\tG$, l'application canonique 
\begin{equation}\label{tfr28b}
\rH^i(\tG,F)\rightarrow \rH^i(Y_\fet, \theta^*F)
\end{equation}
est bijective pour tout $i\geq 0$.
\end{itemize}
\end{cor}

(i) En effet, le diagramme 
\begin{equation}
\xymatrix{
{\Gamma(\tG,F)}\ar[d]_v&{\Gamma(Y_\fet,\lambda_*F)}\ar[d]^u\ar[l]^-(0.5)\sim_-(0.5)w\\
{\Gamma(\tG,\theta_*(\theta^*F))}\ar[rd]_{v'}&{\Gamma(Y_\fet,\lambda_*(\theta_*(\theta^*F)))}\ar[d]^{u'}
\ar[l]_-(0.4){w'}^-(0.4)\sim\\
&{\Gamma(Y_\fet,\theta^*F)}}
\end{equation}
où $w$, $w'$ et $v'$ sont les bijections canoniques, $u$ et $v$ sont 
induits par le morphisme d'adjonction $\id\rightarrow \theta_*\theta^*$ et
$u'$ est induit par l'isomorphisme \eqref{tfr26c}, est  commutatif. Comme $u'\circ u$ est bijectif en vertu de \ref{tfr27}, 
il en est de même de $v'\circ v$, d'où la proposition.  

(ii) En effet, le diagramme 
\begin{equation}
\xymatrix{
{\rH^i(\tG,F)}\ar[d]_v&{\rH^i(Y_\fet,\lambda_*F)}\ar[d]^u\ar[l]_-(0.5)w\\
{\rH^i(\tG,\theta_*(\theta^*F))}\ar[rd]_{v'}&{\rH^i(Y_\fet,\lambda_*(\theta_*(\theta^*F)))}\ar[d]^{u'}\ar[l]_-(0.4){w'}\\
&{\rH^i(Y_\fet,\theta^*F)}}
\end{equation}
où $w$, $w'$ et $v'$ sont induits par la suite spectrale de Cartan-Leray (\cite{sga4} V 5.3), $u$ et $v$ sont 
induits par le morphisme d'adjonction $\id\rightarrow \theta_*\theta^*$ et
$u'$ est induit par l'isomorphisme \eqref{tfr26c}, est  commutatif. D'autre part, $u'\circ u$ est bijectif,
et le foncteur $\lambda_*$ est exact en vertu de \ref{tfr27}; donc $w$ est bijectif, d'où la proposition. 

\subsection{}\label{tfr29}
Soient $\oz$ un point géométrique de $Z$, $\uZ$ le localisé strict de $Z$ en $\oz$, 
$\uX$ le localisé strict de $X$ en $g(\oz)$, $\ug\colon \uZ\rightarrow \uX$ le morphisme induit par $g$, 
$\uY=Y\times_X\uX$, $\uf\colon \uY\rightarrow \uX$ la projection canonique. 
On désigne par $\tuG$ le topos de Faltings relatif associé au couple de morphismes $(\uf,\ug)$ \eqref{tfr1} et par
\begin{equation}\label{tfr29a}
\theta\colon \uY_\fet\rightarrow \tuG
\end{equation}
le morphisme défini dans \eqref{tfr24i}. Le diagramme commutatif de morphismes canoniques 
\begin{equation}
\xymatrix{
\uZ\ar[r]^\ug\ar[d]&\uX\ar[d]&\uY\ar[l]_{\uf}\ar[d]\\
Z\ar[r]^g&X&Y\ar[l]_f}
\end{equation}
induit par fonctorialité un morphisme \eqref{tfr7b}
\begin{equation}\label{tfr29b}
\Phi\colon \tuG\rightarrow \tG.
\end{equation}
On note 
\begin{equation}\label{tfr29c}
\varphi_\oz\colon \tG\rightarrow \uY_\fet
\end{equation}
le foncteur composé $\theta^*\circ \Phi^*$.

\begin{prop}\label{tfr30}
Conservons les hypothèses de \ref{tfr29}, supposons de plus que $f$ soit cohérent. Alors~:
\begin{itemize}
\item[{\rm (i)}] Pour tout faisceau $F$ de $\tG$, on a un isomorphisme canonique fonctoriel 
\begin{equation}\label{tfr30a}
\pi_*(F)_\oz\stackrel{\sim}{\rightarrow}\Gamma(\uY_\fet,\varphi_\oz(F)).
\end{equation}
\item[{\rm (ii)}] Pour tout faisceau abélien $F$ de $\tG$ et tout entier $i\geq 0$, on a un isomorphisme canonique fonctoriel 
\begin{equation}\label{tfr30b}
\rR^i\pi_*(F)_\oz\stackrel{\sim}{\rightarrow}\rH^i(\uY_\fet,\varphi_\oz(F)).
\end{equation}
\item[{\rm (iii)}] Pour toute suite exacte de faisceaux abéliens 
$0\rightarrow F'\rightarrow F\rightarrow F''\rightarrow 0$ de $\tG$ et tout entier $i\geq 0$, le diagramme 
\begin{equation}\label{tfr30c}
\xymatrix{
{\rR^i\pi_*(F'')_\oz} \ar[r]\ar[d]&{\rR^{i+1}\pi_*(F')_\oz}\ar[d]\\
{\rH^i(\uY_\fet,\varphi_\oz(F''))}\ar[r]&{\rH^{i+1}(\uY_\fet,\varphi_\oz(F'))}}
\end{equation}
où les flèches verticales sont les isomorphismes canoniques \eqref{tfr30b} et les flèches horizontales  
sont les bords des suites exactes longues de cohomologie, est commutatif. 
\end{itemize}
\end{prop}

Cette proposition sera démontrée dans \ref{lptfr8}.

\section{Limites projectives de topos de Faltings relatifs}\label{lptfr}

\subsection{}\label{lptfr1}
On désigne par $\fM$ la catégorie des couples de morphismes de schémas de même but
$(Z\rightarrow X\leftarrow Y)$ et par $\fD$ la catégorie des morphismes de $\fM$.
Les objets de $\fD$ sont donc des diagrammes commutatifs de morphismes de schémas 
\begin{equation}\label{lptfr1a}
\xymatrix{
W\ar[r]\ar[d]&U\ar[d]&V\ar[l]\ar[d]\\
Z\ar[r]&X&Y\ar[l]}
\end{equation} 
où l'on considère les flèches horizontales comme des objets de $\fM$
et les flèches verticales comme des morphismes de $\fM$; un tel objet sera noté $(W\rightarrow U\leftarrow V,Z\rightarrow X\leftarrow Y)$.   
On désigne par $\fG$ la sous-catégorie pleine de $\fD$ formée des objets $(W\rightarrow U\leftarrow V,Z\rightarrow X\leftarrow Y)$ tels que les morphismes
$U\rightarrow X$ et $W\rightarrow Z$ soient étales de présentation finie et que le morphisme $V\rightarrow U\times_XY$ soit étale fini. 
Le ``foncteur but''
\begin{equation}\label{lptfr1b}
\fG\rightarrow \fM, \ \ \ (W\rightarrow U\leftarrow V,Z\rightarrow X\leftarrow Y)\mapsto (Z\rightarrow X\leftarrow Y)
\end{equation}
fait de $\fG$ une catégorie fibrée, clivée et normalisée (cf. \ref{tfr2}).  
La catégorie fibre au-dessus d'un objet $(Z\rightarrow X\leftarrow Y)$ de $\fM$ est la catégorie sous-jacente  au site 
$G_\coh$ défini dans \ref{tfr14} relativement au couple de morphismes $(Y\rightarrow X, Z\rightarrow X)$. 
Pour tout morphisme de $\fM$ défini par un diagramme commutatif
\begin{equation}\label{lptfr1c}
\xymatrix{
Z'\ar[r]\ar[d]&X'\ar[d]&Y'\ar[l]\ar[d]\\
Z\ar[r]&X&Y\ar[l]}
\end{equation}
le foncteur image inverse de \eqref{lptfr1b} associé est le foncteur
\begin{equation}\label{lptfr1d}
\Phi^+\colon 
\begin{array}[t]{clcr}\fG_{(Z\rightarrow X\leftarrow Y)}&\rightarrow &\fG_{(Z'\rightarrow X'\leftarrow Y')}, \\
(W\rightarrow U\leftarrow V)&\mapsto & (W\times_ZZ'\rightarrow U\times_XX'\leftarrow V\times_YY').
\end{array}
\end{equation}
Celui-ci est clairement exact à gauche.
On munit chaque fibre de $\fG/\fM$ de la topologie co-évanes\-cente, c'est-à-dire de la topologie engendrée par les recouvrements de type (a), (b) et (c) 
définis dans \ref{tfr1}. On vérifie aussitôt que la caractérisation des faisceaux de \ref{tfr3} vaut encore pour les fibres de $\fG/\fM$. 
On en déduit que pour tout morphisme de $\fM$ défini par un diagramme commutatif \eqref{lptfr1c}, 
le foncteur $\Phi^+$ \eqref{lptfr1d} est continu. 
Par suite, $\fG/\fM$ devient un $\mU$-site fibré (\cite{sga4} VI 7.2.4). On désigne par 
\begin{equation}\label{lptfr1e}
\fF\rightarrow \fM
\end{equation}
le $\mU$-topos fibré associé à $\fG/\fM$ (\cite{sga4} VI 7.2.6)~: 
la catégorie fibre de $\fF$ au-dessus d'un objet $(Z\rightarrow X\leftarrow Y)$
de $\fM$ est le topos $\tfG_{(Z\rightarrow X\leftarrow Y)}$ des faisceaux de $\mU$-ensembles sur le site co-évanescent 
$\fG_{(Z\rightarrow X\leftarrow Y)}$, et le foncteur image inverse relatif à un morphisme de $\fM$ défini par un diagramme commutatif \eqref{lptfr1c} 
est le foncteur $\Phi^*\colon \tfG_{(Z\rightarrow X\leftarrow Y)}\rightarrow \tfG_{(Z'\rightarrow X'\leftarrow Y')}$ 
image inverse par le morphisme de topos 
\begin{equation}\label{lptfr1f}
\Phi\colon \tfG_{(Z'\rightarrow X'\leftarrow Y')}\rightarrow \tfG_{(Z\rightarrow X\leftarrow Y)}
\end{equation} 
associé au morphisme de sites $\Phi^+$ \eqref{lptfr1d}. 

Pour tout objet $(Z\rightarrow X\leftarrow Y)$ de $\fM$ tel que les schémas $X$, $Y$ et $Z$ soient cohérents, la topologie de 
$\fG_{(Z\rightarrow X\leftarrow Y)}$ est engendrée par les familles couvrantes de type (c) et les familles couvrantes {\em finies} de type (a), (b).  
Par suite, le topos $\tfG_{(Z\rightarrow X\leftarrow Y)}$ est canoniquement équivalent au topos de Faltings relatif associé au couple de morphismes
$(Y\rightarrow X, Z\rightarrow X)$ d'après \ref{tfr15}(iii). 

On note 
\begin{equation}\label{lptfr1g}
\fF^\vee\rightarrow \fM^\circ
\end{equation}
la catégorie fibrée obtenue en associant à tout objet $(Z\rightarrow X\leftarrow Y)$ de $\fM$ 
la catégorie $\fF_{(Z\rightarrow X\leftarrow Y)}=\tfG_{(Z\rightarrow X\leftarrow Y)}$, 
et à tout morphisme de $\fM$ défini par un diagramme commutatif \eqref{lptfr1c} le foncteur 
\begin{equation}\label{lptfr1h}
\Phi_*\colon \tfG_{(Z'\rightarrow X'\leftarrow Y')}\rightarrow \tfG_{(Z\rightarrow X\leftarrow Y)}
\end{equation} 
image directe par le morphisme de topos $\Phi$ \eqref{lptfr1f}.

\subsection{}\label{lptfr2}
Soient $I$ une catégorie cofiltrante essentiellement petite (\cite{sga4} I 2.7 et 8.1.8), 
\begin{equation}\label{lptfr2a}
\varphi\colon I\rightarrow \fM, \ \ \ i\mapsto (Z_i\rightarrow X_i\leftarrow Y_i)
\end{equation}
un foncteur tel que pour tout morphisme $j\rightarrow i$ de $I$, 
les morphismes $X_j\rightarrow X_i$, $Y_j\rightarrow Y_i$ et $Z_j\rightarrow Z_i$ soient affines. 
On suppose qu'il existe $i_0\in \ob(I)$ tel que les schémas $X_{i_0}$, $Y_{i_0}$ et $Z_{i_0}$ soient cohérents.
On désigne par 
\begin{eqnarray}
\fG_\varphi&\rightarrow& I\label{lptfr2b}\\
\fF_\varphi&\rightarrow& I\label{lptfr2bb}\\
\fF^\vee_\varphi&\rightarrow& I\label{lptfr2bbb}
\end{eqnarray}
les site, topos et catégorie fibrés déduits de $\fG$ \eqref{lptfr1b},  $\fF$ \eqref{lptfr1e} et  $\fF^\vee$ \eqref{lptfr1g}, respectivement,
par changement de base par le foncteur $\varphi$. On notera que $\fF_\varphi$ est le topos fibré associé à $\fG_\varphi$
(\cite{sga4} VI 7.2.6.8). D'après (\cite{ega4} 8.2.3), les limites projectives 
\begin{equation}\label{lptfr2c}
X=\underset{\underset{i\in I}{\longleftarrow}}{\lim}\ X_i,\ \ 
Y=\underset{\underset{i\in I}{\longleftarrow}}{\lim}\ Y_i  \ \ {\rm et}\ \ Z=\underset{\underset{i\in I}{\longleftarrow}}{\lim}\ Z_i  
\end{equation}
sont représentables dans la catégorie des schémas. Les schémas $X$, $Y$ et $Z$ sont cohérents. 
Les objets $(Z_i\rightarrow X_i\leftarrow Y_i)_{i\in I}$ de $\fM$ induisent 
un objet $(Z\rightarrow X\leftarrow Y)$, qui représente la limite projective du foncteur \eqref{lptfr2a}. 

Pour tout $i\in \ob(I)$, on a un diagramme commutatif canonique
\begin{equation}\label{lptfr2d}
\xymatrix{
Z\ar[r]\ar[d]&X\ar[d]&Y\ar[l]\ar[d]\\
Z_i\ar[r]&X_i&Y_i\ar[l]}
\end{equation}
Il lui correspond un foncteur image inverse \eqref{tfr7d}
\begin{equation}\label{lptfr2e}
\Phi_i^+\colon \fG_{(Z_i\rightarrow X_i\leftarrow Y_i)}\rightarrow \fG_{(Z\rightarrow X\leftarrow Y)},
\end{equation}
qui est continu et exact à gauche, et par suite un morphisme de topos 
\begin{equation}\label{lptfr2f}
\Phi_i\colon \tfG_{(Z\rightarrow X\leftarrow Y)}\rightarrow \tfG_{(Z_i\rightarrow X_i\leftarrow Y_i)}.
\end{equation}

On a un foncteur naturel 
\begin{equation}\label{lptfr2g}
\fG_\varphi \rightarrow \fG_{(Z\rightarrow X\leftarrow Y)},
\end{equation}
dont la restriction à la fibre au-dessus de tout $i\in \ob(I)$ est le foncteur $\Phi_i^+$ \eqref{lptfr2e}. 
Ce foncteur transforme morphisme cartésien en isomorphisme. Il se factorise donc de façon unique à travers 
un foncteur (\cite{sga4} VI 6.3)
\begin{equation}\label{lptfr2h}
\underset{\underset{I^\circ}{\longrightarrow}}{\lim}\ \fG_\varphi \rightarrow \fG_{(Z\rightarrow X\leftarrow Y)}.
\end{equation}
Le $I$-foncteur $\fG_\varphi\rightarrow \fG_{(Z\rightarrow X\leftarrow Y)}\times I$ déduit de \eqref{lptfr2g} est un morphisme cartésien
de sites fibrés (\cite{sga4} VI 7.2.2). Il induit donc un morphisme cartésien de topos fibrés (\cite{sga4} VI 7.2.7)
\begin{equation}\label{lptfr2i}
\tfG_{(Z\rightarrow X\leftarrow Y)}\times I\rightarrow \fF_\varphi.
\end{equation}

\begin{prop}\label{lptfr3}
Le couple formé du topos $\tfG_{(Z\rightarrow X\leftarrow Y)}$ et du morphisme \eqref{lptfr2i} 
est une limite projective du topos fibré $\fF_\varphi/I$ {\rm (\cite{sga4} VI 8.1.1)}.
\end{prop} 

On notera d'abord que le foncteur \eqref{lptfr2h} est une équivalence de catégories en vertu de 
(\cite{ega4} 8.8.2,  8.10.5 et 17.7.8). Soient $T$ un $\mU$-topos, 
\begin{equation}\label{lptfr3a}
h\colon T\times I\rightarrow \fF_\varphi
\end{equation}
un morphisme cartésien de topos fibrés au-dessus de $I$. Notons $\varepsilon_I\colon \fG_\varphi\rightarrow \fF_\varphi$
le foncteur cartésien canonique (\cite{sga4} VI (7.2.6.7)), et posons 
\begin{equation}\label{lptfr3b}
h^+= h^*\circ \varepsilon_I\colon \fG_\varphi\rightarrow T\times I.
\end{equation}
Pour tout $i\in \ob(I)$, on désigne par
\begin{equation}\label{lptfr3c}
h_i^+\colon \fG_{(Z_i\rightarrow X_i\leftarrow Y_i)}\rightarrow T
\end{equation}
la restriction de $h^+$ aux fibres au-dessus de $i$. 
Compte tenu de l'équivalence de catégories \eqref{lptfr2h} et de  (\cite{sga4} VI 6.2), il existe un et essentiellement 
un unique foncteur 
\begin{equation}\label{lptfr3d}
g^+\colon \fG_{(Z\rightarrow X\leftarrow Y)}\rightarrow T
\end{equation}
tel que $h^+$ soit isomorphe au composé 
\begin{equation}\label{lptfr3e}
\xymatrix{
{\fG_\varphi}\ar[r]&{\fG_{(Z\rightarrow X\leftarrow Y)}\times I}\ar[rr]^-(0.5){g^+\times \id_I}&&{T\times I}},
\end{equation}
où la première flèche est le foncteur déduit de \eqref{lptfr2g}. Montrons que $g^+$ est un morphisme de sites. 
Pour tout objet $(W\rightarrow U\leftarrow V)$ de $\fG_{(Z\rightarrow X\leftarrow Y)}$, 
il existe $i\in \ob(I)$, un objet $(W_i\rightarrow U_i\leftarrow V_i)$ de $\fG_{(Z_i\rightarrow X_i\leftarrow Y_i)}$ et un isomorphisme de 
$\fG_{(Z\rightarrow X\leftarrow Y)}$
\begin{equation}
(W\rightarrow U\leftarrow V)\stackrel{\sim}{\rightarrow} \Phi_i^+(W_i\rightarrow U_i\leftarrow V_i).
\end{equation}
Comme les foncteurs $h_i^+$ et $\Phi_i^+$ sont exacts à gauche, on en déduit que $g^+$ est exact à gauche. 
D'autre part, tout recouvrement de type (c), (resp. {\em fini} de type (a), resp. {\em fini} de type (b)) de $\fG_{(Z\rightarrow X\leftarrow Y)}$  
est l'image inverse d'un recouvrement de type (c) (resp. fini de type (a), resp. fini de type (b)) de $\fG_{(Z_i\rightarrow X_i\leftarrow Y_i)}$ 
pour un objet $i$ de $I$,  en vertu de (\cite{ega4} 8.10.5(vi)). 
On en déduit que $g^+$ transforme les recouvrements {\em finis} de type (a), (resp. (b)) de $\fG_{(Z\rightarrow X\leftarrow Y)}$ 
en familles épimorphiques de $T$, et les recouvrements de type (c) en isomorphismes de $T$. 
Comme les schémas $X$, $Y$ et $Z$ sont cohérents, $g^+$ est continu en vertu de \ref{tfr14} et \ref{tfr15}.  
Il définit donc un morphisme de topos 
\begin{equation}\label{lptfr3f}
g\colon T\rightarrow \tfG_{(Z\rightarrow X\leftarrow Y)}
\end{equation}
tel que $h$ soit isomorphe au composé 
\begin{equation}\label{lptfr3g}
\xymatrix{
{T\times I}\ar[rr]^-(0.5){g\times \id_I}&&{\tfG_{(Z\rightarrow X\leftarrow Y)}\times I}\ar[r]&{\fF_\varphi}},
\end{equation}
où la seconde flèche est le morphisme \eqref{lptfr2i}. Un tel morphisme $g$ est essentiellement unique car la ``restriction'' 
$g^+\colon \fG_{(Z\rightarrow X\leftarrow Y)}\rightarrow T$ du foncteur $g^*$ est essentiellement unique d'après ce qui précède, d'où la proposition.

\subsection{}\label{lptfr4}
Munissons $\fG_\varphi$ de la topologie totale (\cite{sga4} VI 7.4.1) 
et notons $\Top(\fG_\varphi)$ le topos des faisceaux de $\mU$-ensembles sur $\fG_\varphi$. 
D'après (\cite{sga4} VI 7.4.7), on a une équivalence canonique de catégories
\begin{equation}\label{lptfr4a}
\Top(\fG_\varphi)\stackrel{\sim}{\rightarrow}\bHom_{I^\circ}(I^\circ, \fF^\vee_\varphi). 
\end{equation}
D'autre part, le foncteur naturel $\fG_\varphi\rightarrow \fG_{(Z\rightarrow X\leftarrow Y)}$ \eqref{lptfr2g} 
est un morphisme de sites (\cite{sga4} VI 7.4.4) et il définit donc un morphisme de topos 
\begin{equation}\label{lptfr4b}
\varpi\colon \tfG_{(Z\rightarrow X\leftarrow Y)}\rightarrow \Top(\fG_\varphi). 
\end{equation}
En vertu de \ref{lptfr3} et (\cite{sga4} VI 8.2.9), il existe une équivalence de catégories $\Theta$
qui s'insère dans un diagramme commutatif 
\begin{equation}\label{lptfr4c}
\xymatrix{
{\tfG_{(Z\rightarrow X\leftarrow Y)}}\ar[r]^-(0.5){\Theta}_-(0.5)\sim\ar[d]_{\varpi_*}&{\bHom_{\cart/I^\circ}(I^\circ, \fF^\vee_\varphi)}\ar@{^(->}[d]\\
{\Top(\fG_\varphi)}\ar[r]^-(0.5)\sim&{\bHom_{I^\circ}(I^\circ, \fF^\vee_\varphi)}}
\end{equation}
où la flèche horizontale inférieure est l'équivalence de catégories \eqref{lptfr4a} et 
la flèche verticale de droite est l'injection canonique. 
Pour tout objet $F$ de $\tfG_{(Z\rightarrow X\leftarrow Y)}$, on a (\cite{sga4} VI 8.5.1)
\begin{equation}\label{lptfr4e}
\varpi_*(F)=\{i\mapsto \Phi_{i*}(F)\}.
\end{equation}
Pour tout objet $F=\{i\mapsto F_i\}$ de $\Top(\fG_\varphi)$, on a un 
isomorphisme canonique fonctoriel (\cite{sga4} VI 8.5.2)
\begin{equation}\label{lptfr4d}
\varpi^*(F)\stackrel{\sim}{\rightarrow} \underset{\underset{i\in I^\circ}{\longrightarrow}}{\lim}\ \Phi_i^*(F_i).
\end{equation}

Le foncteur $\varpi_*$ étant pleinement fidèle, le morphisme d'adjonction $\varpi^*\varpi_*\rightarrow \id$ est un isomorphisme. 
On en déduit que pour tout objet $F$ de $\tfG_{(Z\rightarrow X\leftarrow Y)}$, le morphisme canonique 
\begin{equation}\label{lptfr4f}
\underset{\underset{i\in I^\circ}{\longrightarrow}}{\lim}\ \Phi_i^*(\Phi_{i*}(F)) \stackrel{\sim}{\rightarrow} F
\end{equation}
est un isomorphisme. 

\begin{cor}\label{lptfr5}
Soient $F$ un faisceau de $\Top(\fG_\varphi)$, 
$\{i\mapsto F_i\}$ la section de $\bHom_{I^\circ}(I^\circ, \fF^\vee_\varphi)$ qui lui est 
associée par l'équivalence de catégories \eqref{lptfr4a}.
Alors, on a un isomorphisme canonique fonctoriel
\begin{equation}\label{lptfr5a}
\underset{\underset{i\in I^\circ}{\longrightarrow}}{\lim}\ \Gamma(\tfG_{(Z_i\rightarrow X_i\leftarrow Y_i)},F_i)\stackrel{\sim}{\rightarrow}
\Gamma(\tfG_{(Z\rightarrow X\leftarrow Y)},\underset{\underset{i\in I^\circ}{\longrightarrow}}{\lim}\ \Phi_i^*(F_i)).
\end{equation}
\end{cor}

\begin{cor}\label{lptfr6}
Soit $F$ un faisceau abélien de $\Top(\fG_\varphi)$ et soit
$\{i\mapsto F_i\}$ la section de $\bHom_{I^\circ}(I^\circ, \fF^\vee_\varphi)$ qui lui est 
associée par l'équivalence de catégories \eqref{lptfr4a}.
Alors, pour tout entier $q\geq 0$, on a un isomorphisme canonique fonctoriel
\begin{equation}\label{lptfr6a}
\underset{\underset{i\in I^\circ}{\longrightarrow}}{\lim}\ \rH^q(\tfG_{(Z_i\rightarrow X_i\leftarrow Y_i)},F_i)\stackrel{\sim}{\rightarrow}
\rH^q(\tfG_{(Z\rightarrow X\leftarrow Y)},\underset{\underset{i\in I^\circ}{\longrightarrow}}{\lim}\ \Phi_i^*(F_i)).
\end{equation}
\end{cor}

Les corollaires \ref{lptfr5} et \ref{lptfr6} résultent de \ref{lptfr3} et (\cite{sga4} VI 8.7.7). 
On notera que les conditions requises dans 
(\cite{sga4} VI 8.7.1 et 8.7.7) sont satisfaites en vertu de \ref{tfr16} et (\cite{sga4} VI 3.3, 5.1 et 5.2).

\begin{prop}\label{lptfr7}
Considérons un diagramme commutatif de morphismes de schémas à carré cartésien
\begin{equation}\label{lptfr7a}
\xymatrix{
Z\ar[r]^{g'}\ar[dr]_g&X'\ar[d]^u\ar@{}[rd]|\Box&Y'\ar[l]_{f'}\ar[d]^v\\
&X&\ar[l]_fY}
\end{equation}
et notons $\tG$ (resp. $\tG'$) le topos de Faltings relatif associé au couple de morphismes  $(f,g)$ (resp. $(f',g')$)
et $\Phi\colon \tG'\rightarrow \tG$ le morphisme de fonctorialité induit par le diagramme \eqref{lptfr7a}.
Supposons, de plus, que l'une des hypothèses suivantes soit remplie:
\begin{itemize}
\item[{\rm (i)}] le morphisme $u\colon X'\rightarrow X$ est étale;
\item[{\rm (ii)}] le schéma $X'$ est le localisé strict de $X$ en un point géométrique $\ox$, 
$u \colon X'\rightarrow X$ étant le morphisme canonique, le morphisme $f$ est cohérent et le schéma $Z$ est cohérent. 
\end{itemize}
Alors, $\Phi$ est une équivalence de topos.
\end{prop}

Notons $G$ (resp. $G'$) le site de Faltings relatif associé au couple de morphismes  $(f,g)$ (resp. $(f',g')$) \eqref{tfr1}. 
Dans le cas (i), $(Z\rightarrow X'\leftarrow Y')$ est un objet de $G$ et le morphisme 
\begin{equation}\label{lptfr7b}
(Z\rightarrow X'\leftarrow Y')^a\rightarrow (Z\rightarrow X\leftarrow Y)^a
\end{equation}
est un isomorphisme de $\tG$ en vertu de \ref{tfr3}. Le morphisme $\Phi$ qui s'identifie au morphisme de localisation du topos $\tG$ en 
$(Z\rightarrow X'\leftarrow Y')^a$ d'après \ref{tfr42}, est donc une équivalence de topos. 

Considérons le cas (ii). D'après le cas (i), quitte à remplacer $X$ par un voisinage étale affine $X_1$ de $\ox$ et $Y$ par $X_1\times_XY$,  
on peut supposer $X$ affine et $Y$ cohérent. Notons $I$ la catégorie des $X$-schémas étales affines $\ox$-pointés, 
et pour tout $i\in \ob(I)$, $X_i$ le $X$-schéma étale affine $\ox$-pointé correspondant et posons $Y_i=X_i\times_XY$. 
On a alors dans la catégorie des schémas
\begin{equation}\label{lptfr7c}
X'=\underset{\underset{i\in I}{\longleftarrow}}{\lim}\ X_i\ \ {\rm et}\ \ 
Y'=\underset{\underset{i\in I}{\longleftarrow}}{\lim}\ Y_i.
\end{equation}
Considérons le foncteur
\begin{equation}\label{lptfr7d}
\varphi\colon I\rightarrow \fM_{/(Z\rightarrow X\leftarrow Y)}, \ \ \ i\mapsto (Z\rightarrow X_i\leftarrow Y_i),
\end{equation}
où $\fM$ est la catégorie définie dans \ref{lptfr1}.
Reprenant les notations de \ref{lptfr2} pour ce foncteur $\varphi$, le topos $\tfG_{(Z\rightarrow X'\leftarrow Y')}$ 
est la limite projective du topos fibré $\fF_\varphi/I$ en vertu de \ref{lptfr3}. Par ailleurs, pour tout $i\in \ob(I)$, le morphisme
de fonctorialité 
\begin{equation}
\tfG_{(Z\rightarrow X_i\leftarrow Y_i)}\rightarrow \tfG_{(Z\rightarrow X\leftarrow Y)}
\end{equation}
est une équivalence de topos d'après le cas (i). La proposition s'ensuit.

\begin{lem}\label{lptfr18}
Soient $f\colon Y\rightarrow X$, $g\colon Z\rightarrow X$ deux morphismes tels que $f$ soit cohérent, 
$G$ le site de Faltings relatif associé à $(f,g)$, $G_\coh$ la sous-catégorie pleine définie dans \ref{tfr14}. 
Pour tout point géométrique $\oz$ de $Z$, on désigne par $Z_{(\oz)}$ le localisé strict de $Z$ en $\oz$, 
par $h_\oz\colon Z_{(\oz)}\rightarrow Z$ le morphisme canonique, 
par $G_{\oz}$ le site de Faltings relatif associé à $(f,g\circ h_{\oz})$ et 
par $\Phi^+_{\oz}\colon G\rightarrow G_{\oz}$ le morphisme de fonctorialité induit par $h_{\oz}$ \eqref{tfr7}. 
Soit $(A_n\rightarrow A)_{n\in N}$ une famille de morphismes de $G_\coh$.
Alors, pour que $(A_n\rightarrow A)_{n\in N}$ soit couvrante dans $G$,
il faut et il suffit que pour tout point géométrique $\oz$ de $Z$, la famille $(\Phi^+_\oz(A_n)\rightarrow \Phi^+_\oz(A))_{n\in N}$ 
soit couvrante dans $G_\oz$.
\end{lem}

En effet, la condition est nécessaire puisque les foncteurs $\Phi^+_\oz$ sont continus (\cite{sga4} III 1.6). 
Montrons qu'elle est suffisante. Supposons que pour tout point géométrique $\oz$ de $Z$, la famille 
$(\Phi^+_\oz(A_n)\rightarrow \Phi^+_\oz(A))_{n\in N}$ soit couvrante dans $G_\oz$ et montrons que 
 la famille $(A_n\rightarrow A)_{n\in N}$ est couvrante dans $G$. 
 La question étant locale sur $G$ \eqref{tfr41}, on peut se borner au cas où $X$, $Y$ et $Z$ sont cohérents. 
 
Soient $\oz$ un point géométrique de $Z$, $Z_0$ un $Z$-schéma étale $\oz$-pointé qui soit un schéma affine.
On note $I$ la catégorie des $Z_0$-schémas étales affines $\oz$-pointés, et pour tout $i\in \ob(I)$, 
$Z_i$ le $Z$-schéma étale $\oz$-pointé correspondant. On a alors 
\begin{equation}\label{lptfr18a}
Z_{(\oz)}=\underset{\underset{i\in I}{\longleftarrow}}{\lim}\ Z_i.
\end{equation}
Considérons le foncteur
\begin{equation}\label{lptfr18b}
\varphi\colon I\rightarrow \fM_{/(Z\rightarrow X\leftarrow Y)}, \ \ \ i\mapsto (Z_i\rightarrow X\leftarrow Y),
\end{equation}
où $\fM$ est la catégorie définie dans \ref{lptfr1}. Reprenons les notations de \ref{lptfr2} en les adaptant à ce foncteur. 
Pour tout $i\in \ob(I)$, notons
\begin{eqnarray}
\Psi_i^+\colon G_\coh=\fG_{(Z\rightarrow X\leftarrow Y)}\rightarrow \fG_{(Z_i\rightarrow X\leftarrow Y)},\label{lptfr18c}\\ 
\Phi_i^+\colon \fG_{(Z_i\rightarrow X\leftarrow Y)}\rightarrow \fG_{(Z_{(\oz)}\rightarrow X\leftarrow Y)}, \label{lptfr18d}
\end{eqnarray}
les foncteurs image inverse \eqref{tfr7d}, de sorte que $\Phi_i^+\circ \Psi_i^+$ est la restriction du foncteur $\Phi_{\oz}$. 
On a un foncteur naturel 
\begin{equation}\label{lptfr18e}
\fG_\varphi \rightarrow \fG_{(Z_{(\oz)}\rightarrow X\leftarrow Y)},
\end{equation}
dont la restriction à la fibre au-dessus de tout $i\in \ob(I)$ est le foncteur $\Phi_i^+$. 
Ce foncteur transforme morphisme cartésien en isomorphisme. Il se factorise donc de façon unique à travers 
un foncteur
\begin{equation}\label{lptfr18f}
\underset{\underset{I^\circ}{\longrightarrow}}{\lim}\ \fG_\varphi \rightarrow \fG_{(Z_{(\oz)}\rightarrow X\leftarrow Y)}, 
\end{equation}
qui est une équivalence de catégories en vertu de (\cite{ega4} 8.8.2,  8.10.5 et 17.7.8).
D'après (\cite{ega4} 8.10.5(vi)), tout recouvrement de type (c), (resp. {\em fini} de type (a), resp. {\em fini} de type (b)) de 
$\fG_{(Z_{(\oz)}\rightarrow X\leftarrow Y)}$  est l'image inverse d'un recouvrement de type (c) (resp. fini de type (a), resp. fini de type (b)) de 
$\fG_{(Z_i\rightarrow X\leftarrow Y)}$ pour un objet $i$ de $I$. 

La famille de morphismes $(\Phi^+_\oz(A_n)\rightarrow \Phi^+_\oz(A))_{n\in N}$ de $\fG_{(Z_{(\oz)}\rightarrow X\leftarrow Y)}$
étant couvrante dans $G_\oz$, elle l'est aussi dans $\fG_{(Z_{(\oz)}\rightarrow X\leftarrow Y)}$ en vertu de \ref{tfr15}(iii) et (\cite{sga4} III 3.5). 
Compte tenu de ce qui précède et de \ref{tfr15}(ii), il existe donc un objet $i$ de $I$ et une famille couvrante finie 
$(B_m\rightarrow \Psi_i^+(A))_{m\in M}$ de $\fG_{(Z_i\rightarrow X\leftarrow Y)}$ tels que le recouvrement 
$(\Phi_i^+(B_m)\rightarrow \Phi_\oz^+(A))_{m\in M}$ raffine le recouvrement $(\Phi^+_\oz(A_n)\rightarrow \Phi^+_\oz(A))_{n\in N}$. 
Quitte à remplacer $i$ par un objet de $I$ au-dessus de lui, le recouvrement $(B_m\rightarrow \Psi_i^+(A))_{m\in M}$ 
raffine alors la famille $(\Psi^+_i(A_n)\rightarrow \Psi^+_i(A))_{n\in N}$ qui est donc couvrante. 
On notera que $(Z_i\rightarrow X\leftarrow Y)$ est un objet de $G_\coh$ et que le site $\fG_{(Z_i\rightarrow X\leftarrow Y)}$ s'identifie 
à la catégorie $G_{\coh/(Z_i\rightarrow X\leftarrow Y)}$ munie de la topologie induite par celle de $G_\coh$ \eqref{tfr40}. 

On a ainsi montré qu'il existe un recouvrement $(Z_j\rightarrow X\leftarrow Y)_{j\in J}$ de type (a) de l'objet final de $G$ 
tel que pour tout $j\in J$, la famille $(A_n\times (Z_j\rightarrow X\leftarrow Y) \rightarrow A\times (Z_j\rightarrow X\leftarrow Y))_{n\in N}$ 
soit couvrante dans $G$. Par suite, la famille $(A_n\rightarrow A)_{n\in N}$ est couvrante dans $G$; d'où la proposition.

\begin{prop}\label{lptfr19}
Soient $X$, $Y$, $Z$ trois schémas cohérents, $f\colon Y\rightarrow X$, $g\colon Z\rightarrow X$ deux morphismes.
On désigne par $G$ le site de Faltings relatif associé au couple de morphismes $(f,g)$, par $G_\coh$ la sous-catégorie pleine définie dans 
\ref{tfr14} et par $C$ le site défini dans \ref{topfl1} relativement aux foncteurs de changement de base par $f$ et $g$ respectivement, 
\begin{equation}
f^+\colon \Et_{/X}\rightarrow \Et_{/Y}\ \ \  {\rm et} \ \ \ g^+\colon \Et_{/X}\rightarrow \Et_{/Z}.
\end{equation}
La catégorie $G$ est naturellement une sous-catégorie pleine de $C$ \eqref{tfr6}. 
Alors, pour qu'une famille $((W_i\rightarrow U_i\leftarrow V_i)\rightarrow (W\rightarrow U\leftarrow V))_{i\in I}$ 
de morphismes de $G_\coh$ soit couvrante dans $G$, il faut et il suffit qu'elle le soit dans $C$.
\end{prop}

En effet, la condition est nécessaire puisque le foncteur canonique $\rho^+\colon G\rightarrow C$ \eqref{tfr6a} est continu (\cite{sga4} III 1.6).
Montrons qu'une famille $\cF=((W_i\rightarrow U_i\leftarrow V_i)\rightarrow (W\rightarrow U\leftarrow V))_{i\in I}$
de morphismes de $G_\coh$ qui est couvrante dans $C$ est couvrante dans $G$. On peut supposer que  
$(W\rightarrow U\leftarrow V)$ est l'objet final $(Z\rightarrow X\leftarrow Y)$ de $G$ \eqref{tfr41}. 
En vertu de \ref{lptfr18}, on peut se borner au cas où le schéma $Z$ est strictement local. 

Considérons les sites $\Et_{\coh/X}$, $\Et_{\coh/Y}$ et $\Et_{\coh/Z}$ \eqref{notconv10} qui sont stables par limites projectives finies. 
On désigne par $C_\coh$ le site défini dans \ref{topfl1} relativement aux foncteurs de changement de base par $f$ et $g$, 
\begin{equation}
f^+_\coh\colon \Et_{\coh/X}\rightarrow \Et_{\coh/Y}\ \ \  {\rm et} \ \ \ g^+_\coh\colon \Et_{\coh/X}\rightarrow \Et_{\coh/Z}.
\end{equation}
Tout objet de $C_\coh$ est naturellement un objet de $C$. On définit ainsi
un foncteur pleinement fidèle $C_\coh\rightarrow C$. Il est clairement exact à gauche et il est continu d'après \ref{topfl2}. 
Le morphisme de topos qu'il induit $\tC\rightarrow \tC_\coh$ est une équivalence de catégories d'après la propriété universelle 
des produits orientés de topos \eqref{topfl4}.  
Par suite, une famille de morphismes de même but de $C_\coh$ est couvrante dans $C_\coh$ si et seulement si 
elle l'est dans $C$ en vertu de (\cite{sga4} III 3.5). 

Par ailleurs, tout objet de $G_{\coh}$ est naturellement un objet de $C_\coh$. 
On peut donc considérer $\cF$ comme une famille couvrante de $C_\coh$. 
D'après \ref{topfl15}(ii), la topologie de $C_\coh$ est engendrée par la prétopologie définie pour chaque objet $A$ de $C_\coh$
par la donnée de l'ensemble $\Cov(A)$ des familles de morphismes $(A_j\rightarrow A)_{j\in J}$ obtenues par composition 
d'un nombre fini de familles de type (c) et de familles finies de type (a) et (b); en particulier, l'ensemble $J$ est fini.
Comme $Z$ est strictement local, tout recouvrement de $\Cov(Z\rightarrow X\leftarrow Y)$ peut être raffiné par un 
recouvrement composé d'un nombre fini de familles de type (c) et de familles finies de type (b). 

En vertu de \ref{tfr32}, il existe donc un recouvrement de type (c),  
$(Z\rightarrow X'\leftarrow Y')\rightarrow (Z \rightarrow X\leftarrow Y)$ et un recouvrement fini de type (b),
$((Z\rightarrow X'\leftarrow Y'_j)\rightarrow (Z \rightarrow X'\leftarrow Y'))_{j\in J}$ de $C_\coh$, tels que 
le recouvrement composé $((Z\rightarrow X'\leftarrow Y'_j)\rightarrow (Z \rightarrow X\leftarrow Y))_{j\in J}$ raffine $\cF$. 
Par suite, pour tout $j\in J$, il existe $i\in I$ et un morphisme de $C$, représenté par le diagramme commutatif 
\begin{equation}
\xymatrix{
Z\ar[r]\ar[d]&X'\ar[d]&Y'_j\ar[d]\ar[l]\\
W_i\ar[r]&U_i&V_i\ar[l]}
\end{equation}
Posons $Y''_j=X'\times_{U_i}V_i$. Comme $(W_i\rightarrow U_i\leftarrow V_i)$ est un objet de $G$, 
le morphisme canonique $V_i\rightarrow U_i\times_XY$ est (étale) fini. Il en est donc de même du morphisme $Y''_j\rightarrow X'\times_XY=Y'$. 
Par suite, $(Y''_j\rightarrow Y')_{j\in J}$ est un recouvrement de $\Et_{\rf/Y'}$,
de sorte que $((Z\rightarrow X'\leftarrow Y''_j)\rightarrow (Z \rightarrow X'\leftarrow Y'))_{j\in J}$
est un recouvrement de type (b) de $G_\coh$, ce qui implique que $\cF$ est un recouvrement de $G_\coh$.

\begin{prop}\label{lptfr20}
Soient $X$, $Y$, $Z$ trois schémas cohérents, $f\colon Y\rightarrow X$, $g\colon Z\rightarrow X$ deux morphismes.
On désigne par $\tG$ le topos de Faltings relatif associé au couple de morphismes $(f,g)$ et par 
\begin{equation}
\rho\colon Z_\et\gtimes_{X_\et}Y_\et\rightarrow \tG
\end{equation}
le morphisme canonique \eqref{tfr6b}. 
Alors, lorsque $(\oy \rightsquigarrow \oz)$ décrit la famille des points de $Z_\et\gtimes_{X_\et}Y_\et$ \eqref{topfl17}, 
la famille des foncteurs fibres de $\tG$ associés aux points $\rho(\oy \rightsquigarrow \oz)$ est conservative
{\rm (\cite{sga4} IV 6.4.0)}. 
\end{prop}

On désigne par par $G$ le site de Faltings relatif associé au couple de morphismes $(f,g)$ \eqref{tfr1}, 
par $G_\coh$ la sous-catégorie pleine de $G$ définie dans \ref{tfr14}, 
par $C$ le site défini dans \ref{topfl1} relativement aux foncteurs de changement de base par $f$ et $g$ respectivement, 
\begin{equation}
f^+\colon \Et_{/X}\rightarrow \Et_{/Y}\ \ \  {\rm et} \ \ \ g^+\colon \Et_{/X}\rightarrow \Et_{/Z}.
\end{equation}
On rappelle que le topos de faisceaux de $\mU$-ensembles sur $C$ est canoniquement équivalent à $Z_\et\gtimes_{X_\et}Y_\et$,
que la catégorie $G$ est naturellement une sous-catégorie pleine de $C$ et que le foncteur d'injection canonique $\rho^+\colon G\rightarrow C$
est continu et exact à gauche; il induit le morphisme $\rho$ (cf. \ref{tfr6}).  

Soit $u\colon F\rightarrow F'$ un morphisme de $\tG$ tel que pour tout point $(\oy \rightsquigarrow \oz)$
de $Z_\et\gtimes_{X_\et}Y_\et$, le morphisme correspondant 
$u_{\rho(\oy \rightsquigarrow \oz)}\colon F_{\rho(\oy \rightsquigarrow \oz)}\rightarrow 
F'_{\rho(\oy \rightsquigarrow \oz)}$ soit un monomorphisme. Montrons que $u$ est un monomorphisme.
Il s'agit de montrer que pour tous $(W\rightarrow U\leftarrow V)\in \ob(G_\coh)$ 
et $a,b \in F(W\rightarrow U\leftarrow V)$ tels que $u(a)=u(b)$, on a $a=b$. 
Compte tenu de \ref{tfr42} et \eqref{tfr7e}, on peut supposer $(W\rightarrow U\leftarrow V)=(Z\rightarrow X\leftarrow Y)$. 
Pour tout point $(\oy \rightsquigarrow \oz)$ de $Z_\et\gtimes_{X_\et}Y_\et$, on a 
$a_{\rho(\oy \rightsquigarrow \oz)}=b_{\rho(\oy \rightsquigarrow \oz)}$ puisque 
$u_{\rho(\oy \rightsquigarrow \oz)}$ est un monomorphisme. 
On désigne par $\cP_{\rho(\oy \rightsquigarrow \oz)}$ la catégorie des objets $\rho(\oy \rightsquigarrow \oz)$-pointés de $G$ définie par
\ref{tfr22} et par $\cP^\coh_{\rho(\oy \rightsquigarrow \oz)}$ la sous-catégorie pleine formée des objets 
$((W\rightarrow U\leftarrow V),\kappa,\xi,\zeta)$ tels que $(W\rightarrow U\leftarrow V)$ soit un objet de $G_\coh$. 
Il résulte de \eqref{tfr21c} et (\cite{agt} VI.10.19(i)) que $\cP^\coh_{\rho(\oy \rightsquigarrow \oz)}$ est canoniquement équivalente à 
la catégorie des voisinages de $\rho(\oy \rightsquigarrow \oz)$ dans $G_\coh$ (\cite{sga4} IV 6.8.2).
Il existe donc un objet 
$((W_{(\oy \rightsquigarrow \oz)}\rightarrow  U_{(\oy \rightsquigarrow \oz)}\leftarrow V_{(\oy \rightsquigarrow \oz)}),
\xi_{(\oy \rightsquigarrow \oz)},\kappa_{(\oy \rightsquigarrow \oz)},
\zeta_{(\oy \rightsquigarrow \oz)})$ de $\cP^\coh_{\rho(\oy \rightsquigarrow \oz)}$ tel que $a$ et $b$ aient mêmes images dans
$F(W_{(\oy \rightsquigarrow \oz)}\rightarrow  U_{(\oy \rightsquigarrow \oz)}\leftarrow V_{(\oy \rightsquigarrow \oz)})$. 
Par ailleurs, les topos $X_\et$, $Y_\et$ et $Z_\et$ étant cohérents et 
les morphismes $f\colon Y_\et\rightarrow X_\et$ et $g\colon Z_\et\rightarrow X_\et$ étant cohérents (\cite{sga4} VI 3.3), 
la famille des foncteurs fibres de $Z_\et\gtimes_{X_\et}Y_\et$ associés aux points 
$(\oy \rightsquigarrow \oz)$ est conservative, d'après \ref{topfl16}(i) et \ref{topfl17}.
On en déduit que la famille des morphismes 
\begin{equation}
\{(W_{(\oy \rightsquigarrow \oz)}\rightarrow  U_{(\oy \rightsquigarrow \oz)}\leftarrow V_{(\oy \rightsquigarrow \oz)})
\rightarrow (Z\rightarrow X\leftarrow Y)\}_{(\oy \rightsquigarrow \oz)}
\end{equation}
est couvrante dans $C$ (\cite{sga4} IV 6.5). Elle est donc couvrante dans $G_\coh$ en vertu de \ref{lptfr19}.
Par suite, $a=b$ et $u$ est un monomorphisme. 

Supposons, de plus, que pour tout point $(\oy \rightsquigarrow \oz)$
de $Z_\et\gtimes_{X_\et}Y_\et$, le morphisme 
$u_{\rho(\oy \rightsquigarrow \oz)}$ soit un épimorphisme et montrons qu'il en est de même de $u$. 
Il suffit de montrer que pour tous $(W\rightarrow U\leftarrow V)\in \ob(G_\coh)$ et $b \in F'(W\rightarrow U\leftarrow V)$, 
il existe $a\in F(W\rightarrow U\leftarrow V)$ tel que $b=u(a)$. 
On peut encore supposer $(W\rightarrow U\leftarrow V)=(X\rightarrow X\leftarrow Y)$.
D'après \eqref{tfr22b}, pour tout point $(\oy \rightsquigarrow \oz)$ de $Z_\et\gtimes_{X_\et}Y_\et$, 
il existe un objet $((W_{(\oy \rightsquigarrow \oz)}\rightarrow  U_{(\oy \rightsquigarrow \oz)}\leftarrow V_{(\oy \rightsquigarrow \oz)}),
\kappa_{(\oy \rightsquigarrow \oz)},\xi_{(\oy \rightsquigarrow \oz)}, \zeta_{(\oy \rightsquigarrow \oz)})$ 
de $\cP^\coh_{\rho(\oy \rightsquigarrow \oz)}$ et une section 
$a_{(\oy \rightsquigarrow \oz)}\in F(W_{(\oy \rightsquigarrow \oz)}\rightarrow  U_{(\oy \rightsquigarrow \oz)}\leftarrow V_{(\oy \rightsquigarrow \oz)})$
dont l'image par $u$ dans $F'(W_{(\oy \rightsquigarrow \oz)}\rightarrow  U_{(\oy \rightsquigarrow \oz)}\leftarrow V_{(\oy \rightsquigarrow \oz)})$ 
soit la restriction de $b$. 
Comme $u$ est un monomorphisme, les sections $a_{(\oy \rightsquigarrow \oz)}$
coïncident sur 
\begin{equation}
(W_{(\oy \rightsquigarrow \oz)}\times_ZW_{(\oy' \rightsquigarrow \oz')}\rightarrow  
U_{(\oy \rightsquigarrow \oz)}\times_XU_{(\oy' \rightsquigarrow \oz')}\leftarrow 
V_{(\oy \rightsquigarrow \oz)}\times_YV_{(\oy' \rightsquigarrow \oz')}),
\end{equation}
pour tous les points $(\oy \rightsquigarrow \oz)$
et $(\oy' \rightsquigarrow \oz')$ de $Z_\et\gtimes_{X_\et}Y_\et$ \eqref{tfr2}. 
Ils proviennent donc d'une section $a\in F(X\rightarrow X\leftarrow Y)$, et on a $u(a)=b$ puisque les restrictions 
aux $(W_{(\oy \rightsquigarrow \oz)}\rightarrow  U_{(\oy \rightsquigarrow \oz)}\leftarrow V_{(\oy \rightsquigarrow \oz)})$ coïncident, 
d'où la proposition.

\begin{prop}\label{lptfr12}
Considérons un diagramme commutatif de morphismes de schémas 
\begin{equation}\label{lptfr12c}
\xymatrix{
Z'\ar[r]^{g'}\ar[d]_w&X'\ar[d]_u\ar@{}[rd]|\Box&Y'\ar[l]_{f'}\ar[d]^v\\
Z\ar[r]^g&X&Y\ar[l]_f}
\end{equation}
tel que $f$ soit cohérent, $X'$ soit cohérent, $Z'$ soit le localisé strict de $Z$ en un point géométrique $\oz$, $w$ étant le morphisme canonique,
$u$ soit étale et le carré de droite soit cartésien. 
Notons $\tG$ (resp. $\tG'$) le topos de Faltings relatif associé au couple de morphismes  
$(f,g)$ (resp. $(f',g')$), $\Phi\colon \tG'\rightarrow \tG$ le morphisme de fonctorialité induit par le diagramme \eqref{lptfr12c}
et $\pi\colon \tG\rightarrow Z_\et$ le morphisme canonique \eqref{tfr4c}.
Alors, 
\begin{itemize}
\item[{\rm (i)}] Pour tout faisceau $F$ de $\tG$, on a un isomorphisme canonique fonctoriel 
\begin{equation}\label{lptfr12a}
\pi_*(F)_{\oz}\stackrel{\sim}{\rightarrow}\Gamma(\tG',\Phi^*(F)).
\end{equation}
\item[{\rm (ii)}] Pour tout faisceau abélien $F$ de $\tG$ et tout entier $q\geq 0$, on a un isomorphisme canonique fonctoriel 
\begin{equation}\label{lptfr12b}
\rR^q\pi_*(F)_{\oz}\stackrel{\sim}{\rightarrow}\rH^q(\tG',\Phi^*(F)).
\end{equation}
\end{itemize}
\end{prop}

On désigne par $J$ la catégorie des $Z$-schémas étales $\oz$-pointés 
et par $I$ la sous-catégorie pleine de $J_{/Z\times_XX'}$ formée des objets $(\oz\rightarrow W)$ 
tels que le schéma $W$ soit affine. Ce sont des catégories cofiltrantes, et le foncteur canonique
$I^\circ\rightarrow J^\circ$ est cofinal (\cite{sga4} I 8.1.3(c)). On a dans la catégorie des schémas 
\begin{equation}
Z'=\underset{\underset{(\oz\rightarrow W)\in I}{\longleftarrow}}{\lim}\ W.
\end{equation}
Considérons le foncteur
\begin{equation}
\varphi\colon I\rightarrow \fM, \ \ \ (\oz\rightarrow W)\mapsto (W\rightarrow X'\leftarrow Y'),
\end{equation}
où $\fM$ est la catégorie définie dans \ref{lptfr1}.
Pour tout $i=(\oz\rightarrow W)\in \ob(I)$, notons
\begin{equation}
\Psi_i\colon \tfG_{(W\rightarrow X'\leftarrow Y')}\rightarrow \tfG_{(Z\rightarrow X\leftarrow Y)}=\tG
\end{equation}
le morphisme canonique \eqref{lptfr2f}. 
Reprenons les notations de \ref{lptfr4} en les adaptant; on note
\begin{equation}
\varpi\colon \tG'\rightarrow \Top(\fG_\varphi)
\end{equation}
le morphisme canonique de topos \eqref{lptfr4b}.
Pour tout faisceau $F$ de $\tG$, $\{i\mapsto \Psi_i^*(F)\}$ est naturellement une section de $\bHom_{I^\circ}(I^\circ,\fF_\varphi^\vee)$.
Elle définit donc un faisceau de $\Top(\fG_\varphi)$ \eqref{lptfr4a}.  On a alors un isomorphisme canonique fonctoriel \eqref{lptfr4d}
\begin{equation}
\Phi^*(F)\stackrel{\sim}{\rightarrow} \varpi^*(\{i\mapsto \Psi_i^*(F)\}).
\end{equation}

(i) D'après (\cite{sga4} VIII 3.9), on a un isomorphisme canonique fonctoriel
\begin{equation}
\pi_*(F)_{\oz}\stackrel{\sim}{\rightarrow} \underset{\underset{(\oz\rightarrow W)\in I^\circ}{\longrightarrow}}{\lim}\ 
\Gamma((W\rightarrow X\leftarrow Y)^a,F). 
\end{equation}
Le morphisme canonique $(W\rightarrow X'\leftarrow Y')^a\rightarrow (W\rightarrow X\leftarrow Y)^a$ 
est un isomorphisme d'après \ref{tfr3}.
On en déduit, en vertu de \ref{tfr42} et \ref{lptfr5}, un isomorphisme fonctoriel
\begin{equation}
\pi_*(F)_\oz\stackrel{\sim}{\rightarrow} \Gamma(\tG',\Phi^*(F)).
\end{equation}

(ii) D'après (\cite{sga4} VIII 3.9 et V 5.1), pour tout entier $q\geq 0$, on a un isomorphisme canonique fonctoriel
\begin{equation}
\rR^q\pi_*(F)_{\oz}\stackrel{\sim}{\rightarrow} \underset{\underset{(\oz\rightarrow W)\in I^\circ}{\longrightarrow}}{\lim}\ 
\rH^q((W\rightarrow X\leftarrow Y)^a,F). 
\end{equation}
Le morphisme canonique $(W\rightarrow X'\leftarrow Y')^a\rightarrow (W\rightarrow X\leftarrow Y)^a$ 
étant un isomorphisme, on en déduit, en vertu de \ref{tfr42} et \ref{lptfr6}, un isomorphisme fonctoriel
\begin{equation}
\rR^q\pi_*(F)_\oz\stackrel{\sim}{\rightarrow} \rH^q(\tG',\Phi^*(F)).
\end{equation}

\subsection{}\label{lptfr8}
Nous pouvons maintenant démontrer la proposition \ref{tfr30}. Reprenons les notations de \ref{tfr29}. 

(i) En vertu de \ref{lptfr7} et \ref{lptfr12}(i),  on a un isomorphisme fonctoriel
\begin{equation}\label{lptfr8a}
\pi_*(F)_\oz\stackrel{\sim}{\rightarrow} \Gamma(\tuG,\Phi^*(F)).
\end{equation}
La proposition s'en déduit compte tenu de \ref{tfr28}(i).

(ii) En vertu de \ref{lptfr7}, \ref{lptfr12}(ii), on a un isomorphisme fonctoriel
\begin{equation}\label{lptfr8b}
\rR^q\pi_*(F)_\oz\stackrel{\sim}{\rightarrow} \rH^q(\tuG,\Phi^*(F)).
\end{equation}
La proposition s'en déduit compte tenu de \ref{tfr28}(ii).

(iii) Cela résulte aussitôt de la preuve de (ii).

\begin{prop}\label{lptfr9}
Considérons un diagramme commutatif de morphismes de schémas 
\begin{equation}\label{lptfr9a}
\xymatrix{
Z'\ar[r]^{g'}\ar[d]_{h_Z}&X'\ar[d]_{h_X}&Y'\ar[l]_{f'}\ar[d]^{h_Y}\\
Z\ar[r]^g&X&\ar[l]_fY}
\end{equation}
tel que les flèches verticales soient des morphismes cohérents, 
et notons $\tG$ (resp. $\tG'$) le topos de Faltings relatif associé au couple de morphismes  
$(f,g)$ (resp. $(f',g')$), $\Phi\colon \tG'\rightarrow \tG$ le morphisme de fonctorialité induit par le diagramme \eqref{lptfr9a}
et $\rho\colon Z_\et\gtimes_{X_\et}Y_\et\rightarrow \tG$ le morphisme canonique \eqref{tfr6b}.
Soient $(\oy\rightsquigarrow \oz)$ un point de $Z_\et\gtimes_{X_\et}Y_\et$ \eqref{tfr21}, $\uZ$ le localisé strict de $Z$ en $\oz$, 
$\uX$ le localisé strict de $X$ en $\ox=g(\oz)$, $\uY$ le localisé strict de $Y$ en $\oy$, $\ug\colon \uZ\rightarrow \uX$ le morphisme induit par $g$, 
$\uf\colon \uY\rightarrow \uX$ le morphisme induit par $f$ et le point $(\oy\rightsquigarrow \oz)$ {\rm (\cite{sga4} VIII 7.3)}. 
Posons $\uX'=\uX\times_XX'$, $\uY'=\uY\times_YY'$ et $\uZ'=\uZ\times_ZZ'$. Les morphismes $f'$, $g'$, $\uf$ et $\ug$ induisent 
des morphismes $\uf'\colon \uY'\rightarrow \uX'$ et $\ug'\colon \uZ'\rightarrow \uX'$ qui s'insèrent dans un diagramme commutatif 
\begin{equation}\label{lptfr9b}
\xymatrix{
Z'\ar[r]^{g'}&X'&Y'\ar[l]_{f'}\\
\uZ'\ar[r]^{\ug'}\ar[d]\ar[u]&\uX'\ar[d]\ar[u]&\uY'\ar[l]_{\uf'}\ar[d]\ar[u]\\
\uZ\ar[r]^\ug&\uX&\ar[l]_\uf\uY}
\end{equation}
Notons $\tuG'$ le topos de Faltings relatif associé au couple de morphismes  $(\uf',\ug')$
et $\Psi\colon \tuG'\rightarrow \tG'$ le morphisme de fonctorialité induit par les deux carrés supérieurs du diagramme \eqref{lptfr9b}. 
Alors, 
\begin{itemize}
\item[{\rm (i)}] Pour tout faisceau $F$ de $\tG'$, on a un isomorphisme canonique fonctoriel 
\begin{equation}\label{lptfr9c}
\Phi_*(F)_{\rho(\oy\rightsquigarrow \oz)}\stackrel{\sim}{\rightarrow}\Gamma(\tuG',\Psi^*(F)).
\end{equation}
\item[{\rm (ii)}] Pour tout faisceau abélien $F$ de $\tG'$ et tout entier $q\geq 0$, on a un isomorphisme canonique fonctoriel 
\begin{equation}\label{lptfr9d}
\rR^q\Phi_*(F)_{\rho(\oy\rightsquigarrow \oz)}\stackrel{\sim}{\rightarrow}\rH^q(\tuG',\Psi^*(F)). 
\end{equation}
\end{itemize}
\end{prop}

Soit $\cP_{\rho(\oy \rightsquigarrow \oz)}$ la catégorie des objets $\rho(\oy \rightsquigarrow \oz)$-pointés de $G$ (cf. \ref{tfr22}). 
Pour tout objet $((W\rightarrow U\leftarrow V),\kappa,\xi,\zeta)$ de $\cP_{\rho(\oy \rightsquigarrow \oz)}$, posons 
$(W'\rightarrow U'\leftarrow V')=\Phi^+(W\rightarrow U\leftarrow V)$ \eqref{tfr7d}.   
Fixons un objet $((W_0\rightarrow U_0\leftarrow V_0),\kappa_0,\xi_0,\zeta_0)$ de $\cP_{\rho(\oy \rightsquigarrow \oz)}$
tel que les schémas $U_0$, $V_0$ et $W_0$ soient affines, et notons $I$ la sous-catégorie pleine de 
\begin{equation}
(\cP_{\rho(\oy \rightsquigarrow \oz)})_{/((W_0\rightarrow U_0\leftarrow V_0),\kappa_0,\xi_0,\zeta_0)}
\end{equation}
formée des objets $((W\rightarrow U\leftarrow V),\kappa,\xi,\zeta)$ tels que les schémas $U$, $V$ et $W$ soient affines.
Les catégories $\cP_{\rho(\oy \rightsquigarrow \oz)}$ et $I$ sont cofiltrantes 
et le foncteur canonique $I^\circ\rightarrow \cP^\circ_{\rho(\oy \rightsquigarrow \oz)}$ est cofinal (\cite{sga4} IV 6.8.2 et I 8.1.3(c)).

Pour tout $i\in \ob(I)$, notons $((W_i\rightarrow U_i\leftarrow V_i),\kappa_i,\xi_i,\zeta_i)$ l'objet de
$\cP_{\rho(\oy \rightsquigarrow \oz)}$ correspondant. 
Le foncteur de la catégorie $I^\circ$ dans la catégorie opposée de la catégorie des $Z$-schémas étales $\oz$-pointés 
(resp. $X$-schémas étales $\ox$-pointés, resp. $Y$-schémas étales $\oy$-pointés)
qui à $i$ associe $(W_i,\kappa_i)$ (resp. $(U_i,\xi_i)$, resp. $(V_i,\zeta_i)$) est cofinal. Cela résulte de (\cite{sga4} I 8.1.3(b))
et de la description explicite de la catégorie $\cP_{\rho(\oy \rightsquigarrow \oz)}$ dans \ref{tfr22}. 
On en déduit que l'on a dans la catégorie des schémas
\begin{equation}
\uX=\underset{\underset{i\in I}{\longleftarrow}}{\lim}\ U_i,\ \ 
\uY=\underset{\underset{i\in I}{\longleftarrow}}{\lim}\ V_i  \ \ {\rm et}\ \ \uZ=\underset{\underset{i\in I}{\longleftarrow}}{\lim}\ W_i,
\end{equation} 
et par suite 
\begin{equation}
\uX'=\underset{\underset{i\in I}{\longleftarrow}}{\lim}\ U'_i,\ \ 
\uY'=\underset{\underset{i\in I}{\longleftarrow}}{\lim}\ V'_i  \ \ {\rm et}\ \ \uZ'=\underset{\underset{i\in I}{\longleftarrow}}{\lim}\ W'_i.
\end{equation} 
Par ailleurs, il résulte des hypothèses que les schémas $U'_i$, $V'_i$ et $W'_i$ sont cohérents. 
Considérons le foncteur 
\begin{equation}
\varphi \colon I\rightarrow \fM, \ \ \ i\mapsto (W'_i\rightarrow U'_i\leftarrow V'_i),
\end{equation} 
où $\fM$ est la catégorie définie dans \ref{lptfr1}.
Pour tout $i\in \ob(I)$, notons
\begin{equation}
\Psi_i\colon \tfG_{(W'_i\rightarrow U'_i\leftarrow V'_i)}\rightarrow \tfG_{(Z'\rightarrow X'\leftarrow Y')}=\tG'
\end{equation}
le morphisme canonique \eqref{lptfr2f}. Reprenons les notations de \ref{lptfr4} en les adaptant; on note
\begin{equation}
\varpi\colon \tuG'\rightarrow \Top(\fG_\varphi). 
\end{equation}
le morphisme canonique de topos \eqref{lptfr4b}.
Pour tout faisceau $F$ de $\tG'$, 
$\{i\mapsto \Psi_i^*(F)\}$ est naturellement une section de $\bHom_{I^\circ}(I^\circ,\fF_\varphi^\vee)$.
Elle définit donc un faisceau de $\Top(\fG_\varphi)$ \eqref{lptfr4a}.  On a alors un isomorphisme canonique fonctoriel \eqref{lptfr4d}
\begin{equation}
\Psi^*(F)\stackrel{\sim}{\rightarrow} \varpi^*(\{i\mapsto \Psi_i^*(F)\}).
\end{equation}

(i) D'après (\cite{sga4} IV (6.8.3)), pour tout faisceau $F$ de $\tG'$, on a un isomorphisme canonique
\begin{equation}
\Phi_*(F)_{\rho(\oy\rightsquigarrow \oz)}\stackrel{\sim}{\rightarrow} \underset{\underset{i\in I^\circ}{\longrightarrow}}{\lim}\ 
\Gamma((W'_i\rightarrow U'_i\leftarrow V'_i)^a,F). 
\end{equation}
On en déduit, en vertu de \ref{tfr42} et \ref{lptfr5}, un isomorphisme fonctoriel
\begin{equation}
\Phi_*(F)_{\rho(\oy\rightsquigarrow \oz)}\stackrel{\sim}{\rightarrow} \Gamma(\tuG',\Psi^*(F)).
\end{equation}

(ii) D'après (\cite{sga4} IV (6.8.4) et V 5.1), pour tout tout faisceau abélien $F$ de $\tG'$ et tout entier $q\geq 0$, on a un isomorphisme canonique
\begin{equation}
\rR^q\Phi_*(F)_{\rho(\oy\rightsquigarrow \oz)}\stackrel{\sim}{\rightarrow} \underset{\underset{i\in I^\circ}{\longrightarrow}}{\lim}\ 
\rH^q((W'_i\rightarrow U'_i\leftarrow V'_i)^a,F). 
\end{equation}
On en déduit, en vertu de \ref{tfr42} et \ref{lptfr6}, un isomorphisme fonctoriel
\begin{equation}
\rR^q\Phi_*(F)_{\rho(\oy\rightsquigarrow \oz)}\stackrel{\sim}{\rightarrow} \rH^q(\tuG',\Psi^*(F)).
\end{equation}

\begin{cor}\label{lptfr90}
Conservons les hypothèses et notations de \ref{lptfr9}, supposons, de plus, les morphismes $h_X$, $h_Y$ et $h_Z$ finis. 
Tout point géométrique de $\uX'$ se spécialise donc uniquement en un point de $X'_\ox$, la fibre de $h_X$ au-dessus de $\ox$ .  
On désigne par $Z'_\oz=\{\oz'_i; i\in I\}$ la fibre de $h_Z$ au-dessus de $\oz$, et pour tout $i\in I$, 
par $\{\oy'_{ij}; j\in J_i\}$ l'ensemble des points géométriques $\oy'$ de $Y'$ au-dessus de $\oy$ tels que $\uf'(\oy')$ se spécialise en 
$\ox'_i=g'(\oz'_i)$. Pour tous $i\in I$ et $j\in J_i$, on a donc un point canonique $(\oy'_{ij}\rightsquigarrow \oz'_i)$ 
de $Z'_\et\gtimes_{X'_\et}Y'_\et$. Alors, pour tout faisceau $F$ de $\tG'$, on a un isomorphisme canonique fonctoriel 
\begin{equation}
\Phi_*(F)_{\rho(\oy\rightsquigarrow \oz)}\stackrel{\sim}{\rightarrow}\prod_{i\in I,j\in J_i}F_{\rho'(\oy'_{ij}\rightsquigarrow \oz'_i)},
\end{equation}
où $\rho'\colon Z'_\et\gtimes_{X'_\et}Y'_\et\rightarrow \tG'$ est le morphisme canonique \eqref{tfr6b}.
\end{cor}

Pour tout $i\in I$, on désigne par $\uZ'_i$ le localisé strict de $Z'$ en $\oz'_i$, 
par $\uX'_i$ le localisé strict de $X'$ en $\ox'_i=g'(\oz'_i)$ et par $\ug'_i\colon \uZ'_i\rightarrow \uX'_i$ le morphisme induit par $g'$. 
Comme $h_Z$ est fini, $\uZ'$ est canoniquement isomorphe à la somme disjointe des schémas $(\uZ'_i)_{i\in I}$.  
De même, $h_X$ étant fini, pour tout $i\in I$, le morphisme canonique $\uX'_i\rightarrow \uX'$ est une immersion ouverte et fermée. 
De plus, le diagramme 
\begin{equation}
\xymatrix{
{\uZ'_i}\ar[r]^{\ug'_i}\ar[d]&{\uX'_i}\ar[d]\\
{\uZ}\ar[r]^{\ug}&{\uX}}
\end{equation}
est commutatif. 
On rappelle qu'on a un morphisme canonique $\uf'\colon \uY'\rightarrow \uX'$ qui relève $f'$ et $\uf$. 
Comme $h_Y$ est fini, $\uY'$ est canoniquement isomorphe à la somme disjointe des localisés stricts de $Y'$ 
en les points de $Y'_\oy$. 
Pour tout $j\in J_i$, on désigne par $\uY'_{ij}$ le localisé strict de $Y'$ en $\oy'_{ij}$, par 
$\uf'_{ij}\colon \uY'_{ij}\rightarrow \uX'_i$ le morphisme induit par $\uf'$, par 
$\tuG'_{ij}$ le topos de Faltings relatif associé au couple de morphismes  $(\uf'_{ij},\ug'_i)$
et par $\Psi_{ij}\colon \tuG'_{ij}\rightarrow \tG'$ le morphisme de fonctorialité déduit du diagramme commutatif
\begin{equation}
\xymatrix{
\uZ'_i\ar[r]^{\ug'_i}\ar[d]&\uX'_i\ar[d]&\uY'_{ij}\ar[l]_{\uf'_{ij}}\ar[d]\\
\uZ'\ar[r]^{\ug'}\ar[d]&\uX'\ar[d]&\uY'\ar[l]_{\uf'}\ar[d]\\
Z'\ar[r]^{g'}&X'&\ar[l]_{f'}Y'}
\end{equation}

On voit aussitôt que $\uY'_i=\uY'\times_{\uX'}\uX'_i$ est canoniquement isomorphe 
à la somme disjointe des schémas $(\uY'_{ij})_{j\in J_i}$. 
Par suite, en vertu de \ref{tfr44}, l'objet final $(\uZ'\rightarrow \uX'\leftarrow \uY')^a$ de $\tuG'$
est la somme disjointe des faisceaux $(\uZ'_i\rightarrow \uX'_i\leftarrow \uY'_{ij})^a$ pour $i\in I$ et $j\in J_i$.  
Il résulte alors de \ref{lptfr9}(i) que pour tout faisceau $F$ de $\tG'$, on a un isomorphisme canonique fonctoriel 
\begin{equation}
\Phi_*(F)_{\rho(\oy\rightsquigarrow \oz)}\stackrel{\sim}{\rightarrow}\prod_{i\in I,j\in J_i}\Gamma(\tuG'_{ij},\Psi^*_{ij}(F)).
\end{equation}
Pour tous $i\in I$ et $j\in J_i$, le topos $\tuG'_{ij}$ est local de centre $\rho'_{ij}(\oy'_{ij}\rightsquigarrow \oz'_i)$ d'après \ref{tfr280},
où $\rho'_{ij}\colon \uZ'_{i,\et}\gtimes_{\uX'_{i,\et}}\uY'_{ij,\et}\rightarrow \tG'_{ij}$ est le morphisme canonique. 
La proposition s'ensuit compte tenu de \eqref{tfr7e}.

\begin{cor}\label{lptfr91}
Considérons un diagramme commutatif de morphismes de schémas cohérents
\begin{equation}\label{lptfr91a}
\xymatrix{
Z'\ar[r]^{g'}\ar[d]_{h_Z}&X'\ar[d]_{h_X}&Y'\ar[l]_{f'}\ar[d]^{h_Y}\\
Z\ar[r]^g&X&\ar[l]_fY}
\end{equation}
tel que les flèches verticales soient des morphismes entiers, et notons $\tG$ (resp. $\tG'$) le topos de Faltings relatif associé au couple 
de morphismes $(f,g)$ (resp. $(f',g')$) et $\Phi\colon \tG'\rightarrow \tG$ le morphisme de fonctorialité induit par le diagramme \eqref{lptfr91a}. 
Alors, pour tout faisceau abélien $F$ de $\tG'$ et tout entier $q\geq 1$, on a $\rR^q\Phi_*(F)=0$.
\end{cor}

En effet, si les morphismes $h_X$, $h_Y$ et $h_Z$ sont finis, 
il résulte de \ref{lptfr20} et \ref{lptfr90} que le foncteur $F\mapsto \Phi_*(F)$ est exact sur la catégorie des faisceaux abéliens de $\tG'$, 
d'où la proposition dans ce cas. 

Dans le cas général, en vertu de \ref{lptfr20} et \ref{lptfr9}(ii), on se ramène à prouver que lorsque les schémas 
$X$, $Y$ et $Z$ sont locaux, pour tout faisceau abélien $F$ de $\tG'$ et tout entier $q\geq 1$, on a $\rH^q(\tG',F)=0$.  
On aura alors $X=\Spec(A)$, $Y=\Spec(B)$, $Z=\Spec(C)$, $X'=\Spec(A')$, $Y'=\Spec(B')$ et $Z'=\Spec(C')$, avec
$A'$ (resp. $B'$, resp. $C'$) une algèbre entière sur $A$ (resp. $B$, resp. $C$).
\'Ecrivant $A'$ (resp. $B'$, resp. $C'$) comme limite inductive filtrante de sous-$A$-algèbres 
(resp. sous-$B$-algèbres, resp. sous-$C$-algèbres) de type fini $(A_i)_{i\in I}$ (resp. $(B_i)_{i\in I}$, resp. $(C_i)_{i\in I}$), 
telles que pour tout $i\in I$, l'homomorphisme $A'\rightarrow B'$ (resp. $A'\rightarrow C'$) envoie $A'_i$ dans $B'_i$ (resp. $C'_i$), 
on voit qu'il existe une catégorie cofiltrante essentiellement petite $I$ et un foncteur 
\begin{equation}
\varphi\colon I\rightarrow \fM_{/(Z\rightarrow X\leftarrow Y)}, \ \ \ i\mapsto (Z'_i\rightarrow X'_i\leftarrow Y'_i),
\end{equation}
vérifiant les propriétés énumérées dans \ref{lptfr2}, tels que pour tout $i\in \ob(I)$, 
le $X$-schéma $X'_i$ (resp. $Y$-schéma $Y'_i$, resp. $Z$-schéma $Z'_i$) soit fini, 
et dont la limite projective est représentable par $(Z'\rightarrow X'\leftarrow Y')$. 
La proposition recherchée résulte alors du cas fini traité plus haut, compte tenu de \eqref{lptfr4f} et \ref{lptfr6}

\subsection{}\label{lptfr10}
Soient $f\colon Y\rightarrow X$, $g\colon Z\rightarrow X$, $h\colon Z'\rightarrow Z$ trois morphismes de schémas, $g'=g\circ h$, 
$(\oy\rightsquigarrow \oz)$ un point de $Z_\et\gtimes_{X_\et}Y_\et$ \eqref{tfr21}. 
Notons $i_\oy\colon \Ens\rightarrow Y_\et$ (resp. $i_\oz\colon \Ens\rightarrow Z_\et$) le morphismes de topos défini par le point 
géométrique $\oy$ de $Y$ (resp. $\oz$ de $Z$), $\varepsilon \colon Z'_{\oz,\et}\rightarrow \Ens$ la projection canonique (\cite{sga4} IV 4.3)
et $j_\oz\colon Z'_{\oz}\rightarrow Z'$ le morphisme canonique. 
Le point $(\oy\rightsquigarrow \oz)$ est défini par un $2$-morphisme $f_\et\circ i_\oy\rightarrow g_\et\circ i_\oz$. Composant avec $\varepsilon$
et utilisant l'isomorphisme canonique $i_\oz\varepsilon \stackrel{\sim}{\rightarrow} h_\et j_{\oz,\et}$, on obtient un $2$-morphisme 
\begin{equation}\label{lptfr10a}
f_\et i_\oy \varepsilon \rightarrow g'_\et j_{\oz,\et}.
\end{equation}
Le couple de morphismes $j_{\oz,\et}$ et  $i_\oy\circ \varepsilon$ et le $2$-morphisme ci-dessus définissent un morphisme 
\begin{equation}\label{lptfr10b}
\iota\colon Z'_{\oz,\et}\rightarrow Z'_\et\gtimes_{X_\et}Y_\et
\end{equation}
qui s'insère dans un diagramme non-commutatif
\begin{equation}\label{lptfr10c}
\xymatrix{
&{Z'_{\oz,\et}}\ar[d]^{\iota}\ar[ld]_{j_{\oz,\et}}\ar[r]^\varepsilon&{\Ens}\ar[d]^{i_\oy}\\
{Z'_\et}\ar[rd]_{g'_\et}&{Z'_\et\gtimes_{X_\et}Y_\et}\ar[l]_-(0.4){\rp_1}\ar[r]^-(0.4){\rp_2}&{Y_\et}\ar[ld]^{f_\et}\\
&{X_\et}&}
\end{equation}

On désigne par $\tG$ (resp. $\tG'$) le topos de Faltings relatif associé au couple de morphismes  
$(f,g)$ (resp. $(f,g')$), par $\Phi\colon \tG'\rightarrow \tG$ le morphisme de fonctorialité induit par $h$, 
par $\pi'\colon \tG'\rightarrow Z'_\et$ le morphisme canonique \eqref{tfr4c}, par 
\begin{eqnarray}
\rho\colon Z_\et\gtimes_{X_\et}Y_\et&\rightarrow&\tG,\label{lptfr10d}\\
\rho'\colon Z'_\et\gtimes_{X_\et}Y_\et&\rightarrow&\tG',\label{lptfr10e}
\end{eqnarray}
les morphismes canoniques \eqref{tfr6b} et par 
\begin{equation}\label{lptfr10f}
\jmath\colon Z'_{\oz,\et}\rightarrow \tG'
\end{equation}
le morphisme composé $\rho' \iota$ \eqref{lptfr10b}. 
Il résulte aussitôt de la définition de $\iota$ \eqref{lptfr10c} et de \eqref{tfr6c} que le diagramme 
\begin{equation}\label{lptfr10g}
\xymatrix{
&\tG'\ar[d]^{\pi'}\\
{Z'_{\oz,\et}}\ar[r]_{j_{\oz,\et}}\ar[ru]^{\jmath}&{Z'_\et}}
\end{equation}
est commutatif à isomorphisme canonique près. On en déduit un morphisme de changement de base 
\begin{equation}\label{lptfr10h}
j_{\oz,\et}^*\pi'_* \rightarrow \jmath^*
\end{equation}
adjoint du morphisme composé $\pi'_*\rightarrow \pi'_*\jmath_*\jmath^*\rightarrow j_{\oz,\et*}\jmath^*$, où la première flèche est induite par 
le morphisme d'adjonction
et la seconde flèche par l'isomorphisme canonique $\pi'_*\jmath_*\stackrel{\sim}{\rightarrow} j_{\oz,\et*}$ \eqref{lptfr10g}.

Il résulte de la propriété universelle des produits orientés de topos que le diagramme de morphismes de topos 
\begin{equation}\label{lptfr10i}
\xymatrix{
{Z'_{\oz,\et}}\ar[r]^-(0.5){\iota}\ar[d]_\varepsilon&{Z'_\et\gtimes_{X_\et}Y_\et}\ar[d]^{\lh_\et}\\
\Ens\ar[r]^-(0.5){(\oy\rightsquigarrow \oz)}&{Z_\et\gtimes_{X_\et}Y_\et}}
\end{equation}
où $\lh_\et$ est le morphisme de fonctorialité induit par $h$, est commutatif à isomorphisme canonique près. 
Par ailleurs, le diagramme 
\begin{equation}\label{lptfr10j}
\xymatrix{
{Z'_\et\gtimes_{X_\et}Y_\et}\ar[d]_{\lh_\et}\ar[r]^-(0.5){\rho'}&\tG'\ar[d]^\Phi\\
{Z_\et\gtimes_{X_\et}Y_\et}\ar[r]^-(0.5){\rho}&\tG}
\end{equation}
est commutatif à isomorphisme canonique près. On en déduit que le diagramme 
\begin{equation}\label{lptfr10k}
\xymatrix{
{Z'_{\oz,\et}}\ar[r]^-(0.5){\jmath}\ar[d]_\varepsilon&{\tG'}\ar[d]^{\Phi}\\
\Ens\ar[r]^-(0.5){\rho(\oy\rightsquigarrow \oz)}&{\tG}}
\end{equation}
est commutatif à isomorphisme canonique près.

\begin{rema}\label{lptfr14}
Les formations des morphismes $\iota$ \eqref{lptfr10b} et $\jmath$ \eqref{lptfr10f} sont fonctorielles 
dans un sens évident que nous laissons au lecteur le soin d'expliciter. 
\end{rema}

\begin{prop}\label{lptfr13}
Les hypothèses étant celles de \ref{lptfr10}, supposons de plus les schémas 
$X$ et $Y$ strictement locaux de points fermés $g(\oz)$ et $\oy$, respectivement. Alors, 
\begin{itemize}
\item[{\rm (i)}] Pour tout faisceau $F$ de $\tG'$, le morphisme de changement de base \eqref{lptfr10h}
\begin{equation}\label{lptfr13a}
j_{\oz,\et}^*(\pi'_*(F))\rightarrow \jmath^*(F).
\end{equation}
est un isomorphisme. 
\item[{\rm (ii)}] Pour tout faisceau abélien $F$ de $\tG'$ et tout entier $q\geq 1$, on a $j_{\oz,\et}^*(\rR^q\pi'_*(F))=0$.
\end{itemize}
\end{prop}

En effet, soient $\oz'$ un point géométrique de $Z'$ au-dessus de $\oz$, $\uZ'$ le localisé strict de $Z'$ en $\oz'$, $\ug'\colon \uZ'\rightarrow X$ 
le morphisme induit par $g'$, $\tuG'$ le topos de Faltings relatif associé au couple de morphismes $(f,\ug')$,
$\Psi\colon \tuG'\rightarrow \tG'$ le morphisme de fonctorialité. 

(i) En vertu de \ref{lptfr12}(i), on a un isomorphisme canonique fonctoriel
\begin{equation}\label{lptfr13c}
\Gamma(\tuG',\Psi^*(F))\stackrel{\sim}{\rightarrow}\pi'_*(F)_{\oz'}.
\end{equation}
La fibre du morphisme \eqref{lptfr13a} en $\oz'$ s'identifie à un morphisme
\begin{equation}\label{lptfr13b}
\pi'_*(F)_{\oz'}\rightarrow F_{\rho'(\oy\rightsquigarrow \oz')}.
\end{equation}
On vérifie aussitôt que le diagramme de morphismes de topos 
\begin{equation}
\xymatrix{
{\uZ'_\et\gtimes_{X_\et}Y_\et}\ar[d]\ar[r]^-(0.5){\urho'}&\tuG'\ar[d]^\Psi\\
{Z'_\et\gtimes_{X_\et}Y_\et}\ar[r]^-(0.5){\rho'}&\tG'}
\end{equation}
où $\urho'$ est le morphisme canonique \eqref{tfr6b} et la flèche non libellée est le morphisme de fonctorialité, est commutatif à isomorphisme canonique près. 
Considérant $(\oy\rightsquigarrow \oz')$ comme un point de $\uZ'_\et\gtimes_{X_\et}Y_\et$,
le composé de \eqref{lptfr13c} et \eqref{lptfr13b} s'identifie au morphisme d'évaluation canonique 
\begin{equation}
\Gamma(\tuG',\Psi^*(F))\rightarrow\Psi^*(F)_{\urho'(\oy\rightsquigarrow \oz')}.
\end{equation}
Celui-ci est un isomorphisme puisque le topos $\tuG'$ est local de centre $\urho'(\oy\rightsquigarrow \oz')$ d'après \ref{tfr280}. 
Par suite, \eqref{lptfr13b} est un isomorphisme, d'où la proposition. 

(ii) En vertu de \ref{lptfr12}(ii), on a un isomorphisme canonique fonctoriel
\begin{equation}\label{lptfr13d}
\rR^q\pi'_*(F)_{\oz'}\stackrel{\sim}{\rightarrow}\rH^q(\tuG',\Psi^*(F)).
\end{equation}
Le terme de droite est nul puisque le topos $\tuG'$ est local d'après \ref{tfr280}, d'où la proposition.

\begin{prop}\label{lptfr16}
Les hypothèses étant celles de \ref{lptfr10}, supposons de plus le morphisme $h$ propre.
Alors,  
\begin{itemize}
\item[{\rm (i)}] Pour tout faisceau $F$ de $\tG'$, le morphisme de changement de base relativement au diagramme \eqref{lptfr10k}
\begin{equation}\label{lptfr16a}
\Phi_*(F)_{\rho(\oy\rightsquigarrow \oz)}\rightarrow \Gamma(Z'_{\oz,\et},\jmath^*(F))
\end{equation}
est un isomorphisme. 
\item[{\rm (ii)}] Pour tout faisceau abélien de torsion $F$ de $\tG'$ et tout entier $q\geq 0$, 
le morphisme de changement de base relativement au diagramme \eqref{lptfr10k}
\begin{equation}\label{lptfr16b}
\rR^q\Phi_*(F)_{\rho(\oy\rightsquigarrow \oz)}\rightarrow \rH^q(Z'_{\oz,\et},\jmath^*(F))
\end{equation}
est un isomorphisme. 
\end{itemize}
\end{prop}

Nous démontrons seulement la seconde proposition; la preuve de la première étant similaire et plus simple, est laissée au lecteur. 
Notons $\uZ$ le localisé strict de $Z$ en $\oz$, 
$\uX$ le localisé strict de $X$ en $g(\oz)$, $\uY$ le localisé strict de $Y$ en $\oy$, $\ug\colon \uZ\rightarrow \uX$ le morphisme induit par $g$,
et $\uf\colon \uY\rightarrow \uX$ le morphisme induit par $f$ et le point $(\oy\rightsquigarrow \oz)$ (\cite{sga4} VIII 7.3). 
Posons $\uZ'=Z'\times_Z\uZ$ et notons $\uh\colon \uZ'\rightarrow \uZ$ la projection canonique 
et $\tuG$ (resp. $\tuG'$) le topos de Faltings relatif associé au couple de morphismes $(\uf,\ug)$ (resp. $(\uf,\ug\circ \uh)$).
Le diagramme des morphismes de fonctorialité
\begin{equation}\label{lptfr16c}
\xymatrix{
\tuG'\ar[r]^{\Psi'}\ar[d]_{\uPhi}&\tG'\ar[d]^{\Phi}\\
\tuG\ar[r]^\Psi&\tG}
\end{equation}
est commutatif à isomorphisme canonique près. 
On peut naturellement considérer  $(\oy\rightsquigarrow \oz)$ comme un point de $\uZ_\et\gtimes_{\uX_\et}\uY_\et$; on note
\begin{equation}
\ujmath\colon Z'_{\oz,\et}\rightarrow \tuG'
\end{equation}
le morphisme associé défini dans \eqref{lptfr10f}. On vérifie aussitôt \eqref{lptfr14} qu'on a un isomorphisme canonique 
\begin{equation}
\jmath\stackrel{\sim}{\rightarrow}\Psi' \ujmath.
\end{equation}

En vertu de \ref{lptfr9}(ii), pour tout faisceau abélien $F$ de $\tG'$ et tout entier $q\geq 0$, 
le morphisme de changement de base relativement au diagramme \eqref{lptfr16c} induit un isomorphisme
\begin{equation}
\rR^q\Phi_*(F)_{\rho(\oy\rightsquigarrow \oz)}\stackrel{\sim}{\rightarrow} \rR^q\uPhi_*(\Psi'^*(F))_{\urho(\oy\rightsquigarrow \oz)},
\end{equation}
où $\urho\colon \uZ_\et\gtimes_{\uX_\et}\uY_\et\rightarrow \tuG$ est le morphisme canonique \eqref{tfr6b}.
Compte tenu de (\cite{egr1} 1.2.4(ii)), on peut donc se borner au cas où les schémas $X$, $Y$ et $Z$ sont strictement locaux 
de points fermés $g(\oz)$, $\oy$ et $\oz$ respectivement. 
D'après \ref{lptfr9}(ii), il s'agit alors de montrer que pour tout faisceau abélien de torsion $F$ de $\tG'$ et tout entier $q\geq 0$, 
le morphisme d'image inverse par $\jmath$,
\begin{equation}\label{lptfr16d}
\rH^q(\tG',F) \rightarrow \rH^q(Z'_{\oz,\et},\jmath^*(F)),
\end{equation}
est un isomorphisme. 

Considérons la suite spectrale de Cartan-Leray 
\begin{equation}\label{lptfr16e}
\rE_2^{a,b}=\rH^a(Z'_\et,\rR^b\pi'_*(F))\Rightarrow \rH^{a+b}(\tG',F).
\end{equation}
Pour tous $a\geq 0$ et $b\geq 1$, on a $\rE_2^{a,b}=0$ en vertu de \ref{lptfr13} et (\cite{sga4} XII 5.5).
On en déduit pour tout entier $q\geq 0$, un isomorphisme canonique 
\begin{equation}\label{lptfr16f}
\rH^q(Z'_{\et},\pi'_*(F))\stackrel{\sim}{\rightarrow} \rH^q(\tG',F).
\end{equation}

Considérons le diagramme commutatif
\begin{equation}\label{lptfr16g}
\xymatrix{
{\rH^q(Z'_{\et},\pi'_*(F))}\ar[r]^-(0.5){j_{\oz,\et}^*}\ar[d]_{\pi'^*}&{\rH^q(Z'_{\oz,\et},j_{\oz,\et}^*(\pi'_*(F)))}\ar[d]^w\\
{\rH^q(\tG',\pi'^*(\pi'_*(F)))}\ar[r]^-(0.5){\jmath^*}\ar[d]_u&{\rH^q(Z'_{\oz,\et},\jmath^*(\pi'^*(\pi'_*(F))))}\ar[d]^v\\
{\rH^q(\tG',F)}\ar[r]^-(0.5){\jmath^*}&{\rH^q(Z'_{\oz,\et},\jmath^*(F))}}
\end{equation}
où $u$ et $v$ sont définis par adjonction et $w$ est induit par l'isomorphisme canonique \eqref{lptfr10g}.
Le composé $u\circ \pi'^*$ est l'isomorphisme \eqref{lptfr16f}. 
D'après (\cite{sga4} XVII 2.1.3), $v\circ w$ est induit par le morphisme de changement de base \eqref{lptfr10h};
il est donc bijectif d'après \ref{lptfr13}. Le morphisme d'image inverse $j_{\oz,\et}^*$ est un isomorphisme, en vertu de (\cite{sga4} XII 5.5).
Il s'ensuit que le morphisme \eqref{lptfr16d}, représenté par la flèche horizontale inférieure, est un isomorphisme, d'où la proposition. 

\begin{prop}\label{lptfr15}
Soient $f\colon Y\rightarrow X$, $g\colon Z\rightarrow X$, $h\colon Z'\rightarrow Z$ trois morphismes de schémas tels que $h$ soit propre, 
$\tG$ (resp. $\tG'$) le topos de Faltings relatif associé au couple de morphismes  $(f,g)$ (resp. $(f,g\circ h)$). 
Considérons le diagramme commutatif à isomorphisme canonique près \eqref{tfr7c},
\begin{equation}\label{lptfr15a}
\xymatrix{
\tG'\ar[r]^{\pi'}\ar[d]_{\Phi}&{Z'_\et}\ar[d]^{h_\et}\\
\tG\ar[r]^{\pi}&{Z_\et}}
\end{equation}
où $\Phi$ est le morphisme de fonctorialité défini par $h$, et $\pi$ et $\pi'$ sont les morphismes canoniques \eqref{tfr4c}.
Alors, 
\begin{itemize}
\item[{\rm (i)}] Pour tout faisceau $F$ de $Z'_\et$, le morphisme de changement de base relativement au diagramme \eqref{lptfr15a}
\begin{equation}\label{lptfr15b}
\pi^*(h_{\et*}(F))\rightarrow \Phi_*(\pi'^*(F))
\end{equation}  
est un isomorphisme. 
\item[{\rm (ii)}] Pour tout faisceau abélien de torsion $F$ de $Z'_\et$
et tout entier $q\geq 0$, le morphisme de changement de base relativement au diagramme \eqref{lptfr15a}
\begin{equation}\label{lptfr15c}
\pi^*(\rR^qh_{\et*}(F))\rightarrow \rR^q\Phi_*(\pi'^*(F))
\end{equation}  
est un isomorphisme. 
\end{itemize}
\end{prop}

Nous démontrons seulement la seconde proposition; la preuve de la première étant similaire et plus simple, est laissée au lecteur. 
Notons $\rho\colon Z_\et\gtimes_{X_\et}Y_\et\rightarrow\tG$ le morphisme canonique \eqref{tfr6b}. Compte tenu de \ref{lptfr20},
il suffit de montrer que pour tout point $(\oy\rightsquigarrow \oz)$ de $Z_\et\gtimes_{X_\et}Y_\et$ \eqref{tfr21}, la fibre 
\begin{equation}
\rR^qh_{\et*}(F)_\oz\rightarrow \rR^q\Phi_*(\pi'^*(F))_{\rho(\oy\rightsquigarrow \oz)}
\end{equation}
de \eqref{lptfr15c} en $\rho(\oy\rightsquigarrow \oz)$  est un isomorphisme. 
Reprenons les notations de \ref{lptfr10} et considérons le diagramme 
\begin{equation}
\xymatrix{
{Z'_{\oz,\et}}\ar[r]^-(0.5){\jmath}\ar[d]_\varepsilon&\tG'\ar[r]^{\pi'}\ar[d]^{\Phi}&{Z'_\et}\ar[d]^{h_\et}\\
\Ens\ar[r]^-(0.5){\rho(\oy\rightsquigarrow \oz)}&\tG\ar[r]^{\pi}&{Z_\et}}
\end{equation}
dont les carrés sont commutatifs à isomorphismes canoniques près. 
Comme $\pi(\rho(\oy\rightsquigarrow \oz))$ est le point de $Z_\et$ défini par $\oz$, 
la proposition résulte de \ref{lptfr16}, (\cite{sga4} XII 5.2) et (\cite{egr1} 1.2.4(ii)).

\section{Topos de Faltings relatif annelé}\label{tfra}

\subsection{}\label{tfra1}
Dans cette section, $f\colon Y\rightarrow X$ et $g\colon Z\rightarrow X$ désignent deux morphismes de schémas.
On note $G$ (resp. $\tG$) le site (resp. topos) de Faltings relatif associé aux morphismes $(f,g)$ et on reprend les notations de la section \ref{tfr}.
On suppose, de plus, que la projection canonique $Z\times_XY\rightarrow Z$ est le composé de deux morphismes 
\begin{equation}\label{tfra1a}
\xymatrix{
 Z\times_XY\ar[r]^-(0.5)j&\oZ\ar[r]^-(0.5)\hbar&Z},
 \end{equation} 
où $\oZ$ est un schéma normal et localement irréductible et $j$ est une immersion ouverte quasi-compacte.  
Pour tout objet $(W\rightarrow U\leftarrow V)$ de $G$, on pose $\oW=W\times_{Z}\oZ$
et on note $\oW^V$ la fermeture intégrale de $\oW$ dans $W\times_UV$. 
\begin{equation}\label{tfra1b}
\xymatrix{
&\oW\ar[ld]&{\oW^V}\ar[l]\\
W\ar[d]&{W\times_{X}{Y}}\ar[d]\ar[l]\ar[u]&{W\times_UV}\ar[l]\ar[d]\ar[u]\\
U & {U\times_XY}\ar[l] & V\ar[l]}
\end{equation}
Pour tout morphisme $(W'\rightarrow U'\leftarrow V')\rightarrow (W\rightarrow U\leftarrow V)$ de $G$, on a 
un morphisme canonique $\oW'^{V'}\rightarrow \oW^V$ qui s'insère dans un diagramme commutatif 
\begin{equation}\label{tfra1c}
\xymatrix{
{W'\times_{U'}V'}\ar[r]\ar[d]&{\oW'^{V'}}\ar[r]\ar[d]&\oW'\ar[r]\ar[d]&W'\ar[d]\\
{W\times_UV}\ar[r]&{\oW^V}\ar[r]&\oW\ar[r]&W}
\end{equation}

On désigne par $\ocB$ le préfaisceau défini pour tout objet $(W\rightarrow U\leftarrow V)$ de $G$ par
\begin{equation}\label{tfra1d}
\ocB(W\rightarrow U\leftarrow V)=\Gamma(\oW^V,\co_{\oW^V}).
\end{equation}
Nous montrerons dans \ref{tfra3} ci-dessous que c'est un faisceau pour la topologie co-évanescente de $G$. On a un homomorphisme canonique
\begin{equation}\label{tfra1e}
\hbar_*(\co_\oZ)\rightarrow \pi_*(\ocB),
\end{equation}
défini pour tout $W\in \ob(\Et_{/Z})$ par l'homomorphisme canonique
\begin{equation}\label{tfra1f}
\Gamma(\oW,\co_\oW)\rightarrow \Gamma(\oW^Y,\co_{\oW^Y}),
\end{equation}
où l'on a considéré l'objet $\pi^+(W)=(W\rightarrow X\leftarrow Y)$ de $G$ \eqref{tfr4a}.

\begin{rema}\label{tfra2}
Pour tout objet $(W\rightarrow U\leftarrow V)$ de $G$, le schéma $\oW^V$ est normal et localement irréductible
et la projection canonique $(W\times_XY)\times_{\oW}\oW^V \rightarrow W\times_UV$ est un isomorphisme. 
En effet, $\oW$ et $W\times_UV$ sont normaux et localement irréductibles d'après (\cite{agt} III.3.3).
Soit $W_0$ un ouvert de $\oW$ n'ayant qu'un nombre fini de composantes irréductibles. 
Alors $W_0\times_{\oW}\oW^V$ est la somme finie des fermetures intégrales de $W_0$ dans 
les points génériques de $W\times_UV$ qui sont au-dessus de $W_0$, dont chacune est un schéma intègre et normal 
en vertu de (\cite{ega2} 6.3.7), d'où la première assertion. La seconde assertion résulte du fait $V$ est entier sur $U\times_XY$.
\end{rema}

\begin{rema}\label{tfra8}
Supposons que $Z=X$ et $g=\id_X$ et posons $\oX=\oZ$ de sorte que $f$ se factorise en 
\begin{equation}\label{tfra8a}
\xymatrix{
Y\ar[r]^j&\oX\ar[r]^\hbar&X}. 
\end{equation}
On désigne par $E$ (resp. $\tE$) le site (resp. topos) de Faltings associé à $f$ \eqref{tf1}, 
par $\ocB_E$ l'anneau de $\tE$ défini dans (\cite{agt} III.8.10) par la factorisation $f=\hbar\circ j$ (noté $\ocB$ dans {\em loc. cit.})
et par $\iota \colon \tG\rightarrow \tE$ et $\jmath\colon \tE\rightarrow \tG$ les équivalences de topos \eqref{tfr11a} et  \eqref{tfr11b}.
Il résulte aussitôt des définitions que l'on a $\iota_*(\ocB)=\ocB_E$ et $\jmath_*(\ocB_E)=\ocB$.   
 \end{rema}

\begin{prop}\label{tfra3}
Le préfaisceau $\ocB$ sur $G$ est un faisceau pour la topologie co-évanescente. 
\end{prop}

On désigne par $E'$ (resp. $\tE'$) le site (resp. topos) de Faltings associé au morphisme $Z\times_XY\rightarrow Z$ \eqref{tf1}
et par $\tau\colon \tE'\rightarrow \tG$ le morphisme de fonctorialité \eqref{tfr13c}. Pour tout objet $(V'\rightarrow U')$ de $E'$,
on pose $\oU'=U'\times_Z\oZ$ et on note $\oU'^{V'}$ la fermeture intégrale de $\oU'$ dans $V'$. D'après (\cite{agt} III.8.16), 
la correspondance qui à tout $(V'\rightarrow U')\in \ob(E')$, associe l'anneau 
\begin{equation}
\ocB'(V'\rightarrow U')=\Gamma(\oU'^{V'},\co_{\oU'^{V'}}),
\end{equation}
est un faisceau sur $E'$ pour la topologie co-évanescente. Par ailleurs, compte tenu de \ref{tfr7d} et \ref{tfr9b}, on a un isomorphisme
\begin{equation}
\tau_*(\ocB')\stackrel{\sim}{\rightarrow}\ocB,
\end{equation}
d'où la proposition.

\begin{prop}\label{tfra4}
Supposons les hypothèses suivantes satisfaites: 
\begin{itemize}
\item[{\rm (a)}] les schémas $X$ et $Z$ sont strictement locaux de points fermés $x$ et $z$ respectivement et le morphisme $g$ est local, {\em i.e.}, $g(z)=x$;
\item[{\rm (b)}] l'ensemble des composantes connexes de $Y$ est localement fini et le morphisme $f\colon Y\rightarrow X$ est cohérent;
\item[{\rm (c)}] le schéma $\oZ$ est strictement local;
\item[{\rm (d)}] pour tout revêtement étale $V\rightarrow Y$ tel que $V$ soit connexe, le schéma $Z\times_XV$ est connexe. 
\end{itemize}
Soient $\oy$ un point géométrique de $Y$, $(\oy\rightsquigarrow z)$ le point de $Z_\et\gtimes_{X_\et}Y_\et$ défini par l'unique flèche de spécialisation 
de $f(\oy)$ dans $x$ \eqref{tfr21},  $\rho(\oy\rightsquigarrow z)$ son image par le morphisme $\rho\colon Z_\et\gtimes_{X_\et}Y_\et\rightarrow \tG$ \eqref{tfr6b}. Alors~:
\begin{itemize}
\item[{\rm (i)}] La fibre $\ocB_{\rho(\oy \rightsquigarrow z)}$ de $\ocB$ en $\rho(\oy\rightsquigarrow z)$ est un anneau normal et strictement local. 
\item[{\rm (ii)}] On a un isomorphisme canonique 
\begin{equation}\label{tfra4a}
(\hbar_*(\co_{\oZ}))_{z}\stackrel{\sim}{\rightarrow}\Gamma(\oZ,\co_{\oZ}).
\end{equation}
\item[{\rm (iii)}]  L'homomorphisme 
\begin{equation}\label{tfra4b}
(\hbar_*(\co_{\oZ}))_{z}\rightarrow \ocB_{\rho(\oy \rightsquigarrow z)}
\end{equation} 
induit par l'homomorphisme canonique $\pi^*(\hbar_*(\co_{\oZ}))\rightarrow \ocB$ \eqref{tfra1e} est injectif et local.
\end{itemize}
\end{prop}

(i) D'après (b), les composantes connexes de $Y$ sont ouvertes dans $Y$. 
Soient $Y^\star$ la composante connexe de $Y$ contenant $\oy$, $(V_i)_{i\in I}$ 
le revêtement universel normalisé de $Y^\star$ au point $\oy$ \eqref{notconv11}.
Il résulte de \ref{tfr250}, \ref{tfr27} et (\cite{agt} VI.9.9) qu'on a un isomorphisme canonique
\begin{equation}\label{tfra4c}
\ocB_{\rho(\oy \rightsquigarrow z)}\stackrel{\sim}{\rightarrow}
\underset{\underset{i\in I}{\longrightarrow}}{\lim}\ \ocB(Z\rightarrow X\leftarrow V_i).
\end{equation}
Pour chaque $i\in I$, notons $\oZ^{V_i}$ la fermeture intégrale de $\oZ$ dans $Z\times_X V_i$.   
Le schéma $Z\times_X V_i$ est normal d'après (\cite{agt} III.3.3), et étant connexe par l'hypothèse (d), il est intègre. 
Par suite, le schéma $\oZ^{V_i}$ est normal, intègre et entier sur $\oZ$. Il est donc strictement local d'après (\cite{agt} III.3.5). 
Par ailleurs, pour tous $(i,j)\in I^2$ avec $j\geq i$, le morphisme de transition $\oZ^{V_j}\rightarrow \oZ^{V_i}$ est entier et dominant. 
Par suite, l'homomorphisme de transition
\begin{equation}\label{tfra4d}
\ocB(Z\rightarrow X\leftarrow V_i)\rightarrow \ocB(Z\rightarrow X\leftarrow V_j)
\end{equation}
est local. On en déduit que l'anneau $\ocB_{\rho(\oy \rightsquigarrow z)}$ est local, normal et hensélien 
(\cite{ega1n} 0.6.5.12(ii) et \cite{raynaud1} I § 3 prop.~1). 
Comme l'homomorphisme $\Gamma(\oZ,\co_{\oZ})\rightarrow \ocB_{\rho(\oy \rightsquigarrow z)}$ est entier et donc local, 
le corps résiduel de $\ocB_{\rho(\oy \rightsquigarrow z)}$ est une extension algébrique de celui de 
$\Gamma(\oZ,\co_{\oZ})$. Il est donc séparablement clos. 

(ii) Cela résulte aussitôt du fait que $Z$ est strictement local. 

(iii) Rappelons qu'on a un isomorphisme canonique \eqref{tfr21a} 
\begin{equation}\label{tfra4e}
(\pi^{-1}(\hbar_*(\co_{\oZ})))_{\rho(\oy \rightsquigarrow z)}\stackrel{\sim}{\rightarrow} \hbar_*(\co_{\oZ})_{z}.
\end{equation}
Compte tenu de (ii) et \eqref{tfra4c}, la fibre de l'homomorphisme canonique 
$\pi^{-1}(\hbar_*(\co_{\oZ}))\rightarrow \ocB$ en $\rho(\oy \rightsquigarrow z)$ s'identifie à l'homomorphisme canonique
\begin{equation}\label{tfra4f}
\Gamma(\oZ,\co_{\oZ})\rightarrow \ocB_{\rho(\oy \rightsquigarrow z)}=\underset{\underset{i\in I}{\longrightarrow}}{\lim}\ 
\ocB(Z\rightarrow X\leftarrow V_i), 
\end{equation} 
qui est clairement injectif et entier et donc local.

\subsection{}\label{tfra5}
Considérons un diagramme commutatif de morphismes de schémas 
\begin{equation}\label{tfra5a}
\xymatrix{
{Z'\times_{X'}Y'}\ar[r]^-(0.4){j'}\ar[d]_{w\times_uv}&\oZ'\ar[r]^{\hbar'}\ar[d]_{\ow}&Z'\ar[r]^{g'}\ar[d]^w&X'\ar[d]^u&Y'\ar[l]_{f'}\ar[d]^v\\
{Z\times_XY}\ar[r]^-(0.4)j&\oZ\ar[r]^\hbar&Z\ar[r]^g&X&\ar[l]_fY}
\end{equation}
tel que $\oZ'$ soit un schéma normal et localement irréductible, que $j'$ soit une immersion ouverte quasi-compacte et
que $\hbar'\circ j'$ soit la projection canonique $Z'\times_{X'}Y'\rightarrow Z'$. 
On désigne par $G'$ (resp. $\tG'$) le site (resp. topos) de Faltings relatif associé à $(f',g')$, par $\ocB'$ l'anneau de $\tG'$ associé à  
la factorisation $\hbar'\circ j'$ de la projection canonique \eqref{tfra1} par $\Phi^+\colon G\rightarrow G'$ le foncteur de changement de base \eqref{tfr7d} 
et par  
\begin{equation}\label{tfra5b}
\Phi\colon \tG'\rightarrow \tG.
\end{equation}
le morphisme de fonctorialité \eqref{tfr7b}. Pour tout objet $(W'\rightarrow U'\leftarrow V')$ de $G'$, on pose $\oW'=W'\times_{Z'}\oZ'$
et on note $\oW'^{V'}$ la fermeture intégrale de $\oW'$ dans $W'\times_{U'}V'$ de sorte que 
\begin{equation}\label{tfra5c}
\ocB'((W'\rightarrow U'\leftarrow V'))=\Gamma(\oW'^{V'},\co_{\oW'^{V'}}).
\end{equation}
Pour tout objet $(W\rightarrow U\leftarrow V)$ de $G$, posant $(W'\rightarrow U'\leftarrow V')=\Phi^+(W\rightarrow U\leftarrow V)$, on a 
alors un diagramme commutatif 
\begin{equation}\label{tfra5d}
\xymatrix{
{W'\times_{U'}V'}\ar[r]\ar[d]&{W'\times_{X'}Y'}\ar[r]\ar[d]&{\oW'}\ar[d]\\
{W\times_{U}V}\ar[r]&{W\times_{X}Y}\ar[r]&{\oW}}
\end{equation}
On en déduit un morphisme 
\begin{equation}\label{tfra5e}
\oW'^{V'}\rightarrow \oW^V,
\end{equation}
et par suite un homomorphisme d'anneaux de $\tG$
\begin{equation}\label{tfra5f}
\ocB\rightarrow \Phi_*(\ocB').
\end{equation}
Nous considérons dans la suite $\Phi$ comme un morphisme de topos annelés (respectivement par $\ocB'$ et $\ocB$). 
Nous utilisons pour les modules la notation $\Phi^{-1}$ pour désigner l'image
inverse au sens des faisceaux abéliens et nous réservons la notation 
$\Phi^*$ pour l'image inverse au sens des modules.

\begin{lem}\label{tfra6}
Les hypothèses étant celles de \ref{tfra5}, supposons, de plus, que les morphismes $u\colon X'\rightarrow X$ et
$w\colon Z'\rightarrow Z$ soient étales et que les deux diagrammes
\begin{equation}
\xymatrix{
\oZ'\ar[r]^{\hbar'}\ar[d]_\ow&Z'\ar[d]^w\\
\oZ\ar[r]^{\hbar}&Z}
\ \ \
\xymatrix{
X'\ar[d]_u&Y'\ar[d]^v\ar[l]_{f'}\\
X&Y\ar[l]_{f}}
\end{equation}
soient cartésiens, 
de sorte que $(Z'\rightarrow X'\leftarrow Y')$ est un objet de $G$. Alors~:
\begin{itemize} 
\item[{\rm (i)}] Le morphisme 
\begin{equation}\label{tfra6a}
\Phi^{-1}(\ocB)\rightarrow \ocB'
\end{equation}
adjoint du morphisme \eqref{tfra5f} est un isomorphisme. 
\item[{\rm (ii)}] Le morphisme de topos annelés $\Phi\colon (\tG',\ocB')\rightarrow (\tG,\ocB)$ 
s'identifie au morphisme de localisation de $(\tG,\ocB)$ en $(Z'\rightarrow X'\leftarrow Y')^a$.
\end{itemize}
\end{lem}

En effet, le morphisme de topos $\Phi\colon \tG'\rightarrow \tG$ s'identifie au morphisme de localisation 
de $\tG$ en $(Z'\rightarrow X'\leftarrow Y')^a$ en vertu de \ref{tfr42}. 
Il suffit donc de montrer la première proposition. Le foncteur $\Phi^*\colon \tG\rightarrow \tG'$ s'identifie au 
foncteur de restriction par le foncteur canonique $G'\rightarrow G$. 
Il résulte des hypothèses que le diagramme 
\begin{equation}
\xymatrix{
{Z'\times_{X'}Y'}\ar[r]\ar[d]_{w\times_uv}&Z'\ar[d]^w\\
{Z\times_XY}\ar[r]&Z}
\end{equation}
est cartésien.  
Pour tout $(W\rightarrow U\leftarrow V)\in \ob(G')$, on a des isomorphismes canoniques 
\begin{eqnarray}
W\times_{Z'}\oZ'&\stackrel{\sim}{\rightarrow}&W\times_Z\oZ=\oW,\\
W\times_{X'}Y'&\stackrel{\sim}{\rightarrow}&W\times_XY. 
\end{eqnarray}
On en déduit un isomorphisme (fonctoriel en $(W\rightarrow U\leftarrow V)$)
\begin{equation}\label{tfra6b}
\Phi^{-1}(\ocB)((W\rightarrow U\leftarrow V))\stackrel{\sim}{\rightarrow}\ocB'((W\rightarrow U\leftarrow V)).
\end{equation}
Il reste à montrer que celui-ci est adjoint du morphisme \eqref{tfra5f}. 
Posons $(W'\rightarrow U'\leftarrow V')=\Phi^+(W\rightarrow U\leftarrow V)$ \eqref{tfr7d}.
Le morphisme $(W\rightarrow U\leftarrow V)\rightarrow (Z'\rightarrow X'\leftarrow Y')$ induit une section 
$(W\rightarrow U\leftarrow V)\rightarrow (W'\rightarrow U'\leftarrow V')$ de la projection canonique 
\begin{equation}
(W'\rightarrow U'\leftarrow V')\rightarrow (W\rightarrow U\leftarrow V).
\end{equation}
On en déduit des morphismes $\oW^V\rightarrow \oW'^{V'}\rightarrow \oW^V$ dont le composé est l'identité de $\oW^V$. 
Le composé 
\begin{eqnarray}
\lefteqn{\Phi^{-1}(\ocB)((W\rightarrow U\leftarrow V))\rightarrow \Phi^{-1}(\Phi_*(\ocB'))((W\rightarrow U\leftarrow V))}\\
&=& \ocB'((W'\rightarrow U'\leftarrow V'))\rightarrow \ocB'((W\rightarrow U\leftarrow V))\nonumber
\end{eqnarray}
où la première flèche est induite par \eqref{tfra5f} et la dernière flèche est le morphisme d'adjonction, est donc
l'isomorphisme \eqref{tfra6b}; d'où la première proposition. 

\begin{prop}\label{tfra7}
Les hypothèses étant celles de \ref{tfra5}, supposons, de plus, que $w$ et $\ow$ soient des isomorphismes, que le diagramme
\begin{equation}
\xymatrix{
X'\ar[d]_u&Y'\ar[d]^v\ar[l]_{f'}\\
X&Y\ar[l]_{f}}
\end{equation}
soit cartésien et que l'une des deux conditions suivantes soit remplie 
\begin{itemize}
\item[{\rm (i)}] le morphisme $u\colon X'\rightarrow X$ est étale;
\item[{\rm (ii)}] le schéma $X'$ est le localisé strict de $X$ en un point géométrique $\ox$, 
$u$ étant le morphisme canonique, le morphisme $f$ est cohérent et le schéma $Z$ est cohérent. 
\end{itemize}
Alors, le morphisme $\Phi\colon \tG'\rightarrow \tG$ \eqref{tfra5b} est une équivalence de topos et le morphisme \eqref{tfra5f}
\begin{equation}\label{tfra7a}
\ocB\rightarrow \Phi_*(\ocB')
\end{equation}
est un isomorphisme. 
\end{prop}

La première assertion est démontrée dans \ref{lptfr7} et est mentionnée pour mémoire. 
Soit $(W\rightarrow U\leftarrow V)$ un objet de $G$, et posons $(W\rightarrow U'\leftarrow V')=\Phi^+(W\rightarrow U\leftarrow V)$ \eqref{tfr7d}. 
Les morphismes canoniques 
\begin{eqnarray}
W\times_{X'}Y'&\rightarrow&W\times_XY,\\
V'&\rightarrow& V\times_UU',\\
W\times_{U'}V'&\rightarrow& W\times_UV
\end{eqnarray}
étant des isomorphismes, le morphisme canonique \eqref{tfra5f}
\begin{equation}
\ocB(W\rightarrow U\leftarrow V)\rightarrow \Phi_*(\ocB')(W\rightarrow U\leftarrow V)
\end{equation}
est un isomorphisme, d'où la proposition.

\begin{prop}\label{tfra9}
Les hypothèses étant celles de \ref{tfra5}, supposons, de plus, 
que $Z'$ soit le localisé strict de $Z$ en un point géométrique $\oz$, $w$ étant le morphisme canonique, que les deux diagrammes
\begin{equation}\label{tfra9a}
\xymatrix{
\oZ'\ar[r]^{\hbar'}\ar[d]_\ow&Z'\ar[d]^w\\
\oZ\ar[r]^{\hbar}&Z}
\ \ \
\xymatrix{
X'\ar[d]_u&Y'\ar[d]^v\ar[l]_{f'}\\
X&Y\ar[l]_{f}}
\end{equation}
soient cartésiens,  que les morphismes $f$ et $\hbar$ soient cohérents et que l'une des deux conditions suivantes soit remplie 
\begin{itemize}
\item[{\rm (i)}] le morphisme $u$ est étale;
\item[{\rm (ii)}] le schéma $X'$ est le localisé strict de $X$ en un point géométrique $\ox$, $u$ étant le morphisme canonique.
\end{itemize}
Alors, le morphisme 
\begin{equation}\label{tfra9b}
\Phi^{-1}(\ocB)\rightarrow \ocB'
\end{equation}
adjoint du morphisme \eqref{tfra5f} est un isomorphisme. 
\end{prop}

On notera que dans le cas (ii), le morphisme $Z'\rightarrow X'$ fournit une flèche de spécialisation de $g(\oz)$ dans $\ox$, 
mais $g(\oz)$ n'est pas nécessairement égal à $\ox$. 

Comme le morphisme canonique 
\begin{equation}\label{tfra9c}
Z'\times_{X'}Y'\rightarrow Z'\times_{X}Y
\end{equation}
est un isomorphisme, on peut considérer le diagramme commutatif
\begin{equation}\label{tfra9d}
\xymatrix{
Z'\times_{X'}Y'\ar[r]\ar@{=}[d]&\oZ'\ar[r]\ar@{=}[d]&Z'\ar[r]\ar@{=}[d]&X'\ar[d]\ar@{}[rd]|\Box&Y'\ar[l]\ar[d]\\
Z'\times_{X}Y\ar[r]\ar[d]\ar@{}[rd]|\Box&\oZ'\ar[r]\ar[d]\ar@{}[rd]|\Box&Z'\ar[r]\ar[d]&X\ar@{=}[d]&Y\ar[l]\ar@{=}[d]\\
Z\times_{X}Y\ar[r]&\oZ\ar[r]&Z\ar[r]&X&Y \ar[l] }
\end{equation}
Compte tenu de \ref{tfra7}, on peut se borner au cas où $u$ est un isomorphisme. 
D'après \ref{tfra6}, quitte à remplacer $X$ par un voisinage ouvert de Zariski affine $X_1$ contenant le support de $g(\oz)$, et 
$Z$, $\oZ$ et $Y$ par $Z\times_XX_1$, $\oZ\times_XX_1$ et $Y\times_XX_1$, respectivement, on peut se borner au cas où $X$ et $Y$ sont cohérents. 
Quitte à remplacer $Z$ par un voisinage étale affine $Z_1$ de $\oz$ dans $Z$ et $\oZ$ par $\oZ\times_ZZ_1$, on peut, de plus, supposer $Z$ affine. 

Notons $I$ la catégorie des $Z$-schémas étales affines $\oz$-pointés et pour tout $i\in \ob(I)$, $Z_i$ le $Z$-schéma étale affine $\oz$-pointé 
correspondant. On a alors dans la catégorie des schémas  
\begin{equation}\label{tfra9e}
Z'=\underset{\underset{i\in I}{\longleftarrow}}{\lim}Z_i.
\end{equation}
Avec les notations de \ref{lptfr1}, on obtient un foncteur
\begin{equation}\label{tfra9f}
\varphi\colon I\rightarrow \fM_{/(Z\rightarrow X\leftarrow Y)}, \ \ \ i\mapsto (Z_i\rightarrow X\leftarrow Y). 
\end{equation}
Reprenons les notations de \ref{lptfr2} et \ref{lptfr4} pour ce foncteur.
En vertu de \ref{lptfr3}, le topos $\tG'$ est la limite projective du topos fibré $\fF_\varphi/I$.  

Pour tout $i\in \ob(I)$,  posons $\oZ_i=\oZ\times_ZZ_i$ et notons $\ocB_i$ l'anneau de $\tfG_{(Z_i\rightarrow X\leftarrow Y)}$ 
défini dans \ref{tfra1} par la factorisation $Z_i\times_{X}Y \rightarrow \oZ_i\rightarrow Z_i$ de la projection canonique \eqref{tfra9d}. 
Pour tout morphisme $j\rightarrow i$ de $I$, notons 
\begin{equation}\label{tfra9g}
\Phi_{ji}\colon \tfG_{(Z_j\rightarrow X\leftarrow Y)} \rightarrow \tfG_{(Z_i\rightarrow X\leftarrow Y)}
\end{equation}
le morphisme de fonctorialité \eqref{lptfr1f}. 
On a alors un homomorphisme canonique $\ocB_i\rightarrow \Phi_{ji*}(\ocB_j)$ \eqref{tfra5f}. Ces homomorphismes
vérifient une relation de compatibilité pour la composition des morphismes de $I$ du type (\cite{egr1} (1.1.2.2)).
Ils définissent donc un anneau de $\Top(\fG_\varphi)$ que l'on note $\{i\mapsto \ocB_i\}$ \eqref{lptfr4a}.
Pour tout $i\in \ob(I)$, on a un morphisme canonique $Z'\rightarrow Z_i$ qui induit par fonctorialité un morphisme de topos annelé
\begin{equation}\label{tfra9h}
\Phi_{i}\colon (\tfG_{(Z'\rightarrow X\leftarrow Y)}, \ocB') \rightarrow (\tfG_{(Z_i\rightarrow X\leftarrow Y)}, \ocB_i).
\end{equation}
La collection de ces homomorphismes $(\ocB_i\rightarrow \Phi_{i*}(\ocB'))_{i\in I}$ induit un homomorphisme de $\Top(\fG_\varphi)$
\begin{equation}\label{tfra9i}
\{i\mapsto \ocB_i\}\rightarrow \varpi_*(\ocB').
\end{equation}
Montrons que l'homomorphisme adjoint 
\begin{equation}\label{tfra9j}
\varpi^*(\{i\mapsto \ocB_i\})\rightarrow \ocB'
\end{equation}
est un isomorphisme. Comme le foncteur $\varpi_*$ est pleinement fidèle \eqref{lptfr4c}, il suffit de montrer 
que l'homomorphisme induit 
\begin{equation}\label{tfra9k}
\varpi_*(\varpi^*(\{i\mapsto \ocB_i\})) \rightarrow \varpi_*(\ocB') 
\end{equation}
est un isomorphisme. Compte tenu de (\cite{sga4} VI 8.5.3), cela revient à montrer que pour tout $i\in \ob(I)$, 
l'homomorphisme canonique de $\tfG_{(Z_i\rightarrow X\leftarrow Y)}$
\begin{equation}\label{tfra9l}
\underset{\underset{j\in I_{/i}}{\longrightarrow}}{\lim}\ \Phi_{ji*}(\ocB_j)\rightarrow \Phi_{i*}(\ocB'),
\end{equation}
est un isomorphisme.

Soit $(W\rightarrow U\leftarrow V)$ un objet de $\fG_{(Z_i\rightarrow X\leftarrow Y)}$. On pose $\oW=W\times_Z\oZ$ et on note 
$\oW^V$ la fermeture intégrale de $\oW$ dans $W\times_UV$. Le morphisme canonique
\begin{equation}\label{tfra9m}
Z'\times_{Z_i}\oW^V=\underset{\underset{j\in I_{/i}}{\longleftarrow}}{\lim}Z_j\times_{Z_i}\oW^V
\end{equation}
est un isomorphisme d'après (\cite{ega4} 8.2.5).  Comme $\hbar$ est cohérent, le schéma $\oW^{V}$ est cohérent \eqref{tfra2}. 
On en déduit par (\cite{sga4} VI 5.2) et (\cite{agt} III.8.11(iii) et III.8.22) que l'homomorphisme canonique 
\begin{equation}\label{tfra9n}
\underset{\underset{j\in I_{/i}}{\longrightarrow}}{\lim}\ \ocB_j(W\times_{Z_i}Z_j\rightarrow U\leftarrow V)\rightarrow 
\ocB'(W\times_{Z_i}Z'\rightarrow U\leftarrow V)
\end{equation}
est un isomorphisme, d'où l'isomorphisme \eqref{tfra9l}.

On déduit de l'isomorphisme \eqref{tfra9j} et de \eqref{lptfr4d} que l'homomorphisme canonique 
\begin{equation}\label{tfra9o}
\underset{\underset{i\in I^\circ}{\longrightarrow}}{\lim}\ \Phi_i^{-1}(\ocB_i) \rightarrow \ocB' 
\end{equation}
est un isomorphisme. 
Pour tout $i\in \ob(I)$, le morphisme canonique $Z_i\rightarrow Z$ induit par fonctorialité un 
morphisme de topos annelés 
\begin{equation}
\rho_i\colon (\tfG_{(Z_i\rightarrow X\leftarrow Y)},\ocB_i)\rightarrow (\tG,\ocB).
\end{equation}
Comme l'homomorphisme $\rho_i^{-1}(\ocB)\rightarrow \ocB_i$ est un isomorphisme en vertu de \ref{tfra6}(i), 
la proposition résulte de l'isomorphisme \eqref{tfra9o}.

\chapter{Cohomologie du topos de Faltings}\label{cohtopfal}

\section{Hypothèses et notations; schéma logarithmique adéquat}\label{TFA}

\subsection{}\label{TFA1}
Dans ce chapitre, $K$ désigne un corps de valuation discrète complet de 
caractéristique $0$, à corps résiduel {\em algébriquement clos} $k$ de caractéristique $p>0$,  
$\co_K$ l'anneau de valuation de $K$, $\oK$ une clôture algébrique de $K$, $\co_\oK$ la clôture intégrale de $\co_K$ dans $\oK$,
$\fm_\oK$ l'idéal maximal de $\co_\oK$ et $G_K$ le groupe de Galois de $\oK$ sur $K$.
On note $\co_C$ le séparé complété $p$-adique de $\co_\oK$, $\fm_C$ son idéal maximal,
$C$ son corps des fractions et $v\colon C\rightarrow \mQ\cup\{\infty\}$ la valuation normalisée par $v(p)=1$. 
On désigne par $\hmZ(1)$ et $\mZ_p(1)$ les $\mZ[G_K]$-modules 
\begin{eqnarray}
\hmZ(1)&=&\underset{\underset{n\geq 1}{\longleftarrow}}{\lim}\ \mu_{n}(\co_{\oK}),\label{TFA1aa}\\
\mZ_p(1)&=&\underset{\underset{n\geq 0}{\longleftarrow}}{\lim}\ \mu_{p^n}(\co_{\oK}),\label{TFA1a}
\end{eqnarray}  
où $\mu_n(\co_{\oK})$ désigne le sous-groupe des racines $n$-ièmes de l'unité dans $\co_\oK$. 
Pour tout $\mZ_p[G_K]$-module $M$ et tout entier $n$, on pose $M(n)=M\otimes_{\mZ_p}\mZ_p(1)^{\otimes n}$.

Comme $\co_\oK$ satisfait les conditions requises dans \ref{mptf1}, 
il est loisible de considérer les notions de $\alpha$-algèbre introduites dans les sections \ref{alpha}--\ref{aet} relativement à $\co_\oK$.   

On pose $S=\Spec(\co_K)$ et $\oS=\Spec(\co_\oK)$ et 
on note $s$ (resp.  $\eta$, resp. $\oeta$) le point fermé de $S$ (resp.  générique de $S$, resp. générique de $\oS$).
Pour tout entier $n\geq 1$, on pose $S_n=\Spec(\co_K/p^n\co_K)$. 
On munit $S$ de la structure logarithmique $\cM_S$ définie par son point fermé, 
autrement dit, $\cM_S=u_*(\co_\eta^\times)\cap \co_S$, où $u\colon \eta\rightarrow S$ est l'injection canonique
(cf. \cite{agt} II.5 pour un lexique de géométrie logarithmique). 

Pour tout $S$-schéma $X$, on note $X_s$ (resp. $X_\eta$, resp. $X_\oeta$) 
la fibre fermée (resp. générique, resp. géométrique générique) de $X$ au dessus de $S$, 
et on pose 
\begin{equation}\label{TFA1b}
\oX=X\times_S\oS \ \ \ {\rm et}\ \ \ X_n=X\times_SS_n.
\end{equation}

\subsection{}\label{TFA3}
Dans ce chapitre, $f\colon (X,\cM_X)\rightarrow (S,\cM_S)$ 
désigne un morphisme {\em adéquat} de schémas logarithmiques (\cite{agt} III.4.7).
On désigne par $X^\circ$ le sous-schéma ouvert maximal de $X$
où la structure logarithmique $\cM_X$ est triviale~; c'est un sous-schéma ouvert de $X_\eta$.
On note $j\colon X^\circ\rightarrow X$ l'injection canonique. 
Pour tout $X$-schéma $U$, on pose  
\begin{equation}\label{TFA3a}
U^\circ=U\times_XX^\circ.
\end{equation} 
On note $\hbar\colon \oX\rightarrow X$ et $h\colon \oX^\circ\rightarrow X$ les morphismes canoniques \eqref{TFA1b}, de sorte que 
l'on a $h=\hbar\circ j_\oX$. Pour alléger les notations, on pose
\begin{equation}\label{TFA3b}
\tOmega^1_{X/S}=\Omega^1_{(X,\cM_X)/(S,\cM_S)},
\end{equation}
que l'on considère comme un faisceau de $X_\zar$ ou $X_\et$, selon le contexte (cf. \ref{notconv12}).

\subsection{}\label{TFA5}
Pour tout entier $n\geq 1$, on note $a\colon X_s\rightarrow X$, $a_n\colon X_s\rightarrow X_n$, 
$\iota_n\colon X_n\rightarrow X$ et $\oiota_n\colon \oX_n\rightarrow \oX$ les injections canoniques \eqref{TFA1b}. 
Le corps résiduel de $\co_K$ étant algébriquement clos, il existe un unique $S$-morphisme $s\rightarrow \oS$. 
Celui-ci induit des immersions fermées $\oa\colon X_s\rightarrow \oX$ et $\oa_n\colon X_s\rightarrow \oX_n$
qui relèvent $a$ et $a_n$, respectivement. 
\begin{equation}\label{TFA5a}
\xymatrix{
{X_s}\ar[r]_{\oa_n}\ar@{=}[d]\ar@/^1pc/[rr]^{\oa}&{\oX_n}\ar[r]_{\oiota_n}\ar[d]^{\hbar_n}&\oX\ar[d]^\hbar\\
{X_s}\ar[r]^{a_n}\ar@/_1pc/[rr]_{a}&{X_n}\ar[r]^{\iota_n}&X}
\end{equation}
Comme $\hbar$ est entier et que $\hbar_n$ est un homéomorphisme
universel, pour tout faisceau $\cF$ de $\oX_\et$, le morphisme de changement de base 
\begin{equation}\label{TFA5b}
a^*(\hbar_*(\cF))\rightarrow \oa^*(\cF)
\end{equation}
est un isomorphisme (\cite{sga4} VIII 5.6). Par ailleurs, $\oa_n$ étant un homéomorphisme universel, on peut considérer $\co_{\oX_n}$
comme un faisceau de $X_{s,\zar}$ ou $X_{s,\et}$, selon le contexte (cf. \ref{notconv12}).

\subsection{}\label{TFA6}
On désigne par
\begin{equation}\label{TFA6a}
\pi\colon E\rightarrow \Et_{/X}
\end{equation}
le $\mU$-site fibré de Faltings associé au morphisme $h\colon \oX^\circ\rightarrow X$ (cf. \ref{tf1}).
Pour tout $U\in \ob(\Et_{/X})$, on note 
\begin{equation}\label{TFA6F}
\iota_{U!}\colon \Et_{\rf/\oU^\circ}\rightarrow E
\end{equation} 
le foncteur canonique \eqref{tf1c}. On munit $E$ de la topologie co-évanescente définie par $\pi$ (\cite{agt} VI.5.3)  
et on note $\tE$ le topos des faisceaux de $\mU$-ensembles sur $E$. On désigne par 
\begin{eqnarray}
\sigma\colon \tE \rightarrow X_\et,\label{TFA6c}\\
\beta\colon \tE \rightarrow \oX^\circ_\fet,\label{TFA6b}\\
\psi\colon \oX^\circ_\et\rightarrow \tE,\label{TFA6d}\\
\rho\colon X_\et\gtimes_{X_\et}\oX^\circ_\et\rightarrow \tE,\label{TFA6e}
\end{eqnarray}
les morphismes canoniques \eqref{tf1kk}, \eqref{tf1k}, \eqref{tf1m} et \eqref{tf3b}. 

\subsection{}\label{TFA66}
On note $\tE_s$ le sous-topos fermé de $\tE$ complémentaire de l'ouvert $\sigma^*(X_\eta)$ \eqref{tf20}, 
\begin{equation}\label{TFA66a}
\delta\colon \tE_s\rightarrow \tE
\end{equation} 
le plongement canonique et 
\begin{equation}\label{TFA66b}
\sigma_s\colon \tE_s\rightarrow X_{s,\et}
\end{equation} 
le morphisme de topos induit par $\sigma$ \eqref{TFA6c} (cf. \ref{tf23}). Le diagramme de morphismes de topos 
\begin{equation}\label{TFA66c}
\xymatrix{
{\tE_s}\ar[r]^{\sigma_s}\ar[d]_{\delta}&{X_{s,\et}}\ar[d]^{a}\\
{\tE}\ar[r]^\sigma&{X_\et}}
\end{equation}
est commutatif à isomorphisme près. Les foncteurs $a_*$ et $\delta_*$ étant exacts, 
pour tout groupe abélien $F$ de $\tE_s$ et tout entier $i\geq 0$, on a un isomorphisme canonique 
\begin{equation}\label{TFA66d}
a_*(\rR^i\sigma_{s*}(F))\stackrel{\sim}{\rightarrow}\rR^i\sigma_*(\delta_*F). 
\end{equation}

\subsection{}\label{TFA2}
D'après (\cite{agt} III.4.2(iii)), $\oX$ est normal et localement irréductible (\cite{agt} III.3.1). 
Par ailleurs, l'immersion $j\colon X^\circ\rightarrow X$ est quasi-compacte puisque $X$ est noethérien.  
Pour tout objet $(V\rightarrow U)$ de $E$, on note $\oU^V$ la fermeture intégrale de $\oU$ dans $V$. 
Pour tout morphisme $(V'\rightarrow U')\rightarrow (V\rightarrow U)$ de $E$, on a un morphisme canonique 
$\oU'^{V'}\rightarrow \oU^V$ qui s'insère dans un diagramme commutatif 
\begin{equation}\label{TFA2a}
\xymatrix{
V'\ar[r]\ar[d]&{\oU'^{V'}}\ar[d]\ar[r]&\oU'\ar[r]\ar[d]&U'\ar[d]\\
V\ar[r]&{\oU^V}\ar[r]&\oU\ar[r]&U}
\end{equation} 
On désigne par $\ocB$ le préfaisceau sur $E$ défini pour tout $(V\rightarrow U)\in \ob(E)$, par 
\begin{equation}\label{TFA2b}
\ocB((V\rightarrow U))=\Gamma(\oU^V,\co_{\oU^V}).
\end{equation}
C'est un faisceau pour la topologie co-évanescente de $E$ en vertu de (\cite{agt} III.8.16). 
Pour tout $U\in \ob(\Et_{/X})$, on pose 
\begin{equation}\label{TFA2d}
\ocB_{U}=\ocB\circ \iota_{U!}.
\end{equation} 

D'après (\cite{agt} III.8.17), on a un homomorphisme canonique 
\begin{equation}\label{TFA2c}
\sigma^*(\hbar_*(\co_\oX))\rightarrow \ocB.
\end{equation}
Sauf mention explicite du contraire, on considère $\sigma$ \eqref{TFA6c}
comme un morphisme de topos annelés
\begin{equation}\label{TFA2e}
\sigma\colon (\tE,\ocB)\rightarrow (X_\et,\hbar_*(\co_\oX)).
\end{equation}

Notons encore $\co_\oK$ le faisceau constant de $\oX^\circ_\fet$ de valeur $\co_\oK$. 
Sauf mention explicite du contraire, on considère $\beta$ \eqref{TFA6b} comme un morphisme de topos annelés
\begin{equation}\label{TFA2g}
\beta\colon (\tE,\ocB)\rightarrow (\oX^\circ_\fet,\co_\oK).
\end{equation}

\subsection{}\label{TFA8}
Pour tout entier $n\geq 0$ et tout $U\in \ob(\Et_{/X})$, on pose 
\begin{eqnarray}
\ocB_n&=&\ocB/p^n\ocB,\label{TFA8a}\\
\ocB_{U,n}&=&\ocB_U/p^n\ocB_U.\label{TFA8b}
\end{eqnarray}
Les correspondances $\{U\mapsto p^n\ocB_U\}$ et $\{U\mapsto \ocB_{U,n}\}$
forment naturellement des préfaisceaux sur $E$ \eqref{tf1h}, et les morphismes canoniques 
\begin{eqnarray}
\{U\mapsto p^n\ocB_U\}^\tta&\rightarrow&p^n\ocB,\label{TFA8c}\\
\{U\mapsto \ocB_{U,n}\}^\tta&\rightarrow&\ocB_n,\label{TFA8d}
\end{eqnarray}
où  les termes de gauche désignent les faisceaux associés dans $\tE$, sont des isomorphismes
en vertu de (\cite{agt} VI.8.2 et VI.8.9). D'après (\cite{agt} III.9.7), $\ocB_n$ est un anneau de $\tE_s$. 

Si $n\geq 1$, on désigne par
\begin{equation}\label{TFA8e}
\sigma_n\colon (\tE_s,\ocB_n)\rightarrow (X_{s,\et},\co_{\oX_n})
\end{equation}
le morphisme de topos annelés induit par $\sigma$ \eqref{TFA2e} (cf. \cite{agt} (III.9.9.4)) et par 
\begin{equation}\label{TFA8f}
\tau_n\colon (\tE_s,\ocB_n)\rightarrow (X_{s,\zar},\co_{\oX_n})
\end{equation}
le composé de $\sigma_n$ et du morphisme canonique \eqref{notconv12j} 
\begin{equation}\label{TFA8h}
u_n\colon (X_{s,\et},\co_{\oX_n})\rightarrow (X_{s,\zar},\co_{\oX_n}).
\end{equation}

Notant encore $\co_\oK/p^n\co_\oK$ le faisceau constant de $\oX^\circ_\fet$ de valeur $\co_\oK/p^n\co_\oK$, on désigne par 
\begin{equation}\label{TFA8g}
\beta_n\colon (\tE_s,\ocB_n)\rightarrow (\oX^\circ_\fet,\co_\oK/p^n\co_\oK)
\end{equation}
le morphisme induit par $\beta$ \eqref{TFA2g}.

\begin{rema}\label{TFA81}
Sous les hypothèses de \ref{TFA8}, on notera que l'homomorphisme canonique \eqref{TFA6F}
\begin{equation}\label{TFA81a}
\ocB_{U,n}\rightarrow \ocB_n\circ \iota_{U!}
\end{equation} 
n'est pas en général un isomorphisme~; c'est pourquoi nous n'utiliserons pas la notation $\ocB_{n,U}$. 
Toutefois, sous certaines hypothèses \eqref{tpcg1}, nous montrerons dans \ref{amtF27} que pour 
tout schéma affine $U$ de $\Et_{/X}$, l'homomorphisme \eqref{TFA81a} est un $\alpha$-isomorphisme.
\end{rema}

\subsection{}\label{TFA9}
Soient $U$ un objet de $\Et_{/X}$, $\oy$ un point géométrique de $\oU^\circ$. 
Le schéma $\oU$ étant localement irréductible d'après (\cite{agt} III.3.3 et III.4.2(iii)),  
il est la somme des schémas induits sur ses composantes irréductibles. On note $\oU^\star$
la composante irréductible de $\oU$ (ou ce qui revient au même, sa composante connexe) contenant $\oy$. 
De même, $\oU^\circ$ est la somme des schémas induits sur ses composantes irréductibles
et $\oU^{\star \circ}=\oU^\star\times_{X}X^\circ$ est la composante irréductible de $\oU^\circ$ contenant $\oy$. 
On note $\bB_{\pi_1(\oU^{\star \circ},\oy)}$ le topos classifiant du groupe profini $\pi_1(\oU^{\star \circ},\oy)$ et
\begin{equation}\label{TFA9a}
\nu_\oy\colon \oU^{\star \circ}_\fet \stackrel{\sim}{\rightarrow}\bB_{\pi_1(\oU^{\star \circ},\oy)}
\end{equation}
le foncteur fibre  de $\oU^{\star \circ}_\fet$ en $\oy$ \eqref{notconv11c}. On pose
\begin{equation}\label{TFA9b}
\oR^\oy_U=\nu_\oy(\ocB_U|\oU^{\star \circ}).
\end{equation}
Explicitement, soit $(V_i)_{i\in I}$ le revêtement universel normalisé de $\oU^{\star \circ}$ en $\oy$ \eqref{notconv11}.
Pour chaque $i\in I$, $(V_i\rightarrow U)$ est naturellement un objet de $E$. 
On note $\oU^{V_i}$ la fermeture intégrale de $\oU$ dans $V_i$.
Les schémas $(\oU^{V_i})_{i\in I}$ forment alors un système projectif filtrant, et on a   
\begin{equation}\label{TFA9c}
\oR^\oy_U=\underset{\underset{i\in I}{\longrightarrow}}{\lim}\ \Gamma(\oU^{V_i},\co_{\oU^{V_i}}).
\end{equation} 

\begin{rema}\label{TFA10}
Sous les hypothèses de \ref{TFA9}, si, de plus, $U$ est affine, 
l'anneau $\oR^\oy_U$ est intègre et normal. En effet, le schéma $\oU^\star$ est affine, intègre et normal. 
Pour tout $i\in I$, $V_i$ étant connexe, il est alors intègre et normal. 
Par suite, $\oU^{V_i}$ est intègre, normal et entier sur $\oU^\star$ (\cite{ega2} 6.3.7). 
Par ailleurs, pour tous $(i,j)\in I^2$ avec  
$i\geq j$, le morphisme canonique $\oU^{V_i}\rightarrow \oU^{V_j}$ est entier et dominant. 
L'assertion recherchée s'ensuit d'après (\cite{ega1n} 0.6.1.6(i) et 0.6.5.12(ii)). 
\end{rema}

\subsection{}\label{TFA14}
Soient $\ox$ un point géométrique de $X$, $\uX$ le localisé strict de $X$ en $\ox$. On pose $\uoX=\uX\times_S\oS$ \eqref{TFA1b} et 
$\uoX^\circ=\uoX\times_XX^\circ$ \eqref{TFA3a}. 
On désigne par $\fV_\ox$ la catégorie des $X$-schémas étales $\ox$-pointés, ou ce qui revient au même,
la catégorie des voisinages du point de $X_\et$ associé à $\ox$ dans le site $\Et_{/X}$ (\cite{sga4} IV 6.8.2 et VIII 3.9).
Pour tout objet $(U,\fp\colon \ox\rightarrow U)$ de $\fV_\ox$, 
on note encore $\fp\colon \uX\rightarrow U$ le morphisme déduit de 
$\fp$ (\cite{sga4} VIII 7.3) et on pose
\begin{equation}\label{TFA14e}
\ofp^\circ=\fp\times_X\oX^\circ\colon \uoX^\circ \rightarrow \oU^\circ.
\end{equation}

On note
\begin{equation}\label{TFA14a}
\varphi_\ox\colon \tE\rightarrow \uoX^\circ_\fet
\end{equation}
le foncteur canonique défini dans \eqref{tf10b}.

Le foncteur composé $\varphi_\ox\circ\beta^*$ est canoniquement isomorphe au foncteur image inverse par 
le morphisme canonique $\uoX^\circ_\fet\rightarrow \oX^\circ_\fet$, d'après \eqref{tf8g} et (\cite{agt} (VI.10.12.6)). 
Pour tout objet $F$ de $X_\et$, on a un isomorphisme canonique et fonctoriel 
\begin{equation}\label{TFA14b}
\varphi_\ox(\sigma^*(F))\stackrel{\sim}{\rightarrow} F_\ox,
\end{equation}
de but le faisceau constant de $\uoX^\circ_\fet$ de valeur $F_\ox$, en vertu de (\cite{agt} VI.10.24 et (VI.10.12.6)). 

Pour tout objet $F=\{U\mapsto F_U\}$ de $\hE$ \eqref{tf1},
on note $F^\tta$ le faisceau de $\tE$ associé à $F$ \eqref{TFA6}, et pour tout $U\in \ob(\Et_{/X})$, 
$F_U^\tta$ le faisceau de $\oU^\circ_\fet$ associé à $F_U$. 
D'après (\cite{agt} VI.10.37),  on a un isomorphisme canonique et fonctoriel
\begin{equation}\label{TFA14c}
\varphi_\ox(F^\tta) \stackrel{\sim}{\rightarrow}\underset{\underset{(U,\fp)\in \fV_\ox^\circ}{\longrightarrow}}{\lim}\ (\ofp^\circ)_\fet^*(F^\tta_U).
\end{equation}

En vertu de (\cite{agt} VI.10.30), pour tout groupe abélien $F$ de $\tE$
et tout entier $q\geq 0$, on a un isomorphisme canonique et fonctoriel
\begin{equation}\label{TFA14d}
\rR^q\sigma_*(F)_\ox\stackrel{\sim}{\rightarrow}\rH^q(\uoX^\circ_\fet,\varphi_\ox(F)). 
\end{equation}

\subsection{}\label{TFA11}
Soient $(\oy\rightsquigarrow \ox)$ un point de $X_\et\gtimes_{X_\et}\oX^\circ_\et$ \eqref{topfl17}, 
$\uX$ le localisé strict de $X$ en $\ox$, $\fV_\ox$ la catégorie des $X$-schémas étales $\ox$-pointés. 
On pose $\uoX=\uX\times_S\oS$ \eqref{TFA1b} et $\uoX^\circ=\uoX\times_XX^\circ$ \eqref{TFA3a}.
Le $X$-morphisme $u\colon \oy\rightarrow \uX$
définissant $(\oy\rightsquigarrow \ox)$ se relève en un $\oX^\circ$-morphisme $v\colon \oy\rightarrow \uoX^\circ$ et 
induit donc un point géométrique de $\uoX^\circ$ que l'on note aussi (abusivement) $\oy$.
Pour tout objet $(U,\fp\colon \ox\rightarrow U)$ de $\fV_\ox$, 
on note encore $\fp\colon \uX\rightarrow U$ le morphisme déduit de $\fp$ (\cite{sga4} VIII 7.3) et on pose
\begin{equation}\label{TFA11b}
\ofp^\circ=\fp\times_X\oX^\circ\colon \uoX^\circ \rightarrow \oU^\circ.
\end{equation}
On note aussi (abusivement) $\oy$ le point géométrique $\ofp^\circ(v(\oy))$ de $\oU^\circ$.

Pour tout objet $F=\{U\mapsto F_U\}$ de $\hE$ \eqref{tf1h}, on note
$F^\tta$ le faisceau de $\tE$ associé à $F$, et pour tout $U\in \ob(\Et_{/X})$, 
$F_U^\tta$ le faisceau de $\oU^\circ_\fet$ associé à $F_U$. 
D'après (\cite{agt} VI.10.36 et (VI.9.3.4)), on a un isomorphisme canonique et fonctoriel 
\begin{equation}\label{TFA11a}
(F^\tta)_{\rho(\oy\rightsquigarrow \ox)} \stackrel{\sim}{\rightarrow} 
\underset{\underset{(U,\fp)\in \fV_\ox^\circ}{\longrightarrow}}{\lim}\ (F^\tta_U)_{\rho_{\oU^\circ}(\oy)},
\end{equation}
où $\rho$ est le morphisme \eqref{TFA6e} et $\rho_{\oU^\circ}\colon \oU^\circ_\et\rightarrow \oU^\circ_\fet$ 
est le morphisme canonique \eqref{notconv10a}. 
Compte tenu de \eqref{TFA9b} et (\cite{agt} VI.9.9), on en déduit un isomorphisme canonique de $\Gamma(\uoX,\co_{\uoX})$-algèbres
\begin{equation}\label{TFA11c}
\ocB_{\rho(\oy\rightsquigarrow \ox)}\stackrel{\sim}{\rightarrow} 
\underset{\underset{(U,\fp)\in \fV_\ox^\circ}{\longrightarrow}}{\lim}\ \oR^{\oy}_U,
\end{equation}
où $\oR^{\oy}_U$ est la $\Gamma(\oU,\co_\oU)$-algèbre définie dans \eqref{TFA9c}.

\subsection{}\label{TFA12}
Conservons les hypothèses et notations de \ref{TFA11}; supposons, de plus, que {\em $\ox$ soit au-dessus de $s$}. 
D'après (\cite{agt} III.3.7), $\uoX$ est normal et strictement local (et en particulier intègre). 
Pour tout objet $(U,\fp\colon \ox\rightarrow U)$ de $\fV_\ox$, on désigne par $\oU^\star$  
la composante irréductible de $\oU$ contenant $\oy$ 
et on pose $\oU^{\star\circ}=\oU^\star\times_XX^\circ$ qui est la composante irréductible de $\oU^\circ$ contenant $\oy$ (cf. \ref{TFA9}).
Le morphisme $\ofp^\circ \colon \uoX^\circ\rightarrow \oU^\circ$
\eqref{TFA11b} se factorise donc à travers $\oU^{\star \circ}$. 
On désigne par $\bB_{\pi_1(\uoX^{\circ},\oy)}$ le topos classifiant du groupe profini 
$\pi_1(\uoX^{\circ},\oy)$ et par
\begin{equation}\label{TFA12d}
\nu_\oy\colon \uoX^\circ_\fet \stackrel{\sim}{\rightarrow}\bB_{\pi_1(\uoX^\circ,\oy)}
\end{equation}
le foncteur fibre de $\uoX^\circ_\fet$ en $\oy$ \eqref{notconv11c}. 
D'après (\cite{agt} VI.10.31 et VI.9.9), le foncteur composé
\begin{equation}\label{TFA12e}
\xymatrix{
{\tE}\ar[r]^-(0.5){\varphi_\ox}&{\uoX^\circ_\fet}\ar[r]^-(0.5){\nu_\oy}&{\bB_{\pi_1(\uoX^\circ,\oy)}}\ar[r]&\Ens},
\end{equation}
où $\varphi_\ox$ est le foncteur canonique \eqref{TFA14a} et la dernière flèche est le foncteur d'oubli de l'action de $\pi_1(\uoX^\circ,\oy)$, 
est canoniquement isomorphe au foncteur fibre associé au point $\rho(\oy\rightsquigarrow \ox)$ de $\tE$.  
On observera (\cite{agt} (VI.10.25.6))
que pour tout objet $F$ de $X_\et$, l'isomorphisme \eqref{TFA14b} induit l'isomorphisme canonique (\cite{agt} (VI.10.18.1))
\begin{equation}
(\sigma^*F)_{\rho(\oy\rightsquigarrow \ox)} \stackrel{\sim}{\rightarrow} F_{\ox}.
\end{equation}

On définit la $\Gamma(\uoX,\co_{\uoX})$-algèbre $\oR_{\uX}^\oy$ de $\bB_{\pi_1(\uoX^\circ,\oy)}$, 
que l'on note aussi $\uoR$ lorsqu'il n'y a aucun risque d'ambiguïté, par la formule
\begin{equation}\label{TFA12f}
\oR_{\uX}^\oy=\underset{\underset{(U,\fp)\in \fV_\ox^\circ}{\longrightarrow}}{\lim}\ \oR^{\oy}_U, 
\end{equation}
où l'on considère $\oR^{\oy}_U$ comme une $\Gamma(\oU,\co_\oU)$-algèbre de $\bB_{\pi_1(\oU^{\star\circ},\oy)}$ \eqref{TFA9c}.  
L'isomorphisme \eqref{TFA14c} induit un isomorphisme de $\Gamma(\uoX,\co_{\uoX})$-algèbres de 
$\bB_{\pi_1(\uoX^\circ,\oy)}$
\begin{equation}\label{TFA12g}
\nu_\oy(\varphi_\ox(\ocB))\stackrel{\sim}{\rightarrow} \oR_{\uX}^\oy,
\end{equation} 
dont l'isomorphisme de $\Gamma(\uoX,\co_{\uoX})$-algèbres sous-jacent est \eqref{TFA11c}.

\begin{remas}\label{TFA15}
Soient $\ox$ un point géométrique de $X$ au-dessus de $s$, $\uX$ le localisé strict de $X$ en $\ox$, 
$\varphi_\ox\colon \tE\rightarrow \uoX^\circ_\fet$ le foncteur canonique \eqref{TFA14a}. Alors, 
\begin{itemize}
\item[(i)] Il existe un point $(\oy\rightsquigarrow \ox)$ de $X_\et\gtimes_{X_\et}\oX^\circ_\et$ \eqref{topfl17}. 
En effet, d'après (\cite{agt} III.3.7), $\uoX$ est normal et strictement local (et en particulier intègre). 
Comme $X^\circ$ est schématiquement dense dans $X$ \eqref{TFA3}, $\uoX^\circ$ est intègre et non-vide 
(\cite{ega4} 11.10.5). Soit $v\colon \oy\rightarrow \uoX^\circ$ un point géométrique de $\uoX^\circ$.
On note encore $\oy$ le point géométrique de $\oX^\circ$ et 
$u\colon \oy\rightarrow \uX$ le $X$-morphisme induits par~$v$. On obtient ainsi 
un point $(\oy\rightsquigarrow \ox)$ de $X_\et\gtimes_{X_\et}\oX^\circ_\et$.  
\item[(ii)] Le morphisme d'adjonction $\id\rightarrow \delta_*\delta^*$ \eqref{TFA66a} induit 
un isomorphisme de foncteurs 
\begin{equation}\label{TFA15a}
\varphi_\ox\stackrel{\sim}{\rightarrow}\varphi_\ox\delta_*\delta^*.
\end{equation}
En effet, d'après \ref{tf21}(i), pour tout objet $F$ de $\tE$, le morphisme d'adjonction 
\begin{equation}\label{TFA15b}
F_{\rho(\oy\rightsquigarrow \ox)}\rightarrow(\delta_*(\delta^*(F)))_{\rho(\oy\rightsquigarrow \ox)}
\end{equation}
est un isomorphisme. 
Comme $\uoX^\circ$ est intègre (\cite{agt} III.3.7), la proposition s'ensuit compte tenu de la description \eqref{TFA12e} 
du foncteur fibre en $\rho(\oy\rightsquigarrow \ox)$.
\end{itemize}
\end{remas}

\section{Tour de revêtements associée à une carte logarithmique adéquate}\label{cad}

\subsection{}\label{cad1}
On suppose dans cette section que le morphisme $f\colon (X,\cM_X)\rightarrow (S,\cM_S)$ \eqref{TFA3} admet une carte adéquate 
$((P,\gamma),(\mN,\iota),\vartheta\colon \mN\rightarrow P)$ (\cite{agt} III.4.4), autrement dit qu'il existe une carte $(P,\gamma)$ pour  $(X,\cM_X)$ (\cite{agt} II.5.13), 
une carte $(\mN,\iota)$ pour $(S,\cM_S)$ et un homomorphisme de monoïdes $\vartheta\colon \mN\rightarrow P$ 
tels que les conditions suivantes soient remplies~:  
\begin{itemize}
\item[(i)] Le diagramme d'homomorphismes de monoïdes
\begin{equation}\label{cad1a}
\xymatrix{
P\ar[r]^-(0.5)\gamma&{\Gamma(X,\cM_X)}\\
\mN\ar[r]^-(0.5){\iota}\ar[u]^\vartheta&{\Gamma(S,\cM_S)}\ar[u]_{f^\flat}}
\end{equation}
est commutatif, ou ce qui revient au même (avec les notations de \ref{notconv2}),
le diagramme associé de morphismes de schémas logarithmiques 
\begin{equation}\label{cad1b}
\xymatrix{
{(X,\cM_X)}\ar[r]^-(0.5){\gamma^a}\ar[d]_f&{\bA_P}\ar[d]^{\bA_\vartheta}\\
{(S,\cM_S)}\ar[r]^-(0.5){\iota^a}&{\bA_\mN}}
\end{equation}
est commutatif.
\item[(ii)] Le monoïde $P$ est torique, {\em i.e.}, $P$ est fin et saturé et $P^\gp$, le groupe associé à $P$, 
 est un $\mZ$-module libre (\cite{agt} II.5.1).
\item[(iii)] L'homomorphisme $\vartheta$ est saturé (\cite{agt} II.5.2).
\item[(iv)] L'homomorphisme $\vartheta^\gp\colon \mZ\rightarrow P^\gp$ est injectif, 
le sous-groupe de torsion de $\coker(\vartheta^\gp)$ est d'ordre premier à $p$ et le morphisme de schémas usuels
\begin{equation}\label{cad1c}
X\rightarrow S\times_{\bA_\mN}\bA_P
\end{equation}
déduit de \eqref{cad1b} est étale.  
\item[(v)] Posons $\lambda=\vartheta(1)\in P$, 
\begin{eqnarray}
L&=&\Hom_{\mZ}(P^\gp,\mZ),\label{cad1d}\\
\rH(P)&=&\Hom(P,\mN).\label{cad1e}
\end{eqnarray} 
On notera que $\rH(P)$ est un monoïde fin, saturé et affûté et que l'homomorphisme canonique 
$\rH(P)^\gp\rightarrow \Hom((P^\sharp)^\gp,\mZ)$ est un isomorphisme (\cite{ogus} I 2.2.3). 
On suppose qu'il existe $h_1,\dots,h_r\in \rH(P)$, qui sont $\mZ$-linéairement indépendants dans $L$, tels que  
\begin{equation}\label{cad1f}
\ker(\lambda)\cap \rH(P)=\{\sum_{i=1}^ra_ih_i | \ (a_1,\dots,a_r)\in \mN^r\},
\end{equation}
où l'on considère $\lambda$ comme un homomorphisme $L\rightarrow \mZ$. 
\end{itemize}

\vspace{2mm}
On pose $\pi=\iota(1)$ qui est une uniformisante de $\co_K$. On a alors \eqref{notconv2}
\begin{equation}\label{cad1g}
S\times_{\bA_\mN}\bA_P=\Spec(\co_K[P]/(\pi-e^\lambda)).
\end{equation}

On désigne par $P_\eta$ ou $\lambda^{-1}P$ la localisation de $P$ par $\lambda=\vartheta(1)$ (\cite{ogus} I 1.4.4), 
qui est aussi la somme amalgamée des deux homomorphismes $\vartheta\colon \mN\rightarrow P$
et $\mN\rightarrow \mZ$, de sorte que le diagramme d'homomorphismes canoniques
\begin{equation}\label{cad1h}
\xymatrix{
\mN\ar[r]^\vartheta\ar[d]&P\ar[d]\\
\mZ\ar[r]&{P_\eta}}
\end{equation}
est cocartésien. Comme $P$ est intègre, $P_\eta$ est intègre, et l'homomorphisme canonique $P\rightarrow P_\eta$ est injectif. 
Par ailleurs, il résulte aussitôt des propriétés universelles des localisations de monoïdes et d'anneaux que l'on a 
\begin{equation}\label{cad1i}
\mZ[P_\eta]=\mZ[P][\lambda^{-1}].
\end{equation}
L'élément $e^\lambda$ n'est pas un diviseur de zéro dans $\mZ[P]$ \eqref{notconv2}.

\begin{lem}[\cite{tsuji5} 3.14 et 3.16] \label{cad2}
\
\begin{itemize}
\item[{\rm (i)}] Le groupe $P^\gp/\lambda\mZ$ est libre de type fini.
\item[{\rm (ii)}] Il existe deux entiers $c$ et $d$ tels que $0\leq c\leq d$ et un isomorphisme de monoïdes 
\begin{equation}\label{cad2a}
\mZ^{c+1}\oplus \mN^{d-c}\stackrel{\sim}{\rightarrow} P_\eta
\end{equation}
tel que l'image de $(1,0,\dots,0)$ dans $P_\eta$ soit égale à $\lambda$. 
\end{itemize}
\end{lem}

(i)\ En effet, notant $P^\times$ le groupe des unités de $P$, le conoyau de l'homomorphisme $\mZ\rightarrow P^\gp/P^\times$ défini par $\lambda$ 
est sans torsion d'après (\cite{ogus} I 4.8.11). Par ailleurs, $P^\times$ est contenu dans $P^\gp$ et est donc sans torsion. La suite exacte 
\begin{equation}
0\rightarrow P^\times \rightarrow P^\gp/\lambda \mZ \rightarrow P^\gp/(\lambda\mZ+P^\times)\rightarrow 0
\end{equation}
montre alors que $P^\gp/\lambda \mZ $ est sans torsion.

(ii)\ Soit $F$ la face de $P$ engendrée par $\lambda$, c'est-à-dire l'ensemble des éléments 
$x\in P$ tels qu'il existe $y\in P$ et $n\in \mN$ tels que $x+y=n\lambda$ (\cite{ogus} I 1.4.2).
On note $F^\gp$ le groupe associé à $F$ et 
$P/F$ le conoyau dans la catégorie des monoïdes de l'injection canonique $F\rightarrow P$ (\cite{ogus} I 1.1.6).
Considérons le diagramme commutatif de morphismes de monoïdes 
\begin{equation}
\xymatrix{
&1\ar[r]\ar@{}[dr]|*+[o][F-]{3}& {P/F}\\
\mZ\ar[r]\ar@{}[dr]|*+[o][F-]{1}&F^\gp\ar[r]\ar[u]\ar@{}[dr]|*+[o][F-]{2}&{P_\eta}\ar[u]\\
\mN\ar[r]\ar[u]&F\ar[r]\ar[u]&P\ar[u]}
\end{equation}
Le carré $\xymatrix{\ar@{}|*+[o][F-]{1}}$ est cocartésien. En effet, la somme amalgamée des deux homomorphismes $\mN\rightarrow F$
et $\mN\rightarrow \mZ$ n'est autre que la localisation $\lambda^{-1}F$, c'est-à-dire $F^\gp$. 
Par suite, le carré $\xymatrix{\ar@{}|*+[o][F-]{2}}$ est cocartésien, et il en est alors de même du carré $\xymatrix{\ar@{}|*+[o][F-]{3}}$. 
On en déduit que $P/F$ est le conoyau dans la catégorie des monoïdes de l'homomorphisme canonique $F^\gp\rightarrow P_\eta$.
Par suite, $P/F$ s'identifie au quotient de $P_\eta$ par l'action naturelle de $F^\gp$ (\cite{ogus} I 1.1.5). 
L'homomorphisme $P_\eta\rightarrow P/F$ est donc surjectif. 
Comme $P$ est intègre, l'homomorphisme $F^\gp\rightarrow P_\eta$ est injectif. 
D'après la preuve de (\cite{agt} II.6.3(v)), $P/F$ est un monoïde libre de type fini. 
Choisissons un scindage de $P_\eta\rightarrow P/F$, on obtient un isomorphisme 
\begin{equation}
F^\gp\oplus P/F\stackrel{\sim}{\rightarrow} P_\eta.
\end{equation}
Par ailleurs, $P$ étant intègre,  l'homomorphisme $F^\gp\rightarrow P^\gp$ est injectif, 
et $F^\gp/\lambda \mZ$ est donc libre de type fini d'après (i), d'où la proposition.

\subsection{}\label{cad3}
Pour tout entier $n\geq 1$, on pose
\begin{eqnarray}\label{cad3a}
\co_{K_n}=\co_K[\zeta]/(\zeta^{n}-\pi),
\end{eqnarray}
qui est un anneau de valuation discrète. On note $K_n$ le corps des fractions de $\co_{K_n}$
et $\pi_n$ la classe de $\zeta$ dans $\co_{K_n}$, qui est une uniformisante de $\co_{K_n}$.  
On pose $S^{(n)}=\Spec(\co_{K_n})$
que l'on munit de la structure logarithmique $\cM_{S^{(n)}}$ définie par son point fermé. 
On désigne par $\iota_n\colon \mN\rightarrow \Gamma(S^{(n)},\cM_{S^{(n)}})$
l'homomorphisme défini par $\iota_n(1)=\pi_n$; c'est une carte pour $(S^{(n)},\cM_{S^{(n)}})$.  

Considérons le système inductif de monoïdes $(\mN^{(n)})_{n\geq 1}$, 
indexé par l'ensemble $\mZ_{\geq 1}$ ordonné par la relation de divisibilité, 
défini par $\mN^{(n)}=\mN$ pour tout $n\geq 1$ et dont l'homomorphisme de transition
$\mN^{(n)}\rightarrow \mN^{(mn)}$ (pour $m,n\geq 1$) est l'endomorphisme de Frobenius d'ordre $m$ 
de $\mN$ ({\em i.e.}, l'élévation à la puissance $m$-ième). On notera $\mN^{(1)}$ simplement $\mN$. Les  schémas logarithmiques
$(S^{(n)},\cM_{S^{(n)}})_{n\geq 1}$ forment naturellement un système projectif.
Pour tous entiers $m, n\geq 1$, avec les notations de  \ref{notconv2}, on a un diagramme cartésien de morphismes de schémas logarithmiques
\begin{equation}\label{cad3b}
\xymatrix{
{(S^{(mn)},\cM_{S^{(mn)}})}\ar[r]^-(0.5){\iota^a_{mn}}\ar[d]&{\bA_{\mN^{(mn)}}}\ar[d]\\
{(S^{(n)},\cM_{S^{(n)}})}\ar[r]^-(0.5){\iota^a_n}&{\bA_{\mN^{(n)}}}}
\end{equation}
où $\iota^a_n$ (resp. $\iota^a_{mn}$) est le morphisme associé à $\iota_n$ (resp. $\iota_{mn}$) (\cite{agt} II.5.13). 

\subsection{}\label{cad4}
Considérons le système inductif de monoïdes $(P^{(n)})_{n\geq 1}$, 
indexé par l'ensemble $\mZ_{\geq 1}$ ordonné par la relation de divisibilité, 
défini par $P^{(n)}=P$ pour tout $n\geq 1$ et dont l'homomorphisme de transition
$i_{n,mn}\colon P^{(n)}\rightarrow P^{(mn)}$ (pour $m, n\geq 1$) est l'endomorphisme 
de Frobenius d'ordre $m$ de $P$ ({\em i.e.}, l'élévation à la puissance $m$-ième). Pour tout $n\geq 1$, on note
\begin{equation}\label{cad4a}
P\stackrel{\sim}{\rightarrow}P^{(n)}, \ \ \ t\mapsto t^{(n)},
\end{equation} 
l'isomorphisme canonique. Pour tout $t\in P$ et tous $m, n\geq 1$, on a donc
\begin{equation}\label{cad4b}
i_{n,mn}(t^{(n)})=(t^{(mn)})^m. 
\end{equation}
On notera $P^{(1)}$ simplement $P$. 

On désigne par $P^{(n)}_\eta$ la localisation de $P^{(n)}$ par $\lambda^{(n)}$ (\cite{ogus} I 1.4.4). 
Pour tous $m, n\geq 1$, l'homomorphisme $i_{n,mn}$ induit un homomorphisme que l'on note aussi 
$i_{n,mn}\colon P^{(n)}_\eta\rightarrow P^{(mn)}_\eta$.
L'isomorphisme \eqref{cad4a} induit un isomorphisme 
\begin{equation}\label{cad4c}
P_\eta\stackrel{\sim}{\rightarrow}P^{(n)}_\eta, \ \ \ t\mapsto t^{(n)}.
\end{equation}  
On notera $P^{(1)}_\eta$ simplement $P_\eta$. 

D'après \ref{cad2}(ii), il existe deux entiers $c$ et $d$ tels que $0\leq c\leq d$ et un isomorphisme de monoïdes 
\begin{equation}\label{cad4d}
\mZ^{c+1}\oplus \mN^{d-c}\stackrel{\sim}{\rightarrow} P_\eta
\end{equation}
tel que l'image de $(1,0\dots,0)$ dans $P_\eta$ soit égale à $\lambda$. 
On fixe un tel isomorphisme et pour tout entier $i$ tel que $1\leq i\leq c$ (resp. $c+1\leq i\leq d$), 
on note $\lambda'_i$ l'image de $(\underline{1}_{i+1},0)$ (resp. $(0,\underline{1}_{i-c})$)
dans $P_\eta$, où $\underline{1}_{i+1}$ est l'élément de $\mZ^{c+1}$ dont toutes les composantes sont nulles sauf la $(i+1)$-ième qui vaut $1$
(resp. $\underline{1}_{i-c}$ est l'élément de $\mN^{d-c}$ dont toutes les composantes sont nulles sauf la $(i-c)$-ième qui vaut $1$). 
On pose $\lambda_0=\lambda$. 
Il existe un entier $\alpha\geq 1$ tel que $\lambda_i=\lambda^\alpha\lambda_i'\in P$ pour tout $1\leq i\leq d$.

\begin{lem}[\cite{tsuji5} 4.3]\label{cad5}
Il existe un entier $N\geq 0$ tel que pour tout entier $n\geq 1$, on ait les inclusions 
\begin{equation}\label{cad5a}
\lambda_0^N\mZ[P^{(n)}]\subset \bigoplus_{(\nu_0,\nu_1,\dots,\nu_d)\in (\mN\cap [0,n-1])^{d+1}} \mZ[P] \prod_{0\leq i \leq d}(\lambda_i^{(n)})^{\nu_i}
\subset \mZ[P^{(n)}].
\end{equation}
\end{lem}

En effet, l'homomorphisme canonique $P_\eta\rightarrow P^{(n)}_\eta$  est le composé 
de l'endomorphisme de Frobenius d'ordre $n$ de $P_\eta$ ({\em i.e.}, l'élévation à la puissance $n$-ième) et de l'isomorphisme \eqref{cad4c}.
Par suite, compte tenu de  \eqref{cad4d},
l'algèbre $\mZ[P^{(n)}_\eta]$ est un $\mZ[P_\eta]$-module libre de type fini, de base $\prod_{0\leq i \leq d}(\lambda_i^{(n)})^{\nu_i}$ pour 
$(\nu_0,\nu_1,\dots,\nu_d)\in (\mN\cap [0,n-1])^{d+1}$. 
Comme les homomorphismes canoniques $P\rightarrow P_\eta$ et $P^{(n)}\rightarrow P^{(n)}_\eta$ sont injectifs, 
le morphisme $\mZ[P]$-linéaire canonique 
\begin{equation}
\bigoplus_{(\nu_0,\nu_1,\dots,\nu_d)\in (\mN\cap [0,n-1])^{d+1}} \mZ[P] \prod_{0\leq i \leq d}(\lambda_i^{(n)})^{\nu_i} \rightarrow \mZ[P^{(n)}]
\end{equation}
est injectif. 

Soient $\rho_1,\dots,\rho_r$ des générateurs de $P$. Pour tout $1\leq j\leq r$, il existe $(m_{j0},m_{j1},\dots,m_{jd})\in \mZ^{c+1}\oplus \mN^{d-c}$
tel que 
\begin{equation}
\rho_j=\prod_{0\leq i\leq d}\lambda_i^{m_{ji}}\in P_\eta.
\end{equation}
Pour tout $0\leq i\leq c$, il existe $M_i\geq 1$ tel que $\lambda_0^{M_i}\lambda^{-1}_i\in P$.  
Montrons que 
\begin{equation}
N=\sum_{0\leq i\leq c}M_i\sum_{1\leq j\leq r}|m_{ji}|
\end{equation} 
répond à la question. Le $\mZ[P]$-module $\mZ[P^{(n)}]$ est engendré par les éléments 
\begin{equation}
\rho_{\umu}^{(n)}=\prod_{1\leq j\leq r}(\rho_j^{(n)})^{\mu_j}=\prod_{0\leq i\leq d}(\lambda_i^{(n)})^{\sum_{1\leq j\leq r}m_{ji}\mu_j}
\end{equation}
pour $\umu=(\mu_1,\dots,\mu_r)\in (\mN\cap [0,n-1])^r$. Pour tout entier $0\leq i\leq d$, on pose $\sum_{1\leq j\leq r}m_{ji}\mu_j=\kappa_in +\nu_i$,
où $\kappa_i\in \mZ$ et $\nu_i\in \mN\cap [0,n-1]$. Pour tout $c+1\leq i\leq d$, on a $\kappa_i\in \mN$ et par suite, 
\begin{equation}
(\lambda_i^{(n)})^{\sum_{1\leq j\leq r}m_{ji}\mu_j}=\lambda_i^{\kappa_i}(\lambda_i^{(n)})^{\nu_i}\in P(\lambda_i^{(n)})^{\nu_i}\subset P^{(n)}. 
\end{equation}
Pour tout $0\leq i\leq c$, on a $|\kappa_i|\leq \sum_{1\leq j\leq r}|m_{ji}|$ et par suite, 
\begin{equation}
\lambda_0^{M_i\sum_{1\leq j\leq r}|m_{ji}|}(\lambda_i^{(n)})^{\sum_{1\leq j\leq r}m_{ji}\mu_j}=\lambda_0^{M_i\sum_{1\leq j\leq r}|m_{ji}|}
\lambda_i^{\kappa_i}(\lambda_i^{(n)})^{\nu_i}\in P(\lambda_i^{(n)})^{\nu_i}\subset P^{(n)}.
\end{equation}
Par conséquent, $\lambda_0^N\rho^{(n)}_{\umu}\in P\prod_{0\leq i\leq d}(\lambda_i^{(n)})^{\nu_i}\subset P^{(n)}$.

\subsection{}\label{cad6}
Pour tout entier $n\geq 1$, on pose (avec les notations de \ref{notconv2})
\begin{equation}\label{cad6a}
(X^{(n)},\cM_{X^{(n)}})=(X,\cM_{X})\times_{\bA_P}\bA_{P^{(n)}}.
\end{equation}
On notera que la projection canonique $(X^{(n)},\cM_{X^{(n)}})\rightarrow \bA_{P^{(n)}}$ est stricte. 
Comme le diagramme \eqref{cad3b} est cartésien, il existe un unique morphisme
\begin{equation}\label{cad6b}
f^{(n)}\colon (X^{(n)},\cM_{X^{(n)}})\rightarrow (S^{(n)},\cM_{S^{(n)}}),
\end{equation} 
qui s'insère dans le diagramme commutatif
\begin{equation}\label{cad6c}
\xymatrix{
{(X^{(n)},\cM_{X^{(n)}})}\ar[rrr]\ar[ddd]\ar[rd]_{f^{(n)}}\ar@{}[drrr]|*+[o][F-]{1}&&&
{\bA_{P^{(n)}}}\ar[ld]^{\bA_{\vartheta}}\ar[ddd]\\
&{(S^{(n)},\cM_{S^{(n)}})}\ar[r]^-(0.5){\iota^a_n}\ar[d]&{\bA_{\mN^{(n)}}}\ar[d]&\\
&{(S,\cM_{S})}\ar[r]^-(0.5){\iota^a}&{\bA_\mN}&\\
{(X,\cM_{X})}\ar[rrr]\ar[ru]^f&&&{\bA_P}\ar[lu]_{\bA_\vartheta}}
\end{equation}

\begin{prop}\label{cad7}
Soient $n$ un entier $\geq 1$, $h\colon \oS\rightarrow S^{(n)}$ un $S$-morphisme. Alors,
\begin{itemize}
\item[{\rm (i)}] La face $\xymatrix{\ar@{}|*+[o][F-]{1}}$ du diagramme \eqref{cad6c} est une carte adéquate pour $f^{(n)}$ 
{\rm (\cite{agt} III.4.4)}; en particulier, $f^{(n)}$ est lisse et saturé. 
\item[{\rm (ii)}] Le schéma $X^{(n)}$ est intègre, normal, Cohen-Macaulay et $S^{(n)}$-plat. 
\item[{\rm (iii)}] Le schéma $X^{(n)}\times_{S^{(n)}}\oS$ est normal et localement irréductible {\rm (\cite{agt} III.3.1)}; 
il est donc la somme des schémas induits sur ses composantes irréductibles. 
\item[{\rm (iv)}] Le sous-schéma ouvert maximal de $X^{(n)}$ 
où la structure logarithmique $\cM_{X^{(n)}}$ est triviale est égal à $X^{(n)\circ}=X^{(n)}\times_XX^\circ$ \eqref{TFA3a}. 
\item[{\rm (v)}] Le morphisme $X^{(n)\circ}\rightarrow X^\circ$ est étale, fini et surjectif~; 
de plus $X^{(n)\circ}\times_{S^{(n)}}\oS$ est un espace principal homogène pour la topologie étale 
au-dessus de $\oX^\circ$ sous le groupe 
\[
\Hom_\mZ(P^\gp,\mu_n(\oK)).
\]
\end{itemize}
\end{prop}

En effet, les questions étant locales sur $X$, on peut le supposer affine. Les propositions résultent alors de (\cite{agt} II.6.6).
On notera que l'hypothèse que $X_s$ est non-vide requise dans (\cite{agt} II.6.1) ne sert pas dans les preuves de ces assertions.

\begin{prop}\label{cad8}
Variant l'entier $n\geq 1$ et le $S$-morphisme $h\colon \oS\rightarrow S^{(n)}$,
le nombre de composantes irréductibles de $X^{(n)}\times_{S^{(n)}}\oS$ reste borné. 
\end{prop}

Soient $n$ un entier $\geq 1$, $h\colon \oS\rightarrow S^{(n)}$ un $S$-morphisme. 
Posons $\oX^{(n)}=X^{(n)}\times_{S^{(n)}}\oS$ et $\oX^{(n)\circ}=\oX^{(n)}\times_XX^\circ$.
D'après \eqref{cad6c} et avec les notations de \ref{notconv2} et \ref{cad4}, on a un diagramme cartésien de $S$-morphismes 
\begin{equation}
\xymatrix{
{X^{(n)}}\ar[r]\ar[d]&{\Spec(\co_{K_n}[P^{(n)}]/(\pi_n-e^{\lambda^{(n)}}))}\ar[d]\\
{X}\ar[r]&{\Spec(\co_K[P]/(\pi-e^\lambda))}}
\end{equation}
où les flèches horizontales sont étales. On en déduit un diagramme cartésien de $\oS$-morphismes
\begin{equation}
\xymatrix{
{\oX^{(n)}}\ar[r]\ar[d]&{\Spec(\co_{\oK}[P^{(n)}]/(\pi_n-e^{\lambda^{(n)}}))}\ar[d]\\
{\oX}\ar[r]&{\Spec(\co_{\oK}[P]/(\pi-e^\lambda))}}
\end{equation}
Compte tenu de \ref{cad7}(iv), on en déduit un diagramme cartésien de $\oK$-morphismes
\begin{equation}
\xymatrix{
{\oX^{(n)\circ}}\ar[r]\ar[d]&{\Spec(\oK[(P^{(n)})^\gp]/(\pi_n-e^{\lambda^{(n)}}))}\ar[d]\\
{\oX^\circ}\ar[r]&{\Spec(\oK[P^\gp]/(\pi-e^\lambda))}}
\end{equation}
La flèche verticale de gauche est étale d'après \ref{cad7}(v), ainsi que les flèches horizontales. 
Par suite, les points génériques de $\oX^{(n)\circ}$ s'envoient sur ceux de 
$\Spec(\oK[(P^{(n)})^\gp]/(\pi_n-e^{\lambda^{(n)}}))$ et sur ceux de $\oX^{\circ}$, et donc aussi sur ceux de 
$\Spec(\oK[P^\gp]/(\pi-e^\lambda))$. 

Par ailleurs, comme le $\mZ$-module $P^\gp$ est libre de type fini, pour tout entier $n\geq 1$, on a 
un isomorphisme de $\oK$-algèbres
\begin{equation}
\oK[P^\gp]/(1-e^\lambda)\stackrel{\sim}{\rightarrow} \oK[(P^{(n)})^\gp]/(\pi_n-e^{\lambda^{(n)}}).
\end{equation}
Le nombre de composantes irréductibles de $\Spec(\oK[(P^{(n)})^\gp]/(\pi_n-e^{\lambda^{(n)}}))$ 
(pour $n\geq 1$) est donc constant. Par suite, le nombre de composantes irréductibles de $\oX^{(n)\circ}$ (pour $n\geq 1$) est borné. 
Il en est alors de même de $\oX^{(n)}$ puisque $\oX^{(n)\circ}$ est schématiquement dense dans $\oX^{(n)}$ d'après \ref{cad7}(i)-(iv) et (\cite{agt} III.4.2(vi)).

\subsection{}\label{cad9}
On suppose dans la suite de cette section que $X=\Spec(R)$ est affine.
Pour tout entier $n\geq 1$, le schéma $X^{(n)}$ est affine d'anneau 
\begin{equation}\label{cad9a}
A_n=R\otimes_{\mZ[P]}\mZ[P^{(n)}].
\end{equation}
On désigne par $A_n^\circ$ l'anneau du schéma affine $X^{(n)\circ}=X^{(n)}\times_XX^\circ$ \eqref{TFA3a}. 
Les schémas logarithmiques $(X^{(n)},\cM_{X^{(n)}})_{n\geq 1}$ forment naturellement un système projectif au-dessus 
de $(\bA_{P^{(n)}})_{n\geq 1}$ \eqref{notconv2}, indexé par l'ensemble $\mZ_{\geq 1}$ ordonné par la relation de divisibilité. On pose 
\begin{eqnarray}
A_\infty&=&\underset{\underset{n\geq 1}{\longrightarrow}}{\lim}\ A_n,\label{cad9b}\\
A^\circ_\infty&=&\underset{\underset{n\geq 1}{\longrightarrow}}{\lim}\ A^\circ_n.\label{cad9c}
\end{eqnarray}

Etant donné une suite $(u_n)_{n\geq 0}$ d'entiers $\geq 1$ telle que $u_n$ divise $u_{n+1}$, on pose 
\begin{eqnarray}
A_{u_\infty}&=&\underset{\underset{n\geq 0}{\longrightarrow}}{\lim}\ A_{u_n}, \label{cad9d}\\
A^\circ_{u_\infty}&=&\underset{\underset{n\geq 0}{\longrightarrow}}{\lim}\ A^\circ_{u_n}.\label{cad9e}
\end{eqnarray}
Avec les notations de \ref{cad3}, on pose 
\begin{eqnarray}
K_{u_\infty}&=&\underset{\underset{n\geq 0}{\longrightarrow}}{\lim}\ K_{u_n}, \label{cad9f}\\
\co_{K_{u_\infty}}&=&\underset{\underset{n\geq 0}{\longrightarrow}}{\lim}\ \co_{K_{u_n}},\label{cad9g}
\end{eqnarray}
de sorte que $\co_{K_{u_\infty}}$ est la fermeture intégrale de $\co_K$ dans $K_{u_\infty}$.
Si la suite $(u_n)_{n\geq 0}$ est non-bornée, $\co_{K_{u_\infty}}$ est un anneau de valuation non-discrète de hauteur $1$. 
On peut alors considérer les notions de $\alpha$-algèbre (\ref{alpha}, \ref{finita}, \ref{aet}...) sur les $\co_{K_{u_\infty}}$-algèbres, en particulier 
sur $A_{u_\infty}$.

Si $u_n=n!$ $(n\geq 0)$, alors $A_{u_\infty}=A_\infty$; on notera $K_{u_\infty}$ et $\co_{K_{u_\infty}}$ par $K_\infty$ et $\co_{K_\infty}$. 
On considérera aussi le cas où $u_n=p^n$ $(n\geq 0)$, auquel cas on notera $A_{u_\infty}$, $K_{u_\infty}$  et $\co_{K_{u_\infty}}$ par $A_{p^\infty}$, 
$K_{p^\infty}$ et $\co_{K_{p^\infty}}$.

\begin{lem}\label{cad100}
\
\begin{itemize}
\item[{\rm (i)}] Pour tous entiers $m,n\geq 1$, le morphisme canonique 
\begin{equation}
\bigoplus_{(\nu_0,\nu_1,\dots,\nu_d)\in (\mN\cap [0,m-1])^{d+1}} A_n \prod_{0\leq i \leq d}(t_i^{(mn)})^{\nu_i}
\rightarrow A_{mn},
\end{equation}
où $t_i^{(mn)}$ désigne l'image canonique de $\lambda_i^{(mn)}\in P^{(mn)}$ dans $A_{mn}$ \eqref{cad4}, est injectif. 
\item[{\rm (ii)}] Il existe un entier $N>0$ tel que tous entiers $m,n\geq 1$, on ait
\begin{equation}
\pi^N_nA_{mn}\subset \bigoplus_{(\nu_0,\nu_1,\dots,\nu_d)\in (\mN\cap [0,m-1])^{d+1}} 
A_n \prod_{0\leq i \leq d}(t_i^{(mn)})^{\nu_i} \subset A_{mn}.
\end{equation}
\end{itemize}
\end{lem}

(i) Comme le morphisme $X\rightarrow S\otimes_{\bA_\mN}\bA_P$ \eqref{cad1c} est étale, on peut se réduire au cas où 
$R=\co_K[P]/(\pi-e^\lambda)$ \eqref{cad1g}. Pour tout entier $n\geq 1$, on a alors 
\begin{equation}
A_n=\co_K[P^{(n)}]/(\pi-e^\lambda).
\end{equation}
Compte tenu de \eqref{cad1i}, on en déduit que 
\begin{equation}
A_n[\pi_n^{-1}]=K[P^{(n)}_\eta]/(\pi-e^\lambda).
\end{equation}
Il résulte de \eqref{cad2a} que pour tous entiers $m,n\geq 1$, le morphisme $K[P^{(n)}_\eta]$-linéaire canonique 
\begin{equation}
\bigoplus_{(\nu_0,\nu_1,\dots,\nu_d)\in (\mN\cap [0,m-1])^{d+1}}  K[P^{(n)}_\eta]/(\pi-e^\lambda) \prod_{0\leq i \leq d}(t_i^{(mn)})^{\nu_i}
\rightarrow K[P^{(mn)}_\eta]/(\pi-e^\lambda)
\end{equation}
est bijectif. Comme $\pi_n$ n'est pas diviseur de zéro dans $A_n$ d'après \ref{cad7}(ii), 
l'homomorphisme canonique $A_n\rightarrow A_n[\pi_n^{-1}]$ est injectif, d'où la proposition. 

(ii) Cela résulte aussitôt de (i) et \ref{cad5}.

\begin{prop}\label{cad10}
Il existe un entier $N\geq 1$ tel que les propriétés suivantes soient remplies~:
\begin{itemize}
\item[{\rm (i)}]  Pour tous entiers $n,q\geq 1$ et tout $A_n$-module $M$, on a
\begin{equation}\label{cad10d}
\pi_{n}^N\Tor_q^{A_n}(A_\infty,M)=0,
\end{equation}
où $\pi_n$ est l'uniformisante de $K_n$ définie dans \ref{cad3}. 
\item[{\rm (ii)}] Pour toute suite $(u_n)_{n\geq 0}$ d'entiers $\geq 1$ telle que $u_n$ divise $u_{n+1}$, 
tous entiers $m\geq 0$ et $q\geq 1$ et tout $A_{u_m}$-module $M$, on a
\begin{equation}\label{cad10a}
\pi_{u_m}^N\Tor_q^{A_{u_m}}(A_{u_\infty},M)=0.
\end{equation}
\end{itemize}
\end{prop}

En effet, d'après \ref{cad100}(ii), il existe un entier $N\geq 1$ tel que pour tous entiers $m,n,q\geq 1$ et tout $A_n$-module $M$, on ait 
\begin{equation}\label{cad10c}
\pi^N_n\Tor_q^{A_n}(A_{mn},M)=0.
\end{equation}

\begin{cor}\label{cad11}
Il existe un entier $N\geq 1$ tel que les propriétés suivantes soient remplies:
\begin{itemize}
\item[{\rm (i)}] Pour tout entier $n\geq1$ et toute suite exacte de $A_n$-modules
$M'\rightarrow M\rightarrow M''$, la suite 
\begin{equation}
M'\otimes_{A_n}A_\infty\rightarrow M\otimes_{A_n}A_\infty\rightarrow M''\otimes_{A_n}A_\infty
\end{equation}
est $\pi^N_n$-exacte \eqref{finita3}.
\item[{\rm (ii)}] Pour toute suite $(u_n)_{n\geq 0}$ d'entiers $\geq 1$ telle que $u_n$ divise $u_{n+1}$, 
tout entier $n\geq 0$ et toute suite exacte de $A_{u_n}$-modules $M'\rightarrow M\rightarrow M''$, la suite 
\begin{equation}
M'\otimes_{A_{u_n}}A_{u_\infty}\rightarrow M\otimes_{A_{u_n}}A_{u_\infty}\rightarrow M''\otimes_{A_{u_n}}A_{u_\infty}
\end{equation}
est $\pi^N_{u_n}$-exacte.
\end{itemize}
\end{cor}

En effet, d'après \ref{cad10}, il existe un entier $N\geq 1$ tel que les propriétés suivantes soient remplies:
\begin{itemize}
\item[{\rm (i)}] Pour tout entier $n\geq 1$ et tout morphisme injectif de $A_n$-modules
$\nu$, le noyau de $\nu\otimes_{A_n}A_\infty$ est annulé par $\pi^N_n$.
\item[{\rm (ii)}] Pour toute suite $(u_n)_{n\geq 0}$ d'entiers $\geq 1$ telle que $u_n$ divise $u_{n+1}$, tout entier $m\geq 0$ et 
tout morphisme injectif de $A_{u_m}$-modules $\nu$, le noyau de $\nu\otimes_{A_{u_m}}A_{u_\infty}$ est annulé par $\pi^N_{u_m}$.
\end{itemize}

\begin{prop}\label{cad12}
Pour toute suite non bornée $(u_n)_{n\geq 0}$ d'entiers $\geq 1$ telle que $u_n$ divise $u_{n+1}$,
l'anneau $A_{u_\infty}$ est universellement $\alpha$-cohérent \eqref{afini14}. En particulier, $A_\infty$ est universellement $\alpha$-cohérent.
\end{prop} 

Soit $N$ l'entier $\geq 1$ de la proposition \ref{cad11} et soit $g$ un entier $\geq 0$. 
Posons $H_{\infty}=A_{u_\infty}[t_1,\dots,t_g]$ et pour tout entier $n\geq 0$, 
$H_n=A_{u_n}[t_1,\dots,t_g]$.  Considérons un entier $a\geq 1$ et une forme $H_\infty$-linéaire 
$\varphi_\infty \colon H_\infty^a\rightarrow H_\infty$.
Il existe un entier $n\geq 0$ et une forme $H_n$-linéaire $\varphi_n\colon H_n^a\rightarrow H_n$ tels que $\varphi_\infty=\varphi_n\otimes_{A_{u_n}}A_{u_\infty}$. 
L'anneau $A_{u_n}$ étant noethérien, il existe un entier $b\geq 0$ et une suite exacte de $H_n$-modules 
\begin{equation}\label{cad12b}
H_n^b\longrightarrow H_n^a\stackrel{\varphi_n}{\longrightarrow} H_n.
\end{equation}
La suite induite par extension des scalaires de $A_{u_n}$ à $A_{u_\infty}$,
\begin{equation}\label{cad12c}
H_\infty^b\longrightarrow H_\infty^a\stackrel{\varphi_\infty}{\longrightarrow} H_\infty
\end{equation}
est $\pi^N_{u_n}$-exacte en vertu de \ref{cad11}. Prenant $n$ de plus en plus grand, on en déduit que
$\ker(\varphi_\infty)$ est de type $\alpha$-fini sur $H_\infty$, d'où la proposition.

\subsection{}\label{cad13}
Soit $\oy$ un point géométrique de $\oX^\circ$ que l'on suppose fixé dans la suite de cette section. 
Pour tous entiers $m,n\geq 1$, le morphisme canonique $X^{(mn)}\rightarrow X^{(n)}$ 
est fini et surjectif. D'après (\cite{ega4} 8.3.8(i)), il existe alors un $X$-morphisme 
\begin{equation}\label{cad13a}
\oy\rightarrow \underset{\underset{n\geq 1}{\longleftarrow}}{\lim}\ X^{(n)},
\end{equation} 
la limite inductive étant indexée par l'ensemble $\mZ_{\geq 1}$ ordonné par la relation de divisibilité.
On fixe un tel morphisme dans toute la suite de cette section. Celui-ci induit un $S$-morphisme
\begin{equation}\label{cad13b}
\oS\rightarrow \underset{\underset{n\geq 1}{\longleftarrow}}{\lim}\ S^{(n)},
\end{equation}
et par suite un $\co_K$-homomorphisme $\co_{K_\infty}\rightarrow \co_\oK$ (cf. \ref{cad9}). Pour tout entier $n\geq 1$, on pose
\begin{equation}\label{cad13c}
\oX^{(n)}= X^{(n)}\times_{S^{(n)}}\oS \ \ \ {\rm et}\ \ \ \oX^{(n)\circ}=\oX^{(n)}\times_XX^\circ.
\end{equation} 
On en déduit un $\oX$-morphisme 
\begin{equation}\label{cad13d}
\oy\rightarrow \underset{\underset{n\geq 1}{\longleftarrow}}{\lim}\ \oX^{(n)}.
\end{equation} 

Pour tout entier $n\geq 1$, le schéma $\oX^{(n)}$ étant normal et localement irréductible d'après \ref{cad7}(iii), 
il est la somme des schémas induits sur ses composantes irréductibles. 
On note $\oX^{(n)\star}$ la composante irréductible de $\oX^{(n)}$ contenant l'image de $\oy$ \eqref{cad13d}. 
On pose 
\begin{equation}\label{cad13f}
R_n=\Gamma(\oX^{(n)\star},\co_{\oX^{(n)}}).
\end{equation}
Ces anneaux forment naturellement un système inductif. On pose 
\begin{eqnarray}
R_\infty&=&\underset{\underset{n\geq 1}{\longrightarrow}}{\lim}\ R_n,\label{cad13g}\\
R_{p^\infty}&=&\underset{\underset{n\geq 0}{\longrightarrow}}{\lim}\ R_{p^n}.\label{cad13gg}
\end{eqnarray}
Ce sont des anneaux normaux et intègres d'après (\cite{ega1n} 0.6.1.6(i) et 0.6.5.12(ii)).

\begin{lem}\label{cad14}
On désigne par $\star$ l'un des deux symboles $p^\infty$ ou $\infty$,  
et pour toute extension algébrique $L$ de $K$, par $\co_L$ la clôture intégrale de $\co_K$ dans $L$. Alors,
\begin{itemize}
\item[{\rm (i)}] Pour toute extension algébrique $L$ de $K_\star$ \eqref{cad9}, 
le schéma $\Spec(A_\star\otimes_{\co_{K_\star}}\co_L)$ est normal et localement irréductible. 
\item[{\rm (ii)}] Le schéma $\Spec(R_\star)$ est une composante connexe de $\Spec(A_\star\otimes_{\co_{K_\star}}\co_\oK)$.
\end{itemize}
\end{lem}
(i) Considérons le cas $\star=\infty$, l'autre cas étant similaire. 
On a un isomorphisme canonique (\cite{ega4} 8.2.3)
\begin{equation}
\Spec(A_\infty\otimes_{\co_{K_\infty}}\co_L)\stackrel{\sim}{\rightarrow}
\underset{\underset{n\geq 1}{\longleftarrow}}{\lim}\ X^{(n)}\times_{S^{(n)}}\Spec(\co_L),
\end{equation}
Pour tout entier $n\geq 1$, on démontre, en calquant la preuve de (\cite{agt} II.6.6(iii)), que le schéma 
$X^{(n)}\times_{S^{(n)}}\Spec(\co_L)$ est normal et localement irréductible.
L'immersion $X^{(n)\circ}\rightarrow X^{(n)}$ est schématiquement dominante d'après \ref{cad7}(i)-(iv) et (\cite{agt} III.4.2(vi)). 
Pour tout entier $m\geq 1$, le morphisme canonique $X^{(mn)\circ}\rightarrow X^{(n)\circ}$ est étale, fini et surjectif en vertu de  \ref{cad7}(v).
Toute composante irréductible de $X^{(mn)}\times_{S^{(n)}}\Spec(\co_L)$ domine donc une composante irréductible de $X^{(n)}\times_{S^{(n)}}\Spec(\co_L)$. 
En vertu de \ref{cad8}, il existe un entier $n\geq 1$ tel que pour toute composante 
irréductible $C$ de $X^{(n)}\times_{S^{(n)}}\Spec(\co_L)$ et tout entier $m\geq 1$, $C\times _{X^{(n)}}X^{(mn)}$ soit irréductible. 
La proposition résulte alors de (\cite{ega1n} 0.6.1.6 et 0.6.5.12(ii)).

(ii) Cela résulte aussitôt de la preuve de (i).

\begin{prop}\label{cad15}
Les $\co_{\oK}$-algèbres $R_\infty$ et $R_{p^\infty}$ sont universellement $\alpha$-cohérentes \eqref{afini14}.
\end{prop}

En effet, d'après (\cite{tate} Prop.~9), pour toute extension finie $L$ de $K_\infty$ (resp. $K_{p^\infty}$), la clôture intégrale 
de $\co_{K}$ dans $L$ est $\alpha$-finie-étale sur $\co_{K_\infty}$ (resp. $\co_{K_{p^\infty}}$) \eqref{aet2}. 
Compte tenu de \ref{cad14}, $R_\infty$ est donc la limite inductive 
d'un système inductif filtrant de $A_\infty$-algèbres dont chacune est $\alpha$-finie-étale (\cite{agt} V.7.4). 
De même, $R_{p^\infty}$ est la limite inductive d'un système inductif filtrant de $A_{p^\infty}$-algèbres dont chacune est $\alpha$-finie-étale.
La proposition s'ensuit en vertu de \ref{aet9} et \ref{cad12}. 

\begin{prop}\label{cad16}
Il existe un entier $N\geq 1$ tel que pour tous entiers $n,q\geq 1$ et tout $A_n$-module $M$, on ait
\begin{equation}\label{cad16a}
\pi_{n}^N\Tor_q^{A_n}(R_\infty,M)=0,
\end{equation}
où $\pi_n$ est l'uniformisante de $K_n$ définie dans \ref{cad3}. 
\end{prop}

En effet, on a un isomorphisme canonique 
\begin{equation}
(A_\infty\otimes_{A_n}M)\otimes_{\co_{K_\infty}}\co_\oK\stackrel{\sim}{\rightarrow}
(A_\infty\otimes_{\co_{K_\infty}}\co_\oK)\otimes_{A_n}M.
\end{equation}
Comme $\co_\oK$ est $\co_{K_\infty}$-plat, on en déduit un isomorphisme 
\begin{equation}
\Tor_q^{A_n}(A_\infty,M)\otimes_{\co_{K_\infty}}\co_\oK\stackrel{\sim}{\rightarrow}
\Tor_q^{A_n}(A_\infty\otimes_{\co_{K_\infty}}\co_\oK,M). 
\end{equation}
La proposition résulte alors de \ref{cad10}(i) et \ref{cad14}. 

\begin{cor}\label{cad17}
Il existe un entier $N\geq 1$ tel que pour tout entier $n\geq1$ et toute suite exacte de $A_n$-modules
$M'\rightarrow M\rightarrow M''$, la suite 
\begin{equation}
M'\otimes_{A_n}R_\infty\rightarrow M\otimes_{A_n}R_\infty\rightarrow M''\otimes_{A_n}R_\infty
\end{equation}
soit $\pi^N_n$-exacte \eqref{finita3}.
\end{cor}
Cela résulte aussitôt de \ref{cad16}.

\section{\texorpdfstring{Schémas $K(\pi,1)$, suite}{Schémas K(pi,1), suite}}\label{Kp1s}

\subsection{}\label{Kp1s1}
Les hypothèses et notations de \ref{cad1} sont en vigueur dans cette section. On reprend les notations de \ref{cad6}. 
Pour tout entier $n\geq 1$ et tout $S$-morphisme $\tau\colon \oS\rightarrow S^{(n)}$, on pose \eqref{TFA3a}
\begin{equation}\label{Kp1s1a}
Y^{(n)}_\tau=X^{(n)}\times_{S^{(n)}}\oeta \ \ \ {\rm et}\ \ \ Y^{(n)\circ}_\tau=Y^{(n)}_\tau\times_XX^\circ.
\end{equation}
On laissera tomber l'indice $\tau$ lorsqu'il n'y a aucun risque d'ambiguïté. 
On a $Y^{(1)}=X_\oeta$ et $Y^{(1)\circ}=\oX^\circ$. 
On note $g\colon \oX^\circ\rightarrow X_\oeta$ et $g_n\colon Y^{(n)\circ}\rightarrow Y^{(n)}$ les injections canoniques, et
$\varphi_n\colon Y^{(n)}\rightarrow X_\oeta$ et $\varphi_n^\circ\colon Y^{(n)\circ}\rightarrow \oX^\circ$ les morphismes canoniques.
Pour tout faisceau $\cF$ de $\oX^\circ_\et$ et tout entier $q\geq 0$, 
on a un morphisme de changement de base
\begin{equation}\label{Kp1s1c}
\varphi_n^*(\rR^q g_*\cF)\rightarrow \rR^q g_{n*}(\varphi_n^{\circ *}(\cF)).
\end{equation}

\begin{lem}\label{Kp1s2}
Il existe un entier $r\geq 0$ et un morphisme lisse $t\colon X_\eta\rightarrow \mA_\eta^r$ vérifiant les propriétés suivantes~:
\begin{itemize}
\item[{\rm (i)}] La structure logarithmique $\cM_X|X_\eta$ est définie par le diviseur à croisement normaux $D$ sur $X_\eta$ image inverse 
du diviseur des coordonnées sur $\mA_\eta^r$. En particulier, 
$X^\circ=t^{-1}(\mG_\eta^r)$, où $\mG_\eta^r$ est l'ouvert de $\mA_\eta^r$ où les coordonnées ne s'annulent pas.
\item[{\rm (ii)}] Pour tout entier $n\geq 1$, notons $Z^{(n)}$ le schéma défini par le diagramme cartésien
\begin{equation}
\xymatrix{
Z^{(n)}\ar[r]\ar[d]&X_\eta\ar[d]^t\\
{\mA_\eta^r}\ar[r]&{\mA_\eta^r}}
\end{equation}
où la flèche horizontale inférieure est définie par l'élévation à la puissance $n$-ième des coordonnées de $\mA_\eta^r$. 
Alors, pour tout $S$-morphisme $\tau\colon \oS\rightarrow S^{(n)}$, avec les notations de \ref{Kp1s1}, il existe un $X_\oeta$-morphisme étale et fini 
\begin{equation}
Y^{(n)}_\tau\rightarrow Z^{(n)}_\oeta.
\end{equation}
\end{itemize}
\end{lem}

D'après (\cite{agt} II.6.3(v)), le schéma (usuel) $X_\eta$ est lisse sur $\eta$, $X^\circ$ est l'ouvert 
complémentaire dans $X_\eta$ d'un diviseur à croisements normaux stricts $D$ et  $\cM_X|X_\eta$ 
est la structure logarithmique sur $X_\eta$ définie par $D$. Rappelons brièvement la démonstration de {\em loc. cit.} 
Soit $F$ la face de $P$ engendrée par $\lambda$, c'est-à-dire l'ensemble des éléments 
$x\in P$ tels qu'il existe $y\in P$ et $m\in \mN$ tels que $x+y=m\lambda$ (cf. \cite{ogus} I 1.4.2).
On note $F^{-1}P$ la localisation de $P$ par $F$ (\cite{ogus} I 1.4.4), $Q$ le sous-groupe des unités de $F^{-1}P$,
$(F^{-1}P)^\sharp$ le quotient de $F^{-1}P$ par $Q$ (\cite{agt} II.5.1) et $P/F$ le conoyau dans la catégorie des monoïdes de l'injection canonique 
$F\rightarrow P$ (\cite{ogus} I 1.1.5). Il résulte aussitôt des propriétés universelles des localisations de monoïdes et d'anneaux que 
l'homomorphisme canonique $\mZ[P]\rightarrow \mZ[F^{-1}P]$ induit un isomorphisme
\begin{equation}\label{Kp1s2b}
\mZ[P]_{\lambda}\stackrel{\sim}{\rightarrow} \mZ[F^{-1}P],
\end{equation}
de sorte qu'on a un diagramme cartésien
\begin{equation}\label{Kp1s2c}
\xymatrix{
{X_\eta}\ar[r]\ar[d]&{\bA_{F^{-1}P}}\ar[d]\\
X\ar[r]&{\bA_P}}
\end{equation} 
Par ailleurs, on a un isomorphisme canonique 
\begin{equation}\label{Kp1s2d}
P/F\stackrel{\sim}{\rightarrow}(F^{-1}P)^\sharp.
\end{equation} 
La preuve de (\cite{agt} II.6.3(v)) montre que $P/F$ est un monoïde libre de type fini. On en déduit un isomorphisme 
\begin{equation}\label{Kp1s2e}
Q\times (P/F) \stackrel{\sim}{\rightarrow} F^{-1}P.
\end{equation}

(i) Comme le monoïde $P/F$ est libre de type fini, il existe un entier $r\geq 0$ et un isomorphisme $\Spec(K[P/F])\simeq \mA_\eta^r$.
Les diagrammes \eqref{cad1b} et \eqref{Kp1s2c} et l'isomorphisme \eqref{Kp1s2e} induisent donc un morphisme \eqref{notconv2}
\begin{equation}\label{Kp1s2i}
X_\eta\rightarrow \mA_\eta^r\times_\eta\Spec(K[Q]/(\pi-e^\lambda)),
\end{equation}
et par suite un morphisme $t\colon X_\eta\rightarrow \mA_\eta^r$. La structure logarithmique $\cM_X|X_\eta$ sur $X_\eta$ est l'image inverse par $t$
de la structure logarithmique sur $\mA_\eta^r$ définie par les axes des coordonnées (\cite{agt} II.5.10).
Compte tenu de \ref{cad1}(iv), le morphisme \eqref{Kp1s2i} est étale et le morphisme $t$ est lisse; d'où la proposition.

(ii) Notant $u$ l'endomorphisme de Frobenius d'ordre $n$ de $P$ ({\em i.e.}, l'élévation à la puissance $n$-ième),
le diagramme \eqref{cad6c}  induit un diagramme commutatif à carrés cartésiens 
\begin{equation}
\xymatrix{
{X^{(n)}}\ar[r]\ar[d]&{S^{(n)}\times_{\bA_\mN}\bA_P}\ar[r]\ar[d]&{\bA_P}\ar[d]^{\bA_{u}}\\
X\ar[r]&{S\times_{\bA_\mN}\bA_P}\ar[r]&{\bA_P}}
\end{equation}
On en déduit un diagramme cartésien \eqref{notconv2}
\begin{equation}
\xymatrix{
{X^{(n)}\times_{S^{(n)},\tau}\oS}\ar[r]\ar[d]&{\Spec(\co_{\oK}[P]/(\pi_n-e^\lambda))}\ar[d]\\
{\oX}\ar[r]&{\Spec(\co_{\oK}[P]/(\pi-e^\lambda))}}
\end{equation}
et par suite, compte tenu de \eqref{Kp1s2b}, un diagramme cartésien
\begin{equation}
\xymatrix{
{Y^{(n)}_\tau}\ar[r]\ar[d]_{\varphi_n}&{\Spec(\oK[F^{-1}P]/(\pi_n-e^\lambda))}\ar[d]\\
{X_\oeta}\ar[r]&{\Spec(\oK[F^{-1}P]/(\pi-e^\lambda))}}
\end{equation}
où les flèches verticales de droite sont induites par les homomorphismes de Frobenius d'ordre $n$ de $P$ et $F^{-1}P$. 
On notera que le morphisme 
\begin{equation}
\Spec(\oK[Q]/(\pi_n-e^\lambda))\rightarrow \Spec(\oK[Q]/(\pi-e^\lambda))
\end{equation}
induit par l'homomorphisme de Frobenius d'ordre $n$ de $Q$, est étale et fini. 
Compte tenu de \eqref{Kp1s2e} et des définitions, on en déduit un $X_\oeta$-morphisme étale et fini $Y^{(n)}_\tau\rightarrow Z^{(n)}_\oeta$.

\begin{lem}\label{Kp1s3}
Soient $n$ un entier $\geq 1$, $F$ un faisceau de $(\mZ/n\mZ)$-modules 
localement constant et constructible de $X_{\oeta,\et}$, $\tau\colon \oS\rightarrow S^{(n)}$ un $S$-morphisme, $q$ un entier $\geq 1$. 
Alors, avec les notations de \ref{Kp1s1}, le morphisme de changement de base \eqref{Kp1s1c}
\begin{equation}
\varphi_n^*(\rR^q g_*(g^*F))\rightarrow \rR^q g_{n*}(g_n^*(\varphi_n^*(F)))
\end{equation}
est nul.
\end{lem}

Cela résulte de \ref{Kp1s2} et \ref{Kpun14}.

\begin{lem}\label{Kp1s4}
Pour tout faisceau abélien de torsion localement constant et constructible $F$ de $\oX^\circ_\et$, il existe un entier $n\geq 1$ 
tel que pour tout $S$-morphisme $\tau\colon \oS\rightarrow S^{(n)}$, avec les notations de \ref{Kp1s1}, $\varphi^{\circ*}_n(F)$ se prolonge en un 
faisceau abélien de torsion localement constant et constructible de $Y^{(n)}_\et$. 
\end{lem}

Par descente, $F$ est représentable par un schéma en groupes abéliens étale et fini $V$ au-dessus de $\oX^\circ$.  
Notons $T$ la fermeture intégrale de $X_\oeta$ dans $V$. Soit $n$ un entier $\geq 1$ et divisible par tous les indices de ramification
du morphisme $T\rightarrow X_\oeta$ en les points de codimension un de $T$. Reprenons les notations de \ref{Kp1s2} et notons 
$T^{(n)}$ la fermeture intégrale de $Z_\oeta^{(n)}$ dans $V\times_{X_\oeta}Z^{(n)}_\oeta$. On observera que $Z^{(n)}$ est lisse sur $\eta$.  
D'après le lemme d'Abhyankar (\cite{sga1} X 3.6), $T^{(n)}\rightarrow Z_\oeta^{(n)}$ est étale au-dessus des points de codimension un de $Z_\oeta^{(n)}$. 
Par suite, le morphisme canonique $T^{(n)}\rightarrow Z_\oeta^{(n)}$ est étale en vertu du théorème de pureté de Zariski-Nagata (\cite{sga2} X 3.4). 
Par suite, $T^{(n)}\times_{Z_\oeta^{(n)}}T^{(n)}$ est la fermeture intégrale de $Z_\oeta^{(n)}$ dans 
$(V\times_{\oX^\circ}V)\times_{X_\oeta}Z^{(n)}_\oeta$. On en déduit une structure de schéma en groupes abéliens sur $T^{(n)}$ au-dessus de $Z_\oeta^{(n)}$
qui prolonge celle de $V\times_{X_\oeta}Z^{(n)}_\oeta$ au-dessus de $\oX^\circ\times_{X_\oeta}Z^{(n)}_\oeta$. Le faisceau 
$F|(\oX^\circ\times_{X_\oeta}Z^{(n)}_\oeta)$ se prolonge donc en un faisceau abélien de torsion localement constant et constructible de $(Z^{(n)}_\oeta)_\et$.
La proposition s'ensuit compte tenu de \ref{Kp1s2}(ii).

\begin{prop}\label{Kp1s5}
Soient $F$ un faisceau abélien de torsion, localement constant et constructible de $\oX^\circ_\et$, $V\rightarrow \oX^\circ$ un revêtement étale, 
$q$ un entier $\geq 1$, $\xi\in \rH^q(V,F)$, $\ox$ un point géométrique de $X$.
Alors, il existe un voisinage étale $U$ de $\ox$ dans $X$ et un revêtement étale surjectif $W\rightarrow V\times_{\oX^\circ}\oU^\circ$ tel que l'image canonique
de $\xi$ dans $\rH^q(W,F)$ soit nulle. 
\end{prop}

Par la preuve de l'implication (iv)$\Rightarrow$(v) de \ref{Kpun1}, on peut se borner au cas où $V=\oX^\circ$. 
En vertu de \ref{Kp1s4}, il existe un entier $n\geq 1$ tel que pour tout $S$-morphisme $\tau\colon \oS\rightarrow S^{(n)}$, 
avec les notations de \ref{Kp1s1}, $F|Y^{(n)\circ}_\tau$ se prolonge en un faisceau abélien de torsion localement constant et constructible de $Y^{(n)}_{\tau,\et}$.  
D'après  \ref{cad7}(i), le morphisme $f_n\colon (X^{(n)},\cM_{X^{(n)}})\rightarrow (S^{(n)},\cM_{S^{(n)}})$ est adéquat, 
et il vérifie les conditions de \ref{cad1}. Supposons que pour tout $S$-morphisme $\tau\colon \oS\rightarrow S^{(n)}$, 
la proposition soit démontrée pour l'image de $\xi$ dans  $\rH^q(Y^{(n)\circ}_\tau,F)$ et tout point
géométrique de $X^{(n)}$ au-dessus de $\ox$. Il existe donc un morphisme étale $U^{(n)}_\tau\rightarrow X^{(n)}$ tel que  le morphisme
\begin{equation}
U^{(n)}_\tau\otimes_X\kappa(\ox)\rightarrow X^{(n)}\otimes_X\kappa(\ox)
\end{equation}
soit surjectif, et un revêtement étale surjectif $W_\tau\rightarrow U^{(n)\circ}_\tau\times_{S^{(n)},\tau}\oS$ tel que l'image canonique de $\xi$ dans 
$\rH^q(W_\tau,F)$ soit nulle.
D'après \ref{Kpun15}, il existe un $X$-schéma étale $\ox$-pointé $U$ et pour tout $S$-morphisme $\tau\colon \oS\rightarrow S^{(n)}$,
un $X^{(n)}$-morphisme $U\times_XX^{(n)}\rightarrow U^{(n)}_\tau$. Considérons le diagramme commutatif à carrés cartésiens 
\begin{equation}
\xymatrix{
{W_\tau}\ar[d]&{W_\tau\times_{U^{(n)}_\tau}(U\times_XX^{(n)})}\ar[d]^{v_\tau}\ar[l]&\\
{U^{(n)\circ}_\tau\times_{S^{(n)},\tau}\oS}\ar[d]&{(U^\circ\times_XX^{(n)})\times_{S^{(n)},\tau}\oS}\ar[l]\ar[r]^-(0.5){u_\tau}\ar[d]&{\oU^\circ}\\
{U^{(n)}_\tau}&{U\times_XX^{(n)}}\ar[l]&}
\end{equation}
Le couple formé de $U$ et du revêtement étale surjectif de $\oU^\circ$ somme disjointe des $u_\tau\circ v_\tau$ pour tous les $\tau$, 
répond alors à la question. On peut donc se borner au cas où $F$ se prolonge en un faisceau abélien de torsion, 
localement constant et constructible $G$ de $X_{\oeta,\et}$. L'isomorphisme $g^*(G)\stackrel{\sim}{\rightarrow} F$ \eqref{Kp1s1}
induit par adjonction un morphisme $G\rightarrow g_*(F)$ qui est en fait un isomorphisme; en particulier, $g_*(F)$ est localement constant et constructible.
En effet, la question étant locale pour la topologie étale de $X$, on peut 
supposer $G$ constant, auquel cas l'assertion recherchée résulte de (\cite{sga4} IX 2.14.1). 

On désigne par $g\colon \oX^\circ\rightarrow X_\oeta$ l'injection canonique. Considérons la suite spectrale de Cartan-Leray
\begin{equation}
E_2^{a,b}=\rH^a(X_\oeta,\rR^bg_*(F))\Rightarrow \rH^{a+b}(\oX^\circ,F),
\end{equation} 
et notons $(E^q_i)_{0\leq i\leq q}$ la filtration aboutissement sur $\rH^q(\oX^\circ,F)$, de sorte que 
\begin{equation}
E_i^q/E_{i+1}^q=E_\infty^{i,q-i}
\end{equation} 
pour tout $0\leq i\leq q$. Il existe un entier $0\leq i_\xi\leq q$ tel que $\xi\in E_{i_\xi}^q\backslash E_{i_\xi+1}^q$. 
Soit $n$ un entier $\geq 1$ qui annule $F$. Reprenons de nouveau les notations de \ref{Kp1s1}.  
Il résulte de (\cite{sga4} XII 4.4) qu'on a un morphisme de suites spectrales de Cartan-Leray
\begin{equation}
\xymatrix{
{\rH^a(X_\oeta,\rR^bg_*(F))}\ar[d]_{u_{a,b}}\ar@{=>}[r]&{\rH^{a+b}(\oX^\circ,F)}\ar[d]^{u}\\
{\rH^a(Y^{(n)},\rR^b g_{n*}(\varphi_n^{\circ*}(F)))}\ar@{=>}[r]&{\rH^{a+b}(Y^{(n)\circ},\varphi_n^{\circ*}(F))}}
\end{equation}
où $u_{a,b}$ est induit par le morphisme de changement de base \eqref{Kp1s1c} et 
$u$ est le morphisme canonique. Pour tout $b\geq 1$, le morphisme $u_{a,b}$ est nul, en vertu de \ref{Kp1s3}. Il en est alors de même du morphisme
induit par $u$ entre les gradués des filtrations aboutissement.
Procédant comme plus haut, on peut donc se réduire par une récurrence finie, au cas où $i_\xi=q$. 
Par suite, $\xi$ est l'image canonique d'un élément $\zeta\in \rH^q(X_\oeta,g_*(F))$. En vertu de \ref{Kpun21} et \ref{Kpun11}, 
il existe un voisinage étale $U$ de $\ox$ dans $X$ et un revêtement étale surjectif $W'\rightarrow U_\oeta$ tel que l'image canonique
de $\zeta$ dans $\rH^q(W',g_*(F))$ soit nulle. Par suite, l'image canonique de $\xi$ dans $\rH^q(W'^\circ,F)$ est nulle.

\begin{cor}[\cite{achinger} 9.5]\label{Kp1s6}
Soient $\ox$ un point géométrique de $X$ au-dessus de $s$, 
$\uX$ le localisé strict de $X$ en $\ox$. Alors, $\uoX^\circ$ est un schéma $K(\pi,1)$ \eqref{Kpun2}.
\end{cor}

En effet, $\uoX$ étant normal et strictement local (et en particulier intègre) d'après (\cite{agt} III.3.7), 
la proposition résulte de \ref{Kpun18} et \ref{Kp1s5}. 

\begin{rema}\label{Kp1s7}
La proposition \ref{Kp1s5} lorsque le schéma (usuel) $X$ est lisse sur $S$, est due à Faltings (\cite{faltings1} Lemme 2.3 page 281). 
Le cas général est dû à Achinger (\cite{achinger} 9.5). 
On notera que le cas où $\ox$ est un point géométrique de $X_\eta$ a déjà été traité dans \ref{Kpun16}.  
\end{rema}

\section{Acyclicité locale}\label{acycloc}

\subsection{}\label{acycloc1}
On rappelle que l'on dispose des morphismes canoniques
\begin{eqnarray}
\psi\colon \oX^\circ_\et\rightarrow \tE,\label{acycloc1a}\\
\rho\colon X_\et\gtimes_{X_\et}\oX^\circ_\et\rightarrow \tE,\label{acycloc1b}
\end{eqnarray}
\eqref{tf1m} et \eqref{tf3b}, et du plongement canonique $\delta\colon \tE_s\rightarrow \tE$ \eqref{TFA66a}.

\begin{prop}\label{acycloc2}
Pour tout faisceau abélien de torsion, localement constant et constructible $F$ de $\oX^\circ_\et$ et 
tout entier $q\geq 1$, on a $\rR^q\psi_*(F)=0$.
\end{prop}

En effet, d'après (\cite{sga4} V 5.1), $\rR^q\psi_*(F)$ est le faisceau de $\tE$ associé au préfaisceau sur $E$ défini par 
\begin{equation}
(V\rightarrow U)\mapsto \rH^q(V,F).
\end{equation}
Soient $(V\rightarrow U)$ un objet de $E$, $\xi\in \rH^q(V,F)$, $\ox$ un point géométrique de $U$. 
D'après \ref{Kp1s5}, il existe un voisinage étale $U_1$ de $\ox$ dans $X$
et un objet $(V_1\rightarrow U_1)$ de $E$ au-dessus de $(V\rightarrow U)$ tels que le morphisme canonique
$V_1\rightarrow V\times_{\oU^\circ}\oU^\circ_1$ soit surjectif et que l'image canonique de $\xi$ dans $\rH^q(V_1,F)$ soit nulle. 
Par suite, il existe un recouvrement 
$((V_\alpha\rightarrow U_\alpha)\rightarrow (V\rightarrow U))_{\alpha\in \Sigma}$ de $E$ pour la topologie co-évanescente \eqref{tf1} 
tel que pour tout $\alpha\in \Sigma$, l'image canonique de $\xi$ dans $\rH^q(V_\alpha,F)$ soit nulle; d'où la proposition. 

\begin{prop}
\begin{itemize}
\item[{\rm (i)}] Le morphisme d'adjonction $\id\rightarrow \rho_*\rho^*$ induit un isomorphisme $\delta_*\stackrel{\sim}{\rightarrow} \rho_*\rho^*\delta_*$;
en particulier, le foncteur composé
\begin{equation}
\rho^* \delta_*\colon \tE_s\rightarrow X_\et\gtimes_{X_\et}\oX^\circ_\et
\end{equation}
est pleinement fidèle. 
\item[{\rm (ii)}] Pour tout faisceau abélien $F$ de $\tE_s$ et tout entier $q\geq 1$, on a $\rR^q\rho_*(\rho^*(\delta_*(F)))=0$. 
\end{itemize}
\end{prop}

Soient $\ox$ un point géométrique de $X$ au-dessus de $s$,  $\uX$ le localisé strict de $X$ en $\ox$. 
Notons $\rho_{\uoX^\circ}\colon \uoX^\circ_\et\rightarrow \uoX^\circ_\fet$ le morphisme canonique \eqref{notconv10a}  et 
$\varphi_\ox\colon \tE\rightarrow \uoX^\circ_\fet$ le foncteur canonique  défini dans \eqref{tf10b}. 
D'après (\cite{agt} III.3.7), $\uoX$ est normal et strictement local (et en particulier intègre). 
La preuve de \ref{tf12} montre alors que le morphisme 
d'adjonction $\id\rightarrow \rho_*\rho^*$ induit un isomorphisme $\varphi_\ox \stackrel{\sim}{\rightarrow} \varphi_\ox \rho_*\rho^*$. 
En vertu de \ref{Kp1s6}, $\uoX^\circ$ est un schéma $K(\pi,1)$. 
La preuve de \ref{tf13} établit alors que pour tout faisceau abélien $G$ de $\tE$ et tout entier $q\geq 1$, on a $\varphi_\ox(\rR^q\rho_*(\rho^*(G)))=0$. 

Par ailleurs, pour tout faisceau $F$ de $\tE_s$, $\rho_*(\rho^*(\delta_*(F)))|\sigma^*(X_\eta)$ est l'objet final de $\tE_{/\sigma^*(X_\eta)}$. 
Pour tout faisceau abélien $G$ de $\tE_s$ et tout entier $q\geq 0$, on a $\rR^q\rho_*(\rho^*(\delta_*(G)))|\sigma^*(X_\eta)=0$. La proposition s'ensuit puisque la famille des foncteurs $\varphi_\ox \delta_*$ 
lorsque $\ox$ décrit l'ensemble des points géométriques de $X$ au-dessus de $s$, 
est conservative en vertu de \ref{tf21}(ii) et (\cite{agt} VI.10.31).

\section{Théorème de pureté de Faltings et cohomologie galoisienne}\label{tpcg}

\subsection{}\label{tpcg1}
On suppose dans cette section que le morphisme 
$f\colon (X,\cM_X)\rightarrow (S,\cM_S)$ \eqref{TFA3} admet une carte adéquate 
$((P,\gamma),(\mN,\iota),\vartheta\colon \mN\rightarrow P)$ (cf. \ref{cad1}), que $X=\Spec(R)$ est affine et connexe et que $X_s$ est non-vide. 
Ces conditions correspondent à celles fixées dans (\cite{agt} II.6.2).
On reprend les notations de \ref{cad}, en particulier celles fixées dans \ref{cad6} et \ref{cad9}.

\begin{teo}[Faltings, \cite{faltings2} §~2b]\label{tpcg2}
Pour toute $A^\circ_\infty$-algèbre étale finie $B'$ \eqref{cad9c}, la fermeture intégrale  $B$  de $A_\infty$ dans $B'$
est $\alpha$-finie-étale sur $A_\infty$ \eqref{aet2}.
\end{teo}

Signalons ici que Scholze dispose, dans le cadre de sa théorie des perfectoïdes,  d'une généralisation de ce  résultat  
(\cite{scholze} 1.10 et 7.9).

\subsection{}\label{tpcg3}
Soit $\oy$ un point géométrique de $\oX^\circ$ que l'on suppose fixé dans la suite de cette section. 
Le schéma $\oX$ étant localement irréductible d'après \ref{cad7}(iii),  
il est la somme des schémas induits sur ses composantes irréductibles. On note $\oX^\star$
la composante irréductible de $\oX$ contenant $\oy$. 
De même, $\oX^\circ$ est la somme des schémas induits sur ses composantes irréductibles
et $\oX^{\star \circ}=\oX^\star\times_{X}X^\circ$ est la composante irréductible de $\oX^\circ$ contenant $\oy$. 
On désigne par $\Delta$ le groupe profini $\pi_1(\oX^{\star \circ},\oy)$ 
et pour alléger les notations, par $\oR$ la représentation discrète $\oR^\oy_X$ de $\Delta$ définie dans \eqref{TFA9b}. 
On notera que $\oR$ est intègre et normal \eqref{TFA10}.

\subsection{}\label{tpcg4}
Pour tous entiers $m,n\geq 1$, le morphisme canonique $X^{(mn)}\rightarrow X^{(n)}$ 
est fini et surjectif. D'après (\cite{ega4} 8.3.8(i)), il existe alors un $X$-morphisme
\begin{equation}\label{tpcg4a}
\oy\rightarrow \underset{\underset{n\geq 1}{\longleftarrow}}{\lim}\ X^{(n)},
\end{equation} 
la limite projective étant indexée par l'ensemble $\mZ_{\geq 1}$ ordonné par la relation de divisibilité.
On fixe un tel morphisme dans toute la suite de cette section. Celui-ci induit un $S$-morphisme
\begin{equation}\label{tpcg4b}
\oS\rightarrow \underset{\underset{n\geq 1}{\longleftarrow}}{\lim}\ S^{(n)}. 
\end{equation}
Pour tout entier $n\geq 1$, on pose
\begin{equation}\label{tpcg4c}
\oX^{(n)}= X^{(n)}\times_{S^{(n)}}\oS \ \ \ {\rm et}\ \ \ \oX^{(n)\circ}=\oX^{(n)}\times_XX^\circ.
\end{equation} 
On en déduit un $\oX$-morphisme 
\begin{equation}\label{tpcg4d}
\oy\rightarrow \underset{\underset{n\geq 1}{\longleftarrow}}{\lim}\ \oX^{(n)}.
\end{equation} 

Pour tout entier $n\geq 1$, le schéma $\oX^{(n)}$ étant normal et localement irréductible d'après \ref{cad7}(iii), 
il est la somme des schémas induits sur ses composantes irréductibles. 
On note $\oX^{(n)\star}$ la composante irréductible de $\oX^{(n)}$ contenant l'image de $\oy$ \eqref{tpcg4d}.
De même, $\oX^{(n)\circ}$ est la somme des schémas induits sur ses composantes irréductibles
et $\oX^{(n)\star\circ}=\oX^{(n)\star}\times_XX^\circ$  est la composante irréductible de $\oX^{(n)\circ}$ 
contenant l'image de $\oy$. On notera que $\oX^{(n)}$ étant fini sur $\oX$ \eqref{cad6c},  
$\oX^{(n)\star}$ est la fermeture intégrale de $\oX^{\star}$ dans $\oX^{(n)\star\circ}$. 
On pose 
\begin{equation}\label{tpcg4f}
R_n=\Gamma(\oX^{(n)\star},\co_{\oX^{(n)}}).
\end{equation}
D'après \ref{cad7}(v), le morphisme $\oX^{(n)\star\circ}\rightarrow \oX^{\star\circ}$ 
est étale fini. Il résulte de la preuve de (\cite{agt} II.6.8(iv)) que $\oX^{(n)\star \circ}$ 
est en fait un revêtement étale fini et galoisien de $\oX^{\star \circ}$ de groupe $\Delta_n$ 
canoniquement isomorphe à un sous-groupe de $\Hom_\mZ(P^\gp/\mZ\lambda,\mu_{n}(\oK))$.  
Le groupe $\Delta_n$ agit naturellement sur $R_n$. 

Si $n$ est une puissance de $p$, le morphisme canonique
\begin{equation}
\oX^{(n)\star}\rightarrow\oX^{(n)}\times_\oX\oX^\star
\end{equation}
est un isomorphisme en vertu de (\cite{agt} II.6.6(v)), et on a donc $\Delta_n\simeq \Hom_\mZ(P^\gp/\mZ\lambda,\mu_{n}(\oK))$.

Les anneaux $(R_n)_{n\geq 1}$ forment naturellement un système inductif. On pose 
\begin{eqnarray}
R_\infty&=&\underset{\underset{n\geq 1}{\longrightarrow}}{\lim}\ R_n,\label{tpcg4g}\\
R_{p^\infty}&=&\underset{\underset{n\geq 0}{\longrightarrow}}{\lim}\ R_{p^n}.\label{tpcg4gg}
\end{eqnarray}
Ce sont des anneaux normaux et intègres d'après (\cite{ega1n} 0.6.1.6(i) et 0.6.5.12(ii)).

Le morphisme \eqref{tpcg4d} induit un $\oX^\star$-morphisme 
\begin{equation}\label{tpcg4e}
\Spec(\oR)\rightarrow \underset{\underset{n\geq 1}{\longleftarrow}}{\lim}\ \oX^{(n)\star},
\end{equation} 
et par suite des homomorphismes injectifs de $R_1$-algèbres 
\begin{equation}\label{tpcg4h}
R_{p^\infty}\rightarrow R_\infty\rightarrow \oR.
\end{equation}

Les groupes $(\Delta_n)_{n\geq 1}$ forment naturellement un système projectif. On pose 
\begin{eqnarray}
\Delta_\infty&=&\underset{\underset{n\geq 1}{\longleftarrow}}{\lim}\ \Delta_n,\label{tpcg4k}\\
\Delta_{p^\infty}&=&\underset{\underset{n\geq 0}{\longleftarrow}}{\lim}\ \Delta_{p^n}.\label{tpcg4kk}
\end{eqnarray}
On a des homomorphismes canoniques 
\begin{equation}\label{tpcg4l}
\xymatrix{
{\Delta_{\infty}}\ar@{->>}[d]\ar@{^(->}[r]&{\Hom_{\mZ}(P^\gp/\mZ\lambda,\hmZ(1))}\ar[d]\\
{\Delta_{p^\infty}}\ar[r]^-(0.5){\sim}&{\Hom_{\mZ}(P^\gp/\mZ \lambda,\mZ_p(1))}}
\end{equation}
Le noyau $\Sigma_0$ de l'homomorphisme canonique $\Delta_\infty\rightarrow \Delta_{p^\infty}$
est un groupe profini d'ordre premier à $p$. Par ailleurs, le morphisme \eqref{tpcg4d} détermine un homomorphisme 
surjectif $\Delta\rightarrow \Delta_\infty$. On note $\Sigma$ son noyau. Les homomorphismes \eqref{tpcg4h} sont alors $\Delta$-équivariants. 
\begin{equation}\label{tpcg4m}
\xymatrix{
R_1\ar[rr]^{\Delta_{p^\infty}}\ar@/^2pc/[rrrr]|{\Delta_\infty}\ar@/_2pc/[rrrrrr]|{\Delta}&&
{R_{p^\infty}}\ar[rr]^{\Sigma_0}&&{R_\infty}\ar[rr]^\Sigma&&\oR}
\end{equation}
On note $\oF$ (resp. $F_\infty$, resp. $F_{p^\infty}$) le corps des fractions de $\oR$ (resp. $R_\infty$, resp. $R_{p^\infty}$).

\begin{prop}[\cite{agt} II.6.17] \label{tpcg5}
L'extension $\oF$ de $F_\infty$ \eqref{tpcg4m} est la réunion d'un système inductif filtrant de 
sous-extensions galoisiennes finies $N$ de $F_\infty$ telles que la clôture intégrale
de $R_\infty$ dans $N$ soit $\alpha$-finie-étale sur $R_\infty$ \eqref{aet2}.
\end{prop}

\begin{cor}\label{tpcg17}
La $\co_{\oK}$-algèbre $\oR$ est universellement $\alpha$-cohérente.
\end{cor}

Cela résulte de \ref{aet9}, \ref{cad15} et \ref{tpcg5}.

\begin{cor}\label{tpcg21}
La $R_\infty$-algèbre $\oR$ est $\alpha$-fidèlement plate.
\end{cor}
Cela résulte de \ref{aet11}, \ref{tpcg5}, \ref{aet8}  et (\cite{agt} V.12.4).

\begin{cor}\label{tpcg22}
Il existe un entier $N\geq 1$ tel que pour toute suite exacte de $R$-modules
$M'\rightarrow M\rightarrow M''$, la suite 
\begin{equation}
M'\otimes_{R}\oR\rightarrow M\otimes_{R}\oR\rightarrow M''\otimes_{R}\oR
\end{equation}
soit $p^N$-exacte \eqref{finita3}.
\end{cor}
Cela résulte de \ref{cad17} et \ref{tpcg21}.

\subsection{}\label{tpcg20}
On pose \eqref{TFA3b}
\begin{equation}\label{tpcg20a}
\tOmega^1_{R/\co_K}=\tOmega^1_{X/S}(X).
\end{equation}
D'après \ref{cad1}(iv) et (\cite{kato1} 1.8), on a un isomorphisme $R$-linéaire canonique
\begin{equation}\label{tpcg20b}
\tOmega^1_{R/\co_K}\stackrel{\sim}{\rightarrow} (P^\gp/\mZ\lambda)\otimes_\mZ R.
\end{equation}
En particulier, $\tOmega^1_{R/\co_K}$ est un $R$-module libre de rang $d=\rg(P^\gp)-1=\dim(X/S)$.  
Comme le groupe $P^\gp/\mZ\lambda$ est libre d'après \ref{cad2}(i),  on a un isomorphisme canonique
\begin{equation}\label{tpcg20c}
(P^\gp/\mZ\lambda)\otimes_\mZ\mZ_p(-1)\stackrel{\sim}{\rightarrow}
\Hom_{\mZ_p}(\Delta_{p^\infty}, \mZ_p).
\end{equation}
On en déduit, pour tout élément non nul $a$ de $\co_\oK$, un isomorphisme $R_1$-linéaire
\begin{equation}\label{tpcg20d}
\tOmega^1_{R/\co_K}\otimes_R(R_1/aR_1)(-1)\stackrel{\sim}{\rightarrow} 
\Hom_{\mZ}(\Delta_{p^\infty},R_1/aR_1).
\end{equation}

\begin{prop}[\cite{agt} II.8.17]\label{tpcg6}
Pour tout élément non nul $a$ de $\co_\oK$, il existe un unique homomorphisme de $R_1$-algèbres graduées 
\begin{equation}\label{tpcg6a}
\wedge(\tOmega^1_{R/\co_K}\otimes_R(R_1/aR_1)(-1))\rightarrow \rH^*(\Delta,\oR/a\oR)
\end{equation}
dont la composante en degré un est induite par \eqref{tpcg20d}. 
Celui-ci est $\alpha$-injectif et son conoyau est annulé par $p^{\frac{1}{p-1}}\fm_\oK$.   
\end{prop}

\begin{prop}[\cite{agt} II.9.10]\label{tpcg7} 
L'endomorphisme de Frobenius absolu de $\oR/p\oR$ est surjectif.
\end{prop}

\begin{lem}\label{tpcg18}
Soient $A$ une $\co_{\oK}$-algèbre, $G$ un groupe abélien fini, $M$ un $A$-module $\alpha$-cohérent \eqref{finita9}, 
muni d'une action linéaire de $G$. Alors, pour tout entier $i\geq 0$, le $A$-module $\rH^i(G,M)$ est $\alpha$-cohérent.
\end{lem}

Compte tenu de \ref{finita17}, \ref{finita18} et de la suite spectrale de Hochschild-Serre (\cite{agt} II.3.4), on peut se borner au cas où $G$ est cyclique. 
Soient $\sigma$ un générateur de $G$, $\Tr$ l'endomorphisme de $M$ induit par 
$\sum_{\sigma\in G}\sigma$. Alors, $\rR\Gamma(G,M)$ est quasi-isomorphe au complexe de $A$-modules, 
placé en degrés $\geq 0$, 
\begin{equation}
M\stackrel{\sigma-1}{\longrightarrow} M\stackrel{\Tr}{\longrightarrow} M\stackrel{\sigma-1}{\longrightarrow}
M\stackrel{\Tr}{\longrightarrow}M\stackrel{\sigma-1}{\longrightarrow}  \dots,
\end{equation}
d'où la proposition \eqref{finita18}.

\begin{lem}\label{tpcg8}
Soient $A$ une $\co_{\oK}$-algèbre, $G$ un groupe profini, 
$H$ un sous-groupe fermé et distingué de $G$, d'ordre (comme groupe profini) premier à $p$ et tel que $G/H$ soit isomorphe à $\mZ_p^d$, 
$M$ un $A$-module $\alpha$-cohérent \eqref{finita9} dont tous les éléments sont de torsion $p$-primaire, muni d'une action linéaire discrète de $G$. 
Alors, pour tout entier $i\geq 0$, le $A$-module $\rH^i(G,M)$ est $\alpha$-cohérent et est nul pour tout $i\geq d+1=\rg(P^\gp)$.
\end{lem}

En effet, pour tout $i\geq 0$, on a un isomorphisme canonique $\rH^i(G,M) \stackrel{\sim}{\rightarrow} \rH^i(G/H,M^{H})$.  
Montrons que le $A$-module $M^H$ est $\alpha$-cohérent. Pour tout $\gamma\in \fm_\oK$, il existe un sous-$A$-module de type 
fini $M'$ de $M$, engendré par $u_1,\dots,u_t$, tel que $M/M'$ soit annulé par $\gamma$. Quitte à remplacer $M'$ par le sous-$A$-module de $M$
engendré par les conjugués des $u_i$ par $G$ $(1\leq i\leq t)$, qui sont en nombre fini, on peut supposer $M'$ stable par
l'action de $G$. Comme $M$ est $\alpha$-cohérent, $M'$ est $\alpha$-cohérent. 
L'action de $G$ sur $M'$ se factorise à travers un quotient fini $G_0$. Soit $H_0$ l'image de $H$ dans $G_0$. 
Le $A$-module $M'^{H}=M'^{H_0}$ est alors $\alpha$-cohérent \eqref{finita18}. 
Comme $M^H/M'^H$ est annulé par $\gamma$, on en déduit que $M^H$ est $\alpha$-cohérent \eqref{finita21}. 
On se réduit ainsi au cas où $G\simeq\mZ_p^d$. 
D'après (\cite{agt} II.3.23), $\rR\Gamma(G,M)$ est alors quasi-isomorphe 
à un complexe borné, en degrés compris entre $0$ et $d$, de $A$-modules $\alpha$-cohérents~; d'où la proposition \eqref{finita18}.

\begin{lem}\label{tpcg19}
Soient $A$ une $\co_{\oK}$-algèbre, $G$ un groupe profini isomorphe à $\mZ_p^d$, 
$M$ un $A$-module de type $\alpha$-fini dont tous les éléments sont de torsion $p$-primaire, 
muni d'une action linéaire discrète de $G$, $\gamma\in \fm_\oK$.
Alors, il existe un entier $n_0\geq 0$ vérifiant la propriété suivante: 
pour tout $A$-module $N$ muni de l'action triviale de $G$, 
tout entier $n\geq n_0$ et tout homomorphisme {\em surjectif} $\nu\colon G\rightarrow \mu_{p^n}(\co_\oK)$, 
notant $\co_{\oK}(\nu)$ le $(\co_\oK)$-$G$-module topologique $\co_\oK$, 
muni de la topologie $p$-adique et de l'action de $G$ définie par $\nu$, pour tout entier $i\geq 0$, 
$\rH^i(G,M\otimes_AN\otimes_{\co_\oK}\co_\oK(\nu))$ est annulé par $\gamma$.
\end{lem}

Il existe un sous-$A$-module de type 
fini $M'$ de $M$, engendré par $u_1,\dots,u_\ell$, tel que $M/M'$ soit annulé par $\gamma$. Quitte à remplacer $M'$ par le sous-$A$-module de $M$
engendré par les conjugués des $u_i$ par $G$ $(1\leq i\leq \ell)$, qui sont en nombre fini, on peut supposer $M'$ stable par l'action de $G$. 
Pour tout $A$-module $N$, le morphisme canonique $M'\otimes_A N\rightarrow M\otimes_A N$ est alors un $\gamma^2$-isomorphisme \eqref{alpha3}. 
On peut donc se borner au cas où $M$ est de type fini sur $A$. Soit $e_1,\dots,e_d$ une $\mZ_p$-base de $G$. 
Il existe un entier $t\geq 1$ tel que pour tout $1\leq j\leq d$, $e_j^{p^t}$ agisse trivialement sur $M$. 
Soient $n$ un entier $\geq t$, $\nu\colon G\rightarrow \mu_{p^n}(\co_\oK)$ un homomorphisme surjectif.
Il existe $1\leq j\leq d$ tel que que $\zeta_n=\nu(e_j)$ soit une racine primitive $p^n$-ième de l'unité.
Notons $G_j$ le sous-groupe de $G$ engendré par $e_j$. 
D'après (\cite{agt} II.3.23), pour tout $A$-module $N$ muni de l'action triviale de $G$ et tout entier $i\geq 0$, 
$\rH^i(G_j,M\otimes_A N \otimes_{\co_\oK}\co_\oK(\nu))$ est la cohomologie du complexe, concentré en degrés $0, 1$, suivant 
\begin{equation}
\zeta_n e_j-\id_{M\otimes_AN}\colon M\otimes_AN\rightarrow M\otimes_AN.
\end{equation}
La relation 
\begin{equation}
(\zeta_n e_j-\id_{M\otimes_AN})((\zeta_n e_j)^{p^t-1}+(\zeta_n e_j)^{p^t-2}+\dots+\id_{M\otimes_AN})=(\zeta_n^{p^t}-1)\id_{M\otimes_AN}
\end{equation}
montre que le noyau et le conoyau de $\zeta_n e_j -\id_{M\otimes_AN}$ sont annulés par $\zeta_n^{p^t}-1$. 
Par suite, compte tenu de la suite spectrale (\cite{agt} II.3.4)
\begin{equation}
E_2^{a,b}=\rH^a(G/G_j,\rH^b(G_j,M\otimes_A N \otimes_{\co_\oK}\co_\oK(\nu)))\Rightarrow \rH^{a+b}(G, M\otimes_A N \otimes_{\co_\oK}\co_\oK(\nu)),
\end{equation}
pour tout $i\geq 0$, $\rH^i(G, M\otimes_A N \otimes_{\co_\oK}\co_\oK(\nu))$ est annulé par $(\zeta_n^{p^t}-1)^2$.
La proposition s'ensuit puisqu'il existe $n_0\geq 0$ tel que pour tout $n\geq n_0$, $\gamma\in (\zeta_{n}^{p^t}-1)^2 \co_\oK$.

\begin{lem}\label{tpcg9}
Soit $m$ un entier $\geq 0$. Alors, 
\begin{itemize}
\item[{\rm (i)}] Pour tout entier $n\geq 0$, tout $(R_{p^n}/p^mR_{p^n})$-module de présentation finie $M$, 
muni d'une action semi-liné\-aire de $\Delta_{p^n}$ et tout entier $i\geq 0$, 
le $R_1$-module $\rH^i(\Delta_{p^\infty},M\otimes_{R_{p^n}}R_{p^\infty})$ est de présentation $\alpha$-finie, 
et est nul pour tout $i\geq d+1=\rg(P^\gp)$.
\item[{\rm (ii)}] Pour tout entier $n\geq 1$, tout $(R_{n}/p^mR_{n})$-module de présentation finie $M$, 
muni d'une action semi-linéaire de $\Delta_n$ et tout entier $i\geq 0$, 
le $R_1$-module $\rH^i(\Delta_{\infty},M\otimes_{R_n}R_\infty)$ est de présentation $\alpha$-finie, 
et est nul pour tout $i\geq d+1$.
\end{itemize}
\end{lem}

On notera d'abord que $R_1/p^mR_1$ est un anneau cohérent (\cite{egr1} 1.10.3). Tout $(R_1/p^mR_1)$-module de présentation 
$\alpha$-finie est donc $\alpha$-cohérent \eqref{afini6}.

(i) En effet, la $p$-dimension cohomologique de $\Delta_{p^\infty}$ étant égale à $d$ (\cite{agt} II.3.24),  
il suffit de montrer que pour tout $i\geq 0$, le $(R_1/p^mR_1)$-module $\rH^i(\Delta_{p^\infty},M\otimes_{R_{p^n}}R_{p^\infty})$ est 
$\alpha$-cohérent. Comme $R_{p^n}/p^mR_{p^n}$ est une $(R_1/p^mR_1)$-algèbre cohérente, 
tout $(R_{p^n}/p^mR_{p^n})$-module de présentation $\alpha$-finie est un $(R_1/p^mR_1)$-module  
de présentation $\alpha$-finie et est donc $\alpha$-cohérent.
D'après \ref{tpcg18} et la suite spectrale de Hochschild-Serre (\cite{agt} II.3.4),
en tenant compte de \ref{cad7}(i), il est alors loisible de remplacer $f$ par $f^{(p^n)}$ \eqref{cad6c}.
On peut donc se borner au cas où $n=0$. 

Pour tout entier $j\geq 0$, posons 
\begin{eqnarray}
\Xi_{p^j}&=&\Hom(\Delta_{p^\infty},\mu_{p^j}(\co_\oK)),\\
\Xi_{p^\infty}&=&\Hom(\Delta_{p^\infty},\mu_{p^\infty}(\co_\oK)).
\end{eqnarray}
On identifie $\Xi_{p^j}$ à un sous-groupe de $\Xi_{p^\infty}$. 
Pour tout $\nu\in \Xi_{p^\infty}$, on note $\co_\oK(\nu)$ le $(\co_\oK)$-$\Delta_{p^\infty}$-module topologique $\co_\oK$, 
muni de la topologie $p$-adique et de l'action de $\Delta_{p^\infty}$ définie par $\nu$.
D'après (\cite{agt} II.8.9), il existe une décomposition canonique de $R_{p^\infty}$ en somme directe de $R_1$-modules de présentation finie,
stables sous l'action de $\Delta_{p^\infty}$,
\begin{equation}
R_{p^\infty}=\bigoplus_{\nu\in \Xi_{p^\infty}}R_{p^\infty}^{(\nu)}\otimes_{\co_\oK}\co_\oK(\nu),
\end{equation}
où $\Delta_{p^\infty}$ agit trivialement sur $R_{p^\infty}^{(\nu)}$. Pour tout entier $i\geq 0$, on a donc une décomposition canonique 
\begin{equation}
\rH^i(\Delta_{p^\infty},M\otimes_{R_1}R_{p^\infty})=\bigoplus_{\nu\in \Xi_{p^\infty}}\rH^i(\Delta_{p^\infty},M\otimes_{R_1}R_{p^\infty}^{(\nu)}\otimes_{\co_\oK}\co_\oK(\nu)).
\end{equation}
Pour tout $\nu \in \Xi_{p^\infty}$, le $R_1$-module $M\otimes_{R_1}R^{(\nu)}_{p^\infty}$ est de présentation finie. D'après \ref{tpcg19},
pour tout $\gamma\in \fm_\oK$, il existe donc un entier $j\geq 1$ tel que pour tout $\nu\in \Xi_{p^\infty}-\Xi_{p^j}$, le $R_1$-module
$\rH^i(\Delta_{p^\infty},M\otimes_{R_1}R_{p^\infty}^{(\nu)}\otimes_{\co_\oK}\co_\oK(\nu))$ soit annulé par $\gamma$. 
Par suite, en vertu \ref{tpcg8}, le $R_1$-module $\rH^i(\Delta_{p^\infty},M\otimes_{R_1}R_{p^\infty})$ est de présentation $\alpha$-finie. 

(ii) Soit $n'$ le plus grand diviseur premier à $p$ de $n$.   
Comme $R_{n'}/p^mR_{n'}$ est une $(R_1/p^mR_1)$-algèbre cohérente, 
tout $(R_{n'}/p^mR_{n'})$-module de présentation $\alpha$-finie est un $(R_1/p^mR_1)$-module de présentation $\alpha$-finie 
et est donc $\alpha$-cohérent. Le foncteur $\Gamma(\Delta_{n'}, -)$ est exact,
et il transforme les $(R_1/p^mR_1)$-modules $\alpha$-cohérents en $(R_1/p^mR_1)$-modules $\alpha$-cohérents \eqref{finita18}.  
Il est alors loisible de remplacer $f$ par $f^{(n')}$ \eqref{cad6c}, d'après \ref{cad7}(i).
On peut donc supposer que $n$ est une puissance de $p$.

Comme $\Sigma_0$ est un groupe profini d'ordre premier à $p$ \eqref{tpcg4l}, pour tout $i\geq 0$,
on a un isomorphisme canonique 
\begin{equation}
\rH^i(\Delta_{p^\infty},(M\otimes_{R_n}R_\infty)^{\Sigma_0})\stackrel{\sim}{\rightarrow} \rH^i(\Delta_\infty,M\otimes_{R_n}R_\infty).
\end{equation} 
Par ailleurs, considérant une présentation finie de $M$ sur $R_n$, on déduit de (\cite{agt} II.6.13) que le morphisme canonique
\begin{equation}
M\otimes_{R_n}R_{p^\infty}\rightarrow (M\otimes_{R_n}R_\infty)^{\Sigma_0}
\end{equation}
est un isomorphisme. La proposition résulte donc de (i).

\begin{prop}\label{tpcg10}
Soient $m$ un entier $\geq 1$, $\mL$ un $(\co_\oK/p^m\co_\oK)$-module libre de type fini, muni d'une action linéaire continue de $\Delta$. 
Alors, 
\begin{itemize}
\item[{\rm (i)}] Le $(R_\infty/p^mR_\infty)$-module $(\mL\otimes_{\co_\oK}\oR)^\Sigma$ est $\alpha$-projectif de type $\alpha$-fini \eqref{aet1}, 
et le morphisme canonique 
\begin{equation}\label{tpcg10a}
(\mL\otimes_{\co_\oK}\oR)^\Sigma\otimes_{R_\infty}\oR\rightarrow \mL\otimes_{\co_\oK}\oR
\end{equation}
est un $\alpha$-isomorphisme. 
\item[{\rm (ii)}]  Pour tout $\gamma\in \fm_\oK$, il existe un entier $n\geq 1$, un $(R_n/p^mR_n)$-module de présentation finie $M_n$, 
muni d'une action semi-linéaire de $\Delta_n$, et un morphisme $R_{\infty}$-linéaire et $\Delta_{\infty}$-équivariant 
\begin{equation}\label{tpcg10b}
M_n\otimes_{R_n}R_{\infty}\rightarrow (\mL\otimes_{\co_\oK}\oR)^\Sigma,
\end{equation}
dont le noyau et le conoyau sont annulés par $\gamma$. 
\item[{\rm (iii)}] Pour tout entier $i\geq 0$, le $R_1$-module $\rH^i(\Delta,\mL\otimes_{\co_\oK}\oR)$ est de présentation $\alpha$-finie, 
et est $\alpha$-nul pour tout $i\geq d+1=\rg(P^\gp)$.
\end{itemize}
\end{prop}

(i) Soient $N$ une extension galoisienne finie de $F_\infty$ contenue dans $\oF$ \eqref{tpcg4m}, 
$D$ la clôture intégrale de $R_\infty$ dans $N$, $G=\Gal(N/F_\infty)$, 
$\Tr_G$ l'endomorphisme $R_\infty$-linéaire de $D$ (ou de $D/p^m D$) induit par $\sum_{\sigma\in G}\sigma$.
Comme on a $D=\oR\cap N$ et $R_\infty=D\cap F_\infty$ d'après \ref{cad7}(iii), 
les homomorphismes $R_\infty/p^m R_\infty\rightarrow D/p^m D\rightarrow \oR/p^m \oR$ sont injectifs. 
Supposons que la $R_\infty$-algèbre $D$ soit $\alpha$-finie-étale.  
Alors, $D$ est un $\alpha$-$G$-torseur sur $R_\infty$ en vertu de \ref{aet8}. Par suite, le quotient
\[
\frac{(D/p^m D)^G}{\Tr_G(D/p^m D)}
\]
est $\alpha$-nul en vertu de \ref{aet7}. Comme $\Tr_G(D)\subset R_\infty$, l'homomorphisme 
$R_\infty/p^m R_\infty\rightarrow (D/p^m D)^G$ est un $\alpha$-isomorphisme. 
Compte tenu de \ref{tpcg5}, on en déduit, par passage à la limite inductive, que l'homomorphisme 
$R_\infty/p^mR_\infty\rightarrow (\oR/p^m\oR)^\Sigma$ est un $\alpha$-isomorphisme.
On démontre par le même argument, en tenant compte de (\cite{agt} V.7.8), que l'homomorphisme
\begin{equation}\label{tpcg10c}
D/p^mD\rightarrow (\oR/p^m\oR)^{\Gal(\oF/N)}
\end{equation}
est un $\alpha$-isomorphisme. 

Choisissons $N$ de sorte que l'action de $\Sigma$ sur $\mL$ se factorise à travers $G$ \eqref{tpcg5}.
On en déduit que le morphisme canonique
\begin{equation}\label{tpcg10d}
(\mL\otimes_{\co_\oK} D)^G\rightarrow (\mL\otimes_{\co_\oK}\oR)^\Sigma
\end{equation}
est un $\alpha$-isomorphisme. Par ailleurs, il résulte de \ref{aet5} que le morphisme canonique 
\begin{equation}\label{tpcg10e}
(\mL\otimes_{\co_\oK} D)^G\otimes_{R_\infty}D\rightarrow \mL\otimes_{\co_\oK} D
\end{equation}
est un $\alpha$-isomorphisme~; il en est donc de même de \eqref{tpcg10a}. Comme $D$ est $\alpha$-fidèlement plat sur $R_\infty$ 
d'après (\cite{agt} V.6.7), 
on en déduit que $(\mL\otimes_{\co_\oK}\oR)^\Sigma$ est un $(R_\infty/p^mR_\infty)$-module $\alpha$-projectif de type $\alpha$-fini 
en vertu de (\cite{agt} V.8.7). 

(ii) Il suffit de montrer que pour tout $\gamma\in \fm_\oK$, il existe une suite de $R_\infty$-représentations continues de $\Delta_\infty$  
\begin{equation}
(R_\infty/p^mR_\infty)^a\stackrel{v}{\rightarrow}(R_\infty/p^mR_\infty)^b\stackrel{u}{\rightarrow} (\mL\otimes_{\co_\oK}\oR)^\Sigma\rightarrow 0,
\end{equation}
où $a$ et $b$ sont deux entiers $\geq 0$, telle que la suite des $R_\infty$-modules sous-jacents soit $\gamma$-exacte \eqref{finita3}. 
En effet, il existerait alors un entier $n\geq 1$ tel que $v$ se déduise par extension des scalaires d'un morphisme de $R_n$-représentations 
de $\Delta_n$, $w\colon (R_n/p^mR_n)^a\rightarrow (R_n/p^mR_n)^b$, et on prendrait pour $M_n$ le conoyau de $w$.

D'après (i), il existe deux entiers $r,q\geq 1$ et un morphisme $R_\infty$-linéaire 
\begin{equation}\label{tpcg10f}
R_\infty^r\rightarrow (\mL\otimes_{\co_\oK}\oR)^\Sigma
\end{equation}
dont le conoyau est annulé par $\gamma^{1/12}$ et tels que les images de la base canonique de $R_\infty^r$ soient invariants par le noyau 
de l'homomorphisme $\Delta_\infty\rightarrow \Delta_q$. Notons $R_q\langle \Delta_q\rangle$ l'anneau
non-commutatif de groupe sous-jacent le $R_q$-module libre de base $(e_\sigma)_{\sigma\in \Delta_q}$,
et dont la multiplication est donnée par $(be_\sigma)(b'e_{\sigma'})=b\sigma(b')e_{\sigma\sigma'}$. 
La catégorie des $R_q$-modules munis d'une action semi-linéaire de $\Delta_q$ 
est équivalente à la catégorie des $R_q\langle \Delta_q\rangle$-modules à gauche. 
Le morphisme \eqref{tpcg10f} induit donc un morphisme $R_\infty$-linéaire et $\Delta_\infty$-équivariant
\begin{equation}\label{tpcg10g}
u\colon R_q\langle \Delta_q\rangle^r\otimes_{R_q}(R_\infty/p^mR_\infty)\rightarrow (\mL\otimes_{\co_\oK}\oR)^\Sigma.
\end{equation}
D'après \ref{finita12}(iii), le noyau de $u$ est de type $\gamma$-fini. Recommençant la même construction pour $\ker(u)$, on obtient 
l'assertion recherchée.

(iii) En effet, pour tout $i\geq 0$, le morphisme canonique 
\begin{equation}
\rH^i(\Delta_\infty,(\mL\otimes_{\co_\oK} \oR)^\Sigma)\rightarrow \rH^i(\Delta,\mL\otimes_{\co_\oK} \oR)
\end{equation} 
est un $\alpha$-isomorphisme d'après (\cite{agt} II.6.20). La proposition résulte alors de (ii) et \ref{tpcg9}(ii).

\section{\texorpdfstring{$\alpha$-modules du topos de Faltings}{alpha-modules du topos de Faltings}}

\subsection{}\label{amtF1}
Les hypothèses et notations de \ref{tpcg} sont en vigueur dans cette section.  
Pour tout entier $n\geq 1$, le morphisme $\oX^{(n)\circ}\rightarrow \oX^\circ$ \eqref{tpcg4c} est étale fini
d'après \ref{cad7}(v), de sorte que $(\oX^{(n)\circ}\rightarrow X)$ est un objet de $E$ \eqref{TFA6a}. 
On note $E^{(n)}\rightarrow \Et_{/X}$ (resp. $\tE^{(n)}$) le site fibré (resp. topos) de Faltings associé 
au morphisme $\oX^{(n)\circ}\rightarrow X$ \eqref{tf1}. Tout objet de $E^{(n)}$
est naturellement un objet de $E$. On définit ainsi un foncteur 
\begin{equation}\label{amtF1a}
\Phi^{(n)}\colon E^{(n)}\rightarrow E.
\end{equation}
On vérifie aussitôt que $\Phi^{(n)}$ se factorise à travers une équivalence de catégories 
\begin{equation}\label{amtF1b}
E^{(n)}\stackrel{\sim}{\rightarrow} E_{/(\oX^{(n)\circ}\rightarrow X)}.
\end{equation}
Il résulte alors de (\cite{agt} VI.5.38) que la topologie co-évanescente de $E^{(n)}$ 
est induite par celle de $E$ au moyen du foncteur $\Phi^{(n)}$. 
Par suite, $\Phi^{(n)}$ est continu et cocontinu (\cite{sga4} III 5.2). 
Il définit donc une suite de trois foncteurs adjoints~:
\begin{equation}\label{amtF1c}
(\Phi^{(n)})_!\colon \tE^{(n)}\rightarrow \tE, \ \ \ \Phi^{(n)*}\colon \tE\rightarrow \tE^{(n)}, 
\ \ \ \Phi^{(n)}_*\colon \tE^{(n)}\rightarrow \tE,
\end{equation}
dans le sens que pour deux foncteurs consécutifs de la suite, celui de droite est
adjoint à droite de l'autre. D'après (\cite{sga4} III 5.4), le foncteur $(\Phi^{(n)})_!$ se factorise à travers 
une équivalence de catégories 
\begin{equation}\label{amtF1d}
\tE^{(n)}\stackrel{\sim}{\rightarrow} \tE_{/(\oX^{(n)\circ}\rightarrow X)^\tta},
\end{equation}
où $(\oX^{(n)\circ}\rightarrow X)^\tta$ est l'objet de $\tE$ associé à $(\oX^{(n)\circ}\rightarrow X)$.
Le couple de foncteurs $(\Phi^{(n)*},\Phi^{(n)}_*)$ définit un morphisme de topos que l'on note aussi
\begin{equation}\label{amtF1e}
\Phi^{(n)}\colon \tE^{(n)}\rightarrow \tE
\end{equation}
et qui n'est autre que le composé de \eqref{amtF1d} et du morphisme de localisation de $\tE$ en 
$(\oX^{(n)\circ}\rightarrow X)^\tta$. 
Le foncteur \eqref{amtF1a} est un adjoint à gauche du foncteur 
\begin{equation}\label{amtF1f}
\Phi^{(n)+}\colon E\rightarrow E^{(n)}, \ \ \ (V\rightarrow U)\mapsto (V\times_{\oX^\circ}\oX^{(n)\circ}\rightarrow U).
\end{equation}
D'après (\cite{sga4} III 2.5), le morphisme \eqref{amtF1e} s'identifie donc au morphisme  
défini par fonctorialité (\cite{agt} VI.10.12) à partir du diagramme 
\begin{equation}\label{amtF1g}
\xymatrix{
{\oX^{(n)\circ}}\ar[r]\ar[d]&X\ar@{=}[d]\\
{\oX^\circ}\ar[r]&X}
\end{equation}

\subsection{}\label{amtF2}
Soit $n$ un entier $\geq 1$. 
Pour tout objet $(V\rightarrow U)$ de $E$, on pose $\oU^{(n)}=U\times_X\oX^{(n)}$ et on note 
$\oU^{V(n)}$ la fermeture intégrale de $\oU^{(n)}$ dans $V\times_{\oU}\oU^{(n)}$. 
Pour tout morphisme $(V'\rightarrow U')\rightarrow (V\rightarrow U)$ de $E$, on a un morphisme canonique 
$\oU'^{V'(n)}\rightarrow \oU^{V(n)}$ qui s'insère dans un diagramme commutatif 
\begin{equation}\label{amtF2a}
\xymatrix{
{V'\times_{\oU'}\oU'^{(n)}}\ar[r]\ar[d]&{\oU'^{V'(n)}}\ar[d]\ar[r]&{\oU'^{(n)}}\ar[r]\ar[d]&U'\ar[d]\\
{V\times_{\oU}\oU^{(n)}}\ar[r]&{\oU^{V(n)}}\ar[r]&{\oU^{(n)}}\ar[r]&U}
\end{equation} 
Comme l'injection canonique $X^{(n)\circ}\rightarrow X^{(n)}$ est schématiquement dominante d'après \ref{cad7}(i) et (\cite{agt} III.4.2(iv)), 
$\oU^{(n)\circ}$ est un ouvert schématiquement dominant de $\oU^{(n)}$ (\cite{ega4} 11.10.5).
On a donc $\oU^{\oU^\circ(n)}=\oU^{(n)}$ en vertu de \ref{cad7}(iii).

On désigne par $\ocB^{(n)}$ le préfaisceau sur $E$ défini pour tout $(V\rightarrow U)\in \ob(E)$, par 
\begin{equation}\label{amtF2b}
\ocB^{(n)}(V\rightarrow U)=\Gamma(\oU^{V(n)},\co_{\oU^{V(n)}}).
\end{equation}
Comme $\oX^{(n)}$ est fini sur $\oX$, on a 
\begin{equation}\label{amtF2c}
\ocB^{(n)}(V\rightarrow U)=\ocB(V\times_{\oX^\circ}\oX^{(n)\circ}\rightarrow U).
\end{equation}
Par suite, $\ocB^{(n)}$ est un faisceau sur $E$, et 
on a un isomorphisme canonique d'anneaux sur $\tE$
\begin{equation}\label{amtF2d}
\ocB^{(n)}\stackrel{\sim}{\rightarrow} \Phi^{(n)}_*(\Phi^{(n)*}\ocB),
\end{equation}
où $\Phi^{(n)}$ est le morphisme de localisation \eqref{amtF1e}. Pour tout $U\in \ob(\Et_{/X})$, on pose 
\begin{equation}\label{amtF2e}
\ocB^{(n)}_{U}=\ocB^{(n)}\circ \iota_{U!},
\end{equation} 
où $\iota_{U!}$ est le foncteur \eqref{TFA6F}. 

Les $(\ocB^{(n)})_{n\geq 1}$ forment naturellement un système inductif d'anneaux de $\tE$, indexé par l'ensemble 
$\mZ_{\geq 1}$ ordonné par la relation de divisibilité. On pose 
\begin{equation}\label{amtF2f}
\ocB^{(\infty)}=\underset{\underset{n\geq 1}{\longrightarrow}}{\lim}\ \ocB^{(n)}.
\end{equation}
Pour tout $U\in \ob(\Et_{/X})$, on pose 
\begin{equation}\label{amtF2g}
\ocB^{(\infty)}_{U}=\ocB^{(\infty)}\circ \iota_{U!}.
\end{equation} 

D'après \ref{tf4}, le topos $\tE$ est cohérent. Supposons $U$ cohérent. 
Le faisceau associé à tout objet $(V\rightarrow U)$ de $E$ est alors cohérent d'après {\em loc. cit.}
Par suite, le morphisme canonique
\begin{equation}\label{amtF2h}
\underset{\underset{n\geq 1}{\longrightarrow}}{\lim}\ \ocB^{(n)}_U(V)\rightarrow \ocB^{(\infty)}_U(V)
\end{equation}
est un isomorphisme en vertu de (\cite{sga4} VI 5.2). Comme le schéma $\oU^\circ$ est cohérent, on en déduit 
par (\cite{agt} VI.9.12) que le morphisme canonique
\begin{equation}\label{amtF2i}
\underset{\underset{n\geq 1}{\longrightarrow}}{\lim}\ \ocB^{(n)}_U\rightarrow \ocB^{(\infty)}_U
\end{equation}
est un isomorphisme de $\oU^\circ_\fet$. 

\begin{lem}\label{amtF3}
Soit $(V\rightarrow U)$ un objet de $E$. Alors,
\begin{itemize}
\item[{\rm (i)}] Pour tout entier $n\geq 1$, 
le schéma $\oU^{V(n)}$ est normal et localement irréductible et le morphisme canonique 
$V\times_{\oU}\oU^{(n)}\rightarrow \oU^{V(n)}$ est une immersion ouverte schématiquement dominante. 
\item[{\rm (ii)}] Si $U$ est quasi-compact, le nombre de composantes irréductibles de $\oU^{V(n)}$ lorsque $n$
décrit les entiers $\geq 1$, est borné. 
\end{itemize}
\end{lem}

(i) En effet, $V$ étant entier sur $\oU^\circ$, 
le morphisme canonique $V\times_{\oU}\oU^{(n)}\rightarrow \oU^\circ \times_{\oU}\oU^{V(n)}$ est un isomorphisme.
Le morphisme canonique $V\times_{\oU}\oU^{(n)}\rightarrow \oU^{V(n)}$ est donc une immersion ouverte 
schématiquement dominante.
Le schéma $\oU^{(n)}$ est normal et localement irréductible d'après \ref{cad7}(iii) et (\cite{agt} III.3.3). 
Soit $Q$ un ouvert de $\oU^{(n)}$ n'ayant qu'un nombre fini de composantes irréductibles. 
Alors, $Q\times_{\oU^{(n)}}\oU^{V(n)}$ est la somme finie des fermetures intégrales de $Q$ dans 
les points génériques de $V\times_{\oU}\oU^{(n)}$ qui sont au-dessus de $Q$, 
dont chacune est un schéma intègre et normal en vertu de (\cite{ega2} 6.3.7). 
Le schéma $\oU^{V(n)}$ est donc normal et localement irréductible.

(ii) Il résulte de \ref{cad8}, en remplaçant $X$ par des ouverts affines qui couvrent $U$,  que 
le nombre de composantes irréductibles de $\oU^{(n)\circ}$ (pour $n\geq 1$) est borné.  La proposition s'ensuit 
compte tenu de (i) et du fait que le morphisme canonique $V\rightarrow \oU^\circ$ est étale et fini.

\subsection{}\label{amtF4}
On désigne par $\cC$ la sous-catégorie pleine de $\Et_{/X}$ formée des schémas affines, que l'on munit 
de la topologie induite par celle de $\Et_{/X}$. On voit aussitôt que $\cC$ est 
une famille $\mU$-petite, topologiquement génératrice du site $\Et_{/X}$ et est stable par produits fibrés. 
On désigne par 
\begin{equation}\label{amtF4a}
\pi_\cC\colon E_\cC\rightarrow \cC
\end{equation} 
le site fibré déduit de $\pi$ \eqref{TFA6a}
par changement de base par le foncteur d'injection canonique $\cC\rightarrow \Et_{/X}$. 
On munit $E_\cC$ de 
la topologie co-évanescente définie par $\pi_\cC$ et on note $\tE_\cC$ le topos des faisceaux de $\mU$-ensembles 
sur $E_\cC$. D'après (\cite{agt} VI.5.21 et VI.5.22), 
la topologie de $E_\cC$ est induite par celle de $E$ au moyen du foncteur de projection canonique 
$E_\cC\rightarrow E$, et celui-ci induit par restriction une équivalence de catégories 
\begin{equation}\label{amtF4b}
\tE\stackrel{\sim}{\rightarrow}\tE_\cC.
\end{equation}

\begin{lem}\label{amtF5}
Soit $(V\rightarrow U)$ un objet de $E_\cC$. Alors,
\begin{itemize}
\item[{\rm (i)}] Le schéma $\Spec(\ocB^{(\infty)}_U(V))$ est normal et localement irréductible. 
Il est la fermeture intégrale de $\oU$ dans $\Spec(\ocB_U^{(\infty)}(V))\times_XX^\circ$.
\item[{\rm (ii)}] Le diagramme de morphismes canoniques
\begin{equation}\label{amtF5a}
\xymatrix{
{\Spec(\ocB_U^{(\infty)}(V))\times_XX^\circ}\ar[r]\ar[d]&V\ar[d]\\
{\Spec(\ocB_U^{(\infty)}(\oU^\circ))\times_XX^\circ}\ar[r]&{\oU^\circ}}
\end{equation}
est cartésien. 
\end{itemize}
\end{lem}

(i) Pour tout entier $n\geq 1$, comme $X^{(n)}$ est fini sur $X$, $\oU^{V(n)}$ est la fermeture intégrale de 
$\oU$ dans $V\times_\oU\oU^{(n)}$.  D'après \ref{amtF3}(i), $\oU^{V(n)}$ est normal et localement irréductible.
Pour tout entier $m\geq 1$, le morphisme canonique $X^{(mn)\circ}\rightarrow X^{(n)\circ}$ étant étale, fini et surjectif, 
le morphisme canonique $\oU^{V(mn)}\rightarrow \oU^{V(n)}$ est entier et surjectif,  
et toute composante irréductible de $\oU^{V(mn)}$ domine une composante irréductible de 
$\oU^{V(n)}$. De plus, en vertu de \ref{amtF3}(ii), il existe un entier $n\geq 1$ tel que pour toute composante 
irréductible $C$ de $\oU^{V(n)}$ et tout entier $m\geq 1$, $C\times _{\oU^{V(n)}}\oU^{V(mn)}$ soit irréductible. 
La proposition résulte alors de (\cite{ega1n} 0.6.1.6 et 0.6.5.12(ii)).

(ii) En effet, $\Spec(\ocB_U^{(\infty)}(V))$ est la limite projective des schémas $(\oU^{V(n)})_{n\geq 1}$ (\cite{ega4} 8.2.3). 
Par suite,   $\Spec(\ocB_U^{(\infty)}(V))\times_XX^\circ$ est la limite projective des schémas 
$(V\times_\oU\oU^{(n)})_{n\geq 1}$ (\cite{ega4} 8.2.5); d'où la proposition.

\begin{prop}\label{amtF6}
Soient $(V\rightarrow U)$ un objet de $E_\cC$ tel que $V$ soit connexe, $\ov$ un point géométrique de $V$.  
Alors, il existe un revêtement universel $(V_i)_{i\in I}$ de $V$ en $\ov$ tel que pour tout $i\in I$, $V_i$ soit un revêtement étale galoisien de $V$
et que la $\ocB_U^{(\infty)}(V)$-algèbre $\ocB_U^{(\infty)}(V_i)$ soit $\alpha$-finie-étale \eqref{aet2}. 
\end{prop}

En effet, on peut se borner au cas où $U=X$.  Soient $Y$ un revêtement étale fini de $X^\circ$, $W=V\times_{X^\circ}Y$.
Considérons le diagramme de morphismes canoniques 
\begin{equation}\label{amtF6a}
\xymatrix{
{\Spec(\ocB^{(\infty)}_X(W))\times_XX^\circ}\ar[rd]\ar[ddd]\ar[rrr]
\ar@{}[rddd]|*+[o][F-]{1}&&&{\Spec(A_\infty)\times_XY}\ar[ddd]\ar[ld]\\
&W\ar[r]\ar[d]&Y\ar[d]&\\
&V\ar[r]&X^\circ&\\
{\Spec(\ocB^{(\infty)}_X(V))\times_XX^\circ}\ar[rrr]\ar[ru]&&&{\Spec(A_\infty)\times_XX^\circ}\ar[lu]}
\end{equation}
où $A_\infty$ est l'algèbre définie dans \eqref{cad9b}.
Notons $\Spec(B)$ la fermeture intégrale de $\Spec(A_\infty)$ dans $\Spec(A_\infty)\times_XY$. En vertu de \ref{tpcg2}, 
$B$ est une $A_\infty$-algèbre $\alpha$-finie-étale. 
Par suite, $B\otimes_{A_\infty}\ocB^{(\infty)}_X(V)$ est une $\ocB^{(\infty)}_X(V)$-algèbre $\alpha$-finie-étale (\cite{agt} V.7.4(3)).

Le carré $\xymatrix{\ar@{}|*+[o][F-]{1}}$ du diagramme \eqref{amtF6a} étant cartésien en vertu de 
\ref{amtF5}(ii), il en est de même du grand carré extérieur. On en déduit un diagramme cartésien
\begin{equation}
\xymatrix{
{\Spec(\ocB^{(\infty)}_X(W))\times_XX^\circ}\ar[r]\ar[d]&{\Spec(B\otimes_{A_\infty}\ocB^{(\infty)}_X(V)[\frac 1 p])}\ar[d]\\
{\Spec(\ocB^{(\infty)}_X(V))\times_XX^\circ}\ar[r]&{\Spec(\ocB^{(\infty)}_X(V)[\frac 1 p])}}
\end{equation}
Comme $B[\frac 1 p]$ est une $A_\infty[\frac 1 p]$-algèbre finie étale (\cite{agt} V.7.3), 
on en déduit en vertu de \ref{amtF5}(i) que $\ocB^{(\infty)}_X(W)$ est la fermeture intégrale de 
$\ocB^{(\infty)}_X(V)$ dans $B\otimes_{A_\infty}\ocB^{(\infty)}_X(V)[\frac 1 p]$.
Par suite, $\ocB^{(\infty)}_X(W)$ est une $\ocB^{(\infty)}_X(V)$-algèbre $\alpha$-finie-étale d'après \ref{aet3}.

Comme $\oX$ est normal et localement irréductible (\cite{agt} III.4.2(iii)), il en est de même de $W$ 
(\cite{agt} III.3.3). Par suite, $W$ est la réunion des schémas induits sur ses composantes irréductibles.
Pour toute composante irréductible $W'$ de $W$, 
$\ocB^{(\infty)}_X(W')$ est une $\ocB^{(\infty)}_X(V)$-algèbre $\alpha$-finie-étale d'après (\cite{agt} V.7.4).
On conclut la preuve en observant que l'homomorphisme $\pi_1(V,\ov)\rightarrow \pi_1(X^\circ,\ov)$ est injectif puisque $U=X$.

\subsection{}\label{amtF7}
Soient $U$ un objet de $\Et_{/X}$, $n$ un entier $\geq 1$. 
Avec les notations de \ref{tpcg4}, on pose 
\begin{equation}\label{amtF7a}
\oU^{\star\circ}=U\times_X\oX^{\star\circ} \ \ \ {\rm et}\ \ \ \oU^{(n)\star\circ}=U\times_X\oX^{(n)\star\circ},
\end{equation}
et on note $\pi_U^{(n)\star}\colon \oU^{(n)\star\circ}\rightarrow \oU^{\star\circ}$ le morphisme canonique.
On rappelle que $\oX^{(n)\star\circ}$ est un revêtement étale fini et galoisien de $\oX^{\star\circ}$ 
de groupe $\Delta_n$. On définit le foncteur 
\begin{equation}\label{amtF7b}
\oU^{\circ}_\fet\rightarrow \oU^{\star \circ}_\fet, \ \ \ F\mapsto F^{(n)\star}=(\pi_{U}^{(n)\star})_*(F|\oU^{(n)\star\circ}).
\end{equation}
Celui-ci est exact puisque le foncteur $(\pi_U^{(n)\star})_*$ l'est.

Soit $F$ un objet de $\oU^\circ_\fet$. D'après (\cite{agt} VI.9.4), pour tout $V\in \ob(\Et_{\rf/\oU^{\star \circ}})$, on a 
\begin{equation}\label{amtF7c}
F^{(n)\star}(V)=F(V\times_{\oX^{\star\circ}}\oX^{(n)\star\circ}).
\end{equation}
Cet ensemble est donc naturellement muni d'une action de $\Delta_n$, et on a un isomorphisme canonique 
\begin{equation}\label{amtF7d}
F(V) \stackrel{\sim}{\rightarrow}(F^{(n)\star}(V))^{\Delta_n}.
\end{equation}

Les faisceaux $(F^{(n)\star})_{n\geq 1}$  forment naturellement un système inductif de $\oU^{\star\circ}_\fet$, 
indexé par l'ensemble $\mZ_{\geq 1}$ ordonné par la relation de divisibilité. On pose 
\begin{equation}\label{amtF7e}
F^{(\infty)\star}=\underset{\underset{n\geq 1}{\longrightarrow}}{\lim}\ F^{(n)\star}.
\end{equation}

On voit aussitôt que le morphisme canonique $\oX^{\star\circ}\rightarrow X$ est cohérent \eqref{tpcg3}.
Supposons $U$ cohérent, de sorte que $\oU^{\star\circ}$ est cohérent.
D'après (\cite{agt} VI.9.12) et (\cite{sga4} VI 5.2), pour tout $V\in \ob(\Et_{\rf/\oU^{\star \circ}})$, on a alors un isomorphisme
canonique
\begin{equation}\label{amtF7f}
F^{(\infty)\star}(V)\stackrel{\sim}{\rightarrow}\underset{\underset{n\geq 1}{\longrightarrow}}{\lim}\ 
F^{(n)\star}(V).
\end{equation}
Par suite, l'ensemble discret $F^{(\infty)\star}(V)$ est naturellement muni d'une action continue de $\Delta_\infty$,
et on a un isomorphisme canonique
\begin{equation}\label{amtF7g}
F(V) \stackrel{\sim}{\rightarrow}(F^{(\infty)\star}(V))^{\Delta_\infty}.
\end{equation}

\subsection{}\label{amtF8}
Soient $U$ un objet de $\Et_{/X}$, $m\geq 0$ et $n\geq 1$ deux entiers. Appliquant le foncteur \eqref{amtF7b} aux anneaux 
$\ocB_U$ \eqref{TFA2d} et $\ocB_{U,m}=\ocB_U/p^m\ocB_U$ \eqref{TFA8b} de $\oU^\circ_\fet$, 
on obtient deux anneaux $\ocB^{(n)\star}_U$ et $\ocB^{(n)\star}_{U,m}$ de $\oU^{\star\circ}_\fet$. Ce foncteur étant exact, on a 
un isomorphisme canonique 
\begin{equation}\label{amtF8a}
\ocB^{(n)\star}_{U,m}\stackrel{\sim}{\rightarrow} \ocB^{(n)\star}_U/p^m\ocB^{(n)\star}_U.
\end{equation}
Compte tenu de \eqref{amtF2c}, l'injection canonique $\oX^{(n)\star}\rightarrow \oX^{(n)}$ induit un homomorphisme \eqref{amtF2e}
\begin{equation}\label{amtF8b}
\ocB^{(n)}_U|\oU^{\star \circ}\rightarrow \ocB^{(n)\star}_U.
\end{equation}
On en déduit par passage à la limite inductive un homomorphisme \eqref{amtF2g}
\begin{equation}\label{amtF8c}
\ocB^{(\infty)}_U|\oU^{\star \circ}\rightarrow \ocB^{(\infty)\star}_U.
\end{equation}
Par ailleurs, les isomorphismes \eqref{amtF8a} induisent un isomorphisme 
\begin{equation}\label{amtF8d}
\ocB^{(\infty)\star}_{U,m}\stackrel{\sim}{\rightarrow} \ocB^{(\infty)\star}_U/p^m\ocB^{(\infty)\star}_U.
\end{equation}

\subsection{}\label{amtF9}
On désigne par $E^\star_\cC$ la sous-catégorie pleine de $E_\cC$ \eqref{amtF4} formée des objets $(V\rightarrow U)$ 
tels que le morphisme canonique $V\rightarrow \oU^{\circ}$ se factorise à travers  $\oU^{\star\circ}$ \eqref{amtF7a}.
Le foncteur \eqref{amtF4a} induit un foncteur fibrant
\begin{equation}\label{amtF9a}
E^\star_\cC\rightarrow\cC.
\end{equation}

\begin{lem}\label{amtF10}
Soit $(V\rightarrow U)$ un objet de $E^\star_\cC$. Alors, 
\begin{itemize}
\item[{\rm (i)}] Le morphisme canonique $\Spec(\ocB_U^{(\infty)\star}(V))\rightarrow \Spec(\ocB_U^{(\infty)}(V))$ 
\eqref{amtF8b} est une immersion ouverte et fermée. 
\item[{\rm (ii)}] Le diagramme de morphismes canoniques
\begin{equation}
\xymatrix{
{\Spec(\ocB_U^{(\infty)\star}(V))}\ar[r]\ar[d]&{\Spec(\ocB_U^{(\infty)}(V))}\ar[d]\\
{\Spec(\ocB_U^{(\infty)\star}(\oU^\circ))}\ar[r]&{\Spec(\ocB_U^{(\infty)}(\oU^\circ))}}
\end{equation}
est cartésien. 
\end{itemize}
\end{lem}

(i) Pour tout entier $n\geq 1$, $\Spec(\ocB^{(n)\star}_U(V))$ est un ouvert et fermé 
de $\Spec(\ocB^{(n)}_U(V))=\oU^{V(n)}$. D'après \ref{amtF3}, il existe $n\geq1$ tel que pour tout $m\geq1$, 
le diagramme canonique
\begin{equation}
\xymatrix{
{\oX^{(mn)\star}}\ar[r]\ar[d]&{\oX^{(mn)}}\ar[d]\\ 
{\oX^{(n)\star}}\ar[r]&{\oX^{(n)}}}
\end{equation}
soit cartésien. Il en est donc de même du diagramme canonique
\begin{equation}
\xymatrix{
{\Spec(\ocB^{(mn)\star}_U(V))}\ar[r]\ar[d]&{\Spec(\ocB^{(mn)}_U(V))}\ar[d]\\ 
{\Spec(\ocB^{(n)\star}_U(V))}\ar[r]&{\Spec(\ocB^{(n)}_U(V))}}
\end{equation}
La proposition s'ensuit en vertu de (\cite{ega4} 8.3.12). 

(ii) En effet, pour tout entier $n\geq 1$, le diagramme 
\begin{equation}
\xymatrix{
{\Spec(\ocB^{(n)\star}_U(V))}\ar[r]\ar[d]&{\oU^{V(n)}}\ar[d]\\ 
{\oU^{(n)\star}}\ar[r]&{\oU^{(n)}}}
\end{equation}
est cartésien.

\begin{prop}\label{amtF11}
Soient $(V\rightarrow U)$ un objet de $E^\star_\cC$ tel que $V$ soit connexe,  
$\ov$ un point géométrique de $V$.  
Alors, il existe un revêtement universel $(V_i)_{i\in I}$ de $V$ en $\ov$ tel que pour tout $i\in I$, $V_i$ soit un revêtement étale galoisien de $V$ et 
que $\ocB_U^{(\infty)\star}(V_i)$ soit une $\ocB_U^{(\infty)\star}(V)$-algèbre $\alpha$-finie-étale. 
\end{prop}

Cela résulte de \ref{amtF6}, \ref{amtF10} et (\cite{agt} V.7.4(2) et V.7.8).

\begin{cor}\label{amtF12}
Soient $(V\rightarrow U)$ un objet de $E^\star_\cC$, $\cF$ un  $(\ocB^{(\infty)\star}_U|V)$-module de $V_\fet$, $q$ un entier $\geq 1$.  
Alors, $\rH^q(V_\fet, \cF)$ est $\alpha$-nul. 
\end{cor}

Comme $\oX$ est normal et localement irréductible (\cite{agt} III.4.2(iii)), il en est de même de $V$ (\cite{agt} III.3.3).
Par ailleurs, $j\colon X^\circ\rightarrow X$ étant quasi-compact, $V$ est quasi-compact. On peut donc supposer
$V$ connexe. Soit $\ov$ un point géométrique de $V$.
D'après \ref{amtF11}, il existe un revêtement universel $(V_i)_{i\in I}$ de $V$ en $\ov$ 
tel que pour tout $i\in I$, $V_i$ soit un revêtement étale galoisien de $V$ et 
que $\ocB_U^{(\infty)\star}(V_i)$ soit une $\ocB_U^{(\infty)\star}(V)$-algèbre $\alpha$-finie-étale.
Pour tout $i\in I$, notons $G_i$ le groupe des $V$-automorphismes de $V_i$. Il résulte de \ref{amtF5} et \ref{amtF10} que
$\Spec(\ocB_U^{(\infty)\star}(V))$ est normal et localement irréductible,
de fermeture intégrale $\Spec(\ocB_U^{(\infty)\star}(V_i))$ 
dans $\Spec(\ocB_U^{(\infty)\star}(V_i))\times_XX^\circ$ et que le diagramme 
\begin{equation}\label{amtF12a}
\xymatrix{
{\Spec(\ocB_U^{(\infty)\star}(V_i))\times_XX^\circ}\ar[r]\ar[d]&V_i\ar[d]\\
{\Spec(\ocB_U^{(\infty)\star}(V))\times_XX^\circ}\ar[r]&V}
\end{equation}
soit cartésien. En vertu de \ref{aet8}, $\ocB_U^{(\infty)\star}(V_i)$ est donc un $\alpha$-$G_i$-torseur de $\ocB_U^{(\infty)\star}(V)$. 
Par suite, $\rH^q(G_i,\cF(V_i))$ est $\alpha$-nul d'après \ref{aet7}. 
La proposition s'ensuit par passage à la limite. En effet, notant
\begin{equation}
\nu_\ov\colon V_\fet\stackrel{\sim}{\rightarrow} \bB_{\pi_1(V,\ov)},
\ \ \ M\mapsto \underset{\underset{i\in I}{\longrightarrow}}{\lim}\ M(V_i)
\end{equation}
le foncteur fibre en $\ov$ \eqref{notconv11c}, on a des isomorphismes canoniques 
\begin{equation}
\rH^q(V_\fet,\cF)\stackrel{\sim}{\rightarrow}\rH^q(\pi_1(V,\ov),\nu_\ov(\cF))\stackrel{\sim}{\rightarrow}
\underset{\underset{i\in I}{\longrightarrow}}{\lim}\ \rH^q(G_i,\cF(V_i)).
\end{equation}

\begin{cor}\label{amtF13}
Pour tout objet $(V\rightarrow U)$ de $E^\star_\cC$ et tout entier $m\geq 1$, 
l'homomorphisme canonique
\begin{equation}
\ocB^{(\infty)\star}_U(V)/p^m\ocB^{(\infty)\star}_U(V)\rightarrow \ocB^{(\infty)\star}_{U,m}(V)
\end{equation} 
est un $\alpha$-isomorphisme. 
\end{cor}

En effet, d'après \eqref{amtF8d}, comme $\ocB^{(\infty)\star}_U$ est $\mZ_p$-plat, 
on a une suite exacte de groupes abéliens de $\oU^{\star \circ}_\fet$
\begin{equation}
0\longrightarrow \ocB^{(\infty)\star}_U\stackrel{\cdot p^m}{\longrightarrow} \ocB^{(\infty)\star}_U\longrightarrow \ocB^{(\infty)\star}_{U,m}\longrightarrow 0.
\end{equation}
La proposition résulte alors du fait que $\rH^1(V_\fet, \ocB^{(\infty)\star}_U|V)$ est $\alpha$-nul en vertu de \ref{amtF12}.

\begin{cor}\label{amtF14}
Pour tout morphisme cartésien $(V'\rightarrow U')\rightarrow (V\rightarrow U)$ de $E^\star_\cC$ \eqref{amtF9a} 
et tout entier $m\geq 1$, l'homomorphisme canonique 
\begin{equation}
\ocB_{U,m}^{(\infty)\star}(V)\otimes_{\co_X(U)}\co_X(U')\rightarrow \ocB^{(\infty)\star}_{U',m}(V')
\end{equation}
est un $\alpha$-isomorphisme. 
\end{cor}
En effet, pour tout entier $n\geq 1$, l'homomorphisme canonique 
\begin{equation}
\ocB_U^{(n)\star}(V)\otimes_{\co_X(U)}\co_X(U')\rightarrow \ocB^{(n)\star}_{U'}(V')
\end{equation}
est un isomorphisme en vertu de (\cite{agt} III.8.11(iii)). 
On en déduit que l'homomorphisme canonique 
\begin{equation}
\ocB_U^{(\infty)\star}(V)\otimes_{\co_X(U)}\co_X(U')\rightarrow \ocB^{(\infty)\star}_{U'}(V')
\end{equation}
est un isomorphisme \eqref{amtF7f}. La proposition est donc une conséquence de \ref{amtF13}.

\begin{cor}\label{amtF15}
Soient $m$ un entier $\geq 1$, $\mL$ un $(\co_\oK/p^m\co_\oK)$-module localement libre de type fini de $\oX^\circ_\fet$. 
Pour tout $X$-schéma étale $U$, on note $\mL_U$ l'image inverse de $\mL$ par le morphisme canonique $\oU^\circ\rightarrow \oX^\circ$,
on pose $\cL_U=\mL_U\otimes_{\co_\oK}\ocB_U$ et on désigne par $\cL_U^{(\infty)\star}$ 
le $(\ocB_U^{(\infty)\star})$-module de $\oU^\circ_\fet$ associé à $\cL_U$ \eqref{amtF7e}. 
Alors, pour tout morphisme cartésien $(V'\rightarrow U')\rightarrow (V\rightarrow U)$ de $E^\star_\cC$, le morphisme canonique 
\begin{equation}
\cL_U^{(\infty)\star}(V)\otimes_{\co_X(U)}\co_X(U')\rightarrow \cL_{U'}^{(\infty)\star}(V')
\end{equation}
est un $\alpha$-isomorphisme. 
\end{cor} 

En effet, comme $\oX$ est normal et localement irréductible (\cite{agt} III.4.2(iii)), il en est de même de $V$ (\cite{agt} III.3.3).
Par ailleurs, $j\colon X^\circ\rightarrow X$ étant quasi-compact, $V$ est quasi-compact. 
On peut donc supposer $V$ connexe. 
Pour tout revêtement étale galoisien $V_1$ de $V$ de groupe $G$, posant $V'_1=V'\times_VV_1$, les morphismes canoniques
\begin{equation}
\cL^{(\infty)\star}_U(V)\rightarrow (\cL^{(\infty)\star}_U(V_1))^G\ \ \ {\rm et}\ \ \ \cL^{(\infty)\star}_{U'}(V')\rightarrow (\cL^{(\infty)\star}_{U'}(V'_1))^G
\end{equation} 
sont des isomorphismes. Compte tenu de (\cite{agt} II.3.15 et III.2.11), comme $\co_X(U')$ est plat sur $\co_X(U)$, 
on peut donc se réduire au cas où $\mL_U|V$ est constant de valeur un $(\co_\oK/p^m\co_\oK)$-module libre de type fini $M$.
On a un isomorphisme canonique \eqref{amtF7f}
\begin{equation}
\cL_U^{(\infty)\star}(V)\stackrel{\sim}{\rightarrow}M\otimes_{\co_\oK/p^m\co_\oK}\ocB_{U,m}^{(\infty)\star}(V);
\end{equation}
et de même pour $(V'\rightarrow U')$. La proposition résulte alors de \ref{amtF14}.

\begin{cor}\label{amtF16}
Soient $m$ un entier $\geq 1$, $\mL$ un $(\co_\oK/p^m\co_\oK)$-module localement libre de type fini de $\oX^\circ_\fet$. 
Pour tout $X$-schéma étale $U$, on note $\mL_U$ l'image inverse de $\mL$ par le morphisme canonique $\oU^\circ\rightarrow \oX^\circ$,
et on pose $\cL_U=\mL_U\otimes_{\co_\oK}\ocB_U$. Alors, 
pour tout morphisme cartésien $(V'\rightarrow U')\rightarrow (V\rightarrow U)$ de $E_\cC$ \eqref{amtF4a}, le morphisme canonique 
\begin{equation}
\cL_U(V)\otimes_{\co_X(U)}\co_X(U')\rightarrow \cL_{U'}(V')
\end{equation}
est un $\alpha$-isomorphisme. 
\end{cor}

En effet, on peut se borner au cas où $V$ est connexe (cf. la preuve de \ref{amtF15}). 
Quitte à changer de point géométrique $\oy$ de $\oX^\circ$ \eqref{tpcg3}, 
on peut supposer de plus que $(V\rightarrow U)$ est un objet de $E^\star_\cC$. 
La proposition résulte alors de \ref{amtF15} et \eqref{amtF7g} puisque $\co_X(U')$ est plat sur $\co_X(U)$ (\cite{agt} II.3.15).

\begin{cor}\label{amtF17}
Pour tout morphisme cartésien $(V'\rightarrow U')\rightarrow (V\rightarrow U)$ de $E_\cC$
et tout entier $m\geq 1$, l'homomorphisme canonique 
\begin{equation}
\ocB_{U,m}(V)\otimes_{\co_X(U)}\co_X(U')\rightarrow \ocB_{U',m}(V')
\end{equation}
est un $\alpha$-isomorphisme. 
\end{cor}

\begin{prop}\label{amtF18}
Soient $(V\rightarrow U)$ un objet de $E^\star_\cC$, $\cF$ un  $\ocB_U$-module de $\oU^\circ_\fet$,
$\cF^{(\infty)\star}$ le $(\ocB_U^{(\infty)\star})$-module de $\oU^{\star\circ}_\fet$ associé à $\cF$ \eqref{amtF7e}, 
$q$ un entier $\geq 0$. Alors, il existe un morphisme canonique
\begin{equation}\label{amtF18a}
\rH^q(\Delta_\infty,\cF^{(\infty)\star}(V))\rightarrow \rH^q(V_\fet,\cF|V),
\end{equation}
qui est un $\alpha$-isomorphisme.
\end{prop}

Soit $n$ un entier $\geq 1$. Avec les notations de \ref{amtF7}, on désigne par $\Gamma^{\Delta_n}$  le 
foncteur ``sous-faisceau des sections $\Delta_n$-invariantes'' sur la catégorie des $\mZ[\Delta_n]$-modules de $\oU^{\star\circ}_\fet$. 
Par descente, les foncteurs $(\pi_U^{(n)\star})^*$ et $\Gamma^{\Delta_n}\circ (\pi_U^{(n)\star})_*$ établissent des équivalences de 
catégories quasi-inverses l'une de l'autre entre la catégorie des groupes abéliens de $\oU^{\star\circ}_\fet$ et 
la catégorie des groupes abéliens de $\oU^{(n)\star\circ}_\fet$ sur lesquels $\Delta_n$ agit d'une façon compatible 
avec son action sur $\oU^{(n)\star \circ}$. On a donc un isomorphisme canonique de groupes abéliens de $\oU^{\star\circ}_\fet$
\begin{equation}
\cF|\oU^{\star\circ}\stackrel{\sim}{\rightarrow}\Gamma^{\Delta_n}(\cF^{(n)\star}).
\end{equation}
Le foncteur $(\pi_U^{(n)\star})_*$ étant exact, appliquant (\cite{verdier} II §~2 prop.~3.1), on obtient un isomorphisme de groupes abéliens
\begin{equation}
\cF|\oU^{\star\circ}\stackrel{\sim}{\rightarrow}\rR\Gamma^{\Delta_n}(\cF^{(n)\star}).
\end{equation}
Le foncteur $\Gamma^{\Delta_n}$ transforme les $\mZ[\Delta_n]$-modules injectifs de $V_\fet$ en groupes abéliens injectifs
de $V_\fet$, et le foncteur $\Gamma(V_\fet,-)$ transforme les $\mZ[\Delta_n]$-modules injectifs de $V_\fet$
en $\mZ[\Delta_n]$-modules injectifs (\cite{sga4} V 0.2). D'autre part, on a  
\begin{equation}
\Gamma(V_\fet, \Gamma^{\Delta_n}(-))=\Gamma^{\Delta_n}(\Gamma(V_\fet,-)).
\end{equation}
Prenant les foncteurs dérivés des deux membres et appliquant de nouveau (\cite{verdier} II §~2 prop.~3.1), on en déduit une suite spectrale 
\begin{equation}
{^{(n)}E}_2^{a,b}=\rH^a(\Delta_n,\rH^b(V_\fet,\cF^{(n)\star}|V))
\Rightarrow \rH^{a+b}(V_\fet,\cF|V).
\end{equation}
Comme l'immersion $j\colon X^\circ\rightarrow X$ est quasi-compacte, le schéma $V$ est cohérent.
Par suite, en vertu de (\cite{sga4} VI 5.3) et (\cite{agt} VI.9.12), pour tout entier $b\geq 0$, on a un isomorphisme canonique
\begin{equation}
\underset{\underset{n\geq 1}{\longrightarrow}}{\lim}\ \rH^b(V_\fet,\cF^{(n)\star}|V) 
\stackrel{\sim}{\rightarrow} \rH^b(V_\fet,\cF^{(\infty)\star}|V).
\end{equation}
On en déduit par passage à la limite inductive une suite spectrale 
\begin{equation}\label{amtF18b}
E_2^{a,b}=\rH^a(\Delta_\infty,\rH^b(V_\fet,\cF^{(\infty)\star}|V))\Rightarrow \rH^{a+b}(V_\fet,\cF|V).
\end{equation}
La proposition s'ensuit puisque $\rH^b(V_\fet, \cF^{(\infty)\star}|V)$ est $\alpha$-nul pour tout $b\geq 1$ en vertu de \ref{amtF12}.

\begin{cor}\label{amtF19}
Soient $m$, $q$ deux entiers tels que $m\geq 1$ et $q\geq 0$, $\mL$ un $(\co_\oK/p^m\co_\oK)$-module localement libre de type fini de $\oX^\circ_\fet$. 
Pour tout $X$-schéma étale $U$, on note $\mL_U$ l'image inverse de $\mL$ par le morphisme canonique $\oU^\circ\rightarrow \oX^\circ$,
et on pose $\cL_U=\mL_U\otimes_{\co_\oK}\ocB_U$. Alors, pour tout morphisme cartésien 
$(V'\rightarrow U')\rightarrow (V\rightarrow U)$ de $E_\cC$, le morphisme canonique
\begin{equation}\label{amtF19a}
\rH^q(V_\fet,\cL_U|V)\otimes_{\co_X(U)}\co_X(U')\rightarrow \rH^q(V'_\fet,\cL_{U'}|V')
\end{equation}
est un $\alpha$-isomorphisme.
\end{cor}

Pour tout $X$-schéma étale $U$, on note $\cL_U^{(\infty)\star}$ 
le $(\ocB_U^{(\infty)\star})$-module de $\oU^\circ_\fet$ associé à $\cL_U$ \eqref{amtF7e}. 
Supposons d'abord que $(V\rightarrow U)$ soit un objet de $E^\star_\cC$. D'après \ref{amtF18}, on a un morphisme canonique 
\begin{equation}\label{amtF19b}
\rH^q(\Delta_\infty,\cL_U^{(\infty)\star}(V))\rightarrow \rH^q(V_\fet,\cL_U|V),
\end{equation}
qui est un $\alpha$-isomorphisme. Celui-ci est fonctoriel en $(V\rightarrow U)\in \ob(E^\star_\cC)$ d'après la preuve de {\em loc. cit.}
Par ailleurs, d'après \ref{amtF15}, le morphisme canonique 
\begin{equation}
\cL_{U}^{(\infty)\star}(V)\otimes_{\co_X(U)}\co_X(U')\rightarrow \cL^{(\infty)\star}_{U'}(V')
\end{equation}
est un $\alpha$-isomorphisme. Comme $\co_X(U')$ est plat sur $\co_X(U)$, on en déduit que le morphisme canonique
\begin{equation}\label{amtF19c}
\rH^q(V_\fet,\cL_{U}|V)\otimes_{\co_X(U)}\co_X(U')\rightarrow \rH^q(V'_\fet,\cL_{U'}|V')
\end{equation}
est un $\alpha$-isomorphisme (\cite{agt} II.3.15). La proposition s'ensuit en variant le point géométrique $\oy$ de $\oX^\circ$ \eqref{tpcg3}.

\begin{cor}\label{amtF20}
Pour tout morphisme cartésien $(V'\rightarrow U')\rightarrow (V\rightarrow U)$ de $E_\cC$ 
et tous entiers $m\geq 1$ et $q\geq 0$, le morphisme canonique
\begin{equation}\label{amtF20a}
\rH^q(V_\fet,\ocB_{U,m}|V)\otimes_{\co_X(U)}\co_X(U')\rightarrow \rH^q(V'_\fet,\ocB_{U',m}|V')
\end{equation}
est un $\alpha$-isomorphisme.
\end{cor}

\begin{prop}\label{amtF21}
Soient $m$, $q$ deux entiers tels que $m\geq 1$ et $q\geq 0$, $\mL$ un $(\co_\oK/p^m\co_\oK)$-module localement libre de type fini de $\oX^\circ_\fet$,
$(V\rightarrow X)$ un objet de $E$. 
Pour tout $X$-schéma étale $U$, on note $\mL_U$ l'image inverse de $\mL$ par le morphisme canonique $\oU^\circ\rightarrow \oX^\circ$,
et on pose 
\begin{equation}
\cH^q(U)=\rH^q((V\times_XU)_\fet,\mL_U\otimes_{\co_\oK}\ocB_U|(V\times_XU)).
\end{equation} 
Alors, le préfaisceau de $\alpha$-$(\co_\oK)$-modules sur $\cC$ qui à tout objet $U$ de $\cC$ 
associe le $\alpha$-$(\co_\oK)$-module $\alpha(\cH^q(U))$ \eqref{alpha4} 
est un faisceau pour la topologie étale \eqref{alpha18}. 
\end{prop}

Soit $(U_i\rightarrow U)_{i\in I}$ un recouvrement de $\cC$. Pour tout $(i,j)\in I^2$, posons $U_{ij}=U_i\times_UU_j$. 
Par descente fidèlement plate, la suite de $\co_\oX(\oU)$-modules 
\begin{equation}
0\rightarrow \cH^q(U)\rightarrow \prod_{i\in I}\cH^q(U)\otimes_{\co_X(U)}\co_X(U_i)
\rightarrow \prod_{(i,j)\in I^2}\cH^q(U)\otimes_{\co_X(U)}\co_X(U_{ij}),
\end{equation}
où la dernière flèche est la différence des morphismes induits par les projections de $U_{ij}$ sur les deux facteurs, 
est exacte. On en déduit par \ref{amtF19} que la suite de $\co_\oK$-modules 
\begin{equation}
0\rightarrow \cH^q(U)\rightarrow \prod_{i\in I}\cH^q(U_i)\rightarrow \prod_{(i,j)\in I^2}\cH^q(U_{ij}),
\end{equation}
où la dernière flèche est la différence des morphismes induits par les projections de $U_{ij}$ 
sur les deux facteurs, est $\alpha$-exacte \eqref{finita3}. La proposition s'ensuit en vertu de \ref{alpha19}.

\subsection{}\label{amtF22}
Soit $(V\rightarrow X)$ un objet de $E$. On désigne par
\begin{equation}\label{amtF22a}
\pi^\dagger\colon E^\dagger\rightarrow \Et_{/X}
\end{equation}
le site fibré de Faltings associé au morphisme $V\rightarrow X$ (\cite{agt} VI.10.1). On munit $E^\dagger$ de la topologie co-évanescente associée à $\pi^\dagger$ 
(\cite{agt} VI.5.3) et on note $\tE^\dagger$ le topos des faisceaux de $\mU$-ensembles sur $E^\dagger$.
Tout objet de $E^\dagger$ est naturellement un objet de $E$. 
On définit ainsi un foncteur 
\begin{equation}\label{amtF22b}
\Phi\colon E^\dagger\rightarrow E
\end{equation}
qui induit une équivalence de catégories 
\begin{equation}\label{amtF22c}
E^\dagger\stackrel{\sim}{\rightarrow} E_{/(V\rightarrow X)}. 
\end{equation}
D'après (\cite{agt} VI.5.32), la topologie co-évanescente de $E^\dagger$ 
est induite par celle de $E$ au moyen du foncteur $\Phi$. Par suite, $\Phi$ est continu et cocontinu (\cite{sga4} III 5.2). 
Il définit donc une suite de trois foncteurs adjoints~:
\begin{equation}\label{amtF22d}
\Phi_!\colon \tE^\dagger\rightarrow \tE, \ \ \ \Phi^*\colon \tE\rightarrow \tE^\dagger, \ \ \ \Phi_*\colon \tE^\dagger\rightarrow \tE,
\end{equation}
dans le sens que pour deux foncteurs consécutifs de la suite, celui de droite est
adjoint à droite de l'autre. D'après (\cite{sga4} III 5.4), le foncteur $\Phi_!$ se factorise à travers 
une équivalence de catégories 
\begin{equation}\label{amtF22e}
\tE^\dagger\stackrel{\sim}{\rightarrow} \tE_{/(V\rightarrow X)^\tta},
\end{equation}
où $(V\rightarrow X)^\tta$ est le faisceau de $\tE$ associé à $(V\rightarrow X)$.
Le couple de foncteurs $(\Phi^*,\Phi_*)$ définit un morphisme de topos qu'on note aussi
\begin{equation}\label{amtF22f}
\Phi\colon \tE^\dagger\rightarrow \tE
\end{equation}
et qui n'est autre que le composé de \eqref{amtF22e} et du morphisme de localisation 
de $\tE$ en $(V\rightarrow X)^\tta$. Par ailleurs, le foncteur 
\begin{equation}\label{amtF22g}
\Phi^+\colon E\rightarrow E^\dagger,\ \ \ (W\rightarrow U)\mapsto (W\times_{\oX^\circ}V\rightarrow U), 
\end{equation}
est clairement un adjoint à droite du foncteur $\Phi\colon E^\dagger\rightarrow E$, et est continu et exact à gauche (\cite{agt} VI.10.12). 
Il définit donc le morphisme de topos $\Phi\colon \tE^\dagger\rightarrow \tE$ (\cite{sga4} III 2.5). 

Notons
\begin{equation}\label{amtF22h}
\sigma^\dagger\colon \tE^\dagger\rightarrow X_\et
\end{equation}
le morphisme canonique (\cite{agt} (VI.10.6.4)). Il résulte aussitôt des définitions (\cite{agt} VI.10.6) qu'on a un isomorphisme canonique 
\begin{equation}\label{amtF22i}
\sigma^\dagger\stackrel{\sim}{\rightarrow} \sigma\circ  \Phi,
\end{equation}
où $\sigma$ est le morphisme de topos \eqref{TFA6c}. 
On note $\cA$ la $\co_X$-algèbre de $X_\et$ associée à la $R$-algèbre $\ocB_X(V)$ \eqref{notconv12a}.
On a un morphisme canonique 
\begin{equation}
\cA\rightarrow \sigma^\dagger_*(\Phi^*(\ocB)).
\end{equation}
Sauf mention explicite du contraire, on considère $\sigma^\dagger$ comme un morphisme de topos annelés
\begin{equation}
\sigma^\dagger\colon (\tE^\dagger,\Phi^*(\ocB))\rightarrow (X_\et,\cA).
\end{equation}

\begin{prop}\label{amtF23}
Conservons les hypothèses de \ref{amtF22}, soient, de plus, $m\geq 1$ et $q\geq 0$ deux entiers, 
$\mL$ un $(\co_\oK/p^m\co_\oK)$-module localement libre de type fini de $\oX^\circ_\fet$. 
On note $M$ le $\ocB_X(V)$-module $\rH^q(V_\fet,\mL\otimes_{\co_\oK}\ocB_X|V)$ et $\tM$ le $\cA$-module associé de $X_\et$ \eqref{notconv12a}. 
Alors, on a un morphisme $\cA$-linéaire canonique de $X_\et$ 
\begin{equation}\label{amtF23a}
\tM\rightarrow \rR^q\sigma^\dagger_*(\Phi^*(\beta^*(\mL))),
\end{equation}
qui est un $\alpha$-isomorphisme, où $\beta$ est le morphisme de topos annelés \eqref{TFA2g}.
\end{prop}

Pour tout $X$-schéma étale $U$, on note $\mL_U$ l'image inverse de $\mL$ par le morphisme canonique $\oU^\circ\rightarrow \oX^\circ$.  
La correspondance $\cL=\{U\mapsto  \mL_U\otimes_{\co_\oK}\ocB_U\}$ définit un préfaisceau abélien sur $E$ \eqref{tf1h}. 
Notons $\cL^\tta$ le faisceau associé de $\tE$. D'après (\cite{agt} VI.5.34, VI.8.9 et VI.5.17), on a un isomorphisme canonique 
\begin{equation}
\beta^*(\mL)\stackrel{\sim}{\rightarrow}\cL^\tta. 
\end{equation}
Notons $(\cL\circ \Phi)^\tta$ le faisceau de $\tE^\dagger$ associé au préfaisceau $\cL\circ \Phi$ sur $E^\dagger$ \eqref{amtF22b}. 
En vertu de (\cite{sga4} III 2.3), on a un isomorphisme canonique 
\begin{equation}\label{amtF23i}
\Phi^*(\beta^*(\mL))\stackrel{\sim}{\rightarrow} (\cL\circ \Phi)^\tta.
\end{equation}
Le préfaisceau $\cL\circ \Phi$ sur $E^\dagger$ est défini par la correspondance $\{U\mapsto  \mL_U\otimes_{\co_\oK}\ocB_U|(V\times_XU)\}$ 
($U\in \ob(\Et_{/X})$) (\cite{agt} VI.5.2). 

Suivant (\cite{agt} VI.10.39), notons $\cH^q(\cL\circ \Phi)$ le faisceau de $X_\et$ associé au préfaisceau sur $\Et_{/X}$ 
défini pour tout $U\in \ob(\Et_{/X})$ par le groupe abélien $\rH^q((V\times_XU)_\fet,\mL_U\otimes_{\co_\oK}\ocB_U|(V\times_XU))$.  
En vertu de (\cite{agt} VI.10.40), on a un isomorphisme canonique de $X_\et$
\begin{equation}\label{amtF23j}
\cH^q(\cL\circ \Phi)\stackrel{\sim}{\rightarrow} \rR^q\sigma^\dagger_*(\Phi^*(\beta^*(\mL))).
\end{equation}
On observera que $\cH^q(\cL\circ \Phi)$ est naturellement un $\sigma^\dagger_*(\Phi^*(\ocB))$-module et on vérifie aussitôt que le morphisme 
\eqref{amtF23j} est $\sigma^\dagger_*(\Phi^*(\ocB))$-linéaire (\cite{agt} VI.10.39). 

Compte tenu de (\cite{sga4} II 3.0.4) et de la définition du foncteur ``faisceau associé'' (\cite{sga4} II 3.4),   
$\cH^q(\cL\circ \Phi)$ est le faisceau de $X_\et$ associé au préfaisceau sur $\cC$ \eqref{amtF4}
défini pour tout $U\in \ob(\cC)$ par le groupe abélien $\rH^q((V\times_XU)_\fet,\mL_U\otimes_{\co_\oK}\ocB_U|(V\times_XU))$.
Il résulte alors de \ref{alpha27} et \ref{amtF21}, que pour tout $U\in \ob(\cC)$, le morphisme canonique
\begin{equation}\label{amtF23m}
\rH^q((V\times_XU)_\fet,\mL_U\otimes_{\co_\oK}\ocB_U|(V\times_XU))\rightarrow \cH^q(\cL\circ \Phi)(U)
\end{equation}
est un $\alpha$-isomorphisme. On en déduit d'après \ref{amtF19} que le morphisme canonique
\begin{equation}\label{amtF23n}
M\otimes_R\co_X(U)\rightarrow \cH^q(\cL\circ \Phi)(U)
\end{equation}
est un $\alpha$-isomorphisme~; d'où la proposition.

\begin{cor}\label{amtF24}
Soient $m,q$ deux entiers tels que $m\geq 1$ et $q\geq 0$, $\mL$ un $(\co_\oK/p^m\co_\oK)$-module localement libre de type fini de $\oX^\circ_\fet$. 
On note $M$ le $\co_\oX(\oX)$-module 
$\rH^q(\oX^\circ_\fet,\mL\otimes_{\co_\oK}\ocB_X)$ et $\tM$ le $\hbar_*(\co_\oX)$-module associé de $X_\et$ \eqref{notconv12a}. 
Alors, on a un morphisme $\hbar_*(\co_\oX)$-linéaire canonique de $X_\et$ 
\begin{equation}\label{amtF24a}
\tM\rightarrow \rR^q\sigma_*(\beta^*(\mL)),
\end{equation}
qui est un $\alpha$-isomorphisme, où $\sigma$ et $\beta$ sont les morphismes de topos annelés \eqref{TFA2e} et \eqref{TFA2g}.
\end{cor}

\begin{cor}\label{amtF240}
Soient $m,q$ deux entiers tels que $m\geq 1$ et $q\geq 0$, $\mL$ un $(\co_\oK/p^m\co_\oK)$-module localement libre de type fini de $\oX^\circ_\fet$. 
On note $N$ le $\co_{\oX_m}(\oX_m)$-module $\rH^q(\oX^\circ_\fet,\mL\otimes_{\co_\oK}\ocB_X)$ 
et $\tN$ le $\co_{\oX_m}$-module associé de $\oX_{m,\et}$ \eqref{notconv12a}. 
Alors, on a un morphisme $\co_{\oX_m}$-linéaire canonique de $\oX_{m,\et}$  
\begin{equation}
\tN\rightarrow \rR^q\sigma_{m*}(\beta^*_m(\mL))
\end{equation}
qui est un $\alpha$-isomorphisme, où $\sigma_m$ et $\beta_m$ sont les morphismes de topos annelés \eqref{TFA8e} et \eqref{TFA8g}.
\end{cor}

En effet, on a $\beta^*(\mL)=\beta^{-1}(\mL)\otimes_{(\co_\oK/p^m\co_\oK)}\ocB_m$, où $\beta$ est le morphisme de topos annelés \eqref{TFA2g},
qui est donc un objet de $\tE_s$ d'après (\cite{agt} III.9.6 et III.9.7). 
Le morphisme canonique $\delta^*(\beta^*(\mL))\rightarrow \beta^*_m(\mL)$ est un isomorphisme, 
où $\delta$ est le plongement \eqref{TFA66a}. On en déduit par adjonction un isomorphisme
\begin{equation}
\beta^*(\mL)\stackrel{\sim}{\rightarrow} \delta_*(\beta^*_m(\mL)).
\end{equation} 
Il résulte alors de \eqref{TFA66d} qu'on a un isomorphisme canonique
\begin{equation}
\rR^q\sigma_{m*}(\beta^*_m(\mL))\stackrel{\sim}{\rightarrow}a^{-1}(\rR^q\sigma_*(\beta^*(\mL))),
\end{equation}
où $a^{-1}$ désigne l'image inverse au sens des faisceaux abéliens par l'injection canonique $a\colon X_s\rightarrow X$. 
Par ailleurs, avec les notations de \ref{amtF24}, on a un isomorphisme canonique 
\begin{equation}
\tM\stackrel{\sim}{\rightarrow} a_*(\tN).
\end{equation}
On en déduit par adjonction un isomorphisme $a^{-1}(\tM)\stackrel{\sim}{\rightarrow} \tN$. La proposition résulte donc de \ref{amtF24}.

\begin{rema}\label{amtF241}
Conservons les hypothèses de \ref{amtF240}. 
D'après \ref{aet100} et (\cite{egr1} 1.4.3), $\co_{\oX_m}$ est un faisceau d'anneaux cohérent de $\oX_{m,\zar}$. 
En vertu de \ref{afini104}, \ref{afini6} et \ref{tpcg10}(iii), le $\co_{\oX_m}$-module de $\oX_{m,\zar}$ associé à $N$ 
est $\alpha$-cohérent \eqref{finita9} et est $\alpha$-nul pour tout $i\geq d+1$, où $d=\dim(X/S)$. 
\end{rema}

\subsection{}\label{amtF25}
On désigne par $\hE$ (resp. $\hE_\cC$) la catégorie des préfaisceaux de $\mU$-ensembles sur $E$ (resp. $E_\cC$)
et par 
\begin{equation}\label{amtF25a}
\cP^\vee\rightarrow \Et_{/X}^\circ
\end{equation}
la catégorie fibrée obtenue en associant à tout $U\in \ob(\Et_{/X})$ la catégorie $(\Et_{\rf/\oU^\circ})^\wedge$ 
des préfaisceaux de $\mU$-ensembles
sur $\Et_{\rf/\oU^\circ}$, et à tout morphisme $f\colon U'\rightarrow U$ de $\Et_{/X}$ 
le foncteur 
\begin{equation}\label{amtF25b}
\of^\circ_{\fet*}\colon (\Et_{\rf/\oU'^\circ})^\wedge\rightarrow (\Et_{\rf/\oU^\circ})^\wedge
\end{equation} 
défini par la composition avec le foncteur image inverse $\Et_{\rf/\oU^\circ} \rightarrow \Et_{\rf/\oU'^\circ}$
par le morphisme $\of^\circ\colon \oU'^\circ\rightarrow \oU^\circ$ déduit de $f$. On note 
\begin{equation}\label{amtF25c}
\cP^\vee_\cC\rightarrow \cC^\circ
\end{equation}
la catégorie fibrée au-dessus de $\cC^\circ$ déduite de la catégorie fibrée \eqref{amtF25a} 
par changement de base par le foncteur d'injection canonique $\cC\rightarrow \Et_{/X}$.
On a alors une équivalence de catégories (\cite{agt} VI.5.2)
\begin{eqnarray}
\hE&\stackrel{\sim}{\rightarrow}& \bHom_{(\Et_{/X})^\circ}((\Et_{/X})^\circ,\cP^\vee)\\
F&\mapsto &\{U\mapsto F\circ \iota_{U!}\},\nonumber
\end{eqnarray}
où $\iota_{U!}$ est le foncteur \eqref{TFA6F}. 
On identifiera dans la suite $F$ à la section $\{U\mapsto F\circ \iota_{U!}\}$ qui lui est associée par cette équivalence.
De même, on a une équivalence de catégories
\begin{eqnarray}
\hE_\cC&\stackrel{\sim}{\rightarrow}& \bHom_{\cC^\circ}(\cC^\circ,\cP^\vee_\cC)\label{amtF25d}\\
F&\mapsto &\{U\mapsto F\circ \iota_{U!}\}.\nonumber
\end{eqnarray}

On note 
\begin{equation}\label{amtF25e}
\hE\rightarrow \tE, \ \ \ F\mapsto F^\tta.
\end{equation}
le foncteur ``faisceau associé'' sur $E$. 
Comme $E_\cC$ est une sous-catégorie topologiquement génératrice de $E$, le foncteur ``faisceau associé'' sur $E_\cC$
induit un foncteur que l'on note aussi 
\begin{equation}\label{amtF25f}
\hE_\cC\rightarrow \tE, \ \ \ F\mapsto F^\tta.
\end{equation}

Soient $F=\{W\mapsto G_W\}$ $(W\in \ob(\Et_{/X}))$ un objet de $\hE$, $F_\cC=\{U\mapsto G_U\}$ 
$(U\in \ob(\cC))$ l'objet de $\hE_\cC$ obtenu en restreignant $F$ à $E_\cC$. 
Il résulte aussitôt de (\cite{sga4} II 3.0.4) et de la définition du foncteur ``faisceau associé'' (\cite{sga4} II 3.4)
qu'on a un isomorphisme canonique de $\tE$
\begin{equation}\label{amtF25g}
(F_\cC)^\tta\stackrel{\sim}{\rightarrow} F^\tta.
\end{equation}

On désigne par $\ahE_\cC$ la catégorie des $\alpha$-$\co_\oK$-préfaisceaux sur $E_\cC$ \eqref{alpha14}, par 
$\atE_\cC$ la catégorie des $\alpha$-$\co_\oK$-faisceaux \eqref{alpha18} sur $E_\cC$ et par 
\begin{eqnarray}
\halpha\colon \bMod(\co_\oK,\hE_\cC)&\rightarrow& \ahE_\cC,\label{amtF25h}\\
\talpha\colon \bMod(\co_\oK,\tE_\cC)&\rightarrow& \atE_\cC,\label{amtF25i}
\end{eqnarray}
les foncteurs canoniques \eqref{alpha16a} et \eqref{alpha22a}. 

\begin{prop}\label{amtF26}
Soient $m$ un entier $\geq 1$, $\mL$ un $(\co_\oK/p^m\co_\oK)$-module localement libre de type fini de $\oX^\circ_\fet$. 
Pour tout $X$-schéma étale $U$, on note $\mL_U$ l'image inverse de $\mL$ par le morphisme canonique $\oU^\circ\rightarrow \oX^\circ$,
et on pose $\cL_U=\mL_U\otimes_{\co_\oK}\ocB_U$. Alors, 
le préfaisceau de $\alpha$-$\co_\oK$-modules $\halpha(\{U\mapsto \cL_U\})$ $(U\in \ob(\cC))$ \eqref{amtF25h} est un faisceau.
\end{prop}

Soient $(V\rightarrow U)$ un objet de $E_\cC$, $(U_i\rightarrow U)_{i\in I}$ un recouvrement de $\cC$. 
Pour tout $(i,j)\in I^2$, posons $U_{ij}=U_i\times_UU_j$, $V_i=V\times_UU_i$ et $V_{ij}=V\times_UU_{ij}$. 
Par descente fidèlement plate, la suite de $\co_X(U)$-modules 
\begin{equation}
0\rightarrow \cL_U(V)\rightarrow \prod_{i\in I}\cL_U(V)\otimes_{\co_X(U)}\co_X(U_i)
\rightarrow \prod_{(i,j)\in I^2}\cL_U(V)\otimes_{\co_X(U)}\co_X(U_{ij}),
\end{equation}
où la dernière flèche est la différence des morphismes induits par les projections de $U_{ij}$ sur les deux facteurs, 
est exacte. On en déduit par \ref{amtF16} que la suite de $\co_\oK$-modules 
\begin{equation}
0\rightarrow \cL_U(V)\rightarrow \prod_{i\in I}\cL_{U_i}(V_i)
\rightarrow \prod_{(i,j)\in I^2}\cL_{U_{ij}}(V_{ij}),
\end{equation}
où la dernière flèche est la différence des morphismes induits par les projections de $U_{ij}$ 
sur les deux facteurs, est $\alpha$-exacte \eqref{finita3}. La proposition s'ensuit en vertu de \ref{alpha181}.

\begin{cor}\label{amtF27}
Pour tout entier $m\geq 1$ et tout schéma affine $U$ de $\Et_{/X}$, l'homomorphisme canonique
\begin{equation}
\ocB_{U,m}\rightarrow \ocB_m\circ \iota_{U!},
\end{equation}
où $\ocB_m=\ocB/p^m\ocB$ \eqref{TFA8a} et $\iota_{U!}$ est le foncteur \eqref{TFA6F},  est un $\alpha$-isomorphisme.
\end{cor}
En effet, le préfaisceau de $\alpha$-$\co_\oK$-modules $\halpha(\{U\mapsto \ocB_{U,m}\})$ $(U\in \ob(\cC))$ \eqref{amtF25h}
est un faisceau en vertu de \ref{amtF26}.  La proposition résulte alors de \ref{alpha27}, compte tenu de \eqref{TFA8d} et \eqref{amtF25g}.

\begin{prop}\label{amtF28}
Reprenons les notations de \ref{tpcg3}. 
Soient $(V_i)_{i\in I}$ le revêtement universel normalisé de $\oX^{\star \circ}$ en $\oy$ \eqref{notconv11},
$m$ un entier $\geq 1$. Alors,
\begin{itemize}
\item[{\rm (i)}] Le morphisme canonique  
\begin{equation}\label{amtF28a}
\oR/p^m\oR\rightarrow  \underset{\underset{i\in I}{\longrightarrow}}{\lim}\ \ocB_m(V_i\rightarrow X)
\end{equation}
est un $\alpha$-isomorphisme.
\item[{\rm (ii)}] Pour tout entier $q\geq 1$, le $\co_C$-module $\underset{\underset{i\in I}{\longrightarrow}}{\lim}\ \rH^q((V_i\rightarrow X),\ocB_m)$ est $\alpha$-nul.
\end{itemize}
\end{prop}

(i) D'après \ref{amtF27}, pour tout $i\in I$, l'homomorphisme canonique 
\begin{equation}
\ocB_{X,m}(V_i)\rightarrow \ocB_m(V_i\rightarrow X)
\end{equation}
est un $\alpha$-isomorphisme. Il en est donc de même de l'homomorphisme canonique 
\begin{equation}
\underset{\underset{i\in I}{\longrightarrow}}{\lim}\ \ocB_{X,m}(V_i)\rightarrow \underset{\underset{i\in I}{\longrightarrow}}{\lim}\ \ocB_m(V_i\rightarrow X).
\end{equation}
La proposition s'ensuit en interprétant le terme de gauche comme la fibre du faisceau $\ocB_{X,m}$ de $\oX^{\circ}_\fet$ en $\oy$, et $\oR$ comme celle 
du faisceau $\ocB_X$ \eqref{TFA9b}. 

(ii) Montrons d'abord que la limite inductive
\begin{equation}\label{amtF28b}
\underset{\underset{i\in I}{\longrightarrow}}{\lim}\ \rH^q(V_{i,\fet}, \ocB_{X,m}|V_i)
\end{equation}
est nulle. Posons $G=\pi_1(\oX^{\star\circ},\oy)$ et soit $B_m$ la représentation continue de $G$ définie par $\ocB_{X,m}$. 
La limite inductive \eqref{amtF28b} est alors isomorphe à la limite inductive des groupes 
$\rH^q(H,B_m)$ prise sur les sous-groupes ouverts $H$ de $G$, qui est nulle puisqu'on a 
\begin{equation}
\rH^q(H,B_m) = \underset{\underset{H'\subset H}{\longrightarrow}}{\lim}\ \rH^q(H/H',B_m^{H'}),
\end{equation}
la limite étant prise sur les sous-groupes normaux ouverts $H'$ de $H$, et le composé des morphismes canoniques
\begin{equation}
\rH^q(H/H',B_m^{H'}) \rightarrow \rH^q(H,B_m^{H'}) \rightarrow \rH^q(H',B_m^{H'})
\end{equation}
est nul. 

Pour tout $i\in I$, on désigne par $\tE_i$ le topos de Faltings associé au morphisme $V_i\rightarrow X$, que l'on identifie au localisé du topos $\tE$ 
en $(V_i\rightarrow X)^\ra$ \eqref{amtF22}, et par 
\begin{equation}
\sigma_i\colon \tE_i\rightarrow X_\et
\end{equation}
le morphisme canonique \eqref{amtF22h}. Pour tout entier $q\geq 0$, on note $M^q_i$ le $\ocB_X(V_i)$-module $\rH^q(V_{i,\fet}, \ocB_{X,m}|V_i)$ \eqref{TFA8b}, 
$\cA_i$ la $\co_X$-algèbre de $X_\et$ associée à $\ocB_X(V_i)$ \eqref{notconv12a}, et $\tM^q_i$ le $\cA_i$-module de $X_\et$ associé à $M_i$.
En vertu \ref{amtF23}, on a un morphisme $\cA_i$-linéaire canonique de $X_\et$ 
\begin{equation}\label{amtF28c}
\tM^q_i\rightarrow \rR^q\sigma_{i*}(\ocB_m|(V_i\rightarrow X)^\ra),
\end{equation}
qui est un $\alpha$-isomorphisme. Considérons la suite spectrale de Cartan-Leray 
\begin{equation}
\rE_2^{a,b}=\rH^a(X_\et, \rR^b\sigma_{i*}(\ocB_m|(V_i\rightarrow X)^\ra))\Rightarrow \rH^{a+b}((V_i\rightarrow X),\ocB_m).
\end{equation}
Comme $X$ est affine, on déduit de \eqref{amtF28c} que le morphisme canonique 
\begin{equation}\label{amtF28d}
\rH^q((V_i\rightarrow X),\ocB_m)\rightarrow \Gamma(X, \rR^q\sigma_{i*}(\ocB_m|(V_i\rightarrow X)^\ra))
\end{equation}
est un $\alpha$-isomorphisme. La proposition s'ensuit en passant à la limite inductive sur $I$, 
puisque pour tout $q\geq 1$, $\underset{\underset{i\in I}{\longrightarrow}}{\lim}\ \tM^q_i$ est $\alpha$-nul, d'après \eqref{amtF28b}.

\section{\texorpdfstring{$\alpha$-finitude}{alpha-Finitude}}\label{finitude}

\begin{prop}\label{finitude1}
Soient $n,q$ deux entiers tels que $n\geq 1$ et $q\geq 0$, $\mL$ un $(\co_\oK/p^n\co_\oK)$-module localement libre de type fini de $\oX^\circ_\fet$. 
Alors, le $\co_{\oX_n}$-module $\fm_\oK\otimes_{\co_\oK}\rR^q\sigma_{n*}(\beta_n^*(\mL))$ \eqref{TFA8} est associé à un  $\co_{\oX_n}$-module
quasi-cohérent et $\alpha$-cohérent \eqref{finita9} de $\oX_{n,\zar}$ \eqref{notconv12a}, 
et il est nul pour tout $q\geq d+1$, où $d=\dim(X/S)$.
\end{prop}

En effet, soit $(u_j\colon X_j\rightarrow X)_{j\in J}$ un recouvrement étale tel que $J$ soit fini et que 
pour tout $j\in J$, $(X_j,\cM_X|X_j)$ vérifie les hypothèses de \ref{amtF1}. 
Soit $j\in J$. D'après (\cite{agt} (VI.10.12.6) et (III.9.11.12)), on a un diagramme de morphismes de topos annelés 
\begin{equation}\label{finitude1f}
\xymatrix{
{((\oX_j^\circ)_\fet,\co_\oK/p^n\co_\oK)}\ar[d]_{\ou^\circ_j}&{(\tE_{j,s},\ocB_{j,n})}\ar[l]_-(0.5){\beta_{j,n}}\ar[d]^{\Phi_j}\ar[r]^-(0.5){\sigma_{j,n}}&
{((X_{j,s})_{\et},\co_{\oX_{j,n}})}\ar[d]^{\ou_{j,s}}\\
{(\oX^\circ_\fet,\co_\oK/p^n\co_\oK)}&{(\tE_s,\ocB_n)}\ar[l]_-(0.5){\beta_n}\ar[r]^-(0.5){\sigma_n}&{(X_{s,\et},\co_{\oX_n})}}
\end{equation}
où la ligne supérieure est l'analogue pour $X_j$ de la ligne inférieure pour $X$, les flèches
verticales sont définies par fonctorialité, et dont les carrés sont commutatifs à isomorphismes canoniques près. 
Par ailleurs, le morphisme $\Phi_j$ s'identifie au morphisme de localisation de $(\tE_s,\ocB_n)$ en $\sigma_s^*(X_{j,s})$ 
en vertu de (\cite{agt} III.9.12 et III.9.13). 
En fait, $\sigma_{j,n}$ s'identifie au morphisme $(\sigma_n)_{/X_{j,s}}$, de sorte que le carré de droite de \eqref{finitude1f} est 2-cartésien 
(cf. \cite{sga4} IV 5.10).
Par suite, en vertu de \ref{amtF240} et \ref{amtF241}, il existe un $\co_{\oX_{j,n}}$-module quasi-cohérent et $\alpha$-cohérent
$\cF_j$ de $(\oX_{j,n})_\zar$ et un morphisme $\co_{\oX_{j,n}}$-linéaire
\begin{equation}\label{finitude1a}
\iota_j(\cF_j)\rightarrow  \rR^q\sigma_{n*}(\beta_n^*(\mL))|X_{j,s},
\end{equation}
qui est un $\alpha$-isomorphisme, où $\iota_j(\cF_j)$ est le $\co_{\oX_{j,n}}$-module de $(\oX_{j,n})_\et$ associé à $\cF_j$ \eqref{notconv12a}.
En particulier, le morphisme induit
\begin{equation}\label{finitude1b}
\fm_\oK\otimes_{\co_\oK}\iota_j(\cF_j)\rightarrow  \fm_\oK\otimes_{\co_\oK}\rR^q\sigma_{n*}(\beta_n^*(\mL))|X_{j,s}
\end{equation}
est un isomorphisme \eqref{alpha21}. 
Par suite, $\fm_\oK\otimes_{\co_\oK}\rR^q\sigma_{n*}(\beta_n^*(\mL))$ est nul pour tout $q\geq d+1$, d'après \ref{tpcg10}(iii).
Compte tenu de \eqref{notconv12e}, on a un isomorphisme canonique 
\begin{equation}\label{finitude1c}
\fm_\oK\otimes_{\co_\oK}\iota_j(\cF_j)\stackrel{\sim}{\rightarrow} \iota_j(\fm_\oK\otimes_{\co_\oK}\cF_j).
\end{equation}
Les isomorphismes \eqref{finitude1b} définissent une donnée de descente sur $(\iota_j(\fm_\oK\otimes_{\co_\oK}\cF_j))_{j\in J}$ relativement 
au recouvrement étale $(\oX_{j,n}\rightarrow \oX_n)_{j\in J}$. Comme le foncteur \eqref{notconv12a} est pleinement fidèle 
et compte tenu de \eqref{notconv120c}, 
on en déduit une donnée de descente sur $(\fm_\oK\otimes_{\co_\oK}\cF_j)_{j\in J}$ relativement au même recouvrement. 
D'après (\cite{sga1} VIII 1.1), il existe alors un $\co_{\oX_n}$-module quasi-cohérent $\cG$ et pour tout $j\in J$, 
un isomorphisme
\begin{equation}\label{finitude1d}
\cG\otimes_{\co_{\oX_n}}\co_{\oX_{j,n}}\stackrel{\sim}{\rightarrow} \fm_\oK\otimes_{\co_\oK}\cF_j
\end{equation} 
tels que si l'on munit $(\cG\otimes_{\co_{\oX_n}}\co_{\oX_{j,n}})_{j\in J}$ de la donnée de descente définie par $\cG$, 
les isomorphismes \eqref{finitude1d} soient compatibles aux données de descente.
Les isomorphismes \eqref{finitude1b} se descendent en 
un isomorphisme 
\begin{equation}\label{finitude1e}
\iota(\cG)\rightarrow  \fm_\oK\otimes_{\co_\oK}\rR^q\sigma_{n*}(\beta_n^*(\mL)),
\end{equation}
où $\iota(\cG)$ est le $\co_{\oX_n}$-module de $\oX_{n,\et}$ associé à $\cG$.
Par ailleurs, $\cG$ est $\alpha$-cohérent, d'après \ref{finita16} et \ref{afini12}; d'où la proposition.

\begin{cor}\label{finitude2}
Soient $n,q$ deux entiers tels que $n\geq 1$ et $q\geq 0$, $\mL$ un $(\co_\oK/p^n\co_\oK)$-module localement libre de type fini de $\oX^\circ_\fet$. 
Alors, le $\co_{\oX_n}$-module $\rR^q\tau_{n*}(\beta_n^*(\mL))$ \eqref{TFA8} est $\alpha$-quasi-cohérent \eqref{afini3} et $\alpha$-cohérent \eqref{finita9}, 
et est $\alpha$-nul pour tout $q\geq d+1$, où $d=\dim(X/S)$.
\end{cor}

En effet, pour tout $\co_{\oX_n}$-module $\cF$ de $\oX_{n,\et}$ et pour tout entier $i\geq 0$, les morphismes canoniques \eqref{TFA8h}
\begin{equation}\label{finitude2a}
\fm_\oK\otimes_{\co_\oK}\rR^ju_{n*}(\cF)\rightarrow \rR^iu_{n*}(\fm_{\co_\oK}\otimes_{\co_\oK}\cF)\rightarrow 
\rR^ju_{n*}(\cF)
\end{equation}
sont des $\alpha$-isomorphismes d'après \ref{alpha3}. Par suite, 
$\rR^iu_{n*}(\rR^j\sigma_{n*}(\beta_n^*(\mL)))$ est $\alpha$-nul pour tous $i\geq 1$ et $j\geq 0$, en vertu de \ref{finitude1} et (\cite{sga4} VII 4.3). 
La suite spectrale de Cartan-Leray
\begin{equation}\label{finitude2b}
E_2^{i,j}=\rR^iu_{n*}(\rR^j\sigma_{n*}(\beta_n^*(\mL)))\Rightarrow \rR^{i+j}\tau_{n*}(\beta_n^*(\mL))
\end{equation}
induit donc un isomorphisme de $\alpha$-$\co_{\oX_n}$-modules de $\oX_{n,\zar}$ \eqref{alpha20}
\begin{equation}\label{finitude2c}
\alpha(u_{n*}(\fm_\oK\otimes_{\co_\oK}\rR^q\sigma_{n*}(\beta_n^*(\mL))))\stackrel{\sim}{\rightarrow} \alpha(\rR^q\tau_{n*}(\beta_n^*(\mL))).
\end{equation}
La proposition s'ensuit compte tenu de \ref{finitude1} et \ref{finita16}.

\begin{cor}\label{finitude3}
Supposons $X$ propre sur $S$, et soient $n,q$ deux entiers tels que $n\geq 1$ et $q\geq 0$, 
$\mL$ un $(\co_\oK/p^n\co_\oK)$-module localement libre de type fini de $\oX^\circ_\fet$. 
Alors, le $\co_\oK$-module $\rH^q(\tE_s,\beta_n^*(\mL))$ est de présentation $\alpha$-finie.
\end{cor}

On notera d'abord que l'anneau $\co_\oK/p^n\co_\oK$ est universellement cohérent d'après \ref{aet100}.
Considérons la suite spectrale de Cartan Leray
\begin{equation}\label{finitude3a}
E_2^{i,j}=\rH^i(\oX_n,\rR^j\tau_{n*}(\beta_n^*(\mL)))\Rightarrow \rH^{i+j}(\tE_s,\beta_n^*(\mL)).
\end{equation}
En vertu de \ref{finitude2} et \ref{afini9}, pour tous entiers $i,j\geq 0$, le $(\co_\oK/p^n\co_\oK)$-module $E_2^{i,j}$ 
est $\alpha$-cohérent. La proposition s'ensuit compte tenu de \ref{finita18}.

\section{\texorpdfstring{Principal théorème de comparaison $p$-adique de Faltings: cas absolu}
{Principal théorème de comparaison p-adique de Faltings: cas absolu}}\label{TPCF}

\subsection{}\label{TFA101}
Suivant \eqref{eipo3a}, on désigne par $(\co_\oK)^\flat$ la limite projective du système projectif $(\co_\oK/p\co_\oK)_{\mN}$ 
dont les morphismes de transition sont les itérés de l'endomorphisme de Frobenius de $\co_\oK/p\co_\oK$.
\begin{equation}\label{TFA101a}
(\co_\oK)^\flat= \underset{\underset{\mN}{\longleftarrow}}{\lim}\co_\oK/p\co_\oK.
\end{equation}
Considérons le système projectif de monoïdes $(\co_C)_{\mN}$ 
dont les morphismes de transition sont les itérés de l'élévation à la puissance $p$ de $\co_C$. 
D'après \ref{eip5}(i), l'application 
\begin{equation}\label{TFA101b}
\underset{\underset{\mN}{\longleftarrow}}{\lim}\ \co_C \rightarrow (\co_\oK)^\flat, \ \ \  
(x_n)_{n\geq 0}\mapsto (\ox_n)_{n\geq 0},
\end{equation}
où $\ox_n$ désigne la classe de $x_n$ dans $\co_\oK/p\co_\oK$, est un isomorphisme de monoïdes multiplicatifs. 
Par suite, $(\co_\oK)^\flat$ est un anneau intègre, parfait de caractéristique $p$. 
Pour tout $x\in (\co_\oK)^\flat$, d'image inverse $(x_n)_{n\geq 0}$ par \eqref{TFA101b}, on pose  
\begin{equation}\label{TFA101c}
x^\sharp=x_0.
\end{equation}

On fixe une suite $(p_n)_{n\geq 0}$ d'éléments de $\co_\oK$ telle que $p_0=p$ et $p_{n+1}^p=p_n$ pour tout $n\geq 0$
et on note $\varpi$ l'élément associé de $(\co_{\oK})^\flat$. On pose 
\begin{equation}\label{TFA101e}
\oK^\flat=(\co_{\oK})^\flat[\frac{1}{\varpi}].
\end{equation} 
L'isomorphisme \eqref{TFA101b} se prolonge alors en un isomorphisme de monoïdes multiplicatifs
\begin{equation}\label{TFA101f}
\underset{\underset{\mN}{\longleftarrow}}{\lim}\ C \stackrel{\sim}{\rightarrow} \oK^\flat.
\end{equation}
Par conséquent, $\oK^\flat$ est un corps. C'est le corps des fractions de $(\co_{\oK})^\flat$. 
Il est algébriquement clos en vertu de (\cite{scholze} 3.8).

L'application $(\co_{\oK})^\flat\rightarrow \co_C$, $x\mapsto x^\sharp$ \eqref{TFA101c} se prolonge en une application multiplicative
\begin{equation}
\oK^\flat\rightarrow C, \ \ \ x\mapsto x^\sharp.
\end{equation}
On vérifie aussitôt que l'application 
\begin{equation}\label{TFA101d}
v_{\oK^\flat}\colon \oK^\flat\rightarrow \mQ\cup \{\infty\}, \ \ \ x\mapsto v_{\oK^\flat}(x)=v(x^\sharp),
\end{equation}
où $v$ est la valuation normalisée de $C$ \eqref{TFA1}, est une valuation non discrète, de hauteur $1$ d'anneau de valuation $(\co_\oK)^\flat$. 
Il est alors naturel de noter $(\co_\oK)^\flat$ aussi $\co_{\oK^\flat}$, et son idéal maximal $\fm_{\oK^\flat}$. 

Comme $\co_{\oK^\flat}$ satisfait les conditions requises dans \ref{mptf1}, 
il est loisible de considérer les notions de $\alpha$-algèbre introduites dans les sections \ref{alpha}--\ref{mptf} relativement à $\co_{\oK^\flat}$.

\subsection{}\label{TFA102}
Pour tout entier $i\geq 0$, l'homomorphisme 
\begin{equation}\label{TFA102a}
\pi_i \colon \co_{\oK^\flat}\rightarrow \co_\oK/p\co_\oK, \ \ \ (x_n)_{n\geq 0}\mapsto x_i
\end{equation}
est surjectif. On vérifie aussitôt que $\ker(\pi_0)=\varpi \co_{\oK^\flat}$. On note $\varphi$ l'endomorphisme de Frobenius de $\co_{\oK^\flat}$.
Comme $\pi_0=\pi_i\circ \varphi^i$, on en déduit que $\ker(\pi_i)=\varpi^{p^i}\co_{\oK^\flat}$. En particulier, $\pi_i$ induit un isomorphisme
\begin{equation}\label{TFA102b}
\co_{\oK^\flat}/\varpi^{p^i}\co_{\oK^\flat}\stackrel{\sim}{\rightarrow}\co_\oK/p\co_\oK.
\end{equation}
Par ailleurs, le diagramme 
\begin{equation}\label{TFA102c}
\xymatrix{
{\co_{\oK^\flat}/\varpi^{p^{i+1}}\co_{\oK^\flat}}\ar[r]^-(0.5)\sim\ar[d]&{\co_\oK/p\co_\oK}\ar[d]\\
{\co_{\oK^\flat}/\varpi^{p^i}\co_{\oK^\flat}}\ar[r]^-(0.5)\sim&{\co_\oK/p\co_\oK}}
\end{equation}
où la flèche verticale de gauche est le morphisme canonique et la flèche verticale de droite est l'endomorphisme de Frobenius, est commutatif.
On en déduit que $\co_{\oK^\flat}$ est complet et séparé pour la topologie $\varpi$-adique, et que cette dernière 
coïncide avec la topologie limite projective des topologies discrètes de $\co_\oK/p\co_\oK$.

On observera que pour toute $(\co_\oK/p\co_\oK)$-algèbre, les notions de $\alpha$-algèbre introduites dans les sections \ref{alpha}--\ref{aet}
relativement à $\co_{\oK}$ coïncident avec celles relativement à $\co_{\oK^\flat}$ via l'isomorphisme 
$\co_{\oK^\flat}/\varpi\co_{\oK^\flat}\stackrel{\sim}{\rightarrow}\co_\oK/p\co_\oK$ \eqref{TFA102b}.

\subsection{}\label{TPCF6}
On désigne par $\tE_s^{\mN^\circ}$ le topos des systèmes projectifs de $\tE_s$, indexés par l'ensemble ordonné $\mN$ des entiers naturels \eqref{notconv13}, et par 
\begin{equation}\label{TPCF6a}
\uplambda\colon \tE_s^{\mN^\circ}\rightarrow \tE_s
\end{equation}
le morphisme canonique \eqref{notconv13a}.
Pour tout objet $F=(F_n)_{n\geq 0}$ de $\tE_s^{\mN^\circ}$ et tout entier $i\geq 0$, on note $\varepsilon_{\leq i}(F)$ l'objet $(F_n^i)_{n\geq 0}$  
de $\tE_s^{\mN^\circ}$, défini par 
\begin{equation}\label{TPCF6b}
F_n^i=\left\{
\begin{array}{clcr}
F_n&{\rm si}\ 0\leq n\leq i,\\
F_i&{\rm si}\ n\geq i,
\end{array}
\right.
\end{equation}
le morphisme de transition $F_{n+1}^i\rightarrow F_n^i$ étant égal au morphisme de transition $F_{n+1}\rightarrow F_n$ de $F$ si $0\leq n<i$ et à l'identité 
si $n\geq i$. La correspondance $F\mapsto \varepsilon_{\leq i}(F)$ est clairement fonctorielle, et on a un morphisme fonctoriel canonique
\begin{equation}\label{TPCF6c}
F\rightarrow  \varepsilon_{\leq i}(F).
\end{equation}
Si $F$ est un anneau (resp. module), il en est de même de $\varepsilon_{\leq i}(F)$ et le morphisme \eqref{TPCF6c} est un homomorphisme. 

Si $A$ est un anneau, on note encore $A$ l'anneau constant de valeur $A$ de $\tE_s$ ou de $\tE_s^{\mN^\circ}$.
On a donc un isomorphisme canonique $\uplambda^*(A)\stackrel{\sim}{\rightarrow} A$ de $\tE_s^{\mN^\circ}$.
On note $(\co_\oK)^\wp$ l'anneau $(\co_\oK/p\co_\oK)_{\mN}$ de $\tE_s^{\mN^\circ}$ dont les morphismes de transition sont les itérés
de l'endomorphisme de Frobenius de $\co_\oK/p\co_\oK$. On le considère comme  une $\co_{\oK^\flat}$-algèbre 
via les homomorphismes \eqref{TFA102b}.

On désigne par $\cP$ la catégorie abélienne 
des $\mF_p$-modules de $\tE_s^{\mN^\circ}$ et par $\cP_\AR$ son quotient par la sous-catégorie 
épaisse des $\mF_p$-modules Artin-Rees-nuls (AR-nuls en abrégé) (\cite{sga5} V 2.2.1). On rappelle que $\cP_\AR$
est canoniquement équivalente à la catégorie des systèmes projectifs de $\cP$ à translation près (\cite{sga5} V 2.4.2 et 2.4.4).

\subsection{}\label{TPCF7}
On note $\phi$ l'endomorphisme de Frobenius de $\ocB_1$ \eqref{TFA8a}
et on désigne par $\ocB^\wp$ la $(\co_\oK)^\wp$-algèbre $(\ocB_1)_{\mN}$ de $\tE_s^{\mN^\circ}$  
dont les morphismes de transition sont les itérés de $\phi$ \eqref{TPCF6}. 
On considère $\ocB^\wp$ comme une $\co_{\oK^\flat}$-algèbre.
On note encore $\phi$ l'endomorphisme de Frobenius de $\ocB^\wp$.
Cette notation n'induit aucun risque d'ambiguïté puisque l'endomorphisme de Frobenius de $\ocB^\wp$ est induit par celui de $\ocB_1$. 
Pour tout entier $i\geq 0$, on considère $\varepsilon_{\leq i}(\ocB^\wp)$ comme une $\ocB^\wp$-algèbre via l'homomorphisme canonique \eqref{TPCF6c}.

\begin{lem}\label{TPCF8}
\
\begin{itemize}
\item[{\rm (i)}] Pour tout entier $n\geq 0$, notant $\op_n$ la classe de $p_n$ dans $\co_\oK/p\co_\oK$ \eqref{TFA101}, la suite
\begin{equation}\label{TPCF8a}
\xymatrix{
{\ocB_1}\ar[r]^-(0.5){\cdot \op_n}&{\ocB_1}\ar[r]^{\phi^n}&{\ocB_1}\ar[r]&0}
\end{equation}
est exacte.
\item[{\rm (ii)}] L'endomorphisme de Frobenius $\phi$ de $\ocB^\wp$ induit un isomorphisme de $\cP_\AR$ \eqref{TPCF6}.
\item[{\rm (iii)}] Pour tout entier $i\geq 0$, la suite de $\ocB^\wp$-modules
\begin{equation}\label{TPCF8b}
\xymatrix{
0\ar[r]&{\ocB^\wp}\ar[r]^-(0.5){\cdot \varpi^{p^i}}&{\ocB^\wp}\ar[r]&{\varepsilon_{\leq i}(\ocB^\wp)}\ar[r]&0}
\end{equation}
où la deuxième flèche est le morphisme canonique \eqref{TPCF6c}, est exacte dans $\cP_\AR$. De plus, privée du zéro de gauche, elle est exacte dans $\cP$.
\end{itemize}
\end{lem}

(i) En effet, soient $(\oy\rightsquigarrow \ox)$ un point de $X_\et\gtimes_{X_\et}\oX^\circ_\et$ \eqref{topfl17} 
tel que $\ox$ soit au-dessus de $s$.
On désigne par $\fV_\ox$ la catégorie des $X$-schémas étales $\ox$-pointés, 
et par $\fW_\ox$ la sous-catégorie pleine de $\fV_\ox$ formée des objets $(U,\fp\colon \ox\rightarrow U)$ 
tels que la restriction $(U,\cM_X|U)\rightarrow (S,\cM_S)$ du morphisme $f$ \eqref{TFA3} vérifie les hypothèses de \ref{tpcg1}.
Ce sont des catégories cofiltrantes, et le foncteur d'injection canonique
$\fW^\circ_\ox\rightarrow \fV^\circ_\ox$ est cofinal (\cite{sga4} I 8.1.3(c)).
Avec les notations de \ref{TFA11}, on a donc un isomorphisme canonique \eqref{TFA11c} 
\begin{equation}
\ocB_{\rho(\oy\rightsquigarrow \ox)} \stackrel{\sim}{\rightarrow} 
\underset{\underset{(U,\fp)\in \fW^\circ_\ox}{\longrightarrow}}{\lim}\ \oR^{\oy}_U.
\end{equation} 
La proposition s'ensuit en vertu de \ref{tpcg7} et \ref{tf21}(ii).

(ii) Il résulte de (i) et (\cite{agt} III.7.3(i)) que l'endomorphisme $\phi$ de $\ocB^\wp$ est surjectif, de noyau $(\op_1\ocB_1)_{\mN}$, où les morphismes
de transition sont induits par l'endomorphisme de Frobenius $\phi$ de $\ocB_1$, et sont donc nuls. 

(iii) Pour tout entier $j<0$, posons $p_j=p$ et $\op_j=0$. On notera que $\varpi^{p^i}$ correspond à l'élément $(\op_{n-i})_{n\geq 0}$ via l'égalité \eqref{TFA101a}.
Compte tenu de (i) et (\cite{agt} III.7.3(i)), la suite privée du zéro de gauche est donc exacte dans $\cP$. 
Pour tout entier $n$, posons $q_n=p/p_n\in \co_\oK$. Comme $\ocB$ est $\co_\oK$-plat,
le noyau de la multiplication par $\varpi^{p^i}$ dans $\ocB^\wp$ est  $(q_{n-i}\ocB_1)_{n\geq 0}$, 
où les morphismes de transition sont induits par $\phi$ et sont donc nuls en degrés $\geq i$ puisque $v(q_n)=1-1/p^n$; d'où la proposition.

\begin{prop}\label{TPCF9}
Supposons $X$ propre sur $S$ et soient  $\mL$ un $\mF_p$-module de type fini de $\oX^\circ_\fet$, $q$ un entier $\geq 0$. 
Posons $\cL=\delta^*(\beta^*(\mL))$, où $\beta$ est le morphisme \eqref{TFA6b} et $\delta$ est le plongement \eqref{TFA66a},
et notons $\uphi$ l'endomorphisme $\co_{\oK^\flat}$-semi-linéaire de $\rH^q(\tE_s^{\mN^\circ},\uplambda^*(\cL)\otimes_{\mF_p}\ocB^\wp)$
induit par l'endomorphisme de Frobenius $\phi$ de $\ocB^\wp$. Alors, il existe un entier $r\geq 0$ et un morphisme $\co_{\oK^\flat}$-linéaire 
\begin{equation}\label{TPCF9a}
\fm_{\oK^\flat}\otimes_{\co_{\oK^\flat}}\rH^q(\tE_s^{\mN^\circ},\uplambda^*(\cL)\otimes_{\mF_p}\ocB^\wp)\rightarrow \co_{\oK^\flat}^r,
\end{equation}
compatible à $\uphi$ et à l'endomorphisme de Frobenius $\varphi$ de $\co_{\oK^\flat}$, et qui est un $\alpha$-isomorphisme.
\end{prop}

En effet, considérons le système projectif de $\co_{\oK^\flat}$-modules $M=(M_n)_{n\geq 1}$ défini pour tout entier $n\geq 1$, par 
\begin{equation}\label{TPCF9b}
M_n=\rH^q(\tE_s^{\mN^\circ},\uplambda^*(\cL)\otimes_{\mF_p}(\ocB^\wp/\varpi^n\ocB^\wp)),
\end{equation} 
les morphismes de transition étant induits par les projections canoniques 
\begin{equation}
\ocB^\wp/\varpi^{n+1}\ocB^\wp\rightarrow \ocB^\wp/\varpi^n\ocB^\wp.
\end{equation}
D'après \ref{TPCF8}(iii), pour tout entier $i\geq 0$, on a un isomorphisme canonique 
\begin{equation}
\ocB^\wp/\varpi^{p^i}\ocB^\wp\stackrel{\sim}{\rightarrow}  \varepsilon_{\leq i}(\ocB^\wp).
\end{equation}
On en déduit, compte tenu de (\cite{agt} III.7.11) et (\cite{jannsen} 1.15), un isomorphisme $\co_{\oK^\flat}$-linéaire
\begin{equation}\label{TPCF9c}
M_{p^i}\stackrel{\sim}{\rightarrow} \rH^q(\tE_s,\cL\otimes_{\mF_p} \ocB_1),
\end{equation}
où le but est considéré comme un $\co_{\oK^\flat}$-module via $\pi_i$ \eqref{TFA102a}.
De plus, le diagramme 
\begin{equation}\label{TPCF9d}
\xymatrix{
{M_{p^{i+1}}}\ar[r]^-(0.5){\sim}\ar[d]&{\rH^q(\tE_s,\cL\otimes_{\mF_p} \ocB_1)}\ar[d]\\
{M_{p^i}}\ar[r]^-(0.5){\sim}&{\rH^q(\tE_s,\cL\otimes_{\mF_p} \ocB_1)}}
\end{equation}
où la flèche verticale de gauche est le morphisme canonique et la flèche verticale de droite est induite par l'endomorphisme de Frobenius de $\ocB_1$, 
est commutatif. Les $\co_{\oK^\flat}$-modules $M_{p^i}$ sont donc de type $\alpha$-fini en vertu de \ref{finitude3}.

Reprenons les notations de \ref{mptf10} pour $R=\co_{\oK^\flat}$. 
La multiplication par $\varpi$ dans $\ocB^\wp$ induit un morphisme $\co_{\oK^\flat}$-linéaire 
\begin{equation}
u\colon M[-1]\rightarrow M.
\end{equation}
Il est clair que le composé $M[-1]\stackrel{u}{\rightarrow}M\rightarrow M[-1]$, où la seconde flèche est le morphisme canonique, est induit par 
la multiplication par $\varpi$. Par ailleurs, il résulte de \ref{TPCF8}(iii) (pour $i=0$) que la suite
\begin{equation}
\xymatrix{
0\ar[r]&{\ocB^\wp/\varpi \ocB^\wp}\ar[r]^-(0.5){\cdot\varpi^n}&{\ocB^\wp/\varpi^{n+1}\ocB^\wp}\ar[r]&{\ocB^\wp/\varpi^n\ocB^\wp}\ar[r]&0}
\end{equation}
est exacte dans $\cP_\AR$. On en déduit que la suite 
\begin{equation}
\xymatrix{M_1\ar[r]^-(0.5){u^{\otimes n}}&M_{n+1}\ar[r]&M_n}
\end{equation}
où la seconde flèche est le morphisme canonique, est exacte au centre. 

Il résulte de \ref{TPCF8}(ii) que l'endomorphisme de Frobenius de $\ocB^\wp$ induit pour tout $n\geq 1$, un isomorphisme de $\cP_\AR$ 
\begin{equation}\label{TPCF9e}
\phi\colon \ocB^\wp/\varpi^n\ocB^\wp\rightarrow\ocB^\wp/\varpi^{pn}\ocB^\wp,
\end{equation}
et par suite un isomorphisme $\co_{\oK^\flat}$-semi-linéaire
\begin{equation}\label{TPCF9f}
M\stackrel{\sim}{\rightarrow}M^{(p)}.
\end{equation}
On notera que pour tout $i\geq 0$, le diagramme 
\begin{equation}
\xymatrix{
{M_{p^i}}\ar[r]\ar[d]&{\rH^q(\tE_s,\cL\otimes_{\mF_p} \ocB_1)}\ar@{=}[d]\\
{M_{p^{i+1}}}\ar[r]&{\rH^q(\tE_s,\cL\otimes_{\mF_p} \ocB_1)}}
\end{equation}
où la flèche verticale de gauche est la composante de degré $p^i$ de 
l'isomorphisme \eqref{TPCF9f} et les flèches horizontales sont les isomorphismes \eqref{TPCF9c}, est commutatif. 

En vertu de \ref{mptf11}(i), les morphismes de transition du système projectif $M$ sont $\alpha$-surjectifs.  
Par ailleurs, d'après (\cite{agt} III.7.11), on a une suite exacte
\begin{eqnarray}
\lefteqn{0\rightarrow \rR^1\underset{\underset{\mN}{\longleftarrow}}{\lim}\ \rH^{q-1}(\tE_s,\cL\otimes_{\mF_p}\ocB_1)}\\
&&\rightarrow \rH^q(\tE_s^{\mN^\circ},\uplambda^*(\cL)\otimes_{\mF_p}\ocB^\wp)\rightarrow \underset{\underset{\mN}{\longleftarrow}}{\lim}\ 
\rH^q(\tE_s,\cL\otimes_{\mF_p}\ocB_1)\rightarrow 0,\nonumber
\end{eqnarray}
où les morphismes de transition des systèmes projectifs sont induits par l'endomorphisme de Frobenius $\phi$ de $\ocB_1$.
Compte tenu de \ref{alpha10}, \eqref{TPCF9d}, (\cite{jannsen} 1.15) et (\cite{roos} théo.~1), le terme de gauche de cette suite est donc $\alpha$-nul. 
On en déduit que le morphisme canonique
\begin{equation}
\rH^q(\tE_s^{\mN^\circ},\uplambda^*(\cL)\otimes_{\mF_p}\ocB^\wp)\rightarrow \underset{\longleftarrow}{\lim}\ M_n
\end{equation}
est un $\alpha$-isomorphisme. La proposition s'ensuit en vertu de \ref{mptf11}(iii).

\begin{cor}\label{TPCF90}
Sous les hypothèses de \eqref{TPCF9}, pour tout entier $i\geq 0$, il existe un morphisme canonique $\co_{\oK^\flat}$-linéaire
\begin{equation}
\rH^q(\tE_s^{\mN^\circ},\uplambda^*(\cL)\otimes_{\mF_p}\ocB^\wp)\otimes_{\co_{\oK^\flat}}(\co_{\oK^\flat}/\varpi^{p^i}\co_{\oK^\flat})\rightarrow
\rH^q(\tE_s,\cL\otimes_{\mF_p}\ocB_1),
\end{equation}
où le but est considéré comme un $\co_{\oK^\flat}$-module via $\pi_i$ \eqref{TFA102a}, qui est un $\alpha$-isomorphisme.
\end{cor}

\begin{prop}\label{TPCF10}
Pour tout élément $b\in \co_\oK$ tel que $b^p\not\in p\co_\oK$, la suite 
\begin{equation}\label{TPCF10a}
\xymatrix{
0\ar[r]&{\mF_p \cdot b\oplus (\frac{p}{b^{p-1}})\cdot \ocB_1}\ar[r]&{\ocB_1}\ar[rr]^{\phi-b^{p-1}\id}&&{\ocB_1}\ar[r]&0}
\end{equation}
est exacte.
\end{prop} 

En effet, la question étant locale, on peut se borner au cas où $X$ est affine. 
Soit $(\oy\rightsquigarrow \ox)$ un point de $X_\et\gtimes_{X_\et}\oX^\circ_\et$ \eqref{topfl17} tel que $\ox$ soit au-dessus de $s$.
On désigne par $\fV_\ox$ la catégorie des $X$-schémas étales $\ox$-pointés, 
et par $\fW_\ox$ la sous-catégorie pleine de $\fV_\ox$ formée des objets $(U,\fp\colon \ox\rightarrow U)$
tels que le schéma $U$ soit affine. 
Ce sont des catégories cofiltrantes, et le foncteur d'injection canonique
$\fW^\circ_\ox\rightarrow \fV^\circ_\ox$ est cofinal (\cite{sga4} I 8.1.3(c)). 
Avec les notations de \ref{TFA11}, on a un isomorphisme canonique \eqref{TFA11c} 
\begin{equation}\label{TPCF10b}
\ocB_{\rho(\oy\rightsquigarrow \ox)} \stackrel{\sim}{\rightarrow} 
\underset{\underset{(U,\fp)\in \fW^\circ_\ox}{\longrightarrow}}{\lim}\ \oR^{\oy}_U.
\end{equation} 
Montrons que pour tout $z\in \ocB_{\rho(\oy\rightsquigarrow \ox)}$, il existe $t\in \ocB_{\rho(\oy\rightsquigarrow \ox)}$ tel que 
$t^p-b^{p-1}t=z$. On peut supposer que $z\in \oR^\oy_X$. 
Comme l'anneau $\ocB_{\rho(\oy\rightsquigarrow \ox)}$ est local (\cite{agt} III.10.10(i)), d'idéal maximal contenant $p$, 
\begin{equation}\label{TPCF10c}
(\frac{p}{b^{p-1}})^pz^{p-1}-(1-p)^{p-1}
\end{equation} 
est une unité de $\ocB_{\rho(\oy\rightsquigarrow \ox)}$. 
Il existe alors un objet $(U,\fp\colon \ox\rightarrow U)$ de $\fW_\ox$ tel que l'élément \eqref{TPCF10c} soit une unité de $\oR^\oy_U$. 
Il suffit de montrer que l'extension 
\begin{equation}\label{TPCF10d}
\oR^\oy_U[T]/(T^p-b^{p-1}T-z)
\end{equation}
est étale au-dessus de $\oR^\oy_U[\frac 1 p]$, ou encore que pour tout corps $L$ et tout homomorphisme $u\colon \oR^\oy_U[\frac 1 p]\rightarrow L$, 
les équations $T^p-b^{p-1} T -u(z)$ et $p T^{p-1}-b^{p-1}$ n'ont aucune solution commune dans $L$. 
Si une telle solution $y\in L$ existait, on aurait $y^{p-1}=\frac{b^{p-1}}{p}$ et $u(z)=y^p-b^{p-1} y\in L$, et par suite,
\begin{equation}\label{TPCF10e}
(\frac{p}{b^{p-1}})^pu(z)^{p-1}=(1-p)^{p-1},
\end{equation}
ce qui est absurde puisque $(\frac{p}{b^{p-1}})^pz^{p-1}-(1-p)^{p-1}$ est une unité de $\oR^\oy_U$.

Considérons ensuite, pour tout $t\in \ocB_{\rho(\oy\rightsquigarrow \ox)}$, l'équation
\begin{equation}\label{TPCF10f}
t^p-b^{p-1}t\in p\ocB_{\rho(\oy\rightsquigarrow \ox)},
\end{equation}
qui est équivalente, puisque $\ocB_{\rho(\oy\rightsquigarrow \ox)}$ est $\co_\oK$-plat, à l'équation
\begin{equation}\label{TPCF10g}
(\frac{t}{b})^p-\frac{t}{b}\in \frac{p}{b^p}\ocB_{\rho(\oy\rightsquigarrow \ox)}.
\end{equation}
Comme l'anneau $\ocB_{\rho(\oy\rightsquigarrow \ox)}$ est normal et strictement local (\cite{agt} III.10.10(i)), d'idéal maximal contenant $p$, 
on en déduit d'abord que $\frac t b\in \ocB_{\rho(\oy\rightsquigarrow \ox)}$, puis que 
\begin{equation}
t\in b\cdot (\mZ_p+\frac{p}{b^{p}}\ocB_{\rho(\oy\rightsquigarrow \ox)}).
\end{equation}

Il résulte de ce qui précède et de \ref{tf21}(ii) que la suite \eqref{TPCF10a} est exacte.

\begin{cor}\label{TPCF11}
Pour tout élément non-nul $b$ de $\co_{\oK^\flat}$, la suite  
\begin{equation}\label{TPCF11a}
\xymatrix{
0\ar[r]&{\mF_p}\ar[r]^-(0.4){\cdot b}&{\ocB^\wp}\ar[rr]^{\phi-b^{p-1}\id}&&{\ocB^\wp}\ar[r]&0}
\end{equation}
est exacte dans $\cP_\AR$.
\end{cor}

Soit $(b_n)_{n\geq 0}\in (\co_C)^{\mN}$ l'image inverse de $b$ par l'isomorphisme \eqref{TFA101b}. 
Il existe $n_0\geq 0$ tel que pour tout $n\geq n_0$, $b_n^p\not\in p\co_C$. En vertu de \ref{TPCF10}, la suite 
\begin{equation}\label{TPCF11b}
\xymatrix{
0\ar[r]&{\mF_p \cdot b_n\oplus (\frac{p}{b_n^{p-1}})\cdot \ocB_1}\ar[r]&{\ocB_1}\ar[rr]^{\phi-b_n^{p-1}\id}&&{\ocB_1}\ar[r]&0}
\end{equation}
est exacte. Comme $v(b_n)<\frac 1 p$, on a $pv(\frac{p}{b_n^{p-1}})>1$, d'où la proposition.

\begin{prop}\label{TPCF12}
Supposons $X$ propre sur $S$, et soient  $\mL$ un $\mF_p$-module de type fini de $\oX^\circ_\fet$, $q$ un entier $\geq 0$, 
$b$ un élément non-nul de $\co_{\oK^\flat}$. 
Posons $\cL=\delta^*(\beta^*(\mL))$, où $\beta$ est le morphisme \eqref{TFA6b} et $\delta$ est le plongement \eqref{TFA66a},
et notons $\uphi$  l'endomorphisme $\co_{\oK^\flat}$-semi-linéaire de 
\[
\rH^q(\tE_s^{\mN^\circ},\uplambda^*(\cL)\otimes_{\mF_p}\ocB^\wp)
\]
induit par l'endomorphisme de Frobenius $\phi$ de $\ocB^\wp$. Alors,
\begin{itemize}
\item[{\rm (i)}] La suite  
\begin{equation}\label{TPCF12a}
\xymatrix{ 
0\ar[r]&{\rH^q(\tE_s,\cL)}\ar[r]^-(0.5){\cdot b}&{\rH^q(\tE_s^{\mN^\circ},\uplambda^*(\cL)\otimes_{\mF_p}\ocB^\wp)}\ar[rr]^-(0.5){\uphi-b^{p-1}\id}&&
{\rH^q(\tE_s^{\mN^\circ},\uplambda^*(\cL)\otimes_{\mF_p}\ocB^\wp)}\ar[r]& 0}
\end{equation}
est exacte.
\item[{\rm (ii)}] Le morphisme canonique
\begin{equation}\label{TPCF12b}
\ker(\uphi-b^{p-1}\id|\rH^q(\tE_s^{\mN^\circ},\uplambda^*(\cL)\otimes_{\mF_p}\ocB^\wp))\rightarrow 
\ker(\uphi-b^{p-1}\id|\rH^q(\tE_s^{\mN^\circ},\uplambda^*(\cL)\otimes_{\mF_p}\ocB^\wp)\otimes_{\co_{\oK^\flat}}\oK^\flat)
\end{equation}
est un isomorphisme.
\end{itemize}
\end{prop}
En effet, d'après (\cite{agt} III.7.11) et (\cite{jannsen} 1.15), la suite \eqref{TPCF11a} induit une suite exacte longue de cohomologie
\begin{equation}
\xymatrix{
\rM^{q-1}\ar[r]&{\rH^q(\tE_s,\cL)}\ar[r]&{\rM^q}\ar[rr]^{\uphi-b^{p-1}\id}&&{\rM^q}\ar[r]&{\rH^{q+1}(\tE_s,\cL)}},
\end{equation}
où on a posé $\rM^j=\rH^j(\tE_s^{\mN^\circ},\uplambda^*(\cL)\otimes_{\mF_p}\ocB^\wp)$ pour $j\geq 0$ et $\rM^{-1}=0$. 
En vertu de \ref{TPCF9}, il existe un entier $r\geq 0$ et un morphisme $\co_{\oK^\flat}$-linéaire
\begin{equation}\label{TPCF12c}
u\colon \fm_{\oK^\flat}\otimes_{\co_{\oK^\flat}}\rM^q\rightarrow \co_{\oK^\flat}^r 
\end{equation}
compatible à $\uphi$ et à l'endomorphisme de Frobenius $\varphi$ de $\co_{\oK^\flat}$, et qui est un $\alpha$-isomorphisme. 
Comme $\oK^\flat$ est algébriquement clos, pour tout $c\in \co_{\oK^\flat}$, l'endomorphisme $\varphi-c^{p-1}\id$ de $\co_{\oK^\flat}$ est surjectif.
Pour tous $y\in \rM^q$ et $t_1\in \fm_{\oK^\flat}$, il existe $z\in \co_{\oK^\flat}^r$ tel que 
\begin{equation}\label{TPCF12d}
u(t_1^p\otimes y)=\varphi(z)-c^{p-1}z,
\end{equation}
où on a encore noté $\varphi$ l'endomorphisme de $\co_{\oK^\flat}^r$ induit par l'endomorphisme $\varphi$ de $\co_{\oK^\flat}$. 
Pour tout $t_2\in \fm_{\oK^\flat}$, il existe $x\in \rM^q$ et $t\in \fm_{\oK^\flat}$ tels que $u(t\otimes x)=t_2z$. On en déduit que 
\begin{equation}\label{TPCF12e}
u((t_1t_2)^p\otimes y)=u(t^p\otimes \uphi(x)-(ct_2)^{p-1} t\otimes x)
\end{equation}
Pour tout $t_3\in \fm_{\oK^\flat}$, on a donc 
\begin{equation}\label{TPCF12f}
(t_1t_2t_3)^p y=\uphi(tt_3x)-(ct_2t_3)^{p-1} (tt_3 x).
\end{equation}
Par ailleurs, pour tout élément non-nul $b'$ de $\co_{\oK^\flat}$ on a un diagramme commutatif
\begin{equation}\label{TPCF12g}
\xymatrix{
{\ker(\uphi-b^{p-1}\id|\rM^q)}\ar[r]\ar[d]_{\cdot b'}&{\rM^q}\ar[rr]^{\uphi-b^{p-1}\id}\ar[d]_{\cdot b'}&&{\rM^q}\ar[r]\ar[d]^{\cdot b'^p}&
{\coker(\uphi-b^{p-1}\id|\rM^q)}\ar[d]^{\cdot b'^p}\\
{\ker(\uphi-(bb')^{p-1}\id|\rM^q)}\ar[r]&{\rM^q}\ar[rr]^{\uphi-(bb')^{p-1}\id}&&{\rM^q}\ar[r]&{\coker(\uphi-(bb')^{p-1}\id|\rM^q)}}
\end{equation}
On déduit de ce qui précède (en prenant $c=bt_1 $ et $b'=t_1t_2t_3$) 
que pour tout $b'\in \fm_{\oK^\flat}$, le dernier morphisme vertical du diagramme \eqref{TPCF12g} est nul.

D'après \ref{TPCF11}, pour tout élément non-nul $b'$ de $\co_{\oK^\flat}$, on a un diagramme commutatif à lignes exactes
\begin{equation}\label{TPCF12h}
\xymatrix{
0\ar[r]&{\mF_p}\ar@{=}[d]\ar[r]^-(0.5){\cdot b}&{\ocB^\wp}\ar[d]^{\cdot b'}\ar[rr]^-(0.5){\phi-b^{p-1}\id}&&{\ocB^\wp}\ar[d]^{\cdot b'^p}\ar[r]&0\\
0\ar[r]&{\mF_p}\ar[r]^-(0.5){\cdot bb'}&{\ocB^\wp}\ar[rr]^-(0.5){\phi-(bb')^{p-1}\id}&&{\ocB^\wp}\ar[r]&0}
\end{equation}
Il induit un diagramme commutatif à lignes exactes
\begin{equation}\label{TPCF12i}
\xymatrix{
0\ar[r]&{\coker(\uphi'-b^{p-1}\id|\rM^{q-1})}\ar[r]\ar[d]_{\cdot b'^p}&{\rH^q(\tE_s,\cL)}\ar@{=}[d]\ar[r]&{\ker(\uphi-b^{p-1}\id|\rM^q)}\ar[d]^{\cdot b'}\ar[r]&0\\
0\ar[r]&{\coker(\uphi'-(bb')^{p-1}\id|\rM^{q-1})}\ar[r]&{\rH^q(\tE_s,\cL)}\ar[r]&{\ker(\uphi-(bb')^{p-1}\id|\rM^q)}\ar[r]&0}
\end{equation}
où on a noté $\uphi'$ l'endomorphisme de $\rM^{q-1}$ induit par $\phi$. 
Comme le morphisme vertical de gauche est nul, on en déduit que l'endomorphisme $\uphi'-b^{p-1}\id$ de $\rM^{q-1}$ est surjectif, 
et il en est donc de même de l'endomorphisme $\uphi-b^{p-1}\id$ de $\rM^q$; 
d'où l'exactitude de la suite \eqref{TPCF12a}.
On en déduit  aussi que la flèche verticale de droite est un isomorphisme~:
\begin{equation}\label{TPCF12j}
\cdot b'\colon \ker(\uphi-b^{p-1}\id|\rM^q)\stackrel{\sim}{\rightarrow}\ker(\uphi-(bb')^{p-1}\id|\rM^q).
\end{equation}
Il s'ensuit que l'application \eqref{TPCF12b} est injective. 

Soit $x\in \rM^q\otimes_{\co_{\oK^\flat}}\oK^\flat$ tel que 
$\uphi(x)=b^{p-1}x$. Il existe $z\in \rM^q$ et $b'\in \fm_{\oK^\flat}$ tels que $b'\not=0$ et $z=b'x\in \rM^q\otimes_{\co_{\oK^\flat}}\oK^\flat$.
Par suite, $\uphi(z)-(bb')^{p-1}z$ est annulé par une puissance de $\varpi$. 
Quitte à remplacer $b'$ par un multiple, on peut supposer que $\uphi(z)=(bb')^{p-1}z$.
D'après l'isomorphisme \eqref{TPCF12j}, il existe $x'\in \ker(\uphi-b^{p-1}\id|\rM^q)$ tel que $z=b'x'$. On a alors $x=x'$ dans 
$\rM^q\otimes_{\co_{\oK^\flat}}\oK^\flat$; d'où la surjectivité de l'application \eqref{TPCF12b}. 

\begin{cor}\label{TPCF13}
Supposons $X$ propre sur $S$, et soient  $\mL$ un $\mF_p$-module de type fini de $\oX^\circ_\fet$, $\cL=\delta^*(\beta^*(\mL))$.
Alors, pour tout entier $q\geq 0$, $\rH^q(\tE_s,\cL)$ est un $\mF_p$-espace vectoriel de dimension finie,
et le morphisme canonique
\begin{equation}\label{TPCF13a}
\rH^q(\tE_s,\cL)\otimes_{\mF_p}\co_{\oK^\flat}\rightarrow \rH^q(\tE_s^{\mN^\circ},\uplambda^*(\cL)\otimes_{\mF_p}\ocB^\wp)
\end{equation}
est un $\alpha$-isomorphisme.
\end{cor}

En effet, en vertu de \ref{TPCF12} et avec les notations de sa preuve, on a des isomorphismes canoniques 
\begin{equation}
\rH^q(\tE_s,\cL)\stackrel{\sim}{\rightarrow}\ker(\uphi-\id |\rM^q)\stackrel{\sim}{\rightarrow}
\ker(\uphi-\id |\rM^q\otimes_{\co_{\oK^\flat}}\oK^\flat).
\end{equation}
D'après \ref{TPCF9}, $\rM^q\otimes_{\co_{\oK^\flat}}\oK^\flat$ est un $\oK^\flat$-espace vectoriel de dimension finie,
et $\uphi$ est un isomorphisme $\co_{\oK^\flat}$-semi-linéaire de $\rM^q\otimes_{\co_{\oK^\flat}}\oK^\flat$. 
Il résulte alors de (\cite{katz} 4.1.1) que $\rH^q(\tE_s,\cL)$ est un $\mF_p$-espace vectoriel de dimension finie, et que le morphisme $\oK^\flat$-linéaire
\begin{equation}
\rH^q(\tE_s,\cL)\otimes_{\mF_p}\oK^\flat\rightarrow \rM^q\otimes_{\co_{\oK^\flat}}\oK^\flat
\end{equation}
est bijectif. Par suite, le morphisme \eqref{TPCF13a} est injectif. Montrons qu'il est $\alpha$-surjectif.
Considérons le morphisme $u\colon \fm_{\oK^\flat}\otimes_{\co_{\oK^\flat}}\rM^q\rightarrow \co_{\oK^\flat}^r$ \eqref{TPCF12c} 
et notons $e_1,\dots,e_r$ la base canonique de $\co_{\oK^\flat}^r$. Comme $u$ est $\alpha$-injectif, 
il suffit de montrer que pour tout $1\leq i\leq r$ et tout $\gamma\in \fm_{\oK^\flat}$, il existe $x\in \rM^q$ tel que 
$\uphi(x)=x$ et $u(\gamma\otimes x)=\gamma e_i$. Pour tout $t\in \fm_{\oK^\flat}$,  
il existe $y\in \rM^q$ et $\beta\in \fm_{\oK^\flat}$ tels que $u(\beta\otimes y)=te_i$. 
On a $u(\beta^p\otimes \uphi(y))=t^pe_i$ et donc $u(\beta^p\otimes \uphi(y)-t^{p-1}\beta\otimes y)=0$. Par suite, pour 
tout $t'\in \fm_{\oK^\flat}$, on a 
\begin{equation}
\uphi(t'\beta y)=(tt')^{p-1}t'\beta y.
\end{equation}
Compte tenu de l'isomorphisme \eqref{TPCF12j} (pour $b=1$ et $b'=tt'$), il existe $x\in \rM^q$ tel que $\uphi(x)=x$ et $t'\beta y=tt'x$. 
Par suite, pour tout $t''\in \fm_{\oK^\flat}$, on a $t't''\beta\otimes y=tt't''\otimes x$ dans $\fm_{\oK^\flat}\otimes_{\co_{\oK^\flat}}\rM^q$
et donc $u(tt't''\otimes x)=tt't''e_i$, d'où l'assertion recherchée. 

\begin{cor}\label{TPCF14}
Supposons $X$ propre sur $S$, et soient  $n$ un entier $\geq 1$, $\mL$ un $(\mZ/p^n\mZ)$-module de type fini de 
$\oX^\circ_\fet$, $\cL=\delta^*(\beta^*(\mL))$. Alors, pour tout entier $q\geq 0$, le morphisme canonique
\begin{equation}\label{TPCF14a}
\rH^q(\tE_s,\cL)\otimes_{\mZ_p}\co_\oK\rightarrow \rH^q(\tE_s,\cL\otimes_{\mZ_p}\ocB_n)
\end{equation}
est un $\alpha$-isomorphisme.
\end{cor}

En effet, par dévissage, on peut se réduire au cas où $n=1$. On notera que $\ocB_n$ est $(\mZ/p^n\mZ)$-plat (\cite{agt} III.9.2). 
La proposition résulte alors de \ref{TPCF13}, \ref{TPCF90} et \eqref{TFA102b}.

\begin{teo}[\cite{faltings2} Theorem 8 page 223]\label{TPCF15}
Supposons $X$ propre sur $S$, et soit $F$ un faisceau abélien de torsion, localement constant constructible de $\oX^\circ_{\et}$.
Alors, pour tout entier $q\geq 0$, on a un morphisme canonique 
\begin{equation}
\rH^q(\oX^\circ_{\et},F)\otimes_{\mZ_p}\co_\oK\rightarrow \rH^q(\tE,\psi_*(F)\otimes_{\mZ_p}\ocB),
\end{equation}
où $\psi$ est le morphisme \eqref{acycloc1a}, qui est un $\alpha$-isomorphisme. 
\end{teo}

On peut évidemment se borner au cas où $F$ est annulé par $p^n$, pour un entier $n\geq 1$. 
En vertu de \ref{acycloc2} et de la suite spectrale de Cartan-Leray, le morphisme canonique 
\begin{equation}
\rH^q(\tE,\psi_*(F)) \rightarrow  \rH^q(\oX^\circ_{\et},F)
\end{equation}
est un isomorphisme. 
Notant $\rho_{\oX^\circ}\colon \oX^\circ_\et\rightarrow \oX^\circ_\fet$ le morphisme canonique \eqref{notconv10a}, 
il existe un $(\mZ/p^n\mZ)$-module de type fini canonique $\mL$ de $\oX^\circ_\fet$ tel que $F=\rho^*_{\oX^\circ}(\mL)$ d'après 
(\cite{agt} VI.9.18, VI.9.20 et III.2.11). En vertu de (\cite{agt} VI.10.9), on a un isomorphisme canonique 
\begin{equation}
\beta^*(\mL)\stackrel{\sim}{\rightarrow}\psi_*(\rho_{\oX^\circ}^*(\mL)).
\end{equation}
Considérons le diagramme commutatif de morphismes canoniques
\begin{equation}
\xymatrix{
{\rH^q(\tE,\beta^*(\mL))\otimes_{\mZ_p}\co_\oK} \ar[r]\ar[d]&{\rH^q(\tE,\beta^*(\mL)\otimes_{\mZ_p}\ocB_n)}\ar[d]\\
{\rH^q(\tE_s,\delta^*(\beta^*(\mL)))\otimes_{\mZ_p}\co_\oK} \ar[r]&{\rH^q(\tE_s,\delta^*(\beta^*(\mL))\otimes_{\mZ_p}\ocB_n)}}
\end{equation}
Les flèches verticales sont des isomorphismes d'après \ref{tf27}. La flèche horizontale inférieure est
un isomorphisme en vertu de \ref{TPCF14}. Il en est donc de même de la flèche horizontale supérieure; d'où la proposition.

\begin{cor}\label{TPCF16}
Si $X$ est propre sur $S$, pour tous entiers $n\geq 1$ et $q\geq 0$, on a un morphisme canonique 
\begin{equation}\label{TPCF16a}
\rH^q(\oX^\circ_{\et},\mZ/p^n\mZ)\otimes_{\mZ_p}\co_\oK\rightarrow \rH^q(\tE_s,\ocB_n),
\end{equation}
qui est $\alpha$-isomorphisme. 
\end{cor}

En effet, d'après (\cite{agt} VI.10.9(iii)), on a un isomorphisme canonique 
\begin{equation}
\beta^*(\mZ/p^n\mZ)\stackrel{\sim}{\rightarrow} \psi_*(\rho_{\oX^\circ}^*(\mZ/p^n\mZ)),
\end{equation}
où $\rho_{\oX^\circ}\colon \oX^\circ_\et\rightarrow \oX^\circ_\fet$ est le morphisme canonique \eqref{notconv10a}. 
On en déduit un isomorphisme 
\begin{equation}
\mZ/p^n\mZ\stackrel{\sim}{\rightarrow} \psi_*(\mZ/p^n\mZ),
\end{equation}
qui n'est autre que le morphisme canonique par adjonction. Par ailleurs, comme le foncteur $\delta_*$ \eqref{TFA66a} est exact, 
on a un isomorphisme canonique 
\begin{equation}
\rH^q(\tE_s,\ocB_n) \stackrel{\sim}{\rightarrow} \rH^q(\tE,\ocB_n).
\end{equation}
La proposition résulte alors de \ref{TPCF15}.

\chapter{Cohomologie relative des topos de Faltings}\label{cohtopfalrel}

\section{Hypothèses et notations; morphismes de schémas logarithmiques adéquats}\label{mtfla}

\subsection{}\label{mtfla1}
Dans ce chapitre, $K$ désigne un corps de valuation discrète complet de 
caractéristique $0$, à corps résiduel {\em algébriquement clos} $k$ de caractéristique $p>0$,  
$\co_K$ l'anneau de valuation de $K$, $\oK$ une clôture algébrique de $K$, $\co_\oK$ la clôture intégrale de $\co_K$ dans $\oK$,
$\fm_\oK$ l'idéal maximal de $\co_\oK$ et $G_K$ le groupe de Galois de $\oK$ sur $K$.
On note $\co_C$ le séparé complété $p$-adique de $\co_\oK$, $\fm_C$ son idéal maximal,
$C$ son corps des fractions et $v$ sa valuation, normalisée par $v(p)=1$. 
On désigne par $\hmZ(1)$ et $\mZ_p(1)$ les $\mZ[G_K]$-modules 
\begin{eqnarray}
\hmZ(1)&=&\underset{\underset{n\geq 1}{\longleftarrow}}{\lim}\ \mu_{n}(\co_{\oK}),\label{mtfla1a}\\
\mZ_p(1)&=&\underset{\underset{n\geq 0}{\longleftarrow}}{\lim}\ \mu_{p^n}(\co_{\oK}),\label{mtfla1c}
\end{eqnarray}  
où $\mu_n(\co_{\oK})$ désigne le sous-groupe des racines $n$-ièmes de l'unité dans $\co_\oK$. 
Pour tout $\mZ_p[G_K]$-module $M$ et tout entier $n$, on pose $M(n)=M\otimes_{\mZ_p}\mZ_p(1)^{\otimes n}$.

Comme $\co_\oK$ satisfait les conditions requises dans \ref{mptf1}, 
il est loisible de considérer les notions de $\alpha$-algèbre introduites dans les sections \ref{alpha}--\ref{mptf} relativement à $\co_\oK$.   

On pose $S=\Spec(\co_K)$ et $\oS=\Spec(\co_\oK)$. 
On note $s$ (resp.  $\eta$, resp. $\oeta$) le point fermé de $S$ (resp.  générique de $S$, resp. générique de $\oS$).
Pour tout entier $n\geq 1$, on pose $S_n=\Spec(\co_K/p^n\co_K)$. 
On munit $S$ de la structure logarithmique $\cM_S$ définie par son point fermé. 

Pour tout $S$-schéma $X$, on note $X_s$ (resp. $X_\eta$, resp. $X_\oeta$) 
la fibre fermée (resp. générique, resp. géométrique générique) de $X$ au dessus de $S$, et on pose 
\begin{equation}\label{mtfla1b}
\oX=X\times_S\oS \ \ \ {\rm et}\ \ \  X_n=X\times_SS_n.
\end{equation}

\subsection{}\label{mtfla2}
Dans ce chapitre, $f\colon (X,\cM_X)\rightarrow (S,\cM_S)$ et $f'\colon (X',\cM_{X'})\rightarrow (S,\cM_S)$
désignent des morphismes {\em adéquats} de schémas logarithmiques (\cite{agt} III.4.7) et 
\begin{equation}\label{mtfla2a}
g\colon (X',\cM_{X'})\rightarrow (X,\cM_X)
\end{equation} 
un $(S,\cM_S)$-morphisme lisse et saturé. 
On désigne par $X^\circ$ le sous-schéma ouvert maximal de $X$
où la structure logarithmique $\cM_X$ est triviale; c'est un sous-schéma ouvert de $X_\eta$.
On note $j\colon X^\circ\rightarrow X$ l'injection canonique et $j_\oX\colon \oX^\circ\rightarrow \oX$ le changement de base de $j$. 
Pour tout $X$-schéma $U$, on pose  
\begin{equation}\label{mtfla2b}
U^\circ=U\times_XX^\circ.
\end{equation} 
On note $\hbar\colon \oX\rightarrow X$ et $h\colon \oX^\circ\rightarrow X$ les morphismes canoniques \eqref{TFA1b}, de sorte que 
l'on a $h=\hbar\circ j_\oX$. 

On désigne par $X'^\rhd$ le sous-schéma ouvert maximal de $X'$
où la structure logarithmique $\cM_{X'}$ est triviale; c'est un sous-schéma ouvert de $X'^\circ=X'\times_XX^\circ$.
On note $j'\colon X'^\rhd\rightarrow X'$ l'injection canonique et $j'_{\oX'}\colon \oX'^\rhd\rightarrow \oX'$ le changement de base de $j'$. 
Pour tout $X'$-schéma $U'$, on pose  
\begin{equation}\label{mtfla2c}
U'^\rhd=U'\times_{X'}X'^\rhd.
\end{equation} 
On note $\hbar'\colon \oX'\rightarrow X'$ et $h'\colon \oX'^\rhd\rightarrow X'$ les morphismes canoniques \eqref{TFA1b}, de sorte que 
l'on a $h'=\hbar'\circ j'_{\oX'}$. 

Pour alléger les notations, on pose
\begin{equation}\label{mtfla2d}
\tOmega^1_{X/S}=\Omega^1_{(X,\cM_X)/(S,\cM_S)},
\end{equation}
que l'on considère comme un faisceau de $X_\zar$ ou $X_\et$, selon le contexte \eqref{notconv12}. 
C'est un $\co_X$-module localement libre de type fini.
De même, on pose
\begin{equation}\label{mtfla2e}
\tOmega^1_{X'/S}=\Omega^1_{(X',\cM_{X'})/(S,\cM_S)} \ \ \ {\rm et} \ \ \
\tOmega^1_{X'/X}=\Omega^1_{(X',\cM_{X'})/(X,\cM_X)},
\end{equation}
que l'on considère comme un faisceau de $X'_\zar$ ou $X'_\et$, selon le contexte. 
Ce sont des $\co_{X'}$-modules  localement libres de type fini. 
De plus, on a alors une suite exacte canonique de $\co_{X'}$-modules 
\begin{equation}\label{mtfla2f}
0\rightarrow g^*(\tOmega^1_{X/S})\rightarrow \tOmega^1_{X'/S}\rightarrow \tOmega^1_{X'/X}\rightarrow 0.
\end{equation}

\subsection{}\label{mtfla3}
Pour tout entier $n\geq 1$, on note $a\colon X_s\rightarrow X$, $a_n\colon X_s\rightarrow X_n$, $a'\colon X'_s\rightarrow X'$ et $a'_n\colon X'_s\rightarrow X'_n$ 
les injections canoniques \eqref{mtfla1b}. 
Le corps résiduel de $\co_K$ étant algébriquement clos, il existe un unique $S$-morphisme $s\rightarrow \oS$. 
Celui-ci induit des immersions fermées $\oa\colon X_s\rightarrow \oX$, $\oa_n\colon X_s\rightarrow \oX_n$, 
$\oa'\colon X'_s\rightarrow \oX'$ et $\oa'_n\colon X'_s\rightarrow \oX'_n$ qui relèvent $a$, $a_n$, $a'$ et $a'_n$, respectivement. 
Comme $\oa_n$ (resp. $\oa'_n$) est un homéomorphisme universel, on peut considérer $\co_{\oX_n}$ (resp. $\co_{\oX'_n}$)
comme un faisceau de $X_{s,\zar}$ ou $X_{s,\et}$ (resp. $X'_{s,\zar}$ ou $X'_{s,\et}$), selon le contexte (cf. \ref{notconv12}).

\subsection{}\label{mtfla7}
On désigne par
\begin{equation}\label{mtfla7a}
\pi\colon E\rightarrow \Et_{/X}
\end{equation}
le $\mU$-site fibré de Faltings associé au morphisme $h\colon \oX^\circ\rightarrow X$ (cf. \ref{tf1}).
On munit $E$ de la topologie co-évanescente définie par $\pi$ (\cite{agt} VI.5.3)  
et on note $\tE$ le topos des faisceaux de $\mU$-ensembles sur $E$. On désigne par 
\begin{eqnarray}
\sigma\colon \tE \rightarrow X_\et,\label{mtfla7b}\\
\beta\colon \tE \rightarrow \oX^\circ_\fet,\label{mtfla7c}\\
\psi\colon \oX^\circ_\et\rightarrow \tE,\label{mtfla7d}\\
\rho\colon X_\et\gtimes_{X_\et}\oX^\circ_\et\rightarrow \tE,\label{mtfla7e}
\end{eqnarray}
les morphismes canoniques \eqref{tf1kk}, \eqref{tf1k}, \eqref{tf1m} et \eqref{tf3b}. 

On note $\tE_s$ le sous-topos fermé de $\tE$ complémentaire de l'ouvert $\sigma^*(X_\eta)$ \eqref{tf20}, 
\begin{equation}\label{mtfla7f}
\delta\colon \tE_s\rightarrow \tE
\end{equation} 
le plongement canonique et 
\begin{equation}\label{mtfla7g}
\sigma_s\colon \tE_s\rightarrow X_{s,\et}
\end{equation} 
le morphisme de topos induit par $\sigma$ \eqref{tf23c}. 

\subsection{}\label{mtfla8}
Pour tout $(V\rightarrow U)\in \ob(E)$, on note $\oU^V$ la fermeture intégrale de $\oU$ dans $V$.
On désigne par $\ocB$ le préfaisceau sur $E$ défini pour tout $(V\rightarrow U)\in \ob(E)$, par 
\begin{equation}\label{mtfla8a}
\ocB((V\rightarrow U))=\Gamma(\oU^V,\co_{\oU^V}).
\end{equation}
C'est un anneau de $\tE$ (\cite{agt} III.8.16). 
D'après (\cite{agt} III.8.17), on a un homomorphisme canonique 
\begin{equation}\label{mtfla8b}
\sigma^*(\hbar_*(\co_\oX))\rightarrow \ocB.
\end{equation}
Sauf mention explicite du contraire, on considère $\sigma$ \eqref{mtfla7b}
comme un morphisme de topos annelés
\begin{equation}\label{mtfla8c}
\sigma\colon (\tE,\ocB)\rightarrow (X_\et,\hbar_*(\co_\oX)).
\end{equation}

Notant encore $\co_\oK$ le faisceau constant de $\oX^\circ_\fet$ de valeur $\co_\oK$. 
Sauf mention explicite du contraire, on considère $\beta$ \eqref{mtfla7c} comme un morphisme de topos annelés
\begin{equation}\label{mtfla8d}
\beta\colon (\tE,\ocB)\rightarrow (\oX^\circ_\fet,\co_\oK).
\end{equation}

Pour tout entier $m\geq 1$, on pose 
\begin{equation}\label{mtfla8e}
\ocB_m=\ocB/p^m\ocB.
\end{equation}
C'est un anneau de $\tE_s$ (\cite{agt} III.9.7). On désigne par
\begin{equation}\label{mtfla8f}
\sigma_m\colon (\tE_s,\ocB_m)\rightarrow (X_{s,\et},\co_{\oX_m})
\end{equation}
le morphisme de topos annelés induit par $\sigma$ \eqref{mtfla8c} (cf. \cite{agt} (III.9.9.4)) et par 
\begin{equation}\label{mtfla8g}
\Sigma_m\colon (\tE_s,\ocB_m)\rightarrow (X_{s,\zar},\co_{\oX_m})
\end{equation}
le composé de $\sigma_m$ et du morphisme canonique \eqref{notconv12j} 
\begin{equation}\label{mtfla8h}
u_m\colon (X_{s,\et},\co_{\oX_m})\rightarrow (X_{s,\zar},\co_{\oX_m}).
\end{equation}

Notant encore $\co_\oK/p^m\co_\oK$ le faisceau constant de $\oX^\circ_\fet$ de valeur $\co_\oK/p^m\co_\oK$, on désigne par 
\begin{equation}\label{mtfla8i}
\beta_m\colon (\tE_s,\ocB_m)\rightarrow (\oX^\circ_\fet,\co_\oK/p^m\co_\oK)
\end{equation}
le morphisme induit par $\beta$ \eqref{mtfla8d}.

\subsection{}\label{mtfla9}
On désigne par
\begin{equation}\label{mtfla9a}
\pi'\colon E'\rightarrow \Et_{/X'}
\end{equation}
le $\mU$-site fibré de Faltings associé au morphisme $h'\colon \oX'^\rhd\rightarrow X'$ (cf. \ref{tf1}). 
On munit $E'$ de la topologie co-évanescente définie par $\pi'$ (\cite{agt} VI.5.3)  
et on note $\tE'$ le topos des faisceaux de $\mU$-ensembles sur $E'$. On désigne par 
\begin{eqnarray}
\sigma'\colon \tE' \rightarrow X'_\et,\label{mtfla9b}\\
\beta'\colon \tE' \rightarrow \oX'^\rhd_\fet,\label{mtfla9c}\\
\psi'\colon \oX'^\rhd_\et\rightarrow \tE',\label{mtfla9d}\\
\rho'\colon X'_\et\gtimes_{X'_\et}\oX'^\rhd_\et\rightarrow \tE,\label{mtfla9de}
\end{eqnarray}
les morphismes canoniques \eqref{tf1kk}, \eqref{tf1k}, \eqref{tf1m} et \eqref{tf3b}. 

On note $\tE'_s$ le sous-topos fermé de $\tE'$ complémentaire de l'ouvert $\sigma'^*(X'_\eta)$ \eqref{tf20}, 
\begin{equation}\label{mtfla9e}
\delta'\colon \tE'_s\rightarrow \tE'
\end{equation} 
le plongement canonique et 
\begin{equation}\label{mtfla9f}
\sigma'_s\colon \tE'_s\rightarrow X'_{s,\et}
\end{equation} 
le morphisme de topos induit par $\sigma'$ \eqref{tf23c}.

\subsection{}\label{mtfla10}
Pour tout objet $(V\rightarrow U)$ de $E'$, on note $\oU^V$ la fermeture intégrale de $\oU$ dans $V$. 
On désigne par $\ocB'$ le préfaisceau sur $E'$ défini pour tout $(V\rightarrow U)\in \ob(E')$, par 
\begin{equation}\label{mtfla10a}
\ocB'((V\rightarrow U))=\Gamma(\oU^V,\co_{\oU^V}).
\end{equation}
C'est un anneau de $\tE'$ (\cite{agt} III.8.16). 
D'après (\cite{agt} III.8.17), on a un homomorphisme canonique 
\begin{equation}\label{mtfla10b}
\sigma'^*(\hbar'_*(\co_{\oX'}))\rightarrow \ocB'.
\end{equation}
Sauf mention explicite du contraire, on considère $\sigma'$ \eqref{mtfla9b}
comme un morphisme de topos annelés
\begin{equation}\label{mtfla10c}
\sigma'\colon (\tE',\ocB')\rightarrow (X'_\et,\hbar'_*(\co_{\oX'})).
\end{equation}

Notant encore $\co_\oK$ le faisceau constant de $\oX'^\rhd_\fet$ de valeur $\co_\oK$. 
Sauf mention explicite du contraire, on considère $\beta'$ \eqref{mtfla9c} comme un morphisme de topos annelés
\begin{equation}\label{mtfla10d}
\beta'\colon (\tE',\ocB')\rightarrow (\oX'^\rhd_\fet,\co_\oK).
\end{equation}

Pour tout entier $m\geq 1$, on pose
\begin{equation}\label{mtfla10e}
\ocB'_m=\ocB'/p^m\ocB'.
\end{equation}
C'est un objet de $\tE'_s$ (\cite{agt} III.9.7). On désigne par
\begin{equation}\label{mtfla10f}
\sigma'_m\colon (\tE'_s,\ocB'_m)\rightarrow (X'_{s,\et},\co_{\oX'_m})
\end{equation}
le morphisme de topos annelés induit par $\sigma'$ \eqref{mtfla10c} (cf. \cite{agt} (III.9.9.4)) et par 
\begin{equation}\label{mtfla10g}
\Sigma'_m\colon (\tE'_s,\ocB'_m)\rightarrow (X'_{s,\zar},\co_{\oX'_m})
\end{equation}
le composé de $\sigma'_m$ et du morphisme canonique \eqref{notconv12j} 
\begin{equation}\label{mtfla10h}
u'_m\colon (X'_{s,\et},\co_{\oX'_m})\rightarrow (X'_{s,\zar},\co_{\oX'_m}).
\end{equation}

Notant encore $\co_\oK/p^m\co_\oK$ le faisceau constant de $\oX'^\rhd_\fet$ de valeur $\co_\oK/p^m\co_\oK$, on désigne par 
\begin{equation}\label{mtfla10i}
\beta'_m\colon (\tE'_s,\ocB'_m)\rightarrow (\oX'^\rhd_\fet,\co_\oK/p^m\co_\oK)
\end{equation}
le morphisme induit par $\beta'$ \eqref{mtfla10d}.

\subsection{}\label{mtfla11}
Le foncteur 
\begin{equation}\label{mtfla11a}
\Theta^+\colon E\rightarrow E', \ \ \ (V\rightarrow U)\mapsto (V\times_{X^\circ}X'^\rhd\rightarrow U\times_XX').
\end{equation}
est continu et exact à gauche (\cite{agt} VI.10.12). Il définit donc un morphisme de topos 
\begin{equation}\label{mtfla11b}
\Theta\colon \tE'\rightarrow \tE.
\end{equation}
Il résulte aussitôt des définitions que les carrés du diagramme
\begin{equation}\label{mtfla11c}
\xymatrix{
{X'_\et}\ar[d]_{g}&{\tE'}\ar[l]_-(0.5){\sigma'}\ar[d]^{\Theta}\ar[r]^-(0.5){\beta'}&
{\oX'^\rhd_\fet}\ar[d]^{\upgamma}\\
{X_\et}&{\tE}\ar[l]_{\sigma}\ar[r]^{\beta}&{\oX^\circ_\fet}}
\end{equation}
où $\upgamma\colon \oX'^\rhd\rightarrow \oX^\circ$ est le morphisme induit par $g$, sont commutatifs à isomorphismes canoniques près.  

On a un isomorphisme canonique $\Theta^*(\sigma^*(X_\eta))\simeq \sigma'^*(X'_\eta)$ \eqref{mtfla11c}.
En vertu de (\cite{sga4} IV 9.4.3), il existe donc un morphisme de topos
\begin{equation}\label{mtfla11d}
\Theta_s\colon \tE'_s\rightarrow \tE_s
\end{equation}
unique à isomorphisme canonique près tel que le diagramme 
\begin{equation}\label{mtfla11e}
\xymatrix{
{\tE'_s}\ar[r]^{\Theta_s}\ar[d]_{\delta'}&{\tE_s}\ar[d]^{\delta}\\
{\tE'}\ar[r]^\Theta&{\tE}}
\end{equation}
soit commutatif à isomorphisme près, et même $2$-cartésien. Il résulte de \eqref{mtfla11c} et (\cite{sga4} IV 9.4.3) 
que le diagramme de morphismes de topos 
\begin{equation}\label{mtfla11f}
\xymatrix{
{\tE'_s}\ar[r]^{\Theta_s}\ar[d]_{\sigma'_s}&{\tE_s}\ar[d]^{\sigma_s}\\
{X'_{s,\et}}\ar[r]^{g_s}&{X_{s,\et}}}
\end{equation}
est commutatif à isomorphisme canonique près.

\subsection{}\label{mtfla12}
Pour tout $(V\rightarrow U)\in \ob(E)$, posons $(V'\rightarrow U')=\Theta^+(V\rightarrow U)$ de sorte qu'on a un diagramme commutatif 
\begin{equation}\label{mtfla12a}
\xymatrix{\oX'^\rhd\ar[d]&V'\ar@{}[ld]|{\Box}\ar[l]\ar[r]\ar[d]&\oU'\ar[d]\ar[r]\ar@{}[rd]|{\Box}&\oX'\ar[d]\\
\oX^\circ&V\ar[l]\ar[r]&\oU\ar[r]&\oX}
\end{equation}
On en déduit un morphisme 
\begin{equation}\label{mtfla12b}
\oU'^{V'}\rightarrow \oU^V,
\end{equation}
et par suite un homomorphisme d'anneaux de $\tE$
\begin{equation}\label{mtfla12c}
\ocB\rightarrow \Theta_*(\ocB').
\end{equation}
Nous considérons dans la suite $\Theta$ \eqref{mtfla11b} comme un morphisme de topos annelés (respectivement par $\ocB'$ et $\ocB$). 
Nous utilisons pour les modules la notation $\Theta^{-1}$ pour désigner l'image
inverse au sens des faisceaux abéliens et nous réservons la notation 
$\Theta^*$ pour l'image inverse au sens des modules.

Pour tout entier $m\geq 1$, l'homomorphisme canonique $\Theta^{-1}(\ocB)\rightarrow \ocB'$
induit un homomorphisme $\Theta_s^*(\ocB_m)\rightarrow \ocB'_m$. 
Le morphisme $\Theta_s$ est donc sous-jacent à un morphisme de topos annelés, que l'on note 
\begin{equation}\label{mtfla12d}
\Theta_m\colon (\tE'_s,\ocB'_m)\rightarrow (\tE_s,\ocB_m).
\end{equation}

\begin{lem}\label{mtfla4}
Soient $(P,\gamma)$ une carte fine pour $(X,\cM_X)$ {\rm (\cite{agt} II.5.13)}, $\ox'$ un point géométrique de $X'$, $\ox=g(\ox')$.  
Alors, il existe un voisinage étale $U'$
de $\ox'$ dans $X'$, une carte $(P',\gamma')$ pour $(U',\cM_{X'}|U')$ et un homomorphisme $h\colon P\rightarrow P'$ 
de monoïdes tels que les conditions suivantes soient remplies~:  
\begin{itemize}
\item[{\rm (i)}]  $((P',\gamma'),(P,\gamma),h)$ est une carte pour la restriction
$g_{U'}\colon (U',\cM_{X'}|U')\rightarrow (X,\cM_X)$ de $g$ à $U'$, autrement dit,
le diagramme d'homomorphismes de monoïdes 
\begin{equation}\label{mtfla4a}
\xymatrix{
P'\ar[r]^-(0.5){\gamma'}&{\Gamma(U',\cM_{X'})}\\
P\ar[r]^-(0.5){\gamma}\ar[u]^h&{\Gamma(X,\cM_X)}\ar[u]_{g^\flat_{U'}}}
\end{equation}
est commutatif, ou ce qui revient au même le diagramme associé de morphismes de schémas logarithmiques 
\begin{equation}\label{mtfla4b}
\xymatrix{
{(U',\cM_{X'}|U')}\ar[r]^-(0.5){\gamma'^a}\ar[d]_{g_{U'}}&{\bA_{P'}}\ar[d]^{\bA_h}\\
{(X,\cM_X)}\ar[r]^-(0.5){\gamma^a}&{\bA_P}}
\end{equation}
est commutatif \eqref{notconv2}.
\item[{\rm (ii)}]  $P'$ est fin et saturé {\rm (\cite{agt} II.5.1)}. Si, de plus, $P^\gp$ est libre sur $\mZ$, il en est de même de $P'^\gp$. 
\item[{\rm (iii)}] L'homomorphisme $h^\gp\colon P^\gp\rightarrow P'^\gp$ est injectif, 
le sous-groupe de torsion de $\coker(h^\gp)$ est fini d'ordre inversible dans $\co_{X',\ox'}$ et le morphisme de schémas usuels
\begin{equation}\label{mtfla4c}
U'\rightarrow X\times_{\bA_P}\bA_{P'}
\end{equation}
déduit de \eqref{mtfla4b} est étale.  
\item[{\rm (iv)}] Il existe un sous-groupe $A'$ de $P'$ tel que $\gamma'$ induise un isomorphisme
\begin{equation}\label{mtfla4d}
P'/A'\stackrel{\sim}{\rightarrow} \cM_{X',\ox'}/\co_{X',\ox'}^\times.
\end{equation}
\item[{\rm (v)}] Si, de plus, il existe un sous-groupe $A$ de $P$ tel que $\gamma$ induise un isomorphisme
\begin{equation}\label{mtfla4e}
P/A\stackrel{\sim}{\rightarrow} \cM_{X,\ox}/\co_{X,\ox}^\times,
\end{equation}
alors l'homomorphisme $h$ est saturé {\rm (\cite{agt} II.5.2)}. 
\end{itemize}
\end{lem}

Nous reprenons la preuve de (\cite{fkato} 4.1; cf. § 6) en l'adaptant (cf. aussi \cite{agt} III.4.3). 
Soient $t_1,\dots,t_r\in \cM_{X',\ox'}$ tels que $d\log(t_1),\dots,d\log(t_r)$ forment une base du $\co_{X',\ox'}$-module 
\begin{equation}
\Omega^1_{(X',\cM_{X'})/(X,\cM_X),\ox'}. 
\end{equation}
Posons $H=\mN^r\oplus P$ et considérons l'homomorphisme
\begin{equation}
\varphi\colon H=\mN^r\oplus P\rightarrow \cM_{X',\ox'}, \ \ \ (n_1,\dots,n_r)+t\mapsto \prod_{i=1}^r t^{n_i}_i \cdot g^\flat(\gamma(t)).
\end{equation}
Notons $\alpha\colon \cM_{X',\ox'}^\gp\rightarrow \cM_{X',\ox'}^\gp/\co_{X',\ox'}^\times$ la projection 
canonique et $L$ l'image de l'homomorphisme
\begin{equation}
\alpha\circ \varphi^\gp \colon H^\gp\rightarrow \cM_{X',\ox'}^\gp/\co_{X',\ox'}^\times.
\end{equation}
Comme $\cM_{X'}$ est fin et saturé (\cite{tsuji4} II 2.12), 
$\cM_{X',\ox'}^\gp/\co_{X',\ox'}^\times$ est un $\mZ$-module libre de type fini. 
D'après l'étape 2 de (\cite{fkato} page 331), on voit que le conoyau de $\alpha\circ \varphi^\gp$ 
est annulé par un entier inversible dans $\co_{X',\ox'}$. 
Il existe donc une $\mZ$-base $e_1,\dots,e_d$ de $\cM_{X',\ox'}^\gp/\co_{X',\ox'}^\times$ et des entiers $f_1,\dots,f_d$ tels que 
$e_1^{f_1},\dots,e_d^{f_d}$ forment une $\mZ$-base de $L$, que $f_i$ divise $f_{i+1}$ pour tout $1\leq i\leq d-1$ et que $f_d$ soit inversible dans 
$\co_{X',\ox'}$. Soient $F_1,\dots,F_d\in H^\gp$ et $\te_1,\dots,\te_d\in \cM^\gp_{X',\ox'}$ 
tels que $\alpha(\varphi^\gp(F_i))=e_i^{f_i}$  et $\alpha(\te_i)=e_i$ pour tout $1\leq i\leq d$. 
Il existe alors $u_i\in \co_{X',\ox'}^\times$ tel que 
$\varphi^\gp(F_i)=u_i\te_i^{f_i}$. Comme $\co_{X',\ox'}^\times$ est $f_i$-divisible, il existe $v_i\in \co_{X',\ox'}^\times$
tel que $v_i^{f_i}=u_i$. Remplaçant $\te_i$ par $v_i\te_i$, on peut supposer que $\varphi^\gp(F_i)=\te_i^{f_i}$.
On désigne par $\beta\colon \cM_{X',\ox'}^\gp/\co_{X',\ox'}^\times\rightarrow \cM_{X',\ox'}^\gp$ le scindage de $\alpha$
défini par $\beta(e_i)=\te_i$ pour tout $1\leq i\leq d$,  
par $\rho\colon H^\gp\rightarrow L$ l'homomorphisme surjectif induit par $\alpha\circ \varphi^\gp$, 
par $\sigma\colon L\rightarrow H^\gp$ le scindage de $\rho$ défini par 
$\sigma(e_i^{f_i})=F_i$ pour tout $1\leq i\leq d$, par $M$ le noyau de $\rho$, 
et par $\tau\colon H^\gp\rightarrow M$ l'homomorphisme qui à tout $h\in H^\gp$ associe $h-\sigma(\rho(h))$. 
On pose $G=M\oplus \cM_{X',\ox'}^\gp/\co_{X',\ox'}^\times$ et
\begin{equation}
\phi\colon G=M\oplus \cM_{X',\ox'}^\gp/\co_{X',\ox'}^\times \rightarrow \cM^\gp_{X',\ox'}, 
\ \ \ (m,t)\mapsto \phi(m,t)=\varphi^\gp(m)\cdot \beta(t). 
\end{equation} 
On a $\varphi^\gp = \phi\circ (\tau\oplus \alpha\circ \varphi^\gp)\colon H^\gp\rightarrow \cM_{X',\ox'}^\gp$. On le vérifie immédiatement
pour les éléments de $M$ et pour les éléments $(F_i)_{1\leq i\leq d}$. On pose $P'=\phi^{-1}(\cM_{X',\ox'})$. D'après (\cite{kato1}
2.10), il existe un voisinage étale $U'$ de $\ox'$ dans $X'$ et un homomorphisme $\gamma'\colon P'\rightarrow \Gamma(U',\cM_{X'})$
qui est une carte pour $(U',\cM_{X'}|U')$ et dont la fibre $\gamma'_{\ox'}\colon P'\rightarrow \cM_{X',\ox'}$ en $\ox'$ est induite par $\phi$. 
L'homomorphisme $\tau\oplus \alpha\circ \varphi^\gp\colon H^\gp\rightarrow G$ 
induit un homomorphisme $H\rightarrow P'$ et par suite un homomorphisme
$h\colon P\rightarrow P'$ qui rend commutatif le diagramme \eqref{mtfla4a}. 

Comme l'homomorphisme $G\rightarrow \cM^\gp_{X',\ox'}/\co_{X',\ox'}^\times$ induit par $\phi$ est surjectif, 
$P'^\gp=G$ et $P'$ est intègre.  
Il résulte aussitôt de la définition qu'il existe un sous-groupe $A'$ de $P'$ tel que $\gamma'_{\ox'}$ induise un isomorphisme  
\begin{equation}\label{mtfla4f}
P'/A'\stackrel{\sim}{\rightarrow} \cM_{X',\ox'}/\co_{X',\ox'}^\times.
\end{equation}
Par suite, $A'=P'^\times$. 
Comme $\cM_{X',\ox'}$ est saturé, $\cM_{X',\ox'}/\co_{X',\ox'}^\times$ est saturé et donc $P'$ est saturé (\cite{agt} II.5.1). 
Si $P^\gp$ est libre sur $\mZ$, il en est de même de $H^\gp$ et donc de $G=P'^\gp$.  

L'homomorphisme $\tau\oplus \alpha\circ \varphi^\gp\colon H^\gp\rightarrow G$ est injectif et 
son conoyau est isomorphe à celui de $\alpha\circ \varphi^\gp$. 
On en déduit que l'homomorphisme $h^\gp\colon P^\gp\rightarrow P'^\gp$ est injectif et que le sous-groupe de torsion de 
$\coker(h^\gp)$ est fini d'ordre inversible dans $\co_{X',\ox'}$. L'étape 4 de (\cite{fkato} § 6 page 332) montre que 
le morphisme de schémas usuels $U'\rightarrow X\times_{\bA_P}\bA_{P'}$ déduit de \eqref{mtfla4b} est étale.

Sous les hypothèses de (v), l'homomorphisme $h$ est saturé en vertu de \eqref{mtfla4f} et (\cite{tsuji4} I 3.16).

\begin{defi}\label{mtfla5}
On appelle {\em carte relativement adéquate} pour le morphisme $g$ la donnée d'une carte $(\mN,\iota)$ pour $(S,\cM_S)$, 
d'une carte $(P,\gamma)$ pour $(X,\cM_X)$, d'une carte $(P',\gamma')$ pour $(X',\cM_{X'})$ et de deux homomorphismes de monoïdes
$\vartheta\colon \mN\rightarrow P$ et $h\colon P\rightarrow P'$ tels que les conditions suivantes soient remplies:
\begin{itemize}
\item[{\rm (i)}] $((P',\gamma'),(P,\gamma),h)$ est une carte pour $g$ (\cite{agt} II.5.14). 
\item[{\rm (ii)}] L'homomorphisme $h$ est saturé  (\cite{agt} II.5.2). 
\item[{\rm (iii)}] L'homomorphisme $h^\gp\colon P^\gp\rightarrow P'^\gp$ est injectif, 
le sous-groupe de torsion de $\coker(h^\gp)$ est fini d'ordre premier à $p$ et le morphisme de schémas usuels
\begin{equation}\label{mtfla5a}
X'\rightarrow X\times_{\bA_P}\bA_{P'}
\end{equation}
déduit de la carte (i) est étale.  
\item[{\rm (iv)}] $((P,\gamma),(\mN,\iota),\vartheta)$ est une carte adéquate pour $f$ (\cite{agt} III.4.4).
\item[{\rm (v)}] $((P',\gamma'),(\mN,\iota),h\circ \vartheta)$ est une carte adéquate pour $f'$.
\end{itemize}
\end{defi}

\begin{prop}\label{mtfla6}
Pour tout point géométrique $\ox'$ de $X'$ au-dessus de $s$, 
il existe un diagramme commutatif de morphismes de schémas logarithmiques
\begin{equation}\label{mtfla6a}
\xymatrix{
{(U',\cM_{U'})}\ar[r]\ar[d]_{g^\dagger}&{(X',\cM_{X'})}\ar[d]^g\\
{(U,\cM_U)}\ar[r]&{(X,\cM_X)}}
\end{equation}
tel que les flèches horizontales soient étales et strictes, que $U'$ soit un voisinage étale de $\ox'$ dans $X'$, et que $g^\dagger$
admette une carte relativement adéquate $((P',\gamma'),(P,\gamma),(\mN,\iota),\vartheta\colon \mN\rightarrow P, h\colon P\rightarrow P')$ \eqref{mtfla5}. 
\end{prop}

En effet, d'après (\cite{agt} III.4.6) et la preuve de (c)$\Rightarrow$(a), il existe un voisinage étale $U$ de $\ox=g(\ox')$ dans $X$, une carte adéquate 
$((P,\gamma),(\mN,\iota),\vartheta\colon \mN\rightarrow P)$ de la restriction $f_U\colon (U,\cM_X|U)\rightarrow (S,\cM_S)$ de $f$ à $U$
et un sous-groupe $A$ de $P$ tels que $\gamma$ induise un isomorphisme
\begin{equation}\label{mtfla6b}
P/A\stackrel{\sim}{\rightarrow} \cM_{X,\ox}/\co_{X,\ox}^\times.
\end{equation}
Remplaçant $g$ par le changement de base par le morphisme strict $(U,\cM_X|U)\rightarrow (X,\cM_X)$, on peut supposer qu'il existe 
une carte adéquate $((P,\gamma),(\mN,\iota),\vartheta\colon \mN\rightarrow P)$ pour $f$ 
et un sous-groupe $A$ de $P$ tels que $\gamma$ induise un isomorphisme
\begin{equation}\label{mtfla6c}
P/A\stackrel{\sim}{\rightarrow} \cM_{X,\ox}/\co_{X,\ox}^\times.
\end{equation}
D'après \ref{mtfla4}, il existe alors un voisinage étale $U'$ de $\ox'$ dans $X'$ et une carte 
\[
((P',\gamma'),(P,\gamma),h\colon P\rightarrow P')
\]
pour la restriction $g_{U'}\colon (U',\cM_{X'}|U')\rightarrow (X,\cM_X)$ de $g$ à $U'$ tels que les conditions suivantes soient remplies:
\begin{itemize}
\item[(i)] Le monoïde $P'$ est torique (\cite{agt} II.5.1) et l'homomorphisme $h$ est saturé.
\item[(ii)] L'homomorphisme $h^\gp\colon P^\gp\rightarrow P'^\gp$ est injectif, 
le sous-groupe de torsion de $\coker(h^\gp)$ est fini d'ordre premier à $p$ et le morphisme de schémas usuels
\begin{equation}
U'\rightarrow X\times_{\bA_P}\bA_{P'}
\end{equation}
induit par la carte est étale. 
\item[(iii)] Il existe un sous-groupe $A'$ de $P'$ tel que $\gamma'$ induise un isomorphisme
\begin{equation}\label{mtfla6d}
P'/A'\stackrel{\sim}{\rightarrow} \cM_{X',\ox'}/\co_{X',\ox'}^\times.
\end{equation}
\end{itemize}
Il résulte alors de la preuve de (c)$\Rightarrow$(a) de (\cite{agt} III.4.6) que 
$((P',\gamma'),(\mN,\iota),h\circ \vartheta)$ est une carte adéquate pour $f'_{U'}=f\circ g_{U'}$.

\section{Cohomologie galoisienne dans le cas relatif}\label{eccr}

\subsection{}\label{eccr1}
Dans cette section, on suppose que les schémas $X=\Spec(R)$ et $X'=\Spec(R')$ sont affines et connexes, 
que $X'_s$ est non-vide et que le morphisme $g$ \eqref{mtfla2a} admet une carte relativement adéquate \eqref{mtfla5}
\begin{equation}\label{eccr1a}
((P',\gamma'),(P,\gamma),(\mN,\iota),\vartheta\colon \mN\rightarrow P, h\colon P\rightarrow P').
\end{equation} 
On pose $\pi=\iota(1)$ qui est une uniformisante de $\co_K$, 
\begin{equation}\label{eccr1b}
\tOmega^1_{R/\co_K}=\tOmega^1_{X/S}(X),  \ \ \ \tOmega^1_{R'/\co_K}=\tOmega^1_{X'/S}(X')  \ \ \ {\rm et} \ \ \
\tOmega^1_{R'/R}=\tOmega^1_{X'/X}(X').
\end{equation}
On a alors une suite exacte scindée de $R'$-modules \eqref{mtfla2f}
\begin{equation}\label{eccr1c}
0\rightarrow \tOmega^1_{R/\co_K}\otimes_RR'\rightarrow \tOmega^1_{R'/\co_K}\rightarrow \tOmega^1_{R'/R}\rightarrow 0.
\end{equation}

\subsection{}\label{eccr20}
Pour tout entier $n\geq 1$, on pose
\begin{eqnarray}\label{eccr20a}
\co_{K_n}=\co_K[\zeta]/(\zeta^{n}-\pi),
\end{eqnarray}
qui est un anneau de valuation discrète. On note $K_n$ le corps des fractions de $\co_{K_n}$
et $\pi_n$ la classe de $\zeta$ dans $\co_{K_n}$, qui est une uniformisante de $\co_{K_n}$.  
On pose $S^{(n)}=\Spec(\co_{K_n})$
que l'on munit de la structure logarithmique $\cM_{S^{(n)}}$ définie par son point fermé. 
On désigne par $\iota_n\colon \mN\rightarrow \Gamma(S^{(n)},\cM_{S^{(n)}})$
l'homomorphisme défini par $\iota_n(1)=\pi_n$; c'est une carte pour $(S^{(n)},\cM_{S^{(n)}})$.  

Considérons le système inductif de monoïdes $(\mN^{(n)})_{n\geq 1}$, 
indexé par l'ensemble $\mZ_{\geq 1}$ ordonné par la relation de divisibilité, 
défini par $\mN^{(n)}=\mN$ pour tout $n\geq 1$ et dont l'homomorphisme de transition
$\mN^{(n)}\rightarrow \mN^{(mn)}$ (pour $m,n\geq 1$) est l'endomorphisme de Frobenius d'ordre $m$ 
de $\mN$ ({\em i.e.}, l'élévation à la puissance $m$-ième). On notera $\mN^{(1)}$ simplement $\mN$. Les  schémas logarithmiques
$(S^{(n)},\cM_{S^{(n)}})_{n\geq 1}$ forment naturellement un système projectif.
Pour tous entiers $m, n\geq 1$, avec les notations de  \ref{notconv2}, on a un diagramme cartésien de morphismes de schémas logarithmiques
\begin{equation}\label{eccr20b}
\xymatrix{
{(S^{(mn)},\cM_{S^{(mn)}})}\ar[r]^-(0.5){\iota^a_{mn}}\ar[d]&{\bA_{\mN^{(mn)}}}\ar[d]\\
{(S^{(n)},\cM_{S^{(n)}})}\ar[r]^-(0.5){\iota^a_n}&{\bA_{\mN^{(n)}}}}
\end{equation}
où $\iota^a_n$ (resp. $\iota^a_{mn}$) est le morphisme associé à $\iota_n$ (resp. $\iota_{mn}$) (\cite{agt} II.5.13). 

\subsection{}\label{eccr2}
On désigne par $(P^{(n)})_{n\geq 1}$ le système inductif de monoïdes,
indexé par l'ensemble $\mZ_{\geq 1}$ ordonné par la relation de divisibilité, 
défini par $P^{(n)}=P$ pour tout $n\geq 1$ et dont l'homomorphisme de transition
$i_{n,mn}\colon P^{(n)}\rightarrow P^{(mn)}$ (pour $m, n\geq 1$) est l'endomorphisme 
de Frobenius d'ordre $m$ de $P$ ({\em i.e.}, l'élévation à la puissance $m$-ième). On note $P^{(1)}$ simplement $P$.

Pour tout $n\geq 1$, on pose (avec les notations de \ref{notconv2})
\begin{equation}\label{eccr2a}
(X^{(n)},\cM_{X^{(n)}})=(X,\cM_{X})\times_{\bA_P}\bA_{P^{(n)}}.
\end{equation}
On rappelle \eqref{cad6c} qu'il existe un unique morphisme 
\begin{equation}\label{eccr2b}
f^{(n)}\colon (X^{(n)},\cM_{X^{(n)}})\rightarrow (S^{(n)},\cM_{S^{(n)}})
\end{equation} 
qui s'insère dans le diagramme commutatif
\begin{equation}\label{eccr2c}
\xymatrix{
{(X^{(n)},\cM_{X^{(n)}})}\ar[rrr]\ar[ddd]\ar[rd]_{f^{(n)}}&&&
{\bA_{P^{(n)}}}\ar[ld]^{\bA_{\vartheta}}\ar[ddd]\\
&{(S^{(n)},\cM_{S^{(n)}})}\ar[r]^-(0.5){\iota^a_n}\ar[d]&{\bA_{\mN^{(n)}}}\ar[d]&\\
&{(S,\cM_{S})}\ar[r]^-(0.5){\iota^a}&{\bA_\mN}&\\
{(X,\cM_{X})}\ar[rrr]\ar[ru]^f&&&{\bA_P}\ar[lu]_{\bA_\vartheta}}
\end{equation}
D'après \ref{cad7}, le morphisme $f^{(n)}$ est lisse et saturé et le sous-schéma ouvert maximal de $X^{(n)}$ 
où la structure logarithmique $\cM_{X^{(n)}}$ est triviale est égal à $X^{(n)\circ}=X^{(n)}\times_{X}X^\circ$ \eqref{mtfla2b}.
D'après (\cite{agt} III.4.2), le schéma $X^{(n)}$ est donc normal et $S^{(n)}$-plat, et l'immersion $X^{(n)\circ}\rightarrow X^{(n)}$ est schématiquement dominante.

\subsection{}\label{eccr3}
On désigne par $(P'^{[n]})_{n\geq 1}$ le système inductif de monoïdes 
indexé par l'ensemble $\mZ_{\geq 1}$ ordonné par la relation de divisibilité, 
défini par $P'^{[n]}=P'$ pour tout $n\geq 1$ et dont l'homomorphisme de transition
$P'^{[n]}\rightarrow P'^{[mn]}$ (pour $m, n\geq 1$) est l'endomorphisme 
de Frobenius d'ordre $m$ de $P'$. On note $P'^{[1]}$ simplement $P'$. 

Pour tout $n\geq 1$, on pose (avec les notations de \ref{notconv2})
\begin{equation}\label{eccr3a}
(X'^{[n]},\cM_{X'^{[n]}})=(X',\cM_{X'})\times_{\bA_{P'}}\bA_{P'^{[n]}}.
\end{equation}
Il existe alors un unique morphisme
\begin{equation}\label{eccr3b}
f'^{[n]}\colon (X'^{[n]},\cM_{X'^{[n]}})\rightarrow (S^{(n)},\cM_{S^{(n)}}) 
\end{equation} 
qui s'insère dans le diagramme commutatif 
\begin{equation}\label{eccr3c}
\xymatrix{
{(X'^{[n]},\cM_{X'^{[n]}})}\ar[rrr]\ar[ddd]\ar[rd]_{f'^{[n]}}&&&
{\bA_{P'^{[n]}}}\ar[ld]^{\bA_{h\circ \vartheta}}\ar[ddd]\\
&{(S^{(n)},\cM_{S^{(n)}})}\ar[r]\ar[d]&{\bA_{\mN^{(n)}}}\ar[d]&\\
&{(S,\cM_{S})}\ar[r]&{\bA_\mN}&\\
{(X',\cM_{X'})}\ar[rrr]\ar[ru]^{f'}&&&{\bA_{P'}}\ar[lu]_{\bA_{h\circ \vartheta}}}
\end{equation}
D'après \ref{cad7}, le morphisme $f'^{[n]}$ est lisse et saturé et le sous-schéma ouvert maximal de $X'^{[n]}$ 
où la structure logarithmique $\cM_{X'^{[n]}}$ est triviale est égal à $X'^{[n]\rhd}=X'^{[n]}\times_{X'}X'^\rhd$ \eqref{mtfla2c}.
D'après (\cite{agt} III.4.2), le schéma $X'^{[n]}$ est donc normal et $S^{(n)}$-plat, et l'immersion $X'^{[n]\rhd}\rightarrow X'^{[n]}$ est schématiquement dominante.

\subsection{}\label{eccr4}
Pour tout $n\geq 1$, on pose 
\begin{equation}\label{eccr4a}
(X'^{(n)},\cM_{X'^{(n)}})=(X^{(n)},\cM_{X^{(n)}})\times_{(X,\cM_X)}(X',\cM_{X'}),
\end{equation}
le produit étant indifféremment pris dans la catégorie des schémas logarithmiques ou dans celle des schémas logarithmiques saturés (\cite{agt} II.5.20). 
Le morphisme 
\begin{equation}\label{eccr4b}
f'^{(n)}\colon (X'^{(n)},\cM_{X'^{(n)}})\rightarrow (S^{(n)},\cM_{S^{(n)}}),
\end{equation}
déduit de $f^{(n)}$ \eqref{eccr2b}, est donc lisse et saturé (\cite{tsuji4} II 2.11). 
Il résulte de \ref{cad7}(iv) que $X'^{(n)\rhd}=X'^{(n)}\times_{X'}X'^\rhd$  \eqref{mtfla2c}
est le sous-schéma ouvert maximal de $X'^{(n)}$ où la structure logarithmique $\cM_{X'^{(n)}}$ est triviale (\cite{ogus} III 1.2.8).
D'après (\cite{agt} III.4.2), le schéma $X'^{(n)}$ est normal et $S^{(n)}$-plat, et l'immersion $X'^{(n)\rhd}\rightarrow X'^{(n)}$ est schématiquement dominante.

On désigne par $P'^{(n)}$ la somme amalgamée des homomorphismes de monoïdes $h \colon P\rightarrow P'$ et 
$P\rightarrow P^{(n)}$. Comme $h$ est saturé, le diagramme 
\begin{equation}\label{eccr4c}
\xymatrix{
P\ar[r]\ar[d]&{P^{(n)}}\ar[d]\\
P'\ar[r]&{P'^{(n)}}}
\end{equation}
est cocartésien indifféremment dans la catégorie des monoïdes ou dans celle des monoïdes fins et saturés (\cite{agt} II.5.2). 
Le diagramme
\begin{equation}\label{eccr4d}
\xymatrix{
P^\gp\ar[r]\ar[d]&{P^{(n)\gp}}\ar[d]\\
P'^\gp\ar[r]&{P'^{(n)\gp}}}
\end{equation}
est donc cocartésien dans la catégorie des groupes abéliens (\cite{sga4} I 2.11), 
et $P'^{(n)}$ s'identifie au sous-monoïde fin de $P'^{(n)\gp}$ engendré par les images de $P'$ et $P^{(n)}$. 
Considérant $P'^{[n]}$ comme un monoïde au-dessus de $P^{(n)}$ via l'homomorphisme $h$, on a un homomorphisme canonique 
\begin{equation}\label{eccr4e}
P'^{(n)}\rightarrow P'^{[n]}
\end{equation}
au-dessus de $P'$ et $P^{(n)}$, et un morphisme strict canonique \eqref{notconv2}
\begin{equation}\label{eccr4f}
(X'^{(n)},\cM_{X'^{(n)}})\rightarrow \bA_{P'^{(n)}}.
\end{equation}

Il existe un unique $(S^{(n)},\cM_{S^{(n)}})$-morphisme 
\begin{equation}\label{eccr4g}
g^{(n)}\colon (X'^{[n]},\cM_{X'^{[n]}})\rightarrow (X^{(n)},\cM_{X^{(n)}})
\end{equation}
au-dessus de $g$ et du morphisme $\bA_{h}\colon \bA_{P'^{[n]}}\rightarrow \bA_{P^{(n)}}$. Celui-ci induit un morphisme canonique 
\begin{equation}\label{eccr4h}
(X'^{[n]},\cM_{X'^{[n]}})\rightarrow (X'^{(n)},\cM_{X'^{(n)}})
\end{equation}
au-dessus de $(X',\cM_{X'})$, de $(X^{(n)},\cM_{X^{(n)}})$ et du morphisme $\bA_{P'^{[n]}}\rightarrow \bA_{P'^{(n)}}$ \eqref{eccr4e}.

\begin{prop}\label{eccr5}
Soit $n$ un entier $\geq 1$
On note $C$ le conoyau de l'homomorphisme $h^\gp\colon P^\gp\rightarrow P'^\gp$ et $C[n]$ le sous-groupe de $n$-torsion de $C$. 
\begin{itemize}
\item[{\rm (i)}] Les morphismes $f^{(n)}$ \eqref{eccr2b} et $f'^{[n]}$ \eqref{eccr3b} sont adéquats, le morphisme $g^{(n)}$ \eqref{eccr4g}
est lisse et saturé, et 
\begin{equation}\label{eccr5a}
(P'^{[n]},P^{(n)},\mN^{(n)},\vartheta\colon \mN^{(n)}\rightarrow P^{(n)},h\colon P^{(n)}\rightarrow P'^{[n]})
\end{equation} 
est une carte relativement adéquate pour $g^{(n)}$.
\item[{\rm (ii)}] Le noyau (resp. conoyau) de 
l'homomorphisme canonique $P'^{(n)\gp}\rightarrow P'^{[n]\gp}$ \eqref{eccr4e} 
est canoniquement isomorphisme à $C[n]$ (resp. $C/nC$). 
En particulier, le morphisme de schémas en groupes $\Spec(K[P'^{[n]\gp}])\rightarrow \Spec(K[P'^{(n)\gp}])$ est fini et étale. 
\item[{\rm (iii)}] Les carrés du diagramme d'homomorphismes de monoïdes \eqref{eccr4e}
\begin{equation}\label{eccr5b}
\xymatrix{
{P'}\ar[r]\ar[d]&{P'^{(n)}}\ar[r]\ar[d]&{P'^{[n]}}\ar[d]\\
{P'^\gp}\ar[r]&{P'^{(n)\gp}}\ar[r]&{P'^{[n]\gp}}}
\end{equation}
sont cocartésiens. 
\item[{\rm (iv)}] Les carrés du diagramme canonique
\begin{equation}\label{eccr5c}
\xymatrix{
{X'^{[n]}}\ar[d]\ar[r]&{\bA_{P'^{[n]}}}\ar[d]\\
{X'^{(n)}}\ar[d]\ar[r]&{\bA_{P'^{(n)}}}\ar[d]\\
{X'}\ar[r]&{\bA_{P'}}}
\end{equation}
sont cartésiens.
\item[{\rm (v)}] Les carrés du diagramme canonique
\begin{equation}\label{eccr5d}
\xymatrix{
{X'^{[n]\rhd}}\ar[d]\ar[r]&{\bA_{P'^{[n]\gp}}}\ar[d]\\
{X'^{(n)\rhd}}\ar[r]\ar[d]&{\bA_{P'^{(n)\gp}}}\ar[d]\\
{X'^\rhd}\ar[r]&{\bA_{P'^{\gp}}}}
\end{equation}
sont cartésiens. 
\item[{\rm (vi)}] Il existe un sous-schéma ouvert et fermé $X'^{\{n\}}$ de $X'^{(n)}$ tel que le morphisme canonique 
$X'^{[n]}\rightarrow X'^{(n)}$ se factorise à travers un morphisme fini et surjectif $X'^{[n]}\rightarrow X'^{\{n\}}$.
De plus, posant $X'^{\{n\}\rhd}=X'^{\{n\}}\times_{X'}X'^\rhd$, le morphisme $X'^{[n]\rhd}\rightarrow X'^{\{n\}\rhd}$ est étale. 
\end{itemize}
\end{prop}

(i) En effet, les morphismes $f^{(n)}$ et $f'^{[n]}$ sont adéquats d'après \ref{cad7}(i).
Il est clair que \eqref{eccr5a} est une carte pour le morphisme $g^{(n)}$, et elle est relativement adéquate \eqref{mtfla5} en vertu de \ref{cad7}(i) 
et du diagramme commutatif (de schémas non-logarithmiques) dont les flèches verticales sont des isomorphismes
\begin{equation}
\xymatrix{
{X'^{[n]}}\ar[r]\ar[d]&{X^{(n)}\times_{\bA_{P^{(n)}}}\bA_{P'^{[n]}}}\ar[d]\\
{X'\times_{\bA_{P'}}\bA_{P'^{[n]}}}\ar[r]&{X\times_{\bA_P}\bA_{P'^{[n]}}}}
\end{equation}

(ii) En effet, on a un diagramme commutatif d'homomorphismes canoniques dont les lignes sont exactes
\begin{equation}\label{eccr5e}
\xymatrix{
{P^\gp}\ar@{^(->}[r]\ar[d]&{P'^\gp}\ar@{->>}[r]\ar[d]&C\ar@{=}[d]\\
{P^{(n)\gp}}\ar@{=}[d]\ar@{^(->}[r]&{P'^{(n)\gp}}\ar[d]\ar@{->>}[r]&C\ar[d]^{\cdot n}\\
{P^{(n)\gp}}\ar@{^(->}[r]&{P'^{[n]\gp}}\ar@{->>}[r]&C}
\end{equation}
La première assertion s'ensuit et elle implique la seconde assertion compte tenu de (\cite{ogus} IV 3.1.10).

(iii)  Considérons le diagramme commutatif d'homomorphismes de monoïdes
\begin{equation}
\xymatrix{
{P}\ar[rrr]\ar[ddd]\ar[rd]\ar@{}[drrr]|*+[o][F-]{1}&&&{P^{(n)}}\ar[ld]\ar[ddd]\\
&{P^\gp}\ar[r]\ar[d]\ar@{}[dr]|*+[o][F-]{2}&{P^{(n)\gp}}\ar[d]&\\
&{P'^\gp}\ar[r]&{P'^{(n)\gp}}&\\
{P'}\ar[rrr]\ar[ru]\ar@{}[urrr]|*+[o][F-]{3}&&&{P'^{(n)}}\ar[lu]}
\end{equation}
Le grand carré extérieur est cocartésien. 
La face $\xymatrix{\ar@{}|*+[o][F-]{2}}$ est cocartésienne 
puisque le foncteur oubli de la catégorie des groupes abéliens dans la catégorie des monoïdes préserve les
sommes amalgamées (cf. \cite{kato1} 1.3, Remark).  
La face $\xymatrix{\ar@{}|*+[o][F-]{1}}$ est cocartésienne puisque $P^{(n)\gp}$ est engendré, en tant que monoïde, 
par les images canoniques de $P^\gp$ et $P^{(n)}$ (cf. aussi la preuve de \cite{agt} II.6.6(vi)).
On en déduit que la face $\xymatrix{\ar@{}|*+[o][F-]{3}}$ est cocartésienne. 
De même que pour $\xymatrix{\ar@{}|*+[o][F-]{1}}$, le diagramme 
\begin{equation}
\xymatrix{
P'\ar[r]\ar[d]&P'^{[n]}\ar[d]\\
{P'^\gp}\ar[r]&{P'^{[n]\gp}}}
\end{equation}
est cocartésien. La proposition s'ensuit. 

(iv) En effet, par définition, le grand rectangle extérieur du diagramme commutatif canonique
\begin{equation}
\xymatrix{
{X'^{(n)}}\ar[d]\ar[r]&X'\ar[d]\\
{\bA_{P'^{(n)}}}\ar[r]\ar[d]&{\bA_{P'}}\ar[d]\\
{\bA_{P^{(n)}}}\ar[r]&{\bA_{P}}}
\end{equation}
est cartésien. Par ailleurs, le carré inférieur est cartésien puisque le diagramme de morphismes de monoïdes \eqref{eccr4c} est cocartésien. 
Par suite, le carré supérieur est aussi cartésien. 
La proposition s'ensuit puisque le rectangle extérieur de \eqref{eccr5c} est cartésien par définition. 

(v) Cela résulte de (iii) et (iv) compte tenu des isomorphismes canoniques 
\begin{eqnarray}
X'^{(n)\rhd}&\stackrel{\sim}{\rightarrow}&X'^{(n)}\times_{\bA_{P'^{(n)}}}\bA_{P'^{(n)\gp}},\\
X'^{[n]\rhd}&\stackrel{\sim}{\rightarrow}&X'^{[n]}\times_{\bA_{P'^{[n]}}}\bA_{P'^{[n]\gp}}.
\end{eqnarray}

(vi) En effet, le morphisme $X'^{[n]\rhd}\rightarrow X'^{(n)\rhd}$ est fini et étale d'après (ii) et (v).  
Il se factorise donc à travers un morphisme fini et surjectif $X'^{[n]\rhd}\rightarrow X'^{\{n\}\rhd}$, 
où $X'^{\{n\}\rhd}$ est un sous-schéma ouvert et fermé de $X'^{(n)\rhd}$. 
La proposition s'ensuit puisque le morphisme $X'^{[n]}\rightarrow X'^{(n)}$ est fini, 
que le schéma $X'^{(n)}$ est normal et noethérien 
et que l'immersion $X'^{(n)\rhd}\rightarrow X'^{(n)}$ est schématiquement dominante \eqref{eccr4}.

\subsection{}\label{eccr41}
D'après \eqref{eccr4h}, on a des morphismes canoniques de systèmes projectifs de schémas, indexés par l'ensemble 
$\mZ_{\geq 1}$ ordonné par la relation de divisibilité, 
\begin{equation}\label{eccr41a}
X'^{[n]}\rightarrow X'^{\{n\}} \rightarrow X'^{(n)}\rightarrow X^{(n)}\rightarrow S^{(n)},
\end{equation}
où $X'^{\{n\}}$ est le sous-schéma ouvert et fermé $X'^{(n)}$ défini dans \ref{eccr5}(vi). 
On se donne un $S$-morphisme
\begin{equation}\label{eccr41c}
\oS\rightarrow \underset{\underset{n\geq 1}{\longleftarrow}}{\lim}\ S^{(n)}.
\end{equation}
Pour tout entier $n\geq 1$, on pose
\begin{eqnarray}
\oX'^{[n]}= X'^{[n]}\times_{S^{(n)}}\oS,&& \oX'^{[n]\rhd}=\oX'^{[n]}\times_{X'}X'^\rhd,\label{eccr41d1}\\
\oX'^{\{n\}}=X'^{\{n\}}\times_{S^{(n)}}\oS,&&  \oX'^{\{n\}\rhd}=\oX'^{\{n\}}\times_{X'}X'^\rhd,\label{eccr41d4}\\
\oX'^{(n)}= X'^{(n)}\times_{S^{(n)}}\oS,&& \oX'^{(n)\rhd}=\oX'^{(n)}\times_{X'}X'^\rhd,\label{eccr41d2}\\
\oX^{(n)}= X^{(n)}\times_{S^{(n)}}\oS,&& \oX^{(n)\circ}=\oX^{(n)}\times_{X}X^\circ.\label{eccr41d3}\\
\end{eqnarray} 
Les morphismes \eqref{eccr41a} induisent des morphismes de systèmes projectifs de $\oS$-schémas
\begin{equation}\label{eccr41e}
\oX'^{[n]}\rightarrow \oX'^{\{n\}}\rightarrow \oX'^{(n)}\rightarrow \oX^{(n)}.
\end{equation}

\begin{lem}\label{eccr410}
Soit $n$ un entier $\geq 1$.
On note $C$ le conoyau de l'homomorphisme $h^\gp\colon P^\gp\rightarrow P'^\gp$ et $C[n]$ le sous-groupe de $n$-torsion de $C$, qu'on identifie
au noyau de l'homomorphisme canonique $P'^{(n)\gp}\rightarrow P'^{[n]\gp}$ par {\rm \ref{eccr5}(ii)}.  
\begin{itemize}
\item[{\rm (i)}]   Le schéma  $\oX'^{[n]\rhd}$ est un espace principal homogène pour la topologie étale 
au-dessus de $\oX'^{\{n\}\rhd}$, sous le groupe $\Hom(C,\mu_n(\oK))$.
\item[{\rm (ii)}] Le schéma $\oX'^{(n)\rhd}$ est un espace principal homogène pour la topologie étale 
au-dessus de $\oX'^{\rhd}$ sous le groupe $\Hom(P^{\gp},\mu_n(\oK))$. 
\item[{\rm (iii)}] Le choix d'une section ensembliste $s_n$ de l'homomorphisme surjectif canonique 
\begin{equation}\label{eccr410d}
\Hom(P'^{(n)\gp},\oK^\times)\rightarrow \Hom(C[n],\oK^\times)
\end{equation}
permet d'identifier le schéma $\oX'^{(n)}$ (resp. $\oX'^{(n)\rhd}$) avec la somme disjointe d'un nombre fini 
de copies de $\oX'^{\{n\}}$ (resp. $\oX'^{\{n\}\rhd}$), canoniquement indexées par le groupe fini $\Hom(C[n],\oK^\times)$.
\end{itemize}
\end{lem}

(i) Cela résulte de \ref{eccr5}(ii)-(v)-(vi) et (\cite{ogus} I 3.3.6). 

(ii) Cela résulte de \ref{cad7}(v) compte tenu de la définition \eqref{eccr4a}. 

Comme $\oX'^{(n)}$ est la fermeture intégrale de $\oX'$ dans $\oX'^{(n)\rhd}$, on en déduit une action du groupe $\Hom(P^{\gp},\mu_n(\oK))$
sur $\oX'^{(n)}$ par des $\oX'$-automorphismes.

(iii) On notera d'abord que $\oK^\times$ étant divisible, l'homomorphisme \eqref{eccr410d} est surjectif. 
Comme $\oX'^{(n)}$ est la fermeture intégrale de $\oX'$ dans $\oX'^{(n)\rhd}$, il suffit de montrer la proposition pour $\oX'^{(n)\rhd}$. 

Considérons les morphismes canoniques de schémas en groupes 
\begin{equation}\label{eccr410a}
\Spec(K[P'^{[n]\gp}])\stackrel{u}{\rightarrow} \Spec(K[P'^{(n)\gp}/C[n]]) \stackrel{v}{\rightarrow} \Spec(K[P'^{(n)\gp}]) 
\stackrel{w}{\rightarrow} \Spec(K[C[n]]).
\end{equation}
Le schéma $\Spec(K[P'^{(n)\gp}/C[n]])$ s'identifie à l'image de $v\circ u$, qui est fini. 
Il s'identifie par ailleurs au noyau de $w$ d'après (\cite{ogus} I 3.3.6).
Le $K$-schéma en groupes $\Spec(K[C[n]])$ étant fini étale, on voit aussitôt que $\ow$ induit une décomposition de 
$\Spec(\oK[P'^{(n)\gp}])$ en somme disjointe finie de sous-schémas ouverts et fermés, 
dont chacun est isomorphe à $\Spec(\oK[P'^{(n)\gp}/C[n]])$, la somme étant canoniquement indexées par le groupe fini $\Hom(C[n],\oK^\times)$.   
La donnée d'une section ensembliste $s_n$ de \eqref{eccr410d} permet d'identifier chacun de ces ouverts et fermés de $\Spec(\oK[P'^{(n)\gp}])$ à $\Spec(\oK[P'^{(n)\gp}/C[n]])$. 

L'homomorphisme nul $P'^\gp\rightarrow P^\gp/nP^\gp$ et l'homomorphisme canonique $P^{(n)\gp}\rightarrow P^\gp/nP^\gp$ 
induisent un homomorphisme $P'^{(n)\gp}\rightarrow P^\gp/nP^\gp$ \eqref{eccr5e}. L'homomorphisme composé 
\begin{equation}\label{eccr410b}
C[n]\rightarrow P'^{(n)\gp}\rightarrow P^\gp/nP^\gp
\end{equation}
s'identifie à l'homomorphisme défini par le lemme du serpent relativement au rectangle extérieur du diagramme \eqref{eccr5e}. Il est donc 
injectif. Comme $\oK^\times$ est divisible, \eqref{eccr410b} induit un homomorphisme surjectif 
\begin{equation}\label{eccr410c}
G_n=\Hom(P^\gp/nP^\gp,\oK^\times)\rightarrow \Hom(C[n],\oK^\times).
\end{equation}

Les homomorphismes \eqref{eccr410b} induisent des actions de $G_n$ sur $\Spec(\oK[P'^{(n)\gp}])$ et $\Spec(\oK[C[n]])$ 
qui font de $\ow$ \eqref{eccr410a} un homomorphisme $G_n$-équivariant. 
En particulier, $G_n$ permute transitivement les sous-schémas ouverts et fermés de $\Spec(\oK[P'^{(n)\gp}])$ 
définis par $\ow$. 
Fixons une section ensembliste $s_n$ de \eqref{eccr410d}.  
Comme le morphisme $\Spec(\oK[P'^{(n)\gp}])\rightarrow \Spec(\oK[P'^{\gp}])$ est un revêtement étale fini de groupe $G_n$ (\cite{ogus} I 3.3.6), 
on en déduit, compte tenu de \ref{eccr5}(v), que $\oX'^{(n)\rhd}$ est la somme disjointe finie de copies de $\oX'^{\{n\}\rhd}$, 
canoniquement indexés par le groupe fini $\Hom(C[n],\oK^\times)$.

\begin{remas}\label{eccr411}
Conservons les hypothèses et notations de \ref{eccr410}.
\begin{itemize}
\item[(i)] Pour tous entiers $n,m\geq 1$, le diagramme d'homomorphismes canoniques
\begin{equation}
\xymatrix{
{\Hom(P'^{(n)\gp},\oK^\times)}\ar[r]^-(0.5){q_n}&{\Hom(C[n],\oK^\times)}\\
{\Hom(P'^{(nm)\gp},\oK^\times)}\ar[r]^-(0.5){q_{nm}}\ar[u]^{p_{n,m}}&{\Hom(C[nm],\oK^\times)}\ar[u]_{c_{n,m}}}
\end{equation}
est commutatif. 
\item[(ii)] Le sous-groupe de torsion $C_\tor$ de $C$ étant fini, notons $n_0$ son exposant de sorte que $C_\tor=C[n_0]$.  Pour tous entiers $n,m\geq 1$ tels que $n_0$ divise $n$, 
le morphisme $c_{n,m}$ est donc un isomorphisme. 
\item[(iii)] Pour tous entiers $n,m\geq 1$, l'homomorphisme canonique $P'^{(n)\gp}\rightarrow P'^{(nm)\gp}$ est injectif de conoyau fini. L'homomorphisme $p_{n,m}$ est donc surjectif de noyau fini. 
\item[(iv)] D'après (iii) et (\cite{te} III § 7.4 prop.~5), le morphisme canonique
\begin{equation}
\underset{\underset{n_0|n}{\longleftarrow}}{\lim}\ \Hom(P'^{(n)\gp},\oK^\times)\rightarrow \Hom(P'^{(n_0)\gp},\oK^\times)
\end{equation}
est surjectif. Compte tenu de ce qui précède, il existe donc pour tout entier $n\geq 1$ divisible par $n_0$, une section ensembliste $s_n$ de $p_n$ 
telle que pour tous entiers $n,m\geq 1$ avec $n$ divisible par $n_0$, le diagramme 
\begin{equation}
\xymatrix{
{\Hom(P'^{(n)\gp},\oK^\times)}&{\Hom(C[n],\oK^\times)}\ar[l]_-(0.5){s_n}\\
{\Hom(P'^{(nm)\gp},\oK^\times)}\ar[u]^{p_{n,m}}&{\Hom(C[nm],\oK^\times)}\ar[u]_{c_{n,m}}\ar[l]_-(0.5){s_{nm}}}
\end{equation} 
soit commutatif. 
\end{itemize}
\end{remas}
 
\begin{lem}\label{eccr412}
Notons $C$ le conoyau de l'homomorphisme $h^\gp\colon P^\gp\rightarrow P'^\gp$ et $C_\tor$ le sous-groupe de torsion de $C$, et posons 
\begin{eqnarray}
\oX'^{(\infty)}&=&\underset{\underset{n\geq 1}{\longleftarrow}}{\lim}\ \oX'^{(n)},\\
\oX'^{\{\infty\}}&=&\underset{\underset{n\geq 1}{\longleftarrow}}{\lim}\ \oX'^{\{n\}}.
\end{eqnarray}
Alors, le schéma $\oX'^{(\infty)}$ s'identifie à la somme disjointe d'un nombre fini de copies de $\oX'^{\{\infty\}}$, 
canoniquement indexées par le groupe fini $\Hom(C_\tor,\oK^\times)$. 
\end{lem}
Cela résulte de \ref{eccr410}(iii) et \ref{eccr411}(iv).

\subsection{}\label{eccr40}
Soient $\oy'$ un point géométrique de $\oX'^\rhd$ et $\oy$ son image dans $\oX^\circ$ (cf. \ref{mtfla1} et \ref{mtfla2}). 
Les schémas $\oX$ et $\oX'$ étant localement irréductibles d'après \ref{cad7}(iii),  
ils sont les sommes des schémas induits sur leurs composantes irréductibles. 
On note $\oX^\star$ (resp. $\oX'^\star$)
la composante irréductible de $\oX$ (resp. $\oX'$) contenant $\oy$ (resp. $\oy'$). 
De même, $\oX^\circ$ (resp. $\oX'^\rhd$) est la somme des schémas induits sur ses composantes irréductibles
et $\oX^{\star \circ}=\oX^\star\times_{X}X^\circ$ (resp. $\oX'^{\star \rhd}=\oX'^\star\times_{X'}X'^\rhd$) est la composante irréductible de 
$\oX^\circ$ (resp. $\oX'^\rhd$) contenant $\oy$ (resp. $\oy'$).

On désigne par $\Delta$ le groupe profini $\pi_1(\oX^{\star \circ},\oy)$ et par $\oR$ la représentation discrète $\oR^{\oy}_{X}$ de $\Delta$ 
définie dans \eqref{TFA9b}. Soit $(V_i)_{i\in I}$ le revêtement universel normalisé de $\oX^{\star \circ}$ en $\oy$ \eqref{notconv11}. 
On rappelle que l'on a \eqref{TFA9c}  
\begin{equation}\label{eccr40b}
\oR=\underset{\underset{i\in I}{\longrightarrow}}{\lim}\  \ocB(V_i\rightarrow X).
\end{equation}

Considérons l'anneau
\begin{equation}\label{eccr40g}
\oR^!=\underset{\underset{i\in I}{\longrightarrow}}{\lim}\  \ocB'(V_i\times_{\oX^\circ}\oX'^\rhd\rightarrow X'),
\end{equation}
qui est naturellement muni d'une action discrète de $\Delta$ par des automorphismes d'anneaux. 
On a un homomorphisme canonique $\Delta$-équivariant $\oR\rightarrow \oR^!$. 

On désigne par $\Delta'$ le groupe profini $\pi_1(\oX'^{\star \rhd},\oy')$ et par $\oR'$ la représentation discrète $\oR'^{\oy'}_{X'}$ de $\Delta'$ 
définie dans \eqref{TFA9b}. Soit $(W_j)_{j\in J}$ le revêtement universel normalisé de $\oX'^{\star \rhd}$ en $\oy'$. 
On a donc  
\begin{equation}\label{eccr40a}
\oR'=\underset{\underset{j\in J}{\longrightarrow}}{\lim}\  \ocB'(W_j\rightarrow X').
\end{equation}

Pour tout $i\in I$, on a un $\oX^\circ$-morphisme canonique $\oy\rightarrow V_i$.
On en déduit un $\oX'^\rhd$-morphisme $\oy'\rightarrow V_i\times_{\oX^\circ}\oX'^\rhd$.
Le schéma $V_i\times_{\oX^\circ}\oX'^\rhd$ étant localement irréductible,  
il est la somme des schémas induits sur ses composantes irréductibles. 
On note $V'_i$ la composante irréductible de $V_i\times_{\oX^\circ}\oX'^\rhd$ contenant l'image de $\oy'$. 
Les schémas $(V'_i)_{i\in I}$ forment naturellement un système projectif de revêtements étales finis connexes $\oy'$-pointés de $\oX'^{\star \rhd}$.
Pour tout $i\in I$, on note $\Pi_i$ le sous-groupe ouvert de $\Delta'$ correspondant à $V'_i$, autrement dit le noyau de l'action canonique 
de $\Delta'$ sur la fibre de $V'_i$ au-dessus de $\oy'$ \eqref{notconv11b}. On désigne par $\Pi$ le sous-groupe fermé de $\Delta'$ défini par
\begin{equation}\label{eccr40c}
\Pi=\cap_{i\in I}\Pi_i.
\end{equation}
Posant $\Delta^\intern=\Delta'/\Pi$, on a un homomorphisme canonique $\Delta'\rightarrow \Delta$ qui se factorise à travers un homomorphisme 
injectif $\Delta^\intern\rightarrow \Delta$.

Considérons l'anneau 
\begin{equation}\label{eccr40e}
\oR^\intern=\underset{\underset{i\in I}{\longrightarrow}}{\lim}\  \ocB'(V'_i\rightarrow X'),
\end{equation}
qui est naturellement muni d'une action discrète de $\Delta^\intern$ par des automorphismes d'anneaux. 
On a un homomorphisme canonique surjectif et $\Delta^\intern$-équivariant $\oR^!\rightarrow \oR^\intern$. 
 
Pour tous $i\in I$ et $j\in J$, il existe au plus un morphisme de $\oX'^{\star \rhd}$-schémas pointés $W_j\rightarrow V'_i$. 
De plus, pour tout $i\in I$, il existe $j\in J$ et un morphisme de $\oX'^{\star \rhd}$-schémas pointés $W_j\rightarrow V'_i$. 
On a donc un homomorphisme canonique
\begin{equation}\label{eccr40d}
\underset{\underset{i\in I}{\longrightarrow}}{\lim}\  \ocB'(V'_i\rightarrow X')\rightarrow 
\underset{\underset{j\in J}{\longrightarrow}}{\lim}\  \ocB'(W_j\rightarrow X').
\end{equation}
On en déduit un homomorphisme canonique $\Delta'$-équivariant de $R'_1$-algèbres
\begin{equation}\label{eccr40f}
\oR^\intern\rightarrow \oR'.
\end{equation}

On désigne par $\hoR$, $\hoR^!$, $\hoR^\intern$ et $\hoR'$ les séparés complétés $p$-adiques de $\oR$, $\oR^!$,
$\oR^\intern$ et $\oR'$, respectivement.

\subsection{}\label{eccr42}
Pour tous entiers $m,n\geq 1$, le morphisme canonique $X'^{[mn]}\rightarrow X'^{[n]}$ est fini et surjectif \eqref{eccr3a}.
D'après (\cite{ega4} 8.3.8(i)), le point géométrique $\oy'\rightarrow X'$ fixé dans \ref{eccr40}, se relève donc en un $X'$-morphisme 
\begin{equation}\label{eccr41b}
\oy'\rightarrow \underset{\underset{n\geq 1}{\longleftarrow}}{\lim}\ X'^{[n]},
\end{equation} 
que l'on fixe dans la suite de cette section.
Quitte à changer \eqref{eccr41c}, on peut supposer qu'il est induit par \eqref{eccr41b}. 
On déduit de \eqref{eccr41b}, \eqref{eccr41a} et \eqref{eccr41c} un $\oX'$-morphisme 
\begin{equation}\label{eccr41f}
\oy'\rightarrow \underset{\underset{n\geq 1}{\longleftarrow}}{\lim}\ \oX'^{[n]}.
\end{equation}

Pour tout entier $n\geq 1$, le schéma $\oX'^{[n]}$ est normal et localement irréductible d'après \ref{cad7}(iii). 
Il est donc la somme des schémas induits sur ses composantes irréductibles. 
On note $\oX'^{[n]\star}$ la composante irréductible de $\oX'^{[n]}$ contenant l'image de $\oy'$ \eqref{eccr41f}.
De même, $\oX'^{[n]\rhd}$ est la somme des schémas induits sur ses composantes irréductibles
et $\oX'^{[n]\star\rhd}=\oX'^{[n]\star}\times_{X'}X'^\rhd$ est la composante irréductible de $\oX'^{[n]\rhd}$ contenant l'image de $\oy'$.
On pose 
\begin{equation}\label{eccr42a}
R'_n=\Gamma(\oX'^{[n]\star},\co_{\oX'^{[n]}}).
\end{equation}
Les anneaux $(R'_n)_{n\geq 1}$ forment naturellement un système inductif. On pose 
\begin{eqnarray}
R'_\infty&=&\underset{\underset{n\geq 1}{\longrightarrow}}{\lim}\ R'_n,\label{eccr42b1}\\
R'_{p^\infty}&=&\underset{\underset{n\geq 0}{\longrightarrow}}{\lim}\ R'_{p^n}.\label{eccr42b2}
\end{eqnarray}
Le morphisme \eqref{eccr41f} induit des homomorphismes injectifs de $R'_1$-algèbres 
\begin{equation}\label{eccr42c}
R'_{p^\infty}\rightarrow R'_\infty\rightarrow \oR'.
\end{equation}

De même, pour tout entier $n\geq 1$, le schéma $\oX^{(n)}$ est la somme des schémas induits sur ses composantes irréductibles. 
On note $\oX^{(n)\star}$ la composante irréductible de $\oX^{(n)}$ contenant l'image de $\oy'$ \eqref{eccr41f}.
Le schéma $\oX^{(n)\circ}$ est la somme des schémas induits sur ses composantes irréductibles
et $\oX^{(n)\star\circ}=\oX^{(n)\star}\times_XX^\circ$ est la composante irréductible de $\oX^{(n)\circ}$ contenant l'image de $\oy'$.
On pose 
\begin{equation}\label{eccr42d}
R_n=\Gamma(\oX^{(n)\star},\co_{\oX^{(n)}}).
\end{equation}
Les anneaux $(R_n)_{n\geq 1}$ forment naturellement un système inductif. On pose 
\begin{eqnarray}
R_\infty&=&\underset{\underset{n\geq 1}{\longrightarrow}}{\lim}\ R_n,\label{eccr42e1}\\
R_{p^\infty}&=&\underset{\underset{n\geq 0}{\longrightarrow}}{\lim}\ R_{p^n}.\label{eccr42e2}
\end{eqnarray}
Les morphismes \eqref{eccr41e} et \eqref{eccr41f} induisent des homomorphismes injectifs de $R_1$-algèbres 
\begin{equation}\label{eccr42f}
R_{p^\infty}\rightarrow R_\infty\rightarrow \oR.
\end{equation}

Pour tout entier $n\geq 1$, le schéma $\oX'^{(n)}$ est normal et localement irréductible d'après \eqref{eccr4a} et (\cite{agt} III.4.2(iii)). 
Il est donc la somme des schémas induits sur ses composantes irréductibles. 
Par ailleurs, l'immersion $\oX'^{(n)\rhd}\rightarrow \oX'^{(n)}$ est schématiquement dominante (\cite{ega4} 11.10.5). 
On note $\oX'^{(n)\star}$ la composante irréductible de $\oX'^{(n)}$ contenant l'image de $\oy'$ \eqref{eccr41f}.
Le schéma $\oX'^{(n)\rhd}$ est aussi la somme des schémas induits sur ses composantes irréductibles
et $\oX'^{(n)\star\rhd}=\oX'^{(n)\star}\times_{X'}X'^\rhd$  est la composante irréductible de $\oX'^{(n)\rhd}$ 
contenant l'image de $\oy'$. On pose 
\begin{equation}\label{eccr42g}
R^\intern_n=\Gamma(\oX'^{(n)\star},\co_{\oX'^{(n)}}).
\end{equation}
Les anneaux $(R^\intern_n)_{n\geq 1}$ forment naturellement un système inductif. On pose 
\begin{eqnarray}
R^\intern_\infty&=&\underset{\underset{n\geq 1}{\longrightarrow}}{\lim}\ R^\intern_n,\label{eccr42h1}\\
R^\intern_{p^\infty}&=&\underset{\underset{n\geq 0}{\longrightarrow}}{\lim}\ R^\intern_{p^n}.\label{eccr42h2}
\end{eqnarray}
Ce sont des anneaux normaux et intègres d'après (\cite{ega1n} 0.6.1.6(i) et 0.6.5.12(ii)). 
Compte tenu de la définition \eqref{eccr40e}, les morphismes \eqref{eccr41e} et \eqref{eccr41f} 
induisent des homomorphismes injectifs de $R^\intern_1$-algèbres 
\begin{equation}\label{eccr42i}
R^\intern_{p^\infty}\rightarrow R^\intern_\infty\rightarrow \oR^\intern.
\end{equation}

Les morphismes de \eqref{eccr41e} induisent des homomorphismes de systèmes inductifs d'algèbres 
\begin{equation}\label{eccr42j}
(R_n)_{n\geq 1}\rightarrow (R^\intern_n)_{n\geq 1}\rightarrow (R'_n)_{n\geq 1}.
\end{equation} 
L'homomorphisme $R^\intern_1\rightarrow R'_1$ est un isomorphisme.

\subsection{}\label{eccr43}
Soit $n$ un entier $\geq 1$. Il résulte de  \ref{cad7}(v) que le morphisme $\oX'^{[n]\star\rhd}\rightarrow \oX'^{\star\rhd}$ 
est un revêtement étale fini et galoisien de groupe $\Delta'_n$, 
un sous-groupe de $\Hom_\mZ(P'^\gp/h^\gp(\vartheta^\gp(\mZ)),\mu_{n}(\oK))$ 
(cf. la preuve de \cite{agt} II.6.8(iv)).  
Le groupe $\Delta'_n$ agit naturellement sur $R'_n$. Si $n$ est une puissance de $p$, le morphisme canonique
\begin{equation}\label{eccr43a}
\oX'^{[n]\star}\rightarrow\oX'^{[n]}\times_{\oX'}\oX'^\star
\end{equation}
est un isomorphisme en vertu de (\cite{agt} II.6.6(v)), et on a donc $\Delta'_n\simeq \Hom_\mZ(P'^\gp/h^\gp(\vartheta^\gp(\mZ)),\mu_{n}(\oK))$.

Les groupes $(\Delta'_n)_{n\geq 1}$ forment naturellement un système projectif. 
On pose 
\begin{eqnarray}
\Delta'_\infty&=&\underset{\underset{n\geq 1}{\longleftarrow}}{\lim}\ \Delta'_n,\label{eccr43b1}\\
\Delta'_{p^\infty}&=&\underset{\underset{n\geq 0}{\longleftarrow}}{\lim}\ \Delta'_{p^n}.\label{eccr43b2}
\end{eqnarray}
On identifie $\Delta'_\infty$ (resp. $\Delta'_{p^\infty}$) au groupe de Galois de l'extension des corps de fractions de $R'_\infty$ sur $R'_1$
(resp. $R'_{p^\infty}$ sur $R'_1$). On a des homomorphismes canoniques \eqref{mtfla1a}
\begin{equation}\label{eccr43b3}
\xymatrix{
{\Delta'_{\infty}}\ar@{->>}[d]\ar@{^(->}[r]&{\Hom(P'^\gp/h^\gp(\vartheta^\gp(\mZ)),\hmZ(1))}\ar[d]\\
{\Delta'_{p^\infty}}\ar[r]^-(0.5)\sim&{\Hom(P'^\gp/h^\gp(\vartheta^\gp(\mZ)),\mZ_p(1))}}
\end{equation}
Le noyau $\Sigma'_0$ de l'homomorphisme canonique $\Delta'_\infty\rightarrow \Delta'_{p^\infty}$
est un groupe profini d'ordre premier à $p$. 
Par ailleurs, le morphisme \eqref{eccr41f} détermine un homomorphisme surjectif $\Delta'\rightarrow \Delta'_\infty$ \eqref{eccr40}. 
On note $\Sigma'$ son noyau. Les homomorphismes \eqref{eccr42c} sont $\Delta'$-équivariants.

De même, le morphisme $\oX^{(n)\star\circ}\rightarrow \oX^{\star\circ}$ 
est un revêtement étale fini et galoisien de groupe $\Delta_n$, 
un sous-groupe de $\Hom_\mZ(P^\gp/\vartheta^\gp(\mZ),\mu_{n}(\oK))$.  
Le groupe $\Delta_n$ agit naturellement sur $R_n$. Si $n$ est une puissance de $p$, le morphisme canonique
\begin{equation}\label{eccr43c}
\oX^{(n)\star}\rightarrow\oX^{(n)}\times_{\oX}\oX^\star
\end{equation}
est un isomorphisme, et on a donc $\Delta_n\simeq \Hom_\mZ(P^\gp/\vartheta^\gp(\mZ),\mu_{n}(\oK))$.

Les groupes $(\Delta_n)_{n\geq 1}$ forment naturellement un système projectif. 
On pose 
\begin{eqnarray}
\Delta_\infty&=&\underset{\underset{n\geq 1}{\longleftarrow}}{\lim}\ \Delta_n,\label{eccr43d1}\\
\Delta_{p^\infty}&=&\underset{\underset{n\geq 0}{\longleftarrow}}{\lim}\ \Delta_{p^n}.\label{eccr43d2}
\end{eqnarray}
On identifie $\Delta_\infty$ (resp. $\Delta_{p^\infty}$) au groupe de Galois de l'extension des corps de fractions de $R_\infty$ sur $R_1$
(resp. $R_{p^\infty}$ sur $R_1$). On a des homomorphismes canoniques \eqref{mtfla1a}
\begin{equation}\label{eccr43d3}
\xymatrix{
{\Delta_{\infty}}\ar@{->>}[d]\ar@{^(->}[r]&{\Hom(P^\gp/\vartheta^\gp(\mZ),\hmZ(1))}\ar[d]\\
{\Delta_{p^\infty}}\ar[r]^-(0.5)\sim&{\Hom(P^\gp/\vartheta^\gp(\mZ),\mZ_p(1))}}
\end{equation}
Le noyau $\Sigma_0$ de l'homomorphisme canonique $\Delta_\infty\rightarrow \Delta_{p^\infty}$
est un groupe profini d'ordre premier à $p$. 
Par ailleurs, les morphismes \eqref{eccr41e} et \eqref{eccr41f} déterminent un homomorphisme surjectif $\Delta\rightarrow \Delta_\infty$ \eqref{eccr40}.
On note $\Sigma$ son noyau. Les homomorphismes \eqref{eccr42f} sont $\Delta$-équivariants.

Il résulte de \ref{eccr410}(i) que le morphisme $\oX'^{[n]\star\rhd}\rightarrow \oX'^{(n)\star\rhd}$ est un revêtement étale, 
fini et galoisien de groupe $\fS_n$, un sous-groupe de $\Hom_\mZ(P'^{\gp}/h^\gp(P^{\gp}),\mu_n(\oK))$. 
Il résulte de \ref{tpcg4} et \eqref{eccr4a} que le morphisme $\oX'^{(n)\star\rhd}\rightarrow \oX'^{\star\rhd}$ 
est un revêtement étale, fini et galoisien de groupe $\Delta^\intern_n$, un sous-groupe de $\Delta_n$. 
Le groupe $\Delta^\intern_n$ agit naturellement sur $R^\intern_n$. On a une suite exacte canonique
\begin{equation}\label{eccr43e}
0\rightarrow \fS_n\rightarrow \Delta'_n\rightarrow \Delta^\intern_n \rightarrow 0.
\end{equation}
Par ailleurs, le diagramme 
\begin{equation}\label{eccr43f}
\xymatrix{
\fS_n\ar[r]\ar[d]&{\Hom_\mZ(P'^{\gp}/h^\gp(P^{\gp}),\mu_n(\oK))}\ar[d]\\
{\Delta'_n}\ar[r]&{\Hom_\mZ(P'^\gp/h^\gp(\vartheta^\gp(\mZ)),\mu_{n}(\oK))}}
\end{equation}
où la flèche verticale de droite est l'homomorphisme canonique, est commutatif.
En vertu de \ref{eccr8}(iv) ci-dessous, si $n$ est une puissance de $p$, le morphisme canonique
\begin{equation}\label{eccr43g}
\oX'^{(n)\star}\rightarrow \oX^{(n)}\times_{\oX}\oX'^\star
\end{equation}
est un isomorphisme, et on a 
\begin{eqnarray}
\fS_{n}&=&\Hom(P'^{\gp}/h^\gp(P^{\gp}),\mu_{n}(\oK)),\label{eccr43h1}\\
\Delta^\intern_{n}&=&\Hom(P^{\gp}/\vartheta^\gp(\mZ),\mu_{n}(\oK)).\label{eccr43h2}
\end{eqnarray}

Les groupes $(\fS_n)_{n\geq 1}$ et $(\Delta^\intern_n)_{n\geq 1}$ forment naturellement des systèmes projectifs. On pose 
\begin{eqnarray}
\fS_\infty&=&\underset{\underset{n\geq 1}{\longleftarrow}}{\lim}\ \fS_n,\label{eccr43i1}\\
\fS_{p^\infty}&=&\underset{\underset{n\geq 0}{\longleftarrow}}{\lim}\ \fS_{p^n}, \label{eccr43i2}\\
\Delta^\intern_\infty&=&\underset{\underset{n\geq 1}{\longleftarrow}}{\lim}\ \Delta^\intern_n,\label{eccr43i3}\\
\Delta^\intern_{p^\infty}&=&\underset{\underset{n\geq 0}{\longleftarrow}}{\lim}\ \Delta^\intern_{p^n}.\label{eccr43i4}
\end{eqnarray}
On identifie $\Delta^\intern_\infty$ (resp. $\Delta^\intern_{p^\infty}$) 
au groupe de Galois de l'extension des corps des fractions de $R^\intern_\infty$ (resp. $R^\intern_{p^\infty}$) sur $R'_1$. 
D'après la théorie de Galois, les suites 
\begin{equation}\label{eccr43j}
0\rightarrow \fS_\infty\rightarrow \Delta'_\infty\rightarrow \Delta^\intern_\infty\rightarrow 0,
\end{equation}
\begin{equation}\label{eccr43k}
0\rightarrow \fS_{p^\infty}\rightarrow \Delta'_{p^\infty}\rightarrow \Delta^\intern_{p^\infty}\rightarrow 0,
\end{equation}
déduites de \eqref{eccr43e} sont exactes. Le groupe $\fS_\infty$ (resp. $\fS_{p^\infty}$) s'identifie donc au groupe de Galois
de l'extension des corps de fractions de $R'_\infty$ sur $R^\intern_\infty$ (resp. $R'_{p^\infty}$ sur $R^\intern_{p^\infty}$).  
Comme le noyau de l'homomorphisme canonique $\Delta^\intern_\infty\rightarrow \Delta^\intern_{p^\infty}$ 
est un groupe profini d'ordre premier à $p$, 
on en déduit que l'homomorphisme canonique $\fS_\infty\rightarrow \fS_{p^\infty}$ est surjectif. 
On désigne par $\fN$ son noyau, qui est un groupe profini d'ordre premier à $p$, de sorte qu'on a la suite exacte
\begin{equation}\label{eccr43l}
0\rightarrow \fN\rightarrow \fS_\infty \rightarrow \fS_{p^\infty}\rightarrow 0.
\end{equation}

On désigne par $\fS$ le noyau de l'homomorphisme canonique $\Delta'\rightarrow \Delta^\intern_\infty$. 
On a alors un isomorphisme canonique 
\begin{equation}\label{eccr43m}
\fS\stackrel{\sim}{\rightarrow}\underset{\underset{n\geq 1}{\longleftarrow}}{\lim} \ \pi_1(\oX'^{(n)\star\rhd},\oy').
\end{equation}

\begin{prop}\label{eccr8}
\
\begin{itemize}
\item[{\rm (i)}] Soit $n$ un entier $\geq 1$. Remplaçant $g$ \eqref{mtfla2} et sa carte relativement adéquate \eqref{eccr1a}
par $g^{(n)}$ \eqref{eccr4g} et sa carte relativement adéquate \eqref{eccr5a}, ce qui est loisible d'après \ref{eccr5}{\rm (i)}, 
la tour de schémas $(X^{(m)})_{m\geq 1}$ (resp. $(X'^{[m]})_{m\geq 1}$, resp. $(X'^{(m)})_{m\geq 1}$) 
est remplacée par $(X^{(mn)})_{m\geq 1}$ (resp. $(X'^{[mn]})_{m\geq 1}$, resp. $(X'^{(mn)}\times_{X'^{(n)}}X'^{[n]})_{m\geq 1}$). 
\item[{\rm (ii)}] Pour tous entiers $m,n\geq 1$, le schéma 
$\oX'^{(mn)}\times_{\oX'^{(n)}}\oX'^{[n]}$ est normal, localement irréductible et $\oS$-plat; 
l'immersion $\oX'^{(mn)\rhd}\times_{\oX'^{(n)\rhd}}\oX'^{[n]\rhd}\rightarrow \oX'^{(mn)}\times_{\oX'^{(n)}}\oX'^{[n]}$
est schématiquement dominante.
\item[{\rm (iii)}] Le nombre de composantes irréductibles de $\oX'^{[n]}$ (resp. $\oX'^{(n)}$) lorsque $n$ décrit les entiers $\geq 1$, est borné. 
\item[{\rm (iv)}] Pour tout entier $n\geq 0$, les morphismes canoniques 
\begin{eqnarray}
\oX'^{[p^n]\star}&\rightarrow& \oX'^{[p^n]}\times_{\oX'}\oX'^\star, \label{eccr8b}\\
\oX'^{(p^n)\star}&\rightarrow& \oX^{(p^n)}\times_{\oX}\oX'^\star, \label{eccr8c}
\end{eqnarray}
sont des isomorphismes, et on a 
\begin{eqnarray}
\fS_{p^n}&=&\Hom(P'^{\gp}/h^\gp(P^{\gp}),\mu_{p^n}(\oK)),\label{eccr8d}\\
\Delta^\intern_{p^n}&=&\Delta_{p^n}=\Hom(P^{\gp}/\vartheta^\gp(\mZ),\mu_{p^n}(\oK)).\label{eccr8h}
\end{eqnarray}
\item[{\rm (v)}] Pour tous entiers $n\geq 0$ et $m\geq 1$,  le schéma 
$\oX'^{(mp^n)\star}\times_{\oX'^{(p^n)\star}}\oX'^{[p^n]\star}$ est intègre et normal. 
\end{itemize}
\end{prop}

(i) C'est immédiat. 

(ii) Cela résulte de (i) compte tenu de  \ref{eccr42} et de l'isomorphisme
\begin{equation}
\oX'^{(mn)}\times_{\oX'^{(n)}}\oX'^{[n]}\stackrel{\sim}{\rightarrow}(X'^{(mn)}\times_{X'^{(n)}}X'^{[n]})\times_{S^{(mn)}}\oS.
\end{equation}

(iii) En effet, l'assertion pour les $\oX'^{[n]}$ résulte de \ref{cad8} appliqué au morphisme $f'$. 
L'assertion pour les $\oX'^{(n)}$ résulte de celle pour les $\oX'^{[n]}$ et de \ref{eccr410}(iii), 
puisque le sous-groupe de torsion de $P'^\gp/h^\gp(P^\gp)$ est fini. 

(iv) L'isomorphisme \eqref{eccr8b} résulte de (\cite{agt} II.6.6(v)) appliqué au morphisme $f'$. 
Comme $P'^\gp/h^\gp(P^\gp)$ n'a pas de torsion $p$-primaire non nulle, on a $X'^{\{p^n\}}=X'^{(p^n)}$ d'après \ref{eccr410}(iii).
Le morphisme canonique $\oX'^{[p^n]}\times_{\oX'}\oX'^\star\rightarrow \oX^{(p^n)}\times_{\oX}\oX'^\star$ est donc fini et surjectif. 
On en déduit que le schéma $\oX^{(p^n)}\times_{\oX}\oX'^\star$ est connexe, d'où l'isomorphisme \eqref{eccr8c}. 

L'équation \eqref{eccr8d} résulte de \eqref{eccr8b}, \eqref{eccr8c} et \ref{eccr410}(i).
L'équation \eqref{eccr8h} résulte de \eqref{eccr8c} et \ref{eccr410}(ii) (cf. la preuve de \cite{agt} II.6.8(iv)).

(v) En effet, si $m$ est une puissance de $p$, la proposition résulte de \eqref{eccr8c} appliqué au morphisme $g^{(p^n)}$, compte tenu de (i).
Dans le cas général, posant $mp^n= m'p^{n'}$, où $n'\geq n$ et $m'$ est premier à $p$, on a 
\begin{equation}
\oX'^{(mp^n)\star}\times_{\oX'^{(p^n)\star}}\oX'^{[p^n]\star}\stackrel{\sim}{\rightarrow}
\oX'^{(m'p^{n'})\star}\times_{\oX'^{(p^{n'})\star}}\oX'^{(p^{n'})\star}\times_{\oX'^{(p^n)\star}}\oX'^{[p^n]\star}.
\end{equation}
D'après ce qui précède, $\oX'^{(p^{n'})\star}\times_{\oX'^{(p^n)\star}}\oX'^{[p^n]\star}$ est intègre et normal. 
On rappelle que les morphismes $\oX'^{[p^n]\star\rhd}\rightarrow \oX'^{(p^n)\star\rhd}$ et 
$\oX'^{(m'p^{n'})\star\rhd}\rightarrow \oX'^{(p^{n'})\star\rhd}$ sont finis, étales et galoisiens. 
Notons $G$ (resp. $F$, resp. $H$) le corps de fractions de 
$\oX'^{(p^{n'})\star}\times_{\oX'^{(p^n)\star}}\oX'^{[p^n]\star}$ (resp. $\oX'^{(p^{n'})\star}$, resp.  $\oX'^{(m'p^{n'})\star}$). 
L'extension $G/F$ est galoisienne de degré une puissance de $p$ et l'extension 
$H/F$ est galoisienne de degré premier à $p$. Par suite, $H\otimes_FG$ est un corps, d'où la proposition compte tenu de (ii).

\begin{cor}\label{eccr9}
\
\begin{itemize}
\item[{\rm (i)}] Pour tout entier $n\geq 1$, l'homomorphisme canonique 
\begin{equation}\label{eccr9a}
R_n\otimes_RR'\rightarrow R^\intern_n
\end{equation}
est surjectif, de noyau engendré par un idempotent de $R_n\otimes_RR'$. 
\item[{\rm (ii)}] Le schéma $\Spec(R_{\infty}\otimes_{R}R')$ est normal et localement irréductible,  
et $\Spec(R^\intern_{\infty})$ est une composante irréductible de $\Spec(R_{\infty}\otimes_{R}R')$. 
\item[{\rm (iii)}] L'homomorphisme canonique 
\begin{equation}\label{eccr9b}
R_{p^\infty}\otimes_{R_1}R'_1\rightarrow R^\intern_{p^\infty}
\end{equation}
est un isomorphisme, et on a 
\begin{equation}\label{eccr9c}
\fS_{p^\infty}=\Hom(P'^{\gp}/h^\gp(P^{\gp}),\mZ_p(1)).
\end{equation} 
\item[{\rm (iv)}] Pour tout entier $n\geq 1$, le schéma $\Spec(R^\intern_{\infty}\otimes_{R^\intern_n}R'_n)$ 
est normal et localement irréductible.
\item[{\rm (v)}] Pour tout entier $n\geq 0$, les anneaux $R^\intern_{p^\infty}\otimes_{R^\intern_{p^n}}R'_{p^n}$ et
$R^\intern_{\infty}\otimes_{R^\intern_{p^n}}R'_{p^n}$ sont normaux et intègres.
\item[{\rm (vi)}] L'anneau $R'_{p^\infty}\otimes_{R^\intern_{p^\infty}}R^\intern_\infty$ est normal et intègre, 
et on a un isomorphisme canonique
\eqref{eccr43l}
\begin{equation}\label{eccr9d}
R^\intern_\infty\otimes_{R^\intern_{p^\infty}}R'_{p^\infty}\stackrel{\sim}{\rightarrow}(R'_\infty)^\fN. 
\end{equation}
\item[{\rm (vii)}] Pour tout $a\in \co_\oK$, l'homomorphisme canonique
\begin{equation}\label{eccr9e}
R'_{p^\infty}\otimes_{R^\intern_{p^\infty}}(R^\intern_\infty/aR^\intern_\infty)\stackrel{\sim}{\rightarrow}(R'_\infty/aR'_\infty)^\fN 
\end{equation}
est un isomorphisme. 
\end{itemize}
\end{cor}

(i) En effet, la première assertion résulte de l'isomorphisme canonique \eqref{eccr4a}
\begin{equation}
\oX'^{(n)}\stackrel{\sim}{\rightarrow} \oX^{(n)}\times_XX',
\end{equation}
et du fait que $\oX'^{(n)}$ est la somme des schémas induits sur ses composantes irréductibles. 

(ii) Pour tout entier $n\geq 1$, le schéma $\oX'^{(n)}$ est normal et localement irréductible et l'immersion $\oX'^{(n)\rhd}\rightarrow \oX'^{(n)}$
est schématiquement dominante. Pour tous entiers $m,n\geq 1$, le morphisme $\oX'^{(mn)\rhd}\rightarrow \oX'^{(n)\rhd}$ est étale fini et surjectif. 
Par suite, toute composante irréductible $\oX'^{(mn)}$ domine une composante irréductible de $\oX'^{(n)}$. 
En vertu de \ref{eccr8}(iii), il existe un entier $n\geq 1$ tel que pour toute composante 
irréductible $C$ de $\oX'^{(n)}$ et tout entier $m\geq 1$, $C\times_{\oX'^{(n)}}\oX'^{(mn)}$ soit irréductible. 
Le schéma $\Spec(R_{\infty}\otimes_{R}R')$ est donc normal et localement irréductible en vertu de (\cite{ega1n} 0.6.1.6 et 0.6.5.12(ii)).
La seconde assertion est une conséquence de (i). 

(iii) L'isomorphisme \eqref{eccr9b} résulte de \eqref{eccr8c}, et l'équation \eqref{eccr9c} résulte aussitôt de \eqref{eccr8d}.

(iv) En effet, pour tout entier $m\geq 1$, le schéma 
$\oX'^{(mn)\star}\times_{\oX'^{(n)\star}}\oX'^{[n]\star}$ est normal et localement irréductible, et
l'immersion $\oX'^{(mn)\star\rhd}\times_{\oX'^{(n)\star\rhd}}\oX'^{[n]\star\rhd}\rightarrow \oX'^{(mn)\star}\times_{\oX'^{(n)\star}}\oX'^{[n]\star}$
est schématiquement dominante d'après \ref{eccr8}(ii). Par ailleurs, pour tous entiers $m',m\geq 1$,  le morphisme canonique
\begin{equation}
\varpi^{m'm}_{m}\colon \oX'^{(m'mn)\star}\times_{\oX'^{(n)\star}}\oX'^{[n]\star}\rightarrow \oX'^{(mn)\star}\times_{\oX'^{(n)\star}}\oX'^{[n]\star}
\end{equation}
est étale, fini et surjectif au-dessus de $X'^{[n]\star \rhd}$, et est donc dominant. 
De plus, il existe un entier $m\geq 1$ tel que pour  tout entier $m'\geq 1$ et 
toute composante irréductible $C$ de $\oX'^{(mn)\star}\times_{\oX'^{(n)\star}}\oX'^{[n]\star}$, $(\varpi^{m'm}_{m})^{-1}(C)$ est irréductible
d'après \ref{eccr8}(i)-(iii).  
La proposition s'ensuit par passage à la limite inductive (\cite{ega1n} 0.6.1.6(i) et 0.6.5.12(ii)).

(v) En effet, pour tout entier $m\geq 1$, le schéma $\oX'^{(mp^{n})\star}\times_{\oX'^{(p^{n})\star}}\oX'^{[p^{n}]\star}$ est normal et intègre en vertu
de \ref{eccr8}(v). Par ailleurs, pour tous entiers $m,m'\geq 1$,  le morphisme 
\begin{equation}
\oX'^{(m'mp^{n})\star}\times_{\oX'^{(p^{n})\star}}\oX'^{[p^{n}]\star}\rightarrow \oX'^{(mp^{n})\star}\times_{\oX'^{(p^{n})\star}}\oX'^{[p^{n}]\star}
\end{equation}
est étale et fini au-dessus de $X'^{\star \rhd}$, et est donc dominant. 
La proposition s'ensuit par passage à la limite inductive compte tenu de (\cite{ega1n} 0.6.1.6(i) et 0.6.5.12(ii)).

(vi) En effet, pour tous entiers $n\geq 0$ et $m\geq 1$, 
le schéma $\oX'^{(mp^n)\star}\times_{\oX'^{(p^n)\star}}\oX'^{[p^n]\star}$ est intègre et normal, d'après  \ref{eccr8}(v).
Par ailleurs, pour tous entiers $n'\geq n\geq 0$ et $m,m'\geq 1$, le morphisme canonique
\begin{equation}
\oX'^{(m'mp^{n'})\star}\times_{\oX'^{(p^{n'})\star}} \oX'^{[p^{n'}]\star}\rightarrow  \oX'^{(mp^n)\star}\times_{\oX'^{(p^n)\star}}\oX'^{[p^n]\star}
\end{equation}
est étale et fini au-dessus de $X'^{\star\rhd}$, et est donc dominant.  
On obtient la première assertion par passage à la limite inductive (\cite{ega1n} 0.6.1.6(i) et 0.6.5.12(ii)).
La seconde assertion s'ensuit compte tenu du diagramme d'extensions entières d'anneaux intègres et normaux 
\begin{equation}
\xymatrix{
{R^\intern_\infty}\ar[r]^-(0.5){\fS_{p^\infty}}\ar@/_1pc/[rr]_{\fS_\infty}&{R^\intern_\infty\otimes_{R^\intern_{p^\infty}}R'_{p^\infty}}
\ar[r]^-(0.5)\fN&{R'_\infty}}
\end{equation}
où au-dessus de chaque flèche, on a noté le groupe de Galois de l'extension associée des corps de fractions \eqref{eccr43l}. 

(vii) Comme le groupe profini $\fN$ est la limite inverse de groupes finis d'ordres premiers à $p$, la suite 
\begin{equation}
0\rightarrow (R'_\infty)^\fN\stackrel{a}{\rightarrow} (R'_\infty)^\fN \rightarrow (R'_\infty/aR'_\infty)^\fN\rightarrow 0
\end{equation} 
est exacte. 
L'isomorphisme \eqref{eccr9e} se déduit alors de l'isomorphisme \eqref{eccr9d} par réduction modulo $a$.

\subsection{}\label{eccr34}
On peut résumer les principales constructions de cette section dans le diagramme commutatif suivant d'extensions d'anneaux intègres et normaux
\begin{equation}\label{eccr34a}
\xymatrix{
\oR\ar[r]&\oR^\intern\ar[rr]^-(0.5)\Pi&&\oR'\\
R_\infty\ar[r]\ar[u]^{\Sigma}&{R^\intern_\infty}\ar[r]_-(0.5){\fS_{p^\infty}}\ar[u]^{\Sigma^\intern}
\ar@/^1pc/[rr]^-(0.5){\fS_\infty}\ar@/^1pc/[rru]^-(0.3){\fS}&
{R^\intern_\infty\otimes_{R^\intern_{p^\infty}}R'_{p^\infty}}\ar[r]_-(0.4){\fN}&{R'_\infty}\ar[u]_{\Sigma'}\\
R_{p^\infty}\ar[r]\ar[u]^{\Sigma_0}&R^\intern_{p^\infty}\ar[r]^-(0.5){\fS_{p^\infty}}\ar[u]&R'_{p^\infty}\ar[ru]_{\Sigma'_0}\ar[u]&\\
R_1\ar[u]^{\Delta_{p^\infty}}\ar[r]&R^\intern_1\ar@{=}[r]\ar[u]^{\Delta^\intern_{p^\infty}}&R'_1\ar[u]_{\Delta'_{p^\infty}}&}
\end{equation}
où pour certaines extensions entières, on a noté le groupe de Galois de l'extension associée des corps de fractions. 
On observera que l'homomorphisme $R^\intern_1\rightarrow \oR^\intern$ induit une extension galoisienne des corps de fractions, 
et il en est donc de même de l'homomorphisme $R^\intern_\infty\rightarrow \oR^\intern$. 
On note $\oF^\intern$ (resp. $F^\intern_\infty$) le corps de fractions
de $\oR^\intern$ (resp. $R^\intern_\infty$) et $\Sigma^\intern$ le groupe de Galois de l'extension $\oF^\intern/F^\intern_\infty$. 

\begin{prop}\label{eccr35}
L'extension $\oF^\intern$ de $F^\intern_\infty$ est la réunion d'un système inductif filtrant de 
sous-extensions galoisiennes finies $E$ de $F^\intern_\infty$ telles que la clôture intégrale
de $R^\intern_\infty$ dans $E$ soit $\alpha$-finie-étale sur $R^\intern_\infty$ \eqref{aet2}. 
\end{prop}

Notons $\oF$ (resp. $F_\infty$) le corps de fractions de $\oR$ (resp. $R_\infty$). 
D'après \ref{tpcg5}, l'extension $\oF$ de $F_\infty$ est la réunion d'un système inductif filtrant de 
sous-extensions galoisiennes finies $N$ de $F_\infty$ telles que la clôture intégrale
de $R_\infty$ dans $N$ soit $\alpha$-finie-étale sur $R_\infty$. Soient $N$ une telle extension de $F_\infty$, 
$\cN$ la clôture intégrale de $R_\infty$ dans $N$, $E$ l'image de l'homomorphisme canonique 
$N\otimes_{F_\infty}F^\intern_\infty\rightarrow \oF^\intern$, $\cE$ la clôture intégrale de $R^\intern_\infty$ dans $E$. 
Comme la $R_\infty$-algèbre $\cN$ est $\alpha$-finie-étale, la $R^\intern_\infty$-algèbre $\cN\otimes_{R_\infty}R^\intern_\infty$ est 
$\alpha$-finie-étale (\cite{agt} V.7.4). 
L'anneau $R^\intern_\infty$ étant normal et intègre, il existe donc une $R^\intern_\infty$-algèbre $\cE'$ telle que la clôture intégrale
de $R^\intern_\infty$ dans $\cN\otimes_{R_\infty}R^\intern_\infty\otimes_{\co_K}K$ soit isomorphe à $\cE\times \cE'$.
Par suite, la $R^\intern_\infty$-algèbre $\cE$ est $\alpha$-finie-étale en vertu de (\cite{agt} V.7.11 et V.7.4).
Par ailleurs, l'extension $N/F_\infty$ étant galoisienne, il en est de même de l'extension $E/F^\intern_\infty$, d'où la proposition. 

\begin{cor}\label{eccr37}
Soit $M$ un $\oR^\intern$-module muni d'une action $\oR^\intern$-semi-linéaire continue de $\Sigma^\intern$
pour la topologie discrète de $M$. Alors $\rH^i(\Sigma^\intern,M)$ est $\alpha$-nul pour tout $i\geq 1$, 
et le morphisme canonique $M^{\Sigma^\intern}\otimes_{R^\intern_\infty}\oR^\intern\rightarrow M$ est un $\alpha$-isomorphisme.
\end{cor}

Soient $E$ une extension galoisienne finie de $F^\intern_\infty$ contenue dans $\oF^\intern$, $D$ la clôture intégrale
de $R^\intern_\infty$ dans $E$, $G=\Gal(E/F^\intern_\infty)$, $\Sigma_E=\Gal(\oF^\intern/E)$. Supposons que la $R^\intern_\infty$-algèbre $D$ 
soit $\alpha$-finie-étale. Alors $D$ est un $\alpha$-$G$-torseur sur $R^\intern_\infty$ en vertu de \ref{aet8}. 
On en déduit par (\cite{agt} V.12.5 et V.12.8) que, 
pour tout $i\geq 1$, $\rH^i(G,M^{\Sigma_E})$ est $\alpha$-nul, et le morphisme canonique 
$M^{\Sigma^\intern}\otimes_{R^\intern_\infty}D\rightarrow M^{\Sigma_E}$ est un $\alpha$-isomorphisme.
La proposition s'ensuit par passage à la limite inductive en vertu de \ref{eccr35}. 

\begin{cor}\label{eccr36}
Pour tout $a\in \co_\oK$, l'homomorphisme canonique
\begin{equation}
R^\intern_\infty/a R^\intern_\infty \rightarrow (\oR^\intern/a \oR^\intern)^{\Sigma^\intern}
\end{equation}
est un $\alpha$-isomorphisme.
\end{cor}

Soient $E$ une extension galoisienne finie de $F^\intern_\infty$ contenue dans $\oF^\intern$, $D$ la clôture intégrale
de $R^\intern_\infty$ dans $E$, $G=\Gal(E/F^\intern_\infty)$, 
$\Tr_G$ l'endomorphisme $R^\intern_\infty$-linéaire de $D$ (ou de $D/aD$) induit par 
$\sum_{\sigma\in G}\sigma$. Comme on a $D=\oR^\intern\cap E$ et $R^\intern_\infty=D\cap F^\intern_\infty$, 
les homomorphismes $R^\intern_\infty/aR^\intern_\infty\rightarrow D/aD\rightarrow \oR^\intern/a\oR^\intern$ sont injectifs. 
Supposons que $D$ soit $\alpha$-finie-étale sur $R^\intern_\infty$.  
Alors $D$ est un $\alpha$-$G$-torseur sur $R^\intern_\infty$ en vertu de \ref{aet8}. Par suite, le quotient
\begin{equation}
\frac{(D/aD)^G}{\Tr_G(D/aD)}
\end{equation}
est $\alpha$-nul en vertu de (\cite{agt} V.12.8). Comme $\Tr_G(D)\subset R^\intern_\infty$, l'homomorphisme 
$R^\intern_\infty/aR^\intern_\infty\rightarrow (D/aD)^G$ est un $\alpha$-isomorphisme. 
La proposition s'en déduit par passage à la limite inductive d'après \ref{eccr35}.

\begin{prop}\label{eccr11}
Les $\co_{\oK}$-algèbres $R^\intern_{\infty}$ et $R^\intern_{p^\infty}$ sont universellement $\alpha$-cohérentes \eqref{afini14}.
\end{prop}

En effet, $R^\intern_{\infty}$ (resp. $R^\intern_{p^\infty}$) est une $R_\infty$-algèbre (resp. $R_{p^\infty}$-algèbre) 
de présentation finie d'après \ref{eccr9}(ii)-(iii). La proposition résulte donc de \ref{cad15}, compte tenu de \ref{finita20}.

\begin{rema}\label{eccr10}
Soit $n$ un entier $\geq 1$. D'après \ref{eccr5}(i), il est loisible de remplacer $g$ \eqref{mtfla2} et sa carte relativement adéquate \eqref{eccr1a}
par le morphisme $g^{(n)}$ \eqref{eccr4g} et sa carte relativement adéquate \eqref{eccr5a}. 
Les anneaux $(R_m)_{m\geq 1}$ (resp. $(R'_m)_{m\geq 1}$) sont remplacés par $(R_{mn})_{m\geq 1}$ (resp. $(R'_{mn})_{m\geq 1}$)
en vertu de \ref{eccr8}(i). 
Les anneaux $R_{\infty}$ et $R'_{\infty}$ restent donc inchangés. Si $n$ est une puissance de $p$, 
les anneaux $R_{p^\infty}$ et $R'_{p^\infty}$ restent aussi inchangés.

Pour tout entier $m\geq 1$, l'anneau $R^\intern_{m}$ est remplacé par l'anneau 
d'une composante irréductible (ou ce qui revient au même connexe) de $\Spec(R^\intern_{mn}\otimes_{R^\intern_{n}}R'_{n})$
d'après \ref{eccr8}(i)-(ii). Ce dernier schéma est intègre si $n$ est une puissance de $p$ d'après \ref{eccr8}(v). 
L'anneau $R^\intern_\infty$ est donc remplacé par l'anneau d'une composante irréductible de 
$\Spec(R^\intern_\infty\otimes_{R^\intern_{n}}R'_{n})$, compte tenu de \ref{eccr9}(iv). Si $n$ est une puissance de $p$,  
les anneaux $R^\intern_{p^\infty}$ et $R^\intern_{\infty}$ sont remplacés par 
$R^\intern_{p^\infty}\otimes_{R^\intern_{n}}R'_{n}$ et $R^\intern_{\infty}\otimes_{R^\intern_{n}}R'_{n}$, respectivement, 
compte tenu de \ref{eccr9}(v); l'anneau $R'_{p^\infty}\otimes_{R^\intern_{p^\infty}}R^\intern_\infty$
reste donc inchangé. 
\end{rema}

\begin{prop}\label{eccr100}
\
\begin{itemize}
\item[{\rm (i)}] L'homomorphisme canonique 
\begin{equation}\label{eccr100a}
R_\infty\otimes_{R}R'\rightarrow (\oR^!)^\Sigma
\end{equation}
est un isomorphisme  \eqref{eccr40g}. 
\item[{\rm (ii)}] Pour tout $a\in \co_\oK$, l'homomorphisme canonique 
\begin{equation}\label{eccr100b}
(R_\infty/aR_\infty)\otimes_{R}R'\rightarrow (\oR^!/a\oR^!)^\Sigma
\end{equation}
est un $\alpha$-isomorphisme. 
\item[{\rm (iii)}] Pour tout $a\in \co_\oK$, l'homomorphisme canonique 
\begin{equation}\label{eccr100c}
(\oR/a\oR)\otimes_RR'\rightarrow \oR^!/a\oR^!
\end{equation}
est un $\alpha$-isomorphisme. 
\end{itemize}
\end{prop}

(i) Reprenons les notations de \ref{eccr40}. Pour tout $i\in I$, notons $A_i$ l'anneau du schéma affine normal et intègre $V_i$ 
et posons 
\begin{equation}
A=\underset{\underset{i\in I}{\longrightarrow}}{\lim}\ A_i. 
\end{equation}
Le schéma $\oX'$ étant normal et localement irréductible (\cite{agt} III.4.2(iii)), il en est de même des schémas 
$V_i\times_XX'$ $(i\in I)$ (\cite{agt} III.3.3).   
L'immersion $X'^\rhd\rightarrow X'$ étant schématiquement dominante (\cite{agt} III.4.2(iv)), il en est de même des immersions 
$V_i\times_{X^\circ}X'^\rhd\rightarrow V_i\times_XX'$ $(i\in I)$ (\cite{ega4} 11.10.5).
L'anneau  $\oR^!$ \eqref{eccr40g} est donc la clôture intégrale de $R'_1$ dans $A\otimes_{R}R'$.
Par suite, l'anneau $(\oR^!)^\Sigma$ est la clôture intégrale de $R'_1$ dans $(A\otimes_{R}R')^\Sigma$. 
L'homomorphisme canonique 
\begin{equation}
A^\Sigma\otimes_{R}R'\rightarrow (A\otimes_{R}R')^\Sigma
\end{equation}
est un isomorphisme en vertu de (\cite{agt} II.3.15). 
L'anneau $R_\infty$ est la clôture intégrale de $R_1$ dans $A^\Sigma$ d'après \ref{cad14}.
Comme le schéma $\Spec(R_\infty\otimes_{R}R')$ est normal et localement irréductible d'après \ref{eccr9}(ii), 
on déduit de ce qui précède que l'anneau $R_\infty\otimes_{R}R'$ est la clôture intégrale de $R'_1$ dans $(A\otimes_{R}R')^\Sigma$, 
d'où la proposition. 

(ii) Notons $\oF$ (resp. $F_\infty$) le corps des fractions de $\oR$ (resp. $R_\infty$). 
D'après \ref{tpcg5}, l'extension $\oF$ de $F_\infty$ est la réunion d'un système inductif filtrant de 
sous-extensions galoisiennes finies $(N_i)_{i\in I}$ de $F_\infty$ telles que la clôture intégrale $\cN_i$ 
de $R_\infty$ dans $N_i$ soit $\alpha$-finie-étale sur $R_\infty$. 
Pour tout $i\in I$, la $(R_\infty\otimes_RR')$-algèbre $\cN_i\otimes_{R}R'$ est $\alpha$-finie-étale d'après (\cite{agt} V.7.4). 
Notons $\cN^!_i$ la clôture intégrale de $R'_1$ dans $\cN_i\otimes_{R_\infty}A^\Sigma\otimes_RR'$. 
Le schéma $\Spec(R_\infty\otimes_RR')$ étant normal et localement irréductible d'après \ref{eccr9}(ii), 
$\cN^!_i$ est aussi la clôture intégrale de $R_\infty\otimes_RR'$ dans $\cN_i\otimes_RR'\otimes_{\co_K}K$. 
Par suite, la $(R_\infty\otimes_RR')$-algèbre $\cN^!_i$ est $\alpha$-finie-étale en vertu de \ref{aet3}.
Par ailleurs, l'homomorphisme canonique de $(R_\infty\otimes_RR')$-algèbres
\begin{equation}\label{eccr100d}
\underset{\underset{i\in I}{\longrightarrow}}{\lim}\ \cN^!_i \rightarrow \oR^!
\end{equation}
est un isomorphisme. On en déduit, en vertu de (\cite{agt} V.12.8), que $\rH^1(\Sigma,\oR^!)$ est $\alpha$-nul.  Par suite, le morphisme
canonique 
\begin{equation}
(\oR^!)^\Sigma/a(\oR^!)^\Sigma\rightarrow (\oR^!/a\oR^!)^\Sigma
\end{equation}
est injectif et $\alpha$-surjectif. 
Composant avec la réduction modulo $a$ de l'isomorphisme \eqref{eccr100a}, on en déduit que le morphisme \eqref{eccr100b}
est injectif et $\alpha$-surjectif. 

(iii) On peut se borner au cas où $a=0$ (\cite{agt} V.2.3). 
On notera que l'action de $\Sigma$ sur $\oR^!$ est continue pour la topologie discrète. 
En vertu de (\cite{agt} II.6.19), le morphisme canonique
\begin{equation}
(\oR^!)^\Sigma\otimes_{R_\infty}\oR\rightarrow \oR^!
\end{equation}
est un $\alpha$-isomorphisme. La proposition s'ensuit compte tenu de (i).

\subsection{}\label{eccr12}
Posons 
\begin{eqnarray}
\Xi_{p^n}&=&\Hom(\Delta_{p^\infty},\mu_{p^n}(\co_\oK)),\label{eccr12a}\\
\Xi_{p^\infty}&=&\Hom(\Delta_{p^\infty},\mu_{p^\infty}(\co_\oK)),\label{eccr12b}\\
\Xi'_{p^n}&=&\Hom(\Delta'_{p^\infty},\mu_{p^n}(\co_\oK)),\label{eccr12c}\\
\Xi'_{p^\infty}&=&\Hom(\Delta'_{p^\infty},\mu_{p^\infty}(\co_\oK)).\label{eccr12d}
\end{eqnarray}
On identifie $\Xi_{p^n}$ (resp. $\Xi'_{p^n}$) à un sous-groupe de $\Xi_{p^\infty}$ (resp. $\Xi'_{p^\infty}$). D'après (\cite{agt} II.8.9), 
il existe une décomposition canonique de $R_{p^\infty}$ en somme directe de $R_1$-modules de présentation finie,
stables sous l'action de $\Delta_{p^\infty}$,
\begin{equation}\label{eccr12e}
R_{p^\infty}=\bigoplus_{\lambda\in \Xi_{p^\infty}}R^{(\lambda)}_{p^\infty},
\end{equation}
telle que l'action de $\Delta_{p^\infty}$ sur le facteur $R^{(\lambda)}_{p^\infty}$ soit donnée par le caractère $\lambda$. 
De plus, pour tout $n\geq 0$, on a 
\begin{equation}\label{eccr12f}
R_{p^n}=\bigoplus_{\lambda\in \Xi_{p^n}}R^{(\lambda)}_{p^\infty}. 
\end{equation}
De même, il existe une décomposition canonique de $R'_{p^\infty}$ en somme directe de $R'_1$-modules de présentation finie,
stables sous l'action de $\Delta'_{p^\infty}$,
\begin{equation}\label{eccr12g}
R'_{p^\infty}=\bigoplus_{\lambda\in \Xi'_{p^\infty}}R'^{(\lambda)}_{p^\infty},
\end{equation}
telle que l'action de $\Delta'_{p^\infty}$ sur le facteur $R'^{(\lambda)}_{p^\infty}$ soit donnée par le caractère $\lambda$. 
De plus, pour tout $n\geq 0$, on a 
\begin{equation}\label{eccr12h}
R'_{p^n}=\bigoplus_{\lambda\in \Xi'_{p^n}}R'^{(\lambda)}_{p^\infty}. 
\end{equation}

Posons 
\begin{eqnarray}
\chi_{p^n}&=&\Hom(\fS_{p^\infty},\mu_{p^n}(\co_\oK)),\label{eccr12i}\\
\chi_{p^\infty}&=&\Hom(\fS_{p^\infty},\mu_{p^\infty}(\co_\oK)).\label{eccr12j}
\end{eqnarray}
On identifie $\chi_{p^n}$ à un sous-groupe de $\chi_{p^\infty}$. 
Compte tenu de \eqref{eccr43k} et \eqref{eccr8h}, la suite canonique de $\mZ_p$-modules  
\begin{equation}\label{eccr12k}
0\rightarrow \fS_{p^\infty}\rightarrow \Delta'_{p^\infty}\rightarrow \Delta_{p^\infty} \rightarrow 0
\end{equation}
est exacte et scindée. On en déduit une suite exacte 
\begin{equation}\label{eccr12l}
0\rightarrow \Xi_{p^\infty}\rightarrow \Xi'_{p^\infty}\rightarrow \chi_{p^\infty} \rightarrow 0.
\end{equation}

\begin{lem}\label{eccr13}
L'homomorphisme canonique $R^\intern_{p^\infty}\rightarrow R'_{p^\infty}$ induit un isomorphisme
\begin{equation}\label{eccr13a}
R^\intern_{p^\infty}\stackrel{\sim}{\rightarrow}\bigoplus_{\lambda\in \Xi_{p^\infty}}R'^{(\lambda)}_{p^\infty}. 
\end{equation}
\end{lem}

Comme $R^\intern_{p^\infty}$ est intègre et normal \eqref{eccr42h2} et 
que $R'_{p^\infty}$ est entier sur $R^\intern_{p^\infty}$,  on a d'après \eqref{eccr43k}
\begin{equation}
R^\intern_{p^\infty}=(R'_{p^\infty})^{\fS_{p^\infty}}.
\end{equation}
La proposition s'ensuit compte tenu de \eqref{eccr12l}.

\subsection{}\label{eccr16}
Pour tout $\nu\in \chi_{p^\infty}$, fixons un relèvement $\tnu\in \Xi'_{p^\infty}$ \eqref{eccr12l} et posons 
\begin{equation}\label{eccr16a}
R'^{(\nu)}_{p^\infty}=\bigoplus_{\lambda\in \Xi_{p^\infty}}R'^{(\tnu+\lambda)}_{p^\infty},
\end{equation}
qui est naturellement un $R^\intern_{p^\infty}$-module \eqref{eccr13}. 
La décomposition \eqref{eccr12g} induit alors une décomposition en $R^\intern_{p^\infty}$-modules
\begin{equation}\label{eccr16b}
R'_{p^\infty}=\bigoplus_{\nu\in \chi_{p^\infty}}R'^{(\nu)}_{p^\infty}
\end{equation}
telle que l'action de $\fS_{p^\infty}$ sur le facteur $R'^{(\nu)}_{p^\infty}$ soit donnée par le caractère $\nu$.

\begin{lem}\label{eccr17}
Pour tout $n\geq 0$, l'homomorphisme canonique $R^\intern_{p^\infty}\otimes_{R^\intern_{p^n}}R'_{p^n}\rightarrow R'_{p^\infty}$ induit un isomorphisme
\begin{equation}\label{eccr17a}
R^\intern_{p^\infty}\otimes_{R^\intern_{p^n}}R'_{p^n}\stackrel{\sim}{\rightarrow}\bigoplus_{\nu\in \chi_{p^n}}R'^{(\nu)}_{p^\infty}. 
\end{equation}
En particulier, pour tout $\nu\in \chi_{p^n}$, le $R^\intern_{p^\infty}$-module $R'^{(\nu)}_{p^\infty}$ est de présentation finie.
\end{lem}

En effet, remplaçant $g$ par le morphisme $g^{(p^n)}$, cela résulte de \ref{eccr13} compte tenu de \eqref{eccr8d} et \ref{eccr10}.

\begin{prop}[\cite{agt} II.8.1]\label{eccr31}
Soient $G$ un $\mZ_p$-module libre de type fini, $n$ un entier $\geq 0$, $\nu\colon G\rightarrow \mu_{p^n}(\co_\oK)$ 
un homomorphisme surjectif, $\zeta$ un générateur du groupe $\mu_{p^n}(\co_\oK)$, $a\in \co_\oK$, 
$\fq$ l'idéal de $\co_\oK$ engendré par $a$ et $\zeta-1$. 
Soit $A$ une $\co_\oK$-algèbre complète et séparée pour la topologie $p$-adique et $\co_{\oK}$-plate. 
On note  $A(\nu)$ le $A$-$G$-module topologique $A$, muni de la topologie $p$-adique
et de l'action de $G$ définie par la multiplication par $\nu$. 
\begin{itemize}
\item[{\rm (i)}] Si $\nu=1$ (i.e., $n=0$), alors on a un isomorphisme canonique de $A$-algèbres graduées 
\begin{equation}\label{eccr31a}
\wedge(\Hom_{\mZ_p}(G,A/aA))
\stackrel{\sim}{\rightarrow}  \rH^*_\cont(G,A(\nu)/aA(\nu)).
\end{equation}
\item[{\rm (ii)}] Si $\nu\not=1$ (i.e., $n\not=0$), alors $\rH^i_\cont(G,A(\nu)/aA(\nu))$
est un $A/\fq A$-module libre de type fini pour tout $i\geq 0$ et est nul pour tout $i\geq \rg_{\mZ_p}(G)+1$.
\item[{\rm (iii)}] Le système projectif $(\rH^*(G,A(\nu)/p^rA(\nu)))_{r\geq 0}$ 
vérifie la condition de Mittag-Leffler uniformément en $\nu$, autrement dit, si on note, pour tous entiers $r'\geq r\geq 0$,
\begin{equation}\label{eccr31b}
h_{r,r'}^\nu\colon \rH^*(G,A(\nu)/p^{r'}A(\nu))\rightarrow \rH^*(G,A(\nu)/p^rA(\nu))
\end{equation}
le morphisme canonique, alors pour tout entier $r\geq 1$,
il existe un entier $r'\geq r$, dépendant du $\mZ_p$-rang de $G$ mais pas de $\nu$, tel que pour tout entier $r''\geq r'$, 
les images de $h_{r,r'}^\nu$ et $h_{r,r''}^\nu$ soient égales. 
\end{itemize}
\end{prop}

La preuve est identique à celle de (\cite{agt} II.8.1). 

\begin{prop}\label{eccr30}
Soient $a$ un élément non nul de $\co_\oK$, $\zeta$ une racine primitive $p$-ième de $1$ dans $\co_\oK$.
Posons $\delta=\dim_{\mQ}((P'^\gp/h^\gp(P^\gp))\otimes_{\mZ}\mQ)$. 
Alors~: 
\begin{itemize}
\item[{\rm (i)}] Il existe un et un unique homomorphisme de $R^\intern_\infty$-algèbres graduées
\begin{equation}\label{eccr30a}
\wedge(\Hom_{\mZ}(\fS_{p^\infty},R^\intern_\infty/aR^\intern_\infty))\rightarrow 
\rH^*(\fS_{p^\infty},R'_{p^\infty}\otimes_{R^\intern_{p^\infty}}(R^\intern_\infty/aR^\intern_\infty))
\end{equation}
dont la composante en degré un est le composé des morphismes canoniques
\begin{equation}\label{eccr30b}
\Hom_{\mZ}(\fS_{p^\infty},R^\intern_\infty/aR^\intern_\infty)\stackrel{\sim}{\rightarrow}
\rH^1(\fS_{p^\infty},R^\intern_\infty/aR^\intern_\infty)\rightarrow 
\rH^1(\fS_{p^\infty},R'_{p^\infty}\otimes_{R^\intern_{p^\infty}}(R^\intern_\infty/aR^\intern_\infty)).
\end{equation}
Celui-ci admet, en tant que morphisme de $R^\intern_\infty$-modules gradués,  un inverse à gauche canonique
\begin{equation}\label{eccr30c}
\rH^*(\fS_{p^\infty},R'_{p^\infty}\otimes_{R^\intern_{p^\infty}}(R^\intern_\infty/aR^\intern_\infty))\rightarrow 
\wedge(\Hom_{\mZ}(\fS_{p^\infty},R^\intern_\infty/aR^\intern_\infty)),
\end{equation}
dont le noyau est annulé par $\zeta-1$.  
\item[{\rm (ii)}] Le $R^\intern_\infty$-module $\rH^i(\fS_{p^\infty},R'_{p^\infty}\otimes_{R^\intern_{p^\infty}}(R^\intern_\infty/aR^\intern_\infty))$ 
est de présentation $\alpha$-finie pour tout $i\geq 0$, et est nul pour tout $i\geq \delta+1$. 
\item[{\rm (iii)}] Le système projectif $(\rH^*(\fS_{p^\infty},R'_{p^\infty}\otimes_{R^\intern_{p^\infty}}(R^\intern_\infty/p^rR^\intern_\infty))_{r\geq 0}$ vérifie la condition 
de Mittag-Leffler~; plus précisément, si on note, pour tous entiers $r'\geq r\geq 0$,
\begin{equation}\label{eccr30d}
h_{r,r'}\colon \rH^*(\fS_{p^\infty},R'_{p^\infty}\otimes_{R^\intern_{p^\infty}}(R^\intern_\infty/p^{r'}R^\intern_\infty))\rightarrow 
\rH^*(\fS_{p^\infty},R'_{p^\infty}\otimes_{R^\intern_{p^\infty}}(R^\intern_\infty/p^rR^\intern_\infty))
\end{equation}
le morphisme canonique, alors pour tout entier $r\geq 1$,
il existe un entier $r'\geq r$, dépendant seulement de $\delta$ mais pas des autres données dans \ref{eccr1}, 
tel que pour tout entier $r''\geq r'$, les images de $h_{r,r'}$ et $h_{r,r''}$ soient égales.  
\end{itemize}
\end{prop}

En effet, d'après \eqref{eccr16b}, on a une décomposition canonique de $R'_{p^\infty}\otimes_{R^\intern_{p^\infty}}R^\intern_\infty$ en 
somme directe de $R^\intern_\infty[\fS_{p^\infty}]$-modules
\begin{equation}\label{eccr30e}
R'_{p^\infty}\otimes_{R^\intern_{p^\infty}}R^\intern_\infty=\bigoplus_{\nu\in \chi_{p^\infty}}
R'^{(\nu)}_{p^\infty}\otimes_{R^\intern_{p^\infty}}R^\intern_\infty\otimes_{\co_{\oK}}\co_{\oK}(\nu),
\end{equation}
où $\fS_{p^\infty}$ agit trivialement sur $R'^{(\nu)}_{p^\infty}\otimes_{R^\intern_{p^\infty}}R^\intern_\infty$ et agit sur $\co_\oK(\nu)=\co_\oK$ 
par le caractère $\nu$. Comme $R'_{p^\infty}\otimes_{R^\intern_{p^\infty}}R^\intern_\infty$ est intègre d'après \ref{eccr9}(vi) et est donc $\co_{\oK}$-plat,
les $\co_{\oK}$-modules $R'^{(\nu)}_{p^\infty}\otimes_{R^\intern_{p^\infty}}R^\intern_\infty$ sont plats. On a donc,  
en vertu de (\cite{agt} II.3.15), une décomposition canonique 
en somme directe de $R^\intern_\infty$-modules
\begin{eqnarray}\label{eccr30f}
\lefteqn{\rH^*(\fS_{p^\infty},R'_{p^\infty}\otimes_{R^\intern_{p^\infty}}(R^\intern_\infty/aR^\intern_\infty))}\\
&=&\bigoplus_{\nu\in \chi_{p^\infty}}
\rH^*(\fS_{p^\infty},\co_{\oK}(\nu)/a\co_\oK(\nu))\otimes_{\co_{\oK}}R'^{(\nu)}_{p^\infty}\otimes_{R^\intern_{p^\infty}}R^\intern_\infty. \nonumber
\end{eqnarray}

(i) Si $\nu_0$ est l'unité de $\chi_{p^\infty}$, on a $R'^{(\nu_0)}_{p^\infty}=R^\intern_{p^\infty}$ d'après \ref{eccr13}, 
de sorte que la composante en $\nu_0$ de la décomposition \eqref{eccr30f}
est l'image de l'homomorphisme canonique de $R^\intern_{\infty}$-algèbres graduées 
\begin{equation}\label{eccr30g}
\rH^*(\fS_{p^\infty},\co_{\oK}/a\co_\oK)\otimes_{\co_{\oK}}R^\intern_\infty \rightarrow 
\rH^*(\fS_{p^\infty},R'_{p^\infty}\otimes_{R^\intern_{p^\infty}}(R^\intern_\infty/aR^\intern_\infty)).
\end{equation}
Par ailleurs, l'homomorphisme canonique de $R^\intern_\infty$-algèbres graduées
\begin{equation}\label{eccr30h}
\wedge(\Hom_{\mZ}(\fS_{p^\infty},\co_\oK/a\co_\oK))\otimes_{\co_\oK} R^\intern_\infty
\rightarrow \wedge(\Hom_{\mZ}(\fS_{p^\infty},R^\intern_\infty/aR^\intern_\infty))
\end{equation}
est un isomorphisme. La proposition résulte alors de \eqref{eccr30f}, \ref{eccr31} et \eqref{eccr9c}.

(ii) Soient $i,n$ deux entiers $\geq 0$, $\zeta_n$ une racine primitive $p^n$-ième de l'unité. 
Il résulte de \eqref{eccr30f}, \ref{eccr17}, \ref{eccr31} et \eqref{eccr9c} 
que $\rH^i(\fS_{p^\infty},R'_{p^\infty}\otimes_{R^\intern_{p^\infty}}(R^\intern_\infty/aR^\intern_\infty))$ 
est la somme directe d'un $R^\intern_\infty$-module de présentation finie et d'un $R^\intern_\infty$-module annulé par $\zeta_n-1$. 
Il est donc de présentation $\alpha$-finie sur $R^\intern_\infty$.
La seconde assertion est claire puisque 
la $p$-dimension cohomologique de $\fS_{p^\infty}$ est égale à $\delta$ d'après \eqref{eccr9c}.

(iii) Cela résulte de \eqref{eccr30f}, \ref{eccr31}(iii) et \eqref{eccr9c}. 

\subsection{}\label{eccr32}
Soit $a$ un élément non nul de $\co_\oK$. 
Comme le sous-groupe de torsion de $P'^\gp/h^\gp(P^\gp)$ est d'ordre premier à $p$,  
on a d'après \eqref{eccr9c}, un isomorphisme canonique 
\begin{equation}\label{eccr32a}
(P'^\gp/h^\gp(P^\gp))\otimes_\mZ\mZ_p(-1)\stackrel{\sim}{\rightarrow}
\Hom_{\mZ_p}(\fS_{p^\infty}, \mZ_p).
\end{equation}
On en déduit, compte tenu de (\cite{ogus} IV 1.1.4), un isomorphisme $R^\intern_\infty$-linéaire \eqref{eccr1b}
\begin{equation}\label{eccr32b}
\tOmega^1_{R'/R}\otimes_{R'}(R^\intern_\infty/aR^\intern_\infty)(-1)\stackrel{\sim}{\rightarrow} 
\Hom_{\mZ}(\fS_{p^\infty},R^\intern_\infty/aR^\intern_\infty).
\end{equation}
On interprètera dans la suite le but de ce morphisme comme un groupe de cohomologie. Posons 
\begin{equation}\label{eccr32c}
\delta=\dim_{\mQ}((P'^\gp/h^\gp(P^\gp))\otimes_{\mZ}\mQ).
\end{equation}

\begin{cor}\label{eccr33}
Soit $a$ un élément non nul de $\co_\oK$. 
Il existe un et un unique homomorphisme de $R^\intern_\infty$-algèbres graduées
\begin{equation}\label{eccr33a}
\wedge(\tOmega^1_{R'/R}\otimes_{R'}(R^\intern_\infty/aR^\intern_\infty)(-1))\rightarrow 
\rH^*(\fS_{p^\infty},R'_{p^\infty}\otimes_{R^\intern_{p^\infty}}(R^\intern_\infty/aR^\intern_\infty))
\end{equation}
dont la composante en degré un est induite par \eqref{eccr32b} et \eqref{eccr30b}.
Celui-ci admet, en tant que morphisme de $R^\intern_\infty$-modules gradués, un inverse à gauche canonique 
\begin{equation}\label{eccr33b}
\rH^*(\fS_{p^\infty},R'_{p^\infty}\otimes_{R^\intern_{p^\infty}}(R^\intern_\infty/aR^\intern_\infty))\rightarrow 
\wedge(\tOmega^1_{R'/R}\otimes_{R'}(R^\intern_\infty/aR^\intern_\infty)(-1)),
\end{equation}
dont le noyau est annulé par $p^{\frac{1}{p-1}}$.  
\end{cor}

\begin{prop}\label{eccr50}
Pour tout élément non nul $a$ de $\co_\oK$ et tout entier $i\geq 0$, le morphisme canonique 
\begin{equation}\label{eccr50a}
\rH^i(\fS_\infty,R'_\infty/a R'_\infty)\otimes_{R^\intern_\infty}\oR^\intern\rightarrow \rH^i(\Pi,\oR'/a \oR')
\end{equation}
est un $\alpha$-isomorphisme, et le morphisme canonique
\begin{equation}\label{eccr50b}
\rH^i(\fS_{p^\infty},R'_{p^\infty}\otimes_{R^\intern_{p^\infty}}(R^\intern_\infty/aR^\intern_\infty))\rightarrow
\rH^i(\fS_\infty,R'_\infty/a R'_\infty)
\end{equation}
est un isomorphisme. 
\end{prop}

Nous utiliserons les notations de \ref{eccr34}.
Il résulte de (\cite{agt} II.6.19 et II.6.22) et de la suite spectrale 
\begin{equation}
\rE_2^{ij}=\rH^i(\fS_\infty,\rH^j(\Sigma',\oR'/a\oR'))\Rightarrow \rH^{i+j}(\fS,\oR'/a\oR')
\end{equation}
que pour tout entier $i\geq 0$, le morphisme canonique 
\begin{equation}
\rH^i(\fS_\infty,R'_\infty/aR'_\infty)\rightarrow \rH^i(\fS,\oR'/a\oR')
\end{equation}
est un $\alpha$-isomorphisme. 
Il résulte de \ref{eccr37} et de la suite spectrale 
\begin{equation}
\rE_2^{ij}=\rH^i(\Sigma^\intern,\rH^j(\Pi,\oR'/a\oR'))\Rightarrow \rH^{i+j}(\fS,\oR'/a\oR')
\end{equation}
que pour tout entier $i\geq 0$, le morphisme canonique
\begin{equation}
\rH^i(\fS,\oR'/a\oR')\rightarrow \rH^i(\Pi,\oR'/a\oR')^{\Sigma^\intern}
\end{equation}
est un $\alpha$-isomorphisme, et que le morphisme $\oR^\intern$-linéaire canonique 
\begin{equation}
\rH^i(\Pi,\oR'/a\oR')^{\Sigma^\intern}\otimes_{R^\intern_\infty}\oR^\intern\rightarrow \rH^i(\Pi,\oR'/a\oR')
\end{equation}
est un $\alpha$-isomorphisme.
Le diagramme commutatif de morphismes canoniques
\begin{equation}
\xymatrix{
{\rH^i(\fS,\oR'/a\oR')}\ar[r]&{\rH^i(\Pi,\oR'/a\oR')^{\Sigma^\intern}}\ar[d]\\
{\rH^i(\fS_{\infty},R'_\infty/aR'_\infty)}\ar[u]\ar[r]&{\rH^i(\Pi,\oR'/a\oR')}}
\end{equation}
implique alors que le morphisme \eqref{eccr50a} est un $\alpha$-isomorphisme. 

Par ailleurs, le groupe profini $\fN$ étant d'ordre premier à $p$ \eqref{eccr43l}, il est de $p$-dimension cohomologique nulle.
Il résulte alors de \eqref{eccr9e} et de la suite spectrale de Hochschild-Serre que le morphisme \eqref{eccr50b} est un isomorphisme.

\begin{prop}\label{eccr51}
Soit $a$ un élément non nul de $\co_\oK$. 
\begin{itemize}
\item[{\rm (i)}] Il existe un et un unique homomorphisme de $\oR^\intern$-algèbres graduées
\begin{equation}\label{eccr51a}
\wedge(\tOmega^1_{R'/R}\otimes_{R'}(\oR^\intern/a\oR^\intern)(-1))\rightarrow \rH^*(\Pi,\oR'/a\oR')
\end{equation}
dont la composante en degré un est induite par \eqref{eccr32b} (cf. \ref{eccr34}). Celui-ci 
est $\alpha$-injectif et son conoyau est annulé par $p^{\frac{1}{p-1}}\fm_\oK$.  
\item[{\rm (ii)}] Le $\oR^\intern$-module $\rH^i(\Pi,\oR'/a \oR')$ est de présentation $\alpha$-finie 
pour tout $i\geq 0$, et est $\alpha$-nul pour tout $i\geq \delta+1$ \eqref{eccr32c}. 
\item[{\rm (iii)}] Pour tous entiers $r'\geq r\geq 0$, notons
\begin{equation}\label{eccr51b}
\hbar_{r,r'}\colon \rH^*(\Pi,\oR'/p^{r'} \oR')\rightarrow \rH^*(\Pi,\oR'/p^r \oR')
\end{equation}
le morphisme canonique. Alors pour tout entier $r\geq 1$,
il existe un entier $r'\geq r$, dépendant seulement de $\delta$ mais pas des autres données dans \eqref{eccr1}, 
tel que pour tout entier $r''\geq r'$, les images de $\hbar_{r,r'}$ et $\hbar_{r,r''}$ soient $\alpha$-isomorphes.  
\end{itemize}
\end{prop}
Cela résulte de \ref{eccr30}, \ref{eccr33} et \ref{eccr50}.

\begin{prop}\label{eccr14}
Soit $m$ un entier $\geq 0$ et posons $\delta=\dim_{\mQ}((P'^\gp/h^\gp(P^\gp))\otimes_{\mZ}\mQ)$. Alors, 
\begin{itemize}
\item[{\rm (i)}] Pour tout entier $n\geq 0$, tout $(R'_{p^n}/p^mR'_{p^n})$-module de présentation finie $M$, 
muni d'une action semi-liné\-aire de $\fS_{p^n}$ et tout entier $i\geq 0$, 
le $R^\intern_{p^\infty}$-module $\rH^i(\fS_{p^\infty},M\otimes_{R'_{p^n}}R'_{p^\infty})$ est de présentation $\alpha$-finie, 
et est nul pour tout $i\geq \delta+1$.
\item[{\rm (ii)}] Pour tout entier $n\geq 0$, tout $(R'_{p^n}/p^mR'_{p^n})$-module de présentation finie $M$, 
muni d'une action semi-liné\-aire de $\fS_{p^n}$ et tout entier $i\geq 0$, 
le $R^\intern_\infty$-module 
\[
\rH^i(\fS_{p^\infty},M\otimes_{R'_{p^n}}R'_{p^\infty}\otimes_{R^\intern_{p^\infty}}R^\intern_\infty)
\] 
est de présentation $\alpha$-finie, et est nul pour tout $i\geq \delta+1$.
\item[{\rm (iii)}] Pour tout entier $n\geq 1$,  tout $(R'_{n}/p^mR'_{n})$-module de présentation finie $M$, 
muni d'une action semi-linéaire de $\fS_n$ et tout entier $i\geq 0$, 
le $R^\intern_\infty$-module $\rH^i(\fS_{\infty},M\otimes_{R'_n}R'_\infty)$ est de présentation $\alpha$-finie, 
et est nul pour tout $i\geq \delta+1$.
\end{itemize}
\end{prop}

(i) En effet, la $p$-dimension cohomologique de $\fS_{p^\infty}$ étant égale à $\delta$ (\cite{agt} II.3.24), 
il suffit de montrer que pour tout entier $i\geq 0$, 
le $R^\intern_{p^\infty}$-module $\rH^i(\fS_{p^\infty},M\otimes_{R'_{p^n}}R'_{p^\infty})$ est de présentation $\alpha$-finie. 
La $\co_\oK$-algèbre $R^\intern_{p^\infty}$ est universellement $\alpha$-cohérente en vertu de \ref{eccr11}. 
L'algèbre $R'_{p^n}$ étant finie et de présentation finie sur $R^\intern_{p^n}$ \eqref{eccr5c}, 
tout $(R^\intern_{p^\infty}\otimes_{R^\intern_{p^n}}R'_{p^n})$-module 
de présentation $\alpha$-finie est aussi un $R^\intern_{p^\infty}$-module de présentation 
$\alpha$-finie et est donc $\alpha$-cohérent d'après \ref{finita19}. 
Remplaçant $g$ par le morphisme $g^{(p^n)}$, on peut alors se réduire au cas où $n=0$,
compte tenu de \ref{eccr10}, \ref{tpcg18} et de la suite spectrale de Hochschild-Serre (\cite{agt} II.3.4). 

D'après \eqref{eccr16b}, on a une décomposition canonique de $R'_{p^\infty}$ en $R^\intern_{p^\infty}[\fS_{p^\infty}]$-modules 
\begin{equation}
R'_{p^\infty}=\bigoplus_{\nu\in \chi_{p^\infty}}R'^{(\nu)}_{p^\infty}\otimes_{\co_\oK}\co_\oK(\nu),
\end{equation}
où $\fS_{p^\infty}$ agit trivialement sur $R'^{(\nu)}_{p^\infty}$ et agit sur $\co_\oK(\nu)=\co_\oK$ par le caractère $\nu$. 
Par suite, pour tout entier $i\geq 0$, on a une décomposition canonique en $R^\intern_{p^\infty}$-modules
\begin{equation}\label{eccr14b}
\rH^i(\fS_{p^\infty},M\otimes_{R'_1}R'_{p^\infty})=\bigoplus_{\nu\in \chi_{p^\infty}}
\rH^i(\fS_{p^\infty},M\otimes_{R'_1}R'^{(\nu)}_{p^\infty}\otimes_{\co_\oK}\co_\oK(\nu)).
\end{equation}
D'après \ref{eccr17}, pour tout $\nu \in \chi_{p^\infty}$, le $R^\intern_{p^\infty}$-module $M\otimes_{R'_1}R'^{(\nu)}_{p^\infty}$ est de présentation finie 
et est donc $\alpha$-cohérent. 
En vertu de \ref{tpcg19}, pour tout $\gamma\in \fm_\oK$, il existe un entier $j\geq 1$ tel que pour tout $\nu\in \chi_{p^\infty}-\chi_{p^j}$ \eqref{eccr12j}, le $R_1$-module
$\rH^i(\fS_{p^\infty},M\otimes_{R'_1}R'^{(\nu)}_{p^\infty}\otimes_{\co_\oK}\co_\oK(\nu))$ soit annulé par $\gamma$. 
Par suite, d'après \ref{tpcg8}, le $R^\intern_{p^\infty}$-module $\rH^i(\fS_{p^\infty},M\otimes_{R'_1}R'_{p^\infty})$ est de présentation $\alpha$-finie. 

(ii) En effet, l'anneau  $R'_{p^\infty}\otimes_{R^\intern_{p^\infty}}R^\intern_\infty$ est intègre et normal d'après \ref{eccr9}(vi). 
La $\co_\oK$-algèbre $R^\intern_{\infty}$ est universellement $\alpha$-cohérente en vertu de \ref{eccr11}. 
Par ailleurs, remplaçant $g$ par le morphisme $g^{(p^n)}$, l'anneau $R'_{p^\infty}\otimes_{R^\intern_{p^\infty}}R^\intern_\infty$
reste inchangé alors que l'anneau $R^\intern_\infty$ est remplacé par $R^\intern_\infty\otimes_{R^\intern_{p^n}}R'_{p^n}$ d'après \ref{eccr10}. 
Il suffit donc de calquer la preuve de (i).

(iii) On notera d'abord que la $\co_\oK$-algèbre $R^\intern_{\infty}$ est universellement $\alpha$-cohérente en vertu de \ref{eccr11}. 
Soit $n'$ le plus grand diviseur premier à $p$ de $n$. Remplaçant $g$ par le morphisme $g^{(n')}$, 
l'anneau $R^\intern_{\infty}$ est remplacé par l'anneau $B$ d'une composante irréductible (ou ce qui revient au même connexe)
de $\Spec(R^\intern_{\infty}\otimes_{R^\intern_{n'}}R'_{n'})$ d'après \ref{eccr10}.
L'algèbre $R'_{n'}$ étant finie et de présentation finie sur $R^\intern_{n'}$ \eqref{eccr5c}, 
tout $B$-module de présentation $\alpha$-finie est aussi un $R^\intern_{\infty}$-module de présentation 
$\alpha$-finie et est donc $\alpha$-cohérent d'après \ref{finita19}. L'extension des corps de fractions 
de $B$ sur $R^\intern_{\infty}$ est galoisienne, de groupe $H$, un sous-groupe de $\fS_{n'}$. 
En particulier, $H$ est un groupe abélien d'ordre premier à $p$. 
Le foncteur $\Gamma(H, -)$ est exact, et il transforme les $R^\intern_{\infty}$-modules $\alpha$-cohérents en 
$R^\intern_{\infty}$-modules $\alpha$-cohérents d'après \ref{tpcg18}.  
Compte tenu de \ref{eccr10} et de la suite spectrale de Hochschild-Serre (\cite{agt} II.3.4), 
on peut donc supposer que $n$ est une puissance de $p$.

Considérons le groupe $\fN$ défini dans \eqref{eccr43l} qui est un groupe profini d'ordre premier à $p$. 
Pour tout $i\geq 0$, on a un isomorphisme canonique 
\begin{equation}
\rH^i(\fS_{p^\infty},(M\otimes_{R'_n}R'_\infty)^{\fN})\stackrel{\sim}{\rightarrow} \rH^i(\fS_\infty,M\otimes_{R'_n}R'_\infty).
\end{equation} 
D'après \ref{eccr9}(vii), l'homomorphisme canonique
\begin{equation}\label{eccr14a}
(R'_{p^\infty}\otimes_{R^\intern_{p^\infty}}R^\intern_\infty)/p^m(R'_{p^\infty}\otimes_{R^\intern_{p^\infty}}R^\intern_\infty) \rightarrow 
(R'_\infty/p^mR'_\infty)^\fN
\end{equation}
est un isomorphisme. 
Considérant une présentation finie de $M$ sur $R'_n/p^mR'_n$, 
on déduit de \eqref{eccr14a} que le morphisme canonique
\begin{equation}
M\otimes_{R'_n}R'_{p^\infty}\otimes_{R^\intern_{p^\infty}}R^\intern_\infty\rightarrow (M\otimes_{R'_n}R'_\infty)^{\fN}
\end{equation}
est un isomorphisme. La proposition résulte donc de (ii).

\begin{cor}\label{eccr15}
Soient $m$ un entier $\geq 1$, $\mL$ un $(\co_\oK/p^m\co_\oK)$-module libre de type fini, muni d'une action linéaire discrète de $\Delta'$,
$i$ un entier $\geq 0$. Posons $\delta=\dim_{\mQ}((P'^\gp/h^\gp(P^\gp))\otimes_{\mZ}\mQ)$.
Alors, le $R^\intern_\infty$-module $\rH^i(\fS,\mL\otimes_{\co_\oK}\oR')$ est de présentation $\alpha$-finie, 
et est $\alpha$-nul pour tout $i\geq \delta+1$ \eqref{eccr43m}.
\end{cor}

On rappelle que $\Sigma'$ désigne le noyau de l'homomorphisme canonique $\Delta'\rightarrow \Delta'_\infty$ \eqref{eccr43b3}. 
Compte tenu de \eqref{eccr43j}, on a une suite exacte 
\begin{equation}
0\rightarrow \Sigma'\rightarrow \fS\rightarrow \fS_\infty \rightarrow 0. 
\end{equation}
En vertu de (\cite{agt} II.6.19), $\rH^i(\Sigma',\mL\otimes_{\co_\oK} \oR')$ est $\alpha$-nul pour tout $i\geq 1$. 
Par la suite spectrale de Hochschild-Serre (\cite{agt} II.3.4), pour tout $i\geq 0$, le morphisme canonique 
\begin{equation}
\rH^i(\fS_\infty,(\mL\otimes_{\co_\oK} \oR')^{\Sigma'})\rightarrow \rH^i(\fS,\mL\otimes_{\co_\oK} \oR')
\end{equation} 
est donc un $\alpha$-isomorphisme. En vertu de \ref{tpcg10}(ii), 
pour tout $\gamma\in \fm_\oK$, il existe un entier $n\geq 1$, un $(R'_n/p^mR'_n)$-module de présentation finie $M_n$, 
muni d'une action semi-linéaire de $\Delta'_n$, et un morphisme $R'_{\infty}$-linéaire et $\Delta'_{\infty}$-équivariant 
\begin{equation}
M_n\otimes_{R'_n}R'_{\infty}\rightarrow (\mL\otimes_{\co_\oK}\oR')^{\Sigma'},
\end{equation}
dont le noyau et le conoyau sont annulés par $\gamma$. 
La proposition résulte alors de \ref{eccr14}(iii).

\section{\texorpdfstring{$\alpha$-finitude relative}{alpha-finitude relative}}\label{AFR}

Les hypothèses et notations de \ref{mtfla} sont en vigueur dans cette section.

\begin{prop}\label{AFR6}
Supposons le morphisme $g\colon X'\rightarrow X$ propre et le schéma $X$ affine d'anneau $R$ \eqref{mtfla2}, 
posons $R_1=R\otimes_{\co_K}\co_\oK$, et soient $m\geq 1$, $q\geq 0$ deux entiers, 
$\mL$ un $(\co_\oK/p^m\co_\oK)$-module localement libre de type fini de $\oX'^\rhd_\fet$. 
Alors, le $R_1$-module $\rH^q(\tE'_s,\beta'^*_m(\mL))$ est de présentation $\alpha$-finie \eqref{mtfla10i}.
\end{prop}

On notera d'abord que l'anneau $R_1/p^mR_1$ est universellement cohérent (\cite{egr1} 1.4.1 et 1.12.15).
Considérons la suite spectrale de Cartan Leray
\begin{equation}\label{AFR6a}
\rE_2^{i,j}=\rH^i(\oX'_m,\rR^j\Sigma'_{m*}(\beta'^*_m(\mL)))\Rightarrow \rH^{i+j}(\tE'_s,\beta'^*_m(\mL)).
\end{equation}
En vertu de \ref{finitude2} et \ref{afini9}, pour tous entiers $i,j\geq 0$, le $(R_1/p^mR_1)$-module $E_2^{i,j}$ 
est $\alpha$-cohérent. La proposition s'ensuit compte tenu de \ref{finita18}. 

\subsection{}\label{AFR9}
Dans la suite de cette section, on se donne un point $(\oy\rightsquigarrow \ox)$ de $X_\et\gtimes_{X_\et}\oX^\circ_\et$ \eqref{topfl17} 
tel que $\ox$ soit au-dessus de $s$. On note $\uX$ le localisé strict de $X$ en $\ox$ et on pose 
\begin{equation}\label{AFR9a}
\uoX=\oX\times_X\uX\ \ \ {\rm et} \ \ \ \uoX^\circ=\uoX\times_XX^\circ.
\end{equation}
D'après (\cite{agt} III.3.7), le schéma $\uoX$ est normal et strictement local (et en particulier intègre).
Le $X$-morphisme $u\colon \oy\rightarrow \uX$
définissant le point $(\oy\rightsquigarrow \ox)$ se relève en un $\oX^\circ$-morphisme $v\colon \oy\rightarrow \uoX^\circ$ et 
il induit donc un point géométrique de $\uoX^\circ$ que l'on note aussi (abusivement) $\oy$. On désigne par $\uDelta$
le groupe profini $\pi_1(\uoX^{\circ},\oy)$, par $\bB_{\uDelta}$ son topos classifiant et par
\begin{equation}\label{AFR9b}
\nu_\oy\colon \uoX^\circ_\fet \stackrel{\sim}{\rightarrow}\bB_{\uDelta}
\end{equation}
le foncteur fibre de $\uoX^\circ_\fet$ en $\oy$ \eqref{notconv11c}. On note
\begin{equation}\label{AFR9c}
\varphi_\ox\colon \tE\rightarrow \uoX^\circ_\fet
\end{equation}
le foncteur canonique défini dans \eqref{tf10b}. 
D'après (\cite{agt} VI.10.31 et VI.9.9), le foncteur composé
\begin{equation}\label{AFR9d}
\xymatrix{
{\tE}\ar[r]^-(0.5){\varphi_\ox}&{\uoX^\circ_\fet}\ar[r]^-(0.5){\nu_\oy}&{\bB_{\uDelta}}\ar[r]&\Ens},
\end{equation}
où la dernière flèche est le foncteur d'oubli de l'action de $\uDelta$, 
est canoniquement isomorphe au foncteur fibre associé au point $\rho(\oy\rightsquigarrow \ox)$ de $\tE$ \eqref{mtfla7e}. 

On pose $\uR_1=\Gamma(\uoX,\co_{\uoX})$ et on note $\uoR$ la $\uR_1$-algèbre $\oR_{\uX}^\oy$ de $\bB_{\uDelta}$ définie dans \eqref{TFA12f}.
On rappelle qu'on a un isomorphisme canonique de $\uR_1$-algèbres \eqref{TFA11c}
\begin{equation}\label{AFR9e}
\ocB_{\rho(\oy\rightsquigarrow \ox)}\stackrel{\sim}{\rightarrow} \uoR. 
\end{equation}

\begin{teo}\label{AFR8}
Supposons le morphisme $g\colon X'\rightarrow X$ projectif \eqref{mtfla2a}. Alors, 
\begin{itemize}
\item[{\rm (i)}] Pour tous entiers $m\geq 1$, $q\geq 0$ et tout $(\co_\oK/p^m\co_\oK)$-module localement libre de type fini $\mL'$ de $\oX'^\rhd_\fet$, 
le $\uoR$-module $(\rR^q\Theta_{m*}(\beta'^*_m(\mL')))_{\rho(\oy\rightsquigarrow \ox)}$ est de présentation $\alpha$-finie, 
où $\beta'_m$ et $\Theta_m$ sont les morphismes de topos annelés \eqref{mtfla10i} et \eqref{mtfla12d}.

\item[{\rm (ii)}] Il existe un entier $q_0\geq 0$ tel que pour tous entiers 
$m\geq 1$, $q\geq q_0$ et tout $(\co_\oK/p^m\co_\oK)$-module localement libre de type fini $\mL'$ de $\oX'^\rhd_\fet$, 
le $\uoR$-module $(\rR^q\Theta_{m*}(\beta'^*_m(\mL')))_{\rho(\oy\rightsquigarrow \ox)}$ soit $\alpha$-nul.
\end{itemize}
\end{teo}

Le reste de cette section est consacré à la preuve de ce théorème qui sera donnée dans \ref{AFR29}.
On observera qu'il est loisible de remplacer $X$ par un voisinage étale de $\ox$ (\cite{agt} VI.10.14).

\subsection{}\label{AFR10}
Dans la suite de cette section, on suppose que le morphisme $f\colon (X,\cM_X)\rightarrow (S,\cM_S)$ admet une carte adéquate 
$((P,\gamma),(\mN,\iota),\vartheta \colon \mN\rightarrow P)$ (cf. \ref{cad1}) et on reprend les notations de \ref{eccr20} et \ref{eccr2},
dont la construction est faisable sous cette hypothèse qui est plus faible que celle fixée dans \ref{eccr1}.
On fixe, de plus, un $S$-morphisme
\begin{equation}\label{AFR10a}
\oS\rightarrow \underset{\underset{n\geq 1}{\longleftarrow}}{\lim}\ S^{(n)}.
\end{equation}
Pour tout entier $n\geq 1$, on pose
\begin{eqnarray}
\oX^{(n)}&=& X^{(n)}\times_{S^{(n)}}\oS,\label{AFR10b}\\
\oX^{(n)\circ}&=&\oX^{(n)}\times_XX^\circ.\label{AFR10c}
\end{eqnarray}
Pour tous entiers $m,n\geq 1$ tels que $m$ divise $n$, on note $\xi^{(n,m)}\colon \oX^{(n)}\rightarrow \oX^{(m)}$ le morphisme canonique
qui est fini et surjectif. On pose $\xi^{(n)}=\xi^{(n,1)}\colon \oX^{(n)}\rightarrow \oX$, de sorte qu'on a 
\begin{equation}\label{AFR10d}
\xi^{(n)}=\xi^{(m)}\circ \xi^{(n,m)}.
\end{equation}
Considérons le $\oX$-schéma 
\begin{equation}\label{AFR10e}
\oX^{(\infty)}= \underset{\underset{n\geq 1}{\longleftarrow}}{\lim}\ \oX^{(n)},
\end{equation}
et la $\co_{\oX}$-algèbre quasi-cohérente 
\begin{equation}\label{AFR10f}
\cR_\infty= \underset{\underset{n\geq 1}{\longrightarrow}}{\lim}\ \xi^{(n)}_*(\co_{\oX^{(n)}}),
\end{equation}
de sorte que $\oX^{(\infty)}=\Spec_{\oX}(\cR_\infty)$ (\cite{ega4} 8.2.3). 

\begin{lem}\label{AFR25}
\
\begin{itemize}
\item[{\rm (i)}] Le schéma $\oX^{(\infty)}$ est normal et localement irréductible.
\item[{\rm (ii)}] Pour tout ouvert affine $\oU$ de $\oX$, la $\co_\oK$-algèbre 
$\Gamma(\oU,\cR_\infty)$ est universellement $\alpha$-cohérente \eqref{afini14}. 
\item[{\rm (iii)}] La $\co_\oK$-algèbre $\co_{\oX^{(\infty)}}$ de $\oX^{(\infty)}_\zar$ est $\alpha$-cohérente. 
\item[{\rm (iv)}]  Pour tout entier $g\geq 0$, 
la $\co_\oK$-algèbre $\cR_\infty[t_1,\dots,t_g]$ (des polynômes en $g$ variables à coefficients dans $\cR_\infty$) de $\oX_\zar$ 
est $\alpha$-cohérente \eqref{finita9}. 
\end{itemize}
\end{lem}

(i) En effet, la question étant locale, on peut supposer $X$ affine. 
Pour tout entier $n\geq 1$, le schéma $\oX^{(n)}$ est normal et localement irréductible et l'immersion $\oX^{(n)\circ}\rightarrow \oX^{(n)}$
est schématiquement dominante (\cite{agt} III.4.2). Pour tous entiers $m,n\geq 1$, le morphisme $\oX^{(mn)\circ}\rightarrow \oX^{(n)\circ}$ est étale fini et surjectif. 
Par suite, toute composante irréductible $\oX^{(mn)}$ domine une composante irréductible de $\oX^{(n)}$. 
En vertu de \ref{cad8}, il existe un entier $n\geq 1$ tel que pour toute composante 
irréductible $C$ de $\oX^{(n)}$ et tout entier $m\geq 1$, $C\times_{\oX^{(n)}}\oX^{(mn)}$ soit irréductible. 
On en déduit que le schéma $\oX^{(\infty)}$ est normal et localement irréductible en vertu de (\cite{ega1n} 0.6.1.6 et 0.6.5.12(ii)).

(ii) Cela résulte de \ref{cad15}, compte tenu de la preuve de (i) et (\cite{sga4} VI 5.2). 

(iii) \& (iv) Cela résulte de (ii) et \ref{afini103}.

\subsection{}\label{AFR24}
Reprenons les notations de \ref{AFR9}. Pour tout entier $n\geq 1$, on pose
\begin{eqnarray}
\uoX^{(n)}&=&\oX^{(n)}\times_X\uX,\label{AFR24a}\\
\uoX^{(n)\circ}&=&\uoX^{(n)}\times_XX^\circ.\label{AFR24b}
\end{eqnarray}
Pour tous entiers $m,n\geq 1$ tels que $m$ divise $n$, on désigne par $\uxi^{(n,m)}\colon \uoX^{(n)}\rightarrow \uoX^{(m)}$ 
et $\uxi^{(n)}\colon \uoX^{(n)}\rightarrow \uoX$ les morphismes induits par $\xi^{(n,m)}$ et $\xi^{(n)}$. 
D'après (\cite{ega4} 8.3.8(i)), il existe un $\uoX$-morphisme 
\begin{equation}\label{AFR24e}
\oy\rightarrow \underset{\underset{n\geq 1}{\longleftarrow}}{\lim}\ \uoX^{(n)},
\end{equation} 
que l'on fixe dans la suite de cette section. Comme le schéma $\uoX$ est strictement local (\cite{agt} III.3.7), pour tout 
$n\geq 1$, le schéma $\uoX^{(n)}$ est une union disjointe finie de schémas normaux et strictement locaux (\ref{cad7}(iii)). 
On note $\uoX^{(n)\star}$ la composante irréductible de $\uoX^{(n)}$ contenant l'image de $\oy$. On pose
\begin{eqnarray}
\uR_n&=&\Gamma(\uoX^{(n)\star},\co_{\uoX^{(n)}}), \label{AFR24c}\\
\uR_\infty&=&\underset{\underset{n\geq 1}{\longrightarrow}}{\lim}\ \uR_n.\label{AFR24d}
\end{eqnarray}
Les anneaux $\uR_n$ sont donc strictement henséliens et normaux, et il en de même de l'anneau $\uR_\infty$ 
(\cite{raynaud1} I §3 prop.~1 et \cite{ega1n} 0.6.5.12(ii)).  

On désigne par $\fV_\ox$ la catégorie des $X$-schémas étales $\ox$-pointés $(U,\fp\colon \ox\rightarrow U)$ tels que le schéma $U$ soit affine. 
Pour tout objet $(U,\fp\colon \ox\rightarrow U)$ de $\fV_\ox$, 
on note encore $\fp\colon \uX\rightarrow U$ le $X$-morphisme déduit de $\fp$ (\cite{sga4} VIII 7.3), et on pose, pour tout entier $n\geq 1$, 
\begin{equation}\label{AFR24i}
\oU^{(n)}=\oX^{(n)}\times_XU.
\end{equation}
Le morphisme $\fp$ induit alors des morphismes $\fp^{(n)}_U\colon \uoX^{(n)}\rightarrow \oU^{(n)}$. 
Le morphisme \eqref{AFR24e} induit un morphisme 
\begin{equation}\label{AFR24f}
\oy\rightarrow \underset{\underset{n\geq 1}{\longleftarrow}}{\lim}\ \oU^{(n)}.
\end{equation} 
Notant $\oR^\oy_U$ la $\Gamma(\oU,\co_\oU)$-algèbre définie dans \eqref{TFA9c}, on en déduit, compte tenu de \ref{cad7}(v), un $\oU$-morphisme 
\begin{equation}\label{AFR24g}
\Spec(\oR^\oy_U)\rightarrow \underset{\underset{n\geq 1}{\longleftarrow}}{\lim}\ \oU^{(n)}.
\end{equation} 
Ces morphismes sont clairement compatibles. Compte tenu de la définition de l'algèbre $\uoR$ \eqref{TFA12f},
on obtient par passage à la limite projective sur la catégorie $\fV_\ox$ un $\uoX$-morphisme 
\begin{equation}\label{AFR24j}
\Spec(\uoR)\rightarrow \underset{\underset{n\geq 1}{\longleftarrow}}{\lim}\ \uoX^{(n)}.
\end{equation}
Le schéma $\uoR$ étant normal et strictement local (en particulier intègre) d'après (\cite{agt} III.10.10(i)), on en déduit un $\uoX$-morphisme 
\begin{equation}\label{AFR24k}
\Spec(\uoR)\rightarrow \underset{\underset{n\geq 1}{\longleftarrow}}{\lim}\ \uoX^{(n)\star},
\end{equation}
et par suite un homomorphisme de $\uR_1$-algèbres
\begin{equation}\label{AFR24h}
\uR_\infty\rightarrow \uoR.
\end{equation}

Pour tout objet $(U,\fp)$ de $\fV_\ox$ et tout entier $n\geq 1$, 
le schéma $\oU^{(n)}$ étant localement irréductible d'après \ref{cad7}(iii),    
il est la somme des schémas induits sur ses composantes irréductibles. On désigne par $\oU^{(n)\star}$  la composante irréductible de $\oU^{(n)}$ 
contenant l'image de $\oy$ \eqref{AFR24f}. On pose
\begin{eqnarray}
R_{U,n}&=&\Gamma(\oU^{(n)\star},\co_{\oX^{(n)}}), \label{AFR24l}\\
R_{U,\infty}&=&\underset{\underset{n\geq 1}{\longrightarrow}}{\lim}\ R_{U,n}.\label{AFR24m}
\end{eqnarray}
On a un isomorphisme canonique 
\begin{equation} \label{AFR24n}
\underset{\underset{(U,\fp)\in \fV_\ox^\circ}{\longrightarrow}}{\lim}\ R_{U,\infty}\stackrel{\sim}{\rightarrow}\uR_\infty. 
\end{equation}

\begin{prop}\label{AFR30}
Les $\co_{\oK}$-algèbres $\uR_\infty$ et $\uoR$ sont universellement $\alpha$-cohérentes \eqref{afini14}.
\end{prop}

Conservons les notations de \ref{AFR24}. Pour tout objet $(U,\fp)$ de $\fV_\ox$, 
la $\co_{\oK}$-algèbre $R_{U,\infty}$ est universellement $\alpha$-cohérente d'après \ref{cad15}. 
La $\co_{\oK}$-algèbre $\uR_\infty$ est donc universellement $\alpha$-cohérente en vertu de \ref{aet10} et \eqref{AFR24n}. 

Par ailleurs, compte tenu de \eqref{TFA12f}, l'homomorphisme canonique
\begin{equation}\label{AFR30d}
\underset{\underset{(U,\fp)\in \fV_\ox^\circ}{\longrightarrow}}{\lim}\ \oR^{\oy}_U\otimes_{R_{U,\infty}}\uR_\infty\rightarrow \uoR 
\end{equation}
est un isomorphisme. Il résulte alors de \ref{tpcg5}  et (\cite{agt} V.7.4) que $\uoR$ est une limite inductive filtrante de $\uR_\infty$-algèbres
$\alpha$-finies-étales. Par suite, la $\co_{\oK}$-algèbre $\uoR$ est universellement $\alpha$-cohérente en vertu de \ref{aet9}.

\begin{prop}\label{AFR31}
Il existe un entier $N\geq 1$ qui ne dépend que de la carte adéquate pour $f$ fixée dans \ref{AFR10}, 
tel que pour toute suite exacte de $\co_{X,\ox}$-modules $M'\rightarrow M\rightarrow M''$, la suite 
\begin{equation}
M'\otimes_{\co_{X,\ox}}\uoR\rightarrow M\otimes_{\co_{X,\ox}}\uoR\rightarrow M''\otimes_{\co_{X,\ox}}\uoR
\end{equation}
soit $p^N$-exacte \eqref{finita3}.
\end{prop}
Cela résulte de \ref{tpcg22} et de la définition de $\uoR$ \eqref{TFA12f}. 
On notera dans \ref{tpcg22} que l'entier $N$ ne dépend que de la carte adéquate pour $f$. 

\begin{cor}\label{AFR32}
Il existe un entier $N\geq 1$ qui ne dépend que de la carte adéquate pour $f$ fixée dans \ref{AFR10},
tel que pour tout complexe de $\co_{X,\ox}$-modules $M^\bullet$
et tout entier $q$, le noyau et le conoyau du morphisme canonique
\begin{equation}
\rH^q(M^\bullet)\otimes_{\co_{X,\ox}}\uoR\rightarrow \rH^q(M^\bullet\otimes_{\co_{X,\ox}}\uoR)
\end{equation}
soient annulés par $p^N$.
\end{cor}

En effet, soit $N$ l'entier fourni par \ref{AFR31} et soit $M^\bullet$ un complexe de $\co_{X,\ox}$-modules. 
Pour tout entier $n$, notons $Z^n$ (resp. $\uoZ^n$) le noyau du morphisme $M^n\rightarrow M^{n+1}$ 
(resp. $M^n\otimes_{\co_{X,\ox}}\uoR\rightarrow M^{n+1}\otimes_{\co_{X,\ox}}\uoR$) 
et  $B^n$ (resp. $\uoB^n$) l'image du morphisme $M^{n-1}\rightarrow M^n$ (resp. $M^{n-1}\otimes_{\co_{X,\ox}}\uoR\rightarrow M^n\otimes_{\co_{X,\ox}}\uoR$).
D'après \ref{AFR31}, le noyau et le conoyau du morphisme canonique $Z^n\otimes_{\co_{X,\ox}}\uoR\rightarrow \uoZ^n$ sont annulés par $p^N$.
Par ailleurs, le morphisme canonique $B^n\otimes_{\co_{X,\ox}}\uoR\rightarrow \uoB^n$ est surjectif et son noyau est annulé par $p^N$. 
Notons $\uoC^n$ l'image du morphisme canonique $B^n\otimes_{\co_{X,\ox}}\uoR\rightarrow Z^n\otimes_{\co_{X,\ox}}\uoR$. On a alors un diagramme commutatif 
canonique à lignes exactes 
\begin{equation}
\xymatrix{
0\ar[r]&\uoC^n\ar[r]\ar@{->>}[d]&{Z^n\otimes_{\co_{X,\ox}}\uoR}\ar[r]\ar[d]&{\rH^n(M^\bullet)\otimes_{\co_{X,\ox}}\uoR}\ar[d]\ar[r]&0\\
0\ar[r]&\uoB^n\ar[r]&\uoZ^n\ar[r]&{\rH^n(M^\bullet\otimes_{\co_{X,\ox}}\uoR)}\ar[r]&0}
\end{equation}
La proposition s'ensuit puisque la  flèche verticale de gauche est surjective.

\subsection{}\label{decprim1}
Le schéma $\uoX$ étant strictement local (\cite{agt} III.3.7), il s'identifie au localisé strict de $\oX$ en $\ox$. 
Pour tout entier $n\geq 1$, le schéma $\uoX^{(n)}$ est donc une union disjointe finie de schémas {\em normaux} et strictement locaux (\ref{cad7}(iii)),  
à savoir les localisés stricts de $\oX^{(n)}$ en les points de la fibre $\oX^{(n)}_\ox$ de $\oX^{(n)}$ au-dessus de $\ox$,
\begin{equation}\label{decprim1a}
\uoX^{(n)}=\bigsqcup_{\oz\in \oX^{(n)}_\ox}\oX^{(n)}_{(\oz)}.
\end{equation}
Pour tout $\oz\in \oX^{(n)}_\ox$, on pose 
\begin{equation}\label{decprim1b}
\uR_n^{(\oz)}=\Gamma(\uoX^{(n)}_{(\oz)},\co_{\uoX^{(n)}}).
\end{equation}
On note $\ox_n$ le point de $\oX^{(n)}_\ox$ tel que $\oX^{(n)}_{(\ox_n)}$ soit la composante $\uoX^{(n)\star}$ de $\uoX^{(n)}$ contenant
l'image de $\oy$ \eqref{AFR24}, de sorte que $\uR_n^{(\ox_n)}=\uR_n$ \eqref{AFR24c}. 

En vertu de \ref{cad7}(v), pour tout entier $m\geq 1$, le morphisme canonique $\uxi^{(mn,n)}\colon \uoX^{(mn)}\rightarrow \uoX^{(n)}$  
induit un morphisme étale, fini et surjectif $\uoX^{(mn)\circ}\rightarrow \uoX^{(n)\circ}$ qui fait de $\uoX^{(mn)\circ}$ un espace principal homogène 
pour la topologie étale de $\uoX^{(n)\circ}$ sous le groupe $\Hom_\mZ(P^\gp,\mu_m(\oK))$. 
Par ailleurs, l'ouvert $\uoX^{(mn)\circ}$ de $\uoX^{(mn)}$ est schématiquement dense d'après \ref{cad7}(i)-(iv) et (\cite{agt} III.4.2(iv)). 
On en déduit une action naturelle de $\Hom_\mZ(P^\gp,\mu_m(\oK))$ sur $\uoX^{(mn)}$ par des $\uoX^{(n)}$-automorphismes. 
L'action induite sur la fibre $\oX^{(mn)}_{\ox_n}$ est transitive. 
Comme on a $\ox_n=\uxi^{(mn,n)}(\ox_{mn})$, 
pour tout $\oz\in \oX^{(mn)}_{\ox_n}$, il existe $\iota\in \Hom_\mZ(P^\gp,\mu_m(\oK))$ induisant un $\uR_n$-isomorphisme 
\begin{equation}\label{decprim1c}
\iota^*\colon \uR_{mn}^{(\oz)}\stackrel{\sim}{\rightarrow} \uR_{mn}.
\end{equation}

\begin{prop}\label{decprim2}
Conservons les notations de \ref{decprim1}. 
Il existe un entier $N\geq 1$ qui ne dépend que de la carte adéquate pour $f$ fixée dans \ref{AFR10}, 
tel que pour tous entiers $m,n\geq 1$, on ait un morphisme $\uR_n$-linéaire injectif 
\begin{equation}\label{decprim2a}
\bigoplus_{\unu\in (\mN\cap [0,m-1])^{d+1}}\uR_n\rightarrow \bigoplus_{\oz\in \oX^{(mn)}_{\ox_n}}\uR_{mn}^{(\oz)}
\end{equation}
dont le conoyau est annulé par $\pi_n^N$ \eqref{eccr20} et dont la composante correspondant à $\unu=(0,\dots,0)$
et $\oz\in \oX^{(mn)}_{\ox_n}$ est l'homomorphisme canonique d'anneaux $\uR_n\rightarrow \uR_{mn}^{(\oz)}$. 
\end{prop}

Cela résulte aussitôt de \ref{cad100} par changement de base par le morphisme $\oX^{(n)}_{(\ox_n)}\rightarrow X^{(n)}$.

\begin{cor}\label{decprim3}
Il existe un entier $N\geq 1$ qui ne dépend que de la carte adéquate pour $f$ fixée dans \ref{AFR10}, 
tel que pour tous entiers $m,n\geq 1$ et tout $\uR_n$-module $M_n$, le noyau du morphisme canonique
\begin{equation}\label{decprim3a}
M_n\rightarrow M_n\otimes_{\uR_n}\uR_{mn}
\end{equation}
soit annulé par $\pi_n^N$. 
\end{cor}

Montrons que l'entier $N$ fourni par la proposition \ref{decprim2} convient. En effet, d'après {\em loc. cit.}, le noyau du morphisme  
\begin{equation}\label{decprim3b}
\bigoplus_{\unu\in (\mN\cap [0,m-1])^{d+1}}M_n\rightarrow \bigoplus_{\oz\in \oX^{(mn)}_{\ox_n}}M_n\otimes_{\uR_n}\uR_{mn}^{(\oz)}
\end{equation}
induit par \eqref{decprim2a} est annulé par $\pi_n^N$. On notera que le morphisme canonique \eqref{decprim3a} est la composante de \eqref{decprim3b} 
correspondant à $\unu=(0,\dots,0)$ et $\oz=\ox_{mn}$. Compte tenu de \eqref{decprim1c}, le noyau du morphisme \eqref{decprim3a} coïncide donc avec 
celui du morphisme canonique 
\begin{equation}\label{decprim3c}
M_n\rightarrow \bigoplus_{\oz\in \oX^{(mn)}_{\ox_n}}M_n\otimes_{\uR_n}\uR_{mn}^{(\oz)}
\end{equation}
qui est bien annulé par $\pi_n^N$ d'après ce qui précède. 

\subsection{}\label{decprim4}
Pour tout entier $n\geq 1$ et tout $\uR_n$-module $M$, on désigne par $\Gamma_{\ox_n}(M)$ le sous-$\uR_n$-module de $M$ des sections à support dans $\ox_n$. 
On notera qu'il existe un idéal de type fini $I$ de $\uR_n$ tel que $\rV(I)=\{\ox_n\}$. En effet, l'anneau local $\uR_n\otimes_{\co_\oK}k$ s'identifie au localisé
strict de $\oX^{(n)}_s$ \eqref{mtfla1}; il est donc noethérien. On peut alors prendre pour $I$ l'idéal de $\uR_n$ engendré par $(\gamma,\lambda_1,\dots,\lambda_q)$,
où $\gamma$ est un élément non nul de $\fm_\oK$ et $\lambda_1,\dots,\lambda_q$ relèvent des générateurs de l'idéal maximal de $\uR_n\otimes_{\co_\oK}k$. 
Par suite, l'ouvert $\Spec(\uR_n)-\{\ox_n\}$ est quasi-compact. 

Le schéma $\Spec(\uR_\infty)$ étant strictement local, on désigne par $\ox_\infty$ son point fermé. 
Pour tout $\uR_\infty$-module $M$, on désigne par $\Gamma_{\ox_\infty}(M)$ le sous-$\uR_\infty$-module de $M$ des sections à support dans $\ox_\infty$. 

Comme $\ox_\infty$ est l'unique point de $\Spec(\uR_\infty)$ au-dessus du point $\ox_n\in \oX^{(n)}_{(\ox_n)}$, 
pour tout $\uR_\infty$-module $M$ que l'on considère aussi comme un $\uR_n$-module, on a  $\Gamma_{\ox_\infty}(M)=\Gamma_{\ox_n}(M)$. 
De même, pour tout $m\geq 1$, $\ox_{mn}$ étant l'unique point de $\oX^{(mn)}_{(\ox_{mn})}$ au-dessus du point $\ox_n\in \oX^{(n)}_{(\ox_n)}$, 
pour tout $\uR_{mn}$-module $M$ que l'on considère aussi comme un $\uR_n$-module, on a  $\Gamma_{\ox_{mn}}(M)=\Gamma_{\ox_n}(M)$. 

Pour tout entier $n\geq 1$ et tout $\uR_n$-module $M_n$, le morphisme canonique
\begin{equation}\label{decprim4a}
\underset{\underset{m\geq 1}{\longrightarrow}}{\lim}\ \Gamma_{\ox_{mn}}(M_n\otimes_{\uR_n}\uR_{mn})\rightarrow 
\Gamma_{\ox_{\infty}}(M_n\otimes_{\uR_n}\uR_\infty)
\end{equation}
est un isomorphisme. En effet, comme $U_n=\Spec(\uR_n)-\{\ox_n\}$ est quasi-compact, le morphisme canonique
\begin{equation}\label{decprim4b}
\underset{\underset{m\geq 1}{\longrightarrow}}{\lim}\ \Gamma(U_n\otimes_{\uR_n}\uR_{mn},(M_n\otimes_{\uR_n}\uR_{mn})^\sim)\rightarrow 
\Gamma(U_n\otimes_{\uR_n}\uR_\infty,(M_n\otimes_{\uR_n}\uR_\infty)^\sim),
\end{equation}
où $(-)^\sim$ désigne les faisceaux associés, est un isomorphisme.

\begin{prop}\label{decprim5}
Il existe un entier $N\geq 1$ qui ne dépend que de la carte adéquate pour $f$ fixée dans \ref{AFR10}, 
tel que pour tous entiers $m,n\geq 1$ et tout $\uR_n$-module $M_n$, le noyau et le conoyau du morphisme $\uR_{mn}$-linéaire canonique
\begin{equation}\label{decprim5a}
\Gamma_{\ox_n}(M_n)\otimes_{\uR_n}\uR_{mn}\rightarrow \Gamma_{\ox_{mn}}(M_n\otimes_{\uR_n}\uR_{mn})
\end{equation}
soient annulés par $\pi_n^N$.
\end{prop}

On notera d'abord que pour tout $\oz\in \oX^{(mn)}_{\ox_n}$ et tout $\uR^{(\oz)}_{mn}$-module $F$ que l'on considère aussi comme un $\uR_n$-module, 
on a  $\Gamma_{\oz}(F)=\Gamma_{\ox_n}(F)$
puisque $\oz$ est le seul point de $\oX_{(\oz)}^{(mn)}$ au-dessus du point $\ox_n\in \oX_{(\ox_{n})}^{(n)}$. 
On a donc un morphisme $\uR^{(\oz)}_{mn}$-linéaire canonique
\begin{equation}\label{decprim5b}
v_\oz\colon \Gamma_{\ox_n}(M_n)\otimes_{\uR_n}\uR^{(\oz)}_{mn}\rightarrow \Gamma_{\oz}(M_n\otimes_{\uR_n}\uR^{(\oz)}_{mn}).
\end{equation}
Le morphisme \eqref{decprim5a} n'est autre que $v_{\ox_{mn}}$.

Montrons que si $N$ est l'entier fourni par la proposition \ref{decprim2}, l'entier $2N$ convient.
On a un diagramme commutatif
\begin{equation}\label{decprim5c}
\xymatrix{
{\bigoplus_{\unu\in (\mN\cap [0,m-1])^{d+1}}\Gamma_{\ox_n}(M_n)}\ar[r]^-(0.5){h_1}\ar@{=}[d]&
{\bigoplus_{\oz\in \oX^{(mn)}_{\ox_n}}\Gamma_{\ox_n}(M_n)\otimes_{\uR_n}\uR_{mn}^{(\oz)}}\ar[d]^{v=\bigoplus_{\oz\in \oX^{(mn)}_{\ox_n}}v_\oz}\\
{\bigoplus_{\unu\in (\mN\cap [0,m-1])^{d+1}}\Gamma_{\ox_n}(M_n)}\ar[r]^-(0.5){h_2}&
{\bigoplus_{\oz\in \oX^{(mn)}_{\ox_n}}\Gamma_{\oz}(M_n\otimes_{\uR_n}\uR_{mn}^{(\oz)})}}
\end{equation}
où les flèches horizontales sont induites par \eqref{decprim2a}.  
Il résulte de \ref{decprim2} que le noyau et le conoyau de $h_1$ et le noyau de $h_2$ sont annulés par $\pi_n^N$. 
Montrons que le conoyau de $h_2$ est annulé par $\pi_n^{2N}$. Soit 
$\zeta$ une section de $\bigoplus_{\oz\in \oX^{(mn)}_{\ox_n}}\Gamma_{\oz}(M_n\otimes_{\uR_n}\uR_{mn}^{(\oz)})$. 
D'après \ref{decprim2}, $\pi_n^N\zeta$ provient d'une section $\zeta'$ de $\bigoplus_{\unu\in (\mN\cap [0,d-1])^{d+1}}M_n$. 
Pour tout ouvert affine $U=\Spec(A)$ de $\oX^{(n)}_{\ox_n}-\{\ox_n\}$, le noyau du morphisme 
\begin{equation}
\bigoplus_{\unu\in (\mN\cap [0,m-1])^{d+1}}M_n\otimes_{\uR_n}A\rightarrow \bigoplus_{\oz\in \oX^{(mn)}_{\ox_n}}(M_n\otimes_{\uR_n}A)\otimes_{\uR_n}\uR_{mn}^{(\oz)}
\end{equation}
induit par \eqref{decprim2a}, est annulé par $\pi_n^N$. Par suite, $\pi_n^N\zeta'$ est à support dans $\ox_n$.  Le conoyau de $h_2$ est donc annulé 
par $\pi_n^{2N}$. La proposition s'ensuit par une chasse au diagramme \eqref{decprim5c}.

\begin{cor}\label{decprim7}
Il existe un entier $N\geq 1$ qui ne dépend que de la carte adéquate pour $f$ fixée dans \ref{AFR10}, 
tel que pour tout entier $n\geq 1$ et tout $\uR_n$-module $M_n$, le noyau et le conoyau du morphisme $\uR_\infty$-linéaire canonique
\begin{equation}\label{decprim7a}
\Gamma_{\ox_n}(M_n)\otimes_{\uR_n}\uR_\infty\rightarrow \Gamma_{\ox_\infty}(M_n\otimes_{\uR_n}\uR_\infty)
\end{equation}
soient annulés par $\pi_n^N$.
\end{cor}

Cela résulte de \ref{decprim5} par passage à la limite inductive en tenant compte de l'isomorphisme \eqref{decprim4a}.

\begin{prop}\label{decprim8}
Pour tout entier $n\geq 1$ et tout $\uR_n$-module de présentation finie $M$, le $\uR_n$-module $\Gamma_{\ox_n}(M)$ est de présentation finie. 
\end{prop}

Pour toute extension finie $L$ de $K_n$ contenue dans $\oK$, on désigne par $\uR_{L,n}$ l'anneau du localisé strict de 
$X^{(n)}\times_{S^{(n)}}\Spec(\co_L)$ en $\ox_n$. Pour toutes extensions finies $L\subset L'$ de $K_n$ contenues dans $\oK$, 
comme le morphisme $X^{(n)}\times_{S^{(n)}}\Spec(\co_{L'})\rightarrow X^{(n)}\times_{S^{(n)}}\Spec(\co_L)$ induit un isomorphisme entre
les fibres spéciales, l'homomorphisme canonique $\uR_{L,n}\otimes_{\co_L}\co_{L'}\rightarrow \uR_{L',n}$ est un isomorphisme 
d'après (\cite{ega4} 18.8.10). Passant à la limite inductive sur les extensions $L'$ de $L$ contenues dans $\oK$, on en déduit que le 
morphisme canonique $\uR_{L,n}\otimes_{\co_L}\co_{\oK}\rightarrow \uR_n$ est un isomorphisme (\cite{ega4} 18.8.18(ii)). 
Notons $\fm_{L,n}$ (resp. $\fm_n$) l'idéal maximal de $\uR_{L,n}$ (resp. $\uR_n$). 
On en déduit que l'homomorphisme canonique $\uR_{L,n}\rightarrow \uR_n$ est plat et que $\fm_n$ est le radical de l'idéal $\fm_{L,n}\uR_n$. 

Soient $L$ une extension finie $L$ de $K_n$ contenue dans $\oK$, $M_L$ un $\uR_{L,n}$-module de présentation finie tels qu'il existe
un isomorphisme $\uR_n$-linéaire $M_L\otimes_{\uR_{L,n}}\uR_n\stackrel{\sim}{\rightarrow}M$. Soient $h_1,\dots, h_e$ des générateurs de 
$\fm_{L,n}$. La suite de morphismes canoniques
\begin{equation}
0\rightarrow \Gamma_{\fm_{L,n}}(M_L)\rightarrow M_L\rightarrow \oplus_{1\leq i\leq e}(M_L)_{h_i}
\end{equation}
est alors exacte. Comme $\rV((h_1,\dots,h_e)\uR_n)=\{\fm_n\}$ et que $\uR_n$ est $\uR_{L,n}$-plat, on en déduit un isomorphisme 
\begin{equation}
\Gamma_{\fm_{L,n}}(M_L)\otimes_{\uR_{L,n}}\uR_n\stackrel{\sim}{\rightarrow}\Gamma_{\fm_n}(M).
\end{equation}
La proposition s'ensuit puisque l'anneau $\uR_{L,n}$ est noethérien. 

\begin{cor}\label{decprim9}
Pour tout $\uR_\infty$-module de présentation $\alpha$-finie $M$ \eqref{finita2}, le $\uR_\infty$-module $\Gamma_{\ox_\infty}(M)$ est de présentation $\alpha$-finie.
\end{cor}

En effet, pour tout $\gamma\in \fm_\oK$, il existe un $\uR_\infty$-module de présentation finie $M'$ et un $\gamma$-isomorphisme $u\colon M'\rightarrow M$ \eqref{finita3}. 
Il existe donc un morphisme $\uR_\infty$-linéaire $v\colon M\rightarrow M'$ tel que $u\circ v=\gamma^2\id_M$ et $v\circ u=\gamma^2\id_{M'}$ \eqref{alpha3}.
On en déduit que le morphisme $\Gamma_{\ox_\infty}(M')\rightarrow  \Gamma_{\ox_\infty}(M)$ est un $\gamma^2$-isomorphisme. 
On peut donc se réduire au cas où $M$ est un $\uR_\infty$-module de présentation finie. Il existe alors un entier $n\geq 1$, un $\uR_n$-module de présentation finie
$M_n$ et un isomorphisme $M_n\otimes_{\uR_n}\uR_\infty\stackrel{\sim}{\rightarrow} M$. 
En vertu de \ref{decprim7},  le morphisme $\uR_\infty$-linéaire canonique
\begin{equation}
\Gamma_{\ox_n}(M_n)\otimes_{\uR_n}\uR_\infty\rightarrow \Gamma_{\ox_\infty}(M)
\end{equation}
est un $\pi_n^N$-isomorphisme, où l'entier $N$ ne dépend que de la carte adéquate pour $f$ fixée dans \ref{AFR10}.
La proposition résulte alors de \ref{decprim8} en choisissant $n$ assez grand.

\subsection{}\label{AFR11}
Pour tout entier $n\geq 1$, 
le morphisme $\xi^{(n)\circ}\colon \oX^{(n)\circ}\rightarrow \oX^\circ$ est étale et fini d'après \ref{cad7}(v), 
de sorte que $(\oX^{(n)\circ}\rightarrow X)$ est un objet de $E$. 
On note $E^{(n)}\rightarrow \Et_{/X}$ (resp. $\tE^{(n)}$) le site fibré (resp. topos) de Faltings associé 
au morphisme $\oX^{(n)\circ}\rightarrow X$ \eqref{tf1}. Tout objet de $E^{(n)}$
est naturellement un objet de $E$. On définit ainsi un foncteur 
\begin{equation}\label{AFR11a}
\Phi^{(n)}\colon E^{(n)}\rightarrow E
\end{equation}
qui se factorise à travers une équivalence de catégories 
\begin{equation}\label{AFR11b}
E^{(n)}\stackrel{\sim}{\rightarrow} E_{/(\oX^{(n)\circ}\rightarrow X)}.
\end{equation}
Le foncteur $\Phi^{(n)}$ est continu et cocontinu (cf. \ref{amtF1}). 
Il définit donc une suite de trois foncteurs adjoints~:
\begin{equation}\label{AFR11c}
(\Phi^{(n)})_!\colon \tE^{(n)}\rightarrow \tE, \ \ \ \Phi^{(n)*}\colon \tE\rightarrow \tE^{(n)}, 
\ \ \ \Phi^{(n)}_*\colon \tE^{(n)}\rightarrow \tE,
\end{equation}
dans le sens que pour deux foncteurs consécutifs de la suite, celui de droite est
adjoint à droite de l'autre. Le foncteur $(\Phi^{(n)})_!$ se factorise à travers 
une équivalence de catégories 
\begin{equation}\label{AFR11d}
\tE^{(n)}\stackrel{\sim}{\rightarrow} \tE_{/(\oX^{(n)\circ}\rightarrow X)^\tta},
\end{equation}
où $(\oX^{(n)\circ}\rightarrow X)^\tta$ est l'objet de $\tE$ associé à $(\oX^{(n)\circ}\rightarrow X)$.
Le couple de foncteurs $(\Phi^{(n)*},\Phi^{(n)}_*)$ définit un morphisme de topos que l'on note aussi
\begin{equation}\label{AFR11e}
\Phi^{(n)}\colon \tE^{(n)}\rightarrow \tE,
\end{equation}
et qui n'est autre que le composé de l'équivalence \eqref{AFR11d} et du morphisme de localisation de $\tE$ en $(\oX^{(n)\circ}\rightarrow X)^\tta$. 
Celui-ci se déduit aussi par fonctorialité du diagramme commutatif (\cite{agt} VI.10.12) 
\begin{equation}\label{AFR11f}
\xymatrix{
{\oX^{(n)\circ}}\ar[r]\ar[d]_{\xi^{(n)\circ}}&X\ar@{=}[d]\\
{\oX^\circ}\ar[r]&X}
\end{equation}

On désigne par 
\begin{equation}\label{AFR11h}
\varphi_\ox^{(n)}\colon \tE^{(n)}\rightarrow \uoX^{(n)\circ}_\fet
\end{equation}
le foncteur canonique défini dans \eqref{tf10b}.

D'après (\cite{sga4} IV 5.5), pour tous entiers $m,n\geq 1$ tels que $m$ divise $n$, on a un morphisme canonique de topos 
\begin{equation}\label{AFR11i}
\Phi^{(n,m)}\colon \tE^{(n)}\rightarrow \tE^{(m)}
\end{equation}
qui s'insère dans un diagramme commutatif à isomorphisme canonique près
\begin{equation}\label{AFR11j}
\xymatrix{
{\tE^{(n)}}\ar[rd]^{\Phi^{(n)}}\ar[d]_{\Phi^{(n,m)}}&\\
{\tE^{(m)}}\ar[r]_{\Phi^{(m)}}&{\tE}}
\end{equation}

\begin{lem}\label{AFR12}
Pour tous entiers $m,n\geq 1$ tels que $m$ divise $n$, le diagramme de foncteurs 
\begin{equation}\label{AFR12a}
\xymatrix{
{\tE^{(n)}}\ar[r]^{\varphi_\ox^{(n)}}&{\uoX^{(n)\circ}_\fet}\\
{\tE^{(m)}}\ar[r]^{\varphi_\ox^{(m)}}\ar[u]^{\Phi^{(n,m)*}}&{\uoX^{(m)\circ}_\fet}\ar[u]_{\uxi^{(n,m)\circ*}}}
\end{equation}
est commutatif à isomorphisme canonique près
\begin{equation}\label{AFR12b}
\varphi_\ox^{(n)}\circ \Phi^{(n,m)*}\stackrel{\sim}{\rightarrow} \uxi^{(n,m)\circ*}\circ \varphi_\ox^{(m)}.
\end{equation}
De plus, ces isomorphismes vérifient une relation de cocycle que nous n'expliciterons pas.
\end{lem}

On désigne par $\tuE^{(n)}$ le topos de Faltings associé au morphisme $\uoX^{(n)\circ}\rightarrow \uX$, par 
\begin{eqnarray}
\Xi^{(n)}\colon \tuE^{(n)}&\rightarrow &\tE^{(n)}\\
\uPhi^{(n,m)}\colon \tuE^{(n)}&\rightarrow &\tuE^{(m)} 
\end{eqnarray}
les morphismes définis par fonctorialité (\cite{agt} VI.10.12), et par
\begin{equation}
\theta^{(n)}\colon \uoX^{(n)\circ}_\fet\rightarrow \tuE^{(n)}
\end{equation}
le morphisme défini dans (\cite{agt} VI.10.23(i)). On a alors $\varphi^{(n)}_\ox=\theta^{(n)*}\circ \Xi^{(n)*}$ \eqref{tf10}. 
Il résulte aussitôt des définitions (\cite{agt} VI.10.22) que les carrés du diagramme 
\begin{equation}
\xymatrix{
{\uoX^{(n)}_\fet}\ar[r]^-(0.5){\theta^{(n)}}\ar[d]_{\uxi^{(n,m)\circ}}&{\tuE^{(n)}}\ar[d]^{\uPhi^{(n,m)}}\ar[r]^{\Xi^{(n)}}&{\tE^{(n)}}\ar[d]^{\Phi^{(n,m)}}\\
{\uoX^{(m)}_\fet}\ar[r]^-(0.5){\theta^{(m)}}&{\tuE^{(m)}}\ar[r]^{\Xi^{(m)}}&{\tE^{(m)}}}
\end{equation}
sont commutatifs à isomorphismes près; d'où la proposition.

\subsection{}\label{AFR13}
Pour tout entier $n\geq 1$, on note 
\begin{equation}\label{AFR13a}
\sigma^{(n)}\colon \tE^{(n)}\rightarrow X_\et
\end{equation}
le morphisme canonique (\cite{agt} (VI.10.6.4)). Il résulte aussitôt des définitions (\cite{agt} VI.10.6) qu'on a un isomorphisme canonique 
\begin{equation}\label{AFR13b}
\sigma^{(n)}\stackrel{\sim}{\rightarrow} \sigma\circ  \Phi^{(n)},
\end{equation}
où $\sigma$ est le morphisme \eqref{mtfla7b}.

Le schéma $\oX^{(n)}$ est normal et localement irréductible et  
l'immersion $\oX^{(n)\circ}\rightarrow \oX^{(n)}$ est schématiquement dominante d'après \ref{cad7} et (\cite{agt} III.4.2(iv)).
Comme le morphisme $\xi^{(n)}\colon \oX^{(n)}\rightarrow \oX$ est fini, on en déduit un isomorphisme canonique 
\begin{equation}\label{AFR13c}
\sigma^{(n)}_*(\Phi^{(n)*}(\ocB))\stackrel{\sim}{\rightarrow} \hbar_*(\xi^{(n)}_*(\co_{\oX^{(n)}})).
\end{equation}
Dans la suite, on considère $\sigma^{(n)}$ comme un morphisme de topos annelés 
\begin{equation}\label{AFR13d}
\sigma^{(n)}\colon (\tE^{(n)},\Phi^{(n)*}(\ocB))\rightarrow (X_\et,\hbar_*(\xi^{(n)}_*(\co_{\oX^{(n)}}))). 
\end{equation}

Pour tout objet $\cF$ de $\tE$, les faisceaux $(\sigma^{(n)}_*(\Phi^{(n)*}(\cF)))_{n\geq 1}$ 
forment naturellement un système inductif de $X_{\et}$ indexé par l'ensemble $\mZ_{\geq 1}$ ordonné par la relation de divisibilité. 
On pose abusivement
\begin{equation}\label{AFR13e}
\sigma^{(\infty)}_*(\cF)=\underset{\underset{n\geq 1}{\longrightarrow}}{\lim}\ \sigma^{(n)}_*(\Phi^{(n)*}(\cF)).
\end{equation}
D'après (\cite{sga4} VII 5.14), pour tout $X$-schéma étale de présentation finie $U$, on a 
\begin{equation}\label{AFR13f}
\sigma^{(\infty)}_*(\cF)(U)=\underset{\underset{n\geq 1}{\longrightarrow}}{\lim}\ \cF(U\times_X\oX^{(n)\circ}\rightarrow U).
\end{equation}

L'algèbre $\sigma^{(\infty)}_*(\ocB)$ de $X_\et$ est associée à l'algèbre quasi-cohérente $\hbar_*(\cR_\infty)$ de $X_\zar$ 
\eqref{AFR10f}, \eqref{notconv12a}. Ceci résulte aussitôt de \eqref{AFR13c} et du fait que les foncteurs \eqref{notconv12a} et $\hbar_*$ 
commutent aux limites inductives filtrantes compte tenu de \eqref{notconv12e} et (\cite{sga4} VI 5.1). 
D'après \eqref{AFR13c} et \eqref{notconv12g}, on a un homomorphisme canonique 
\begin{equation}\label{AFR13k}
\sigma^{(\infty)}_*(\ocB)_\ox\rightarrow \uR_\infty.
\end{equation}
Composant avec l'homomorphisme \eqref{AFR24h}, on obtient un homomorphisme 
\begin{equation}\label{AFR13l}
\sigma^{(\infty)}_*(\ocB)_\ox\rightarrow \uoR.
\end{equation}

Pour tout entier $q\geq 0$ et tout groupe abélien $\cF$ de $\tE$, les faisceaux $(\rR^q\sigma^{(n)}_*(\Phi^{(n)*}(\cF)))_{n\geq 1}$ 
forment naturellement un système inductif de groupes abéliens de $X_{\et}$. On pose abusivement
\begin{equation}\label{AFR13g}
\rR^q\sigma^{(\infty)}_*(\cF)=\underset{\underset{n\geq 1}{\longrightarrow}}{\lim}\ \rR^q \sigma^{(n)}_*(\Phi^{(n)*}(\cF)).
\end{equation}

Soient $\cF$ un groupe abélien de $\tE$, $n\geq 1$ et $q\geq 0$ deux entiers. 
En vertu de (\cite{agt} VI.10.30), on a un isomorphisme  canonique et fonctoriel \eqref{AFR11h}
\begin{equation}\label{AFR13h}
\rR^q\sigma^{(n)}_*(\Phi^{(n)*}(\cF))_\ox\stackrel{\sim}{\rightarrow}\rH^q(\uoX^{(n)\circ}_\fet,\varphi^{(n)}_\ox(\Phi^{(n)*}(\cF))).
\end{equation}
Par passage à la limite inductive \eqref{AFR12}, on en déduit un isomorphisme canonique et fonctoriel
\begin{equation}\label{AFR13i}
\rR^q\sigma^{(\infty)}_*(\cF)_\ox\stackrel{\sim}{\rightarrow}\underset{\underset{n\geq 1}{\longrightarrow}}{\lim}\  
\rH^q(\uoX^{(n)\circ}_\fet,\varphi^{(n)}_\ox(\Phi^{(n)*}(\cF))).
\end{equation}

\begin{prop}\label{AFR15}
Soient $m$ un entier $\geq 1$, $\cF$ un $\ocB_m$-module de $\tE$. Alors, 
\begin{itemize}
\item[{\rm (i)}] On a un morphisme $\uoR$-linéaire canonique et fonctoriel
\begin{equation}\label{AFR15a}
\sigma^{(\infty)}_*(\cF)_\ox\otimes_{\sigma^{(\infty)}_*(\ocB)_\ox}\uoR\rightarrow \cF_{\rho(\oy\rightsquigarrow \ox)},
\end{equation}
où $\rho$ est le morphisme \eqref{mtfla7e}, cf. \eqref{AFR9e} et \eqref{AFR13l}. C'est, de plus, un $\alpha$-isomor\-phisme.
\item[{\rm (ii)}] Pour tout entier $q\geq 1$, le $\co_\oK$-module $\rR^q\sigma^{(\infty)}_*(\cF)_\ox$ est $\alpha$-nul. 
\end{itemize}
\end{prop}

Reprenons les notations de \ref{AFR24} et \ref{AFR11}. 
Pour tout objet $(U,\fp\colon \ox\rightarrow U)$ de $\fV_\ox$, on note 
\begin{equation}
\iota_{U!}\colon \Et_{\rf/\oU^\circ}\rightarrow E
\end{equation} 
le foncteur canonique \eqref{tf1c}, et on pose $\cF_U=\cF\circ \iota_{U!}$ qui est un objet de $\oU^\circ_\fet$. 
D'après (\cite{agt} VI.10.37), pour tout entier $n\geq 1$, on a un isomorphisme canonique et fonctoriel
\begin{equation}\label{AFR15b}
\varphi^{(n)}_\ox(\Phi^{(n)*}(\cF))\stackrel{\sim}{\rightarrow} 
\underset{\underset{(U,\fp)\in \fV_\ox^\circ}{\longrightarrow}}{\lim} \ (\fp_U^{(n)})^{\circ*}(\cF_U|\oU^{(n)\circ}),
\end{equation}
où $\fp_U^{(n)}\colon \uoX^{(n)}\rightarrow \oU^{(n)}$ est le morphisme induit par $\fp$ (cf. \ref{AFR24}).
En vertu de (\cite{agt} VI.11.10), on en déduit, pour tout entier $q\geq 0$, un isomorphisme fonctoriel
\begin{equation}\label{AFR15c}
\rH^q(\uoX^{(n)\circ}_\fet,\varphi^{(n)}_\ox(\Phi^{(n)*}(\cF)))\stackrel{\sim}{\rightarrow} 
\underset{\underset{(U,\fp)\in \fV_\ox^\circ}{\longrightarrow}}{\lim} \ \rH^q(\oU^{(n)\circ}_\fet,\cF_U|\oU^{(n)\circ}).
\end{equation}
Compte tenu de \eqref{AFR13i}, celui-ci induit un isomorphisme fonctoriel
\begin{equation}\label{AFR15d}
\rR^q\sigma^{(\infty)}_*(\cF)_\ox\stackrel{\sim}{\rightarrow} 
\underset{\underset{n\geq 1}{\longrightarrow}}{\lim} \
\underset{\underset{(U,\fp)\in \fV_\ox^\circ}{\longrightarrow}}{\lim} \ \rH^q(\oU^{(n)\circ}_\fet,\cF_U|\oU^{(n)\circ}).
\end{equation}

Soit $(U,\fp)$ un objet de $\fV_\ox$. On pose $\oU=\oU^{(1)}$, $\oU^\star=\oU^{(1)\star}$ (cf. \ref{AFR24}), $\Delta_U=\pi_1(\oU^{\star\circ},\oy)$ et 
\begin{equation}\label{AFR15g}
\Sigma_U=\underset{\underset{n\geq 1}{\longleftarrow}}{\lim} \ \pi_1(\oU^{(n)\star\circ},\oy).
\end{equation}
D'après (\cite{ega4} 8.8.2, 8.10.5 et 17.7.8), on a un isomorphisme canonique
\begin{equation}\label{AFR15k}
\pi_1(\underset{\underset{n\geq 1}{\longleftarrow}}{\lim} \ \oU^{(n)\star\circ},\oy) \stackrel{\sim}{\rightarrow}\Sigma_U.
\end{equation}
On désigne par
\begin{equation}\label{AFR15h}
\nu_{U,\oy}\colon \oU^{\star\circ}_\fet \stackrel{\sim}{\rightarrow}\bB_{\Delta_U}
\end{equation}
le foncteur fibre de $\oU^{\star \circ}_\fet$ en $\oy$ \eqref{notconv11c}. On rappelle qu'on a $\oR^\oy_U=\nu_{U,\oy}(\ocB_U|\oU^{\star\circ})$ \eqref{TFA9b}.
On pose $M_U=\nu_{U,\oy}(\cF_U|\oU^{\star\circ})$.
En vertu de (\cite{agt} II.6.19), pour tout entier $q\geq 1$, le $\co_\oK$-module $\rH^q(\Sigma_U,M_U)$ est $\alpha$-nul, et le morphisme canonique 
\begin{equation}\label{AFR15i}
M_U^{\Sigma_U}\otimes_{R_{U,\infty}}\oR^\oy_U\rightarrow M_U
\end{equation}
est un $\alpha$-isomorphisme. 

Compte tenu de l'isomorphisme \eqref{AFR24n} et de la définition de l'algèbre $\uoR$ \eqref{TFA12f},
la limite inductive des homomorphismes canoniques $R_{\infty,U}\rightarrow \oR^\oy_U$, prise sur la catégorie $\fV^\circ_\ox$,
s'identifie à l'homomorphisme $\uR_\infty\rightarrow \uoR$ \eqref{AFR24h}.

Soit $q$ un entier $\geq 0$. 
D'après \eqref{AFR15k} et (\cite{agt} VI.11.10), le morphisme canonique 
\begin{equation}\label{AFR15j}
\underset{\underset{n\geq 1}{\longrightarrow}}{\lim} \ \rH^q(\oU^{(n)\star\circ}_\fet,\cF_U|\oU^{(n)\star\circ})\rightarrow\rH^q(\Sigma_U,M_U)
\end{equation}
est un isomorphisme. Par ailleurs, les isomorphismes \eqref{AFR24n} et \eqref{AFR15d} induisent un isomorphisme fonctoriel
\begin{equation}\label{AFR15n}
\rR^q\sigma^{(\infty)}_*(\cF)_\ox\otimes_{\sigma^{(\infty)}_*(\ocB)_\ox}\uR_\infty\stackrel{\sim}{\rightarrow} 
\underset{\underset{n\geq 1}{\longrightarrow}}{\lim} \
\underset{\underset{(U,\fp)\in \fV_\ox^\circ}{\longrightarrow}}{\lim} \ \rH^q(\oU^{(n)\star\circ}_\fet,\cF_U|\oU^{(n)\star\circ}).
\end{equation}
On en déduit un isomorphisme fonctoriel
\begin{equation}\label{AFR15p}
\rR^q\sigma^{(\infty)}_*(\cF)_\ox\otimes_{\sigma^{(\infty)}_*(\ocB)_\ox}\uR_\infty \stackrel{\sim}{\rightarrow} 
\underset{\underset{(U,\fp)\in \fV_\ox^\circ}{\longrightarrow}}{\lim} \ \rH^q(\Sigma_U,M_U).
\end{equation}

D'après \ref{cad8}, le nombre de composantes connexes de $\oU^{(n)}$ reste borné quand $n$ varie. 
Considérant des morphismes  
\begin{equation}
\oy'\rightarrow \underset{\underset{n\geq 1}{\longleftarrow}}{\lim}\ \oU^{(n)\circ}
\end{equation} 
pas nécessairement induits par le morphisme fixé dans \eqref{AFR24e}, 
on déduit de ce qui précède que pour tout $q\geq 1$, le $\co_\oK$-module 
\begin{equation}\label{AFR15l}
\underset{\underset{n\geq 1}{\longrightarrow}}{\lim} \ \rH^q(\oU^{(n)\circ}_\fet,\cF_U)
\end{equation}
est $\alpha$-nul. Intervertissant les limites inductives dans \eqref{AFR15d}, on en déduit 
que pour tout $q\geq 1$, le $\co_\oK$-module $\rR^q\sigma^{(\infty)}_*(\cF)_\ox$ est $\alpha$-nul. 

Reprenons de nouveau le point géométrique fixé dans \eqref{AFR24f}. En vertu de (\cite{agt} VI.10.36), on a un isomorphisme canonique 
\begin{equation}\label{AFR15m}
\cF_{\rho(\oy\rightsquigarrow \ox)}\stackrel{\sim}{\rightarrow} \underset{\underset{(U,\fp)\in \fV_\ox^\circ}{\longrightarrow}}{\lim} \ M_U.
\end{equation}
L'analogue de cet isomorphisme pour le faisceau $\ocB$ n'est autre que l'isomorphisme \eqref{TFA11c}. 
Compte tenu de \eqref{AFR15p} et \eqref{AFR15m}, 
on prend pour morphisme \eqref{AFR15a} la limite inductive des morphismes \eqref{AFR15i}, qui est un $\alpha$-isomorphisme.

\begin{cor}\label{AFR150}
Soient $m, q$ deux entiers $\geq 1$, $\cF$ un $\ocB_m$-module de $\tE$. Alors, le $\hbar_*(\co_{\oX})$-module $\rR^q\sigma^{(\infty)}_*(\cF)$ est $\alpha$-nul.
\end{cor}

En effet, le $\hbar_*(\co_\oX)$-module $\rR^q\sigma^{(\infty)}_*(\cF)$ est $\alpha$-nul en tout point géométrique de $X_s$ en vertu de \ref{AFR15}(ii).
On notera que pour tout point géométrique $\ox'$ de $X_s$, il existe un point $(\oy'\rightsquigarrow \ox')$ de $X_\et\gtimes_{X_\et}\oX^\circ_\et$
en vertu de (\cite{agt} III.3.7), ce qui permet d'appliquer \ref{AFR15}(ii) en $\ox'$. 
Par ailleurs, le $\hbar_*(\co_\oX)$-module $\rR^q\sigma^{(\infty)}_*(\cF)$ étant de $p^m$-torsion, il est nul sur $X_\eta$, d'où la proposition.

\subsection{}\label{AFR16}
Pour tout entier $n\geq 1$, on pose 
\begin{equation}\label{AFR16a}
(X'^{(n)},\cM_{X'^{(n)}})=(X^{(n)},\cM_{X^{(n)}})\times_{(X,\cM_X)}(X',\cM_{X'}),
\end{equation}
le produit étant indifféremment pris dans la catégorie des schémas logarithmiques ou dans celle des schémas logarithmiques saturés (\cite{agt} II.5.20). 
Le morphisme 
\begin{equation}\label{AFR16b}
f'^{(n)}\colon (X'^{(n)},\cM_{X'^{(n)}})\rightarrow (S^{(n)},\cM_{S^{(n)}}),
\end{equation}
déduit de $f^{(n)}$ \eqref{eccr2b}, est donc lisse et saturé (\cite{tsuji4} II 2.11). 
Il résulte de \ref{cad7}(iv) que $X'^{(n)\rhd}=X'^{(n)}\times_{X'}X'^\rhd$  \eqref{mtfla2c} 
est le sous-schéma ouvert maximal de $X'^{(n)}$ où la structure logarithmique $\cM_{X'^{(n)}}$ est triviale (\cite{ogus} III 1.2.8).
D'après (\cite{agt} III.4.2), le schéma $X'^{(n)}$ est normal et $S^{(n)}$-plat, et l'immersion $X'^{(n)\rhd}\rightarrow X'^{(n)}$ est schématiquement dominante.

On pose \eqref{AFR10a}
\begin{equation}\label{AFR16c}
\oX'^{(n)}= X'^{(n)}\times_{S^{(n)}}\oS  \ \ \ {\rm et}\ \ \ \oX'^{(n)\rhd}=\oX'^{(n)}\times_{X'}X'^\rhd.
\end{equation} 
D'après (\cite{agt} III.4.2(iii)), le schéma $\oX'^{(n)}$ est normal et localement irréductible. 
Par ailleurs, l'immersion $\oX'^{(n)\rhd}\rightarrow \oX'^{(n)}$ est schématiquement dominante (\cite{ega4} 11.10.5). 

Pour tous entiers $m,n\geq 1$ tels que $m$ divise $n$, on désigne par $\xi'^{(n,m)}\colon \oX'^{(n)}\rightarrow \oX'^{(m)}$ 
le morphisme induit par $\xi^{(n,m)}$ \eqref{AFR10} et on pose $\xi'^{(n)}=\xi'^{(n,1)}\colon \oX'^{(n)}\rightarrow \oX'$, 
de sorte qu'on a 
\begin{equation}\label{AFR16d}
\xi'^{(n)}=\xi'^{(m)}\circ \xi'^{(n,m)}.
\end{equation}

Considérons le $\oX'$-schéma 
\begin{equation}\label{AFR16e}
\oX'^{(\infty)}= \underset{\underset{n\geq 1}{\longleftarrow}}{\lim}\ \oX'^{(n)},
\end{equation}
et la $\co_{\oX'}$-algèbre quasi-cohérente 
\begin{equation}\label{AFR16f}
\cR^!_\infty= \underset{\underset{n\geq 1}{\longrightarrow}}{\lim}\ \xi'^{(n)}_*(\co_{\oX'^{(n)}}),
\end{equation}
de sorte que $\oX'^{(\infty)}=\Spec_{\oX'}(\cR^!_\infty)$ (\cite{ega4} 8.2.3).

\begin{lem}\label{AFR26}
\
\begin{itemize}
\item[{\rm (i)}] Pour tout ouvert affine $\oU'$ de $\oX'$ au-dessus d'un ouvert affine $\oU$ de $\oX$, 
la $\co_\oK$-algèbre $\Gamma(\oU',\cR^!_\infty)$ est universellement $\alpha$-cohérente \eqref{afini14}. 
\item[{\rm (ii)}] La $\co_\oK$-algèbre $\co_{\oX'^{(\infty)}}$ de $\oX'^{(\infty)}_{\zar}$ est $\alpha$-cohérente \eqref{finita9}. 
\item[{\rm (iii)}]  Pour tout entier $g\geq 0$, 
la $\co_\oK$-algèbre $\cR^!_\infty[t_1,\dots,t_g]$ (des polynômes en $g$ variables à coefficients dans $\cR^!_\infty$) de $\oX'_\zar$ 
est $\alpha$-cohérente \eqref{finita9}. 
\end{itemize}
\end{lem}

(i) En effet, pour tout entier $n\geq 1$, on a un isomorphisme canonique
\begin{equation}
\Gamma(\oU'\times_{\oX'}\oX'^{(n)},\co_{\oX'^{(n)}})\stackrel{\sim}{\rightarrow} \Gamma(\oU\times_{\oX}\oX^{(n)},\co_{\oX^{(n)}})\otimes_{\Gamma(\oU,\co_\oX)}
\Gamma(\oU',\co_{\oX'}).
\end{equation}
On en déduit par passage à la limite inductive un isomorphisme \eqref{AFR10f}
\begin{equation}
\Gamma(\oU',\cR^!_\infty)\stackrel{\sim}{\rightarrow} \Gamma(\oU,\cR_\infty)\otimes_{\Gamma(\oU,\co_\oX)}
\Gamma(\oU',\co_{\oX'}).
\end{equation}
Comme le morphisme canonique $\oX'\rightarrow \oX$ est de présentation finie et que la $\co_\oK$-algèbre $\Gamma(\oU,\cR_\infty)$ est 
universellement $\alpha$-cohérente d'après \ref{AFR25}(ii), il s'ensuit que la $\co_\oK$-algèbre $\Gamma(\oU',\cR^!_\infty)$ est 
universellement $\alpha$-cohérente \eqref{finita20}. 

(ii) \& (iii) Cela résulte de (i) et \ref{afini103}.

\subsection{}\label{AFR17}
On a \eqref{mtfla11a}
\begin{equation}\label{AFR17a}
\Theta^+(\oX^{(n)\circ}\rightarrow X)= (\oX'^{(n)\rhd}\rightarrow X').
\end{equation}
On désigne par $E'^{(n)}\rightarrow \Et_{/X'}$ (resp. $\tE'^{(n)}$) le site fibré (resp. topos) de Faltings associé 
au morphisme $\oX'^{(n)\rhd}\rightarrow X'$ \eqref{tf1}. Tout objet de $E'^{(n)}$
est naturellement un objet de $E'$. On définit ainsi un foncteur 
\begin{equation}\label{AFR17b}
\Phi'^{(n)}\colon E'^{(n)}\rightarrow E',
\end{equation}
qui se factorise à travers une équivalence de catégories 
\begin{equation}\label{AFR17c}
E'^{(n)}\stackrel{\sim}{\rightarrow} E'_{/(\oX'^{(n)\rhd}\rightarrow X')}.
\end{equation}
Le foncteur $\Phi'^{(n)}$ est continu et cocontinu (cf. \ref{amtF1}). 
Il définit donc une suite de trois foncteurs adjoints~:
\begin{equation}\label{AFR17d}
\Phi'^{(n)}_!\colon \tE'^{(n)}\rightarrow \tE', \ \ \ \Phi'^{(n)*}\colon \tE'\rightarrow \tE'^{(n)}, 
\ \ \ \Phi'^{(n)}_*\colon \tE'^{(n)}\rightarrow \tE',
\end{equation}
dans le sens que pour deux foncteurs consécutifs de la suite, celui de droite est
adjoint à droite de l'autre. D'après (\cite{sga4} III 5.4), le foncteur $\Phi'^{(n)}_!$ se factorise à travers 
une équivalence de catégories 
\begin{equation}\label{AFR17e}
\tE'^{(n)}\stackrel{\sim}{\rightarrow} \tE'_{/(\oX'^{(n)\rhd}\rightarrow X')^\tta},
\end{equation}
où $(\oX'^{(n)\rhd}\rightarrow X')^\tta$ est l'objet de $\tE'$ associé à $(\oX'^{(n)\rhd}\rightarrow X')$.
Le couple de foncteurs $(\Phi'^{(n)*},\Phi'^{(n)}_*)$ définit un morphisme de topos que l'on note aussi
\begin{equation}\label{AFR17f}
\Phi'^{(n)}\colon \tE'^{(n)}\rightarrow \tE'.
\end{equation}
Celui-ci est égal au composé de l'équivalence \eqref{AFR17e} et du morphisme de localisation de $\tE'$ en l'objet 
$(\oX'^{(n)\rhd}\rightarrow X')^\tta$. Il se déduit aussi par fonctorialité du diagramme commutatif (\cite{agt} VI.10.12) 
\begin{equation}\label{AFR17g}
\xymatrix{
{\oX'^{(n)\rhd}}\ar[r]\ar[d]&X'\ar@{=}[d]\\
{\oX'^\rhd}\ar[r]&X'}
\end{equation}

D'après (\cite{sga4} IV 5.5),  pour tous entiers $m,n\geq 1$ tels que $m$ divise $n$, on a un morphisme canonique de topos 
\begin{equation}\label{AFR17i}
\Phi'^{(n,m)}\colon \tE'^{(n)}\rightarrow \tE'^{(m)}
\end{equation}
qui s'insère dans un diagramme commutatif à isomorphisme canonique près
\begin{equation}\label{AFR17j}
\xymatrix{
{\tE'^{(n)}}\ar[rd]^{\Phi'^{(n)}}\ar[d]_{\Phi'^{(n,m)}}&\\
{\tE'^{(m)}}\ar[r]_{\Phi'^{(m)}}&{\tE'}}
\end{equation}

\subsection{}\label{AFR18}
On note 
\begin{equation}\label{AFR18a}
\sigma'^{(n)}\colon \tE'^{(n)}\rightarrow X'_\et
\end{equation}
le morphisme canonique (\cite{agt} (VI.10.6.4)). Il résulte aussitôt des définitions (\cite{agt} VI.10.6) qu'on a un isomorphisme canonique 
\begin{equation}\label{AFR18b}
\sigma'^{(n)}\stackrel{\sim}{\rightarrow} \sigma'\circ  \Phi'^{(n)},
\end{equation}
où $\sigma'$ est le morphisme \eqref{mtfla9b}.  

Le schéma $\oX'^{(n)}$ est normal et localement irréductible et 
l'immersion $\oX'^{(n)\rhd}\rightarrow \oX'^{(n)}$ est schématiquement dominante \eqref{AFR16}.
Comme le morphisme $\xi'^{(n)}\colon \oX'^{(n)}\rightarrow \oX'$ est fini, on en déduit un isomorphisme canonique 
\begin{equation}\label{AFR18c}
\sigma'^{(n)}_*(\Phi'^{(n)*}(\ocB'))\stackrel{\sim}{\rightarrow} \hbar'_*(\xi'^{(n)}_*(\co_{\oX'^{(n)}})),
\end{equation}
où $\hbar'\colon \oX'\rightarrow X'$ est le morphisme canonique \eqref{mtfla2}. 
Dans la suite, on considère $\sigma'^{(n)}$ comme un morphisme de topos annelés 
\begin{equation}\label{AFR18d}
\sigma'^{(n)}\colon (\tE'^{(n)},\Phi'^{(n)*}(\ocB'))\rightarrow (X'_\et,\hbar'_*(\xi'^{(n)}_*(\co_{\oX'^{(n)}}))). 
\end{equation}

Pour tout objet $\cF$ de $\tE'$, les faisceaux $(\sigma'^{(n)}_*(\Phi'^{(n)*}(\cF)))_{n\geq 1}$ 
forment naturellement un système inductif de $X'_{\et}$ indexé par l'ensemble $\mZ_{\geq 1}$ ordonné par la relation de divisibilité. 
On pose abusivement
\begin{equation}\label{AFR18e}
\sigma'^{(\infty)}_*(\cF)=\underset{\underset{n\geq 1}{\longrightarrow}}{\lim}\ \sigma'^{(n)}_*(\Phi'^{(n)*}(\cF)).
\end{equation}
D'après (\cite{sga4} VII 5.14), pour tout $X'$-schéma étale de présentation finie $U'$, on a 
\begin{equation}\label{AFR18f}
\sigma'^{(\infty)}_*(\cF)(U')=\underset{\underset{n\geq 1}{\longrightarrow}}{\lim}\ \cF(U'\times_{X'}\oX'^{(n)\rhd}\rightarrow U').
\end{equation}

L'algèbre $\sigma'^{(\infty)}_*(\ocB')$ de $X'_\et$ est associée à l'algèbre quasi-cohérente $\hbar'_*(\cR^!_\infty)$ de $X'_\zar$ 
\eqref{AFR16f}, \eqref{notconv12a}. Ceci résulte aussitôt de \eqref{AFR18c} et du fait que les foncteurs \eqref{notconv12a} et $\hbar'_*$ 
commutent aux limites inductives filtrantes compte tenu de \eqref{notconv12e} et (\cite{sga4} VI 5.1). 

Pour tout entier $q\geq 0$ et tout groupe abélien $\cF$ de $\tE'$, les faisceaux $(\rR^q\sigma'^{(n)}_*(\Phi'^{(n)*}(\cF)))_{n\geq 1}$ 
forment naturellement un système inductif de groupes abéliens de $X'_{s,\et}$. On pose abusivement
\begin{equation}\label{AFR18g}
\rR^q\sigma'^{(\infty)}_*(\cF)=\underset{\underset{n\geq 1}{\longrightarrow}}{\lim}\ \rR^q \sigma'^{(n)}_*(\Phi'^{(n)*}(\cF)).
\end{equation}

\begin{prop}\label{AFR19}
Supposons que le morphisme 
$f'\colon (X',\cM_{X'})\rightarrow (S,\cM_S)$ \eqref{mtfla2} admette une carte adéquate, 
que $X'$ soit affine et connexe et que $X'_s$ soit non-vide. Soient $m\geq 1$, $n\geq 1$, $q\geq 0$ des entiers, 
$\mL'$ un $(\co_\oK/p^m\co_\oK)$-module localement libre de type fini de $\oX'^\rhd_\fet$. 
On note $M'^{(n)}$ le $\ocB'_{X'}(\oX'^{(n)\rhd})$-module 
\begin{equation}\label{AFR19a}
M'^{(n)}=\rH^q(\oX'^{(n)\rhd}_\fet,\mL'\otimes_{\co_\oK}\ocB'_{X'}),
\end{equation} 
et $\tM'^{(n)}$ le $\ocB'_{X'}(\oX'^{(n)\rhd})$-module associé de $X'_\et$ \eqref{notconv12a}. 
Alors, on a un morphisme $\ocB'_{X'}(\oX'^{(n)\rhd})$-linéaire canonique de $X'_\et$ 
\begin{equation}\label{AFR19b}
\tM'^{(n)}\rightarrow \rR^q\sigma'^{(n)}_*(\Phi'^{(n)*}(\beta'^*(\mL'))),
\end{equation}
qui est un $\alpha$-isomorphisme, où $\beta'$ est le morphisme de topos annelés \eqref{mtfla10d}. 
\end{prop}

C'est un cas particulier de \ref{amtF23}. 
 
\begin{cor}\label{AFR20}
Supposons que le morphisme 
$f'\colon (X',\cM_{X'})\rightarrow (S,\cM_S)$ \eqref{mtfla2} admette une carte adéquate, 
que $X'$ soit affine et connexe et que $X'_s$ soit non-vide. Soient $m \geq 1$ et $q\geq 0$ deux entiers et 
$\mL'$ un $(\co_\oK/p^m\co_\oK)$-module localement libre de type fini de $\oX'^\rhd_\fet$. 
On note $M'^{(\infty)}$ le $\sigma'^{(\infty)}_*(\ocB')(X')$-module 
\begin{equation}\label{AFR20a}
M'^{(\infty)}=\underset{\underset{n\geq 1}{\longrightarrow}}{\lim}\
\rH^q(\oX'^{(n)\rhd}_\fet,\mL'\otimes_{\co_\oK}\ocB'_{X'}),
\end{equation} 
et $\tM'^{(\infty)}$ le $\sigma'^{(\infty)}_*(\ocB')(X')$-module associé de $X'_\et$. 
Alors, il existe un morphisme $\sigma'^{(\infty)}_*(\ocB')(X')$-linéaire canonique de $X'_\et$ 
\begin{equation}\label{AFR20b}
\tM'^{(\infty)}\rightarrow \rR^q\sigma'^{(\infty)}_*(\beta'^*(\mL')),
\end{equation}
qui est un $\alpha$-isomorphisme. 
\end{cor}

Cela résulte aussitôt de \ref{AFR19} puisque le foncteur \eqref{notconv12a} commute aux limites inductives compte tenu de \eqref{notconv12e}.

\begin{lem}\label{AFR201}
Conservons les hypothèses de \ref{AFR20} et supposons de plus que le morphisme $g$ \eqref{mtfla2a} 
admette une carte relativement adéquate \eqref{mtfla5}
\begin{equation}\label{AFR201a}
((P',\gamma'),(P,\gamma),(\mN,\iota),\vartheta\colon \mN\rightarrow P, h\colon P\rightarrow P').
\end{equation} 
Alors, le $\cR^!_\infty(\oX')$-module $M'^{(\infty)}$ est $\alpha$-cohérent.
De plus, il est $\alpha$-nul si $q\geq \dim(X'/X)+1$.
\end{lem}

En effet, reprenant les notations de \ref{eccr41}, on observera que pour tout entier $n\geq 1$, $(\oX'^{\{n\}\rhd}\rightarrow X')$ est naturellement
un objet de $E'$. Posons alors 
\begin{eqnarray}
\oX'^{\{\infty\}}&=&\underset{\underset{n\geq 1}{\longleftarrow}}{\lim}\ \oX'^{\{n\}},\\
\sigma'^{\{\infty\}}_*(\ocB')(X')&=&\underset{\underset{n\geq 1}{\longrightarrow}}{\lim}\ \ocB'_{X'}(\oX'^{\{n\}\rhd}),\\
M'^{\{\infty\}}&=&\underset{\underset{n\geq 1}{\longrightarrow}}{\lim}\
\rH^q(\oX'^{\{n\}\rhd}_\fet,\mL'\otimes_{\co_\oK}\ocB'_{X'}).
\end{eqnarray}
En vertu de \ref{eccr412}, $\oX'^{(\infty)}$ est la somme disjointe d'un nombre
fini de copies de $\oX'^{\{\infty\}}$, ou de façon équivalente, $\sigma'^{(\infty)}_*(\ocB')(X')$ est le produit d'un nombre
fini de copies de $\sigma'^{\{\infty\}}_*(\ocB')(X')$. Par ailleurs, le $\sigma'^{\{\infty\}}_*(\ocB')(X')$-module $M'^{\{\infty\}}$ 
est de présentation $\alpha$-finie en vertu de \ref{eccr15}, \ref{eccr8}(iii) et (\cite{agt} VI.11.10). 
Il est, de plus, $\alpha$-nul si $q\geq \dim(X'/X)+1$. 
La $\co_\oK$-algèbre $\cR^!_\infty(\oX')=\sigma'^{(\infty)}_*(\ocB')(X')$ étant $\alpha$-cohérente d'après \ref{AFR26}(i), 
on en déduit que le $\cR^!_\infty(\oX')$-module $M'^{(\infty)}$ est $\alpha$-cohérent \eqref{finita19}. 
De plus, il est $\alpha$-nul si $q\geq \dim(X'/X)+1$.

\begin{prop}\label{AFR200}
Soient $m\geq 1$, $q\geq 0$ deux entiers,
$\mL'$ un $(\co_\oK/p^m\co_\oK)$-module localement libre de type fini de $\oX'^\rhd_\fet$. 
Alors, le $\sigma'^{(\infty)}_*(\ocB')$-module $\fm_\oK\otimes_{\co_\oK}\rR^q \sigma'^{(\infty)}_*(\beta'^*(\mL'))$ est associé à 
un $\hbar'_*(\cR^!_\infty)$-module quasi-cohérent et $\alpha$-cohérent de $X'_\zar$ pour tout $q\geq 0$ \eqref{AFR16f}, 
et il est nul pour tout $q\geq \dim(X'/X)+1$, où $\beta'$ est le morphisme de topos annelés \eqref{mtfla10d}.
\end{prop}

En effet, d'après \ref{mtfla4} (cf. la preuve de \ref{mtfla6}), 
il existe un recouvrement étale $(v'_j\colon X'_j\rightarrow X')_{j\in J}$, un recouvrement de Zariski $(v_j\colon X_j\rightarrow X)_{j\in J}$ 
et pour tout $j\in J$, un diagramme commutatif de schémas logarithmiques 
\begin{equation}\label{AFR200b}
\xymatrix{
{(X'_j,\cM_{X'}|X'_j)}\ar[r]^{v'_j}\ar[d]_{g_j}&{(X',\cM_{X'})}\ar[d]^{g}\\
{(X_j,\cM_X|X_j)}\ar[r]^{v_j}&{(X,\cM_X)}}
\end{equation}
tels que $J$ soit fini et que pour tout $j\in J$, l'une des deux conditions suivantes soit remplie:
\begin{itemize}
\item[(a)] les schémas $X_j$ et $X'_j$ sont affines et connexes, $(X'_{j})_s$ est non-vide et
le morphisme $g_j$ admet une carte relativement adéquate \eqref{mtfla5}
\begin{equation}\label{AFR200c}
((P'_j,\gamma'_j),(P,\gamma_j),(\mN,\iota),\vartheta\colon \mN\rightarrow P, h_j\colon P\rightarrow P'_j),
\end{equation} 
où $((P,\gamma_j),(\mN,\iota),\vartheta)$ est la carte adéquate pour $f\circ v_j$ induite par celle pour $f$ fixée dans \ref{AFR9}; 
\item[(b)] $(X'_{j})_s$ est vide.
\end{itemize}

Pour tout $j\in J$, posons $\cA_j=\hbar'_*(\cR^!_\infty)\otimes_{\co_{X'}}\co_{X'_j}$ et notons 
\begin{equation}\label{AFR200d}
\iota'_j \colon \bMod^\qcoh(\cA_j,X'_{j,\zar})\rightarrow \bMod(\sigma'^{(\infty)}_*(\ocB')|X'_j,X'_{j,\et})
\end{equation}
le foncteur canonique \eqref{notconv12a}. Il existe un $\cA_j$-module quasi-cohérent $\cF^q_j$ de $X'_{j,\zar}$
et un morphisme $\sigma'^{(\infty)}_*(\ocB')(X'_j)$-linéaire 
\begin{equation}\label{AFR200e}
\iota'_j(\cF^q_j)\rightarrow \rR^q \sigma'^{(\infty)}_*(\beta'^*(\mL'))|X'_j
\end{equation}
qui est un $\alpha$-isomorphisme.  Dans le cas (a), c'est une conséquence de \ref{AFR20}, compte tenu de (\cite{agt} VI 10.14).
Dans le cas (b), on prend $\cF^q_j=0$. 
Par ailleurs, en vertu de \ref{afini103} et \ref{AFR201}, 
le $\cA_j$-module $\cF^q_j$ est $\alpha$-cohérent pour tout $q\geq 0$, et il est $\alpha$-nul pour tout $q\geq \dim(X'/X)+1$. 

Compte tenu de \ref{alpha21} et \eqref{notconv12e}, le morphisme \eqref{AFR200e} induit un $\sigma'^{(\infty)}_*(\ocB')(X'_j)$-isomor\-phisme 
\begin{equation}\label{AFR200f}
\iota'_j(\fm_\oK\otimes_{\co_\oK}\cF^q_j)\stackrel{\sim}{\rightarrow}  \fm_\oK\otimes_{\co_\oK}\rR^q \sigma'^{(\infty)}_*(\beta'^*(\mL'))|X'_j.
\end{equation}
On en déduit une donnée de descente sur 
$\iota'_j(\fm_\oK\otimes_{\co_\oK}\cF^q_j)_{j\in J}$ relativement au recouvrement étale $(X'_j\rightarrow X')_{j\in J}$. 
Comme les foncteurs \eqref{AFR200d} sont pleinement fidèles \eqref{notconv12}, 
il en résulte une donnée de descente sur 
$(\fm_\oK\otimes_{\co_\oK}\cF^q_j)_{j\in J}$ relativement au recouvrement étale $(X'_j\rightarrow X')_{j\in J}$. 
D'après (\cite{sga1} VIII 1.1), il existe alors un $\hbar'_*(\cR^!_\infty)$-module quasi-cohérent $\cG^q$ et pour tout $j\in J$, 
un isomorphisme
\begin{equation}\label{AFR200g}
v'^*_j(\cG^q)\stackrel{\sim}{\rightarrow} \fm_\oK\otimes_{\co_\oK}\cF^q_j
\end{equation} 
tels que si l'on munit $(v'^*_j(\cG^q))_{j\in J}$ de la donnée de descente définie par $\cG^q$, 
les isomorphismes \eqref{AFR200g} soient compatibles aux données de descente.
Le $\hbar'_*(\cR^!_\infty)$-module $\cG^q$ est $\alpha$-cohérent d'après \ref{afini103}, \ref{afini12} et \ref{AFR26}(ii). 
De plus, il est nul si $q\geq \dim(X'/X)+1$. Par ailleurs, notant 
\begin{equation}
\iota' \colon \bMod^\qcoh(\hbar'_*(\cR^!_\infty),X'_{\zar})\rightarrow \bMod(\sigma'^{(\infty)}_*(\ocB'),X'_{\et})
\end{equation}
le foncteur canonique \eqref{notconv12a}, les isomorphismes \eqref{AFR200f} se descendent en un isomorphisme $\sigma'^{(\infty)}_*(\ocB')$-linéaire
\begin{equation}\label{AFR200h}
\iota'(\cG^q)\rightarrow  \fm_\oK\otimes_{\co_\oK}\rR^q \sigma'^{(\infty)}_*(\beta'^*(\mL'));
\end{equation}
d'où la proposition.

\begin{prop}\label{AFR23}
Supposons le morphisme $g$ projectif. Alors,
\begin{itemize}
\item[{\rm (i)}] Pour tous entiers $m\geq 1$, $i,j\geq 0$
et tout $(\co_\oK/p^m\co_\oK)$-module localement libre de type fini $\mL'$ de $\oX'^\rhd_\fet$, 
le $\sigma^{(\infty)}_*(\ocB)$-module $\fm_\oK\otimes_{\co_\oK}\rR^ig_*(\rR^j\sigma'^{(\infty)}_*(\beta'^*(\mL')))$ est associé à 
un  $\hbar_*(\cR_\infty)$-module quasi-cohérent et $\alpha$-cohérent de $X_\zar$ \eqref{AFR10f}, 
où $\beta'$ est le morphisme de topos annelés \eqref{mtfla10d}.
\item[{\rm (ii)}] Il existe un entier $q_0\geq 0$ tel que pour tous entiers $m\geq 1$, $i,j\geq 0$ avec $i+j\geq q_0$,
et tout $(\co_\oK/p^m\co_\oK)$-module localement libre de type fini $\mL'$ de $\oX'^\rhd_\fet$, 
le $\sigma^{(\infty)}_*(\ocB)$-module $\fm_\oK\otimes_{\co_\oK}\rR^ig_*(\rR^j\sigma'^{(\infty)}_*(\beta'^*(\mL')))$ soit nul.
\end{itemize}
\end{prop}

En effet, le morphisme canonique
\begin{equation}\label{AFR23a}
\rR^ig_*(\fm_{\oK}\otimes_{\co_\oK} \rR^j\sigma'^{(\infty)}_*(\beta'^*(\mL')))\rightarrow \rR^ig_*(\rR^j\sigma'^{(\infty)}_*(\beta'^*(\mL')))
\end{equation}
est un $\alpha$-isomorphisme \eqref{alpha3}. 
D'après \ref{AFR200}, le $\sigma'^{(\infty)}_*(\ocB')$-module
$\fm_{\oK}\otimes_{\co_\oK} \rR^j\sigma'^{(\infty)}_*(\beta'^*(\mL'))$ est associé à un 
$\hbar'_*(\cR^!_\infty)$-module quasi-cohérent et $\alpha$-cohérent $\cF$. 
Le $\sigma^{(\infty)}_*(\ocB)$-module $\rR^ig_*(\fm_{\oK}\otimes_{\co_\oK} \rR^j\sigma'^{(\infty)}_*(\beta'^*(\mL')))$ 
est donc associé au $\hbar_*(\cR_\infty)$-module $\rR^ig_*(\cF)$ \eqref{notconv120e}.

Le diagramme de morphismes canoniques 
\begin{equation}\label{AFR23b}
\xymatrix{
{\oX'^{(\infty)}}\ar[r]\ar[d]_{g^{(\infty)}}&{X'}\ar[d]^g\\
{\oX^{(\infty)}}\ar[r]&X}
\end{equation}
est cartésien, et ses flèches horizontales sont affines \eqref{AFR10e}, \eqref{AFR16e}. 
Comme $g$ est plat (\cite{kato1} 4.5), $X'$ et $\oX^{(\infty)}$ sont tor-indépendants sur $X$ dans le sens de (\cite{sga6} III 1.5). 
Comme $g$ est pseudo-cohérent (\cite{sga6} III 1.2 et 1.3), il en est de même de $g^{(\infty)}$ d'après (\cite{sga6} III 1.10(a)). 
Le $\hbar'_*(\cR^!_\infty)$-module $\cF$ est l'image directe 
d'un $\co_{\oX'^{(\infty)}}$-module quasi-cohérent et $\alpha$-cohérent $\cG$ \eqref{afini103}.
Par suite, le $\co_{\oX^{(\infty)}}$-module quasi-cohérent $\rR^ig^{(\infty)}_*(\cG)$ est $\alpha$-cohérent en vertu de \ref{afini15} et \ref{AFR25}(ii). 
Le $\hbar_*(\cR_\infty)$-module quasi-cohérent $\rR^ig_*(\cF)$ est donc $\alpha$-cohérent \eqref{afini103}.
De plus, il existe un entier $q_0\geq 0$ qui ne dépend que de $g$, tel que $\rR^ig_*(\cF)=0$ pour $i+j\geq q_0$ 
en vertu de \ref{AFR200} et (\cite{ega3} 1.4.12), d'où la proposition \eqref{alpha21}.

\subsection{}\label{AFR21}
Soit $n$ un entier $\geq 1$. D'après \eqref{AFR17a} et (\cite{sga4} IV 5.10), le morphisme $\Theta$ \eqref{mtfla11b}
induit un morphisme de topos 
\begin{equation}\label{AFR21a}
\Theta^{(n)}\colon \tE'^{(n)}\rightarrow \tE^{(n)}
\end{equation}
qui s'insère dans un diagramme commutatif à isomorphisme canonique près
\begin{equation}\label{AFR21b}
\xymatrix{
{\tE'^{(n)}}\ar[d]_{\Theta^{(n)}}\ar[r]^{\Phi'^{(n)}}&{\tE'}\ar[d]^{\Theta}\\
{\tE^{(n)}}\ar[r]^{\Phi^{(n)}}&{\tE}}
\end{equation}
Pour tout entier $q\geq 0$, le morphisme de changement de base relatif à ce diagramme (\cite{egr1} (1.2.2.3)) 
\begin{equation}\label{AFR21c}
\Phi^{(n)*}\circ \rR^q\Theta_*\rightarrow \rR^q\Theta^{(n)}_*\circ \Phi'^{(n)*}
\end{equation}
est un isomorphisme. 

D'après \eqref{mtfla11c} et \eqref{AFR21b}, le diagramme de morphismes de topos
\begin{equation}\label{AFR21d}
\xymatrix{
{\tE'^{(n)}}\ar[r]^{\sigma'^{(n)}}\ar[d]_{\Theta^{(n)}}&{X'_\et}\ar[d]^{g}\\
{\tE^{(n)}}\ar[r]^{\sigma^{(n)}}&{X_\et}}
\end{equation}
est commutatif à isomorphisme canonique près
\begin{equation}\label{AFR21e}
\sigma^{(n)}_*\circ \Theta^{(n)}_*\stackrel{\sim}{\rightarrow}g_*\circ \sigma'^{(n)}_*.
\end{equation}
Le morphisme $g$ étant cohérent, on en déduit, d'après (\cite{sga4} VI 5.1) et compte tenu de \eqref{AFR21c}, que le diagramme de foncteurs
\begin{equation}\label{AFR21f}
\xymatrix{
{\tE'}\ar[r]^{\sigma'^{(\infty)}_*}\ar[d]_{\Theta_*}&{X'_\et}\ar[d]^{g_*}\\
{\tE}\ar[r]^{\sigma^{(\infty)}_*}&{X_\et}}
\end{equation}
est commutatif à isomorphisme canonique près
\begin{equation}\label{AFR21g}
\sigma^{(\infty)}_*\circ \Theta_*\stackrel{\sim}{\rightarrow}g_*\circ \sigma'^{(\infty)}_*.
\end{equation}

\begin{prop}\label{AFR22}
Soient $m$ un entier $\geq 1$, $\cF'$ un $\ocB'_m$-module de $\tE'$. 
Alors, on a une suite spectrale canonique et fonctorielle
\begin{equation}\label{AFR22a}
\rE_2^{i,j}=\fm_{\oK}\otimes_{\co_\oK}\rR^ig_*(\rR^j\sigma'^{(\infty)}_*(\cF'))\Rightarrow \fm_{\oK}\otimes_{\co_\oK}\sigma^{(\infty)}_*(\rR^{i+j}\Theta_*(\cF')).
\end{equation}
\end{prop}

En effet, pour tout entier $n\geq 1$, on a la suite spectrale de Cartan-Leray (\cite{sga4} V 5.4)
\begin{equation}\label{AFR22b}
\rR^ig_*(\rR^j\sigma'^{(n)}_*(\Phi'^{(n)*}(\cF')))\Rightarrow \rR^{i+j}(g\circ \sigma'^{(n)})_*(\Phi'^{(n)*}(\cF')).
\end{equation}
Ces suites forment un système inductif indexé par l'ensemble $\mZ_{\geq 1}$ ordonné par la relation de divisibilité. 
D'après (\cite{sga4} VI 5.1), le morphisme $g$ étant cohérent, 
on en déduit par passage à la limite inductive une suite spectrale 
\begin{equation}\label{AFR22c}
\rR^ig_*(\rR^j\sigma'^{(\infty)}_*(\cF'))\Rightarrow \underset{\underset{n\geq 1}{\longrightarrow}}{\lim}\ 
\rR^{i+j}(g\circ \sigma'^{(n)})_*(\Phi'^{(n)*}(\cF')).
\end{equation}
De même, compte tenu de \eqref{AFR21c}, on a une suite spectrale 
\begin{equation}\label{AFR22d}
\rR^i\sigma^{(\infty)}_*(\rR^j\Theta_*(\cF'))\Rightarrow \underset{\underset{n\geq 1}{\longrightarrow}}{\lim}\ 
\rR^{i+j}(\sigma^{(n)}\circ \Theta^{(n)})_*(\Phi'^{(n)*}(\cF')).
\end{equation}
On en déduit, pour tout entier $q\geq 0$,  un morphisme 
\begin{equation}\label{AFR22e}
\underset{\underset{n\geq 1}{\longrightarrow}}{\lim}\ 
\rR^q(\sigma^{(n)}\circ \Theta^{(n)})_*(\Phi'^{(n)*}(\cF'))
\rightarrow \sigma^{(\infty)}_*(\rR^q\Theta_*(\cF')).
\end{equation}
Celui-ci est un $\alpha$-isomorphisme d'après \ref{AFR150}. La suite spectrale recherchée \eqref{AFR22a} se déduit de \eqref{AFR22c} et \eqref{AFR22e},
compte tenu des isomorphismes \eqref{AFR21e}.

\begin{prop}\label{AFR27}
Supposons le morphisme $g$ projectif. Alors, 
\begin{itemize}
\item[{\rm (i)}] Pour tous entiers $m\geq 1$, $q\geq 0$ et tout $(\co_\oK/p^m\co_\oK)$-module localement libre de type fini $\mL'$ de $\oX'^\rhd_\fet$, 
le $\sigma^{(\infty)}_*(\ocB)_\ox$-module $\sigma^{(\infty)}_*(\rR^q\Theta_*(\beta'^*(\mL')))_\ox$
est $\alpha$-cohérent, où $\beta'$ est le morphisme de topos annelés \eqref{mtfla10d}.
\item[{\rm (ii)}] Il existe un entier $q_0\geq 0$ tel que pour tous entiers 
$m\geq 1$, $q\geq q_0$ et tout $(\co_\oK/p^m\co_\oK)$-module localement libre de type fini $\mL'$ de $\oX'^\rhd_\fet$, 
le $\sigma^{(\infty)}_*(\ocB)_\ox$-module $\sigma^{(\infty)}_*(\rR^q\Theta_*(\beta'^*(\mL')))_\ox$ soit $\alpha$-nul.
\end{itemize}
\end{prop}

Cela résulte de \ref{AFR23} et \ref{AFR22}, compte tenu de \ref{finita17}, \ref{finita18}, \ref{afini17}(ii) et \ref{AFR25}(ii).

\begin{cor}\label{AFR28}
Supposons le morphisme $g$ projectif. Alors, 
\begin{itemize}
\item[{\rm (i)}] Pour tous entiers $m\geq 1$, $q\geq 0$ et tout $(\co_\oK/p^m\co_\oK)$-module localement libre de type fini $\mL'$ de $\oX'^\rhd_\fet$, 
le $\uoR$-module $\rR^q\Theta_*(\beta'^*(\mL'))_{\rho(\oy\rightsquigarrow \ox)}$ est $\alpha$-cohérent, où $\rho$ est le morphisme \eqref{mtfla7e}
et $\beta'$ est le morphisme de topos annelés \eqref{mtfla10d}. 
\item[{\rm (ii)}] Il existe un entier $q_0\geq 0$ tel que pour tous entiers $m\geq 1$, $q\geq q_0$ 
et tout $(\co_\oK/p^m\co_\oK)$-module localement libre de type fini $\mL'$ de $\oX'^\rhd_\fet$, 
le $\uoR$-module $\rR^q\Theta_*(\beta'^*(\mL'))_{\rho(\oy\rightsquigarrow \ox)}$ soit $\alpha$-nul.
\end{itemize} 
\end{cor}

Cela résulte de  \ref{AFR30}, \ref{AFR15}(i) et \ref{AFR27}, compte tenu de \ref{finita14}(iii) et \ref{finita19}.

\subsection{}\label{AFR29}
On peut maintenant démontrer le théorème \ref{AFR8}. Comme les foncteurs $\delta_*$ et $\delta'_*$ sont exacts \eqref{mtfla11e}, 
pour tout groupe abélien $F$ de $\tE_s$ et tout entier $q\geq 0$, on a un isomorphisme canonique 
\begin{equation}
\delta_*(\rR^q\Theta_{s*}(F))\stackrel{\sim}{\rightarrow}\rR^q\Theta_*(\delta'_*F). 
\end{equation}
Par ailleurs, $\beta'^*(\mL')$ étant un objet de $\tE'_s$ (\cite{agt} III.9.6), on a un isomorphisme canonique 
\begin{equation}
\delta'_*(\beta'^*_m(\mL'))\stackrel{\sim}{\rightarrow}\beta'^*(\mL').
\end{equation}
Le théorème \ref{AFR8} résulte alors de \ref{AFR28}, compte tenu de (\cite{agt} III.9.5(i)).

\section{Longueurs normalisées}\label{longnm}

\subsection{}\label{longnm1}
Dans cette section, on suppose que le morphisme $f\colon (X,\cM_X)\rightarrow (S,\cM_S)$ admet une carte adéquate 
$((P,\gamma),(\mN,\iota),\vartheta \colon \mN\rightarrow P)$ (cf. \ref{cad1}) et on reprend les notations de \ref{eccr2}.
On fixe, de plus, un $S$-morphisme
\begin{equation}\label{longnm1a}
\oS\rightarrow \underset{\underset{n\geq 1}{\longleftarrow}}{\lim}\ S^{(n)}.
\end{equation}
Pour tout entier $n\geq 1$, on pose
\begin{equation}\label{longnm1b}
\oX^{(n)}= X^{(n)}\times_{S^{(n)}}\oS.
\end{equation}
Pour tous entiers $m,n\geq 1$ tels que $m$ divise $n$, on note $\xi^{(n,m)}\colon \oX^{(n)}\rightarrow \oX^{(m)}$ le morphisme canonique
qui est fini et surjectif. On pose $\xi^{(n)}=\xi^{(n,1)}\colon \oX^{(n)}\rightarrow \oX$, de sorte qu'on a 
$\xi^{(n)}=\xi^{(m)}\circ \xi^{(n,m)}$.

\subsection{}\label{longnm3}
Soient $(\oy\rightsquigarrow \ox)$ un point de $X_\et\gtimes_{X_\et}\oX^\circ_\et$ \eqref{topfl17} 
tel que $\ox$ soit au-dessus de $s$, $x$ le support de $\ox$ dans $X_s$, 
$\uX$ le localisé strict de $X$ en $\ox$, $\uoX=\uX\times_S\oS$ \eqref{mtfla1b}. 
Le corps résiduel de $\co_K$ étant algébriquement clos, on peut identifier $x$ (resp. $\ox$) à un point (resp. point géométrique) de $\oX$.
Le schéma $\uoX$ s'identifie donc au localisé strict de $\oX$ en $\ox$ (\cite{agt} III.3.7). 
Le $X$-morphisme $u\colon \oy\rightarrow \uX$
définissant le point $(\oy\rightsquigarrow \ox)$ se relève en un $\oX^\circ$-morphisme 
$v\colon \oy\rightarrow \uoX^\circ=\uoX\times_XX^\circ$ \eqref{mtfla2b}. 
Pour tout entier $n\geq 1$, on pose
\begin{equation}\label{longnm3a}
\uoX^{(n)}=\oX^{(n)}\times_X\uX.
\end{equation}
Pour tous entiers $m,n\geq 1$ tels que $m$ divise $n$, on désigne par $\uxi^{(n,m)}\colon \uoX^{(n)}\rightarrow \uoX^{(m)}$ et 
$\uxi^{(n)}\colon \uoX^{(n)}\rightarrow \uoX$ les morphismes induits par $\xi^{(n,m)}$ et $\xi^{(n)}$. 
D'après (\cite{ega4} 8.3.8(i)), il existe un $\uoX$-morphisme 
\begin{equation}\label{longnm3b}
\oy\rightarrow \underset{\underset{n\geq 1}{\longleftarrow}}{\lim}\ \uoX^{(n)},
\end{equation} 
la limite inductive étant indexée par l'ensemble $\mZ_{\geq 1}$ ordonné par la relation de divisibilité, 
que l'on fixe dans la suite de cette section.

D'après (\cite{agt} III.3.7), le schéma $\uoX$ est strictement local et normal (et en particulier intègre).
Le schéma $\uoX^{(n)}$ est donc une union disjointe finie de schémas normaux et strictement locaux (\ref{cad7}(iii)). 
On note $c(n)$ de nombre de ses composantes irréductibles. 
On désigne par $\uoX^{(n)\star}$ la composante irréductible de $\uoX^{(n)}$ contenant l'image de $\oy$. On pose
\begin{eqnarray}
\uR_n&=&\Gamma(\uoX^{(n)\star},\co_{\uoX^{(n)}}), \label{longnm3c}\\
\uR_\infty&=&\underset{\underset{n\geq 1}{\longrightarrow}}{\lim}\ \uR_n,\label{longnm3d}
\end{eqnarray}
la limite inductive étant indexée par l'ensemble $\mZ_{\geq 1}$ ordonné par la relation de divisibilité. L'anneau $\uR_n$ est donc un facteur direct direct de $\uB_n$. 
Les anneaux $\uR_n$ sont locaux, strictement henséliens et normaux, et il en est de même de l'anneau $\uR_\infty$ 
(\cite{raynaud1} I §~3 prop.~1 et \cite{ega1n} 0.6.5.12(ii)).  On note $\ox_\infty$ le point fermé de $\Spec(\uR_\infty)$.
Pour tout $\uR_\infty$-module $M$, on désigne par $\Gamma_{\ox_\infty}(M)$ le sous-$\uR_\infty$-module de $M$ des sections à support dans $\ox_\infty$ \eqref{decprim4}. 

Pour tout $n\geq 1$, le morphisme $\xi^{(n)\circ}\colon \oX^{(n)\circ}\rightarrow \oX^{\circ}$ est un revêtement étale fini et galoisien d'après \ref{cad7}(v).
Munissant le schéma intègre $\uoX^\circ$ du point géométrique défini par le morphisme $v\colon \oy\rightarrow \uoX^\circ$,
le groupe profini $\pi_1(\uoX^{\circ},\oy)$ agit alors sur l'anneau $\uR_\infty$ et l'action est discrète.

\begin{rema}\label{longnm4}
Soient $(\oy^\iota\rightsquigarrow \ox)$ un autre point de $X_\et\gtimes_{X_\et}\oX^\circ_\et$,
\begin{equation}\label{longnm4a}
\oy^\iota\rightarrow \underset{\underset{n\geq 1}{\longleftarrow}}{\lim}\ \uoX^{(n)},
\end{equation} 
un $\uoX$-morphisme. Nous affectons d'un exposant $\iota$ les objets associés à ces données,
analogues à ceux associés à $(\oy\rightsquigarrow \ox)$ et \eqref{longnm3b} dans \ref{longnm3}. 
Comme $X^{(n)\circ}$ est un espace principal homogène pour la topologie étale 
au-dessus de $X^\circ$ sous le groupe $\Hom_\mZ(P^\gp,\mu_n(\oK))$ d'après \ref{cad7}(v),
il existe un isomorphisme de $\uR_1$-algèbres
\begin{equation}\label{longnm4c}
\iota_\infty\colon \uR_\infty\stackrel{\sim}{\rightarrow} \uR_\infty^\iota.
\end{equation} 
\end{rema}

\begin{prop}\label{longnm5}
La catégorie des $\uR_\infty$-modules de présentation $\alpha$-finie est abélienne, autrement dit, elle est stable par noyaux et conoyaux. 
Elle est de plus stable par les foncteurs $\Tor^{\uR_\infty}_i(-,-)$ et $\Ext^i_{\uR_\infty}(-,-)$ $(i\geq 0)$. 
\end{prop}

En effet, la $\co_\oK$-algèbre $\uR_\infty$ étant $\alpha$-cohérente d'après \ref{AFR30}, pour qu'un $\uR_\infty$-module 
soit de présentation $\alpha$-finie, il faut et il suffit qu'il soit $\alpha$-cohérent \eqref{finita19}.  
Par ailleurs, la catégorie des $\uR_\infty$-modules $\alpha$-cohérents est abélienne en vertu de \ref{finita18}, d'où la première assertion. 

Soient $M, N$ deux $\uR_\infty$-modules $\alpha$-cohérents, $\gamma\in \fm_{\oK}$. Il existe un entier $n\geq 0$
et un morphisme $\uR_\infty$-linéaire $u\colon (\uR_\infty)^n\rightarrow M$ dont le conoyau est annulé par $\gamma$. 
Posons $M_0=(\uR_\infty)^n$, $M_1=\ker(u)$ et $M_{-1}=\im(u)$ de sorte qu'on a une suite exacte 
$0\rightarrow M_1\rightarrow M_0\rightarrow M_{-1}\rightarrow 0$ et un $\gamma$-isomorphisme $M_{-1}\rightarrow M$. 
On en déduit des suites exactes 
\begin{eqnarray*}
0\rightarrow \Tor^{\uR_\infty}_1(M_{-1},N)\rightarrow M_1\otimes_{\uR_\infty}N\rightarrow M_0\otimes_{\uR_\infty}N\rightarrow 
M_{-1}\otimes_{\uR_\infty}N\rightarrow 0,\\
0\rightarrow \Hom_{\uR_\infty}(M_{-1},N) \rightarrow \Hom_{\uR_\infty}(M_0,N)\rightarrow \Hom_{\uR_\infty}(M_{1},N) \rightarrow \Ext^1_{\uR_\infty}(M_{-1},N)
\rightarrow 0,
\end{eqnarray*}
et pour tout $i\geq 1$, des isomorphismes 
\begin{eqnarray}
\Tor^{\uR_\infty}_{i+1}(M_{-1},N)&\stackrel{\sim}{\rightarrow}&\Tor^{\uR_\infty}_i(M_1,N),\label{longnm5a}\\
\Ext^i_{\uR_\infty}(M_1,N)&\stackrel{\sim}{\rightarrow}& \Ext^{i+1}_{\uR_\infty}(M_{-1},N).\label{longnm5b}
\end{eqnarray}
Les $\uR_\infty$-modules $M\otimes_{\uR_\infty}N$ et $\Hom_{\uR_\infty}(M,N)$ sont $\alpha$-cohérents  d'après \ref{finita180}. 
Pour tout $j\in \{-1,0,1\}$, le $\uR_\infty$-module $M_j$ étant $\alpha$-cohérent, il en est de même de $M_j\otimes_{\uR_\infty}N$
et $\Hom_{\uR_\infty}(M_j,N)$. Par suite, les $\uR_\infty$-modules $\Tor^{\uR_\infty}_1(M_{-1},N)$
et $\Ext^1_{\uR_\infty}(M_{-1},N)$ sont $\alpha$-cohérents \eqref{finita18}. Il en est alors de même des $\uR_\infty$-modules $\Tor^{\uR_\infty}_1(M,N)$
et $\Ext^1_{\uR_\infty}(M,N)$, compte tenu de \ref{alpha3} et \ref{finita21}.

On démontre par récurrence que les $\uR_\infty$-modules $\Tor^{\uR_\infty}_i(M,N)$
et $\Ext^i_{\uR_\infty}(M,N)$, pour tout $i\geq 1$, sont $\alpha$-cohérents. 
En effet, supposons que pour tous $\uR_\infty$-modules $\alpha$-cohérents $M'$ et $N'$,
les $\uR_\infty$-modules $\Tor^{\uR_\infty}_i(M',N')$ et $\Ext^i_{\uR_\infty}(M',N')$ sont $\alpha$-cohérents. 
On en déduit par \eqref{longnm5a} et \eqref{longnm5b} que les $\uR_\infty$-modules $\Tor^{\uR_\infty}_{i+1}(M_{-1},N)$
et $\Ext^{i+1}_{\uR_\infty}(M_{-1},N)$ sont $\alpha$-cohérents.
Il en est alors de même des $\uR_\infty$-modules $\Tor^{\uR_\infty}_{i+1}(M,N)$ et $\Ext^{i+1}_{\uR_\infty}(M,N)$, 
compte tenu encore de \ref{alpha3} et \ref{finita21}.

\begin{lem}\label{longnm50}
Soient $A$ une $\uR_\infty$-algèbre plate, $M$, $N$ deux $\uR_\infty$-modules. 
On a alors un morphisme de bifoncteurs cohomologiques 
\begin{equation}\label{longnm50a}
\Ext^i_{\uR_\infty}(M,N)\otimes_{\uR_\infty}A\rightarrow \Ext^i_A(M\otimes_{\uR_\infty}A,N\otimes_{\uR_\infty}A)
\end{equation}
se réduisant en degré $0$ au morphisme canonique 
\begin{equation}\label{longnm50b}
\Hom_{\uR_\infty}(M,N)\otimes_{\uR_\infty}A\rightarrow \Hom_A(M\otimes_{\uR_\infty}A,N\otimes_{\uR_\infty}A).
\end{equation}
Si, de plus, $M$ est $\alpha$-cohérent, les morphismes \eqref{longnm50a} sont des $\alpha$-isomorphismes. 
\end{lem}

Comme les foncteurs $M\mapsto \Hom_{\uR_\infty}(M,N)\otimes_{\uR_\infty}A$ et $M\mapsto \Hom_A(M\otimes_{\uR_\infty}A,N\otimes_{\uR_\infty}A)$
de la catégorie opposée de la catégorie des $\uR_\infty$-modules dans la catégorie des $A$-modules, sont exacts à gauche, 
le morphisme \eqref{longnm50b} induit un morphisme de leurs foncteurs dérivés. 
Considérant une résolution $M_\bullet$ de $M$ par des $\uR_\infty$-modules projectifs, on a un morphisme des modules de cohomologie 
\begin{equation}\label{longnm50c}
\rH^i(\Hom_{\uR_\infty}(M_\bullet,N)\otimes_{\uR_\infty}A)\rightarrow \rH^i(\Hom_A(M_\bullet\otimes_{\uR_\infty}A,N\otimes_{\uR_\infty}A)).
\end{equation} 
Comme $A$ est $\uR_\infty$-plat, la source de ce morphisme est canoniquement isomorphe à $\Ext^i_{\uR_\infty}(M,N)\otimes_{\uR_\infty}A$. 
D'autre part, $M_\bullet\otimes_{\uR_\infty}A$ étant une résolution de $M\otimes_{\uR_\infty}A$ par des $A$-modules projectifs,
le but de ce morphisme est canoniquement isomorphe à $\Ext^i_A(M\otimes_{\uR_\infty}A,N\otimes_{\uR_\infty}A)$. 
On prend alors pour \eqref{longnm50a} le morphisme \eqref{longnm50c}.

Supposons que $M$ soit $\alpha$-cohérent et montrons que  les morphismes \eqref{longnm50a} sont des $\alpha$-isomorphismes.
Soit $\gamma\in \fm_{\oK}$. D'après \ref{finita2}(i), il existe un $\uR_\infty$-module de présentation finie $M_{-1}$ 
et un $\gamma$-isomorphisme $v\colon M_{-1}\rightarrow M$. Soient $n$ un entier $\geq 0$,
$u\colon (\uR_\infty)^n\rightarrow M_{-1}$ un morphisme $\uR_\infty$-linéaire surjectif. 
Posons $M_0=(\uR_\infty)^n$ et $M_1=\ker(u)$ de sorte qu'on a une suite exacte $0\rightarrow M_1\rightarrow M_0\rightarrow M_{-1}\rightarrow 0$. 
On en déduit des suites exactes 
\begin{equation}\label{longnm50f}
0\rightarrow \Hom_{\uR_\infty}(M_{-1},N) \rightarrow \Hom_{\uR_\infty}(M_0,N)\rightarrow \Hom_{\uR_\infty}(M_{1},N) \rightarrow \Ext^1_{\uR_\infty}(M_{-1},N)
\rightarrow 0,
\end{equation}
\begin{eqnarray}\label{longnm50g}
\lefteqn{
0\rightarrow \Hom_A(M_{-1}\otimes_{\uR_\infty}A,N\otimes_{\uR_\infty}A) \rightarrow 
\Hom_{\uR_\infty}(M_0\otimes_{\uR_\infty}A,N\otimes_{\uR_\infty}A)\rightarrow}\\
&&\Hom_A(M_{1}\otimes_{\uR_\infty}A,N\otimes_{\uR_\infty}A) \rightarrow \Ext^1_A(M_{-1}\otimes_{\uR_\infty}A,N\otimes_{\uR_\infty}A)\rightarrow 0, \nonumber
\end{eqnarray}
et pour tout $i\geq 1$, des isomorphismes 
\begin{eqnarray}
\Ext^i_{\uR_\infty}(M_1,N)&\stackrel{\sim}{\rightarrow}& \Ext^{i+1}_{\uR_\infty}(M_{-1},N),\label{longnm50d}\\
\Ext^i_A(M_1\otimes_{\uR_\infty}A,N\otimes_{\uR_\infty}A)&\stackrel{\sim}{\rightarrow}& \Ext^{i+1}_A(M_{-1}\otimes_{\uR_\infty}A,N\otimes_{\uR_\infty}A).\label{longnm50e}
\end{eqnarray}
De plus, le morphisme \eqref{longnm50a} induit un morphisme du complexe déduit de la suite \eqref{longnm50f} en étendant les scalaires de $\uR_\infty$ vers $A$, 
dans le complexe défini par la suite \eqref{longnm50g}. Par ailleurs, le diagramme 
\begin{equation}\label{longnm50h}
\xymatrix{
{\Ext^i_{\uR_\infty}(M_1,N)\otimes_{\uR_\infty}A}\ar[r]^-(0.5)\sim\ar[d]_{w^i_1}& {\Ext^{i+1}_{\uR_\infty}(M_{-1},N)\otimes_{\uR_\infty}A}\ar[d]^{w^{i+1}_{-1}}\\
{\Ext^i_A(M_1\otimes_{\uR_\infty}A,N\otimes_{\uR_\infty}A)}\ar[r]^-(0.5)\sim&{\Ext^{i+1}_A(M_{-1}\otimes_{\uR_\infty}A,N\otimes_{\uR_\infty}A)}}
\end{equation}
où les flèches horizontales sont induites par \eqref{longnm50d} et \eqref{longnm50e} et les flèches verticales par \eqref{longnm50a}, est commutatif.  

D'après \ref{alpha3}, il existe un morphisme $\uR_\infty$-linéaire $v'\colon M\rightarrow M_{-1}$ tel que $v\circ v'=\gamma^2\id_{M}$ et 
$v'\circ v=\gamma^2\id_{M_{-1}}$. Pour tout entier $i\geq 0$, considérons le diagramme commutatif 
\begin{equation}\label{longnm50i}
\xymatrix{
{\Ext^i_{\uR_\infty}(M,N)\otimes_{\uR_\infty}A}\ar[r]^-(0.5){a^i}\ar[d]_{w^i}& {\Ext^i_{\uR_\infty}(M_{-1},N)\otimes_{\uR_\infty}A}\ar[d]^{w^i_{-1}}\\
{\Ext^i_A(M\otimes_{\uR_\infty}A,N\otimes_{\uR_\infty}A)}\ar[r]^-(0.5){b^i}&{\Ext^i_A(M_{-1}\otimes_{\uR_\infty}A,N\otimes_{\uR_\infty}A)}}
\end{equation}
où les flèches $a^i$ et $b^i$ sont induites par $v$, et les flèches $w^i$ et $w^i_{-1}$ sont les morphismes canoniques \eqref{longnm50a}.
On notera que $a^i$ et $b^i$ sont des $\gamma^2$-isomorphismes.

Montrons par récurrence sur $i\geq 0$ que le morphisme $w^i$ est un $\alpha$-isomorphisme. 
Comme $M_{-1}$ est de présentation finie, $w^0_{-1}$ est un isomorphisme. Par ailleurs, $a^0$ et $b^0$ étant des $\gamma^2$-isomorphismes, 
on en déduit que $w^0$ est un $\gamma^4$-isomorphisme \eqref{longnm50i}. 
Par suite, $w^0$ est un $\alpha$-isomorphisme, et ceci vaut pour tout $\uR_\infty$-module $\alpha$-cohérent $M$. 
Comme le $\uR_\infty$-module $M_1$ est $\alpha$-cohérent, $w^0_1$ est un $\alpha$-isomorphisme. 
Les suites exactes \eqref{longnm50f} et \eqref{longnm50g} montrent alors que $w^1_{-1}$ est un 
$\alpha$-isomorphisme. On en déduit que $w^1$ est un $\gamma^5$-isomorphisme \eqref{longnm50i}. Par suite, $w^1$ est un $\alpha$-isomorphisme. 

Supposons que $w^i$ soit un $\alpha$-isomorphisme pour un entier $i\geq 1$ et tout $\uR_\infty$-module $\alpha$-cohérent $M$, et montrons que 
$w^{i+1}$ est un $\alpha$-isomorphisme. L'hypothèse de récurrence implique que $w_1^i$ est un $\alpha$-isomorphisme, et il en est alors de même de 
$w_{-1}^{i+1}$ compte tenu de \eqref{longnm50h}. Le diagramme  \eqref{longnm50i} montre alors que $w^{i+1}$ est un $\gamma^5$-isomorphisme. 
Par suite, $w^{i+1}$ est un $\alpha$-isomorphisme, d'où la proposition.

\subsection{}\label{longnm7}
Soit $(A_i,\varphi_{ji})$ un système inductif d'anneaux de valuation dont l'ensemble d'indices est préordonné filtrant à droite tel que 
pour tout $i\leq j$, l'homomorphisme $\varphi_{ji}$ est injectif. Pour tout $i\in I$, notons $K_i$ le corps des fractions de $A_i$. 
Posons $A=\underset{\longrightarrow}{\lim}\ A_i$ et $K=\underset{\longrightarrow}{\lim}\ K_i$. 
Alors, $A$ est un anneau intègre de corps des fractions canoniquement isomorphe à $K$ (\cite{ega1n} 0.6.1.5). 
Il résulte aussitôt de la définition (\cite{ac} VI §~1.2 déf.~2) que $A$ est un anneau de valuation. Le groupe des valeurs $K^\times/A^\times$ de la valuation
canonique de $K$ est la limite inductive des groupes des valeurs $K_i^\times/A_i^\times$ des valuations canoniques des $K_i$.

\begin{prop}\label{longnm8}
Si $x$ est un point générique de $X_s$ \eqref{longnm3}, alors $\uR_\infty$ est un anneau de valuation non discrète de hauteur un,
ayant une valuation qui prolonge la valuation $v$ de $\co_\oK$ \eqref{mtfla1}.
\end{prop}

Soient $n$ un entier $\geq 1$, $L$ une extension finie de $K_n$ \eqref{eccr20}, $\co_L$ la clôture intégrale de $\co_K$ dans $L$.
L'anneau $\uR_n$ \eqref{longnm3c} est le localisé strict de $\oX^{(n)}$ en un point géométrique $\ox^{(n)}$ au-dessus de $\ox$. 
Comme le morphisme $\oX^{(n)}\rightarrow \oX$ est fini et surjectif, $\ox^{(n)}$ est localisé en 
un point générique de la fibre spéciale de $\oX^{(n)}\rightarrow \oS$.
Le schéma $X^{(n)}_L=X^{(n)}\otimes_{\co_{K_n}}\co_L$ \eqref{eccr2a} est normal et $\co_L$-plat 
en vertu de \ref{cad7}(i), (\cite{kato1} 4.5) et (\cite{kato2} 8.2 et 4.1); cf. aussi (\cite{tsuji1} 1.5.1). De plus, 
la fibre spéciale de $X^{(n)}_L\rightarrow \Spec(\co_L)$ est réduite d'après (\cite{tsuji4} II.4.2).  
Notons $R_{L,n}$ le localisé strict de $X^{(n)}_L$ en $\ox^{(n)}$. 
Il résulte de ce qui précède que $R_{L,n}$ est un anneau de valuation discrète dont les uniformisantes de $\co_L$ sont des uniformisantes. 

D'après (\cite{raynaud1} VIII §~2 prop.~3), $\uR_n$  est la limite inductive 
des anneaux $R_{L,n}$ lorsque $L$ décrit les extensions finies de $K_n$ contenues dans $\oK$. 
Il s'ensuit, compte tenu de \ref{longnm7}, que $\uR_n$ est un anneau de valuation et qu'il est muni 
d'une valuation surjective $\uR_n\rightarrow \mQ_{\geq 0}\cup \{\infty\}$ qui prolonge la valuation
$v\colon \co_{\oK}\rightarrow \mQ_{\geq 0}\cup \{\infty\}$ \eqref{mtfla1}.

Pour tout entier $m\geq 1$, le morphisme canonique $\oX^{(mn)}\rightarrow \oX^{(n)}$ est dominant. 
L'homomorphisme canonique $\uR_n\rightarrow \uR_{mn}$ est donc injectif. D'après ce qui précède, celui-ci induit un isomorphisme
entre les groupes des valeurs des valuations canoniques des corps de fractions. 
La proposition s'ensuit compte tenu de \ref{longnm7}.

\begin{defi}[\cite{gr1} 14.5.3]\label{longnm16}
Soit $A$ une $\co_\oK$-algèbre locale, d'idéal maximal contenant $\fm_\oK$. On dit que $A$ est {\em mesurable} s'il existe une 
$\co_\oK$-algèbre locale $B$, essentiellement de présentation finie sur $\co_\oK$, d'idéal maximal contenant $\fm_\oK$, 
et un homomorphisme local et ind-étale $B\rightarrow A$ (\cite{raynaud1} VIII défi.~3).
\end{defi}

\begin{prop}[\cite{gr1} 14.5.51, 14.5.54] \label{longnm17}
Soit $A$ une $\co_\oK$-algèbre mesurable \eqref{longnm16} 
d'idéal maximal $\fm_A$ et de corps résiduel $\kappa_A$. On peut associer canoniquement à tout $A$-module $M$ 
à support dans le point fermé $s_A$ de $\Spec(A)$, une {\em longueur} $\lambda_A(M)\in \mR_{\geq 0}\cup\{\infty\}$ 
vérifiant les propriétés suivantes:
\begin{itemize}
\item[{\rm (i)}] Pour tout $A$-module de présentation finie $M$ à support dans $s_A$, $\lambda_A(M)$ est fini.  
\item[{\rm (ii)}] Pour toute suite exacte de $A$-modules à supports dans $s_A$, $0\rightarrow M'\rightarrow M\rightarrow M''\rightarrow 0$, on~a 
\begin{equation}
\lambda_A(M)= \lambda_A(M')+\lambda_A(M'').
\end{equation}
\item[{\rm (iii)}] Soient $\psi\colon A\rightarrow B$ un homomorphisme local de $\co_\oK$-algèbres mesurables induisant un morphisme fini 
d'extensions résiduelles $\kappa_A\rightarrow \kappa_B$, $M$ un $B$-module à support dans le point fermé $s_B$ de $\Spec(B)$. Alors, 
on a 
\begin{equation}
\lambda_B(M)=\frac{\lambda_A(\psi_*(M))}{[\kappa_B:\kappa_A]}.
\end{equation}
\item[{\rm (iv)}] Soient $\psi\colon A\rightarrow B$ un homomorphisme local plat de $\co_\oK$-algèbres mesurables induisant un morphisme entier 
$\kappa_A\rightarrow B/\fm_A B$, $M$ un $A$-module à support dans  $s_A$. Alors, 
on a 
\begin{equation}
\lambda_B(M\otimes_AB)=\longueur(B/\fm_A B)\cdot \lambda_A(M).
\end{equation}
\end{itemize}
\end{prop}

\subsection{}\label{longnm19}
Soit $A$ une $\co_\oK$-algèbre plate et mesurable qui est un anneau de valuation.  
Notons $\Gamma_A$ (resp. $\Gamma_{\co_\oK}$) le groupe des valeurs de la valuation de $A$ (resp. $\co_\oK$). 
D'après (\cite{gr1} 14.5.5), la valuation de $A$ est de rang un  
et l'indice de ramification $(\Gamma_A:\Gamma_{\co_\oK})$ est fini. 
Comme $\Gamma_{\co_\oK}\simeq \mQ$ est un $\mZ$-module injectif, la suite exacte 
$0\rightarrow \Gamma_{\co_\oK}\rightarrow \Gamma_A\rightarrow \Gamma_A/\Gamma_{\co_\oK} \rightarrow 0$ est scindée. 
Par suite, $\Gamma_A=\Gamma_{\co_\oK}$ puisque $\Gamma_A$ est sans torsion. 
En vertu de (\cite{gr1} 14.5.23 et 14.5.18), pour tout $A$-module
de type $\alpha$-fini, $\lambda_A(M)$ est donc la longueur définie dans \eqref{mptf9b}.

\subsection{}\label{longnm9}
On désigne par $\Lambda$ le conoyau dans la catégorie des monoïdes 
de l'homomorphisme $\vartheta\colon \mN\rightarrow P$ \eqref{longnm1} 
et par $q\colon P\rightarrow \Lambda$ l'homomorphisme canonique. Posons $\lambda=\vartheta(1)$. 
D'après (\cite{ogus} I 1.1.5), $\Lambda$ est le quotient de $P$ par la relation de congruence $E$ formée des éléments $(x,y)\in P\times P$ 
pour lesquels il existe $a,b\in \mN$ tels que $x+a\lambda=y+b\lambda$. Dire que $E$ est une relation 
de congruence signifie qu'elle est une relation d'équivalence et que $E$ est un sous-monoïde de $P\times P$.  
Le groupe $\Lambda^\gp$ s'identifie canoniquement à $P^\gp/\mZ\lambda$. 
Comme $P$ est intègre, $\Lambda$ est intègre~; on peut donc l'identifier à l'image 
de $P$ dans $P^\gp/\mZ\lambda$.

On désigne par $P_\eta$ la localisation de $P$ par $\lambda$ (\cite{ogus} I 1.4.4), 
qui est aussi la somme amalgamée des deux homomorphismes $\vartheta\colon \mN\rightarrow P$
et $\mN\rightarrow \mZ$, de sorte que le diagramme 
\begin{equation}\label{longnm9a}
\xymatrix{
\mN\ar[r]^\vartheta\ar[d]&P\ar[d]\\
\mZ\ar[r]&{P_\eta}}
\end{equation}
est cocartésien. Par suite, $\Lambda$ est canoniquement isomorphe au conoyau dans la catégorie des monoïdes 
de l'homomorphisme canonique $\mZ\rightarrow P_\eta$. D'après \ref{cad2}(ii), il existe donc deux entiers $0\leq c\leq d$
et un isomorphisme de monoïdes 
\begin{equation}\label{longnm9b}
\mZ^{c}\oplus \mN^{d-c}\stackrel{\sim}{\rightarrow} \Lambda.
\end{equation}
En particulier, le monoïde $\Lambda$ est saturé.

\begin{lem}[\cite{agt} II.8.6]\label{longnm10}
Conservons les notations de \ref{longnm9}. Alors~:
\begin{itemize}
\item[{\rm (i)}] Pour tout $x\in \Lambda$, l'ensemble $q^{-1}(x)$ admet un unique élément minimal $\tx$ 
pour la relation de préordre de $P$ définie par sa structure de monoïde.
\item[{\rm (ii)}] Pour tout $x\in \Lambda$ et tout entier $n\geq 0$, on a $\widetilde{nx}=n\tx$. 
\end{itemize}
\end{lem}

\subsection{}\label{longnm11}
Soit $M$ l'un des monoïdes $\mN$, $P$ ou $\Lambda$. 
Considérons le système inductif de monoïdes $(M^{(n)})_{n\geq 1}$, 
indexé par l'ensemble $\mZ_{\geq 1}$ ordonné par la relation de divisibilité, 
défini par $M^{(n)}=M$ pour tout $n\geq 1$ et dont l'homomorphisme de transition
$i_{n,mn}\colon M^{(n)}\rightarrow M^{(mn)}$ (pour $m, n\geq 1$) est l'endomorphisme 
de Frobenius d'ordre $m$ de $M$ ({\em i.e.}, l'élévation à la puissance $m$-ième). 
On désigne par $M_\infty$ la limite inductive
\begin{equation}\label{longnm11a}
M_\infty=\underset{\underset{n\geq 1}{\longrightarrow}}{\lim}\ M^{(n)}.
\end{equation} 

Pour tout $n\geq 1$, on note
\begin{equation}\label{longnm11b}
P\stackrel{\sim}{\rightarrow}P^{(n)}, \ \ \ t\mapsto t^{(n)},
\end{equation} 
l'isomorphisme canonique. Pour tout $t\in P$ et tous $m, n\geq 1$, on a donc
\begin{equation}\label{longnm11c}
i_{n,mn}(t^{(n)})=(t^{(mn)})^m. 
\end{equation}

En tant que foncteur adjoint à gauche, le foncteur ``groupe associé'' commute aux limites inductives. 
On a donc un isomorphisme canonique 
\begin{equation}\label{longnm11d}
\Lambda^\gp_\infty\stackrel{\sim}{\rightarrow} \Lambda^\gp_\mQ=\Lambda^\gp\otimes_\mZ\mQ. 
\end{equation}

On note aussi $q\colon P_{\infty}\rightarrow \Lambda_{\infty}$ la limite inductive des homomorphismes canoniques 
$q\colon P^{(n)}\rightarrow \Lambda^{(n)}$ \eqref{longnm9}.
En vertu de \ref{longnm10}(ii), les applications $\Lambda^{(n)}\rightarrow P^{(n)}, x\mapsto \tx$ sont compatibles
aux morphismes de transition et définissent donc par passage à la limite inductive une application que l'on note aussi 
\begin{equation}\label{longnm11e}
\Lambda_{\infty}\rightarrow P_{\infty}, \ \ \ x\mapsto \tx.
\end{equation}
On a clairement $q(\tx)=x$.

\subsection{}\label{longnm12}
Avec les notations de \ref{notconv2}, on a un diagramme commutatif de morphismes de schémas logarithmiques
\begin{equation}\label{longnm12a}
\xymatrix{
{(X,\cM_X)}\ar[r]^-(0.5){\gamma^a}\ar[d]_f&{\bA_P}\ar[d]^{\bA_\vartheta}\\
{(S,\cM_S)}\ar[r]^-(0.5){\iota^a}&{\bA_\mN}}
\end{equation}
où les flèches horizontales sont induites par les cartes $\gamma:P\rightarrow \Gamma(X,\cM_X)$ et $\iota\colon \mN\rightarrow \Gamma(S,\cM_S)$ 
(cf. \ref{cad1}). Le morphisme induit de schémas usuels
\begin{equation}\label{longnm12b}
X\rightarrow S\times_{\bA_\mN}\bA_P
\end{equation}
est étale. Il résulte aussitôt des propriétés universelles des localisations de monoïdes et d'anneaux que 
l'homomorphisme canonique $\mZ[P]\rightarrow \mZ[P_\eta]$ \eqref{longnm9a} induit un isomorphisme
\begin{equation}\label{longnm12c}
\mZ[P]_{\lambda}\stackrel{\sim}{\rightarrow} \mZ[P_\eta].
\end{equation}
Compte tenu de \ref{cad2}(ii), on en déduit un $K$-isomorphisme 
\begin{equation}\label{longnm12d}
\mZ[P]\otimes_{\mZ[\mN]}K\stackrel{\sim}{\rightarrow} K[\Lambda].
\end{equation}
Par passage à la limite inductive, on en déduit un $\oK$-isomorphisme 
\begin{equation}\label{longnm12e}
\mZ[P_\infty]\otimes_{\mZ[\mN_\infty]}\oK\stackrel{\sim}{\rightarrow} \oK[\Lambda_\infty],
\end{equation}
où l'on considère $\oK$ comme un $\mZ[\mN_\infty]$-module via le composé des homomorphismes canoniques 
\begin{equation}
\mZ[\mN_\infty]\rightarrow \underset{\underset{n\geq 1}{\longrightarrow}}{\lim}\ \Gamma(S^{(n)},\co_{S^{(n)}})\rightarrow \Gamma(\oS,\co_\oS)=\co_\oK,
\end{equation}
induit par le choix du point géométrique $\oy$ \eqref{longnm1a}.

\subsection{}\label{longnm13}
Soit $n$ un entier $\geq 1$.
D'après \ref{cad7}(i), on a un morphisme étale canonique $X_n\rightarrow S^{(n)}\times_{\bA_{\mN^{(n)}}}\bA_{P^{(n)}}$ 
(cf. \ref{cad3} et \ref{cad6} pour les notations). 
En particulier, $\uR_n$ s'identifie au localisé strict de 
$\oS\times_{\bA_{\mN^{(n)}}}\bA_{P^{(n)}}$ en le point géométrique image canonique de $\oy$ \eqref{longnm3b}.

On désigne par $C_n$ la $\co_\oK$-algèbre du schéma affine $\oS\times_{\bA_{\mN^{(n)}}}\bA_{P^{(n)}}$.
Notant $\lambda^{(n)}$ l'image de $\lambda$ dans $P^{(n)}$ par l'isomorphisme \eqref{longnm11b},  on a 
\begin{equation}\label{longnm13a}
C_n=\co_{\oK}[P^{(n)}]/(\pi_n-e^{\lambda^{(n)}}).
\end{equation}
On désigne par $\alpha_n\colon P^{(n)}\rightarrow C_n$ l'homomorphisme de monoïdes canonique. 
Suivant \ref{longnm10}, pour tout $x\in \Lambda^{(n)}$, on note $\tx\in P^{(n)}$ 
l'unique élément minimal de $q^{-1}(x)$ pour la relation de préordre de $P^{(n)}$ 
définie par sa structure de monoïde. On a un isomorphisme de $\co_{K_n}$-modules
\begin{equation}\label{longnm13b}
C_{n}=\bigoplus_{x\in \Lambda^{(n)}}\co_{\oK}\cdot \alpha_n(\tx).
\end{equation}  
Les algèbres $(C_n)_{n\geq 1}$ forment un système inductif indexé par l'ensemble $\mZ_{\geq 1}$ ordonné par la relation de divisibilité. 
Posons 
\begin{equation}\label{longnm13c}
C_\infty=\underset{\underset{n\geq 1}{\longrightarrow}}{\lim}\ C_{n},
\end{equation}
et notons $\alpha_\infty\colon P_\infty\rightarrow C_\infty$ l'homomorphisme de monoïdes obtenu par passage à la limite inductive 
des homomorphismes $\alpha_n$. Les décompositions \eqref{longnm13b} induisent une décomposition 
\begin{equation}\label{longnm13d}
C_{\infty}=\bigoplus_{x\in \Lambda_\infty}\co_{\oK}\cdot \alpha_\infty(\tx),
\end{equation}  
où $\tx$ désigne l'image de $x$ dans $P_\infty$ par l'application \eqref{longnm11e}.
L'isomorphisme \eqref{longnm12e} induit un isomorphisme de $\oK$-algèbres 
\begin{equation}\label{longnm13e}
C_\infty\otimes_{\co_\oK}\oK\stackrel{\sim}{\rightarrow} \oK[\Lambda_\infty].
\end{equation}

Pour tout $x\in \Lambda_\infty$, il existe $a_x\in \oK$ tel que $\alpha_\infty(\tx)\otimes 1=a_x \cdot x$ via l'isomorphisme \eqref{longnm13e}. 
On définit alors l'application 
\begin{equation}\label{longnm13f}
\psi\colon \Lambda_\infty\rightarrow \mR, \ \ \ x\mapsto v(a_x),
\end{equation}
où $v$ est la valuation de $\oK$ fixée dans \ref{mtfla1}. Il résulte aussitôt de la définition que pour tout $x,y\in \Lambda_\infty$, on a 
\begin{equation}\label{longnm13g}
\psi(x)+\psi(y)\geq \psi(x+y).
\end{equation}
D'après \ref{longnm10}(ii), pour tout $x\in  \Lambda_\infty$ et tout $n\geq 1$, on a 
\begin{equation}\label{longnm13h}
\psi(nx)=n\psi(x).
\end{equation}

\subsection{}\label{longnm15}
Soit $n$ un entier $\geq 1$. 
On identifie le groupe $\Lambda^{(n)\gp}$ au sous-groupe $\frac 1 n \Lambda^\gp$ de $\Lambda^\gp_\mQ$ (cf. \ref{longnm11}). 
Pour tout $\gamma\in \Lambda^\gp_\mQ$, on désigne par $[\gamma]$ la classe de $\gamma$ dans $\Lambda^\gp_\mQ/\Lambda^{(n)\gp}$ et 
par $C_{[\gamma]}$ le $C_n$-module
\begin{equation}\label{longnm15a}
C_{[\gamma]}=\bigoplus_{x\in [\gamma]\cap \Lambda_\infty}\co_{\oK}\cdot \alpha_\infty(\tx).
\end{equation}
La multiplication par $\gamma^{-1}$ dans $\oK[\Lambda^{\gp}_\mQ]$ induit un isomorphisme 
\begin{equation}\label{longnm15b}
C_{[\gamma]}\stackrel{\sim}{\rightarrow}
\bigoplus_{x\in \Lambda^{(n)\gp}\cap (\Lambda_\infty-\gamma)} \{x\otimes a\in \oK[\Lambda^{(n)\gp}]\ |\ v(a)\geq \psi(\gamma+x)\}.
\end{equation}
Pour tout entier $m\geq 1$, comme $\Lambda$ est saturé, 
on a $\Lambda^{(mn)}=\Lambda^{(mn)\gp}\cap \Lambda_\infty$. On en déduit une décomposition canonique
\begin{equation}\label{longnm15c}
C_{mn}=\oplus_{\gamma\in \Lambda^{(mn)\gp}/\Lambda^{(n)\gp}}C_{[\gamma]}.
\end{equation}

La $\co_\oK$-algèbre $\uR_n$ étant mesurable \eqref{longnm16}, on désigne par $\lambda_n=\lambda_{\uR_n}$ la longueur normalisée associée \eqref{longnm17}.
Soit $M$ un $\uR_n$-module de présentation finie à support dans le point fermé de $\Spec(\uR_n)$. On définit l'application 
\begin{equation}\label{longnm15d}
l_M\colon \Lambda_\mQ^\gp\rightarrow \mR, \ \ \ \gamma \mapsto \lambda_n(C_{[\gamma]}\otimes_{C_n}M).
\end{equation}
Posons $\Lambda_\mR^\gp=\Lambda^\gp\otimes_\mZ\mR$. On démontre que $l_M$ s'étend uniquement en une fonction  
$\Lambda_\mR^\gp\rightarrow \mR$ qui se descend en une fonction 
\begin{equation}\label{longnm15e}
\ol_M\colon \Lambda_\mR^\gp/\Lambda^{(n)\gp}\rightarrow \mR.
\end{equation}
La démonstration est une étape de la preuve de (\cite{gr1} 17.1.43) qui ne nécessite aucun changement dans notre situation. 

Fixons une $\mZ$-base $e_1,\dots,e_d$ de $\Lambda^\gp$ \eqref{longnm9b} et posons 
\begin{equation}\label{longnm15f}
D_n=\{ \sum_{i=1}^dt_ie_i \in \Lambda^\gp_\mR \ | \ - \frac 1 {2n} \leq t_i<\frac 1 {2n} \ \forall 1\leq i\leq d\}
\end{equation} 
qui est un domaine fondamental du réseau $\Lambda^{(n)\gp}$ dans $\Lambda^{\gp}_\mR$. 

Soit $m$ un entier $\geq 1$. 
D'après \ref{cad7}(v), $\xi^{(mn)\circ}\colon \oX^{(mn)\circ}\rightarrow \oX^{\circ}$ est un revêtement étale fini galoisien. 
Son groupe de Galois agit naturellement sur $\uoX^{(mn)}$ et il permute transitivement ses composantes irréductibles. 
Compte tenu de \ref{longnm17} et de \eqref{longnm15c}, on a donc
\begin{equation}\label{longnm15g}
\lambda_{mn}(\uR_{mn}\otimes_{\uR_n} M)=\frac{c(n)}{c(mn)}\sum_{\gamma\in \Lambda^{(mn)\gp}\cap D_n}l_M(\gamma),
\end{equation}
où $c(mn)$ désigne le nombre de composantes irréductibles de $\uoX^{(mn)}$ \eqref{longnm3}. 
On notera que les anneaux locaux $\uR_n$ et $\uR_{mn}$ ont même corps résiduel. 

\begin{prop}[\cite{gr1} 17.1.43]\label{longnm18}
Conservons les hypothèses et notations de \ref{longnm15}. Alors,
\begin{itemize}
\item[{\rm (i)}] La fonction $\ol_M$ est intégrable par rapport à la mesure invariante $d\mu$ de $\Lambda_\mR^\gp/\Lambda^{(n)\gp}$ de volume 
total $1$. De plus, on a  
\begin{equation}
\underset{m\to +\infty}{\lim}\ \frac{c(n (m!))}{(n (m!))^d} \lambda_{n (m!)}(\uR_{n (m!)}\otimes_{\uR_n} M)=\frac{c(n)}{n^d}\int_{D_n}\ol_M d\mu.
\end{equation} 
Notons $\lambda_\infty(\uR_\infty\otimes_{\uR_n} M)$ cette limite.
\item[{\rm (ii)}] Pour tout entier $r\geq 1$,  tout nombre réel $\varepsilon >0$ et 
tout idéal de type fini $I$ de $\uR_{1}$ tel que le seul point associé de $\uR_{1}/I$ soit le point fermé de $\Spec(\uR_1)$, 
il existe un nombre réel $\delta>0$
vérifiant la propriété suivante. Pour tout entier $n\geq 1$ et tout morphisme surjectif de $(\uR_n/I\uR_n)$-modules $N\rightarrow N'$, 
où $N$ est un $(\uR_n/I\uR_n)$-module de présentation finie, engendré par au plus $r$ éléments, tel que 
\begin{equation}
|\lambda_\infty(\uR_\infty\otimes_{\uR_n} N)-\lambda_\infty(\uR_\infty\otimes_{\uR_n} N')|\leq \delta,
\end{equation}
on a 
\begin{equation}
\frac{c(n)}{n^d}|\lambda_n(N)-\lambda_n(N')|\leq \varepsilon.
\end{equation}
\end{itemize}
\end{prop}

La preuve est identique à celle de (\cite{gr1} 17.1.43), une fois modifiés les coefficients de normalisation $d_n$ de {\em loc. cit.} en $n^d/c(n)$, 
compte tenu de \eqref{longnm15g}.

\begin{prop}\label{longnm20}
On peut associer canoniquement à tout $\uR_\infty$-module $M$ à support dans le point fermé $\ox_\infty$ de $\Spec(\uR_\infty)$, 
une {\em longueur} $\lambda(M)\in \mR_{\geq 0}\cup\{\infty\}$ 
vérifiant les propriétés suivantes:
\begin{itemize}
\item[{\rm (i)}] Pour toute suite exacte de $\uR_\infty$-modules à supports dans $\ox_\infty$, $0\rightarrow M'\rightarrow M\rightarrow M''\rightarrow 0$, on~a 
\begin{equation}
\lambda(M)= \lambda(M')+\lambda(M'').
\end{equation}
\item[{\rm (ii)}] Pour tout $\uR_\infty$-module $\alpha$-nul $M$, à support dans $\ox_\infty$, on a $\lambda(M)=0$.
\item[{\rm (ii')}] Pour tout morphisme de $\uR_\infty$-modules à supports dans $\ox_\infty$, $M\rightarrow M'$, qui est un $\alpha$-isomor\-phisme,
on a $\lambda(M)=\lambda(M')$. 
\item[{\rm (iii)}] Pour tout $\uR_\infty$-module de type $\alpha$-fini $M$, à support dans $\ox_\infty$, et tout $\gamma\in \fm_\oK$, la longueur $\lambda(\gamma M)$ est finie. 
\item[{\rm (iv)}] Soit $M$ un $\uR_\infty$-module à support dans $\ox_\infty$ et contenu dans un $\uR_\infty$-module de présentation $\alpha$-finie. 
Pour que $M$ soit $\alpha$-nul, il faut et il suffit que $\lambda(M)=0$. 
\end{itemize}
\end{prop}

Cela résulte de \ref{longnm18} et (\cite{gr1} 14.5.75(ii), 14.5.77, 14.5.80(ii), 14.5.83(ii)).

\begin{rema}\label{longnm21}
Posons 
\begin{equation}\label{longnm21a}
\uR_{p^\infty}=\underset{\underset{n\geq 0}{\longrightarrow}}{\lim}\ \uR_{p^n},
\end{equation}
la limite inductive étant indexée par l'ensemble $\mN$ ordonné par la relation d'ordre habituelle. 
Gabber et Ramero montrent dans (\cite{gr1} 17.1.43) l'existence d'une longueur 
normalisée $\lambda_{p^\infty}$ pour les $\uR_{p^\infty}$-modules à support dans le point fermé de $\Spec(\uR_{p^\infty})$. En fait, leurs arguments s'adaptent 
immédiatement au cas des $\uR_\infty$-modules comme expliqué plus haut. 
Notons momentanément $\lambda_\infty$ la longueur normalisée pour les $\uR_\infty$-modules à support dans le point fermé de $\Spec(\uR_{\infty})$ \eqref{longnm20}.

Pour tout $\uR_{p^\infty}$-module de présentation finie $M$ à support dans le point fermé de $\Spec(\uR_{p^\infty})$,  on a 
\begin{equation}\label{longnm21b}
\lambda_\infty(M\otimes_{\uR_{p^\infty}}\uR_\infty)=\lambda_{p^\infty}(M).
\end{equation}
En effet, il existe un entier $n\geq 0$ et un $\uR_{p^n}$-module de présentation finie $N$ à support dans le point fermé de $\Spec(\uR_{p^n})$
tels que $M=N\otimes_{\uR_{p^n}}\uR_{p^\infty}$.
L'assertion résulte alors de \ref{longnm18}(i) et de la formule analogue pour $\uR_{p^\infty}$ (\cite{gr1} 17.1.43). 
\end{rema}

\section{\texorpdfstring{\'Etude locale de certains $\varphi$-modules $\alpha$-étales}{Etude locale de certains phi-modules alpha-etales}}\label{elcpam}

\subsection{}\label{elcpam2}
Les hypothèses et notations de §~\ref{longnm} sont en vigueur dans cette section.
On note, de plus, $\phi_\oK$ l'endomorphisme de Frobenius absolu de $\co_\oK/p\co_\oK$. 
On désigne par $\co_{\oK^\flat}$ la limite projective du système projectif $(\co_\oK/p\co_\oK)_{\mN}$ 
dont les morphismes de transition sont les itérés de $\phi_\oK$; 
\begin{equation}\label{elcpam2a}
\co_{\oK^\flat}= \underset{\underset{\mN}{\longleftarrow}}{\lim}\ \co_\oK/p\co_\oK.
\end{equation}
C'est un anneau de valuation non-discrète, de hauteur $1$, complet et 
parfait de caractéristique $p$ (cf. \ref{TFA101} et \ref{TFA102}). On note $\oK^\flat$ son corps des fractions
et $\fm_{\oK^\flat}$ son idéal maximal. La valuation $p$-adique normalisée $v$ de $C$ \eqref{TFA1}
induit une valuation $v_{\oK^\flat}$ de $\oK^\flat$ \eqref{TFA101d}. 
On fixe une suite $(p_n)_{n\geq 0}$ d'éléments de $\co_\oK$ 
telle que $p_0=p$ et $p_{n+1}^p=p_n$ pour tout $n\geq 0$ et on note $\varpi$ l'élément associé de $\co_{\oK^\flat}$. 
Pour tout entier $i\geq 0$, l'homomorphisme
\begin{equation}\label{elcpam2b}
\co_{\oK^\flat}\rightarrow \co_\oK/p\co_\oK, \ \ \ (x_n)_{n\geq 0}\mapsto x_i
\end{equation}
induit un isomorphisme \eqref{TFA102b}
\begin{equation}\label{elcpam2c}
\co_{\oK^\flat}/\varpi^{p^i}\co_{\oK^\flat}\stackrel{\sim}{\rightarrow}\co_\oK/p\co_\oK.
\end{equation}

Comme $\co_{\oK^\flat}$ satisfait les conditions requises dans \ref{mptf1}, 
il est loisible de considérer les notions de $\alpha$-algèbre introduites dans les sections \ref{alpha}--\ref{aet} relativement à $\co_{\oK^\flat}$.
On observera que pour toute $(\co_\oK/p\co_\oK)$-algèbre, les notions de $\alpha$-algèbre
relativement à $\co_{\oK}$ coïncident avec celles relativement à $\co_{\oK^\flat}$ via l'isomorphisme 
$\co_{\oK^\flat}/\varpi\co_{\oK^\flat}\stackrel{\sim}{\rightarrow}\co_\oK/p\co_\oK$ \eqref{elcpam2c}.

\subsection{}\label{elcpam3}
On note $\phi$ l'endomorphisme de Frobenius absolu de $\uR_\infty/p\uR_\infty$. 
Suivant \eqref{eipo3a}, on désigne par $\uR_\infty^\flat$ la limite projective du système projectif $(\uR_\infty/p\uR_\infty)_{\mN}$ 
dont les morphismes de transition sont les itérés de $\phi$; 
\begin{equation}\label{elcpam3a}
\uR_\infty^\flat= \underset{\underset{\mN}{\longleftarrow}}{\lim}\ \uR_\infty/p\uR_\infty.
\end{equation}
C'est une $\co_{\oK^\flat}$-algèbre parfaite. On note encore $\phi$ l'endomorphisme de Frobenius de $\uR^\flat_\infty$; 
cet abus de notation n'induit aucun risque d'ambiguïté. 
Pour tout entier $i\geq 0$, on désigne par $\pi_i$ l'homomorphisme
\begin{equation}\label{elcpam3b}
\pi_i\colon \uR_\infty^\flat\rightarrow \uR_\infty/p\uR_\infty, \ \ \ (x_n)_{n\geq 0}\mapsto x_i.
\end{equation}

\begin{prop}\label{elcpam4} 
L'endomorphisme de Frobenius absolu $\phi$ de $\uR_\infty/p\uR_\infty$ est surjectif.
\end{prop}

En effet, avec les notations de \ref{AFR24}, on a un isomorphisme canonique \eqref{AFR24n}
\begin{equation} \label{}
\underset{\underset{(U,\fp)\in \fV_\ox^\circ}{\longrightarrow}}{\lim}\ R_{U,\infty}\stackrel{\sim}{\rightarrow}\uR_\infty. 
\end{equation}
Pour tout objet $(U,\fp)$ de $\fV_\ox^\circ$, l'endomorphisme de Frobenius absolu de $R_{U,\infty}/pR_{U,\infty}$ 
est surjectif d'après (\cite{agt} II.9.10); d'où la proposition.

\begin{cor}\label{elcpam6}
Pour tout entier $n\geq 1$ et tout $\gamma\in \fm_\oK$ tels que $p^nv(\gamma)\leq 1$ \eqref{mtfla1}, l'homomorphisme 
\begin{equation}\label{elcpam6a}
\uR_\infty/\gamma\uR_\infty\rightarrow \uR_\infty/\gamma^{p^n}\uR_\infty
\end{equation}
induit par l'endomorphisme $\phi^n$ de $\uR_\infty/p\uR_\infty$ est un isomorphisme.
\end{cor}

En effet, la surjectivité de \eqref{elcpam6a} résulte de \ref{elcpam4}. L'injectivité résulte du fait que l'anneau 
$\uR_\infty$ est local, strictement hensélien normal et en particulier intègre \eqref{longnm3}.

\begin{cor}\label{elcpam5}
Pour tout entier $i\geq 0$, la suite
\begin{equation}\label{elcpam5a}
\xymatrix{
0\ar[r]&{\uR_\infty^\flat}\ar[r]^-(0.5){\cdot \varpi^{p^i}}&{\uR_\infty^\flat}\ar[r]^-(0.5){\pi_i}&{\uR_\infty/p\uR_\infty}\ar[r]&0}
\end{equation}
est exacte. En particulier, $\pi_0$ induit un isomorphisme 
\begin{equation}\label{elcpam5b}
\uR_\infty^\flat/\phi^{-i}(\varpi)\uR_\infty^\flat\stackrel{\sim}{\rightarrow}\uR_\infty/p_i\uR_\infty,
\end{equation}
où $p_i$ est l'élément de $\co_\oK$ fixé dans \ref{elcpam2} tel que $p_i^{p^i}=p$.
\end{cor}

On notera d'abord que $\varpi^{p^i}$ correspond à l'élément $(\op_{n-i})_{n\geq 0}$ via l'égalité \eqref{elcpam2a}, où 
$\op_j$ est la classe de $p_j$ dans $\co_\oK/p\co_\oK$ pour tout $j\geq 0$ et où l'on a posé $\op_j=0$ pour tout $j<0$.
Pour tout entier $n\geq 0$, on a le diagramme commutatif 
\begin{equation}
\xymatrix{
0\ar[r]&{\uR_\infty/pp_{n+1}^{-1}\uR_\infty}\ar[r]^-(0.5){\cdot p_{n+1}}\ar[d]_{\phi}&{\uR_\infty/p\uR_\infty}\ar[r]^{\phi^{n+1}}\ar[d]^\phi&
{\uR_\infty/p\uR_\infty}\ar[r]\ar@{=}[d]&0\\
0\ar[r]&{\uR_\infty/pp_n^{-1}\uR_\infty}\ar[r]^-(0.5){\cdot p_{n}}&{\uR_\infty/p\uR_\infty}\ar[r]^{\phi^{n}}&{\uR_\infty/p\uR_\infty}\ar[r]&0}
\end{equation}
dont les lignes sont exactes d'après \ref{elcpam6}. On notera que le morphisme $\phi$ est bien défini car $pv(pp_{n+1}^{-1})=p-1/p^n\geq 1$. 
On en déduit, compte tenu de \ref{elcpam4}, que la suite \eqref{elcpam5a} est exacte. 
La seconde assertion s'ensuit compte tenu du fait que $\phi^{-i}(\varpi)\in \co_{\oK^\flat}$ correspond à l'élément $(\op_{n+i})_{n\geq 0}$ via l'égalité \eqref{elcpam2a}.

\begin{cor}\label{elcpam500}
L'anneau  $\uR^\flat_\infty$ est complet et séparé pour la topologie $\varpi$-adique, 
et celle-ci coïncide avec la limite des topologies discrètes sur $\uR_\infty/p\uR_\infty$ \eqref{elcpam3a}.
\end{cor}

En effet, pour tout entier $i\geq 0$, on a un diagramme commutatif 
\begin{equation}\label{elcpam500a}
\xymatrix{
{\uR^\flat_\infty/\varpi^{p^{i+1}}\uR^\flat_\infty}\ar[r]^-(0.5){\pi_{i+1}}\ar[d]&{\uR_\infty/p\uR_\infty}\ar[d]^{\phi}\\
{\uR^\flat_\infty/\varpi^{p^{i}}\uR^\flat_\infty}\ar[r]^-(0.5){\pi_i}&{\uR_\infty/p\uR_\infty}}
\end{equation}
où la flèche verticale de gauche est le morphisme canonique. Les flèches horizontales sont des isomorphismes d'après \ref{elcpam5},  
d'où la proposition. On notera que la première assertion résulte aussi de \ref{mdc3}.

\begin{lem}\label{elcpam25}
\
\begin{itemize}
\item[{\rm (i)}] L'anneau $\uR^\flat_\infty$ est intégralement clos dans $\uR^\flat_\infty[\varpi^{-1}]$. 
\item[{\rm (ii)}] Les seuls idempotents de $\uR^\flat_\infty[\varpi^{-1}]$ sont $0$ et $1$.
\end{itemize}
\end{lem}

(i) Soient $a\in \uR_\infty^\flat$, $n$ un entier $\geq 0$ tels que la section $b=\varpi^{-n}a$ de $\uR^\flat_\infty[\varpi^{-1}]$ soit entière sur $\uR^\flat_\infty$. 
Il existe donc un entier $N\geq 1$ et des sections $c_0,\dots,c_{N-1}\in \uR^\flat_\infty$ tels que $b^N+\sum_{i=0}^{N-1}c_ib^i=0$ dans $\uR^\flat_\infty[\varpi^{-1}]$. 
Multipliant par $\varpi^{p^nN}$, on obtient $a^N+\sum_{i=0}^{N-1}\varpi^{p^n(N-i)}c_ia^i=0$ dans $\uR^\flat_\infty$.
Considérons un entier $m\geq 0$ tel que $p^{m-n}\geq N$ de sorte que $v(p_m^{p^nN})\leq v(p)=1$, et des relèvements $\ta_m$ et $\tc_{im}$ dans $\uR_\infty$ 
de $\pi_m(a)$ et $\pi_m(c_i)$ dans $\uR_\infty/p\uR_\infty$ \eqref{elcpam3b}.
Il existe alors $r\in \uR_\infty$ tel que $\ta_m^N+\sum_{i=0}^{N-1}\varpi^{p^n(N-i)}\tc_{im}\ta_m^i+pr=0$ dans $\uR_\infty$.
Multipliant par $p_m^{-p^nN}$ et posant $\tb_m=p_m^{-p^n}\ta_m\in \uR_\infty[\frac 1 p]$, on voit que 
$\tb_m^N+\sum_{i=0}^{N-1}\tc_{im}\tb_m^i+p_m^{-p^nN}pr=0$ dans $\uR_\infty[\frac 1 p]$. Comme $\uR_\infty$ est intégralement clos dans $\uR_\infty[\frac 1 p]$,
on en déduit que $\tb_m\in \uR_\infty$ et par suite que $\ta_m=p_m^{p^n}  \tb_m\in p_m^{p^n} \uR_\infty$. Il s'ensuit que $\pi_n(a)=0$ \eqref{elcpam500a}.
Par suite, on a $a\in \varpi^{p^n}\uR_\infty^\flat$ d'après \ref{elcpam5} et donc $b\in \uR_\infty^\flat$, ce qui achève la preuve de l'assertion.

(ii) En effet, si $e$ est un idempotent de $\uR^\flat_\infty[\varpi^{-1}]$, alors $e\in \uR_\infty^\flat$ d'après (i). 
Compte tenu de \ref{elcpam5} et du fait que l'anneau $\uR_\infty/p\uR_\infty$ est local,
on a $e\in \varpi\uR_\infty^\flat$ ou $e\in 1+ \varpi\uR_\infty^\flat$. 
Dans le premier cas, $e-1$ est inversible d'après \ref{elcpam500} et donc $e=0$.
Dans le second cas, $e$ est inversible et donc $e=1$.

\subsection{}\label{elcpam7}
Pour tout $\gamma\in \fm_\oK$ tel que $pv(\gamma)\leq 1$ \eqref{mtfla1}, on note
\begin{equation}\label{elcpam7a}
\ophi_\gamma\colon \uR_\infty/\gamma\uR_\infty\stackrel{\sim}{\rightarrow} \uR_\infty/\gamma^p\uR_\infty
\end{equation}
l'isomorphisme induit par $\phi$ \eqref{elcpam6}. 

Si $\theta\colon A\rightarrow B$ est un homomorphisme d'anneaux et $M$ est un $A$-module, il peut être commode de noter aussi $\theta^*(M)$
le $B$-module $M\otimes_AB$.

\begin{prop}\label{elcpam8}
Soient $\gamma\in \fm_\oK$ tel que $pv(\gamma)\leq 1$ \eqref{mtfla1}, 
$M$ un $(\uR_\infty/\gamma\uR_\infty)$-module à support dans le point fermé $\ox_\infty$ de $\Spec(\uR_\infty)$. 
Considérant $\ophi^*_\gamma(M)$ comme un $\uR_\infty$-module via la projection canonique
$\uR_\infty\rightarrow \uR_\infty/\gamma^p\uR_\infty$  \eqref{elcpam7}, on a 
\begin{equation}\label{elcpam8b}
\lambda(\ophi^*_\gamma(M))=p^{d+1}\lambda(M),
\end{equation}
où $\lambda$ est la longueur normalisée définie dans \ref{longnm20} et $d=\dim(\uX_s)$. 
\end{prop}

Supposons d'abord que $M$ soit de présentation finie sur $\uR_\infty/\gamma\uR_\infty$. Il existe alors un entier $n\geq 1$
et un $(\uR_n/\gamma\uR_n)$-module de présentation finie $M_n$ à support dans le point fermé de $\Spec(\uR_n)$ tels que $M=M_n\otimes_{\uR_n}\uR_\infty$.
Remplaçant $X$ par $X^{(n)}$ \eqref{cad7}, on peut supposer que $n=1$. Compte tenu de \ref{longnm21}, la proposition résulte alors de l'assertion analogue pour le
$(\uR_{p^\infty}/\gamma\uR_{p^\infty})$-module $M_1\otimes_{\uR_1}\uR_{p^\infty}$ et la longueur normalisée $\lambda_{p^\infty}$, démontrée dans (\cite{gr0} 8.4.15);
on notera que les hypothèses fixées dans \ref{longnm} correspondent à celles de (\cite{gr0} 8.3.19).  

Si $M$ est de type fini sur $\uR_\infty/\gamma\uR_\infty$, on l'écrit comme limite inductive filtrante 
de $(\uR_\infty/\gamma\uR_\infty)$-modules de présentation finie, avec des morphismes de transition surjectifs. 
La proposition se déduit alors du cas des modules de présentation finie traité plus haut, compte tenu de (\cite{gr0} 7.2.62(i)). 

Dans le cas général, on écrit $M$ comme limite inductive filtrante de ses sous-$(\uR_\infty/\gamma\uR_\infty)$-modules de type fini $(M_i)_{i\in I}$. 
Comme $\ophi_\gamma$ est un isomorphisme, les morphismes de transition du système inductif $(\ophi^*_\gamma(M_i))_{i\in I}$ sont injectifs. 
La proposition pour $M$ se déduit alors de l'assertion correspondante pour les $M_i$ démontrée plus haut, compte tenu de (\cite{gr0} 7.2.62(i)). 

\begin{prop}\label{elcpam26}
Si $x$ est un point générique de $X_s$ \eqref{longnm3}, alors $\uR_\infty^\flat$ est un anneau de valuation non discrète de hauteur un.
\end{prop}

Comme $\uR_\infty$ est un anneau de valuation non discrète de hauteur $1$ d'après \ref{longnm8}, 
il en est de même de son séparé complété $p$-adique $\huRi$. 
Soit $\nu\colon \huRi[\frac 1 p]\rightarrow \mQ\cup \{\infty\}$ la valuation normalisée par $\nu(p)=1$.
Considérons le système projectif de monoïdes $(\huRi)_{\mN}$ 
dont les morphismes de transition sont les itérés de l'élévation à la puissance $p$ de $\huRi$. 
D'après \ref{eip5}(i), la réduction modulo $p$ de $\huRi$ induit un isomorphisme de monoïdes multiplicatifs
\begin{equation}\label{elcpam26d}
\underset{\underset{\mN}{\longleftarrow}}{\lim}\ \huRi \stackrel{\sim}{\rightarrow} \uR_\infty^\flat, \ \ \ (x_n)_{n\geq 0}\mapsto (\ox_n)_{n\geq 0}.
\end{equation}
Par suite, $\uR_\infty^\flat$ est un anneau intègre. 
L'isomorphisme \eqref{elcpam26d} se prolonge en un isomorphisme de monoïdes multiplicatifs 
\begin{equation}\label{elcpam26a}
\underset{\underset{\mN}{\longleftarrow}}{\lim}\ \huRi[\frac 1 p] \stackrel{\sim}{\rightarrow} \uR_\infty^\flat[\frac{1}{\varpi}].
\end{equation}
Par conséquent, $\uR_\infty^\flat[\frac{1}{\varpi}]$ est le corps des fractions de $\uR_\infty^\flat$. Notons 
\begin{equation}\label{elcpam26b}
\uR_\infty^\flat[\frac{1}{\varpi}]\rightarrow  \huRi[\frac 1 p],\ \ \ x\mapsto x^\sharp
\end{equation}
l'application composée de l'inverse de \eqref{elcpam26a} et de la projection $(x_n)_{n\geq 0}\mapsto x_0$. On vérifie aussitôt que l'application
\begin{equation}\label{elcpam26c}
\nu^\flat\colon \uR_\infty^\flat[\frac{1}{\varpi}]\rightarrow  \mQ\cup \{\infty\},\ \ \ x\mapsto \nu(x^\sharp)
\end{equation}
est une valuation d'anneau $\uR_\infty^\flat$. Elle est surjective puisqu'elle prolonge la valuation $v_{\oK^\flat}$ de $\oK^\flat$ \eqref{TFA101d}.

\begin{lem}\label{elcpam13}
Tout $\uR^\flat_\infty$-module $\cM$ séparé pour la topologie $\varpi$-adique 
tel que $\cM/\varpi \cM$ soit de type $\alpha$-fini (resp. $\alpha$-nul), est de type $\alpha$-fini (resp. $\alpha$-nul).
\end{lem}

Soient $r$ un entier $\geq 0$, $\psi\colon (\uR_\infty^\flat)^r\rightarrow \cM$
un morphisme $\uR_\infty^\flat$-linéaire, $\psi_1\colon (\uR_\infty^\flat/\varpi \uR_\infty^\flat)^r\rightarrow \cM/\varpi \cM$ le morphisme induit par $\psi$. 
Supposons que $\coker(\psi_1)$ soit annulé par un élément $\gamma\in \fm_{\oK^\flat}$ tel que $v_{\oK^\flat}(\gamma)<v_{\oK^\flat}(\varpi)$. Alors, la suite 
\begin{equation}\label{elcpam13a}
(\uR_\infty^\flat)^r\stackrel{\psi}{\rightarrow} \cM \stackrel{u}{\rightarrow}\coker(\psi_1)\rightarrow 0,
\end{equation}
où $u$ est le morphisme canonique, est exacte.
En effet, $u$ est clairement surjectif. Soit $x\in \cM$ tel que $u(x)=0$.
On montre par récurrence que pour tout entier $n\geq 1$, qu'il existe $y_n\in (\uR_\infty^\flat)^d$ et $x_n\in \cM$ tels que $x=\psi(y_n)+\gamma^{1-n}\varpi^n x_n$
et $y_{n+1}-y_n\in \gamma^{-n}\varpi^n (\uR_\infty^\flat)^r$. Comme $\uR^\flat_\infty$ est complet et séparé pour la topologie $\varpi$-adique 
\eqref{elcpam500} et $\cM$
est séparé pour la topologie $\varpi$-adique, on en déduit qu'il existe $y\in (\uR_\infty^\flat)^r$ tel que $x=\psi(y)$; d'où l'assertion recherchée.
Par suite, $\coker(\psi)$ est annulé par $\gamma$, ce qui implique la proposition.

\begin{lem}\label{elcpam20}
Soient $\cM, \cN$ deux $\uR^\flat_\infty$-modules tels que $\cM$ soit complet et séparé pour la topologie $\varpi$-adique et 
que $\cN$ soit $\alpha$-projectif de type $\alpha$-fini. Notons $\hcN$ et $\cM\hotimes_{\uR^\flat_\infty}\cN$
les séparés complétés $\varpi$-adiques de $\cN$ et $\cM\otimes_{\uR^\flat_\infty}\cN$. 
Alors, le morphisme canonique $\iota \colon \cM\otimes_{\uR^\flat_\infty}\cN\rightarrow  
\cM\hotimes_{\uR^\flat_\infty}\cN$ est un $\alpha$-isomorphisme. 
En particulier, le morphisme canonique $\cN\rightarrow \hcN$ est un $\alpha$-isomorphisme. 
\end{lem}

En effet, pour tout $\gamma \in \fm_{\oK^\flat}$, il existe $a\geq 1$ et deux morphismes $\uR^\flat_\infty$-linéaires 
\begin{equation}
\cN\stackrel{u}{\longrightarrow} (\uR^\flat_\infty)^a  \stackrel{v}{\longrightarrow} \cN
\end{equation}
tels que $v\circ u$ soit égal à la multiplication par $\gamma$. Notons  
\begin{equation}
\cM\otimes_{\uR^\flat_\infty}\cN\stackrel{\mu}{\longrightarrow} \cM^{\oplus a}  \stackrel{\nu}{\longrightarrow} \cM\otimes_{\uR^\flat_\infty}\cN
\end{equation}
les morphismes déduits de $u$ et $v$. Comme  $\cM$ est complet et séparé pour la topologie $\varpi$-adique, on a 
$\gamma\cdot (\cap_{n\geq 0}\varpi^n (\cM\otimes_{\uR^\flat_\infty}\cN))=0$. De plus, pour toute suite de Cauchy $(x_n)_{n\geq 0}$ de 
$\cM\otimes_{\uR^\flat_\infty}\cN$,
la suite $(\mu(x_n))_{n\geq 0}$ converge vers un élément $y$ de $\cM^{\oplus a}$. La suite $(\gamma x_n)_{n\geq 0}$ converge donc vers 
$\nu(y)\in \cM\otimes_{\uR^\flat_\infty}\cN$. 
Par suite, $\coker(\iota)$ est annulé par $\gamma$. Il est donc $\alpha$-nul, d'où la première assertion de la proposition. 
La seconde assertion s'ensuit puisque $\uR^\flat_\infty$ est complet et séparé pour la topologie $\varpi$-adique \eqref{elcpam500}.

\begin{prop}\label{elcpam9}
Soient $\gamma\in \fm_\oK$ tel que $pv(\gamma)\leq 1$  \eqref{mtfla1}, $M$ un $(\uR_\infty/\gamma^p\uR_\infty)$-module, 
\begin{equation}\label{elcpam9a}
\varphi\colon \ophi^*_\gamma(M/\gamma M)\rightarrow M
\end{equation}
un morphisme $\uR_\infty$-linéaire \eqref{elcpam7}. Supposons les conditions suivantes satisfaites:
\begin{itemize}
\item[{\rm (i)}] $\dim(\uX_s)\geq 1$;
\item[{\rm (ii)}] $\varphi$ est un $\alpha$-isomorphisme;
\item[{\rm (iii)}] $M$ est contenu dans un $\uR_\infty$-module de présentation $\alpha$-finie $N$; 
\item[{\rm (iv)}]  le support de $N$ est contenu dans le point fermé $\ox_\infty$ de $\Spec(\uR_\infty)$.
\end{itemize}
Alors, le $\uR_\infty$-module $M$ est $\alpha$-nul.
\end{prop}

En effet,  comme $\ophi^*_\gamma$ est un isomorphisme \eqref{elcpam6},
pour tout $\delta\in \fm_\oK$ tel que $v(\delta)\leq v(\gamma)$, le morphisme $\varphi$ \eqref{elcpam9a} induit un morphisme 
$\uR_\infty$-linéaire
\begin{equation}\label{elcpam9b}
\ophi^*_\gamma(\delta M/\gamma M)\rightarrow \delta^p M
\end{equation}
qui est un $\alpha$-isomorphisme. Compte tenu de \ref{elcpam8} et \ref{longnm20}(ii'), on en déduit que 
\begin{equation}\label{elcpam9c}
p^{d+1}\lambda(\delta M/\gamma M)= \lambda(\delta^p M),
\end{equation}
où $d=\dim(\uX_s)$. 
Filtrant $\delta^p M$ par les sous-$\uR_\infty$-modules $(\gamma^i\delta^{p-i} M)_{0\leq i\leq p}$,
pour tout $0\leq i\leq p-1$, la multiplication par $\gamma^i\delta^{p-i-1}$ dans $M$ induit un morphisme surjectif 
\begin{equation}\label{elcpam9d}
\delta M/\gamma M\rightarrow (\gamma^i\delta^{p-i}M)/(\gamma^{i+1}\delta^{p-i-1}M).
\end{equation}
Par suite, on a, en vertu de \ref{longnm20},
\begin{equation}\label{elcpam9e}
p^{d+1}\lambda(\delta M/\gamma M)=\lambda(\delta^p M)\leq p \lambda(\delta M/\gamma M).
\end{equation}
On notera que la longueur $\lambda(\delta^p M)$ est finie d'après \ref{longnm20}(i)-(iii).
Comme $p^{d+1}>p$, il s'ensuit que $\lambda(\delta M/\gamma M)=0$ et donc que $\lambda(\delta^p M)=0$.
Par suite, $\delta^p M$ est $\alpha$-nul en vertu de \ref{longnm20}(iv) et il en est alors de même de $M$. 

\begin{rema}[\cite{gr1} 14.5.80(i)]\label{elcpam900}
Soient $A$ un anneau, $0\rightarrow M'\rightarrow M\rightarrow M''\rightarrow 0$ une suite exacte de $A$-modules, $a,b\in A$. 
Alors, $abM$ est une extension d'un sous-module de $aM'$ par un quotient de $bM''$. En effet, la multiplication par $a$ dans $M$
induit un morphisme $A$-linéaire surjectif 
\begin{equation}
bM''=bM/(M'\cap bM)\twoheadrightarrow abM/a(M'\cap bM),
\end{equation}
et on a $a(M'\cap bM)\subset aM'$. 
\end{rema}

\begin{prop}\label{elcpam90}
Soient $\gamma\in \fm_\oK$ tel que $pv(\gamma)\leq 1$  \eqref{mtfla1}, $M$ un $(\uR_\infty/\gamma^p\uR_\infty)$-module, 
\begin{equation}\label{elcpam90a}
\varphi\colon \ophi^*_\gamma(M/\gamma M)\rightarrow M
\end{equation}
un morphisme $\uR_\infty$-linéaire \eqref{elcpam7}. Supposons les conditions suivantes satisfaites:
\begin{itemize}
\item[{\rm (i)}] $\dim(\uX_s)\geq 1$;
\item[{\rm (ii)}] $\varphi$ est un $\alpha$-isomorphisme;
\item[{\rm (iii)}] $M$ est de présentation $\alpha$-finie sur $\uR_\infty$;
\item[{\rm (iv)}] le $(\co_\oK/\gamma^p\co_\oK)$-module $M$ est $\alpha$-plat en dehors du point fermé ${\ox_\infty}$ de $\Spec(\uR_\infty/p\uR_\infty)$ \eqref{decprim4}, 
{\em i.e.}, pour tout $\fq\in \Spec(\uR_\infty/p\uR_\infty)$ 
tel que $\fq\not=\ox_\infty$, le $(\co_\oK/\gamma^p\co_\oK)$-module $M_\fq$ est $\alpha$-plat {\rm (\cite{agt} V.6.1)}.
\end{itemize}
Alors, le $\uR_\infty$-module $\Gamma_{\ox_\infty}(M)$ des sections de $M$ à support dans $\ox_\infty$, est $\alpha$-nul.
\end{prop}

En effet, comme $\ophi_\gamma$ est un isomorphisme \eqref{elcpam6}, pour tout $\delta\in \fm_\oK$ tel que $v(\delta)\leq v(\gamma)$, 
$\varphi$ \eqref{elcpam90a} induit un morphisme $\uR_\infty$-linéaire
\begin{equation}\label{elcpam90b}
\ophi^*_\gamma(\delta \Gamma_{\ox_\infty}(M/\gamma M))\rightarrow \delta^p \Gamma_{\ox_\infty}(M)
\end{equation} 
qui est un $\alpha$-isomorphisme. On en déduit, compte tenu de \ref{elcpam8} et \ref{longnm20}(ii'), que
\begin{equation}\label{elcpam90c}
p^{d+1}\lambda(\delta \Gamma_{\ox_\infty}(M/\gamma M))= \lambda(\delta^p \Gamma_{\ox_\infty}(M)),
\end{equation}
où $d=\dim(\uX_s)$. 
Filtrant $M$ par les sous-$\uR_\infty$-modules $(\gamma^i M)_{0\leq i\leq p}$, pour tout $0\leq i\leq p-1$, on a une suite exacte 
\begin{equation}\label{elcpam90d}
0\rightarrow \Gamma_{\ox_\infty}(\gamma^{i+1}M)\rightarrow \Gamma_{\ox_\infty}(\gamma^iM)\rightarrow \Gamma_{\ox_\infty}(\gamma^iM/\gamma^{i+1}M).
\end{equation}
Notant $\oGamma_{\ox_\infty}(\gamma^iM)$ l'image de la dernière flèche, on en déduit compte tenu de \ref{elcpam900} et  \ref{longnm20},
\begin{equation}\label{elcpam90e}
\lambda(\delta^p \Gamma_{\ox_\infty}(M)) \leq 
\sum_{0\leq i\leq p-1}\lambda (\delta \oGamma_{\ox_\infty}(\gamma^iM)) \leq 
\sum_{0\leq i\leq p-1}\lambda (\delta \Gamma_{\ox_\infty}(\gamma^iM/\gamma^{i+1}M)).
\end{equation}
Pour tout $0\leq i\leq p-1$, la multiplication par $\gamma^i$ dans $M$ induit un morphisme surjectif 
\begin{equation}\label{elcpam90f}
M/\gamma M\rightarrow \gamma^iM/\gamma^{i+1}M.
\end{equation}
Compte tenu de (iv), celui-ci est $\alpha$-injectif au-dessus de $\Spec(\uR_\infty/p\uR_\infty)-\{\ox_\infty\}$. 
Par suite, le morphisme induit
\begin{equation}\label{elcpam90g}
\Gamma_{\ox_\infty}(M/\gamma M)\rightarrow \Gamma_{\ox_\infty}(\gamma^iM/\gamma^{i+1}M)
\end{equation}
est $\alpha$-surjectif. 
En effet, pour tout $u\in \Gamma_{\ox_\infty}(\gamma^iM/\gamma^{i+1}M)$, il existe $v \in M/\gamma M$ qui le relève \eqref{elcpam90f}. 
Pour tout $t\in \fm_\oK$, $t v$ est nul au-dessus $\Spec(\uR_\infty/p\uR_\infty)-\{\ox_\infty\}$. 
Par suite, $t v\in \Gamma_{\ox_\infty}(M/\gamma M)$, d'où l'assertion. On en déduit, compte tenu de \eqref{elcpam90e} et \ref{longnm20}, que
\begin{equation}\label{elcpam90h}
p^{d+1}\lambda(\delta \Gamma_{\ox_\infty}(M/\gamma M))=\lambda(\delta^p \Gamma_{\ox_\infty}(M))\leq p \lambda(\delta \Gamma_{\ox_\infty}(M/\gamma M)).
\end{equation}
En vertu de \ref{decprim9}, le $\uR_\infty$-module $\Gamma_{\ox_\infty}(M)$ est de présentation $\alpha$-finie. 
La longueur $\lambda(\delta^p \Gamma_{\ox_\infty}(M))$ est donc finie d'après \ref{longnm20}(iii).
Comme $p^{d+1}>p$, il s'ensuit que $\lambda(\delta^p \Gamma_{\ox_\infty}(M))=0$.
Par suite, $\delta^p \Gamma_{\ox_\infty}(M)$ est $\alpha$-nul en vertu de \ref{longnm20}(iv),
et il en est alors de même de $\Gamma_{\ox_\infty}(M)$.

\begin{lem}\label{elcpam27}
Supposons que $x$ soit un point générique de $X_s$ \eqref{longnm3}. 
Soient $\gamma\in \fm_\oK$ tel que $pv(\gamma)\leq 1$  \eqref{mtfla1}, 
$M$ un $(\uR_\infty/\gamma^p\uR_\infty)$-module de type $\alpha$-fini, 
\begin{equation}\label{elcpam27a}
\varphi\colon \ophi^*_\gamma(M/\gamma M)\rightarrow M
\end{equation}
un morphisme $\uR_\infty$-linéaire \eqref{elcpam7} qui est un $\alpha$-isomorphisme. 
Alors, le $(\uR_\infty/\gamma^p\uR_\infty)$-module $M$ est $\alpha$-projectif. 
\end{lem}

D'après \ref{longnm8}, $\uR_\infty$ est un anneau de valuation non discrète de hauteur $1$,
ayant une valuation $\nu\colon \uR_\infty\rightarrow \mR_{\geq 0}\cup\{\infty\}$ qui prolonge la valuation $v$ de $\co_{\oK}$. 
Pour tout $\varepsilon \in\mR_{\geq 0}$, 
notons $I_\varepsilon$ l'idéal de $\uR_\infty$ formé des éléments $x\in \uR_\infty$ tels que $\nu(x)>\varepsilon$. 
En vertu de \ref{mptf7}, il existe une unique suite décroissante de nombres réels positifs $(\varepsilon_i)_{i\geq 1}$ tendant vers $0$ telle que 
\begin{equation}\label{elcpam27b}
M\approx \oplus_{i\geq 1}\uR_\infty/I_{\varepsilon_i}. 
\end{equation}
Soit $n$ le sup de $0$ et des entiers $i\geq 1$ tel que $\varepsilon_i\geq v(\gamma)$. On a alors 
\begin{equation}\label{elcpam27c}
M/\gamma M\approx (\uR_\infty/\gamma \uR_\infty)^n\oplus \oplus_{i> n}\uR_\infty/I_{\varepsilon_i}. 
\end{equation}
D'après \ref{mptf7} et compte tenu de \eqref{elcpam27a}, 
on a $\varepsilon_i=pv(\gamma)$ pour tout $1\leq i\leq n$ et $p \varepsilon_i = \varepsilon_i$ pour tout $i >n$. 
Par suite, $M\approx (\uR_\infty/\gamma^p \uR_\infty)^n$, d'où la proposition (\cite{agt} V.3.1).

\begin{defi}\label{elcpam12}
On appelle {\em $\varphi$-$\uR_\infty^\flat$-module étale} la donnée d'un $\uR_\infty^\flat$-module $\cM$ et d'un isomorphisme $\uR_\infty^\flat$-linéaire
\begin{equation}\label{elcpam12a}
\phi^*(\cM)\stackrel{\sim}{\rightarrow} \cM,
\end{equation}
dit structure de Frobenius \eqref{elcpam3}. 
\end{defi}

Soient $\cM_1$ et $\cM_2$ deux $\varphi$-$\uR_\infty^\flat$-modules étales. Un morphisme de $\cM_1$ dans $\cM_2$
est un morphisme $\uR_\infty^\flat$-linéaire $\cM_1\rightarrow \cM_2$ compatible aux structures de Frobenius. 

Comme l'endomorphisme de Frobenius $\phi$ de $\uR_\infty^\flat$ est un isomorphisme, 
si $u\colon \cM_1\rightarrow \cM_2$ est un morphisme de $\varphi$-$\uR_\infty^\flat$-modules étales, 
les $\uR_\infty^\flat$-modules $\ker(u)$ et $\coker(u)$ 
sont naturellement munis de structures de Frobenius qui en font des $\varphi$-$\uR_\infty^\flat$-modules étales. 

Si $\cM_1$ et $\cM_2$ sont deux $\varphi$-$\uR_\infty^\flat$-modules étales, les $\uR_\infty^\flat$-modules $\cM_1\otimes_{\uR_\infty^\flat}\cM_2$
et $\Ext^i_{\uR_\infty^\flat}(\cM_1,\cM_2)$ $(i\geq 0$) sont naturellement munis de structures de Frobenius qui en font des $\varphi$-$\uR_\infty^\flat$-modules étales. 

\begin{remas}\label{elcpam18}
Soient $\cM$ un $\varphi$-$\uR_\infty^\flat$-module étale \eqref{elcpam12}, $M=\cM/\varpi \cM$ que l'on considère comme un 
$(\uR_\infty/p\uR_\infty)$-module \eqref{elcpam5a}, $\gamma\in \fm_\oK$ tel que $pv(\gamma)\leq 1$ \eqref{mtfla1}.
\begin{itemize}
\item[(i)] La structure de Frobenius de $\cM$ induit un isomorphisme $\uR_\infty$-linéaire \eqref{elcpam7}
\begin{equation}\label{elcpam18a}
\ophi_\gamma^*(M/\gamma M)\stackrel{\sim}{\rightarrow}M/\gamma^pM.
\end{equation}
En effet, il suffit d'établir l'assertion lorsque $\gamma=p_1$ est la racine $p$-ième de $p$ fixée dans \ref{elcpam2}. 
Ceci résulte du diagramme commutatif
\begin{equation}\label{elcpam18b}
\xymatrix{
{\uR^\flat_\infty/\phi^{-1}(\varpi)\uR^\flat_\infty}\ar[r]\ar[d]&{\uR^\flat_\infty/\varpi \uR^\flat_\infty}\ar[d]\\
{\uR_\infty/p_1\uR_\infty}\ar[r]&{\uR_\infty/p\uR_\infty}}
\end{equation}
où les flèches verticales sont les isomorphismes \eqref{elcpam5b} et les flèches horizontales sont les isomorphismes induits par les endomorphismes 
de Frobenius. 
\item[(ii)] L'élévation à la puissance $p$ dans $\fm_\oK$ induit un isomorphisme $(\co_\oK/p\co_\oK)$-semi-linéaire
\begin{equation}\label{elcpam18c}
\fm_\oK/\gamma\fm_\oK\stackrel{\sim}{\rightarrow}\fm_\oK/\gamma^p\fm_\oK.
\end{equation}
\item[(iii)] On déduit de (i) et (ii) un isomorphisme $\uR_\infty$-linéaire canonique
\begin{equation}\label{elcpam18d}
\ophi_\gamma^*(\fm_\oK\otimes_{\co_\oK}(M/\gamma M))\stackrel{\sim}{\rightarrow} (\fm_\oK\otimes_{\co_\oK}(M/\gamma^pM)).
\end{equation}
\end{itemize}
\end{remas}

\begin{lem}\label{elcpam19}
Soient $\cM$ un $\varphi$-$\uR_\infty^\flat$-module étale \eqref{elcpam12}, $M=\cM/\varpi \cM$, $\cM[\varpi]$ le noyau de la multiplication par 
$\varpi$ dans $\cM$. On suppose les conditions suivantes satisfaites:
\begin{itemize}
\item[{\rm (i)}] $\dim(\uX_s)\geq 1$;
\item[{\rm (ii)}] le $\uR_\infty$-module $\cM[\varpi]$ est $\alpha$-nul;
\item[{\rm (iii)}] le $\uR_\infty$-module $M$ est de présentation $\alpha$-finie; 
\item[{\rm (iv)}] le $(\uR_\infty/p\uR_\infty)$-module $M$ est $\alpha$-projectif en dehors du point fermé $\ox_\infty$ de $\Spec(\uR_\infty/p\uR_\infty)$, 
{\em i.e.}, pour tout $\fq\in \Spec(\uR_\infty/p\uR_\infty)$ 
tel que $\fq\not=\ox_\infty$, le $(\uR_\infty/p\uR_\infty)_\fq$-module $M_\fq$ est $\alpha$-projectif {\rm (\cite{agt} V.3.2)}. 
\end{itemize}
Alors, le $\uR_\infty^\flat$-module $\cM^\vee=\Hom_{\uR_\infty^\flat}(\cM,\uR_\infty^\flat)$ est complet et séparé pour
la topologie $\varpi$-adique, $\varpi$ n'est pas diviseur de zéro dans $\cM^\vee$ et pour tout entier $n\geq 1$, le morphisme canonique 
\begin{equation}\label{elcpam19a}
\cM^\vee/\varpi^n\cM^\vee\rightarrow \Hom_{\uR_\infty^\flat}(\cM,\uR^\flat_\infty/\varpi^n\uR^\flat_\infty)
\end{equation}
est un $\alpha$-isomorphisme. 
\end{lem} 

On notera d'abord que le $\uR_\infty^\flat$-module $\cM^\vee$ est naturellement muni d'une structure de Frobenius qui en fait un 
$\varphi$-$\uR_\infty^\flat$-module étale. 
Comme $\varpi$ n'est pas diviseur de zéro dans $\uR^\flat_\infty$ \eqref{elcpam5}, pour tout entier $n\geq 1$, la suite canonique
\begin{equation}\label{elcpam19b}
0\longrightarrow \cM^\vee \stackrel{\varpi^n}{\longrightarrow} \cM^\vee\longrightarrow \Hom_{\uR^\flat_\infty}(\cM,\uR^\flat_\infty/\varpi^n \uR^\flat_\infty)
\end{equation}
est exacte. Par ailleurs, $\uR_\infty^\flat$ étant complet et séparé pour la topologie $\varpi$-adique, 
le morphisme canonique 
\begin{equation}\label{elcpam19c}
\Hom_{\uR_\infty^\flat}(\cM,\uR_\infty^\flat)\rightarrow  \underset{\underset{n\geq 1}{\longleftarrow}}{\lim}\ 
\Hom _{\uR_\infty^\flat}(\cM,\uR_\infty^\flat/\varpi^n \uR_\infty^\flat)
\end{equation}
est bijectif. Par suite, le morphisme canonique 
\begin{equation}\label{elcpam19d}
\cM^\vee \rightarrow \underset{\underset{n\geq 1}{\longleftarrow}}{\lim}\ \cM^\vee/\varpi^n \cM^\vee
\end{equation}
est un isomorphisme, autrement dit, $\cM^\vee$ est complet et séparé pour la topologie $\varpi$-adique.

Comme $\varpi$ n'est pas diviseur de zéro dans $\uR_\infty^\flat$, on a $\Tor_1^{\uR^\flat_\infty}(\cM,\uR_\infty^\flat/\varpi \uR_\infty^\flat)=\cM[\varpi]$
et pour tout $j\geq 2$, $\Tor_j^{\uR^\flat_\infty}(\cM,\uR_\infty^\flat/\varpi \uR_\infty^\flat)=0$. Considérons la suite spectrale
\begin{equation}
\rE_2^{i,j}=\Ext^i_{\uR_\infty/p\uR_\infty}(\Tor_j^{\uR^\flat_\infty}(\cM,\uR_\infty^\flat/\varpi \uR_\infty^\flat),\uR_\infty/p\uR_\infty)\Rightarrow
\Ext^{i+j}_{\uR_\infty^\flat}(\cM,\uR_\infty^\flat/\varpi \uR_\infty^\flat).
\end{equation}
Le $\uR_\infty$-module $E_2^{0,1}$ est $\alpha$-nul, compte tenu de (ii). On en déduit que le morphisme canonique
\begin{equation}
\Ext^1_{\uR_\infty/p\uR_\infty}(M,\uR_\infty/p \uR_\infty)\rightarrow \Ext^1_{\uR_\infty^\flat}(\cM,\uR_\infty^\flat/\varpi \uR_\infty^\flat)
\end{equation}
est un $\alpha$-isomorphisme.  De même, on a une suite exacte 
\begin{equation}
0\rightarrow \Ext^1_{\uR_\infty/p\uR_\infty}(M,\uR_\infty/p \uR_\infty)\rightarrow \Ext^1_{\uR_\infty}(M,\uR_\infty/p \uR_\infty)
\rightarrow \Hom_{\uR_\infty}(M,\uR_\infty/p \uR_\infty).
\end{equation}
Par suite, le $\uR_\infty$-module $\Ext^1_{\uR_\infty^\flat}(\cM,\uR_\infty^\flat/\varpi \uR_\infty^\flat)$ est de présentation 
$\alpha$-finie compte tenu de (iii) et en vertu de \ref{longnm5}. 
Par ailleurs,  compte tenu de (iv) et \ref{longnm50}, le $\uR_\infty$-module
$\fm_\oK\otimes_{\co_\oK}\Ext^1_{\uR_\infty^\flat}(\cM,\uR_\infty^\flat/\varpi \uR_\infty^\flat)$ est à support dans $\ox_\infty$.

Le $\uR_\infty^\flat$-module $\cE=\Ext^1_{\uR_\infty^\flat}(\cM,\uR_\infty^\flat)$ est naturellement muni d'une structure de Frobenius qui en fait un 
$\varphi$-$\uR_\infty^\flat$-module étale.
Le $\uR_\infty^\flat$-module $\cE/\varpi\cE$ s'injecte dans $\Ext^1_{\uR_\infty^\flat}(\cM,\uR_\infty^\flat/\varpi \uR_\infty^\flat)$. 
Par suite, $\fm_\oK\otimes_{\co_\oK}(\cE/\varpi\cE)$ est $\alpha$-nul en vertu de \ref{elcpam9} et \ref{elcpam18}(iii). 
Il en est alors de même de $\cE/\varpi\cE$. 
Pour tout entier $n\geq 1$, la suite exacte canonique 
\[
0\longrightarrow \uR_\infty^\flat \stackrel{\varpi^n}{\longrightarrow} \uR_\infty^\flat
\longrightarrow \uR_\infty^\flat/\varpi^n \uR_\infty^\flat\rightarrow 0
\] 
induit un diagramme commutatif à lignes exactes  
\[
\xymatrix{
{\Hom_{\uR_\infty^\flat}(\cM,\uR_\infty^\flat)}\ar[r]\ar[d]&{\Hom_{\uR_\infty^\flat}(\cM,\uR_\infty^\flat/\varpi^n \uR_\infty^\flat)}\ar[r]\ar@{=}[d]&
{\Ext^1_{\uR_\infty^\flat}(\cM,\uR_\infty^\flat)}\ar[d]\\
{\Hom_{\uR_\infty^\flat}(\cM,\uR_\infty^\flat/\varpi^{n+1} \uR_\infty^\flat)}\ar[r]^-(0.5)u&{\Hom_{\uR_\infty^\flat}(\cM,\uR_\infty^\flat/\varpi^n \uR_\infty^\flat)}
\ar[r]^-(0.5)\partial&{\Ext^1_{\uR_\infty^\flat}(\cM,\uR_\infty^\flat/\varpi \uR_\infty^\flat)}}
\]
où les flèches verticales sont les morphismes canoniques. Par suite, le morphisme $\partial$ se factorise à travers $\cE/\varpi \cE$.
Il est donc $\alpha$-nul, et le morphisme $u$ est donc $\alpha$-surjectif. Compte tenu de l'isomorphisme 
\eqref{elcpam19c}, on en déduit que le morphisme canonique 
\begin{equation}\label{elcpam19g}
\cM^\vee/\varpi^n \cM^\vee\rightarrow \Hom_{\uR^\flat_\infty}(\cM,\uR^\flat_\infty/\varpi^n \uR^\flat_\infty)
\end{equation}
est $\alpha$-surjectif. 

\begin{rema}
Il ressort de la preuve de \ref{elcpam19} que pour tout $\uR_\infty^\flat$-module $\cM$, 
le $\uR_\infty^\flat$-module $\cM^\vee=\Hom_{\uR_\infty^\flat}(\cM,\uR_\infty^\flat)$ est complet et séparé pour
la topologie $\varpi$-adique.
\end{rema}

\begin{prop}\label{elcpam10}
Soient $\cM$ un $\varphi$-$\uR_\infty^\flat$-module étale \eqref{elcpam12}, $M=\cM/\varpi \cM$, $\cM[\varpi]$ le noyau de la multiplication par 
$\varpi$ dans $\cM$. On suppose les conditions suivantes satisfaites:
\begin{itemize}
\item[{\rm (i)}] $\dim(\uX_s)\geq 1$;
\item[{\rm (ii)}] le $\uR_\infty$-module $\cM[\varpi]$ est $\alpha$-nul;
\item[{\rm (iii)}] le $\uR_\infty^\flat$-module $\cM$ est séparé pour la topologie $\varpi$-adique;
\item[{\rm (iv)}] le $\uR_\infty$-module $M$ est de présentation $\alpha$-finie;
\item[{\rm (v)}] le $(\uR_\infty/p\uR_\infty)$-module $M$ est $\alpha$-projectif en dehors du point fermé $\ox_\infty$ de $\Spec(\uR_\infty/p\uR_\infty)$, 
{\em i.e.}, pour tout $\fq\in \Spec(\uR_\infty/p\uR_\infty)$ 
tel que $\fq\not=\ox_\infty$, le $(\uR_\infty/p\uR_\infty)_\fq$-module $M_\fq$ est $\alpha$-projectif {\rm (\cite{agt} V.3.2)}.
\end{itemize}
Alors, le $\uR_\infty^\flat$-module $\cM$ est $\alpha$-projectif de type $\alpha$-fini.
\end{prop} 

Compte tenu de (ii), pour tout entier $n\geq 1$, le morphisme canonique
\begin{equation}\label{elcpam10b}
\End_{\uR^\flat_\infty}(\cM)/\varpi^n \End_{\uR^\flat_\infty}(\cM)\rightarrow \End_{\uR^\flat_\infty}(\cM/\varpi^n\cM)
\end{equation}
est $\alpha$-injectif. Par ailleurs, d'après (iii), le morphisme canonique 
\begin{equation}\label{elcpam10c}
\End_{\uR^\flat_\infty}(\cM)\rightarrow  \underset{\underset{n\geq 1}{\longleftarrow}}{\lim}\ 
\End_{\uR^\flat_\infty}(\cM/\varpi^n\cM)
\end{equation}
est injectif. Par suite, le morphisme canonique 
\begin{equation}\label{elcpam10d}
\End_{\uR^\flat_\infty}(\cM) \rightarrow \underset{\underset{n\geq 1}{\longleftarrow}}{\lim}\ \End_{\uR^\flat_\infty}(\cM)/\varpi^n \End_{\uR^\flat_\infty}(\cM)
\end{equation}
est injectif, {\em i.e.}, $\End_{\uR^\flat_\infty}(\cM)$ est séparé pour la topologie $\varpi$-adique. 

Les $\uR_\infty^\flat$-modules $\cM^\vee=\Hom_{\uR_\infty^\flat}(\cM,\uR_\infty^\flat)$ et $\End_{\uR^\flat_\infty}(\cM)$ 
sont naturellement munis de structures de Frobenius qui en font des $\varphi$-$\uR_\infty^\flat$-modules étales. 
Il en est donc de même de $\cM\otimes_{\uR_\infty^\flat}\cM^\vee$. Le morphisme canonique 
\begin{equation}\label{elcpam10e}
\psi\colon \cM\otimes_{\uR_\infty^\flat}\cM^\vee \rightarrow \End_{\uR_\infty^\flat}(\cM)
\end{equation}
est compatible aux structures de Frobenius. 
On désigne par $\cN$ son conoyau, qui est donc naturellement muni d'une structure de Frobenius qui en fait un $\varphi$-$\uR_\infty^\flat$-module étale.
On note $N$ le conoyau du morphisme canonique
\begin{equation}\label{elcpam10f}
M\otimes_{\uR_\infty}\Hom_{\uR_\infty}(M,\uR_\infty/p \uR_\infty) \rightarrow \End_{\uR_\infty}(M). 
\end{equation}
D'après (iv) et \ref{longnm5}, les $\uR_\infty$-modules $\Hom_{\uR_\infty}(M,\uR_\infty/p \uR_\infty)$ et 
$\End_{\uR_\infty}(M)$ sont de présentation $\alpha$-finie et il en est alors de même de $N$.
Il résulte de (v), \ref{longnm50} et (\cite{agt} V.4.1) que le morphisme \eqref{elcpam10f} est un $\alpha$-isomorphisme en dehors de $\ox_\infty$.
Par suite, le $\uR_\infty$-module $\fm_{\co_\oK}\otimes_{\co_\oK}N$ est à support dans $\ox_\infty$.  
D'après \ref{elcpam19},  le morphisme canonique 
\begin{equation}\label{elcpam10a}
\cM^\vee/\varpi\cM^\vee\rightarrow \Hom_{\uR_\infty}(M,\uR_\infty/p\uR_\infty)
\end{equation}
est $\alpha$-surjectif. 
Compte tenu de la $\alpha$-injectivité du morphisme canonique \eqref{elcpam10b}
\begin{equation}\label{elcpam10g}
\End_{\uR_\infty^\flat}(\cM)/\varpi \End_{\uR_\infty^\flat}(\cM)\rightarrow \End_{\uR_\infty}(M),
\end{equation}
on en déduit que le morphisme canonique 
\begin{equation}\label{elcpam10h}
\cN/\varpi\cN\rightarrow N
\end{equation}
est $\alpha$-injectif. Par suite, $\fm_{\co_\oK}\otimes_{\co_\oK}(\cN/\varpi\cN)$ est contenu dans  $\fm_{\co_\oK}\otimes_{\co_\oK}N$.
On en déduit que $\fm_{\co_\oK}\otimes_{\co_\oK}(\cN/\varpi\cN)$ 
est $\alpha$-nul en vertu de \ref{elcpam9} et \ref{elcpam18}(iii), et il en est de même de $\cN/\varpi \cN$. 

Montrons que la suite 
\begin{equation}\label{elcpam10i}
\cM\otimes_{\uR_\infty^\flat}\cM^\vee \stackrel{\psi}{\rightarrow} \End_{\uR_\infty^\flat}(\cM) \stackrel{\tau}{\rightarrow}  \cN/\varpi \cN\rightarrow 0,
\end{equation}
où $\tau$ est le morphisme canonique, est $\alpha$-exacte. Il est clair que $\tau$ est surjectif.  
Soient $\gamma\in \fm_{\oK^\flat}$ tel que $v_{\oK^\flat}(\gamma)<v_{\oK^\flat}(\varpi)$, $u\in \End_{\uR_\infty^\flat}(\cM)$ tel que $\tau(u)=0$.
Le $\uR_\infty^\flat$-module $\cM$ étant de type $\alpha$-fini d'après \ref{elcpam13},  
il existe un entier $r\geq 0$ et un morphisme $h\colon (\uR_\infty^\flat)^r\rightarrow \cM$  dont le conoyau est annulé par $\gamma$. 
Comme $\gamma (\cN/\varpi \cN)=0$, on montre par récurrence que pour tout entier $n\geq 1$, 
il existe $\zeta_n\in (\cM^\vee)^{\oplus r}$ et $u_n\in \End_{\uR_\infty^\flat}(\cM)$ tels que 
$\gamma u=\psi\circ (h\otimes \id_{\cM^\vee})(\zeta_n)+\gamma^{2-n}\varpi^n u_n$ et $\zeta_{n+1}-\zeta_n\in \gamma^{-n}\varpi^n (\cM^\vee)^{\oplus r}$.
Le $\uR_\infty^\flat$-module $\cM^\vee$ étant complet et séparé pour la topologie $\varpi$-adique d'après \ref{elcpam19}, 
la suite $(\zeta_n)_n$ converge vers un élément $\zeta\in (\cM^\vee)^{\oplus r}$. 
Comme $\End_{\uR^\flat_\infty}(\cM)$ est séparé pour la topologie $\varpi$-adique \eqref{elcpam10d}, 
on a $\gamma u=\psi\circ(h\otimes \id_{\cM^\vee})(\zeta)$. La suite \eqref{elcpam10i} est donc $\alpha$-exacte.

Le module $\cN/\varpi \cN$ étant $\alpha$-nul, $\psi$ est donc $\alpha$-surjectif \eqref{elcpam10i}. 
Pour tout $\delta\in \fm_{\oK^\flat}$, il existe alors $\zeta_1,\dots,\zeta_n\in \cM$
et $\rho_1,\dots,\rho_n\in \cM^\vee$ tels que $\psi(\sum_{i=1}^n \zeta_i\otimes \rho_i)=\delta \id_\cM$, d'où la proposition (\cite{agt} V.3.1(iv)).

\begin{teo}\label{elcpam15}
Soient $M$ un $(\uR_\infty/p\uR_\infty)$-module de présentation $\alpha$-finie, 
$(p_n)_{n\geq 0}$ la suite d'éléments de $\co_\oK$ fixée dans \ref{elcpam2}, 
\begin{equation}\label{elcpam15a}
\varphi\colon M/p_1 M\rightarrow M
\end{equation}
un morphisme $\ophi_{p_1}$-linéaire tel que le morphisme $\uR_\infty$-linéaire induit 
\begin{equation}\label{elcpam15b}
\tvarphi\colon \ophi^*_{p_1}(M/p_1M)\rightarrow M
\end{equation}
soit un $\alpha$-isomorphisme \eqref{elcpam7}.
Alors, il existe une $(\uR_\infty/p\uR_\infty)$-algèbre $\alpha$-finie-étale et $\alpha$-fidè\-lement plate $A$,  
un entier $r\geq 0$ et un morphisme $A$-linéaire 
\begin{equation}\label{elcpam15c}
\fm_{\oK}\otimes_{\co_{\oK}}M\otimes_{\uR_\infty}A\rightarrow A^r
\end{equation}
compatible à $\varphi$ et aux endomorphismes de Frobenius de $\co_{\oK}/p\co_{\oK}$ et $A$, qui est un $\alpha$-isomorphisme.
\end{teo}

En effet, l'endomorphisme de Frobenius $\phi_\oK$ de $\co_\oK/p\co_\oK$ induit un isomorphisme
\begin{equation}\label{elcpam15aa}
\ophi_\oK\colon \co_\oK/p_1\co_\oK\stackrel{\sim}{\rightarrow} \co_\oK/p\co_\oK,
\end{equation}
et un isomorphisme $\ophi_\oK$-linéaire
\begin{equation}\label{elcpam15ba}
\opsi_\oK\colon \fm_\oK/p_1\fm_\oK\stackrel{\sim}{\rightarrow} \fm_\oK/p\fm_\oK.
\end{equation}
On en déduit un isomorphisme $\ophi_{p_1}$-linéaire \eqref{elcpam7}
\begin{equation}\label{elcpam15ca}
(\fm_\oK/p_1\fm_\oK)\otimes_{\co_\oK}(M/p_1M)\stackrel{\sim}{\rightarrow} 
(\fm_\oK/p\fm_\oK)\otimes_{\co_\oK}\ophi_{p_1}^*(M/p_1M),
\end{equation}
et par suite un isomorphisme $\uR_\infty$-linéaire
\begin{equation}\label{elcpam15da}
\ophi_{p_1}^*((\fm_\oK\otimes_{\co_\oK}M)/p_1(\fm_\oK\otimes_{\co_\oK}M))\stackrel{\sim}{\rightarrow} 
\fm_\oK\otimes_{\co_\oK}\ophi_{p_1}^*(M/p_1M).
\end{equation}
Remplaçant $M$ par $\fm_\oK\otimes_{\co_\oK}M$, on peut donc supposer que $\tvarphi$ \eqref{elcpam15b} est un isomorphisme,
de sorte que $\varphi$ \eqref{elcpam15a} est un isomorphisme $\ophi_{p_1}$-linéaire. 

Procédant par récurrence sur $d=\dim(\uX_s)$, 
on peut supposer la proposition établie en tout point géométrique de $X$ au-dessus de $s$ dont le support dans $X_s$ 
est de codimension $\leq d-1$, cette condition étant clairement satisfaite si $d=0$. 

On rappelle que $\uoX$ est canoniquement isomorphe au localisé strict de $\oX$ en $\ox$ \eqref{longnm3}.
Les points géométriques de $\uoX$ correspondent donc aux points géométriques de $\oX$ qui se spécialisent en $\ox$. 
De plus, les localisés stricts de $\uoX$ et de $\oX$ en des points géométriques qui se correspondent 
sont canoniquement isomorphes (\cite{sga4} VIII 7.3 et 7.4). Soient $\ox^\iota$ un point géométrique de $\oX$ au-dessus de $s$,
$\xi\colon \ox^\iota\rightarrow \uoX$ un $\oX$-morphisme. Nous affecterons d'un exposant $\iota$ les objets associés à $\ox^\iota$,
analogues à ceux associés à $\ox$ dans \ref{longnm3}. Le morphisme $\xi$ induit un $\oX$-morphisme 
que l'on note encore $\xi \colon \uoX^\iota\rightarrow \uoX$. 
Pour tout entier $n\geq 1$, on a un diagramme cartésien
\begin{equation}\label{elcpam15d}
\xymatrix{
{\uoX^{\iota(n)}}\ar[r]\ar[d]&{\uoX^{(n)}}\ar[d]\\
{\uoX^\iota}\ar[r]^{\xi}&{\uoX}}
\end{equation}
Compte tenu de \ref{longnm4}, on peut supposer qu'il existe un point $(\oy\rightsquigarrow \ox^\iota)$ de $X_\et\gtimes_{X_\et}\oX^\circ_\et$
tel que le composé de $\oy\rightsquigarrow \ox^\iota$ et $\xi\colon \ox^\iota \rightsquigarrow \ox$ soit le point $(\oy\rightsquigarrow \ox)$
donné dans \ref{longnm3}. Le morphisme \eqref{longnm3b} se relève donc en un morphisme 
\begin{equation}\label{elcpam15e}
\oy\rightarrow \underset{\underset{n\geq 1}{\longleftarrow}}{\lim}\ \uoX^{\iota (n)}.
\end{equation} 
On en déduit un morphisme plat 
\begin{equation}
\xi_\infty\colon \Spec(\uR_\infty^\iota)\rightarrow \Spec(\uR_\infty).
\end{equation}
Rappelons que le groupe profini $\pi_1(\uoX^\circ,\oy)$ agit sur l'anneau $\uR_\infty$ \eqref{longnm3}. 
Pour tout entier $n\geq 1$, $\uoX^{(n)\star \circ}\rightarrow \uoX^\circ$ est un revêtement étale galoisien connexe d'après \ref{cad7}(v).
Par ailleurs, le point fermé $\ox_\infty$ de $\Spec(\uR_\infty)$ est le seul point de $\Spec(\uR_\infty)$ au-dessus du point fermé $\ox$ de $\uoX$.  
On en déduit que la collection des morphismes $g^*\circ \xi_\infty\colon \Spec(\uR_\infty^\iota)\rightarrow \Spec(\uR_\infty)$ 
pour $g\in \pi_1(\uoX^\circ,\oy)$ et $\xi\colon \ox^\iota\rightarrow \uoX$ décrivant  
les points géométriques de $\uoX-\{\ox\}$ au-dessus de $s$, couvre $\Spec(\uR_\infty/p\uR_\infty)-\{\ox_\infty\}$. 

Par hypothèse de récurrence et descente $\alpha$-fidèlement plate (\cite{agt} V.8.7), on déduit de ce qui précède 
que le $(\uR_\infty/p\uR_\infty)$-module $M$ est $\alpha$-projectif de type $\alpha$-fini en dehors de $\ox_\infty$, 
{\em i.e.}, pour tout $\fq\in \Spec(\uR_\infty/p\uR_\infty)-\{\ox_\infty\}$, 
le $(\uR_\infty/p\uR_\infty)_\fq$-module $M_\fq$ est $\alpha$-projectif de type $\alpha$-fini.

Pour tout $\gamma\in \fm_\oK$ tel que $pv(\gamma)\leq 1$, $\varphi$ \eqref{elcpam15a} induit un isomorphisme $\ophi_\gamma$-linéaire 
\begin{equation}\label{elcpam15h}
M/\gamma M\rightarrow M/\gamma^p M.
\end{equation}
Pour tout entier $n\geq 1$, on en déduit un isomorphisme $(\ophi_{p_1}\circ\dots\circ \ophi_{p_n})$-linéaire 
\begin{equation}\label{elcpam15l}
M/p_n M\rightarrow M.
\end{equation}
Par conséquent, la suite 
\begin{equation}\label{elcpam15m}
\xymatrix{
0\ar[r]&{M/pp_n^{-1}M}\ar[r]^-(0.5){p_n}&{M}\ar[r]^-(0.5){\varphi^{(n)}}&M\ar[r]&0},
\end{equation}
où $\varphi^{(n)}$ est induit par le morphisme \eqref{elcpam15l}, est exacte au centre et à droite. 
Montrons qu'elle est $\alpha$-exacte. 
La multiplication par $p_n$ dans $\uR_\infty$ induit un morphisme injectif
\begin{equation}
\uR_\infty/pp_n^{-1}\uR_\infty\rightarrow \uR_\infty/p\uR_\infty.
\end{equation}
Supposons d'abord $\dim(\uX_s)=0$, autrement dit que $x$ est un point générique de $X_s$ \eqref{longnm3}. 
En vertu de \ref{elcpam27}, le $(\uR_\infty/p\uR_\infty)$-module $M$ est $\alpha$-plat. 
Par suite, la première flèche de \eqref{elcpam15m} est $\alpha$-injective. 
Supposons ensuite $\dim(\uX_s)\geq 1$.
Il résulte de l'hypothèse de récurrence que la première flèche de \eqref{elcpam15m} 
est $\alpha$-injective au-dessus de $\Spec(\uR_\infty/p\uR_\infty)-\{\ox_\infty\}$. 
Par ailleurs, $\Gamma_{\ox_\infty}(M/pp_n^{-1}M)$ est $\alpha$-nul en vertu de \ref{elcpam90}. 
Par suite, la première flèche de \eqref{elcpam15m} est $\alpha$-injective. 
Dans les deux cas, la suite \eqref{elcpam15m} est donc $\alpha$-exacte.

On désigne par $\cM$ le $\uR_\infty^\flat$-module limite projective du système $(M)_\mN$ dont les morphismes de transition sont
les itérés de l'endomorphisme $\phi$-semi-linéaire de $M$ induit par $\varphi$ \eqref{elcpam15a};
\begin{equation}\label{elcpam15f}
\cM= \underset{\underset{\mN}{\longleftarrow}}{\lim}\ M.
\end{equation}
D'après \ref{mdc3}, $\cM$ est complet et séparé pour la topologie $\varpi$-adique.  
Le morphisme $\varphi$ induit un endomorphisme $\phi$-semi-linéaire  
\begin{equation}\label{elcpam15g}
\Phi\colon \cM\rightarrow \cM,
\end{equation}
qui est clairement bijectif.

Pour tout entier $i\geq 0$, notons $\tau_i\colon \cM\rightarrow M$ le morphisme qui associe à $(x_i)_{i\in \mN}\in \cM$ l'élément $x_i\in M$.
Celui-ci est $\uR_\infty^\flat$-linéaire lorsque l'on considère $M$ comme un $\uR^\flat_\infty$-module via l'homomorphisme $\pi_i$ \eqref{elcpam3b}. 
Il est aussi surjectif puisque $\varphi$ \eqref{elcpam15a} est surjectif. Par ailleurs, il résulte par passage à la limite de \eqref{elcpam15m} que la suite 
\begin{equation}\label{elcpam15i}
0\longrightarrow \cM\stackrel{\varpi^{p^i}}{\longrightarrow} \cM\stackrel{\tau_i}{\longrightarrow} M\longrightarrow 0
\end{equation}
est $\alpha$-exacte.
Notons que $\varpi^{p^i}$ correspond à l'élément $(\op_{n-i})_{n\geq 0}$ via l'égalité \eqref{elcpam2a}, où on a posé $\op_j=0$ pour tout $j<0$. 

Montrons que le $\uR^\flat_\infty$-module $\cM$ est $\alpha$-projectif de type $\alpha$-fini.
Supposons d'abord $\dim(\uX_s)=0$, autrement dit que $x$ est un point générique de $X_s$. 
Le $\uR^\flat_\infty$-module $\cM$ est de type $\alpha$-fini d'après \ref{elcpam13} et \eqref{elcpam15i}. 
D'après  \ref{elcpam26},  $\uR_\infty^\flat$ est un anneau de valuation non discrète de hauteur un.
Notons $\cN$ le sous-module des éléments de torsion de $\cM$ ({\em i.e.}, des éléments annulés par une puissance de $\varpi$). 
Celui-ci est $\alpha$-nul d'après \eqref{elcpam15i}.
En vertu de (\cite{scholze1} 2.8(ii)), il existe un entier $r\geq 0$ tel que $\cM/\cN\approx (\uR^\flat_\infty)^r$ (cf. \ref{mptf3}). 
En particulier, le $\uR^\flat_\infty$-module $\cM/\cN$ est $\alpha$-projectif de type $\alpha$-fini (\cite{agt} V.3.1), et il en est alors de même de $\cM$.
Supposons ensuite $\dim(\uX_s)\geq 1$. En vertu de \ref{elcpam10} et compte tenu de \eqref{elcpam15i} et de l'hypothèse de récurrence, 
le $\uR^\flat_\infty$-module $\cM$ est $\alpha$-projectif de type $\alpha$-fini.

D'après \ref{elcpam25}(ii), le $\uR_\infty^\flat[\varpi^{-1}]$-module $\cM[\varpi^{-1}]=\cM\otimes_{\uR_\infty^\flat}\uR_\infty^\flat[\varpi^{-1}]$ 
est localement libre de rang $r\geq 0$. 
Il existe alors un morphisme étale fini et fidèlement plat $\uR_\infty^\flat[\varpi^{-1}]\rightarrow \cC$ et un isomorphisme $\cC$-linéaire 
\begin{equation}\label{elcpam15o}
u\colon \cM\otimes_{\uR_\infty^\flat}\cC\stackrel{\sim}{\rightarrow}\cC^r
\end{equation}
compatible à $\Phi$ \eqref{elcpam15g} et à l'endomorphisme de Frobenius $\phi_{\cC}$ de $\cC$. 
La preuve est la partie principale de celle de (\cite{katz} 4.1.1); 
on notera que les hypothèses dans {\em loc. cit.} que le schéma est intègre, réduit et normal sont superflues. 

Comme $\cC$ est étale sur $\uR_\infty^\flat[\varpi^{-1}]$, le diagramme d'homomorphismes d'algèbres 
\begin{equation}\label{elcpam15p}
\xymatrix{
{\cC}\ar[r]^{\phi_{\cC}}&{\cC}\\
{\uR_\infty^\flat[\varpi^{-1}]}\ar[u]\ar[r]^{\phi}&{\uR_\infty^\flat[\varpi^{-1}]}\ar[u]}
\end{equation}
est cocartésien. Cela résulte du fait qu'un morphisme étale radiciel et surjectif de schémas est un isomorphisme (\cite{ega4} IV 17.9.1).  
Comme $\phi$ est un isomorphisme de $\uR_\infty^\flat[\varpi^{-1}]$, $\phi_{\cC}$ est un isomorphisme, autrement dit, $\cC$ est une 
$\mF_p$-algèbre parfaite.

Soit $\cA$ la clôture intégrale de $\uR_\infty^\flat$ dans $\cC$. 
En vertu de (\cite{gr} 3.5.28 et 3.5.38(i)), la $\uR_\infty^\flat$-algèbre $\cA$ est $\alpha$-finie-étale. 
On notera qu'une algèbre $\alpha$-étale et uniformément presque finie projective dans le sens de (\cite{gr} 3.1.1) est $\alpha$-finie-étale \eqref{aet2}. 
En particulier, la $\uR_\infty^\flat$-algère $\cA$ est $\alpha$-projective de type $\alpha$-fini et de rang fini.
D'après (\cite{agt} V.5.16), il existe un entier $r\geq 0$ et une décomposition de la $\uR_\infty^\flat$-algèbre
\begin{equation}
\sigma_*((\uR_\infty^\flat)^\alpha)=R_0\times R_1\times \dots \times R_r
\end{equation}
tels que pour tout $0\leq i\leq r$, la $R_i$-algèbre $\cA\otimes_{\uR_\infty^\flat}R_i$ soit de rang $i$; cf. \eqref{alpha5c} et (\cite{agt} V.5.14). 
Inversant $\varpi$, on obtient une décomposition $\uR_\infty^\flat[\varpi^{-1}]=\prod_{i=0}^rR_i[\varpi^{-1}]$ telle que 
$\cA\otimes_{\uR_\infty^\flat}R_i[\varpi^{-1}]$ soit de rang $i$. D'après \ref{elcpam25}(ii), il existe alors $0\leq \rho\leq r$ tel que 
$R_i[\varpi^{-1}]=0$ pour tout $i\not=\rho$. Comme $\varpi$ n'est pas un diviseur de zéro dans $\uR_\infty^\flat$ \eqref{elcpam5a}, la suite 
\begin{equation}
\xymatrix{
0\ar[r]&{\sigma_*((\uR_\infty^\flat)^\alpha)}\ar[r]^-(0.5){\cdot \varpi}&{\sigma_*((\uR_\infty^\flat)^\alpha)}\ar[r]&{\sigma_*((\uR_\infty^\flat/\varpi \uR_\infty^\flat)^\alpha)}}
\end{equation}
est exacte. On en déduit que $R_i=0$ pour tout $i\not=\rho$. 
Appliquant la dernière remarque au-dessus de (\cite{agt} V.5.15) au $\alpha$-isomorphisme $\uR^\flat_\infty\rightarrow R_\rho$,  
on voit que la $\uR_\infty^\flat$-algèbre $\cA$ est $\alpha$-projective de type $\alpha$-fini et de rang $\rho$. 
Comme l'homomorphisme $\uR_\infty^\flat[\varpi^{-1}]\rightarrow \cC$ est fidèlement plat, on a $\rho>0$. 
La $\uR_\infty^\flat$-algèbre $\cA$ est donc $\alpha$-fidèlement plate en vertu de (\cite{agt} V.6.7). 

Pour tout $a\in \cA$, il existe $c\in \cC$ tel que $c^p=a$.
Comme $a$ est entier sur $\uR_\infty^\flat$, il en est de même de $c$. Par suite, $c\in \cA$. La $\mF_p$-algèbre $\cA$ est donc parfaite.
Posons $\cP$ l'image canonique de $\cM\otimes_{\uR_\infty^\flat}\cA$ dans $\cM\otimes_{\uR_\infty^\flat}\cC$. 
Comme le $\uR_\infty^\flat$-module $\cM$ est $\alpha$-plat, le morphisme canonique 
$\cM\otimes_{\uR_\infty^\flat}\cA\rightarrow \cP$ est un $\alpha$-isomorphisme. 
Le $\cA$-module $u(\cP)$ \eqref{elcpam15o} étant de type $\alpha$-fini, il existe $\gamma\in \fm_{\oK^\flat}$ tel que 
\begin{equation}\label{elcpam15q}
\gamma \cA^r\subset u(\cP)\subset \gamma^{-1} \cA^r.
\end{equation} 
Par ailleurs, comme $\Phi$ \eqref{elcpam15g} est un isomorphisme et 
que l'algèbre $\cA$ est parfaite, on a $\phi_{\cC^r}(u(\cP))=u(\cP)$. 
Appliquant les puissances négatives de $\phi_{\cC^r}$ aux inclusions \eqref{elcpam15q}, on en déduit que 
\begin{equation}\label{elcpam15r}
u(\fm_{\oK^\flat}\otimes_{\co_{\oK^\flat}}\cP)\subset \cA^r
\end{equation}
et que le morphisme 
\begin{equation}\label{elcpam15s}
\fm_{\oK^\flat}\otimes_{\co_{\oK^\flat}}\cP\rightarrow \cA^r
\end{equation}
induit par $u$ est un $\alpha$-isomorphisme. Par suite, l'isomorphisme $u$ induit un morphisme $\cA$-linéaire 
\begin{equation}\label{elcpam15t}
\fm_{\oK^\flat}\otimes_{\co_{\oK^\flat}}\cM\otimes_{\uR_\infty^\flat}\cA\rightarrow \cA^r
\end{equation}
qui est un $\alpha$-isomorphisme et est compatible à $\Phi$ et à l'endomorphisme de Frobenius de $\cA$; d'où la proposition.

\begin{cor}\label{elcpam23}
Conservons les hypothèses et notations de \ref{elcpam15}, notons, de plus,  
$\cM$ le $\uR_\infty^\flat$-module limite projective du système projectif $(M)_\mN$ dont les morphismes de transition sont
les itérés de l'endomorphisme $\phi$-semi-linéaire de $M$ induit par $\varphi$ \eqref{elcpam15a}, 
et $\Phi$ l'endomorphisme $\phi$-semi-linéaire de $\cM$ induit par $\varphi$. 
Alors, 
\begin{itemize}
\item[{\rm (i)}] Pour tout entier $i\geq 0$, la suite 
\begin{equation}\label{elcpam23b}
0\longrightarrow \cM\stackrel{\varpi^{p^i}}{\longrightarrow} \cM\stackrel{\tau_i}{\longrightarrow} M\longrightarrow 0,
\end{equation}
où $\tau_i$ est le morphisme qui associe à $(x_i)_{i\in \mN}\in \cM$ l'élément $x_i\in M$, est $\alpha$-exacte. 
\item[{\rm (ii)}] Il existe une $\uR^\flat_\infty$-algèbre $\alpha$-finie-étale et $\alpha$-fidèlement plate $\cA$,  
un entier $r\geq 0$ et un morphisme $\cA$-linéaire 
\begin{equation}\label{elcpam23a}
\fm_{\oK^\flat}\otimes_{\co_{\oK^\flat}}\cM\otimes_{\uR^\flat_\infty}\cA\rightarrow \cA^r
\end{equation}
compatible à $\Phi$ et à l'endomorphisme de Frobenius de $\cA$, qui est un $\alpha$-isomorphisme.
\end{itemize}
\end{cor}

Cela résulte de la preuve de \ref{elcpam15}.

\section{Suite exacte d'Artin-Schreier du topos de Faltings annelé}\label{TCFR}

\subsection{}\label{TCFR1}
Les hypothèses et notations de \ref{longnm1} et \ref{longnm3} sont en vigueur dans cette section.
On reprend aussi les notations de \ref{elcpam2} et \ref{elcpam3}.
On désigne simplement par $\uoR$ la $\uR_1$-algèbre $\oR_{\uX}^\oy$ définie dans \eqref{TFA12f}. 
On rappelle qu'on a un isomorphisme canonique de $\uR_1$-algèbres \eqref{TFA11c}
\begin{equation}\label{TCFR1a}
\ocB_{\rho(\oy\rightsquigarrow \ox)}\stackrel{\sim}{\rightarrow} \uoR, 
\end{equation}  
où $\ocB$ est l'anneau de $\tE$ défini dans \ref{mtfla8} et $\rho$ est le morphisme \eqref{mtfla7e}.

\begin{lem}\label{TCFR2}
La $\uR_\infty$-algèbre $\uoR$ est $\alpha$-fidèlement plate.
\end{lem}
Cela résulte de \ref{aet11}, \ref{tpcg21} et (\cite{agt} V.6.5(2)), compte tenu des définitions \eqref{TFA12f} et \eqref{longnm3d}.

\subsection{}\label{TCFR3}
On note encore $\phi$ l'endomorphisme de Frobenius absolu de $\uoR/p\uoR$. Suivant \eqref{eipo3a},
on désigne par $\uoR^\flat$ la limite projective du système projectif $(\uoR/p\uoR)_{\mN}$ 
dont les morphismes de transition sont les itérés de $\phi$; 
\begin{equation}\label{TCFR3a}
\uoR^\flat= \underset{\underset{\mN}{\longleftarrow}}{\lim}\ \uoR/p\uoR.
\end{equation}
C'est une $\uR^\flat_{\infty}$-algèbre parfaite \eqref{elcpam3a}. 
On note encore $\phi$ l'endomorphisme de Frobenius de $\uoR^\flat$. 
Ces abus de notation n'induisent aucun risque d'ambiguïté. 

Pour tout entier $i\geq 0$, on désigne par $\pi_i$ l'homomorphisme
\begin{equation}\label{TCFR3b}
\pi_i\colon \uoR^\flat\rightarrow \uoR/p\uoR, \ \ \ (x_n)_{n\geq 0}\mapsto x_i.
\end{equation}

\begin{prop}\label{TCFR4} \
\begin{itemize}
\item[{\rm (i)}] Pour tout $\gamma\in \fm_\oK$ tel que $p v(\gamma)\leq 1$ \eqref{mtfla1}, l'homomorphisme 
\begin{equation}\label{TCFR4a}
\ophi_\gamma\colon \uoR/\gamma\uoR\rightarrow \uoR/\gamma^{p}\uoR
\end{equation}
induit par l'endomorphisme $\phi$ de $\uoR/p\uoR$ est un isomorphisme.
\item[{\rm (ii)}] Pour tout entier $i\geq 0$, la suite
\begin{equation}\label{TCFR4b}
\xymatrix{
0\ar[r]&{\uoR^\flat}\ar[r]^-(0.5){\cdot \varpi^{p^i}}&{\uoR^\flat}\ar[r]^-(0.5){\pi_i}&{\uoR/p\uoR}\ar[r]&0}
\end{equation}
est exacte. 
\item[{\rm (iii)}] Le morphisme $\pi_0$ induit un isomorphisme 
\begin{equation}\label{TCFR4c}
\uoR^\flat/\phi^{-i}(\varpi)\uoR^\flat\stackrel{\sim}{\rightarrow}\uoR/p_i\uoR,
\end{equation}
où $p_i$ est l'élément de $\co_\oK$ fixé dans \ref{elcpam2} tel que $p_i^{p^i}=p$.
\item[{\rm (iv)}] L'anneau  $\uoR^\flat$ est complet et séparé pour la topologie $\varpi$-adique,
et celle-ci coïncide avec la limite des topologies discrètes sur $\uoR/p\uoR$ \eqref{TCFR3a}.
\end{itemize}
\end{prop}

En effet, l'endomorphisme de Frobenius de $\uoR/p\uoR$ est surjectif en vertu de \eqref{TFA12f} et (\cite{agt} II.9.10). 
Par suite, $\ophi_\gamma$ est surjectif. L'injectivité résulte du fait que l'anneau 
$\uoR$ est local, strictement hensélien normal et en particulier intègre d'après (\cite{agt} III.10.10(i)). 
Ceci démontre l'assertion (i). Les autres assertions s'en déduisent comme pour $\uR_\infty$ (cf. \ref{elcpam5} et \ref{elcpam500}).  

\begin{lem}\label{TCFR15}
Pour tout element non-nul $b$ de $\co_{\oK^\flat}$, la suite 
\begin{equation}
\xymatrix{
0\ar[r]&{\mF_p}\ar[r]^-(0.5){b}&{\uoR^\flat}\ar[rr]^{\phi-b^{p-1}\id}&&{\uoR^\flat}\ar[r]&0}
\end{equation}
est exacte.  
\end{lem}

Soit $(b_n)_{n\geq 0}\in (\co_C)^{\mN}$ l'image inverse de $b$ par l'isomorphisme \eqref{TFA101b}. 
Il existe $n_0\geq 0$ tel que pour tout $n\geq n_0$, $b_n^p\not\in p\co_C$. 
D'après \ref{TPCF10} et \eqref{TCFR1a}, la suite 
\begin{equation}
\xymatrix{
0\ar[r]&{\mF_p \cdot b_n\oplus (\frac{p}{b_n^{p-1}})\cdot \uoR/p\uoR}\ar[r]&{\uoR/p\uoR}\ar[rr]^{\phi-b_n^{p-1}\id}&&{\uoR/p\uoR}\ar[r]&0}
\end{equation}
est exacte. Considérons $(\mF_p \cdot b_n\oplus (\frac{p}{b_n^{p-1}})\cdot \uoR/p\uoR)_{n\geq n_0}$ comme un sous-système projectif
du système projectif $(\uoR/p\uoR)_{n\geq n_0}$ dont les morphismes de transition sont les itérés de $\phi$. 
Comme $v(b_n)<\frac 1 p$, on a $pv(\frac{p}{b_n^{p-1}})>1$. 
Par suite, $(\mF_p \cdot b_n\oplus (\frac{p}{b_n^{p-1}})\cdot \uoR/p\uoR)_{n\geq n_0}$ 
vérifie la condition de Mittag-Leffler. La proposition s'ensuit par passage à la limite projective.

\subsection{}\label{TCFR5}
On désigne par $\tE'^{\mN^\circ}$ le topos des systèmes projectifs de $\tE'$ \eqref{mtfla9}, indexés par l'ensemble ordonné $\mN$ des entiers naturels 
\eqref{notconv13}, et par 
\begin{equation}\label{TCFR5a}
\uplambda'\colon \tE'^{\mN^\circ}\rightarrow \tE'
\end{equation}
le morphisme canonique \eqref{notconv13a}. Pour tout objet $F=(F_n)_{n\geq 0}$ de $\tE'^{\mN^\circ}$ 
et tout entier $i\geq 0$, on note $\varepsilon_{\leq i}(F)$ l'objet $(F_n^i)_{n\geq 0}$ de $\tE'^{\mN^\circ}$, défini par 
\begin{equation}\label{TCFR5b}
F_n^i=\left\{
\begin{array}{clcr}
F_n&{\rm si}\ 0\leq n\leq i,\\
F_i&{\rm si}\ n\geq i,
\end{array}
\right.
\end{equation}
le morphisme de transition $F_{n+1}^i\rightarrow F_n^i$ étant égal au morphisme de transition $F_{n+1}\rightarrow F_n$ de $F$ si $0\leq n<i$ et à l'identité 
si $n\geq i$. La correspondance $F\mapsto \varepsilon_{\leq i}(F)$ est clairement fonctorielle, et on a un morphisme fonctoriel canonique
\begin{equation}\label{TCFR5c}
F\rightarrow  \varepsilon_{\leq i}(F).
\end{equation}
Si $F$ est un anneau (resp. module), il en est de même de $\varepsilon_{\leq i}(F)$ et le morphisme \eqref{TCFR5c} est un homomorphisme. 

Si $A$ est un anneau, on note encore $A$ l'anneau constant de valeur $A$ de $\tE'$ ou de $\tE'^{\mN^\circ}$.
On a donc un isomorphisme canonique $\uplambda'^*(A)\stackrel{\sim}{\rightarrow} A$ de $\tE'^{\mN^\circ}$.

On note $(\co_\oK)^\wp$ l'anneau $(\co_\oK/p\co_\oK)_{\mN}$ de $\tE'^{\mN^\circ}$ dont les morphismes de transition sont les itérés
de l'endomorphisme de Frobenius de $\co_\oK/p\co_\oK$. On le considère comme  une $\co_{\oK^\flat}$-algèbre 
via les homomorphismes \eqref{elcpam2b}.

On désigne par $\cP'$ la catégorie abélienne 
des $\mF_p$-modules de $\tE'^{\mN^\circ}$ et par $\cP'_\AR$ son quotient par la sous-catégorie 
épaisse des $\mF_p$-modules Artin-Rees-nuls (AR-nuls en abrégé) (\cite{sga5} V 2.2.1). On rappelle que $\cP'_\AR$
est canoniquement équivalente à la catégorie des systèmes projectifs de $\cP'$ à translation près (\cite{sga5} V 2.4.2 et 2.4.4).

\subsection{}\label{TCFR6}
On note $\phi'$ l'endomorphisme de Frobenius de $\ocB'_1$ \eqref{mtfla8e}
et on désigne par $\ocB'^\wp$ la $(\co_\oK)^\wp$-algèbre $(\ocB'_1)_{\mN}$ de $\tE'^{\mN^\circ}$  
dont les morphismes de transition sont les itérés de $\phi'$ \eqref{TCFR5}. 
On considère $\ocB'^\wp$ comme une $\co_{\oK^\flat}$-algèbre.
On note encore $\phi'$ l'endomorphisme de Frobenius de $\ocB'^\wp$.
Cette notation n'induit aucun risque d'ambiguïté puisque l'endomorphisme de Frobenius de $\ocB'^\wp$ est induit par celui de $\ocB'_1$.

\subsection{}\label{TCFR7} 
Pour tout faisceau abélien $F=(F_n)_{n\geq 0}$ de $\tE'^{\mN^\circ}$, on désigne par $\Theta_*(F)_{\rho(\oy\rightsquigarrow \ox)}$ 
la limite projective des fibres des $\Theta_*(F_n)$ en $\rho(\oy\rightsquigarrow \ox)$ \eqref{mtfla7e},
\begin{equation}\label{TCFR7a} 
\Theta_*(F)_{\rho(\oy\rightsquigarrow \ox)}=\underset{\underset{\mN}{\longleftarrow}}{\lim}\ \Theta_*(F_n)_{\rho(\oy\rightsquigarrow \ox)},
\end{equation}
où $\Theta$ est le morphisme de topos \eqref{mtfla11b}.
On définit ainsi un foncteur de la catégorie des faisceaux abéliens de $\tE'^{\mN^\circ}$ dans celle des faisceaux abéliens.
Celui-ci étant exact à gauche, 
on désigne abusivement par $\rR^q\Theta_*(F)_{\rho(\oy\rightsquigarrow \ox)}$ ($q\geq 0$) ses foncteurs dérivés à droite. 

D'après (\cite{jannsen} 1.6), on a une suite exacte canonique 
\begin{equation}\label{TCFR7b}
0\rightarrow \rR^1 \underset{\underset{\mN}{\longleftarrow}}{\lim}\ \rR^{q-1}\Theta_*(F_i)_{\rho(\oy\rightsquigarrow \ox)}\rightarrow
\rR^q\Theta_*(F)_{\rho(\oy\rightsquigarrow \ox)}\rightarrow 
\underset{\underset{\mN}{\longleftarrow}}{\lim}\ \rR^q\Theta_*(F_i)_{\rho(\oy\rightsquigarrow \ox)}\rightarrow 0. 
\end{equation}

\subsection{}\label{TCFR8} 
Pour tout $\mF_p$-module de type fini $\mL'$ de $\oX'^\rhd_\fet$ et tout entier $q\geq 0$,  
on associe le $(\uoR/p\uoR)$-module 
\begin{equation}\label{TCFR8a} 
M^q(\mL')=\rR^q\Theta_*(\beta'^*(\mL')\otimes_{\mF_p}\ocB'_1)_{\rho(\oy\rightsquigarrow \ox)},
\end{equation}
où $\beta'$ est le morphisme \eqref{mtfla9c}. 
On le munit de l'endomorphisme $(\uoR/p\uoR)$-semi-linéaire $\varphi$ 
induit par l'endomorphisme de Frobenius $\phi'$ de $\ocB'_1$. 
On leur associe aussi le $\uoR^\flat$ module 
\begin{equation}\label{TCFR8b} 
\cM^q(\mL')=\rR^q\Theta_*(\uplambda'^*(\beta'^*(\mL'))\otimes_{\mF_p}\ocB'^\wp)_{\rho(\oy\rightsquigarrow \ox)}, 
\end{equation}
que l'on munit de l'endomorphisme $\uoR^\flat$-semi-linéaire $\Phi$ 
induit par l'endomorphisme de Frobenius $\phi'$ de $\ocB'^\wp$ \eqref{TCFR6}. 

D'après \ref{TPCF8}, pour tout entier $i\geq 0$, la suite de $\ocB^\wp$-modules
\begin{equation}\label{TCFR8c}
\xymatrix{
0\ar[r]&{\ocB'^\wp}\ar[r]^-(0.5){\cdot \varpi^{p^i}}&{\ocB'^\wp}\ar[r]&{\varepsilon_{\leq i}(\ocB'^\wp)}\ar[r]&0}
\end{equation}
où la deuxième flèche est le morphisme canonique \eqref{TCFR5c}, est exacte dans $\cP_\AR$. 
Compte tenu de \eqref{TCFR7b} et (\cite{jannsen} 1.15), 
on en déduit une suite exacte longue de $\uoR^\flat$-modules
\begin{equation}\label{TCFR8d}
\xymatrix{
{\cM^q(\mL')}\ar[r]^-(0.5){\cdot \varpi^{p^i}}&{\cM^q(\mL')}\ar[r]^{\fp^q_i(\mL')}&{M^q(\mL')}\ar[r]&{\cM^{q+1}(\mL')}},
\end{equation}
où $M^q(\mL')$ est considéré comme un $\uoR^\flat$-module via l'homomorphisme $\pi_i$ \eqref{TCFR3b}.

\begin{prop}\label{TCFR9} 
Supposons le morphisme $g\colon X'\rightarrow X$ projectif \eqref{mtfla2a}.
Soient $\mL'$ un $\mF_p$-module de type fini de $\oX'^\rhd_\fet$, $q$ un entier $\geq 0$ et reprenons les notations de  \ref{TCFR8}.
Alors,
\begin{itemize}
\item[{\rm (i)}] Pour tout entier $i\geq 0$, la suite de $\uoR^\flat$-modules
\begin{equation}\label{TCFR9a} 
\xymatrix{
0\ar[r]&{\cM^q(\mL')}\ar[r]^{\varpi^{p^i}}&{\cM^q(\mL')}\ar[r]^{\fp^q_i(\mL')}&{M^q(\mL')}\ar[r]&0}
\end{equation}
est $\alpha$-exacte. 
\item[{\rm (ii)}] Le morphisme de $\uoR^\flat$-modules 
\begin{equation}\label{TCFR9b} 
\cM^q(\mL')\rightarrow \underset{\underset{\mN}{\longleftarrow}}{\lim} \ M^q(\mL'),
\end{equation}
induit par les morphismes $\fp^q_i(\mL')$ $(i\geq 0)$,
où les morphismes de transition de la limite projective sont les itérés de $\varphi$, est un $\alpha$-isomorphisme. 
\item[{\rm (iii)}] Il existe une $(\uoR/p\uoR)$-algèbre $\alpha$-finie-étale et $\alpha$-fidèlement plate $A$, un entier $r\geq 0$
et un morphisme $A$-linéaire 
\begin{equation}\label{TCFR9c} 
\fm_\oK\otimes_{\co_\oK}M^q(\mL')\otimes_{\uoR}A\rightarrow A^r,
\end{equation}
compatible à $\varphi$ et aux endomorphismes de Frobenius de $\co_\oK/p\co_\oK$ et $A$, qui est un $\alpha$-isomorphisme. 
\item[{\rm (iv)}] Il existe une $\uoR^\flat$-algèbre $\alpha$-finie-étale et $\alpha$-fidèlement plate $\cA$, un entier $r\geq 0$ 
et un morphisme $\cA$-linéaire 
\begin{equation}\label{TCFR9d} 
\fm_{\oK^\flat}\otimes_{\co_\oK^\flat}\cM^q(\mL')\otimes_{\uoR^\flat}\cA\rightarrow \cA^r,
\end{equation}
compatible à $\Phi$ et aux endomorphismes de Frobenius de $\co_{\oK^\flat}$ et $\cA$, qui est un $\alpha$-isomor\-phisme. 
\end{itemize}
\end{prop}

Notons d'abord que compte tenu de \eqref{TCFR8d}, la proposition (i) est équivalente à la suivante:
\begin{itemize}
\item[{\rm (i')}] Pour tout entier $i\geq 0$, le morphisme 
\begin{equation}\label{TCFR9e} 
\fp^q_i(\mL')\colon \cM^q(\mL')\rightarrow M^q(\mL')
\end{equation}
est $\alpha$-surjectif. 
\end{itemize}

Nous omettrons $\mL'$ des notations $M^q(\mL')$, $\cM^q(\mL')$ et $\fp^q_i(\mL')$ pour les alléger. 
Nous équiperons les propositions (i'), (ii), (iii) et (iv) d'un indice $q$ pour rappeler la dépendance en cet indice.  

Nous procédons par récurrence descendante sur $q$. D'après \ref{AFR28},  
$M^q$ est $\alpha$-nul pour $q$ assez grand, et il en est donc de même de $\cM^q$ 
compte tenu de \eqref{TCFR7b} et (\cite{gr} 2.4.2(ii)), d'où la proposition pour ces valeurs de $q$. 
Supposons la proposition (iv)$_{q+1}$ établie et montrons les propositions (i')$_q$, (ii)$_q$, (iii)$_q$ et (iv)$_q$. 
Comme $\varpi$ n'est pas un diviseur de zéro dans $\uoR^\flat$ \eqref{TCFR4b},  
on déduit de (iv)$_{q+1}$, par descente $\alpha$-fidèlement plate, que la multiplication par $\varpi$ dans $\cM^{q+1}$ est $\alpha$-injective.  
La suite exacte longue \eqref{TCFR8d} implique alors que pour tout entier $i\geq 0$, le morphisme 
\begin{equation}\label{TCFR9f} 
\cM^q/\varpi^{p^i}\cM^q\rightarrow M^q
\end{equation}
induit par $\fp^q_i$ est un $\alpha$-isomorphisme, d'où la proposition (i')$_q$.   

Par ailleurs, l'endomorphisme de Frobenius $\phi'$ de $\ocB'^\wp$ induit un isomorphisme de $\cP'_\AR$ d'après \ref{TPCF8}(ii).
L'endomorphisme $\Phi$ de $\cM^q$ est donc un isomorphisme. 
D'après \ref{TCFR4}(iii), le morphisme $\pi_0$ induit 
un isomorphisme $\uoR^\flat/\phi^{-1}(\varpi)\uoR^\flat\stackrel{\sim}{\rightarrow}\uoR/p_1\uoR$. 
Compte tenu de \eqref{TCFR9f}, on en déduit que le morphisme $\uoR$-linéaire
\begin{equation}\label{TCFR9g} 
\ophi_{p_1}^*(M^q/p_1 M^q)\rightarrow M^q
\end{equation}
induit par l'endomorphisme $\varphi$ de $M^q$  est un $\alpha$-isomorphisme \eqref{TCFR4a}. 

Observons que les hypothèses de \ref{AFR10} sont satisfaites et reprenons les notations de \ref{AFR13}.
En vertu de \ref{AFR15}, on a un morphisme $\uoR$-linéaire canonique 
\begin{equation}\label{TCFR9h} 
\sigma^{(\infty)}_*(\rR^q\Theta_*(\beta'^*(\mL')\otimes_{\mF_p}\ocB'_1))_\ox\otimes_{\sigma^{(\infty)}_*(\ocB)_\ox}\uoR\rightarrow M^q,
\end{equation}
qui est un $\alpha$-isomorphisme. Tenant compte de \eqref{AFR13k}, on pose
\begin{equation}\label{TCFR9i} 
N^q=\sigma^{(\infty)}_*(\rR^q\Theta_*(\beta'^*(\mL')\otimes_{\mF_p}\ocB'_1))_\ox\otimes_{\sigma^{(\infty)}_*(\ocB)_\ox}\uR_\infty,
\end{equation}
que l'on munit de l'endomorphisme $(\uR_\infty/p\uR_\infty)$-semi-linéaire $\psi$ induit par l'endomorphisme de Frobenius $\phi'$ de $\ocB'_1$.
Le morphisme \eqref{TCFR9h} induit un morphisme $\uoR$-linéaire
\begin{equation}\label{TCFR9j} 
N^q\otimes_{\uR_\infty}\uoR\rightarrow M^q,
\end{equation}
qui est un $\alpha$-isomorphisme. Celui-ci est compatible à $\psi$, $\varphi$
et à l'endomorphisme de Frobenius $\phi$ de $\uoR/p\uoR$, d'après la preuve de de \ref{AFR15}.

Comme la $\uR_\infty$-algèbre $\uoR$ est $\alpha$-fidèlement plate \eqref{TCFR2}
et que le morphisme \eqref{TCFR9g} est un $\alpha$-isomorphisme, le morphisme $\uR_\infty$-linéaire
\begin{equation}\label{TCFR9k} 
\ophi_{p_1}^*(N^q/p_1 N^q)\rightarrow N^q
\end{equation}
induit par l'endomorphisme $\psi$ de $N^q$, est un $\alpha$-isomorphisme (cf. \ref{elcpam7}). 
Par ailleurs, le $\uR_\infty$-module $N^q$ est de présentation $\alpha$-finie, en vertu de \ref{AFR27}. 
La proposition (iii)$_q$ résulte alors de \ref{elcpam15} et \eqref{TCFR9j}.

Notons $\cM^q_\proj$ le $\uoR^\flat$-module limite projective du système projectif $(M^q)_\mN$ dont les morphismes de transition sont
les itérés de $\varphi$, et $\Phi$ l'endomorphisme $\phi$-semi-linéaire de $\cM^q_\proj$ induit par $\varphi$. 
De même, notons $\cN^q_\proj$ le $\uR_\infty^\flat$-module limite projective du système projectif $(N^q)_\mN$ dont les morphismes de transition sont
les itérés de $\psi$, et $\Psi$ l'endomorphisme $\phi$-semi-linéaire de $\cN^q_\proj$ induit par $\psi$. 
En vertu de \ref{elcpam23}, on a les propriétés suivantes:
\begin{itemize} 
\item[(a)] pour tout entier $i\geq 0$, la suite 
\begin{equation}\label{TCFR9l}
0\longrightarrow \cN^q_\proj\stackrel{\varpi^{p^i}}{\longrightarrow} \cN^q_\proj\stackrel{\tau_i}{\longrightarrow} N^q\longrightarrow 0,
\end{equation}
où $\tau_i$ est le morphisme qui associe à $(x_i)_{i\in \mN}\in \cN^q_\proj$ l'élément $x_i\in N^q$, est $\alpha$-exacte;
\item[(b)] il existe une $\uR^\flat_\infty$-algèbre $\alpha$-finie-étale et $\alpha$-fidèlement plate $\cA$,  
un entier $r\geq 0$ et un morphisme $\cA$-linéaire 
\begin{equation}\label{TCFR9m}
\fm_{\oK^\flat}\otimes_{\co_{\oK^\flat}}\cN^q_\proj\otimes_{\uR^\flat_\infty}\cA\rightarrow \cA^r
\end{equation}
compatible à $\Psi$ et à l'endomorphisme de Frobenius de $\cA$, qui est un $\alpha$-isomorphisme.
\end{itemize}

\vspace{2mm}

On observera que $\cM^q_\proj$ et $\cN^q_\proj$ sont complets et séparés pour les topologies $\varpi$-adiques d'après \ref{mdc3}. 
Notons $\cN^q_\proj\hotimes_{\uR_\infty^\flat}\uoR^\flat$ le séparé complété $\varpi$-adique de
$\cN^q_\proj\otimes_{\uR_\infty^\flat}\uoR^\flat$. 
Compte tenu de \eqref{elcpam5a}, \eqref{TCFR4b} et \eqref{TCFR9l}, le morphisme \eqref{TCFR9j} induit un morphisme $\uoR^\flat$-linéaire
\begin{equation}\label{TCFR9n}
\cN^q_\proj\hotimes_{\uR_\infty^\flat}\uoR^\flat\rightarrow \cM^q_\proj, 
\end{equation}
qui est un $\alpha$-isomorphisme d'après (\cite{gr} 2.4.2(ii)). Montrons que le morphisme induit
\begin{equation}\label{TCFR9p}
u\colon \cN^q_\proj\otimes_{\uR_\infty^\flat}\uoR^\flat\rightarrow \cM^q_\proj
\end{equation}
est un $\alpha$-isomorphisme.
On note  $\uoR^\flat\hotimes_{\uR_\infty^\flat}\cA$ et $\cM^q_\proj\hotimes_{\uR_\infty^\flat}\cA$ les séparés complétés $\varpi$-adiques de 
$\uoR^\flat\otimes_{\uR_\infty^\flat}\cA$ et $\cM^q_\proj\otimes_{\uR_\infty^\flat}\cA$, et  
\begin{eqnarray}
v\colon \uoR^\flat\otimes_{\uR_\infty^\flat}\cA&\rightarrow& \uoR^\flat\hotimes_{\uR_\infty^\flat}\cA\label{TCFR9q}\\
w\colon \cM^q_\proj\otimes_{\uR_\infty^\flat}\cA&\rightarrow& \cM^q_\proj\hotimes_{\uR_\infty^\flat}\cA\label{TCFR9r}
\end{eqnarray}
les morphismes canoniques.
Comme le $\uR_\infty^\flat$-module $\cA$ est $\alpha$-projectif de type $\alpha$-fini \eqref{aet2}, 
$v$ et $w$ sont des $\alpha$-isomorphismes d'après \ref{elcpam20}. 
Considérons le diagramme commutatif 
\begin{equation}\label{TCFR9s}
\xymatrix{
{\cN^q_\proj\otimes_{\uR^\flat_\infty}\uoR^\flat\otimes_{\uR^\flat_\infty}\cA}\ar[r]\ar[rd]_{u\otimes\id}\ar@/^2pc/[rr]^t&
{\cN^q_\proj\hotimes_{\uR^\flat_\infty}\uoR^\flat\otimes_{\uR^\flat_\infty}\cA}\ar[r]\ar[d]&
{\cN^q_\proj\hotimes_{\uR^\flat_\infty}\uoR^\flat\hotimes_{\uR^\flat_\infty}\cA}\ar[d]\\
&{\cM^q_\proj\otimes_{\uR^\flat_\infty}\cA}\ar[r]^-(0.5)w&{\cM^q_\proj\hotimes_{\uR^\flat_\infty}\cA}}
\end{equation}
où $\cN^q_\proj\hotimes_{\uR^\flat_\infty}\uoR^\flat\hotimes_{\uR^\flat_\infty}\cA$ est le séparé complété $\varpi$-adique de
$\cN^q_\proj\otimes_{\uR^\flat_\infty}\uoR^\flat\otimes_{\uR^\flat_\infty}\cA$. 
Les flèches verticales sont des $\alpha$-isomorphismes d'après \eqref{TCFR9n}, 
et $t$ est un $\alpha$-isomorphisme d'après \eqref{TCFR9q}. Il s'ensuit que $u\otimes\id$ est un $\alpha$-isomorphisme.
Comme la $\uR_\infty^\flat$-algèbre $\cA$ est $\alpha$-fidèlement plate, on en déduit que $u$ est un $\alpha$-isomorphisme.

Pour tout entier $i\geq 0$, la suite 
\begin{equation}\label{TCFR9o}
0\longrightarrow \cM^q_\proj\stackrel{\varpi^{p^i}}{\longrightarrow} \cM^q_\proj\stackrel{\kappa_i}{\longrightarrow} M^q\longrightarrow 0,
\end{equation}
où $\kappa_i$ est le morphisme qui associe à $(x_i)_{i\in \mN}\in \cM^q_\proj$ l'élément $x_i\in M^q$, est $\alpha$-exacte. 
Cela résulte de \eqref{TCFR9l} et \eqref{TCFR9p}, sauf pour l'exactitude à gauche qui découle de \eqref{TCFR9m} et \eqref{TCFR9p}. 
En effet, la multiplication par $\varpi$ dans $\uoR^\flat$ est injective \eqref{TCFR4b} et la $\uR_\infty^\flat$-algèbre $\cA$ est $\alpha$-fidèlement plate.

Considérons de nouveau le système projectif $(M^q)_\mN$ dont les morphismes de transition sont les itérés de $\varphi$ et posons
\begin{equation}\label{TCFR9t}
\cT^q=\rR^1 \underset{\underset{\mN}{\longleftarrow}}{\lim}\ M^q. 
\end{equation}
D'après \eqref{TCFR7b} et (\cite{agt} 7.12), on a une suite exacte canonique
\begin{equation}
0\rightarrow \cT^{q-1}\rightarrow \cM^q \rightarrow \cM^q_\proj\rightarrow 0. 
\end{equation}
Il résulte alors de \eqref{TCFR9f} et \eqref{TCFR9o} que $\cT^{q-1}/\varpi \cT^{q-1}$ est $\alpha$-nul. 
Comme le $\co_{\oK^\flat}$-module $\cT^{q-1}$ est dérivé-$\varpi$-complet d'après \ref{mdc5}(iii), on en déduit que $\cT^{q-1}$ 
est $\alpha$-nul en vertu de \ref{mdc8}. Les propositions (ii)$_q$ et (iv)$_q$ s'ensuivent \eqref{TCFR9m}.

\begin{cor}\label{TCFR10} 
Sous les hypothèses de \ref{TCFR9}, il existe un entier $r\geq 0$ et un morphisme $\uoR^\flat$-linéaire 
\begin{equation}\label{TCFR10a} 
\fm_{\oK^\flat}\otimes_{\co_\oK^\flat}\cM^q(\mL')\rightarrow (\uoR^\flat)^r,
\end{equation}
compatible à $\Phi$ et aux endomorphismes de Frobenius de $\co_{\oK^\flat}$ et $\uoR^\flat$, qui est un $\alpha$-isomorphisme
\end{cor}

En effet, en vertu de \ref{TCFR9}(iv), il existe une $\uoR^\flat$-algèbre $\alpha$-finie-étale et $\alpha$-fidèlement plate $\cA$, un entier $r\geq 0$ 
et un morphisme $\cA$-linéaire 
\begin{equation}\label{TCFR10b} 
\fm_{\oK^\flat}\otimes_{\co_\oK^\flat}\cM^q(\mL')\otimes_{\uoR^\flat}\cA\rightarrow \cA^r,
\end{equation}
compatible à $\Phi$ et aux endomorphismes de Frobenius de $\co_{\oK^\flat}$ et $\cA$, qui est un $\alpha$-isomorphisme.  
D'après \ref{alpha6}, il suffit donc de montrer qu'il existe un morphisme de $\alpha$-$\uoR^\flat$-algèbres $\cA^\alpha\rightarrow (\uoR^\flat)^\alpha$
(cf. \ref{alpha5}), ou ce qui revient au même, un morphisme de $\uoR^\flat$-algèbres $\cA\rightarrow \sigma_*((\uoR^\flat)^\alpha)$,
où $\sigma_*((\uoR^\flat)^\alpha)=\Hom_{\co_{\oK^\flat}}(\fm_{\oK^\flat},\uoR^\flat)$. 

Pour toute $\co_\oK$-algèbre (resp. $\co_{\oK^\flat}$-algèbre) $A$, notons $A\eEt_{\apf}$ la catégorie des 
$\alpha$-$A$-algèbres $\alpha$-étales, de présentation $\alpha$-finie \eqref{aet2}. 

Comme l'anneau $\uoR$ est local et strictement hensélien d'après (\cite{agt} III.10.10.2(i)) et \eqref{TCFR1a} et que $p$ appartient à l'idéal
maximal de $\uoR$, $(\uoR,p\uoR)$ est un couple hensélien au sens de (\cite{raynaud1} XI défi.~3). Par suite, en vertu de (\cite{gr} 5.5.7(iii)), 
le foncteur 
\begin{equation}\label{TCFR10c}
\uoR\eEt_\apf\rightarrow (\uoR/p\uoR)\eEt_\apf, \ \ \ C\mapsto C/pC,
\end{equation}
est une équivalence de catégories. 
De même, comme l'anneau $\uoR^\flat$ est complet et séparé pour la topologie $\varpi$-adique \eqref{TCFR4}, le foncteur 
\begin{equation}\label{TCFR10d}
\uoR^\flat\eEt_\apf\rightarrow (\uoR^\flat/\varpi\uoR^\flat)\eEt_\apf, \ \ \ C\mapsto C/\varpi C,
\end{equation}
est une équivalence de catégories. Compte tenu de l'isomorphisme $\uoR^\flat/\varpi\uoR^\flat\stackrel{\sim}{\rightarrow} \uoR/p\uoR$ \eqref{TCFR4b},
on en déduit une équivalence de catégories 
\begin{equation}\label{TCFR10e}
\uoR^\flat\eEt_\apf\stackrel{\sim}{\rightarrow} \uoR\eEt_\apf.
\end{equation}
Notons $\tcA$ l'image de $\cA^\alpha$ par cette équivalence. Il suffit donc de montrer qu'il existe un morphisme de $\alpha$-$\uoR$-algèbres
$\tcA\rightarrow \uoR^\alpha$.

Comme la $\uoR^\flat$-algèbre $\cA$ est $\alpha$-projective de type $\alpha$-fini et de rang fini \eqref{aet2}, 
la $\alpha$-$\uoR^\alpha$-algèbre $\tcA$ est $\alpha$-projective de type $\alpha$-fini et  de rang fini d'après (\cite{gr} 2.4.18 \& 5.1.7). 
Par ailleurs, la $\uoR^\flat$-algèbre $\cA$ étant $\alpha$-fidèlement plate, elle n'est pas $\alpha$-nulle. 
Par suite, la $\alpha$-$\uoR$-algèbre $\tcA$ n'est pas nulle \eqref{TCFR10e}. 

La $\uoR[\frac 1 p]$-algèbre $C=\sigma_*(\tcA)[\frac 1 p]$ \eqref{alpha5b} est finie étale d'après \ref{alpha6}(ii).
Notons $\cC$ la clôture intégrale de $\uoR$ dans $C$. 
En vertu de (\cite{gr} 8.2.31(i)), les $\alpha$-$\uoR$-algèbres $\cC^\alpha$ et $\tcA$ sont isomorphes. 
En particulier, $C$ n'est pas nulle. 
Reprenant la définition \eqref{TFA12f} de $\uoR=\oR_\uX^\oy$, il existe un $X$-schéma étale $\ox$-pointé $U$ et une 
$\oR^\oy_U[\frac 1 p]$-algèbre finie étale et non nulle $D$ tels que $C=D\otimes_{\oR^\oy_U}\uoR$. Compte tenu de \eqref{TFA9c},
il existe un homomorphisme de $\oR^\oy_U[\frac 1 p]$-algèbres $D\rightarrow \oR^\oy_U[\frac 1 p]$ et par suite un 
homomorphisme de $\uoR[\frac 1 p]$-algèbres $C\rightarrow \uoR[\frac 1 p]$. Comme $\uoR$ est normal, on en déduit 
un homomorphisme de $\uoR$-algèbres $\cC\rightarrow \uoR$, d'où l'assertion recherchée.

\begin{cor}\label{TCFR11} 
Sous les hypothèses de \ref{TCFR9}, il existe un entier $r\geq 0$ et un morphisme $(\uoR/p\uoR)$-linéaire 
\begin{equation}\label{TCFR11a} 
\fm_\oK\otimes_{\co_\oK}M^q(\mL')\rightarrow (\uoR/p\uoR)^r,
\end{equation}
compatible à $\varphi$ et aux endomorphismes de Frobenius de $\co_\oK/p\co_\oK$ et $\uoR/p\uoR$, qui est un $\alpha$-isomor\-phisme. 
\end{cor}

Cela résulte de \ref{TCFR10}, \ref{TCFR9}(i) et \ref{TCFR4}(ii).

\begin{prop}\label{TCFR12}
Supposons le morphisme $g\colon X'\rightarrow X$ projectif \eqref{mtfla2a}.
Soient $\mL'$ un $\mF_p$-module de type fini de $\oX'^\rhd_\fet$, $b$ un élément non-nul de $\co_{\oK^\flat}$, 
$q$ un entier $\geq 0$. Reprenons les notations de  \ref{TCFR8}. Alors,
\begin{itemize}
\item[{\rm (i)}] La suite  
\begin{equation}\label{TCFR12a}
\xymatrix{ 
0\ar[r]&{\rR^q\Theta_*(\beta'^*(\mL'))_{\rho(\oy\rightsquigarrow \ox)}}\ar[r]^-(0.5){\cdot b}&{\cM^q(\mL')}\ar[rr]^-(0.5){\Phi-b^{p-1}\id}&&
{\cM^q(\mL')}\ar[r]& 0}
\end{equation}
est exacte.
\item[{\rm (ii)}] Le morphisme canonique
\begin{equation}\label{TCFR12aa}
\ker(\Phi-b^{p-1}\id|\cM^q(\mL'))\rightarrow \ker(\Phi-b^{p-1}\id|\cM^q(\mL')\otimes_{\co_{\oK^\flat}}\oK^\flat)
\end{equation}
est un isomorphisme. 
\end{itemize}
\end{prop}

En effet, d'après \ref{TPCF11}, la suite 
\begin{equation}\label{TCFR12b}
\xymatrix{
0\ar[r]&{\mF_p}\ar[r]^-(0.4){\cdot b}&{\ocB'^\wp}\ar[rr]^{\phi'-b^{p-1}\id}&&{\ocB'^\wp}\ar[r]&0}
\end{equation}
est exacte dans $\cP'_\AR$ \eqref{TCFR5}. Compte tenu de \eqref{TCFR7b} et (\cite{jannsen} 1.15),
celle-ci induit une suite exacte longue de cohomologie
\[
\xymatrix{ 
{\cM^{q-1}}\ar[r]&{\rR^q\Theta_*(\beta'^*(\mL'))_{\rho(\oy\rightsquigarrow \ox)}}\ar[r]^-(0.5){\cdot b}&{\cM^q}\ar[rr]^-(0.5){\Phi-b^{p-1}\id}&&
{\cM^q}},
\]
où on a posé $\cM^q=\cM^q(\mL')$ pour tout $q\geq 0$ et $\cM^{-1}=0$. 
En vertu de \ref{TCFR10}, il existe un entier $r\geq 0$ et un morphisme $\uoR^\flat$-linéaire
\begin{equation}\label{TCFR12c}
u\colon \fm_{\oK^\flat}\otimes_{\co_{\oK^\flat}}\cM^q\rightarrow (\uoR^\flat)^r 
\end{equation}
compatible à $\Phi$ et aux endomorphismes de Frobenius de $\co_{\oK^\flat}$ et $\uoR^\flat$, et qui est un $\alpha$-isomorphisme. 

En vertu de \ref{TCFR15}, pour tout $c\in \co_{\oK^\flat}$, l'endomorphisme $\phi-c^{p-1}\id$ de $\uoR^\flat$ est surjectif. 
Pour tous $y\in \cM^q$ et $t_1\in \fm_{\oK^\flat}$, il existe donc $z\in (\uoR^\flat)^r$ tel que 
\begin{equation}\label{TCFR12d}
u(t_1^p\otimes y)=\phi(z)-c^{p-1}z,
\end{equation}
où on a encore noté $\phi$ l'endomorphisme de $(\uoR^\flat)^r$ induit par l'endomorphisme de Frobenius $\phi$ de $\uoR^\flat$. 
Pour tout $t_2\in \fm_{\oK^\flat}$, il existe $x\in \cM^q$ et $t\in \fm_{\oK^\flat}$ tels que $u(t\otimes x)=t_2z$. On en déduit que 
\begin{equation}\label{TCFR12e}
u((t_1t_2)^p\otimes y)=u(t^p\otimes \Phi(x)-(ct_2)^{p-1} t\otimes x)
\end{equation}
Pour tout $t_3\in \fm_{\oK^\flat}$, on a donc 
\begin{equation}\label{TCFR12f}
(t_1t_2t_3)^p y=\Phi(tt_3x)-(ct_2t_3)^{p-1} (tt_3 x).
\end{equation}
Par ailleurs, pour tout élément non-nul $b'$ de $\co_{\oK^\flat}$, on a un diagramme commutatif
\begin{equation}\label{TCFR12g}
\xymatrix{
{\ker(\Phi-b^{p-1}\id|\cM^q)}\ar[r]\ar[d]_{\cdot b'}&{\cM^q}\ar[rr]^{\Phi-b^{p-1}\id}\ar[d]_{\cdot b'}&&{\cM^q}\ar[r]\ar[d]^{\cdot b'^p}&
{\coker(\Phi-b^{p-1}\id|\cM^q)}\ar[d]^{\cdot b'^p}\\
{\ker(\Phi-(bb')^{p-1}\id|\cM^q)}\ar[r]&{\cM^q}\ar[rr]^{\Phi-(bb')^{p-1}\id}&&{\cM^q}\ar[r]&{\coker(\Phi-(bb')^{p-1}\id|\cM^q)}}
\end{equation}
On déduit de ce qui précède (en prenant $c=bt_1 $ et $b'=t_1t_2t_3$) 
que pour tout $b'\in \fm_{\oK^\flat}$, le dernier morphisme vertical du diagramme \eqref{TCFR12g} est nul.

D'après \ref{TPCF11}, pour tout élément non-nul $b'$ de $\fm_{\oK^\flat}$, on a un diagramme commutatif à lignes exactes
\begin{equation}\label{TCFR12h}
\xymatrix{
0\ar[r]&{\mF_p}\ar@{=}[d]\ar[r]^-(0.5){\cdot b}&{\ocB'^\wp}\ar[d]^{\cdot b'}\ar[rr]^-(0.5){\phi'-b^{p-1}\id}&&{\ocB'^\wp}\ar[d]^{\cdot b'^p}\ar[r]&0\\
0\ar[r]&{\mF_p}\ar[r]^-(0.5){\cdot bb'}&{\ocB'^\wp}\ar[rr]^-(0.5){\phi'-(bb')^{p-1}\id}&&{\ocB'^\wp}\ar[r]&0}
\end{equation}
Il induit un diagramme commutatif à lignes exactes
\begin{equation}\label{TCFR12i}
\xymatrix{
{\coker(\Phi'-b^{p-1}\id|\cM^{q-1})}\ar@{^(->}[r]\ar[d]_{\cdot b'^p}&{\rR^q\Theta_*(\beta'^*(\mL'))_{\rho(\oy\rightsquigarrow \ox)}}
\ar@{=}[d]\ar@{->>}[r]&{\ker(\Phi-b^{p-1}\id|\cM^q)}\ar[d]^{\cdot b'}\\
{\coker(\Phi'-(bb')^{p-1}\id|\cM^{q-1})}\ar@{^(->}[r]&{\rR^q\Theta_*(\beta'^*(\mL'))_{\rho(\oy\rightsquigarrow \ox)}}\ar@{->>}[r]&{\ker(\Phi-(bb')^{p-1}\id|\cM^q)}}
\end{equation}
où on a noté $\Phi'$ l'endomorphisme de $\cM^{q-1}$ induit par l'endomorphisme $\phi'$ de $\ocB'^\wp$ 
(pour le distinguer de l'endomorphisme $\Phi$ de $\cM^q$). 
Comme le morphisme vertical de gauche est nul, on en déduit que l'endomorphisme $\Phi'-b^{p-1}\id$ de $\cM^{q-1}$ est surjectif, 
et il en est donc de même de l'endomorphisme $\Phi-b^{p-1}\id$ de $\cM^q$; 
d'où l'exactitude de la suite \eqref{TCFR12a}. On en déduit aussi que la flèche verticale de droite est un isomorphisme 
\begin{equation}\label{TCFR12j}
\cdot b'\colon \ker(\Phi-b^{p-1}\id|\cM^q)\stackrel{\sim}{\rightarrow} \ker(\Phi-(bb')^{p-1}\id|\cM^q).
\end{equation}
Il s'ensuit que l'application \eqref{TCFR12aa} est injective. 

Soit $x\in \cM^q\otimes_{\co_{\oK^\flat}}\oK^\flat$ tel que 
$\Phi(x)=b^{p-1}x$. Il existe $z\in \cM^q$ et $b'\in \fm_{\oK^\flat}$ tels que $b'\not=0$ et $z=b'x\in \cM^q\otimes_{\co_{\oK^\flat}}\oK^\flat$.
Par suite, $\Phi(z)-(bb')^{p-1}z$ est annulé par une puissance de $\varpi$. 
Quitte à remplacer $b'$ par un multiple, on peut supposer que $\Phi(z)=(bb')^{p-1}z$.
D'après l'isomorphisme \eqref{TCFR12j}, il existe $x'\in \ker(\Phi-b^{p-1}\id|\cM^q)$ tel que $z=b'x'$. 
On a alors $x=x'$ dans $\cM^q\otimes_{\co_{\oK^\flat}}\oK^\flat$; d'où la surjectivité de l'application \eqref{TCFR12aa}. 

\begin{cor}\label{TCFR13}
Supposons le morphisme $g\colon X'\rightarrow X$ projectif \eqref{mtfla2a}.
Soient $\mL'$ un $\mF_p$-module de type fini de $\oX'^\rhd_\fet$, 
$q$ un entier $\geq 0$. Alors, reprenant les notations de  \ref{TCFR8}, le $\mF_p$-espace vectoriel 
$\rR^q\Theta_*(\beta'^*(\mL'))_{\rho(\oy\rightsquigarrow \ox)}$ est de dimension finie,
et le morphisme canonique
\begin{equation}\label{TCFR13a}
\rR^q\Theta_*(\beta'^*(\mL'))_{\rho(\oy\rightsquigarrow \ox)}\otimes_{\mF_p}\uoR^\flat
\rightarrow \cM^q(\mL')
\end{equation}
est un $\alpha$-isomorphisme.
\end{cor}

En effet, en vertu de \ref{TCFR12} et avec les notations de sa preuve, on a des isomorphismes canoniques 
\begin{equation}\label{TCFR13b}
\rR^q\Theta_*(\beta'^*(\mL'))_{\rho(\oy\rightsquigarrow \ox)}\stackrel{\sim}{\rightarrow}\ker(\Phi-\id |\cM^q)\stackrel{\sim}{\rightarrow}
\ker(\Phi-\id |\cM^q\otimes_{\co_{\oK^\flat}}\oK^\flat).
\end{equation}
D'après \ref{TCFR10}, il existe un entier $r\geq 0$ et un morphisme $\uoR^\flat$-linéaire
\begin{equation}\label{TCFR13c}
u\colon \fm_{\oK^\flat}\otimes_{\co_{\oK^\flat}}\cM^q\rightarrow (\uoR^\flat)^r 
\end{equation}
compatible à $\Phi$ et aux endomorphismes de Frobenius de $\co_{\oK^\flat}$ et $\uoR^\flat$, et qui est un $\alpha$-isomorphisme. 
Il induit donc un isomorphisme $\uoR^\flat\otimes_{\co_{\oK^\flat}}\oK^\flat$-linéaire
\begin{equation}\label{TCFR13d}
\cM^q\otimes_{\co_{\oK^\flat}}\oK^\flat\stackrel{\sim}{\rightarrow} (\uoR^\flat\otimes_{\co_{\oK^\flat}}\oK^\flat)^r.
\end{equation}
On a la suite exacte d'Artin-Schreier
\begin{equation}
\xymatrix{
0\ar[r]&{\mF_p}\ar[r]&{\uoR^\flat\otimes_{\co_{\oK^\flat}}\oK^\flat}\ar[rr]^{\phi-\id}&&{\uoR^\flat\otimes_{\co_{\oK^\flat}}\oK^\flat}}.
\end{equation}
On en déduit que $\rR^q\Theta_*(\beta'^*(\mL'))_{\rho(\oy\rightsquigarrow \ox)}$ est un $\mF_p$-espace vectoriel de dimension finie, 
et que le morphisme $\uoR^\flat\otimes_{\co_{\oK^\flat}}\oK^\flat$-linéaire
\begin{equation}
\rR^q\Theta_*(\beta'^*(\mL'))_{\rho(\oy\rightsquigarrow \ox)}\otimes_{\mF_p}\uoR^\flat\otimes_{\co_{\oK^\flat}}\oK^\flat\rightarrow 
\cM^q\otimes_{\co_{\oK^\flat}}\oK^\flat
\end{equation}
est bijectif. 
D'après \ref{TCFR4}(ii), $\varpi$ n'est pas diviseur de zéro dans $\uoR^\flat$.
Par suite, le morphisme \eqref{TCFR13a} est injectif. Montrons qu'il est $\alpha$-surjectif.
Notons $e_1,\dots,e_r$ la base canonique de $(\uoR^\flat)^r$. Comme $u$ est $\alpha$-injectif, 
il suffit de montrer que pour tout $1\leq i\leq r$ et tout $\gamma\in \fm_{\oK^\flat}$, il existe $x\in \cM^q$ tel que 
$\Phi(x)=x$ et $u(\gamma\otimes x)=\gamma e_i$. Pour tout $t\in \fm_{\oK^\flat}$,  
il existe $y\in \cM^q$ et $\beta\in \fm_{\oK^\flat}$ tels que $u(\beta\otimes y)=te_i$. 
On a $u(\beta^p\otimes \Phi(y))=t^pe_i$ et donc $u(\beta^p\otimes \Phi(y)-t^{p-1}\beta\otimes y)=0$. Par suite, pour 
tout $t'\in \fm_{\oK^\flat}$, on a 
\begin{equation}
\Phi(t'\beta y)=(tt')^{p-1}t'\beta y.
\end{equation}
Compte tenu de l'isomorphisme \eqref{TCFR12j} (pour $b=1$ et $b'=tt'$), il existe $x\in \cM^q$ tel que $\Phi(x)=x$ et $t'\beta y=tt'x$. 
Par suite, pour tout $t''\in \fm_{\oK^\flat}$, on a $t't''\beta\otimes y=tt't''\otimes x$ dans $\fm_{\oK^\flat}\otimes_{\co_{\oK^\flat}}\cM^q$
et donc $u(tt't''\otimes x)=tt't''e_i$, d'où l'assertion recherchée.

\begin{cor}\label{TCFR14}
Supposons le morphisme $g\colon X'\rightarrow X$ projectif \eqref{mtfla2a}.
Soient $n$ un entier $\geq 1$, $\mL'$ un $(\mZ/p^n\mZ)$-module de type fini de $\oX'^\rhd_\fet$.
Alors, pour tout entier $q\geq 0$, le morphisme canonique
\begin{equation}\label{TCFR14a}
\rR^q\Theta_*(\beta'^*(\mL'))_{\rho(\oy\rightsquigarrow \ox)}\otimes_{\mZ_p}\uoR
\rightarrow \rR^q\Theta_*(\beta'^*(\mL')\otimes_{\mZ_p}\ocB')_{\rho(\oy\rightsquigarrow \ox)}
\end{equation}
est un $\alpha$-isomorphisme.
\end{cor}

En effet, par dévissage, on peut se réduire au cas où $n=1$. On notera que $\ocB'$ et $\uoR$ sont $\mZ_p$-plats (\cite{agt} III.9.2). 
La proposition résulte alors de \ref{TCFR13}, \ref{TCFR9}(i) et \eqref{TCFR4b}.

\section{\texorpdfstring{Principal théorème de comparaison $p$-adique de Faltings: cas relatif}
{Principal théorème de comparaison p-adique de Faltings: cas relatif}}

\begin{prop}\label{TCFR16}
Supposons le morphisme $g\colon X'\rightarrow X$ projectif \eqref{mtfla2a}.
Soient $n$ un entier $\geq 1$, $\mL'$ un $(\mZ/p^n\mZ)$-module de type fini de $\oX'^\rhd_\fet$.
Alors, pour tout entier $q\geq 0$, le morphisme canonique
\begin{equation}\label{TCFR16a}
\rR^q\Theta_*(\beta'^*(\mL'))\otimes_{\mZ_p}\ocB
\rightarrow \rR^q\Theta_*(\beta'^*(\mL')\otimes_{\mZ_p}\ocB')
\end{equation}
est un $\alpha$-isomorphisme.
\end{prop}

En effet, la fibre de ce morphisme en tout point de $\tE$ de la forme $\rho(\oy\rightsquigarrow \ox)$,
où $(\oy\rightsquigarrow \ox)$ est un point de $X_\et\gtimes_{X_\et}\oX^\circ_\et$ \eqref{topfl17} 
tel que $\ox$ soit au-dessus de $s$, est un $\alpha$-isomorphisme, en vertu de \ref{TCFR14} et \eqref{TCFR1a}.
Par ailleurs, $\ocB'|\sigma'^*(X'_\eta)$ et $\ocB|\sigma^*(X_\eta)$ étant des $\mQ_p$-algèbres, 
le morphisme \eqref{TCFR16a} est évidemment un isomorphisme au-dessus de l'ouvert $\sigma^*(X_\eta)$ de $\tE$. 
La proposition s'ensuit en vertu de \ref{tf21}(ii).

\begin{lem}\label{TCFR160}
Supposons le morphisme $g_\eta\colon X'_\eta\rightarrow X_\eta$ propre \eqref{mtfla2a}. 
Alors, le morphisme $g^\circ\colon X'^\circ\rightarrow X^\circ$ est lisse, et 
le schéma $X'^\rhd$ est le complémentaire dans $X'^\circ$ d'un diviseur à croisements normaux 
sur $X'^\circ$ relativement à $X^\circ$.
\end{lem}

En effet, tout point géométrique $\ox'$ de $X'_\eta$ étant une générisation d'un point géométrique de $X'$ au-dessus de $s$, 
il existe un entier $r\geq 0$ et un isomorphisme de monoïdes $\cM_{X',\ox'}/\co_{X',\ox'}^\times\stackrel{\sim}{\rightarrow}\mN^r$ d'après (\cite{agt} III.4.6). 
Supposons le support de $\ox'$ contenu dans $X'^\circ$. On a alors $\cM_{X'/X,\ox'}=\cM_{X',\ox'}/\co_{X',\ox'}^\times$ (\cite{ogus} II 1.1.10). 
En vertu de (\cite{ogus} IV 3.3.1), il existe un voisinage étale $U'$ de $\ox'$ dans $X'^\circ$, un voisinage étale $U$ de $g(\ox')$ dans $X^\circ$,
un morphisme $h\colon U'\rightarrow U$ au-dessus de $g$, un monoïde $P'$ et une carte $\gamma\colon P'\rightarrow \Gamma(U',\cM_{X'})$ pour $(U',\cM_{X'}|U')$ 
tels que les conditions suivantes soient satisfaites:
\begin{itemize}
\item[(i)] la carte $\gamma$ est ``neat'' dans la terminologie de (\cite{ogus} II 2.3.1), autrement dit le composé
\begin{equation}
\xymatrix{
P'\ar[r]^-(0.5){\gamma}&{\Gamma(U',\cM_{X'})}\ar[r]&{\cM_{X',\ox'}/\co^\times_{X',\ox'}}}
\end{equation} 
où la seconde flèche est l'homomorphisme canonique, est un isomorphisme; 
\item[(ii)] le morphisme de schémas usuels 
\begin{equation}
U'\rightarrow U\times_{\Spec(\mZ)}\bA_{P'}
\end{equation}
induit par $h$ et $\gamma$, est lisse (voir \cite{agt} II.5.12-II.5.14). 
\end{itemize}
La proposition s'ensuit puisque $U\times_{\Spec(\mZ)}\bA_{P'}$ est isomorphe à l'espace affine $\mA^r_U$ muni de la structure logarithmique associée 
au diviseur défini par le produit de toutes les coordonnées, 
qui est évidemment un diviseur à croisements normaux stricts sur $\mA^r_U$ relativement à $U$.

\subsection{}\label{TCFR17}
Notons $\upgamma\colon \oX'^\rhd\rightarrow \oX^\circ$ le morphisme induit par $g$ \eqref{mtfla2a}.
D'après (\cite{agt} (VI.10.12.7)), le diagramme
\begin{equation}\label{TCFR17a}
\xymatrix{
{\oX'^\rhd_\et}\ar[r]^-(0.5){\psi'}\ar[d]_{\upgamma}&{\tE'}\ar[d]^{\Theta}\\
{\oX^\circ_\et}\ar[r]^{\psi}&{\tE}}
\end{equation}
où $\psi$ et $\psi'$ sont les morphismes \eqref{mtfla7d} et \eqref{mtfla9d} respectivement, est commutatif à isomorphisme canonique près.

\begin{teo}[\cite{faltings2}, Theorem 6, page 266]\label{TCFR18}
Supposons le morphisme $g\colon X'\rightarrow X$ projectif \eqref{mtfla2a}.
Soient $n$ un entier $\geq 1$, $F'$ un $(\mZ/p^n\mZ)$-module localement constant constructible de $\oX'^\rhd_\et$.
Alors, pour tout entier $q\geq 0$, on a un morphisme canonique
\begin{equation}\label{TCFR18a}
\psi_*(\rR^q\upgamma_*(F'))\otimes_{\mZ_p}\ocB
\rightarrow \rR^q\Theta_*(\psi'_*(F')\otimes_{\mZ_p}\ocB')
\end{equation}
qui est un $\alpha$-isomorphisme.
\end{teo}

On peut évidemment se borner au cas où $F'$ est annulé par $p^n$, pour un entier $n\geq 1$. 
En vertu de \ref{acycloc2}, \ref{Kpun34}, \ref{TCFR160} et de la suite spectrale de Cartan-Leray, les morphismes canoniques 
\begin{eqnarray}
\rR^q\Theta_* (\psi'_*(F'))&\rightarrow& \rR^q(\Theta \circ \psi')_*(F'),\label{TCFR18b}\\
\rR^q(\psi \circ \upgamma)_*(F')&\rightarrow& \psi_*(\rR^q\upgamma_*(F')),\label{TCFR18c}
\end{eqnarray}
sont des isomorphismes. Compte tenu de l'isomorphisme de commutativité de \eqref{TCFR17a}
\begin{equation}\label{TCFR18d}
\Theta \circ \psi'\stackrel{\sim}{\rightarrow}\psi \circ \upgamma,
\end{equation}
on en déduit un isomorphisme 
\begin{equation}\label{TCFR18e}
\rR^q\Theta_* (\psi'_*(F'))\stackrel{\sim}{\rightarrow}\psi_*(\rR^q\upgamma_*(F')).
\end{equation}
Par ailleurs, notant $\rho_{\oX'^\rhd}\colon \oX'^\rhd_\et\rightarrow \oX'^\rhd_\fet$ le morphisme canonique \eqref{notconv10a}, 
il existe un $(\mZ/p^n\mZ)$-module de type fini canonique $\mL'$ de $\oX'^\rhd_\fet$ tel que $F'=\rho^*_{\oX'^\rhd}(\mL')$ d'après 
(\cite{agt} VI.9.18, VI.9.20 et III.2.11). En vertu de (\cite{agt} VI.10.9), on a un isomorphisme canonique 
\begin{equation}\label{TCFR18f}
\beta'^*(\mL')\stackrel{\sim}{\rightarrow}\psi'_*(\rho_{\oX'^\rhd}^*(\mL')).
\end{equation}
La proposition résulte alors de \ref{TCFR16}.

\begin{cor}\label{TCFR19}
Si le morphisme $g\colon X'\rightarrow X$ est projectif, pour tous entiers $n\geq 1$ et $q\geq 0$, on a un morphisme $\ocB$-linéaire canonique 
\begin{equation}\label{TCFR19a}
\psi_*(\rR^q\upgamma_*(\mZ/p^n\mZ))\otimes_{\mZ_p}\ocB\rightarrow \rR^q\Theta_*(\ocB'_n),
\end{equation}
qui est un $\alpha$-isomorphisme. 
\end{cor}

En effet, d'après (\cite{agt} VI.10.9(iii)), on a un isomorphisme canonique 
\begin{equation}\label{TCFR19b}
\beta'^*(\mZ/p^n\mZ)\stackrel{\sim}{\rightarrow} \psi'_*(\rho_{\oX'^\rhd}^*(\mZ/p^n\mZ)),
\end{equation}
où $\rho_{\oX'^\rhd}\colon \oX'^\rhd_\et\rightarrow \oX'^\rhd_\fet$ est le morphisme canonique \eqref{notconv10a}. 
On en déduit un isomorphisme 
\begin{equation}\label{TCFR19c}
\mZ/p^n\mZ\stackrel{\sim}{\rightarrow} \psi'_*(\mZ/p^n\mZ),
\end{equation}
qui n'est autre que le morphisme canonique par adjonction. La proposition résulte alors de \ref{TCFR18}. 

\begin{rema}[Note ajoutée le 2 décembre 2022]\label{TCFR180}
Le théorème \ref{TCFR18} vaut en fait sous l'hypothèse plus générale que {\em $g$ soit propre}. 
En effet, la condition de projectivité sur $g$ est utilisée pour prouver le résultat de $\alpha$-cohérence \ref{AFR23} que nous déduisons de \ref{afini15}. 
Alors que dans ce texte, nous nous basons sur les résultats de finitude de \cite{sga6} (plutôt que ceux de \cite{kiehl}),  
He \cite{tongmu3} vient de reprendre la preuve de Kiehl (\cite{kiehl} 2.9'a) et  est parvenu à étendre l'énoncé \ref{afini15} aux morphismes propres.
\end{rema}

\begin{rema}\label{TCFR181}
Dans (\cite{faltings2}, Theorem 6, page 265), Faltings requiert une condition de géométrie logarithmique, formulée en termes de géométrie torique, plus forte que 
celle requise dans \ref{TCFR18}, qui revient essentiellement à demander que $g$ admette localement une carte relativement adéquate  \eqref{mtfla5}
\begin{eqnarray}
((P',\gamma'),(P,\gamma),(\mN,\iota),\vartheta\colon \mN\rightarrow P, h\colon P\rightarrow P')
\end{eqnarray} 
telle que $P'^\gp/h^\gp(P^\gp)$ soit sans torsion. 
\end{rema}

\chapter{Les suites spectrales de Hodge-Tate}\label{suitesspecht}

\section{Hypothèses et notations}\label{hght}

\subsection{}\label{hght1}
Dans ce chapitre, $K$ désigne un corps de valuation discrète complet de 
caractéristique $0$, à corps résiduel {\em algébriquement clos} $k$ de caractéristique $p>0$,  
$\co_K$ l'anneau de valuation de $K$, $\oK$ une clôture algébrique de $K$, $\co_\oK$ la clôture intégrale de $\co_K$ dans $\oK$,
$\fm_\oK$ l'idéal maximal de $\co_\oK$ et $G_K$ le groupe de Galois de $\oK$ sur $K$.
On note $\co_C$ le séparé complété $p$-adique de $\co_\oK$, $\fm_C$ son idéal maximal et $C$ son corps des fractions. 
On désigne par $\mZ_p(1)$ le $\mZ[G_K]$-module 
\begin{equation}\label{hght1a}
\mZ_p(1)=\underset{\underset{n\geq 1}{\longleftarrow}}{\lim}\ \mu_{p^n}(\co_\oK),
\end{equation}  
où $\mu_{p^n}(\co_\oK)$ désigne le sous-groupe des racines $p^n$-ièmes de l'unité dans $\co_\oK$. 
Pour tout $\mZ_p[G_K]$-module $M$ et tout entier $n$, on pose $M(n)=M\otimes_{\mZ_p}\mZ_p(1)^{\otimes n}$.

On pose $S=\Spec(\co_K)$, $\oS=\Spec(\co_\oK)$ et $\coS=\Spec(\co_C)$. 
On note $s$ (resp.  $\eta$, resp. $\oeta$) le point fermé de $S$ (resp.  générique de $S$, resp. générique de $\oS$).
Pour tout entier $n\geq 1$, on pose $S_n=\Spec(\co_K/p^n\co_K)$. 
On munit $S$ de la structure logarithmique $\cM_S$ définie par son point fermé, 
et $\oS$ et $\coS$ des structures logarithmiques $\cM_\oS$ et $\cM_\coS$ images inverses de $\cM_S$.

Pour tout $S$-schéma $X$, on note $X_s$ (resp. $X_\eta$, resp. $X_\oeta$) 
la fibre fermée (resp. générique, resp. géométrique générique) de $X$ au dessus de $S$, et on pose 
\begin{equation}\label{hght1c}
\oX=X\times_S\oS,  \ \ \ \coX=X\times_S\coS \ \ \ {\rm et}\ \ \  X_n=X\times_SS_n.
\end{equation}

\subsection{}\label{hght3}
Pour toute extension finie $L$ de $K$, on note $\co_L$ la fermeture intégrale
de $\co_K$ dans $L$ et on pose $S_L=\Spec(\co_L)$ que l'on munit de la structure logarithmique $\cM_{S_L}$
définie par son point fermé. Les schémas logarithmiques $(S_L,\cM_{S_L})$ forment naturellement un système projectif. 
On désigne par $(\oS,\cL_\oS)$ la limite projective des schémas logarithmiques $(S_L,\cM_{S_L})$, 
indexée par l'ensemble filtrant des sous-$K$-extensions finies $L$ de $\oK$. 
Notant $\ou\colon \oeta\rightarrow \oS$ l'injection canonique, on a $\cL_\oS=\ou_*(\co_\oeta^\times)\cap \co_\oS$,
à ne pas confondre avec la structure logarithmique $\cM_\oS$ sur $\oS$ définie dans \ref{hght1}.

\subsection{}\label{tfkum4}
Pour tout $\mU$-topos $T$ et tout groupe abélien $M$, on note $M_T$ 
(ou simplement $M$ lorsqu'il n'y a aucun risque d'ambiguïté) le groupe constant de $T$ de valeur $M$. 
Pour tout entier $n\geq 1$, on désigne par $\mu_{n,T}$ 
(ou simplement $\mu_n$ lorsqu'il n'y a aucun risque d'ambiguïté) le groupe constant de $T$ de valeur 
le groupe $\mu_n(\co_\oK)$ des racines $n$-ièmes de l'unité dans $\co_\oK$.  
Pour tout $(\mZ/n\mZ)$-module $F$ de $T$ et tout entier $m$, on pose 
\begin{equation}\label{tfkum4a}
F(m)=\left\{
\begin{array}{clcr}
F\otimes_{\mZ}\mu_{n,T}^{\otimes m},&{\rm si}\ m\geq 0,\\
\Hom_{\mZ}(\mu_{n,T}^{\otimes -m},F),&{\rm si}\ m< 0.
\end{array}
\right.
\end{equation}
Cet usage n'étant pas classique pour le topos $X_\et$, 
on prendra garde de ne pas confondre le faisceau $\mu_{n,X_\et}$ ainsi défini avec 
le faisceau des racines $n$-ièmes de l'unité de $X_\et$, c'est-à-dire le noyau de 
la puissance $n$-ième sur $\mG_{m,X}$ (qui ne sera pas utilisé dans ce chapitre).

\begin{lem}\label{tfkum5}
Soient $A$ un anneau intègre, $M$ un monoïde intègre, 
$u$ un homomorphisme injectif de $M$ dans le monoïde multiplicatif $A$, $n$ un entier 
$\geq 1$, $\phi\colon A\rightarrow A$ l'application d'élévation à la puissance $n$-ième. 
Supposons que pour tout $t\in M$, l'équation $x^n=u(t)$ admette 
$n$ racines distinctes dans $A$. Alors, $\phi^{-1}(M)$ est un monoïde intègre et l'homomorphisme  
$\phi^{-1}(M)\rightarrow M$ déduit de $\phi$ identifie $M$ au quotient de $\phi^{-1}(M)$ par le sous-monoïde $\mu_n(A)$ 
des racines $n$-ièmes de l'unité dans $A$  \eqref{notconv7}. 
\end{lem}

En effet, pour tout $t\in M$, on a $u(t)\not=0$ car $u$ est injectif et $M$ est intègre. 
Comme $A$ est intègre, $\phi^{-1}(M)$ est donc intègre. Par ailleurs, l'application canonique 
$\phi^{-1}(M)\rightarrow M$ est surjective et fait de $\phi^{-1}(M)$ un torseur sur $M$ sous le groupe $\mu_n(A)$.
Par suite, $M$ est le quotient de $\phi^{-1}(M)$ par le sous-monoïde $\mu_n(A)$  (\cite{ogus} I 1.1.6).

\subsection{}\label{hght2}
Dans ce chapitre, $f\colon (X,\cM_X)\rightarrow (S,\cM_S)$ 
désigne un morphisme {\em adéquat} de schémas logarithmiques (\cite{agt} III.4.7). 
On désigne par $X^\circ$ le sous-schéma ouvert maximal de $X$
où la structure logarithmique $\cM_X$ est triviale~; c'est un sous-schéma ouvert de $X_\eta$.
On note $j\colon X^\circ\rightarrow X$ l'injection canonique. D'après (\cite{agt} III.4.2(iv)), l'immersion $j$ 
est schématiquement dominante et on a un isomorphisme canonique
\begin{equation}\label{hght2c}
\cM_X\stackrel{\sim}{\rightarrow}j_*(\co_{X^\circ}^\times)\cap \co_X.
\end{equation} 
En particulier, l'homomorphisme structural $\alpha\colon \cM_X\rightarrow \co_X$ est injectif.

Pour tout $X$-schéma $U$, on pose  
\begin{equation}\label{hght2a}
U^\circ=U\times_XX^\circ.
\end{equation} 
On note $\hbar\colon \oX\rightarrow X$ et $h\colon \oX^\circ\rightarrow X$ les morphismes canoniques \eqref{hght1c}, de sorte que 
l'on a $h=\hbar\circ j_\oX$. Pour alléger les notations, on pose
\begin{equation}\label{hght2b}
\tOmega^1_{X/S}=\Omega^1_{(X,\cM_X)/(S,\cM_S)},
\end{equation}
que l'on considère comme un faisceau de $X_\zar$ ou $X_\et$, selon le contexte (cf. \ref{notconv12}). 

On reprend les notations des objets associés à $f$ introduits dans §\ref{TFA}, en particulier ceux relatifs au topos de Faltings: $E$, $\tE$, 
$\ocB$...

\section{Torseur et extension de Higgs-Tate d'un schéma logarithmique affine}\label{taht}

\subsection{}\label{definf3}
On rappelle dans cette section la principale construction introduite dans (\cite{agt} II.10). 
Tous les anneaux des vecteurs de Witt considérés dans ce chapitre sont relatifs à $p$ \eqref{notconv1}. 
Suivant \eqref{eipo3a}, on désigne par $\co_{\oK^\flat}$ la limite projective du système projectif $(\co_\oK/p\co_\oK)_{\mN}$ 
dont les morphismes de transition sont les itérés de l'endomorphisme de Frobenius absolu de $\co_\oK/p\co_\oK$; 
\begin{equation}\label{definf3a}
\co_{\oK^\flat}= \underset{\underset{\mN}{\longleftarrow}}{\lim}\ \co_\oK/p\co_\oK.
\end{equation}
C'est un anneau de valuation non-discrète, de hauteur $1$, complet et 
parfait de caractéristique $p$ (cf. \ref{TFA101} et \ref{TFA102}). 
On fixe une suite $(p_n)_{n\geq 0}$ d'éléments de $\co_\oK$ telle que $p_0=p$ et $p_{n+1}^p=p_n$ (pour tout $n\geq 0$) et 
on note $\varpi$ l'élément associé de $\co_{\oK^\flat}$. On pose  
\begin{equation}\label{definf3b}
\xi=[\varpi]-p \in \rW(\co_{\oK^\flat}),
\end{equation}
où $[\ ]$ est le représentant multiplicatif. D'après \ref{eip4}, la suite 
\begin{equation}\label{definf3c}
0\longrightarrow \rW(\co_{\oK^\flat})\stackrel{\cdot \xi}{\longrightarrow} \rW(\co_{\oK^\flat})
\stackrel{\theta}{\longrightarrow} \co_C \longrightarrow 0,
\end{equation}
où $\theta$ est l'homomorphisme de Fontaine \eqref{eipo3d}, est exacte. Suivant \ref{eipo3}, on pose
\begin{equation}\label{definf3d}
\cA_2(\co_\oK)=\rW(\co_{\oK^\flat})/\ker(\theta)^2,
\end{equation}
et on note $\theta_2\colon \cA_2(\co_{\oK})\rightarrow \co_C$ l'homomorphisme induit par $\theta$.
On en déduit une suite exacte 
\begin{equation}\label{definf3e}
0\longrightarrow \co_C\stackrel{\cdot \xi}{\longrightarrow} \cA_2(\co_\oK)
\stackrel{\theta_2}{\longrightarrow} \co_C \longrightarrow 0,
\end{equation}
où on a encore noté $\cdot \xi$ le morphisme induit par la multiplication par $\xi$ dans $\cA_2(\co_\oK)$. 

L'idéal de carré nul $\ker(\theta_2)$ de $\cA_2(\co_\oK)$ est un $\co_C$-module libre de base $\xi$. Il sera noté $\xi\co_C$. 
Contrairement à $\xi$, il ne dépend pas du choix de la suite $(p_n)_{n\geq 0}$. 
On note $\xi^{-1}\co_C$ le $\co_C$-module dual de $\xi\co_C$. 
Pour tout $\co_C$-module $M$, on désigne les $\co_C$-modules $M\otimes_{\co_C}(\xi \co_C)$ 
et $M\otimes_{\co_C}(\xi^{-1} \co_C)$ simplement par $\xi M$ et $\xi^{-1} M$, respectivement.

Le groupe de Galois $G_K$ agit naturellement sur $\rW(\co_{\oK^\flat})$ par des automorphismes d'anneaux,
et l'homomorphisme $\theta$ est $G_K$-équivariant. On en déduit une action de $G_K$ 
sur $\cA_2(\co_{\oK})$ par des automorphismes d'anneaux tel que l'homomorphisme $\theta_2$ soit $G_K$-équivariant. 

\subsection{}\label{definf17}
On a un homomorphisme canonique 
\begin{equation}\label{definf17a}
\mZ_p(1)\rightarrow \co_{\oK^\flat}^\times.
\end{equation} 
Pour tout $\zeta\in \mZ_p(1)$, on note encore $\zeta$ son image dans $\co_{\oK^\flat}^\times$.  
Comme $\theta([\zeta]-1)=0$, on obtient un homomorphisme de groupes
\begin{equation}\label{definf17b}
\mZ_p(1)\rightarrow \cA_2(\co_\oK),\ \ \ 
\zeta\mapsto\log([\zeta])=[\zeta]-1,
\end{equation}
dont l'image est contenue dans $\ker(\theta_2)=\xi\co_C$. Celui-ci est clairement $\mZ_p$-linéaire. D'après (\cite{agt} II.9.18),
son image engendre l'idéal $p^{\frac{1}{p-1}}\xi\co_C$ de $\cA_2(\co_\oK)$, et le morphisme $\co_C$-linéaire induit
\begin{equation}\label{definf17c}
\co_C(1)\rightarrow p^{\frac{1}{p-1}}\xi \co_C
\end{equation}
est un isomorphisme.

\subsection{}\label{definf5}
Considérons le système projectif de monoïdes multiplicatifs $(\co_\oK)_{n\in \mN}$, 
où les morphismes de transition sont tous égaux à l'élévation à la puissance $p$-ième.
On note $Q_S$ le produit fibré du diagramme d'homomorphismes de monoïdes 
\begin{equation}\label{definf5a}
\xymatrix{
&{\co_K-\{0\}}\ar[d]\\
{\underset{\underset{x\mapsto x^p}{\longleftarrow}}{\lim}\ \co_\oK}\ar[r]&\co_\oK}
\end{equation}
où  la flèche verticale est l'homomorphisme canonique
et la flèche horizontale est la projection sur la première composante ({\em i.e.}, d'indice $0$).
On désigne par $\tau_S$ l'homomorphisme composé
\begin{equation}\label{definf5b}
\tau_S \colon Q_S\longrightarrow
\underset{\underset{x\mapsto x^p}{\longleftarrow}}{\lim}\ \co_\oK \longrightarrow \co_{\oK^\flat} \stackrel{[\ ]}{\longrightarrow} 
\rW(\co_{\oK^\flat}),
\end{equation} 
où $[\ ]$ est le représentant multiplicatif et les autres flèches sont les morphismes canoniques. 
Il résulte aussitôt des définitions que le diagramme 
\begin{equation}\label{definf5c}
\xymatrix{
Q_S\ar[r]\ar[d]_{\tau_S}&{\co_K-\{0\}}\ar[d]\\
{\rW(\co_{\oK^\flat})}\ar[r]^-(0.4)\theta&{\co_C}}
\end{equation}
où les flèches non libellées sont les morphismes canoniques, est commutatif. 
Par ailleurs, le groupe de Galois $G_K$ agit naturellement sur le monoïde $Q_S$,
et l'homomorphisme $\tau_S$ est $G_K$-équivariant.

\subsection{}\label{definf6}
On pose 
\begin{equation}\label{definf6a}
\tS=\Spec(\cA_2(\co_\oK)),
\end{equation} 
que l'on munit de la structure logarithmique $\cM_{\tS}$
associée à la structure pré-logarithmique définie par l'homomorphisme $Q_S\rightarrow \cA_2(\co_\oK)$ induit par $\tau_S$ \eqref{definf5b}. 

On fixe une uniformisante $\pi$ de $\co_K$ et une suite $(\pi_n)_{n\geq 0}$ 
d'éléments de $\co_\oK$ telle que $\pi_0=\pi$ et $\pi_{n+1}^p=\pi_n$ (pour tout $n\geq 0$) et 
on note $\upi$ l'élément associé de $\co_{\oK^\flat}$ \eqref{definf3a}. 
En vertu de (\cite{agt} II.9.7), $\cM_{\tS}$ est la structure logarithmique sur 
$\tS$ associée à la structure pré-logarithmique définie par l'homomorphisme 
\begin{equation}\label{definf6b}
\mN\rightarrow \cA_2(\co_\oK),\ \ \ 1\mapsto [\upi].
\end{equation}
En effet, notant $\tpi$ l'élément de $Q_S$ défini par ses projections \eqref{definf5a}
\begin{equation}\label{definf6c}
(\pi_n)_{n\in \mN}\in \underset{\underset{\mN}{\longleftarrow}}{\lim}\ \co_\oK \ \ \ {\rm et}\ \ \ 
\pi \in \co_K-\{0\},
\end{equation}
on a $\tau_S(\tpi)=[\upi]\in \rW(\co_{\oK^\flat})$. 
Le schéma logarithmique $(\tS,\cM_{\tS})$ est donc fin et saturé, 
et $\theta_2$ \eqref{definf3d} induit une immersion fermée exacte 
\begin{equation}\label{definf6d}
i_S\colon (\coS,\cM_\coS)\rightarrow (\tS,\cM_{\tS}).
\end{equation}

\subsection{}\label{taht1}
On suppose dans la suite de cette section que le morphisme $f\colon (X,\cM_X)\rightarrow (S,\cM_S)$ \eqref{hght2} admet une carte adéquate 
\eqref{cad1} (\cite{agt} III.4.4), que $X=\Spec(R)$ est affine et connexe et que $X_s$ est non-vide. 
Ces conditions correspondent à celles fixées dans \ref{tpcg1}. 
On suppose, de plus, qu'il existe une carte fine et saturée $M\rightarrow \Gamma(X,\cM_X)$ pour $(X,\cM_X)$
induisant un isomorphisme 
\begin{equation}\label{taht1a}
M\stackrel{\sim}{\rightarrow} \Gamma(X,\cM_X)/\Gamma(X,\co^\times_X).
\end{equation}
Celle-ci est a priori indépendante de la carte adéquate donnée. On notera que le monoïde $\Gamma(X,\cM_X)$ est saturé (\cite{agt} II.5.6).
On pose $d=\dim(X/S)$ et \eqref{hght2b}
\begin{equation}\label{taht1b}
\tOmega^1_{R/\co_K}=\tOmega^1_{X/S}(X),
\end{equation}
qui est un $R$-module libre de rang $d$ \eqref{tpcg20b}. 

\subsection{}\label{definf12}
On munit  $\coX=X\times_S\coS$ \eqref{hght1c} de la structure logarithmique $\cM_\coX$ image inverse de $\cM_X$. 
On a alors un isomorphisme canonique
\begin{equation}\label{definf12a}
(\coX,\cM_{\coX})\stackrel{\sim}{\rightarrow}(X,\cM_X)\times_{(S,\cM_S)}(\coS,\cM_\coS),
\end{equation} 
le produit étant indifféremment pris dans la catégorie des schémas logarithmiques ou 
dans celle des schémas logarithmiques fins. 
On se donne une $(\tS,\cM_{\tS})$-déformation lisse $(\tX,\cM_\tX)$ de $(\coX,\cM_{\coX})$,
c'est-à-dire un morphisme lisse de schémas logarithmiques fins $(\tX,\cM_\tX)\rightarrow (\tS, \cM_{\tS})$ 
et un $(\coS,\cM_{\coS})$-isomorphisme 
\begin{equation}\label{definf12b}
(\coX,\cM_{\coX})\stackrel{\sim}{\rightarrow}
(\tX,\cM_\tX)\times_{(\tS, \cM_{\tS})}(\coS,\cM_{\coS}).
\end{equation}
Comme $f$ est lisse et que $\coX$ est affine, une telle déformation existe et est unique à isomorphisme près en vertu de (\cite{kato1} 3.14).

\subsection{}\label{taht100}
Soit $\oy$ un point géométrique de $\oX^\circ$ (cf. \eqref{hght1c} et \eqref{hght2a}). 
Le schéma $\oX$ étant localement irréductible d'après \ref{cad7}(iii),  
il est la somme des schémas induits sur ses composantes irréductibles. On note $\oX^\star$
la composante irréductible de $\oX$ contenant $\oy$. 
De même, $\oX^\circ$ est la somme des schémas induits sur ses composantes irréductibles
et $\oX^{\star \circ}=\oX^\star\times_{X}X^\circ$ est la composante irréductible de $\oX^\circ$ contenant $\oy$. 
On pose 
\begin{equation}\label{taht100c}
R_1=\Gamma(\oX^\star,\co_\oX).
\end{equation}
On désigne par $\Delta$ le groupe profini $\pi_1(\oX^{\star \circ},\oy)$ et par $(V_i)_{i\in I}$ le revêtement universel normalisé de $\oX^{\star \circ}$ en $\oy$ \eqref{notconv11}.
Pour chaque $i\in I$, on note $\oX_i$ la fermeture intégrale de $\oX$ dans $V_i$.
Les schémas $(\oX_i)_{i\in I}$ forment alors un système projectif filtrant. On pose   
\begin{equation}\label{taht100d}
\oR=\underset{\underset{i\in I}{\longrightarrow}}{\lim}\ \Gamma(\oX_i,\co_{\oX_i}).
\end{equation} 
C'est un anneau intègre et normal \eqref{TFA10}, sur lequel agit naturellement $\Delta$ par des homomorphismes d'anneaux et l'action est discrète. 
On note $\hoR$ (resp. $\hRun$) le séparé complété $p$-adique de $\oR$ (resp. $R_1$).  
Les $\co_C$-modules $\hRun$ et $\hoR$ sont plats et l'homomorphisme canonique $\hRun\rightarrow \hoR$ est injectif (\cite{agt} II.6.14). 

On pose 
\begin{equation}\label{taht100e}
\mX=\Spec(\oR)\ \ \ {\rm et} \ \ \ \hmX=\Spec(\hoR)
\end{equation}
que l'on munit des structures logarithmiques images inverses de $\cM_X$, 
notées respectivement $\cM_\mX$ et $\cM_\hmX$. 
Les actions de $\Delta$ sur $\oR$ et $\hoR$ induisent des actions à gauche sur 
les schémas logarithmiques $(\mX,\cM_\mX)$ et $(\hmX,\cM_\hmX)$. 
Munissant $(\coX,\cM_\coX)$ de l'action triviale de $\Delta$ (cf. \ref{definf12}), on a un morphisme canonique $\Delta$-équivariant 
\begin{equation}\label{taht100f}
(\hmX,\cM_\hmX)\rightarrow (\coX,\cM_\coX).
\end{equation}

\subsection{}\label{taht2}
Suivant \eqref{eipo3a}, on désigne par $\oR^\flat$ la limite projective du système projectif $(\oR/p\uoR)_{\mN}$ 
dont les morphismes de transition sont les itérés de l'endomorphisme de Frobenius absolu de $\oR/p\oR$; 
\begin{equation}
\oR^\flat= \underset{\underset{\mN}{\longleftarrow}}{\lim}\ \oR/p\oR.
\end{equation}
D'après \ref{eip4} et \ref{tpcg7}, la suite 
\begin{equation}\label{taht2b}
0\longrightarrow \rW(\oR^\flat)\stackrel{\cdot \xi}{\longrightarrow} \rW(\oR^\flat)
\stackrel{\theta}{\longrightarrow} \hoR \longrightarrow 0,
\end{equation}
où $\theta$ est l'homomorphisme de Fontaine \eqref{eipo3d} et $\xi\in \rW(\co_{\oK^\flat})$ est l'élément défini dans \eqref{definf3b}, est exacte. 
Suivant \ref{eipo3}, on pose
\begin{equation}\label{taht2a}
\cA_2(\oR)=\rW(\oR^\flat)/\ker (\theta)^2,
\end{equation}
et on note $\theta_2\colon \cA_2(\oR)\rightarrow \hoR$ l'homomorphisme induit par $\theta$.
On en déduit une suite exacte 
\begin{equation}\label{taht2c}
0\longrightarrow \hoR\stackrel{\cdot \xi}{\longrightarrow} \cA_2(\oR)
\stackrel{\theta_2}{\longrightarrow} \hoR \longrightarrow 0,
\end{equation}
où on a encore noté $\cdot \xi$ le morphisme induit par la multiplication par $\xi$ dans $\cA_2(\oR)$. 

Le groupe $\Delta$ agit naturellement sur $\rW(\oR^\flat)$ par des automorphismes d'anneaux,
et l'homomorphisme $\theta$ est $\Delta$-équivariant. On en déduit une action de $\Delta$ 
sur $\cA_2(\oR)$ par des automorphismes d'anneaux telle que l'homomorphisme $\theta_2$ soit 
$\Delta$-équivariant.

\subsection{}\label{taht20}
Considérons le système projectif de monoïdes multiplicatifs $(\oR)_{n\in \mN}$, 
où les morphismes de transition sont tous égaux à l'élévation à la puissance $p$-ième.
On note $Q_X$ le produit fibré du diagramme d'homomorphismes de monoïdes 
\begin{equation}\label{taht20a}
\xymatrix{
&{\Gamma(X,\cM_X)}\ar[d]\\
{\underset{\underset{\mN}{\longleftarrow}}{\lim}\ \oR}\ar[r]&\oR}
\end{equation}
où  la flèche verticale est l'homomorphisme canonique
et la flèche horizontale est la projection sur la première composante ({\em i.e.}, d'indice $0$).
On rappelle que les limites projectives sont représentables dans la catégorie des monoïdes commutatifs 
et qu'elles commutent au foncteur d'oubli à valeurs dans la catégorie des ensembles. 
Il existe un unique homomorphisme 
\begin{equation}\label{taht20e}
\mZ_p(1)\rightarrow Q_X
\end{equation}
tel que l'homomorphisme induit $\mZ_p(1)\rightarrow \Gamma(X,\cM_X)$ soit constant (de valeur $1$)
et l'homomorphisme induit
\begin{equation}\label{taht20f}
\mZ_p(1)\rightarrow \underset{\underset{\mN}{\longleftarrow}}{\lim}\ \oR
\end{equation} 
soit l'homomorphisme canonique. Le groupe $\Delta$ agit naturellement sur le monoïde $Q_X$, et l'homomorphisme
\eqref{taht20e} est $\Delta$-équivariant.  

On désigne par $\tau_X$ l'homomorphisme composé
\begin{equation}\label{taht20b}
\tau_X \colon Q_X\longrightarrow
\underset{\underset{\mN}{\longleftarrow}}{\lim}\ \oR \longrightarrow \oR^\flat \stackrel{[\ ]}{\longrightarrow} 
\rW(\oR^\flat),
\end{equation} 
où $[\ ]$ est le représentant multiplicatif et les autres flèches sont les morphismes canoniques. 
Il est clairement $\Delta$-équivariant.  

Il résulte aussitôt des définitions que le diagramme 
\begin{equation}\label{taht20c}
\xymatrix{
Q_X\ar[r]\ar[d]_{\tau_X}&{\Gamma(X,\cM_X)}\ar[d]\\
{\rW(\oR^\flat)}\ar[r]^-(0.4)\theta&{\hoR}}
\end{equation}
où les flèches non libellées sont les morphismes canoniques, est commutatif. On a un homomorphisme canonique 
$Q_S\rightarrow Q_X$ \eqref{definf5} qui s'insère dans un diagramme commutatif 
\begin{equation}\label{taht20d}
\xymatrix{
Q_S\ar[r]\ar[d]_{\tau_S}&Q_X\ar[d]^{\tau_X}\\
{\rW(\co_{\oK^\flat})}\ar[r]&{\rW(\oR^\flat)}}
\end{equation}

\begin{prop}\label{taht22}
\
\begin{itemize}
\item[{\rm (i)}] Le monoïde $Q_X$ est intègre. 
\item[{\rm (ii)}] L'homomorphisme canonique $\mZ_p(1)\rightarrow Q_X$  \eqref{taht20e} est injectif.
\item[{\rm (iii)}] La projection canonique $Q_X\rightarrow \Gamma(X,\cM_X)$ 
identifie $ \Gamma(X,\cM_X)$ au quotient de $Q_X$ par le sous-monoïde $\mZ_p(1)$.
\end{itemize}
\end{prop}

On désigne par $\nu$ l'application d'élévation à la puissance $p$-ième de $\oR$,  
et pour tout entier $n\geq 0$, par $Q_{X,n}$ le produit fibré du diagramme d'homomorphismes de monoïdes 
\begin{equation}\label{taht22a}
\xymatrix{
&{\Gamma(X,\cM_X)}\ar[d]\\
\oR\ar[r]^{\nu^n}&\oR}
\end{equation}
L'homomorphisme canonique $\Gamma(X,\cM_X)\rightarrow \oR$ vérifie les hypothèses de \ref{tfkum5}. 
En effet, l'homomorphisme $\Gamma(X,\cM_X)\rightarrow R$ est injectif d'après \eqref{hght2c},
et comme $X$ est connexe,  l'homomorphisme $R\rightarrow \oR$ est injectif. 
Par suite, le monoïde $Q_{X,n}$ est intègre et la projection canonique $Q_{X,n}\rightarrow \Gamma(X,\cM_X)$ 
identifie $\Gamma(X,\cM_X)$ au quotient de $Q_{X,n}$ par le sous-monoïde $\mu_{p^n}(\co_\oK)$ d'après  \ref{tfkum5}. 
On a donc une suite exacte de groupes abéliens 
\begin{equation}\label{taht22b}
0\rightarrow \mu_{p^n}(\co_\oK)\rightarrow Q_{X,n}^\gp\rightarrow \Gamma(X,\cM_X)^\gp\rightarrow 0.
\end{equation}
On observera que le diagramme canonique
\begin{equation}\label{taht22c}
\xymatrix{
{Q_{X,n}}\ar[d]\ar[r]&{\Gamma(X,\cM_X)}\ar[d]\\
{Q_{X,n}^\gp}\ar[r]&{\Gamma(X,\cM_X)^\gp}}
\end{equation}
est cartésien (\cite{ogus} I 1.1.6).

Comme les limites projectives commutent aux limites projectives, l'homomorphisme canonique 
\begin{equation}\label{taht22d}
Q_X\rightarrow  \underset{\underset{\mN}{\longleftarrow}}{\lim}\ Q_{X,n}
\end{equation}
est un isomorphisme. Les propositions (i) et (ii) s'ensuivent aussitôt.  

Considérons le diagramme commutatif
\begin{equation}\label{taht22e}
\xymatrix{
0\ar[r]&{\mZ_p(1)}\ar@{=}[d]\ar[r]&{Q_X^\gp}\ar[r]\ar[d]^u&{\Gamma(X,\cM_X)^\gp}\ar@{=}[d]\ar[r]&0\\
0\ar[r]&{\mZ_p(1)}\ar[r]&{\underset{\underset{\mN}{\longleftarrow}}{\lim}\ Q_{X,n}^\gp}\ar[r]&{\Gamma(X,\cM_X)^\gp}\ar[r]&0}
\end{equation}
où la suite inférieure est exacte d'après \eqref{taht22b}. Le foncteur ``groupe associé à un monoïde'' étant adjoint à gauche du foncteur d'injection
canonique de la catégorie des monoïdes commutatifs dans la catégorie des groupes abéliens, 
il commute naturellement aux limites inductives; mais il ne commute pas à toutes les limites projectives. 
Toutefois, il résulte de \eqref{taht22d} et du fait que les monoïdes $Q_X$ et $Q_{X,n}$ ($n\geq 0$) sont intègres que $u$ est injectif. 
Comme l'homomorphisme $Q_X\rightarrow \Gamma(X,\cM_X)$ est surjectif, il s'ensuit que la suite supérieure de \eqref{taht22e} est exacte,
et par suite que $u$ est un isomorphisme.  Prenant alors la limite projective des diagrammes cartésiens \eqref{taht22c}, 
on déduit que le diagramme 
\begin{equation}\label{taht22f}
\xymatrix{
{Q_X}\ar[d]\ar[r]&{\Gamma(X,\cM_X)}\ar[d]\\
{Q_X^\gp}\ar[r]&{\Gamma(X,\cM_X)^\gp}}
\end{equation}
est cartésien. La proposition (iii) s'ensuit compte tenu de \eqref{taht22e} et (\cite{ogus} I 1.1.6).

\subsection{}\label{taht21}
On pose 
\begin{equation}\label{taht21a}
\tmX=\Spec(\cA_2(\oR))
\end{equation} 
que l'on munit de la structure logarithmique $\cM_{\tmX}$ associée à la structure pré-logarithmique définie 
par l'homomorphisme $Q_X\rightarrow \cA_2(\oR)$ induit par $\tau_X$ \eqref{taht20b}. 
D'après (\cite{agt} II.9.7) et compte tenu des hypothèses de \ref{taht1}, le schéma logarithmique $(\tmX,\cM_{\tmX})$ est fin et saturé, et 
l'homomorphisme $\theta_2$ induit une immersion fermée exacte 
\begin{equation}\label{taht21b}
i_X\colon (\hmX,\cM_\hmX)\rightarrow (\tmX,\cM_{\tmX}).
\end{equation}

Les actions de $\Delta$ sur $\cA_2(\oR)$ et $Q_X$ induisent une action à gauche sur
le schéma logarithmique $(\tmX,\cM_{\tmX})$. L'immersion fermée $i_X$ est $\Delta$-équivariante.

D'après \eqref{taht20d}, on a un morphisme canonique
\begin{equation}\label{taht21c}
(\tmX,\cM_{\tmX})\rightarrow (\tS,\cM_\tS)
\end{equation}
qui s'insère dans un diagramme commutatif \eqref{definf6d}
\begin{equation}\label{taht21d}
\xymatrix{
{(\hmX,\cM_\hmX)}\ar[r]^{i_X}\ar[d]&{(\tmX,\cM_{\tmX})}\ar[d]\\
{(\coS,\cM_{\coS})}\ar[r]^{i_S}&{(\tS,\cM_\tS)}}
\end{equation}

\subsection{}\label{taht5}
Rappelons les définitions du torseur et de l'extension de Higgs-Tate associés à $(\tX,\cM_\tX)$ (\cite{agt} II.10.3 et II.10.4). 
On désigne par $\hmX_\zar$ le topos de Zariski de $\hmX$ \eqref{taht100e} et par $\trT$ le $\co_\hmX$-module associé au $\hoR$-module 
\eqref{taht1b}
\begin{equation}\label{taht5a}
\rT=\Hom_{\hoR}(\tOmega^1_{R/\co_K}\otimes_R\hoR,\xi\hoR).
\end{equation} 
Soient $U$ un ouvert de Zariski de $\hmX$, $\tU$ l'ouvert correspondant de $\tmX$ \eqref{taht21b}. 
On désigne par $\cL(U)$ 
l'ensemble des morphismes représentés par des flèches pointillées qui complètent  le diagramme canonique \eqref{definf12b}
\begin{equation}\label{taht5d}
\xymatrix{
{(U,\cM_\hmX|U)}\ar[rr]^-(0.5){i_X\times_{\tmX}\tU}\ar[d]&&{(\tU,\cM_{\tmX}|\tU)}\ar@{.>}[d]\ar@/^2pc/[dd]\\
{(\coX,\cM_{\coX})}\ar[rr]\ar[d]&&{(\tX,\cM_\tX)}\ar[d]\\
{(\coS,\cM_\coS)}\ar[rr]^-(0.5){i_S}&&{(\tS,\cM_{\tS})}}
\end{equation}
de façon à le laisser commutatif. D'après (\cite{agt} II.5.23),
le foncteur $U\mapsto \cL(U)$ est un $\trT$-torseur de $\hmX_\zar$.  
On l'appelle {\em torseur de Higgs-Tate} associé à $(\tX,\cM_\tX)$.  
On désigne par $\cF$ le $\hoR$-module des fonctions affines sur $\cL$ (cf. \cite{agt} II.4.9). 
Celui-ci s'insère dans une suite exacte canonique 
\begin{equation}\label{taht5e}
0\rightarrow \hoR\rightarrow \cF\rightarrow \xi^{-1}\tOmega^1_{R/\co_K} \otimes_R \hoR\rightarrow 0,
\end{equation} 
dite {\em extension de Higgs-Tate} associée à $(\tX,\cM_\tX)$.

\subsection{}\label{taht6}
On munit $\hmX$ de l'action naturelle à gauche de $\Delta$~; pour tout $g\in \Delta$, 
l'automorphisme de $\hmX$ défini par $g$, que l'on note aussi $g$, est induit par l'automorphisme $g^{-1}$ de $\hoR$. 
On considère $\trT$ comme un $\co_\hmX$-module $\Delta$-équivariant au moyen 
de la donnée de descente correspondant au $\hRun$-module 
$\Hom_{\hRun}(\tOmega^1_{R/\co_K}\otimes_R\hRun,\xi\hRun)$ (cf. \cite{agt} II.4.18). Pour tout $g\in \Delta$, on a donc 
un isomorphisme canonique de $\co_\hmX$-modules
\begin{equation}\label{taht6a}
\tau_g^\trT\colon \trT\stackrel{\sim}{\rightarrow} g^*(\trT).
\end{equation}

L'action naturelle à gauche de $\Delta$ sur le schéma logarithmique $(\tmX,\cM_{\tmX})$ \eqref{taht21}
induit sur le $\trT$-torseur $\cL$ une structure $\Delta$-équivariante (cf. \cite{agt} II.4.18), 
autrement dit, elle induit pour tout $g\in \Delta$, un isomorphisme $\tau_g^\trT$-équivariant
\begin{equation}\label{taht6c}
\tau^{\cL}_g\colon \cL\stackrel{\sim}{\rightarrow} g^*(\cL);
\end{equation}
ces isomorphismes étant soumis à des relations de compatibilité (cf. \cite{agt} II.4.16). 
En effet, pour tout ouvert de Zariski $U$ de $\hmX$, on prend pour 
\begin{equation}\label{taht6d}
\tau^{\cL}_g(U)\colon \cL(U)\stackrel{\sim}{\rightarrow} \cL(g(U))
\end{equation}
l'isomorphisme défini de la façon suivante. 
Soient $\tU$ l'ouvert de $\tmX$ correspondant à $U$, $\mu\in \cL(U)$ que l'on considère comme un morphisme
\begin{equation}\label{taht6e}
\mu\colon (\tU,\cM_{\tmX}|\tU)\rightarrow (\tX,\cM_\tX).
\end{equation}
Comme $i_X$ \eqref{taht21b} et le morphisme \eqref{taht100f} sont $\Delta$-équivariants, le morphisme composé
\begin{equation}\label{taht6f}
(g(\tU),\cM_{\tmX}|g(\tU))\stackrel{g^{-1}}{\longrightarrow} (\tU,\cM_{\tmX}|\tU)
\stackrel{\mu}{\longrightarrow} (\tX,\cM_\tX)
\end{equation}
prolonge le morphisme canonique $(g(U),\cM_\hmX|g(U))\rightarrow (\tX,\cM_\tX)$. 
Il correspond à l'image de $\mu$ par $\tau^{\cL}_g(U)$. On vérifie aussitôt que 
le morphisme $\tau^{\cL}_g$ ainsi défini est un isomorphisme $\tau_g^\trT$-équivariant et que 
ces isomorphismes vérifient les relations de compatibilité requises dans (\cite{agt} II.4.16).

D'après (\cite{agt} II.4.21), les structures $\Delta$-équivariantes sur $\trT$ et $\cL$ induisent une structure
$\Delta$-équivariante sur le $\co_\hmX$-module associé à $\cF$, ou, ce qui revient au même,
une action $\hoR$-semi-linéaire de $\Delta$ sur $\cF$, telle que les morphismes de la suite \eqref{taht5e} soient 
$\Delta$-équivariants. 

D'après (\cite{agt} (II.10.4.11)), pour tous $\mu\in \cL(\hmX)$ et $\beta\in \cF$, on a 
\begin{equation}\label{taht6j}
(g^{-1}(\beta))(\mu)=g^{-1}(\beta({^g\mu})),
\end{equation}
où ${^g\mu}$ est défini par le morphisme composé 
\begin{equation}\label{taht6k}
(\tmX,\cM_{\tmX})\stackrel{g^{-1}}{\longrightarrow} (\tmX,\cM_{\tmX})
\stackrel{\mu}{\longrightarrow} (\tX,\cM_\tX).
\end{equation}

\subsection{}\label{taht15}
On désigne par $\tQ_X$ le produit fibré du diagramme d'homomorphismes canoniques de monoïdes \eqref{taht20a}
\begin{equation}\label{taht15a}
\xymatrix{
&&{\Gamma(\tX,\cM_\tX)}\ar[d]\\
Q_X\ar[r]&{\Gamma(X,\cM_X)}\ar[r]&{\Gamma(\coX,\cM_\coX)}}
\end{equation}
Il existe un unique homomorphisme 
\begin{equation}\label{taht15b}
\mZ_p(1)\rightarrow \tQ_X
\end{equation}
tel que l'homomorphisme composé $\mZ_p(1)\rightarrow \Gamma(\tX,\cM_\tX)$ soit constant (de valeur $1$)
et l'homomorphisme composé $\mZ_p(1)\rightarrow Q_X$ soit l'homomorphisme \eqref{taht20e}. 
L'action naturelle de $\Delta$ sur $Q_X$ et son action triviale sur $\Gamma(\tX,\cM_\tX)$ induisent une action sur $\tQ_X$ (cf. \ref{taht20}).

Soit $T'$ un élément de $\tQ_X$, d'images canoniques $T\in Q_X$, $t\in \Gamma(X,\cM_X)$ et $t'\in \Gamma(\tX,\cM_\tX)$. 
Tout élément $\mu\in \cL(\hmX)$ \eqref{taht5} détermine un morphisme de schémas logarithmiques que l'on note encore 
\begin{equation}\label{taht15d}
\mu\colon (\tmX,\cM_{\tmX})\rightarrow (\tX,\cM_{\tX}),
\end{equation}
et en particulier un homomorphisme de monoïdes $\mu^\flat\colon \mu^{-1}(\cM_{\tX})\rightarrow \cM_{\tmX}$.
D'après (\cite{ogus} IV 2.1.2),  il existe un et un unique élément $b_\mu \in  \xi\hoR$ tel que 
\begin{equation}\label{taht15e}
\mu^{\flat}(\mu^{-1}(t'))=(1+b_\mu) T \in \Gamma(\tmX,\cM_{\tmX}).
\end{equation}
L'application $\mu\mapsto b_\mu$ est une fonction affine sur $\cL$ à valeurs dans $\xi \hoR$ 
de terme linéaire $d\log(t)\in \tOmega^1_{R/\co_K}\otimes_R\hoR$ (cf. \cite{agt} II.4.7). Elle définit donc naturellement 
un élément de $\xi\cF$ \eqref{taht5e}, que l'on note $d\log(T')\in \xi\cF$, au-dessus de $d\log(t)$.
On voit aussitôt que l'application 
\begin{equation}\label{taht15f}
d\log\colon \tQ_X\rightarrow \xi\cF
\end{equation}
ainsi définie est un homomorphisme de monoïdes. Il est $\Delta$-équivariant d'après \eqref{taht6j}.
Le diagramme
\begin{equation}\label{taht15g}
\xymatrix{
&\mZ_p(1)\ar[r]\ar[d]_{-\log([\ ])}&{\tQ_X}\ar[r]\ar[d]^{d\log}&{\Gamma(X,\cM_X)}\ar[d]^{d\log}&\\
0\ar[r]&{\xi\hoR}\ar[r]&{\xi\cF}\ar[r]&{\tOmega^1_{R/\co_K}\otimes_R\hoR}\ar[r]&0}
\end{equation}
où l'homomorphisme $\log([\ ])$ est induit par \eqref{definf17b}, est clairement commutatif.

\subsection{}\label{taht16}
Soit $t\in \Gamma(X,\cM_X)$. 
On désigne par $Q_X(t)$ l'image inverse de $t$ par la projection canonique $Q_X\rightarrow \Gamma(X,\cM_X)$ \eqref{taht20a},  
qui est un $\mZ_p(1)$-torseur de $\bB_\Delta$ d'après \ref{taht22}, et par $\xi\cF(t)$ la fibre au-dessus de $d\log(t)\otimes 1$ dans l'extension \eqref{taht5e}
\begin{equation}\label{taht16a}
0\rightarrow \xi\hoR\rightarrow \xi\cF\rightarrow \tOmega^1_{R/\co_K} \otimes_R \hoR\rightarrow 0,
\end{equation} 
qui est un $\xi\hoR$-torseur de $\bB_\Delta$. 
La donnée d'une section $t'\in \Gamma(\tX,\cM_\tX)$ ayant même image que $t$ dans $\Gamma(\coX,\cM_\coX)$
définit une application $Q_X(t)\rightarrow \tQ_X$ dont le composé avec la projection canonique sur $Q_X$ est l'injection canonique et 
le composé avec la projection canonique sur $\Gamma(\tX,\cM_\tX)$ est constant de valeur $t'$. 
Cette application est équivariante pour les actions naturelles de $\mZ_p(1)$ (resp. $\Delta$). 
L'application \eqref{taht15f}, induit alors une application $\Delta$-équivariante 
\begin{equation}\label{taht16b}
Q_X(t)\rightarrow \xi\cF(t).
\end{equation}
Cette application est aussi équivariante relativement à l'homomorphisme 
\begin{equation}\label{taht16c}
-\log([\ ])\colon \mZ_p(1)\rightarrow \xi\hoR
\end{equation}
induit par l'homomorphisme $\log([\ ])$ défini dans  \eqref{definf17b}.

\section{Théorie de Kummer du topos de Faltings annelé}\label{tfkum}

\subsection{}\label{tfkum2}
Avec les notations de \ref{hght3} et \ref{hght2}, pour toute extension finie $L$ de $K$, on pose 
\begin{eqnarray}
(X_L,\cM_{X_L})&=&(X,\cM_X)\times_{(S,\cM_S)}(S_L,\cM_{S_L}),\label{tfkum2a}\\
(\oX,\cL_{\oX})&=&(X,\cM_X)\times_{(S,\cM_S)}(\oS,\cL_{\oS}),\label{tfkum2b}
\end{eqnarray}
les produits étant pris dans la catégorie des schémas logarithmiques, de sorte que 
$X_L=X\times_SS_L$ et $\oX=X\times_S\oS$ \eqref{hght1c}. Pour alléger les notations, on pose
\begin{equation}\label{tfkum2c}
\tOmega^1_{\oX/\oS}=\Omega^1_{(\oX,\cL_\oX)/(\oS,\cL_\oS)},
\end{equation}
que l'on considère comme un faisceau de $\oX_\zar$ ou $\oX_\et$, selon le contexte (cf. \ref{notconv12}). 
On a un isomorphisme canonique \eqref{hght2b}
\begin{equation}\label{tfkum2d}
\tOmega^1_{X/S}\otimes_{\co_X}\co_{\oX}\stackrel{\sim}{\rightarrow}\tOmega^1_{\oX/\oS}.
\end{equation}
On note $\cL_\oX^\gp$ le groupe de $\oX_\et$ associé au monoïde $\cL_\oX$
(\cite{agt} II.5.3) et 
\begin{equation}\label{tfkum2e}
d\log \colon \cL_\oX\rightarrow \tOmega^1_{\oX/\oS}
\end{equation}
la dérivation logarithmique universelle. Celle-ci induit un morphisme $\co_\oX$-linéaire et surjectif (\cite{ogus} IV 1.2.11)
\begin{equation}\label{tfkum2f}
\cL_\oX^\gp\otimes_\mZ\co_\oX\rightarrow \tOmega^1_{\oX/\oS}.
\end{equation}

Pour tout entier $n\geq 1$, on pose 
\begin{eqnarray}
\tOmega^1_{X_n/S_n}=\tOmega^1_{X/S}\otimes_{\co_X}\co_{X_n},\label{tfkum2fg}\\
\tOmega^1_{\oX_n/\oS_n}=\tOmega^1_{X/S}\otimes_{\co_X}\co_{\oX_n},\label{tfkum2g}
\end{eqnarray}
que l'on considère aussi comme des faisceaux de $X_{s,\zar}$ ou $X_{s,\et}$, selon le contexte (cf. \ref{TFA5}). 
On note $\oa\colon X_s\rightarrow \oX$ l'immersion fermée canonique \eqref{TFA5a}.
Compte tenu de \eqref{notconv12e}, le morphisme \eqref{tfkum2f}
induit un morphisme $\co_{\oX_n}$-linéaire et surjectif de $X_{s,\et}$
\begin{equation}\label{tfkum2h}
\oa^*(\cL_\oX^\gp)\otimes_\mZ\co_{\oX_n}\rightarrow \tOmega^1_{\oX_n/\oS_n}.
\end{equation}

\begin{lem}\label{tfkum3}
\
\begin{itemize}
\item[{\rm (i)}] Le monoïde $\cL_\oX$ est intègre. 
\item[{\rm (ii)}] L'homomorphisme structural $\oalpha\colon \cL_\oX\rightarrow \co_\oX$ est injectif.
\item[{\rm (iii)}] L'homomorphisme $j^*_{\oX}(\oalpha)\colon j^*_{\oX}(\cL_\oX)\rightarrow \co_{\oX^\circ}$ se factorise
à travers un isomorphisme 
\begin{equation}
j^*_{\oX}(\cL_\oX)\stackrel{\sim}{\rightarrow} \co_{\oX^\circ}^\times.
\end{equation} 
\end{itemize}
\end{lem}

En effet, le schéma logarithmique $(\oX,\cL_\oX)$ est la limite projective des schémas logarithmiques $(X_L,\cM_{X_L})$,
indexée par l'ensemble filtrant des sous-$K$-extensions finies $L$ de $\oK$. 
Notant $\hbar_L\colon \oX\rightarrow X_L$ le morphisme canonique, l'homomorphisme 
\begin{equation}\label{tfkum3a}
\underset{\underset{L\subset \oK}{\longrightarrow}}{\lim}\ \hbar_L^*(\cM_{X_L}) \rightarrow \cL_{\oX},
\end{equation}
où la limite est prise sur l'ensemble filtrant des sous-$K$-extensions finies de $\oK$, est un isomorphisme. 
Par ailleurs, l'homomorphisme canonique  
\begin{equation}\label{tfkum3b}
\underset{\underset{L\subset \oK}{\longrightarrow}}{\lim}\ \hbar_L^{-1}(\co_{X_L})\rightarrow \co_\oX,
\end{equation}
où la limite est prise sur l'ensemble filtrant des sous-$K$-extensions finies de $\oK$, 
est un isomorphisme dans $\oX_\et$. En effet, c'est un isomorphisme dans $\oX_\zar$ d'après (\cite{ega4} 8.2.12),
et il l'est donc dans $\oX_\et$ en vertu de (\cite{ega1n} 0.6.1.6) et (\cite{raynaud1} I §~3 prop.~1).

(i) Cela résulte de \eqref{tfkum3a}. 

(ii) Pour toute sous-$K$-extension finie $L$ de $\oK$,
le morphisme canonique $(X_L,\cM_{X_L})\rightarrow (S_L,\cM_{S_L})$ est adéquat d'après (\cite{agt} III.4.8).
Par suite, l'homomorphisme structural $\cM_{X_L}\rightarrow \co_{X_L}$ est un monomorphisme d'après (\cite{agt} III.4.2(iv)). 
La proposition s'ensuit compte tenu de \eqref{tfkum3a} et \eqref{tfkum3b}. 

(iii) Pour toute sous-$K$-extension finie $L$ de $\oK$,
on voit aussitôt que $X_L^\circ=X_L\times_XX^\circ$ est le sous-schéma ouvert maximal de 
$X_L$ où la structure logarithmique $\cM_{X_L}$ est triviale.
L'homomorphisme $j^*_{X_L}(\cM_{X_L})\rightarrow \co_{X_L^\circ}$ se factorise donc
à travers un isomorphisme $ j^*_{X_L}(\cM_{X_L})\stackrel{\sim}{\rightarrow} \co_{X_L^\circ}^\times$. 
La proposition s'ensuit compte tenu de \eqref{tfkum3a} et \eqref{tfkum3b}.

\subsection{}\label{tfkum6}
Pour tout entier $n\geq 0$,  on désigne par $\cQ_n$ le monoïde de $\tE$ défini par le diagramme cartésien 
de la catégorie des monoïdes de $\tE$ \eqref{notconv7}
\begin{equation}\label{tfkum6a}
\xymatrix{
{\cQ_n}\ar[d]\ar[r]&{\sigma^*(\hbar_*(\cL_\oX))}\ar[d]^\ell\\
{\ocB}\ar[r]^{\nu^n}&{\ocB}}
\end{equation}
où $\nu$ est l'endomorphisme d'élévation à la puissance $p$-ième de $\ocB$ \eqref{TFA2}
et $\ell$ est le composé des homomorphismes canoniques 
\begin{equation}\label{tfkum6b}
\xymatrix{
{\sigma^*(\hbar_*(\cL_\oX))}\ar[rr]^-(0.5){\sigma^*(\hbar_*(\oalpha))}&&{\sigma^*(\hbar_*(\co_\oX))}\ar[r]&{\ocB}},
\end{equation}
le second étant l'homomorphisme \eqref{TFA2c}.
Les monoïdes $(\cQ_n)_{n\in \mN}$ forment naturellement un système projectif.
Compte tenu de \eqref{TFA5b} et \eqref{TFA66c}, le diagramme \eqref{tfkum6a} induit par image inverse par le plongement
$\delta$ \eqref{TFA66a} un diagramme cartésien 
de la catégorie des monoïdes de $\tE_s$
\begin{equation}\label{tfkum6c}
\xymatrix{
{\delta^*(\cQ_n)}\ar[d]\ar[r]&{\sigma^*_s(\oa^*(\cL_\oX))}\ar[d]^{\delta^*(\ell)}\\
{\delta^*(\ocB)}\ar[r]^{\delta^*(\nu^n)}&{\delta^*(\ocB)}}
\end{equation}

\begin{prop}\label{tfkum7}
Soit $n$ un entier $\geq 0$. 
\begin{itemize}
\item[{\rm (i)}] Le monoïde $\delta^*(\cQ_n)$ est intègre. 
\item[{\rm (ii)}] L'homomorphisme canonique $\mu_{p^n,\tE_s}\rightarrow \delta^*(\cQ_n)$ est injectif \eqref{tfkum4}.
\item[{\rm (iii)}] La projection canonique $\delta^*(\cQ_n)\rightarrow \sigma^*_s(\oa^*(\cL_\oX))$ 
identifie $\sigma^*_s(\oa^*(\cL_\oX))$ au quotient de $\delta^*(\cQ_n)$ par le sous-monoïde $\mu_{p^n,\tE_s}$ \eqref{notconv7}.
\end{itemize}
\end{prop} 

Soient $(\oy\rightsquigarrow \ox)$ un point de $X_\et\gtimes_{X_\et}\oX^\circ_\et$ \eqref{topfl17} tel que $\ox$ soit
au-dessus de $s$, $\rho(\oy\rightsquigarrow \ox)$ son image dans $\tE$ par le morphisme $\rho$ \eqref{TFA6e}. 
On note encore $\ox$ le point géométrique $\oa(\ox)$ de $\oX$ \eqref{TFA5a}. 
Compte tenu de (\cite{agt} (VI.10.18.1)) et  \eqref{TFA5b}, la fibre de $\ell$ \eqref{tfkum6b} 
en $\rho(\oy\rightsquigarrow \ox)$ s'identifie à un homomorphisme
\begin{equation}
\ell_{\rho(\oy\rightsquigarrow \ox)} \colon \cL_{\oX,\ox}\rightarrow \ocB_{\rho(\oy \rightsquigarrow \ox)}.
\end{equation}
Montrons que celui-ci remplit les conditions de \ref{tfkum5}.
D'après (\cite{agt} III.10.10), $\ocB_{\rho(\oy \rightsquigarrow \ox)}$ 
est un anneau normal et strictement local (et en particulier intègre) et l'homomorphisme
\begin{equation}
\co_{\oX,\ox}\rightarrow \ocB_{\rho(\oy \rightsquigarrow \ox)}
\end{equation}
fibre de l'homomorphisme canonique $\sigma^*(\hbar_*(\co_\oX))\rightarrow \ocB$ \eqref{TFA2c} en $\rho(\oy \rightsquigarrow \ox)$,
est injectif et local. Comme l'homomorphisme structural 
$\oalpha\colon \cL_\oX\rightarrow \co_\oX$ est injectif d'après \ref{tfkum3}(ii), 
$\ell_{\rho(\oy\rightsquigarrow \ox)}$ est injectif.

On désigne par $\fV_\ox$ la catégorie des $X$-schémas étales $\ox$-pointés, 
et par $\fW_\ox$ la sous-catégorie pleine de $\fV_\ox$ formée des objets $(U,\fp\colon \ox\rightarrow U)$ 
tels que le schéma $U$ soit affine. Ce sont des catégories cofiltrantes, et le foncteur d'injection canonique
$\fW^\circ_\ox\rightarrow \fV^\circ_\ox$ est cofinal (\cite{sga4} I 8.1.3(c)).
Avec les notations de \ref{TFA11}, on a un isomorphisme canonique \eqref{TFA11c} 
\begin{equation}
\ocB_{\rho(\oy\rightsquigarrow \ox)} \stackrel{\sim}{\rightarrow} 
\underset{\underset{(U,\fp)\in \fW_\ox^\circ}{\longrightarrow}}{\lim}\ \oR^{\oy}_U.
\end{equation} 
On a clairement \eqref{TFA5b}
\begin{equation}
\cL_{\oX,\ox} \stackrel{\sim}{\rightarrow} 
\underset{\underset{(U,\fp)\in \fW_\ox^\circ}{\longrightarrow}}{\lim}\ \Gamma(\oU,\cL_\oX).
\end{equation} 
De plus, l'homomorphisme $\ell_{\rho(\oy\rightsquigarrow \ox)}$ s'identifie à la limite inductive des homomorphismes composés
\begin{equation}
\xymatrix{
{\Gamma(\oU,\cL_\oX)}\ar[r]^{\oalpha}&{\Gamma(\oU,\co_{\oX})}\ar[r]&{\oR^\oy_U}},
\end{equation}
où la seconde flèche est l'homomorphisme canonique \eqref{TFA2c}. 
Pour tout $(U,\fp)\in \ob(\fW_\ox)$ et tout $t\in \Gamma(\oU,\cL_\oX)$, 
comme $\oalpha(t)$ est inversible sur $\oU^\circ$ d'après \ref{tfkum3}(iii), le $\oU$-schéma 
\begin{equation}
\Spec(\co_\oU[T]/(T^{p^n}-\oalpha(t)))
\end{equation}
est fini et étale au-dessus de $\oU^\circ$. L'anneau $\oR^{\oy}_U$  étant intègre et normal \eqref{TFA10}, 
il résulte aussitôt de sa définition \eqref{TFA9c}  que l'équation $T^{p^n}=\oalpha(t)$ admet une racine et donc 
$p^n$ racines distinctes dans $\oR_U^\oy$. 
On en déduit que pour tout $m\in \cL_{\oX,\ox}$, l'équation $T^{p^n}=\ell_{\rho(\oy\rightsquigarrow \ox)}(m)$ 
admet $p^n$ racines distinctes dans $\ocB_{\rho(\oy\rightsquigarrow \ox)}$. La proposition s'ensuit en vertu de 
\ref{tfkum5} et \ref{tf21}.

\begin{cor}\label{tfkum8}
Pour tout entier $n\geq 0$, la suite d'homomorphismes canoniques
\begin{equation}\label{tfkum8a}
0\rightarrow \mu_{p^n,\tE_s}\rightarrow \delta^*(\cQ^\gp_n)\rightarrow \sigma_s^*(\oa^*(\cL^\gp_\oX))\rightarrow 0
\end{equation}
est exacte. 
\end{cor} 

En effet, le foncteur d'injection canonique de la catégorie des groupes abéliens de $\tE_s$ 
dans celle des monoïdes de $\tE_s$ admet pour adjoint à gauche le foncteur ``groupe associé'' 
$\cN\mapsto \cN^\gp$. Ce dernier commute donc aux limites inductives.
L'exactitude de la suite \eqref{tfkum8a} au centre et à droite résulte alors de \ref{tfkum7}(iii) et (\cite{agt} (II.5.7.1)). 
L'exactitude de cette suite à gauche est une conséquence de \ref{tfkum7}(i)-(ii).

\begin{cor}\label{tfkum9}
Soient $\ox$ un point géométrique de $X$ au-dessus de $s$, $\uX$ le localisé strict de $X$ en $\ox$, 
$\varphi_\ox\colon \tE\rightarrow \uoX^\circ_\fet$ le foncteur canonique \eqref{TFA14a}, $n$ un entier $\geq 0$. 
On désigne encore par $\ox$ le point géométrique $\oa(\ox)$ de $\oX$ \eqref{TFA5a}. 
Alors, l'image du diagramme \eqref{tfkum6a} par le foncteur $\varphi_\ox$ s'identifie à un diagramme cartésien 
de la catégorie des monoïdes de $\uoX^\circ_\fet$
\begin{equation}\label{tfkum9a}
\xymatrix{
{\varphi_\ox(\cQ_n)}\ar[r]\ar[d]&{\cL_{\oX,\ox}}\ar[d]^{\varphi_\ox({\ell})}\\
{\varphi_\ox(\ocB)}\ar[r]^{\nu^n}&{\varphi_\ox(\ocB)}}
\end{equation}
où on a encore noté $\nu$ l'endomorphisme d'élévation à la puissance $p$-ième de $\varphi_\ox(\ocB)$.
De plus, le monoïde $\varphi_\ox(\cQ_n)$ est intègre, et  la projection canonique $\varphi_\ox(\cQ_n)\rightarrow \cL_{\oX,\ox}$
identifie $\cL_{\oX,\ox}$ au quotient de $\varphi_\ox(\cQ_n)$ par le sous-monoïde $\mu_{p^n,\uoX^\circ_\fet}$ \eqref{notconv7}. 
\end{cor}

En effet, on a un isomorphisme canonique 
$\varphi_\ox(\sigma^*(\hbar_*(\cL_\oX)))\stackrel{\sim}{\rightarrow}\cL_{\oX,\ox}$
d'après \eqref{TFA14b} et \eqref{TFA5b}; d'où la première assertion. 
Comme $\uoX^\circ$ est intègre et non-vide \eqref{TFA15},  
considérant un point $(\oy\rightsquigarrow \ox)$ de $X_\et\gtimes_{X_\et}\oX^\circ_\et$,
la seconde assertion résulte de la preuve de \ref{tfkum7} compte tenu de la description \eqref{TFA12e} 
du foncteur fibre de $\tE$ associé au point $\rho(\oy\rightsquigarrow \ox)$. 

\subsection{}\label{tfkum10}
Reprenons les notations de \ref{definf3}. 
Pour tout entier $n\geq 1$, on note 
\begin{equation}\label{tfkum10a}
\partial_n\colon \oa^*(\cL^\gp_\oX)\rightarrow \rR^1\sigma_{s*}(\mu_{p^n,\tE_s})
\end{equation}
l'homomorphisme composé du morphisme d'adjonction 
$\oa^*(\cL^\gp_\oX)\rightarrow \sigma_{s*}(\sigma^*_s(\oa^*(\cL^\gp_\oX)))$
et du cobord de la suite exacte \eqref{tfkum8a}. Celui-ci induit deux morphismes  $\co_{\oX_n}$-linéaires
\begin{eqnarray}
\oa^*(\cL^\gp_\oX)\otimes_\mZ\co_{\oX_n}&\rightarrow& \rR^1\sigma_{n*}(\ocB_n(1)),\label{tfkum10b}\\
\oa^*(\cL^\gp_\oX)\otimes_\mZ\co_{\oX_n}&\rightarrow& \rR^1\sigma_{n*}(\xi\ocB_n),\label{tfkum10c}
\end{eqnarray}
le second étant défini au moyen du composé des morphismes canoniques \eqref{definf17c}
\begin{equation}\label{tfkum10d}
\co_C(1)\stackrel{\sim}{\rightarrow}p^{\frac{1}{p-1}}\xi\co_C\rightarrow \xi \co_C.
\end{equation}

\begin{teo}\label{tfkum11}
Soit $n$ un entier $\geq 1$.  
\begin{itemize}
\item[{\rm (i)}] Le morphisme \eqref{tfkum10c} se factorise à travers le morphisme
surjectif \eqref{tfkum2h} et induit un morphisme $\co_{\oX_n}$-linéaire de $X_{s,\et}$ 
\begin{equation}\label{tfkum11a}
\xi^{-1}\tOmega^1_{\oX_n/\oS_n}\rightarrow \rR^1\sigma_{n*}(\ocB_n).
\end{equation}
\item[{\rm (ii)}] Il existe un et un unique homomorphisme de $\co_{\oX_n}$-algèbres graduées de $X_{s,\et}$
\begin{equation}\label{tfkum11b}
\wedge (\xi^{-1}\tOmega^1_{\oX_n/\oS_n})\rightarrow \oplus_{i\geq 0}\rR^i\sigma_{n*}(\ocB_n)
\end{equation}
dont la composante de degré un est le morphisme \eqref{tfkum11a}. De plus, 
son noyau est annulé par $p^{\frac{2d}{p-1}}\fm_\oK$ et son conoyau est annulé par $p^{\frac{2d+1}{p-1}}\fm_\oK$, 
où $d=\dim(X/S)$.  
\end{itemize}
\end{teo}

La preuve de cet énoncé sera donnée dans \ref{tfkum19} après quelques résultats préliminaires.

\subsection{}\label{tfkum12}
Soient $(\oy\rightsquigarrow \ox)$ un point de $X_\et\gtimes_{X_\et}\oX^\circ_\et$ \eqref{topfl17} tel que $\ox$ soit
au-dessus de $s$, $\uX$ le localisé strict de $X$ en $\ox$. 
On note encore $\ox$ le point géométrique $\oa(\ox)$ de $\oX$ \eqref{TFA5a}. 
D'après (\cite{agt} III.3.7), $\uoX=\uX\times_S\oS$ est normal et strictement local (et en particulier intègre); on peut donc l'identifier 
au localisé strict de $\oX$ en $\ox$. Le $X$-morphisme $\oy\rightarrow \uX$
définissant $(\oy\rightsquigarrow \ox)$ se relève en un $\oX^\circ$-morphisme $\oy\rightarrow \uoX^\circ$ et 
induit donc un point géométrique de $\uoX^\circ$ que l'on note aussi (abusivement) $\oy$.
On désigne par 
\begin{equation}\label{tfkum12a}
\varphi_\ox\colon \tE\rightarrow \uoX^\circ_\fet
\end{equation} 
le foncteur canonique \eqref{TFA14a}, par $\bB_{\pi_1(\uoX^\circ,\oy)}$ 
le topos classifiant du groupe profini $\pi_1(\uoX^\circ,\oy)$ et par 
\begin{equation}\label{tfkum12b}
\psi_{\oy}\colon \uoX^\circ_\fet \stackrel{\sim}{\rightarrow}\bB_{\pi_1(\uoX^\circ,\oy)}
\end{equation}
le foncteur fibre de $\uoX^\circ_\fet$ en $\oy$ \eqref{notconv11c}. 
Compte tenu de \eqref{TFA14b}, \eqref{TFA5b} et \eqref{TFA12g}, 
l'homomorphisme $\psi_\oy(\varphi_\ox(\ell))$ \eqref{tfkum6b} 
s'identifie à  l'homomorphisme composé
\begin{equation}\label{tfkum12c}
\xymatrix{
{\cL_{\oX,\ox}}\ar[r]^{\oalpha_{\ox}}&{\co_{\oX,\ox}}\ar[r]&{\oR^\oy_{\uX}}},
\end{equation}
où $\oR^\oy_{\uX}$ est l'algèbre définie dans \eqref{TFA12f}  et  
la seconde flèche est l'homomorphisme canonique.
Pour tout entier $n\geq 0$, l'image du diagramme \eqref{tfkum6a} par le foncteur $\psi_\oy\circ \varphi_\ox$ 
s'identifie donc à un diagramme cartésien de la catégorie des monoïdes de  $\bB_{\pi_1(\uoX^\circ,\oy)}$
\begin{equation}\label{tfkum12d}
\xymatrix{
{\psi_\oy(\varphi_\ox(\cQ_n))}\ar[r]\ar[d]&{\cL_{\oX,\ox}}\ar[d]^{\psi_\oy(\varphi_\ox({\ell}))}\\
{\oR^\oy_{\uX}}\ar[r]^{\nu^n}&{\oR^\oy_{\uX}}}
\end{equation}
où on a encore noté $\nu$ l'endomorphisme d'élévation à la puissance $p$-ième de $\oR^\oy_{\uX}$.

Pour tout $t\in \cL_{\oX,\ox}$, la section $\oalpha_\ox(t)$  est inversible sur $\uoX^\circ$ d'après \ref{tfkum3}(iii). Le $\uoX$-schéma 
\begin{equation}\label{tfkum12g}
\cT_n(t)=\Spec(\co_{\uoX}[T]/(T^{p^n}-\oalpha_\ox(t)))
\end{equation}
est donc fini et étale sur $\uoX^\circ$. En particulier, $\cT^\circ_n(t)=\cT_n(t)\times_XX^\circ$ induit un $\mu_{p^n}(\co_\oK)$-torseur 
de $\uoX^\circ_\fet$. On désigne par $\cT_n(t)_\oy$ la fibre de $\cT_n(t)$ au-dessus de $\oy$ et par $\cQ_n(t)$ la fibre 
de la projection canonique $\psi_\oy(\varphi_\ox(\cQ_n))\rightarrow \cL_{\oX,\ox}$ \eqref{tfkum12d} au-dessus de $t$.
On a un isomorphisme canonique de $\mu_{p^n}(\co_\oK)$-torseurs de $\bB_{\pi_1(\uoX^\circ,\oy)}$
\begin{equation}\label{tfkum12h}
\cT_n(t)_\oy\stackrel{\sim}{\rightarrow}\cQ_n(t).
\end{equation}

\begin{lem}\label{tfkum13}
Sous les hypothèses de \ref{tfkum12}, pour tout entier $n\geq 0$, l'application 
\begin{equation}\label{tfkum13a}
\ttt_n\colon \cL_{\oX,\ox}\rightarrow \rH^1(\uoX^\circ_\fet,\mu_{p^n,\uoX^\circ_\fet}), \ \ \ t\mapsto [\cT^\circ_n(t)],
\end{equation}
est un homomorphisme de monoïdes \eqref{tfkum4}, et le diagramme 
\begin{equation}\label{tfkum13b}
\xymatrix{
{\cL_{\oX,\ox}}\ar[r]^-(0.5){\ttt_n}\ar[d]&{\rH^1(\uoX^\circ_\fet,\mu_{p^n,\uoX^\circ_\fet})}\ar[r]_-(0.5)\sim^-(0.5)v&
{\rH^1(\uoX^\circ_\fet,\varphi_\ox(\delta_*(\mu_{p^n,\tE_s})))}\\
{\cL_{\oX,\ox}^\gp}\ar[r]^-(0.5){\partial_{n,\ox}}&{\rR^1\sigma_{s*}(\mu_{p^n,\tE_s})_\ox}\ar[r]_-(0.5)\sim^-(0.5)w&
{\rR^1\sigma_*(\delta_*(\mu_{p^n,\tE_s}))_\ox}\ar[u]_u}
\end{equation}
où $\partial_{n,\ox}$ est la fibre de $\partial_n$ \eqref{tfkum10a} en $\ox$, $u$ est l'isomorphisme \eqref{TFA14d},  
$v$ est induit par l'isomorphisme canonique $\varphi_\ox\stackrel{\sim}{\rightarrow} \varphi_\ox\delta_*\delta^*$ \eqref{TFA15a},
$w$ est induit par l'isomorphisme \eqref{TFA66d} et la flèche non libellée est l'homomorphisme canonique, est commutatif. 
\end{lem}

Considérons la suite exacte de groupes abéliens de $\tE$ 
\begin{equation}\label{tfkum13e}
0\rightarrow \delta_*(\mu_{p^n,\tE_s})\rightarrow \delta_*(\delta^*(\cQ_n^\gp))\rightarrow 
\delta_*(\sigma^*_s(\oa^*(\cL_\oX^\gp)))\rightarrow 0
\end{equation}
déduite de la suite exacte \eqref{tfkum8a} par image directe par le plongement $\delta$. 
D'après (\cite{agt} VI.10.30(iii)), le diagramme 
\begin{equation}
\xymatrix{
{\sigma_*(\delta_*(\sigma_s^*(\oa^*\cL^\gp_\oX)))_\ox}\ar[r]^-(0.5){u'}_-(0.5)\sim\ar[d]_a&
{\rH^0(\uoX^\circ_\fet,\varphi_\ox(\delta_*(\sigma_s^*(\oa^*\cL^\gp_\oX))))}\ar[d]^b\\
{\rR^1\sigma_*(\delta_*(\mu_{p^n,\tE_s}))_\ox}
\ar[r]^-(0.5)u_-(0.5)\sim&{\rH^1(\uoX^\circ_\fet,\varphi_\ox(\delta_*(\mu_{p^n,\tE_s})))}}
\end{equation}
où les flèches horizontales sont les isomorphismes \eqref{TFA14d} et les flèches verticales sont les bords 
des suites exactes longues de cohomologie, est commutatif. Par ailleurs, le diagramme 
\begin{equation}
\xymatrix{
{\cL^\gp_{\oX,\ox}}\ar[r]^-(0.5)i\ar[rd]_{\partial_{n,\ox}}&{\sigma_{s*}(\sigma_s^*(\oa^*\cL^\gp_\oX))_\ox}
\ar[r]^-(0.5){w'}_-(0.5)\sim\ar[d]^c&{\sigma_*(\delta_*(\sigma_s^*(\oa^*\cL^\gp_\oX)))_\ox}\ar[d]^a\\
&{\rR^1\sigma_{s*}(\mu_{p^n,\tE_s})_\ox}\ar[r]^-(0.5){w}_-(0.5)\sim&{\rR^1\sigma_*(\delta_*(\mu_{p^n,\tE_s}))_\ox}}
\end{equation}
où $i$ est induit par le morphisme d'adjonction $\id\rightarrow \sigma_{s*}\sigma_s^*$, 
$c$ est le bord de la suite exacte longue de cohomologie déduite de la suite \eqref{tfkum8a}
et $w$ et $w'$ sont induits par l'isomorphisme \eqref{TFA66d}, est commutatif. 

D'après \eqref{TFA15a}, \eqref{TFA14b} et \eqref{TFA66c}, l'image de \eqref{tfkum13e} par le foncteur $\varphi_\ox$ s'identifie à une suite exacte de 
groupes abéliens de $\uoX^\circ_\fet$
\begin{equation}\label{tfkum13c}
0\rightarrow \mu_{p^n,\uoX^\circ_\fet}\rightarrow \varphi_\ox(\cQ_n^\gp)\rightarrow \cL_{\oX,\ox}^\gp\rightarrow 0. 
\end{equation}
Celle-ci induit un homomorphisme 
\begin{equation}\label{tfkum13d}
\ttt'_n\colon \cL^\gp_{\oX,\ox}\rightarrow\rH^1(\uoX^\circ_\fet,\mu_{p^n,\uoX^\circ_\fet}),
\end{equation}
qui s'identifie donc à l'homomorphisme $b$
\begin{equation}
\xymatrix{
{\cL^\gp_{\oX,\ox}}\ar[d]_{\ttt'_n}\ar[r]^-(0.5){v'}_-(0.5)\sim&
{\rH^0(\uoX^\circ_\fet,\varphi_\ox(\delta_*(\sigma_s^*(\oa^*\cL^\gp_\oX))))}\ar[d]^b\\
{\rH^1(\uoX^\circ_\fet,\mu_{p^n,\uoX^\circ_\fet})}\ar[r]^-(0.5)v_-(0.5)\sim&
{\rH^1(\uoX^\circ_\fet,\varphi_\ox(\delta_*(\mu_{p^n,\tE_s})))}}
\end{equation}
On vérifie que $v'^{-1}\circ u'\circ w'\circ i$ est l'identité. 

Il résulte de \ref{tfkum9} et (\cite{ogus} I.1.1.6) que le diagramme canonique
\begin{equation}\label{tfkum13f}
\xymatrix{
{\varphi_\ox(\cQ_n)}\ar[d]\ar[r]&{\cL_{\oX,\ox}}\ar[d]\\
{\varphi_\ox(\cQ^\gp_n)}\ar[r]&{\cL^\gp_{\oX,\ox}}}
\end{equation}
est cartésien. Compte tenu de \eqref{tfkum12h}, on en déduit  
que $\ttt_n$ est le composé de $\ttt'_n$ et de l'homomorphisme canonique $\cL_{\oX,\ox}\rightarrow \cL_{\oX,\ox}^\gp$;
d'où la proposition.

\subsection{}\label{tfkum14} 
Conservons les hypothèses et notations de \ref{tfkum12}; supposons, de plus, les conditions suivantes remplies~:
\begin{itemize}
\item[(i)] Le schéma $X$ est affine et connexe et le morphisme $f$ \eqref{hght2} admet une carte adéquate \eqref{cad1} (\cite{agt} III.4.4).
\item[(ii)] Il existe une carte fine et saturée $M\rightarrow \Gamma(X,\cM_X)$ pour $(X,\cM_X)$
induisant un isomorphisme 
\begin{equation}\label{TFSLA31a}
M\stackrel{\sim}{\rightarrow} \Gamma(X,\cM_X)/\Gamma(X,\co^\times_X).
\end{equation}
\end{itemize}
Ces conditions correspondent à celles fixées dans \ref{taht1}. 
Notons $R$ l'anneau de $X$ et posons $\tOmega^1_{R/\co_K}=\tOmega^1_{X/S}(X)$ \eqref{taht1b}.
Le schéma $\oX$ étant localement irréductible, 
il est la somme des schémas induits sur ses composantes irréductibles. 
On désigne par $\oX^\star$  la composante irréductible de $\oX$ contenant $\oy$. De même, $\oX^\circ$
est la somme des schémas induits sur ses composantes irréductibles, et $\oX^{\star\circ}=\oX^\star\times_XX^\circ$ 
est la composante irréductible de $\oX^\circ$ contenant $\oy$.  On note $\bB_{\pi_1(\oX^{\star \circ},\oy)}$
le topos classifiant du groupe profini $\pi_1(\oX^{\star \circ},\oy)$.

Reprenons les notations de \ref{taht}.
On munit  $\coX=X\times_S\coS$ \eqref{hght1c} de la structure logarithmique $\cM_\coX$ image inverse de $\cM_X$
et on se donne une $(\tS,\cM_{\tS})$-déformation lisse $(\tX,\cM_\tX)$ de $(\coX,\cM_{\coX})$ (cf. \ref{definf12}). 
Comme $f$ est lisse et que $\coX$ est affine, une telle déformation existe et est unique à isomorphisme près en vertu de (\cite{kato1} 3.14).
On désigne par $\hoR^\oy_X$ le complété $p$-adique de l'anneau $\oR^\oy_X$ \eqref{TFA9c} et par
\begin{equation}\label{tfkum14b} 
0\rightarrow \hoR^\oy_X\rightarrow \cF\rightarrow \xi^{-1} \tOmega^1_{R/\co_K}\otimes_R\hoR^\oy_X\rightarrow 0
\end{equation}
l'extension de Higgs-Tate associée à $(\tX,\cM_\tX)$ \eqref{taht5e};  
c'est une extension de $\hoR^\oy_X$-représentations continues de $\pi_1(\oX^{\star\circ},\oy)$. 
Pour tout entier $n\geq 0$, on désigne par  
\begin{equation}\label{tfkum14c} 
0\rightarrow \xi(\oR^\oy_X/p^n\oR^\oy_X)\rightarrow \xi(\cF/p^n\cF)\rightarrow \tOmega^1_{R/\co_K}\otimes_R
(\oR^\oy_X/p^n\oR^\oy_X)\rightarrow 0
\end{equation}
la suite exacte déduite de \eqref{tfkum14b}, par
\begin{equation}\label{tfkum14e}
A_n\colon \Gamma(\oX^\star,\tOmega^1_{\oX_n/\oS_n})
\rightarrow \rH^1(\pi_1(\oX^{\star\circ},\oy),\xi(\oR^\oy_X/p^n\oR^\oy_X))
\end{equation}
le morphisme $\Gamma(\oX^\star,\co_\oX)$-linéaire induit, et par $Q_n$ le produit fibré du diagramme d'homomorphismes de monoïdes de $\bB_{\pi_1(\oX^{\star \circ},\oy)}$
\begin{equation}\label{tfkum14d} 
\xymatrix{
&{\Gamma(X,\cM_X)}\ar[d]\\
{\oR^\oy_X}\ar[r]^{\nu^n}&{\oR^\oy_X}}
\end{equation}
où on a encore noté $\nu$ l'endomorphisme d'élévation à la puissance $p$-ième de $\oR^\oy_X$.
On observera que $A_n$ ne dépend pas du choix de $(\tX,\cM_\tX)$ (\cite{agt} II.10.10). 
L'homomorphisme canonique $\Gamma(X,\cM_X)\rightarrow \oR^\oy_X$ vérifie les hypothèses de \ref{tfkum5}. 
En effet, l'homomorphisme $\Gamma(X,\cM_X)\rightarrow R$ est injectif d'après \eqref{hght2c},
et comme $X$ est connexe,  l'homomorphisme $R\rightarrow \oR^\oy_X$ est injectif. 
Par suite, le monoïde $Q_n$ est intègre et la projection canonique $Q_n\rightarrow \Gamma(X,\cM_X)$ 
identifie $\Gamma(X,\cM_X)$ au quotient de $Q_n$ par le sous-monoïde $\mu_{p^n}(\co_\oK)$ d'après  \ref{tfkum5}. 
On a donc une suite exacte de groupes abéliens de $\bB_{\pi_1(\oX^{\star\circ},\oy)}$
\begin{equation}\label{tfkum14f}
0\rightarrow \mu_{p^n}(\co_\oK)\rightarrow Q_n^\gp\rightarrow \Gamma(X,\cM_X)^\gp\rightarrow 0.
\end{equation}
Elle induit un homomorphisme 
\begin{equation}\label{tfkum14g}
B_n\colon \Gamma(X,\cM_X)\rightarrow \rH^1(\pi_1(\oX^{\star\circ},\oy),\mu_{p^n}(\co_\oK)).
\end{equation}
On observera que le diagramme canonique
\begin{equation}\label{tfkum14h}
\xymatrix{
{Q_n}\ar[d]\ar[r]&{\Gamma(X,\cM_X)}\ar[d]\\
{Q_n^\gp}\ar[r]&{\Gamma(X,\cM_X)^\gp}}
\end{equation}
est cartésien (\cite{ogus} I 1.1.6). 
On note encore 
\begin{equation}\label{tfkum14hi}
d\log \colon \cM_X\rightarrow \tOmega^1_{X/S}
\end{equation}
la dérivation logarithmique universelle et 
\begin{equation}\label{tfkum14i}
\log([\ ])\colon \mZ_p(1)\rightarrow \xi\hoR^\oy_X
\end{equation}
l'homomorphisme induit par l'homomorphisme $\log([\ ])$ défini dans \eqref{definf17b}.

Certains des objets considérés dans ce numéro ont déjà été introduits dans la preuve de \ref{taht22}. Nous les avons rappelés pour la commodité du lecteur.

\begin{lem}\label{tfkum15}
Sous les hypothèses de \ref{tfkum14}, pour tout entier $n\geq 0$, le diagramme 
\begin{equation}\label{tfkum15a}
\xymatrix{
{\Gamma(X,\cM_X)}\ar[r]^-(0.5){B_n}\ar[d]_{d\log}&
{\rH^1(\pi_1(\oX^{\star\circ},\oy),\mu_{p^n}(\co_\oK))}\ar[d]^{-\log([\ ])}\\
{\Gamma(\oX^\star,\tOmega^1_{\oX_n/\oS_n})}\ar[r]^-(0.5){A_n}&
{\rH^1(\pi_1(\oX^{\star\circ},\oy),\xi(\oR^\oy_X/p^n\oR^\oy_X))}}
\end{equation}
est commutatif. 
\end{lem}

On notera d'abord que le morphisme $\coX\rightarrow \tX$ étant un homéomorphisme universel, on peut identifier 
les topos étales de $\coX$ et $\tX$. D'après (\cite{fkato} 3.6 ou \cite{ogus} IV 2.1.2), 
$\cM_{\coX}$ est le quotient de $\cM_\tX$ par le sous-monoïde $1+\xi\co_{\coX}$ (cf. \cite{ogus} I 1.1.5). 
Comme $\coX$ est affine, on en déduit que l'homomorphisme canonique 
\begin{equation}\label{tfkum15c}
\Gamma(\tX,\cM_\tX)\rightarrow \Gamma(\coX,\cM_{\coX})
\end{equation}
est surjectif. 
Soit $t\in \Gamma(X,\cM_X)$. On désigne par $Q_n(t)$ l'image inverse de $t$ par la projection canonique 
$Q_n\rightarrow \Gamma(X,\cM_X)$, qui est un $\mu_{p^n}(\co_\oK)$-torseur de $\bB_{\pi_1(\oX^{\star\circ},\oy)}$, 
et par $\xi\cF_n(t)$ la fibre au-dessus de $d\log(t)\otimes 1$ dans l'extension \eqref{tfkum14c},
qui est un $\xi(\oR_X^\oy/p^n\oR_X^\oy)$-torseur de $\bB_{\pi_1(\oX^{\star\circ},\oy)}$. 
Compte tenu de \eqref{tfkum14h}, $B_n(t)$ est la classe du $\mu_{p^n}(\co_\oK)$-torseur $Q_n(t)$ de 
$\bB_{\pi_1(\oX^{\star\circ},\oy)}$. D'autre part, $A_n(d\log(t)\otimes 1)$
est la classe du $\xi(\oR_X^\oy/p^n\oR_X^\oy)$-torseur $\xi\cF_n(t)$ de $\bB_{\pi_1(\oX^{\star\circ},\oy)}$.

Avec les notations de \ref{taht20}, on a un homomorphisme canonique $\pi_1(\oX^{\star\circ},\oy)$-équivariant 
$Q_X\rightarrow Q_n$. Celui-ci induit une application $\pi_1(\oX^{\star\circ},\oy)$-équivariante $Q_X(t)\rightarrow Q_n(t)$ \eqref{taht16}. 
Cette application est aussi équivariante relativement à l'homomorphisme canonique $\mZ_p(1)\rightarrow \mu_{p^n}(\co_\oK)$. 
On en déduit un isomorphisme 
\begin{equation}\label{tfkum15d}
Q_X(t)\wedge^{\mZ_p(1)}\mu_{p^n}(\co_\oK)\stackrel{\sim}{\rightarrow} Q_n(t),
\end{equation}
où la source est le $\mu_{p^n}(\co_\oK)$-torseur de $\bB_{\pi_1(\oX^{\star\circ},\oy)}$ déduit de $Q_X(t)$
par extension de son groupe structural, autrement dit, le quotient de 
$Q_X(t)\times\mu_{p^n}(\co_\oK)$ par l'action diagonale de $\mZ_p(1)$ (\cite{giraud2} III 1.4.6). 
De même, le morphisme canonique $\xi\cF\rightarrow \xi(\cF/p^n\cF)$ induit une application 
$\pi_1(\oX^{\star\circ},\oy)$-équivariante $\xi\cF(t)\rightarrow \xi\cF_n(t)$ \eqref{taht16}. 
Cette application est aussi équivariante relativement à l'homomorphisme canonique 
$\xi\hoR^{\oy}_X\rightarrow \xi(\hoR^{\oy}_X/p^n\hoR^{\oy}_X)$. On en déduit un isomorphisme 
\begin{equation}\label{tfkum15e}
\xi\cF(t)\wedge^{\xi\hoR^{\oy}_X}\xi(\hoR^{\oy}_X/p^n\hoR^{\oy}_X)\stackrel{\sim}{\rightarrow} \xi\cF_n(t),
\end{equation}
où la source est le $\xi(\hoR^{\oy}_X/p^n\hoR^{\oy}_X)$-torseur de $\bB_{\pi_1(\oX^{\star\circ},\oy)}$ déduit de $\xi\cF(t)$
par extension de son groupe structural. 

Choisissant un élément $t'\in \Gamma(\tX,\cM_\tX)$ ayant même image que $t$ dans $\Gamma(\coX,\cM_\coX)$ \eqref{tfkum15c}, 
l'application \eqref{taht16b} et les isomorphismes \eqref{tfkum15d} et \eqref{tfkum15e} 
induisent une application $\pi_1(\oX^{\star\circ},\oy)$-équivariante
\begin{equation}\label{tfkum15b} 
Q_n(t)\rightarrow \xi\cF_n(t).
\end{equation}
Cette application est aussi équivariante relativement à l'homomorphisme $-\log([\ ])$ \eqref{tfkum14i}; d'où la proposition.

\begin{prop}\label{tfkum16} 
Sous les hypothèses de \ref{tfkum12}, pour tout entier $n\geq 1$, il existe un unique morphisme $\co_{\oX,\ox}$-linéaire
\begin{equation}\label{tfkum16a} 
\phi_n\colon \tOmega^1_{\oX_n/\oS_n,\ox}\rightarrow \rH^1(\pi_1(\uoX^\circ,\oy),\xi(\oR^\oy_{\uX}/p^n\oR^\oy_{\uX}))
\end{equation}
qui s'insère dans le diagramme commutatif 
\begin{equation}\label{tfkum16b} 
\xymatrix{
{\cL_{\oX,\ox}}\ar[r]^-(0.5){\ttt_n}\ar[d]_{d\log}&{\rH^1(\pi_1(\uoX^\circ,\oy),\mu_{p^n}(\co_\oK))}\ar[d]\ar[d]^{-\log([\ ])}\\
{\tOmega^1_{\oX_n/\oS_n,\ox}}\ar[r]^-(0.5){\phi_n}&{\rH^1(\pi_1(\uoX^\circ,\oy),\xi(\oR^\oy_{\uX}/p^n\oR^\oy_{\uX}))}}
\end{equation}
où $\ttt_n$ est l'application définie dans \eqref{tfkum13a} en tenant compte de \eqref{tfkum12b}.
\end{prop}

En effet, l'unicité de $\phi_n$ est claire \eqref{tfkum2f}. Montrons son existence. 
Pour toute extension finie $L$ de $K$ contenue dans $\oK$,
on rappelle que l'on a $X_L=X\times_SS_L$ \eqref{tfkum2a}.
On note encore $\ox$ l'image du point géométrique $\ox$ par la projection canonique $\oX\rightarrow X_L$. 
D'après \eqref{tfkum3a} et \eqref{tfkum3b}, on a des isomorphismes canoniques 
\begin{eqnarray}
\cL_{\oX,\ox}&\stackrel{\sim}{\rightarrow}&\underset{\underset{L\subset \oK}{\longrightarrow}}{\lim}\  \cM_{X_L,\ox},\\
\tOmega^1_{\oX_n/\oS_n,\ox}&\stackrel{\sim}{\rightarrow}& \underset{\underset{L\subset \oK}{\longrightarrow}}{\lim}\  \tOmega^1_{X_n/S_n,\ox}\otimes_{\co_{X,\ox}}\co_{X_L,\ox}.
\end{eqnarray}
Pour toute extension finie $L/K$, la projection canonique $(X_L,\cM_{X_L})\rightarrow (S_L,\cM_{S_L})$ est adéquate (\cite{agt} III.4.8). 
Par ailleurs, le schéma $\uoX^\circ$ et l'algèbre $\oR^\oy_{\uX}$ ne changent pas en remplaçant $(X,\cM_X)$  par $(X_L,\cM_{X_L})$. 
On est donc réduit à montrer que la restriction de $-\log([\ ])\circ \ttt_n$ à $\cM_{X,\ox}$ se factorise à travers l'homomorphisme 
\begin{equation}
\cM_{X,\ox}\rightarrow \tOmega^1_{X_n/S_n,\ox}
\end{equation}
induit par le morphisme $d\log$ \eqref{tfkum14hi}.
On notera que celui-ci est compatible avec le morphisme $d\log$ \eqref{tfkum2e}.
En vertu de (\cite{agt} II.5.17), on peut supposer les hypothèses de \ref{tfkum14} remplies. 
Il suffit de montrer que  la restriction de l'homomorphisme $-\log([\ ])\circ \ttt_n$ à $\Gamma(X,\cM_X)$ 
se factorise à travers l'homomorphisme 
\begin{equation}
\Gamma(X,\cM_X)\rightarrow \Gamma(X,\tOmega^1_{X_n/S_n})
\end{equation}
induit par $d\log$ \eqref{tfkum14hi}.
Identifions les topos $\uoX^\circ_\fet$ et $\bB_{\pi_1(\uoX^\circ,\oy)}$ via \eqref{tfkum12b}.
Compte tenu de \eqref{tfkum12d}, 
on a un morphisme canonique $Q_n\rightarrow \varphi_\ox(\cQ_n)$ de $\bB_{\pi_1(\uoX^\circ,\oy)}$ 
dont le morphisme induit entre les groupes associés s'insère dans un diagramme commutatif
\begin{equation}
\xymatrix{
{\mu_{p^n}(\co_\oK)}\ar@{^(->}[r]\ar[d]&{Q_n^\gp}\ar@{->>}[r]\ar[d]&{\Gamma(X,\cM_X)^\gp}\ar[d]\\
{\mu_{p^n,\uoX^\circ_\fet}}\ar@{^(->}[r]&{\varphi_\ox(\cQ_n^\gp)}\ar@{->>}[r]&{\cL_{\oX,\ox}^\gp}}
\end{equation}
où les lignes sont les suites exactes \eqref{tfkum13c} et \eqref{tfkum14f}. 
Compte tenu de \eqref{tfkum12h} et \eqref{tfkum13f}, 
$\ttt_n$ est le composé de $\ttt'_n$ \eqref{tfkum13d} et de l'homomorphisme canonique $\cL_{\oX,\ox}\rightarrow \cL_{\oX,\ox}^\gp$.
On en déduit que le diagramme 
\begin{equation}
\xymatrix{
{\Gamma(X,\cM_X)}\ar[r]^-(0.5){B_n}\ar[d]&
{\rH^1(\pi_1(\oX^{\star\circ},\oy),\mu_{p^n}(\co_\oK))}\ar@{=}[d]\\
{\cL_{\oX,\ox}}\ar[r]^-(0.5){\ttt_n}&{\rH^1(\uoX^\circ_\fet,\mu_{p^n,\uoX^\circ_\fet})}}
\end{equation} 
est commutatif. 
L'assertion recherchée résulte alors de \ref{tfkum15}.

\begin{cor}\label{tfkum17}
Sous les hypothèses de \ref{tfkum14}, pour tout entier $n\geq 1$, le diagramme 
\begin{equation}\label{tfkum17a}
\xymatrix{
{\Gamma(\oX^\star,\tOmega^1_{\oX_n/\oS_n})}\ar[r]^-(0.5){A_n}\ar[d]&
{\rH^1(\pi_1(\oX^{\star\circ},\oy),\xi(\oR^\oy_X/p^n\oR^\oy_X))}\ar[d]\\
{\tOmega^1_{\oX_n/\oS_n,\ox}}\ar[r]^-(0.5){\phi_n}&{\rH^1(\pi_1(\uoX^\circ,\oy),\xi(\oR^\oy_{\uX}/p^n\oR^\oy_{\uX}))}}
\end{equation}
où $A_n$ est le morphisme \eqref{tfkum14e} et $\phi_n$ est le morphisme \eqref{tfkum16a}, est commutatif. 
\end{cor}

Cela résulte aussitôt de la preuve de \ref{tfkum16}.

\begin{cor}\label{tfkum18}
Sous les hypothèses de \ref{tfkum12}, pour tout entier $n\geq 1$, 
il existe un et un unique homomorphisme de $\co_{\oX,\ox}$-algèbres graduées 
\begin{equation}\label{tfkum18a}
\wedge (\xi^{-1}\tOmega^1_{\oX_n/\oS_n,\ox})\rightarrow \oplus_{i\geq 0}\rH^i(\pi_1(\uoX^\circ,\oy),\oR^\oy_{\uX}/p^n\oR^\oy_{\uX})
\end{equation}
dont la composante en degré un est le morphisme $\xi^{-1}\phi_n$ \eqref{tfkum16a}. De plus, 
notant $d=\dim(X/S)$ la dimension relative de $X$ sur $S$, 
le noyau de \eqref{tfkum18a} est annulé par $p^{\frac{2d}{p-1}}\fm_\oK$ et 
son conoyau est annulé par $p^{\frac{2d+1}{p-1}}\fm_\oK$.
\end{cor}

On désigne par $\fV_\ox$ la catégorie des $X$-schémas étales $\ox$-pointés, par
$\fW_\ox$ (resp. $\fW'_\ox$) la sous-catégorie pleine de $\fV_\ox$ formée des objets $(U,\fp\colon \ox\rightarrow U)$ 
tels que le schéma $U$ soit affine
(resp. vérifie les conditions (i) et (ii) de \ref{tfkum14}). 
Ce sont des catégories cofiltrantes, et les foncteurs d'injection canonique
$\fW'^\circ_\ox\rightarrow \fW^\circ_\ox\rightarrow \fV^\circ_\ox$ sont cofinaux d'après (\cite{agt} II.5.17) et (\cite{sga4} I 8.1.3(c)).
Reprenons les notations de \ref{TFA12}. 
Pour tout objet $(U,\fp)$ de $\fW'_\ox$, on note 
\begin{equation}\label{tfkum18b}
A_{(U,\fp),n}\colon \Gamma(\oU^\star,\tOmega^1_{\oX_n/\oS_n})
\rightarrow \rH^1(\pi_1(\oU^{\star\circ},\oy),\xi(\oR^\oy_U/p^n\oR^\oy_U))
\end{equation} 
le morphisme canonique \eqref{tfkum14e}; on rappelle que celui-ci ne dépend pas de 
la déformation choisie dans \ref{tfkum14} pour le définir (\cite{agt} II.10.10). D'après \ref{tfkum17}, le diagramme 
\begin{equation}\label{tfkum18c}
\xymatrix{
{\Gamma(\oU^\star,\tOmega^1_{\oX_n/\oS_n})}\ar[rr]^-(0.5){A_{(U,\fp),n}}\ar[d]&&
{\rH^1(\pi_1(\oU^{\star\circ},\oy),\xi(\oR^\oy_U/p^n\oR^\oy_U))}\ar[d]\\
{\tOmega^1_{\oX_n/\oS_n,\ox}}\ar[rr]^-(0.5){\phi_n}&&{\rH^1(\pi_1(\uoX^\circ,\oy),\xi(\oR^\oy_{\uX}/p^n\oR^\oy_{\uX}))}}
\end{equation}
est commutatif. D'après (\cite{agt} VI.11.10) et \eqref{TFA12f}, pour tout entier $i\geq 0$, on a un isomorphisme canonique
\begin{equation}\label{tfkum18d}
\rH^i(\pi_1(\uoX^\circ,\oy),\oR^\oy_{\uX}/p^n\oR^\oy_{\uX})
\stackrel{\sim}{\rightarrow} \underset{\underset{(U,\fp)\in \fW'^\circ_\ox}{\longrightarrow}}{\lim}\
\rH^i(\pi_1(\oU^{\star\circ},\oy),\oR^\oy_U/p^n\oR^\oy_U). 
\end{equation} 
Par ailleurs, $\uoX$ étant strictement local d'après (\cite{agt} III.3.7), il s'identifie au localisé strict de $\oX$ en $\ox$.
On a donc un isomorphisme canonique
\begin{equation}\label{tfkum18e}
\tOmega^1_{\oX_n/\oS_n,\ox}
\stackrel{\sim}{\rightarrow} \underset{\underset{(U,\fp)\in \fW'^\circ_\ox}{\longrightarrow}}{\lim}\
\Gamma(\oU^\star,\tOmega^1_{\oX_n/\oS_n}). 
\end{equation} 
Les homomorphismes $A_{(U,\fp),n}$, pour $(U,\fp)\in \ob(\fW'_\ox)$, forment un homomorphisme 
de systèmes inductifs d'après (\cite{agt} III.10.16). Leur limite inductive s'identifie donc à $\phi_n$ \eqref{tfkum18c}. 
Il suffit alors de montrer que pour tout objet $(U,\fp)$ de $\fW'_\ox$, 
il existe un et un unique homomorphisme de $\co_{\oX_n}(\oU^\star)$-algèbres graduées 
\begin{equation}\label{tfkum18f}
\wedge \Gamma(\oU^\star,\xi^{-1}\tOmega^1_{\oX_n/\oS_n})\rightarrow 
\oplus_{i\geq 0}\rH^i(\pi_1(\oU^{\star\circ},\oy),\oR^\oy_U/p^n\oR^\oy_U)
\end{equation}
dont la composante en degré un est le morphisme $\xi^{-1} A_{(U,\fp),n}$, que 
son noyau est annulé par $p^{\frac{2d}{p-1}}\fm_\oK$ et que son conoyau est annulé par $p^{\frac{2d+1}{p-1}}\fm_\oK$.

D'après (\cite{agt} II.10.16), on a un diagramme commutatif
\begin{equation}\label{tfkum18g}
\xymatrix{
{\Gamma(\oU^\star,\xi^{-1}\tOmega^1_{\oX_n/\oS_n})}\ar[d]_a\ar[r]^-(0.5){\xi^{-1}A_{(U,\fp),n}}&
{\rH^1(\pi_1(\oU^{\star \circ},\oy),\oR^{\oy}_U/p^n\oR^{\oy}_U)}\\
{\rH^1(\pi_1(\oU^{\star \circ},\oy),\xi^{-1}\oR^{\oy}_U(1)/p^n\xi^{-1}\oR^{\oy}_U(1))}\ar[r]^{-b}_\sim&
{\rH^1(\pi_1(\oU^{\star \circ},\oy),p^{\frac{1}{p-1}}\oR^{\oy}_U/p^{n+\frac{1}{p-1}}\oR^{\oy}_U)}\ar[u]_c}
\end{equation}
où $a$ est induit par le morphisme défini dans \eqref{tpcg20d}, 
$b$ est induit par l'isomorphisme canonique $\hoR^{\oy}_U(1)\stackrel{\sim}{\rightarrow} p^{\frac{1}{p-1}}\xi \hoR^{\oy}_U$ \eqref{definf17c}
et $c$ est induit par l'injection canonique $p^{\frac{1}{p-1}}\oR^{\oy}_U \rightarrow \oR^{\oy}_U$.
En vertu de \ref{tpcg6}, il existe un et un unique homomorphisme de $\co_{\oX_n}(\oU^\star)$-algèbres graduées
\begin{equation}\label{tfkum18h}
\wedge (\xi^{-1}\tOmega^1_{\oX_n/\oS_n}(\oU^\star))\rightarrow 
\oplus_{i\geq 0}\rH^i(\pi_1(\oU^{\star \circ},\oy),\xi^{-i}\oR^{\oy}_U(i)/p^n\xi^{-i}\oR^{\oy}_U(i))
\end{equation}
dont la composante en degré un est le morphisme $a$. Son noyau est $\alpha$-nul et son conoyau est annulé par 
$p^{\frac{1}{p-1}}\fm_\oK$. On en déduit qu'il existe un et un unique homomorphisme de 
$\co_{\oX_n}(\oU^\star)$-algèbres graduées
\begin{equation}\label{tfkum18i}
\wedge (\xi^{-1}\tOmega^1_{\oX_n/\oS_n}(\oU^\star))\rightarrow 
\oplus_{i\geq 0}\rH^i(\pi_1(\oU^{\star \circ},\oy),\oR^{\oy}_U/p^n\oR^{\oy}_U)
\end{equation}
dont la composante en degré $1$ est $\xi^{-1}A_{(U,\fp),n}$. 
Une chasse au diagramme \eqref{tfkum18g} 
montre que le noyau de \eqref{tfkum18i} est annulé par $p^{\frac{2d}{p-1}}\fm_\oK$.
Comme $\rH^i(\pi_1(\oU^{\star \circ},\oy),\oR^{\oy}_U/p^n\oR^{\oy}_U)$ 
est $\alpha$-nul pour tout $i\geq d+1$  en vertu de \ref{tpcg10}(iii), 
le conoyau de \eqref{tfkum18i} est annulé par $p^{\frac{2d+1}{p-1}}\fm_\oK$.

\subsection{}\label{tfkum19} 
Nous pouvons maintenant démontrer le théorème \ref{tfkum11}. La proposition étant locale pour la topologie étale de $X_s$, 
soient $\ox$ un point géométrique de $X_s$, $\uX$ le localisé strict de $X$ en $\ox$. 
D'après \ref{TFA15}(i), il existe un point  $(\oy\rightsquigarrow \ox)$ de $X_\et\gtimes_{X_\et}\oX^\circ_\et$ 
\eqref{topfl17}. 
En vertu de \eqref{TFA14d} et \eqref{TFA12g}, pour tout entier $i\geq 0$, on a un isomorphisme canonique 
\begin{equation}\label{tfkum19a}
\rR^i\sigma_{n*}(\xi\ocB_n)_\ox\stackrel{\sim}{\rightarrow}\rH^i(\uoX^\circ_\fet,\xi(\oR^\oy_{\uX}/p^n\oR^\oy_{\uX})). 
\end{equation}
La proposition résulte alors de \ref{tfkum13}, \ref{tfkum16}  et \ref{tfkum18}.

\section{La suite spectrale de Hodge-Tate absolue}\label{ssht}

\subsection{}\label{ssht1}
Avec les notations de \ref{hght2}, 
on désigne par $\tE_s^{\mN^\circ}$ (resp. $X_{s,\et}^{\mN^\circ}$) le topos des systèmes projectifs de $\tE_s$ (resp. $X_{s,\et}$), 
indexés par l'ensemble ordonné $\mN$ des entiers naturels \eqref{notconv13}, par $\bvocB$ l'anneau 
$(\ocB_{n+1})_{n\in \mN}$ de $\tE_s^{\mN^\circ}$  \eqref{TFA8a} et
par $\co_{\bvoX}$ l'anneau $(\co_{\oX_{n+1}})_{n\in \mN}$ de $X_{s,\et}^{\mN^\circ}$. 
On note 
\begin{equation}\label{ssht1a}
\bvsigma\colon (\tE_s^{\mN^\circ},\bvocB)\rightarrow(X_{s,\et}^{\mN^\circ},\co_{\bvoX})
\end{equation}
le morphisme de topos annelés induit par les morphismes $(\sigma_{n+1})_{n\in \mN}$ \eqref{TFA8e}. 
Considérons la suite spectrale de Cartan-Leray (\cite{sga4} V 5.3)
\begin{equation}\label{ssht1b}
\rE_2^{i,j}=\rH^i(X_s^{\mN^\circ},\rR^j \bvsigma_*\bvocB)\Rightarrow \rH^{i+j}(\tE_s^{\mN^\circ},\bvocB),
\end{equation}
qui est aussi la deuxième suite spectrale d'hypercohomologie du foncteur $\Gamma(X^{\mN^\circ}_{s,\et},-)$ 
par rapport au complexe $\rR\bvsigma_*(\bvocB)$ (\cite{ega3} 0.11.4.3 ou \cite{deligne2} 1.4.5 et 1.4.6). 
On rappelle que la deuxième suite spectrale d'hypercohomologie définit un foncteur 
de la catégorie dérivée des complexes de groupes abéliens de $X^{\mN^\circ}_{s,\et}$ 
dans celle des suites spectrales de groupes abéliens (\cite{ega3} 0.11.1.2).

On note $\xi^{-1}\tOmega^1_{\bvoX/\bvoS}$ le $\co_{\bvoX}$-module 
$(\xi^{-1}\tOmega^1_{\oX_{n+1}/\oS_{n+1}})_{n\in \mN}$ de $X_{s,\et}^{\mN^\circ}$ \eqref{tfkum2g} 
et pour tout entier $j\geq 1$, on pose $\xi^{-j}\tOmega^j_{\bvoX/\bvoS}=\wedge^j(\xi^{-1}\tOmega^1_{\bvoX/\bvoS})$. 
On a alors un isomorphisme canonique (\cite{agt} (III.7.12.4))
\begin{equation}\label{ssht1c}
\xi^{-j}\tOmega^j_{\bvoX/\bvoS}\stackrel{\sim}{\rightarrow}(\xi^{-j}\tOmega^j_{\oX_{n+1}/\oS_{n+1}})_{n\in \mN}.
\end{equation}
En vertu de \ref{tfkum11}, pour tout entier $j\geq 0$, on a un morphisme $\co_{\bvoX}$-linéaire canonique
\begin{equation}\label{ssht1d}
\xi^{-j}\tOmega^j_{\bvoX/\bvoS}\rightarrow \rR^j \bvsigma_*(\bvocB),
\end{equation}
dont le noyau et le conoyau sont annulés par  $p^{\frac{2d+1}{p-1}}\fm_\oK$ (\cite{agt} (III.7.5.5) et II.8.2).

\begin{prop}\label{ssht2}
Supposons le $S$-schéma $X$ propre. Alors, pour tout entier $q\geq 0$, le morphisme canonique 
\begin{equation}\label{ssht2a}
\rH^q(\tE_s^{\mN^\circ},\bvocB)\rightarrow \underset{\underset{n\geq 1}{\longleftarrow}}\lim\ \rH^q(\tE_s,\ocB_n)
\end{equation}
est un $\alpha$-isomorphisme. 
\end{prop}

En effet, d'après (\cite{agt} VI.7.10), on a une suite exacte
\begin{equation}\label{ssht2b}
0\rightarrow \rR^1\underset{\underset{n\geq 1}{\longleftarrow}}\lim\ \rH^{q-1}(\tE_s,\ocB_n)\rightarrow 
\rH^q(\tE_s^{\mN^\circ},\bvocB)\rightarrow \underset{\underset{n\geq 1}{\longleftarrow}}\lim\ \rH^q(\tE_s,\ocB_n)\rightarrow 0.
\end{equation}
En vertu de \ref{TPCF16}, pour tout $n\geq 1$, on a un morphisme canonique 
\begin{equation}\label{ssht2c}
u^q_n\colon \rH^q(\oX^\circ_{\et},\mZ/p^n\mZ)\otimes_{\mZ_p}\co_\oK\rightarrow \rH^q(\tE_s,\ocB_n),
\end{equation}
qui est un $\alpha$-isomorphisme. Les morphismes 
\begin{equation}\label{ssht2d}
\underset{\underset{n\geq 1}{\longleftarrow}}\lim\ u^q_n \ \ \ {\rm et}\ \ \ \rR^1\underset{\underset{n\geq 1}{\longleftarrow}}\lim\ u^q_n
\end{equation}
sont donc des $\alpha$-isomorphismes (\cite{gr} 2.4.2(ii)). 
Les groupes $\rH^q(\oX^\circ_{\et},\mZ/p^n\mZ)$ étant finis en vertu de (\cite{sga45} Th.finitude 1.1), 
le système projectif $(\rH^q(\oX^\circ_{\et},\mZ/p^n\mZ))_{n\geq 1}$ vérifie la condition de Mittag-Leffler. 
On en déduit que 
\begin{equation}
\rR^1\underset{\underset{n\geq 1}{\longleftarrow}}\lim\ \rH^q(\tE_s,\ocB_n)
\end{equation}
est $\alpha$-nul d'après (\cite{jannsen} 1.15) et (\cite{roos} théo.~1). La proposition s'ensuit. 

\begin{rema}\label{ssht300}
D'après (\cite{sga5} VI 2.2.2 et la remarque après 2.2.3), le système projectif $\rH^q=(\rH^q(\oX^\circ_{\et},\mZ/p^n\mZ))_{n\geq 0}$ est AR-$p$-adique constructible, 
autrement dit, il existe un système projectif $p$-adique noethérien de groupes abéliens $M^q=(M^q_n)_{n\geq 0}$ et un AR-isomorphisme 
$M^q\rightarrow \rH^q$ (\cite{sga5} V 3.2.2). On en déduit un isomorphisme 
\begin{equation}\label{sshtrl300a}
\underset{\underset{n\geq 0}{\longleftarrow}}\lim\ M^q_n\stackrel{\sim}{\rightarrow} \rH^q(\oX^\circ_{\et},\mZ_p)= \underset{\underset{n\geq 0}{\longleftarrow}}\lim\ 
\rH^q(\oX^\circ_{\et},\mZ/p^n\mZ).
\end{equation}
En particulier, le $\mZ_p$-module $\rH^q(\oX^\circ_{\et},\mZ_p)$ est de type fini. 
Le morphisme de systèmes projectifs $(M^q_n\otimes_{\mZ_p}\co_C)_{n\geq 0}\rightarrow (\rH^q(\oX^\circ_{\et},\mZ/p^n\mZ)\otimes_{\mZ_p}\co_C)_{n\geq 0}$ est aussi 
un AR-isomorphisme. Comme la source est un système $p$-adique, on en déduit que le morphisme canonique 
\begin{equation}\label{sshtrl300b}
\rH^q(\oX^\circ_{\et},\mZ_p)\otimes_{\mZ_p}\co_C \rightarrow \underset{\underset{n\geq 0}{\longleftarrow}}\lim\  
\rH^q(\oX^\circ_{\et},\mZ/p^n\mZ) \otimes_{\mZ_p}\co_C
\end{equation}
est un isomorphisme. 
\end{rema}

\begin{cor}\label{ssht3}
Supposons le $S$-schéma $X$ propre. Alors, pour tout entier $q\geq 0$, on a un $\alpha$-isomorphisme canonique 
\begin{equation}\label{ssht3a}
\rH^q(\tE_s^{\mN^\circ},\bvocB)\stackrel{\approx}{\longrightarrow} \rH^q(\oX^\circ_\et,\mZ_p) \otimes_{\mZ_p}\co_C.
\end{equation}
\end{cor}

Cela résulte de \ref{TPCF16}, \ref{ssht2} et \ref{ssht300}. 
On prendra garde que le $\alpha$-isomorphisme \eqref{ssht3a} n'est pas en général induit par un vrai morphisme.

\begin{prop}\label{ssht4}
Supposons le $S$-schéma $X$ propre et soit $\cF$ un $\co_X$-module cohérent.
Pour tout $n\geq 1$, posons $\ocF_n=\cF\otimes_{\co_X}\co_{\oX_n}$, que l'on considère 
comme un faisceau de $X_{s,\zar}$ ou $X_{s,\et}$, selon le contexte \eqref{notconv12}, et $\bvocF=(\ocF_{n+1})_{n\in \mN}$. 
Alors, pour tout entier $q\geq 0$, les morphismes canoniques
\begin{equation}\label{ssht4a}
\rH^q(X,\cF)\otimes_{\co_K}\co_C\rightarrow \rH^q(X_{s,\et}^{\mN^\circ},\bvocF)\rightarrow 
\underset{\underset{n\geq 1}{\longleftarrow}}\lim\ \rH^q(X_{s,\et},\ocF_n)
\end{equation}
sont des isomorphismes. 
\end{prop}

En effet, d'après (\cite{agt} VI.7.10), on a une suite exacte
\begin{equation}\label{ssht4b}
0\rightarrow \rR^1\underset{\underset{n\geq 1}{\longleftarrow}}\lim\ \rH^{q-1}(X_{s,\et},\ocF_n)\rightarrow 
\rH^q(X_{s,\et}^{\mN^\circ},\bvocF)\rightarrow \underset{\underset{n\geq 1}{\longleftarrow}}\lim\ \rH^q(X_{s,\et},\ocF_n)\rightarrow 0.
\end{equation}
Pour tout $n\geq 1$, le morphisme canonique 
\begin{equation}\label{ssht4c}
\rH^q(X_{s,\zar},\cF\otimes_{\co_X}\co_{X_n})\otimes_{\co_K}\co_\oK\rightarrow  \rH^q(X_{s,\zar},\ocF_n)
\end{equation}
est un isomorphisme. Il résulte alors de (\cite{ega3} 4.1.7) et (\cite{sga4} VII 4.3) que 
le système projectif $(\rH^q(X_{s,\et},\ocF_n))_{n\geq 1}$ vérifie la condition de Mittag-Leffler et que le morphisme canonique 
\begin{equation}\label{ssht4d}
\rH^q(X,\cF)\hotimes_{\co_K}\co_C\rightarrow \underset{\underset{n\geq 1}{\longleftarrow}}\lim\ \rH^q(X_{s,\et},\ocF_n)
\end{equation}
est un isomorphisme. La proposition s'ensuit compte tenu de \eqref{ssht4b}, (\cite{jannsen} 1.15) et du fait que le $\co_K$-module 
$\rH^q(X,\cF)$ est de type fini. 

\begin{teo}\label{ssht5}
Si le $S$-schéma $X$ est propre, on a une suite spectrale canonique 
\begin{equation}\label{ssht5a}
\rE_2^{i,j}=\rH^i(X,\tOmega^j_{X/S})\otimes_{\co_K}C(-j)\Rightarrow \rH^{i+j}(\oX^\circ_\et,\mQ_p)\otimes_{\mQ_p}C.
\end{equation}
\end{teo}

Cela résulte de la suite spectrale de Cartan-Leray \eqref{ssht1b}, compte tenu de \ref{ssht3}, \ref{ssht4} et \ref{tfkum11}.

\begin{defi}\label{ssht11}
Si le $S$-schéma $X$ est propre, la suite spectrale \eqref{ssht5a} est appelée {\em suite spectrale de Hodge-Tate} de $f$ (ou de $X$ sur $S$).
\end{defi}

\subsection{}\label{ssht12}
On munit le schéma $\oS$ de l'action naturelle à gauche du groupe de Galois absolu $G_K$ de $K$ \eqref{hght1}; 
pour tout $u\in G_K$, l'automorphisme de $\oS$ défini par $u$, que l'on note aussi $u$, 
est induit par l'automorphisme $u^{-1}$ de $\co_\oK$. 
Pour tout $S$-schéma $Z$, on en déduit par changement de base une action à gauche de $G_K$ sur $\oZ=Z\times_S\oS$
et par suite des actions à gauche sur les topos $\oZ_\zar$ et $\oZ_\et$ (cf. \ref{notconv17}). 
On notera que $\co_\oZ$ est un anneau $G_K$-équivariant de $\oZ_\zar$ ou $\oZ_\et$ \eqref{notconv12}.

Pour tout $u\in G_K$, le diagramme commutatif
\begin{equation}\label{ssht12a}
\xymatrix{
{\oX^\circ}\ar[r]^{u}\ar[d]_h&{\oX^\circ}\ar[d]^h\\
X\ar[r]^{\id}&X}
\end{equation}
induit par fonctorialité (\cite{agt} VI.10.12) un morphisme de topos que l'on note encore 
\begin{equation}\label{ssht12b}
\gamma_u\colon \tE\rightarrow \tE.
\end{equation}
Pour tout $(u,v)\in G_K^2$, on a un isomorphisme canonique $c_{u,v}\colon \gamma_{v}^*\gamma_u^*\stackrel{\sim}{\rightarrow}\gamma_{uv}^*$.
On voit aussitôt que ces isomorphismes vérifient les relations de compatibilité \ref{notconv17}(i)--(iii). 
On obtient ainsi une action à gauche de $G_K$ sur le topos $\tE$ dans le sens de  {\em loc. cit.} Suivant les conventions générales, 
pour tout $u\in G_K$, le morphisme $\gamma_u$ sera encore noté $u\colon \tE\rightarrow \tE$, ce qui n'induit aucun risque d'ambiguïté.

Soit $u\in G_K$. Pour tout $\oS$-schéma $V$, notons
\begin{equation}\label{ssht12c}
{^uV}=V\times_{\oS,u}\oS,
\end{equation}
le changement de base de $V$ par l'automorphisme $u$ de $\oS$.  Pour tout $S$-schéma $U$, on a un $\oS$-isomorphisme
canonique  $\oU\stackrel{\sim}{\rightarrow} {^u\oU}$ dont le composé avec la projection canonique ${^u\oU}\rightarrow \oU$ est l'automorphisme
$\oU$ induit par $u$. Pour tout objet $(V\rightarrow U)$ de $E$, le morphisme ${^uV}\rightarrow U$ composé de la projection 
canonique ${^uV}\rightarrow V$ et du morphisme $V\rightarrow U$ est clairement un objet de $E$ et on a 
\begin{equation}\label{ssht12d}
u^*((V\rightarrow U)^a)\stackrel{\sim}{\rightarrow}({^uV}\rightarrow U)^a.
\end{equation}
Par ailleurs, le diagramme commutatif à carrés cartésiens
\begin{equation}\label{ssht12e}
\xymatrix{
{^u V}\ar[r]\ar[d]&V\ar[d]\\
{\oU^\circ}\ar[r]^u\ar[d]&{\oU^\circ}\ar[d]\\
{\oU}\ar[r]^u&{\oU}}
\end{equation}
où la flèche horizontale supérieure est la projection canonique, induit un isomorphisme 
\begin{equation}\label{ssht12f}
\oU^{^uV}\stackrel{\sim}{\rightarrow}\oU^V
\end{equation} 
au-dessus de l'automorphisme $u$ de $\oU$ \eqref{TFA2}. 

Pour tout $u\in G_K$, on a un isomorphisme d'anneaux 
\begin{equation}\label{ssht12g}
\tau^\ocB_u\colon \ocB \stackrel{\sim}{\rightarrow}u^*(\ocB),
\end{equation}
dont l'adjoint de son inverse est défini pour tout $(V\rightarrow U)\in \ob(E)$ par l'isomorphisme 
\begin{equation}
\Gamma(\oU^V,\co_{\oU^V})\stackrel{\sim}{\rightarrow} \Gamma(\oU^{^uV},\co_{\oU^{^uV}})
\end{equation}
induit par \eqref{ssht12f}.
Ces isomorphismes vérifient clairement la relation de compatibilité \eqref{notconv17d} et font donc 
de $\ocB$ un anneau $G_K$-équivariant de $\tE$ (cf. \ref{notconv17}). 

Munissant le topos $X_\et$ de l'action triviale de $G_K$, le morphisme canonique $\sigma\colon \tE\rightarrow X_\et$ \eqref{TFA6c} 
est $G_K$-équivariant dans le sens de \ref{notconv18}. 
Il revient au même de dire que pour tout $u\in G_K$, le diagramme 
\begin{equation}\label{ssht12h}
\xymatrix{
\tE\ar[r]^-(0.5)u\ar[dr]_{\sigma}&\tE\ar[d]^{\sigma}\\
&{X_\et}}
\end{equation}
est commutatif à isomorphisme canonique près, et que ces isomorphismes vérifient des relations de compatibilité. 
Le faisceau $\sigma^*(X_\eta)$ de $\tE$ est donc naturellement muni d'une structure de faisceau $G_K$-équivariant. 
Par suite, l'action à gauche de $G_K$ sur $\tE$ induit une action à gauche  sur le sous-topos fermé $\tE_s$
de sorte que le plongement $\delta\colon \tE_s\rightarrow E$ est $G_K$-équivariant. 
Concrètement, pour tout $u\in G_K$, le diagramme \eqref{ssht12h} induit un isomorphisme 
\begin{equation}
\tau_u^{\sigma^*(X_\eta)}\colon \sigma^*(X_\eta)\stackrel{\sim}{\rightarrow}u^*(\sigma^*(X_\eta)). 
\end{equation}
En vertu de (\cite{sga4} IV 9.4.3), il existe donc un morphisme de topos que l'on note encore
\begin{equation}\label{ssht12i}
u\colon \tE_s\rightarrow \tE_s,
\end{equation}
unique à isomorphisme canonique près tel que le diagramme 
\begin{equation}\label{ssht12j}
\xymatrix{
{\tE_s}\ar[r]^{u}\ar[d]_{\delta}&{\tE_s}\ar[d]^{\delta}\\
{\tE}\ar[r]^u&{\tE}}
\end{equation}
soit commutatif à isomorphisme canonique près. 

Compte tenu de (\cite{sga4} IV 9.4.3), munissant le topos $X_{s,\et}$ de l'action triviale de $G_K$, le morphisme 
$\sigma_s\colon \tE_s\rightarrow X_{s,\et}$ \eqref{TFA66b} est $G_K$-équivariant.

Pour tout entier $n\geq 1$, la structure d'anneau $G_K$-équivariant sur $\ocB$ induit sur $\ocB_n$ 
une structure d'anneau $G_K$-équivariant de $\tE_s$ (ou de $\tE$). 
En particulier, pour tout $u\in G_K$, le morphisme \eqref{ssht12i} est sous-jacent à un morphisme de topos annelés que l'on note encore 
\begin{equation}\label{ssht12k}
u\colon (\tE_s,\ocB_n)\rightarrow (\tE_s,\ocB_n).
\end{equation}
Par ailleurs, $\co_{\oX_n}$ est naturellement un anneau $G_K$-équivariant de $X_{s,\et}$, et l'homomorphisme 
$\co_{\oX_n}\rightarrow \sigma_{s*}(\ocB_n)$ est un morphisme d'anneaux $G_K$-équivariants  (cf. \ref{notconv18}).
Le morphisme de topos annelés \eqref{TFA8e}
\begin{equation}\label{ssht12l}
\sigma_n\colon (\tE_s,\ocB_n)\rightarrow (X_{s,\et},\co_{\oX_n})
\end{equation} 
est donc $G_K$-équivariant.

\subsection{}\label{ssht120}
L'action de $G_K$ sur $\tE_s$ \eqref{ssht12} induit une action à gauche sur $\tE^{\mN^\circ}_s$. 
Concrètement, pour tout $u\in G_K$, le morphisme \eqref{ssht12i} induit un morphisme de topos que l'on note encore 
\begin{equation}\label{ssht120a}
u\colon \tE^{\mN^\circ}_s\rightarrow \tE^{\mN^\circ}_s. 
\end{equation}
La structure d'anneau $G_K$-équivariant sur $\ocB$ induit sur $\bvocB$ 
une structure d'anneau $G_K$-équivariant de $\tE^{\mN^\circ}_s$. 
Munissant le topos $X^{\mN^\circ}_{s,\et}$ de l'action triviale de $G_K$, le morphisme 
$\sigma_s^{\mN^\circ}\colon \tE^{\mN^\circ}_s\rightarrow X^{\mN^\circ}_{s,\et}$ est $G_K$-équivariant et 
l'homomorphisme canonique 
\begin{equation}\label{ssht120b}
\co_{\bvoX}\rightarrow  (\sigma_s^{\mN^\circ})_*(\bvocB)
\end{equation}
est un morphisme d'anneaux $G_K$-équivariants (cf. \ref{notconv18}). Le morphismes de topos annelés \eqref{ssht1a}
\begin{equation}\label{ssht120c}
\bvsigma \colon (\tE^{\mN^\circ}_s,\bvocB)\rightarrow (X^{\mN^\circ}_{s,\et},\co_{\bvoX})
\end{equation}
est donc $G_K$-équivariant.

Soit $\cF$ un $\bvocB$-module $G_K$-équivariant de $\tE^{\mN^\circ}_s$ \eqref{notconv17}. 
D'après \ref{notconv19}, pour tout entier $q\geq 0$,
le $\co_C$-module $\rH^q(\tE^{\mN^\circ}_s,\bvocB)$ est naturellement muni d'une action $\co_C$-semi-linéaire à gauche de $G_K$. 
D'après \ref{notconv18}, pour tout entier $j\geq 0$,
$\rR^j\bvsigma_*(\cF)$ est naturellement muni d'une structure de $\co_{\bvoX}$-module $G_K$-équivariant. 
Comme l'action de $G_K$ sur $X^{\mN^\circ}_{s,\et}$ est triviale, cette structure se résume en une action $\co_{\bvoX}$-semi-linéaire 
à gauche de $G_K$ sur $\rR^j\bvsigma_*(\cF)$. On en déduit, pour tout $i\geq 0$,  une action $\co_C$-semi-linéaire à gauche de $G_K$ sur 
le $\co_C$-module $\rH^i(X^{\mN^\circ}_{s,\et},\rR^j\bvsigma_*(\cF))$.

On rappelle que le morphisme $\theta\colon \rW(\co_{\oK^\flat})\rightarrow \co_C$ de Fontaine \eqref{eipo3d} est $G_K$-équivariant.  
Par suite, $G_K$ agit naturellement sur $\cA_2(\co_\oK)$ \eqref{eipo3e} et cette action stabilise l'idéal $\ker(\theta_2)=\xi\co_C$.  
Par ailleurs, pour tout $u\in G_K$, les actions de $u$ sur $\oX$ et $\oS$ induisent un isomorphisme $\co_\oX$-linéaire
\begin{equation}
\tOmega^1_{\oX/\oS}\stackrel{\sim}{\rightarrow} u^*(\tOmega^1_{\oX/\oS}).
\end{equation}
Ces isomorphismes munissent $\tOmega^1_{\oX/\oS}$ d'une structure de $\co_\oX$-module $G_K$-équivariant de $\oX_\et$ 
\eqref{notconv17}. 
Comme l'action de $G_K$ sur $X_{s}$ est triviale, on en déduit, pour tout entier $j\geq 0$, 
une action à gauche $\co_{\bvoX}$-semi-linéaire de $G_K$ sur le $\co_{\bvoX}$-module $\xi^{-j}\tOmega^j_{\bvoX/\bvoS}$.

\begin{prop}\label{ssht13}
La suite spectrale de Cartan-Leray \eqref{ssht1b} est $G_K$-équivariante \eqref{ssht120}. 
\end{prop}

Pour tout $u\in G_K$, considérons l'isomorphisme composé 
\begin{equation}\label{ssht13a}
\rR\bvsigma_*(\bvocB)\stackrel{\sim}{\rightarrow}  \rR\bvsigma_*(u_*(\bvocB))\stackrel{\sim}{\rightarrow} \rR\bvsigma_*(\bvocB),
\end{equation}
où la première flèche est induite par l'isomorphisme canonique $\bvocB\stackrel{\sim}{\rightarrow} u_*(\bvocB)$,
et la seconde flèche est induite par la $G_K$-équivariance de $\bvsigma$. 
On vérifie aussitôt que ces morphismes induisent les structures $G_K$-équivariantes sur 
$\rH^q(\tE^{\mN^\circ}_s, \bvocB)$ $(q\geq 0)$ et $\rR^j\bvsigma_*(\bvocB)$ $(j\geq 0)$, définies dans \ref{ssht120}.
Par fonctorialité de la deuxième suite spectrale d'hypercohomologie du foncteur $\Gamma(X^{\mN^\circ}_{s,\et},-)$
sur la catégorie dérivée des groupes abéliens de $X^{\mN^\circ}_{s,\et}$, les morphismes \eqref{ssht13a} induisent   
une action de $G_K$ sur la suite spectrale de Cartan-Leray \eqref{ssht1b}.

\begin{prop}\label{ssht14}
Pour tout entier $j\geq 0$, le morphisme \eqref{ssht1d}
\begin{equation}\label{ssht14a}
\xi^{-j}\tOmega^j_{\bvoX/\bvoS}\rightarrow \rR^j \bvsigma_*(\bvocB)
\end{equation}
est $G_K$-équivariant pour les structures $G_K$-équivariantes naturelles \eqref{ssht120}. 
\end{prop}

Par ailleurs, le groupe $G_K$ agit naturellement sur le schéma logarithmique $(\oS,\cL_\oS)$ \eqref{hght3} et par suite sur 
le schéma logarithmique  $(\oX,\cL_\oX)$ \eqref{tfkum2b}. Pour tout entier $n\geq 0$, le monoïde $\cQ_n$ \eqref{tfkum6a} est donc
muni d'une structure de monoïde $G_K$-équivariant de $\tE$, de sorte que les morphismes du diagramme \eqref{tfkum6a} sont 
des morphismes de monoïdes $G_K$-équivariants. 
Munissant le groupe $\mu_{p^n}(\co_\oK)$ de l'action canonique de $G_K$, 
la  suite exacte \eqref{tfkum8a}  est donc une suite de groupe abéliens $G_K$-équivariants de $\tE_s$. 
Par suite, $\partial_n$ \eqref{tfkum10a} est un morphisme de groupes abéliens $G_K$-équivariants de $X_{s,\et}$, 
et il en est alors de même du morphisme \eqref{tfkum10c}.
Le caractère universel de la dérivation logarithmique $d\log$ \eqref{tfkum2e}
implique que le morphisme \eqref{tfkum2f} 
\begin{equation}
\cL_\oX^\gp\otimes_\mZ\co_\oX\rightarrow \tOmega^1_{\oX/\oS}
\end{equation}
est un morphisme de $\co_\oX$-modules $G_K$-équivariants de $\oX_\et$. 
Le morphisme \eqref{tfkum11a} est donc $G_K$-équivariant, d'où la proposition. 

\begin{prop}\label{ssht17}
Supposons le $S$-schéma $X$ propre. Alors, pour tout entier $q\geq 0$, l'isomorphisme 
\begin{equation}\label{ssht17a}
\rH^q(\tE_s^{\mN^\circ},\bvocB)\otimes_\mZ\mQ\stackrel{\sim}{\longrightarrow} \rH^q(\oX^\circ_\et,\mZ_p) \otimes_{\mZ_p}C,
\end{equation}
induit par \eqref{ssht3a}, est $G_K$-équivariant. 
\end{prop}

Compte tenu de \ref{ssht2} et de sa preuve, il suffit de montrer que pour tout entier $n\geq 1$,
le morphisme canonique $u_n^q$ \eqref{ssht2c} est $G_K$-équivariant. 
Le morphisme $\psi\colon \oX^\circ_\et\rightarrow \tE$ \eqref{TFA6d} est $G_K$-équivariant  \eqref{notconv18}. 
Le faisceau constant $\mZ/p^n\mZ$ de $\oX^\circ_\et$ est canoniquement muni d'une structure d'algèbre $G_K$-équivariante. 
Par suite, $\psi_*(\mZ/p^n\mZ)$ est naturellement muni d'une structure d'algèbre $G_K$-équivariante de $\tE$,
et l'homomorphisme canonique $\mZ/p^n\mZ\rightarrow \psi_*(\mZ/p^n\mZ)$ est un morphisme d'algèbres $G_K$-équivariantes. 
Cet homomorphisme est en fait un isomorphisme d'après (\cite{agt} VI.10.9(iii)). 
En vertu de \ref{acycloc2}, pour tout entier $q\geq 0$, le morphisme canonique
\begin{equation}
\rH^q(\tE, \psi_*(\mZ/p^n\mZ))\rightarrow \rH^q(\oX^\circ_\et,\mZ/p^n\mZ)
\end{equation}
est un isomorphisme. Il est de plus $G_K$-équivariant lorsque l'on munit les deux groupes des actions naturelles de $G_K$ 
(cf. \ref{notconv19}). Par ailleurs, les morphismes canoniques
\begin{equation}
\rH^q(\tE,\mZ/p^n\mZ)\otimes_{\mZ_p}\co_C\rightarrow \rH^q(\tE,\ocB_n)\stackrel{\sim}{\rightarrow} \rH^q(\tE_s,\ocB_n)
\end{equation}
sont clairement $G_K$-équivariants, d'où la proposition.

\begin{cor}\label{ssht15}
Si le $S$-schéma $X$ est propre, la suite spectrale de Hodge-Tate \eqref{ssht5a} est $G_K$-équivariante. 
\end{cor}

Cela résulte de \ref{ssht13}, \ref{ssht14} et \ref{ssht17}.

\begin{cor}\label{ssht16}
Si le $S$-schéma $X$ est propre, la suite spectrale de Hodge-Tate \eqref{ssht5a} dégénère en $\rE_2$ et est décomposée: 
elle induit, pour tout entier $n\geq 0$, une décomposition canonique $G_K$-équivariante
\begin{equation}
\rH^n(\oX^\circ_\et,\mQ_p)\otimes_{\mQ_p}C\stackrel{\sim}{\rightarrow}\bigoplus_{0\leq i\leq n}\rH^i(X,\tOmega^{n-i}_{X/S})\otimes_{\co_K}C(i-n).
\end{equation}
\end{cor}

Cela résulte de \ref{ssht15} et du fait que $\rH^i(G_K,C(j))$ est nul pour $i=0,1$ et tout entier non-nul $j$ (\cite{tate} prop.~8).

\subsection{}\label{ssht20}
Soient $u\colon \cT'\rightarrow \cT$, $v\colon \cT\rightarrow \cS$ deux morphismes de topos, $v'=v\circ u$,
$\rK^\bullet$ un complexe de faisceaux abéliens de $\cT'$ borné inférieurement. 
Considérons la seconde suite spectrale d'hypercohomologie du foncteur $v'_*$ 
par rapport au complexe $\rK^\bullet$ (\cite{ega3} 0.11.4.3), 
\begin{equation}\label{ssht20a}
{^\prime\rE}_2^{i,j}=\rR^iv'_*(\rH^j(\rK^\bullet))\Rightarrow \rR^{i+j}v'_*(\rK^\bullet),
\end{equation}
et la seconde suite spectrale d'hypercohomologie du foncteur $v_*$ par rapport au complexe $\rR u_*(\rK^\bullet)$, 
\begin{equation}\label{ssht20b}
\rE_2^{i,j}=\rR^iv_*(\rR^ju_*(\rK^\bullet))\Rightarrow \rR^{i+j}v_*(\rR u_*(\rK^\bullet)).
\end{equation}
Pour tout entier $j\geq 0$, considérons l'edge-homomorphisme de la 
seconde suite spectrale d'hypercohomologie du foncteur $u_*$ par rapport au complexe $\rK^\bullet$
\begin{equation}\label{ssht20c}
\rR^ju_*(\rK^\bullet)\rightarrow u_*(\rH^j(\rK^\bullet)).
\end{equation}
Pour tout entier $i\geq 0$, considérons le morphisme 
\begin{equation}\label{ssht20d}
\rR^iv_*(\rR^ju_*(\rK^\bullet))\rightarrow \rR^iv'_*(\rH^j(\rK^\bullet)),
\end{equation}
composé des morphismes
\begin{equation}\label{ssht20e}
\rR^iv_*(\rR^ju_*(\rK^\bullet))\rightarrow \rR^iv_*(u_*(\rH^j(\rK^\bullet))) \rightarrow \rR^iv'_*(\rH^j(\rK^\bullet)),
\end{equation}
où la première flèche est induite par \eqref{ssht20c} et la seconde flèche est l'edge-homomorphisme de la 
suite spectrale de Cartan-Leray (\cite{sga4} V 5.4). 
Pour tout entier $q\geq 0$, considérons l'isomorphisme induit par la même suite spectrale 
\begin{equation}\label{ssht20f}
\rR^qv_*(\rR u_*(\rK^\bullet))\stackrel{\sim}{\rightarrow} \rR^qv'_*(\rK^\bullet). 
\end{equation}

\begin{prop}\label{ssht21}
Sous les hypothèses de \ref{ssht20}, les morphismes \eqref{ssht20d} pour les entiers $i,j\geq 0$ et les morphismes 
\eqref{ssht20f} pour les entiers $q\geq 0$, définissent un morphisme de suites spectrales $\rE\rightarrow {^\prime \rE}$ \eqref{ssht20a},
\eqref{ssht20b} {\rm (cf. \cite{ega3} 0.11.1.2)}.
\end{prop}

On désigne par $\bD(\cT',\mZ)$ (resp. $\bD(\cT,\mZ)$, resp. $\bD(\cS,\mZ)$) la catégorie dérivée des faisceaux abéliens de $\cT'$ (resp. $\cT$, resp. $\cS$). 
Nous utiliserons la notion d'objet spectral d'une catégorie triangulée introduite dans (\cite{verdier2} II 4.1). 
Soit $\cK$ le deuxième objet spectral canonique de $\bD(\cT',\mZ)$ associé au complexe $K^\bullet$ (\cite{verdier2} III 4.3). 
Pour tous éléments $i\leq j\leq \ell$ de $\tmZ=\mZ\cup\{-\infty,+\infty\}$, on a 
\begin{equation}\label{ssht21a}
\cK(i,j)=\tau_{[i-1,j-1]}K^\bullet
\end{equation}
et
\begin{equation}\label{ssht21b}
\delta(i,j,\ell)\colon \cK(j,\ell)\rightarrow \cK(i,j)[+1]
\end{equation}
est le morphisme canonique (cf. \cite{verdier2} II 4.1.2 et II 4.3 pour les notations). 
L'image $\rR u_*(\cK)$ de $\cK$ par le foncteur exact $\rR u_*$ est un objet spectral de  $\bD(\cT,\mZ)$. 
Soit $\cL$ le deuxième objet spectral canonique de $\bD(\cT,\mZ)$ associé au complexe $\rR u_*(K^\bullet)$ (\cite{verdier2} III 4.3.3). 
Pour tout entier $j$, la suite spectrale de Cartan-Leray montre que le morphisme canonique 
\begin{equation}\label{ssht21c}
\tau_{\leq j}(\rR u_*(\tau_{\leq j} K^\bullet))\rightarrow \tau_{\leq j}(\rR u_*(K^\bullet))
\end{equation}
est un isomorphisme de $\bD(\cT,\mZ)$. On en déduit un système compatible de morphismes 
\begin{equation}\label{ssht21d}
\tau_{\leq j}(\rR u_*(K^\bullet))\rightarrow \rR u_*(\tau_{\leq j}(K^\bullet)) \ \ \ (j\in \tmZ), 
\end{equation}
et des morphismes d'objets spectraux $\cL\rightarrow \rR u_*(\cK)$ et 
\begin{equation}\label{ssht21e}
\rR v_*(\cL)\rightarrow \rR v'_*(\cK).
\end{equation}
Pour tout entier $j$, on a 
\begin{eqnarray}
\cK(j,j+1)&=&\rH^j(K^\bullet)[-j],\\
\cL(j,j+1)&=&\rR^ju_*(K^\bullet)[-j].
\end{eqnarray}
Le morphisme $\cL(j,j+1)\rightarrow \rR u_*(\cK(j,j+1))$ s'identifie au morphisme composé
\begin{equation}
\rR^ju_*(K^\bullet)[-j]\rightarrow u_*(\rH^j(K^\bullet)[-j]) \rightarrow \rR u_*(\rH^j(K^\bullet)[-j]),
\end{equation}
où la première flèche est induite par \eqref{ssht20c} et la second flèche est le morphisme canonique. 
Par suite, la suite spectrale associée à $\rR v_*(\cL)$ (resp. $\rR v'_*(\cK)$) n'est autre que $\rE$ (resp. ${^\prime \rE}$) (\cite{verdier2} III 4.4.6).
Le morphisme $\rE\rightarrow {^\prime \rE}$ déduit du morphisme \eqref{ssht21e} induit sur les termes initiaux les morphismes \eqref{ssht20d}  
et sur les aboutissements les morphismes \eqref{ssht20f}.

\subsection{}\label{ssht6}
Soit $f'\colon (X',\cM_{X'})\rightarrow (S,\cM_S)$ un morphisme adéquat de schémas logarithmiques (\cite{agt} III.4.7). 
On désigne par $X'^\rhd$ le sous-schéma ouvert maximal de $X'$
où la structure logarithmique $\cM_{X'}$ est triviale. Pour tout $X'$-schéma $U'$, on pose  
\begin{equation}
U'^\rhd=U'\times_{X'}X'^\rhd.
\end{equation} 
On associe à $f'$ des objets analogues à ceux associés à $f$ dans que l'on note avec les mêmes symboles équipés d'un prime $^\prime$. 
On dispose en particulier d'un morphisme de topos annelés 
\begin{equation}\label{ssht6a}
\bvsigma'\colon (\tE'^{\mN^\circ}_s,\bvocB')\rightarrow(X'^{\mN^\circ}_{s,\et},\co_{\bvoX'}),
\end{equation}
et par suite d'une suite spectrale de Cartan-Leray 
\begin{equation}\label{ssht6b}
{^\prime\rE}^{i,j}_2=\rH^i(X'^{\mN^\circ}_s,\rR^j \bvsigma'_*\bvocB')\Rightarrow \rH^{i+j}(\tE'^{\mN^\circ}_s,\bvocB').
\end{equation}

Soit $g\colon (X',\cM_{X'})\rightarrow (X,\cM_X)$ un $(S,\cM_S)$-morphisme. 
Nous reprenons dans la suite les constructions introduites dans \ref{mtfla11} et \ref{mtfla12} que nous rappelons pour la commodité du lecteur. 
On désigne par $\upgamma\colon \oX'^\rhd\rightarrow \oX^\circ$ le morphisme de schémas et par 
\begin{equation}\label{ssht6c}
\Theta\colon \tE'\rightarrow \tE
\end{equation}
le morphisme de topos induits par $g$ (\cite{agt} VI.10.12). Le diagramme 
\begin{equation}\label{ssht6d}
\xymatrix{
{\tE'}\ar[r]^-(0.5){\sigma'}\ar[d]_{\Theta}&{X'_\et}\ar[d]^{g_\et}\\
{\tE}\ar[r]^{\sigma}&{X_\et}&}
\end{equation}
est commutatif à isomorphisme canonique près. 
On en déduit un isomorphisme $\Theta^*(\sigma^*(X_\eta))\simeq\sigma'^*(X'_\eta)$.
En vertu de (\cite{sga4} IV 9.4.3), il existe donc un morphisme de topos
\begin{equation}\label{ssht6e}
\Theta_s\colon \tE'_s\rightarrow \tE_s
\end{equation}
unique à isomorphisme canonique près tel que le diagramme 
\begin{equation}\label{ssht6f}
\xymatrix{
{\tE'_s}\ar[r]^{\Theta_s}\ar[d]_{\delta'}&{\tE_s}\ar[d]^{\delta}\\
{\tE'}\ar[r]^\Theta&{\tE}}
\end{equation}
soit commutatif à isomorphisme près. Il résulte de \eqref{ssht6d} et (\cite{sga4} IV 9.4.3) 
que le diagramme de morphismes de topos 
\begin{equation}\label{ssht6g}
\xymatrix{
{\tE'_s}\ar[d]_{\Theta_s}\ar[r]^-(0.5){\sigma'_s}&{X'_{s,\et}}\ar[d]^{g_s}\\
{\tE_s}\ar[r]^-(0.5){\sigma_s}&{X_{s,\et}}}
\end{equation}
est commutatif à isomorphisme canonique près.

On a un homomorphisme canonique $\Theta^{-1}(\ocB)\rightarrow \ocB'$ (cf. \ref{mtfla12}).
Pour tout entier $n\geq 1$, le morphisme \eqref{ssht6e} est donc sous-jacent à un morphisme de topos annelés 
\begin{equation}\label{ssht6gh}
\Theta_n\colon (\tE'_s,\ocB'_n)\rightarrow (\tE_s,\ocB_n).
\end{equation}
Le diagramme de morphismes de topos annelés
\begin{equation}\label{ssht6hh}
\xymatrix{
{(\tE'_s,\ocB'_n)}\ar[r]^{\Theta_n}\ar[d]_{\sigma'_n}&{(\tE_s,\ocB_n)}\ar[d]^{\sigma_n}\\
{(X'_{s,\et},\co_{\oX'_n})}\ar[r]^{\ogg_n}&{(X_{s,\et},\co_{\oX_n})}}
\end{equation}
est commutatif à isomorphisme canonique près.

Les morphismes \eqref{ssht6gh} induisent un morphisme de topos annelés 
\begin{equation}\label{ssht6h}
\bvTheta\colon (\tE'^{\mN^\circ}_s,\bvocB')\rightarrow (\tE^{\mN^\circ}_s,\bvocB).
\end{equation}
Le diagramme de morphismes de topos annelés 
\begin{equation}\label{ssht6i}
\xymatrix{
{(\tE'^{\mN^\circ}_s,\bvocB')}\ar[r]^{\bvTheta}\ar[d]_{\bvsigma'}&{(\tE^{\mN^\circ}_s,\bvocB)}\ar[d]^{\bvsigma}\\
{(X'^{\mN^\circ}_{s,\et},\co_{\bvoX'})}\ar[r]^{\bvogg}&{(X^{\mN^\circ}_{s,\et},\co_{\bvoX})}}
\end{equation}
où $\bvogg$ est induit par le morphisme de schémas $\ogg\colon \oX'\rightarrow \oX$, 
est commutatif à isomorphisme canonique près. 

Pour tout entier $j\geq 0$, considérons les morphismes 
\begin{equation}\label{ssht6j}
\rR^j\bvsigma_*(\bvTheta_*(\bvocB'))\rightarrow \rR^j(\bvsigma\circ \bvTheta)_*(\bvocB') \stackrel{\sim}{\rightarrow} 
\rR^j(\bvogg\circ \bvsigma')_*(\bvocB') \rightarrow \bvogg_*(\rR^j\bvsigma'_*(\bvocB')),
\end{equation}
le premier et le dernier étant les edge-homomorphismes de la suite spectrale de Cartan-Leray (\cite{sga4} V 5.4) et l'isomorphisme
central étant induit par \eqref{ssht6i}. 
Composant avec l'image directe supérieure $\rR^j\bvsigma_*$ de l'homomorphisme canonique $\bvocB\rightarrow \bvTheta_*(\bvocB')$, 
on obtient un morphisme
\begin{equation}\label{ssht6k}
\rR^j\bvsigma_*(\bvocB)\rightarrow \bvogg_*(\rR^j\bvsigma'_*(\bvocB')). 
\end{equation}
Pour tous entiers $i\geq 0$ et $j\geq 0$, considérons le morphisme composé 
\begin{equation}\label{ssht6m}
\rH^i(X^{\mN^\circ}_{s,\et},\rR^j\bvsigma_*(\bvocB))\rightarrow 
\rH^i(X^{\mN^\circ}_{s,\et},\bvogg_*(\rR^j\bvsigma'_*(\bvocB')))\rightarrow \rH^i(X'^{\mN^\circ}_{s,\et},\rR^j\bvsigma'_*(\bvocB')),
\end{equation}
où la première flèche est induite par \eqref{ssht6k} et la seconde flèche est l'edge-homomorphisme de 
la suite spectrale de Cartan-Leray (\cite{sga4} V 5.3). 

Pour tout entier $q\geq 0$, considérons le morphisme composé 
\begin{equation}\label{ssht6l}
\rH^q(\tE^{\mN^\circ}_s,\bvocB)\rightarrow \rH^q(\tE^{\mN^\circ}_s,\bvTheta_*(\bvocB'))\rightarrow \rH^q(\tE'^{\mN^\circ}_s,\bvocB'),
\end{equation}
où la première flèche est induite par l'homomorphisme canonique $\bvocB\rightarrow \bvTheta_*(\bvocB')$
et la seconde flèche est l'edge-homomorphisme de la suite spectrale de Cartan-Leray. 

\begin{prop}\label{ssht7}
Sous les hypothèses de \ref{ssht6}, les morphismes \eqref{ssht6m} 
\begin{equation}\label{ssht7a}
\rH^i(X^{\mN^\circ}_{s,\et},\rR^j\bvsigma_*(\bvocB))\rightarrow \rH^i(X'^{\mN^\circ}_{s,\et},\rR^j\bvsigma'_*(\bvocB')),
\end{equation}
pour les entiers $i,j\geq 0$, et les morphismes \eqref{ssht6l}
\begin{equation}\label{ssht7b}
\rH^q(\tE^{\mN^\circ}_s,\bvocB)\rightarrow \rH^q(\tE'^{\mN^\circ}_s,\bvocB'),
\end{equation}
pour les entiers $q\geq 0$, définissent un morphisme de suites spectrales $\rE\rightarrow {^\prime \rE}$ \eqref{ssht1b}, \eqref{ssht6b}
{\rm (cf. \cite{ega3} 0.11.1.2)}. 
\end{prop}

Notons 
\begin{equation}
{^\prime{^\prime \rE}}_2^{i,j}=\rH^i(X^{\mN^\circ}_{s,\et},\rR^j(\bvogg\circ \bvsigma')_*(\bvocB'))\Rightarrow \rH^{i+j}(\tE'^{\mN^\circ}_s,\bvocB').
\end{equation}
la seconde suite spectrale d'hypercohomologie du foncteur $\Gamma(X^{\mN^\circ}_{s,\et},-)$ par rapport au complexe $\rR(\bvogg\circ \bvsigma')_*(\bvocB')$.
D'après \ref{ssht21}, appliqué au composé du morphisme $\bvogg$ et du morphisme canonique $X^{\mN^\circ}_{s,\et}\rightarrow \Ens$ et au
complexe $\rR\bvsigma'_*(\bvocB')$, on a un morphisme de suites spectrales 
\begin{equation}
{^\prime {^\prime \rE}}\rightarrow {^\prime \rE}.
\end{equation}
Par ailleurs, l'homomorphisme canonique $\bvocB\rightarrow \bvTheta_*(\bvocB')$ induit un morphisme $\bvocB\rightarrow \rR\bvTheta_*(\bvocB')$
de la catégorie dérivée des groupes abéliens de $\tE^{\mN^\circ}_s$. Tenant compte de \eqref{ssht6i}, on en déduit un morphisme
\begin{equation}\label{ssht7c}
\rR\bvsigma_*(\bvocB)\rightarrow \rR(\bvogg\circ \bvsigma')_*(\bvocB').
\end{equation}
Par fonctorialité de la seconde suite spectrale d'hypercohomologie du foncteur $\Gamma(X^{\mN^\circ}_{s,\et},-)$, 
on obtient un morphisme de suites spectrales 
\begin{equation}
\rE\rightarrow {^\prime {^\prime \rE}}.
\end{equation}
On vérifie aussitôt que le morphisme composé de suites spectrales $\rE\rightarrow {^\prime \rE}$ est défini sur les termes initiaux par 
les morphismes \eqref{ssht7a} et sur les aboutissements par les morphismes \eqref{ssht7b}. 

\begin{prop}\label{ssht8}
Sous les hypothèses de \ref{ssht6}, pour tous entiers $n,j\geq 0$, le diagramme 
\begin{equation}\label{ssht8a}
\xymatrix{
{g_s^{-1}(\xi^{-j}\tOmega^j_{\oX_n/\oS_n})}\ar[r]^-(0.4){u^j}\ar[d]
&{g_s^{-1}(\rR^j\sigma_{s*}(\ocB_n))}\ar[r]^-(0.4){v^j}&{\rR^j\sigma'_{s*}(\Theta_s^{-1}(\ocB_n))}\ar[d]\\
{\xi^{-j}\tOmega^j_{\oX'_n/\oS'_n}}\ar[rr]^-(0.5){u'^j}&&{\rR^j\sigma'_{s*}(\ocB'_n)}}
\end{equation}
où les flèches $u^j$ et $u'^j$ sont induites par le morphisme \eqref{tfkum11a}, 
$v^j$ est le morphisme de changement de base relativement au diagramme commutatif \eqref{ssht6g},
et les flèches verticales sont les morphismes canoniques, est commutatif. 
\end{prop}

D'après (\cite{sga4} XVII 4.1.5), pour tout faisceau abélien $F$ de $\tE'_s$, 
on a un morphisme de changement de base relativement au diagramme commutatif \eqref{ssht6g}
dans la catégorie dérivée des faisceaux abéliens de $X'_{s,\et}$,
\begin{equation}
g_s^*(\rR\sigma_{s*}(F))\rightarrow \rR\sigma'_*(\Theta_s^*(F)). 
\end{equation}
Celui-ci étant fonctoriel en $F$, on en déduit un digramme commutatif 
\begin{equation}\label{ssht8d}
\xymatrix{
{g_s^*(\oa^*(\cL^\gp_\oX))}\ar[r]\ar[rd]&{g_s^*(\rR^1\sigma_{s*}(\mu_{p^n,\tE_s}))}\ar[d]\\
&{\rR^1\sigma'_{s*}(\mu_{p^n,\tE'_s}))}}
\end{equation}
où la flèche verticale est le morphisme de changement de base relativement au diagramme commutatif \eqref{ssht6g}, 
la flèche horizontale est induite par la suite exacte de faisceaux abéliens de $\tE_s$ \eqref{tfkum8a}
\begin{equation}\label{ssht8b}
0\rightarrow \mu_{p^n,\tE_s}\rightarrow \delta^*(\cQ^\gp_n)\rightarrow \sigma_s^*(\oa^*(\cL^\gp_\oX))\rightarrow 0,
\end{equation}
et la flèche oblique est induite par l'image inverse de cette dernière sur $\tE'_s$, réécrite en tenant compte de \eqref{ssht6g},
\begin{equation}\label{ssht8c}
0\rightarrow \mu_{p^n,\tE'_s}\rightarrow \Theta_s^*(\delta^*(\cQ^\gp_n))\rightarrow \sigma'^*_s(g_s^*(\oa^*(\cL^\gp_\oX)))\rightarrow 0.
\end{equation}

On a un morphisme canonique de monoïdes de $\oX'_\et$, 
\begin{equation}
\nu\colon \ogg^*(\cL_\oX)\rightarrow \cL_{\oX'},
\end{equation}
où le monoïde $\cL_\oX$ de $\oX_\et$ est défini dans \eqref{tfkum2b} et $\cL_{\oX'}$ est l'analogue de $\oX'_\et$. 
Par ailleurs, le diagramme d'homomorphismes de monoïdes 
\begin{equation}\label{ssht8e}
\xymatrix{
{\sigma'^*(g^*(\hbar_*(\cL_\oX)))}\ar@{=}[d]\ar[r]&{\sigma'^*(\hbar'_*(\ogg^*(\cL_\oX)))}\ar[rr]^-(0.5){\sigma'^*(\hbar'_*(\nu))}
&&{\sigma'^*(\hbar'_*(\cL_{\oX'}))}\ar[d]^{\ell'}\\
{\Theta^*(\sigma^*(\hbar_*(\cL_\oX)))}\ar[r]^-(0.5){\Theta^*(\ell)}&{\Theta^*(\ocB)}\ar[rr]&&{\ocB'}}
\end{equation}
où $\ell$ et $\ell'$ sont les homomorphismes \eqref{tfkum6b}, la flèche non libellée en haut est le morphisme de changement de base et celle en bas est
l'homomorphisme canonique, est commutatif. Comme le foncteur $\Theta^*$ est exact à gauche, 
on en déduit un homomorphisme de monoïdes \eqref{tfkum6a}
\begin{equation}\label{ssht8f}
\fq\colon \Theta^*(\cQ_n)\rightarrow \cQ'_n.
\end{equation}
On a un diagramme commutatif de suites exactes de $\tE'_s$ 
\begin{equation}\label{ssht8g}
\xymatrix{
0\ar[r]&{\mu_{p^n,\tE'_s}}\ar[r]\ar@{=}[d]&{\delta'^*(\Theta^*(\cQ^\gp_n))}\ar[r]\ar[d]^{\delta'^*(\fq^\gp)}&
{\sigma'^*_s(\oa'^*(\ogg^*(\cL^\gp_\oX)))}\ar[r]\ar[d]^{\sigma'^*_s(\oa'^*(\nu^\gp))}&0\\
0\ar[r]&{\mu_{p^n,\tE'_s}}\ar[r]&{\delta'^*(\cQ'^\gp_n))}\ar[r]&{\sigma'^*_s(\oa'^*(\cL^\gp_{\oX'}))}\ar[r]&0}
\end{equation}
où la suite inférieure est la suite \eqref{tfkum8a} pour $X'$ et la suite supérieure est l'image inverse par $\Theta$ de la suite \eqref{tfkum8a} pour $X$ réécrite 
en tenant compte de \eqref{ssht6f} et \eqref{ssht6g}. 

On déduit de \eqref{ssht8d} et \eqref{ssht8g} que le diagramme 
\begin{equation}\label{ssht8h}
\xymatrix{
{g_s^*(\oa^*(\cL_\oX^\gp))}\ar[r]^-(0.5){g_s^*(\partial_n)}\ar[d]_{\oa'^*(\nu^\gp)}&{g_s^*(\rR^1\sigma_{s*}(\mu_{p^n,\tE_s}))}\ar[d]\\
{\oa'^*(\cL_{\oX'}^\gp)}\ar[r]^-(0.5){\partial'_n}&{\rR^1\sigma'_{s*}(\mu_{p^n,\tE'_s})}}
\end{equation}
où $\partial_n$ et $\partial'_n$ sont les morphismes \eqref{tfkum10a} et la flèche non libellée est le morphisme de changement de base relativement au diagramme commutatif \eqref{ssht6g}, est commutatif. Compte tenu de \ref{tfkum11}, on en déduit que le diagramme 
\begin{equation}\label{ssht8i}
\xymatrix{
{g_s^{-1}(\xi^{-1}\tOmega^1_{\oX_n/\oS_n})}\ar[r]^-(0.4){u^1}\ar[d]
&{g_s^{-1}(\rR^1\sigma_{s*}(\ocB_n))}\ar[r]^-(0.4){v^1}&{\rR^1\sigma'_{s*}(\Theta_s^{-1}(\ocB_n))}\ar[d]\\
{\xi^{-1}\tOmega^1_{\oX'_n/\oS'_n}}\ar[rr]^-(0.5){u'^1}&&{\rR^1\sigma'_{s*}(\ocB'_n)}}
\end{equation}
où les flèches $u^1$ et $u'^1$ sont induites par le morphisme \eqref{tfkum11a}, 
$v^1$ est le morphisme de changement de base relativement au diagramme commutatif \eqref{ssht6g},
et les flèches verticales sont les morphismes canoniques, est commutatif. 

Le diagramme \eqref{ssht8a} est commutatif pour $j=0$ puisque $v^0$ est un homomorphisme de $g_s^{-1}(\co_{\oX_n})$-algèbres. 
Compte tenu de \ref{tfkum11} et \eqref{ssht8i}, 
pour montrer que le diagramme \eqref{ssht8a} est commutatif pour tout $j\geq 1$, il suffit de montrer que le morphisme 
$v^j$  est compatible au cup-produit, 
ou ce qui revient au même, que le morphisme adjoint 
$\rR^j\sigma_{s*}(\ocB_n)\rightarrow g_{s*}(\rR^j\sigma'_{s*}(\Theta_s^{-1}(\ocB_n)))$ est compatible au cup-produit. 
Ce dernier est composé des morphismes 
\begin{eqnarray}\label{ssht8j}
\lefteqn{
{\rR^j\sigma_{s*}(\ocB_n)}\rightarrow {\rR^j\sigma_{s*}(\Theta_{s*}(\Theta_s^{-1}(\ocB_n)))}\rightarrow
{\rR^j(\sigma_{s}\Theta_{s})_*(\Theta_s^{-1}(\ocB_n))}}\\
&&\stackrel{\sim}{\rightarrow}{\rR^j(g_s\sigma'_{s})_*(\Theta_s^{-1}(\ocB_n))}\rightarrow {g_{s*}(\rR^j\sigma'_{s*}(\Theta_s^{-1}(\ocB_n)))},\nonumber
\end{eqnarray}
où la première flèche est induite par le morphisme d'adjonction $\ocB_n \rightarrow \Theta_{s*}(\Theta_s^{-1}(\ocB_n))$, qui est un homomorphisme 
d'algèbres, et est donc compatible au cup-produit, la troisième flèche est induite par l'isomorphisme de commutativité du diagramme \eqref{ssht6g}
et est clairement compatible au cup-produit, et la deuxième et la quatrième flèche sont les edge-homomorphismes de 
la suite spectrale de Cartan-Leray (\cite{sga4} V 5.4).  

Soit $A'$ un anneau commutatif de $\tE'_s$. Considérons le préfaisceau $P^j$ sur $\Et_{/X'_s}$ défini par 
\begin{equation}
U'\in \ob(\Et_{/X'_s})\mapsto \rH^j(\sigma'^*_s(U'),A'), 
\end{equation}
et le préfaisceau $Q^j$ sur $\Et_{/X_s}$ défini par 
\begin{equation}
U\in \ob(\Et_{/X_s})\mapsto P^j(U\times_XX'). 
\end{equation}
Le faisceau $\rR^j\sigma'_{s*}(A')$ (resp. $\rR^j(g_s\sigma'_{s})_*(A')$) est canoniquement isomorphe au faisceau associé à $P^j$ (resp. $Q^j$) 
(\cite{sga4} V 5.1). 
Le morphisme canonique $P^j\rightarrow \rR^j\sigma'_{s*}(A')$ induit un morphisme de préfaisceaux $Q^j\rightarrow g_{s*}(\rR^j\sigma'_{s*}(A'))$ dont le morphisme de faisceaux associés $\rR^j(g_s\sigma'_{s})_*(A')\rightarrow g_{s*}(\rR^j\sigma'_{s*}(A'))$ n'est autre que l'edge-homomorphisme de 
la suite spectrale de Cartan-Leray. Celui-ci est donc compatible au cup-produit. 

Appliquant le foncteur $\rR^j\sigma_{s*}$ à l'homomorphisme canonique de la catégorie dérivée des faisceaux abéliens de $\tE_s$, 
\begin{equation}\label{ssht8k}
\Theta_{s*}(A')\rightarrow \rR \Theta_{s*}(A'),
\end{equation}
on obtient l'edge-homomorphisme de la suite spectrale de Cartan-Leray
\begin{equation}\label{ssht8l}
\rR^j\sigma_{s*}(\Theta_{s*}(A'))\rightarrow \rR^j(\sigma_{s} \Theta_{s})_*(A').
\end{equation}
Soit $C$ un anneau (constant) et supposons que $A'$ soit une $C$-algèbre. D'après (\cite{sp} \href{https://stacks.math.columbia.edu/tag/0B68}{Tag 0B68}),
il existe un morphisme de la catégorie dérivée des $C$-modules de $\tE_s$, 
\begin{equation}\label{ssht8m}
\rR \Theta_{s*}(A')\otimes_C^\rL \rR \Theta_{s*}(A')\rightarrow \rR \Theta_{s*}(A'\otimes_C^\rL A'),
\end{equation}
où le produit tensoriel dérivé est défini dans (\cite{sp} \href{https://stacks.math.columbia.edu/tag/06YH}{Tag 06YH}). 
Par ailleurs, le diagramme de morphisme canoniques
\begin{equation}\label{ssht8n}
\xymatrix{
{\rR \Theta_{s*}(A')\otimes_C^\rL \rR \Theta_{s*}(A')}\ar[rr]&&{\rR \Theta_{s*}(A'\otimes_C^\rL A')}\ar[d]\\
{\Theta_{s*}(A')\otimes_C^\rL \Theta_{s*}(A')}\ar[u]\ar[r]&{\Theta_{s*}(A')}\ar[r]&{\rR\Theta_{s*}(A')}}
\end{equation}
est commutatif. Ceci résulte de la propriété universelle de \eqref{ssht8m} en observant que pour tout $C$-module $F$ de $\tE_s$, 
on a $\rL \Theta_s^*(F)=\Theta_s^*(F)$. On en déduit que le diagramme
\begin{equation}
\xymatrix{
{\rR\sigma_{s*}(\Theta_{s*}(A'))\otimes_C^\rL \rR\sigma_{s*}(\Theta_{s*}(A'))}\ar[r]\ar[d]&
{\rR(\sigma_{s}\Theta_{s})_*(A')\otimes_C^\rL \rR(\sigma_{s}\Theta_{s})_*(A')}\ar[d]\\
{\rR\sigma_{s*}(\Theta_{s*}(A')\otimes_C^\rL\Theta_{s*}(A'))}\ar[r]\ar[dd]
&{\rR\sigma_{s*}(\rR\Theta_{s*}(A')\otimes_C^\rL\rR\Theta_{s*}(A'))}\ar[d]\\
&{\rR(\sigma_{s}\Theta_{s})_*(A'\otimes_C^\rL A')}\ar[d]\\
{\rR\sigma_{s*}(\Theta_{s*}(A'))}\ar[r]&{\rR(\sigma_{s}\Theta_{s})_*(A')}}
\end{equation}
où le carré supérieur traduit la fonctorialité de \eqref{ssht8m} et le rectangle inférieur est induit par le diagramme \eqref{ssht8n}. 
Il s'ensuit par fonctorialité que le morphisme \eqref{ssht8l} est compatible au cup-produit. Il en est alors de même de $v^j$,
d'où la commutativité de \eqref{ssht8a}.

\begin{lem}\label{ssht9}
Sous les hypothèses de \ref{ssht6}, pour tous entiers $n,q\geq 0$, le diagramme 
\begin{equation}\label{ssht9a}
\xymatrix{
{\rH^q(\oX^\circ_\et,\mZ/p^n\mZ)}\ar[r]\ar[d]&{\rH^q(\tE_s,\ocB_n)}\ar[d]\\
{\rH^q(\oX'^\rhd_\et,\mZ/p^n\mZ)}\ar[r]&{\rH^q(\tE'_s,\ocB'_n)}}
\end{equation}
où les flèches horizontales sont induites par le morphisme \eqref{TPCF16a} et les flèches verticales sont induites par les foncteurs 
``image inverse'' et l'homomorphisme canonique $\Theta^{-1}(\ocB)\rightarrow \ocB'$, est commutatif.
\end{lem}

En effet, considérons un diagramme de morphismes de topos, commutatif à isomorphisme près 
\begin{equation}
\xymatrix{
{(T'_2,A'_2)}\ar[r]^{v'}\ar[d]_{u'}&{(T_2,A_2)}\ar[d]^u\\
{(T'_1,A'_1)}\ar[r]^v&{(T_1,A_1)}}
\end{equation}
Le diagramme commutatif d'homomorphismes canoniques d'anneaux de $T_1$ 
\begin{equation}
\xymatrix{
A_1\ar[r]\ar@{=}[d]&{u_*(A_2)}\ar[r]&{u_*(v'_*(A'_2))}\ar@{=}[d]\\ 
A_1\ar[r]&{v_*(A'_1)}\ar[r]&{v_*(u'_*(A'_2))}}
\end{equation}
et le diagramme commutatif de morphismes canoniques de $\bD^+(T_1,A_1)$
\begin{equation}
\xymatrix{
{u_*(v'_*(A'_2))}\ar[r]\ar@{=}[d]&{\rR u_*(v'_*(A'_2))}\ar[r]&{\rR u_*(\rR v'_*(A'_2))}\ar@{=}[d]\\ 
{v_*(u'_*(A'_2))}\ar[r]&{\rR v_*(u'_*(A'_2))}\ar[r]&{\rR v_*(\rR u'_*(A'_2))}}
\end{equation}
induisent un diagramme commutatif de $\bD^+(T_1,A_1)$
\begin{equation}
\xymatrix{
A_1\ar[r]\ar@{=}[d]&{\rR u_*(A_2)}\ar[r]&{\rR u_*(\rR v'_*(A'_2))}\ar@{=}[d]\\ 
A_1\ar[r]&{\rR v_*(A'_1)}\ar[r]&{\rR v_*(\rR u'_*(A'_2))}}
\end{equation}
Prenant la cohomologie $\rR\Gamma(T_1,-)$, on en déduit que le diagramme de morphismes canoniques de $\bD^+(\bMod(\mZ))$
\begin{equation}\label{ssht9b}
\xymatrix{
{\rR\Gamma(T_1,A_1)}\ar[r]\ar[d]&{\rR\Gamma(T_2,A_2)}\ar[d]\\
{\rR\Gamma(T'_1,A'_1)}\ar[r]&{\rR\Gamma(T'_2,A'_2)}}
\end{equation}
est commutatif. 

Considérons maintenant le diagramme de morphismes de topos annelés
\begin{equation}
\xymatrix{
{(\oX'^\rhd_\et,\mZ/p^n\mZ)}\ar[d]_{\upgamma}\ar[r]^{\psi'}&{(\tE',\mZ/p^n\mZ)}\ar[d]^{\Theta}&{(\tE',\ocB'_n)}\ar[d]\ar[l]&
{(\tE'_s,\ocB'_n)}\ar[d]^{\Theta_n}\ar[l]_-(0.5){\delta}\\
{(\oX^\circ_\et,\mZ/p^n\mZ)}\ar[r]^{\psi}&{(\tE,\mZ/p^n\mZ)}&{(\tE,\ocB_n)}\ar[l]&{(\tE_s,\ocB_n)}\ar[l]_-(0.5){\delta'}}
\end{equation}
dont les carrés sont commutatifs à isomorphismes canoniques près. 
Il résulte de ce qui précède que le diagramme
\begin{equation}\label{ssht9c}
\xymatrix{
{\rR\Gamma(\oX^\circ_\et,\mZ/p^n\mZ)}\ar[d]_{\upgamma^*}&{\rR\Gamma(\tE,\mZ/p^n\mZ)}\ar[r]\ar[d]\ar[l]_-(0.5){\psi^*}&{\rR\Gamma(\tE,\ocB_n)}\ar[r]^-(0.5){\delta^*}\ar[d]&
{(\tE_s,\ocB_n)}\ar[d]^{\Theta_n^*}\\
{\rR\Gamma(\oX'^\rhd_\et,\mZ/p^n\mZ)}&{\rR\Gamma(\tE',\mZ/p^n\mZ)}\ar[r]\ar[l]_-(0.5){\psi'^*}&{\rR\Gamma(\tE',\ocB'_n)}\ar[r]^-(0.5){\delta'^*}&
{\rR\Gamma(\tE'_s,\ocB'_n)}}
\end{equation}
est commutatif. 
Les morphismes canoniques $\mZ/p^n\mZ\rightarrow \psi_*(\mZ/p^n\mZ)$ et $\mZ/p^n\mZ\rightarrow \psi'_*(\mZ/p^n\mZ)$ sont des isomorphismes
d'après la preuve de \ref{TPCF16}. Par suite, d'après \ref{acycloc2}, les morphismes $\psi^*$ et $\psi'^*$ de \eqref{ssht9c} sont des isomorphismes.  
La proposition résulte alors de \eqref{ssht9c}.

\begin{prop}\label{ssht10}
Pour tout morphisme $(X',\cM_{X'})\rightarrow (X,\cM_X)$ de $(S,\cM_S)$-schémas logarithmiques adéquats,  
il existe un morphisme canonique de la suite spectrale de Hodge Tate de $X$ sur $S$ vers celle de $X'$ sur $S$,
défini sur les termes initiaux (resp. l'aboutissement) par l'image inverse des différentielles logarithmiques 
(resp. cohomologie étale $p$-adique).
\end{prop}

Cela résulte de \ref{ssht7}, \ref{ssht8} et \ref{ssht9}.

\section{Topos de Faltings relatif, II}\label{ktfr}

\subsection{}\label{ktfr1}
Dans la suite de ce chapitre, en plus des hypothèses de § \ref{hght}, on se donne un morphisme lisse et saturé de schémas logarithmiques 
\begin{equation}\label{ktfr1a}
g\colon (X',\cM_{X'})\rightarrow (X,\cM_X)
\end{equation} 
tel que le morphisme $f'=f\circ g\colon (X',\cM_{X'})\rightarrow (S,\cM_S)$ soit adéquat (\cite{agt} III.4.7).  
On reprend les notations de \ref{mtfla}. 
De plus, on désigne par $G$ (resp. $\tG$) le site (resp. topos) de Faltings relatif associé au couple de morphismes $(h\colon \oX^\circ \rightarrow X,g\colon X'\rightarrow X)$ (cf. \ref{tfr1}). 
D'après \ref{tfr13}, le diagramme commutatif canonique
\begin{equation}\label{ktfr1b}
\xymatrix{
X'\ar[d]_g&{\oX'^\rhd}\ar[d]^{\upgamma}\ar[l]_{h'}\\
X&{\oX^\circ}\ar[l]_h}
\end{equation}
induit des morphismes de topos 
\begin{equation}\label{ktfr1c}
\xymatrix{
\tE'\ar[r]^{\tau}&\tG\ar[r]^{\lgg}&\tE}
\end{equation}
dont le composé est le morphisme de fonctorialité $\Theta\colon \tE'\rightarrow \tE$ \eqref{mtfla11b}.
On rappelle \eqref{tfr13d} que les triangles et carrés du diagramme de morphismes de topos 
\begin{equation}\label{ktfr1d}
\xymatrix{
&\tE'\ar[d]^{\tau}\ar[r]^{\beta'}\ar[ld]_{\sigma'} & \oX'^\rhd_\fet\ar[d]^{\upgamma_\fet}\\
X'_\et\ar[d]_{g_\et}&\tG\ar[d]_-(0.5){\lgg}\ar[r]^-(0.4){\lambda}\ar[l]_-(0.4){\pi}&\oX^\circ_\fet\\
X_\et&\tE\ar[ur]_{\beta}\ar[l]_{\sigma}&}
\end{equation}
où $\pi$ \eqref{tfr4c} et $\lambda$ \eqref{tfr4d}, 
$\sigma$, $\sigma'$ \eqref{tf1kk}, $\beta$ et $\beta'$ \eqref{tf1k} sont les morphismes canoniques, 
sont commutatifs à isomorphismes canoniques près. 

On désigne par 
\begin{equation}\label{ktfr1e}
\varrho\colon X'_\et\gtimes_{X_\et}\oX^\circ_\et\rightarrow \tG
\end{equation}
le morphisme canonique de topos \eqref{tfr6b}. On vérifie aussitôt que les carrés du diagramme de morphismes de topos 
\begin{equation}\label{ktfr1f}
\xymatrix{
{X'_\et\gtimes_{X'_\et}\oX'^\rhd_\et}\ar[d]\ar[r]^-(0.5){\rho'}&\tE'\ar[d]^\tau\\
{X'_\et\gtimes_{X_\et}\oX^\circ_\et}\ar[d]\ar[r]^-(0.5){\varrho}&\tG\ar[d]^\lgg\\
{X_\et\gtimes_{X_\et}\oX^\circ_\et}\ar[r]^-(0.5){\rho}&\tE}
\end{equation}
où $\rho$ et $\rho'$ sont les morphismes canoniques \eqref{mtfla7e} et \eqref{mtfla9de}, 
sont commutatifs à isomorphismes canoniques près (\ref{tfr12} et \cite{agt} VI.4.10).

\subsection{}\label{ktfr20}
Comme $X'_\eta$ est un ouvert de $X'_\et$, {\em i.e.}, un sous-objet de l'objet final, 
$\pi^*(X'_\eta)=(X'_\eta\rightarrow X\leftarrow \oX^\circ)^a$ est un ouvert de $\tG$ \eqref{tfr4a}. 
En vertu de \ref{tfr42}, le topos $\tG_{/\pi^*(X'_\eta)}$ est équivalent au topos de Faltings relatif associé au couple de morphismes 
$(\oX^\circ\rightarrow X, X'_\eta\rightarrow X)$. On note 
\begin{equation}\label{ktfr20a}
\jmath\colon \tG_{/\pi^*(X'_\eta)}\rightarrow \tG
\end{equation}
le morphisme de localisation de $\tG$ en $\pi^*(X'_\eta)$, que l'on identifie au morphisme de fonctorialité relativement 
au diagramme commutatif
\begin{equation}\label{ktfr20b}
\xymatrix{
{X'_\eta}\ar[r]\ar[d]&X\ar[d]&{\oX^\circ}\ar[l]\ar[d]\\
X'\ar[r]&X&\ar[l]{\oX^\circ}}
\end{equation}

On désigne par $\tG_s$ le sous-topos fermé de $\tG$ complémentaire de l'ouvert $\pi^*(X'_\eta)$, 
c'est-à-dire la sous-catégorie pleine de $\tG$ formée des faisceaux $F$ tels que $\jmath^*(F)$
soit un objet final de $\tG_{/\pi^*(X'_\eta)}$  (\cite{sga4} IV 9.3.5), et par 
\begin{equation}\label{ktfr20c}
\kappa\colon \tG_s\rightarrow \tG
\end{equation} 
le plongement canonique (\cite{sga4} IV 9.3.5).

On désigne par $\Pt(\tG)$, $\Pt(\tG_{/\pi^*(X'_\eta)})$ et $\Pt(\tG_s)$ les catégories des points de $\tG$,
$\tG_{/\pi^*(X'_\eta)}$ et $\tG_s$, respectivement, et par 
\begin{equation}\label{ktfr20d}
u\colon \Pt(\tG_{/\pi^*(X'_\eta)})\rightarrow \Pt(\tG) \ \ \ {\rm et}\ \ \ v\colon \Pt(\tG_s)\rightarrow \Pt(\tG)
\end{equation}
les foncteurs induits par $\jmath$ et $\kappa$, respectivement. Ces foncteurs sont pleinement fidèles,
et tout point de $\tG$ appartient à l'image essentielle de l'un ou l'autre de ces foncteurs exclusivement 
(\cite{sga4} IV 9.7.2).

\begin{prop}\label{ktfr22}
\
\begin{itemize}
\item[{\rm (i)}] Soit $(\oy\rightsquigarrow \ox')$ un point de $X'_\et\gtimes_{X_\et}\oX^\circ_\et$ \eqref{tfr21}. 
Pour que $\varrho(\oy\rightsquigarrow \ox')$ \eqref{ktfr1e} appartienne à l'image essentielle de $u$ (resp. $v$) \eqref{ktfr20d}, 
il faut et il suffit que $\ox'$ soit à support dans $X'_\eta$ (resp. $X'_s$). 
\item[{\rm (ii)}] La famille des points de $\tG_{/\pi^*(X'_\eta)}$ (resp. $\tG_s$) définie par la famille des points 
$\varrho(\oy\rightsquigarrow \ox')$ de $\tG$ tels que $\ox'$ soit à support dans $X'_\eta$ (resp. $X'_s$)
est conservative. 
\end{itemize}
\end{prop}

(i) En effet, pour que $\varrho(\oy\rightsquigarrow \ox')$ appartienne à l'image essentielle de $u$ (resp. $v$), 
il faut et il suffit que $(\pi^*(X'_\eta))_{\varrho(\oy\rightsquigarrow \ox')}$ soit un singleton (resp. vide). 
Par ailleurs, on a un isomorphisme canonique \eqref{tfr21a}
\begin{equation}\label{tktfr22a}
(\pi^*(X'_\eta))_{\varrho(\oy\rightsquigarrow \ox')}\stackrel{\sim}{\rightarrow} (X'_\eta)_{\ox'},
\end{equation}
d'où la proposition.  

(ii) Cela résulte de (i), \ref{lptfr20} et (\cite{sga4} IV 9.7.3).

\subsection{}\label{ktfr21}
En vertu de (\cite{sga4} IV 9.4.3), il existe un morphisme de topos
\begin{equation}\label{ktfr21a}
\pi_s\colon \tG_s\rightarrow X'_{s,\et}
\end{equation} 
unique à isomorphisme canonique près tel que le diagramme
\begin{equation}\label{ktfr21b}
\xymatrix{
{\tG_s}\ar[r]^{\pi_s}\ar[d]_{\kappa}&{X'_{s,\et}}\ar[d]^{a'}\\
{\tG}\ar[r]^{\pi}&{X'_\et}}
\end{equation}
où $a'$ est l'injection canonique, soit commutatif à isomorphisme près, et même $2$-cartésien.
Par définition, on a un isomorphisme canonique 
\begin{equation}\label{ktfr21j}
\pi^*\circ a'_*\stackrel{\sim}{\rightarrow}\kappa_*\circ \pi_s^*.
\end{equation}
Les foncteurs $a'_*$ et $\kappa_*$ étant exacts, 
pour tout groupe abélien $F$ de $\tG_s$ et tout entier $i\geq 0$, on a un isomorphisme canonique 
\begin{equation}\label{ktfr21k}
a'_*(\rR^i\pi_{s*}(F))\stackrel{\sim}{\rightarrow}\rR^i\pi_*(\kappa_*F). 
\end{equation}

On a un isomorphisme canonique $\tau^*(\pi^*(X'_\eta))\simeq \sigma'^*(X'_\eta)$ \eqref{ktfr1d}.
En vertu de (\cite{sga4} IV 9.4.3), il existe donc un morphisme de topos
\begin{equation}\label{ktfr21c}
\tau_s\colon \tE'_s\rightarrow \tG_s
\end{equation}
unique à isomorphisme canonique près tel que le diagramme 
\begin{equation}\label{ktfr21d}
\xymatrix{
{\tE'_s}\ar[r]^{\tau_s}\ar[d]_{\delta'}&{\tG_s}\ar[d]^{\kappa}\\
{\tE'}\ar[r]^\tau&{\tG}}
\end{equation}
soit commutatif à isomorphisme près (cf. \ref{mtfla9}). 
Les foncteurs $\delta'_*$ et $\kappa_*$ étant exacts, 
pour tout groupe abélien $F$ de $\tE'_s$ et tout entier $i\geq 0$, on a un isomorphisme canonique 
\begin{equation}\label{ktfr21l}
\kappa_*(\rR^i\tau_{s*}(F))\stackrel{\sim}{\rightarrow}\rR^i\tau_*(\delta'_*F). 
\end{equation}
Il résulte de \eqref{ktfr1d} et (\cite{sga4} IV 9.4.3) 
que le diagramme de morphismes de topos 
\begin{equation}\label{ktfr21e}
\xymatrix{
{\tE'_s}\ar[r]^{\tau_s}\ar[d]_{\sigma'_s}&{\tG_s}\ar[dl]^{\pi_s}\\
{X'_{s,\et}}&}
\end{equation}
est commutatif à isomorphisme canonique près. 

On a un isomorphisme canonique $\lgg^*(\sigma^*(X_\eta))\simeq \pi^*(X'_\eta)$ \eqref{ktfr1d}.
En vertu de (\cite{sga4} IV 9.4.3), il existe donc un morphisme de topos
\begin{equation}\label{ktfr21f}
\lgg_s\colon \tG_s\rightarrow \tE_s
\end{equation}
unique à isomorphisme canonique près tel que le diagramme 
\begin{equation}\label{ktfr21g}
\xymatrix{
{\tG_s}\ar[r]^{\lgg_s}\ar[d]_{\kappa}&{\tE_s}\ar[d]^{\delta}\\
{\tG}\ar[r]^\lgg&{\tE}}
\end{equation}
soit commutatif à isomorphisme près (cf. \ref{mtfla7}). Il résulte de \eqref{ktfr1d} et (\cite{sga4} IV 9.4.3) 
que le diagramme de morphismes de topos 
\begin{equation}\label{ktfr21h}
\xymatrix{
{\tG_s}\ar[r]^{\lgg_s}\ar[d]_{\pi_s}&{\tE_s}\ar[d]^{\sigma_s}\\
{X'_{s,\et}}\ar[r]^{g_{s,\et}}&{X_{s,\et}}}
\end{equation}
est commutatif à isomorphisme canonique près.

Il résulte encore de (\cite{sga4} IV 9.4.3) que le composé $\lgg_s\circ \tau_s$ est le morphisme \eqref{mtfla11d}
\begin{equation}\label{ktfr21i}
\Theta_s\colon \tE'_s\rightarrow \tE_s.
\end{equation}

\begin{prop}\label{ktfr25}
Supposons $g$ propre, et considérons le diagramme de morphismes de topos, commutatif à isomorphisme canonique près \eqref{ktfr1d},
\begin{equation}\label{ktfr25a}
\xymatrix{
{\tG}\ar[r]^{\lgg}\ar[d]_{\pi}&{\tE}\ar[d]^{\sigma}\\
{X'_{\et}}\ar[r]^{g_{\et}}&{X_{\et}}}
\end{equation}
Alors, 
\begin{itemize}
\item[{\rm (i)}] Pour tout faisceau $F$ de $X'_\et$, le morphisme de changement de base relativement au diagramme \eqref{ktfr25a}
\begin{equation}\label{ktfr25b}
\sigma^*(g_{\et*}(F))\rightarrow \lgg_*(\pi^*(F))
\end{equation}  
est un isomorphisme. 
\item[{\rm (ii)}] Pour tout faisceau abélien de torsion $F$ de $X'_\et$
et tout entier $q\geq 0$, le morphisme de changement de base relativement au diagramme \eqref{ktfr25a}
\begin{equation}\label{ktfr25c}
\sigma^*(\rR^qg_{\et*}(F))\rightarrow \rR^q\lgg_*(\pi^*(F))
\end{equation}  
est un isomorphisme. 
\end{itemize}
\end{prop}

C'est un cas particulier de \ref{lptfr15}. 

\begin{cor}\label{ktfr26}
Supposons $g$ propre. Alors, 
\begin{itemize}
\item[{\rm (i)}] Pour tout faisceau $F$ de $X'_{s,\et}$, le morphisme de changement de base relativement au diagramme \eqref{ktfr21h}
\begin{equation}\label{ktfr26a}
\sigma^*_s(g_{s,\et*}(F))\rightarrow \lgg_{s*}(\pi^*_s(F))
\end{equation}  
est un isomorphisme. 
\item[{\rm (ii)}] Pour tout faisceau abélien de torsion $F$ de $X'_{s,\et}$
et tout entier $q\geq 0$, le morphisme de changement de base relativement au diagramme \eqref{ktfr21h}
\begin{equation}\label{ktfr26b}
\sigma^*_s(\rR^qg_{s,\et*}(F))\rightarrow \rR^q\lgg_{s*}(\pi^*_s(F))
\end{equation}  
est un isomorphisme. 
\end{itemize}
\end{cor}

Cela résulte de \ref{ktfr25} et (\cite{egr1} 1.2.4). 
On notera pour (ii) que les foncteurs $a_*$, $a'_*$, $\delta_*$ et $\kappa_*$ étant exacts, pour tout 
$q\geq 0$, on a des isomorphismes canoniques $a_*\circ \rR^qg_{s,\et*}\simeq \rR^qg_{\et*}\circ a'_*$ et 
$\delta_*\circ \rR^q\lgg_{s*}\simeq \rR^q\lgg_*\circ \kappa_*$.

\subsection{}\label{ktfr12}
Considérons un diagramme commutatif de morphismes de schémas
\begin{equation}\label{ktfr12a}
\xymatrix{
U'\ar[r]\ar[d]&U\ar[d]\\
X'\ar[r]&X}
\end{equation}
tel que les flèches verticales soient des morphismes étales. 
Pour tout préfaisceau $F$ sur $G$, on définit le préfaisceau $F_{U'\rightarrow U}$ sur $\Et_{\rf/\oU^\circ}$ en posant pour tout $V\in \ob(\Et_{\rf/\oU^\circ})$,
\begin{equation}\label{ktfr12b}
F_{U'\rightarrow U}(V)=F(U'\rightarrow U\leftarrow V).
\end{equation}
Si $F$ est un faisceau de $\tG$, alors $F_{U'\rightarrow U}$ est un faisceau de $\oU^\circ_\fet$.

\subsection{}\label{ktfr13}
A tout point géométrique $\ox'$ de $X'$, on associe une catégorie $\cC_{\ox'}$ de la façon suivante. 
Les objets de $\cC_{\ox'}$ sont les diagrammes commutatifs de morphismes de schémas 
\begin{equation}\label{ktfr13a}
\xymatrix{
&U'\ar[r]\ar[d]&U\ar[d]\\
\ox'\ar[r]\ar[ru]&X'\ar[r]&X}
\end{equation}
tels que les morphismes $U'\rightarrow X'$ et $U\rightarrow X$ soient étales. 
Un tel objet sera noté $(\ox'\rightarrow U'\rightarrow U)$. Soient $(\ox'\rightarrow U'\rightarrow U)$, $(\ox'\rightarrow U'_1\rightarrow U_1)$ deux objets de
$\cC_{\ox'}$. Un morphisme de $(\ox'\rightarrow U'_1\rightarrow U_1)$ vers $(\ox'\rightarrow U'\rightarrow U)$ est la donnée d'un $X'$-morphisme
$U'_1\rightarrow U'$ et d'un $X$-morphisme $U_1\rightarrow U$ tels que le diagramme 
\begin{equation}\label{ktfr13b}
\xymatrix{
&U'_1\ar[r]\ar[d]&U_1\ar[d]\\
\ox'\ar[r]\ar[ru]&U'\ar[r]&U}
\end{equation}
soit commutatif. On observera que les produits fibrés sont représentables dans $\cC_{\ox'}$ (cf. la preuve de \cite{agt} VI.10.3).
Les limites projectives finies sont donc représentables dans $\cC_{\ox'}$ (cf. \cite{sga4} I 2.3).
Par suite, la catégorie $\cC_{\ox'}$ est cofiltrante (\cite{sga4} I 2.7.1).

\subsection{}\label{ktfr14}
Soient $(\oy\rightsquigarrow \ox')$ un point de $X'_\et\gtimes_{X_\et}\oX^\circ_\et$ \eqref{tfr21}, 
$\uX'$ le localisé strict de $X'$ en $\ox'$, $\uX$ le localisé strict de $X$ en $g(\ox')$, 
$\cC_{\ox'}$ la catégorie associée à $\ox'$ dans \ref{ktfr13}. 
On désigne par $u\colon \oy\rightarrow \uX$ le $X$-morphisme qui définit le point $(\oy\rightsquigarrow \ox')$, 
et par $v\colon \oy\rightarrow \uoX^\circ$ le $\oX^\circ$-morphisme induit \eqref{mtfla2b}. 
Pour chaque objet $(\ox'\rightarrow U'\rightarrow U)$ de $\cC_{\ox'}$ \eqref{ktfr13}, on a un $X'$-morphisme canonique 
$\uX'\rightarrow U'$ et un $X$-morphisme canonique $\uX\rightarrow U$ qui s'insèrent dans un diagramme commutatif  
\begin{equation}\label{ktfr14a}
\xymatrix{
&\uX'\ar[r]^g\ar[d]&\uX\ar[d]\\
\ox'\ar[r]\ar[ru]&U'\ar[r]&U}
\end{equation}
où on a encore noté $g$ le morphisme induit par $g$. On en déduit un morphisme $\uoX\rightarrow \oU$. 
Le morphisme $v\colon \oy\rightarrow \uoX^\circ$ induit alors un point géométrique de $\oU^\circ$ que l'on note encore $\oy$. 
Le diagramme 
\begin{equation}\label{ktfr14b}
\xymatrix{
\uX\ar[d]&\oy\ar[l]_-(0.5)u\ar[d]\\
U&\oU^\circ\ar[l]}
\end{equation}
est commutatif. 

On désigne par $\varrho(\oy\rightsquigarrow \ox')$ l'image de $(\oy\rightsquigarrow \ox')$ par  le morphisme $\varrho$ \eqref{ktfr1e}, 
qui est donc un point de $\tG$, et par $\cP_{\varrho(\oy\rightsquigarrow \ox')}$ 
la catégorie des objets $\varrho(\oy \rightsquigarrow \ox')$-pointés de $G$. On rappelle \eqref{tfr22} que 
les objets de $\cP_{\varrho(\oy \rightsquigarrow \ox')}$ sont des quadruplets 
formés d'un objet $(U'\rightarrow U\leftarrow V)$ de $G$, d'un $X'$-morphisme $\uX'\rightarrow U'$,
d'un $X$-morphisme $\uX\rightarrow U$ et 
d'un $\oX^\circ$-morphisme $\oy\rightarrow V$ tels que les carrés du diagramme
\begin{equation}\label{ktfr14c}
\xymatrix{
{\uX'}\ar[d]\ar[r]^-(0.5)g&{\uX}\ar[d]&\oy\ar[l]_-(0.5)u\ar[d]\\
U'\ar[r]&U&V\ar[l]}
\end{equation}
soient commutatifs. On dispose donc du foncteur 
\begin{equation}\label{ktfr14d}
\begin{array}{clcr}
\cP_{\varrho(\oy\rightsquigarrow \ox')}&\rightarrow& \cC_{\ox'}\\
\xymatrix{
{\uX'}\ar[d]\ar[r]^-(0.5)g&{\uX}\ar[d]&\oy\ar[l]_-(0.5)u\ar[d]\\
U'\ar[r]&U&V\ar[l]}&\mapsto &
\xymatrix{
{\uX'}\ar[d]\ar[r]^-(0.5)g&{\uX}\ar[d]\\
U'\ar[r]&U}
\end{array}
\end{equation}
dont la catégorie fibre au-dessus d'un objet $(\ox'\rightarrow U'\rightarrow U)$ de $\cC_{\ox'}$ s'identifie à la catégorie 
des voisinages de $\oy$ dans $\Et_{\rf/\oU^\circ}$, autrement dit à la catégorie des revêtements étales $\oy$-pointés de $\oU^\circ$ \eqref{ktfr14b}.

D'après \eqref{tfr22b} et (\cite{sga4} IV (6.8.4)), pour tout préfaisceau $F$ sur $G$, on a un isomorphisme canonique
\begin{equation}\label{ktfr14e}
F^a_{\varrho(\oy\rightsquigarrow \ox')}\stackrel{\sim}{\rightarrow} \underset{\underset{(\ox'\rightarrow U'\rightarrow U)\in \cC^\circ_{\ox'}}{\longrightarrow}}{\lim}\ 
(F^a_{U'\rightarrow U})_{\rho_{\oU^\circ}(\oy)},
\end{equation}
où $F_{U'\rightarrow U}$ est le préfaisceau sur $\Et_{\rf/\oU^\circ}$ défini dans \eqref{ktfr12b} et 
$\rho_{\oU^\circ}\colon \oU^\circ_\et \rightarrow \oU^\circ_\fet$ est le morphisme canonique \eqref{notconv10a}.

\subsection{}\label{ktfr4}
Soient $\ox'$ un point géométrique de $X'$, $\uX'$ le localisé strict de $X'$ en $\ox'$, $\uX$ le localisé strict de $X$ en $g(\ox')$. 
On désigne par $\uG$ (resp. $\tuG$) le site (resp. topos) de Faltings relatif associé au couple de morphismes 
$(\uh\colon \uoX^\circ\rightarrow \uX, \ug\colon \uX'\rightarrow \uX)$ induits par $h$ et $g$ \eqref{mtfla2},  par 
\begin{equation}\label{ktfr4a}
\Phi\colon \tuG\rightarrow \tG
\end{equation}
le morphisme de fonctorialité \eqref{tfr7b} et par 
\begin{equation}\label{ktfr4b}
\vartheta\colon \uoX^\circ_\fet\rightarrow \tuG
\end{equation}
le morphisme défini dans \eqref{tfr24i}. On pose  
\begin{equation}\label{ktfr4c}
\phi_{\ox'}=\vartheta^*\circ \Phi^*\colon \tG\rightarrow \uoX^\circ_\fet.
\end{equation}
On note 
\begin{equation}\label{ktfr4g}
\uupgamma\colon \uoX'^\rhd\rightarrow \uoX^\circ
\end{equation} 
le morphisme induit par $\ug$. 

D'après \ref{tfr30}, pour tout groupe abélien $F$ de $\tG$
et tout entier $q\geq 0$, on a un isomorphisme canonique et fonctoriel
\begin{equation}\label{ktfr4d}
\rR^q\pi_*(F)_{\ox'}\stackrel{\sim}{\rightarrow}\rH^q(\uoX^\circ_\fet,\phi_{\ox'}(F)). 
\end{equation}

\subsection{}\label{ktfr28}
Soit $(\oy\rightsquigarrow \ox')$ un point de $X'_\et\gtimes_{X_\et}\oX^\circ_\et$ \eqref{tfr21}. 
Reprenons les notations de \ref{ktfr4} et notons $\rho_{\uoX^\circ}\colon \uoX^\circ_\et\rightarrow \uoX^\circ_\fet$ 
le morphisme canonique \eqref{notconv10a}. 
Le point $(\oy\rightsquigarrow \ox')$ est défini par un $X$-morphisme $\oy\rightarrow \uX$ 
et il induit donc un $\oX^\circ$-morphisme $\oy\rightarrow \uoX^\circ$.
Pour tout objet $F$ de $\tG$, on a un isomorphisme canonique 
\begin{equation}\label{ktfr28a}
F_{\varrho(\oy\rightsquigarrow \ox')}\stackrel{\sim}{\rightarrow}(\phi_{\ox'}(F))_{\rho_{\uoX^\circ}(\oy)},
\end{equation}
où $\varrho$ est le morphisme \eqref{ktfr1e}.  
En effet, on vérifie aussitôt que le diagramme de morphismes de topos 
\begin{equation}\label{ktfr28b}
\xymatrix{
{\uX'_\et\gtimes_{\uX_\et}\uoX^\circ_\et}\ar[d]\ar[r]^-(0.5){\uvarrho}&\tuG\ar[d]^\Phi\\
{X'_\et\gtimes_{X_\et}\oX^\circ_\et}\ar[r]^-(0.5){\varrho}&\tG}
\end{equation}
où $\uvarrho$ est le morphisme canonique \eqref{tfr6b} et la flèche non libellée est le morphisme de fonctorialité, est commutatif. 
L'isomorphisme \eqref{ktfr28a} résulte alors de \ref{tfr250} en considérant $(\oy\rightsquigarrow \ox')$ comme un point de 
$\uX'_\et\gtimes_{\uX_\et}\uoX^\circ_\et$. 

\begin{lem}\label{ktfr29}
Soient $(\oy_1\rightsquigarrow \ox')$ et $(\oy_2\rightsquigarrow \ox')$ deux points de $X'_\et\gtimes_{X_\et}\oX^\circ_\et$, 
$F$ un objet de $\tG$. Alors, il existe un isomorphisme fonctoriel 
\begin{equation}\label{ktfr29a}
F_{\varrho(\oy_1\rightsquigarrow \ox')}\stackrel{\sim}{\rightarrow}F_{\varrho(\oy_2\rightsquigarrow \ox')}.
\end{equation}
\end{lem}

En effet, notant $\uX$ le localisé strict de $X$ en $g(\ox')$,  
le schéma $\uoX$ est normal et strictement local d'après (\cite{agt} III.3.7). Par suite, le schéma $\uoX^\circ$ est intègre. 
La proposition résulte alors de \eqref{ktfr28a} et du fait que tous les foncteurs fibres de $\uoX^\circ_\fet$ sont isomorphes (\cite{agt} VI.9.11).

\subsection{}\label{ktfr15}
Reprenons les hypothèses et notations de \ref{ktfr4},  
et désignons, de plus, par $\uE'$ (resp. $\tuE'$) le site (resp. topos) de Faltings associé au morphisme 
$\uh'\colon \uoX'^\rhd\rightarrow \uX'$ induit par $h'$ \eqref{mtfla2}  (cf. \ref{tf1}), par 
\begin{equation}\label{ktfr15a}
\Phi'\colon \tuE'\rightarrow \tE'
\end{equation}
le morphisme de fonctorialité (\cite{agt} VI.10.12), par 
\begin{equation}\label{ktfr15b}
\utau\colon \tuE'\rightarrow \tuG
\end{equation}
le morphisme canonique \eqref{tfr13c} et par 
\begin{equation}\label{ktfr15c}
\theta'\colon \uoX'^\rhd_\fet\rightarrow \tuE'
\end{equation}
le morphisme défini dans (\cite{agt} (VI.10.23.1)).  On pose  
\begin{equation}\label{ktfr15d}
\varphi'_{\ox'}=\theta'^*\circ \Phi'^*\colon \tE'\rightarrow \uoX'^\rhd_\fet.
\end{equation}

\begin{prop}\label{ktfr16}
Conservons les hypothèses et notations de \ref{ktfr4} et \ref{ktfr15}. 
\begin{itemize} 
\item[{\rm (i)}] Le diagramme de morphismes de topos 
\begin{equation}\label{ktfr16a}
\xymatrix{
{\tuE'}\ar[r]^{\Phi'}\ar[d]_\utau&{\tE'}\ar[d]^{\tau}\\
{\tuG}\ar[r]^{\Phi}&{\tG}}
\end{equation}
est commutatif à isomorphisme canonique près. 
\item[{\rm (ii)}] Pour tout groupe abélien $F$ de $\tE'$ et tout entier $q\geq 0$, 
le morphisme de changement de base relativement au diagramme \eqref{ktfr16a}
\begin{equation}\label{ktfr16b}
\Phi^*(\rR^q\tau_*(F))\rightarrow \rR^q\utau_*(\Phi'^*(F))
\end{equation}
est un isomorphisme. 
\end{itemize}
\end{prop}

(i) Cela résulte aussitôt des définitions \eqref{tfr7b} et \eqref{tfr13c}.

(ii) Soient $\cC_{\ox'}$ la catégorie définie dans \ref{ktfr13}, $(\ox'\rightarrow U'\rightarrow U)$ un objet de $\cC_{\ox'}$. 
On a un isomorphisme canonique
\begin{equation}\label{ktfr16e}
\tau^*((U'\rightarrow U\leftarrow \oU^\circ)^a)\stackrel{\sim}{\rightarrow} (\oU'^\rhd\rightarrow U')^a. 
\end{equation}
On désigne par $\tG_{U'\rightarrow U}$ le topos de Faltings relatif associé au couple de morphismes $(\oU^\circ\rightarrow U, U'\rightarrow U)$ \eqref{tfr1} et
par $\tE'_{U'}$ le topos de Faltings associé au morphisme $\oU'^\rhd\rightarrow U'$ \eqref{tf1}. Le diagramme de morphismes de fonctorialité de topos
\begin{equation}\label{ktfr16f}
\xymatrix{
{\tuE'}\ar[r]^{\Psi'_{U'}}\ar[d]_{\utau}&{\tE'_{U'}}\ar[r]^{\Phi'_{U'}}\ar[d]_{\tau_{U'\rightarrow U}}&{\tE'}\ar[d]^{\tau}\\
{\tuG}\ar[r]^-(0.45){\Psi_{U'\rightarrow U}}&{\tG_{U'\rightarrow U}}\ar[r]^-(0.45){\Phi_{U'\rightarrow U}}&{\tG}}
\end{equation}
est commutatif à isomorphisme canonique près. D'après \ref{tfr42}, on a une équivalence canonique de topos
\begin{equation}\label{ktfr16g}
\tG_{U'\rightarrow U}\stackrel{\sim}{\rightarrow} \tG_{/(U'\rightarrow U\leftarrow \oU^\circ)^a},
\end{equation}
et le morphisme $\Phi_{U'\rightarrow U}$ s'identifie au morphisme de localisation de $\tG$ en $(U'\rightarrow U\leftarrow \oU^\circ)^a$. 
D'après (\cite{agt} VI.10.14), on a une équivalence canonique de topos
\begin{equation}\label{ktfr16h}
\tE'_{U'}\stackrel{\sim}{\rightarrow} \tE'_{/(\oU'^\rhd\rightarrow U')^a},
\end{equation}
et le morphisme $\Phi'_{U'}$ s'identifie au morphisme de localisation de $\tE'$ en $(\oU'^\rhd\rightarrow U')^a$. 
De plus, $\tau_{U'\rightarrow U}$ s'identifie au morphisme déduit de $\tau$ par localisation en vertu de  \eqref{ktfr16e} et (\cite{sga4} IV 5.11).

Soit $(\ox'\rightarrow U'_0\rightarrow U_0)$ un objet de $\cC_{\ox'}$ tel que les schémas $U_0$ et $U'_0$ soient affines. 
On désigne par $I$ la sous-catégorie pleine de $(\cC_{\ox'})_{/(\ox'\rightarrow U'_0\rightarrow U_0)}$ formée des objets 
$(\ox'\rightarrow U'\rightarrow U)$ tels que les morphismes $U'\rightarrow U'_0$ et $U\rightarrow U_0$ soient affines. 
Le foncteur canonique $I\rightarrow \cC_{\ox'}$ est alors cofinal et la catégorie $I$ est cofiltrante d'après (\cite{sga4} I 8.1.3).

Suivant \ref{lptfr2}, on désigne par 
\begin{equation}\label{ktfr16i}
\fA\rightarrow I
\end{equation}
le topos fibré obtenu en associant à tout objet $(\ox'\rightarrow U'\rightarrow U)$ de $I$ le topos $\tG_{U'\rightarrow U}$,
et à tout morphisme de $I$ donné par le diagramme commutatif
\begin{equation}\label{ktfr16j}
\xymatrix{
&U'_1\ar[r]\ar[d]&U_1\ar[d]\\
\ox'\ar[r]\ar[ru]&U'\ar[r]&U}
\end{equation}
le foncteur image inverse $\tG_{U'\rightarrow U}\rightarrow \tG_{U'_1\rightarrow U_1}$ par le morphisme de fonctorialité. 
En vertu de \ref{lptfr3}, les morphismes $\Psi_{U'\rightarrow U}$ identifient le topos $\tuG$ à la limite projective du topos fibré $\fA$.

Suivant (\cite{agt} VI.11.2), on désigne par 
\begin{equation}\label{ktfr16k}
\fB\rightarrow I
\end{equation}
le topos fibré obtenu en associant à tout objet $(\ox'\rightarrow U'\rightarrow U)$ de $I$ le topos $\tE'_{U'}$,
et à tout morphisme de $I$ donné par le diagramme commutatif \eqref{ktfr16j}
le foncteur image inverse $\tE'_{U'}\rightarrow \tE'_{U'_1}$ par le morphisme de fonctorialité. 
En vertu de (\cite{agt} VI.11.3), les morphismes $\Psi_{U'}$ identifient le topos $\tuE'$ à la limite projective du topos fibré $\fB$.

Les morphismes $\tau_{U'\rightarrow U}$, pour $(\ox'\rightarrow U'\rightarrow U)\in \ob(I)$, 
définissent un morphisme de topos fibrés (\cite{sga4} VI 7.1.6)
\begin{equation}\label{ktfr16l}
\Lambda\colon \fB\rightarrow \fA.
\end{equation}
Le morphisme $\utau$ se déduit de $\Lambda$ par passage à la limite projective (\cite{sga4} VI 8.1.4). 

D'après (\cite{sga4} VI 8.7.5), pour tout groupe abélien $F$ de $\tE'$, on a un isomorphisme canonique
\begin{equation}
\rR^q\utau_*(\Phi'^*(F))\stackrel{\sim}{\rightarrow} \underset{\underset{(\ox'\rightarrow U'\rightarrow U)\in \ob(I)}{\longrightarrow}}{\lim}\ 
\Psi_{U'\rightarrow U}^*(\rR^q\tau_{U'\rightarrow U*}(\Phi'^*_{U'}(F))).
\end{equation}
On notera que les conditions requises dans (\cite{sga4} VI 8.7.1) sont satisfaites en vertu de \ref{tf4}, \ref{tfr16} et (\cite{sga4} VI 3.3 et 5.1).
Par ailleurs, compte tenu de (\cite{egr1} 1.2.4(ii)), le morphisme \eqref{ktfr16b} s'identifie à la limite inductive des morphismes 
\begin{equation}
\Psi_{U'\rightarrow U}^*(\Phi_{U'\rightarrow U}^*(\rR^q\tau_*(F)))\rightarrow \Psi_{U'\rightarrow U}^*(\rR^q\tau_{U'\rightarrow U*}(\Phi'^*_{U'}(F)))
\end{equation}
déduits des morphismes de changement de base relativement au carré de droite de \eqref{ktfr16f}. 
Ces derniers sont des isomorphismes d'après (\cite{sga4} V 5.1(3)); d'où la proposition.

\begin{prop}\label{ktfr17}
Conservons les hypothèses et notations de \ref{ktfr4} et \ref{ktfr15}. 
\begin{itemize} 
\item[{\rm (i)}] Le diagramme de morphismes de topos 
\begin{equation}\label{ktfr17a}
\xymatrix{
{\uoX'^\rhd_\fet}\ar[r]^{\theta'}\ar[d]_{\uupgamma_\fet}&{\tuE'}\ar[d]^{\utau}\\
{\uoX^\circ_\fet}\ar[r]^{\vartheta}&{\tuG}}
\end{equation}
est commutatif à isomorphisme canonique près. 
\item[{\rm (ii)}] Pour tout groupe abélien $F$ de $\tuE'$ et tout entier $q\geq 0$, 
le morphisme de changement de base relativement au diagramme \eqref{ktfr17a}
\begin{equation}\label{ktfr17b}
\vartheta^*(\rR^q\utau_*(F))\rightarrow \rR^q\uupgamma_{\fet*}(\theta'^*(F))
\end{equation}
est un isomorphisme. 
\end{itemize}
\end{prop}

(i) On désigne par $\uG_\rf$ la sous-catégorie pleine de $\uG$ \eqref{ktfr4} formée des objets $(U'\rightarrow U\leftarrow V)$ 
tels que les morphismes $U\rightarrow \uX$ et $U'\rightarrow \uX'$ soient finis étales. Il s'ensuit que le morphisme $V\rightarrow \uoX^\circ$ est aussi fini étale. 
On munit $\uG_\rf$ de la topologie engendrée par les recouvrements de type (a), (b) et (c) \eqref{tfr1}.  
Le topos des faisceaux de $\mU$-ensembles sur $\uG_\rf$ est le produit orienté $\uX'_\fet\gtimes_{\uX_\fet} \uoX^\circ_\fet$ 
des morphismes de topos $\ug_\fet\colon \uX'_\fet\rightarrow \uX_\fet$ et $\uh_\fet\colon \uoX^\circ_\fet\rightarrow \uX_\fet$  
(cf. \cite{agt} VI.3.7 et VI.3.10). 

Pour tout $U\in \ob(\Et_{\scoh/\uX})$ \eqref{notconv10}, 
on désigne par $U^\rf$ la somme disjointe des localisés stricts de $U$ en les points de $U_{g(\ox')}$;
c'est un sous-schéma ouvert et fermé de $U$, qui est fini sur $\uX$ (\cite{ega4} 18.5.11). 
De même, pour tout $U'\in \ob(\Et_{\scoh/\uX'})$, 
on désigne par $U'^\rf$ la somme disjointe des localisés stricts de $U'$ en les points de $U'_{\ox'}$. 

D'après \ref{tfr24}, le foncteur 
\begin{equation}\label{ktfr17c}
\nu^+\colon 
\begin{array}[t]{clcr}
\uG_\scoh&\rightarrow& \uG_\rf\\
(U'\rightarrow U\leftarrow V)&\mapsto&  (U'^\rf\rightarrow U^\rf\leftarrow U^\rf\times_UV).
\end{array}
\end{equation}
est exact à gauche et continu. Compte tenu de \ref{tfr15}(iii), il définit donc un morphisme de topos 
\begin{equation}\label{ktfr17d}
\nu\colon \uX'_\fet\gtimes_{\uX_\fet} \uoX^\circ_\fet\rightarrow \tuG.
\end{equation}

Le $2$-isomorphisme canonique $\uh_\fet \uupgamma_\fet\stackrel{\sim}{\rightarrow} \ug_\fet \uh'_\fet$ induit un morphisme 
\begin{equation}\label{ktfr17e}
t\colon \uoX'^\rhd_\fet \rightarrow \uX'_\fet\gtimes_{\uX_\fet} \uoX^\circ_\fet
\end{equation}
qui s'insère dans un diagramme (non-commutatif)
\begin{equation}\label{ktfr17f}
\xymatrix{
&{\uoX'^\rhd_\fet}\ar[d]_{t}\ar[dl]_{\uh'_\fet}\ar[rd]^{\uupgamma_\fet}&\\
{\uX'_\fet}\ar[rd]_{\ug_\fet}&{\uX'_\fet\gtimes_{\uX_\fet}\uoX^\circ_\fet}\ar[l]_-(0.4){\rp_1}\ar[r]^-(0.4){\rp_2}&
{\uoX^\circ_\fet}\ar[ld]^{\uh_\fet}\\
&{\uX_\fet}&}
\end{equation} 
D'après (\cite{agt} VI.3.6), pour tout objet $(U'\rightarrow U\leftarrow V)$ de $\uG_\rf$, on a un isomorphisme canonique 
\begin{equation}\label{ktfr17g}
t^*((U'\rightarrow U\leftarrow V)^a)\stackrel{\sim}{\rightarrow} \oU'^\rhd\times_{\oU^\circ}V.
\end{equation}

Pour tout objet $(U'\rightarrow U\leftarrow V)$ de $\uG$, on a un isomorphisme canonique (cf. \eqref{tfr7d} et \eqref{tfr9b})
\begin{equation}\label{ktfr17h}
\utau^*((U'\rightarrow U\leftarrow V)^a)\stackrel{\sim}{\rightarrow} (\oU'^\rhd\times_{\oU^\circ}V\rightarrow U')^a.
\end{equation}

Pour tout objet $(W\rightarrow U')$ de $E'_\scoh$, on a un isomorphisme canonique 
\begin{equation}\label{ktfr17i}
\theta'^*((W\rightarrow U')^a)\stackrel{\sim}{\rightarrow}U'^\rf\times_{U'}W.
\end{equation}

Pour tout objet $(U'\rightarrow U\leftarrow V)$ de $\uG_\scoh$, posant $\oU'^{\rf\rhd}=U'^\rf\times_{U'}\oU'^\rhd$ et 
$\oU^{\rf\circ}=U^\rf\times_{U}\oU^\circ$, on a un isomorphisme canonique 
\begin{equation}\label{ktfr17j}
U'^\rf\times_{U'}\oU'^\rhd\times_{\oU^\circ}V \stackrel{\sim}{\rightarrow} \oU'^{\rf\rhd}\times_{\oU^{\rf\circ}} U^\rf\times_UV.
\end{equation}
On en déduit que le diagramme 
\begin{equation}\label{ktfr17k}
\xymatrix{
{\uoX'^\rhd_\fet}\ar[r]^{\theta'}\ar[d]_t&{\tuE'}\ar[d]^{\utau}\\
{\uX'_\fet\gtimes_{\uX_\fet}\uoX^\circ_\fet}\ar[r]^-(0.45){\nu}&{\tuG}}
\end{equation}
est commutatif à isomorphisme canonique près.

On désigne par $\rho_\uX\colon \uX_\et\rightarrow \uX_\fet$ et $\rho_{\uX'}\colon \uX'_\et\rightarrow \uX'_\fet$
les morphismes canoniques  \eqref{notconv10a}, par  $\epsilon \colon \uX_\fet\rightarrow \Ens$ et $\epsilon' \colon \uX'_\fet\rightarrow \Ens$
les projections canoniques et par
$\iota\colon \Ens\rightarrow \uX_\fet$ et $\iota'\colon \Ens\rightarrow \uX'_\fet$ les points $\rho_{\uX}(g(\ox'))$ et $\rho_{\uX'}(\ox')$, respectivement. 
On notera que $\epsilon$, $\epsilon'$, $\iota$ et $\iota'$ sont des  équivalences de topos.
Les morphismes de topos $\id_{\uoX^\circ_\fet}$ et $\iota' \epsilon \uh_\fet\colon \uoX^\circ_\fet\rightarrow \uX'_\fet$ et le $2$-isomorphisme 
\begin{equation}\label{ktfr17l}
\uh_\fet\stackrel{\sim}{\rightarrow} \ug_\fet \iota' \epsilon \uh_\fet
\end{equation}
induit par l'isomorphisme canonique $\id_{\uX_\fet}\stackrel{\sim}{\rightarrow} \iota \epsilon$ et la relation $\iota=\ug_\fet \iota'$, 
définissent un morphisme \eqref{tfr23b}
\begin{equation}\label{ktfr17m}
\mu\colon \uoX^\circ_\fet \rightarrow \uX'_\fet\gtimes_{\uX_\fet} \uoX^\circ_\fet
\end{equation}
qui s'insère dans un diagramme (non-commutatif)
\begin{equation}\label{ktfr17n}
\xymatrix{
{\Ens}\ar[d]_{\iota'}&{\uoX^\circ_\fet}\ar[d]_{\mu}\ar[l]_{\epsilon \uh_\fet}\ar[rd]^\id&\\
{\uX'_\fet}\ar[rd]_{\ug_\fet}&{\uX'_\fet\gtimes_{\uX_\fet}\uoX^\circ_\fet}\ar[l]_-(0.4){\rp_1}\ar[r]^-(0.4){\rp_2}&{\uoX^\circ_\fet}\ar[ld]^{\uh_\fet}\\
&{\uX_\fet}&}
\end{equation} 

Le carré et le triangle du diagramme 
\begin{equation}\label{ktfr17o}
\xymatrix{
&{\uX'_\fet}\ar[ld]_{\epsilon'}\ar[d]^{\ug_\fet}&{\uoX'^\rhd_\fet}\ar[d]^{\uupgamma_\fet}\ar[l]_{\uh'_\fet}\\
\Ens&{\uX_\fet}\ar[l]_{\epsilon}&{\uoX^\circ_\fet}\ar[l]_{\uh_\fet}}
\end{equation}
sont commutatifs à isomorphismes canoniques près. 
D'après la propriété universelle des produits orientés (\cite{agt} VI.3.7), 
l'isomorphisme canonique $\id_{\uX'_\fet}\stackrel{\sim}{\rightarrow} \iota' \epsilon'$ induit donc un isomorphisme 
\begin{equation}\label{ktfr17p}
t\stackrel{\sim}{\rightarrow} \mu \uupgamma_\fet.
\end{equation}
La proposition résulte alors de \eqref{ktfr17k} et \eqref{ktfr17p}.

(ii) Considérons le diagramme commutatif 
\begin{equation}\label{ktfr17iia}
\xymatrix{
{\uoX'^\rhd_\fet}\ar[r]^{\theta'}\ar[d]_{\uupgamma_\fet}&{\tuE'}\ar[d]^{\utau}\\
{\uoX^\circ_\fet}\ar[r]^{\vartheta}\ar@{=}[d]&{\tuG}\ar[d]^{\ulambda}\\
{\uoX^\circ_\fet}\ar@{=}[r]&{\uoX^\circ_\fet}}
\end{equation}
où $\ulambda$ est le morphisme canonique \eqref{tfr4d}. En vertu de \ref{tfr27}, le morphisme de changement de base 
\begin{equation}\label{ktfr17iib}
\ulambda_*\rightarrow \vartheta^*
\end{equation}
est un isomorphisme; en particulier, le foncteur $\ulambda_*$ est exact. Par suite, d'après (\cite{egr1}  1.2.4(v)), le diagramme 
\begin{equation}\label{ktfr17iic}
\xymatrix{
{\rR^q(\ulambda\utau)_*}\ar[rr]\ar@{=}[d]&&{\rR^q\uupgamma_{\fet*}\theta'^*}\ar@{=}[d]\\
{\ulambda_*\rR^q \utau_*}\ar[r]&{\vartheta^*\rR^q \utau_*}\ar[r]& {\rR^q\uupgamma_{\fet*}\theta'^*}}
\end{equation}
où les flèches horizontales sont les morphismes de changement de base relativement aux deux carrés et au rectangle extérieur du diagramme 
\eqref{ktfr17iia}, est commutatif. 

Considérons le diagramme commutatif 
\begin{equation}\label{ktfr17iid}
\xymatrix{
{\uoX'^\rhd_\fet}\ar[r]^{\theta'}\ar@{=}[d]&{\tuE'}\ar[d]^{\ubeta'}\\
{\uoX'^\rhd_\fet}\ar[d]_{\uupgamma_\fet}\ar@{=}[r]&{\uoX'^\rhd_\fet}\ar[d]^{\uupgamma_\fet}\\
{\uoX^\circ_\fet}\ar@{=}[r]&{\uoX^\circ_\fet}}
\end{equation}
où $\ubeta'$ est le morphisme canonique (\cite{agt} (VI.10.6.3)). En vertu de (\cite{agt} VI.10.27), le morphisme de changement de base 
\begin{equation}\label{ktfr17iie}
\ubeta'_*\rightarrow \theta'^*
\end{equation}
est un isomorphisme; en particulier, le foncteur $\ubeta'_*$ est exact. Par suite, d'après (\cite{egr1}  1.2.4(v)), le diagramme 
\begin{equation}\label{ktfr17iif}
\xymatrix{
{\rR^q(\uupgamma_\fet\ubeta')_*}\ar[r]\ar@{=}[d]&{(\rR^q\uupgamma_{\fet*})\theta'^*}\ar@{=}[d]\\
{(\rR^q\uupgamma_{\fet*})\ubeta'_*}\ar[r]&{(\rR^q\uupgamma_{\fet*})\theta'^*}}
\end{equation}
où les flèches horizontales sont les morphismes de changement de base relativement au carré supérieur et au rectangle extérieur du diagramme 
\eqref{ktfr17iid}, est commutatif. La flèche horizontale supérieure est donc un isomorphisme. 

Comme le diagramme 
\begin{equation}\label{ktfr17iig}
\xymatrix{
{\tuE'}\ar[r]^{\ubeta'}\ar[d]_{\utau}&{\uoX'^\rhd_\fet}\ar[d]^{\uupgamma_\fet}\\
{\tuG}\ar[r]^{\ulambda}&{\uoX^\circ_\fet}}
\end{equation}
est commutatif à isomorphisme canonique près \eqref{tfr13d}, on déduit de ce qui précède 
que la flèche horizontale supérieure du diagramme \eqref{ktfr17iic} est un isomorphisme.
La proposition s'ensuit puisque \eqref{ktfr17iib} est un isomorphisme.

\begin{cor}\label{ktfr18}
Conservons les hypothèses et notations de \ref{ktfr4} et \ref{ktfr15}. 
\begin{itemize}
\item[{\rm (i)}] Pour tout groupe abélien $F$ de $\tE'$ et tout entier $q\geq 0$, on a un isomorphisme canonique fonctoriel
\begin{equation}\label{ktfr18a}
\phi_{\ox'}(\rR^q\tau_*(F))\stackrel{\sim}{\rightarrow}\rR^q\uupgamma_{\fet*}(\varphi'_{\ox'}(F)).
\end{equation}
\item[{\rm (ii)}] Pour toute suite exacte de faisceaux abéliens 
$0\rightarrow F'\rightarrow F\rightarrow F''\rightarrow 0$ de $\tE'$ et tout entier $q\geq 0$, le diagramme 
\begin{equation}\label{ktfr18b}
\xymatrix{
{\phi_{\ox'}(\rR^q\tau_*(F''))} \ar[r]\ar[d]&{\phi_{\ox'}(\rR^{q+1}\tau_*(F'))}\ar[d]\\
{\rR^q\uupgamma_{\fet*}(\varphi_{\ox'}(F''))}\ar[r]&{\rR^{q+1}\uupgamma_{\fet*}(\varphi_{\ox'}(F'))}}
\end{equation}
où les flèches verticales sont les isomorphismes canoniques \eqref{ktfr18a} et les flèches horizontales  
sont les bords des suites exactes longues de cohomologie, est commutatif. 
\end{itemize}
\end{cor}

(i) Cela résulte aussitôt de \ref{ktfr16} et \ref{ktfr17}. 

(ii) Cela résulte de la fonctorialité de la version dérivée du morphisme de changement de base.  

\begin{cor}\label{ktfr19}
Soient $(\oy\rightsquigarrow \ox')$ un point de $X'_\et\gtimes_{X_\et}\oX^\circ_\et$ \eqref{tfr21},
$\uX'$ le localisé strict de $X'$ en $\ox'$, $\uX$ le localisé strict de $X$ en $g(\ox')$,
$\uupgamma\colon \uoX'^\rhd\rightarrow \uoX^\circ$ le morphisme induit par $g$, 
$\rho_{\uoX^\circ}\colon \uoX^\circ_\et\rightarrow \uoX^\circ_\fet$ le morphisme canonique \eqref{notconv10a}. 
On considère $\oy$ comme un point géométrique de $\uoX^\circ$ par  le $\oX^\circ$-morphisme
$\oy\rightarrow \uoX^\circ$ induit $(\oy\rightsquigarrow \ox')$.
\begin{itemize}
\item[{\rm (i)}] Pour tout groupe abélien $F$ de $\tE'$ et tout entier $q\geq 0$, on a un isomorphisme canonique fonctoriel
\begin{equation}\label{ktfr19a}
(\rR^q\tau_*(F))_{\varrho(\oy\rightsquigarrow \ox')}\stackrel{\sim}{\rightarrow}
\rR^q\uupgamma_{\fet*}(\varphi'_{\ox'}(F))_{\rho_{\uoX^\circ}(\oy)}.
\end{equation}
\item[{\rm (ii)}] Pour toute suite exacte de faisceaux abéliens 
$0\rightarrow F'\rightarrow F\rightarrow F''\rightarrow 0$ de $\tE'$ et tout entier $q\geq 0$, le diagramme 
\begin{equation}\label{ktfr19b}
\xymatrix{
{(\rR^q\tau_*(F''))_{\varrho(\oy\rightsquigarrow \ox')}} \ar[r]\ar[d]&{(\rR^{q+1}\tau_*(F'))_{\varrho(\oy\rightsquigarrow \ox')}}\ar[d]\\
{\rR^q\uupgamma_{\fet*}(\varphi_{\ox'}(F''))_{\rho_{\uoX^\circ}(\oy)}}\ar[r]&
{\rR^{q+1}\uupgamma_{\fet*}(\varphi_{\ox'}(F'))_{\rho_{\uoX^\circ}(\oy)}}}
\end{equation}
où les flèches verticales sont les isomorphismes canoniques \eqref{ktfr19a} et les flèches horizontales  
sont les bords des suites exactes longues de cohomologie, est commutatif. 
\end{itemize}
\end{cor}

Cela résulte de \ref{ktfr18} et \eqref{ktfr28a}.

\subsection{}\label{ktfr2}
Pour tout objet $(U'\rightarrow U\leftarrow V)$ de $G$, on note $\oU'^V$ la fermeture intégrale de $\oU'$ dans $U'\times_UV$. 
\begin{equation}\label{ktfr2a}
\xymatrix{
&\oU'&{\oU'^V}\ar[l]\\
U'\ar[d]&{\oU'^\circ}\ar[d]\ar[l]\ar[u]&{U'\times_UV}\ar[l]\ar[d]\ar[u]\\
U & {\oU^\circ}\ar[l] & V\ar[l]}
\end{equation}
On désigne par $\ocB^!$ le préfaisceau sur $G$ défini pour tout $(U'\rightarrow U\leftarrow V)\in \ob(G)$ par
\begin{equation}\label{ktfr2b}
\ocB^!(U'\rightarrow U\leftarrow V)=\Gamma(\oU'^V,\co_{\oU'^V}).
\end{equation}
C'est un faisceau pour la topologie co-évanescente de $G$ d'après \ref{tfra3}. 
On notera que la définition de $\ocB^!$ correspond à celle de $\ocB$ dans \ref{tfra1} en prenant pour \eqref{tfra1a}
le composé des morphismes canoniques $\oX'^\circ=X'\times_X\oX^\circ\rightarrow \oX'\rightarrow X'$; 
et on observera que $\oX'$ est normal et localement irréductible d'après (\cite{agt} III.4.2(iii)).

Compte tenu de \ref{tfra8}, on a un homomorphisme canonique \eqref{tfra5f} 
\begin{equation}\label{ktfr2c}
\ocB\rightarrow \lgg_*(\ocB^!),
\end{equation}
où $\ocB$ est l'anneau de $\tE$ défini dans \eqref{mtfla8a}.
De même, on a un homomorphisme canonique 
\begin{equation}\label{ktfr2d}
\ocB^!\rightarrow\tau_*(\ocB'),
\end{equation}
où $\ocB'$ est l'anneau de $\tE'$ défini dans \eqref{mtfla10a}. 
On a enfin un homomorphisme canonique \eqref{tfra1e}
\begin{equation}\label{ktfr2e}
\hbar'_*(\co_{\oX'})\rightarrow \pi_*(\ocB^!).
\end{equation}

\begin{lem}\label{ktfr3}
L'homomorphisme canonique $\ocB^!\rightarrow\tau_*(\ocB')$ \eqref{ktfr2d} est un isomorphisme.
\end{lem}

En effet, comme $\oX'$ est normal et localement irréductible  (\cite{agt} III.4.2(iii)), pour tout objet $(U'\rightarrow U\leftarrow V)$ de $G$,
les schémas $\oU'$ et $U'\times_UV$ sont normaux et localement irréductibles (\cite{agt} III.3.3) et il en est de même de $\oU'^V$ \eqref{tfra2}. 
Par ailleurs, l'immersion $X'^\rhd\rightarrow X'$ étant schématiquement dominante (\cite{agt} III.4.2(iv)), il en est de même des immersions 
$V\times_{X^\circ}X'^\rhd\rightarrow V\times_XX'$ et $U'\times_{(U\times_XX')}(V\times_{X^\circ}X'^\rhd)\rightarrow U'\times_UV$ (\cite{ega4} 11.10.5).
Par suite, $\oU'^V$ est la fermeture intégrale de $\oU'$ dans $U'\times_{(U\times_XX')}(V\times_{X^\circ}X'^\rhd)$, d'où la proposition. 
On notera que le morphisme canonique 
\begin{equation}\label{ktfr3a}
\oU'^\rhd\times_{\oU^\circ}V\rightarrow U'\times_{(U\times_XX')}(V\times_{X^\circ}X'^\rhd)
\end{equation}
est un isomorphisme. 

\begin{prop}\label{ktfr5}
Soient $(\oy\rightsquigarrow \ox')$ un point de $X'_\et\gtimes_{X_\et}\oX^\circ_\et$ \eqref{tfr21} tel que $\ox'$ soit au-dessus de $s$, 
$\varrho(\oy\rightsquigarrow \ox')$ son image par  le morphisme $\varrho$ \eqref{ktfr1e}, qui est donc un point de $\tG$,
$\uX'$ le localisé strict de $X'$ en $\ox'$. Alors, 
\begin{itemize}
\item[{\rm (i)}] La fibre $\ocB^!_{\varrho(\oy \rightsquigarrow \ox')}$ de $\ocB^!$ en $\varrho(\oy\rightsquigarrow \ox')$ 
est un anneau normal et strictement local. 
\item[{\rm (ii)}] On a un isomorphisme canonique
\begin{equation}\label{ktfr5a}
(\hbar'_*(\co_{\oX'}))_{\ox'}\stackrel{\sim}{\rightarrow}\Gamma(\uoX',\co_{\uoX'}).
\end{equation}
\item[{\rm (iii)}]  L'homomorphisme 
\begin{equation}\label{ktfr5b}
(\hbar'_*(\co_{\oX'}))_{\ox'}\rightarrow \ocB^!_{\varrho(\oy \rightsquigarrow \ox')}
\end{equation} 
induit par l'homomorphisme canonique $\pi^*(\hbar'_*(\co_{\oX'}))\rightarrow \ocB^!$ \eqref{ktfr2e} est injectif et local.
\end{itemize}
\end{prop}

On rappelle que les schémas $X$ et $\oX'$ sont normaux et localement irréductibles d'après (\cite{agt} III.4.2(iii)).
Le corps résiduel de $\co_K$ étant algébriquement clos \eqref{mtfla1}, on peut identifier $\ox'$ à un point géométrique de $\oX'$.
D'après (\cite{agt} III.3.7), le schéma $\uoX'$ est normal et strictement local. Il s'identifie donc au localisé strict de $\oX'$ en $\ox'$.  
De même, notant $\uX$ le localisé strict de $X$ en $g(\ox')$, 
le schéma $\uoX$ est normal et strictement local, et il s'identifie au localisé strict de $\oX$ en $g(\ox')$. 

On désigne par $\tuG$ le topos de Faltings relatif associé au couple de morphismes $(\uoX^\circ\rightarrow \uX, \uX'\rightarrow \uX)$,  par 
\begin{equation}\label{ktfr5c}
\Phi\colon \tuG\rightarrow \tG
\end{equation}
le morphisme de fonctorialité \eqref{tfr7b}, et par 
\begin{eqnarray}
\upi\colon \tuG&\rightarrow& \uX'_\et,\label{ktfr5d}\\
\uvarrho\colon \uX'_\et\gtimes_{\uX_\et}\uoX^\circ_\et &\rightarrow& \tuG,\label{ktfr5dd}
\end{eqnarray}
les morphismes canoniques \eqref{tfr4c} et \eqref{tfr6b}. 
On note $\uocB^!$ l'anneau de $\tuG$ défini dans \ref{tfra1} par la factorisation canonique
\begin{equation}\label{ktfr5e}
\uoX'^\circ\longrightarrow \uoX'\stackrel{\uhbar'}{\longrightarrow} \uX'.
\end{equation}
D'après \ref{tfra9}, on a un isomorphisme canonique
\begin{equation}\label{ktfr5f}
\Phi^{-1}(\ocB^!)\stackrel{\sim}{\rightarrow} \uocB^!. 
\end{equation}

En vertu de (\cite{kato1} 4.5, \cite{tsuji4} II.4.2 et \cite{gr1} 12.7.8(iii)), le morphisme $g\colon X'\rightarrow X$ est plat à fibres géométriquement réduites. 
Par suite, d'après \ref{csc2} appliqué au morphisme $\oX'\rightarrow \oX$ induit par $g$,
pour tout revêtement étale $V\rightarrow \uoX^\circ$ tel que $V$ soit connexe, le schéma $V\times_{\uX}\uX'$ est connexe. 
Considérons $(\oy\rightsquigarrow \ox')$ aussi comme un point de $\uX'_\et\gtimes_{\uX_\et}\uoX^\circ_\et$.
En vertu de  \ref{tfra4}, la fibre $\uocB^!_{\uvarrho(\oy \rightsquigarrow \ox')}$ de $\uocB^!$ en $\uvarrho(\oy\rightsquigarrow \ox')$ 
est donc un anneau normal et strictement local et l'homomorphisme 
\begin{equation}\label{ktfr5g}
(\uhbar'_*(\co_{\uoX'}))_{\ox'}\rightarrow \uocB^!_{\uvarrho(\oy \rightsquigarrow \ox')}
\end{equation} 
induit par l'homomorphisme canonique $\upi^*(\uhbar'_*(\co_{\uoX'}))\rightarrow \uocB^!$ \eqref{ktfr2e} est injectif et local. 
Les propositions (i) et (iii) s'ensuivent compte tenu de \eqref{ktfr28b} et \eqref{ktfr5f}. 
La proposition (ii) résulte aussitôt de \eqref{notconv12g}.

\subsection{}\label{ktfr6}
Considérons un diagramme commutatif de morphismes de schémas
\begin{equation}\label{ktfr6a}
\xymatrix{
U'\ar[r]\ar[d]&U\ar[d]\\
X'\ar[r]&X}
\end{equation}
tel que les flèches verticales soient des morphismes étales et soit $\oy$ un point géométrique de $\oU^\circ$. 
Le schéma $\oU$ étant localement irréductible d'après (\cite{agt} III.3.3 et III.4.2(iii)),  
il est la somme des schémas induits sur ses composantes irréductibles. On note $\oU^\star$
la composante irréductible de $\oU$ contenant $\oy$. 
De même, $\oU^\circ$ est la somme des schémas induits sur ses composantes irréductibles
et $\oU^{\star \circ}=\oU^\star\times_{X}X^\circ$ est la composante irréductible de $\oU^\circ$ contenant $\oy$. 
On note $\bB_{\pi_1(\oU^{\star \circ},\oy)}$ le topos classifiant du groupe profini $\pi_1(\oU^{\star \circ},\oy)$ et
\begin{equation}\label{ktfr6b}
\psi_\oy\colon \oU^{\star \circ}_\fet \stackrel{\sim}{\rightarrow}\bB_{\pi_1(\oU^{\star \circ},\oy)}
\end{equation}
le foncteur fibre  de $\oU^{\star \circ}_\fet$ en $\oy$ \eqref{notconv11c}.

Considérons l'anneau $\ocB^!_{U'\rightarrow U}$ de $\oU^{\circ}_\fet $ défini dans \eqref{ktfr12b} et posons
\begin{equation}\label{ktfr6d}
\oR^{!\oy}_{U'\rightarrow U}=\psi_\oy(\ocB^!_{U'\rightarrow U}|\oU^{\star \circ}).
\end{equation}
Explicitement, soit $(V_i)_{i\in I}$ le revêtement universel normalisé de $\oU^{\star \circ}$ en $\oy$ \eqref{notconv11}.
Pour chaque $i\in I$, $(U'\rightarrow U\leftarrow V_i)$ est naturellement un objet de $G$. On a  alors
\begin{equation}\label{ktfr6e}
\oR^{!\oy}_{U'\rightarrow U}=\underset{\underset{i\in I}{\longrightarrow}}{\lim}\ \ocB^!(U'\rightarrow U\leftarrow V_i).
\end{equation} 

Considérons l'anneau $\ocB_{U}$ de $\oU^{\circ}_\fet $ défini dans \eqref{TFA2d} et posons
\begin{equation}\label{ktfr6f}
\oR^{\oy}_{U}=\psi_\oy(\ocB_{U}|\oU^{\star \circ}).
\end{equation}
Pour chaque $i\in I$, $(V_i\rightarrow U)$ est naturellement un objet de $E$. On a  alors
\begin{equation}\label{ktfr6g}
\oR^{\oy}_{U}=\underset{\underset{i\in I}{\longrightarrow}}{\lim}\ \ocB(V_i\rightarrow U).
\end{equation} 

L'homomorphisme canonique $\lgg^*(\ocB)\rightarrow \ocB^!$ \eqref{ktfr2c} induit pour tout $i\in I$, un morphisme (fonctoriel en $i$)
\begin{equation}\label{ktfr6h}
\ocB(V_i\rightarrow U)\rightarrow \ocB^!(U'\rightarrow U\leftarrow V_i).
\end{equation}
On en déduit par passage à la limite inductive un homomorphisme 
\begin{equation}\label{ktfr6i}
\oR^{\oy}_{U}\rightarrow \oR^{!\oy}_{U'\rightarrow U}.
\end{equation}

\subsection{}\label{ktfr7}
Conservons les hypothèses et notations de \ref{ktfr6}, supposons de plus qu'il existe un point géométrique $\oy'$ de $\oU'^\rhd$ au-dessus 
du point géométrique $\oy$ de $\oU^\circ$. 
Le schéma $\oU'$ étant normal et localement irréductible, il est la somme des schémas induits sur ses composantes irréductibles. 
On note $\oU'^\star$ la composante irréductible de $\oU'$ contenant $\oy'$. 
De même, $\oU'^\rhd$ est la somme des schémas induits sur ses composantes irréductibles
et $\oU'^{\star \rhd}=\oU'^\star\times_{X'}X'^\rhd$ est la composante irréductible de $\oU'^\rhd$ contenant $\oy'$. 
On note $\bB_{\pi_1(\oU'^{\rhd \circ},\oy')}$ le topos classifiant du groupe profini $\pi_1(\oU'^{\star \rhd},\oy')$ et
\begin{equation}\label{ktfr7a}
\psi'_{\oy'}\colon \oU'^{\star \rhd}_\fet \stackrel{\sim}{\rightarrow}\bB_{\pi_1(\oU'^{\star \rhd},\oy')}
\end{equation}
le foncteur fibre  de $\oU'^{\star \rhd}_\fet$ en $\oy'$ \eqref{notconv11c}. On pose 
\begin{equation}\label{ktfr7b}
\oR'^{\oy'}_{U'}=\psi'_{\oy'}(\ocB'_{U'}|\oU'^{\star \rhd}).
\end{equation}
Explicitement, soit $(W_j)_{j\in J}$ le revêtement universel normalisé de $\oU'^{\star \rhd}$ en $\oy'$ \eqref{notconv11}. 
Pour chaque $j\in J$, $(W_j\rightarrow U')$ est naturellement un objet de $E'$. On a alors 
\begin{equation}\label{ktfr7c}
\oR'^{\oy'}_{U'}=\underset{\underset{j\in J}{\longrightarrow}}{\lim}\  \ocB'(W_j\rightarrow U').
\end{equation}

Pour tout $i\in I$, on a un $\oU^\circ$-morphisme canonique $\oy\rightarrow V_i$.
On en déduit un $\oU'^\rhd$-morphisme $\oy'\rightarrow V_i\times_{\oU^\circ}\oU'^\rhd$.
Le schéma $V_i\times_{\oU^\circ}\oU'^\rhd$ étant localement irréductible,  
il est la somme des schémas induits sur ses composantes irréductibles. 
On note $V'_i$ la composante irréductible de $V_i\times_{\oU^\circ}\oU'^\rhd$ contenant l'image de $\oy'$. 
Les schémas $(V'_i)_{i\in I}$ forment naturellement un système projectif de revêtements étales finis connexes $\oy'$-pointés de $\oU'^{\star \rhd}$.
Compte tenu de l'isomorphisme \eqref{ktfr3a}, on a un isomorphisme canonique de $\tE'$ 
\begin{equation}\label{ktfr7d}
\tau^*((U'\rightarrow U\leftarrow V_i)^a)\stackrel{\sim}{\rightarrow} (V_i\times_{\oU^\circ}\oU'^\rhd\rightarrow U')^a.
\end{equation}
D'après \ref{ktfr3}, on en déduit un homomorphisme canonique (fonctoriel en $i$)
\begin{equation}\label{ktfr7e}
\ocB^!(U'\rightarrow U\leftarrow V_i)=\ocB'(V_i\times_{\oU^\circ}\oU'^\rhd\rightarrow U')\rightarrow \ocB'(V'_i\rightarrow U').
\end{equation}
Posant
\begin{equation}\label{ktfr7f}
\oR^{\intern\oy'}_{U'\rightarrow U}=\underset{\underset{i\in I}{\longrightarrow}}{\lim}\  \ocB'(V'_i\rightarrow U'),
\end{equation}
on a donc un homomorphisme canonique
\begin{equation}\label{ktfr7g}
\oR^{!\oy}_{U'\rightarrow U}\rightarrow \oR^{\intern\oy'}_{U'\rightarrow U}.
\end{equation}

Pour tous $i\in I$ et $j\in J$, il existe au plus un morphisme de $\oU'^{\star \rhd}$-schémas pointés $W_j\rightarrow V'_i$. 
De plus, pour tout $i\in I$, il existe $j\in J$ et un morphisme de $\oU'^{\star \rhd}$-schémas pointés $W_j\rightarrow V'_i$. 
On a donc un homomorphisme
\begin{equation}\label{ktfr7h}
\oR^{\intern\oy'}_{U'\rightarrow U}=\underset{\underset{i\in I}{\longrightarrow}}{\lim}\  \ocB'(V'_i\rightarrow U')\rightarrow 
\underset{\underset{j\in J}{\longrightarrow}}{\lim}\  \ocB'(W_j\rightarrow U')=\oR'^{\oy'}_{U'}. 
\end{equation}

On dispose donc de trois homomorphismes canoniques 
\begin{equation}\label{ktfr7i}
\oR^\oy_U\rightarrow \oR^{!\oy}_{U'\rightarrow U} \rightarrow \oR^{\intern\oy'}_{U'\rightarrow U}\rightarrow \oR'^{\oy'}_{U'}.
\end{equation}

\begin{rema}\label{ktfr8}
Conservons les hypothèses et notations de \ref{ktfr6} et \ref{ktfr7}. 
Pour tout $i\in I$, le schéma $\oU'^{V_i}$ est normal et localement irréductible \eqref{tfra2}. 
Il est donc la somme des schémas induits sur ses composantes irréductibles, 
et $\oU'^{V'_i}$ est la composante irréductible contenant l'image canonique de $\oy'$. 
Si, de plus, le schéma $U'$ est affine, l'homomorphisme \eqref{ktfr7e} 
\begin{equation}\label{ktfr8a}
\ocB^!(U'\rightarrow U\leftarrow V_i)\rightarrow \ocB'(V'_i\rightarrow U')
\end{equation}
est donc surjectif, et son noyau est engendré par un idempotent de $\ocB^!(U'\rightarrow U\leftarrow V_i)$. 
Par ailleurs, pour tous $(i,j)\in I^2$ avec  
$i\geq j$, le morphisme canonique $\oU'^{V'_i}\rightarrow \oU'^{V'_j}$ est entier et dominant. 
On en déduit que l'anneau $\oR^{\intern\oy'}_{U'\rightarrow U}$ est intègre et normal (\cite{ega1n} 0.6.1.6(i) et 0.6.5.12(ii)),  
et que l'homomorphisme \eqref{ktfr7g} 
\begin{equation}\label{ktfr8b}
\oR^{!\oy}_{U'\rightarrow U}\rightarrow \oR^{\intern\oy'}_{U'\rightarrow U}
\end{equation}
est surjectif et son noyau engendré par des idempotents de $\oR^{!\oy}_{U'\rightarrow U}$.

On observera que l'anneau $\oR^{!\oy}_{U'\rightarrow U}$ \eqref{ktfr6e} (resp.  $\oR^{\intern\oy'}_{U'\rightarrow U}$ \eqref{ktfr7f})
correspond à l'anneau noté $\oR^!$ dans \eqref{eccr40g} (resp. $\oR^\intern$ dans \eqref{eccr40e}).
\end{rema}

\begin{rema} \label{ktfr9}
D'après \eqref{ktfr14e}, pour tout point $(\oy\rightsquigarrow \ox')$ de $X'_\et\gtimes_{X_\et}\oX^\circ_\et$, on a un isomorphisme canonique
\begin{equation}\label{ktfr9a}
\ocB^!_{\varrho(\oy\rightsquigarrow \ox')}\stackrel{\sim}{\rightarrow} \underset{\underset{(\ox'\rightarrow U'\rightarrow U)\in \cC^\circ_{\ox'}}{\longrightarrow}}{\lim}\ \oR^{!\oy}_{U'\rightarrow U},
\end{equation}
où $\cC_{\ox'}$ est la catégorie définie dans \ref{ktfr13} et $\oR^{!\oy}_{U'\rightarrow U}$ est l'anneau \eqref{ktfr6d} qui est clairement fonctoriel 
sur $\cC^\circ_{\ox'}$.
\end{rema}

\subsection{}\label{ktfr10}
Soient $(\oy'\rightsquigarrow \ox')$ un point de $X'_\et\gtimes_{X'_\et}\oX'^\rhd_\et$, $\uX'$ le localisé strict de $X'$ en $\ox'$.
La donnée de $(\oy'\rightsquigarrow \ox')$ est équivalente à la donnée 
d'un $X'$-morphisme $u'\colon \oy'\rightarrow \uX'$. 
On en déduit un $\oX'^\rhd$-morphisme $v'\colon \oy'\rightarrow \uoX'^\rhd$  \eqref{mtfla2c}. 
On désigne par $\fV'_{\ox'}$ la catégorie des $X'$-schémas étales $\ox'$-pointés. 
Pour tout objet $(\ox'\rightarrow U')$ de $\fV'_{\ox'}$, on a un $X'$-morphisme canonique $\uX'\rightarrow U'$.
On en déduit un morphisme $\uoX'\rightarrow \oU'$. Le morphisme $v'\colon \oy'\rightarrow \uoX'^\rhd$ induit 
alors un point géométrique de $\oU'^\rhd$ que l'on note encore $\oy'$. 

La correspondance qui à tout objet $(\ox'\rightarrow U')$ de $\fV'^\circ_{\ox'}$ associe l'anneau $\oR'^{\oy'}_{U'}$ \eqref{ktfr7b} est fonctorielle. 
On pose 
\begin{equation}\label{ktfr10a}
\oR'^{\oy'}_{\uX'}=\underset{\underset{(\ox'\rightarrow U')\in \fV'^\circ_{\ox'}}{\longrightarrow}}{\lim}\ \oR'^{\oy'}_{U'}. 
\end{equation}

Posons $\ox=g(\ox')$ et $\oy=\upgamma(\oy')$ \eqref{ktfr1b} qui sont donc des points géométriques de $X$ et $\oX^\circ$ respectivement,
et notons $(\oy\rightsquigarrow \ox')$ l'image de $(\oy'\rightsquigarrow \ox')$ par le morphisme canonique 
\begin{equation}\label{ktfr10b}
X'_\et\gtimes_{X'_\et}\oX'^\rhd_\et\rightarrow X'_\et\gtimes_{X_\et}\oX^\circ_\et.
\end{equation}
On désigne par $\uX$ le localisé strict de $X$ en $\ox$ et par $\fV_{\ox}$ la catégorie des $X$-schémas étales $\ox$-pointés.
Le point $(\oy\rightsquigarrow \ox')$ définit pour tout objet $(\ox\rightarrow U)$ de $\fV^\circ_{\ox}$, un point géométrique 
de $\oU^\circ$ que l'on note encore $\oy$. 
La correspondance qui à tout objet $(\ox\rightarrow U)$ de $\fV^\circ_{\ox}$ associe l'anneau $\oR^{\oy}_{U}$ \eqref{ktfr6g} est fonctorielle. 
On pose 
\begin{equation}\label{ktfr10c}
\oR^{\oy}_{\uX}=\underset{\underset{(\ox\rightarrow U)\in \fV^\circ_{\ox}}{\longrightarrow}}{\lim}\ \oR^{\oy}_U. 
\end{equation}

Considérons la catégorie $\cC_{\ox'}$ définie dans \ref{ktfr13}.  
La correspondance qui à tout objet $(\ox'\rightarrow U'\rightarrow U)$ de $\cC^\circ_{\ox'}$ associe l'anneau 
$\oR^{!\oy}_{U'\rightarrow U}$ \eqref{ktfr6d} est fonctorielle.  On pose
\begin{equation}\label{ktfr10d}
\oR^{!\oy}_{\uX'\rightarrow \uX}=\underset{\underset{(\ox'\rightarrow U'\rightarrow U)\in \cC^\circ_{\ox'}}{\longrightarrow}}{\lim}\ \oR^{!\oy}_{U'\rightarrow U}.
\end{equation}
De même, la correspondance qui à tout objet $(\ox'\rightarrow U'\rightarrow U)$ de $\cC^\circ_{\ox'}$ associe l'anneau 
$\oR^{\intern\oy'}_{U'\rightarrow U}$ \eqref{ktfr7h} est fonctorielle.  On pose
\begin{equation}\label{ktfr10e}
\oR^{\intern\oy'}_{\uX'\rightarrow \uX}=\underset{\underset{(\ox'\rightarrow U'\rightarrow U)\in \cC^\circ_{\ox'}}{\longrightarrow}}{\lim}\ \oR^{\intern\oy'}_{U'\rightarrow U}.
\end{equation}

Pour tout objet $(\ox'\rightarrow U'\rightarrow U)$ de $\cC_{\ox'}$, le diagramme 
\begin{equation}\label{ktfr10f}
\xymatrix{
\oy'\ar[r]\ar[d]&\oy\ar[d]\\
\oU'^\rhd\ar[r]&\oU^\circ}
\end{equation}
est commutatif. On a donc trois homomorphismes canoniques fonctoriels \eqref{ktfr7i}
\begin{equation}\label{ktfr10g}
\oR^\oy_U\rightarrow \oR^{!\oy}_{U'\rightarrow U} \rightarrow \oR^{\intern\oy'}_{U'\rightarrow U}\rightarrow \oR'^{\oy'}_{U'}.
\end{equation}
On en déduit par passage à la limite trois homomorphismes 
\begin{equation}\label{ktfr10h}
\oR^\oy_{\uX}\rightarrow \oR^{!\oy}_{\uX'\rightarrow \uX} \rightarrow \oR^{\intern\oy'}_{\uX'\rightarrow \uX}\rightarrow \oR'^{\oy'}_{\uX'}.
\end{equation}
Compte tenu de \eqref{TFA11c}, \eqref{ktfr9a} et \eqref{ktfr1f}, le premier homomorphisme et le composé des deux autres s'identifient aux homomorphismes 
\begin{equation}
\ocB_{\rho(\oy\rightsquigarrow \ox)}\rightarrow \ocB^!_{\varrho(\oy\rightsquigarrow \ox')}\rightarrow \ocB'_{\rho'(\oy'\rightsquigarrow \ox')}
\end{equation}
induits par les adjoints des homomorphismes \eqref{ktfr2c} et \eqref{ktfr2d}, respectivement.

\begin{prop}\label{ktfr11}
Conservons les hypothèses de \ref{ktfr10}, supposons de plus que $\ox'$ soit au-dessus de $s$. Alors, l'homomorphisme \eqref{ktfr10h}
\begin{equation}
\oR^{!\oy}_{\uX'\rightarrow \uX} \rightarrow \oR^{\intern\oy'}_{\uX'\rightarrow \uX}
\end{equation}
est un isomorphisme. 
\end{prop}

En effet, pour tout objet $(\ox'\rightarrow U'\rightarrow U)$ de $\cC_{\ox'}$, on a un diagramme commutatif d'homomorphismes d'anneaux 
(sans la flèche pointillée) 
\begin{equation}
\xymatrix{
{\oR^{!\oy}_{U'\rightarrow U}}\ar@{->>}[r]^-(0.4)u\ar[d]_\iota&{\oR^{\intern\oy'}_{U'\rightarrow U}}\ar[r]\ar[d]^\jmath\ar@{.>}[ld]_w&{\oR'^{\oy'}_{U'}}\ar[d]\\
{\oR^{!\oy}_{\uX'\rightarrow \uX}}\ar@{->>}[r]^-(0.4)v\ar@{=}[d]&{\oR^{\intern\oy'}_{\uX'\rightarrow \uX}}\ar[r]&{\oR'^{\oy'}_{\uX'}}\ar@{=}[d]\\
{\ocB^!_{\varrho(\oy\rightsquigarrow \ox')}}\ar[rr]&&{\ocB'_{\rho'(\oy'\rightsquigarrow \ox')}}}
\end{equation}
D'après \ref{ktfr8},  l'homomorphisme $u$ est surjectif et son noyau est engendré par des idempotents de $\oR^{!\oy}_{U'\rightarrow U}$.
L'homomorphisme $v$ est donc surjectif. Montrons que $\iota$ se factorise à travers $u$. 
Il suffit de montrer que pour tout idempotent $e$ de $\oR^{!\oy}_{U'\rightarrow U}$ tel que $u(e)=0$, on a $\iota(e)=0$. 
En vertu de \ref{ktfr5}, l'anneau $\oR^{!\oy}_{\uX'\rightarrow \uX}$ est normal et strictement local (et donc intègre).
Si $\iota(e)$ n'était pas nul, on aurait $\iota(e)=1$. 
Par suite, l'image de $e$ dans $\oR'^{\oy'}_{\uX'}$ vaudrait également $1$. Mais elle vaut aussi $0$ car $u(e)=0$. 
Par suite, l'anneau $\oR'^{\oy'}_{\uX'}$ serait nul, ce qui est absurde puisqu'il est intègre d'après (\cite{agt} III.10.10(i)). 
Il existe donc une flèche pointillée $w$ telle que $\iota=w\circ u$. Comme $u$ est surjectif, on en déduit que $\jmath=v\circ w$. 
De plus, les homomorphismes $w$ pour $(\ox'\rightarrow U'\rightarrow U)\in \ob(\cC^\circ_{\ox'})$ forment un système compatible, 
dont la limite inductive est une section $v'$ de $v$ qui est surjective compte tenu de la relation  $\iota=w\circ u$; d'où la proposition.

\begin{lem}\label{ktfr23}
Pour tout $n\geq 0$, l'anneau $\ocB^!/p^n\ocB^!$ est un objet de $\tG_s$.
\end{lem}

En effet, pour tout point $(\oy\rightsquigarrow \ox')$ de $X'_\et\gtimes_{X_\et}\oX^\circ_\et$ \eqref{tfr21}
tel que $\ox'$ soit au-dessus de $\eta$, l'homomorphisme canonique $\pi^{-1}(\co_{X'})\rightarrow \ocB^!$ \eqref{ktfr2e}
induit un homomorphisme $\co_{X',\ox'}\rightarrow \ocB^!_{\varrho(\oy\rightsquigarrow \ox')}$ \eqref{tfr21a}.
Par suite, $p$ est inversible dans $\ocB^!_{\varrho(\oy\rightsquigarrow \ox')}$, d'où la proposition en vertu de \ref{ktfr22}.

\subsection{}\label{ktfr24} 
Pour tout entier $n\geq 1$, on pose 
\begin{equation}\label{ktfr24a}
\ocB^!_n=\ocB^!/p^n\ocB^!.
\end{equation}
On note $a'\colon X'_s\rightarrow X'$, $a'_n\colon X'_s\rightarrow X'_n$ \eqref{TFA1b}, 
$\iota'_n\colon X'_n\rightarrow X'$
et $\oiota'_n\colon \oX'_n\rightarrow \oX'$ les injections canoniques. 
Le corps résiduel de $\co_K$ étant algébriquement clos, il existe un unique $S$-morphisme $s\rightarrow \oS$. 
Celui-ci induit des immersions fermées $\oa'\colon X'_s\rightarrow \oX'$ et $\oa'_n\colon X'_s\rightarrow \oX'_n$
qui relèvent $a'$ et $a'_n$, respectivement. 
\begin{equation}\label{ktfr24b}
\xymatrix{
{X'_s}\ar[r]_{\oa_n}\ar@{=}[d]\ar@/^1pc/[rr]^{\oa'}&{\oX'_n}\ar[r]_{\oiota'_n}\ar[d]^{\hbar'_n}&\oX'\ar[d]^{\hbar'}\\
{X'_s}\ar[r]^{a'_n}\ar@/_1pc/[rr]_{a'}&{X'_n}\ar[r]^{\iota'_n}&X'}
\end{equation}
L'homomorphisme canonique $\pi^{-1}(\hbar_*(\co_{\oX'}))\rightarrow \ocB^!$ \eqref{ktfr2e} induit un homomorphisme 
\begin{equation}\label{ktfr24c}
\pi^{-1}(\iota'_{n*}(\hbar'_{n*}(\co_{\oX'_n})))\rightarrow \ocB^!_n.
\end{equation}
Comme $\oa'_n$ est un homéomorphisme universel, on peut considérer $\co_{\oX'_n}$ comme un faisceau de $X'_{s,\et}$. 
On peut alors identifier les anneaux $\iota'_{n*}(\hbar'_{n*}(\co_{\oX'_n}))$ et $a'_*(\co_{\oX'_n})$ 
de $X'_\et$ et par suite les anneaux $\pi^{-1}(\iota'_{n*}(\hbar'_{n*}(\co_{\oX'_n})))$ et 
$\kappa_*(\pi^*_s(\co_{\oX'_n}))$ de $\tG_s$ \eqref{ktfr21j}. 
Comme $\ocB^!_n$ est un objet de $\tG_s$ \eqref{ktfr23}, on peut considérer 
\eqref{ktfr24c} comme un homomorphisme de $\tG_s$ 
\begin{equation}\label{ktfr24d}
\pi_s^*(\co_{\oX'_n})\rightarrow \ocB^!_n.
\end{equation}
Le morphisme $\pi_s$ \eqref{ktfr21a} est donc sous-jacent à un morphisme de topos annelés, que l'on note 
\begin{equation}\label{ktfr24e}
\pi_n\colon (\tG_s,\ocB^!_n)\rightarrow (\oX'_{s,\et},\co_{\oX'_n}).
\end{equation}
Nous utilisons pour les $\co_{\oX'_n}$ modules la notation $\pi_s^{-1}$ pour désigner l'image
inverse au sens des faisceaux abéliens et nous réservons la notation 
$\pi_n^*$ pour l'image inverse au sens des modules. 

L'homomorphisme canonique $\tau^{-1}(\ocB^!)\rightarrow \ocB'$ \eqref{ktfr2d}
induit un homomorphisme $\tau_s^*(\ocB^!_n)\rightarrow \ocB'_n$. 
Le morphisme $\tau_s$ \eqref{ktfr21c} est donc sous-jacent à un morphisme de topos annelés, que l'on note 
\begin{equation}\label{ktfr24f}
\tau_n\colon (\tE'_s,\ocB'_n)\rightarrow (\tG_s,\ocB^!_n).
\end{equation}

L'homomorphisme canonique $\lgg^{-1}(\ocB)\rightarrow \ocB^!$ \eqref{ktfr2c}
induit un homomorphisme $\lgg_s^*(\ocB_n)\rightarrow \ocB^!_n$. 
Le morphisme $\lgg_s$ \eqref{ktfr21f} est donc sous-jacent à un morphisme de topos annelés, que l'on note 
\begin{equation}\label{ktfr24g}
\lgg_n\colon (\tG_s,\ocB^!_n)\rightarrow (\tE_s,\ocB_n).
\end{equation}

On vérifie aussitôt que le composé $\lgg_n\circ \tau_n$ est le morphisme \eqref{mtfla12d}
\begin{equation}\label{ktfr24h}
\Theta_n\colon (\tE'_s,\ocB'_n)\rightarrow (\tE_s,\ocB_n).
\end{equation}

Avec les notations et conventions de \ref{TFA5}, le diagramme d'homomorphismes d'anneaux 
\begin{equation}\label{ktfr24i}
\xymatrix{
{g_{s*}(\co_{\oX'_n})}\ar[d]&{\co_{\oX_n}}\ar[r]\ar[l]&{\sigma_{s*}(\ocB_n)}\ar[d]\\
{g_{s*}(\pi_{s*}(\ocB^!_n))}\ar[rr]^-(0.5)\sim&&{\sigma_{s*}(\lgg_{s*}(\ocB^!_n))}}
\end{equation}
où l'isomorphisme horizontal inférieur est induit par le diagramme commutatif \eqref{ktfr21h}, 
est clairement commutatif. Le diagramme de morphismes de topos annelés
\begin{equation}\label{ktfr24j}
\xymatrix{
{(\tG_s,\ocB_n^!)}\ar[r]^-(0.5){\pi_n}\ar[d]_-(0.5){\lgg_n}&{(X'_{s,\et},\co_{\oX'_n})}\ar[d]^{\ogg_n}\\
{(\tE_s,\ocB_n)}\ar[r]^-(0.5){\sigma_n}&{(X_{s,\et},\co_{\oX_n})}}
\end{equation}
est donc commutatif à isomorphisme canonique près. 
Pour tout $\co_{\oX'_n}$-module $\cF'$ de $X'_{s,\et}$ et tout entier $q\geq 0$, on a un morphisme canonique  de changement de base (\cite{egr1} (1.2.3.3))
\begin{equation}\label{ktfr24k}
\sigma_n^*(\rR^q\ogg_{n*}(\cF'))\rightarrow \rR^q\lgg_{n*}(\pi_n^*(\cF')),
\end{equation} 
où $\sigma_n^*$ et $\pi_n^*$ désignent les images inverses au sens des topos annelés.

\begin{prop}\label{ktfr27}
Pour tout entier $n\geq 1$, l'homomorphisme canonique
\begin{equation}\label{ktfr27a}
\lgg_s^{-1}(\ocB_n)\otimes_{\pi_s^{-1}(g_s^{-1}(\co_{\oX_n}))} \pi_s^{-1}(\co_{\oX'_n}) \rightarrow \ocB^!_n
\end{equation}
est un $\alpha$-isomorphisme \eqref{ktfr21h}. 
\end{prop}

La question étant locale pour les topologies étales de $X$ et $X'$ (cf. \ref{tfr42} et \ref{tfra6}), on peut supposer 
que le morphisme $g$ admet une carte relativement adéquate \eqref{mtfla5}. 
En vertu de \ref{ktfr22}(ii), il suffit de montrer que pour tout point $(\oy\rightsquigarrow \ox')$ de $X'_\et\gtimes_{X_\et}\oX^\circ_\et$ \eqref{tfr21}
tel que $\ox'$ soit au-dessus de $s$, 
la fibre de l'homomorphisme \eqref{ktfr27a} en $\varrho(\oy\rightsquigarrow \ox')$ est un $\alpha$-isomorphisme. 
D'après \ref{TFA15}(i) appliqué à $X'$, il existe un point $(\oy'\rightsquigarrow \ox')$ de $X'_\et\gtimes_{X'_\et}\oX'^\rhd_\et$.
Compte tenu de \ref{ktfr29}, on peut supposer que $(\oy\rightsquigarrow \ox')$ est l'image de  $(\oy'\rightsquigarrow \ox')$ par le morphisme
canonique
\begin{equation}\label{ktfr27b}
X'_\et\gtimes_{X'_\et}\oX'^\rhd_\et\rightarrow X'_\et\gtimes_{X_\et}\oX^\circ_\et.
\end{equation}
Reprenons les notations de \ref{ktfr10}. Compte tenu de \eqref{tfr21a}, \eqref{TFA11c}, \eqref{ktfr9a} et \eqref{ktfr1f}, 
la fibre de l'homomorphisme \eqref{ktfr27a} en $\varrho(\oy\rightsquigarrow \ox')$ s'identifie à l'homomorphisme canonique 
\begin{equation}\label{ktfr27c}
(\oR^{\oy}_\uX/p^n\oR^{\oy}_\uX)\otimes_{\Gamma(\uoX,\co_{\uoX})} \Gamma(\uoX',\co_{\uoX'})\rightarrow
\oR^{!\oy}_{\uX'\rightarrow \uX}/p^n\oR^{!\oy}_{\uX'\rightarrow \uX}.
\end{equation}

Considérons la catégorie $\cC_{\ox'}$ définie dans \ref{ktfr13} et notons $\cC^\aff_{\ox'}$ la sous-catégorie pleine formée des 
objets $(\ox'\rightarrow U'\rightarrow U)$ de $\cC_{\ox'}$ tels que $U$ et $U'$ soient affines et connexes. 
La catégorie $\cC^\aff_{\ox'}$ est cofiltrante et le foncteur d'injection canonique $\cC^{\aff\circ}_{\ox'}\rightarrow \cC^\circ_{\ox'}$ 
est cofinal d'après (\cite{sga4} I 8.1.3(c)). Pour tout objet $(\ox'\rightarrow U'\rightarrow U)$ de $\cC^\aff_{\ox'}$,  l'homomorphisme canonique 
\begin{equation}\label{ktfr27d}
(\oR^{\oy}_U/p^n\oR^{\oy}_U)\otimes_{\Gamma(U,\co_{U})} \Gamma(U',\co_{U'})\rightarrow
\oR^{!\oy}_{U'\rightarrow U}/p^n\oR^{!\oy}_{U'\rightarrow U}
\end{equation}
est un $\alpha$-isomorphisme en vertu de \ref{eccr100}(iii). On en déduit par passage à la limite inductive sur $\cC^{\aff\circ}_{\ox'}$
que \eqref{ktfr27c} est un $\alpha$-isomorphisme, d'où la proposition.

\begin{lem}\label{ktfr300}
Supposons le morphisme $g\colon X'\rightarrow X$ séparé. 
Il existe alors un entier $N\geq 0$ tel que pour tout point $(\oy\rightsquigarrow \ox)$ de $X_\et\gtimes_{X_\et}\oX^\circ_\et$, avec $\ox$ au-dessus de $s$, 
notant $\uX$ le localisé strict de $X$ en $\ox$ et posant $\uX'=X'\times_X\uX$, $\uR_1=\Gamma(\uoX,\co_\uoX)$  et $\uoR=\ocB_{\rho(\oy\rightsquigarrow \ox)}$, 
pour tout $\co_{\uoX'}$-module quasi-cohérent $\cM$ \eqref{hght1c} et tout entier $q\geq 0$,
le noyau et le conoyau du morphisme canonique
\begin{equation}
\rH^q(\uoX',\cM)\otimes_{\uR_1}\uoR\rightarrow \rH^q(\uoX',\cM\otimes_{\co_{\uoX}}\uoR)
\end{equation}
soient annulés par $p^N$.
\end{lem}

On notera d'abord que le schéma $\uX'$ n'est pas local; on n'a d'ailleurs fixé aucun point géométrique de $X'$ au-dessus de $\ox$. 
On prendra garde de ne pas confondre $\uR_1$ et $\uoR$.

Comme $X$ est quasi-compact, on peut supposer que le morphisme $f$ \eqref{hght2} admet une carte adéquate \eqref{cad1}. 
Montrons alors que l'entier $N$ fourni par la proposition \ref{AFR32} convient. 
Soit $(\oy\rightsquigarrow \ox)$ un point de $X_\et\gtimes_{X_\et}\oX^\circ_\et$, avec $\ox$ au-dessus de $s$, pour lequel nous reprenons les notations de l'énoncé.
Soient $\cM$ un $\co_{\uoX'}$-module quasi-cohérent, $\cU'=(U'_i)_{1\leq i\leq r}$ un recouvrement ouvert affine de $\uoX'$.  
Comme $g$ est séparé, pour tout entier $q$, on a un isomorphisme canonique, fonctoriel en $\cM$,
\begin{equation}
\rH^q(\rC^\bullet(\cU',\cM))\stackrel{\sim}{\rightarrow} \rH^q(\uoX',\cM).
\end{equation}

Le schéma $\uoX$ étant strictement local d'après (\cite{agt} III.3.7), il s'identifie au localisé strict de $\oX$ en $\ox$.
Comme le morphisme canonique de complexes
\begin{equation}
\rC^\bullet(\cU',\cM)\otimes_{\uR_1}\uoR\rightarrow \rC^\bullet(\cU',\cM\otimes_{\co_{\uoX}}\uoR)
\end{equation}
est un isomorphisme, la proposition résulte alors de \ref{AFR32}.

\begin{teo}\label{ktfr31}
Supposons le morphisme $g\colon X'\rightarrow X$ propre. 
Il existe alors un entier $N\geq 0$ tel que pour tous entiers $n\geq 1$ et $q\geq 0$ et 
tout $\co_{\oX'_n}$-module quasi-cohérent $\ocF'$ de $X'_{s,\zar}$, 
que l'on considère aussi comme un $\co_{\oX'_n}$-module de $X'_{s,\et}$ \eqref{notconv12}, 
le noyau et le conoyau du morphisme de changement de base \eqref{ktfr24k}
\begin{equation}\label{ktfr31a}
\sigma_n^*(\rR^q\ogg_{n*}(\ocF'))\rightarrow \rR^q\lgg_{n*}(\pi_n^*(\ocF'))
\end{equation}
soient annulés par $p^N$.
\end{teo}

En effet, la question étant locale pour la topologie étale sur $X$ \eqref{tfr42}, on peut supposer $X=\Spec(R)$ affine. 
Montrons que si $N$ est l'entier fourni par la proposition \ref{ktfr300}, alors $N+1$ convient. 
Soient $n$ un entier $\geq 1$, $\ocF'$ un $\co_{\oX'_n}$-module quasi-cohérent, 
$(\oy\rightsquigarrow \ox)$ un point de $X_\et\gtimes_{X_\et}\oX^\circ_\et$ tel que $\ox$ soit au-dessus de $s$. 
Notons $\uX$ le localisé strict de $X$ en $\ox$ et posons $\uX'=X'\times_X\uX$, 
$\uR_1=\Gamma(\uoX,\co_\uoX)$ et $\uoR=\ocB_{\rho(\oy\rightsquigarrow \ox)}$. 
En vertu de \ref{ktfr300}, pour tout entier $q\geq 0$, le noyau et le conoyau du morphisme canonique
\begin{equation}\label{ktfr31b}
\rH^q(\uoX',\ocF'\otimes_{\co_{\oX}}\co_{\uoX})\otimes_{\uR_1}\uoR
\rightarrow \rH^q(\uoX',\ocF'\otimes_{\co_{\oX}}\uoR)
\end{equation}
sont annulés par $p^N$. 
On notera que les groupes de cohomologie ne changent pas si l'on remplace la topologie de Zariski par la topologie étale. 
Par ailleurs, le schéma $\uoX$ étant strictement local d'après (\cite{agt} III.3.7), il s'identifie au localisé strict de $\oX$ en $\ox$. 
Pour tout entier $q\geq 0$, on a donc un isomorphisme canonique 
\begin{equation}\label{hmdf60c}
\rR^q\ogg_{n*}(\ocF')_\ox\stackrel{\sim}{\rightarrow} \rH^q(\uoX',\ocF'\otimes_{\co_{\oX}}\co_{\uoX}).
\end{equation}
 
Notons $j_{\ox}\colon X'_{\ox}\rightarrow X'_s$ le morphisme canonique et $\jmath\colon X'_{\ox,\et}\rightarrow \tG$
le morphisme de topos associé à $(\oy\rightsquigarrow \ox)$ défini dans \eqref{lptfr10f}. D'après \eqref{lptfr10k} et \eqref{ktfr21h}, les carrés du diagramme 
\begin{equation}\label{ktfr31e}
\xymatrix{
{X'_{\ox,\et}}\ar[rr]^-(0.5){\jmath}\ar[d]_\varepsilon&&{\tG}\ar[d]_{\lgg_s}\ar[r]^-(0.5){\pi_s}&{X'_{s,\et}}\ar[d]^{g_s}\\
\Ens\ar[rr]^-(0.5){\rho(\oy\rightsquigarrow \ox)}&&{\tE}\ar[r]^-(0.5){\sigma_s}&{X_{s,\et}}}
\end{equation}
où $\Ens$ est le topos ponctuel et $\varepsilon$ est la projection canonique (\cite{sga4} IV 4.3), sont commutatifs à isomorphismes canoniques près.

Par le théorème de changement de base propre (\cite{sga4} XII 5.5), 
on déduit de \eqref{ktfr31b} que le noyau et le conoyau du morphisme canonique
\begin{equation}\label{ktfr31c}
\rH^q(X'_{\ox,\et},j_\ox^{-1}(\ocF'))\otimes_{\uR_1}\uoR\rightarrow \rH^q(X'_{\ox,\et},j_\ox^{-1}(\ocF')\otimes_{\uR_1}\uoR)
\end{equation}
sont annulés par $p^N$. En vertu de \ref{ktfr27} et \eqref{ktfr31e}, le morphisme canonique 
\begin{equation}
j_\ox^{-1}(\co_{\oX'_n}) \otimes_{\uR_1}\uoR \rightarrow \jmath^{-1}(\ocB^!_n)
\end{equation}
est un $\alpha$-isomorphisme de $\co_\oK$-modules de $X'_{\ox,\et}$. Par suite, le morphisme canonique 
\begin{equation}\label{ktfr31d}
\rH^q(X'_{\ox,\et},j_\ox^{-1}(\ocF')\otimes_{\uR_1}\uoR) \rightarrow \rH^q(X'_{\ox,\et},j_\ox^{-1}(\ocF')\otimes_{j_\ox^{-1}(\co_{\oX'_n})} 
\jmath^{-1}(\ocB^!_n))
\end{equation}
est un $\alpha$-isomorphisme. 
Composant \eqref{ktfr31c} et \eqref{ktfr31d}, on en déduit que le noyau et le conoyau du morphisme canonique  
\begin{equation}
\rR^q\ogg_{n*}(\ocF')_\ox\otimes_{\uR_1}\uoR \rightarrow \rH^q(X'_{\ox,\et},j_\ox^{-1}(\ocF')\otimes_{j_\ox^{-1}(\co_{\oX'_n})} 
\jmath^{-1}(\ocB^!_n))
\end{equation}
sont annulés par $p^{N+1}$. 
Celui-ci s'identifie au morphisme de changement de base (\cite{egr1} (1.2.3.3)) par rapport au rectangle extérieur du diagramme \eqref{ktfr31e}. 
Par suite, en vertu de \ref{lptfr16} et (\cite{egr1} 1.2.4(ii)), le noyau et le conoyau de la fibre du morphisme 
\eqref{ktfr31a} en le point $\rho(\oy\rightsquigarrow \ox)$ de $\tE$
sont annulés par $p^{N+1}$. La proposition s'ensuit compte tenu de (\cite{agt} III.9.5).

\begin{prop}\label{ktfr30}
Soient $n, q$ deux entiers tels que $n\geq 1$ et $q\geq 0$, $\cF'$ un $\co_{X'_n}$-module cohérent de $X'_{s,\zar}$
qui soit $X_n$-plat. Supposons que $g\colon X'\rightarrow X$ soit propre et que pour tout entier $i\geq 0$, 
le $\co_{X_n}$-module $\rR^ig_{n*}(\cF')$ soit localement libre (de type fini). 
Posons $\ocF'=\cF'\otimes_{\co_{S}}\co_{\oS}$, que l'on considère aussi comme un $\co_{\oX'_n}$-module de $X'_{s,\et}$ \eqref{notconv12}. 
Alors, le morphisme de changement de base \eqref{ktfr24k}
\begin{equation}\label{ktfr30a}
\sigma_n^*(\rR^q\ogg_{n*}(\ocF'))\rightarrow \rR^q\lgg_{n*}(\pi_n^*(\ocF'))
\end{equation}
est un $\alpha$-isomorphisme.
\end{prop}

En effet, la question étant locale pour la topologie étale sur $X$ \eqref{tfr42}, on peut supposer $X=\Spec(R)$ affine.  
D'après (\cite{sp}, \href{https://stacks.math.columbia.edu/tag/07VK}{07VK}), $\rR\Gamma(X',\cF')$ est un complexe 
parfait de la catégorie dérivée $\bD(\bMod(R/p^nR))$, 
et pour toute $R$-algèbre $A$, posons $X'\otimes_{X}A=X'\times_{X}\Spec(A)$, le morphisme canonique 
\begin{equation}
\rR\Gamma(X',\cF')\otimes^\rL_{R/p^nR}(A/p^nA)
\rightarrow \rR\Gamma(X'\otimes_{X}A,\cF'\otimes_{R}A)
\end{equation}
est un isomorphisme. Compte tenu des hypothèses, on en déduit que le morphisme canonique 
\begin{equation}\label{ktfr30e}
\rH^q(X',\cF')\otimes_{R}A \rightarrow \rH^q(X'\otimes_{X}A,\cF'\otimes_{R}A)
\end{equation}
est un isomorphisme.

Soient $(\oy\rightsquigarrow \ox)$ un point de $X_\et\gtimes_{X_\et}\oX^\circ_\et$ tel que $\ox$ soit au-dessus de $s$,
$\uX$ le localisé strict de $X$ en $\ox$, $\uX'=X'\times_X\uX$. Posons $\uR_1=\Gamma(\uoX,\co_\uoX)$ et  
$\uoR=\ocB_{\rho(\oy\rightsquigarrow \ox)}$. D'après \eqref{ktfr30e}, le  morphisme canonique 
\begin{equation}\label{ktfr30b}
\rH^q(\uoX',\ocF'\otimes_{\co_{\oX}}\co_{\uoX})\otimes_{\uR_1}\uoR
\rightarrow \rH^q(\uoX',\ocF'\otimes_{\co_{\oX}}\uoR)
\end{equation}
est un isomorphisme. 

Notons $j_{\ox}\colon X'_{\ox}\rightarrow X'_s$ le morphisme canonique et $\jmath\colon X'_{\ox,\et}\rightarrow \tG$
le morphisme de topos associé à $(\oy\rightsquigarrow \ox)$ défini dans \eqref{lptfr10f} (cf. la preuve de \ref{ktfr31}). 
Par le théorème de changement de base propre (\cite{sga4} XII 5.5), on déduit de \eqref{ktfr30b} que le morphisme canonique
\begin{equation}\label{ktfr30c}
\rH^q(X'_{\ox,\et},j_\ox^{-1}(\ocF'))\otimes_{\uR_1}\uoR\rightarrow \rH^q(X'_{\ox,\et},j_\ox^{-1}(\ocF')\otimes_{\uR_1}\uoR)
\end{equation}
est un isomorphisme. En vertu de \ref{ktfr27} et \eqref{ktfr31e}, le morphisme canonique 
\begin{equation}
j_\ox^{-1}(\co_{\oX'_n}) \otimes_{\uR_1}\uoR \rightarrow \jmath^{-1}(\ocB^!_n)
\end{equation}
est un $\alpha$-isomorphisme de $\co_\oK$-modules de $X'_{\ox,\et}$. Par suite, le morphisme canonique 
\begin{equation}\label{ktfr30d}
\rH^q(X'_{\ox,\et},j_\ox^{-1}(\ocF')\otimes_{\uR_1}\uoR) \rightarrow \rH^q(X'_{\ox,\et},j_\ox^{-1}(\ocF')\otimes_{j_\ox^{-1}(\co_{\oX'_n})} 
\jmath^{-1}(\ocB^!_n))
\end{equation}
est un $\alpha$-isomorphisme. 
Composant \eqref{ktfr30c} et \eqref{ktfr30d}, on en déduit que morphisme canonique  
\begin{equation}
\rR^q\ogg_{n*}(\ocF')_\ox\otimes_{\uR_1}\uoR \rightarrow \rH^q(X'_{\ox,\et},j_\ox^{-1}(\ocF')\otimes_{j_\ox^{-1}(\co_{\oX'_n})} 
\jmath^{-1}(\ocB^!_n))
\end{equation}
est un $\alpha$-isomorphisme. Par suite, en vertu de \ref{lptfr16} et (\cite{egr1} 1.2.4(ii)),
la fibre du morphisme \eqref{ktfr30a} en le point $\rho(\oy\rightsquigarrow \ox)$ de $\tE$
est un $\alpha$-isomorphisme (cf. la preuve de \ref{ktfr31}). La proposition s'ensuit compte tenu de (\cite{agt} III.9.5).

\subsection{}\label{ktfr32}
Pour tout $\mU$-topos $T$, on note $T^{\mN^\circ}$ le topos des systèmes projectifs de $T$, 
indexés par l'ensemble ordonné $\mN$ des entiers naturels \eqref{notconv13}. 
En plus des notations introduites dans \ref{ssht1}, 
on désigne par $\bvocB'$ l'anneau $(\ocB'_{n+1})_{n\in \mN}$ de $\tE'^{\mN^\circ}_s$  \eqref{mtfla10e}, 
par $\bvocB^!$ l'anneau $(\ocB^!_{n+1})_{n\in \mN}$ de $\tG_s^{\mN^\circ}$ \eqref{ktfr24a} et
par $\co_{\bvoX'}$ l'anneau $(\co_{\oX'_{n+1}})_{n\in \mN}$ de $X'^{\mN^\circ}_{s,\et}$. 
Considérons le diagramme de morphismes de topos annelés 
\begin{equation}\label{ktfr32a}
\xymatrix{
{(\tE'^{\mN^\circ}_s,\bvocB')}\ar[rd]^-(0.5){\bvsigma'}\ar[d]_{\bvtau}\ar@/_3pc/[dd]_{\bvTheta}&\\
{(\tG^{\mN^\circ}_s,\bvocB^!)}\ar[r]^-(0.5){\bvpi}\ar[d]_{\bvlgg}&{(X'^{\mN^\circ}_{s,\et},\co_{\bvoX'})}\ar[d]^{\bvogg}\\
{(\tE^{\mN^\circ}_s,\bvocB)}\ar[r]^-(0.5){\bvsigma}&{(X^{\mN^\circ}_{s,\et},\co_{\bvoX})}}
\end{equation}
induits par $(\sigma'_{n+1})_{n\in \mN}$ \eqref{mtfla10f}, $(\tau_{n+1})_{n\in \mN}$ \eqref{ktfr24f}, 
$(\pi_{n+1})_{n\in \mN}$ \eqref{ktfr24e}, $(\lgg_{n+1})_{n\in \mN}$ \eqref{ktfr24g} et $(\Theta_{n+1})_{n\in \mN}$ \eqref{ktfr24h} (cf. \cite{agt} III.7.5). 
Les deux triangles et le carré sont commutatifs à isomorphismes canoniques près \eqref{ktfr24j}. 

Pour tout $\co_{\bvoX'}$-module $\cF'$ de $X'^{\mN^\circ}_{s,\et}$ et tout entier $q\geq 0$, 
on a un morphisme canonique de changement de base (\cite{egr1} (1.2.3.3))
\begin{equation}\label{ktfr32b}
\bvsigma^*(\rR^q\bvogg_*(\cF'))\rightarrow \rR^q\bvlgg_*(\bvpi^*(\cF')),
\end{equation}
où $\bvsigma^*$ et $\bvpi^*$ désignent les images inverses au sens des topos annelés.

\begin{lem}\label{ktfr33}
Pour tout $\co_{\bvoX'}$-module $\cF'=(\cF'_n)_{n\geq 0}$ de $X'^{\mN^\circ}_{s,\et}$ et tout entier $q\geq 0$, le morphisme de changement de base
\eqref{ktfr32b}
\begin{equation}\label{ktfr33a}
\bvsigma^*(\rR^q\bvogg_*(\cF'))\rightarrow \rR^q\bvlgg_*(\bvpi^*(\cF'))
\end{equation}
est induit par les morphismes de changement de base \eqref{ktfr24k}
\begin{equation}\label{ktfr33b}
\sigma_n^*(\rR^q\ogg_{n*}(\cF'_n))\rightarrow \rR^q\lgg_{n*}(\pi_n^*(\cF'_n)), \ \ \ (n\in \mN). 
\end{equation}
\end{lem}
En effet, d'après (\cite{agt} III.7.1), pour tout entier $n\geq 0$, il existe deux morphismes de topos 
$a_n\colon \tE_s\rightarrow \tE^{\mN^\circ}_s$ et $b_n\colon \tG_s\rightarrow \tG^{\mN^\circ}_s$ tels que pour tous objets 
$M=(M_n)_{n\geq 0}$ de $\tE^{\mN^\circ}_s$ et $N=(N_n)_{n\geq 0}$ de $\tG^{\mN^\circ}_s$, on ait $a_n^*(M)=M_n$ et $b_n^*(N)=N_n$. 
D'après (\cite{agt} (III.7.5.4)), le diagramme de morphisme de topos
\begin{equation}
\xymatrix{
{\tG_s}\ar[r]^{b_n}\ar[d]_{\lgg_s}&{\tG^{\mN^\circ}_s}\ar[d]^{\bvlgg}\\
{\tE_s}\ar[r]^{a_n}&{\tE^{\mN^\circ}_s}}
\end{equation}
est commutatif à isomorphisme canonique près. Par ailleurs, pour tout $\bvocB^!$-module $N$ de $\tG^{\mN^\circ}_s$, le morphisme de changement de base 
\begin{equation}
a_n^*(\rR^q\bvlgg_*(N))\rightarrow\rR^q\lgg_{n*}(b_n^*(N))
\end{equation}
est un isomorphisme (\cite{agt} (III.7.5.5)). La proposition résulte alors de (\cite{egr1} 1.2.4(ii)).

\subsection{}\label{ktfr34}
On pose $\cS=\Spf(\co_C)$ et on désigne par $\fX$ (resp. $\fX'$) le schéma formel complété $p$-adique de $\oX$ (resp. $\oX'$), 
et par $\fgg\colon \fX'\rightarrow \fX$ le morphisme induit par $g$ \eqref{ktfr1a}. 
Pour tout entier $n\geq 1$, on note
\begin{equation}\label{ktfr34a}
u_n\colon (X_{s,\et},\co_{\oX_n})\rightarrow (X_{s,\zar},\co_{\oX_n})
\end{equation}
le  morphisme canonique \eqref{notconv12}. 
On désigne par 
\begin{equation}\label{ktfr34b}
\bvu\colon (X_{s,\et}^{\mN^\circ},\co_{\bvoX})\rightarrow (X_{s,\zar}^{\mN^\circ},\co_{\bvoX})
\end{equation}
le morphisme de topos annelés défini par les $(u_{n+1})_{n\in \mN}$ et par
\begin{equation}\label{ktfr34c}
\uplambda\colon (X_{s,\zar}^{\mN^\circ},\co_{\bvoX})\rightarrow (X_{s,\zar}, \co_\fX)
\end{equation}
le morphisme de topos annelés pour lequel  le foncteur $\uplambda_*$ est le foncteur limite projective \eqref{notconv13a}. 
On considère aussi les notations analogues pour $f'$, que l'on munit d'un exposant $^\prime$.

Considérons le diagramme de morphismes de topos annelés 
\begin{equation}\label{ktfr34d}
\xymatrix{
{(X'^{\mN^\circ}_{s,\et},\co_{\bvoX'})}\ar[r]^-(0.5){\bvu'}\ar[d]_{\bvogg_\et}&{(X'^{\mN^\circ}_{s,\zar},\co_{\bvoX'})}\ar[r]^-(0.5){\uplambda'}\ar[d]_{\bvogg_\zar}&
{(X'_{s,\zar},\co_{\fX'})}\ar[d]_{\fgg}\ar[r]^-(0.5){i_{\oX'}}&{(\oX'_\zar,\co_{\oX'})}\ar[d]^g\\
{(X_{s,\et}^{\mN^\circ},\co_{\bvoX})}\ar[r]^-(0.5){\bvu}&{(X_{s,\zar}^{\mN^\circ},\co_{\bvoX})}\ar[r]^-(0.5){\uplambda}&
{(X_{s,\zar},\co_{\fX})}\ar[r]^-(0.5){i_\oX}&{(\oX_\zar,\co_\oX)}}
\end{equation}
où $i_\oX$ et $i_{\oX'}$ sont les morphismes canoniques. 
On vérifie aussitôt que les trois carrés sont commutatifs à isomorphismes canoniques près (\cite{egr1} 1.2.3).

\begin{lem}\label{ktfr35}
Supposons le morphisme $g\colon X'\rightarrow X$ propre. Soit $\cF'$ un $\co_{\oX'}$-module cohérent de $\oX'_\zar$, $q$ un entier $\geq 0$. 
Alors, le morphisme de changement de base relativement au troisième carré du diagramme \eqref{ktfr34d}
\begin{equation}\label{ktfr35a}
i_\oX^*(\rR^qg_*(\cF'))\rightarrow \rR^q\fgg(i^*_{\oX'}(\cF'))
\end{equation}
est un isomorphisme.
\end{lem}

Cela résulte de (\cite{egr1} 2.5.5 et  2.12.2).

\begin{lem}\label{ktfr37}
Supposons le morphisme $g\colon X'\rightarrow X$ séparé et quasi-compact. 
Soient $\cF'=(\cF'_n)_{n\in \mN}$ un $\co_{\bvoX'}$-module de $X'^{\mN^\circ}_{s,\zar}$ tel que pour tout entier $n\geq 0$, 
le $\co_{\oX'_n}$-module $\cF'_n$ soit quasi-cohérent, 
$q$ un entier $\geq 0$. Alors, le morphisme de changement de base relativement au premier carré du diagramme \eqref{ktfr34d}
\begin{equation}\label{ktfr37a}
\bvu^*(\rR^q\bvogg_{\zar*}(\cF'))\rightarrow \rR^q\bvogg_{\et*}(\bvu'^*(\cF'))
\end{equation}
est un isomorphisme.
\end{lem}

En effet, pour tout entier $n\geq 0$, le diagramme de morphismes de topos annelés 
\begin{equation}\label{ktfr37b}
\xymatrix{
{(X'_{s,\et},\co_{\oX'_n})}\ar[r]^{u'_n}\ar[d]_{\ogg_{n,\et}}&{(X'_{s,\zar},\co_{\oX'_n})}\ar[d]^{\ogg_{n,\zar}}\\
{(X_{s,\et},\co_{\oX_n})}\ar[r]^{u_n}&{(X_{s,\zar},\co_{\oX_n})}}
\end{equation}
est commutatif à isomorphisme canonique près. Compte tenu de (\cite{agt} III.7.5; cf. la preuve de \ref{ktfr33}), il suffit de montrer que 
le morphisme de changement de base 
\begin{equation}\label{ktfr37c}
u^*_n(\rR^q\ogg_{n,\zar*}(\cF'_n))\rightarrow \rR^q\ogg_{n,\et*}(u'^*_n(\cF'_n))
\end{equation}
est un isomorphisme.  

D'après \eqref{notconv12e}, on a un isomorphisme canonique 
\begin{equation}\label{ktfr37d}
u'^*_n(\cF'_n)\stackrel{\sim}{\rightarrow} \cF'_{n,\et},
\end{equation}
où $\cF'_{n,\et}$ désigne le $\co_{\oX'_n}$-module de $X'_{s,\et}$ associé à $\cF'_n$ \eqref{notconv12a}. 
De même, le $\co_{\oX_n}$-module $\rR^q\ogg_{n,\zar*}(\cF'_n)$ étant quasi-cohérent, on a un isomorphisme canonique 
\begin{equation}\label{ktfr37e}
u^*_n(\rR^q\ogg_{n,\zar*}(\cF'_n))\stackrel{\sim}{\rightarrow} (\rR^q\ogg_{n,\zar*}(\cF'_n))_\et.
\end{equation}

Pour tout $X$-schéma $U$, on pose $U'=U\times_XX'$. D'après (\cite{sga4} V 5.1),
le faisceau $\rR^q\ogg_{n,\et*}(\cF'_{n,\et})$ s'identifie au faisceau associé au préfaisceau défini par
\begin{equation}\label{ktfr37f}
U\in \ob(\Et_{/\oX_n})\mapsto \rH^q(U'_\et, (\cF'_n\otimes_{\co_{\oX'_n}}\co_{U'})_\et).
\end{equation}
Par ailleurs, pour tout $\oX_n$-schéma étale $U$, le morphisme de changement de base 
\begin{equation}\label{ktfr37h}
\rR^q\ogg_{n,\zar*}(\cF'_n)\otimes_{\co_{\oX_n}}\co_U\rightarrow \rR^q\ogg_{U,\zar*}(\cF'_n\otimes_{\co_{\oX'_n}}\co_{U'})
\end{equation}
est un isomorphisme (\cite{ega3} 1.4.15). On en déduit par un calcul de fibres similaire à \eqref{notconv12f} 
que le faisceau $(\rR^q\ogg_{n,\zar*}(\cF'_n))_\et$ s'identifie au faisceau associé au préfaisceau défini par
\begin{equation}\label{ktfr37g}
U\in \ob(\Et_{/\oX_n})\mapsto \rH^q(U'_\zar, \cF'_n\otimes_{\co_{\oX'_n}}\co_{U'}).
\end{equation}
 
Pour tout $\oX_n$-schéma étale $U$, le diagramme 
\begin{equation}\label{ktfr37i}
\xymatrix{
{\rH^q(U'_\zar, \cF'_n\otimes_{\co_{\oX'_n}}\co_{U'})}\ar[r]^-(0.5){a_U}\ar[d]_{c_U}&{\rH^q(U'_\et, (\cF'_n\otimes_{\co_{\oX'_n}}\co_{U'})_\et)}\ar[d]^{d_U}\\
{\Gamma(U,u_n^*(\rR^q\ogg_{n,\zar*}(\cF'_n)))}\ar[r]^-(0.5){b_U}&{\Gamma(U,\rR^q\ogg_{n,\et*}(u'^*_n(\cF'_{n})))}}
\end{equation}
où $a_U$ et $d_U$ sont les morphismes canoniques et $b_U$ (resp. $c_U$) est induit par \eqref{ktfr37c} (resp. \eqref{ktfr37h}), est commutatif.
En effet, on peut se réduire facilement au cas où $U=\oX_n$ (\cite{egr1} 1.2.4(ii)). Considérons alors le diagramme
\begin{equation}\label{ktfr37j}
{\tiny 
\xymatrix{
{\rH^q(\oX'_{n},\cF'_n)}\ar[r]\ar[d]\ar@{}[rd]|{(1)}&{\rH^q(X'_{n},u'_{n*}(u'^*_n(\cF'_n)))}\ar@{}[rd]|{(2)}\ar[r]\ar[d]&{\rH^q(X'_{n},u'^*_n(\cF'_n))}\ar[d]\\
{\rH^0(\oX_{n},\rR^q\ogg_{n,\zar*}(\cF'_n))}\ar[r]\ar[d]\ar@{}[rrd]|{(3)}&{\rH^0(\oX_{n},\rR^q\ogg_{n,\zar*}(u'_{n*}(u'^*_n(\cF'_n))))}\ar[r]&
{\rH^0(\oX_{n},\rR^q(\ogg_{n,\zar}\circ u'_n)_*(u'^*_n(\cF'_n)))}\ar[d]\\
{\rH^0(\oX_{n},u_n^*(\rR^q\ogg_{n,\zar*}(\cF'_n)))}\ar[r]&{\rH^0(\oX_{n},\rR^q\ogg_{n,\et*}(u'^*_n(\cF'_n)))}&
{\rH^0(\oX_{n},u_{n*}(\rR^q\ogg_{n,\et*}(u'^*_n(\cF'_n))))}\ar@{=}[l]}}
\end{equation}
Le carré $(1)$ et le rectangle $(3)$ sont clairement commutatifs. Pour tout $\co_{\oX'_n}$-module $\cG'$ de $X'_{s,\et}$
et tout ouvert de Zariski $U$ de $\oX_n$, les edge-homomorphismes des suites spectrales de Cartan-Leray 
\begin{eqnarray}
\rH^q(U',u'_{n*}(\cG'))&\rightarrow& \rH^q(U',\cG'), \\
\rR^q\ogg_{n,\zar*}(u'_{n*}(\cG'))&\rightarrow& \rR^q(\ogg_{n,\zar}\circ u'_n)_*(\cG'),
\end{eqnarray}
sont compatibles, dans le sens où le second se déduit du premier en passant aux faisceaux associés. Il s'ensuit que le carré $(2)$ de \eqref{ktfr37j} est commutatif.
Par ailleurs, il résulte de \ref{ssht21} que $d_{\oX_n}$ \eqref{ktfr37i} 
s'identifie au composé des flèches verticales droite et de l'égalité de \eqref{ktfr37j}. 
Ceci achève la preuve de la commutativité du diagramme \eqref{ktfr37i}. 
On en déduit que le morphisme de préfaisceaux défini par les morphismes $a_U$, pour $U\in \ob(\Et_{/\oX_n})$, 
induit le morphisme \eqref{ktfr37c} entre les faisceaux associés.
La proposition s'ensuit puisque les $a_U$ sont des isomorphismes d'après (\cite{sga4} VII 4.3).

\begin{prop}\label{ktfr36}
Supposons le morphisme $g\colon X'\rightarrow X$ propre. Soit $\cF'$ un $\co_{\fX'}$-module cohérent de $X'_{s,\zar}$, $q$ un entier $\geq 0$. 
Alors, il existe un entier $N\geq 0$ tel que le noyau et le conoyau du morphisme de changement de base relativement au deuxième carré du diagramme \eqref{ktfr34d}
\begin{equation}\label{ktfr36a}
\uplambda^*(\rR^q\fgg_*(\cF'))\rightarrow \rR^q\bvogg_{\zar*}(\uplambda'^*(\cF'))
\end{equation}
soient annulés par $p^N$.
\end{prop}

Pour tout entier $n\geq 0$, posons $\cF'_n=\cF'/p^n\cF'$. Il s'agit de montrer que le noyau et le conoyau du morphisme canonique
\begin{equation}\label{ktfr36b}
\rR^q\fgg_*(\cF')/p^n \rR^q\fgg_*(\cF')\rightarrow \rR^q\ogg_{n*}(\cF'_n)
\end{equation}
sont annulés par $p^N$ pour un entier $N\geq 0$ indépendant de $n$ (\cite{agt} (III.7.5.5)). 
La question étant locale sur $X_\zar$, on peut supposer $X=\Spec(R)$ affine. On note $\hR$ le séparé complété $p$-adique de $R$. 
Comme les $\co_\fX$-modules $\rR^q\fgg_*(\cF')$ et $\rR^q\ogg_{n*}(\cF'_n)$ sont cohérents (\cite{egr1} 2.11.5),
il suffit de montrer que le noyau et le conoyau du morphisme canonique  
\begin{equation}\label{ktfr36c}
\rH^q(\fX',\cF')/p^n \rH^q(\fX',\cF')\rightarrow \rH^q(\oX'_n,\cF'_n)
\end{equation}
sont annulés par $p^N$ pour un entier $N\geq 0$ indépendant de $n$ (\cite{egr1} 2.7.2). 
Notons $N_n^q$ (resp. $Q_n^q$) le noyau (resp. conoyau) du morphisme canonique
\begin{equation}\label{ktfr36d}
\rH^q(\fX',\cF')\rightarrow \rH^q(\oX'_n,\cF'_n).
\end{equation}
Soient $\fU'$ un recouvrement ouvert affine fini de $\fX'$, $\rC^\bullet=\rC^\bullet(\fU',\cF')$ le complexe de \v{C}ech de $\cF'$ relativement à $\fU'$,
$\rC^\bullet_n=\rC^\bullet\otimes_R(R/p^nR)$. Pour tout entier $i\geq 0$, $\rH^q(\fX',\cF')=\rH^q(\rC^\bullet)$, $\rH^q(\fX',p^n\cF')=\rH^q(p^n\rC^\bullet)$ et 
$\rH^q(\fX',\cF'_n)=\rH^q(\rC^\bullet_n)$ (\cite{egr1} 2.11.4). 
Les $\hR$-modules $\rH^i(\rC^\bullet)$, $\rH^i(p^n\rC^\bullet)$ et $\rH^i(\rC^\bullet_n)$ sont cohérents en vertu de (\cite{egr1} 2.11.7).
Par ailleurs, le $\hR$-module $\rC^i$ est complet et séparé pour la topologie $p$-adique et il existe un entier $r_i\geq 0$ tel que 
$p^{r_i}\rC^i_{\tor}=0$ (\cite{egr1} 1.10.2). Il résulte alors de (\cite{egr1} 1.17.7) que les $(N^q_n)_{n\geq 0}$ définissent une filtration $p$-bonne sur $\rH^q(\fX',\cF')$;
autrement dit, il existe un entier $m\geq 0$ tel que $N^q_{n+1}=pN^q_n$ pour tout $n\geq m$. 

Le noyau du morphisme \eqref{ktfr36c} est $N_n^q/p^n\rH^q(\fX',\cF')$. Pour tout $n\geq m$, on a 
\begin{equation}
N_n^q/p^n\rH^q(\fX',\cF')=p^{n-m}N_m^q/p^n\rH^q(\fX',\cF'),
\end{equation}
qui est donc annulé par $p^m$. 

Pour tous entiers $n,i\geq 0$, le $\hR$-module $\rH^i(\fX',p^n\cF')$ étant cohérent d'après (\cite{egr1} 2.11.7), il existe un entier $t_{i,n}\geq 0$ 
tel que $p^{t_{i,n}}\rH^i(\fX',p^n\cF')_{\tor}=0$ (\cite{egr1} 1.10.2). Il s'ensuit en vertu de (\cite{egr1} 1.17.3) que le système projectif $(Q_n^q)_{n\geq 0}$ est AR-nul, 
autrement dit, il existe un entier $r\geq 0$ tel que pour tout entier $n\geq 0$, le morphisme canonique $Q^q_{n+r}\rightarrow Q^q_n$ soit nul. 
Soient $\zeta\in Q_n^q$, $\zeta'\in \rH^q(\oX'_n,\cF'_n)$ un relèvement. La suite exacte canonique $\cF'_r\rightarrow \cF'_{n+r}\rightarrow \cF'_n\rightarrow 0$
montre que $p^r\zeta'$ se relève en une section de  $\rH^q(\oX'_{n+r},\cF'_{n+r})$. 
Comme le morphisme $Q^q_{n+r}\rightarrow Q^q_n$ est nul, on en déduit que $p^r\zeta$ est nul. Par suite, $p^rQ_n^q=0$, ce qui achève la preuve de la proposition.

\subsection{}\label{ktfr38}
Les diagrammes \eqref{ktfr32a} et \eqref{ktfr34d} induisent un diagramme de morphismes de topos annelés, commutatif à isomorphisme canonique près,
\begin{equation}\label{ktfr38a}
\xymatrix{
{(\tG^{\mN^\circ}_s,\bvocB^!)}\ar[rr]^-(0.5){i_{\oX'}\circ \uplambda'\circ \bvu'\circ \bvpi}\ar[d]_{\bvlgg}&&{(\oX'_\zar,\co_{\oX'})}\ar[d]^\ogg\\
{(\tE^{\mN^\circ}_s,\bvocB)}\ar[rr]^-(0.5){i_{\oX}\circ \uplambda\circ \bvu\circ \bvsigma}&&{(\oX_\zar,\co_\oX)}}
\end{equation}
Pour tout $\co_\oX$-module $\cF$ de $\oX_\zar$, on note $\cF\otimes_{\co_\oX}\co_{\bvoX}$ son image inverse par $i_{\oX}\circ \uplambda\circ \bvu$ \eqref{ktfr34d},
qui n'est autre que le $\co_\bvoX$-module $(u_{n+1}^*(\cF\otimes_{\co_\oX}\co_{\oX_{n+1}}))_{n\in \mN}$ de  $X^{\mN^\circ}_{s,\et}$.
De même, pour tout $\co_{\oX'}$-module $\cF'$ de $\oX'_\zar$, on note $\cF'\otimes_{\co_{\oX'}}\co_{\bvoX'}$ son image inverse par $i_{\oX'}\circ \uplambda'\circ \bvu'$. 
Pour tout entier $q\geq 0$, on a donc un morphisme canonique de changement de base (\cite{egr1} (1.2.3.3))
\begin{equation}\label{ktfr38b}
\bvsigma^*(\rR^q\ogg_*(\cF')\otimes_{\co_\oX}\co_{\bvoX})\rightarrow \rR^q\bvlgg_*(\bvpi^*(\cF'\otimes_{\co_{\oX'}}\co_{\bvoX'})).
\end{equation}

\begin{prop}\label{ktfr39}
Supposons le morphisme $g\colon X'\rightarrow X$ propre. 
Il existe alors un entier $N\geq 0$ tel que pour tout $\co_{\oX'}$-module cohérent $\cF'$ de $\oX'_\zar$ et tout entier $q\geq 0$,  
le noyau et le conoyau du morphisme de changement de base \eqref{ktfr38b}
\begin{equation}\label{ktfr39a}
\bvsigma^*(\rR^q\ogg_*(\cF')\otimes_{\co_\oX}\co_{\bvoX})\rightarrow \rR^q\bvlgg_*(\bvpi^*(\cF'\otimes_{\co_{\oX'}}\co_{\bvoX'}))
\end{equation}
soient annulés par $p^N$.
\end{prop}

Cela résulte aussitôt de \ref{ktfr31}, \ref{ktfr33}, \ref{ktfr35}, \ref{ktfr37}, \ref{ktfr36} et (\cite{egr1} 1.2.4(ii)). 
On notera que le faisceau d'anneaux $\co_{\oX'}$ de $\oX'_\zar$ (resp. $\co_{\fX'}$ de $X'_{s,\zar}$ ) est cohérent (\cite{gr1} 9.1.27) (resp. (\cite{egr1} 2.8.1)).

\begin{prop}\label{ktfr40}
Soient $\cF$ un $\co_\oX$-module cohérent de $\oX_\zar$ tel que $\cF|X_\oeta$ soit un $\co_{X_\oeta}$-module localement libre, 
$\hcF$ son complété $p$-adique, qui est un $\co_\fX$-module. 
Alors, l'homomorphisme canonique
\begin{equation}\label{ktfr40a}
\Gamma(\fX,\hcF)\otimes_{\mZ_p}\mQ_p\rightarrow 
\Gamma(\tE^{\mN^\circ}_s,\bvsigma^*(\cF\otimes_{\co_\oX}\co_\bvoX))\otimes_{\mZ_p}\mQ_p
\end{equation}
est un isomorphisme. 
\end{prop}

Avec les notations de \ref{ktfr34}, on désigne par
\begin{equation}\label{ktfr40b}
\top\colon (\tE_s^{\mN^\circ},\bvocB)\rightarrow (X_{s,\zar},\co_{\fX})
\end{equation}
le morphisme composé 
\begin{equation}\label{ktfr40c}
(\tE_s^{\mN^\circ},\bvocB)\stackrel{\bvsigma}{\longrightarrow} (X_{s,\et}^{\mN^\circ},\co_{\coX})
\stackrel{\bvu}{\longrightarrow} (X_{s,\zar}^{\mN^\circ},\co_{\coX})
\stackrel{\uplambda}{\longrightarrow} (X_{s,\zar},\co_{\fX}).
\end{equation} 
Notant $\bMod_{\mQ}(\bvocB)$ la catégorie des $\bvocB$-modules à isogénie près (\cite{agt} III.6.1),
le foncteur $\top_*$ induit un foncteur additif et exact à gauche que l'on note encore 
\begin{equation}
\top_*\colon \bMod_{\mQ}(\bvocB) \rightarrow \bMod(\co_{\fX}[\frac 1 p]).
\end{equation} 
D'après (\cite{agt} III.6.16), le foncteur $\top^*$ induit un foncteur additif que l'on note encore
\begin{equation}
\top^*\colon \bMod^\coh(\co_{\fX}[\frac 1 p]) \rightarrow \bMod_{\mQ}^\atf(\bvocB). 
\end{equation}

Montrons que le morphisme $\top_*(\bvocB)$-linéaire 
\begin{equation}\label{ktfr40e}
\hcF\otimes_{\co_\fX}\top_*(\bvocB)[\frac 1 p]\rightarrow \top_*(\top^*(\hcF))[\frac 1 p]
\end{equation}
induit par le morphisme d'adjonction $\id\rightarrow \top_*\top^*$, est un isomorphisme. 
La question étant locale pour la topologie de Zariski de $\oX$, et le faisceau d'anneaux $\co_{\oX}$ de $\oX_\zar$ étant cohérent (\cite{gr1} 9.1.27),
on peut supposer qu'il existe un $\co_\oX$-module cohérent $\cF'$,
un entier $n\geq 1$ et un morphisme $\co_\oX$-linéaire $u\colon \cF\oplus\cF'\rightarrow \co_\oX^n$ induisant un isomorphisme sur $X_\oeta$. 
Il existe un entier $m\geq 0$, tel que le noyau et le conoyau de $u$ sont annulés par $p^m$. 
D'après \ref{alpha3}, il existe alors un morphisme $\co_X$-linéaire $v\colon \co_\oX^n\rightarrow \cF\oplus\cF'$
tel que $u\circ v= p^{2m}\id_{\co_\oX^n}$ et $v\circ u=p^{2m}\id_{\cF\oplus\cF'}$. 
Par suite, le morphisme 
\begin{equation}
\hcF[\frac 1 p]\oplus\hcF'[\frac 1 p]\rightarrow \co_\fX^n[\frac 1 p]
\end{equation}
induit par $u$, est un isomorphisme. On se réduit ainsi au cas où $\cF=\co_\oX$ auquel cas l'assertion recherchée est triviale. 

En vertu de \ref{tfkum11}, le noyau et le conoyau de l'homomorphisme canonique
\begin{equation}\label{sshtr10f}
\co_{\bvoX}\rightarrow \bvsigma_*(\bvocB)
\end{equation}
sont annulés par  $p^{\frac{2d+2}{p-1}}$, où $d=\dim(X/S)$. 
Par suite, l'homomorphisme canonique
\begin{equation}\label{sshtr10g}
\co_\fX[\frac 1 p]\rightarrow \top_*(\bvocB)[\frac 1 p]
\end{equation}
est un isomorphisme. On déduit de ce qui précède que l'homomorphisme \eqref{ktfr40a}
est un isomorphisme (\cite{egr1} 2.10.5).

\section{Cohomologie relative}\label{crtf}

\subsection{}\label{crtf1}
Rappelons que les hypothèses de \ref{ktfr1} sont en vigueur dans cette section. Posons 
\begin{equation}\label{crtf1a}
(\oX',\cL_{\oX'})=(X',\cM_{X'})\times_{(S,\cM_S)}(\oS,\cL_{\oS}),
\end{equation}
le produit étant pris dans la catégorie des schémas logarithmiques, de sorte que 
$\oX'=X'\times_S\oS$ \eqref{mtfla1b}.
Pour alléger les notations, on pose
\begin{equation}\label{crtf1b}
\tOmega^1_{\oX'/\oS}=\Omega^1_{(\oX',\cL_{\oX'})/(\oS,\cL_\oS)} \ \ \ {\rm et}\ \ \ 
\tOmega^1_{\oX'/\oX}=\Omega^1_{(\oX',\cL_{\oX'})/(\oX,\cL_\oX)},
\end{equation}
que l'on considère comme des faisceaux de $\oX'_\zar$ ou $\oX'_\et$, selon le contexte (cf. \ref{notconv12}). 
On a des isomorphismes canoniques \eqref{mtfla2e}
\begin{eqnarray}
\tOmega^1_{\oX'/\oS}&\stackrel{\sim}{\rightarrow}&\tOmega^1_{X'/S}\otimes_{\co_{X'}}\co_{\oX'},\label{crtf1d}\\
\tOmega^1_{\oX'/\oX}&\stackrel{\sim}{\rightarrow}&\tOmega^1_{X'/X}\otimes_{\co_{X'}}\co_{\oX'}.\label{crtf1e}
\end{eqnarray}
On note $\cL_{\oX'}^\gp$ le groupe de $\oX'_\et$ associé au monoïde $\cL_{\oX'}$ et 
\begin{equation}\label{crtf1f}
d\log \colon \cL_{\oX'}\rightarrow \tOmega^1_{\oX'/\oX}
\end{equation}
la dérivation logarithmique universelle. Celle-ci induit un morphisme $\co_{\oX'}$-linéaire et surjectif
\begin{equation}\label{crtf1g}
\cL_{\oX'}^\gp\otimes_\mZ\co_{\oX'}\rightarrow \tOmega^1_{\oX'/\oX}.
\end{equation}

Pour tout entier $n\geq 1$, on pose 
\begin{equation}\label{crtf1h}
\tOmega^1_{\oX'_n/\oS_n}=\tOmega^1_{X'/S}\otimes_{\co_{X'}}\co_{\oX'_n} \ \ \ {\rm et}\ \ \ 
\tOmega^1_{\oX'_n/\oX_n}=\tOmega^1_{X'/X}\otimes_{\co_{X'}}\co_{\oX'_n},
\end{equation}
que l'on considère aussi comme des faisceaux de $X'_{s,\zar}$ ou $X'_{s,\et}$, selon le contexte (cf. \ref{mtfla3}). 
On note $\oa'\colon X'_s\rightarrow \oX'$ l'immersion fermée canonique.
Compte tenu de \eqref{notconv12e}, le morphisme \eqref{crtf1g}
induit un morphisme $\co_{\oX'_n}$-linéaire et surjectif de $X'_{s,\et}$
\begin{equation}\label{crtf1j}
\oa'^*(\cL_{\oX'}^\gp)\otimes_\mZ\co_{\oX'_n}\rightarrow \tOmega^1_{\oX'_n/\oX_n}.
\end{equation}

\subsection{}\label{crtf2}
Nous considérons pour $(X',\cM_{X'})$ les objets analogues à ceux associés à $(X,\cM_{X})$ introduits dans \ref{tfkum6}, 
que nous équipons d'un exposant $^\prime$; nous les explicitons ici pour la commodité du lecteur.
Pour tout entier $n\geq 0$,  on désigne par $\cQ'_n$ le monoïde de $\tE'$ \eqref{mtfla9} défini par le diagramme cartésien 
de la catégorie des monoïdes de $\tE'$ 
\begin{equation}\label{crtf2a}
\xymatrix{
{\cQ'_n}\ar[d]\ar[r]&{\sigma'^*(\hbar'_*(\cL_{\oX'}))}\ar[d]^\ell\\
{\ocB'}\ar[r]^{\nu^n}&{\ocB'}}
\end{equation}
où $\nu$ est l'endomorphisme d'élévation à la puissance $p$-ième de $\ocB'$ \eqref{mtfla10a}
et $\ell$ est le composé des homomorphismes canoniques 
\begin{equation}\label{crtf2b}
\xymatrix{
{\sigma'^*(\hbar'_*(\cL_{\oX'}))}\ar[r]&{\sigma'^*(\hbar'_*(\co_{\oX'}))}\ar[r]&{\ocB'}},
\end{equation}
où la première flèche est induite par le morphisme structural 
$\oalpha'\colon\cL_{\oX'}\rightarrow \co_{\oX'}$ et la seconde flèche est l'homomorphisme \eqref{mtfla10b}.
Les monoïdes $(\cQ'_n)_{n\in \mN}$ forment naturellement un système projectif.
Le diagramme \eqref{crtf2a} induit par image inverse par le plongement
$\delta'$ \eqref{mtfla9e} un diagramme cartésien 
de la catégorie des monoïdes de $\tE'_s$
\begin{equation}\label{crtf2c}
\xymatrix{
{\delta'^*(\cQ'_n)}\ar[d]\ar[r]&{\sigma'^*_s(\oa'^*(\cL_{\oX'}))}\ar[d]^{\delta'^*(\ell)}\\
{\delta'^*(\ocB')}\ar[r]^{\delta'^*(\nu^n)}&{\delta'^*(\ocB')}}
\end{equation}

D'après \ref{tfkum8}, la suite d'homomorphismes canoniques
\begin{equation}\label{crtf2d}
0\rightarrow \mu_{p^n,\tE'_s}\rightarrow \delta'^*(\cQ'^\gp_n)\rightarrow \sigma'^*_s(\oa'^*(\cL^\gp_{\oX'}))\rightarrow 0
\end{equation}
est exacte.

\subsection{}\label{crtf3}
Reprenons les notations de \ref{definf3} et \ref{ktfr24}. 
Pour tout entier $n\geq 1$, on note 
\begin{equation}\label{crtf3a}
\partial'_n\colon \pi_s^*(\oa'^*(\cL^\gp_{\oX'}))\rightarrow \rR^1\tau_{s*}(\mu_{p^n,\tE'_s})
\end{equation}
l'homomorphisme composé du morphisme d'adjonction 
$\pi_s^*(\oa'^*(\cL^\gp_{\oX'}))\rightarrow \tau_{s*}(\tau^*_s(\pi_s^*(\oa'^*(\cL^\gp_{\oX'}))))$
et du cobord de la suite exacte \eqref{crtf2d} en tenant compte de l'isomorphisme $\sigma'_s\simeq \pi_s \tau_s$ \eqref{ktfr21e}. 
Celui-ci induit deux morphismes  $\ocB^!_n$-linéaires
\begin{eqnarray}
\pi_s^*(\oa'^*(\cL^\gp_{\oX'}))\otimes_\mZ\ocB^!_n&\rightarrow& \rR^1\tau_{n*}(\ocB'_n(1)),\label{crtf3b}\\
\pi_s^*(\oa'^*(\cL^\gp_{\oX'}))\otimes_\mZ\ocB^!_n&\rightarrow& \rR^1\tau_{n*}(\xi\ocB'_n),\label{crtf3c}
\end{eqnarray}
le second étant défini au moyen du composé des morphismes canoniques \eqref{definf17c}
\begin{equation}\label{crtf3d}
\co_C(1)\stackrel{\sim}{\rightarrow}p^{\frac{1}{p-1}}\xi\co_C\rightarrow \xi \co_C.
\end{equation}

\begin{teo}\label{crtf4}
Soit $n$ un entier $\geq 1$.  
\begin{itemize}
\item[{\rm (i)}] Le morphisme \eqref{crtf3c} se factorise à travers le morphisme
surjectif 
\begin{equation}\label{crtf4a}
\pi_s^*(\oa'^*(\cL_{\oX'}^\gp))\otimes_\mZ\ocB^!_n\rightarrow \pi_n^*(\tOmega^1_{\oX'_n/\oX_n})
\end{equation}
déduit de \eqref{crtf1j}, et il induit un morphisme $\ocB^!_n$-linéaire de $\tG_s$ 
\begin{equation}\label{crtf4b}
\pi_n^*(\xi^{-1}\tOmega^1_{\oX'_n/\oX_n})\rightarrow \rR^1\tau_{n*}(\ocB'_n),
\end{equation}
où $\pi_n^*$ désigne l'image inverse par le morphisme de topos annelés \eqref{ktfr24e}.
\item[{\rm (ii)}] Il existe un et un unique homomorphisme de $\ocB^!_n$-algèbres graduées de $\tG_s$
\begin{equation}\label{crtf4c}
\wedge (\pi_n^*(\xi^{-1}\tOmega^1_{\oX'_n/\oX_n}))\rightarrow \oplus_{i\geq 0}\rR^i\tau_{n*}(\ocB'_n)
\end{equation}
dont la composante en degré un est le morphisme \eqref{crtf4b}. De plus, 
son noyau est annulé par $p^{\frac{2r}{p-1}}\fm_\oK$ et son conoyau est annulé par $p^{\frac{2r+1}{p-1}}\fm_\oK$, 
où $r=\dim(X'/X)$.  
\end{itemize}
\end{teo}

La preuve de ce théorème sera donnée dans \ref{crtf18} après quelques énoncés préliminaires.

\subsection{}\label{crtf5}
Soient $(\oy'\rightsquigarrow \ox')$ un point de $X'_\et\gtimes_{X'_\et}\oX'^\rhd_\et$ \eqref{tfr21} 
tel que $\ox'$ soit au-dessus de $s$, $\uX'$ le localisé strict de $X'$ en $\ox'$.
On note encore $\ox'$ le point géométrique $\oa'(\ox')$ de $\oX'$ \eqref{mtfla3}. 
D'après (\cite{agt} III.3.7), $\uoX'$ est normal et strictement local (et en particulier intègre); 
on peut donc l'identifier au localisé strict de $\oX'$ en $\ox'$. 
Le $X'$-morphisme $\oy'\rightarrow \uX'$
définissant $(\oy'\rightsquigarrow \ox')$ se relève en un $\oX'^\rhd$-morphisme $\oy'\rightarrow \uoX'^\rhd$ et 
induit donc un point géométrique de $\uoX'^\rhd$ que l'on note aussi (abusivement) $\oy'$.
On pose $\uDelta'=\pi_1(\uoX'^\rhd,\oy')$. 

Posons $\ox=g(\ox')$ et $\oy=\upgamma(\oy')$ \eqref{ktfr1b} qui sont donc des points géométriques de $X$ et 
$\oX^\circ$ respectivement,
et notons $(\oy\rightsquigarrow \ox')$ l'image de $(\oy'\rightsquigarrow \ox')$ par le morphisme canonique 
\begin{equation}\label{crtf5b}
X'_\et\gtimes_{X'_\et}\oX'^\rhd_\et\rightarrow X'_\et\gtimes_{X_\et}\oX^\circ_\et.
\end{equation}
On désigne par $\uX$ le localisé strict de $X$ en $\ox$. 
On note encore $\ox$ le point géométrique $\oa(\ox)$ de $\oX$ \eqref{mtfla3}. 
D'après (\cite{agt} III.3.7), $\uoX$ est normal et strictement local (et en particulier intègre); 
on peut donc l'identifier au localisé strict de $\oX$ en $\ox$. 
Le $X$-morphisme $\oy\rightarrow \uX$
définissant $(\oy\rightsquigarrow \ox')$ se relève en un $\oX^\circ$-morphisme $\oy\rightarrow \uoX^\circ$ et 
induit donc un point géométrique de $\uoX^\circ$ que l'on note aussi (abusivement) $\oy$. On pose $\uDelta=\pi_1(\uoX^\circ,\oy)$.

On désigne par 
\begin{eqnarray}
\varphi'_{\ox'}\colon \tE'&\rightarrow& \uoX'^\rhd_\fet,\\
\phi_{\ox'}\colon \tG&\rightarrow& \uoX^\circ_\fet, 
\end{eqnarray} 
les foncteurs canoniques définis dans \eqref{tf10b} et \eqref{tfr29c}, 
par $\bB_{\uDelta'}$ (resp. $\bB_{\uDelta}$) le topos classifiant du groupe profini $\uDelta'$ (resp. $\uDelta$) et par 
\begin{eqnarray}
\psi'_{\oy'}\colon \uoX'^\rhd_\fet&\stackrel{\sim}{\rightarrow}&\bB_{\uDelta'},\\
\psi_{\oy}\colon \uoX^\circ_\fet&\stackrel{\sim}{\rightarrow}&\bB_{\uDelta},
\end{eqnarray}
les foncteurs fibres \eqref{notconv11c}. On note
\begin{equation}\label{crtf5c}
\uupgamma\colon \uoX'^\rhd\rightarrow \uoX^\circ
\end{equation} 
le morphisme induit par $g$. Comme $\uupgamma(\oy')=\oy$, celui-ci induit un homomorphisme
\begin{equation}
\uDelta'\rightarrow \uDelta.
\end{equation}
On désigne par $\uPi$ son noyau.

\begin{prop}\label{crtf6}
Sous les hypothèses de \ref{crtf5}, 
pour tout groupe abélien $F$ de $\uoX'^\rhd_\fet$ et tout entier $i\geq 0$, 
on a un isomorphisme canonique, fonctoriel et $\uDelta$-équivariant
\begin{equation}\label{crtf6a}
\psi_\oy(\rR^i\uupgamma_{\fet*}(F))\stackrel{\sim}{\rightarrow} \rH^i(\uPi,\psi'_{\oy'}(F)).
\end{equation}
\end{prop}

En vertu de (\cite{kato1} 4.5, \cite{tsuji4} II.4.2 et \cite{gr1} 12.7.8(iii)), 
le morphisme $g\colon X'\rightarrow X$ est plat à fibres géométriquement réduites. 
Par suite, d'après \ref{csc2} appliqué au morphisme $\ogg\colon \oX'\rightarrow \oX$,
pour tout revêtement étale $V\rightarrow \uoX^\circ$ tel que $V$ soit connexe, le schéma $V\times_{\uX}\uX'$ est connexe. 
Comme $\uoX'$ est normal, $V\times_{\uX}\uX'$ est normal et est donc intègre. 
Par suite, $V\times_{\uoX^\circ}\uoX'^\rhd$ est connexe. 
On en déduit que pour tout groupe abélien $F$ de $\uoX'^\rhd_\fet$, 
on a un isomorphisme canonique, fonctoriel et $\uDelta$-équivariant,
\begin{equation}\label{crtf6b}
\psi_\oy(\uupgamma_{\fet*}(F))\stackrel{\sim}{\rightarrow} \Gamma(\uPi,\psi'_{\oy'}(F)),
\end{equation}
d'où la proposition. 

\begin{cor}\label{crtf7}
Sous les hypothèses de \ref{crtf5}, 
pour tout groupe abélien $F$ de $\tE'$ et tout entier $i\geq 0$, on a un isomorphisme canonique, fonctoriel et $\uDelta$-équivariant
\begin{equation}\label{crtf7a}
\psi_\oy(\phi_{\ox'}(\rR^i\tau_*(F)))\stackrel{\sim}{\rightarrow} \rH^i(\uPi,\psi'_{\oy'}(\varphi'_{\ox'}(F))).
\end{equation}
\end{cor}

Cela résulte de \ref{ktfr18} et \ref{crtf6}.

\begin{cor}\label{crtf8}
Sous les hypothèses de \ref{crtf5}, 
pour tout groupe abélien $F$ de $\tE'$ et tout entier $i\geq 0$, on a un isomorphisme canonique et fonctoriel 
\begin{equation}\label{crtf8a}
(\rR^i\tau_*(F))_{\varrho(\oy\rightsquigarrow \ox')}\stackrel{\sim}{\rightarrow} \rH^i(\uPi,\psi'_{\oy'}(\varphi'_{\ox'}(F))).
\end{equation}
\end{cor}

Cela résulte de \ref{crtf7} et \eqref{ktfr28a}.

\subsection{}\label{crtf9}
Conservons les hypothèses et notations de \ref{crtf5}. 
Compte tenu de \eqref{TFA14b}, \eqref{TFA5b} et \eqref{TFA12g}, l'homomorphisme $\psi'_{\oy'}(\varphi'_{\ox'}(\ell))$ \eqref{crtf2b} 
s'identifie à  l'homomorphisme composé
\begin{equation}\label{crtf9a}
\xymatrix{
{\cL_{\oX',\ox'}}\ar[r]^{\oalpha'_{\ox'}}&{\co_{\oX',\ox'}}\ar[r]&{\oR'^{\oy'}_{\uX'}}},
\end{equation}
où la première flèche est la fibre en $\ox'$ du morphisme structural $\oalpha'\colon \cL_{\oX'}\rightarrow \co_{\oX'}$, 
$\oR'^{\oy'}_{\uX'}$ est l'algèbre définie dans \eqref{ktfr10a} et la seconde flèche est l'homomorphisme canonique.
Pour tout entier $n\geq 0$, l'image du diagramme \eqref{crtf2a} par le foncteur $\psi'_{\oy'}\circ \varphi'_{\ox'}$ 
s'identifie donc à un diagramme cartésien de la catégorie des monoïdes de  $\bB_{\uDelta'}$
\begin{equation}\label{crtf9b}
\xymatrix{
{\psi'_{\oy'}(\varphi_{\ox'}(\cQ'_n))}\ar[r]\ar[d]&{\cL_{\oX',\ox'}}\ar[d]^{\psi'_{\oy'}(\varphi'_{\ox'}({\ell}))}\\
{\oR'^{\oy'}_{\uX'}}\ar[r]^{\nu^n}&{\oR'^{\oy'}_{\uX'}}}
\end{equation}
où on a encore noté $\nu$ l'endomorphisme d'élévation à la puissance $p$-ième de $\oR'^{\oy'}_{\uX'}$.

Pour tout $t\in \cL_{\oX',\ox'}$, la section $\oalpha'_{\ox'}(t)$  est inversible sur $\uoX'^\rhd$ d'après \ref{tfkum3}(iii). 
Le $\uoX'$-schéma 
\begin{equation}\label{crtf9c}
\cT'_n(t)=\Spec(\co_{\uoX'}[T]/(T^{p^n}-\oalpha'_{\ox'}(t)))
\end{equation}
est donc fini et étale sur $\uoX'^\rhd$. En particulier, $\cT'^\rhd_n(t)=\cT'_n(t)\times_{X'}X'^\rhd$ 
induit un $\mu_{p^n}(\co_\oK)$-torseur 
de $\uoX'^\rhd_\fet$. On désigne par $\cT'_n(t)_{\oy'}$ la fibre de $\cT'_n(t)$ au-dessus de $\oy'$ et par $\cQ'_n(t)$ la fibre 
de la projection canonique $\psi'_{\oy'}(\varphi'_{\ox'}(\cQ'_n))\rightarrow \cL_{\oX',\ox'}$ \eqref{crtf9b} au-dessus de $t$.
On a un isomorphisme canonique de $\mu_{p^n}(\co_\oK)$-torseurs de $\bB_{\uDelta'}$
\begin{equation}\label{crtf9d}
\cT'_n(t)_{\oy'}\stackrel{\sim}{\rightarrow}\cQ'_n(t).
\end{equation}

\begin{lem}\label{crtf10}
Sous les hypothèses de \ref{crtf9}, pour tout entier $n\geq 0$, l'application 
\begin{equation}\label{crtf10a}
\ttt_n\colon \cL_{\oX',\ox'}\rightarrow \rH^1(\uPi,\mu_{p^n}(\co_\oK)), \ \ \ t\mapsto [\cT'^\rhd_n(t)],
\end{equation}
est un homomorphisme de monoïdes, et le diagramme \eqref{tfkum4}
\begin{equation}\label{crtf10b}
\xymatrix{
{\cL_{\oX',\ox'}}\ar[rr]^-(0.5){\ttt_n}\ar[d]&&{\rH^1(\uPi,\mu_{p^n}(\co_\oK))}\ar[r]_-(0.5)\sim^-(0.5)u&
{\rH^1(\uPi,\psi'_{\oy'}(\varphi'_{\ox'}(\delta'_*(\mu_{p^n,\tE'_s}))))}\\
{\cL_{\oX',\ox'}^\gp}\ar[rr]^-(0.5){\partial'_{n,\varrho(\oy\rightsquigarrow \ox')}}&&
{\rR^1\tau_{s*}(\mu_{p^n,\tE'_s})_{\varrho(\oy\rightsquigarrow \ox')}}\ar[r]_-(0.5)\sim^-(0.5)w&
{\rR^1\tau_*(\delta'_*(\mu_{p^n,\tE'_s}))_{\varrho(\oy\rightsquigarrow \ox')}}\ar[u]_v}
\end{equation}
où $\partial'_{n,\varrho(\oy\rightsquigarrow \ox')}$ est la fibre de  $\partial'_n$ \eqref{crtf3a} en 
$\varrho(\oy\rightsquigarrow \ox')$, 
$u$ est induit par l'isomorphisme canonique $\varphi'_{\ox'}\stackrel{\sim}{\rightarrow} \varphi'_{\ox'}\delta'_*\delta'^*$ \eqref{TFA15a},
$v$ est l'isomorphisme \eqref{crtf8a}, 
$w$ est induit par l'isomorphisme \eqref{ktfr21l} et la flèche non libellée est l'homomorphisme canonique, est commutatif. 
\end{lem}

Considérons la suite exacte de groupes abéliens de $\tE'$ 
\begin{equation}\label{crtf10c}
0\rightarrow \delta'_*(\mu_{p^n,\tE'_s})\rightarrow \delta'_*(\delta'^*(\cQ'^\gp_n))\rightarrow 
\delta'_*(\sigma'^*_s(\oa'^*(\cL_{\oX'}^\gp)))\rightarrow 0
\end{equation}
déduite de la suite exacte \eqref{crtf2d} par image directe par le plongement $\delta'$. 
D'après \ref{ktfr19}(ii), le diagramme 
\begin{equation}\label{crtf10d}
\xymatrix{
{\tau_*(\delta'_*(\sigma'^*_s(\oa'^*(\cL^\gp_{\oX'}))))_{\varrho(\oy\rightsquigarrow \ox')}}\ar[r]^-(0.5){v'}_-(0.5)\sim\ar[d]_a&
{\rH^0(\uPi,\psi'_{\oy'}(\varphi'_{\ox'}(\delta'_*(\sigma'^*_s(\oa'^*(\cL^\gp_{\oX'}))))))}\ar[d]^b\\
{\rR^1\tau_*(\delta'_*(\mu_{p^n,\tE'_s}))_{\varrho(\oy\rightsquigarrow \ox')}}
\ar[r]^-(0.5)v_-(0.5)\sim&{\rH^1(\uPi,\psi'_{\oy'}(\varphi'_{\ox'}(\delta'_*(\mu_{p^n,\tE'_s}))))}}
\end{equation}
où les flèches horizontales sont les isomorphismes \eqref{crtf8a} et les flèches verticales sont les bords 
des suites exactes longues de cohomologie, est commutatif. Par ailleurs, d'après \eqref{tfr21a} et \eqref{ktfr21b},
on a un isomorphisme canonique
\begin{equation}\label{crtf10e}
\pi_s^*(\oa'^*(\cL^\gp_{\oX'}))_{\varrho(\oy\rightsquigarrow \ox')}\stackrel{\sim}{\rightarrow} \cL^\gp_{\oX',\ox'}.
\end{equation}
Le diagramme 
\begin{equation}\label{crtf10f}
\xymatrix{
{\cL^\gp_{\oX',\ox'}}\ar[r]^-(0.5)i\ar[rd]_{\partial'_{n,\varrho(\oy\rightsquigarrow \ox')}}
&{\tau_{s*}(\sigma'^*_s(\oa'^*(\cL^\gp_{\oX'})))_{\varrho(\oy\rightsquigarrow \ox')}}
\ar[r]^-(0.5){w'}_-(0.5)\sim\ar[d]^c&{\tau_*(\delta'_*(\sigma'^*_s(\oa'^*(\cL^\gp_{\oX'}))))_{\varrho(\oy\rightsquigarrow \ox')}}\ar[d]^a\\
&{\rR^1\tau_{s*}(\mu_{p^n,\tE'_s})_{\varrho(\oy\rightsquigarrow \ox')}}\ar[r]^-(0.5){w}_-(0.5)\sim&
{\rR^1\tau_*(\delta'_*(\mu_{p^n,\tE'_s}))_{\varrho(\oy\rightsquigarrow \ox')}}}
\end{equation}
où $i$ est induit par le morphisme d'adjonction $\id\rightarrow \tau_{s*}\tau_s^*$ et l'isomorphisme $\sigma'^*_s\simeq \tau_s^*\pi_s^*$
\eqref{ktfr21e}, $c$ est le bord de la suite exacte longue de cohomologie déduite de la suite \eqref{crtf2d}
et $w$ et $w'$ sont induits par l'isomorphisme \eqref{ktfr21l}, est commutatif. 

D'après \eqref{TFA15a}, \eqref{TFA14b} et \eqref{TFA66c}, 
l'image de \eqref{crtf10c} par le foncteur $\psi'_{\oy'}\circ \varphi'_{\ox'}$ s'identifie à une suite exacte de 
groupes abéliens de $\bB_{\uDelta'}$
\begin{equation}\label{crtf10g}
0\rightarrow \mu_{p^n}(\co_\oK)\rightarrow \psi'_{\oy'}(\varphi'_{\ox'}(\cQ'^\gp_n))\rightarrow \cL_{\oX',\ox'}^\gp\rightarrow 0. 
\end{equation}
Celle-ci induit un homomorphisme 
\begin{equation}\label{crtf10h}
\ttt'_n\colon \cL^\gp_{\oX',\ox'}\rightarrow\rH^1(\uPi,\mu_{p^n}(\co_\oK)),
\end{equation}
qui s'identifie donc à l'homomorphisme $b$
\begin{equation}\label{crtf10i}
\xymatrix{
{\cL^\gp_{\oX',\ox'}}\ar[d]_{\ttt'_n}\ar[r]^-(0.5){u'}_-(0.5)\sim&
{\rH^0(\uPi,\psi'_{\oy'}(\varphi'_{\ox'}(\delta'_*(\sigma'^*_s(\oa'^*(\cL^\gp_{\oX'}))))))}\ar[d]^b\\
{\rH^1(\uPi,\mu_{p^n}(\co_\oK))}\ar[r]^-(0.5)u_-(0.5)\sim&
{\rH^1(\uPi,\psi'_{\oy'}(\varphi'_{\ox'}(\delta'_*(\mu_{p^n,\tE'_s}))))}}
\end{equation}
On vérifie que $u'^{-1}\circ v'\circ w'\circ i$ est l'identité. 

Il résulte de \ref{tfkum9} et (\cite{ogus} I.1.1.6) que le diagramme canonique
\begin{equation}\label{crtf10j}
\xymatrix{
{\varphi'_{\ox'}(\cQ'_n)}\ar[d]\ar[r]&{\cL_{\oX',\ox'}}\ar[d]\\
{\varphi'_{\ox'}(\cQ'^\gp_n)}\ar[r]&{\cL^\gp_{\oX',\ox'}}}
\end{equation}
est cartésien. Compte tenu de \eqref{tfkum12h}, on en déduit  
que $\ttt_n$ est le composé de $\ttt'_n$ et de l'homomorphisme canonique $\cL_{\oX',\ox'}\rightarrow \cL_{\oX',\ox'}^\gp$;
d'où la proposition.

\subsection{}\label{crtf11}
Nous développons pour $f'$ des constructions analogues à celles développées pour $f$ dans \ref{tfkum14}, ainsi que des variantes relatives.
Conservons les hypothèses et notations de \ref{crtf5}; supposons, de plus, les conditions suivantes remplies~:
\begin{itemize}
\item[(i)] Les schémas $X=\Spec(R)$ et $X'=\Spec(R')$ sont affines et connexes, 
et le morphisme $g$ \eqref{ktfr1a} admet une carte relativement adéquate \eqref{mtfla5}.
Ces conditions correspondent à celles fixées dans \ref{eccr1}.
\item[(ii)] Il existe une carte fine et saturée $M\rightarrow \Gamma(X',\cM_{X'})$ pour $(X',\cM_{X'})$
induisant un isomorphisme 
\begin{equation}\label{crtf11a}
M\stackrel{\sim}{\rightarrow} \Gamma(X',\cM_{X'})/\Gamma(X',\co^\times_{X'}).
\end{equation}
\end{itemize}
Posons \eqref{mtfla2e}
\begin{equation}\label{crtf11b}
\tOmega^1_{R'/\co_K}=\tOmega^1_{X'/S}(X') \ \ \ {\rm et}\ \ \ \ \tOmega^1_{R'/R}=\tOmega^1_{X'/X}(X').
\end{equation}

Les schémas $\oX$ et $\oX'$ étant localement irréductibles d'après \ref{cad7}(iii),  
ils sont les sommes des schémas induits sur leurs composantes irréductibles. 
On note $\oX^\star$ (resp. $\oX'^\star$)
la composante irréductible de $\oX$ (resp. $\oX'$) contenant $\oy$ (resp. $\oy'$). 
De même, $\oX^\circ$ (resp. $\oX'^\rhd$) est la somme des schémas induits sur ses composantes irréductibles
et $\oX^{\star \circ}=\oX^\star\times_{X}X^\circ$ (resp. $\oX'^{\star \rhd}=\oX'^\star\times_{X'}X'^\rhd$) 
est la composante irréductible de $\oX^\circ$ (resp. $\oX'^\rhd$) contenant $\oy$ (resp. $\oy'$).
On pose $\Delta=\pi_1(\oX^{\star\circ},\oy)$ et $\Delta'=\pi_1(\oX'^{\star\rhd},\oy')$. Comme 
$\upgamma(\oy')=\oy$ \eqref{ktfr1b}, $\upgamma$ induit un homomorphisme $\Delta'\rightarrow \Delta$. 
On note $\Pi$ son noyau.

On désigne par $\oR^\oy_X$ et $\oR^{!\oy}_{X'\rightarrow X}$ les représentations de $\Delta$ définies dans 
\eqref{ktfr6f} et \eqref{ktfr6d}, respectivement, 
et par $\oR'^{\oy'}_{X'}$ et $\oR^{\intern\oy'}_{X'\rightarrow X}$ les représentations de $\Delta'$ définies dans 
\eqref{ktfr7b} et \eqref{ktfr7f}, respectivement. 
On a alors des homomorphismes canoniques $\Delta'$-équivariants \eqref{ktfr7i}
\begin{equation}\label{crtf11c}
\oR^\oy_X\rightarrow \oR^{!\oy}_{X'\rightarrow X} \rightarrow \oR^{\intern\oy'}_{X'\rightarrow X} \rightarrow \oR'^{\oy'}_{X'}.
\end{equation}

Reprenons les notations de \ref{taht}.
On munit  $\coX'=X'\times_S\coS$ \eqref{hght1c} de la structure logarithmique $\cM_{\coX'}$ image inverse de $\cM_{X'}$
et on se donne une $(\tS,\cM_{\tS})$-déformation lisse $(\tX',\cM_{\tX'})$ de $(\coX',\cM_{\coX'})$ (cf. \ref{definf12}). 
Comme $f'$ est lisse et que $\coX'$ est affine, 
une telle déformation existe et est unique à isomorphisme près en vertu de (\cite{kato1} 3.14).
On désigne par $\hoR'^{\oy'}_{X'}$ le séparé complété $p$-adique de $\oR'^{\oy'}_{X'}$ et par
\begin{equation}\label{crtf11d}
0\rightarrow \hoR'^{\oy'}_{X'}\rightarrow \cF'\rightarrow \xi^{-1} \tOmega^1_{R'/\co_K}\otimes_{R'}\hoR'^{\oy'}_{X'}\rightarrow 0
\end{equation}
l'extension de Higgs-Tate associée à $(\tX',\cM_{\tX'})$ \eqref{taht5e};  
c'est une extension de $\hoR'^{\oy'}_{X'}$-représentations continues de $\Delta'$. 
Pour tout entier $n\geq 0$, on désigne par  
\begin{equation}\label{crtf11e} 
0\rightarrow \xi(\oR'^{\oy'}_{X'}/p^n\oR'^{\oy'}_{X'})\rightarrow \xi(\cF'/p^n\cF')\rightarrow \tOmega^1_{R'/\co_K}\otimes_{R'}
(\oR'^{\oy'}_{X'}/p^n\oR'^{\oy'}_{X'})\rightarrow 0
\end{equation}
la suite exacte déduite de \eqref{crtf11d}, par
\begin{equation}\label{crtf11f}
A_n\colon \Gamma(\oX'^\star,\tOmega^1_{\oX'_n/\oS_n})
\rightarrow \rH^1(\Pi,\xi(\oR'^{\oy'}_{X'}/p^n\oR'^{\oy'}_{X'}))
\end{equation}
le morphisme $\Gamma(\oX'^\star,\co_{\oX'})$-linéaire induit, 
et par $Q'_n$ le produit fibré du diagramme d'homomorphismes de monoïdes de $\Delta'$
\begin{equation}\label{crtf11g}
\xymatrix{
&{\Gamma(X',\cM_{X'})}\ar[d]\\
{\oR'^{\oy'}_{X'}}\ar[r]^{\nu^n}&{\oR'^{\oy'}_{X'}}}
\end{equation}
où on a encore noté $\nu$ l'endomorphisme d'élévation à la puissance $p$-ième de $\oR'^{\oy'}_{X'}$.
On observera que $A_n$ ne dépend pas du choix de $(\tX',\cM_{\tX'})$ (\cite{agt} II.10.10). 
L'homomorphisme canonique $\Gamma(X',\cM_{X'})\rightarrow \oR'^{\oy'}_{X'}$ vérifie les hypothèses de \ref{tfkum5}. 
En effet, l'homomorphisme $\Gamma(X',\cM_{X'})\rightarrow R'$ est injectif d'après (\cite{agt} III.4.2(iv)),
et comme $X'$ est connexe,  l'homomorphisme $R'\rightarrow \oR'^{\oy'}_{X'}$ est injectif. 
Par suite, le monoïde $Q'_n$ est intègre et la projection canonique $Q'_n\rightarrow \Gamma(X',\cM_{X'})$ 
identifie $\Gamma(X',\cM_{X'})$ au quotient de $Q'_n$ par le sous-monoïde $\mu_{p^n}(\co_\oK)$ d'après  \ref{tfkum5}. 
On a donc une suite exacte de $\mZ[\Delta']$-modules
\begin{equation}\label{crtf11h}
0\rightarrow \mu_{p^n}(\co_\oK)\rightarrow Q'^\gp_n\rightarrow \Gamma(X',\cM_{X'})^\gp\rightarrow 0.
\end{equation}
Elle induit un homomorphisme 
\begin{equation}\label{crtf11i}
B_n\colon \Gamma(X',\cM_{X'})\rightarrow \rH^1(\Pi,\mu_{p^n}(\co_\oK)).
\end{equation}
On observera que le diagramme canonique
\begin{equation}\label{crtf11j}
\xymatrix{
{Q'_n}\ar[d]\ar[r]&{\Gamma(X',\cM_{X'})}\ar[d]\\
{Q'^\gp_n}\ar[r]&{\Gamma(X',\cM_{X'})^\gp}}
\end{equation}
est cartésien (\cite{ogus} I 1.1.6). On note encore 
\begin{equation}\label{crtf11k}
d\log \colon \cM_{X'}\rightarrow \tOmega^1_{X'/S}
\end{equation}
la dérivation logarithmique universelle et 
\begin{equation}\label{crtf11l}
\log([\ ])\colon \mZ_p(1)\rightarrow \xi\hoR'^{\oy'}_{X'}
\end{equation}
l'homomorphisme induit par l'homomorphisme $\log([\ ])$ défini dans  \eqref{definf17b}.

\begin{prop}\label{crtf12}
Sous les hypothèses de \ref{crtf11}, pour tout entier $n\geq 0$, le diagramme 
\begin{equation}\label{crtf12a}
\xymatrix{
{\Gamma(\oX'^{\star},\tOmega^1_{\oX'_n/\oS_n})}\ar[r]^-(0.5){A_n}\ar[d]_u&{\rH^1(\Pi,\xi(\oR'^{\oy'}_{X'}/p^n\oR'^{\oy'}_{X'}))}\\
{\tOmega^1_{R'/R}\otimes_{R'}(\oR^{\intern\oy'}_{X'\rightarrow X}/p^n\oR^{\intern\oy'}_{X'\rightarrow X})}\ar[r]^-(0.5)v&
{\rH^1(\Pi,(\oR'^{\oy'}_{X'}/p^n\oR'^{\oy'}_{X'})(1))}\ar[u]_{-w}}
\end{equation}
où $u$ est le morphisme canonique, $v$ est la composante en degré un du morphisme \eqref{eccr51a} et 
$w$ est induit par le morphisme \eqref{crtf11l}, est commutatif. 
\end{prop}

Pour alléger les notations, nous adoptons celles introduites dans la section \ref{eccr} (cf. \ref{eccr34} pour un résumé). 
On observera que $\oR'^{\oy'}_{X'}$ (resp. $\oR^{\intern\oy'}_{X'\rightarrow X}$) correspond à l'algèbre $\oR'$ (resp. $\oR^\intern$).
Il suffit alors de montrer que le diagramme 
\begin{equation}\label{crtf12b}
\xymatrix{
{\tOmega^1_{R'/\co_K}\otimes_{R'}(R'_1/p^nR'_1)}\ar[r]^-(0.5){A_n}\ar[d]_{u'}&{\rH^1(\Pi,\xi(\oR'/p^n\oR'))}\\
{\tOmega^1_{R'/R}\otimes_{R'}(R^\intern_\infty/p^nR^\intern_\infty)}\ar[r]^-(0.5){v'}&
{\Hom(\fS_{\infty},(R^\intern_\infty/p^nR^\intern_\infty)(1))}\ar[u]_{-w'}}
\end{equation}
où $u'$ est le morphisme canonique, $v'$ est induit par l'isomorphisme \eqref{eccr32b} et 
$w'$ est induit par le morphisme \eqref{crtf3d}, est commutatif. D'après (\cite{agt} II.10.10), 
$A_n$ ne dépend pas de la déformation $(\tX',\cM_{\tX'})$. 
On peut donc se borner au cas où $(\tX',\cM_{\tX'})$
est la déformation définie par la carte adéquate de $f'$ fixée dans \ref{crtf11}(i) (cf. \cite{agt} II.10.13).
La même carte détermine une section du torseur de Higgs-Tate sur $\Spec(\hoR')$ \eqref{taht5},
et par suite un scindage de l'extension \eqref{crtf11d} (cf. \cite{agt} II.10.11). 
L'assertion résulte alors de (\cite{agt} II.10.14).

\begin{cor}\label{crtf13}
Sous les hypothèses de \ref{crtf11}, pour tout entier $n\geq 0$, le morphisme $A_n$ \eqref{crtf11f} se factorise à travers le morphisme
surjectif
\begin{equation}
\Gamma(\oX'^{\star},\tOmega^1_{\oX'_n/\oS_n})\rightarrow \Gamma(\oX'^{\star},\tOmega^1_{\oX'_n/\oX_n}),
\end{equation}
et il induit un morphisme que l'on note encore 
\begin{equation}\label{crtf13a}
A_n\colon \Gamma(\oX'^\star,\tOmega^1_{\oX'_n/\oX_n})
\rightarrow \rH^1(\Pi,\xi(\oR'^{\oy'}_{X'}/p^n\oR'^{\oy'}_{X'})).
\end{equation}
\end{cor}

\begin{lem}\label{crtf14}
Sous les hypothèses de \ref{crtf11}, pour tout entier $n\geq 0$, le diagramme 
\begin{equation}\label{crtf14a}
\xymatrix{
{\Gamma(X',\cM_{X'})}\ar[r]^-(0.5){B_n}\ar[d]_{d\log}&
{\rH^1(\Pi,\mu_{p^n}(\co_\oK))}\ar[d]^{-\log([\ ])}\\
{\Gamma(\oX'^\star,\tOmega^1_{\oX'_n/\oS_n})}\ar[r]^-(0.5){A_n}&
{\rH^1(\Pi,\xi(\oR'^{\oy'}_{X'}/p^n\oR'^{\oy'}_{X'}))}}
\end{equation}
est commutatif. 
\end{lem}

Cela résulte aussitôt de \ref{tfkum15} et de la fonctorialité du morphisme de restriction 
\begin{equation}
\rH^1(\pi_1(\oX'^{\star\rhd},\oy'),-)\rightarrow \rH^1(\Pi,-).
\end{equation}

\begin{prop}\label{crtf15}
Sous les hypothèses de \ref{crtf5}, pour tout entier $n\geq 1$, il existe un unique morphisme $\co_{\oX',\ox'}$-linéaire
\begin{equation}\label{crtf15a}
\phi_n\colon \tOmega^1_{\oX'_n/\oX_n,\ox'}\rightarrow \rH^1(\uPi,\xi(\oR'^{\oy'}_{\uX'}/p^n\oR'^{\oy'}_{\uX'}))
\end{equation}
où $\oR'^{\oy'}_{\uX'}$ est l'algèbre définie dans \eqref{ktfr10a},
qui s'insère dans le diagramme commutatif 
\begin{equation}\label{crtf15b}
\xymatrix{
{\cL_{\oX',\ox'}}\ar[r]^-(0.5){\ttt_n}\ar[d]_{d\log}&{\rH^1(\uPi,\mu_{p^n}(\co_\oK))}\ar[d]\ar[d]^{-\log([\ ])}\\
{\tOmega^1_{\oX'_n/\oX_n,\ox'}}\ar[r]^-(0.5){\phi_n}&{\rH^1(\uPi,\xi(\oR'^{\oy'}_{\uX'}/p^n\oR'^{\oy'}_{\uX'}))}}
\end{equation}
où $\ttt_n$ est l'homomorphisme \eqref{crtf10a}, $d\log$ est induit par la dérivation logarithmique universelle \eqref{crtf1f} et $\log([\ ])$
est induit par l'homomorphisme \eqref{crtf11l}. 
\end{prop}

En effet, l'unicité de $\phi_n$ est claire \eqref{crtf1g}. Montrons son existence. 
Pour toute extension finie $L$ de $K$ contenue dans $\oK$, on pose \eqref{hght3}
\begin{equation}
(X'_L,\cM_{X'_L})=(X',\cM_{X'})\times_{(S,\cM_S)}(S_L,\cM_{S_L})
\end{equation}
les produits étant pris dans la catégorie des schémas logarithmiques, de sorte que 
$X'_L=X'\times_SS_L$.
On note encore $\ox'$ l'image du point géométrique $\ox'$ par la projection canonique $\oX'\rightarrow X'_L$. 
D'après \eqref{tfkum3a} et \eqref{tfkum3b}, on a des isomorphismes canoniques 
\begin{eqnarray}
\cL_{\oX',\ox'}&\stackrel{\sim}{\rightarrow}&\underset{\underset{L\subset \oK}{\longrightarrow}}{\lim}\  \cM_{X'_L,\ox'},\\
\tOmega^1_{\oX'_n/\oX_n,\ox'}&\stackrel{\sim}{\rightarrow}& \underset{\underset{L\subset \oK}{\longrightarrow}}{\lim}\  
\tOmega^1_{X'_n/X_n,\ox'}\otimes_{\co_{X',\ox'}}\co_{X'_L,\ox'}.
\end{eqnarray}
Pour toute extension finie $L/K$, la projection canonique $(X'_L,\cM_{X'_L})\rightarrow (S_L,\cM_{S_L})$ est adéquate (\cite{agt} III.4.8). 
Par ailleurs, le groupe $\uPi$ et l'algèbre $\oR'^{\oy'}_{\uX'}$ ne changent pas en remplaçant $(X',\cM_{X'})$  par $(X'_L,\cM_{X'_L})$. 
On est donc réduit à montrer que la restriction de $-\log([\ ])\circ \ttt_n$ à $\cM_{X',\ox'}$ se factorise à travers l'homomorphisme 
\begin{equation}\label{crtf15c}
\cM_{X',\ox'}\rightarrow \tOmega^1_{X'_n/X_n,\ox'}
\end{equation}
induit par la dérivation logarithmique universelle
\begin{equation}\label{crtf15d}
d\log\colon \cM_{X'}\rightarrow \tOmega^1_{X'/X}.
\end{equation}
On notera que cette dernière est compatible avec l'homomorphisme $d\log$ \eqref{crtf1f}.
En vertu de (\cite{agt} II.5.17), on peut supposer les hypothèses de \ref{crtf11} remplies. 
Il suffit de montrer que  la restriction de l'homomorphisme $-\log([\ ])\circ \ttt_n$ à $\Gamma(X',\cM_{X'})$ 
se factorise à travers l'homomorphisme 
\begin{equation}\label{crtf15e}
\Gamma(X',\cM_{X'})\rightarrow \Gamma(X',\tOmega^1_{X'_n/X_n})
\end{equation}
induit par $d\log$ \eqref{crtf15d}, ce qui résulte de \ref{crtf13} et \ref{crtf14}.

\begin{cor}\label{crtf16}
Sous les hypothèses de \ref{crtf11}, pour tout entier $n\geq 1$, le diagramme 
\begin{equation}\label{crtf16a}
\xymatrix{
{\Gamma(\oX'^\star,\tOmega^1_{\oX'_n/\oX_n})}\ar[r]^-(0.5){A_n}\ar[d]&
{\rH^1(\Pi,\xi(\oR'^{\oy'}_{X'}/p^n\oR'^{\oy'}_{X'}))}\ar[d]\\
{\tOmega^1_{\oX'_n/\oX_n,\ox'}}\ar[r]^-(0.5){\phi_n}&{\rH^1(\uPi,\xi(\oR'^{\oy'}_{\uX'}/p^n\oR'^{\oy'}_{\uX'}))}}
\end{equation}
où $A_n$ est le morphisme \eqref{crtf13a} et $\phi_n$ est le morphisme \eqref{crtf15a}, est commutatif. 
\end{cor}

Cela résulte aussitôt de la preuve de \ref{crtf15}.

\begin{cor}\label{crtf17}
Les hypothèses étant celles de \ref{crtf5}, notons $\oR'^{\oy'}_{\uX'}$  et $\oR^{\intern\oy'}_{\uX'\rightarrow \uX}$ 
les anneaux définis dans \eqref{ktfr10a} et \eqref{ktfr10e} respectivement. 
Alors, pour tout entier $n\geq 1$, il existe un et un unique homomorphisme de $\oR^{\intern\oy'}_{\uX'\rightarrow \uX}$-algèbres graduées 
\begin{equation}\label{crtf17a}
\wedge \left(\xi^{-1}\tOmega^1_{\oX'_n/\oX_n,\ox'}\otimes_{\co_{\oX',\ox'}}\oR^{\intern\oy'}_{\uX'\rightarrow \uX}\right)\rightarrow 
\oplus_{i\geq 0}\rH^i(\uPi,\oR'^{\oy'}_{\uX'}/p^n\oR'^{\oy'}_{\uX'})
\end{equation}
dont la composante en degré un est induite par le morphisme $\phi_n$ \eqref{crtf15a}. De plus, 
notant $r=\dim(X'/X)$ la dimension relative de $X'$ sur $X$, 
le noyau de \eqref{crtf17a} est annulé par $p^{\frac{2r}{p-1}}\fm_\oK$ et 
son conoyau est annulé par $p^{\frac{2r+1}{p-1}}\fm_\oK$.
\end{cor}

On désigne par $\cC_{\ox'}$ la catégorie associée à $\ox'$ dans \ref{ktfr13} et par $\cC'_{\ox'}$ la sous-catégorie pleine 
formée des objets $(\ox'\rightarrow U'\rightarrow U)$ tels que la restriction 
$(U',\cM_{X'}|U')\rightarrow (U,\cM_X|U)$ du morphisme $g$ vérifie les conditions \ref{crtf11}(i)-(ii). 
Les catégories $\cC_{\ox'}$ et $\cC'_{\ox'}$ sont cofiltrantes, et le foncteur d'injection canonique
$\cC'^\circ_{\ox'}\rightarrow \cC^\circ_{\ox'}$ est cofinal d'après (\cite{agt} II.5.17) et (\cite{sga4} I 8.1.3(c)).

Soit $(\ox'\rightarrow U'\rightarrow U)$ un objet de $\cC'_{\ox'}$; on omettra $\ox'$ pour alléger les notations. 
On a un $X'$-morphisme canonique 
$\uX'\rightarrow U'$ et un $X$-morphisme canonique $\uX\rightarrow U$ qui s'insèrent dans un diagramme commutatif  
\begin{equation}\label{crtf17b}
\xymatrix{
&\uX'\ar[r]^g\ar[d]&\uX\ar[d]\\
\ox'\ar[r]\ar[ru]&U'\ar[r]&U}
\end{equation}
où on a encore noté $g$ le morphisme induit par $g$. On en déduit des morphismes $\uoX'\rightarrow \oU'$ et $\uoX\rightarrow \oU$. 
Le morphisme $\oy'\rightarrow \uoX'^\rhd$ qui définit le point $(\oy'\rightsquigarrow \ox')$ \eqref{crtf5}
induit un point géométrique de $\oU'^\rhd$ que l'on note encore $\oy'$. 
De même, le morphisme $\oy\rightarrow \uoX^\circ$ qui définit le point $(\oy\rightsquigarrow \ox)$ \eqref{crtf5}
induit un point géométrique de $\oU^\circ$ que l'on note encore $\oy$. On observera que le diagramme 
\begin{equation}\label{crtf17c}
\xymatrix{
\oy'\ar[r]\ar[d]&\oy\ar[d]\\
\oU'^\rhd\ar[r]&\oU^\circ}
\end{equation}
est commutatif. 

Les schémas $\oU$ et $\oU'$ étant localement irréductibles d'après (\cite{agt} III.3.3 et III.4.2(iii)),  
ils sont les sommes des schémas induits sur leurs composantes irréductibles. 
On note $\oU^\star$ (resp. $\oU'^\star$)
la composante irréductible de $\oU$ (resp. $\oU'$) contenant $\oy$ (resp. $\oy'$). 
De même, $\oU^\circ$ (resp. $\oU'^\rhd$) est la somme des schémas induits sur ses composantes irréductibles
et $\oU^{\star \circ}=\oU^\star\times_{X}X^\circ$ (resp. $\oU'^{\star \rhd}=\oU'^\star\times_{X'}X'^\rhd$) 
est la composante irréductible de $\oU^\circ$ (resp. $\oU'^\rhd$) contenant $\oy$ (resp. $\oy'$).
On pose $\Delta_U=\pi_1(\oU^{\star\circ},\oy)$ et $\Delta'_{U'}=\pi_1(\oU'^{\star\rhd},\oy')$. 
D'après \eqref{crtf17c}, on a un homomorphisme canonique $\Delta'_{U'}\rightarrow \Delta_U$. 
On note $\Pi_{U'\rightarrow U}$ son noyau et  
\begin{equation}\label{crtf17d}
A_{U'\rightarrow U,n}\colon \tOmega^1_{\oX'_n/\oX_n}(\oU'^\star)
\rightarrow \rH^1(\Pi_{U'\rightarrow U},\xi(\oR'^{\oy'}_{U'}/p^n\oR^{\oy'}_{U'}))
\end{equation} 
le morphisme canonique \eqref{crtf13a}; on rappelle que celui-ci ne dépend pas de 
la déformation choisie pour le définir \eqref{crtf11}  (\cite{agt} II.10.10). D'après \ref{crtf16}, le diagramme 
\begin{equation}\label{crtf17e}
\xymatrix{
{\tOmega^1_{\oX'_n/\oX_n}(\oU'^\star)}\ar[rr]^-(0.5){A_{U'\rightarrow U,n}}\ar[d]&&
{\rH^1(\Pi_{U'\rightarrow U},\xi(\oR'^{\oy'}_{U'}/p^n\oR'^{\oy'}_{U'}))}\ar[d]\\
{\tOmega^1_{\oX'_n/\oX_n,\ox'}}\ar[rr]^-(0.5){\phi_n}&&{\rH^1(\uPi,\xi(\oR'^{\oy'}_{\uX'}/p^n\oR'^{\oy'}_{\uX'}))}}
\end{equation}
est commutatif. D'après (\cite{agt} VI.11.8), les morphismes canoniques 
\begin{eqnarray}
\uDelta&\rightarrow& \underset{\underset{(U'\rightarrow U)\in \cC'_{\ox'}}{\longleftarrow}}{\lim}\ \Delta_U\\
\uDelta'&\rightarrow& \underset{\underset{(U'\rightarrow U)\in \cC'_{\ox'}}{\longleftarrow}}{\lim}\ \Delta'_{U'}.
\end{eqnarray}
sont des isomorphismes. Par suite, le morphisme canonique 
\begin{equation}
\uPi\rightarrow \underset{\underset{(U'\rightarrow U)\in \cC'_{\ox'}}{\longleftarrow}}{\lim}\ \Pi_{U'\rightarrow U}
\end{equation}
est un isomorphisme. On en déduit, compte tenu de \eqref{ktfr10a} et (\cite{serre1} I prop.~8), 
que pour tout entier $i\geq 0$, le morphisme canonique
\begin{equation}\label{crtf17f}
\underset{\underset{(U'\rightarrow U)\in \cC'^\circ_{\ox'}}{\longrightarrow}}{\lim}\
\rH^i(\Pi_{U'\rightarrow U},\oR'^{\oy'}_{U'}/p^n\oR'^{\oy'}_{U'}) \rightarrow \rH^i(\uPi,\oR'^{\oy'}_{\uX'}/p^n\oR'^{\oy'}_{\uX'})
\end{equation} 
est un isomorphisme. 

Par ailleurs, $\uoX'$ étant strictement local d'après (\cite{agt} III.3.7), il s'identifie au localisé strict de $\oX'$ en $\ox'$.
On a donc un isomorphisme canonique
\begin{equation}\label{crtf17g}
\tOmega^1_{\oX'_n/\oX_n,\ox'}
\stackrel{\sim}{\rightarrow} \underset{\underset{(U'\rightarrow U)\in \cC'^\circ_{\ox'}}{\longrightarrow}}{\lim}\
\tOmega^1_{\oX'_n/\oX_n}(\oU'^\star). 
\end{equation} 
Les homomorphismes $A_{U'\rightarrow U,n}$, pour $(\ox'\rightarrow U'\rightarrow U)\in \ob(\cC'_{\ox'})$, 
forment un homomorphisme 
de systèmes inductifs d'après (\cite{agt} III.10.16). Leur limite inductive s'identifie donc à $\phi_n$ \eqref{crtf17e}. 
Il suffit alors de montrer que pour tout objet $(\ox'\rightarrow U'\rightarrow U)$ de $\cC'_{\ox'}$, 
il existe un et un unique homomorphisme de $\oR^{\intern\oy'}_{U'\rightarrow U}$-algèbres graduées 
\begin{equation}\label{crtf17h}
\wedge \left(\xi^{-1}\tOmega^1_{\oX'_n/\oX_n}(\oU'^\star)\otimes_{\co_{\oX'}(\oU'^\star)}\oR^{\intern\oy'}_{U'\rightarrow U}\right)
\rightarrow \oplus_{i\geq 0}\rH^i(\Pi_{U'\rightarrow U},\oR'^{\oy'}_{U'}/p^n\oR'^{\oy'}_{U'})
\end{equation}
dont la composante en degré un est induite par le morphisme $A_{U'\rightarrow U,n}$, que 
son noyau est annulé par $p^{\frac{2r}{p-1}}\fm_\oK$ et que son conoyau est annulé par $p^{\frac{2r+1}{p-1}}\fm_\oK$.

D'après \ref{crtf12}, on a un diagramme commutatif
\begin{equation}\label{crtf17i}
\xymatrix{
{\xi^{-1}\tOmega^1_{\oX'_n/\oX_n}(\oU'^\star)\otimes_{\co_{\oX'}(\oU'^\star)}\oR^{\intern\oy'}_{U'\rightarrow U}}
\ar[d]_a\ar[r]^-(0.5)A&
{\rH^1(\Pi_{U'\rightarrow U},\oR'^{\oy'}_{U'}/p^n\oR'^{\oy'}_{U'})}\\
{\rH^1(\Pi_{U'\rightarrow U},\xi^{-1}(\oR'^{\oy'}_{U'}/p^n\oR'^{\oy'}_{U'})(1))}\ar[r]^{-b}_\sim&
{\rH^1(\Pi_{U'\rightarrow U},p^{\frac{1}{p-1}}\oR'^{\oy'}_{U'}/p^{n+\frac{1}{p-1}}\oR'^{\oy'}_{U'})}\ar[u]_c}
\end{equation}
où  $A$ est induit par le morphisme $A_{U'\rightarrow U,n}$, 
$a$ est induit par la composante de degré un du morphisme \eqref{eccr51a},
$b$ est induit par l'isomorphisme \eqref{definf17c} et 
$c$ est induit par l'injection canonique $p^{\frac{1}{p-1}}\co_C \rightarrow \co_C$.
En vertu de \ref{eccr51}(i), il existe un et un unique homomorphisme de $\oR^{\intern\oy'}_{U'\rightarrow U}$-algèbres graduées
\begin{equation}\label{crtf17j}
\wedge (\xi^{-1}\tOmega^1_{\oX'_n/\oX_n}(\oU'^\star)\otimes_{\co_{\oX'}(\oU'^\star)}\oR^{\intern\oy'}_{U'\rightarrow U})
\rightarrow \oplus_{i\geq 0}\rH^i(\Pi_{U'\rightarrow U},\xi^{-i}(\oR'^{\oy'}_{U'}/p^n\oR'^{\oy'}_{U'})(i))
\end{equation}
dont la composante en degré un est le morphisme $a$. Son noyau est $\alpha$-nul et son conoyau est annulé par 
$p^{\frac{1}{p-1}}\fm_\oK$. On en déduit qu'il existe un et un unique homomorphisme de 
$\oR^{\intern\oy'}_{U'\rightarrow U}$-algèbres graduées
\begin{equation}\label{crtf17k}
\wedge \left(\xi^{-1}\tOmega^1_{\oX'_n/\oX_n}(\oU'^\star)\otimes_{\co_{\oX'}(\oU'^\star)}\oR^{\intern\oy'}_{U'\rightarrow U}\right)
\rightarrow \oplus_{i\geq 0}\rH^i(\Pi_{U'\rightarrow U},\oR'^{\oy'}_{U'}/p^n\oR'^{\oy'}_{U'})
\end{equation}
dont la composante en degré $1$ est induite par $A_{U'\rightarrow U,n}$. 
Une chasse au diagramme \eqref{crtf17i} 
montre que le noyau de \eqref{crtf17k} est annulé par $p^{\frac{2r}{p-1}}\fm_\oK$.
Comme $\rH^i(\Pi_{U'\rightarrow U},\oR'^{\oy'}_{U'}/p^n\oR'^{\oy'}_{U'})$ 
est $\alpha$-nul pour tout $i\geq r+1$  en vertu de \ref{eccr51}(ii), 
le conoyau de \eqref{crtf17k} est annulé par $p^{\frac{2r+1}{p-1}}\fm_\oK$.

\subsection{}\label{crtf18}
Nous pouvons maintenant démontrer le théorème \ref{crtf4}. 
Soient $(\oy\rightsquigarrow \ox')$ un point de $X'_\et\gtimes_{X_\et}\oX^\circ_\et$ \eqref{tfr21} 
tel que $\ox'$ soit au-dessus de $s$, $\uX'$ le localisé strict de $X'$ en $\ox'$, $\ox=g(\ox')$,  $\uX$ le localisé strict de $X$ en $\ox$.
Le morphisme $\uoX'^\rhd\rightarrow \uoX^\circ$ est fidèlement plat d'après \ref{cgl1}. 
Il existe donc un point $(\oy'\rightsquigarrow \ox')$ de $X'_\et\gtimes_{X'_\et}\oX'^\rhd_\et$ 
qui relève le point $(\oy\rightsquigarrow \ox')$ \eqref{crtf5b}. Reprenons les notations de \ref{crtf5}. 
En vertu de \ref{crtf8} et \eqref{TFA12g}, pour tout entier $i\geq 0$, on a un isomorphisme canonique 
\begin{equation}\label{crtf18a}
\rR^i\tau_{n*}(\ocB_n)_{\varrho(\oy\rightsquigarrow \ox')}
\stackrel{\sim}{\rightarrow}\rH^i(\uPi,\oR'^{\oy'}_{\uX'}/p^n\oR'^{\oy'}_{\uX'}). 
\end{equation}
Par ailleurs, d'après \eqref{tfr21a}, \eqref{ktfr21b} et \eqref{ktfr9a}, on a des isomorphismes canoniques
\begin{eqnarray}
\pi_s^*(\oa'^*(\cL_{\oX'}^\gp))_{\varrho(\oy\rightsquigarrow \ox')}&\stackrel{\sim}{\rightarrow}& \cL^\gp_{\oX',\ox'},\label{crtf18b}\\
\pi_n^*(\xi^{-1}\tOmega^1_{\oX'_n/\oX_n})_{\varrho(\oy\rightsquigarrow \ox')}&\stackrel{\sim}{\rightarrow}&
\xi^{-1}\tOmega^1_{\oX'_n/\oX_n,\ox'}\otimes_{\co_{\oX',\ox'}}\oR^{\intern\oy'}_{\uX'\rightarrow \uX}.\label{crtf18c}
\end{eqnarray}
La proposition résulte alors de \ref{ktfr22}, \ref{crtf10}, \ref{crtf15}  et \ref{crtf17}.

\section{La suite spectrale de Hodge-Tate relative}\label{sshtr}

\subsection{}\label{sshtr1}
Rappelons que les hypothèses de \ref{ktfr1} sont en vigueur dans cette section. 
Considérons le morphisme composé $\bvTheta=\bvlgg\circ \bvtau$ \eqref{ktfr32a} et la suite spectrale de Cartan-Leray (\cite{sga4} V 5.4)
\begin{equation}\label{sshtr1b}
\rE_2^{i,j}=\rR^i\bvlgg_*(\rR^j \bvtau_*(\bvocB'))\Rightarrow \rR^{i+j}\bvTheta_*(\bvocB'),
\end{equation}
qui est aussi la deuxième suite spectrale d'hypercohomologie du foncteur $\bvlgg_*$ par rapport au complexe 
$\rR\bvtau_*(\bvocB')$ (\cite{ega3} 0.11.4.3). 
On rappelle que la deuxième suite spectrale d'hypercohomologie définit un foncteur 
de la catégorie dérivée des complexes de groupes abéliens de $\tG^{\mN^\circ}_s$ 
dans celle des suites spectrales de groupes abéliens de $\tE^{\mN^\circ}_s$ (\cite{ega3} 0.11.1.2).

On note $\xi^{-1}\tOmega^1_{\bvoX'/\bvoX}$ le $\co_{\bvoX'}$-module 
$(\xi^{-1}\tOmega^1_{\oX'_{n+1}/\oX_{n+1}})_{n\in \mN}$ de $X'^{\mN^\circ}_{s,\et}$ \eqref{crtf1h} 
et pour tout entier $j\geq 1$, on pose $\xi^{-j}\tOmega^j_{\bvoX'/\bvoX}=\wedge^j(\xi^{-1}\tOmega^1_{\bvoX'/\bvoX})$. 
On a alors un isomorphisme canonique (\cite{agt} (III.7.12.4))
\begin{equation}\label{sshtr1c}
\xi^{-j}\tOmega^j_{\bvoX'/\bvoX}\stackrel{\sim}{\rightarrow}(\xi^{-j}\tOmega^j_{\oX'_{n+1}/\oX_{n+1}})_{n\in \mN}.
\end{equation}
En vertu de \ref{crtf4}, pour tout entier $j\geq 0$, on a un morphisme $\bvocB^!$-linéaire canonique
\begin{equation}\label{sshtr1d}
\bvpi^*(\xi^{-j}\tOmega^j_{\bvoX'/\bvoX})\rightarrow \rR^j \bvtau_*(\bvocB'),
\end{equation}
dont le noyau et le conoyau sont annulés par  $p^{\frac{2r+1}{p-1}}\fm_\oK$ où $r=\dim(X'/X)$ (\cite{agt} (III.7.5.5) et VI.8.2).

\begin{prop}\label{sshtr2}
Supposons le morphisme $g\colon X'\rightarrow X$ propre. 
Il existe alors un entier $N\geq 0$ tel que 
pour tous entiers $i,j\geq 0$, le noyau et le conoyau du morphisme de changement de base \eqref{ktfr38b}
\begin{equation}\label{sshtr2a}
\bvsigma^*(\rR^ig_*(\tOmega^j_{X'/X})\otimes_{\co_X}\co_{\bvoX})\rightarrow \rR^i\bvlgg_*(\bvpi^*(\tOmega^j_{\bvoX'/\bvoX}))
\end{equation}
soient annulés par $p^N$. 
\end{prop}

C'est un cas particulier de \ref{ktfr39}.

\subsection{}\label{sshtr3}
On désigne par 
\begin{equation}\label{sshtr3a}
\bvdelta\colon \tE^{\mN^\circ}_s \rightarrow \tE^{\mN^\circ}
\end{equation}
le foncteur induit par le plongement $\delta$ \eqref{mtfla7f}. On a alors $\bvdelta_*(\bvocB)=(\ocB_{n+1})_{n\in \mN}$, 
où l'on considère les $\ocB_{n+1}$ comme des anneaux de $\tE$. On peut donc indifféremment 
considérer $\bvocB$ comme un anneau de $\tE^{\mN^\circ}_s$ ou de $\tE^{\mN^\circ}$.
Le foncteur $\bvdelta_*$ est exact sur la catégorie des groupes abéliens (\cite{agt} III.7.3(i)). 
Pour tout $\mZ_p$-module $\cF=(\cF_{n+1})_{n\in \mN}$ de $\tE^{\mN^\circ}$, on a un isomorphisme $\bvocB$-linéaire canonique
\begin{equation}\label{sshtr3b}
\cF\otimes_{\mZ_p}\bvocB\stackrel{\sim}{\rightarrow} \bvdelta_*(\bvdelta^*(\cF)\otimes_{\mZ_p}\bvocB). 
\end{equation}
On peut donc indifféremment  considérer $\cF\otimes_{\mZ_p}\bvocB$ comme un $\bvocB$-module de $\tE^{\mN^\circ}$ ou de $\tE^{\mN^\circ}_s$.

On désigne par $\bMod(\bvocB)$ la catégorie des $\bvocB$-modules de $\tE_s^{\mN^\circ}$,
par $\bMod_{\mQ}(\bvocB)$ la catégorie des objets de $\bMod(\bvocB)$ à isogénie près et par
\begin{equation}\label{sshtr3c}
\bMod(\bvocB)\rightarrow \bMod_\mQ(\bvocB), \ \ \ \cF\mapsto \cF_\mQ,
\end{equation}
le foncteur de localisation (\cite{agt} 6.1.1). On notera que la catégorie $\bMod_{\mQ}(\bvocB)$ est abélienne.

\subsection{}\label{sshtr4}
On désigne par $\upgamma\colon \oX'^\rhd\rightarrow \oX^\circ$ le morphisme induit par $g$, par 
\begin{equation}\label{sshtr4a}
(\oX'^\rhd_\et)^{\mN^\circ}\stackrel{\bvupgamma}{\longrightarrow} (\oX^\circ_\et)^{\mN^\circ} \stackrel{\bvpsi}{\longrightarrow} \tE^{\mN^\circ}
\end{equation}
les morphismes induits par $\upgamma$ et $\psi$ \eqref{mtfla7d}, et par $\bvmZ_p$ la $\mZ_p$-algèbre 
$(\mZ/p^{n+1})_{n\in \mN}$ de $(\oX'^\rhd_\et)^{\mN^\circ}$, à ne pas confondre avec le système projectif constant $\mZ_p$. 
 
En vertu de \ref{TCFR19}, si le morphisme $g\colon X'\rightarrow X$ est projectif \eqref{ktfr1a},
pour tout entier $q\geq 0$, on a un morphisme $\bvocB$-linéaire canonique de $\tE^{\mN^\circ}_s$,
\begin{equation}\label{sshtr4b}
\bvpsi_*(\rR^q\bvupgamma_*(\bvmZ_p))\otimes_{\mZ_p}\bvocB\rightarrow \rR^q\bvTheta_*(\bvocB'),
\end{equation}
qui est un $\alpha$-isomorphisme. 

\begin{teo}\label{sshtr5}
Si le morphisme $g\colon X'\rightarrow X$ est projectif, il existe une suite spectrale canonique de $\bvocB_\mQ$-modules
\begin{equation}\label{sshtr5a}
\rE_2^{i,j}=\bvsigma^*(\rR^ig_*(\tOmega^j_{X'/X})\otimes_{\co_X}\co_{\bvoX})_\mQ(-j)\Rightarrow \bvpsi_*(\rR^{i+j}\bvupgamma_*(\bvmZ_p))\otimes_{\mZ_p}\bvocB_\mQ.
\end{equation}
\end{teo}

La suite spectrale en question est induite par la suite spectrale de Cartan-Leray \eqref{sshtr1b} compte tenu de \eqref{sshtr1d},
\ref{sshtr2} et \eqref{sshtr4b}.

\begin{defi}\label{sshtr50}
La suite spectrale \eqref{sshtr5a} est appelée {\em suite spectrale de Hodge-Tate relative}.
\end{defi}

Le théorème \ref{sshtr5}, ainsi que les énoncés \ref{sshtrl6} et \ref{sshtrsg6} ci-dessous, valent en fait sous l'hypothèse plus générale que {\em $g$ soit propre}. 
En effet, l'hypothèse de projectivité sur $g$ est utilisée dans la preuve de \ref{TCFR18} qui s'étend également aux morphismes propres (voir \ref{TCFR180}).

\subsection{}\label{sshtr6}
On rappelle \eqref{ssht120} que le groupe de Galois $G_K$ agit naturellement à gauche sur le topos $\tE^{\mN^\circ}_s$
et que $\bvocB$ est un anneau $G_K$-équivariant de $\tE^{\mN^\circ}_s$.
De même, $G_K$ agit naturellement à gauche sur les topos annelés $\tE'^{\mN^\circ}_s$ et $\tG^{\mN^\circ}_s$,
et $\bvocB'$ (resp. $\bvocB^!$) est un anneau $G_K$-équivariant de $\tE'^{\mN^\circ}_s$ (resp. $\tG^{\mN^\circ}_s$). 
Les morphismes $\bvtau$, $\bvlgg$, $\bvTheta$, $\bvsigma'$ et $\bvpi$ sont $G_K$-équivariants \eqref{ktfr32a}. 

Soit $\cF$ un $\bvocB'$-module $G_K$-équivariant. D'après \ref{notconv18}, pour tout $j\geq 0$,
$\rR^j\bvtau_*(\cF)$ est naturellement muni d'une structure de $\bvocB^!$-module $G_K$-équivariant. 
Donc pour tout $i\geq 0$, $\rR^i\bvlgg_*(\rR^j\bvtau_*(\cF))$ est naturellement muni d'une structure de $\bvocB$-module $G_K$-équivariant. 
De même, pour tout $q\geq 0$, $\rR^q\bvTheta_*(\cF)$ est naturellement muni d'une structure de $\bvocB$-module $G_K$-équivariant.

Par ailleurs, pour tout $u\in G_K$, les actions de $u$ sur $\oX'$ et $\oX$ induisent un isomorphisme $\co_{\oX'}$-linéaire
\begin{equation}
\tOmega^1_{\oX'/\oX}\stackrel{\sim}{\rightarrow} u^*(\tOmega^1_{\oX'/\oX}).
\end{equation}
Ces isomorphismes munissent $\tOmega^1_{\oX'/\oX}$ d'une structure de $\co_{\oX'}$-module $G_K$-équivariant \eqref{notconv17}. 
On rappelle que $G_K$ agit naturellement sur $\xi\co_C$ \eqref{ssht120}.  
On en déduit pour tout entier $j\geq 0$, 
une structure de $\co_{\bvoX'}$-module $G_K$-équivariant sur $\xi^{-j}\tOmega^j_{\bvoX'/\bvoX}$, 
ou ce qui revient au même puisque l'action de $G_K$ sur le topos $X'^{\mN^\circ}_{s,\et}$ est triviale,
une action à gauche $\co_{\bvoX'}$-semi-linéaire de $G_K$ sur le $\co_{\bvoX'}$-module $\xi^{-j}\tOmega^j_{\bvoX'/\bvoX}$. 
Par suite, $\bvpi^*(\xi^{-j}\tOmega^j_{\bvoX'/\bvoX})$ est naturellement muni d'une structure de $\bvocB^!$-module $G_K$-équivariant 
(cf. \ref{notconv18}).

\begin{prop}\label{sshtr7}
La suite spectrale de Cartan-Leray \eqref{sshtr1b} est $G_K$-équivariante \eqref{sshtr6}. 
\end{prop}

Pour tout $u\in G_K$, considérons l'isomorphisme composé 
\begin{equation}\label{sshtr7a}
\rR\bvtau_*(\bvocB')\stackrel{\sim}{\rightarrow}  \rR\bvtau_*(u_*(\bvocB'))\stackrel{\sim}{\rightarrow}
\rR(\bvtau u)_*(\bvocB') \stackrel{\sim}{\rightarrow}
\rR(u\bvtau)_*(\bvocB') \stackrel{\sim}{\rightarrow} \rR u_* (\rR\bvtau_*(\bvocB')),
\end{equation}
où la première flèche est induite par l'isomorphisme canonique $\bvocB'\stackrel{\sim}{\rightarrow} u_*(\bvocB')$,
la deuxième et la quatrième flèches sont les isomorphismes canoniques, 
et la troisième flèche est l'isomorphisme qui traduit la $G_K$-équivariance de $\bvtau$. 
On en déduit un morphisme de la suite spectrale $\rE$ \eqref{sshtr1b} vers la deuxième suite spectrale d'hypercohomologie du foncteur 
$\bvlgg_*$ relativement au complexe $\rR u_* (\rR\bvtau_*(\bvocB'))$. Par ailleurs, l'isomorphisme 
\begin{equation}
\rR\bvlgg_*(\rR u_* (\rR\bvtau_*(\bvocB')))  \stackrel{\sim}{\rightarrow} \rR u_* (\rR\bvlgg_*(\rR\bvtau_*(\bvocB')))
\end{equation}
déduit de la $G_K$-équivariance de $\bvlgg$, induit un isomorphisme de la deuxième suite spectrale d'hypercohomologie du foncteur 
$\bvlgg_*$ relativement au complexe $\rR u_* (\rR\bvtau_*(\bvocB'))$ vers la suite spectrale $u_*(\rE)$. 
On vérifie aussitôt que le morphisme composé de suites spectrales $\rE\rightarrow u_*(\rE)$ est donné sur les termes $\rE_2$
et sur les aboutissements par les structures $G_K$-équivariantes définies dans \ref{sshtr6}

\begin{prop}\label{sshtr8}
Pour tout entier $j\geq 0$, le morphisme \eqref{sshtr1d}
\begin{equation}\label{sshtr8a}
\bvpi^*(\xi^{-j}\tOmega^j_{\bvoX'/\bvoX})\rightarrow \rR^j \bvtau_*(\bvocB')
\end{equation}
est $G_K$-équivariant pour les structures $G_K$-équivariantes naturelles \eqref{sshtr6}. 
\end{prop}

Il suffit de calquer la preuve de \ref{ssht14}. 
Le groupe $G_K$ agit naturellement sur le schéma logarithmique $(\oS,\cL_\oS)$ \eqref{hght3} et par suite sur 
le schéma logarithmique  $(\oX',\cL_{\oX'})$ \eqref{crtf1a}. Pour tout entier $n\geq 0$, le monoïde $\cQ'_n$ \eqref{crtf2a} est donc
muni d'une structure de monoïde $G_K$-équivariant de $\tE'$, de sorte que les morphismes du diagramme \eqref{crtf2a} sont 
des morphismes de monoïdes $G_K$-équivariants. 
Munissant le groupe $\mu_{p^n}(\co_\oK)$ de l'action canonique de $G_K$, 
la  suite exacte \eqref{crtf2d}  est donc une suite de groupes abéliens $G_K$-équivariants de $\tE'_s$. 
Par suite, $\partial'_n$ \eqref{crtf3a} est un morphisme de groupes abéliens $G_K$-équivariants de $\tG_s$, 
et il en est alors de même du morphisme \eqref{crtf3c}.
Le caractère universel de la dérivation logarithmique $d\log$ \eqref{crtf1f}
implique que \eqref{crtf1g} 
\begin{equation}
\cL_{\oX'}^\gp\otimes_\mZ\co_{\oX'}\rightarrow \tOmega^1_{\oX'/\oX}
\end{equation}
est un morphisme de $\co_{\oX'}$-modules $G_K$-équivariants de $\oX'_\et$. 
Le morphisme \eqref{crtf4b} est donc $G_K$-équivariant, d'où la proposition.

\begin{cor}\label{sshtr9}
Si le morphisme $g\colon X'\rightarrow X$ est projectif, la suite spectrale de Hodge-Tate relative \eqref{sshtr5a} est $G_K$-équivariante.
\end{cor}

On rappelle d'abord que $G_K$ agit naturellement à gauche sur les topos $(\oX^\circ_\et)^{\mN^\circ}$, $(\oX'^\rhd_\et)^{\mN^\circ}$
et $\tE^{\mN^\circ}$, et que les morphismes $\bvupgamma$ et $\bvpsi$ \eqref{sshtr4a} sont $G_K$-équivariants \eqref{notconv18}. 
Par suite, pour tout entier $q\geq 0$, 
$\bvpsi_*(\rR^q\bvupgamma_*(\bvmZ_p))\otimes_{\mZ_p}\bvocB$ est naturellement muni d'une structure de $\bvocB$-module 
$G_K$-équivariant. Le morphisme \eqref{sshtr4b}
\begin{equation}
\bvpsi_*(\rR^q\bvupgamma_*(\bvmZ_p))\otimes_{\mZ_p}\bvocB\rightarrow \rR^q\bvTheta_*(\bvocB'),
\end{equation}
est un morphisme de $\bvocB$-modules $G_K$-équivariants. En effet, avec les notations des preuves de \ref{TCFR18} et \ref{TCFR19},
les morphismes \eqref{TCFR18e} et \eqref{TCFR19c} sont $G_K$-équivariants. 
On vérifie aussi que le morphisme de changement de base \eqref{sshtr2a} est un morphisme de $\bvocB$-modules $G_K$-équivariants. 
La proposition résulte alors de \ref{sshtr7}, \ref{sshtr8} et \ref{TCFR19}.

\begin{lem}\label{sshtr33}
Soient $A$ une $\co_K$-algèbre noethérienne, $M$ un $A$-module de type fini, $j$ un entier non-nul. 
Alors, on a 
\begin{equation}\label{sshtr33a}
((M\hotimes_{\co_K}\co_C)\otimes_{\mZ_p}\mQ_p(j))^{G_K}=0,
\end{equation}
où le produit tensoriel $\hotimes$ est séparé complété pour la topologie $p$-adique.
\end{lem}

Observant que $(M\hotimes_{\co_K}\co_C)\otimes_{\mZ_p}\mQ_p(j)$ est la limite du système inductif 
$(M\hotimes_{\co_K}\co_C(j))_{\mN}$ où les morphismes de transition sont induits par la multiplication par $p$, on voit que 
\begin{equation}\label{sshtr33b}
((M\hotimes_{\co_K}\co_C)\otimes_{\mZ_p}\mQ_p(j))^{G_K}=(M\hotimes_{\co_K}\co_C(j))^{G_K}\otimes_{\mZ_p}\mQ_p.
\end{equation}

Soient $M_\tor$ le sous-$A$-module des sections $x$ de $M$ pour lesquelles il existe un entier $i\geq 1$ tel que $\fm_K^i x=0$,
$M'$ le quotient de $M$ par $M_\tor$, $u\colon M\rightarrow M'$ le morphisme canonique. 
Il existe un entier $m\geq 0$ tel que $M_\tor$ soit annulé par $p^m$.
D'après \ref{alpha3}, il existe alors un morphisme $A$-linéaire $v\colon M'\rightarrow M$
tel que $u\circ v= p^{2m}\id_{M'}$ et $v\circ u=p^{2m}\id_{M}$. Compte tenu de \eqref{sshtr33b} appliqué à $M$ et à $M'$, 
on peut se borner au cas $M$ est $\co_K$-plat et il suffit alors de montrer que 
\begin{equation}\label{sshtr33c}
(M\hotimes_{\co_K}\co_C(j))^{G_K}=0.
\end{equation}

On a des isomorphismes canoniques
\begin{eqnarray}
(M\hotimes_{\co_K}\co_C(j))^{G_K}&\stackrel{\sim}{\rightarrow}&
\underset{\underset{n\in \mN}{\longleftarrow}}{\lim}\ (M\otimes_{\co_K}(\co_C/p^n\co_C)(j))^{G_K},\label{sshtr33d}\\
(\co_C(j))^{G_K}&\stackrel{\sim}{\rightarrow}&
\underset{\underset{n\in \mN}{\longleftarrow}}{\lim}\ ((\co_C/p^n\co_C)(j))^{G_K}.\label{sshtr33e}
\end{eqnarray}
En vertu de (\cite{agt} II.3.15), on a un isomorphisme canonique
\begin{equation}\label{sshtr33f}
(M\otimes_{\co_K}(\co_C/p^n\co_C)(j))^{G_K}
\stackrel{\sim}{\rightarrow}
M\otimes_{\co_K}((\co_C/p^n\co_C)(j))^{G_K}.
\end{equation}
Soient $x\in (M\hotimes_{\co_K}\co_C(j))^{G_K}$, $x_n$ son image dans $(M\otimes_{\co_K}(\co_C/p^n\co_C)(j))^{G_K}$ pour tout $n\geq 0$.
Choisissons une $k$-base $(\oee_i)_{i\in I}$ de $M/\fm_KM$, et pour tout $i\in I$, un relèvement $e_i\in M$ de $\oee_i$.
Comme $M$ est $\co_K$-plat, pour tout $n\geq 0$, les classes $(e_{i,n})_{i\in I}$ des $(e_i)_{i\in I}$ modulo $p^n$ forment 
une $(\co_K/p^n\co_K)$-base de $M/p^nM$. 
L'image de $x_n$ par l'isomorphisme \eqref{sshtr33f} s'écrit sous la forme $\sum_{i\in I}e_{i,n} \otimes x_{i,n}$, 
où les $x_{i,n}$ sont des sections de $((\co_C/p^n\co_C)(j))^{G_K}$ presque toutes nulles. 
Pour tout $i\in I$, le système projectif $(x_{i,n})_{n\geq 0}$ définit un élément
de $(\co_C(j))^{G_K}$ compte tenu de \eqref{sshtr33e}. 
Comme $(\co_C(j))^{G_K}\subset (C(j))^{G_K}=0$ d'après le théorème de Tate (\cite{tate} theo.~2), 
on en déduit que les $x_{i,n}$ sont nuls pour tous $i$ et $n$. Par suite, $x$ est nul d'après \eqref{sshtr33d}, d'où l'assertion recherchée \eqref{sshtr33c}.

\begin{prop}\label{sshtr10}
Pour tout $\co_X$-module cohérent $\cF$ tel que $\cF|X_\eta$ soit un $\co_{X_\eta}$-module localement libre 
et pour tout entier non nul $j$, on a 
\begin{equation}\label{sshtr10a}
(\Gamma(\tE^{\mN^\circ}_s,\bvsigma^*(\cF\otimes_{\co_X}\co_\bvoX))\otimes_{\mZ_p}\mQ_p(j))^{G_K}=0.
\end{equation}
\end{prop}

Reprenons les notations de \ref{ktfr34}. Posons $\ocF=\cF\otimes_{\co_X}\co_\oX$ et notons $\hocF$ son complété $p$-adique qui est un $\co_\fX$-module. 
D'après \ref{ktfr40}, le morphisme canonique
\begin{equation}\label{sshtr10h}
\Gamma(\fX,\hocF)\otimes_{\mZ_p}\mQ_p\rightarrow 
\Gamma(\tE^{\mN^\circ}_s,\bvsigma^*(\cF\otimes_{\co_X}\co_\bvoX))\otimes_{\mZ_p}\mQ_p
\end{equation}
est un isomorphisme. Celui-ci est en fait $G_K$-équivariant pour les structures $G_K$-équivariantes canoniques \eqref{notconv19}. 
Cela résulte de la preuve de \ref{ktfr40} en observant que 
l'homomorphisme canonique $\co_{\fX}\rightarrow \top_*(\bvocB)$ \eqref{ktfr40b} est $G_K$-équivariant (cf. \ref{notconv18}). 
Il suffit donc de montrer que pour tout entier non nul $j$, on a 
\begin{equation}\label{sshtr10i}
(\Gamma(\fX,\hocF)\otimes_{\mZ_p}\mQ_p(j))^{G_K}=0.
\end{equation}

On peut se borner au cas où $X=\Spec(R)$ est affine et $\cF=\tM$ où $M$ est un $R$-module de type fini. 
L'assertion recherchée résulte alors de \ref{sshtr33}.

\begin{cor}\label{sshtr11}
Si le morphisme $g\colon X'\rightarrow X$ est projectif, la suite spectrale de Hodge-Tate relative \eqref{sshtr5a} dégénère en $\rE_2$.
\end{cor}

En effet, pour tout $i,j\geq 0$, posons
\begin{equation}\label{sshtr11a}
\cF^{i,j}=\rR^ig_*(\tOmega^j_{X'/X}).
\end{equation}
D'après \ref{sshtr9}, la différentielle $d_2^{i,j}\colon \rE_2^{i,j}\rightarrow \rE_2^{i+2,j-1}$ s'identifie à une section $G_K$-équivariante de 
\begin{equation}\label{sshtr11b}
\Gamma(\tE^{\mN^\circ}_s,\cHom_{\bvocB}(\bvsigma^*(\cF^{i,j}\otimes_{\co_X}\co_\bvoX),\bvsigma^*(\cF^{i+2,j-1}\otimes_{\co_X}\co_\bvoX)))
\otimes_{\mZ_p}\mQ_p(1).
\end{equation}

Montrons qu'il existe un entier $n\geq 0$ tel que le noyau et le conoyau du morphisme canonique
\begin{eqnarray}\label{sshtr11c}
\lefteqn{\bvsigma^*(\cHom_{\co_X}(\cF^{i,j},\cF^{i+2,j-1})\otimes_{\co_X}\co_{\bvoX})\rightarrow}\nonumber \\
&& \cHom_{\bvocB}(\bvsigma^*(\cF^{i,j}\otimes_{\co_X}\co_\bvoX),\bvsigma^*(\cF^{i+2,j-1}\otimes_{\co_X}\co_\bvoX))
\end{eqnarray}
soient annulés par $p^n$. La question étant locale sur $X$, on peut se borner au cas où $X$ est affine. 
Le morphisme $g$ étant saturé, est exact (\cite{ogus} III 2.5.2), et les fibres du monoïde $(\cM_X/\co^\times_X)|X_\eta$
sont libres puisque $\cM_X|X_\eta$ est défini par un diviseur à croisements normaux (\cite{agt} III.4.7). 
Par suite, le $\co_{X_\eta}$-module $\cF^{i,j}|X_\eta$ est localement libre d'après (\cite{ikn} 7.2)
(cf. \cite{deligne1} 5.5 pour le cas lisse sans structures logarithmiques). Il existe donc un $\co_X$-module cohérent $\cG^{i,j}$,
un entier $n\geq 1$ et un morphisme $\co_X$-linéaire $u\colon \cF^{i,j}\oplus\cG^{i,j}\rightarrow \co_X^n$
induisant un isomorphisme sur $X_\eta$. Il existe un entier $m\geq 0$ tel que $p^m$ annule le noyau et le conoyau de $u$. 
D'après \ref{alpha3}, il existe alors un morphisme $\co_X$-linéaire $v\colon \co_X^n\rightarrow \cF^{i,j}\oplus\cG^{i,j}$
tel que $u\circ v= p^{2m}\id_{\co_X^n}$ et $v\circ u=p^{2m}\id_{\cF^{i,j}\oplus\cG^{i,j}}$, autrement dit, $u$ est une isogénie (\cite{agt} III.6.1.1).
On en déduit aussitôt que le morphisme \eqref{sshtr11c} est une isogénie, d'où l'assertion recherchée.  
Compte tenu de \ref{alpha3}, on peut alors s'identifier $d_2^{i,j}$ à une section $G_K$-équivariante de 
\begin{equation}\label{sshtr11d}
\Gamma(\tE^{\mN^\circ}_s,\bvsigma^*(\cHom_{\co_X}(\cF^{i,j},\cF^{i+2,j-1})\otimes_{\co_X}\co_{\bvoX}))
\otimes_{\mZ_p}\mQ_p(1).
\end{equation}
La proposition s'ensuit en vertu de \ref{sshtr10}.

\section{La suite spectrale de Hodge-Tate relative: localisation}\label{sshtrl}

\subsection{}\label{sshtrl1}
Rappelons que les hypothèses de \ref{ktfr1} sont en vigueur dans cette section. 
On se donne, de plus, un point  $(\oy\rightsquigarrow \ox)$ de $X_\et\gtimes_{X_\et}\oX^\circ_\et$ \eqref{topfl17} tel que $\ox$ soit au-dessus de $s$.
On désigne par $\uX$ le localisé strict de $X$ en $\ox$. On pose $\uX'=X'\times_X\uX$,
\begin{eqnarray}
\uR&=&\Gamma(\uX,\co_\uX),\label{sshtrl1ee}\\
\uoR&=&\ocB_{\rho(\oy\rightsquigarrow \ox)},\label{sshtrl1e}
\end{eqnarray}
et on note $\uhoR$ le séparé complété $p$-adique de $\uoR$. 

Pour tout faisceau abélien $F=(F_n)_{n\geq 0}$ de $\tE'^{\mN^\circ}$, on désigne par $\Theta_*(F)_{\rho(\oy\rightsquigarrow \ox)}$ 
la limite projective des fibres des $\Theta_*(F_n)$ en $\rho(\oy\rightsquigarrow \ox)$,
\begin{equation}\label{sshtrl1a}
\Theta_*(F)_{\rho(\oy\rightsquigarrow \ox)}=\underset{\underset{n\geq 0}{\longleftarrow}}{\lim}\ \Theta_*(F_n)_{\rho(\oy\rightsquigarrow \ox)},
\end{equation}
où $\Theta\colon \tE'\rightarrow \tE$ est le morphisme de fonctorialité \eqref{mtfla11b}.
On définit ainsi un foncteur de la catégorie des faisceaux abéliens de $\tE'^{\mN^\circ}$ dans celle des groupes abéliens.
Celui-ci étant exact à gauche, 
on désigne abusivement par $\rR^q\Theta_*(F)_{\rho(\oy\rightsquigarrow \ox)}$ ($q\geq 0$) ses foncteurs dérivés à droite. 
D'après (\cite{jannsen} 1.6), on a une suite exacte canonique 
\begin{equation}\label{sshtrl1b}
0\rightarrow \rR^1 \underset{\underset{n\geq 0}{\longleftarrow}}{\lim}\ \rR^{q-1}\Theta_*(F_n)_{\rho(\oy\rightsquigarrow \ox)}\rightarrow
\rR^q\Theta_*(F)_{\rho(\oy\rightsquigarrow \ox)}\rightarrow 
\underset{\underset{n\geq 0}{\longleftarrow}}{\lim}\ \rR^q\Theta_*(F_n)_{\rho(\oy\rightsquigarrow \ox)}\rightarrow 0,
\end{equation}
où on a posé $\rR^{-1}\Theta_*(F_n)_{\rho(\oy\rightsquigarrow \ox)}=0$ pour tout $n\geq 0$. 

De même, pour tout faisceau abélien $H=(H_n)_{n\geq 0}$ de $\tG^{\mN^\circ}$, on désigne par $\lgg_*(H)_{\rho(\oy\rightsquigarrow \ox)}$ 
la limite projective des fibres des $\lgg_*(H_n)$ en $\rho(\oy\rightsquigarrow \ox)$,
\begin{equation}\label{sshtrl1c}
\lgg_*(H)_{\rho(\oy\rightsquigarrow \ox)}=\underset{\underset{n\geq 0}{\longleftarrow}}{\lim}\ \lgg_*(H_n)_{\rho(\oy\rightsquigarrow \ox)}, 
\end{equation}
où $\lgg\colon \tG\rightarrow \tE$ est le morphisme défini dans \eqref{ktfr1c}.
On définit ainsi un foncteur de la catégorie des faisceaux abéliens de $\tG^{\mN^\circ}$ dans celle des groupes abéliens.
Celui-ci étant exact à gauche, on désigne abusivement par $\rR^q\lgg_*(H)_{\rho(\oy\rightsquigarrow \ox)}$ ($q\geq 0$) ses foncteurs dérivés à droite. 
On a une suite exacte canonique 
\begin{equation}\label{sshtrl1d}
0\rightarrow \rR^1 \underset{\underset{n\geq 0}{\longleftarrow}}{\lim}\ \rR^{q-1}\lgg_*(H_n)_{\rho(\oy\rightsquigarrow \ox)}\rightarrow
\rR^q\lgg_*(H)_{\rho(\oy\rightsquigarrow \ox)}\rightarrow 
\underset{\underset{n\geq 0}{\longleftarrow}}{\lim}\ \rR^q\lgg_*(H_n)_{\rho(\oy\rightsquigarrow \ox)}\rightarrow 0,
\end{equation}
où on a posé $\rR^{-1}\lgg_*(H_n)_{\rho(\oy\rightsquigarrow \ox)}=0$ pour tout $n\geq 0$. 

Il résulte de \ref{alpha3} et \ref{crtf4} que pour tous entiers $i,j\geq 0$, le morphisme \eqref{sshtr1d} induit un isomorphisme 
\begin{equation}\label{sshtrl1j}
\rR^i\lgg_*(\bvpi^*(\xi^{-j}\tOmega^j_{\bvoX'/\bvoX}))_{\rho(\oy\rightsquigarrow \ox)}\otimes_{\mZ_p}\mQ_p
\stackrel{\sim}{\rightarrow} \rR^i\lgg_*(\rR^j \bvtau_*(\bvocB'))_{\rho(\oy\rightsquigarrow \ox)}\otimes_{\mZ_p}\mQ_p.
\end{equation}

Pour tout faisceau abélien $F$ de $\tE'^{\mN^\circ}$, on a un isomorphisme canonique \eqref{ktfr1c} 
\begin{equation}\label{sshtrl1f}
\Theta_*(F)_{\rho(\oy\rightsquigarrow \ox)}\stackrel{\sim}{\rightarrow}\lgg_*(\bvtau_*(F))_{\rho(\oy\rightsquigarrow \ox)}.
\end{equation}
D'après (\cite{sga4} V 4.9(3)), le foncteur $\bvtau_*$ transforme les faisceaux abéliens injectifs en faisceaux abéliens injectifs. 
Par suite, en vertu de (\cite{sga4} V 0.3), on a une suite spectrale 
\begin{equation}\label{sshtrl1g}
\rE_2^{i,j}=\rR^i\lgg_*(\rR^j\bvtau_*(F))_{\rho(\oy\rightsquigarrow \ox)}\Rightarrow \rR^{i+j}\Theta_*(F)_{\rho(\oy\rightsquigarrow \ox)}.
\end{equation}
Par ailleurs, l'isomorphisme \eqref{sshtrl1f} induit, pour tout $q\geq 0$, un isomorphisme canonique 
\begin{equation}\label{sshtrl1h}
\rR^q\Theta_*(F)_{\rho(\oy\rightsquigarrow \ox)}\stackrel{\sim}{\rightarrow} \rR^q\lgg_*(\rR\bvtau_*(F))_{\rho(\oy\rightsquigarrow \ox)}.
\end{equation}
Celui-ci identifie la suite spectrale \eqref{sshtrl1g} avec la deuxième suite spectrale d'hypercohomologie du foncteur 
$\lgg_*(-)_{\rho(\oy\rightsquigarrow \ox)}$ par rapport au complexe $\rR\bvtau_*(F)$ (\cite{ega3} 0.11.4.3 ou \cite{deligne2} 1.4.5 et 1.4.6),
\begin{equation}\label{sshtrl1i}
\rE_2^{i,j}=\rR^i\lgg_*(\rR^j\bvtau_*(F))_{\rho(\oy\rightsquigarrow \ox)}\Rightarrow \rR^{i+j}\lgg_*(\rR\bvtau_*(F))_{\rho(\oy\rightsquigarrow \ox)}. 
\end{equation}

\begin{prop}\label{sshtrl2}
Supposons le morphisme $g\colon X'\rightarrow X$ projectif. Alors, pour tout entier $q\geq 0$, le morphisme canonique 
\begin{equation}\label{sshtrl2a}
\rR^q\Theta_*(\bvocB')_{\rho(\oy\rightsquigarrow \ox)}\rightarrow \underset{\underset{n\geq 0}{\longleftarrow}}\lim\ 
\rR^q\Theta_*(\ocB'_n)_{\rho(\oy\rightsquigarrow \ox)}
\end{equation}
et surjectif et son noyau est $\alpha$-nul. 
\end{prop}

En effet, le morphisme \eqref{sshtrl2a} est surjectif d'après \eqref{sshtrl1b}.  
Par ailleurs, en vertu de \ref{TCFR19}, pour tous entiers $n,q\geq 0$, on a un morphisme $\ocB$-linéaire canonique 
\begin{equation}\label{sshtrl2b}
\psi_*(\rR^q\upgamma_*(\mZ/p^n\mZ))\otimes_{\mZ_p}\ocB
\rightarrow \rR^q\Theta_*(\ocB'_n),
\end{equation}
qui est un $\alpha$-isomorphisme. D'après \ref{Kpun34} et \ref{TCFR160},
le faisceau abélien $\rR^q\upgamma_*(\mZ/p^n\mZ)$ est localement constant constructible et on a un isomorphisme canonique
\begin{equation}\label{sshtrl2c}
\rR^q\upgamma_*(\mZ/p^n\mZ)_\oy\stackrel{\sim}{\rightarrow}\rH^q(\oX'^\rhd_{\oy,\et},\mZ/p^n\mZ).
\end{equation}
En vertu de (\cite{agt} VI.9.18 et VI.9.20), il existe donc un $(\mZ/p^n\mZ)$-module canonique 
$\mL^q_n$ de $\oX^\circ_\fet$ et un isomorphisme canonique
$\rR^q\upgamma_*(\mZ/p^n\mZ)\stackrel{\sim}{\rightarrow} \rho_{\oX^\circ}^*(\mL^q_n)$,
où $\rho_{\oX^\circ}\colon \oX^\circ_\et\rightarrow \oX^\circ_\fet$ est le morphisme canonique \eqref{notconv10a}.
D'après (\cite{agt} VI.10.9(iii)), on a un isomorphisme canonique 
\begin{equation}\label{sshtrl2d}
\beta^*(\mL^q_n)\stackrel{\sim}{\rightarrow} \psi_*(\rR^q\upgamma_*(\mZ/p^n\mZ)).
\end{equation}
On en déduit, compte tenu de \eqref{tfr21b} et (\cite{agt} VI.9.9), un morphisme $\uoR$-linéaire canonique 
\begin{equation}\label{sshtrl2e}
u_n^q\colon \rH^q(\oX'^\rhd_{\oy,\et},\mZ/p^n\mZ)\otimes_{\mZ_p}\uoR
\rightarrow \rR^q\Theta_*(\ocB'_n)_{\rho(\oy\rightsquigarrow \ox)},
\end{equation}
qui est un $\alpha$-isomorphisme \eqref{sshtrl1e}. Les morphismes 
\begin{equation}\label{sshtrl2f}
\underset{\underset{n\geq 0}{\longleftarrow}}\lim\ u^q_n \ \ \ {\rm et}\ \ \ \rR^1\underset{\underset{n\geq 1}{\longleftarrow}}\lim\ u^q_n
\end{equation}
sont donc des $\alpha$-isomorphismes (\cite{gr} 2.4.2(ii)). Les groupes $\rH^q(\oX'^\rhd_{\oy,\et},\mZ/p^n\mZ)$ étant finis, 
le système projectif $(\rH^q(\oX'^\rhd_{\oy,\et},\mZ/p^n\mZ))_{n\geq 1}$ vérifie la condition de Mittag-Leffler. 
On en déduit que 
\begin{equation}\label{sshtrl2g}
\rR^1\underset{\underset{n\geq 0}{\longleftarrow}}\lim\ \rR^{q-1}\Theta_*(\ocB'_n)_{\rho(\oy\rightsquigarrow \ox)}
\end{equation}
est $\alpha$-nul d'après (\cite{jannsen} 1.15) et (\cite{roos} théo.~1). La proposition s'ensuit compte tenu de \eqref{sshtrl1b}. 

\begin{rema}\label{sshtrl200}
Calquant la preuve de \ref{ssht300}, on montre que le $\mZ_p$-module $\rH^q(\oX'^\rhd_{\oy,\et},\mZ_p)$ est de type fini, et que le morphisme canonique 
\begin{equation}\label{sshtrl200a}
\rH^q(\oX'^\rhd_{\oy,\et},\mZ_p)\otimes_{\mZ_p}\uhoR \rightarrow \underset{\underset{n\geq 0}{\longleftarrow}}\lim\  
\rH^q(\oX'^\rhd_{\oy,\et},\mZ/p^n\mZ) \otimes_{\mZ_p}\uoR
\end{equation}
est un isomorphisme. 
\end{rema}

\begin{cor}\label{sshtrl3}
Supposons le morphisme $g\colon X'\rightarrow X$ projectif. Alors, pour tout entier $q\geq 0$, on a un $\alpha$-isomorphisme canonique 
\begin{equation}\label{sshtrl3a}
\rR^q\Theta_*(\bvocB')_{\rho(\oy\rightsquigarrow \ox)}\stackrel{\approx}{\longrightarrow} 
\rH^q(\oX'^\rhd_{\oy,\et},\mZ_p)\otimes_{\mZ_p}\uhoR. 
\end{equation}
\end{cor}

Cela résulte de \ref{sshtrl2}, \eqref{sshtrl2e} et \ref{sshtrl200}. 
On prendra garde que le $\alpha$-isomorphisme \eqref{sshtrl3a} n'est pas en général induit par un vrai morphisme.

\begin{prop}\label{sshtrl4}
Supposons le morphisme $g\colon X'\rightarrow X$ propre. 
Il existe alors un entier $N\geq 0$ vérifiant la propriété suivante. Soit $\cF'$ un $\co_{X'}$-module cohérent de $X'_\zar$.
Pour tout $n\geq 1$, posons $\ocF'_n=\cF'\otimes_{\co_S}\co_{\oS_n}$, que l'on considère comme un $\co_{\oX'_n}$-module de $X'_{s,\et}$ \eqref{notconv12}. 
Notons $\bvocF'=(\ocF'_{n+1})_{n\in \mN}$ et $\bvpi^*(\bvocF')$ son image inverse par le morphisme de topos annelés $\bvpi$ \eqref{ktfr32a}.
Alors, pour tout entier $q\geq 0$, le morphisme canonique 
\begin{equation}\label{sshtrl4a}
\rR^q\lgg_*(\bvpi^*(\bvocF'))_{\rho(\oy\rightsquigarrow \ox)}\rightarrow 
\underset{\underset{n\geq 0}{\longleftarrow}}\lim\ \rR^q\lgg_{n*}(\pi^*_n(\ocF'_n))_{\rho(\oy\rightsquigarrow \ox)}
\end{equation}
est surjectif et son noyau est annulé par $p^N$.  
\end{prop}

En effet, le morphisme \eqref{sshtrl4a} est surjectif d'après \eqref{sshtrl1d}.
Par ailleurs, en vertu de \ref{lptfr16}(ii), \ref{ktfr27}, \eqref{lptfr10k} et (\cite{sga4} XII 5.5), on a un morphisme canonique 
\begin{equation}\label{sshtrl4e}
\rH^q(\uoX', \ocF'_n\otimes_{\co_\oX}\uoR)\rightarrow 
\rR^q\lgg_{n*}(\pi^*_n(\ocF'_n))_{\rho(\oy\rightsquigarrow \ox)},
\end{equation}
qui est un $\alpha$-isomorphisme. 

D'après \ref{ktfr300} il existe un entier $N\geq 0$ tel que pour tous entiers $n\geq 1$ et $q\geq 0$, 
posant $\cF'_n=\cF'/p^n\cF'$, le noyau et le conoyau du morphisme canonique
\begin{equation}\label{sshtrl4c}
\rH^q(\uX',\cF'_n\otimes_{\co_{X}}\co_{\uX})\otimes_{\co_{\uX}(\uX)}\uoR
\rightarrow \rH^q(\uoX',\ocF'_n\otimes_{\co_{\oX}}\uoR)
\end{equation}
soient annulés par $p^N$. Notons $I^q_n$ (resp. $C^q_n$) l'image (resp. le conoyau) de ce morphisme. 

En vertu de (\cite{ega3} 4.1.7), 
le système projectif $(\rH^q(\uX',\cF'_n\otimes_{\co_{X}}\co_{\uX}))_{n\geq 1}$ vérifie la condition de Mittag-Leffler. 
Par suite, le système projectif $(I^{q}_n)_{n\geq 1}$ vérifie la condition de Mittag-Leffler et le morphisme canonique 
\begin{equation}
\rR^1\underset{\underset{n\geq 0}{\longleftarrow}}\lim\ \rH^q(\uoX',\cF'_n\otimes_{\co_{X}}\uoR)\rightarrow \rR^1\underset{\underset{n\geq 0}{\longleftarrow}}\lim\ C^{q}_n
\end{equation}
est un isomorphisme. La proposition s'ensuit compte tenu de \eqref{sshtrl1d}.

\begin{cor}\label{sshtrl5}
Supposons le morphisme $g\colon X'\rightarrow X$ propre. Soit $\cF'$ un $\co_{X'}$-module cohérent.
Pour tout $n\geq 1$, posons $\ocF'_n=\cF'\otimes_{\co_S}\co_{\oS_n}$, que l'on considère comme un $\co_{\oX'_n}$-module de $X'_{s,\et}$ \eqref{notconv12}. 
Notons $\bvocF'=(\ocF'_{n+1})_{n\in \mN}$ et $\bvpi^*(\bvocF')$ son image inverse par le morphisme de topos annelés $\bvpi$ \eqref{ktfr32a}.
Alors, pour tout entier $q\geq 0$, on a un isomorphisme canonique
\begin{equation}\label{sshtrl5a}
\rR^q\lgg_*(\bvpi^*(\bvocF'))_{\rho(\oy\rightsquigarrow \ox)}\otimes_{\mZ_p}\mQ_p\stackrel{\sim}{\rightarrow} 
(\rH^q(\uX',\cF'\otimes_{\co_{X}}\co_{\uX})\hotimes_{\uR}\uoR)\otimes_{\mZ_p}\mQ_p,
\end{equation}
où le produit tensoriel $\hotimes$ est séparé complété pour la topologie $p$-adique.
\end{cor}

Reprenons les notations de la preuve de \ref{sshtrl4}. 
Pour tous entiers $n\geq 1$ et $q\geq 0$, on note $L^q_n$ (resp. $Q_n^q$) le noyau (resp. conoyau) du morphisme canonique 
\begin{equation}
\rH^q(\uX',\cF'\otimes_{\co_{X}}\co_{\uX})\rightarrow  \rH^q(\uX',\cF'_n\otimes_{\co_{X}}\co_{\uX}).
\end{equation}
En vertu de (\cite{ega3} 4.1.7 et 4.1.7.20), 
la filtration $(L^q_n)_{n\geq 1}$ sur $\rH^q(\uX',\cF'\otimes_{\co_{X}}\co_{\uX})$ est $p$-bonne et le système projectif $(Q_n^q)_{n\geq 1}$ est AR-nul.
Par suite, la filtration définie sur $\rH^q(\uX',\cF'\otimes_{\co_{X}}\co_{\uX})\otimes_{\uR}\uoR$ par les images 
des $L_n^q\otimes_{\uR}\uoR$ ($n\geq 1$) est $p$-bonne; et le morphisme canonique 
\begin{equation}
\rH^q(\uX',\cF'\otimes_{\co_{X}}\co_{\uX})\hotimes_{\uR}\uoR\rightarrow 
\underset{\underset{n\geq 0}{\longleftarrow}}\lim\ \rH^q(\uX',\cF'_n\otimes_{\co_{X}}\co_{\uX})\otimes_{\uR}\uoR,
\end{equation}
où le produit tensoriel $\hotimes$ est séparé complété pour la topologie $p$-adique, est un isomorphisme. 
On notera que $L_n^q\supset p^n \rH^q(\uX',\cF'\otimes_{\co_{X}}\co_{\uX})$. 

Comme les noyaux et conoyaux des morphismes \eqref{sshtrl4c} sont annulés par $p^N$, on en déduit que le noyau du morphisme canonique
\begin{equation}
\rH^q(\uX',\cF'\otimes_{\co_{X}}\co_{\uX})\hotimes_{\uR}\uoR\rightarrow 
\underset{\underset{n\geq 0}{\longleftarrow}}\lim\ \rH^q(\uoX',\ocF'_n\otimes_{\co_{\oX}}\uoR)
\end{equation}
est annulé par $p^N$ et que son conoyau est annulé par $p^{2N}$. La proposition s'ensuit compte tenu de \ref{sshtrl4} et \eqref{sshtrl4e}.

\begin{teo}\label{sshtrl6}
Supposons le morphisme $g\colon X'\rightarrow X$ projectif. 
Alors, il existe une suite spectrale canonique de $\uhoR$-modules 
\begin{equation}\label{sshtrl6a}
\rE_2^{i,j}=(\rH^i(\uX',\tOmega^j_{\uX'/\uX})\hotimes_{\uR}\uoR)\otimes_{\mZ_p}\mQ_p(-j)
\Rightarrow \rH^{i+j}(\oX'^\rhd_{\oy,\et},\mZ_p)\otimes_{\mZ_p}\uhoR[\frac 1 p],
\end{equation}
où le produit tensoriel $\hotimes$ est séparé complété pour la topologie $p$-adique.
\end{teo}

Cela résulte de la suite spectrale \eqref{sshtrl1g} en prenant pour faisceau $F=\bvocB'$,
compte tenu de \ref{sshtrl3}, \eqref{sshtrl1j} et \ref{sshtrl5}. 

On peut donner deux autres constructions de la suite spectrale \eqref{sshtrl6a}. On peut la voir comme une suite spectrale de Cartan-Leray 
relativement à un morphisme de topos (cf. \ref{sshtrl8}). On peut aussi la déduire directement de la suite spectrale de Hodge-Tate relative globale \eqref{sshtr5a} (voir \ref{sshtrl20}). 

\subsection{}\label{sshtrl7} 
On note $\uY$ le localisé strict de $\oX^\circ$ en $\oy$, $X_{(\oy)}$ le localisé strict de $X$ en $h(\oy)$ \eqref{mtfla2}, 
$u\colon \uX\rightarrow X$ et $v\colon \uY\rightarrow \oX^\circ$ les morphismes canoniques 
et $\uh\colon \uY\rightarrow \uX$ le composé du morphisme $\uY\rightarrow X_{(\oy)}$ induit par $h$ 
et du morphisme $X_{(\oy)}\rightarrow \uX$ induit par la flèche de spécialisation de $h(\oy)$ dans $\ox$ (\cite{sga4} VIII 7.3).  
On pose $\uY'=\uY\times_{\oX^\circ}\oX'^\rhd$ et on considère le diagramme commutatif
\begin{equation}\label{sshtrl7a}
\xymatrix{
X'\ar[ddd]_g\ar@{}[rddd]|\Box&&&{\oX'^\rhd}\ar[lll]_-(0.5){h'}\ar[ddd]^\upgamma\ar@{}[lddd]|\Box\\
&{\uX'}\ar[d]_{\ug}\ar[lu]_{u'}&{\uY'}\ar[ru]^{v'}\ar[d]^{\uupgamma}\ar[l]_-(0.5){\uh'}&\\
&{\uX}\ar[ld]^u&{\uY}\ar[l]_-(0.5){\uh}\ar[rd]_v&\\
X&&&{\oX^\circ}\ar[lll]_-(0.5){h}}
\end{equation}
où $u'$, $\ug$, $v'$ et $\uupgamma$ sont les projections canoniques et $\uh'$ est induit par les morphismes $\uh\circ \uupgamma$ et $h'\circ v'$. 
On désigne par $\tuE$ (resp. $\tuE'$) le topos de Faltings du morphisme $\uh$ (resp. $\uh'$) \eqref{tf1} et par $\tuG$
le topos de Faltings relatif associé au couple de morphismes $(\uh,\ug)$ \eqref{tfr1}. 
On prendra garde de ne pas confondre les topos $\tuG$ et $\tuE'$ avec ceux définis dans \ref{ktfr4} et \ref{ktfr15}.  
Le diagramme \eqref{sshtrl7a} induit par fonctorialité du topos de Faltings relatif \eqref{tfr7} un diagramme de morphismes de topos 
à carrés commutatifs à isomorphismes canoniques près
\begin{equation}\label{sshtrl7b}
\xymatrix{
\tuE'\ar[r]^{\Phi'}\ar[d]_{\utau}&\tE'\ar[d]^\tau\\
\tuG\ar[r]^\Psi\ar[d]_{\ulgg}&\tG\ar[d]^\lgg\\
\tuE\ar[r]^\Phi&\tE}
\end{equation}
Celui-ci induit un diagramme de morphismes de topos à carrés commutatifs à isomorphismes canoniques près
\begin{equation}\label{sshtrl7c}
\xymatrix{
{\tuE'^{\mN^\circ}}\ar[r]^-(0.5){\bvPhi'}\ar[d]_{\bvutau}&{\tE'^{\mN^\circ}}\ar[d]^{\bvtau}\\
{\tuG^{\mN^\circ}}\ar[r]^-(0.5){\bvPsi}\ar[d]_{\bvulgg}&{\tG^{\mN^\circ}}\ar[d]^{\bvlgg}\\
{\tuE^{\mN^\circ}}\ar[r]^-(0.5){\bvPhi}&{\tE^{\mN^\circ}}}
\end{equation}
Ces notations n'induisent aucune confusion avec celles introduites dans \eqref{ktfr32a}. En effet, on peut indifféremment 
considérer $\bvocB$ comme un anneau de $\tE^{\mN^\circ}_s$ ou de $\tE^{\mN^\circ}$ \eqref{sshtr3}.
De même, on peut indifféremment considérer $\bvocB'$ (resp. $\bvocB^!$) comme un anneau de $\tE'^{\mN^\circ}_s$ ou de $\tE'^{\mN^\circ}$
(resp. $\tG^{\mN^\circ}_s$ ou de $\tG^{\mN^\circ}$): notant $\bvdelta'\colon \tE'^{\mN^\circ}_s \rightarrow \tE'^{\mN^\circ}$ et 
$\bvkappa\colon \tG^{\mN^\circ}_s \rightarrow \tG^{\mN^\circ}$ les foncteurs induits par les plongements $\delta'$ \eqref{mtfla9e} et $\kappa$ \eqref{ktfr20c}, respectivement,  
on a $\bvdelta'_*(\bvocB')=(\ocB'_{n+1})_{n\in \mN}$ et $\bvkappa_*(\bvocB^!)=(\ocB^!_{n+1})_{n\in \mN}$, 
où l'on considère les $\ocB'_{n+1}$ (resp. $\ocB^!_{n+1}$) comme des anneaux de $\tE'$ (resp. $\tG$).

En vertu de \ref{lptfr9}(i), pour tout  faisceau abélien $F$ de $\tE'^{\mN^\circ}$, on a un isomorphisme canonique fonctoriel 
\begin{equation}\label{sshtrl7d}
\Theta_*(F)_{\rho(\oy\rightsquigarrow \ox)}\stackrel{\sim}{\rightarrow} \Gamma(\tuE'^{\mN^\circ},\bvPhi'^*(F)).
\end{equation}
De même, pour tout  faisceau abélien $H$ de $\tG^{\mN^\circ}$, on a un isomorphisme canonique fonctoriel 
\begin{equation}\label{sshtrl7f}
\lgg_*(H)_{\rho(\oy\rightsquigarrow \ox)}\stackrel{\sim}{\rightarrow} \Gamma(\tuG^{\mN^\circ},\bvPsi^*(H)).
\end{equation}

Par ailleurs, il résulte de \ref{lptfr9} et (\cite{sga4} VIII 7.3) que pour tout faisceau abélien $F$ de $\tE'^{\mN^\circ}$ et tout entier $q\geq 0$, 
le morphisme de changement de base relativement au carré supérieur de \eqref{sshtrl7c} 
\begin{equation}\label{sshtrl7h}
\bvPsi^*(\rR^q\bvtau_*(F))\rightarrow \rR^q \bvutau_*(\bvPhi'^*(F))
\end{equation}
est un isomorphisme.

\begin{prop}\label{sshtrl8}
Conservons les notations de \ref{sshtrl7}.
\begin{itemize}
\item[{\rm (i)}] Pour tout  faisceau abélien $F$ de $\tE'^{\mN^\circ}$ et tout entier $q\geq 0$, 
l'isomorphisme \eqref{sshtrl7d} induit un isomorphisme canonique fonctoriel 
\begin{equation}\label{sshtrl8a}
\rR^q\Theta_*(F)_{\rho(\oy\rightsquigarrow \ox)}\stackrel{\sim}{\rightarrow} \rH^q(\tuE'^{\mN^\circ},\bvPhi'^*(F)).
\end{equation}
\item[{\rm (ii)}] Pour tout  faisceau abélien $H$ de $\tG^{\mN^\circ}$ et 
tout entier $q\geq 0$, l'isomorphisme \eqref{sshtrl7f} induit un isomorphisme canonique fonctoriel 
\begin{equation}\label{sshtrl8b}
\rR^q\lgg_*(H)_{\rho(\oy\rightsquigarrow \ox)}\stackrel{\sim}{\rightarrow} \rH^q(\tuG^{\mN^\circ},\bvPsi^*(H)).
\end{equation}
\item[{\rm (iii)}]  Pour tout faisceau abélien $F$ de $\tE'^{\mN^\circ}$, l'isomorphisme composé
\begin{equation}\label{sshtrl8c}
\rR^i\lgg_*(\rR^j\bvtau_*(F))_{\rho(\oy\rightsquigarrow \ox)}\stackrel{\sim}{\rightarrow} \rH^i(\tuG^{\mN^\circ},\bvPsi^*(\rR^j\bvtau_*(F)))
\stackrel{\sim}{\rightarrow} \rH^i(\tuG^{\mN^\circ},\rR^j\bvtau_*(\bvPhi'^*(F))),
\end{equation}
induit par \eqref{sshtrl8b} et \eqref{sshtrl7h}, et l'isomorphisme \eqref{sshtrl8a} induisent un isomorphisme de la suite spectrale \eqref{sshtrl1g} vers 
la suite spectrale de Cartan-Leray {\rm (\cite{sga4} V 5.3)}
\begin{equation}\label{sshtrl8d}
\rE_2^{i,j}=\rH^i(\tuG^{\mN^\circ},\rR^j \bvutau_*(\bvPhi'^{-1}(F)))\Rightarrow \rH^{i+j}(\tuE'^{\mN^\circ},\bvPhi'^{-1}(F)).
\end{equation}
\end{itemize}
\end{prop}

(i) Il revient au même de dire que pour tout  faisceau abélien injectif $F=(F_n)_{n\geq 0}$ de $\tE'^{\mN^\circ}$ et tout entier $q\geq 1$, 
on a $\rH^q(\tuE'^{\mN^\circ},\bvPhi'^*(F))=0$. 
En effet, les faisceaux abéliens $F_n$ de $\tE'$ sont injectifs et les morphismes $F_{n+1}\rightarrow F_n$ sont des surjections scindées d'après (\cite{jannsen} 1.1). 
On a donc $\rR^j\Theta_*(F_n)_{\rho(\oy\rightsquigarrow \ox)}=0$ pour tous $j\geq 1$ et $n\geq 0$, et par suite $\rH^j(\tuE',\Phi'^*(F_n))=0$ en vertu de \eqref{lptfr9d}. 
L'assertion recherchée résulte alors de la suite exacte canonique (\cite{agt} VI.7.10)
\begin{equation}\label{sshtrl8e}
0\rightarrow \rR^1\underset{\underset{n\geq 0}{\longleftarrow}}\lim\ \rH^{q-1}(\tuE',\Phi'^*(F_n))\rightarrow 
\rH^q(\tuE'^{\mN^\circ},\bvPhi'^*(F))\rightarrow \underset{\underset{n\geq 0}{\longleftarrow}}\lim\ \rH^q(\tuE',\Phi'^*(F_n))\rightarrow 0.
\end{equation}

(ii) La preuve est identique à celle de (i). 

(iii) En effet, d'après (ii), le foncteur $\bvPsi^*$ transforme les faisceaux abéliens injectifs de $\tG^{\mN^\circ}$ 
en faisceaux abéliens de $\tuG^{\mN^\circ}$ acycliques pour le foncteur $\Gamma(\tuG^{\mN^\circ},-)$. 
Compte tenu de \eqref{sshtrl7h}, l'isomorphisme \eqref{sshtrl8b} 
induit donc un isomorphisme de la suite spectrale \eqref{sshtrl1i} vers la suite spectrale d'hypercohomologie du foncteur 
$\Gamma(\tuG^{\mN^\circ}, -)$ par rapport au complexe $\rR\bvutau_*(\bvPhi'^*(F))$,
\begin{equation}\label{sshtrl8h}
\rE_2^{i,j}=\rH^i(\tuG^{\mN^\circ}, \rR^j\bvutau_*(\bvPhi'^*(F)))\Rightarrow 
\rH^{i+j}(\tuG^{\mN^\circ}, \rR\bvutau_*(\bvPhi'^*(F))).
\end{equation}
Celle-ci est canoniquement isomorphe à la suite spectrale de Cartan-Leray \eqref{sshtrl8d}.

\subsection{}\label{sshtrl9}
On rappelle que $\uoX$ est normal et strictement local (et en particulier intègre) (\cite{agt} III.3.7). 
Le $X$-morphisme $\oy\rightarrow \uX$ définissant le point $(\oy\rightsquigarrow \ox)$ se relève en un $\oX^\circ$-morphisme $\oy\rightarrow \uoX^\circ$ et 
il induit donc un point géométrique de $\uoX^\circ$ que l'on note aussi (abusivement) $\oy$. 
On pose $\uDelta=\pi_1(\uoX^\circ,\oy)$ et on note $\mZ[\uDelta]$ son algèbre de groupe. 
On désigne par $\bB_{\uDelta}$ le topos classifiant du groupe profini $\uDelta$ \eqref{notconv11}, par 
\begin{equation}\label{sshtrl9aa}
\nu_\oy\colon \uoX^\circ_\fet\rightarrow \bB_{\uDelta},
\end{equation}
le foncteur fibre en $\oy$ \eqref{notconv11c}, par 
\begin{equation}\label{sshtrl9a}
\varphi_\ox\colon \tE\rightarrow \uoX^\circ_\fet
\end{equation}
le foncteur canonique défini dans \eqref{tf10b}, et par 
\begin{equation}\label{sshtrl9b}
\tE\rightarrow \bB_{\uDelta}, \ \ \ F\mapsto F_{[\oy\rightsquigarrow \ox]}=\nu_\oy(\varphi_\ox(F))
\end{equation}
le foncteur composé $\nu_\oy\circ\varphi_\ox$. Ce sont deux foncteurs exacts. 
D'après (\cite{agt} VI.10.31 et VI.9.9), pour tout objet $F$ de $\tE$, oubliant l'action de $\uDelta$, 
on a un isomorphisme canonique fonctoriel d'ensembles 
\begin{equation}\label{sshtrl9c}
F_{\rho(\oy\rightsquigarrow \ox)}\stackrel{\sim}{\rightarrow}F_{[\oy\rightsquigarrow \ox]}.
\end{equation}

L'anneau $\uoR=\ocB_{[\oy\rightsquigarrow \ox]}$ \eqref{sshtrl1e}
est donc naturellement muni d'une action de $\uDelta$ par des automorphismes d'anneaux (cf. \ref{TFA12}). Il en est donc de même de $\uhoR$. 

Pour tout faisceau abélien $F=(F_n)_{n\geq 0}$ de $\tE'^{\mN^\circ}$, on désigne par $\Theta_*(F)_{[\oy\rightsquigarrow \ox]}$ 
la limite projective des $\mZ[\uDelta]$-modules $\Theta_*(F_n)_{[\oy\rightsquigarrow \ox]}$ \eqref{sshtrl9b},
\begin{equation}\label{sshtrl9d}
\Theta_*(F)_{[\oy\rightsquigarrow \ox]}=\underset{\underset{n\geq 0}{\longleftarrow}}{\lim}\ \Theta_*(F_n)_{[\oy\rightsquigarrow \ox]}.
\end{equation}
On définit ainsi un foncteur de la catégorie des faisceaux abéliens de $\tE'^{\mN^\circ}$ dans celle des $\mZ[\uDelta]$-modules.
Celui-ci étant exact à gauche, 
on désigne abusivement par $\rR^q\Theta_*(F)_{[\oy\rightsquigarrow \ox]}$ ($q\geq 0$) ses foncteurs dérivés à droite. 
Oubliant l'action de $\uDelta$, l'isomorphisme \eqref{sshtrl9d} induit un isomorphisme canonique de groupes abéliens 
\begin{equation}\label{sshtrl9e}
\rR^q\Theta_*(F)_{\rho(\oy\rightsquigarrow \ox)}\stackrel{\sim}{\rightarrow}\rR^q\Theta_*(F)_{[\oy\rightsquigarrow \ox]},
\end{equation}
où la source est définie dans \ref{sshtrl1}.
D'après (\cite{jannsen} 1.6), on a une suite exacte canonique 
\begin{equation}\label{sshtrl9f}
0\rightarrow \rR^1 \underset{\underset{n\geq 0}{\longleftarrow}}{\lim}\ \rR^{q-1}\Theta_*(F_n)_{[\oy\rightsquigarrow \ox]}\rightarrow
\rR^q\Theta_*(F)_{[\oy\rightsquigarrow \ox]}\rightarrow 
\underset{\underset{n\geq 0}{\longleftarrow}}{\lim}\ \rR^q\Theta_*(F_n)_{[\oy\rightsquigarrow \ox]}\rightarrow 0,
\end{equation}
où on a posé $\rR^{-1}\Theta_*(F_n)_{[\oy\rightsquigarrow \ox]}=0$ pour tout $n\geq 0$. 

De même, pour tout faisceau abélien $H=(H_n)_{n\geq 0}$ de $\tG^{\mN^\circ}$, on désigne par $\lgg_*(H)_{[\oy\rightsquigarrow \ox]}$ 
la limite projective des $\mZ[\uDelta]$-modules $\lgg_*(H_n)_{[\oy\rightsquigarrow \ox]}$ \eqref{sshtrl9b},
\begin{equation}\label{sshtrl9g}
\lgg_*(H)_{[\oy\rightsquigarrow \ox]}=\underset{\underset{n\geq 0}{\longleftarrow}}{\lim}\ \lgg_*(H_n)_{[\oy\rightsquigarrow \ox]}.
\end{equation}
On définit ainsi un foncteur de la catégorie des faisceaux abéliens de $\tG^{\mN^\circ}$ dans celle des $\mZ[\uDelta]$-modules.
Celui-ci étant exact à gauche, 
on désigne abusivement par $\rR^q\lgg_*(H)_{[\oy\rightsquigarrow \ox]}$ ($q\geq 0$) ses foncteurs dérivés à droite. 
Oubliant l'action de $\uDelta$, l'isomorphisme \eqref{sshtrl9g} induit un isomorphisme canonique de groupes abéliens 
\begin{equation}\label{sshtrl9h}
\rR^q\lgg_*(H)_{[\oy\rightsquigarrow \ox]}\stackrel{\sim}{\rightarrow}\rR^q\lgg_*(H)_{\rho(\oy\rightsquigarrow \ox)},
\end{equation}
où la source est définie dans \ref{sshtrl1}.
D'après (\cite{jannsen} 1.6), on a une suite exacte canonique 
\begin{equation}\label{sshtrl9i}
0\rightarrow \rR^1 \underset{\underset{n\geq 0}{\longleftarrow}}{\lim}\ \rR^{q-1}\lgg_*(H_n)_{[\oy\rightsquigarrow \ox]}\rightarrow
\rR^q\lgg_*(H)_{[\oy\rightsquigarrow \ox]}\rightarrow 
\underset{\underset{n\geq 0}{\longleftarrow}}{\lim}\ \rR^q\lgg_*(H_n)_{[\oy\rightsquigarrow \ox]}\rightarrow 0,
\end{equation}
où on a posé $\rR^{-1}\lgg_*(H_n)_{[\oy\rightsquigarrow \ox]}=0$ pour tout $n\geq 0$. 

D'après (\cite{sga4} V 4.9(3)), le foncteur $\bvtau_*$ transforme les faisceaux abéliens injectifs en faisceaux abéliens injectifs. 
Par suite, en vertu de (\cite{sga4} V 0.3), pour tout faisceau abélien $F$ de $\tE'^{\mN^\circ}$, on a une suite spectrale de $\mZ[\uDelta]$-modules
\begin{equation}\label{sshtrl9n}
\rE_2^{i,j}=\rR^i\lgg_*(\rR^j\bvtau_*(F))_{[\oy\rightsquigarrow \ox]}\Rightarrow \rR^{i+j}\Theta_*(F)_{[\oy\rightsquigarrow \ox]}. 
\end{equation}
Celle-ci raffine la suite spectrale \eqref{sshtrl1g}.

\subsection{}\label{sshtrl90}
Supposons le morphisme $g\colon X'\rightarrow X$ propre. Soient $n,q$ deux entiers $\geq 0$.
En vertu de \ref{Kpun34} et \ref{TCFR160}, le faisceau $\rR^q\upgamma_*(\mZ/p^n\mZ)$ est localement constant constructible 
et on a un isomorphisme canonique
\begin{equation}\label{sshtrl9j}
\rR^q\upgamma_*(\mZ/p^n\mZ)_\oy\stackrel{\sim}{\rightarrow}\rH^q(\oX'^\rhd_{\oy,\et},\mZ/p^n\mZ).
\end{equation}
D'après  (\cite{agt} VI.9.18 et VI.9.20), 
il existe donc un $(\mZ/p^n\mZ)$-module canonique $\mL^q_n$ de $\oX^\circ_\fet$ et un isomorphisme canonique
$\rR^q\upgamma_*(\mZ/p^n\mZ)\stackrel{\sim}{\rightarrow} \rho_{\oX^\circ}^*(\mL^q_n)$,
où $\rho_{\oX^\circ}\colon \oX^\circ_\et\rightarrow \oX^\circ_\fet$ est le morphisme canonique \eqref{notconv10a}.
On désigne par $\umL^q_n$ l'image inverse de $\mL^q_n$ sur $\uoX^\circ_\fet$. 
D'après (\cite{agt} VI.9.9), oubliant l'action de $\uDelta$, on a un isomorphisme canonique \eqref{sshtrl9aa}
\begin{equation}
\nu_\oy(\umL^q_n)\stackrel{\sim}{\rightarrow}\rR^q\upgamma_*(\mZ/p^n\mZ)_\oy.
\end{equation}
On en déduit une action naturelle de $\uDelta$ sur 
$\rH^q(\oX'^\rhd_{\oy,\et},\mZ/p^n\mZ)$ et par suite une action naturelle de $\uDelta$ sur $\rH^q(\oX'^\rhd_{\oy,\et},\mZ_p)$. 

D'après (\cite{agt} VI.10.9(iii)), on a un isomorphisme canonique 
\begin{equation}\label{sshtrl9k}
\beta^*(\mL^q_n)\stackrel{\sim}{\rightarrow} \psi_*(\rR^q\upgamma_*(\mZ/p^n\mZ)).
\end{equation}
Compte tenu de (\cite{agt} (VI.10.24.3)), celui-ci induit un isomorphisme de $\uoX^\circ_\fet$,
\begin{equation}\label{sshtrl9m}
\umL^q_n\stackrel{\sim}{\rightarrow} \varphi_\ox(\psi_*(\rR^q\upgamma_*(\mZ/p^n\mZ))),
\end{equation}
et par suite, en prenant les fibres, un isomorphisme de $\uDelta$-modules
\begin{equation}\label{sshtrl9l}
\nu_\oy(\umL^q_n)\stackrel{\sim}{\rightarrow} \psi_*(\rR^q\upgamma_*(\mZ/p^n\mZ))_{[\oy\rightsquigarrow \ox]}.
\end{equation}

\subsection{}\label{sshtrl11}
Reprenons les hypothèses et notations de \ref{sshtrl4} et de sa preuve. Pour tous entiers $n\geq 1$ et $q\geq 0$, le morphisme 
\begin{equation}\label{sshtrl11a}
\rH^q(\uX',\cF'_n\otimes_{\co_{X}}\co_{\uX})\otimes_{\co_{\uX}(\uX)}\uoR  
\rightarrow \rR^q\lgg_{n*}(\pi^*_n(\ocF'_n))_{\rho(\oy\rightsquigarrow \ox)}
\end{equation}
composé de \eqref{sshtrl4e} et \eqref{sshtrl4c}, induit, compte tenu de \eqref{sshtrl9c}, un morphisme $\uDelta$-équivariant 
\begin{equation}\label{sshtrl11b}
\rH^q(\uX',\cF'_n\otimes_{\co_{X}}\co_{\uX})\otimes_{\co_{\uX}(\uX)}\uoR  
\rightarrow \rR^q\lgg_{n*}(\pi^*_n(\ocF'_n))_{[\oy\rightsquigarrow \ox]}.
\end{equation}
En effet, le morphisme \eqref{sshtrl11a} n'est autre que la fibre en $\rho(\oy\rightsquigarrow \ox)$ du morphisme de changement de base \eqref{ktfr24k}
\begin{equation}\label{sshtrl11c}
\sigma_n^*(\rR^q\ogg_{n*}(\ocF'_n))\rightarrow \rR^q\lgg_{n*}(\pi_n^*(\ocF'_n)). 
\end{equation}
Celui-ci est le composé 
\begin{eqnarray}\label{sshtrl11d}
\lefteqn{\sigma_s^{-1}(\rR^q\ogg_{s*}(\ocF'_n))\otimes_{\sigma_s^{-1}(\co_{\oX_n})}\ocB_n\rightarrow }\\
&&\rR^q\lgg_*(\pi_s^{-1}(\ocF'_n)) \otimes_{\sigma_s^{-1}(\co_{\oX_n})}\ocB_n
\rightarrow \rR^q\lgg_*(\pi_s^{-1}(\ocF'_n) \otimes_{\pi_s^{-1}(\co_{\oX'_n})}\ocB^!_n),\nonumber
\end{eqnarray}
où la première flèche est induite par le morphisme de changement de base pour les faisceaux abéliens relativement au diagramme commutatif \eqref{ktfr21h} et 
la seconde flèche est définie par fonctorialité et $\ocB_n$-linéarité.  
Pour tout faisceau $L$ de $X_\et$, $\uDelta$ agit trivialement sur $\sigma^*(L)_{[\oy\rightsquigarrow \ox]}$ en vertu de (\cite{agt} VI.10.24).
L'assertion recherchée résulte donc de \eqref{sshtrl11d}.   

\begin{prop}\label{sshtrl10}
Supposons le morphisme $g\colon X'\rightarrow X$ projectif. 
Alors, avec les notations de \ref{sshtrl9}, la suite spectrale \eqref{sshtrl6a} est $\uDelta$-équivariante. 
\end{prop}

Rappelons que la suite spectrale \eqref{sshtrl6a} a été définie à partir de \eqref{sshtrl1g} qui est raffinée par \eqref{sshtrl9n}. 
Il suffit donc de montrer que les isomorphismes envisagés dans la preuve de \ref{sshtrl6} sont $\uDelta$-équivariants.

Pour tous entiers $i,j\geq 0$, l'isomorphisme \eqref{sshtrl1j}
\begin{equation}
\rR^i\lgg_*(\bvpi^*(\xi^{-j}\tOmega^j_{\bvoX'/\bvoX}))_{[\oy\rightsquigarrow \ox]}\otimes_{\mZ_p}\mQ_p
\stackrel{\sim}{\rightarrow} \rR^i\lgg_*(\rR^j \bvtau_*(\bvocB'))_{[\oy\rightsquigarrow \ox]}\otimes_{\mZ_p}\mQ_p
\end{equation}
est clairement $\uDelta$-équivariant. 
Par ailleurs, il résulte de \ref{sshtrl11} et de la preuve de \ref{sshtrl5} que l'isomorphisme \eqref{sshtrl5a}, pour $\cF'=\tOmega^j_{X'/X}$,
est sous-jacent à un isomorphisme $\uDelta$-équivariant
\begin{equation}
\rR^i\lgg_*(\bvpi^*(\tOmega^j_{\bvoX'/\bvoX}))_{[\oy\rightsquigarrow \ox]}\otimes_{\mZ_p}\mQ_p
\stackrel{\sim}{\rightarrow}
(\rH^i(\uX',\tOmega^j_{X'/X}\otimes_{\co_{X}}\co_{\uX})\hotimes_{\co_\uX}\uoR)\otimes_{\mZ_p}\mQ_p.
\end{equation}

Pour tout entier $q\geq 0$, le morphisme \eqref{sshtrl2b} induit un morphisme $\uDelta$-équivariant
\begin{equation}
\psi_*(\rR^q\upgamma_*(\mZ/p^n\mZ))_{[\oy\rightsquigarrow \ox]}\otimes_{\mZ_p}\ocB_{[\oy\rightsquigarrow \ox]}
\rightarrow \rR^q\Theta_*(\ocB'_n)_{[\oy\rightsquigarrow \ox]}.
\end{equation}
D'après \eqref{sshtrl9l}, le morphisme $\uoR$-linéaire canonique  
\begin{equation}
\rH^q(\oX'^\rhd_{\oy,\et},\mZ/p^n\mZ)\otimes_{\mZ_p}\uoR
\rightarrow \rR^q\Theta_*(\ocB'_n)_{[\oy\rightsquigarrow \ox]}
\end{equation}
est donc $\uDelta$-équivariant. Il résulte alors de la preuve de \ref{sshtrl3} que le $\alpha$-isomorphisme \eqref{sshtrl3a} induit un isomorphisme $\uDelta$-équivariant
\begin{equation}
\rR^q\Theta_*(\bvocB')_{[\oy\rightsquigarrow \ox]}\otimes_{\mZ_p}\mQ_p\stackrel{\sim}{\rightarrow}
\rH^q(\oX'^\rhd_{\oy,\et},\mZ_p)\otimes_{\mZ_p}\uhoR[\frac 1 p],
\end{equation}
ce qui achève la preuve.

\subsection{}\label{sshtrl12}
Reprenons les notations de \ref{sshtrl9} et notons $\varpi\colon \uoX^\circ\rightarrow \uX^\circ$ la projection canonique et
$\Pi(\uX^\circ)$ (resp. $\Pi(\uoX^\circ)$) le groupoïde fondamental de $\uX^\circ$ (resp. $\uoX^\circ$) \eqref{notconv20}. 
On pose $\uGamma=\pi_1(\uX^\circ,\varpi(\oy))$, de sorte qu'on a une suite exacte de groupes profinis
\begin{equation}\label{sshtrl12a}
1\rightarrow \uDelta\rightarrow \uGamma\rightarrow G_K \rightarrow 1.
\end{equation}
On note $\mZ[\uGamma]$ l'algèbre de groupe associée à $\uGamma$. 
On fait agir $\uGamma$ sur $\mZ_p(1)$ \eqref{hght1a} à travers son quotient $G_K$. 
On désigne par $\bB_{\uGamma}$ le topos classifiant du groupe profini $\uGamma$ \eqref{notconv11} et par 
\begin{equation}\label{sshtrl12c}
\upnu_\oy\colon \uX^\circ_\fet\rightarrow \bB_{\uGamma},
\end{equation}
le foncteur fibre en $\varpi(\oy)$.

Le groupe $G_K$ agit naturellement à gauche sur $\uoX^\circ$ \eqref{ssht12} et le morphisme $\varpi\colon \uoX^\circ\rightarrow \uX^\circ$ 
est $G_K$-équivariant lorsque l'on fait agir $G_K$ trivialement sur $\uX^\circ$.  
L'action de $G_K$ sur $\uoX^\circ$ induit une action à gauche de $G_K$ sur le topos $\uoX^\circ_\fet$ \eqref{notconv17},  
une action à gauche de $G_K$ sur le groupoïde $\Pi(\uoX^\circ)$ et une action à droite de $G_K$ sur la catégorie  
$\bHom(\Pi(\uoX^\circ),\Ens)$ des préfaisceaux de $\mU$-ensembles sur $\Pi(\uoX^\circ)^\circ$ \eqref{notconv21}. 
Le foncteur canonique \eqref{notconv21f}
\begin{equation}\label{sshtrl12b}
\uoX^\circ_\fet\rightarrow \bHom(\Pi(\uoX^\circ),\Ens), \ \ \ F\mapsto (\oz\mapsto v_\oz(F))
\end{equation}
est $G_K$-équivariant lorsque l'on fait agir $G_K$ à droite sur la catégorie $\uoX^\circ_\fet$ par image inverse.
Il induit une équivalence entre la catégorie des faisceaux $G_K$-équivariants de $\uoX^\circ_\fet$ 
et la catégorie des préfaisceaux $G_K$-équivariants de
$\mU$-ensembles $\varphi$ sur $\Pi(\uoX^\circ)^\circ$ tels que pour tout point géométrique $\oz$ de $\uoX^\circ$, 
l'action de $\pi_1(\uoX^\circ,\oz)$ sur l'ensemble $\varphi(\oz)$ 
induite par son action canonique sur l'objet $\oz$ de $\Pi(\uoX^\circ)$, soit continue pour la topologie discrète.

\subsection{}\label{sshtrl14}
D'après \eqref{notconv20b}, on associe à tout faisceau abélien $F$ de $\tE$ le préfaisceau abélien
$F_{[\bullet\rightsquigarrow \ox]}$ sur $\Pi(\uoX^\circ)^\circ$ \eqref{sshtrl12} défini par 
la correspondance qui à tout point géométrique $\oz$ de $\uoX^\circ$ associe le groupe abélien $F_{\rho(\oz\rightsquigarrow \ox)}$, 
de sorte que l'action de $\pi_1(\uoX^\circ,\oz)$ sur $F_{\rho(\oz\rightsquigarrow \ox)}$ induite par son action sur l'objet 
$\oz$ de $\Pi(\uoX^\circ)$ coïncide avec l'action sous-jacente au $\mZ[\pi_1(\uoX^\circ,\oz)]$-module 
$F_{[\oz\rightsquigarrow \ox]}$ \eqref{sshtrl9b}.

On associe à tout faisceau abélien $F$ de $\tE'^{\mN^\circ}$ le préfaisceau abélien
$\Theta_*(F)_{[\bullet\rightsquigarrow \ox]}$ sur $\Pi(\uoX^\circ)^\circ$ défini par 
la correspondance qui à tout point géométrique $\oz$ de $\uoX^\circ$ associe le groupe abélien 
$\Theta_*(F)_{\rho(\oz\rightsquigarrow \ox)}$ \eqref{sshtrl1a} 
de sorte que l'action de $\pi_1(\uoX^\circ,\oz)$ sur $\Theta_*(F)_{\rho(\oz\rightsquigarrow \ox)}$ induite par son action sur l'objet 
$\oz$ de $\Pi(\uoX^\circ)$ coïncide avec l'action sous-jacente au $\mZ[\pi_1(\uoX^\circ,\oz)]$-module 
$\Theta_*(F)_{[\oz\rightsquigarrow \ox]}$ \eqref{sshtrl9d}. On définit ainsi un foncteur de la catégorie des 
faisceaux abéliens de $\tE'^{\mN^\circ}$ dans celle des préfaisceaux abéliens sur $\Pi(\uoX^\circ)^\circ$. Celui-ci étant exact à gauche, 
on désigne abusivement par $\rR^q\Theta_*(F)_{[\bullet\rightsquigarrow \ox]}$ ($q\geq 0$) ses foncteurs dérivés à droite. 
Pour tout point géométrique $\oz$ de $\uoX^\circ$, l'évaluation de $\rR^q\Theta_*(F)_{[\bullet\rightsquigarrow \ox]}$ en $\oz$,
munie de l'action de $\pi_1(\uoX^\circ,\oz)$ induite par son action sur l'objet $\oz$ de $\Pi(\uoX^\circ)$, 
n'est autre que le $\mZ[\pi_1(\uoX^\circ,\oz)]$-module $\rR^q\Theta_*(F)_{[\oz\rightsquigarrow \ox]}$ 
défini dans \ref{sshtrl9} en vertu de (\cite{sga4} I 3.1).

De même, on associe à tout faisceau abélien $H$ de $\tG^{\mN^\circ}$ et tout entier $q\geq 0$, le préfaisceau 
$\rR^q\lgg_*(H)_{[\bullet\rightsquigarrow \ox]}$ sur $\Pi(\uoX^\circ)^\circ$ défini par 
la correspondance qui à tout point géométrique $\oz$ de $\uoX^\circ$ associe le groupe abélien 
$\rR^q\lgg_*(H)_{\rho(\oz\rightsquigarrow \ox)}$ défini dans \ref{sshtrl1}, de sorte que l'action de $\pi_1(\uoX^\circ,\oz)$ sur 
$\rR^q\lgg_*(H)_{\rho(\oz\rightsquigarrow \ox)}$ induite par son action sur l'objet 
$\oz$ de $\Pi(\uoX^\circ)$ coïncide avec l'action sous-jacente au $\mZ[\pi_1(\uoX^\circ,\oz)]$-module 
$\rR^q\lgg_*(H)_{[\oz\rightsquigarrow \ox]}$ défini sous \eqref{sshtrl9g}. 

On vérifie aussitôt que pour tout faisceau abélien $F$ de $\tE'^{\mN^\circ}$, la suite spectrale \eqref{sshtrl9n} induit 
une suite spectrale de préfaisceaux abéliens sur $\Pi(\uoX^\circ)^\circ$,
\begin{equation}\label{sshtrl14a}
\rE_2^{i,j}=\rR^i\lgg_*(\rR^j\bvtau_*(F))_{[\bullet\rightsquigarrow \ox]}\Rightarrow \rR^{i+j}\Theta_*(F)_{[\bullet\rightsquigarrow \ox]}. 
\end{equation}

\subsection{}\label{sshtrl15}
Le groupe $G_K$ agit naturellement à gauche sur les topos $\tE'$, $\tG$ et $\tE$ (cf. \ref{notconv17} et \ref{ssht12}).
Les morphismes $\Theta$ et $\lgg$ \eqref{ktfr1c} sont clairement $G_K$-équivariants \eqref{notconv18}.

Pour tout faisceau $G_K$-équivariant $F$ de $\tE$, $\varphi_\ox(F)$ \eqref{sshtrl9a} est canoniquement muni 
d'une structure $G_K$-équivariante  \eqref{sshtrl12}. 
Cela résulte immédiatement des définitions \ref{tf10}, les morphismes $\Phi$ et $\theta$ de \eqref{tf10a} 
étant $G_K$-équivariants \eqref{notconv18}. Par suite, pour tout faisceau abélien $G_K$-équivariant $F$ de $\tE$, le préfaisceau abélien
$F_{[\bullet\rightsquigarrow \ox]}$ sur $\Pi(\uoX^\circ)^\circ$ \eqref{sshtrl14} est canoniquement muni 
d'une structure $G_K$-équivariante \eqref{sshtrl12}. D'après \ref{notconv23}, celle-ci définit par descente un préfaisceau abélien 
sur $\Pi(\uX^\circ)^\circ$ que l'on note $F_{\langle\bullet\rightsquigarrow \ox\rangle}$.

Comme $\ocB$ est un anneau $G_K$-équivariant de $\tE$ \eqref{ssht12g}, il définit un préfaisceau d'anneaux 
$\ocB_{\langle\bullet\rightsquigarrow \ox\rangle}$ sur $\Pi(\uX^\circ)^\circ$.
Par suite, l'anneau $\uoR=\ocB_{\langle\oy\rightsquigarrow \ox\rangle}$ \eqref{sshtrl1e}
est naturellement muni d'une action de $\uGamma$ par des automorphismes d'anneaux \eqref{sshtrl12}. Il en est donc de même de $\uhoR$. 

Pour tout faisceau abélien $G_K$-équivariant $F$ de $\tE'^{\mN^\circ}$,
le préfaisceau abélien $\Theta_*(F)_{[\bullet\rightsquigarrow \ox]}$ sur $\Pi(\uoX^\circ)^\circ$ \eqref{sshtrl14} est canoniquement muni 
d'une structure $G_K$-équivariante. D'après \ref{notconv23}, celle-ci définit par descente un préfaisceau abélien 
sur $\Pi(\uX^\circ)^\circ$ que l'on note $\Theta_*(F)_{\langle\bullet\rightsquigarrow \ox\rangle}$.
Pour tout entier $q\geq 0$, tout point géométrique $\oz$ de $\uoX^\circ$ et tout $u\in G_K$, on considère l'isomorphisme 
\begin{equation}\label{sshtrl15c}
\rR^q\Theta_*(F)_{[\oz\rightsquigarrow \ox]}\stackrel{\sim}{\rightarrow} \rR^q\Theta_*(F)_{[u(\oz)\rightsquigarrow \ox]}
\end{equation}
composé de 
\begin{equation}\label{sshtrl15d}
\rR^q\Theta_*(F)_{[\oz\rightsquigarrow \ox]}\rightarrow \rR^q\Theta_*(u^*(F))_{[\oz\rightsquigarrow \ox]}
\rightarrow \rR^q\Theta_*(F)_{[u(\oz)\rightsquigarrow \ox]},
\end{equation}
où la première flèche est induite par l'isomorphisme $\tau_u^F\colon F\stackrel{\sim}{\rightarrow} u^*(F)$ et la seconde flèche est induite 
par l'isomorphisme canonique fonctoriel en $F$,
\begin{equation}\label{sshtrl15e}
\Theta_*(u^*(F))_{[\oz\rightsquigarrow \ox]}\stackrel{\sim}{\rightarrow} \Theta_*(F)_{[u(\oz)\rightsquigarrow \ox]}. 
\end{equation}
On vérifie aussitôt que ces morphismes munissent $\rR^q\Theta_*(F)_{[\bullet\rightsquigarrow \ox]}$ d'une structure de 
préfaisceau abélien $G_K$-équivariant sur $\Pi(\uoX^\circ)^\circ$. 
D'après \ref{notconv23}, celle-ci définit par descente un préfaisceau abélien 
sur $\Pi(\uX^\circ)^\circ$ que l'on note $\rR^q\Theta_*(F)_{\langle\bullet\rightsquigarrow \ox\rangle}$.
Pour $q=0$, on retrouve le préfaisceau $\Theta_*(F)_{\langle\bullet\rightsquigarrow \ox\rangle}$ introduit plus haut. 

De même, pour tout faisceau abélien $G_K$-équivariant $H$ de $\tG^{\mN^\circ}$ et pour tout entier $q\geq 0$,
le préfaisceau abélien $\rR^q\lgg_*(H)_{[\bullet\rightsquigarrow \ox]}$ sur $\Pi(\uoX^\circ)^\circ$ \eqref{sshtrl14} 
est naturellement muni d'une structure $G_K$-équivariante. 
D'après \ref{notconv23}, celle-ci définit par descente un préfaisceau abélien 
sur $\Pi(\uX^\circ)^\circ$ que l'on note $\rR^q\lgg_*(H)_{\langle\bullet\rightsquigarrow \ox\rangle}$.

\begin{lem}\label{sshtrl16}
Avec les notations de \ref{sshtrl15},
pour tout faisceau abélien $G_K$-équivariant $F$ de $\tE'^{\mN^\circ}$, la suite spectrale \eqref{sshtrl14a} est $G_K$-équivariante. 
\end{lem}

Pour tout $u\in G_K$, considérons l'isomorphisme composé 
\begin{equation}\label{sshtrl16a}
\rR\bvtau_*(F)\stackrel{\sim}{\rightarrow}  \rR\bvtau_*(u^*(F))\stackrel{\sim}{\rightarrow} u^* (\rR\bvtau_*(F)),
\end{equation}
où la première flèche est induite par l'isomorphisme $\tau_u^F$ et la deuxième flèche est l'isomorphisme qui traduit la $G_K$-équivariance de $\bvtau$. 
On en déduit un morphisme de la suite spectrale $\rE$ \eqref{sshtrl14a} vers la deuxième suite spectrale d'hypercohomologie du foncteur 
$\lgg_*(-)_{[\bullet\rightsquigarrow \ox]}$ relativement au complexe $u^*(\rR\bvtau_*(F))$. 
Par ailleurs, pour tout faisceau abélien $H$ de $\tG^{\mN^\circ}$ et tout point géométrique $\oz$ de $\uoX^\circ$,
l'isomorphisme canonique fonctoriel en $H$,
\begin{equation}
\lgg_*(u^* (H))_{[\oz\rightsquigarrow \ox]}\stackrel{\sim}{\rightarrow} \lgg_*(H)_{[u(\oz)\rightsquigarrow \ox]}
\end{equation}
induit un isomorphisme 
\begin{equation}
\rR \lgg_*(u^* (\rR\bvtau_*(F)))_{[\oz\rightsquigarrow \ox]}\stackrel{\sim}{\rightarrow} \rR \lgg_*(\rR\bvtau_*(F))_{[u(\oz)\rightsquigarrow \ox]}.
\end{equation}
Notant $\gamma_u\colon \Pi(\uoX^\circ)\rightarrow \Pi(\uoX^\circ)$ l'action de $G_K$ sur $\Pi(\uoX^\circ)$ \eqref{notconv21a},
on en déduit un isomorphisme de la deuxième suite spectrale d'hypercohomologie du foncteur 
$\lgg_*(-)_{[\bullet\rightsquigarrow \ox]}$ relativement au complexe $u^* (\rR\bvtau_*(F))$ vers la suite spectrale $\rE\circ \gamma_u$. 
On vérifie aussitôt que le morphisme composé de suites spectrales $\rE\rightarrow \rE\circ \gamma_u$ est donné sur les termes $\rE_2$
et sur les aboutissements par les structures $G_K$-équivariantes définies dans \ref{sshtrl15}

\subsection{}\label{sshtrl17}
Supposons le morphisme $g\colon X'\rightarrow X$ propre. Notons $\updelta\colon X'^\rhd\rightarrow X^\circ$ le morphisme induit par $g\colon X'\rightarrow X$, de sorte que 
$\upgamma=\oupdelta$ \eqref{ktfr1b}. Soient $n, q$ deux entiers $\geq 0$. 
En vertu de \ref{Kpun34} et \ref{TCFR160}, le faisceau $\rR^q\updelta_*(\mZ/p^n\mZ)$ est localement constant constructible sur $\oX^\circ$ 
et on a un isomorphisme canonique
\begin{equation}\label{sshtrl17b}
\rR^q\updelta_*(\mZ/p^n\mZ)_\oy\stackrel{\sim}{\rightarrow}\rH^q(\oX'^\rhd_{\oy,\et},\mZ/p^n\mZ).
\end{equation}
D'après  (\cite{agt} VI.9.18 et VI.9.20), 
il existe donc un $(\mZ/p^n\mZ)$-module canonique $\cL^q_n$ de $X^\circ_\fet$ et un isomorphisme canonique
$\rR^q\updelta_*(\mZ/p^n\mZ)\stackrel{\sim}{\rightarrow} \rho_{X^\circ}^*(\cL^q_n)$,
où $\rho_{X^\circ}\colon X^\circ_\et\rightarrow X^\circ_\fet$ est le morphisme canonique \eqref{notconv10a}.
On désigne par $\ucL^q_n$ l'image inverse de $\cL^q_n$ sur $\uX^\circ_\fet$.  
D'après (\cite{agt} VI.9.9), oubliant l'action de $\uGamma$, on a un isomorphisme canonique \eqref{sshtrl12c}
\begin{equation}
\upnu_\oy(\ucL^q_n)\stackrel{\sim}{\rightarrow}\rR^q\updelta_*(\mZ/p^n\mZ)_\oy.
\end{equation}
On en déduit une action naturelle de $\uGamma$ sur 
$\rH^q(\oX'^\rhd_{\oy,\et},\mZ/p^n\mZ)$ et par suite une action naturelle de $\uGamma$ sur $\rH^q(\oX'^\rhd_{\oy,\et},\mZ_p)$. 

En vertu \ref{Kpun34} et \ref{TCFR160}, la formation du faisceau $\rR^q\updelta_*(\mZ/p^n\mZ)$ commute 
à tout changement de la base $X^\circ$. Il s'ensuit compte tenu de (\cite{agt} VI.9.18) que l'image inverse de $\cL^q_n$ sur $\oX^\circ_\fet$ 
s'identifie canoniquement au faisceau $\mL^q_n$ défini après \eqref{sshtrl9j}, et l'image inverse de $\ucL^q_n$ sur $\uoX^\circ_\fet$ 
s'identifie canoniquement au faisceau $\umL^q_n$. On obtient ainsi des structures $G_K$-équivariantes canoniques sur $\rR^q\upgamma_*(\mZ/p^n\mZ)$, 
$\mL^q_n$ et $\umL^q_n$. L'isomorphisme \eqref{sshtrl9k} est clairement $G_K$-équivariant 
et il en est donc de même de \eqref{sshtrl9m}. D'après \ref{notconv23}, ce dernier induit donc un isomorphisme $\uGamma$-équivariant
\begin{equation}\label{sshtrl17d}
\upnu_\oy(\ucL^q_n)\stackrel{\sim}{\rightarrow} \psi_*(\rR^q\upgamma_*(\mZ/p^n\mZ))_{\langle\oy\rightsquigarrow \ox\rangle}.
\end{equation}

\subsection{}\label{sshtrl18}
Les morphismes $\sigma_s$ \eqref{TFA66b} et $\pi_s$ \eqref{ktfr21a} sont $G_K$-équivariants lorsque l'on fait agir $G_K$ trivialement sur 
$X_{s,\et}$ et $X'_{s,\et}$. Pour tout entier $n\geq 1$, le groupe $G_K$ agit naturellement sur les anneaux $\co_{\oX_n}$ et $\co_{\oX'_n}$,
et les morphismes $\sigma_n$ \eqref{TFA8e} et $\pi_n$ \eqref{ktfr24e} sont $G_K$-équivariants. 

Reprenons les hypothèses et notations de \ref{sshtrl4} et de sa preuve.  
Le faisceau $\ocF'_n$ est canoniquement muni d'une structure $G_K$-équivariante. 
On en déduit une structure $G_K$-équivariante canonique sur $\pi^*_n(\ocF'_n)$.
Pour tout entier $q\geq 0$, le morphisme de changement de base \eqref{sshtrl11c} est $G_K$-équivariant. 
Par suite, le morphisme \eqref{sshtrl11b} est sous-jacent à un morphisme $\uGamma$-équivariant 
\begin{equation}\label{sshtrl18a}
\rH^q(\uX',\cF'_n\otimes_{\co_{X}}\co_{\uX})\otimes_{\co_{\uX}(\uX)}\uoR  
\rightarrow \rR^q\lgg_{n*}(\pi^*_n(\ocF'_n))_{\langle\oy\rightsquigarrow \ox\rangle}.
\end{equation}
On notera que $G_K$ et donc $\uGamma$ agissent trivialement sur $\uX$ puisque $k$ est algébriquement clos.

\begin{prop}\label{sshtr30}
Supposons le morphisme $g\colon X'\rightarrow X$ projectif.
Alors, avec les notations de \ref{sshtrl17}, la suite spectrale \eqref{sshtrl6a} est $\uGamma$-équivariante. 
\end{prop} 

En effet, d'après \ref{sshtrl16}, la suite spectrale de préfaisceaux abéliens sur $\Pi(\uoX^\circ)^\circ$ \eqref{sshtrl14a} 
\begin{equation}\label{sshtr30a}
\rE_2^{i,j}=\rR^i\lgg_*(\rR^j\bvtau_*(\bvocB'))_{[\bullet\rightsquigarrow \ox]}\Rightarrow \rR^{i+j}\Theta_*(\bvocB')_{[\bullet\rightsquigarrow \ox]}
\end{equation}
est $G_K$-équivariante. En vertu de \ref{notconv23}, elle induit donc une suite spectrale de préfaisceaux abéliens sur $\Pi(\uX^\circ)^\circ$,
\begin{equation}
\rE_2^{i,j}=\rR^i\lgg_*(\rR^j\bvtau_*(F))_{\langle \bullet\rightsquigarrow \ox\rangle}\Rightarrow \rR^{i+j}\Theta_*(F)_{\langle\bullet\rightsquigarrow \ox\rangle}. 
\end{equation}
On en déduit une suite spectrale de $\mZ[\uGamma]$-modules 
\begin{equation}\label{sshtr30b}
\rE_2^{i,j}=\rR^i\lgg_*(\rR^j\bvtau_*(\bvocB'))_{\langle\oy\rightsquigarrow \ox\rangle}\Rightarrow \rR^{i+j}\Theta_*(\bvocB')_{\langle\oy\rightsquigarrow \ox\rangle}. 
\end{equation}

Pour tout entier $j\geq 0$, le morphisme $\bvocB^!$-linéaire canonique \eqref{sshtr1d}
\begin{equation}\label{sshtr30d}
\bvpi^*(\xi^{-j}\tOmega^j_{\bvoX'/\bvoX})\rightarrow \rR^j \bvtau_*(\bvocB')
\end{equation}
est $G_K$-équivariant. Par suite, pour tout entier $i\geq 0$, l'isomorphisme \eqref{sshtrl1j} est sous-jacent à un isomorphisme $\uGamma$-équivariant
\begin{equation}
\rR^i\lgg_*(\bvpi^*(\xi^{-j}\tOmega^j_{\bvoX'/\bvoX}))_{\langle\oy\rightsquigarrow \ox\rangle}\otimes_{\mZ_p}\mQ_p
\stackrel{\sim}{\rightarrow} \rR^i\lgg_*(\rR^j \bvtau_*(\bvocB'))_{\langle\oy\rightsquigarrow \ox\rangle}\otimes_{\mZ_p}\mQ_p.
\end{equation}
Par ailleurs, il résulte de \ref{sshtrl18} et de la preuve de \ref{sshtrl5} que l'isomorphisme \eqref{sshtrl5a}, pour $\cF'=\tOmega^j_{X'/X}$, 
est sous-jacent à un isomorphisme $\uGamma$-équivariant
\begin{equation}
\rR^i\lgg_*(\bvpi^*(\tOmega^j_{\bvoX'/\bvoX}))_{\langle\oy\rightsquigarrow \ox\rangle}\otimes_{\mZ_p}\mQ_p
\stackrel{\sim}{\rightarrow}
(\rH^i(\uX',\tOmega^j_{X'/X}\otimes_{\co_{X}}\co_{\uX})\hotimes_{\co_\uX}\uoR)\otimes_{\mZ_p}\mQ_p.
\end{equation}

Les morphismes $\psi\colon \oX^\circ_\et\rightarrow \tE$ \eqref{mtfla7d} et 
$\psi'\colon \oX'^\rhd_\et\rightarrow \tE'$ \eqref{mtfla9d} sont clairement $G_K$-équivariants \eqref{ssht12d}. 
Pour tous entiers $q,n\geq 0$, l'isomorphisme canonique \eqref{TCFR18e}
\begin{equation}
\rR^q\Theta_*(\psi'_*(\mZ/p^n\mZ))\stackrel{\sim}{\rightarrow}\psi_*(\rR^q\upgamma_*(\mZ/p^n\mZ))
\end{equation}
et le morphisme $\ocB$-linéaire canonique
\begin{equation}
\rR^q\Theta_*(\psi'_*(\mZ/p^n\mZ))\otimes_{\mZ_p}\ocB\rightarrow \rR^q\Theta_*(\psi'_*(\mZ/p^n\mZ)\otimes_{\mZ_p}\ocB')
\end{equation}
sont $G_K$-équivariants. Compte tenu de l'isomorphisme canonique $\mZ/p^n\mZ\stackrel{\sim}{\rightarrow} \psi'_*(\mZ/p^n\mZ)$ \eqref{TCFR19c}, 
on en déduit que le morphisme $\ocB$-linéaire canonique \eqref{TCFR19a}
\begin{equation}
\psi_*(\rR^q\upgamma_*(\mZ/p^n\mZ))\otimes_{\mZ_p}\ocB\rightarrow \rR^q\Theta_*(\ocB'_n)
\end{equation}
est aussi $G_K$-équivariant. D'après \ref{notconv23}, ce dernier induit un morphisme $\uGamma$-équivariant
\begin{equation}
\psi_*(\rR^q\upgamma_*(\mZ/p^n\mZ))_{\langle\oy\rightsquigarrow \ox\rangle}\otimes_{\mZ_p}\ocB_{\langle\oy\rightsquigarrow \ox\rangle}
\rightarrow \rR^q\Theta_*(\ocB'_n)_{\langle\oy\rightsquigarrow \ox\rangle}.
\end{equation}
D'après \eqref{sshtrl17d}, le morphisme $\uoR$-linéaire canonique \eqref{sshtrl2e}
\begin{equation}\label{sshtr30c}
\rH^q(\oX'^\rhd_{\oy,\et},\mZ/p^n\mZ)\otimes_{\mZ_p}\uoR
\rightarrow \rR^q\Theta_*(\ocB'_n)_{\langle\oy\rightsquigarrow \ox\rangle}
\end{equation}
est donc $\uGamma$-équivariant. Il résulte alors de la preuve de \ref{sshtrl3} que le $\alpha$-isomorphisme \eqref{sshtrl3a} induit un isomorphisme $\uGamma$-équivariant
\begin{equation}
\rR^q\Theta_*(\bvocB')_{\langle\oy\rightsquigarrow \ox\rangle}\otimes_{\mZ_p}\mQ_p\stackrel{\sim}{\rightarrow}
\rH^q(\oX'^\rhd_{\oy,\et},\mZ_p)\otimes_{\mZ_p}\uhoR[\frac 1 p],
\end{equation}
ce qui achève la preuve.

\begin{prop}\label{sshtr31}
Posant $\uR_1=\Gamma(\uoX,\co_\uoX)$ et notant $\huRun$ son séparé complété $p$-adique, on a 
\begin{equation}\label{sshtr31a}
\huRun[\frac 1 p] = \uhoR[\frac 1 p]^\uDelta.
\end{equation}
\end{prop}

La question étant locale pour la topologie étale de $X$, on peut supposer que 
le morphisme $f$ \eqref{hght2} admet une carte adéquate \eqref{cad1}.
On désigne par $\fV_\ox$ la catégorie des $X$-schémas étales $\ox$-pointés.
Pour tout objet $(U,\fp\colon \ox\rightarrow U)$ de $\fV_\ox$, le morphisme $\uX\rightarrow U$ déduit de $\fp$ (\cite{sga4} VIII 7.3)
et le $\oX^\circ$-morphisme donné $\oy\rightarrow \uoX^\circ$ \eqref{sshtrl9} induisent un $U$-morphisme $\oy\rightarrow \oU^\circ$. 
On peut donc considérer $\oy$ comme un point géométrique de $\oU^\circ$.

Pour tout objet $(U,\fp)$ de $\fV_\ox$, le schéma $\oU$ étant localement irréductible d'après (\cite{agt} III.3.3 et III.4.2(iii)),  
il est la somme des schémas induits sur ses composantes irréductibles. On note $\oU^\star$
la composante irréductible de $\oU$ contenant $\oy$. 
De même, $\oU^\circ$ est la somme des schémas induits sur ses composantes irréductibles
et $\oU^{\star \circ}=\oU^\star\times_{X}X^\circ$ est la composante irréductible de $\oU^\circ$ contenant $\oy$. 
On pose $\Delta_U=\pi_1(\oU^{\star \circ},\oy)$.  On a \eqref{TFA12f}
\begin{eqnarray}
\uoX&=&\underset{\underset{(U,\fp)\in \fV_\ox}{\longleftarrow}}{\lim}\oU^\star,\label{sshtr31c}\\
\uoR&=&\underset{\underset{(U,\fp)\in \fV_\ox^\circ}{\longrightarrow}}{\lim}\ \oR_U, \label{sshtr31e}
\end{eqnarray}
où $\oR_U$ est la $\Gamma(\oU^\star,\co_\oU)$-algèbre $\oR^\oy_U$ définie dans \eqref{TFA9c}, 
munie de l'action naturelle de $\Delta_U$ par des automorphismes d'anneaux.

On désigne par $\fW_\ox$ la sous-catégorie pleine de $\fV_\ox$ formée des objets $(U,\fp)$ tels que le schéma $U$ soit affine et connexe. 
Le foncteur d'injection canonique $\fW_\ox^\circ\rightarrow \fV_\ox^\circ$ est cofinal en vertu de (\cite{sga4} I 8.1.3(c)).
En vertu de \ref{tpcg6}, pour tout objet $(U,\fp)$ de $\fW_\ox$ et tout entier $m\geq 0$, le morphisme canonique \eqref{hght1c}
\begin{equation}
\Gamma(\oU^\star_m,\co_{\oU_m})\rightarrow (\oR_U/p^m\oR_U)^{\Delta_U}
\end{equation}
est injectif et son conoyau est annulé par $p^{\frac{1}{p-1}}\fm_\oK$. 
Par ailleurs, d'après \eqref{sshtr31c} et (\cite{agt} VI.11.10), on a un isomorphisme canonique 
\begin{equation}
\underset{\underset{(U,\fp)\in \fW_\ox^\circ}{\longrightarrow}}{\lim}\ (\oR_U/p^m\oR_U)^{\Delta_U}
\stackrel{\sim}{\rightarrow} (\uoR/p^m\uoR)^{\uDelta}.
\end{equation}
On en déduit que le morphisme canonique 
\begin{equation}
\uR_1/p^m\uR_1\rightarrow (\uoR/p^m\uoR)^{\uDelta}
\end{equation}
est injectif et que son conoyau est annulé par $p^{\frac{1}{p-1}}\fm_\oK$. 
Passant à la limite projective, on voit que le morphisme canonique 
\begin{equation}
\huRun\rightarrow \uhoR^{\uDelta}
\end{equation}
est injectif et que son conoyau est annulé par $p^{\frac{1}{p-1}}\fm_\oK$. La proposition s'ensuit en observant que $\uhoR[\frac 1 p]$
est la limite du système inductif $(\uhoR)_{\mN}$ où les morphismes de transition sont induits par la multiplication par $p$.

\begin{prop}\label{sshtr32}
Supposons $X$ affine d'anneau $R$. 
Soient $M$ un $R$-module de type fini tel que $M[\frac 1 p]$ soit projectif sur $R[\frac 1 p]$, $j$ un entier non nul. On a alors 
\begin{equation}\label{sshtr32a}
((M\hotimes_{R}\uoR)\otimes_{\mZ_p}\mQ_p(j))^{\uGamma}=0,
\end{equation}
où le produit tensoriel $\hotimes$ est séparé complété pour la topologie $p$-adique.
\end{prop}

Posons $\uR_1=\Gamma(\uoX,\co_\uoX)$ et notons $\huRun$ son séparé complété $p$-adique.
D'après \ref{sshtr33}, on a 
\begin{equation}\label{sshtr32c}
((M\hotimes_R\uR_1)\otimes_{\mZ_p}\mQ_p(j))^{G_K}=0.
\end{equation}
Il suffit alors de montrer que le morphisme canonique
\begin{equation}\label{sshtr32b}
(M\hotimes_R\uR_1)\otimes_{\mZ_p}\mQ_p\rightarrow ((M\hotimes_{R}\uoR)\otimes_{\mZ_p}\mQ_p)^{\uDelta}
\end{equation}
est un isomorphisme. Il existe un $R$-module de type fini $M'$, deux entiers $m,n\geq 0$ et un morphisme $R$-linéaire 
\begin{equation}
u\colon M\oplus M'\rightarrow R^n
\end{equation}
dont le noyau et le conoyau sont annulés par $p^m$. D'après \ref{alpha3}, il existe un morphisme $R$-linéaire $v\colon R^n\rightarrow M\oplus M'$
tel que $u\circ v= p^{2m}\id_{R^n}$ et $v\circ u=p^{2m}\id_{M\oplus M'}$. On en déduit que les noyaux et conoyaux des morphismes  
\begin{eqnarray}
(M\oplus M')\hotimes_R\uR_1&\rightarrow& \huRun^n,\\
(M\oplus M')\hotimes_{R}\uoR&\rightarrow& \uhoR^n,
\end{eqnarray}
induits par $u$, sont annulés par $p^{2m}$. Considérons le diagramme commutatif 
\begin{equation}
\xymatrix{
{((M\oplus M')\hotimes_R\uR_1)\otimes_{\mZ_p}\mQ_p}\ar[r]\ar[d]_{w_1}&
{(((M\oplus M')\hotimes_R\uoR)\otimes_{\mZ_p}\mQ_p)^{\uDelta}}\ar[d]^{w_2}\\
{\huRun[\frac 1 p]^n}\ar[r]&{(\uhoR[\frac 1 p]^n)^{\uDelta}}}
\end{equation}
Comme $w_1$ est un isomorphisme et que $w_2$ est injectif, on peut se réduire au cas où $M=R$, 
auquel cas l'assertion recherchée résulte de \ref{sshtr31}.

\begin{cor}\label{sshtr34}
Supposons le morphisme $g\colon X'\rightarrow X$ projectif. Alors, la suite spectrale \eqref{sshtrl6a} dégénère en $\rE_2$.
\end{cor}

En effet, la question étant locale sur $X$ pour la topologie de Zariski, on peut supposer $X$ affine d'anneau $R$. 
Pour tous $i,j\geq 0$, posons
\begin{equation}
F^{i,j}=\rH^i(X',\tOmega^j_{X'/X}).
\end{equation}
La différentielle $d_2^{i,j}\colon \rE_2^{i,j}\rightarrow \rE_2^{i+2,j-1}$ est une section $\uGamma$-équivariante de 
\begin{equation}
\Hom_{\uhoR}(F^{i,j}\hotimes_{R}\uoR,F^{i+2,j-1}\hotimes_{R}\uoR)\otimes_{\mZ_p}\mQ_p(1),
\end{equation}
où le produit tensoriel $\hotimes$ est séparé complété pour la topologie $p$-adique.

Nous avons établi dans la preuve de \ref{sshtr11} que le $R[\frac 1 p]$-module $F^{i,j}[\frac 1 p]$ est projectif. 
Il existe donc un $R$-module de type fini $G^{i,j}$,
un entier $n\geq 1$ et un morphisme $R$-linéaire $u\colon F^{i,j}\oplus G^{i,j}\rightarrow R^n$
induisant un isomorphisme après inversion de $p$. Il existe un entier $m\geq 0$ tel que $p^m$ annule le noyau et le conoyau de $u$. 
D'après \ref{alpha3}, il existe alors un morphisme $R$-linéaire $v\colon R^n\rightarrow F^{i,j}\oplus G^{i,j}$
tel que $u\circ v= p^{2m}\id_{R^n}$ et $v\circ u=p^{2m}\id_{F^{i,j}\oplus G^{i,j}}$. 

Comme $u$ est une isogénie (\cite{agt} 6.1.1), on vérifie aisément que le morphisme canonique
\begin{equation}
\Hom_{R}(F^{i,j},F^{i+2,j-1})\hotimes_{R}\uoR
\rightarrow \Hom_{\uhoR}(F^{i,j}\hotimes_R\uoR,F^{i+2,j-1}\hotimes_R\uoR)
\end{equation}
est une isogénie. La proposition résulte alors de \ref{sshtr32} appliqué à $\Hom_{R}(F^{i,j},F^{i+2,j-1})$.

\subsection{}\label{sshtrl19}
On pose $\bvuoR=(\uoR/p^{n+1}\uoR)_{n\geq 0}$ que l'on considère comme un anneau du topos $\Ens^{\mN^\circ}$ \eqref{notconv13}.
On désigne par $\bRep_\bvuoR(\uDelta)$ la catégorie des $\bvuoR$-représentations de $\uDelta$, c'est-à-dire la catégorie des représentations 
semi-linéaires de $\uDelta$ dans des $\bvuoR$-modules de $\Ens^{\mN^\circ}$ \eqref{notconv16}, par $\bRep_\bvuoR(\uDelta)_\mQ$
la catégorie des $\bvuoR$-représentations de $\uDelta$ à isogénie près (\cite{agt} III.6.1), par $\bRep_\uhoR(\uDelta)$ 
la catégorie des $\uhoR$-représentations de $\uDelta$ et par $\bRep_{\uhoR[\frac 1 p]}(\uDelta)$ 
la catégorie des ${\uhoR[\frac 1 p]}$-représentations de $\uDelta$. Le foncteur  
\begin{equation}\label{sshtrl19a}
\bRep_\bvuoR(\uDelta)\rightarrow \bRep_{\uhoR[\frac 1 p]}(\uDelta), \ \ (M_n)_{n\geq 0}\mapsto 
(\underset{\underset{n\geq 0}{\longleftarrow}}{\lim}\ M_n)\otimes_{\mZ_p}\mQ_p,
\end{equation}
transforme les isogénies en isomorphismes. Il induit donc un foncteur que l'on note
\begin{equation}\label{sshtrl19b}
{\lim}_\mQ \colon \bRep_{\bvuoR}(\uDelta)_\mQ\rightarrow \bRep_{\uhoR[\frac 1 p]}(\uDelta).
\end{equation}
De même, le foncteur 
\begin{equation}\label{sshtrl19c}
\bRep_\bvuoR(\uDelta)\rightarrow \bRep_{\uhoR[\frac 1 p]}(\uDelta), 
\ \ (M_n)_{n\geq 0}\mapsto (\rR^1\underset{\underset{n\geq 0}{\longleftarrow}}{\lim}\ M_n)\otimes_{\mZ_p}\mQ_p,
\end{equation}
induit un foncteur que l'on note 
\begin{equation}\label{sshtrl19d}
{\lim}^1_\mQ \colon \bRep_\bvuoR(\uDelta)_\mQ\rightarrow \bRep_{\uhoR[\frac 1 p]}(\uDelta).
\end{equation}
Compte tenu de (\cite{agt} III.6.1.4) et (\cite{gabriel} III cor.~1 à prop.~1), 
toute suite exacte $0\rightarrow M'\rightarrow M\rightarrow M''\rightarrow 0$ de $\bRep_\bvuoR(\uDelta)_\mQ$
induit une suite exacte longue canonique fonctorielle de $\uhoR[\frac 1 p]$-représentations de $\uDelta$, 
\begin{equation}\label{sshtrl19e}
0\rightarrow {\lim}_\mQ(M') \rightarrow {\lim}_\mQ(M)
\rightarrow {\lim}_\mQ(M'') \rightarrow {\lim}^1_\mQ(M')\rightarrow
{\lim}^1_\mQ(M)\rightarrow {\lim}^1_\mQ(M'')\rightarrow 0.
\end{equation}

Par ailleurs, le foncteur composé $\nu_\oy\circ\varphi_\ox$ \eqref{sshtrl9} induit un foncteur exact
\begin{equation}\label{sshtrl19f}
\bMod(\bvocB)\rightarrow \bRep_\bvuoR(\uDelta), \ \ \cF=(\cF_n)_{n\geq 0}\mapsto \cF_{[\oy\rightsquigarrow \ox]}=((\cF_n)_{[\oy\rightsquigarrow \ox]}),
\end{equation}
et par suite un foncteur exact
\begin{equation}\label{sshtrl19g}
\bMod_\mQ(\bvocB)\rightarrow \bRep_\bvuoR(\uDelta)_\mQ.
\end{equation}

\begin{rema}\label{sshtrl190}
Conservons les notations de \ref{sshtrl90}. Pour tout $\co_X$-module $\cF$ de $X_\et$, on a un isomorphisme canonique de $\bvuoR$-representations de $\uDelta$
\begin{equation}
(\bvsigma^*(\cF\otimes_{\co_X}\co_{\bvoX}))_{[\oy\rightsquigarrow \ox]}\stackrel{\sim}{\rightarrow} \cF_\ox\otimes_{\uR}\bvuoR,
\end{equation}
où $\cF_\ox$ et $\uR$ \eqref{sshtrl1ee} sont considérés comme des représentations triviales de $\uDelta$. 
En effet, d'après \eqref{TFA14b}, pour tout entier $n\geq 0$, on a un isomorphisme canonique 
\begin{equation}
\varphi_\ox(\sigma^*_n(\cF\otimes_{\co_X}\co_{\oX_n}))\stackrel{\sim}{\rightarrow} \cF_\ox\otimes_{\uR}\varphi_\ox(\ocB_n),
\end{equation}
où $\sigma_n$ est le morphisme de topos annelés \eqref{ssht12l} et $\cF_\ox$ et $\uR$ sont considérés comme des faisceaux constants de $\uoX^\circ_\fet$. 
\end{rema}

\subsection{}\label{sshtrl20}
Supposons le morphisme $g\colon X'\rightarrow X$ projectif. 
Avec les notations de \ref{sshtrl90} et \ref{sshtrl19}, considérons, pour tout entier $i\geq 0$, la $\bvuoR$-representation de $\uDelta$
\begin{equation}
\rH^i=(\nu_\oy(\umL^i_{n+1})\otimes_{\mZ_p}\uoR)_{n\geq 0}.
\end{equation}
Pour tous entiers $i,j\geq 0$, on a un isomorphisme $\uR$-linéaire canonique
\begin{equation}
\rR^ig_*(\tOmega^j_{X'/X})_\ox\stackrel{\sim}{\rightarrow} \rH^i(\uX',\tOmega^j_{\uX'/\uX}).
\end{equation}
Par suite, compte tenu de \eqref{sshtrl9l} et \ref{sshtrl190}, l'image de la suite spectrale de Hodge-Tate relative \eqref{sshtr5a} 
par le foncteur \eqref{sshtrl19g} induit une suite spectrale de $\bRep_\bvuoR(\uDelta)_\mQ$
\begin{equation}
\rE_2^{i,j}=\rH^i(\uX',\tOmega^j_{\uX'/\uX})\otimes_{\uR}\bvuoR_\mQ(-j)\Rightarrow \rH^{i+j}_\mQ.
\end{equation}
Celle-ci dégénère en $\rE_2$ d'après \ref{sshtr11}. Elle définit donc pour tout $q\geq 0$, une filtration décroissante exhaustive 
$(\fil^q_r)_{0\leq r\leq q+1}$ de $\rH^q_\mQ$ telle que $\fil_{q+1}^q=0$ et 
que pour tout $0\leq r\leq q$, on ait une suite exacte canonique
\begin{equation}\label{sshtrl20f}
0\rightarrow \fil_{r+1}^q \rightarrow \fil_r^q \rightarrow \rH^r(\uX',\tOmega^{q-r}_{\uX'/\uX})\otimes_{\uR}\bvuoR_\mQ(r-q)\rightarrow 0. 
\end{equation}
Comme le système $\rH^r(\uX',\tOmega^{q-r}_{\uX'/\uX})\otimes_{\uR}\bvuoR$ vérifie clairement la condition de Mittag-Leffler, 
on en déduit par récurrence que pour tout $0\leq r\leq q$, on a \eqref{sshtrl19d} 
\begin{equation}\label{sshtrl20g}
{\lim}^1_\mQ(\fil_r^q)=0.
\end{equation}
Par suite, d'après \ref{sshtrl200} et \eqref{sshtrl19e}, $({\lim}_\mQ(\fil_r^q))_{0\leq r\leq q+1}$ \eqref{sshtrl19b} définit une filtration décroissante exhaustive 
de $\rH^q(\oX'^\rhd_{\oy,\et},\mZ_p)\otimes_{\mZ_p}\uhoR[\frac 1 p]$ par des $\uhoR[\frac 1 p]$-représentations de 
$\uDelta$ telle que ${\lim}_\mQ(\fil_{q+1}^q)=0$ et que pour tout $0\leq r\leq q$, 
on ait une suite exacte canonique
\begin{equation}\label{sshtrl20h}
0\rightarrow {\lim}_\mQ(\fil_{r+1}^q) \rightarrow {\lim}_\mQ(\fil_r^q) \rightarrow 
\rH^{r}(\uX',\tOmega^{q-r}_{\uX'/\uX})\hotimes_{\uR}\uoR[\frac 1 p](r-q)\rightarrow 0,
\end{equation}
où le produit tensoriel $\hotimes$ est séparé et complété pour la topologie $p$-adique. Cette filtration n'est autre que la filtration aboutissement de la suite spectrale \eqref{sshtrl6a}.
On peut déduire directement de \ref{sshtr9} qu'elle est $\uGamma$-équivariante en procédant comme dans le preuve de \ref{sshtr30}.

\section[La suite spectrale de Hodge-Tate relative: sections globales]{La suite spectrale de Hodge-Tate relative: sections globales au-dessus d'un petit schéma affine}\label{sshtrsg}

\subsection{}\label{sshtrsg1}
Rappelons que les hypothèses de \ref{ktfr1} sont en vigueur dans cette section. 
On suppose de plus que le morphisme 
$f\colon (X,\cM_X)\rightarrow (S,\cM_S)$ \eqref{hght2} admet une carte adéquate 
\eqref{cad1}, que $X=\Spec(R)$ est affine et connexe, et que $X_s$ est non-vide. 

Soit $\oy$ un point géométrique de $\oX^\circ$. 
Le schéma $\oX$ étant localement irréductible d'après \ref{cad7}(iii),  
il est la somme des schémas induits sur ses composantes irréductibles. On note $\oX^\star$
la composante irréductible de $\oX$ contenant $\oy$. 
De même, $\oX^\circ$ est la somme des schémas induits sur ses composantes irréductibles
et $\oX^{\star \circ}=\oX^\star\times_{X}X^\circ$ est la composante irréductible de $\oX^\circ$ contenant $\oy$. 

On pose $\Delta=\pi_1(\oX^{\star\circ},\oy)$ et on note $\mZ[\Delta]$ son algèbre de groupe. 
On désigne par $\bB_{\Delta}$ le topos classifiant du groupe profini $\Delta$ \eqref{notconv11}, par 
\begin{equation}\label{sshtrsg1a}
\nu_\oy\colon \oX^{\star\circ}_\fet\rightarrow \bB_{\Delta},
\end{equation}
le foncteur fibre en $\oy$ \eqref{notconv11c}, et par $\upbeta_\oy$ le foncteur composé 
\begin{equation}\label{sshtrsg1b}
\upbeta_\oy \colon \tE\rightarrow \bB_{\Delta}, \ \ \ F\mapsto \nu_\oy\circ (\beta_*(F)|\oX^{\star\circ}), 
\end{equation}
où $\beta$ est le morphisme de topos \eqref{TFA6b}. 

Soit $(V_i)_{i\in I}$ le revêtement universel normalisé de $\oX^{\star \circ}$ en $\oy$ \eqref{notconv11}.  
Pour tout $i\in I$, $(V_i\rightarrow X)$ est un objet de $E$. 
D'après (\cite{agt} VI.9.9), pour tout objet $F$ de $\tE$, oubliant l'action de $\Delta$, on a un isomorphisme canonique fonctoriel
\begin{equation}\label{sshtrsg1c}
\upbeta_\oy(F)\stackrel{\sim}{\rightarrow} \underset{\underset{i\in I}{\longrightarrow}}{\lim}\ F(V_i\rightarrow X).
\end{equation}
On définit ainsi un foncteur de la catégorie des groupes abéliens dans celle des $\mZ[\Delta]$-modules. Celui-ci étant exact à gauche, 
on désigne par $\rR^q \upbeta_\oy$ $(q\geq 0)$ ses foncteurs dérivés à droite. 
Pour tout faisceau abélien $F$ de $\tE$ et tout entier $q\geq 0$, on a un isomorphisme canonique fonctoriel
\begin{equation}\label{sshtrsg1d}
\rR^q\upbeta_\oy(F)\stackrel{\sim}{\rightarrow} \nu_\oy\circ (\rR^q\beta_*(F)|\oX^{\star\circ}).
\end{equation}
Oubliant l'action de $\Delta$, on a un isomorphisme canonique fonctoriel
\begin{equation}\label{sshtrsg1e}
\rR^q\upbeta_\oy(F)\stackrel{\sim}{\rightarrow} \underset{\underset{i\in I}{\longrightarrow}}{\lim}\ \rH^q((V_i\rightarrow X),F).
\end{equation}

On pose $\oR=\upbeta_\oy(\ocB)$ qui n'est autre que la représentation discrète $\oR^\oy_X$ de $\Delta$ définie dans \eqref{TFA9b}.
On note $\hoR$ le complété séparé $p$-adique de $\oR$.

\subsection{}\label{sshtrsg2}
On reprend les notations de \ref{ssht1}.
Pour tout faisceau abélien $F=(F_n)_{n\geq 0}$ de $\tE^{\mN^\circ}$, on désigne abusivement par $\upbeta_\oy(F)$ 
la limite projective des $\upbeta_\oy(F_n)$,
\begin{equation}\label{sshtrsg2a}
\upbeta_\oy(F)=\underset{\underset{n\geq 0}{\longleftarrow}}{\lim}\ \upbeta_\oy(F_n),
\end{equation}
où $\upbeta_\oy(F_n)$ est le foncteur défini dans \eqref{sshtrsg1b}.
On définit ainsi un foncteur de la catégorie des faisceaux abéliens de $\tE^{\mN^\circ}$ dans celle des $\mZ[\Delta]$-modules.
Celui-ci étant exact à gauche, on désigne abusivement par $\rR^q\upbeta_\oy(F)$ ($q\geq 0$) ses foncteurs dérivés à droite. 
D'après (\cite{jannsen} 1.6), on a une suite exacte canonique 
\begin{equation}\label{sshtrsg2b}
0\rightarrow \rR^1 \underset{\underset{n\geq 0}{\longleftarrow}}{\lim}\ \rR^{q-1}\upbeta_\oy(F_n)\rightarrow
\rR^q\upbeta_\oy(F)\rightarrow 
\underset{\underset{n\geq 0}{\longleftarrow}}{\lim}\ \rR^q\upbeta_\oy(F_n)\rightarrow 0,
\end{equation}
où on a posé $\rR^{-1}\upbeta_\oy(F_n)=0$ pour tout $n\geq 0$. 

Pour tout entier $q\geq 0$, le foncteur $\rR^q\upbeta_\oy$ induit un foncteur qu'on note aussi 
\begin{equation}\label{sshtrsg2c}
\rR^q\upbeta_\oy\colon \bMod(\bvocB)\rightarrow \bRep_{\hoR}(\Delta),
\end{equation}
dans la catégorie des $\hoR$-représentations de $\Delta$ \eqref{notconv16}. 
Celui-ci induit un foncteur qu'on note aussi 
\begin{equation}\label{sshtrsg2d}
\rR^q\upbeta_\oy\colon \bMod_\mQ(\bvocB)\rightarrow \bRep_{\hoR[\frac 1 p]}(\Delta). 
\end{equation}

Il résulte facilement de \eqref{sshtrsg2b}, \eqref{sshtrsg1e} et (\cite{sga4} V 3.5) que les foncteurs $\rR^q\upbeta_\oy$ $(q\geq 0)$ \eqref{sshtrsg2c}
sont les foncteurs dérivés à droite du foncteur $\upbeta_\oy=\rR^0\upbeta_\oy$ sur la catégorie des $\bvocB$-modules. 
Par suite, les foncteurs $\rR^q\upbeta_\oy$ $(q\geq 0)$ \eqref{sshtrsg2d}
sont les foncteurs dérivés à droite du foncteur $\upbeta_\oy=\rR^0\upbeta_\oy$ sur la catégorie des $\bvocB_\mQ$-modules (\cite{agt} III.6.1.6). 

\begin{prop}\label{sshtrsg3}
Soient $\cF$ un $\co_X$-module cohérent tel que le $\co_{X_\eta}$-module $\cF|X_\eta$ soit localement libre, $M=\Gamma(X,\cF)$,
$\bvocF=(\ocF_n)_{n\geq 1}=\bvsigma^*(\cF\otimes_{\co_X}\co_{\bvoX})$, où $\bvsigma$ est le morphisme de topos annelés \eqref{ssht1a} 
(cf. \ref{ktfr38}), que l'on considère aussi comme un $\bvocB$-module de $\tE^{\mN^\circ}$. Alors, 
\begin{itemize}
\item[{\rm (i)}] Le morphisme canonique
\begin{equation}\label{sshtrsg3a}
M\otimes_R\hoR[\frac 1 p]\rightarrow \upbeta_\oy(\bvocF)\otimes_{\mZ}\mQ
\end{equation}
est un isomorphisme \eqref{sshtrsg2a}. 
\item[{\rm (ii)}] Pour tout entier $q\geq 1$, $\rR^q\upbeta_\oy(\bvocF)\otimes_{\mZ}\mQ$ est nul. 
\end{itemize}
\end{prop}

En effet, il existe un $\co_X$-module cohérent $\cG$, deux entiers $m,N\geq 1$ et un morphisme $\co_X$-linéaire 
$u\colon \cF\oplus \cG\rightarrow \co_X^m$ dont le noyau et le conoyau sont annulés par $p^N$. D'après \ref{alpha3}, il existe alors un 
morphisme $\co_X$-linéaire $v\colon \co_X^m\rightarrow \cF\oplus \cG$ tel que $u\circ v=p^{2N}\id_{\co_X^m}$ et 
$v\circ u= p^{2N}\id_{\cF\oplus \cG}$. On peut donc se réduire au cas où $\cF=\co_X$. 
Il résulte de \ref{amtF28}(i) et \eqref{sshtrsg1c} que pour tout entier $n\geq 1$, le morphisme canonique 
\begin{equation}\label{sshtrsg3b}
\oR/p^n\oR \rightarrow \upbeta_\oy(\ocB_n)
\end{equation}
est un $\alpha$-isomorphisme, d'où la proposition (i).  Il s'ensuit que 
$\rR^1\underset{\underset{n\geq 1}{\longleftarrow}}{\lim}\ \upbeta_\oy(\bvocB_n)$
est $\alpha$-nul. D'autre part, d'après \ref{amtF28}(ii) et \eqref{sshtrsg1e}, pour tous entiers $n,q\geq 1$, $\rR^q\upbeta_\oy(\ocB_n)$ est 
$\alpha$-nul.  La proposition (ii) s'ensuit compte tenu de \eqref{sshtrsg2b}.

\begin{lem}\label{sshtrsg4}
Soient $n$ un entier $\geq 1$, $\mL$ un $(\mZ/p^n\mZ)$-module de type fini $\oX^{\circ}_\fet$. Alors, le morphisme $\oR$-linéaire
\begin{equation}\label{sshtrsg4a}
\nu_\oy(\mL|\oX^{\star\circ})\otimes_\mZ\oR\rightarrow \upbeta_\oy(\beta^*(\mL)\otimes_\mZ\ocB),
\end{equation}
induit par le morphisme d'adjonction $\mL\rightarrow \beta_*(\beta^*(\mL))$, est un $\alpha$-isomorphisme. 
\end{lem}
 
En effet, le morphisme \eqref{sshtrsg4a} est induit par le morphisme canonique 
\begin{equation}
\underset{\underset{i\in I}{\longrightarrow}}{\lim}\ \mL(V_i)\rightarrow 
\underset{\underset{i\in I}{\longrightarrow}}{\lim}\ \Gamma((V_i\rightarrow X),\beta^*(\mL)\otimes_\mZ\ocB).
\end{equation}
Pour tout $i\in I$, on désigne par $\tE_i$ le topos de Faltings associé au morphisme $V_i\rightarrow X$, 
que l'on identifie au localisé du topos $\tE$ en $(V_i\rightarrow X)^\ra$ (\cite{agt} VI.10.14), et par 
\begin{equation}\label{sshtrsg4b}
\beta_i\colon \tE_i\rightarrow V_{i,\fet}
\end{equation}
le morphisme canonique (\cite{agt} (VI.10.6.3)). D'après (\cite{agt} (VI.10.12.6)), le diagramme 
\begin{equation}\label{sshtrsg4c}
\xymatrix{
{\tE_i}\ar[r]\ar[d]_{\beta_i}&{\tE}\ar[d]^{\beta}\\
{V_{i,\fet}}\ar[r]&{\oX^\circ_\fet}}
\end{equation}
est commutatif à isomorphisme canonique près. 

Par ailleurs, il existe une partie cofinale $J\subset I$ telle que pour tout $j\in J$, le faisceau $\mL|V_j$ soit constant. 
On peut donc se réduire au cas où $\mL$ est constant, de valeur $\mZ/p^n\mZ$, auquel cas la proposition résulte de \ref{amtF28}(i).

\subsection{}\label{sshtrsg5}
Supposons, le morphisme $g\colon X'\rightarrow X$ propre. Soient $n,q$ deux entiers $\geq 0$.
En vertu de \ref{Kpun34} et \ref{TCFR160}, le faisceau $\rR^q\upgamma_*(\mZ/p^n\mZ)$ est localement constant constructible sur 
$\oX^\circ$ et on a un isomorphisme canonique
\begin{equation}\label{sshtrsg5a}
\rR^q\upgamma_*(\mZ/p^n\mZ)_\oy\stackrel{\sim}{\rightarrow}\rH^q(\oX'^\rhd_{\oy,\et},\mZ/p^n\mZ).
\end{equation}
D'après (\cite{agt} VI.9.18 et VI.9.20), 
il existe donc un $(\mZ/p^n\mZ)$-module canonique $\mL^q_n$ de $\oX^\circ_\fet$ et un isomorphisme canonique
$\rR^q\upgamma_*(\mZ/p^n\mZ)\stackrel{\sim}{\rightarrow} \rho_{\oX^\circ}^*(\mL^q_n)$,
où $\rho_{\oX^\circ}\colon \oX^\circ_\et\rightarrow \oX^\circ_\fet$ est le morphisme canonique \eqref{notconv10a}.
D'après (\cite{agt} VI.9.9), oubliant l'action de $\Delta$, on a un isomorphisme canonique \eqref{sshtrsg1a}
\begin{equation}\label{sshtrsg5b}
\nu_\oy(\mL^q_n|\oX^{\star\circ})\stackrel{\sim}{\rightarrow}\rR^q\upgamma_*(\mZ/p^n\mZ)_\oy.
\end{equation}
On en déduit une action naturelle de $\Delta$ sur 
$\rH^q(\oX'^\rhd_{\oy,\et},\mZ/p^n\mZ)$ et par suite une action naturelle de $\Delta$ sur $\rH^q(\oX'^\rhd_{\oy,\et},\mZ_p)$. 

D'après (\cite{agt} VI.10.9(iii)), on a un isomorphisme canonique 
\begin{equation}\label{sshtrsg5c}
\beta^*(\mL^q_n)\stackrel{\sim}{\rightarrow} \psi_*(\rR^q\upgamma_*(\mZ/p^n\mZ)).
\end{equation}
On en déduit par adjonction un morphisme canonique $\oR$-linéaire et $\Delta$-équivariant 
\begin{equation}\label{sshtrsg5d}
\rH^q(\oX'^\rhd_{\oy,\et},\mZ/p^n\mZ)\otimes_\mZ\oR\rightarrow \upbeta_\oy(\psi_*(\rR^q\upgamma_*(\mZ/p^n\mZ))\otimes_\mZ\ocB),
\end{equation}
qui est un $\alpha$-isomorphisme d'après \ref{sshtrsg4}. Calquant la preuve de \ref{ssht300}, on en déduit un morphisme canonique
$\hoR$-linéaire et $\Delta$-équivariant \eqref{sshtrsg2c}
\begin{equation}\label{sshtrsg5e}
\rH^q(\oX'^\rhd_{\oy,\et},\mZ_p)\otimes_{\mZ_p}\hoR\rightarrow \upbeta_\oy(\bvpsi_*(\rR^q\bvupgamma_*(\bvmZ_p))
\otimes_{\mZ_p}\bvocB),
\end{equation}
qui est un $\alpha$-isomorphisme.

\begin{teo}\label{sshtrsg6}
Supposons le morphisme $g\colon X'\rightarrow X$ projectif. 
Alors, pour tout entier $q\geq 0$, 
il existe $(\fil^q_r)_{0\leq r\leq q+1}$, une filtration décroissante exhaustive canonique 
de $\rH^q(\oX'^\rhd_{\oy,\et},\mZ_p)\otimes_{\mZ_p}\hoR[\frac 1 p]$ par des $\hoR[\frac 1 p]$-représentations de $\Delta$, 
telle que $\fil^q_{q+1}$ soit nul et que pour tout entier $0\leq r\leq q$, on ait une suite exacte $\Delta$-équivariante canonique 
\begin{equation}\label{sshtrsg6e}
0\rightarrow \fil^q_{r+1}\rightarrow \fil^q_r\rightarrow \rH^r(X',\tOmega^{q-r}_{X'/X})\otimes_R\hoR[\frac 1 p](r-q)
\rightarrow 0.
\end{equation}
\end{teo}

En effet, notons $(\Fil^q_r)_{0\leq r\leq q+1}$ la filtration décroissante exhaustive de $\bvpsi_*(\rR^q\bvupgamma_*(\bvmZ_p))
\otimes_{\mZ_p}\bvocB_\mQ$ définie par la suite spectrale \eqref{sshtr5a}, où $\Fil^q_{q+1}=0$. 
Compte tenu de \eqref{sshtrsg5e}, celle-ci induit $(\upbeta_\oy(\Fil^q_r))_{0\leq r\leq q+1}$ \eqref{sshtrsg2d},  
une filtration décroissante exhaustive par des $\hoR[\frac 1 p]$-représentations de $\Delta$, 
de $\rH^q(\oX'^\rhd_{\oy,\et},\mZ_p)\otimes_{\mZ_p}\hoR[\frac 1 p]$. 

Pour tous entiers $i,j\geq 0$, posons
\begin{equation}\label{sshtrsg6a}
\cF^{i,j}=\rR^ig_*(\tOmega^j_{X'/X}).
\end{equation}
On notera que le $\co_{X_\eta}$-module $\cF^{i,j}|X_\eta$ est localement libre (cf. la preuve de \ref{sshtr11}).  
En vertu de \ref{sshtrsg3}, on a un $\hoR$-isomorphisme $\Delta$-équivariant canonique 
\begin{equation}\label{sshtrsg6b}
\rH^i(X',\tOmega^j_{X'/X})\otimes_R\hoR[\frac 1 p]
\stackrel{\sim}{\rightarrow}\upbeta_\oy(\bvsigma^*(\cF^{i,j}\otimes_{\co_X}\co_{\bvoX}))\otimes_{\mZ}\mQ,
\end{equation}
et pour tout $a\geq 1$, $\rR^a\upbeta_\oy(\bvsigma^*(\cF^{i,j}\otimes_{\co_X}\co_{\bvoX}))\otimes_{\mZ}\mQ$ est nul. 

D'après \ref{sshtr11}, pour tous entiers $0\leq r\leq q$, on a une suite exacte canonique
\begin{equation}\label{sshtrsg6c}
0\rightarrow \Fil^q_{r+1}\rightarrow \Fil^q_r\rightarrow \bvsigma^*(\cF^{r,q-r}\otimes_{\co_X}\co_{\bvoX})_\mQ(r-q)
\rightarrow 0.
\end{equation}
On en déduit, par une récurrence croissante sur $r$, que pour tout $a\geq 1$, $\rR^a\upbeta_\oy(\Fil^q_r)=0$,
et qu'on a une suite exacte $\Delta$-équivariante canonique 
\begin{equation}\label{sshtrsg6d}
0\rightarrow \upbeta_\oy(\Fil^q_{r+1})\rightarrow \upbeta_\oy(\Fil^q_r)\rightarrow \rH^r(X',\tOmega^{q-r}_{X'/X})\otimes_R\hoR[\frac 1 p](r-q)
\rightarrow 0.
\end{equation}
La proposition s'ensuit en prenant $\fil^q_r=\upbeta_\oy(\Fil^q_r)$.

\begin{rema}\label{sshtrsg7}
La construction de la filtration $(\fil^q_r)_{0\leq r\leq q+1}$ dans \ref{sshtrsg6} est parallèle à celle donnée dans \ref{sshtrl20} pour la localisation.
On peut donner deux autres constructions de cette filtration comme filtration aboutissement de deux suites spectrales de Cartan-Leray,
similaires à \ref{sshtrl6} et \ref{sshtrl8} pour la localisation.
\end{rema}

\begin{rema}\label{sshtrsg8}
Notant $\varpi\colon \oX^\circ\rightarrow X^\circ$ la projection canonique et posant $\Gamma=\pi_1(X^\circ,\varpi(\oy))$, 
on peut munir $\rH^q(\oX'^\rhd_{\oy,\et},\mZ_p)$ et $\hoR$ d'actions canoniques de $\Gamma$ qui prolongent celles de $\Delta$.
Alors, la filtration $(\fil^q_r)_{0\leq r\leq q+1}$ de $\rH^q(\oX'^\rhd_{\oy,\et},\mZ_p)\otimes_{\mZ_p}\hoR[\frac 1 p]$ définie dans 
\ref{sshtrsg6} est $\Gamma$-équivariante, et les morphismes des suites exactes \eqref{sshtrsg6e} sont $\Gamma$-équivariants. 
En effet, notant $L$ la clôture algébrique de $K$ dans le corps des fonctions de $X^\circ$ et $G_L$ le sous-groupe ouvert de $G_K$
défini par un plongement de $L$ dans $\oK$, on a une suite exacte de groupes profinis
\begin{equation}
1\rightarrow \Delta\rightarrow \Gamma\rightarrow G_L \rightarrow 1.
\end{equation}
L'assertion résulte alors de \ref{sshtr9} par descente galoisienne \eqref{notconv23}, de façon similaire à la preuve de \ref{sshtr30}. 
\end{rema}

\end{document}